\documentclass[12pt,a4paper]{amsart}
\usepackage{latexsym,eucal}
\usepackage[all,cmtip,arrow]{xy}
\usepackage{hyperref}

\calclayout

\let\le\undefined
\DeclareMathSymbol{\le}{\mathrel}{AMSa}{"36}         
\let\ge\undefined
\DeclareMathSymbol{\ge}{\mathrel}{AMSa}{"3E}         
\DeclareMathSymbol{\birarrow}{\mathrel}{AMSa}{"13}   
\DeclareMathSymbol{\empt}{\mathord}{AMSb}{"3F}       
\DeclareMathSymbol{\bethl}{\mathalpha}{AMSb}{"69}    
\DeclareMathSymbol{\subneq}{\mathrel}{AMSb}{"20}     
\DeclareMathSymbol{\rightthreetimes}{\mathbin}{AMSa}{"69}

\renewcommand{\:}{\colon}
\newcommand{\+}{\protect\nobreakdash-}
\renewcommand{\.}{\mskip 0.5\thinmuskip}
\renewcommand{\;}{,\>}

\newcommand{\rarrow}{\longrightarrow}
\newcommand{\larrow}{\longleftarrow}
\newcommand{\mpsto}{\longmapsto}
\newcommand{\ot}{\otimes}
\newcommand{\sub}{\subset}
\newcommand{\ocn}{\odot}
\newcommand{\oc}{\mathbin{\text{\smaller$\square$}}}
\newcommand{\tim}{\.{\rightthreetimes}\.}

\newcommand{\beth}{\text{\smaller{$\bethl$}}}
\newcommand{\lrarrow}{\.\relbar\joinrel\relbar\joinrel\rightarrow\.}

\newcommand{\bu}{{\text{\smaller\smaller$\scriptstyle\bullet$}}}
\newcommand{\subbu}{{\text{\smaller\smaller
                                  $\scriptscriptstyle\bullet$}}}
\newcommand{\st}{\star} \newcommand{\bst}{\star}

\DeclareFontFamily{U}{mathx}{\hyphenchar\font45}
\DeclareFontShape{U}{mathx}{m}{n}{
      <5> <6> <7> <8> <9> <10>
      <10.95> <12> <14.4> <17.28> <20.74> <24.88>
      mathx10
      }{}
\DeclareSymbolFont{mathx}{U}{mathx}{m}{n}
\DeclareFontSubstitution{U}{mathx}{m}{n}
\DeclareMathAccent{\widecheck}{0}{mathx}{"71}

\DeclareMathOperator{\Spec}{Spec}

\DeclareMathOperator{\PSupp}{PSupp}

\DeclareMathOperator{\coker}{coker}
\DeclareMathOperator{\cone}{cone}
\DeclareMathOperator{\Hom}{Hom}
\DeclareMathOperator{\Ext}{Ext}
\DeclareMathOperator{\Tor}{Tor}

\DeclareMathOperator{\qHom}{\mathcal H \mskip-.3\thinmuskip
  \text{\rmfamily\mdseries\fontshape{ui}\selectfont om}}
\DeclareMathOperator{\fHom}{\mathfrak{Hom}}
\DeclareMathOperator{\Cohom}{\mathfrak{Cohom}}
\DeclareMathOperator{\cohom}{Cohom}
\DeclareMathOperator{\Cr}{\mathfrak C}
\DeclareMathOperator{\Qr}{\mathcal Q}
\DeclareMathOperator{\Pro}{\mathsf{Pro}}
\DeclareMathOperator{\pro}{\mathsf{pro}}
\DeclareMathOperator{\FC}{\mathrm{FC}}

\DeclareMathOperator{\rsd}{rsd}
\DeclareMathOperator{\crd}{crd}
\newcommand{\id}{\mathrm{id}}

\newcommand{\Hot}{\mathsf{Hot}}
\newcommand{\Com}{\mathsf{Com}}
\newcommand{\Acycl}{\mathsf{Acycl}}
\newcommand{\Fil}{\mathsf{Fil}}

\newcommand{\qcoh}{{\operatorname{\mathsf{--qcoh}}}}
\newcommand{\coh}{{\operatorname{\mathsf{--coh}}}}
\renewcommand{\cosh}{{\operatorname{\mathsf{--cosh}}}}
\newcommand{\ctrh}{{\operatorname{\mathsf{--ctrh}}}}
\newcommand{\lcth}{{\operatorname{\mathsf{--lcth}}}}
\newcommand{\vect}{{\operatorname{\mathsf{--vect}}}}
\newcommand{\modl}{{\operatorname{\mathsf{--mod}}}}
\newcommand{\modr}{{\operatorname{\mathsf{mod--}}}}
\newcommand{\modrfp}{{\operatorname{\mathsf{mod_{fp}--}}}}
\newcommand{\comodl}{{\operatorname{\mathsf{--comod}}}}
\newcommand{\comodr}{{\operatorname{\mathsf{comod--}}}}
\newcommand{\contra}{{\operatorname{\mathsf{--contra}}}}
\newcommand{\tors}{{\operatorname{\mathsf{--tors}}}}
\newcommand{\discr}{{\operatorname{\mathsf{--discr}}}}
\newcommand{\discrR}{{\operatorname{\mathsf{discr--}}}}
\newcommand{\fadj}{{f{\operatorname{\mathsf{--adj}}}}}
\newcommand{\fac}{{f{\operatorname{\mathsf{--ac}}}}}
\newcommand{\fWac}{{f/\bW{\operatorname{\mathsf{--ac}}}}}
\newcommand{\dashac}{{{\operatorname{\mathsf{--ac}}}}}
\newcommand{\ffdd}{{{\operatorname{\mathsf{ffd--}}}d}}
\newcommand{\ffD}{{{\operatorname{\mathsf{ffd--}}}D}}
\newcommand{\ffdD}{{{\operatorname{\mathsf{ffd--}}}(d+D)}}
\newcommand{\flidd}{{{\operatorname{\mathsf{flid--}}}d}}
\newcommand{\fliD}{{{\operatorname{\mathsf{flid--}}}D}}
\newcommand{\flidD}{{{\operatorname{\mathsf{flid--}}}(d+D)}}
\newcommand{\fvfdd}{{{\operatorname{\mathsf{fvfd--}}}d}}
\newcommand{\fvfdD}{{{\operatorname{\mathsf{fvfd--}}}(d+D)}}
\newcommand{\fidd}{{{\operatorname{\mathsf{fid--}}}d}}
\newcommand{\fiD}{{{\operatorname{\mathsf{fid--}}}D}}
\newcommand{\fpdd}{{{\operatorname{\mathsf{fpd--}}}d}}
\newcommand{\fpD}{{{\operatorname{\mathsf{fpd--}}}D}}
\newcommand{\fpdD}{{{\operatorname{\mathsf{fpd--}}}(d+D)}}
\newcommand{\alfdd}{{{\operatorname{\mathsf{alfd--}}}d}}

\newcommand{\CApr}{{\C/A{\operatorname{\mathsf{--pr}}}}}
\newcommand{\Acot}{{A{\operatorname{\mathsf{--cot}}}}}
\newcommand{\Bcot}{{B{\operatorname{\mathsf{--cot}}}}}

\newcommand{\Xinj}{{X{\operatorname{\mathsf{--inj}}}}}
\newcommand{\Xlin}{{X{\operatorname{\mathsf{--lin}}}}}
\newcommand{\Xlctin}{{X{\operatorname{\mathsf{--lctin}}}}}
\newcommand{\Xprj}{{X{\operatorname{\mathsf{--prj}}}}}
\newcommand{\Xfl}{{X{\operatorname{\mathsf{--fl}}}}}
\newcommand{\Xvfl}{{X{\operatorname{\mathsf{--vfl}}}}}
\newcommand{\Xalf}{{X{\operatorname{\mathsf{--alf}}}}}
\newcommand{\ct}{{\operatorname{\mathrm{-ct}}}}
\newcommand{\qc}{{\operatorname{\mathrm{-qc}}}}

\newcommand{\co}{{\mathsf{co}}}
\newcommand{\ctr}{{\mathsf{ctr}}}
\newcommand{\abs}{{\mathsf{abs}}}
\newcommand{\bco}{{\mathsf{bco}}}
\newcommand{\bctr}{{\mathsf{bctr}}}
\renewcommand{\b}{{\mathsf{b}}}
\newcommand{\si}{{\mathsf{si}}}
\newcommand{\lh}{{\mathsf{lh}}}
\newcommand{\rh}{{\mathsf{rh}}}
\newcommand{\fd}{{\mathsf{ffd}}}
\newcommand{\fpd}{{\mathsf{fpd}}}
\newcommand{\fid}{{\mathsf{fid}}}
\newcommand{\lid}{{\mathsf{flid}}}
\newcommand{\vfd}{{\mathsf{fvfd}}}

\newcommand{\alfd}{{\mathsf{alfd}}}
\newcommand{\ac}{{\mathsf{ac}}}

\newcommand{\op}{{\mathsf{op}}}
\newcommand{\rop}{{\mathrm{op}}}
\newcommand{\fl}{{\mathsf{fl}}}
\newcommand{\inj}{{\mathsf{inj}}}
\newcommand{\fpi}{{\mathsf{fpi}}}
\newcommand{\lin}{{\mathsf{lin}}}
\newcommand{\lct}{{\mathsf{lct}}}

\newcommand{\al}{{\mathsf{al}}}

\newcommand{\alf}{{\mathsf{alf}}}
\newcommand{\prj}{{\mathsf{prj}}}
\newcommand{\cta}{{\mathsf{cta}}}
\newcommand{\vfl}{{\mathsf{vfl}}}
\renewcommand{\cot}{{\mathsf{cot}}}
\newcommand{\fp}{{\mathsf{fp}}}
\newcommand{\coa}{{\mathsf{coa}}}
\newcommand{\cfq}{{\mathsf{cfq}}}
\newcommand{\fq}{{\mathsf{fq}}}
\newcommand{\dil}{{\mathsf{dil}}}

\renewcommand{\O}{{\mathcal O}}
\newcommand{\cA}{{\mathcal A}}
\newcommand{\cB}{{\mathcal B}}
\newcommand{\C}{{\mathcal C}}
\newcommand{\D}{{\mathcal D}}
\newcommand{\E}{{\mathcal E}}
\newcommand{\F}{{\mathcal F}}
\newcommand{\G}{{\mathcal G}}
\newcommand{\cH}{{\mathcal H}}
\newcommand{\cP}{{\mathcal P}}
\newcommand{\cQ}{{\mathcal Q}}
\newcommand{\cR}{{\mathcal R}}
\newcommand{\K}{{\mathcal K}}
\renewcommand{\L}{{\mathcal L}}
\newcommand{\M}{{\mathcal M}}
\newcommand{\N}{{\mathcal N}}
\newcommand{\J}{{\mathcal J}}
\newcommand{\I}{{\mathcal I}}

\newcommand{\bcD}{{\boldsymbol{\mathcal D}}}
\newcommand{\bOmega}{{\boldsymbol{\Omega}}}

\newcommand{\p}{{\mathfrak p}}
\newcommand{\q}{{\mathfrak q}}
\newcommand{\m}{{\mathfrak m}}
\newcommand{\mm}{{\text{\rmfamily\mdseries
                \fontshape{ui}\selectfont m}}}

\renewcommand{\P}{{\mathfrak P}}
\newcommand{\Q}{{\mathfrak Q}}
\newcommand{\R}{{\mathfrak R}}
\renewcommand{\S}{{\mathfrak S}}
\newcommand{\gA}{{\mathfrak A}}
\newcommand{\gB}{{\mathfrak B}}
\newcommand{\gC}{{\mathfrak C}}
\newcommand{\gJ}{{\mathfrak J}}
\newcommand{\gI}{{\mathfrak I}}
\newcommand{\gE}{{\mathfrak E}}
\newcommand{\gF}{{\mathfrak F}}
\newcommand{\gG}{{\mathfrak G}}
\newcommand{\gH}{{\mathfrak H}}
\newcommand{\gK}{{\mathfrak K}}
\newcommand{\gT}{{\mathfrak T}}
\newcommand{\gL}{{\mathfrak L}}
\newcommand{\gM}{{\mathfrak M}}
\newcommand{\gN}{{\mathfrak N}}
\newcommand{\gU}{{\mathfrak U}}
\newcommand{\gV}{{\mathfrak V}}
\newcommand{\gW}{{\mathfrak W}}

\newcommand{\sA}{{\mathsf A}}
\newcommand{\sB}{{\mathsf B}}
\newcommand{\sC}{{\mathsf C}}
\newcommand{\sD}{{\mathsf D}}
\newcommand{\sE}{{\mathsf E}}
\newcommand{\sF}{{\mathsf F}}
\newcommand{\sG}{{\mathsf G}}
\newcommand{\sH}{{\mathsf H}}
\newcommand{\sK}{{\mathsf K}}
\newcommand{\sL}{{\mathsf L}}
\newcommand{\sR}{{\mathsf R}}
\newcommand{\sS}{{\mathsf S}}
\newcommand{\sT}{{\mathsf T}}
\newcommand{\sX}{{\mathsf X}}

\newcommand{\boR}{{\mathbb R}}
\newcommand{\boL}{{\mathbb L}}
\newcommand{\boZ}{{\mathbb Z}}
\newcommand{\boQ}{{\mathbb Q}}
\newcommand{\boA}{{\mathbb A}}

\newcommand{\bB}{{\mathbf B}}
\newcommand{\bD}{{\mathbf D}}
\newcommand{\bW}{{\mathbf W}}
\newcommand{\bT}{{\mathbf T}}
\newcommand{\bV}{{\mathbf V}}
\newcommand{\bw}{{\mathbf w}}
\newcommand{\bt}{{\mathbf t}}
\newcommand{\bS}{{\mathbf S}}

\theoremstyle{plain}
\newtheorem{thm}{Theorem}[subsection]

\newtheorem{lem}[thm]{Lemma}
\newtheorem{prop}[thm]{Proposition}
\newtheorem{cor}[thm]{Corollary}
\theoremstyle{definition}
\newtheorem{rem}[thm]{Remark}
\newtheorem{qst}[thm]{Question}
\newtheorem{ex}[thm]{Example}
\newtheorem{exs}[thm]{Examples}

\newcommand{\Section}[1]{\bigskip\addtocontents{toc}{\medskip}
 \section{#1}\addtocontents{toc}{\protect\nopagebreak\smallskip}\medskip}
\numberwithin{equation}{section}

\begin{document}

\title{Contraherent cosheaves on schemes}
\author{Leonid Positselski}

\address{Institute of Mathematics, Czech Academy of Sciences \\
\v Zitn\'a~25, 115~67 Praha~1 \\ Czech Republic}

\email{positselski@math.cas.cz}

\begin{abstract}
 Contraherent cosheaves are globalizations of contraadjusted or
cotorsion modules over commutative rings obtained by gluing together
over a scheme, with the colocalization functors $\Hom_R(S,{-})$
used for the gluing (where $S$ is the ring of functions on an affine
open subscheme in $\Spec R$).
 The category of contraherent cosheaves over a scheme is a Quillen
exact category with exact functors of infinite product.
 Over a quasi-compact semi-separated scheme or a Noetherian scheme of
finite Krull dimension (in a different version---over any locally
Noetherian scheme), it also has enough projectives.
 We construct the derived co-contra correspondence over a scheme
in two forms.
 The ``na\"\i ve'' one is an equivalence of the conventional derived
categories of quasi-coherent sheaves and contraherent cosheaves, valid
over any quasi-compact semi-separated scheme.
 The more sophisticated version is an equivalence between the coderived
category of quasi-coherent sheaves and the contraderived category of
contraherent cosheaves over a Noetherian scheme with a dualizing complex.
 The former point of view allows us to obtain an explicit
construction of the Lipman--Neeman extraordinary inverse image
functor~$f^!$ for a morphism of quasi-compact semi-separated schemes
$f\:Y\rarrow X$.
 The latter approach provides an expanded version of the covariant
Serre--Grothendieck duality theory and leads to the Hartshorne--Deligne
extraordinary inverse image functor~$f^!$ (which we denote by~$f^+$)
for a morphism of finite type~$f$ between Noetherian schemes.
 We also construct a derived semico-semicontra correspondence, mounting
the ``na\"\i ve'' version of the correspondence along the fibers on top
of the one depending on a dualizing complex on the base of
a flat fibration.
 Semi-separated Noetherian stacks, affine Noetherian formal schemes,
and ind-affine ind-schemes (together with the noncommutative
analogues) are briefly discussed in the appendices.
\end{abstract}

\maketitle

\tableofcontents

\section*{Preface}
\medskip

 This book is the product of work spanning a time period of more than
fifteen years, from 2009 to~2025.
 It solves a problem which I~posed to myself in Spring~2009, when
my typing of the previous thick monograph~\cite{Psemi} was almost
finished and the bulk of the subsequent memoir~\cite{Pkoszul} was
written up.
 An autobiographical account of my work on~\cite{Pkoszul} can be found
in the survey~\cite{Pksurv}.

 The original 2009~problem was to \emph{define a category}.
 Given a scheme $X$, I~wanted to have an abelian category
dual-analogous to the abelian category $X\qcoh$ of quasi-coherent
sheaves on~$X$.
 In particular, while the category $X\qcoh$ has exact functors of
infinite \emph{coproducts}, my desired category was supposed to
have exact functors of infinite \emph{products}.
 The far-reaching aim was to plug the sought-for abelian category into
the construction of the triangulated \emph{contraderived category},
just as the abelian category $X\qcoh$ can be plugged into
the construction of the \emph{coderived category}.

 What does ``dual-analogous'' mean?
 Let me elaborate on that, in order to clarify a possible
misconception.

 The most abstract category theory is self-dual.
 One can always reverse the arrows by passing from a category $\sC$
to the opposite category~$\sC^\op$.
 The passage from $\sC$ to $\sC^\op$ switches products with
coproducts, projective limits with inductive limits, kernels with
cokernels, projective objects with injective objects, etc.

 The theory of modules over associative or commutative rings $R$ in
\emph{not} self-dual.
 The opposite category $R\modl^\op$ to the abelian category of
(arbitrary) $R$\+modules $R\modl$ is \emph{not} equivalent to
the category of modules over any other ring.
 In particular, filtered inductive limits are exact in $R\modl$, while
filtered projective limits aren't.
 The natural morphism from the direct sum of a family of modules to
their direct product is injective, but not surjective.
 The category $R\modl$ is locally (finitely) presentable, but not
locally copresentable.

 However, the theory of modules over rings is
\emph{self-dual-analogous}.
 This means that most module-theoretic constructions have their
\emph{dual analogues}.
 For example, projective modules are dual-analogous to injective
modules, in the obvious sense.
 Flat modules are also dual-analogous to injective modules, in
the sense that a right $R$\+module $F$ is flat if and only if
the left $R$\+module $\Hom_\boZ(F,\boQ/\boZ)$ is injective.

 The module-theoretic duality-analogy is \emph{not} to be confused
with the abstract category-theoretic duality.
 For example, in introductory undergraduate courses of homological
algebra the students learn that projective modules are easy to
describe, while injective modules are more complicated.
 In more advanced studies, the situation is often reversed: some
questions about projective modules are hard (such as, e.~g.,
the Zariski or flat descent of projectivity of modules over
commutative rings, which is a celebrated theorem of Raynaud
and Gruson~\cite{RG,Pe}); while dual-analogous questions about
injective modules may be easier.

 The point is that, whatever these difficulties and differences are,
they certainly \emph{cannot} be resolved or bridged by reversing
the arrows in the category-theoretic sense of the word.
 The assertion that the injective objects of $R\modl$ are nothing
other than the projective objects of $R\modl^\op$ does not tell one
much about injective modules.
 Neither is the assertion that the projective objects of $R\modl$
are nothing other than the injective objects of $R\modl^\op$ very
helpful in studying projective modules.

 Outside of the module theory proper, one prime example of
a duality-analogy is the duality-analogy between \emph{comodules}
and \emph{contramodules}~\cite{EM,Prev}.
 Let me emphasize that, given a coalgebra $\C$ over a field~$k$
(or a coring $\C$ over a ring~$A$), the category of left
$\C$\+contramodules $\C\contra$ is completely different from
the opposite category $\C\comodl^\op$ to the category of left
$\C$\+comodules $\C\comodl$ or the opposite category $(\comodr\C)^\op$
to the category of right $\C$\+comodules $\comodr\C$.
 For example, when $A=k$ and $\C=k$ is just a field, \emph{both}
the categories $\C\comodl$ and $\C\contra$ are equivalent to
the category of (infinite-dimensional) $k$\+vector spaces $k\vect$,
which is, of course, quite different from its opposite category
$k\vect^\op$.

 In order to stress this point, the present author sometimes says
and writes that the category $\C\contra$ (or rather, the class or family
of all such categories ranging over all coalgebras or corings~$\C$)
is ``covariantly dual'' to the category $\C\comodl$ or $\comodr\C$,
as opposed to the ``contravariantly dual'' opposite categories
$\C\comodl^\op$ or $(\comodr\C)^\op$.
 Moreover, the categories $\C\comodl$ and $\C\contra$ are indeed
connected by the covariant functors and derived equivalences of
the \emph{comodule-contramodule
correspondence}~\cite{Psemi,Pkoszul,Pmgm,Pmc,Prev,Pksurv,PS1}.
 Certainly \emph{no} covariant derived equivalence is possible between
$k\vect$ and $k\vect^\op$.

 From category-theoretic perspective, we would suggest
the duality-analogy between Grothendieck abelian categories and
locally presentable abelian categories with enough projective
objects~\cite[Sections~6\+-9]{PS4} as a illustration of what
``covariant duality'' means.
 From the perspective of the abelian group theory, the duality
between torsion abelian groups and reduced cotorsion abelian
groups~\cite{Harrison,Pcta} plays this role.
 We refer to the papers~\cite{PR,PS3} and the preprint~\cite{Pper}
for a further discussion of the class of abelian categories to which
the contramodule categories belong.

 We would argue that the theory of modules over rings \emph{is}
self-dual-analogous (as stated above), but the theory of
quasi-coherent sheaves over schemes is \emph{not}.
 For instance, the infinite coproduct functors are exact in $X\qcoh$,
but the infinite products are not.
 There are enough injective quasi-coherent sheaves on any scheme $X$,
but projective quasi-coherent sheaves usually do not exist.

 One can use flat quasi-coherent sheaves instead of projective ones
(as in~\cite{M-th}), but still a duality-analogy between injective and
flat quasi-coherent sheaves is weak, if one can speak of such
an analogy at all.
 Let's say, the flatness property of quasi-coherent sheaves is local,
but their injectivity property is only local on locally Noetherian
schemes~\cite[Section~II.7]{Har}, \cite[Example~2.8]{Pal}.
 The homotopy flatness property of complexes of quasi-coherent sheaves
is local, but the homotopy injectivity property of such complexes is
\emph{not} local even on affine schemes of finite type over
a field~\cite[Example~6.5]{N-bb}, \cite{Bel}, \cite[Example~2.9]{Pal}).
 We would argue that a better duality-analogy for quasi-coherent
sheaves is needed and actually exists.

 With this preparatory discussion in mind, we can now tell the reader
that the category $X\ctrh$ of \emph{contraherent cosheaves} over
a scheme $X$ is dual-analogous or ``covariantly dual'' to the category
$X\qcoh$ of quasi-coherent sheaves on~$X$.
 So, the category $X\ctrh$ stands in the same relation to the category
$X\qcoh$ as contramodules over a coalgebra or coring stand to comodules,
or as projective/flat modules over a ring stand to injective modules.
 As in the case of comodules and contramodules, there are covariant
functors and various covariant derived equivalences connecting
$X\qcoh$ with $X\ctrh$.
 These are discussed at length in this book.

 To repeat the same point again, the category of contraherent cosheaves
$X\ctrh$ is completely different from the opposite category
$X\qcoh^\op$ to the category of quasi-coherent sheaves $X\qcoh$.
 For example, when $X=\Spec k$ is just the spectrum of a field,
\emph{both} the categories $X\qcoh$ and $X\ctrh$ are equivalent to
the category of $k$\+vector spaces $k\vect$.

 We suggest that the lack of a satisfactory self-duality-analogy on
quasi-coherent sheaves should be interpreted not as a bug,
but as a feature.
 It points to the existence of the category of contraherent cosheaves,
opening the door to nontrivial triangulated equivalences of
co-contra correspondence connecting quasi-coherent sheaves with
contraherent cosheaves.
 Among these, there are the na\"\i ve co-contra correspondence,
the covariant Serre--Grothendieck duality (which becomes a quadrality
when the contraherent cosheaves are included into it), and
the semico-semicontra correspondence.
 We refer to the Introduction below for a more substantial discussion.

 As all duality-analogies tend to be, the duality-analogy between
quasi-coherent sheaves and contraherent cosheaves is not perfect.
 Most importantly, the category of quasi-coherent sheaves $X\qcoh$
is \emph{abelian}, while the category of contraherent cosheaves
$X\ctrh$ is only \emph{exact} (in the sense of Quillen~\cite{Bueh}).
 The duality-analogy does capture many important features: in
particular, the exact category of contraherent cosheaves $X\ctrh$
does have exact functors of infinite direct product (just as I~wanted
from the very beginning).
 However, while the abelian category $X\qcoh$ has enough injective
objects for any scheme $X$, one needs to make various (not too
restrictive) assumptions about $X$ in order to prove that one or
another version of the exact category $X\ctrh$ has enough
projective objects.

 Let me wrap up the historical discussion and finish this preface with
a few additional references.
 After several unsuccessful attempts, I~was able to come up with
my definition of a contraherent cosheaf in Spring~2012.
 In the immediately preceding weeks, my work on Appendix~B to what
became the memoir~\cite{Pweak} convinced me that the problem of
defining a dual-analogous category to quasi-coherent sheaves can and
must be solved, after all; and so I solved it.

 Specifically, I~learned from~\cite[Appendix~B]{Pweak} that
the notion of a contramodule over the adic completion of a Noetherian
commutative ring at an arbitrary (nonmaximal) ideal was perfectly
well-behaved (see also~\cite{Pmgm,Pcta,Pper,Pdc}).
 Contraherent cosheaves were (and, from today's perspective, are)
expected to solve the problem of globalizing such contramodules
over schemes.

 The first early draft of this book manuscript was offered to the public
in September~2012.
 Subsequent updated and expanded versions were posted to
the \texttt{arXiv} between 2012--2017; they are still available
as~\cite{Pcosh}.
 My work on developing broad context surrounding the theory of
contraherent cosheaves continued in the subsequent years in such papers
as~\cite{Pmgm,PR,Pcta,PSl1,Pdc,Pal,PS6}.

 Partly based on this accumulated knowledge, and even more so on
interrelated work of other authors such as~\cite{BCE} and~\cite{ES},
this book manuscript was reworked and greatly expanded in~2024--25.
 A fairly recent, elaborate discussion of what problems the theory of
contraherent cosheaves is supposed to solve, what kind of applications
it should have, and why and how it is possible for it to work as
expected is also available in the long announcement~\cite{Pphil}.

\bigskip
\addtocontents{toc}{\smallskip}
\section*{Introduction}
\addtocontents{toc}{\protect\nopagebreak\smallskip}
\medskip

\subsection{Schemes as corings, sheaves as comodules}
\label{introd-qcoh-sheaves-as-comodules}
 Quasi-coherent sheaves resemble comodules.
 Both form abelian categories with exact functors of infinite direct sum
(and in fact, even of filtered inductive limit) and with enough
injective objects.
 Here we presume comodules over a coalgebra over a field or left
comodules over a right flat coring over a ring.
 When one restricts to quasi-coherent sheaves over Noetherian schemes
and comodules over (flat) corings over Noetherian rings, both
the abelian categories are locally Noetherian.
 Neither has projective objects or exact functors of infinite product,
in general.

 In fact, quasi-coherent sheaves \emph{are} comodules.
 Let $X$ be a quasi-compact semi-separated scheme and $\{U_\alpha\}$ be
its finite affine open covering.
 Denote by $T$ the disconnected union of the schemes~$U_\alpha$; so
$T$ is also an affine scheme and the natural morphism $T\rarrow X$
is affine.
 Then quasi-coherent sheaves over $X$ can be described as quasi-coherent
sheaves $\F$ over $T$ endowed with an isomorphism $\phi\: p_1^*(\F)
\simeq p_2^*(\F)$ between the two inverse images under the natural maps
$p_1$, $p_2\:T\times_X T \birarrow T$.
 The isomorphism~$\phi$ has to satisfy a natural associativity
constraint.

 In other words, this means that the ring of functions
$\C=\O(T\times_X T)$ has a natural structure of a coring over
the ring $A=\O(T)$.
 The quasi-coherent sheaves over $X$ are the same thing as (left or
right) comodules over this coring.
 The quasi-coherent sheaves over a (good enough) stack can be also
described in such way~\cite{KR,KR2,Pflcc}.

\subsection{Comodules and contramodules}
 There are two kinds of module categories over a coalgebra or coring:
in addition to the more familiar comodules, there are also
\emph{contramodules}~\cite{Prev}.
 Introduced originally by Eilenberg and Moore~\cite{EM} in 1965 (see
also the notable paper~\cite{Bar}), they were all but forgotten for
four decades, until the author's preprint and then
monograph~\cite{Psemi} attracted some new attention to them
towards the end of 2000's.
 See~\cite{Prev} for a survey or~\cite[Section~8]{Pksurv} for
an introductory discussion.

 Assuming that a coring $\C$ over a ring $A$ is a projective left
$A$\+module, the category of left $\C$\+contramodules is abelian
with exact functors of infinite product and enough projectives.
 Generally, contramodules are ``dual-analogous'' to comodules in most
respects, i.~e., they behave as though they formed two opposite
categories---which in fact they don't (or otherwise it wouldn't be
interesting).

 On the other hand, there is an important homological phenomenon of
\emph{comodule-contramodule correspondence}, or a \emph{covariant}
equivalence between appropriately defined (``exotic'') derived categories
of left comodules and left contramodules over the same coring.
 This equivalence is often obtained by deriving certain adjoint
functors which act between the abelian categories of comodules and
contramodules and induce a covariant equivalence between their
appropriately picked exact or additive subcategories (e.~g.,
the equivalence of \emph{Kleisli categories}, which was empasized in
application to comodules and contramodules in the paper~\cite{BBW}).
 A detailed discussion of the philosophy and various versions of
the co-contra correspondence can be found in the introductions to
the papers~\cite{Pmgm,Pps}.

\subsection{Contraherent cosheaves as contramodules over corings}
\label{introd-contrah-as-contram}
 \emph{Contraherent cosheaves} are geometric module objects over
a scheme that are similar to (and, sometimes, particular cases of)
contramodules in the same way as quasi-coherent sheaves are similar
to (or particular cases of) comodules.
 Thus the simplest way to define contraherent cosheaves would be to
assume one's scheme $X$ to be quasi-compact and semi-separated, pick
its finite affine open covering $\{U_\alpha\}$, and consider
contramodules over the related coring $\C$ over the ring $A$ as
constructed above in Section~\ref{introd-qcoh-sheaves-as-comodules}.

 This indeed largely agrees with our approach in this book, but there
are several problems to be dealt with.
 First of all, $\C$~is not a projective $A$\+module, but only a flat
one.
 The most immediate consequence is that one cannot hope for an abelian
category of $\C$\+contramodules, but at best for an exact category.
 This is where the \emph{cotorsion modules}~\cite{Xu,En,EJ} (or their
generalizations which we call the \emph{contraadjusted}
modules~\cite{ST,Pcta,PSl1,Pal}) come into play.
 Secondly, it turns out that the exact category of $\C$\+contramodules,
however defined, depends on the choice of an affine open
covering~$\{U_\alpha\}$.

 In the exposition below, we strive to make our theory as similar (or
rather, dual-analogous) to the classical theory of quasi-coherent
sheaves as possible, while refraining from the choice of a covering
to the extent that it remains practicable.
 We start with defining cosheaves of modules over a sheaf of rings on
a topological space, and proceed to introduce the exact subcategory
of contraherent cosheaves in the exact category of cosheaves of
$\O_X$\+modules on an arbitrary scheme~$X$.

\subsection{Cosheaves and costalks}
 Several attempts to develop a theory of cosheaves have been made
in the literature over the years (see, e.~g., \cite{Bre,Schn,BF}).
 The main difficulty arising in this connection is that the conventional
theory of sheaves depends on the exactness property of filtered
inductive limits in the categories of abelian groups or sets
in an essential way.
 E.~g., the most popular approach is based on wide use of
the construction of the stalks, which are defined as filtered inductive
limits.
 It is important that the functors of stalks are exact on both
the categories of presheaves and sheaves.

 The problem is that the costalks of a co(pre)sheaf would be constructed
as filtered projective limits, and filtered projective limits of
abelian groups are not exact.
 One possible way around this difficulty is to restrict oneself to
(co)constructible cosheaves, for which the projective limits defining
the costalks may be stabilizing and hence exact.
 This is apparently the approach taken in dissertation~\cite{Cur} (see
the discussion with further references in~\cite[Section~2.5]{Cur} and
the subsequent publications~\cite{CP,Cur2}).
 The present work is essentially based on the observation that one
does not really need the (co)stalks in the quasi-coherent/contraherent
(co)sheaf theory, as \emph{the functors of (co)sections over affine open
subschemes are already exact on the (co)sheaf categories}.

 On some of our contraherent cosheaf categories (viz., the ``locally
cotorsion locally contraherent cosheaves''), the functors of costalks
are actually exact; but this comes as an afterthought.
 The point is that, in order to dualize the classical approach to
sheaf theory, one would need the functors of costalks to be exact
on copresheaves.
 That certainly does not hold in our context.

\subsection{So, what is a contraherent cosheaf?}
 Let us explain the main definition in some detail.
 A quasi-coherent sheaf $\F$ on a scheme $X$ can be simply defined
as a correspondence assigning to every affine open subscheme
$U\sub X$ an $\O_X(U)$\+module $\F(U)$ and to every pair of
embedded affine open subschemes $V\sub U \sub X$ an isomorphism
of $\O_X(V)$\+modules
$$
 \F(V) \simeq \O_X(V)\ot_{\O_X(U)} \F(U).
$$
 The obvious compatibility condition for three embedded affine
open subschemes $W\sub V\sub U\sub X$ needs to be imposed~\cite{EE}.

 Analogously, a contraherent cosheaf $\P$ on a scheme $X$ is
a correspondence assigning to every affine open subscheme
$U\sub X$ an $\O_X(U)$\+module $\P[U]$ and to every pair of
embedded affine open subschemes $V\sub U\sub X$ an isomorphism
of $\O_X(V)$\+modules
$$
 \P[V]\simeq\Hom_{\O_X(U)}(\O_X(V),\P[U]).
$$
 The difference with the quasi-coherent case is that
the $\O_X(U)$\+module $\O_X(V)$ is always flat, but not necessarily
projective.
 So to make one's contraherent cosheaves well-behaved, one has
to impose the additional Ext-vanishing condition
$$
 \Ext_{\O_X(U)}^1(\O_X(V),\P[U]) = 0
$$
for all affine open subschemes $V\sub U\sub X$.
 Notice that the $\O_X(U)$\+module $\O_X(V)$ always has projective
dimension not exceeding~$1$, so the condition on $\Ext^1$ is
sufficient.

 This nonprojectivity problem is the reason why one does not have
an abelian category of contraherent cosheaves.
 Imposing the Ext-vanishing requirement allows to obtain, at least,
an exact category.
 Given a rule assigning to affine open subschemes $U\sub X$
the $\O_X(U)$\+modules $\P[U]$ together with the isomorphisms for
embedded affine open subschemes $V\sub U$ as above, and assuming
that the $\Ext^1$\+vanishing condition holds, one can show
that $\P$ satisfies the cosheaf axiom for coverings of affine open
subschemes of $X$ by other affine open subschemes.
 Then a general result from~\cite{Groth} says that $\P$ extends
uniquely from the base of affine open subschemes to a cosheaf of
$\O_X$\+modules defined, as it should be, on all the open subsets of~$X$.

\subsection{Contraadjusted and cotorsion modules}
 A module $P$ over a commutative ring $R$ is said to be
\emph{contraadjusted} if $\Ext^1_R(R[s^{-1}],P)=0$ for all $s\in R$.
 This is equivalent to the vanishing of $\Ext^1_R(S,P)$ for all
the $R$\+algebras $S$ of functions on the affine open subschemes of
$\Spec R$.
 Furthermore, a left module $P$ over a (not necessarily commutative)
ring $R$ is said to be \emph{cotorsion} if $\Ext^1_R(F,P)=0$ 
(or equivalently, $\Ext^{>0}_R(F,P)=0$) for all flat left
$R$\+modules~$F$.
 It has been proved that cotorsion modules are ``numerous enough''
(see~\cite{ET,BBE}); one can prove the same for contraadjusted modules
in the similar way.

 The theory of contraherent cosheaves, as developed in this book, has
two main branches, the \emph{locally contraadjusted} and \emph{locally
cotorsion} contraherent cosheaves, corresponding to the two classes of
modules over commutative rings defined in the previous paragraph.
 Each one of these two versions of the theory comes with its own set of
technical problems to be dealt with.

\subsection{Very flat modules and Very Flat Conjecture}
 Any quotient module of a contraadjusted module is contraadjusted, so
the ``contraadjusted dimension'' of any module does not exceed~$1$;
while the cotorsion dimension of a module may be infinite if
the projective dimensions of flat modules are.
 This is one main advantage of the class of contraadjusted modules.

 Modules of the complementary class to the contraadjusted ones (in
the same way as the flats are complementary to the cotorsion modules) we
call \emph{very flat}.
 All very flat modules have projective dimensions not exceeding~$1$.
 Contraadjusted modules are abundant, while very flat modules are
relatively rare.
 Still, one needs one's morphisms of schemes and one's quasi-coherent
sheaves to be very flat in order to perform some basic constructions
with locally contraadjusted contraherent cosheaves, such as the inverse
image and $\Cohom$ functors.
 We refer to~\cite[Section~6.4]{Pphil} for a discussion of this 
phenomenon.
 Very flat quasi-coherent sheaves are also applicable outside of
the contraherent cosheaf theory~\cite[Remark~2.6]{EP}, \cite{ES}.

 The \emph{Very Flat Conjecture}, formulated in an early
February~2014 preliminary version of this book 
manuscript~\cite[Conjecture~1.7.2]{Pcosh} and proved in 2017 in
the paper~\cite{PSl1}, comes to one's rescue here.
 This conjecture (now theorem) tells us that very flat modules, and
particularly very flat morphisms of commutative rings or schemes,
are not so rare after all.
 Specifically, if $R$ is a commutative ring and $S$ is a finitely
presented $R$\+algebra, and if $S$ is also a flat $R$\+module, then
$S$ is a very flat $R$\+module~\cite[Main Theorem~1.1]{PSl1}.

\subsection{Periodicity theorems}
 On the other hand, one needs to know that restricting one's class
of modules from all modules to contraadjusted or cotorsion ones does
not affect the derived category.
 For contraadjusted modules, this holds simply because all modules
have finite contraadjusted dimensions.
 For cotorsion modules, this is a much more nontrivial result, known
as the \emph{cotorsion periodicity theorem}~\cite{BCE}.
 The cotorsion periodicity theorem claims that in any acyclic complex
of cotorsion modules, the modules of cocycles are also
cotorsion~\cite[Theorem~5.1(2)]{BCE}.
 We refer to~\cite[Section~7]{Pphil} for a discussion of periodicity
theorems and their relevance for the contraherent cosheaf theory.

\subsection{Locally contraherent cosheaves}
 The phenomenon of dependence of the category of contramodules over
the coring $\C=\O(T\times_XT)$ over the ring $A=\O(T)$ on the affine
open covering $T\rarrow X$ used to construct it (mentioned in
Section~\ref{introd-contrah-as-contram} above) manifests itself in our
approach in an unexpected predicament of the \emph{contraherence
property} of a cosheaf of $\O_X$\+modules \emph{being not local}.
 So we have to deal with the \emph{locally contraherent cosheaves},
and the necessity to control the extension of this locality brings
the coverings back.

 Once an open covering in restriction to which our cosheaf becomes
contraherent is safely fixed, though, many other cosheaf properties
that we consider in this book become indeed local.
 And any locally contraherent cosheaf on a quasi-compact
semi-separated scheme has a finite \v Cech resolution by
contraherent cosheaves.

\subsection{Partially defined functors} \label{introd-partially-defined}
 One difference between homological theories developed in the settings
of exact and abelian categories is that whenever a functor between
abelian categories isn't exact, a similar functor between exact
categories will tend to have a shrinked domain.
 Any functor between abelian categories that has an everywhere defined
left or right derived functor will tend to be everywhere defined itself,
if only because one can always pass to the degree-zero cohomology of
derived category objects.
 Not so with exact categories, in which complexes may have no cohomology
objects in general.

 Hence the (sometimes annoying) necessity to deal with multitudes of
domains of definitions of various functors in our exposition.
 On the other hand, a functor with the domain consisting of adjusted
objects is typically exact on the exact subcategory where it is defined
(see~\cite[Section~5.6]{Pphil} for a discussion).

\subsection{Na\"\i ve co-contra correspondence}
 Another difference between the theory of comodules and contramodules
over corings as developed in~\cite{Psemi} and our present setting is
that in \emph{loc.\ cit.} we considered corings $\C$ over base rings $A$
of finite homological dimension.
 On the other hand, the ring $A=\O(T)$ constructed above has infinite
homological dimension in most cases, while the coring
$\C=\O(T\times_XT)$ can be said to have ``finite homological dimension
relative to~$A$''.
 In other words, one can say that the gluing procedure producing
(e.~g., quasi-compact semi-separated) schemes from affine ones has
finite homological dimension.
 For this reason, while the comodule-contramodule correspondence
theorem~\cite[Theorem~5.4]{Psemi} was stated for the derived categories
of the second kind, one of our most general co-contra correspondence
results in this book features an equivalence of the conventional
derived categories.
 We call it a ``na\"\i ve co-contra correspondence'' (with
an understanding that a ``non-na\"\i ve'' version would have
the coderived and contraderived categories involved, as per
Sections~\ref{introd-second-kind}\+-%
\ref{introd-covariant-serre-grothendieck-global} below).

\subsection{Application to the Lipman--Neeman extraordinary
inverse image} \label{introd-lipman-neeman}
 One application of the na\"\i ve co-contra correspondence is
a new construction of the extraordinary inverse image functor~$f^!$
on the derived categories of quasi-coherent sheaves for any morphism
of quasi-compact semi-separated schemes $f\:Y\rarrow X$.
 To be more precise, this is a construction of what we call
\emph{the Lipman--Neeman extraordinary inverse image functor}, that is
the right adjoint functor~$f^!$ to the derived direct image functor
$\boR f_*\:\sD(Y\qcoh)\rarrow\sD(X\qcoh)$ on the derived categories
of quasi-coherent sheaves (see the discussion below
in Section~\ref{introd-two-extraord-inverse-images}).
 In fact, we construct a \emph{right derived} functor $\boR f^!$, rather
than just a triangulated functor~$f^!$, as it was usually done
before~\cite[Example~4.2]{N-bb}, \cite[Chapter~4]{Lip}.
 See~\cite[Remark~6.1.4]{N-s} for a discussion; and notice that
the functor we denote by~$f^!$ is denoted by~$f^!$
in~\cite{N-bb} but by~$f^\times$ in~\cite{Lip,N-s}.

 To a morphism $f\:Y\rarrow X$ one assigns the direct and inverse image
functors $f_*\:Y\qcoh\rarrow X\qcoh$ and $f^*\:X\qcoh\rarrow Y\qcoh$
between the abelian categories of quasi-coherent sheaves on $X$ and~$Y$;
the functor~$f^*$ is left adjoint to the functor~$f_*$.
 To the same morphism, we also assign the direct image functor
$f_!\:Y\ctrh_\al\rarrow X\ctrh_\al$ between the exact categories of
antilocal contraherent cosheaves and the inverse image functor
$f^!\:X\lcth^\lin \rarrow Y\lcth^\lin$ between the exact categories
of locally injective locally contraherent cosheaves on $X$ and~$Y$.
 The functor $f^!$ is ``partially'' right adjoint to the functor~$f_!$.

 Passing to the derived functors, one obtains the adjoint
functors $\boR f_*\: \sD(Y\qcoh)\allowbreak\rarrow \sD(X\qcoh)$ and
$\boL f^*\: \sD(X\qcoh)\rarrow \sD(Y\qcoh)$ between the (conventional
unbounded) derived categories of quasi-coherent sheaves.
 We also obtain the adjoint functors $\boL f_!\: \sD(Y\ctrh)\rarrow
\sD(X\ctrh)$ and $\boR f^!\: \sD(X\ctrh)\rarrow\sD(Y\ctrh)$
between the derived categories of contraherent cosheaves on $X$ and~$Y$.

 The derived co-contra correspondence (for the conventional derived
categories) provides equivalences of triangulated categories
$\sD(X\qcoh)\simeq\sD(X\ctrh)$ and $\sD(Y\qcoh)\simeq\sD(Y\ctrh)$
transforming the direct image functor $\boR f_*$ into the direct image
functor $\boL f_!$.
 So the two inverse image functors $\boL f^*$ and $\boR f^!$ can be
viewed as the adjoints on the two sides to one and the same triangulated
functor of direct image.
 This finishes our construction of the triangulated functor
$f^!\:\sD(X\qcoh)\rarrow\sD(Y\qcoh)$ right adjoint to~$\boR f_*$.

\subsection{The calculus of tensor and Hom operations}
 As usually in homological algebra, the tensor product and Hom-type
operations on quasi-coherent sheaves and contraherent cosheaves
play an important role in our theory.
 First of all, under appropriate adjustedness assumptions one can assign
a contraherent cosheaf $\Cohom_X(\F,\P)$ to a quasi-coherent sheaf $\F$
and a contraherent cosheaf $\P$ over a scheme~$X$.
 This operation is the analogue of the tensor product of quasi-coherent
sheaves in the contraherent world.

 Secondly, to a quasi-coherent sheaf $\F$ and a cosheaf of
$\O_X$\+modules $\P$ over a scheme $X$ one can assign a cosheaf of
$\O_X$\+modules $\F\ot_{\O_X}\P$.
 Under our duality-analogy, this corresponds to taking the sheaf of
$\qHom$ from a quasi-coherent sheaf to a sheaf of $\O_X$\+modules.
 When the scheme $X$ is Noetherian, the sheaf $\F$ is coherent,
and the cosheaf $\P$ is contraherent, the cosheaf $\F\ot_{\O_X}\P$
is contraherent, too.
 Under some other assumptions, one can apply the (derived or underived)
\emph{contraherator} functor $\boL\Cr$ or $\Cr$ to the cosheaf
$\F\ot_{\O_X}\P$ to obtain the (complex of) contraherent cosheaves
$\F\ot_{X\ct}^\boL \P$ or $\F\ot_{X\ct}\P$.
 These are the analogues of the quasi-coherent internal Hom functor
$\qHom_{X\qc}$ on quasi-coherent sheaves, which can be obtained
by applying the coherator functor $\Qr$ \cite[Appendix~B]{TT} to
the $\qHom_{\O_X}({-},{-})$ sheaf.

 The remaining two operations are harder to come by.
 Modelled after the comodule-contramodule correspondence functors
$\Phi_\C$ and $\Psi_\C$ from~\cite{Psemi}, they play a similarly
crucial role in our present co-contra correspondence theory.
 Given a quasi-coherent sheaf $\F$ and a cosheaf of $\O_X$\+modules
$\P$ on a quasi-separated scheme $X$, one constructs a quasi-coherent
sheaf $\F\ocn_X\P$ on~$X$.
 Given two quasi-coherent sheaves $\F$ and $\cP$ on a quasi-separated
scheme, under certain adjustedness assumptions one can construct
a contraherent cosheaf $\fHom_X(\F,\cP)$.

\subsection{Coderived and contraderived categories}
\label{introd-second-kind}
 \emph{Derived categories of the second kind}, whose roots go back to
the work of Husemoller, Moore, and Stasheff on two kinds of
differential derived functors~\cite{HMS} and Hinich's paper about
DG\+coalgebras~\cite{Hin}, were introduced in fully developed form
in the author's monograph~\cite{Psemi} and memoir~\cite{Pkoszul}
(see also~\cite{PP2,Orl,EP,Pweak,Pfp,Pctrl,Psemten}).
 The most important representatives of this class of derived category
constructions are known as the \emph{coderived} and
the \emph{contraderived} categories; the difference between them
consists in the use of the closure with respect to infinite direct sums
in one case and with respect to infinite products in the other.

 The previous sentence describes what have been recently called
the derived categories of the second kind \emph{in the sense of 
Positselski}.
 Another approach, known under the name of the coderived and
contraderived categories \emph{in the sense of Becker}~\cite{Bec},
goes back to J\o rgensen~\cite{Jorg} and Krause~\cite{K-st}.
 It is still an open question whether the Positselski and Becker
approaches \emph{ever} disagree with each other within their common
domain of definition.
 They are known to agree under somewhat restrictive assumptions,
such as Noetherianity and finiteness of the Krull dimension.
 For all we know, the advantage of the Becker approach is that it
gives the ``correct'' categories, while the Positselski definitions
are technically easier to work with.
 Both the approaches are developed in this book, both in abstract
theory and in application to various categories of quasi-coherent
sheaves and contraherent cosheaves.
 We refer to~\cite[Examples~2.5(3) and~2.6(3)]{Pps},
\cite[Remark~9.2]{PS4}, and~\cite[Section~7]{Pksurv} for
a historical and philosophical discussion.

\subsection{Covariant Serre--Grothendieck duality over rings}
\label{introd-covariant-serre-grothendieck-affine}
 Here is a typical example of how derived categories of the second
kind occur.
 According to Iyengar and Krause~\cite{IK}, the homotopy category of
complexes of projective modules over a Noetherian commutative ring
with a dualizing complex is equivalent to the homotopy category of
complexes of injective modules.
 This theorem was extended to semi-separated Noetherian schemes with
dualizing complexes by Neeman~\cite{N-f} and Murfet~\cite{M-th} in
the following form: the derived category $\sD(X\qcoh_\fl)$ of
the exact category of flat quasi-coherent sheaves on such a scheme $X$
is equivalent to the homotopy category of injective quasi-coherent
sheaves $\Hot(X\qcoh^\inj)$.
 These results are known as the \emph{covariant Serre--Grothendieck
duality} theory.

 One would like to reformulate this equivalence so that it connects
certain derived categories of modules/sheaves, rather than just
subcategories of (co)resolutions.
 In other words, it would be nice to have some procedure assigning
complexes of projective, flat, and/or injective modules/sheaves to
arbitrary complexes.

 In the case of modules, the homotopy category of projectives is
identified with the contraderived category of the abelian category of
modules, while the homotopy category of injectives is equivalent to
the coderived category of the same abelian category.
 Hence the Iyengar--Krause result is interpreted as an instance of
the ``co-contra correspondence''---in this case, an equivalence between
the coderived and contraderived categories of the same abelian category.

\subsection{Covariant Serre-Grothendieck duality over schemes} \label{introd-covariant-serre-grothendieck-global}
 In the case of quasi-coherent sheaves, however, only a half of
the above picture remains true.
 The homotopy category of injectives $\Hot(X\qcoh^\inj)$ is still
equivalent to the coderived category of quasi-coherent sheaves
$\sD^\co(X\qcoh)$.
 But the attempt to similarly describe the derived category of flats
in terms of the whole abelian category of quasi-coherent sheaves runs
into the problem that the infinite products of quasi-coherent sheaves
are not exact, so the contraderived category construction does not
make sense for them (a detailed discussion can be found
in~\cite[Section~3]{Pphil}).

 This is where contraherent cosheaves come into play.
 The covariant Serre--Grothendieck duality for a nonaffine (but
semi-separated) scheme with a dualizing complex is an equivalence
of \emph{four} triangulated categories, rather than just two.
 In addition to the derived category of flat quasi-coherent sheaves
$\sD(X\qcoh_\fl)$ and the homotopy category of injective quasi-coherent
sheaves $\Hot(X\qcoh^\inj)$, there are also the homotopy category of
projective contraherent cosheaves $\Hot(X\ctrh_\prj)$ and
the derived category of locally injective contraherent cosheaves
$\sD(X\ctrh^\lin)$.

 Just as the homotopy category of injectives $\Hot(X\qcoh^\inj)$
is equivalent to the coderived category $\sD^\co(X\qcoh)$, the homotopy
category of projectives $\Hot(X\ctrh_\prj)$ is identified with
the contraderived category $\sD^\ctr(X\ctrh)$ of contraherent cosheaves.
 The equivalence between the two ``injective'' categories
$\sD^\co(X\qcoh)$ and $\sD(X\ctrh^\lin)$ does not depend on
the dualizing complex, and neither does the equivalence
between the two ``projective'' (or ``flat'') categories
$\sD(X\qcoh_\fl)$ and $\sD^\ctr(X\ctrh)$.
 The equivalences connecting the ``injective'' categories with
the ``projective'' ones do.

\subsection{Nonseparated schemes}
 So far, our discussion in this introduction was essentially limited
to semi-separated schemes.
 Such a restriction of generality is not really that natural or
desirable.
 Indeed, the affine plane with a double point (which is not
semi-separated) is no less a worthy object of study than the line
with a double point (which is).
 There is a problem, however, that appears to stand in the way of
a substantial development of the theory of contraherent cosheaves on
the more general kinds of schemes.

 The problem is that even when the quasi-coherent sheaves on schemes
remain well-behaved objects, the complexes of such sheaves may be
no longer so.
 More precisely, the trouble is with the derived functors of direct
image of complexes of quasi-coherent sheaves, which do not always have
good properties (such as locality along the base, etc.)
 Hence the common wisdom that for the more complicated schemes $X$ one
is supposed to consider complexes of sheaves of $\O_X$\+modules with
quasi-coherent cohomology sheaves, rather than complexes of
quasi-coherent sheaves as such (see, e.~g., \cite[Sections~2
and~6.2]{Rou}).
 As we do not know what is supposed to be either a second kind or
a contraherent analogue of the construction of the derived category of
complexes of sheaves of $\O_X$\+modules with quasi-coherent cohomology
sheaves, we have to restrict our exposition to, approximately, those
situations where the derived category of the abelian category of
quasi-coherent sheaves on $X$ is still a good category to work in.

 There are, basically, two such situations~\cite[Appendix~B]{TT}:
(1)~the quasi-compact semi-separated schemes and (2)~the Noetherian
or, sometimes, locally Noetherian schemes (which, while always
quasi-separated, do not have to be semi-separated).
 Accordingly, our exposition largely splits into two streams
corresponding to the situations~(1) and~(2), where different
techniques are applicable.
 The main difference in the generality level with the quasi-coherent
case is that in the contraherent context one often also needs to
assume one's schemes to have finite Krull dimension, in order to
use Raynaud and Gruson's homological dimension
results~\cite[Corollaire~II.3.2.7]{RG}.

 In particular, we prove the equivalence of the conventional derived
categories of quasi-coherent sheaves and contraherent cosheaves
not only for quasi-compact semi-separated schemes, but also, separately,
for all Noetherian schemes of finite Krull dimension.
 As to the covariant Serre--Grothendieck duality theorem, it remains
valid in the case of a non-semi-separated Noetherian scheme with
a dualizing complex in the form of an equivalence between two derived
categories of the second kind: the coderived category of quasi-coherent
sheaves $\sD^\co(X\qcoh)\simeq\Hot(X\qcoh^\inj)$ and the contraderived
category of contraherent cosheaves $\sD^\ctr(X\ctrh)\simeq
\Hot(X\ctrh_\prj)$.

\subsection{Projective contraherent cosheaves}
 Hartshorne's theory of injective quasi-coherent sheaves on locally
Noetherian schemes~\cite[Section~II.7]{Har}, including the assertions
that injectivity of quasi-coherent sheaves on locally Noetherian schemes
is a local property and such sheaves are flasque, as well as their
classification as direct sums of direct images from the embeddings of
the scheme points, is an important technical tool of the quasi-coherent
sheaf theory.

 We obtain dual-analogous results for projective locally cotorsion
contraherent cosheaves on locally Noetherian schemes, proving that
projectivity of locally cotorsion contraherent cosheaves on such
schemes is a local property, that such projective cosheaves are
coflasque, and classifying them as products of direct images from
the scheme points.
 Just as Hartshorne's theory is based on Matlis' classification
of injective modules over Noetherian rings~\cite{Mat} (cf., however,
\cite[Section~A.3]{EP}), our results on projective locally cotorsion
contraherent cosheaves are based on Enochs' classification of flat
cotorsion modules over commutative Noetherian rings~\cite{En}.

 We also prove that there are enough projective locally cotorsion
contraherent cosheaves on any locally Noetherian scheme.
 For locally contraadjusted contraherent cosheaves, the similar
assertion holds for all Noetherian schemes of finite Krull
dimension.
 On a quasi-compact semi-separated scheme, there are both enough
projective locally contraadjusted and projective locally cotorsion
contraherent cosheaves.

\subsection{Two extraordinary inverse images}
\label{introd-two-extraord-inverse-images}
 Before describing another application of our theory, let us have
a more detailed discussion of the extraordinary inverse images of
quasi-coherent sheaves.
 There are, in fact, \emph{two} different functors sometimes going by
the name of ``the functor~$f^!$'' in the literature.
 One of them, which we name after Lipman and Neeman, is simply
the functor right adjoint to the derived direct image $\boR f_*\:
\sD(Y\qcoh)\rarrow\sD(X\qcoh)$, which we were discussing above in
Section~\ref{introd-lipman-neeman}.
 Neeman proved its existence for an arbitrary morphism of quasi-compact
semi-separated schemes, using the techniques of compact generators
and Brown representability~\cite[Example~4.2]{N-bb}.
 Similar arguments apply in the case of an arbitrary morphism of
Noetherian schemes $f\:Y\rarrow X$ (see~\cite[Theorem~4.1.1]{Lip}
for full generality; and the discussion in~\cite{N-s}).

 The other functor, which we call \emph{the Hartshorne--Deligne 
extraordinary inverse image} and denote (to avoid ambiguity) by~$f^+$,
coincides with Lipman--Neeman's functor in the case of a proper
morphism~$f$.
 In the case when $f$~is an open embedding, on the other hand,
the functor~$f^+$ coincides with the conventional restriction
(inverse image) functor~$f^*$ (which is left adjoint to~$\boR f_*$,
rather than right adjoint).
 More generally, in the case of a smooth morphism~$f$ the functor~$f^+$
only differs from~$f^*$ by a dimensional shift and a top form line
bundle twist.
 This is the functor that was constructed in Hartshorne's
book~\cite{Har} and Deligne's appendix to it~\cite{Del} (hence the name).
 It is Hartshorne--Deligne's, rather than Lipman--Neeman's,
extraordinary inverse image functor that takes a dualizing complex
on $X$ to a dualizing complex on~$Y$, i.~e., $\D_Y^\bu=f^+\D_X^\bu$.

\subsection{The Hartshorne--Deligne extraordinary inverse image}
 It is an important idea apparently due to Gaitsgory~\cite{Gai,GR} that
the Hartshorne--Deligne extraordinary inverse image functor actually
acts between the coderived categories of quasi-coherent sheaves, rather
than between their conventional derived categories.
 To be sure, Hartshorne and Deligne only construct their functor for
bounded below complexes (for which there is no difference between
the coderived and derived categories).
 Gaitsgory and Rozenblyum's ``ind-coherent sheaves'' are closely
related to our coderived category of quasi-coherent sheaves.

 Concerning the Lipman--Neeman right adjoint functor~$f^!$, it can be
shown to exist on both the derived and (in the case of Noetherian
schemes) the coderived categories. 
 The problem arises when one attempts to define the functor~$f^+$ by
decomposing a morphism $f\:Y\rarrow X$ into proper morphisms~$g$ and
open embeddings~$h$ and subsequently composing the functors~$g^!$ for
the former with the functors~$h^*$ for the latter.
 It just so happens that the functor $\sD(X\qcoh)\rarrow\sD(Y\qcoh)$
obtained in this way depends on the chosen decomposition of
the morphism~$f$.
 This was demonstrated by Neeman in his
counterexample~\cite[Example~6.5]{N-bb}.

 The arguments of Deligne~\cite{Del}, if extended to the conventional
unbounded derived categories, break down on a rather subtle point:
while it is true that the restriction of an injective quasi-coherent
sheaf to an open subscheme of a Noetherian scheme remains injective,
the restriction of a homotopy injective complex of quasi-coherent
sheaves to such a subscheme may no longer be homotopy
injective~\cite{Bel}, \cite[Example~2.9]{Pal}.
 On the other hand, Deligne computes the Hom into the object produced
by his extraordinary inverse image functor from an arbitrary bounded
complex of coherent sheaves, which is essentially sufficient to make
a functor between the coderived categories well-defined, as these
are compactly generated by bounded complexes of coherents.

\subsection{Comparison with the conventional derived inverse image}
 The above discussion of the functor~$f^+$ is to be compared with
the remark that the conventional derived inverse image
functor~$\boL f^*$ is not defined on the coderived categories of
quasi-coherent sheaves (but only on their derived categories),
except in the case of a morphism~$f$ of finite flat
dimension~\cite{EP,Gai}.
 On the other hand, the conventional inverse image~$f^*$ is perfectly
well defined on the derived categories of flat quasi-coherent sheaves
(where one does not even need to resolve anything in order to
construct a triangulated functor).
 We show that the functor $f^*\:\sD(X\qcoh_\fl)\rarrow\sD(Y\qcoh_\fl)$
is transformed into the functor $f^+\:\sD^\co(X\qcoh)\rarrow
\sD^\co(Y\qcoh)$ by the equivalence of triangulated categories described
in Sections~\ref{introd-covariant-serre-grothendieck-affine}\+-%
\ref{introd-covariant-serre-grothendieck-global}.

\subsection{The Hartshorne--Deligne functor and contraherent cosheaves}
 Let us explain the connection between the Hartshorne--Deligne
extrardinary inverse image functor and our contraherent cosheaf theory.
 Given a morphism of Noetherian schemes $f\:Y\rarrow X$, just as
the direct image functor $\boR f_*\:\sD^\co(Y\qcoh)\rarrow
\sD^\co(X\qcoh)$ has a right adjoint functor $f^!\:\sD^\co(X\qcoh)
\rarrow\sD^\co(Y\qcoh)$, so does the direct image functor
$\boL f_!\:\sD^\ctr(Y\ctrh)\rarrow\sD^\ctr(X\ctrh)$ has a left
adjoint functor $f^*\:\sD^\ctr(X\ctrh)\rarrow\sD^\ctr(Y\ctrh)$.
 This is a contraherent analogue of the Lipman--Neeman extraordinary
inverse image functor for quasi-coherent sheaves.

 Now assume that~$f$ is a morphism of finite type and the scheme
$X$ has a dualizing complex~$\D_X^\bu$; set $\D_Y^\bu=f^+\D_X^\bu$.
 As we mentioned above in
Section~\ref{introd-covariant-serre-grothendieck-global}, the choice
of the dualizing complexes induces equivalences of triangulated
categories $\sD^\co(X\qcoh)\simeq\sD^\ctr(X\ctrh)$ and similarly
for~$Y$.
 Those equivalences of categories transform the functor $f^*\:
\sD^\ctr(X\ctrh)\rarrow\sD^\ctr(Y\ctrh)$ into a certain functor
$\sD^\co(X\qcoh)\rarrow\sD^\co(Y\qcoh)$.
 It is the latter functor that turns out to be isomorphic to
the Hartshorne--Deligne extraordinary inverse image functor
(which we denote here by~$f^+$).
 
 Two proofs of this result, working on slightly different generality
levels, are given in this book.
 One applies to ``compactifiable'' separated morphisms of finite type
between semi-separated Noetherian schemes and presumes comparison
with a Deligne-style construction of the functor~$f^+$, involving
a decomposition of the morphism~$f$ into an open embedding followed
by a proper morphism~\cite{Del,Conr}.
 The other one is designed for ``embeddable'' morphisms of finite
type between Noetherian schemes and the comparison with
a Hartshorne-style construction of~$f^+$ based on a factorization
of~$f$ into a finite morphism followed by a smooth one~\cite{Har}.

\subsection{Semiderived categories and semico-semicontra correspondence}
 In the final chapter, we work out a ``semi-infinite'' version of
the homological formalism of quasi-coherent sheaves and
contraherent cosheaves (see Preface to~\cite{Psemi} and
the introduction to~\cite{Pfp} for the related speculations,
and~\cite{Psemten} for a partial realization of the program).
 The setting is that of a quasi-compact semi-separated scheme $Y$
flatly fibered over a semi-separated Noetherian scheme $X$ with
a dualizing complex.
 So $\pi\:Y\rarrow X$ is a flat morphism of schemes.
 We construct the semiderived categories of quasi-coherent sheaves
and contraherent cosheaves on $Y$ relative to $X$ as mixtures of
the co/contraderived categories along the base ind-scheme and
the conventional derived categories along the fibers.
 Let us emphasize that the scheme $X$ is Noetherian, but the scheme $Y$
need not be.

 The idea is to put our two co-contra correspondence theorems
``on top of'' one another.
 The ``na\"\i ve co-contra correspondence'' along the fibers of~$\pi$
is mounted on top of the co-contra correspondence depending on
the dualizing complex on the base scheme~$X$.
 In other words, we join the equivalence of the conventional derived
categories of quasi-coherent sheaves and contraherent cosheaves on
a quasi-compact semi-separated scheme together with the equivalence
between the coderived category of quasi-coherent sheaves and
the contraderived category of contraherent cosheaves on
a semi-separated Noetherian scheme with a dualizing complex into
a single ``semimodule-semicontramodule correspondence'' theorem
claiming an equivalence between two semiderived categories.

 In the semi-infinite context, our approach in this book is both more
and less general than that of the book~\cite{Psemten}.
 It is much less general in that we only work with schemes, while
ind-schemes were the setting in~\cite{Psemten}; and it is much more
general in that only quasi-coherent sheaves were considered
in~\cite{Psemten}, while in this book we bring contraherent cosheaves
into the picture.
 But the approach is also more general in that the morphism~$\pi$ is
\emph{not} assumed to be affine in this book.
 The main results of~\cite{Psemten} were proved for flat affine
morphisms of ind-schemes, and an attempt to define the semiderived
category for a nonaffine flat morphism of ind-schemes was outlined
in~\cite[Appendix~A]{Psemten}.
 In this book, we take up the approach of~\cite[Appendix~A]{Psemten}
and develop it into a full-blown theory of semiderived categories
of quasi-coherent sheaves and contraherent cosheaves for a flat,
nonaffine morphism $\pi\:Y\rarrow X$, with a derived semico-semicontra
correspondence theorem as the main result.
 However, we only do it for schemes and not for ind-schemes.

 The discussion in the previous two paragraphs concerns an equivalence
between the semiderived categories of quasi-coherent sheaves and
contraherent cosheaves on the scheme~$Y$ (relative to~$X$), which is
a generalization of the equivalence between the coderived category of
quasi-coherent sheaves and the contraderived category of contraherent
cosheaves on~$X$.
 The latter equivalence forms the diagonal line in a big square diagram
of triangulated equivalences, involving also the (absolute) derived
categories of flat quasi-coherent sheaves and locally injective
contraherent cosheaves on~$X$.
 We work out an extension of the big square diagram to
the semi-infinite context as well, featuring the derived categories of
$X$\+flat quasi-coherent sheaves and $X$\+locally injective contraherent
cosheaves on $Y$, in addition to the semiderived categories.
 But that works only for affine morphisms~$\pi$.

\subsection{Future prospects and potential applications}
 To end, let us describe some prospects for future research and
applications of contraherent cosheaves.
 One of such applications is related to the $\bcD$\+-$\bOmega$
duality theory.
 Here $\bcD$ stands for the sheaf of rings of differential operators on
a smooth scheme and $\bOmega$ denotes the de~Rham DG\+algebra.
 The derived $\bcD$\+-$\bOmega$ duality, as formulated
in~\cite[Appendix~B]{Pkoszul} (see~\cite{Ryb,Prel} for further
developments), happens on two sides.
 On the ``co'' side, the functor $\qHom_{\O_X}(\bOmega,{-})$ takes
right $\bcD$\+modules to right DG\+modules over $\bOmega$, and there is
the adjoint functor ${-}\ot_{\O_X}\bcD$.
 These functors induce an equivalence between the derived category of
$\bcD$\+modules and the coderived category of DG\+modules
over~$\bOmega$.

 On the ``contra'' side, over an affine scheme $X$, the functor
$\bOmega(X)\ot_{\O(X)}{-}$ takes left $\bcD$\+modules to
left DG\+modules over $\bOmega$, and the adjoint functor is
$\Hom_{\O(X)}(\bcD(X),{-})$.
 The induced equivalence is between the derived category of
$\bcD(X)$\+modules and the contraderived category of DG\+modules
over $\bOmega(X)$.
 One would like to extend the ``contra'' side of the story to nonaffine
schemes using the contraherent cosheafification, as opposed to
the quasi-coherent sheafification on the ``co'' side.
 The need for the contraherent cosheaves arises, once again, because
the contraderived category construction does not make sense for
quasi-coherent sheaves.
 The contraherent cosheaves are designed for being plugged into it.
 This prospect for an application of the contraherent cosheaf theory
to $\bcD$\+-$\bOmega$ duality was actually realized in the recent
preprint~\cite{Pdomc}.

 Another direction in which we would like to extend the theory
presented in this book is that of Noetherian formal schemes, and more
generally, ind-schemes of ind-finite or even ind-infinite type.
 The idea is to define an exact category of contraherent cosheaves
of contramodules over a formal scheme that would serve as
the natural ``contra''-side counterpart to the abelian category
of quasi-coherent \emph{torsion} sheaves on the ``co'' side.
 (For a taste of the contramodule theory over complete Noetherian
rings, the reader is referred to~\cite[Appendix~B]{Pweak};
see also~\cite{Pmgm,Pcta,Prev}.)

 In particular, developing the na\"\i ve co-contra correspondence
theory in this context should lead one to a global (nonaffine
scheme) generalization of the MGM duality/equivalence
of~\cite{Mat2,GM,DG,PSY} as interpreted in~\cite{Pmgm}.
 We also expect generalizations of the covariant Serre--Grothendieck
duality and the semico-semicontra correspondence involving
contraherent cosheaves of contramodules.

 Furthermore, one might hope to join the de~Rham DG\+module and
the contramodule stories together in a single theory by considering
contraherent cosheaves of DG\+contramodules over the de~Rham--Witt
complex, with an eye to applications to crystalline sheaves and
the $p$\+adic Hodge theory.
 As compared so such high hopes, our real advances in this book
are quite modest.

\subsection{Preparatory work: Cotorsion pairs in contramodule
categories}
 The very possibility of having a meaningful theory of contraherent
cosheaves on a scheme rests on there being enough contraadjusted
or cotorsion modules over a commutative or associative ring.
 Thus the first problem one encounters when trying to build up
the theories of contraherent cosheaves of contramodules is
the necessity of developing applicable techniques for constructing
flat, cotorsion, contraadjusted, and very flat contramodules
to be used in the (co)resolutions.
 In other words, in order to develop the theory of contraherent
cosheaves of contramodules, one needs to have certain \emph{complete
cotorsion pairs} in contramodule categories, such as the flat and/or
the very flat complete cotorsion pair, at one's disposal.

 In this book, we partially solve this problem by constructing enough
flat, cotorsion, contraadjusted, and very flat contramodules over
a Noetherian ring in the adic topology, enough flat and cotorsion
contramodules over a pro-Noetherian topological ring of totally finite
Krull dimension, and also enough contraadjusted and very flat
contramodules over the projective limit of a sequence of commutative
rings and surjective morphisms between them with finitely generated
kernel ideals.
 The constructions are based on the existence of enough objects
of the respective classes in the conventional categories of modules,
which is used as a black box, and also on the older, more explicit
construction~\cite{Xu} of cotorsion coresolutions of modules over
commutative Noetherian rings of finite Krull dimension.

 This provides the necessary background for possible definitions of
locally contraadjusted or locally cotorsion contraherent cosheaves
of contramodules over Noetherian formal schemes or ind-Noetherian
ind-schemes of totally finite Krull dimension, and also of
locally contraadjusted contraherent cosheaves over ind-schemes of
a more general nature.
 An even more general contemporary-style set-theoretical approach to
constructing complete cotorsion pairs in contramodules categories,
generalizing the proofs in~\cite{ET,BBE} and based on the small object
argument, is worked out in the paper~\cite{PR}.
 But the more explicit elementary approaches to constructing complete
cotorsion pairs tend to provide more information
(cf.~\cite[Section~5]{PSl1} and~\cite{Pctrl,Pal}).

 This kind of work, intended to prepare ground for future theory
building, is relegated to appendices in this book.
 There we also construct the derived co-contra correspondence
(an equivalence between the coderived category of torsion modules
and the contraderived category of contramodules) over
affine Noetherian formal schemes and ind-affine ind-Noetherian
ind-schemes with dualizing complexes.
 Another appendix is devoted to the co-contra correspondence
over noncommutative quasi-compact semi-separated stacks (otherwise
known as flat corings~\cite{KR,KR2,Pflcc}).

\subsection*{Acknowledgement}
 Among the many people whose questions, remarks, and suggestions
contributed to the present research in the early stages of its
development, I~should mention Sergey Arkhipov, Roman Bezrukavnikov,
Alexander Polishchuk, Sergey Rybakov, Alexander Efimov, Amnon Neeman,
Henning Krause, Mikhail Bondarko, Alexey Gorodentsev, Paul Bressler,
Vladimir Hinich, Greg Stevenson, Michael Finkelberg, and Pierre Deligne.
 Parts of an early version of book were written when I was visiting
the University of Bielefeld, and I~want to thank Collaborative
Research Center~701 and Prof.\ Krause for the invitation.
 In its early stages, this work was supported in part by RFBR grants
in Moscow and by the Grant Agency of the Czech Republic under the grant
P201/12/G028 in Brno.

 In later years of my work on this project, I~benefited a lot from
conversations with Vladimir Ivchenko, Alexander Sl\'avik, Michal Hrbek,
Simone Virili, Jan \v St\!'ov\'\i\v cek, Jan Trlifaj, and
Silvana Bazzoni.
 In particular, the conversation with Jan \v St\!'ov\'\i\v cek was
very helpful when I~was working on
Section~\ref{finite-dim-morphisms-derived-inverse-subsect}.
 Jan Trlifaj also always encouraged me to finish my work on this 
manuscript and publish it as a book.
 The author was supported by the GA\v CR grant 23-05148S and
the Institute of Mathematics, Czech Academy of Sciences (RVO:~67985840)
on the final stages of this project's development. {\uchyph=0\par}

\Section{Contraadjusted and Cotorsion Modules}

\subsection{Contraadjusted and very flat modules}
\label{very-eklof-trlifaj-subsect}
 Let $R$ be a commutative ring.
 We will say that an $R$\+module $P$ is \emph{contraadjusted} if
the $R$\+module $\Ext_R^1(R[r^{-1}],P)$ vanishes for every element
$r\in R$.
 An $R$\+module $F$ is called \emph{very flat} if one has
$\Ext_R^1(F,P)=0$ for every contraadjusted $R$\+module~$P$.

 By the definition, any injective $R$\+module is contraadjusted and
any projective $R$\+module is very flat.
 Notice that the projective dimension of the $R$\+module $R[r^{-1}]$
never exceeds~$1$, as it has a natural two-term free resolution
$0\rarrow\bigoplus_{n=0}^\infty R\rarrow\bigoplus_{n=0}^\infty R
\rarrow R[r^{-1}]\rarrow 0$,
$$
 0\lrarrow\bigoplus\nolimits_{n=0}^\infty Rf_n\lrarrow
 \bigoplus\nolimits_{n=0}^\infty Re_n\lrarrow R[r^{-1}]\lrarrow0,
$$
with the differentials taking the basis vector~$f_n$ to $e_n-re_{n+1}$
and the basis vector~$e_n$ to~$r^{-n}$ for every $n\ge0$.
 It follows that any quotient module of a contraadjusted module is
contraadjusted, and one has $\Ext_R^{>0}(F,P)=0$ for any very flat
$R$\+module $F$ and contraadjusted $R$\+module~$P$.

 Computing the $\Ext^1$ in terms of the above resolution, one can
characterize contraadjusted $R$\+modules more explicitly as follows.
 An $R$\+module $P$ is contraadjusted if and only if for any
sequence of elements $p_0$, $p_1$, $p_2$,~$\dotsc\in P$ and $r\in R$
there exists a (not necessarily unique) sequence of elements
$q_0$, $q_1$, $q_2$,~$\dotsc\in P$ such that $q_i=p_i+rq_{i+1}$
for all $i\ge0$ \,\cite[Lemma~5.1]{ST}, \cite[Lemma~2.1(a)]{Pcta}.

 Furthermore, the projective dimension of any very flat module
is equal to~$1$ or less.
 Indeed, any $R$\+module $M$ has a two-term coresolution by
contraadjusted modules, which can be used to compute $\Ext^*_R(F,M)$
for a very flat $R$\+module~$F$.
 (The converse assertion is \emph{not} true, however: a flat
$R$\+module of projective dimension~$1$ need not be very flat.
 See Examples~\ref{q-over-z-not-very-flat-ex}\+-%
\ref{baer-specker-not-very-flat-ex} below.)

 It is also clear that the classes of contraadjusted and very flat
modules are closed under extensions.
 Besides, the class of very flat modules is closed under the passages
to the kernel of a surjective morphism (i.~e., the kernel of
a surjective morphism of very flat $R$\+modules is very flat).
 In addition, we notice that the class of contraadjusted $R$\+modules
is closed under infinite products, while the class of very flat
$R$\+modules is closed under infinite direct sums.
 We refer to~\cite[paragraph after Lemma~2.4]{Pal}
or~\cite[Theorem~4.5]{Pphil} for some further details.

\begin{thm}  \label{eklof-trlifaj-very}
\textup{(a)} Any $R$\+module can be embedded into a contraadjusted
$R$\+module in such a way that the quotient module is very flat. \par
\textup{(b)} Any $R$\+module admits a surjective map onto it from
a very flat $R$\+module such that the kernel is contraadjusted.
\end{thm}

\begin{proof}
 Both assertions follow from the results of Eklof and
Trlifaj~\cite[Theorem~10]{ET}
(see Section~\ref{small-object-argument-subsect} for generalizations).
 It suffices to point out that all the $R$\+modules of the form
$R[r^{-1}]$ form a set rather than a proper class.
 For the reader's convenience and our future use, parts of
the argument from~\cite{ET} are reproduced below.
\end{proof}

 The following definition will be used in the sequel.
 Let $\sA$ be an abelian category with exact functors of filtered
inductive limit, and let $\sC\sub\sA$ be a class of objects.
 An object $X\in\sA$ is said to be a \emph{transfinitely iterated
extension} of objects from $\sC$ if there exist a well-ordered set
$\Gamma$ and a family of subobjects $X_\gamma\sub X$, \ $\gamma\in
\Gamma$, such that $X_\delta\sub X_\gamma$ whenever $\delta<\gamma$,
the union (inductive limit) $\varinjlim_{\gamma\in\Gamma} X_\gamma$
of all $X_\gamma$ coincides with~$X$, and the quotient objects
$X_\gamma/\varinjlim_{\delta<\gamma}X_\delta$ belong to~$\sC$ for all
$\gamma\in\Gamma$ (cf.\ Section~\ref{quasi-compact-quasi-coherent}).

 By~\cite[Lemma~1]{ET} or Lemma~\ref{eklof-lemma-general}, any
transfinitely iterated extension of the $R$\+modules $R[r^{-1}]$,
with arbitrary $r\in R$, is a very flat $R$\+module.
 Proving a converse assertion will be one of our goals.
 The following lemma is a particular case of~\cite[Theorem~2]{ET}.

\begin{lem}
 Any $R$\+module can be embedded into a contraadjusted $R$\+module
in such a way that the quotient module is a transfinitely iterated
extension of the $R$\+modules $R[r^{-1}]$.
\end{lem}

\begin{proof}
 The proof is a set-theoretic argument based on the fact that
the Cartesian square of any infinite cardinality $\lambda$ is
equicardinal to~$\lambda$.
 In our case, let $\lambda$ be any infinite cardinality no smaller
than the cardinality of the ring~$R$.
 For an $R$\+module $L$ of the cardinality not exceeding~$\lambda$
and an $R$\+module $M$ of the cardinality~$\mu$, the set
$\Ext^1_R(L,M)$ has the cardinality at most $\mu^\lambda$, as one
can see by computing the $\Ext^1$ in terms of a projective resolution
of the first argument.
 In particular; set $\aleph=2^\lambda$; then for any $R$\+module $M$
of the cardinality not exceeding~$\aleph$ and any $r\in R$
the cardinality of the set $\Ext^1_R(R[r^{-1}],M)$ does not
exceed~$\aleph$, either.

 Let $\beth$ be the smallest cardinality that is larger than~$\aleph$
and let $\Delta$ be the smallest ordinal of the cardinality~$\beth$.
 Notice that the natural map $\varinjlim_{\delta\in\Delta}
\Ext^*_R(L,Q_\delta)\rarrow\Ext^*_R(L\;\varinjlim_{\delta\in\Delta}
Q_\delta)$ is an isomorphism for any $R$\+module $L$ of the cardinality
not exceeding~$\lambda$ (or even~$\aleph$) and any inductive system
of $R$\+modules $Q_\delta$ indexed by~$\Delta$.
 Indeed, the functor of filtered inductive limit of abelian groups
is exact and the natural map $\varinjlim_{\delta\in\Delta}
\Hom_R(L,Q_\delta)\rarrow \Hom_R(L\;\varinjlim_{\delta\in\Delta}
Q_\delta)$ is an isomorphism for any $R$\+module $L$ of
the cardinality not exceeding~$\aleph$.
 The latter assertion holds because the image of any map of sets
$L\rarrow\Delta$ is contained in a proper initial segment
$\{\delta'\mid\delta'<\delta\}\sub\Delta$ for some $\delta\in\Delta$.

 We proceed by transfinite induction on $\Delta$ constructing for every
element $\delta\in\Delta$ an $R$\+module $P_\delta$ and a well-ordered
set $\Gamma_\delta$.
 For any $\delta'<\delta$, we will have an embedding of $R$\+modules
$P_{\delta'} \rarrow P_\delta$ (such that the three embeddings form
a commutative diagram for any three elements $\delta''<\delta'<\delta$)
and an ordered embedding $\Gamma_{\delta'}\sub\Gamma_\delta$ (making
$\Gamma_{\delta'}$ an initial segment of $\Gamma_\delta$).
 Furthermore, for every element $\gamma\in\Gamma_\delta$, a particular
extension class $c(\gamma,\delta)\in\Ext^1_R(R[r(\gamma)^{-1}],
P_\delta)$, where $r(\gamma)\in R$, will be defined.
 For every $\delta'<\delta$ and $\gamma\in\Gamma_{\delta'}$, the class
$c(\gamma,\delta)$ will be equal to the image of the class
$c(\gamma,\delta')$ with respect to the natural map
$\Ext^1_R(R[r(\gamma)^{-1}],P_{\delta'})\rarrow
\Ext^1_R(R[r(\gamma)^{-1}],P_\delta)$ induced by the embedding
$P_{\delta'}\rarrow P_\delta$.

 At the starting point $0\in\Delta$, the module $P_0$ is our original
$R$\+module $M$ and the set $\Gamma_0$ is the disjoint union of all
the sets $\Ext^1_R(R[r^{-1}],M)$ with $r\in R$, endowed with an arbitrary
well-ordering.
 The elements $r(\gamma)$ and the classes $c(\gamma,0)$ for $\gamma\in
\Gamma_0$ are defined in the obvious way.

 Given $\delta=\delta'+1\in\Delta$ and assuming that the $R$\+module
$P_{\delta'}$ and the set $\Gamma_{\delta'}$ have been constructed
already, we produce the module $P_\delta$ and the set $\Gamma_\delta$
as follows.
 It will be clear from the construction below that $\delta'$ is always
smaller than or equal to the well-ordering type of the set
$\Gamma_{\delta'}$.
 If the two well-ordering types are equal to each other, we simply put
$P_\delta=P_{\delta'}$ and $\Gamma_\delta=\Gamma_{\delta'}$.
 In the case of a strict inequality, there is a unique element
$\gamma_{\delta'}\in\Gamma_{\delta'}$ corresponding to
the ordinal~$\delta'$ (i.~e., such that the well-ordering type of
the subset all the elements in $\Gamma_{\delta'}$ that are smaller
than~$\gamma_{\delta'}$ is equivalent to~$\delta'$).
 Then the construction works as follows.

 Define the $R$\+module $P_\delta$ as the middle term of the extension
corresponding to the class $c(\gamma_{\delta'},\delta')\in
\Ext^1_R(R[r(\gamma_{\delta'})^{-1}],P_{\delta'})$.
 There is a natural embedding $P_{\delta'}\rarrow P_\delta$, as required.
 Set $\Gamma_\delta$ to be the disjoint union of $\Gamma_{\delta'}$
and the sets $\Ext^1_R(R[r^{-1}],P_{\delta})$ with $r\in R$,
well-ordered so that $\Gamma_{\delta'}$ is an initial segment,
while the well-ordering of the remaining elements is chosen arbitrarily.
 The elements $r(\gamma)$ for $\gamma\in\Gamma_{\delta'}$ have been
defined already on the previous steps and the classes $c(\gamma,\delta)$
for such~$\gamma$ are defined in the unique way consistent with
the previous step, while for the remaining $\gamma\in
\Gamma_\delta\setminus\Gamma_{\delta'}$ these elements and
classes are defined in the obvious way.

 When $\delta$ is a limit ordinal, set $P_\delta=\varinjlim_{\delta'
<\delta}P_\delta'$.
 Let $\Gamma_\delta$ be the disjoint union of the set
$\bigcup_{\delta'<\delta}\Gamma_{\delta'}$ and the sets
$\Ext^1_R(R[r^{-1}],P_{\delta})$ with $r\in R$, well-ordered so that
$\bigcup_{\delta'<\delta}\Gamma_{\delta'}$ is an initial segment.
 The elements $r(\gamma)$ and $c(\gamma,\delta)$ for
$\gamma\in\Gamma_\delta$ are defined as above.

 Arguing by transfinite induction, one easily concludes that
the cardinality of the $R$\+module $P_\delta$ never exceeds~$\aleph$
for $\delta\in\Delta$, and neither does the cardinality of
the set $\Gamma_\delta$.
 It follows that the well-ordering type of the set $\Gamma =
\bigcup_{\delta\in\Delta}\Gamma_\delta$ is equal to~$\Delta$.
 So for every $\gamma\in\Gamma$ there exists $\delta\in\Delta$ such
that $\gamma=\gamma_\delta$.

 Set $P=\varinjlim_{\delta\in\Delta}P_\delta$.
 By construction, there is a natural embedding of $R$\+modules
$M\rarrow P$ and the cokernel is a transfinitely iterated extension
of the $R$\+modules $R[r^{-1}]$.
 As every class $c\in\Ext^1_R(R[r^{-1}],P_\delta)$ corresponds to
an element $\gamma\in\Gamma_\delta$, has the corresponding ordinal
$\delta'\in\Delta$ such that $\gamma=\gamma_{\delta'}$, and dies in
$\Ext^1_R(R[r^{-1}],P_{\delta'+1})$, we conclude that
$\Ext^1_R(R[r^{-1}],P)=0$.
\end{proof}

\begin{lem} \label{eklof-trlifaj-cta}
 Any $R$\+module admits a surjective map onto it from a transfinitely
iterated extension of the $R$\+modules $R[r^{-1}]$ such that
the kernel is contraadjusted.
\end{lem}

\begin{proof}
 The proof follows the second half of the proof of Theorem~10
in~\cite{ET}.
 Specifically, given an $R$\+module $M$, pick a surjective map onto
it from a free $R$\+module~$L$.
 Denote the kernel by $K$ and embed it into a contraadjusted
$R$\+module $P$ so that the quotient module $Q$ is a transfinitely
iterated extension of the $R$\+modules $R[r^{-1}]$.
 Then the fibered coproduct $F$ of the $R$\+modules $L$ and $P$
over $K$ is an extension of the $R$\+modules $Q$ and~$L$.
 It also maps onto $M$ surjectively with the kernel~$P$.
 This argument is known as the \emph{Salce lemma}~\cite{Sal}; cf.\
Lemma~\ref{salce-lemma}(a).
\end{proof}

 Both assertions of Theorem~\ref{eklof-trlifaj-very} are now proved.

\begin{cor}  \label{very-flat-transfinite}
 An $R$\+module is very flat if and only if it is a direct summand
of a transfinitely iterated extension of the $R$\+modules $R[r^{-1}]$.
\end{cor}

\begin{proof}
 The ``if'' part has been explained already (by references
to~\cite[Lemma~1]{ET} or Lemma~\ref{eklof-lemma-general}).
 Let us prove the ``only if''.
 Given a very flat $R$\+module $F$, present it as the quotient module
of a transfinitely iterated extension $E$ of the $R$\+modules
$R[r^{-1}]$ by a contraadjusted $R$\+module~$P$.
 Since $\Ext^1_R(F,P)=0$, we can conclude that $F$ is a direct summand
of~$E$.
\end{proof}

 In particular, we have proved that any very flat $R$\+module is flat.

\begin{cor}  \label{very-rel-proj-inj}
\textup{(a)} Any very flat $R$\+module can be embedded into
a contraadjusted very flat $R$\+module in such a way that
the quotient module is very flat. \par
\textup{(b)} Any contraadjusted $R$\+module admits a surjective map
onto it from a very flat contraadjusted $R$\+module such that
the kernel is contraadjusted.
\end{cor}

\begin{proof}
 Follows from Theorem~\ref{eklof-trlifaj-very} and the fact that
the classes of contraadjusted and very flat $R$\+modules are closed
under extensions.
\end{proof}

\subsection{Affine geometry of contraadjusted and very flat modules}
\label{affine-geometry-subsect}
 The results of this section form the module-theoretic background of
our main definitions and constructions in
Chapters~\ref{contraherent-sect}\+-\ref{loc-contra-sect}.

 The following ``homological formula'' will be used in the proofs of
Lemmas~\ref{very-tensor-hom}\+-\ref{very-scalars-veryflat-case}
and~\ref{cotors-hom}\+-\ref{cotors-coexten}.
 Let $R$ and $S$ be two associative rings.
 Let $F$ be an $R$\+$S$\+bimodule, $G$ be a left $S$\+module, and
$P$ be a left $R$\+module.
 Then, given an integer $n\ge0$, there is a natural (adjunction)
isomorphism of abelian groups {\hbadness=1500
\begin{equation} \label{ext-tor-adjunction}
 \Ext^n_S(G,\Hom_R(F,P))\simeq\Ext_R^n(F\ot_SG\;P)
\end{equation}
provided} that $\Ext^i_R(F,P)=0=\Tor_i^S(F,G)$ for all $0<i\le n$
(cf.~\cite[formula~(4) in Section~XVI.4]{CaE}.
 We will mostly use the isomorphism~\eqref{ext-tor-adjunction}
for $n=1$.

\begin{lem} \label{very-tensor-hom}
 \textup{(a)} The class of very flat $R$\+modules is closed with
respect to the tensor products over~$R$. \par
 \textup{(b)} For any very flat $R$\+module $F$ and contraadjusted
$R$\+module $P$, the $R$\+module\/ $\Hom_R(F,P)$ is contraadjusted.
\end{lem}

\begin{proof}
 One approach is to prove both assertions simultaneously using
the adjunction isomorphism $\Ext_R^1(F\ot_R G\;P)\simeq
\Ext_R^1(G,\Hom_R(F,P))$, which clearly holds for any $R$\+module $G$,
any very flat $R$\+module $F$, and contraadjusted $R$\+module~$P$,
and raising the generality step by step.
 Since $R[r^{-1}]\ot_R R[s^{-1}]\simeq R[(rs)^{-1}]$, it follows
that the $R$\+module $\Hom_R(R[s^{-1}],P)$ is contraadjusted for any
contraadjusted $R$\+module $P$ and $s\in R$.
 Using the same adjunction isomorphism, one then concludes that
the $R$\+module $R[s^{-1}]\ot_R G$ is very flat for any very flat
$R$\+module~$G$.
 From this one can deduce in full generality the assertion~(b), and
then the assertion~(a).

 Alternatively, one can use the full strength of
Corollary~\ref{very-flat-transfinite} and check that the tensor
product of two transfinitely iterated extensions of flat modules is
a transfinitely iterated extension of the pairwise tensor products.
 Then deduce (b) from~(a) using the same adjunction isomorphism.
\end{proof}

\begin{lem} \label{very-scalars-always}
 Let $f\:R\rarrow S$ be a homomorphism of commutative rings.  Then \par
\textup{(a)} any contraadjusted $S$\+module is also a contraadjusted
$R$\+module in the $R$\+module structure obtained by the restriction
of scalars via~$f$; \par
\textup{(b)} if $F$ is a very flat $R$\+module, then the $S$\+module
$S\ot_R F$ obtained by the extension of scalars via~$f$
is also very flat; \par
\textup{(c)} if $F$ is a very flat $R$\+module and $Q$ is
a contraadjusted $S$\+module, then\/ $\Hom_R(F,Q)$ is also
a contraadjusted $S$\+module; \par
\textup{(d)} if $F$ is a very flat $R$\+module and $G$ is a very flat
$S$\+module, then $F\ot_R G$ is also a very flat $S$\+module.
\end{lem}

\begin{proof}
 Part~(a): one has $\Ext_R^*(R[r^{-1}],Q) \simeq
\Ext_S^*(S[f(r)^{-1}],Q)$ for any $S$\+module $Q$ and $r\in R$.
 Part~(b) follows from part~(a) in view of the isomorphism
$\Ext^1_S(S\ot_RF\;Q)\simeq\Ext^1_R(F,Q)$ for any flat $R$\+module $F$
and any $S$\+module~$Q$; or alternatively, (b)~follows from
Corollary~\ref{very-flat-transfinite}.
 To prove part~(c), notice that $\Hom_R(F,Q)\simeq\Hom_S(S\ot_RF\;Q)$
and use part~(b) together with Lemma~\ref{very-tensor-hom}(b)
(applied to the ring~$S$).
 Similarly, part~(d) follows from part~(b) and
Lemma~\ref{very-tensor-hom}(a).
\end{proof}

\begin{lem}  \label{very-scalars-veryflat-case}
 Let $f\:R\rarrow S$ be a homomorphism of commutative rings such that
the localization $S[s^{-1}]$ is a very flat $R$\+module for
any element $s\in S$.
 Then \par
\textup{(a)} the $S$\+module\/ $\Hom_R(S,P)$ obtained by the coextension
of scalars via~$f$ is contraadjusted for any contraadjusted
$R$\+module~$P$; \par
\textup{(b)} any very flat $S$\+module is also a very flat $R$\+module
(in the $R$\+module structure obtained by the restriction of scalars
via~$f$); \par
\textup{(c)} the $S$\+module\/ $\Hom_R(G,P)$ is contraadjusted for any
very flat $S$\+module $G$ and contraadjusted $R$\+module~$P$.
\end{lem}

\begin{proof} \emergencystretch=0em \hfuzz=1.1pt
 Part~(a): one has $\Ext_S^1(S[s^{-1}],\Hom_R(S,P)) \simeq
\Ext_R^1(S[s^{-1}],P)$ for any $R$\+mod\-ule $P$ such that
$\Ext_R^1(S,P)=0$ and any $s\in S$.
 Part~(b) follows from part~(a) in view of the isomorphism
$\Ext^1_R(G,P)\simeq\Ext^1_S(G,\Hom_R(S,P))$ for any $S$\+module $G$
and any contraadjusted $R$\+module~$P$; or alternatively,
(b)~follows from Corollary~\ref{very-flat-transfinite}.
 Part~(c): one has $\Ext_S^1(S[s^{-1}],\Hom_R(G,P)) \simeq
\Ext_R^1(G\ot_S S[s^{-1}]\; P)$.
 By Lemma~\ref{very-tensor-hom}(a), the $S$\+module $G\ot_S S[s^{-1}]$
is very flat; by part~(b), it is also a very flat $R$\+module; so
the desired vanishing follows.
\end{proof}

 A homomorphism of commutative rings $f\:R\rarrow S$ is said to be
\emph{very flat} if it satisfies the assumption of
Lemma~\ref{very-scalars-veryflat-case}, i.~e., if the $R$\+module
$S[s^{-1}]$ is very flat for every element $s\in S$.
 In this case, one also says that $S$ is a \emph{very flat
commutative $R$\+algebra}.
 This terminology was introduced in~\cite[Sections~0.6 and~9]{PSl1}.

 Notice that the condition that $S$ be a very flat $R$\+module is
\emph{not} sufficient for the validity of
Lemma~\ref{very-scalars-veryflat-case}.
 A converse assertion to Lemma~\ref{very-scalars-veryflat-case}(b)
can be found in~\cite[Lemma~9.3(a)]{PSl1} or~\cite[Lemma~6.1]{Pphil};
for a counterexample, see~\cite[Example~9.7]{PSl1}.
 Nevertheless, it obviously follows from the next lemma that $S$ is
a very flat $R$\+algebra in its assumptions.

\begin{lem}  \label{very-open-embedding}
 Let $R\rarrow S$ be a homomorphism of commutative rings such that
the related morphism of affine schemes\/ $\Spec S\rarrow\Spec R$
is an open embedding.
 Then $S$ is a very flat $R$\+module.
\end{lem}

\begin{proof}
 The open subset $\Spec S\sub\Spec R$, being quasi-compact, can be
covered by a finite number of principal affine open subsets
$\Spec R[r_\alpha^{-1}] \sub \Spec R$, where $\alpha=1$,~\dots,~$N$.
 The \v Cech sequence
\begin{multline} \label{cech-principal} \textstyle
 0\lrarrow S \lrarrow \bigoplus_\alpha R[r_\alpha^{-1}] \lrarrow
 \bigoplus_{\alpha<\beta} R[(r_\alpha r_\beta)^{-1}] \\ \lrarrow
 \dotsb\lrarrow R[(r_1\dotsm r_N)^{-1}] \lrarrow 0
\end{multline}
is an exact sequence of $S$\+modules, since is its localization by every
element~$r_\alpha$ is a contractible complex of $S$\+modules.
 It remains to recall that the class of very flat $R$\+modules is closed
under the passage to the kernels of surjective morphisms.
\end{proof}

\begin{cor} \label{very-open-cor}
 The following assertions hold in the assumptions of
Lemma~\textup{\ref{very-open-embedding}}. \par
\textup{(a)} The $S$\+module $S\ot_R F$ is very flat for any very flat
$R$\+module~$F$. \par
\textup{(b)} An $S$\+module $G$ is very flat if and only if it is
very flat as an $R$\+module. \par
\textup{(c)} The $S$\+module\/ $\Hom_R(S,P)$ is contraadjusted for any
contraadjusted $R$\+mod\-ule~$P$. \par
\textup{(d)} An $S$\+module $Q$ is contraadjusted if and only if it is
contraadjusted as an $R$\+mod\-ule.
\end{cor}

\begin{proof}
 Part~(a) is a particular case of Lemma~\ref{very-scalars-always}(b).
 Part~(b): if $G$ is a very flat $S$\+module, then it is also very flat
as an $R$\+module by Lemma~\ref{very-scalars-veryflat-case}(b), which
is applicable in view of Lemma~\ref{very-open-embedding}.
 Conversely, if $G$ is very flat as an $R$\+module, then
$G\simeq S\ot_R G$ is also a very flat $S$\+module
by part~(a).

 Part~(c) follows from Lemma~\ref{very-scalars-veryflat-case}(a),
which is applicable by Lemma~\ref{very-open-embedding}.
 Part~(d): if $Q$ is a contraadjusted $S$\+module, then it is also
contraadjusted as an $R$\+module by Lemma~\ref{very-scalars-always}(a).
 Conversely, for any $S$\+module $Q$ there are natural isomorphisms
of $S$\+modules $Q\simeq\Hom_S(S,Q)\simeq\Hom_S(S\ot_RS\;Q)\simeq
\Hom_R(S,Q)$; and if $Q$ is contraadjusted as an $R$\+module, then
it is also a contraadjusted $S$\+module by part~(c).
\end{proof}

\begin{lem} \label{very-open-covering}
 Let $R\rarrow S_\alpha$, \ $\alpha=1$,~\dots,~$N$, be a collection
of homomorphisms of commutative rings for which the corresponding
collection of morphisms of affine schemes\/ $\Spec S_\alpha\rarrow
\Spec R$ is a finite open covering.
 Then \par
\textup{(a)} an $R$\+module $F$ is very flat if and only if all
the $S_\alpha$\+modules $S_\alpha\ot_R F$ are very flat; \par
\textup{(b)} for any contraadjusted $R$\+module $P$, the \v Cech
sequence
\begin{multline} \label{cech-contra}
 0\lrarrow\Hom_R(S_1\ot_R\dotsb\ot_R S_N\;P)\lrarrow\dotsb \\
 \textstyle\lrarrow\bigoplus_{\alpha<\beta}
 \Hom_R(S_\alpha\ot_R S_\beta\;P)\lrarrow
 \bigoplus_\alpha \Hom_R(S_\alpha\;P)\lrarrow P\lrarrow 0
\end{multline}
is an exact sequence of $R$\+modules.
\end{lem}

\begin{proof}
 Part~(a): by Corollary~\ref{very-open-cor}(a\+b), all the $R$\+modules
$S_{\alpha_1}\ot_R\dotsb\ot_R S_{\alpha_k}\ot_R F$ are very flat
whenever the $S_\alpha$\+modules $S_\alpha\ot_R F$ are very flat.
 For any $R$\+module $F$, the \v Cech sequence
\begin{multline} \label{cech-modules} \textstyle
 0\lrarrow F\lrarrow \bigoplus_\alpha S_\alpha\ot_R F \lrarrow
 \bigoplus_{\alpha<\beta} S_\alpha\ot_R S_\beta\ot_R F \\
 \lrarrow\dotsb\lrarrow S_1\ot_R\dotsb\ot_R S_N\ot_R F \lrarrow0
\end{multline}
is an exact sequence of $R$\+modules (since its localization at any
prime ideal of $R$ is contractible).
 It remains to recall that the class of very flat $R$\+modules is
closed with respect to the passage to the kernels of surjections.

 Part~(b): the exact sequence of $R$\+modules
\begin{multline} \label{cech-rings} \textstyle
 0\lrarrow R\lrarrow \bigoplus_\alpha S_\alpha \lrarrow
 \bigoplus_{\alpha<\beta} S_\alpha\ot_R S_\beta \\
 \lrarrow\dotsb\lrarrow S_1\ot_R\dotsb\ot_R S_N \lrarrow0
\end{multline}
is composed from short exact sequences of very flat $R$\+modules, so
the functor $\Hom_R({-},P)$ into a contraadjusted $R$\+module $P$
preserves its exactness.
\end{proof}

\subsection{Cotorsion modules}  \label{cotorsion-modules}
 Let $R$ be an associative ring.
 A left $R$\+module $P$ is said to be \emph{cotorsion}~\cite{En,Xu,EJ}
if $\Ext^1_R(F,P)=0$ for any flat left $R$\+module $F$, or
equivalently, $\Ext^{>0}_R(F,P)=0$ for any flat left $R$\+module~$F$.
 Clearly, the class of cotorsion left $R$\+modules is closed under
extensions and the passage to the cokernels of embeddings, and
also under infinite products.

 The following theorem, previously known essentially as the ``flat
cover conjecture'', was proved by Eklof--Trlifaj~\cite{ET} and
Bican--Bashir--Enochs~\cite{BBE} (cf.\ our
Theorem~\ref{eklof-trlifaj-very}).
 The case of a Noetherian commutative ring~$R$ of finite Krull
dimension was previously treated by Xu~\cite[Theorem~4.3.5]{Xu}
(cf.\ Lemma~\ref{noetherian-flat-cotorsion-embedding} below).

\begin{thm}  \label{flat-cover-thm}
\textup{(a)} Any $R$\+module can be embedded into a cotorsion
$R$\+module in such a way that the quotient module is flat. \par
\textup{(b)} Any $R$\+module admits a surjective map onto it from
a flat $R$\+module such that the kernel is cotorsion. \qed
\end{thm}

 The following lemmas concerning cotorsion (and injective) modules
are similar to the results about contraadjusted modules presented in
Section~\ref{affine-geometry-subsect}.
 With a possible exception of the last lemma, all of these are very
well known.

\begin{lem}  \label{cotors-hom}
 Let $R$ be a commutative ring.  Then \par
\textup{(a)} for any flat $R$\+module $F$ and cotorsion $R$\+module $P$,
the $R$\+module\/ $\Hom_R(F,P)$ is cotorsion; \par
\textup{(b)} for any $R$\+module $M$ and any injective $R$\+module $J$,
the $R$\+module\/ $\Hom_R(M,J)$ is cotorsion; \par
\textup{(c)} for any flat $R$\+module $M$ and any injective $R$\+module
$J$, the $R$\+module\/ $\Hom_R(F,J)$ is injective.
\end{lem}

\begin{proof}
 One has $\Ext_R^1(G,\Hom_R(F,P))\simeq\Ext_R^1(F\ot_R G\;P)$ for any
$R$\+modules $F$, $G$, and $P$ such that $\Ext_R^1(F,P)=0=\Tor^R_1(F,G)$.
 All the three assertions follow from this simple observation (which
is a particular case of formula~\eqref{ext-tor-adjunction}).
\end{proof}

 Our next lemma is a generalization of Lemma~\ref{cotors-hom}
to the noncommutative case.

\begin{lem}  \label{cotors-hom-nc}
 Let $R$ and $S$ be associative rings.  Then \par
\textup{(a)} for any $R$\+flat $R$\+$S$\+bimodule $F$ and any cotorsion
left $R$\+module $P$, the left $S$\+module\/ $\Hom_R(F,P)$ is cotorsion;
\par
\textup{(b)} for any $R$\+$S$\+bimodule $M$ and injective left
$R$\+module $J$, the left $S$\+module\/ $\Hom_R(M,J)$ is cotorsion; \par
\textup{(c)} for any $S$\+flat $R$\+$S$\+bimodule $F$ and any
injective left $R$\+module $J$, the left $S$\+module\/ $\Hom_R(F,J)$
is injective.
\end{lem}

\begin{proof}
 One has $\Ext_S^1(G,\Hom_R(F,P))\simeq\Ext_R^1(F\ot_S G\;P)$ for any
$R$\+$S$\+bimodule $F$, left $S$\+module $G$, and left $R$\+module $P$
such that $\Ext_R^1(F,P)=0=\Tor^S_1(F,G)$.
 Besides, the tensor product $F\ot_S G$ is flat over $R$ if $F$ is
flat over $R$ and $G$ is flat over~$S$.
 This proves~(a); and (b\+c)~are even easier.
\end{proof}

\begin{lem} \label{cotors-restrict}
 Let $f\:R\rarrow S$ be a homomorphism of associative rings.
 Then \par
\textup{(a)} any cotorsion left $S$\+module is also a cotorsion
left $R$\+module in the $R$\+module structure obtained by
the restriction of scalars via~$f$; \par
\textup{(b)} the left $S$\+module\/ $\Hom_R(S,J)$ obtained by
coextension of scalars via~$f$ is injective for any injective
left $R$\+module~$J$.
\end{lem}

\begin{proof}
 Part~(a): one has $\Ext_R^1(F,Q)\simeq\Ext_S^1(S\ot_R F\;Q)$ for any
flat left $R$\+module $F$ and any left $S$\+module~$Q$.
 Part~(b) is left to reader.
\end{proof}

\begin{lem} \label{cotors-coexten}
 Let $f\:R\rarrow S$ be an associative ring homomorphism such that
$S$ is a flat left $R$\+module in the induced $R$\+module structure.
 Then \par
\textup{(a)} the left $S$\+module\/ $\Hom_R(S,P)$ obtained by coextension
of scalars via~$f$ is cotorsion for any cotorsion left $R$\+module~$P$;
\par
\textup{(b)} any injective right $S$\+module is also an injective
right $R$\+module in the $R$\+module structure obtained by
the restriction of scalars via~$f$.
\end{lem}

\begin{proof}
 Part~(a): one has $\Ext_S^1(F,\Hom_R(S,P))\simeq \Ext_R^1(F,P)$ for
any left $R$\+module $P$ such that $\Ext_R^1(S,P)=0$ and any left
$S$\+module~$F$.
 In addition, in the assumptions of Lemma any flat left $S$\+module $F$
is also a flat left $R$\+module.
\end{proof}

\begin{lem}  \label{cotors-inj-covering}
 Let $R\rarrow S_\alpha$ be a collection of commutative ring
homomorphisms such that the corresponding collection of morphisms
of affine schemes $\Spec S_\alpha\rarrow\Spec R$ is an open covering.
 Then \par
\textup{(a)} a contraadjusted $R$\+module $P$ is cotorsion if and
only if all the contraadjusted $S_\alpha$\+modules\/
$\Hom_R(S_\alpha,P)$ are cotorsion; \par
\textup{(b)} a contraadjusted $R$\+module $J$ is injective if and
only if all the contraadjusted $S_\alpha$\+modules\/
$\Hom_R(S_\alpha,J)$ are injective.
\end{lem}

\begin{proof}
 Part~(a): the assertion ``only if'' follows from
Lemma~\ref{cotors-coexten}(a).
 To prove ``if'', use the \v Cech exact sequence~\eqref{cech-contra}
from Lemma~\ref{very-open-covering}(b).
 By Lemmas~\ref{cotors-restrict}(a) and~\ref{cotors-coexten}(a), all
the terms of the sequence, except perhaps the rightmost one,
are cotorsion $R$\+modules, and since the class of cotorsion
$R$\+modules is closed under the cokernels of embeddings, it follows
that the rightmost term is cotorsion as well.

 Part~(b) is provable in the similar way using parts~(b) of
Lemmas~\ref{cotors-restrict}\+-\ref{cotors-coexten}.
\end{proof}

 Let $R$ be a Noetherian commutative ring and $\p\sub R$ be a prime
ideal.
 Let $R_\p$ denote the localization of $R$ at~$\p$, and let
$\widehat R_\p$ be the completion of the local ring $R_\p$.
 We will use the notion of a \emph{contramodule} over a topological
ring, and in particular over a complete Noetherian local ring,
defined in~\cite[Section~1 and Appendix~B]{Pweak}, and denote
the abelian category of contramodules over a topological ring
$T$ by $T\contra$.
 (See also~\cite[Section~2.2]{Prev}, and the paper~\cite{Pcta} for
an alternative point of view.)

 The restriction of scalars with respect to the natural ring
homomorphism $R\rarrow\widehat R_\p$ provides an exact forgetful
functor $\widehat R_\p\contra\rarrow R\modl$.
 It follows from~\cite[Theorem~B.1.1(1)]{Pweak} that this functor is
fully faithful.
 The following proposition is a particular case of the assertions
of~\cite[Propositions~B.10.1 and~B.9.1]{Pweak} (see
also~\cite[Theorems~9.3 and~10.5]{Pcta}).

\begin{prop}  \label{contramodules-cotorsion}
\textup{(a)} Any $\widehat R_\p$\+contramodule is a cotorsion\/
$R$\+module. \par
\textup{(b)} Any free/projective $\widehat R_\p$\+contramodule is
a flat cotorsion $R$\+module.
\end{prop}

\begin{proof}
 In addition to the cited results from~\cite{Pweak} or~\cite{Pcta},
take into account Lemma~\ref{cotors-restrict}(a) and the fact that
any flat $R_\p$\+module is also a flat $R$\+module.
\end{proof}

 The following theorem is a restatement of the main result
of Enochs' paper~\cite{En}.

\begin{thm}  \label{flat-cotorsion-classification}
 Let $R$ be a Noetherian commutative ring.
 Then an $R$\+module is flat and cotorsion if and only if it is
isomorphic to an infinite product\/ $\prod_\p F_\p$ of
free contramodules\/ $F_\p=\widehat R_\p[[X_\p]]=\varprojlim_n
\widehat R_\p/\p^n[X_\p]$ over the complete local
rings\/~$\widehat R_\p$.
 Here $X_\p$ are some sets, and the direct product is taken over
all prime ideals\/ $\p\sub R$.  \qed
\end{thm}

 The following explicit construction of an injective morphism
with flat cokernel from a flat $R$\+module $G$ to a flat cotorsion
$R$\+module $\FC_R(G)$ was given in the book~\cite{Xu}.
 For any prime ideal $\p\sub R$, consider the localization
$G_\p=R_\p\ot_R G$ of the $R$\+module $G$ at $\p$, and take
its $\p$\+adic completion $\widehat G_\p=\varprojlim_n G_\p/\p^n G_\p$.
 By \cite[Lemma~1.3.2 and Corollary~B.8.2(b)]{Pweak}, the $R$\+module
$\widehat G_\p$ is a free $\widehat R_\p$\+contramodule.
 Furthermore, there is a natural $R$\+module map $G\rarrow
\widehat G_\p$.
 Set $\FC_R(G)=\prod_\p \widehat G_\p$.

\begin{lem}  \label{noetherian-flat-cotorsion-embedding}
 The natural $R$\+module morphism $G\rarrow\FC_R(G)$ is injective,
and its cokernel\/ $\FC_R(G)/G$ is a flat $R$\+module.
\end{lem}

\begin{proof} \hbadness=1900
 The following argument can be found in~\cite[Proposition~4.2.2
and Lemma~3.1.6]{Xu}.
 Any $R$\+module morphism from $G$ to an $R_\p/\p^n$\+module, and
hence also to a projective limit of such modules, factorizes
uniquely through the morphism $G\rarrow\widehat G_\p$.
 Consequently, in view of Theorem~\ref{flat-cotorsion-classification}
any morphism from $G$ to a flat cotorsion $R$\+module factorizes
through the morphism $G\rarrow\FC_R(G)$.

 Now one could apply Theorem~\ref{flat-cover-thm}(a), but it is
more instructive to argue directly as follows.
 Let $E$ be an injective cogenerator of the abelian category of
$R$\+modules; e.~g., $E=\Hom_\boZ(R,\boQ/\boZ)$.
 The $R$\+module $\Hom_R(G,E)$ being injective by
Lemma~\ref{cotors-hom}(c), the $R$\+module $\Hom_R(\Hom_R(G,E),E)$
is flat and cotorsion by part~(b) of the same lemma and by
Lemma~\ref{coherent-tensor-hom-lemma}(b) below.

 For any $R$\+module $N$, the natural $R$\+module morphism
$N\rarrow\Hom_R(\Hom_R(N,E),\allowbreak E)$ is injective.
 In particular, the natural morphism $G\rarrow\Hom_R(\Hom_R(G,E),E)$
is injective, and moreover, for any finitely generated/presented
$R$\+module $M$ the induced map $G\ot_R M\rarrow\Hom_R(\Hom_R(G,E),E)
\ot_R M\simeq\Hom_R(\Hom_R(G\ot_R\nobreak M\;E)\;E)$ is injective, too.
 The map $G\rarrow\Hom_R(\Hom_R(G,E),E)$ factorizes through
the map $G\rarrow\FC_R(G)$, and it follows that the map
$G\ot_RM\rarrow\FC_R(G)\ot_RM$ is injective as well.
 Thus the $R$\+module $FC_R(G)/G$ is flat.  \hfuzz=2.5pt
\end{proof}

\subsection{Exact categories of contraadjusted and cotorsion modules}
\label{contraadjusted-exact-cat}
 Let $R$ be a commutative ring.
 As full subcategories of the abelian category of $R$\+modules closed
under extensions, the categories of contraadjusted and very flat
$R$\+modules have natural exact category structures.
 In the exact category of contraadjusted $R$\+modules every morphism
has a cokernel, which is, in addition, an admissible epimorphism.

 In the exact category of contraadjusted $R$\+modules the functors of
infinite product are everywhere defined and exact; they also agree
with the infinite products in the abelian category of $R$\+modules.
 In the exact category of very flat $R$\+modules, the functors of
infinite direct sum are everywhere defined and exact, and agree with
the infinite direct sums in the abelian category of $R$\+modules.

 It is clear from Corollary~\ref{very-rel-proj-inj}(b) that there are
enough projective objects in the exact category of contraadjusted
$R$\+modules; these are precisely the very flat contraadjusted
$R$\+modules.
 Similarly, by Corollary~\ref{very-rel-proj-inj}(a) in the exact
category of very flat $R$\+modules there are enough injective objects;
these are also precisely the very flat contraadjusted modules.

 Denote the exact category of contraadjusted $R$\+modules by
$R\modl^\cta$ and the exact category of very flat $R$\+modules by
$R\modl_\vfl$.
 The tensor product of two very flat $R$\+modules is an exact functor
of two arguments $R\modl_\vfl\times R\modl_\vfl\rarrow R\modl_\vfl$.
 The $\Hom_R$ from a very flat $R$\+module into a contraadjusted
$R$\+module is an exact functor of two arguments
$(R\modl_\vfl)^\op\times R\modl^\cta\rarrow R\modl^\cta$
(where $\sC^\op$ denotes the opposite category to a category~$\sC$).

 For any homomorphism of commutative rings $f\:R\rarrow S$,
the restriction of scalars with respect to~$f$ is an exact functor
$S\modl^\cta\rarrow R\modl^\cta$.
 The extension of scalars $F\mpsto S\ot_R F$ is an exact functor
$R\modl_\vfl\rarrow S\modl_\vfl$.

 For any very flat homomorphism of commutative rings $f\:R\rarrow S$
(see Lemma~\ref{very-scalars-veryflat-case}), the restriction of scalars
with respect to~$f$  is an exact functor
$S\modl_\vfl\rarrow R\modl_\vfl$.
 The coextension of scalars $P\mpsto\Hom_R(S,P)$ is an exact functor
$R\modl^\cta\lrarrow S\modl^\cta$.
 In particular, these assertions hold for any homomorphism of
commutative rings $R\rarrow S$ such that the related morphism of affine 
schemes $\Spec S\rarrow\Spec R$ is an open embedding
(see Lemma~\ref{very-open-embedding}).

\begin{lem}  \label{very-exact-local}
 Let $R\rarrow S_\alpha$ be a collection of homomorphisms of commutative
rings for which the corresponding collection of morphisms of affine
schemes\/ $\Spec S_\alpha\rarrow\Spec R$ is an open covering.
 Then \par
\textup{(a)} a pair of homomorphisms of contraadjusted $R$\+modules
$K\rarrow L\rarrow M$ is a short exact sequence if and only if
such are the induced sequences of contraadjusted $S_\alpha$\+modules\/
$\Hom_R(S_\alpha,K)\rarrow\Hom_R(S_\alpha,L)\rarrow\Hom_R(S_\alpha,M)$
for all\/~$\alpha$; \par
\textup{(b)} a homomorphism of contraadjusted $R$\+modules $P\rarrow Q$
is an admissible epimorphism in $R\modl^\cta$ if and only if
the induced homomorphisms of contraadjusted $S_\alpha$\+modules\/
$\Hom_R(S_\alpha,P)\rarrow\Hom_R(S_\alpha,Q)$ are admissible
epimorphisms in $S_\alpha\modl^\cta$ for all\/~$\alpha$.
\end{lem}

\begin{proof}
 Part~(a): the ``only if'' assertion follows from
Lemma~\ref{very-open-embedding}.
 For the same reason, if the sequences $0\rarrow\Hom_R(S_\alpha,K)
\rarrow\Hom_R(S_\alpha,L)\rarrow\Hom_R(S_\alpha,M)\rarrow0$ are exact,
then so are the sequences obtained by applying the functors
$\Hom_R(S_{\alpha_1}\ot_R\dotsb\ot_R S_{\alpha_k}\;{-})$, \ $k\ge1$,
to the sequence $K\rarrow L\rarrow M$.
 Now it remains to make use of Lemma~\ref{very-open-covering}(b)
in order to deduce exactness of the original sequence
$0\rarrow K\rarrow L\rarrow M\rarrow 0$.

 Part~(b): it is clear from the very right segment of the exact
sequence~\eqref{cech-contra} that surjectivity of the maps
$\Hom_R(S_\alpha,P)\rarrow\Hom_R(S_\alpha,Q)$ implies surjectivity
of the map $P\rarrow Q$.
 It remains to check that the kernel of the latter morphism
is a contraadjusted $R$\+module.
 Denote this kernel by~$K$.
 Since the morphisms $\Hom_R(S_\alpha,P)\rarrow\Hom_R(S_\alpha,Q)$
are admissible epimorphisms, so are all the morphisms obtained
by applying the coextension of scalars with respect to
the ring homomorphisms $R\rarrow S_{\alpha_1}\ot_R\dotsb\ot_R
S_{\alpha_k}$, \ $k\ge1$, to the morphism $P\rarrow Q$.

 Now Lemma~\ref{very-open-covering}(b) applied to both sides of
the morphism $P\rarrow Q$ provides a termwise surjective morphism
of finite exact sequences of $R$\+modules.
 The corresponding exact sequence of kernels has $K$ as its rightmost
nontrivial term, while all the other terms are contraadjusted modules
over the rings $S_{\alpha_1}\ot_R\dotsb\ot_RS_{\alpha_k}$, \ $k\ge1$.
 By Lemma~\ref{very-scalars-always}(a), these are also contraadjusted
as $R$\+modules.
 It follows that the $R$\+module $K$ is contraadjusted as well.
\end{proof}

 Let $R$ be an associative ring.
 As a full subcategory of the abelian category of $R$\+modules closed
under extensions, the category of cotorsion left $R$\+modules has
a natural exact category structure.

 The functors of infinite product are everywhere defined and exact in
this exact category, and agree with the infinite products in the abelian
category of $R$\+modules.
 Similarly, the category of flat $R$\+modules has a natural exact 
category structure with exact functors of infinite direct sum.

 It follows from Theorem~\ref{flat-cover-thm} that there are enough
projective objects in the exact category of cotorsion $R$\+modules;
these are precisely the flat cotorsion $R$\+modules.
 Similarly, there are enough injective objects in the exact category
of flat $R$\+modules, and these are also precisely the flat
cotorsion $R$\+modules.

 Denote the exact category of cotorsion left $R$\+modules by
$R\modl^\cot$ and the exact category of flat left $R$\+modules by
$R\modl_\fl$.
 The abelian category of left $R$\+modules will be denoted simply
by $R\modl$, and the additive category of injective $R$\+modules
(with the trivial exact category structure) by $R\modl^\inj$.

 For any commutative ring $R$, the $\Hom_R$ from a flat $R$\+module
into a cotorsion $R$\+module is an exact functor of two arguments
$(R\modl_\fl)^\op\times R\modl^\cot\rarrow R\modl^\cot$.
 Analogously, the $\Hom_R$ from an arbitrary $R$\+module into
an injective $R$\+module is an exact functor
$(R\modl)^\op\times R\modl^\inj\rarrow R\modl^\cot$.
 The functors $\Hom$ over a noncommutative ring~$R$ mentioned in
Lemma~\ref{cotors-hom-nc} have similar exactness properties.

 For any associative ring homomorphism $f\:R\rarrow S$, the restriction
of scalars via~$f$ is an exact functor $S\modl^\cot\rarrow R\modl^\cot$.
 For any associative ring homomorphism $f\:R\rarrow S$ making $S$
a flat left $R$\+module, the coextension of scalars $P\mpsto
\Hom_R(S,P)$ is an exact functor $R\modl^\cot\rarrow S\modl^\cot$.

\begin{lem}  \label{cotors-exact-local}
 Let $R\rarrow S_\alpha$ be a collection of homomorphisms of commutative
rings for which the corresponding collection of morphisms of affine
schemes\/ $\Spec S_\alpha\rarrow\Spec R$ is an open covering.
 Then \par
\textup{(a)} a pair of morphisms of cotorsion $R$\+modules
$K\rarrow L\rarrow M$ is a short exact sequence if and only if
such are the sequences of cotorsion $S_\alpha$\+modules\/
$\Hom(S_\alpha,K)\rarrow\Hom(S_\alpha,L)\rarrow\Hom(S_\alpha,M)$
for all\/~$\alpha$; \par
\textup{(b)} a morphism of cotorsion $R$\+modules $P\rarrow Q$ is
an admissible epimorphism if and only if such are the morphisms
of cotorsion $S_\alpha$\+modules\/ $\Hom_R(S_\alpha,P)\rarrow
\Hom_R(S_\alpha,Q)$ for all\/~$\alpha$.
\end{lem}

\begin{proof}
 Part~(a) follows from Lemma~\ref{very-exact-local}(a); part~(b)
can be proved in the way similar to Lemma~\ref{very-exact-local}(b).
\end{proof}

\begin{lem}  \label{injective-exact-local}
 Let $R\rarrow S_\alpha$ be a collection of homomorphisms of commutative
rings for which the corresponding collection of morphisms of affine
schemes\/ $\Spec S_\alpha\rarrow\Spec R$ is an open covering.
 Then \par
\textup{(a)} a pair of morphisms of injective $R$\+modules
$I\rarrow J\rarrow K$ is a (split) short exact sequence if and only if
such are the sequences of injective $S_\alpha$\+modules\/
$\Hom(S_\alpha,I)\rarrow\Hom(S_\alpha,J)\rarrow\Hom(S_\alpha,K)$
for all\/~$\alpha$; \par
\textup{(b)} a morphism of injective $R$\+modules $J\rarrow K$ is
a split epimorphism if and only if such are the morphisms
of injective $S_\alpha$\+modules\/ $\Hom_R(S_\alpha,J)\rarrow
\Hom_R(S_\alpha,K)$ for all\/~$\alpha$.
\end{lem}

\begin{proof}
 Similar to the proof of Lemma~\ref{cotors-exact-local}.
\end{proof}

 The definition of an \emph{acyclic complex} in an exact category $\sE$
can be found in Section~\ref{derived-second-kind}.
 The assertions of the following lemma are mentioned
in~\cite[Examples~7.8, 7.9, and~7.14]{Pal}.

\begin{lem} \label{cta-cot-acyclicity-colocal}
 Let $R\rarrow S_\alpha$ be a collection of homomorphisms of commutative
rings for which the corresponding collection of morphisms of affine
schemes\/ $\Spec S_\alpha\rarrow\Spec R$ is an open covering.
 Then \par
\textup{(a)} a complex $C^\bu$ in the category $R\modl^\cta$ is
acyclic in $R\modl^\cta$ if and only if the complex\/
$\Hom_R(S_\alpha,C^\bu)$ is acyclic in $S_\alpha\modl^\cta$ for
every~$\alpha$; \par
\textup{(b)} a complex $C^\bu$ in the category $R\modl^\cot$ is
acyclic in $R\modl^\cot$ if and only if the complex\/
$\Hom_R(S_\alpha,C^\bu)$ is acyclic in $S_\alpha\modl^\cot$ for
every~$\alpha$.
\end{lem}

\begin{proof}
 Similar to the proof of Lemma~\ref{very-exact-local}(a).
\end{proof}

\subsection{Very flat and cotorsion dimensions}
\label{veryflat-cotors-dim-subsect}
 Let $R$ be a commutative ring.
 By analogy with the definition of the flat dimension of a module,
we define the \emph{very flat dimension} of an $R$\+module $M$ as
the minimal length of its resolution by very flat modules.

 Clearly, the very flat dimension of an $R$\+module $M$ is equal
to the supremum of the set of all integers~$d$ for which there exists
a contraadjusted $R$\+module $P$ such that $\Ext_R^d(M,P)\ne0$.
 The very flat dimension of a module cannot differ from its projective
dimension by more than~$1$.

 Similarly, the \emph{cotorsion dimension} of a left module $M$ over
an associative ring $R$ is conventionally defined as the minimal
length of its coresolution by cotorsion $R$\+modules.
 The cotorsion dimension of a left $R$\+module $M$ is equal to
the supremum of the set of all integers~$d$ for which there exists
a flat left $R$\+module $F$ such that $\Ext_R^d(F,M)\ne0$.

 Both the very flat and the cotorsion dimensions of a module do not
depend on the choice of a particular very flat/cotorsion (co)resolution
in the same sense as the familiar projective, flat, and injective
dimensions do not (see Corollary~\ref{fdim-cor}, \cite[Lemma~2.1]{Zhu},
or~\cite[Proposition~2.3(1)]{Sto0} for the general assertion
of this kind).

\begin{lem}  \label{cotorsion-dim-lemma}
 Let $R\rarrow S$ be a morphism of associative rings.
 Then any left $S$\+module $Q$ of cotorsion dimension~$\le d$ over $S$
has cotorsion dimension~$\le d$ over~$R$.
\end{lem}

\begin{proof}
 Follows from Lemma~\ref{cotors-restrict}(a).
\end{proof}

\begin{lem} \label{affine-d+D-lemma}
 Let $R\rarrow S$ be a morphism of associative rings such that
$S$ is a left $R$\+module of flat dimension\/~$\le D$.
 Then \par
\textup{(a)} any left $S$\+module $G$ of flat dimension\/~$\le d$
over $S$ has flat dimension\/~$\le d+D$ over\/~$R$; \par
\textup{(b)} any right $S$\+module $Q$ of injective dimension\/~$\le d$
over $S$ has injective dimension\/~$\le\nobreak d+D$ over~$R$.
\end{lem}

\begin{proof}
 Part~(a) follows from the spectral sequence
$E^2_{pq}=\Tor^S_p(\Tor^R_q(M,S),G)\Longrightarrow\Tor^R_{p+q}(M,G)$,
which holds for any right $R$\+module~$M$ \,\cite[Case~1 in
Section~XVI.5]{CaE}.
 Part~(b) follows from the spectral sequence
$E_2^{pq}=\Ext_{S^\rop}^p(\Tor^R_q(N,S),Q)\Longrightarrow
\Ext_{R^\rop}^{p+q}(N,Q)$, which holds for any right $R$\+module $N$
(where $S^\rop$ and $R^\rop$ denote the rings opposite to $S$ and~$R$)
\cite[Case~3 in Section~XVI.5]{CaE}.
 Alternatively, one can refer to Corollary~\ref{fdim-resolution}(a)
and its dual version.
\end{proof}

\begin{lem}  \label{affine-very-d+D-lemma}
\textup{(a)} Let $R\rarrow S$ be a morphism of commutative rings
such that $S$ is an $R$\+module of very flat dimension\/~$\le D$.
 Then any $S$\+module $G$ of very flat dimension\/~$\le d$ over $S$
has very flat dimension $\le d+1+D$ over~$R$. \par
\textup{(b)}  Let $R\rarrow S$ be a morphism of commutative rings
such that $S[s^{-1}]$ is an $R$\+module of very flat
dimension\/~$\le D$ for any element $s\in S$.
 Then any $S$\+module $G$ of very flat dimension\/~$\le d$
over $S$ has very flat dimension $\le d+D$ over~$R$.
\end{lem}

\begin{proof}
 Part~(a): the $S$\+module $G$ has a projective resolution of
length~$\le d+1$, and any projective $S$\+module has very flat
dimension~$\le D$ over $R$, which implies the desired assertion
(see Corollary~\ref{fdim-resolution}(a)).

 Part~(b): by Corollary~\ref{very-flat-transfinite}, the $S$\+module
$G$ has a resolution of length~$\le d$ by direct summands of
transfinitely iterated extensions of the $S$\+modules $S[s^{-1}]$.
 Hence it suffices to show that the very flat dimension of $R$\+modules
is not raised by the transfinitely iterated extension.

 More generally, we claim that one has $\Ext^n_R(M,P)=0$ whenever
a module $M$ over an associative ring $R$ is a transfinitely iterated
extension of $R$\+modules $M_\alpha$ and $\Ext^n_R(M_\alpha,P)=0$
for all~$\alpha$.
 The case $n=0$ is easy; the case $n=1$ is the result
of~\cite[Lemma~1]{ET}; and the case $n>1$ is reduced to $n=1$ by
replacing the $R$\+module $P$ with an $R$\+module $Q$ occuring at
the rightmost end of a coresolution $0\rarrow P\rarrow J^0\rarrow\dotsb
\rarrow J^{n-2}\rarrow Q\rarrow0$ with injective $R$\+modules~$J^i$.
\end{proof}

\begin{lem}  \label{flat-veryflat-cotors-inj-dim-local}
 Let $R\rarrow S_\alpha$ be a collection of homomorphisms of commutative
rings for which the corresponding collection of morphisms of affine
schemes is a finite open covering.  Then \par
\textup{(a)} the flat dimension of an $R$\+module $F$ is equal to
the supremum of the flat dimensions of the $S_\alpha$\+modules
$S_\alpha\ot_R F$; \par
\textup{(b)} the very flat dimension of an $R$\+module $F$ is equal to
the supremum of the very flat dimensions of the $S_\alpha$\+modules
$S_\alpha\ot_R F$; \par
\textup{(c)} the cotorsion dimension of a contraadjusted $R$\+module $P$
is equal to the supremum of the cotorsion dimensions of
the contraadjusted $S_\alpha$\+modules\/ $\Hom_R(S_\alpha,P)$; \par
\textup{(d)} the injective dimension of a contraadjusted $R$\+module $P$
is equal to the supremum of the injective dimensions of
the contraadjusted $S_\alpha$\+modules\/ $\Hom_R(S_\alpha,P)$.
\end{lem}

\begin{proof}
 Part~(b) follows easily from Lemma~\ref{very-open-covering}(a),
and the proof of part~(a) is similar.
 Parts~(c\+d) analogously follow from
Lemma~\ref{cotors-inj-covering}(a\+b).
\end{proof}

 The following lemma will be needed in
Sections~\ref{finite-dim-morphisms-I}
and~\ref{finite-dim-morphisms-derived-inverse-subsect}.

\begin{lem}  \label{adjusted-scalars-contra}
\textup{(a)} Let $f\:R\rarrow S$ be a homomorphism of commutative rings
and $P$ be an $R$\+module such that\/ $\Ext^1_R(S[s^{-1}],P)=0$ for all
elements $s\in S$.
 Then the $S$\+module\/ $\Hom_R(S,P)$ is contraadjusted. \par
\textup{(b)} Let $f\:R\rarrow S$ be a homomorphism of associative rings
and $P$ be a left $R$\+module such that $\Ext^1_R(G,P)=0$ for all
flat left $S$\+modules~$G$.
 Then the $S$\+module\/ $\Hom_R(S,P)$ is cotorsion.
\end{lem}

\begin{proof}
 See the proofs of Lemmas~\ref{very-scalars-veryflat-case}(a)
and~\ref{cotors-coexten}(a).
\end{proof}

 The following theorem is due to Raynaud and
Gruson~\cite[Corollaire~II.3.2.7]{RG}.

\begin{thm}  \label{raynaud-gruson-flat-thm}
 Let $R$ be a commutative Noetherian ring of Krull dimension~$D$.
 Then the projective dimension of any flat $R$\+module does not
exceed~$D$.
 Consequently, the very flat dimension of any flat $R$\+module
also does not exceed~$D$. \qed
\end{thm}

\begin{cor}  \label{raynaud-gruson-cotors-cor}
 Let $R$ be a commutative Noetherian ring of Krull dimension~$D$.
 Then the cotorsion dimension of any $R$\+module does not exceed~$D$.
\end{cor}

\begin{proof}
 For any associative ring $R$, the supremum of the projective dimensions
of flat left $R$\+modules and the supremum of the cotorsion dimensions
of arbitrary left $R$\+modules are equal to each other.
 Indeed, both numbers are equal to the supremum of the set of all
integers $d$ for which there exist a flat left $R$\+module $F$ and
a left $R$\+module $P$ such that $\Ext_R^d(F,P)\ne0$.
\end{proof}

\subsection{Periodicity theorems}
 The term \emph{periodicity theorems} refers to a class of results
concerning modules of cocycles in acyclic complexes.
 We refer to the introduction to the paper~\cite{BHP} for a survey and
to~\cite[Section~7]{Pphil} for a discussion.
 The most important ones for our purposes are the \emph{flat/projective}
and \emph{cotorsion} periodicity theorems.

 The following \emph{flat/projective periodicity theorem} is due to
Benson--Goodearl~\cite{BG} and Neeman~\cite{N-f}; see
also~\cite[Proposition~7.6]{CH}.

\begin{thm} \label{flat-projective-periodicity}
 Let $R$ be an associative ring. \par
\textup{(a)} If all the modules of cocycles in an acyclic complex of
projective $R$\+modules $P^\bu$ are flat, then the modules of cocycles
are actually projective and the complex $P^\bu$ is contractible. \par
\textup{(b)} If $P^\bu$ is a complex of projective left $R$\+modules
and $F^\bu$ is an acyclic complex of flat left $R$\+modules with flat
modules of cocycles, then any morphism of complexes $P^\bu\rarrow F^\bu$
is homotopic to zero.
\end{thm}

\begin{proof}
 To deduce part~(a) from part~(b), put $F^\bu=P^\bu$.
 Part~(a) is~\cite[Remark~2.15]{N-f}, or alternatively,
\cite[Theorem~2.5]{BG} with~\cite[Proposition~7.6]{CH}.
 Part~(b) is~\cite[Theorem~8.6\,(iii)\,$\Rightarrow$\,(i)]{N-f}. 
\end{proof}

 The following proposition complements
Theorem~\ref{flat-projective-periodicity}(b).

\begin{prop} \label{flat-projective-periodicity-complements}
 Let $R$ be an associative ring. \par
\textup{(a)} Let $F^\bu$ be a complex of flat left $R$\+modules
such that, for every complex of projective left $R$\+modules $P^\bu$,
any morphism of complexes $P^\bu\rarrow F^\bu$ is homotopic to zero.
 Then $F^\bu$ is a acyclic complex of flat $R$\+modules with
flat $R$\+modules of cocycles. \par
\textup{(b)} For every complex of flat left $R$\+modules $F^\bu$
there exists a complex of projective left $R$\+modules $P^\bu$
together with a morphism of complexes $P^\bu\rarrow F^\bu$ whose
cone is an acyclic complex with flat modules of cocycles.
\end{prop}

\begin{proof}
 Part~(a) is~\cite[Lemma~8.5 or
Theorem~8.6\,(i)\,$\Rightarrow$\,(iii)]{N-f}.
 Part~(b) is part~(a) together with~\cite[Proposition~8.1]{N-f}.
 For a contramodule version of Theorem~\ref{flat-projective-periodicity}
and Proposition~\ref{flat-projective-periodicity-complements}, see
the preprint~\cite[Theorems~5.1 and~6.1]{Pbc}.
\end{proof}

 The next \emph{cotorsion periodicity theorem} is due to
Bazzoni, Cort\'es-Izurdiaga, and Estrada~\cite{BCE}.

\begin{thm} \label{cotorsion-periodicity}
 Let $R$ be an associative ring.
 If $C^\bu$ is a acyclic complex of cotorsion $R$\+modules, then all
the $R$\+modules of cocycles in the complex $C^\bu$ are also cotorsion.
\end{thm}

\begin{proof}
 This is~\cite[Theorem~5.1(2)]{BCE}.
 An additional discussion in category-theoretic and scheme-theoretic
contexts can be found in~\cite[Sections~9\+-10]{PS6}.
\end{proof}

 The following corollary can be summarized by saying that
\emph{any cotorsion resolved module is cotorsion}.
 It is the cotorsion version of the much simpler assertion that
any quotient module of a contraadjusted module is contraadjusted.

\begin{cor} \label{cotorsion-resolved}
 Let $R$ be an associative ring and $M$ be an $R$\+module.
 Assume that $M$ admits a resolution $\dotsb\rarrow C_2\rarrow C_1
\rarrow C_0\rarrow M\rarrow0$ by cotorsion $R$\+modules~$C_n$,
\,$n\ge0$.
 Then $M$ is a cotorsion $R$\+module, too.
\end{cor}

\begin{proof}
 Obviously, any $R$\+module $M$ admits a cotorsion (or injective)
coresolution, so the assertion follows from
Theorem~\ref{cotorsion-periodicity}.
\end{proof}

 The following important proposition is essentially an equivalent
version of Theorem~\ref{cotorsion-periodicity}.

\begin{prop} \label{dw-cotorsion-are-dg-cotorsion}
 Let $R$ be an associative ring.
 If $F^\bu$ is an acyclic complex of flat left $R$\+modules with
flat modules of cocycles and $C^\bu$ is a complex of cotorsion
left $R$\+modules, then any morphism of complexes $F^\bu\rarrow
C^\bu$ is homotopic to zero.
\end{prop}

\begin{proof}
 This is~\cite[Theorem~5.3]{BCE}.
\end{proof}

\subsection{Coherent rings, finite morphisms, and coadjusted modules}
\label{coherent-coadjusted}
 Recall that an associative ring $R$ is called \emph{left coherent}
if all its finitely generated left ideals are finitely presented.
 Finitely presented left modules over a left coherent ring $R$
form an abelian subcategory in $R\modl$ closed under kernels,
cokernels, and extensions.
 The definition of a right coherent associative ring is similar.

 The following lemma is a particular case of~\cite[Lemma~4.1]{Pfp}
and~\cite[Lemma~7.3]{Pps}.

\begin{lem}  \label{coherent-tensor-hom-lemma}
 Let $R$ and $S$ be associative rings.  Then \par
\textup{(a)} Assuming that the ring $R$ is left Noetherian, for any
$R$\+injective $R$\+$S$\+bimodule $J$ and any flat left $S$\+module $F$
the left $R$\+module $J\ot_S F$ is injective. \par
\textup{(b)} Assuming that the ring $S$ is right coherent, for any
$S$\+injective $R$\+$S$\+bimodule $I$ and any injective left
$R$\+module $J$ the left $S$\+module $\Hom_R(I,J)$ is flat.
\end{lem}

\begin{proof}
 Part~(a) holds due to the natural isomorphism $\Hom_R(M\;J\ot_SF)
\simeq\Hom_R(M,J)\ot_S F$ for any finitely presented left module $M$
over an associative ring $R$, any $R$\+$S$\+bimodule $J$, and
any flat left $S$\+module~$F$.
 Part~(b) follows from the natural isomorphism $N\ot_S\Hom_R(I,J)
\simeq\Hom_R(\Hom_{S^\rop}(N,I),J)$ for any finitely presented
right module $N$ over an associative ring $S$, any $R$\+$S$\+bimodule
$I$, and any injective left $R$\+module~$J$ (where $S^\rop$ denotes
the ring opposite to~$S$).

 Here we use the facts that injectivity of a left module $I$ over
a left Noetherian ring $R$ is equivalent to exactness of the functor
$\Hom_R({-},I)$ on the category of finitely generated left
$R$\+modules, while flatness of a left module $F$ over a right
coherent ring $S$ is equivalent to exactness of the functor
${-}\ot_S F$ on the category of finitely presented right
$S$\+modules (cf.\ the proof of the next Lemma~\ref{coherent-lemma}).
\end{proof}

\begin{lem}  \label{coherent-lemma}
 Let $R$ and $S$ be associative rings such that $R$ is left coherent.
 Let $F$ be a left $S$\+module of finite projective dimension,
$P$ be an $R$\+flat $S$\+$R$\+bimodule such that $\Ext_S^{>0}(F,P)=0$,
and $M$ be a finitely presented left $R$\+module.
 Then one has\/ $\Ext_S^{>0}(F\;P\ot_R M)=0$, the natural map
of abelian groups\/ $\Hom_S(F,P)\ot_R M\rarrow\Hom_S(F\;P\ot_R M)$
is an isomorphism, and the right $R$\+module\/ $\Hom_S(F,P)$ is flat.
\end{lem}

\begin{proof}
 Let $L_\bu\rarrow M$ be a resolution of $M$ by finitely generated
projective $R$\+modules.
 Then $P\ot_R L_\bu\rarrow P\ot_R M$ is a resolution of the $S$\+module
$P\ot_RM$ by $S$\+modules annihilated by $\Ext_S^{>0}(F,{-})$.
 Since the $S$\+module $F$ has finite projective dimension, it follows
that $\Ext_S^{>0}(F\;P\ot_R M)=0$.

 Consequently, the functor $M\mpsto\Hom_S(F\;P\ot_R M)$ is exact
on the abelian category of finitely presented left $R$\+modules~$M$.
 Obviously, the functor $M\mpsto\Hom_S(F,P)\ot_S M$ is right exact.
 Since the morphism of functors $\Hom_S(F,P)\ot_R M\rarrow
\Hom_S(F\;P\ot_R M)$ is an isomorphism for finitely generated
projective $R$\+modules $M$, we can conclude that it is an isomorphism
for all finitely presented left $R$\+modules.

 Now we have proved that the functor $M\mpsto\Hom_S(F,P)\ot_R M$
is exact on the abelian category of finitely presented left
$R$\+modules.
 Since any left $R$\+module is a filtered inductive limit of finitely
presented ones and the inductive limits commute with tensor products,
it follows that the $R$\+module $\Hom_S(F,P)$ is flat.
\end{proof}

 The following lemma is a version of the preceding one using
the cotorsion periodicity theorem instead of the finite projective
dimension assumption.

\begin{lem}  \label{coherent-cotorsper-lemma}
 Let $R$ and $S$ be associative rings such that $R$ is left coherent.
 Let $P$ be an $R$\+flat $S$\+cotorsion $S$\+$R$\+bimodule and $M$ be
a finitely presented left $R$\+module.
 Then the left $S$\+module $P\ot_RM$ is cotorsion.
 For any flat left $S$\+module $F$, the natural map of abelian groups\/
$\Hom_S(F,P)\ot_R M\rarrow\Hom_S(F\;P\ot_R M)$ is an isomorphism, and
the right $R$\+module\/ $\Hom_S(F,P)$ is flat.
\end{lem}

\begin{proof}
 Let $L_\bu\rarrow M$ be a resolution of $M$ by finitely generated
projective $R$\+modules.
 Then $P\ot_R L_\bu\rarrow P\ot_R M$ is a resolution of the $S$\+module
$P\ot_RM$ by cotorsion $S$\+modules.
 Hence $P\ot_RM$ is a cotorsion $S$\+module by
Corollary~\ref{cotorsion-resolved}.
 Thus $\Ext_S^{>0}(F\;P\ot_R M)=0$, and the argument proceeds from this
point on similarly to the proof of Lemma~\ref{coherent-lemma}.
\end{proof}

\begin{cor}  \label{coherent-very}
 Let $R$ be a commutative ring.  Then \par
\textup{(a)} for any finitely generated $R$\+module $M$ and any
contraadjusted $R$\+module $P$, the $R$\+module $M\ot_R P$
is contraadjusted; \par
\textup{(b)} if the ring $R$ is coherent, then for any very flat
$R$\+module $F$ and any flat contraadjusted $R$\+module $P$,
the $R$\+module\/ $\Hom_R(F,P)$ is flat and contraadjusted; \par
\textup{(c)} in the situation of~\textup{(b)}, for any
finitely presented $R$\+module $M$ the natural morphism of
$R$\+modules\/ $\Hom_R(F,P)\ot_R M\rarrow\Hom_R(F\;P\ot_R M)$
is an isomorphism.
\end{cor}

\begin{proof}
 Part~(a) immediately follows from the facts that the class of
contraadjusted $R$\+modules is closed under finite direct sums
and quotients.
 Part~(b) is provided Lemma~\ref{coherent-lemma} (for a commutative
ring $R=S$) together with Lemma~\ref{very-tensor-hom}(b), and part~(c)
is also Lemma~\ref{coherent-lemma}.
\end{proof}

\begin{cor}  \label{coherent-cotors}
 Let $R$ be a coherent commutative ring.  Then \par
\textup{(a)} for any finitely presented $R$\+module $M$ and any flat
cotorsion $R$\+module $P$, the $R$\+module $M\ot_R P$ is cotorsion; \par
\textup{(b)} for any flat $R$\+module $F$ and flat cotorsion $R$\+module
$P$, the $R$\+module $\Hom_R(F,P)$ is flat and cotorsion; \par
\textup{(c)} in the situation of~\textup{(a)} and~\textup{(b)},
the natural morphism of $R$\+modules\/ $\Hom_R(F,P)\allowbreak
\ot_R M\rarrow \Hom_R(F\;P\ot_R M)$ is an isomorphism.
\end{cor}

\begin{proof}
 This is Lemma~\ref{coherent-cotorsper-lemma} specialized to the case
of a commutative ring $R=S$.

 Alternatively, for a Noetherian commutative ring $R$, there is
an argument based on Proposition~\ref{contramodules-cotorsion} and
Theorem~\ref{flat-cotorsion-classification}.
 Part~(a): since the functor $M\ot_R{-}$ preserves infinite products,
it suffices to show that the $R$\+module $M\ot_R Q_\p$ has
an $\widehat R_\p$\+contramodule structure for any
$\widehat R_\p$\+contramodule~$Q_\p$.
 Indeed, the full subcategory $\widehat R_\p\contra\sub R\modl$ is
closed with respect to finite direct sums and cokernels.
 
 Part~(c): the morphism in question is clearly an isomorphism for
a finitely generated projective $R$\+module~$M$.
 Hence it suffices to show that the functor $M\mpsto\Hom_R(F\;
P\ot_R M)$ is right exact.
 In fact, it is exact, since the functor $\Hom_R(F,{-})$
preserves exactness of short sequences of
$\widehat R_\p$\+contramodules.
 Therefore, the functor $M\mpsto\Hom_R(F,P)\ot_R M$ is also exact,
and we have proved part~(b) as well.
\end{proof}

\begin{cor} \label{coherent-flat-local}
\textup{(a)} Let $R\rarrow S$ be a homomorphism of commutative rings
such that the related morphism of affine schemes\/
$\Spec S\rarrow\Spec R$ is an open embedding.
 Assume that the ring $R$ is coherent.
 Then the $S$\+module\/ $\Hom_R(S,P)$ is flat and contraadjusted for
any flat contraadjusted $R$\+module~$P$. \par
\textup{(b)} Let $R\rarrow S_\alpha$ be a collection of homomorphisms of
commutative rings for which the corresponding collection of morphisms of
affine schemes $\Spec S_\alpha\rarrow\Spec R$ is an open covering.
 Assume that either the ring $R$ is Noetherian and an $R$\+module $P$
is cotorsion, or the ring $R$ is Noetherian of finite Krull dimension
and an $R$\+module $P$ is contraadjusted.
 Then the $R$\+module $P$ is flat if and only if all
the $S_\alpha$\+modules $\Hom_R(S_\alpha,P)$ are flat.
\end{cor}

\begin{proof}
 Part~(a): the $S$\+module $\Hom_R(S,P)$ is contraadjusted by
Corollary~\ref{very-open-cor}(c).
 The $R$\+module $\Hom_R(S,P)$ is flat by
Lemma~\ref{very-open-embedding} and Corollary~\ref{coherent-very}(b).
 Since $\Spec S\rarrow\Spec R$ is an open embedding, it follows that
$\Hom_R(S,P)$ is also flat as an $S$\+module
(cf.\ Corollary~\ref{very-open-cor}(b)).

 The ``only if'' assertion in part~(b) is provided by part~(a).
 The proof of the ``if'' is postponed to
Chapter~\ref{dualizing-complex-sect}.
 The cotorsion case will follow from Corollary~\ref{lct-prj-flat},
while the finite Krull dimension case will be covered by
Corollary~\ref{finite-krull-flat-contraherent}(b).
\end{proof}

 In the terminology of Section~\ref{colocal-classes-subsect},
Corollary~\ref{coherent-flat-local}(b) tells us that the class of
flat contraadjusted modules is \emph{colocal} over Noetherian rings of
finite Krull dimension (while the class of flat cotorsion modules
is colocal over all Noetherian rings, cf.\
Lemma~\ref{cotors-inj-covering}(a)).

\medskip

 Let $R\modl_\fp$ denote the abelian category of finitely presented
left modules over a left coherent ring~$R$.
 For a coherent commutative ring $R$, the tensor product of
a finitely presented $R$\+module with a flat contraadjusted $R$\+module
is an exact functor $R\modl_\fp\times R\modl_\fl^\cta
\rarrow R\modl^\cta$.
 Here the notation is $R\modl_\fl^\cta=R\modl_\fl\cap R\modl^\cta$;
the exact structure on $R\modl_\fl^\cta$ is inherited from $R\modl$.
 The $\Hom_R$ from a very flat $R$\+module into a flat
contraadjusted $R$\+module is an exact functor
$(R\modl_\vfl)^\op\times R\modl_\fl^\cta\rarrow
R\modl_\fl^\cta$.

 Let $R$ be a coherent commutative ring.
 Then the tensor product of a finitely presented $R$\+module with
a flat cotorsion $R$\+module is an exact functor
$R\modl_\fp\times R\modl_\fl^\cot\rarrow R\modl^\cot$.
 The $\Hom_R$ from a flat $R$\+module to a flat cotorsion $R$\+module
is an exact functor $(R\modl_\fl)^\op\times R\modl_\fl^\cot
\rarrow R\modl_\fl^\cot$.
 Here the additive category of flat cotorsion $R$\+modules
$R\modl_\fl^\cot=R\modl_\fl\cap R\modl^\cot$ is endowed with
a trivial exact category structure.

\begin{lem}  \label{quotient-scalars-contraadjusted}
 Let $R$ be a commutative ring and $I\sub R$ be an ideal.
 Then \par
\textup{(a)} an $R/I$\+module $Q$ is a contraadjusted $R/I$\+module
if and only if it is a contraadjusted $R$\+module; \par
\textup{(b)} the $R/I$\+module $P/IP$ is contraadjusted for any
contraadjusted $R$\+module~$P$; \par
\textup{(c)} assuming that the ring $R$ is coherent and the ideal $I$ is
finitely generated, for any very flat $R$\+module $F$ and flat
contraadjusted $R$\+module $P$ the natural morphism of $R/I$\+modules\/
$\Hom_R(F,P)/I\Hom_R(F,P)\rarrow\Hom_{R/I}(F/IF\;P/IP)$ is
an isomorphism.
\end{lem}

\begin{proof}
 Part~(a): the characterization of contraadjusted modules given in
the beginning of Section~\ref{very-eklof-trlifaj-subsect} shows that
the contraadjustedness property of a module depends only on its abelian
group structure and the operators by which the ring acts in it
(rather than on the ring indexing such operators).
 Alternatively, the ``only if'' assertion is a particular case of
Lemma~\ref{very-scalars-always}(a), and one can deduce the ``if'' from
the observation that any element $\bar r\in R/I$ can be lifted to
an element $r\in R$ so that one has an isomorphism of $R/I$\+modules
$R/I[\bar r^{-1}] \simeq R/I\ot_R R[r^{-1}]$.

 Part~(b): the $R$\+module $P/IP$ is contraadjusted as a quotient module
of a contraadjusted $R$\+module.
 By part~(a), $P/IP$ is also a contraadjusted $R/I$\+module.
 Part~(c): by Corollary~\ref{coherent-very}(c), the natural morphism
of $R$\+modules $\Hom_R(F,P)/I\Hom_R(F,P)\allowbreak\rarrow
\Hom_R(F\;P/IP)$ is an isomorphism.
\end{proof}

 The following lemma establishes locality of the coherence property of
commutative rings.

\begin{lem} \label{coherence-of-rings-local}
\textup{(a)} Let $R\rarrow S$ be a homomorphism of commutative rings
such that the related morphism of affine schemes\/ $\Spec S\rarrow
\Spec R$ is an open embedding.
 Assume that the ring $R$ is coherent.
 Then the ring $S$ is coherent. \par
\textup{(b)} Let $f_\alpha\:R\rarrow S_\alpha$ be a collection of
homomorphisms of commutative rings for which the corresponding
collection of morphisms of affine schemes\/ $\Spec S_\alpha\rarrow
\Spec R$ is a finite open covering.
 Assume that all the rings $S_\alpha$ are coherent.
 Then the ring $R$ is coherent.
\end{lem}

\begin{proof}
 This is~\cite[Lemma~9.1]{PS5}.
 Part~(a) holds because, more generally, for any flat epimorphism of
commutative rings $R\rarrow S$ (in the sense
of~\cite[Sections~XI.1--XI.2]{Sten}), coherence of $R$ implies
coherence of~$S$ \,\cite[Proposition~3.7]{CEI}.
 Part~(b) holds because, more generally, for any faithfully flat
homomorphism of commutative rings $R\rarrow S$, coherence of $S$
implies coherence of~$R$ \,\cite[Corollary~2.1]{Harr}, \cite[Propositions~I.5.9 and~I.6.11]{Bour}.
\end{proof}

 Recall that a morphism of Noetherian rings $R\rarrow S$ is called
\emph{finite} if $S$ is a finitely generated $R$\+module in
the induced $R$\+module structure.

\begin{lem}  \label{quotient-scalars-cotorsion}
 Let $R\rarrow S$ be a finite morphism of Noetherian commutative rings.
 Then \par
\textup{(a)} the $S$\+module $S\ot_R P$ is flat and cotorsion for
any flat cotorsion $R$\+module~$P$; \par
\textup{(b)} assuming that the ring $S$ has finite Krull dimension,
an $S$\+module $Q$ is a cotorsion $S$\+module if and only if it is
a cotorsion $R$\+module; \par
\textup{(c)} for any flat $R$\+module $F$ and flat cotorsion
$R$\+module $P$, the natural morphism of $S$\+modules
$S\ot_R\Hom_R(F,P)\rarrow\Hom_S(S\ot_R F\;S\ot_R P)$ is an isomorphism.
\end{lem}

\begin{proof}
 The proof of parts~(a\+b) is based on
Theorem~\ref{flat-cotorsion-classification}.
 Given a prime ideal $\p\sub R$, consider all the prime ideals
$\q\sub S$ whose preimage in $R$ coincides with~$\p$.
 Such ideals form a nonempty finite set, and there are no inclusions
between them~\cite[Theorems~9.1 and~9.3, and Exercise~9.3]{Mats}.
 Let us denote these ideals by $\q_1$,~\dots,~$\q_m$.
 By~\cite[Theorems~9.4(i), 8.7, and~8.15]{Mats}, we have
$S\ot_R \widehat R_\p \simeq
\widehat S_{\q_1}\oplus\dotsb\oplus\widehat S_{\q_m}$.

 Since the functor $S\ot_R{-}$ preserves infinite products, in order
to prove~(a) it suffices to show that the $S$\+module $S\ot_R F_\p$
is a finite direct sum of certain free
$\widehat S_{\q_i}$\+contramodules~$F_{\q_i}$.
 This can be done either by noticing that $F_\p$ is a direct
summand of an infinite product of copies of $\widehat R_\p$
(see~\cite[Section~1.3]{Pweak}), or by showing that the natural
map $S\ot_R \widehat R_\p[[X]]\rarrow\bigoplus_{i=1}^m
\widehat S_{\q_i}[[X]]$ is an isomorphism for any set~$X$ (see
\cite[proof of the first assertion of Proposition~B.9.1]{Pweak}).
 In addition to the assertion of part~(a), we have also proved that
any flat cotorsion $S$\+module $Q$ is a direct summand of
an $S$\+module $S\ot_R P$ for a certain flat cotorsion $R$\+module~$P$.

 The ``only if'' assertion in part~(b) is a particular case of
Lemma~\ref{cotors-restrict}(a).
 Let us prove the ``if''.
 According to Corollary~\ref{raynaud-gruson-cotors-cor},
the cotorsion dimension of any $R$\+module is finite.
 By Theorem~\ref{flat-cover-thm}(a), it follows that any flat
$S$\+module admits a finite coresolution by flat cotorsion
$S$\+modules (cf.\ the dual version of
Corollary~\ref{fdim-subcategory-cor}).
 Hence it suffices to prove that $\Ext_S^{>0}(G,Q)=0$ for a flat
cotorsion $S$\+module~$G$.
 This allows us to assume that $G=S\ot_R F$, where $F$ is
a flat (cotorsion) $R$\+module.
 It remains to recall the $\Ext$ isomorphism from the proof of
Lemma~\ref{cotors-restrict}(a).

 Part~(c): by Corollary~\ref{coherent-cotors}(c), the natural
morphism of $R$\+modules $S\ot_R\Hom_R(F,P)\rarrow\Hom_R(F\;S\ot_RP)$
is an isomorphism.
\end{proof}

\begin{lem}  \label{fin-gen-nilpotent-quotient-scalars}
 Let $R$ be a commutative ring and $I\sub R$ be a finitely generated
nilpotent ideal.
 Then \par
\textup{(a)} an $R$\+module $P$ is contraadjusted if and only if
the $R/I$\+module $P/IP$ is contraadjusted; \par
\textup{(b)} a flat $R$\+module $F$ is very flat if and only if
the $R/I$\+module $F/IF$ is very flat; \par
\textup{(c)} assuming that the ring $R$ is Noetherian, a flat
$R$\+module $P$ is cotorsion if and only if the $R/I$\+module
$P/IP$ is cotorsion; \par
\textup{(d)} assuming that the ring $R$ is coherent, a flat $R$\+module
$P$ is cotorsion whenever the $R/I$\+module $P/IP$ is cotorsion.
\end{lem}

\begin{proof}
 Part~(a): the ``only if'' assertion is a particular case of
Lemma~\ref{quotient-scalars-contraadjusted}(b).
 To prove the ``if'', notice that in our assumptions about $I$
the $R$\+module $P$ has a finite decreasing filtration by
its submodules~$I^nP$.
 Furthermore, the successive quotients $I^nP/I^{n+1}P$ are
the targets of the natural surjective homomorphisms of
$R/I$\+modules $I^n/I^{n+1}\ot_{R/I}P/IP\rarrow I^nP/I^{n+1}P$.
 Since the $R/I$\+module $I^n/I^{n+1}$ is finitely generated,
the $R/I$\+module $I^nP/I^{n+1}P$ is a quotient module of 
a finite direct sum of copies of the $R/I$\+module $P/IP$.
 It remains to use the facts that the class of contraadjusted
modules over a given commutative ring is closed under extensions
and quotients, together with the result of
Lemma~\ref{quotient-scalars-contraadjusted}(a).

 Part~(b): the ``only if'' is a particular case of
Lemma~\ref{very-scalars-always}(b); let us prove the ``if''.
 Let $P$ be a contraadjusted $R$\+module; we have to show that
$\Ext^1_R(F,P)=0$.
 According to the above proof of part~(a), the $R$\+module $P$
has a finite filtration whose successive quotients are
contraadjusted $R/I$\+modules with the $R$\+module structures
obtained by restriction of scalars.
 So it suffices to check that $\Ext^1_R(F,Q)=0$ for any
contraadjusted $R/I$\+module~$Q$.
 Now, the $R$\+module $F$ being flat, one has $\Ext^1_R(F,Q) =
\Ext^1_{R/I}(F/IF\;Q)$ for any $R/I$\+module~$Q$
(see the proof of Lemma~\ref{cotors-restrict}(a)).

 Parts~(c\+d): the ``only if'' assertion in~(c) is a particular case
of Lemma~\ref{quotient-scalars-cotorsion}(a).
 To prove the ``if'', notice the isomorphisms of $R/I$\+modules
$I^nP/I^{n+1}P\simeq I^n/I^{n+1}\ot_RP\simeq I^n/I^{n+1}\ot_{R/I}P/IP$
(the former of which holds since the $R$\+module $P$ is flat).
 The $R/I$\+module $P/IP$ being flat and cotorsion,
the $R/I$\+modules $I^n/I^{n+1}\ot_{R/I}P/IP$ are cotorsion by
Corollary~\ref{coherent-cotors}(a), the $R$\+modules
$I^n/I^{n+1}\ot_{R/I}P/IP$ are cotorsion by
Lemma~\ref{cotors-restrict}(a), and the $R$\+module $P$
is cotorsion because the class of cotorsion $R$\+modules is
closed under extensions.
\end{proof}

 Let $R$ be a commutative ring and $I\sub R$ be an ideal.
 Then the reduction $P\mpsto P/IP$ of a flat contraadjusted
$R$\+module $P$ modulo $I$ is an exact functor $R\modl_\fl^\cta
\rarrow R/I\modl_\fl^\cta$.

 Let $f\:R\rarrow S$ be a finite morphism of Noetherian commutative
rings.
 Then the extension of scalars $P\mpsto S\ot_R P$ of a flat cotorsion
$R$\+module $P$ with respect to~$f$ is an additive functor
$R\modl_\fl^\cot\rarrow S\modl_\fl^\cot$ (between additive categories
naturally endowed with trivial exact category structures).

\medskip

 Let $R$ be a commutative ring.
 We will say that an $R$\+module $K$ is \emph{coadjusted} if
the functor of tensor product with $K$ over $R$ preserves the class
of contraadjusted $R$\+modules.
 By Corollary~\ref{coherent-very}(a), any finitely generated
$R$\+module is coadjusted.

 An $R$\+module $K$ is coadjusted if and only if the $R$\+module
$K\ot_R P$ is contraadjusted for every flat (or very flat)
contraadjusted $R$\+module~$P$.
 Indeed, by Corollary~\ref{very-rel-proj-inj}(b), any contraajusted
$R$\+module is a quotient module of a very flat contraadjusted
$R$\+module; so it remains to recall that any quotient module of
a contraadjusted $R$\+module is contraadjusted.

 Clearly, any quotient module of a coadjusted $R$\+module is coadjusted.
 Furthermore, the class of coadjusted $R$\+modules is closed under
extensions.
 One can see this either by applying the above criterion of
coadjustedness in terms of tensor products with flat contraadjusted
$R$\+modules, or straightforwardly from the right exactness property
of the functor of tensor product together with the facts that the class
of contraadjusted $R$\+modules is closed under quotients and extensions.

 Consequently, there is the induced exact category structure on
the full subcategory of coadjusted $R$\+modules in the abelian
category $R\modl$.
 We denote this exact category by $R\modl^\coa$.
 The tensor product of a coadjusted $R$\+module with a flat
contraadjusted $R$\+module is an exact functor $R\modl^\coa\times
R\modl_\fl^\cta\rarrow R\modl^\cta$.

 Over a Noetherian commutative ring $R$, any injective module $J$ is
coadjusted.
 Indeed, for any $R$\+module $P$, the tensor product $J\ot_R P$ is
a quotient module of an infinite direct sum of copies of $J$, which
means a quotient module of an injective module, which is contraadjusted.
 Hence any quotient module of an injective module is coadjusted, too,
as is any extension of such modules.

\begin{lem}  \label{coadjusted-local}
\textup{(a)} Let $f\:R\rarrow S$ be a homomorphism of commutative rings
such that the related morphism of affine schemes\/
$\Spec S\rarrow\Spec R$ is an open embedding.
 Then the $S$\+module $S\ot_R K$ obtained by the extension of scalars
via~$f$ is coadjusted for any coadjusted $R$\+module~$K$. \par
\textup{(b)} Let $f_\alpha\:R\rarrow S_\alpha$ be a collection of
homomorphisms of commutative rings for which the corresponding
collection of morphisms of affine schemes\/ $\Spec S_\alpha\rarrow
\Spec R$ is a finite open covering.
 Then an $R$\+module $K$ is coadjusted if and only if all
the $S_\alpha$\+modules $S_\alpha\ot_R K$ are coadjusted.
\end{lem}

\begin{proof}
 Part~(a): any contraadjusted $S$\+module $Q$ is also contraadjusted as
an $R$\+module, so the tensor product $(S\ot_R K)\ot_S Q\simeq
K\ot_R Q$ is a contraadjusted $R$\+module.
 By Corollary~\ref{very-open-cor}(d), it is also a contraadjusted
$S$\+module.

 Part~(b): the ``only if'' assertion is provided by part~(a); let us
prove the ``if''.
 Let $P$ be a contraadjusted $R$\+module.
 Applying the functor $K\ot_R{-}$ to the \v Cech exact
sequence~\eqref{cech-contra} from Lemma~\ref{very-open-covering}(b),
we obtain a sequence of $R$\+modules that is exact at its rightmost
nontrivial term.
 So it suffices to show that the $R$\+modules $K\ot_R\Hom_R(S_\alpha,P)$
are contraadjusted.

 Now one has $K\ot_R\Hom_R(S_\alpha,P)\simeq (S_\alpha\ot_R K)
\ot_{S_\alpha}\Hom_R(S_\alpha,P)$, the $S_\alpha$\+module
$\Hom_R(S_\alpha,P)$ is contraadjusted by
Corollary~\ref{very-open-cor}(c), and the restriction of scalars
from $S_\alpha$ to $R$ preserves contraadjustedness by
Lemma~\ref{very-scalars-always}(a).
\end{proof}

\begin{rem} \label{coadjusted-over-non-Noetherian-remark}
 Coadjusted modules over non-Noetherian commutative rings may be scarce
(this becomes a part of a problem which we encounter in
Chapter~\ref{semiderived-sect}; see Remark~\ref{semicosheaves-remark}).
 Still, one can show that any \emph{fp\+injective} $S$\+module (and
consequently any quotient $S$\+module of an fp\+injective module) is
coadjusted over an arbitrary commutative integral domain~$S$.
 Indeed, the direct sums of fp\+injective modules are always
fp\+injective~\cite[Corollary~2.4]{Sten0}, and it remains to observe
that the $S$\+modules $S[s^{-1}]$, \,$s\in S$, are fp\+projective
(in the sense of~\cite{MD,Pfp}) in the case of an integral domain.
 The obvious argument is that one has $S[s^{-1}]=0$ for $s=0$, and
for $s\ne0$ the $S$\+module $S[s^{-1}]$ is an infinitely iterated
extension of one copy of the free $S$\+module $S$ and a countable
sequence of copies of the cyclically presented $S$\+module $S/sS$.
 Both $S$ and $S/sS$ are finitely presented $S$\+modules.
 Therefore, all very flat $S$\+modules are fp\+projective, and all
fp\+injective $S$\+modules are contraadjusted.
 Together with the preservation of fp\+injectivity by infinite direct
sums, this implies the coadjustedness.
\end{rem}

\subsection{Modules adjusted over a base ring}
 In this section, whose purpose is to prepare ground for the discussion
in Sections~\ref{sheaves-injective-over-base-subsect}\+-%
\ref{cosheaves-projective-over-base-subsect}, we consider a homomorphism
of commutative rings $R\rarrow S$ and the properties of an $S$\+module
to be flat, very flat, or injective over~$R$.

 The following lemma is obvious.

\begin{lem} \label{flat-over-base-ring-tensor-product}
 Let $R\rarrow S$ be a homomorphism of commutative rings.
 Let $F$ be an $S$\+module that is flat over $R$, and let $G$ be a flat
$S$\+module.
 Then $F\ot_SG$ is a flat $R$\+module.
\end{lem}

\begin{proof}
 More generally, for any pair of associative rings $R$ and $S$, any
$R$\+flat $R$\+$S$\+bimodule $F$, and any flat $S$\+module $G$,
the $R$\+module $F\ot_SG$ is flat.
\end{proof}

\begin{lem} \label{flatness-over-base-is-local-for-rings}
 Let $R\rarrow S$ be a homomorphism of commutative rings.
 Let $S\rarrow T_\alpha$ be a collection of homomorphisms of commutative
rings for which the corresponding collection of morphisms of affine
schemes\/ $\Spec T_\alpha\rarrow\Spec S$ is a finite open covering.
 Then the $R$\+module $S$ is flat if and only if the $R$\+modules
$T_\alpha$ are flat for all indices~$\alpha$.
\end{lem}

\begin{proof}
 The ``only if'' implication holds because $T_\alpha$ is a flat
$S$\+module.
 To prove the ``if'', use the \v Cech exact sequence~\eqref{cech-rings}
for the collection of commutative ring homomorphisms
$S\rarrow T_\alpha$,
\begin{multline} \label{cech-rings-S-T} \textstyle
 0\lrarrow S\lrarrow \bigoplus_\alpha T_\alpha \lrarrow
 \bigoplus_{\alpha<\beta} T_\alpha\ot_S T_\beta \\
 \lrarrow\dotsb\lrarrow T_1\ot_S\dotsb\ot_S T_N \lrarrow 0.
\end{multline}
 The point is that, for any nonempty subset of indices
$\alpha_1$,~\dots, $\alpha_k$, \ $k\ge1$, the $R$\+module
$T_{\alpha_1}\ot_S\dotsb\ot_S T_{\alpha_k}$ is flat by
Lemma~\ref{flat-over-base-ring-tensor-product} (or by the ``only if''
implication applied to the ring homomorphisms $R\rarrow T_{\alpha_1}
\rarrow T_{\alpha_1}\ot_S\dotsb\ot_S T_{\alpha_k}$), and the kernel of
a surjective map of flat $R$\+modules is a flat $R$\+module.
\end{proof}

\begin{lem} \label{flatness-over-base-is-local}
 Let $S\rarrow T_\alpha$ be a collection of homomorphisms of commutative
rings for which the corresponding collection of morphisms of affine
schemes\/ $\Spec T_\alpha\rarrow\Spec S$ is a finite open covering.
 Let $F$ be an $S$\+module.
 Then $F$ is a flat $R$\+module if and only if $T_\alpha\ot_SF$ are
flat $R$\+modules for all indices~$\alpha$.
\end{lem}

\begin{proof}
 The ``only if'' implication holds by
Lemma~\ref{flat-over-base-ring-tensor-product}.
 The ``if'' follows from the \v Cech exact sequence~\eqref{cech-modules}
constructed for the collection of ring homomorphisms $S\rarrow T_\alpha$
and the $S$\+module~$F$ (just as in the previous lemma).
\end{proof}

 In the terminology of Section~\ref{local-classes-subsect},
Lemma~\ref{flatness-over-base-is-local} says that the property
of $S$\+modules to be flat over $R$ is \emph{local} with respect to
affine open coverings of $\Spec S$.
 This observation actually goes back to~\cite[Section~2.1]{Groth3};
see also~\cite[Section Tag~01U2]{SP}.

 The next six lemmas constitute a version of the preceding ones for
very flat rather than flat modules.

\begin{lem} \label{very-covering-rel-cech}
 Let $f_\alpha\:S\rarrow T_\alpha$, \ $\alpha=1$,~\dots,~$N$, be
a collection of homomorphisms of commutative rings for which
the corresponding collection of morphisms of affine schemes\/
$\Spec T_\alpha\rarrow \Spec S$ is a finite open covering, and let
$R\rarrow S$ be a homomorphism of commutative rings.
 Assume that all the $R$\+modules $T_{\alpha_1}\ot_S\dotsb\ot_S 
T_{\alpha_k}$, \ $k\ge1$, are very flat.
 Then the $R$\+module $S$ is very flat.
\end{lem}

\begin{proof}
 Follows from the \v Cech exact sequence~\eqref{cech-rings-S-T}, as
the kernel of a surjective map of very flat $R$\+modules is very flat.
\end{proof}

\begin{lem} \label{very-covering-rel-ascent}
 Let $S\rarrow T$ be a homomorphism of commutative rings such that
the corresponding morphism of affine schemes\/ $\Spec T\rarrow\Spec S$
is an open embedding, and let $R\rarrow S$ be a homomorphism of
commutative rings.
 Assume that the $R$\+module $S[s^{-1}]$ is very flat for all $s\in S$.
 Then the $R$\+module $T[t^{-1}]$ is very flat for all $t\in T$.
\end{lem}

\begin{proof}
 The $S$\+module $T[t^{-1}]$ is very flat by
Lemma~\ref{very-open-embedding}.
 By Lemma~\ref{very-scalars-veryflat-case}(b), it follows that
the $R$\+module $T[t^{-1}]$ is very flat.
 Alternatively, one can observe that $\Spec T[t^{-1}]$ as
an open subscheme in $\Spec S$ can be covered by a finite number of
principal affine open subschemes $\Spec S[s^{-1}]$, and
the intersections of these are also principal open affines.
 So Lemma~\ref{very-covering-rel-cech} is applicable.
\end{proof}

\begin{lem} \label{very-covering-rel}
 Let $f_\alpha\:S\rarrow T_\alpha$, \ $\alpha=1$,~\dots,~$N$, be
a collection of homomorphisms of commutative rings for which
the corresponding collection of morphisms of affine schemes\/
$\Spec T_\alpha\rarrow \Spec S$ is a finite open covering, and let
$R\rarrow S$ be a homomorphism of commutative rings.
 Then the $R$\+modules $T_\alpha[t_\alpha^{-1}]$ are very
flat for all\/ $t_\alpha\in T_\alpha$, \ $1\le\alpha\le N$, if and only
if the $R$\+module $S[s^{-1}]$ is very flat for all $s\in S$.
\end{lem}

\begin{proof}
 The ``if'' implication is provided by
Lemma~\ref{very-covering-rel-ascent}.
 To prove the ``only if'', notice that the ring of functions on any
affine open subscheme $U$ in $\Spec T_\alpha$ is a very flat $R$\+module
by the same Lemma~\ref{very-covering-rel-ascent} applied to
the ring homomorphisms $R\rarrow T_\alpha\rarrow\O(U)$.
 It remains to apply Lemma~\ref{very-covering-rel-cech} to the covering
of the affine scheme $\Spec S[s^{-1}]$ by its affine open subschemes
$\Spec T_\alpha[f_\alpha(s)^{-1}]$.
\end{proof}

\begin{lem} \label{very-covering-mod-cech}
 Let $f_\alpha\:S\rarrow T_\alpha$, \ $\alpha=1$,~\dots,~$N$, be
a collection of homomorphisms of commutative rings for which
the corresponding collection of morphisms of affine schemes\/
$\Spec T_\alpha\rarrow \Spec S$ is a finite open covering, and let
$R\rarrow S$ be a homomorphism of commutative rings.
 Let $F$ be an $S$\+module.
 Assume that all the $R$\+modules $T_{\alpha_1}\ot_S\dotsb\ot_S 
T_{\alpha_k}\ot_S F$, \ $k\ge1$, are very flat.
 Then the $R$\+module $F$ is very flat.
\end{lem}

\begin{proof}
 Just as in Lemma~\ref{very-covering-rel-cech}, the assertion follows
from the \v Cech exact sequence~\eqref{cech-modules} constructed for
the collection of commutative ring homomorphisms $S\rarrow T_\alpha$
and the $S$\+module~$F$.
\end{proof}

\begin{lem} \label{very-covering-mod-ascent}
 Let $S\rarrow T$ be a homomorphism of commutative rings such that
the corresponding morphism of affine schemes\/ $\Spec T\rarrow\Spec S$
is an open embedding, and let $R\rarrow S$ be a homomorphism of
commutative rings.
 Let $F$ be an $S$\+module.
 Assume that the $R$\+module $F[s^{-1}]$ is very flat for all $s\in S$.
 Then the $R$\+module $T[t^{-1}]\ot_SF$ is very flat for all $t\in T$.
\end{lem}

\begin{proof}
 Similar to the alternative proof of
Lemma~\ref{very-covering-rel-ascent}.
 Notice that $\Spec T[t^{-1}]$ as an open subscheme in $\Spec S$
can be covered by a finite number of principal affine open subschemes
$\Spec S[s^{-1}]$, and the intersections of these are also
principal open affines in $\Spec S$.
 So Lemma~\ref{very-covering-mod-cech} is applicable to
the $T[t^{-1}]$\+module $T[t^{-1}]\ot_SF$.
\end{proof}

\begin{lem} \label{very-covering-mod}
 Let $f_\alpha\:S\rarrow T_\alpha$, \ $\alpha=1$,~\dots,~$N$, be
a collection of homomorphisms of commutative rings for which
the corresponding collection of morphisms of affine schemes\/
$\Spec T_\alpha\rarrow \Spec S$ is a finite open covering, and let
$R\rarrow S$ be a homomorphism of commutative rings.
 Let $F$ be an $S$\+module.
 Then the $R$\+modules $T_\alpha[t_\alpha^{-1}]\ot_S F$ are very
flat for all\/ $t_\alpha\in T_\alpha$, \ $1\le\alpha\le N$, if and only
if the $R$\+module $F[s^{-1}]$ is very flat for all $s\in S$.
\end{lem}

\begin{proof}
 The ``if'' implication is provided by
Lemma~\ref{very-covering-mod-ascent}.
 To prove the ``only if'', notice that, for any affine open subscheme
$U$ in $\Spec T_\alpha$, the $\O(U)$\+module $\O(U)\ot_SF$ is a very
flat $R$\+module by the same Lemma~\ref{very-covering-mod-ascent}
applied to the ring homomorphisms $R\rarrow T_\alpha\rarrow\O(U)$
and the $T_\alpha$\+module $T_\alpha\ot_SF$.
 It remains to apply Lemma~\ref{very-covering-mod-cech} to
the covering of the affine scheme $\Spec S[s^{-1}]$ by its affine
open subschemes $\Spec T_\alpha[f_\alpha(s)^{-1}]$ and
the $S[s^{-1}]$\+module $F[s^{-1}]$.
\end{proof}

\begin{ex} \label{very-covering-mod-counterex}
 Notice that, even if $S$ is a very flat commutative $R$\+algebra
(i.~e., the $R$\+module $S[s^{-1}]$ is very flat for all $s\in S$),
the condition in Lemma~\ref{very-covering-mod-ascent} \emph{cannot} be
replaced by the condition that the $R$\+module $F$ is very flat.
 Indeed, the following example from~\cite[Example~9.6]{PSl1} shows
that for a very flat commutative $R$\+algebra $S$, an $S$\+module
$F$ very flat over $R$, and an element $s\in S$, the $R$\+module
$F[s^{-1}]$ \emph{need not} be very flat over~$R$.

 Put $R=\boZ$ and $S=\boZ[x]$; so $S$ the algebra of polynomials in
one variable with integer coefficients.
 Then $S$ is a very flat commutative $R$\+algebra
by Corollary~\ref{very-flat-algebra-cor} below.
 Let $P$ be the free abelian group with a countable set of generators
$p_1$, $p_2$, $p_3$,~\dots, viewed as an $S$\+module with the element
$x\in S$ acting by the operator $x(p_n)=(n+1)p_{n+1}$, \,$n\ge1$.
 Then $P$ is an $R$\+very flat (in fact, $R$\+free) $S$\+module.
 The $S$\+module $P[x^{-1}]$, viewed as an abelian group, is
isomorphic to the direct sum of countably many copies of the abelian
group of rational numbers~$\boQ$.
 According to Example~\ref{q-over-z-not-very-flat-ex} below,
$P[x^{-1}]$ is \emph{not} a very flat $R$\+module.

 In the terminology of Section~\ref{local-classes-subsect}, this
counterexample shows that the class of all $S$\+modules that are
very flat over $R$ is \emph{not} local, or more specifically, does
not satisfy ascent with respect to restrictions to (principal) affine
open subschemes in $\Spec S$.
 On the other hand, one can force the ascent to hold by including it
into the definition of the property we are interested in, i.~e.,
consider the class of all $S$\+modules $F$ such that the $R$\+module
$F[s^{-1}]$ is very flat for all $s\in S$.
 Then Lemma~\ref{very-covering-mod} says that the resulting property
of $S$\+modules $F$ is local.

 A dual-analogous difficulty with nonpreservation of $R$\+injectivity
by $S$\+module colocalizations is discussed in
Example~\ref{not-contrainjective-example} below.
\end{ex}

 Now we pass to the $R$\+injective $S$\+modules, under the assumption
that $R$ is a Noetherian ring.

\begin{lem} \label{injective-over-base-ring-tensor-product}
 Let $R$ be a Noetherian commutative ring and $R\rarrow S$ be
a commutative ring homomorphism.
 Let $K$ be an $S$\+module that is injective over $R$, and let $G$ be
a flat $S$\+module.
 Then $K\ot_SG$ is an injective $R$\+module.
\end{lem}

\begin{proof}
 This is a particular case of Lemma~\ref{coherent-tensor-hom-lemma}(a).
\end{proof}

\begin{lem} \label{pure-submodule-is-injective-lemma}
 Let $R$ be a Noetherian commutative ring and $R\rarrow S$ be
a commutative ring homomorphism.
 Let $K$ be an $S$\+module, and let\/ $0\rarrow H\rarrow G\rarrow F
\rarrow0$ be a short exact sequence of $S$\+modules with a flat
$S$\+module~$F$.
 Assume that the $R$\+module $G\ot_SK$ is injective.
 Then the $R$\+module $H\ot_SK$ is injective, too.
\end{lem}

\begin{proof}
 The point is that $0\rarrow H\ot_SK\rarrow G\ot_SK\rarrow F\ot_SK
\rarrow0$ is a \emph{pure} short exact sequence of $S$\+modules
(hence also a pure short exact sequence of $R$\+modules) in the sense
of, e.~g., \cite[Section~2.1]{GT}.
 The purity means that the sequence remains exact after
taking the tensor product over $S$ with any $S$\+module (hence also
after taking the tensor product over $R$ with any $R$\+module).
 The pure submodules of injective modules are precisely all
the fp\+injective (``absolutely pure'') modules.
 Over a Noetherian ring $R$, all fp\+injective modules are injective.
 Essentially, this argument is based on the equivalence of the two
definitions of a pure short exact sequence of
modules~\cite[Definition~2.6 and
Lemma~2.19\,(a)\,$\Leftrightarrow$\,(d)]{GT}.
 To be more specific, the purity of a short exact sequence of
$R$\+modules $0\rarrow H'\rarrow G'\rarrow F'\rarrow0$ also means
surjectivity of the map $\Hom_R(M,G')\rarrow\Hom_R(M,F')$ for any
finitely presented $R$\+module $M$.
 This implies injectivity of the map $\Ext^1_R(M,H')\rarrow
\Ext^1_R(M\;G')$.
 In the situation at hand, take $H'=H\ot_SK$ and $G'=G\ot_SK$.
\end{proof}

\begin{cor} \label{injectivity-over-noetherian-base-is-local}
 Let $R$ be a Noetherian commutative ring and $R\rarrow S$ be
a commutative ring homomorphism.
 Let $S\rarrow T_\alpha$ be a collection of homomorphisms of commutative
rings for which the corresponding collection of morphisms of affine
schemes\/ $\Spec T_\alpha\rarrow\Spec S$ is a finite open covering.
 Let $K$ be an $S$\+module.
 Then $K$ is an injective $R$\+module if and only if $T_\alpha\ot_SK$
are injective $R$\+modules for all indices~$\alpha$.
\end{cor}

\begin{proof}
 The ``only if'' implication holds by
Lemma~\ref{injective-over-base-ring-tensor-product}.
 To prove the ``only if'', consider the \v Cech exact
sequence~\eqref{cech-rings-S-T} from the proof of
Lemma~\ref{flatness-over-base-is-local-for-rings}.
 This is a finite exact sequence of flat $S$\+modules.
 Hence the cokernel $F$ of the injective map
$S\rarrow\bigoplus_\alpha T_\alpha$ is a flat $S$\+module.
 It remains to apply Lemma~\ref{pure-submodule-is-injective-lemma}
to the short exact sequence of $S$\+modules
$0\rarrow S\rarrow\bigoplus_\alpha T_\alpha\rarrow F\rarrow0$ and
the $S$\+module~$K$.
\end{proof}

 In the terminology of Section~\ref{local-classes-subsect},
Corollary~\ref{injectivity-over-noetherian-base-is-local} tells us that
the property of $S$\+modules to be injective over $R$ is local with
respect to affine open coverings of $\Spec S$ (when the ring $R$ is
Noetherian).

 Finally, we are interested in the question whether the property of
$S$\+modules to be flat over $R$ is \emph{colocal}.
 Here we can only prove the coascent.

\begin{lem} \label{commutative-coherent-cotorsper-lemma}
 Let $R\rarrow S$ be a homomorphism of commutative rings, where
the ring $R$ is coherent.
 Let $P$ be an $S$\+module whose underlying $R$\+module is flat,
and let $F$ be an $S$\+module. \par
\textup{(a)} If the $S$\+module $F$ is very flat and the $S$\+module $P$
is contraadjusted, then the underlying $R$\+module of
the (contraadjusted) $S$\+module\/ $\Hom_S(F,P)$ is flat. \par
\textup{(b)} If the $S$\+module $F$ is flat and the $S$\+module $P$
is cotorsion, then the underlying $R$\+module of the (cotorsion)
$S$\+module\/ $\Hom_S(F,P)$ is flat.
\end{lem}

\begin{proof}
 Part~(a) is a particular case of
Lemma~\ref{coherent-lemma}.
 Part~(b) is a particular case of
Lemma~\ref{coherent-cotorsper-lemma}.
\end{proof}

\begin{cor} \label{flatness-over-coherent-base-coascent}
 Let $R$ be a coherent commutative ring and $R\rarrow S$ be
a commutative ring homomorphism.
 Let $S\rarrow T$ be a commutative ring homomorphism for which
the corresponding morphism of affine schemes\/ $\Spec T\rarrow\Spec S$
is an open embedding.
 Let $P$ be a contraadjusted $S$\+module that is flat as
an $R$\+module.
 Then the (contraadjusted) $T$\+module\/ $\Hom_S(T,P)$ is also flat
as an $R$\+module.
\end{cor}

\begin{proof}
 Follows from Lemmas~\ref{commutative-coherent-cotorsper-lemma}(a)
and~\ref{very-open-embedding}.
\end{proof}

\begin{qst} \label{codescent-of-flatness-question}
 We do \emph{not} know how to prove a dual-analogous version of
Lemma~\ref{pure-submodule-is-injective-lemma} for the Hom modules
instead of the tensor products and flatness instead of the injectivity.
 Specifically, let $R\rarrow S$ be a homomorphism of commutative rings
with a Noetherian ring~$R$.
 Let us assume additionally that the ring $R$ has finite Krull
dimension (or even admits a dualizing complex), and that $S$ is
a flat $R$\+module (or even a very flat commutative $R$\+algebra).
 Let $0\rarrow H\rarrow G\rarrow F\rarrow0$ be a short exact sequence
of flat $S$\+modules, and let $P$ be a cotorsion $S$\+module.
 Then we have a short exact sequence of $S$\+modules $0\rarrow
\Hom_S(F,P)\rarrow\Hom_S(G,P)\rarrow\Hom_S(H,P)\rarrow0$.
 Let us even assume that the $S$\+modules $F$, $G$, and $H$ are
countably presented and very flat.
 Suppose that the underlying $R$\+modules of the $S$\+modules
$\Hom_S(G,P)$ and $\Hom_S(F,P)$ are flat.
 Does it follow that the underlying $R$\+module of the $S$\+module
$\Hom_S(H,P)$ is flat?
 In view of this difficulty, we \emph{cannot} say whether the property
of (say, cotorsion) $S$\+modules to be flat over $R$ is colocal (in
the sense of Section~\ref{colocal-classes-subsect}) with respect to
affine open coverings of $\Spec S$.
 In the context of the latter question, the case of $S=R$ is covered
by Corollary~\ref{coherent-flat-local}.
\end{qst}

\subsection{Contrainjective and cotorsinjective modules}
\label{contrainjective-subsect}
 In this section, whose purpose is to prepare ground for the discussion
in Section~\ref{cosheaves-loc-injective-over-base-subsect},
we continue to consider a commutative ring homomorphism $R\rarrow S$ and
the property of an $S$\+module to be injective over~$R$.
 We discuss the (non)colocality of this property.

\begin{lem} \label{contrainjective-lemma}
 The following four conditions on an $S$\+module $K$ are equivalent:
\begin{enumerate}
\item $K$ is injective as an $R$\+module, and for every (finitely
generated) $R$\+module $M$, the $S$\+module\/ $\Hom_R(M,K)$ is
contraadjusted;
\item $K$ is injective as an $R$\+module, and for any short exact
sequence of very flat $S$\+modules\/ $0\rarrow H\rarrow G\rarrow F
\rarrow 0$ and any (finitely generated) $R$\+module $M$, the short
sequence of $S$\+modules\/ $0\rarrow\Hom_S(F,\Hom_R(M,K))\allowbreak
\rarrow\Hom_S(G,\Hom_R(M,K))\rarrow\Hom_S(H,\Hom_R(M,K))\rarrow0$
is exact;
\item $K$ is contraadjusted as an $S$\+module, and for every very
flat $S$\+module $F$ (of the form $F=S[s^{-1}]$, \,$s\in S$),
the $S$\+module\/ $\Hom_S(F,K)$ is injective as an $R$\+module;
\item $K$ is contraadjusted as an $S$\+module, and for every very
flat $S$\+module $F$ (of the form $F=S[s^{-1}]$, \,$s\in S$) and any
short exact sequence of $R$\+modules\/ $0\rarrow L\rarrow N\rarrow M
\rarrow0$, the short sequence of $S$\+modules\/ $0\rarrow\Hom_R(M,
\Hom_S(F,K))\rarrow\Hom_R(N,\Hom_S(F,K))\rarrow\Hom_R(L,\Hom_S(F,K))
\rarrow0$ is exact.
\end{enumerate}
 When $S$ is a flat $R$\+module, the conditions~\textup{(1\+-4)} are
also equivalent to the following condition
\begin{enumerate}
\setcounter{enumi}{4}
\item for every very flat $S$\+module $F$ (of the form $F=S[s^{-1}]$,
\,$s\in S$) and any (finitely generated) $R$\+module $M$, one has\/
$\Ext^1_S(F\ot_RM\;K)=0$.
\end{enumerate}
\end{lem}

\begin{proof}
 The equivalence (1)\,$\Longleftrightarrow$\,(2) holds because
an $S$\+module $C$ is contraadjusted if and only if, for every short
exact sequence of very flat $S$\+modules $0\rarrow H\rarrow G\rarrow F
\rarrow 0$, the short sequence of $S$\+modules $0\rarrow\Hom_S(F,C)
\rarrow\Hom_S(G,C)\rarrow\Hom_S(H,C)\rarrow0$ is exact.
 Moreover, the condition~(2) does not change when restricted to
the special case $F=S[s^{-1}]$.
 
 Similarly, the equivalence (3)\,$\Longleftrightarrow$\,(4) holds
because an $R$\+module $J$ is injective if and only if the functor
$\Hom_R({-},J)$ takes short exact sequences of $R$\+modules to short
exact sequences.
 Moreover, the condition~(4) does not change when one assumes
the $R$\+module $M$ to be finitely generated.

 (1)\,$\Longrightarrow$\,(4) Taking $M=R$ in~(1), we conclude
that the $S$\+module $K$ is contraadjusted.
 Forthermore, the short sequence of $S$\+modules
$0\rarrow\Hom_R(M,K)\rarrow\Hom_R(N,K)\rarrow\Hom_R(L,K)\rarrow0$
is exact, since the $R$\+module $K$ is injective by~(1).
 Now contraadjustedness of the $S$\+module $\Hom_R(M,K)$ stated in~(1)
implies exactness of the short sequence $0\rarrow\Hom_S(F,\Hom_R(M,K))
\rarrow\Hom_S(F,\Hom_R(N,K))\rarrow\Hom_S(F,\Hom_R(L,K))\rarrow0$,
as desired in~(4).
 Here we use the natural isomorphism $\Hom_S(F,\Hom_R(M,K))\simeq
\Hom_R(M,\Hom_S(F,K))$.

 (3)\,$\Longrightarrow$\,(2) Taking $F=S$ in~(3), we conclude
that the $S$\+module $K$ is injective as an $R$\+module.
 Furthermore, the short sequence of $S$\+modules
$0\rarrow\Hom_S(F,K)\rarrow\Hom_S(G,K)\rarrow\Hom_S(H,K)\rarrow0$
is exact, since the $S$\+module $K$ is contraadjusted by~(3).
 Now injectivity of the $R$\+module $\Hom_S(F,K)$ stated in~(3)
implies exactness of the short sequence $0\rarrow\Hom_R(M,\Hom_S(F,K))
\rarrow\Hom_R(M,\Hom_S(G,K)\rarrow\Hom_R(M,\Hom_S(H,K))\rarrow0$,
as desired in~(2).
 Once again, we use the natural isomorphism of abelian groups or
$S$\+modules $\Hom_S(F,\Hom_R(M,K))\simeq\Hom_R(M,\Hom_S(F,K))$.

 (5)\,$\Longleftrightarrow$\,(2) If $S$ is a flat $R$\+module,
then (very) flat $S$\+modules are flat over~$R$; so the short sequence
of $S$\+modules $0\rarrow H\ot_RM\rarrow G\ot_RM\rarrow F\ot_RM\rarrow0$
is exact.
 Now the condition $\Ext^1_S(F\ot_RM\;K)=0$ implies exactness of
the short sequence $0\rarrow\Hom_S(F\ot_RM\;K)\rarrow\Hom_S(G\ot_RM\;K)
\rarrow\Hom_S(H\ot_RM\;K)\rarrow0$, as desired in~(2).
 Here we use the natural isomorphism $\Hom_S(F\ot_RM\;K)\simeq
\Hom_S(F,\Hom_R(M,K))$.
 Furthermore, the natural isomorphism $\Ext^1_R(M,K)\simeq
\Ext^1_S(S\ot_RM\;K)$, which holds whenever $S$ is a flat $R$\+module,
proves injectivity of $K$ as an $R$\+module (take $F=S$ in~(5)).
 Conversely, to deduce~(5) from~(2), take $G$ to be a free $S$\+module;
then injectivity of the $R$\+module~$K$ (stated in~(2)) implies
$\Ext^1_S(G\ot_RM\;K)=0$, and therefore surjectivity of the map
$\Hom_S(G\ot_RM\;K)\rarrow\Hom_S(H\ot_RM\;K)$ implies
$\Ext^1_S(F\ot_RM\;K)=0$.

 (5)\,$\Longleftrightarrow$\,(4) Once again, if $S$ is a flat
$R$\+module, the so is~$F$; hence the short sequence of $S$\+modules
$0\rarrow F\ot_RL\rarrow F\ot_RN\rarrow F\ot_RM\rarrow0$ is exact.
  Now the condition $\Ext^1_S(F\ot_RM\;K)=0$ implies exactness of
the short sequence $0\rarrow\Hom_S(F\ot_RM\;K)\rarrow\Hom_S(F\ot_RN\;K)
\rarrow\Hom_S(F\ot_RL\;K)\rarrow0$, as desired in~(4).
 Here we use the natural isomorphism $\Hom_S(F\ot_RM\;K)\simeq
\Hom_R(M,\Hom_S(F,K))$.
 Furthermore, taking $M=R$ in~(5) shows that the $S$\+module $K$ is
contraadjusted.
 Conversely, to deduce~(5) from~(4), take $N$ to be a free $R$\+module;
then contraadjustedness of the $S$\+module~$K$ (stated in~(4)) implies
$\Ext^1_S(F\ot_RN\;K)=0$, and therefore surjectivity of the map
$\Hom_S(F\ot_RN\;K)\rarrow\Hom_S(F\ot_RL\;K)$ implies
$\Ext^1_S(F\ot_RM\;K)=0$.
\end{proof}

 We will say that an $S$\+module $K$ is \emph{$(S,R)$\+contrainjective}
if it satisfies the equivalent conditions of
Lemma~\ref{contrainjective-lemma}(1\+-4).

\begin{cor} \label{contrainjective-cor}
\textup{(a)} The class of $(S,R)$\+contrainjective modules is closed
under extensions and cokernels of injective morphisms in $S\modl$. \par
\textup{(b)} If\/ $0\rarrow E\rarrow D\rarrow C\rarrow0$ is a short
exact sequence of $S$\+modules, the $S$\+module $D$ is
$(S,R)$\+contrainjective, and the $S$\+module $E$ is injective as
an $R$\+module, then the $S$\+module $C$ is $(S,R)$\+contrainjective.
\par
\textup{(c)} If the ring $S$ is a flat $R$\+module, then any injective
$S$\+module is $(S,R)$\+con\-tra\-in\-jec\-tive.
\end{cor}

\begin{proof}
 Part~(a) follows from the criterion of
Lemma~\ref{contrainjective-lemma}(1) or~(3).
 Part~(b) follows from the criterion of
Lemma~\ref{contrainjective-lemma}(1)
(recall that any quotient module of a contraadjusted module is
contraadjusted).
 Part~(c) follows from Lemma~\ref{contrainjective-lemma}(5).
\end{proof}

 The following lemma is a version of Lemma~\ref{contrainjective-lemma}
with contraadjusted modules replaced by cotorsion ones.

\begin{lem} \label{cotorsinjective-lemma}
 The following four conditions on an $S$\+module $K$ are equivalent:
\begin{enumerate}
\item $K$ is injective as an $R$\+module, and for every (finitely
generated) $R$\+module $M$, the $S$\+module\/ $\Hom_R(M,K)$ is
cotorsion;
\item $K$ is injective as an $R$\+module, and for any short exact
sequence of flat $S$\+modules\/ $0\rarrow H\rarrow G\rarrow F
\rarrow 0$ and any (finitely generated) $R$\+module $M$, the short
sequence of $S$\+modules\/ $0\rarrow\Hom_S(F,\Hom_R(M,K))\allowbreak
\rarrow\Hom_S(G,\Hom_R(M,K))\rarrow\Hom_S(H,\Hom_R(M,K))\rarrow0$
is exact;
\item $K$ is cotorsion as an $S$\+module, and for every flat $S$\+module
$F$, the $S$\+module\/ $\Hom_S(F,K)$ is injective as an $R$\+module;
\item $K$ is cotorsion as an $S$\+module, and for every flat $S$\+module
$F$ and any short exact sequence of $R$\+modules\/ $0\rarrow L\rarrow N
\rarrow M\rarrow0$, the short sequence of $S$\+modules\/ $0\rarrow
\Hom_R(M,\Hom_S(F,K))\rarrow\Hom_R(N,\Hom_S(F,K))\rarrow
\Hom_R(L,\Hom_S(F,K))\rarrow0$ is exact.
\end{enumerate}
 When $S$ is a flat $R$\+module, the conditions~\textup{(1\+-4)} are
also equivalent to the following condition
\begin{enumerate}
\setcounter{enumi}{4}
\item for every flat $S$\+module $F$ and any (finitely generated)
$R$\+module $M$, one has\/ $\Ext^1_S(F\ot_RM\;K)=0$.
\end{enumerate}
\end{lem}

\begin{proof}
 Similar to the proof of Lemma~\ref{contrainjective-lemma}.
\end{proof}

 We will say that an $S$\+module $K$ is \emph{$(S,R)$\+cotorsinjective}
if it satisfies the equivalent conditions of
Lemma~\ref{cotorsinjective-lemma}(1\+-4).
 Obviously, any $(S,R)$\+cotorsinjective $S$\+module is
$(S,R)$\+contrainjective.

\begin{cor} \label{cotorsinjective-cor}
\textup{(a)} The class of $(S,R)$\+cotorsinjective modules is closed
under extensions and cokernels of injective morphisms in $S\modl$. \par
\textup{(b)} Let $K^\bu$ be an acyclic complex of
$(S,R)$\+cotorsinjective $S$\+modules.
 Assume that the $S$\+modules of cocycles of the complex $K^\bu$ are
injective as $R$\+modules.
 Then the $S$\+modules of cocycles of the complex $K^\bu$ are actually
$(S,R)$\+cotorsinjective. \par
\textup{(c)} If the ring $S$ is a flat $R$\+module, then any injective
$S$\+module is $(S,R)$\+co\-tors\-in\-jec\-tive.
\end{cor}

\begin{proof}
 Part~(a) follows from the criterion of
Lemma~\ref{cotorsinjective-lemma}(1) or~(3).
 To deduce part~(b) from the criterion of
Lemma~\ref{cotorsinjective-lemma}(1), notice that, for any
$R$\+module $M$, the complex of $S$\+modules $\Hom_R(M,K^\bu)$ is
acyclic, and its modules of cocycles are obtained by applying
the functor $\Hom_R(M,{-})$ to the modules of cocycles of
the complex~$K^\bu$.
 It remains to refer to the cotorsion periodicity theorem for
$S$\+modules (Theorem~\ref{cotorsion-periodicity}).
 Part~(c) follows immediately from Lemma~\ref{cotorsinjective-lemma}(5).
\end{proof}

\begin{lem} \label{contra-cotors-injective-Hom}
\textup{(a)} For any very flat $S$\+module $F$ and any
$(S,R)$\+contrainjective $S$\+module $K$, the $S$\+module\/
$\Hom_S(F,K)$ is $(S,R)$\+contrainjective. \par
\textup{(b)} For any flat $S$\+module $F$ and any
$(S,R)$\+cotorsinjective $S$\+module $K$, the $S$\+module\/
$\Hom_S(F,K)$ is $(S,R)$\+cotorsinjective.
\end{lem}

\begin{proof}
 Part~(a) follows from the criterion of
Lemma~\ref{contrainjective-lemma}(3) and
Lemma~\ref{very-tensor-hom}(a\+b) (for the ring~$S$).
 Part~(b) follows from the criterion of
Lemma~\ref{cotorsinjective-lemma}(3) and
Lemma~\ref{cotors-hom}(a) (for the ring~$S$).
\end{proof}

\begin{lem} \label{contra-cotors-injective-restriction-of-scalars}
 Let $R\rarrow S\rarrow T$ be two homomorphisms of commutative rings.
 Then \par
\textup{(a)} any $(T,R)$\+contrainjective $T$\+module is
$(S,R)$\+contrainjective as an $S$\+module; \par
\textup{(b)} any $(T,R)$\+cotorsinjective $T$\+module is
$(S,R)$\+cotorsinjective as an $S$\+module.
\end{lem}

\begin{proof}
 Part~(a): follows from the criterion of
Lemma~\ref{contrainjective-lemma}(1) or~(3) together with
Lemma~\ref{very-scalars-always}(a\+b) (applied
to the commutative ring homomorphism $S\rarrow T$).
 Part~(b): follows from the criterion of
Lemma~\ref{cotorsinjective-lemma}(1) or~(3) together with
Lemma~\ref{cotors-restrict}(a) (applied to the ring
homomorphism $S\rarrow T$).
\end{proof}

\begin{lem} \label{contra-cotors-injective-coextension-of-scalars}
 Let $R\rarrow S\rarrow T$ be two homomorphisms of commutative rings.
 In this context: \par
\textup{(a)} Assume that $T$ is a very flat $S$\+algebra.
 Then, for any $(S,R)$\+contrainjective $S$\+module $K$,
the $T$\+module\/ $\Hom_S(T,K)$ is $(T,R)$\+contrainjective. \par
\textup{(b)} Assume that $T$ is a flat $S$\+module.
 Then, for any $(S,R)$\+cotorsinjective $S$\+module $K$,
the $T$\+module\/ $\Hom_S(T,K)$ is $(T,R)$\+cotorsinjective.
\end{lem}

\begin{proof}
 Part~(a): use the criterion of Lemma~\ref{contrainjective-lemma}(3)
together with Lemma~\ref{very-scalars-veryflat-case}(a\+b) (applied
to the commutative ring homomorphism $S\rarrow T$).
 Part~(b): use the criterion of Lemma~\ref{cotorsinjective-lemma}(3)
together with Lemma~\ref{cotors-coexten}(a) (applied to the ring
homomorphism $S\rarrow T$).
\end{proof}

\begin{lem} \label{contra-cotors-injective-covering-mod}
 Let $f_\alpha\:S\rarrow T_\alpha$, \ $\alpha=1$,~\dots,~$N$, be
a collection of homomorphisms of commutative rings for which
the corresponding collection of morphisms of affine schemes\/
$\Spec T_\alpha\rarrow \Spec S$ is a finite open covering, and let
$R\rarrow S$ be  a homomorphism of commutative rings.
 Let $K$ be a contraadjusted $S$\+module.
 In this context: \par
\textup{(a)} The $S$\+module $K$ is $(S,R)$\+contrainjective if and
only if, for every\/ $1\le\alpha\le N$, the $T_\alpha$\+module\/
$\Hom_S(T_\alpha,K)$ is $(T_\alpha,R)$\+contrainjective. \par
\textup{(b)} The $S$\+module $K$ is $(S,R)$\+cotorsinjective if and
only if, for every\/ $1\le\alpha\le N$, the $T_\alpha$\+module\/
$\Hom_S(T_\alpha,K)$ is $(T_\alpha,R)$\+cotorsinjective.
\end{lem}

\begin{proof}
 Part~(a): the ``only if'' implication is provided by
Lemma~\ref{contra-cotors-injective-coextension-of-scalars}(a).
 To prove the ``if'', one can use the \v Cech
sequence~\eqref{cech-contra} from Lemma~\ref{very-open-covering}(b)
together with
Lemmas~\ref{contra-cotors-injective-restriction-of-scalars}(a)
and~\ref{contra-cotors-injective-coextension-of-scalars}(a).
 This argument works because the cokernel of any injective morphism
of $(S,R)$\+contrainjective $S$\+modules is $(S,R)$\+contrainjective
(Corollary~\ref{contrainjective-cor}(a)).
 The proof of part~(b) is similar.
\end{proof}

\begin{ex} \label{not-contrainjective-example}
 The following example proves that an $S$\+cotorsion $R$\+injective
$S$\+module need not be $(S,R)$\+contrainjective, even when $S$ is
a very flat commutative $R$\+algebra (in the sense of
Section~\ref{affine-geometry-subsect}).
 We use the terminology and results of
Sections~\ref{cotorsion-prelim-subsect}\+-%
\ref{small-object-argument-subsect} concerning cotorsion pairs
and ordinal-indexed filtrations.

 Assume that $S$ is a flat $R$\+module.
 Then Lemma~\ref{contrainjective-lemma}(5) tells us that the class of
$(S,R)$\+contrainjective $S$\+modules is the right-hand class in
the cotorsion pair in $S\modl$ generated by $S$\+modules of the form
$S[s^{-1}]\ot_RM$, where $s\in S$ and $M$ ranges over finitely
generated $R$\+modules.
 On the other hand, the class of $S$\+cotorsion $R$\+injective
$S$\+modules is the right-hand class in the cotorsion pair in $S\modl$
generated by $S$\+modules of the form $S\ot_RM$ together with
a suitable set of flat $S$\+modules (as in
Example~\ref{deconstructible-examples}).
 By the Eklof--Trlifaj theorem (\cite[Theorems~2 and~10]{ET} or
Theorem~\ref{eklof-trlifaj-general}(b) and
Corollary~\ref{eklof-trlifaj-corollary}), it follows that the left
$\Ext^1_S$\+orthogonal class to the class of $(S,R)$\+contrainjective
$S$\+modules is the class of all direct summands of modules filtered
by $S[s^{-1}]\ot_RM$, while the left $\Ext^1_S$\+orthogonal class to
the class of $S$\+cotorsion $R$\+injective $S$\+modules is the class
of all direct summands of modules filtered by $S\ot_RM$ and flat
$S$\+modules.
 It remains to present an example of $S$\+module of the form
$S[s^{-1}]\ot_RM$ that is \emph{not} a direct summand of any $S$\+module
filtered by the $S$\+modules $S\ot_RN$ (with $N\in R\modl$) and flat
$S$\+modules.

 Let $k$~be a field.
 Consider the algebra of polynomials in one variable $R=k[x]$ and
the algebra of polynomials in two variables $S=k[x,y]$; so $R$ is
a subring in~$S$.
 By Corollary~\ref{very-flat-algebra-cor} below, $S$ is a very flat
$R$\+algebra.
 Put $s=y\in S$ and $M=k$, with the zero action of~$x$ in~$M$.
 Then $S[s^{-1}]\ot_RM$ is the $S$\+module $k[y,y^{-1}]$, with a zero
action of~$x$.
 Let us show that $S[s^{-1}]\ot_RM$ is \emph{not} a submodule of
an $S$\+module filtered by $S$\+modules $S\ot_RN$ and flat $S$\+modules.
 Moreover, let us check that $S[s^{-1}]\ot_RM$ is not filtered by
submodules of the $S$\+modules $S\ot_RN$ and submodules of
flat $S$\+modules.

 Firstly, no nonzero subquotient module of the $S$\+module
$S[s^{-1}]\ot_RM$ is a submodule of a flat $S$\+module, because
the action of~$x$ is injective in submodules of flat $S$\+modules
and zero in subquotient modules of $S[s^{-1}]\ot_RM$.
 Consequently, if $S[s^{-1}]\ot_RM$ is filtered by submodules of
the $S$\+modules $S\ot_RN$ and submodules of flat $S$\+modules,
then $S[s^{-1}]\ot_RM$ is actually filtered by submodules of
the $S$\+modules $S\ot_RN$ only.
 Secondly, the underlying $k[y]$\+module of $S[s^{-1}]\ot_RM$ is
$k[y,y^{-1}]$, while the underlying $k[y]$\+modules of $S\ot_RN$
are free $k[y]$\+modules.
 So it remains to show that the $k[y]$\+module $k[y,y^{-1}]$ is not
filtered by submodules of free $k[y]$\+modules.
 Indeed, the $k[y]$\+module $k[y,y^{-1}]$ cannot be filtered by
$k[y]$\+modules containing \emph{neither} nonzero $y$\+torsion elements
\emph{nor} nonzero infinitely $y$\+divisible elements, as one can easily
see by considering the bottom nonzero submodule in such a filtration
of~$k[y,y^{-1}]$.

 In other words, we have shown that there exist an $R$\+injective
$S$\+cotorsion $S$\+module $P$ and an element $s\in S$ such that
the $S$\+module $\Hom_S(S[s^{-1}],P)$ is \emph{not} $R$\+injective
(cf.\ Lemma~\ref{contrainjective-lemma}(3)).
 There is also a finitely generated $R$\+module $M$ such that
the $S$\+module $\Hom_R(M,P)$ is \emph{not} contraadjusted
(see Lemma~\ref{contrainjective-lemma}(1)).

 In the terminology of Section~\ref{colocal-classes-subsect}, this
counterexample proves that the class of all contraadjusted $S$\+modules
that are injective over $R$ is \emph{not} colocal, or more specifically,
does not satisfy coascent with respect to restrictions to (principal)
affine open subschemes in $\Spec S$.
 On the other hand, one can force the coascent to hold by including it
into the definition of the property we are interested in, i.~e.,
consider the class of all contraadjusted $S$\+modules $K$ such that
the $S$\+module $\Hom_S(S[s^{-1}],K)$ is injective as an $R$\+module
for all $s\in S$.
 According to the criterion of Lemma~\ref{contrainjective-lemma}(3),
this gives precisely the class of all $(S,R)$\+contrainjective
$S$\+modules.
 Then Lemma~\ref{contra-cotors-injective-covering-mod}(a) says that
the resulting (contrainjectivity) property of $S$\+modules $K$
is colocal.

 For a discussion of the dual-analogous difficulty with nonpreservation
of $R$\+very flatness by $S$\+module localizations, see
Example~\ref{very-covering-mod-counterex} above.
\end{ex}

 On the other hand, for any commutative ring homomorphism $R\rarrow S$
with a Noetherian ring $R$, the property of $S$\+modules to be injective
over $R$ is \emph{local} (in the sense of
Section~\ref{local-classes-subsect}) with respect to affine open
coverings of $\Spec S$, as we have seen in
Corollary~\ref{injectivity-over-noetherian-base-is-local}.
 The dual-analogous question, about the colocality of flatness over $R$
for (contraadjusted or cotorsion) $S$\+modules, seems to be open:
the coascent holds by
Corollary~\ref{flatness-over-coherent-base-coascent}, but the codescent
is not clear (see Question~\ref{codescent-of-flatness-question}).

\subsection{Very flat morphisms of schemes}
\label{very-flat-morphisms-subsect}
 A morphism of schemes $f\:Y\rarrow X$ is called \emph{very flat} if,
for any two affine open subschemes $V\sub Y$ and $U\sub X$ such that
$f(V)\sub U$, the ring of regular functions $\O(V)$ is a very flat
module over the ring $\O(U)$.
 By Lemma~\ref{very-open-embedding}, any embedding of an open subscheme
is a very flat morphism.

 According to Lemmas~\ref{very-open-covering}(a)
and~\ref{very-covering-rel}, the property of a morphism to be
very flat is local in both the source and the target schemes.
 A morphism of affine schemes $\Spec S\rarrow\Spec R$ is very flat
if and only if the morphism of commutative rings $R\rarrow S$ is very
flat, or in other words, if and only if $S$ is a very flat commutative
$R$\+algebra in the sense of the definition in
Section~\ref{affine-geometry-subsect}.
 By Lemma~\ref{very-scalars-veryflat-case}(b), the composition of
very flat morphisms of schemes is a very flat morphism.

 It does not seem to follow from anything that the base change of
a very flat morphism of schemes should be a very flat morphism.
 Here is a partial result in this direction.

\begin{lem} \label{very-flat-injective-base-change}
 The base change of a very flat morphism with respect to any locally
closed embedding or universal homeomorphism of schemes is
a very flat morphism.
\end{lem}

\begin{proof}
 Essentially, given a very flat morphism $f\:Y\rarrow X$ and 
a morphism of schemes $g\:x\rarrow X$, in order to conclude
that the morphism $f'\:y=x\times_x Y\rarrow x$ is very flat it suffices
to know that the morphism $g'\:y\rarrow Y$ is injective and
the topology of~$y$ is induced from the topology of~$Y$.

 Indeed, the very flatness condition being local, one can assume all
the four schemes to be affine.
 Now any affine open subscheme $v\sub y$ is the preimage of
a certain open subscheme $V\sub Y$.
 Covering $V$ with open affines if necessary and using the locality
again, one can assume that $V$ is affine, too.
 Finally, if the $\O(X)$\+module $\O(V)$ is very flat, then
by Lemma~\ref{very-scalars-always}(b) so is the $\O(x)$\+module
$\O(v)=\O(x\times_XV)$.
\end{proof}

 A quasi-coherent sheaf $\F$ on a scheme $X$ is called
\emph{very flat} if the $\O_X(U)$\+module $\F(U)$ is very flat for
any affine open subscheme $U\sub X$.
 According to Lemma~\ref{very-open-covering}(a), very flatness of
a quasi-coherent sheaf on a scheme is a local property, so it suffices
to check the condition above for open subschemes belonging to any
given affine open covering of the scheme~$X$.
 By Lemma~\ref{very-scalars-always}(b), the inverse image of a very
flat quasi-coherent sheaf under any morphism of schemes is very flat.
 By Lemma~\ref{very-scalars-veryflat-case}(b), the direct image of
a very flat quasi-coherent sheaf under a very flat affine morphism of
schemes is very flat.

 More generally, given a morphism of schemes $f\:Y\rarrow X$,
a quasi-coherent sheaf $\F$ on $Y$ is said to be \emph{very flat
over\/~$X$} if for any affine open subschemes $U\sub X$ and
$V\sub Y$ such that $f(V)\sub U$ the module of sections $\F(V)$
is very flat over the ring $\O_X(U)$.
 According to Lemmas~\ref{very-open-covering}(a)
and~\ref{very-covering-mod}, the property of very flatness of $\F$
over $X$ is local in both $X$ and~$Y$; in other words, it suffices to 
check this condition for affine open subschemes $U\sub X$ and $V\sub Y$
subordinate to some given open coverings of the schemes $X$ and~$Y$.
 (Here an open subscheme $U$ in $X$ is said to be \emph{subordinate} to
an open covering $\bW$ of $X$ if $U$ is contained in one of the open
subsets of $X$ belonging to $\bW$; cf.
Section~\ref{locally-contraherent} below.)
 By Lemma~\ref{very-scalars-veryflat-case}(b), if the scheme $Y$ is
very flat over~$X$ (i.~e., the morphism~$f$ is very flat) and
a quasi-coherent sheaf $\F$ is very flat on $Y$, then $\F$ is also very
flat over~$X$.

 Let us issue the following \emph{warning}, however.
 Even when the morphism~$f$ is affine and very flat, the condition
that the $\O_X(U)$\+module $\F(f^{-1}(U))$ is very flat for all
open subschemes $U\sub X$ is \emph{not} sufficient for a quasi-coherent
sheaf $\F$ on $Y$ to be very flat over~$X$.
 This is clear from Example~\ref{very-covering-mod-counterex}.

\medskip

 A particular case of the following theorem was stated as a conjecture
in an early February~2014 version of this book
manuscript~\cite[Conjecture~1.7.2]{Pcosh}.
 The proof was given in the paper~\cite{PSl1}.

\begin{thm} \label{very-flat-theorem}
\textup{(a)} Let $R$ be a commutative ring and $S$ be a finitely 
presented commutative $R$\+algebra.
 Assume that $S$ is a flat $R$\+module.
 Then $S$ is a very flat $R$\+module. \par
\textup{(b)} Let $R$ be a commutative ring, $S$ be a finitely presented
commutative $R$\+algebra, and $F$ be a finitely presented $S$\+module.
 Assume that $F$ is a flat $R$\+module.
 Then $F$ is a very flat $R$\+module.
\end{thm}

\begin{proof}
 Part~(a) is~\cite[Main Theorem~1.1]{PSl1}.
 Part~(b) is~\cite[Main Theorem~1.2]{PSl1}.
\end{proof}

 In other words, Theorem~\ref{very-flat-theorem}(a) tells us that any
flat morphism of finite presentation between schemes is very flat.
 The following corollary explains the situation.

\begin{cor} \label{very-flat-algebra-cor}
 Let $R\rarrow S$ be a homomorphism of commutative rings making $S$
a finitely presented, flat $R$\+algebra.
 Then $S$ is a very flat $R$\+algebra.
\end{cor}

\begin{proof}
 This is~\cite[Corollary~9.1]{PSl1}.
 The point is that, for any element $s\in S$, the $R$\+algebra
$S[s^{-1}]$ is also flat and finitely presented.
 Hence the assertion follows from Theorem~\ref{very-flat-theorem}(a).
\end{proof}

 For any scheme $X$, we denote by $\boA^n_X$ the $n$\+dimensional
relative affine space over~$X$; so if $X=\Spec R$ then
$\boA^n_X=\Spec R[x_1,\dotsc,x_n]=\boA^n_R$.
 Similarly, denote by $\boA^\infty_X$ the infinite-dimensional
relative affine space (of countable relative dimension) over~$X$.
 So if $X=\Spec R$, then $\boA^\infty_X=\Spec R[x_1,x_2,x_3,\dotsc]$.

\begin{cor}  \label{relative-affine-space-very-flat}
 For any scheme $X$, the natural projection $\mathbb A^\infty_X
\rarrow X$ is a very flat morphism.
\end{cor}

\begin{proof}
 One can assume $X$ to be affine.
 Then any principal affine open subscheme in $\boA_X^\infty$ is
the complement to the subscheme of zeroes of an equation depending
on a finite number of variables~$x_i$ only.
 Thus it suffices to show that $\boA_X^n\rarrow X$ is a very flat
morphism for every integer~$n\ge1$, which is a particular case of
Corollary~\ref{very-flat-algebra-cor}.
\end{proof}

\begin{thm}
 Let $X$ and $Y$ be two schemes over a field~$k$.
 Then the natural projection morphism $X\times_k Y\rarrow X$ is very
flat.
\end{thm}

\begin{proof}
 This is~\cite[Corollary~9.8]{PSl1}; see property~(VF16)
in~\cite[Section~0.7]{PSl1} for a discussion.
\end{proof}

 The following lemma, claiming that the very flatness property is local
with respect a certain special class of very flat coverings, is to be
compared with Lemma~\ref{very-flat-injective-base-change}.

\begin{lem}  \label{very-flat-local-for-very-flat}
 Let $g\:x\rarrow X$ be a very flat affine morphism of schemes such that
for any (small enough) affine open subscheme $U\sub X$ the ring\/
$\O(g^{-1}(U))$ is a projective module over\/ $\O(U)$ and a projective
generator of the abelian category of\/ $\O(U)$\+modules.
 Then a morphism of schemes $f\:Y\rarrow X$ is very flat whenever
the morphism $f'\:y=x\times_XY\rarrow x$ is very flat.
\end{lem}

\begin{proof}
 Let $V\sub Y$ and $U\sub X$ be affine open subschemes such that
$f(V)\sub U$.
 Set $u=g^{-1}(U)$ and assume that the $\O(u)$\+module $\O(u)
\ot_{\O(U)}\O(V)$ is very flat.
 Then, the morphism~$g$ being very flat, the ring $\O(u)\ot_{\O(U)}\O(V)$
is very flat as an $\O(U)$\+module, too.
 Finally, since $\O(u)$ is a projective generator of the category
of $\O(U)$\+modules, we can conclude that the ring $\O(V)$ is
a very flat $\O(U)$\+module.
\end{proof}

\begin{rem}
 Notice that any surjective finite flat morphism of Noetherian schemes
$g\:x\rarrow X$ satisfies the assumptions of
Lemma~\ref{very-flat-local-for-very-flat}.
 Indeed, such a morphism is very flat by
Theorem~\ref{very-flat-theorem}(a).
 Any finite morphism of schemes is affine; and it remains to check that,
for any surjective finite flat morphism of affine schemes
$\Spec T\rarrow\Spec R$, the $R$\+module $T$ is a projective generator
of $R\modl$.
 Any finitely presented flat module is projective, so we only need to
show that there are no nonzero $R$\+modules $M$ for which every
morphism of $R$\+modules $T\rarrow M$ vanishes.
 Since the functor Hom from a finitely presented module over
a commutative ring commutes with localizations, and finitely generated
projective modules are locally free in the Zariski topology, the desired
property follows from the surjectivity assumption on the morphism
$\Spec T\rarrow\Spec R$.
 The point is that any nonzero free module is a generator of
the module category (e.~g., over a local ring).
\end{rem}

 It is well-known that any flat morphism of finite type between
Noetherian schemes (or of finite presentation between arbitrary schemes)
is an open map~\cite[Th\'eor\`eme~2.4.6]{Groth3}.
 We will see below that the similar result about very flat morphisms
does not require any finiteness conditions at all.

 Let $R$ be a commutative ring and $M$ be an $R$\+module.
 Following the terminology suggested in~\cite[Definition~2.7]{ST} and
taken up in~\cite[Section~0.3]{PSl1}, we define
the \emph{P\+support} $\PSupp M=\PSupp_RM\sub\Spec R$ of an $R$\+module
$M$ as the set of all prime ideals $\p\sub R$ for which the tensor
product $k(\p)\ot_RM$ does not vanish, where $k(\p)$ denotes
the residue field of the ideal~$\p$.

 Notice that, for any commutative ring homomophism $f\:R\rarrow S$,
the image of the induced map of the spectra $\Spec S\rarrow\Spec R$
coincides with the P\+support of the $R$\+module~$S$ (because
a commutative ring is nonzero if and only if its spectrum is nonempty).

\begin{lem}  \label{veryflat-dense-support}
 Let $R$ be a commutative ring without nilpotent elements and $F$ be
a nonzero very flat $R$\+module.
 Then the P\+support of $F$ contains a nonempty open subset in\/
$\Spec R$.
\end{lem}

\begin{proof}
 By Corollary~\ref{very-flat-transfinite}, any very flat $R$\+module $F$
is a direct summand of a transfinitely iterated extension $M =
\varinjlim_\alpha M_\alpha$ of certain $R$\+modules of the form
$R[s_\alpha^{-1}]$, where $s_\alpha\in R$.
 In particular, $F$ is an $R$\+submodule in $M$; consider the minimal
index~$\alpha$ for which the intersection $F\cap M_\alpha\sub M$
is nonzero.
 Denote this intersection by $G$ and set $s=s_\alpha$; then $G$ is
a nonzero $R$\+submodule in $F$ and in $R[s^{-1}]$ simultaneously.

 Hence the localization $G[s^{-1}]$ is a nonzero ideal in
the ring $R[s^{-1}]$.
 The P\+support $\PSupp_{R[s^{-1}]}G[s^{-1}]\sub \Spec R[s^{-1}]$ of
the $R[s^{-1}]$\+module $G[s^{-1}]$ is equal to the intersection of
$\PSupp_RG\sub\Spec R$ with $\Spec R[s^{-1}]\sub\Spec R$.
 By right exactness of the tensor product functor
$k(\p)\ot_{R[s^{-1}]}{-}$ for any prime ideal $\p\sub R[s^{-1}]$,
applied to the short exact sequence of $R[s^{-1}]$\+modules
$0\rarrow G[s^{-1}]\rarrow R[s^{-1}]\rarrow R[s^{-1}]/G[s^{-1}]
\rarrow0$, one has $k(\p)\ot_{R[s^{-1}]}G[s^{-1}]\ne0$ whenever
$k(\p)\ot_{R[s^{-1}]}R[s^{-1}]/G[s^{-1}]=0$.
 Hence the subset $\PSupp_{R[s^{-1}]}G[s^{-1}]$ contains the complement
$V(s,G)\sub\Spec R[s^{-1}]$ to the P\+support of the quotient module
$R[s^{-1}]/G[s^{-1}]$ in $\Spec R[s^{-1}]$.
 The latter P\+support is a closed subset in $\Spec R[s^{-1}]$
corresponding to the ideal $G[s^{-1}]$.
 If there are no nilpotents in $R$, then there exists a prime ideal
in $R[s^{-1}]$ not containing $G[s^{-1}]$; so the open subset
$V(s,G)\sub\Spec R[s^{-1}]$ is nonempty.

 Let us show that the P\+support of the $R$\+module $F$ contains
$V(s,G)$.
 Let $\p\in V(s,G)\sub\Spec R$ be a prime ideal in $R$ and $k(\p)$ be
its residue field.
 Then the map $k(\p)\ot_R G= k(\p)\ot_RG[s^{-1}]\rarrow
k(\p)\ot_R R[s^{-1}]=k(\p)$ is surjective by the argument above,
the map $k(\p)\ot_RM_\alpha\rarrow k(\p)\ot_RM$ is injective because
the $R$\+module $M/M_\alpha$ is flat, and the map $k(\p)\ot_R F
\rarrow k(\p)\ot_R M$ is injective since $F$ is a direct summand in~$M$
(we do not seem to really use the latter observation).
 Finally, both the composition $G\rarrow F\rarrow M$ and the map
$G\rarrow R[s^{-1}]$ factorize through the same $R$\+module morphism
$G\rarrow M_\alpha$.

 Now it follows from the commutativity of the triangle
$G\rarrow M_\alpha\rarrow R[s^{-1}]$ that the map $k(\p)\ot_R G\rarrow
k(\p)\ot_R M_\alpha$ is nonzero, it is clear from the diagram
$G\rarrow M_\alpha\rarrow M$ that the composition $k(\p)\ot_R G\rarrow
k(\p)\ot_R M$ is nonzero, and it follows from commutativity of
the diagram $G\rarrow F\rarrow M$ that the map $k(\p)\ot_RG\rarrow
k(\p)\ot_RF$ is nonzero.
 The assertion of Lemma is proved.
\end{proof}

\begin{prop}  \label{veryflat-open-intersection}
 Let $F$ be a very flat module over a commutative ring $R$ and
$Z\sub\Spec R$ be a closed subset.
 Suppose that the intersection\/ $\PSupp F\cap Z\sub\Spec R$ is nonempty.
 Then\/ $\PSupp F$ contains a nonempty open subset of~$Z$.
\end{prop}

\begin{proof}
 Endow $Z$ with the structure of a reduced closed subscheme in $\Spec R$
and set $S=\O(Z)$.
 Then the intersection $Z\cap\PSupp F$ coincides with the P\+support
$\PSupp_S(S\ot_RF)\sub\Spec S$ of the $S$\+module $S\ot_R F$ in
$\Spec S = Z\sub\Spec R$.
 By Lemma~\ref{very-scalars-always}(b), the $S$\+module $S\ot_R F$ is
very flat.
 Now if this $S$\+module vanishes, then the intersection $Z\cap\PSupp F$
is empty, while otherwise it contains a nonempty open subset of
$\Spec S$ by Lemma~\ref{veryflat-dense-support}.
\end{proof}

\begin{rem}  \label{nilradical-reduction-remark}
 It is clear from the argument above that one can replace
the condition that the ring $R$ has no nilpotent elements in
Lemma~\ref{veryflat-dense-support} by the condition of nonvanishing
of the tensor product $S\ot_RF$ of the $R$\+module $F$ with
the quotient ring $S=R/J$ of the ring $R$ by its nilradical~$J$.
 In particular, this condition holds automatically if
(the $R$\+module $F$ does not vanish and) the nilradical $J\sub R$
is a nilpotent ideal, i.~e., there exists an integer $N\ge1$
such that $J^N=0$.
 This includes all Noetherian commutative rings~$R$.

 On the other hand, it is not difficult to demonstrate an example
of a flat module over a commutative ring that is annihilated by
the reduction modulo the nilradical.
 E.~g., take $R$ to be the ring of polynomials in $x$, $x^{1/2}$,
$x^{1/4}$,~\dots, $x^{1/2^N}$,~\dots\ over a field $k=S$ with
the imposed relation $x^r=0$ for $r>1$, and $F=J$ to be
the nilradical (i.~e., the kernel of the augmentation morphism
to~$k$) in~$R$.
 The $R$\+module $F$ is flat as the inductive limit of free
$R$\+modules $R\ot_k kx^{1/2^N}$ for $N\to\infty$, and one
clearly has $S\ot_RF = R/J\ot_R J = J/J^2=0$.
\end{rem}

\begin{thm} \label{veryflat-open-support-thm}
 The P\+support of any very flat module over a commutative ring $R$ is
an open subset in\/ $\Spec R$.
\end{thm}

\begin{proof}
 Let $F$ be a very flat $R$\+module.
 Denote by $Z$ the closure of the complement $\Spec R\setminus\PSupp F$
in $\Spec R$.
 Then, by the definition, $\PSupp F$ contains the complement $\Spec R
\setminus Z$, and no nonempty open subset in (the induced topology of)
$Z$ is contained in $\PSupp F$.
 By Proposition~\ref{veryflat-open-intersection}, it follows that
$\PSupp F$ does not intersect $Z$, i.~e., $\PSupp F=\Spec R\setminus Z$
is an open subset in $\Spec R$.
\end{proof}

\begin{ex}  \label{q-over-z-not-very-flat-ex}
 In particular, we have shown that the $\boZ$\+module $\boQ$ is
\emph{not} very flat (by Theorem~\ref{veryflat-open-support-thm}),
even though it is a flat module of projective dimension~$1$.
\end{ex}

\begin{ex} \label{baer-specker-not-very-flat-ex}
 Another example of a flat (i.~e., torsion-free) abelian group that is
\emph{not} very flat is the countable product of copies of the group
of integers, $\prod_{n=0}^\infty\boZ$.
 This is explained in the paper~\cite[Example~3.2]{ST}.
 Notice that, for any right coherent associative ring $R$, the class
of flat left $R$\+modules is closed under infinite products in
$R\modl$.
 This example from~\cite{ST} shows that the class of very flat
$\boZ$\+modules is \emph{not} closed under infinite products.
 (See Example~\ref{X-ctrh-prj-not-closed-under-products-ex} below
for a further discussion.)
\end{ex}

\begin{cor}
 Any very flat morphism of schemes is an open map of their underlying
topological spaces.
\end{cor}

\begin{proof}
 Given a very flat morphism $f\:Y\rarrow X$ and an open subset $W\sub Y$,
cover $W$ with open affines $V\sub Y$ for which there exist open affines
$U\sub X$ such that $f(V)\sub U$, and apply
Theorem~\ref{veryflat-open-support-thm} to the very flat
$\O(U)$\+modules $\O(V)$.
\end{proof}

\Section{Contraherent Cosheaves over a Scheme} \label{contraherent-sect}

\subsection{Cosheaves of modules over a sheaf of rings}
\label{cosheaf-mod-subsect}
 Let $X$ be a topological space.
 A \emph{copresheaf of abelian groups} on $X$ is a covariant functor
from the category of open subsets of $X$ (with the identity embeddings
as morphisms) to the category of abelian groups.

 Given a copresheaf of abelian groups $\P$ on $X$, we will denote
the abelian group it assigns to an open subset $U\sub X$ by $\P[U]$
and call it the group of \emph{cosections} of $\P$ over~$U$.
 For a pair of embedded open subsets $V\sub U\sub X$, the map
$\P[V]\rarrow\P[U]$ that the copresheaf $\P$ assigns to $V\sub U$
will be called the \emph{corestriction} map.

 A copresheaf of abelian groups $\P$ on $X$ is called a \emph{cosheaf}
if for any open subset $U\sub X$ and its open covering $U=\bigcup_\alpha
U_\alpha$ the following sequence of abelian groups is (right) exact:
\begin{equation} \label{cosheaf-definition} \textstyle
 \bigoplus_{\alpha,\beta}\P[U_\alpha\cap U_\beta] \lrarrow
 \bigoplus_\alpha \P[U_\alpha] \lrarrow \P[U] \lrarrow 0.
\end{equation}

 Let $\O$ be a sheaf of associative rings on~$X$.
 A copresheaf of abelian groups $\P$ on $X$ is said to be
a \emph{copresheaf of \textup{(}left\textup{)} $\O$\+modules} if
for each open subset $U\sub X$ the abelian group $\P[U]$ is endowed with
the structure of a (left) module over the ring $\O(U)$, and for each pair
of embedded open subsets $V\sub U\sub X$ the map of corestriction of
cosections $\P[V]\rarrow\P[U]$ in the copresheaf $\P$ is
a homomorphism of $\O(U)$\+modules.
 Here the $\O(U)$\+module structure on $\P[V]$ is obtained from
the $\O(V)$\+module structure by the restriction of scalars via
the ring homomorphism $\O(U)\rarrow\O(V)$.

 A copresheaf of $\O$\+modules on $X$ is called a \emph{cosheaf of\/
$\O$\+modules} if its underlying copresheaf of abelian groups is
a cosheaf of abelian groups.

\begin{rem}
 One can define copresheaves with values in any category, and cosheaves
with values in any category that has coproducts.
 In particular, one can speak of cosheaves of sets, etc.
 Notice, however, that unlike for (pre)sheaves, the underlying
copresheaf of sets of a cosheaf of abelian groups is \emph{not}
a cosheaf of sets in general, as the forgetful functor from
the abelian groups to sets preserves products, but not coproducts.
 Thus cosheaves of sets (as developed, e.~g., in~\cite{BF}) and
cosheaves of abelian groups or modules are two quite distinct
theories.
\end{rem}

 Let $\bB$ be a base of open subsets of $X$.
 We will consider covariant functors from $\bB$ (viewed as a full
subcategory of the category of open subsets in~$X$) to the category
of abelian groups.
 We say that such a functor $\Q$ is \emph{endowed with an $\O$\+module
structure} if the abelian group $\Q[U]$ is endowed with
an $\O(U)$\+module structure for each $U\in\bB$ and the above
compatibility condition holds for the corestriction maps
$\Q[V]\rarrow\Q[U]$ assigned by the functor $\Q$ to all $V$, $U\in\bB$
such that $V\sub U$.

 The following result is essentially contained
in~\cite[Section~0.3.2]{Groth}, as is its (more familiar) sheaf
version, to which we will turn in due order.

\begin{thm} \label{cosheaf-base-thm}
 A covariant functor\/ $\Q$ with an\/ $\O$\+module structure on a base\/
$\bB$ of open subsets of $X$ can be extended to a cosheaf of\/
$\O$\+modules\/ $\P$ on $X$ if and only if the following condition holds.
 For any open subset $V\in\bB$, any its covering
$V=\bigcup_\alpha V_\alpha$ by open subsets $V_\alpha\in\bB$, and any
(or, equivalently, some particular) covering $V_\alpha\cap V_\beta =
\bigcup_\gamma W_{\alpha\beta\gamma}$ of the intersections
$V_\alpha\cap V_\beta$ by open subsets $W_{\alpha\beta\gamma}\in\bB$,
the sequence of abelian groups (or\/ $\O(V)$\+modules)
\begin{equation} \label{base-cosheaf} \textstyle
\bigoplus_{\alpha,\beta,\gamma} \Q[W_{\alpha\beta\gamma}]\lrarrow
\bigoplus_\alpha \Q[V_\alpha] \lrarrow \Q[V]\lrarrow 0
\end{equation}
must be exact.
 The functor of restriction of cosheaves of\/ $\O$\+modules to a base\/
$\bB$ is an equivalence between the category of cosheaves of\/
$\O$\+modules on $X$ and the category of covariant functors on\/ $\bB$,
endowed with\/ $\O$\+module structures and
satisfying~\textup{\eqref{base-cosheaf}}.
\end{thm}

\begin{proof}
 The elementary approach taken in the exposition below is to pick
an appropriate stage at which one can dualize and pass to (pre)sheaves,
where our intuitions work better.
 First we notice that if the functor $\Q$ (with its $\O$\+module
structure) has been extended to a cosheaf of $\O$\+modules $\P$ on $X$,
then for any open subset $U\sub X$ there is an exact sequence of
$\O(U)$\+modules
\begin{equation} \label{cosheaf-recover} \textstyle
 \bigoplus_{W,V',V''} \Q[W] \lrarrow \bigoplus_V \Q[V]\lrarrow
 \P[U]\lrarrow 0,
\end{equation}
where the summation in the middle term runs over all open subsets
$V\in\bB$, \ $V\sub U$, while the summation in the leftmost term is
done over all triples of open subsets $W$, $V'$, $V''\in\bB$, \ 
$W\sub V'\sub U$, \ $W\sub V''\sub U$.
 Conversely, given a functor $\Q$ with an $\O$\+module structure one can
recover the $\O(U)$\+module $\P[U]$ as the cokernel of the leftmost
arrow in~\eqref{cosheaf-recover}.

 Clearly, the modules $\P[U]$ constructed in this way naturally form
a copresheaf of $\O$\+modules on $X$.
 Before proving that it is a cosheaf, one needs to show that for any
open covering $U=\bigcup_\alpha V_\alpha$ of an open subset $U\sub X$
by open subsets $V_\alpha\in\bB$ and any open coverings $V_\alpha\cap
V_\beta = \bigcup_\gamma  W_{\alpha\beta\gamma}$ of the intersections
$V_\alpha\cap V_\beta$ by open subsets $W_{\alpha\beta\gamma}\in\bB$
the natural map from the cokernel of the morphism
\begin{equation} \label{cosheaf-alt-recover} \textstyle
 \bigoplus_{\alpha,\beta,\gamma} \Q[W_{\alpha\beta\gamma}]\lrarrow
 \bigoplus_\alpha \Q[V_\alpha]
\end{equation}
to the (above-defined) $\O(U)$\+module $\P[U]$ is an isomorphism.
 In particular, it will follow that $\P[V]=\Q[V]$ for $V\in\bB$.

 Notice that it suffices to check both assertions for co(pre)sheaves
of abelian groups (though it will not matter in the subsequent argument).
 Notice also that a copresheaf of $\O$\+modules $\P$ is a cosheaf if and
only if the dual presheaf of $\O$\+modules $U\mpsto\Hom_\boZ(\P[U],I)$
is a sheaf on $X$ for every abelian group $I$ (or specifically for
$I=\boQ/\boZ$).
 Similarly, the condition~\eqref{base-cosheaf} holds for a covariant
functor $\Q$ on a base $\bB$ if and only if the dual
condition~\eqref{base-sheaf} below holds for the contravariant
functor $V\mpsto\Hom_\boZ(\Q[V],I)$ on~$\bB$.
 So it remains to prove the following Proposition~\ref{sheaf-base-prop}.
\end{proof}

 Now we will consider contravariant functors $\G$ from $\bB$ to
the category of abelian groups, and say that such a functor is endowed
with an $\O$\+module structure if the abelian group $\G(U)$ is
an $\O(U)$\+module for every $U\in\bB$ and the restriction morphisms
$\G(U)\rarrow\G(V)$ are morphisms of $\O(U)$\+modules for all
$V$,~$U\in\bB$ such that $V\sub U$.

\begin{prop} \label{sheaf-base-prop}
 A contravariant functor\/ $\G$ with an\/ $\O$\+module structure on
a base\/ $\bB$ of open subsets of $X$ can be extended to a sheaf of\/
$\O$\+modules\/ $\F$ on $X$ if and only if the following condition holds.
 For any open subset $V\in\bB$, any its covering
$V=\bigcup_\alpha V_\alpha$ by open subsets $V_\alpha\in\bB$, and any
(or, equivalently, some particular) covering $V_\alpha\cap V_\beta =
\bigcup_\gamma W_{\alpha\beta\gamma}$ of the intersections
$V_\alpha\cap V_\beta$ by open subsets $W_{\alpha\beta\gamma}\in\bB$,
the sequence of abelian groups (or\/ $\O(V)$\+modules)
\begin{equation} \label{base-sheaf} \textstyle
 0\lrarrow\G(V)\lrarrow\prod_\alpha \G(V_\alpha)\lrarrow
 \prod_{\alpha,\beta,\gamma} \G(W_{\alpha\beta\gamma})
\end{equation}
must be exact.
 The functor of restriction of sheaves of\/ $\O$\+modules to a base\/
$\bB$ is an equivalence between the category of sheaves of\/
$\O$\+modules on $X$ and the category of contravariant functors on\/
$\bB$, endowed with\/ $\O$\+module structures and
satisfying~\textup{\eqref{base-sheaf}}.
\end{prop}

\begin{proof}[Sketch of proof]
 As above, we notice that if the functor $\G$ (with its $\O$\+module
structure) has been extended to a sheaf of $\O$\+modules $\F$ on $X$,
then for any open subset $U\sub X$ there is an exact sequence of
$\O(U)$\+modules
$$ \textstyle \label{sheaf-recover} 
 0\lrarrow\F(U) \lrarrow \prod_V \G(V) \lrarrow
 \prod_{W,V',V''}\G(W),
$$
the summation rules being as in~\eqref{cosheaf-recover}.
 Conversely, given a functor $\G$ with an $\O$\+module structure 
one can recover the $\O(U)$\+module $\F(U)$ as the kernel of
the rightmost arrow.

 The rest is a conventional argument with (pre)sheaves and coverings.
 Recall that a presheaf $\F$ on $X$ is called \emph{separated} if
the map $\F(U)\rarrow\prod_\alpha \F(U_\alpha)$ is injective for any
open covering $U=\bigcup_\alpha U_\alpha$ of an open subset $U\sub X$.
 Similarly, a contravariant functor $\G$ on a base $\bB$ is said to be
separated if its sequences~\eqref{base-sheaf} are exact at the leftmost
nontrivial term.

 For any open covering $U=\bigcup_\alpha V_\alpha$ of an open subset
$U\sub X$ by open subsets $V_\alpha\in\bB$ and any open coverings
$V_\alpha\cap V_\beta = \bigcup_\gamma W_{\alpha\beta\gamma}$ of
the intersections $V_\alpha\cap V_\beta$ by open subsets
$W_{\alpha\beta\gamma}\in\bB$ there is a natural map from
the (above-defined) $\O(U)$\+module $\F(U)$ to the kernel of
the morphism
\begin{equation} \label{sheaf-kernel}
 \textstyle
 \prod_\alpha \G(V_\alpha)\lrarrow\prod_{\alpha,\beta,\gamma}
 \G(W_{\alpha\beta\gamma}).
\end{equation}
 Let us show that this map is an isomorphism provided that $\G$
satisfies~\eqref{base-sheaf}.
 In particular, it will follow that $\F(V)=\G(V)$ for $V\in\bB$.

 We start with the observation that, when $\G$ is separated, the kernel 
of~\eqref{sheaf-kernel} does not depend on the choice of the open
subsets $W_{\alpha\beta\gamma}$.
 Indeed, let $\phi_\alpha\in\G(V_\alpha)$ and $\phi_\beta\in\G(V_\beta)$
be two sections such that $\phi_\alpha|_{W_{\alpha\beta\gamma}}=
\phi_\beta|_{W_{\alpha\beta\gamma}}$ for all~$\gamma$.
 Let $W\in\bB$ be an open subset such that $W\sub U_\alpha\cap U_\beta$.
 Choose open subsets $T_{\gamma\delta}\in\bB$ such that
$W\cap W_{\alpha\beta\gamma}=\bigcup_\delta T_{\gamma\delta}$
for every~$\gamma$; so $W=\bigcup_{\gamma,\delta}T_{\gamma\delta}$.
 Then one has $\phi_\alpha|_{T_{\gamma\delta}}=
\phi_\beta|_{T_{\gamma\delta}}$ for all~$\gamma$,~$\delta$; and
injectivity of the map $\G(W)\rarrow\prod_{\gamma,\delta}
\G(T_{\gamma\delta})$ implies the desired equation
$\phi_\alpha|_W=\phi_\beta|_W$.
 Thus we can assume that the collection $\{W_{\alpha\beta\gamma}\}$
for fixed $\alpha$ and~$\beta$ consists of all open subsets $W\in\bB$
such that $W\sub V_\alpha\cap V_\beta$.

 Furthermore, it is easy to see that the map from $\F(U)$ to the kernel
of~\eqref{sheaf-kernel} is injective whenever $\G$ is separated.
 Indeed, let $\phi_V\in\G(V)$ be a collection of sections, defined for
all $V\in\bB$, \ $V\sub U$, such that $\phi_{V'}|_W=\phi_{V''}|_W$ for
all $W\in\bB$, \ $W\sub V'\cap V''$, and $\phi_{V_\alpha}=0$
for all~$\alpha$.
 Then, for every fixed $V$, one has $\phi_V|_W=\phi_{V_\alpha}|_W=0$ for
all~$\alpha$ and all $W\in\bB$, \ $W\sub V\cap V_\alpha$.
 Clearly, $V$ is the union of all such open subsets $W$ (where
$\alpha$~varies), so it follows that $\phi_V=0$.

 To check surjectivity, suppose that we are given a collection of 
sections $\phi_\alpha\in\G(V_\alpha)$ representing an element of
the kernel.
 Fix an open subset $V\in\bB$, \ $V\sub U$, and consider its covering
by all the open subsets $W\in\bB$ such that $W\sub V\cap V_\alpha$ for 
some~$\alpha$.
 Set $\psi_W=\phi_\alpha|_W\in\G(W)$ for every such $W$; by assumption,
if $W\sub V_\alpha\cap V_\beta$, then $\phi_\alpha|_W =\phi_\beta|_W$,
so the element~$\psi_W$ is well-defined.
 Applying~\eqref{base-sheaf}, we conclude that there exists
a unique element $\phi_V\in\G(V)$ such that $\phi_V|_W=\psi_W$
for any $W\sub V\cap V_\alpha$.
 The collection of sections~$\phi_V$ represents an element of $\F(U)$
that is a preimage of our original element of the kernel
of~\eqref{sheaf-kernel}.

 Now let us show that $\F$ is a sheaf.
 Let $U=\bigcup_\alpha U_\alpha$ be an open covering of an open subset
$U\sub X$.
 First let us see that $\F$ is separated provided that $\G$~is.
 Let $s\in\F(U)$ be a section whose restriction to all the open
subsets $U_\alpha$ vanishes.
 The element~$s$ is represented by a collection of sections $\phi_V\in
\G(V)$ defined for all open subsets $V\sub U$, \ $V\in\bB$.
 The condition $s|_{U_\alpha}=0$ means that $\phi_W=0$ whenever
$W\sub U_\alpha$, \ $W\in\bB$.
 To check that $\phi_V=0$ for all $V$, we notice that open subsets
$W\sub V$, \ $W\in\bB$ for which there exists~$\alpha$ such that
$W\sub U_\alpha$ form an open covering of~$V$.

 Finally, let $s_\alpha\in \F(U_\alpha)$ be a collection of sections
such that $s_\alpha|_{U_\alpha\cap U_\beta} = s_\beta|_{U_\alpha\cap
U_\beta}$ for all $\alpha$ and~$\beta$.
 Every element~$s_\alpha$ is represented by a collection of sections
$\phi_V\in\G(V)$ defined for all open subsets $V\sub U_\alpha$, \
$V\in\bB$.
 Clearly, the element $\phi_V$ does not depend on the choice of
a particular~$\alpha$ for which $V\sub U_\alpha$, so our notation is
consistent.
 All the open subsets $V\sub U$, \ $V\in\bB$ for which there exists
some $\alpha$ such that $V\sub U_\alpha$ form an open covering of
the open subset $U\sub X$.
 The collection of sections~$\phi_V$ represents an element of the kernel
of the morphism~\eqref{sheaf-kernel} for this covering, hence it
corresponds to an element of $\F(U)$.
\end{proof}

\begin{rem}  \label{scheme-topology}
 Let $X$ be a topological space with a topology base $\bB$ consisting
of quasi-compact open subsets (in the induced topology) for which
the intersection of any two open subsets from $\bB$ that are
contained in a third open subset from $\bB$ is quasi-compact as well.
 E.~g., any scheme $X$ with the base of all affine (or quasi-compact
quasi-separated) open subschemes has these properties.
 Then it suffices to check both the conditions~\eqref{base-cosheaf}
and~\eqref{base-sheaf} for \emph{finite} coverings $V_\alpha$ and
$W_{\alpha\beta\gamma}$ only.

 Indeed, let us explain the sheaf case.
 Obviously, injectivity of the left arrow in~\eqref{base-sheaf}
for any given covering $V=\bigcup_\alpha V_\alpha$ follows from such
injectivity for a subcovering $V=\bigcup_i V_i$, \ 
$\{V_i\}\sub\{V_\alpha\}$.
 Assuming $\G$ is separated, one checks that exactness of
the sequence~\eqref{base-sheaf} for any given covering follows
from the same exactness for a subcovering.

 It follows that for a topological space $X$ with a fixed topology
base $\bB$ satisfying the above condition there is another duality
construction relating sheaves to cosheaves in addition to the one
we used in the proof of Theorem~\ref{cosheaf-base-thm}.
 Given a sheaf of $\O$\+modules $\F$ on $X$, one restricts it to
the base $\bB$, obtaining a contravariant functor $\G$ with
an $\O$\+module structure, defines the dual covariant functor $\Q$
with an $\O$\+module structure on $\bB$ by the rule $\Q[V] =
\Hom_\boZ(\G(V),I)$, where $I$ is an injective abelian group, and
extends the functor $\Q$ to a cosheaf of $\O$\+modules $\P$ on~$X$.

 It is this second duality functor, rather than the one from
the proof of Theorem~\ref{cosheaf-base-thm}, that will play a role
in the sequel (Section~\ref{coflasque} being a rare exception).
\end{rem}

\subsection{Exact category of contraherent cosheaves}
\label{contraherent-definition}
 Let $X$ be a scheme and $\O=\O_X$ be its structure sheaf.
 A cosheaf of $\O_X$\+modules $\P$ is called \emph{contraherent}
if for any pair of embedded affine open subschemes $V\sub U\sub X$
\begin{enumerate}
\renewcommand{\theenumi}{\roman{enumi}}
\item the morphism of $\O_X(V)$\+modules $\P[V]\rarrow
\Hom_{\O_X(U)}(\O_X(V),\P[U])$ induced by the corestriction morphism
$\P[V]\rarrow\P[U]$ is an isomorphism; and
\item one has $\Ext_{\O_X(U)}^{>0}(\O_X(V),\P[U])=0$.
\end{enumerate}
 It follows from Lemma~\ref{very-open-embedding} that
the $\O_X(U)$\+module $\O_X(V)$ has projective dimension at most~$1$,
so it suffices to require the vanishing of $\Ext^1$ in
the condition~(ii).
 We will call~(ii) the \emph{contraadjustedness condition}, and
(i) the \emph{contraherence condition}.

\begin{thm}  \label{contraherent-base}
 The restriction of cosheaves of\/ $\O_X$\+modules to the base of
all affine open subschemes of $X$ induces an equivalence between
the category of contraherent cosheaves on $X$ and the category
of covariant functors\/ $\Q$ with\/ $\O_X$\+module structures on
the category of affine open subschemes of $X$, satisfying
the conditions~\textup{(i\+ii)} for any pair of embedded affine
open subschemes $V\sub U\sub X$.
\end{thm}

\begin{proof}
 According to Theorem~\ref{cosheaf-base-thm}, a cosheaf of
$\O_X$\+modules is determined by its restriction to the base of
affine open subsets of~$X$.
 The contraadjustedness and contraherence conditions depend only on
this restriction.
 By Lemma~\ref{very-open-embedding}, given any affine scheme $U$,
a module $P$ over $\O(U)$ is contraadjusted if and only if
$\Ext^1_{\O(U)}(\O(V),P)=0$ for all affine open subschemes
$V\sub U$.
 Finally, the key observation is that the contraadjustedness and
contraherence conditions~(i\+ii) for a covariant functor with
an $\O_X$\+module structure on the category of affine open
subschemes of $X$ imply the cosheaf condition~\eqref{base-cosheaf}.
 This follows from Lemma~\ref{very-open-covering}(b)
and Remark~\ref{scheme-topology}.
\end{proof}

\begin{rem} \label{quasi-coherence-rem}
 Of course, one can similarly define quasi-coherent sheaves $\F$ on $X$
as sheaves of $\O_X$\+modules such that for any pair of embedded affine
open subschemes $V\sub U\sub X$ the morphism of $\O_X(V)$\+modules
$\O_X(V)\ot_{\O_X(U)}\F(U)\rarrow\F(V)$ induced by the restriction
morphism $\F(U)\rarrow\F(V)$ is an isomorphism.
 Since the $\O_X(U)$\+module $\O_X(V)$ is always flat, no version of
the condition~(ii) is needed in this context.
 The analogue of Theorem~\ref{contraherent-base} for quasi-coherent
sheaves is due to Enochs and Estrada~\cite[Section~2]{EE}; it can be
proved in the same way.
\end{rem}

 A short sequence of contraherent cosheaves $0\rarrow\P\rarrow\Q\rarrow
\R\rarrow0$ is said to be exact if the sequence of cosection modules
$0\rarrow\P[U]\rarrow\Q[U]\rarrow\R[U]\rarrow0$ is exact for every
affine open subscheme $U\sub X$.
 Notice that if $\{U_\alpha\}$ is an affine open covering of an affine
scheme $U$ and $\P\rarrow\Q\rarrow\R$ is a sequence of contraherent
cosheaves on $U$, then the sequence of $\O(U)$\+modules
$0\rarrow \P[U]\rarrow \Q[U]\rarrow \R[U]\rarrow 0$ is exact if
and only if all the sequences of $\O(U_\alpha)$\+modules $0\rarrow
\P[U_\alpha]\rarrow \Q[U_\alpha]\rarrow \R[U_\alpha]\rarrow0$ are.
 This follows from Lemma~\ref{very-exact-local}(a).

 We denote the exact category of contraherent cosheaves on a scheme $X$
by $X\ctrh$.
 By the definition, the functors of cosections over affine open
subschemes are exact on this exact category.
 It is also has exact functors of infinite product, which commute with
cosections over affine open subschemes (and in fact, over all
quasi-compact quasi-separated open subschemes as well).
 For a more detailed discussion of this exact category structure, we
refer the reader to Section~\ref{locally-contraherent}.

\begin{cor} \label{contraherent-on-affine}
 The functor assigning the\/ $\O(U)$\+module\/ $\P[U]$ to a contraherent
cosheaf\/ $\P$ on an affine scheme $U$ is an equivalence between
the exact category $U\ctrh$ of contraherent cosheaves on $U$ and
the exact category\/ $\O(U)\modl^\cta$ of contraadjusted modules
over the commutative ring\/~$\O(U)$.
\end{cor}

\begin{proof}
 Clear from the arguments above together with
Lemmas~\ref{very-tensor-hom}(b) and~\ref{very-open-embedding}.
\end{proof}

\begin{lem}  \label{hom-into-contraherent}
 Let\/ $\P$ be a cosheaf of\/ $\O_U$\+modules and\/ $\Q$ be
a contraherent cosheaf on an affine scheme~$U$.
 Then the group of morphisms of cosheaves of\/ $\O_U$\+modules\/
$\P\rarrow\Q$ is isomorphic to the group of morphisms of\/
$\O(U)$\+modules\/ $\P[U]\rarrow\Q[U]$.
\end{lem}

\begin{proof}
 Any morphism of cosheaves of $\O_U$\+modules $\P\rarrow\Q$ induces
a morphism of the $\O(U)$\+modules of global cosections.
 Conversely, given a morphism of $\O(U)$\+modules $\P[U]\rarrow\Q[U]$
and an affine open subscheme $V\sub U$, the composition
$\P[V]\rarrow\P[U]\rarrow\Q[U]$ is a morphism of $\O(U)$\+modules
from an $\O(V)$\+module $\P[V]$ to an $\O(U)$\+module $\Q[U]$.
 It induces, therefore, a morphism of $\O(V)$\+modules
$\P[V]\rarrow\Hom_{\O(U)}(\O(V),\Q[U])\simeq\Q[V]$.
 Now a morphism between the restrictions of two cosheaves of
$\O_U$\+modules to the base of affine open subschemes of $U$
extends uniquely to a morphism between the whole cosheaves.
\end{proof}

 A contraherent cosheaf $\P$ on a scheme $X$ is said to be
\emph{locally cotorsion} if for any affine open subscheme $U\sub X$
the $\O_X(U)$\+module $\P[U]$ is cotorsion.
 By Lemma~\ref{cotors-inj-covering}(a), the property of a contraherent
cosheaf on an affine scheme to be cotorsion is indeed local, so our
terminology is constent.

 A contraherent cosheaf $\gJ$ on a scheme $X$ is called
\emph{locally injective} if for any affine open subscheme $U\sub X$
the $\O_X(U)$\+module $\gJ[U]$ is injective.
 By Lemma~\ref{cotors-inj-covering}(b), local injectivity of
a contraherent cosheaf is indeed a local property.

 Just as above, one defines the exact categories $X\ctrh^\lct$ and
$X\ctrh^\lin$ of locally cotorsion and locally injective contraherent
cosheaves on $X$.
 These are full subcategories closed under extensions, infinite
products, and cokernels of admissible monomorphisms in $X\ctrh$,
with the induced exact category structures.

 The exact category $U\ctrh^\lct$ of locally cotorsion contraherent
cosheaves on an affine scheme $U$ is equivalent to the exact category
$\O(U)\modl^\cot$ of cotorsion $\O(U)$\+modules.
 The exact category $U\ctrh^\lin$ of locally injective contraherent
cosheaves on $U$ is equivalent to the additive category
$\O(U)\modl^\inj$ of injective $\O(U)$\+modules endowed with
the trivial exact category structure.
 These assertions follow from Lemmas~\ref{cotors-coexten}(a)
and~\ref{cotors-restrict}(b).

\begin{rem}
 Notice that a morphism of contraherent cosheaves on $X$ is
an admissible monomorphism if and only if it acts injectively
on the cosection modules over all the affine open subschemes on $X$.
 At the same time, the property of a morphism of contraherent cosheaves
on $X$ to be an admissible monomorphism is \emph{not} local in $X$,
and \emph{neither} is the property of a cosheaf of $\O_X$\+modules
to be contraherent (see Section~\ref{counterex-subsect} below).
 The property of a morphism of contraherent cosheaves to be
an admissible epimorphism is local, though (see
Lemma~\ref{very-exact-local}(b)).
 All of the above applies to locally cotorsion and locally injective
contraherent cosheaves as well (see Lemmas~\ref{cotors-exact-local}(b)
and~\ref{injective-exact-local}(b)).

 Notice also that a morphism of locally injective or locally
cotorsion contraherent cosheaves that is an admissible epimorphism
in $X\ctrh$ may \emph{not} be an admissible epimorphism in
$X\ctrh^\lct$ or $X\ctrh^\lin$.
 On the other hand, if a morphism of locally injective or locally
cotorsion contraherent cosheaves is an admissible monomorphism in
$X\ctrh$, then it is also an admissible monomorphism in $X\ctrh^\lct$
or $X\ctrh^\lin$, as it is clear from the discussion above.
\end{rem}

\subsection{Direct and inverse images of contraherent cosheaves}
\label{contra-direct-inverse}
 Let $\O_X$ be a sheaf of associative rings on a topological space $X$
and $\O_Y$ be such a sheaf on a topological space~$Y$.
 Furthermore, let $f\:Y\rarrow X$ be a morphism of ringed spaces, i.~e.,
a continuous map $Y\rarrow X$ together with a morphism
$\O_X\rarrow f_*\O_Y$ of sheaves of rings over~$X$.
 Then for any cosheaf of $\O_Y$\+modules $\Q$ the rule
$(f_!\Q)[W] = \Q[f^{-1}(W)]$ for all open subsets $W\sub X$ defines
a cosheaf of $\O_X$\+modules~$f_!\Q$.

 Let $\O_X$ be a sheaf of associative rings on a topological space $X$
and $Y\sub X$ be an open subspace.
 Denote by $\O_Y=\O_X|_Y$ the restriction of the sheaf of rings $\O_X$
onto $Y$, and by $j\:Y\rarrow X$ the corresponding morphism
(open embedding) of ringed spaces. 
 Given a cosheaf of $\O_X$\+modules $\P$ on $X$, the restriction
$\P|_Y$ of $\P$ onto $Y$ is a cosheaf of $\O_Y$\+modules defined by
the rule $\P|_Y[V]=\P[V]$ for any open subset $V\sub Y$.
 One can easily see that the restriction functor $\P\mpsto\P|_Y$
is right adjoint to the direct image functor $\Q\mpsto j_!\Q$ between
the categories of cosheaves of $\O_Y$\+ and $\O_X$\+modules,
that is, the adjunction isomorphism
\begin{equation}  \label{open-cosheaf-direct-inverse-adjunction}
 \Hom^{\O_X}(j_!\Q,\P)\simeq\Hom^{\O_Y}(\Q,\P|_Y)
\end{equation}
holds for any cosheaf of $\O_X$\+modules $\P$ and cosheaf of
$\O_Y$\+modules~$\Q$, where $\Hom^{\O_X}$ and $\Hom^{\O_Y}$ denote
the abelian groups of morphisms in the categories of cosheaves of
modules over the sheaves of rings $\O_X$ and~$\O_Y$.
 Indeed, for any compatible system of $\O_X(W)$\+module maps
$\Q[Y\cap W]\rarrow\P[W]$, \ $W\sub X$, the map
$\Q[Y\cap W]\rarrow\P[W]$ is equal to the composition $\Q[Y\cap W]
\rarrow\P[Y\cap W]\rarrow\P[W]$; so such a system of maps is determined
by the $\O_Y(V)$\+module maps $\Q[V]\rarrow\P[V]$, \ $V\sub Y$.
 Since one has $(j_!\Q)|_Y\simeq\Q$, it follows, in particular,
that the functor~$j_!$ is fully faithful.

 Let $f\:Y\rarrow X$ be an affine morphism of schemes, and let $\Q$ be
a contraherent cosheaf on~$Y$.
 Then $f_!\Q$ is a contraherent cosheaf on~$X$.
 Indeed, for any affine open subscheme $U\sub X$ the $\O_X(U)$\+module
$(f_!\Q)[U] = \Q[(f^{-1}(U)]$ is contraadjusted according to
Lemma~\ref{very-scalars-always}(a) applied to the morphism of commutative
rings $\O_X(U)\rarrow\O_Y(f^{-1}(U))$.
 For any pair of embedded affine open subschemes $V\sub U\sub X$ we have
natural isomorphisms of $\O_X(U)$\+modules
\begin{multline*}
 (f_!\Q)[V] = \Q[f^{-1}(V)] \simeq
 \Hom_{\O_Y(f^{-1}(U))}(\O_Y(f^{-1}(V)),\Q[f^{-1}(U)]) \\
 \simeq \Hom_{\O_X(U)}(\O_X(V),\Q[f^{-1}(U)]) =
 \Hom_{\O_X(U)}(\O_X(V),(f_!\Q)[U]),
\end{multline*}
since $\O_Y(f^{-1}(V))\simeq\O_Y(f^{-1}(U))\ot_{\O_X(U)}\O_X(V)$.

\begin{rem} \label{noncontraherent-direct-image-remark}
 For any quasi-compact quasi-separated morphism of schemes $f\:Y\rarrow
X$, the functor of direct image of sheaves of $\O$\+modules takes
quasi-coherent sheaves on $Y$ to quasi-coherent sheaves on~$X$.
 This assertion holds, essentially, because the module of sections
$\G(T)$ of a quasi-coherent sheaf $\G$ on $Y$ over a quasi-compact
quasi-separated open subscheme $T\sub Y$ is computed as the kernel
of a map between finite direct sums of the modules of sections of $\G$
over affine open subschemes of~$T$.
 For any affine open subschemes $V\sub U\sub X$, the functor of tensor
product with $\O_X(V)$ over $\O_X(U)$ preserves kernels and finite
direct sums.

 The dual analogue of this property does \emph{not} hold for
contraherent cosheaves in general:
 Given a contraherent cosheaf $\P$ on $Y$, the cosheaf of
$\O_X$\+modules $f_!\P$ need not be contraherent (not even for
a morphism of Noetherian schemes~$f$; not even for
an open embedding $f\:Y\rarrow X$ of smooth schemes of finite
type over a field~$k$ with a smooth complement $X\setminus Y$).
 In fact, the contraadjustedness condition~(ii) from
Section~\ref{contraherent-definition} is preserved by the direct
images with respect to quasi-compact morphisms (by
Lemma~\ref{very-scalars-always}(a) and because the class of
contraadjusted modules is closed under quotients).
 However, the contraherence condition~(i) is not preserved by~$f_!$
when $f$ is not affine.
 The problem is that, even though the functor
$\Hom_{\O_X(U)}(\O_X(V),{-})$ preserves short exact sequences of
contraadjusted $\O_X(U)$\+modules, it does \emph{not} preserve
the cokernels of morphisms of contraadjusted $\O_X(U)$\+modules.

 See the next Example~\ref{noncontraherent-direct-image-example} for
a specific counterexample.
 Consequently, one needs to impose various adjustedness conditions
on a contraherent cosheaf $\P$, such as coflasqueness or antilocality,
in order to guarantee that $f_!\P$ is contraherent.
 See Corollaries~\ref{coflasque-direct} and~\ref{clp-direct} below.
 This is an instance of the general phenomenon of partially defined
functors between exact categories, as per
Section~\ref{introd-partially-defined} of the Introduction
and~\cite[Section~5.6]{Pphil}.
\end{rem}

\begin{ex} \label{noncontraherent-direct-image-example}
 Let $R$ be a commutative ring and $f$, $g\in R$ be two regular
elements (nonzero-divisors).
 Put $U=\Spec R$ and $T=\Spec R[f^{-1}]\cup\Spec R[g^{-1}]\sub U$,
and denote by $j\:T\rarrow U$ the open embedding morphism.
 Let $P$ be a contraadjusted (or even cotorsion) $R$\+module and
$\widecheck P$ be the related contraherent cosheaf on~$U$
(as per Corollary~\ref{contraherent-on-affine}).
 Put $\P=\widecheck P|_T$; so $\P$ is a contraherent cosheaf on~$T$.

 By assumption, the $R$\+modules $R[f^{-1}]$ and $R[g^{-1}]$ are
submodules in $R[f^{-1},g^{-1}]$.
 Denote by $E$ the $R$\+module $R[f^{-1},g^{-1}]/(R[f^{-1}]+R[g^{-1}])$.
 Then we have a four-term exact sequence of $R$\+modules
\begin{multline} \label{cosections-sequence}
 0\lrarrow\Hom_R(E,P)\lrarrow\Hom_R(R[f^{-1},g^{-1}],P) \\ \lrarrow
 \Hom_R(R[f^{-1}],P)\oplus\Hom_R(R[g^{-1}],P)\lrarrow(j_!\P)[U]
 \lrarrow0.
\end{multline}
 Put $M=\Hom_R(R[f^{-1},g^{-1}],P)$ and
$N=\Hom_R(R[f^{-1}],P)\oplus\Hom_R(R[g^{-1}],P)$, and denote by $L$
the image of the map $M\rarrow N$ in
the sequence~\eqref{cosections-sequence}.

 Let $h\in R$ be another element.
 Then we have left exact sequences of $R$\+modules
\begin{gather*}
 0\lrarrow\Hom_R(R[h^{-1}],K)\lrarrow\Hom_R(R[h^{-1}],M)
 \lrarrow\Hom_R(R[h^{-1}],L) \\
 0\lrarrow\Hom_R(R[h^{-1}],L)\lrarrow\Hom_R(R[h^{-1}],N)
 \lrarrow\Hom_R(R[h^{-1}],N/L),
\end{gather*}
where $K=\Hom_R(E,P)$.
 It follows that the functor $\Hom_R(R[h^{-1}],{-})$ does \emph{not}
preserve the cokernel of the map $M\rarrow N$ whenever the map
$\Hom_R(R[h^{-1}],M)\rarrow\Hom_R(R[h^{-1}],L)$ is not surjective.
 In other words, the natural map
\begin{multline*}
 \coker(\Hom_R(R[h^{-1}],M)\to\Hom_R(R[h^{-1}],N)) \\
 \lrarrow\Hom_R(R[h^{-1}]\;\coker(M\to N))
\end{multline*}
is \emph{not} an isomorphism (in fact, not injective) in this case.

 Notice that the $R$\+module $R[f^{-1},g^{-1}]$ is very flat, while
the $R$\+module $P$ is contraadjusted by assumption.
 By Lemma~\ref{very-tensor-hom}(b), it follows that the $R$\+module
$M=\Hom_R(R[f^{-1},g^{-1}],P)$ is contraadjusted.
 In particular, $\Ext_R^1(R[h^{-1}],M)=0$.
 Therefore, it is clear from the short exact sequence of $R$\+modules
$0\rarrow K\rarrow M\rarrow L\rarrow0$ that the map
$\Hom_R(R[h^{-1}],M)\rarrow\Hom_R(R[h^{-1}],L)$ is \emph{not} surjective
whenever $\Ext^1_R(R[h^{-1}],K)\ne0$.

 Assume that $\Ext^1_R(E,P)=0$.
 Then we compute
$$
 \Ext^1_R(R[h^{-1}],K)=
 \Ext^1_R(R[h^{-1}],\Hom_R(E,P))\simeq
 \Ext^1_R(R[h^{-1}]\ot_RE\;P)
$$
using formula~\eqref{ext-tor-adjunction}.
 It remains to come up with an example of a cotorsion $R$\+module $P$
for which $\Ext^1_R(E,P)=0$ but $\Ext^1_R(E[h^{-1}],P)\ne0$.
 Then the arguments above will show that the functor
$\Hom_R(R[h^{-1}],{-})$ does not preserve the cokernel of
the $R$\+module map $M\rarrow N$.
 This will mean that the cosheaf of $\O_U$\+modules $j_!\P$ does
\emph{not} satisfy the contraherence condition for the pair of affine
open subschemes $V=\Spec R[h^{-1}]\sub U=\Spec R$ in $U=\Spec R$.

 We argue similarly to Example~\ref{not-contrainjective-example},
using the terminology and results of
Sections~\ref{cotorsion-prelim-subsect}\+-%
\ref{small-object-argument-subsect} concerning cotorsion pairs
and ordinal-indexed filtrations.
 Let $(\sF,\sC)$ be the cotorsion pair in $R\modl$ generated by
the flat $R$\+modules and the $R$\+module~$E$.
 If $E[h^{-1}]\notin\sF$, then the desired $R$\+module $P\in\sC$
exists.
 If $E[h^{-1}]\in\sF$, then $E[h^{-1}]$ is a direct summand of
an $R$\+module filtered by the flat $R$\+modules and copies of
the $R$\+module~$E$.

 Notice that $E$ is an $(f)$\+torsion and $(g)$\+torsion $R$\+module
(or equivalently, an $(f,g)$\+torsion $R$\+module) in the sense of
the definition in Section~\ref{torsion-subsect}.
 In other words, the actions of~$f$ and~$g$ in $E$ are locally
nilpotent.

 Let us find an example where $E[h^{-1}]$ is not a submodule
of any $R$\+module filtered by flat $R$\+modules and~$E$; moreover,
$E[h^{-1}]$ is not even filtered by submodules of flat $R$\+modules
and of the $R$\+module~$E$.
 Firstly, no nonzero subquotient module of the $R$\+module $E[h^{-1}]$
is a submodule of a flat $R$\+module, because the action of~$f$
(or~$g$) is locally nilpotent in subquotients of $E[h^{-1}]$, while
the action of~$f$ in submodules of flat $R$\+modules is injective.
 Therefore, if $E[h^{-1}]$ is filtered by submodules of flat
$R$\+modules and of the $R$\+module $E$, then $E[h^{-1}]$ is filtered
by submodules of the $R$\+module $E$ only.
 Secondly, if there are neither any nonzero $(h)$\+torsion nor any
nonzero infinitely $h$\+divisible elements in $E$, then the uniquely
$h$\+divisible $R$\+module $E[h^{-1}]$ (assuming that it is nonzero)
cannot be filtered by submodules of the $R$\+module~$E$.
 Let us spell out this implication in some more detail.

 Indeed, let $(F_\beta)_{0\le\beta\le\alpha}$ be a filtration of
the $R$\+module $F_\alpha=E[h^{-1}]$ by $R$\+submodules $F_\beta$
indexed by an ordinal~$\alpha$ such that the successive quotients
$F_\beta/\bigcup_{\gamma<\beta}F_\gamma$ are isomorphic to submodules
of~$E$.
 Let $\delta$~be the minimal ordinal for which $F_\delta\ne0$.
 So, in particular, $F_\delta$ is isomorphic to a submodule of~$E$.
 Let $e\in F_\delta$ be a nonzero element, and let $n$~be the maximal
integer for which the element $h^{-n}e\in E[h^{-1}]$ belongs
to~$F_\delta$.
 In there is no such maximal integer~$n$, then $e$~is a nonzero
infinitely $h$\+divisible element in $F_\delta$, contrary to our
assumption.
 Otherwise, the coset $h^{-n-1}e+F_\delta$ is a nonzero element
annihilated by~$h$ in $E[h^{-1}]/F_\delta$.
 Let $\beta>\delta$ be the minimal ordinal such that the element
$h^{-n-1}e\in E[h^{-1}]$ belongs to~$F_\beta$.
 Then the coset of $h^{-n-1}e$ is a nonzero $(h)$\+torsion element in
$F_\beta/\bigcup_{\gamma<\beta}F_\gamma$, also contradicting our
assumption.

 For a specific counterexample, let $k$~be a field and $k[x,y,z]$ be
the ring of polynomials in three variables.
 Take $f=x$ and $g=y$, and put
$$
 E_1=k[x,y,z,x^{-1},y^{-1}]/(k[x,y,z,x^{-1}]+k[x,y,z,y^{-1}]).
$$
 For every element $t\in k[x,y,z]$ not belonging to the ideal
$(x,y)\sub k[x,y,z]$, consider the ring $R_t=k[x,y,z][t^{-1}]$
and the $R_t$\+module $E_t=E_1[t^{-1}]$.
 Pick an element $h_t\in R_t$ not belonging to the ideal $(x,y)\sub
k[x,y,z][t^{-1}]$ such that the image of~$h_t$ is not invertible in
the quotient ring $k[x,y,z][t^{-1}]/(x,y)$.
 Then there are neither any nonzero $(h_t)$\+torsion nor any nonzero
infinitely $h_t$\+divisible elements in~$E_t$, and $E_t[h_t^{-1}]\ne0$.
 According to the above, we can find a cotorsion $R_t$\+module $Q_t$
such that $\Ext^1_{R_t}(E_t,Q_t)=0$ but
$\Ext^1_{R_t}(E_t[h_t^{-1}],Q_t)\ne0$.
 Consequently, $\Ext^1_{R_t}(R_t[h_t^{-1}],\Hom_{R_t}(E_t,Q_t))\ne0$.

 Put $Q=\prod_{t\in k[x,y,z]\setminus(x,y)}Q_t$.
 So $Q$ is a cotorsion $R$\+module (by Lemma~\ref{cotors-restrict}(a)).
 Let $\widecheck Q$ be the related locally cotorsion contraherent
cosheaf on $X=\Spec k[x,y,z]$ and $\Q$ be the restriction of
$\widecheck Q$ to the open complement $Y$ to the closed subscheme
$\Spec k[x,y,z]/(x,y)$ in $\Spec k[x,y,z]$.
 Let $j\:Y\rarrow X$ be the open embedding morphism.
 Then we have seen that the cosheaf of $\O_X$\+modules $j_!\Q$ is
not contraherent.
 In fact, for any principal affine open subscheme
$U=\Spec R_t=\Spec k[x,y,z][t^{-1}]\sub X$ containing the point
$\p=(x,y)\in X$ there exists a smaller principal affine open subscheme
$V=\Spec R_t[h_t^{-1}]\sub U$ such that the cosheaf of
$\O_X$\+modules $j_!\Q$ does \emph{not} satisfy the contraherence
condition for the pair of affine open subschemes $V\sub U\sub X$.
 Therefore, for any open subscheme $W\sub X$ containing
the point $\p\in X$, the restriction of $j_!\Q$ to $W$ is not
contraherent.
 So the cosheaf $j_!\Q$ is not even locally contraherent on $X$ in
the sense of the definition in Section~\ref{locally-contraherent} below.
\end{ex}

 Recall that a scheme $X$ is called \emph{semi-separated}
\cite[Appendix~B]{TT}, if it admits an affine open covering with
affine pairwise intersections of the open subsets belonging to
the covering.
 Equivalently, a scheme $X$ is semi-separated if and only if the diagonal
morphism $X\rarrow X\times_{\Spec\boZ}X$ is affine, and if and only if
the intersection of any two affine open subschemes of $X$ is affine.
 Any morphism from an affine scheme to a semi-separated scheme is
affine, and the fibered product of any two affine schemes over
a semi-separated base scheme is an affine scheme.

 We will say that a morphism of schemes $f\:Y\rarrow X$ is
\emph{coaffine} if for any affine open subscheme $V\sub Y$ there exists
an affine open subscheme $U\sub X$ such that $f(V)\sub U$, and for any
two such affine open subschemes $f(V)\sub U'\sub X$, \
$f(V)\sub U''\sub X$ there exists a third affine open subscheme
$U\sub X$ such that $f(V)\sub U\sub U'\cap U''$.
 If the scheme $X$ is semi-separated, then the second condition
is trivial.
 (We will see below in Section~\ref{direct-inverse-loc-contra} that
the second condition is not actually needed for our constructions.)

 Any morphism into an affine scheme is coaffine.
 Any embedding of an open subscheme is coaffine.
 The composition of two coaffine morphisms between semi-separated schemes
is a coaffine morphism.

 Let $f\:Y\rarrow X$ be a very flat coaffine morphism of schemes
(see Section~\ref{very-flat-morphisms-subsect} for the definition
and discussion of the former property), and let $\P$ be
a contraherent cosheaf on~$X$. 
 Define a contraherent cosheaf $f^!\P$ on $Y$ as follows.

 Let $V\sub Y$ be an affine open subscheme.
 Pick an affine open subscheme $U\sub X$ such that $f(V)\sub U$,
and set $(f^!\P)[V]=\Hom_{\O_X(U)}(\O_Y(V),\P[U])$.
 Due to the contraherence condition on $\P$, this definition of
the $\O_Y(V)$\+module $(f^!\P)[V]$ does not depend on the choice
of an affine open subscheme $U\sub X$.
 Since $f$~is a very flat morphism, the $\O_Y(V)$\+module
$(f^!\P)[V]$ is contraadjusted by
Lemma~\ref{very-scalars-veryflat-case}(a).
 The contraherence condition obviously holds for~$f^!\P$.

 Let $f\:Y\rarrow X$ be a flat coaffine morphism of schemes, and
$\P$ be a locally cotorsion contraherent cosheaf on~$X$.
 Then the same rule as above defines a locally cotorsion contraherent
cosheaf $f^!\P$ on~$Y$.
 One just uses Lemma~\ref{cotors-coexten}(a) in place of
Lemma~\ref{very-scalars-veryflat-case}(a).
 For any coaffine morphism of schemes $f\:Y\rarrow X$ and a locally
injective contraherent cosheaf $\gJ$ on $X$ the very same rule
defines a locally injective contraherent cosheaf $f^!\gJ$ on~$Y$
(use Lemma~\ref{cotors-restrict}(b)).

 For an open embedding of schemes $j\:Y\rarrow X$ and a contraherent
cosheaf $\P$ on $X$ one clearly has $j^!\P\simeq\P|_Y$.

 If $f\:Y\rarrow X$ is an affine morphism of schemes and $\Q$ is
a locally cotorsion contraherent cosheaf on $Y$, then $f_!\Q$ is
a locally cotorsion contraherent cosheaf on~$X$.
 If $f\:Y\rarrow X$ is a flat affine morphism and $\gI$ is a locally
injective contraherent cosheaf on $Y$, then $f_!\gI$ is a locally
injective contraherent cosheaf on~$X$.
 (Use Lemmas~\ref{cotors-restrict}(a) and~\ref{cotors-coexten}(b).)

 Let $f\:Y\rarrow X$ be an affine coaffine morphism of schemes.
 Then for any contraherent cosheaf $\Q$ on $Y$ and any locally
injective contraherent cosheaf $\P$ on $X$ there is a natural
adjunction isomorphism $\Hom^X(f_!\Q,\P)\simeq\Hom^Y(\Q,f^!\P)$,
where $\Hom^X$ and $\Hom^Y$ denote the abelian groups of morphisms
in the categories of contraherent cosheaves on $X$ and~$Y$
(see the explanation two paragraphs below).

 If, in addition, the morphism~$f$ is flat, then such an isomorphism
holds for any contraherent cosheaf $\Q$ on $Y$ and any locally
cotorsion contraherent cosheaf $\P$ on~$X$; in particular,
$f_!$ and~$f^!$ form an adjoint pair of functors between
the exact categories of locally cotorsion contraherent cosheaves
$X\ctrh^\lct$ and $Y\ctrh^\lct$.
 Their restrictions also act as adjoint functors between
the exact categories of locally injective contraherent cosheaves
$X\ctrh^\lin$ and $Y\ctrh^\lin$.

 If the morphism~$f$ is very flat, then the functor
$f^!\:X\ctrh\rarrow Y\ctrh$ is right adjoint to the functor
$f_!\:Y\ctrh\rarrow X\ctrh$.
 In all the mentioned cases, both the abelian groups $\Hom^X(f_!\Q,\P)$
and $\Hom^Y(\Q,f^!\P)$ are identified with the group whose elements
are the collections of homomorphisms of $\O_X(U)$\+modules
$\Q[V]\rarrow\P[U]$, defined for all affine open subschemes
$U\sub X$ and $V\sub Y$ such that $f(V)\sub U$ and compatible
with the corestriction maps.

 All the functors between exact categories of contraherent cosheaves
constructed in the section above are exact and preserve
infinite products.
 For constructions of the direct image functor~$f_!$ (acting between
appropriate exact subcategories of adjusted objects in the exact
categories of contraherent cosheaves) for a nonaffine morphism of
schemes~$f$, see Sections~\ref{coflasque}
and~\ref{homology-subsection} below.

\subsection{$\Cohom$ from a quasi-coherent sheaf to a contraherent
cosheaf}  \label{cohom-subsection}
 Let $X$ be a scheme over an affine scheme $\Spec R$.
 Let $\M$ be a quasi-coherent sheaf on $X$ and $J$ be an injective
$R$\+module.
 Then the rule $U\mpsto\Hom_R(\M(U),J)$ for affine (and more generally,
quasi-compact quasi-separated) open subschemes $U\sub X$ defines
a contraherent cosheaf on~$X$ (cf.\ Remark~\ref{scheme-topology}).
 We will denote it by $\Cohom_R(\M,J)$.
 Since the $\O_X(U)$\+module $\Hom_R(\M(U),J)$ is cotorsion by
Lemma~\ref{cotors-hom-nc}(b), it is even a locally cotorsion contraherent
cosheaf.
 When $\F$ is a flat quasi-coherent sheaf on $X$ and $J$ is
an injective $R$\+module, the contraherent cosheaf $\Cohom_R(\F,J)$
is locally injective (by Lemma~\ref{cotors-hom-nc}(c)).

 We refer to the definitions of a very flat morphism of schemes and
a very flat quasi-coherent sheaf on a scheme given in
Section~\ref{very-flat-morphisms-subsect}.
 If $X\rarrow\Spec R$ is a very flat morphism of schemes and $\F$ is
a very flat quasi-coherent sheaf on $X$, then for any contraadjusted
$R$\+module $P$ the rule $U\mpsto\Hom_R(\F(U),P)$ for affine open
subschemes $U\sub X$ defines a contraherent cosheaf on~$X$.
 The contraadjustedness condition on the $\O_X(U)$\+modules
$\Hom_R(\F(U),P)$ holds by
Lemma~\ref{very-scalars-veryflat-case}(c), while the contraherence
condition follows from the quasi-coherence of~$\F$.
 We will denote the cosheaf so constructed by $\Cohom_R(\F,P)$.

 Analogously, if a scheme $X$ is flat over $\Spec R$ and
a quasi-coherent sheaf $\F$ on $X$ is flat (or, more generally,
the quasi-coherent sheaf $\F$ on $X$ is flat over $\Spec R$, in
the obvious sense based on Lemma~\ref{flatness-over-base-is-local}),
then for any cotorsion $R$\+module $P$ the rule $U\mpsto\Hom_R(\F(U),P)$
for affine open subschemes $U\sub X$ defines a contraherent cosheaf
on~$X$.
 In fact, the $\O_X(U)$\+modules $\Hom_R(\F(U),P)$ are cotorsion by
Lemma~\ref{cotors-hom-nc}(a), hence the contraherent cosheaf
$\Cohom_R(\F,P)$ constructed in this way is locally cotorsion.

 Let $\F$ be a very flat quasi-coherent sheaf on a scheme $X$ and $\P$
be a contraherent cosheaf on~$X$.
 Then the contraherent cosheaf $\Cohom_X(\F,\P)$ is defined by the rule
$U\mpsto\Hom_{\O_X(U)}(\F(U),\P[U])$ for all affine open subschemes
$U\sub X$.
 For any two embedded affine open subschemes $V\sub U\sub X$ one has
\begin{multline*}
 \Hom_{\O_X(V)}(\F(V),\P[V])\\ \simeq
 \Hom_{\O_X(V)}(\O_X(V)\ot_{\O_X(U)}\F(U)\;\Hom_{\O_X(U)}(\O_X(V),\P[U]))
 \\ \simeq \Hom_{\O_X(U)}(\O_X(V),\Hom_{\O_X(U)}(\F(U),\P[U])),
\end{multline*}
so the contraherence condition holds.
 The contraadjustedness condition follows from
Lemma~\ref{very-tensor-hom}(b).

 Similarly, if $\F$ is a flat quasi-coherent sheaf and $\P$ is 
locally cotorsion contraherent cosheaf on $X$, then the contraherent
cosheaf $\Cohom_X(\F,\P)$ is defined by the same rule
$U\mpsto\Hom_{\O_X(U)}(\F(U),\P[U])$ for all affine open subschemes
$U\sub X$.
 By Lemma~\ref{cotors-hom}(a), $\Cohom_X(\F,\P)$ is a locally cotorsion
contraherent cosheaf on~$X$.

 Finally, if $\M$ is a quasi-coherent sheaf on $X$ and $\gJ$ is
a locally injective contraherent cosheaf, then the contraherent cosheaf
$\Cohom_X(\M,\gJ)$ is defined by the very same rule.
 One checks the contraherence condition in the same way as above.
 By Lemma~\ref{cotors-hom}(b), $\Cohom_X(\M,\gJ)$ is a locally cotorsion
contraherent cosheaf on~$X$.
 If $\F$ is a flat quasi-coherent sheaf and $\gJ$ is a locally
injective contraherent cosheaf on $X$, then the contraherent
cosheaf $\Cohom_X(\F,\gJ)$ is locally injective (by
Lemma~\ref{cotors-hom}(c)).

 For any contraadjusted module $P$ over a commutative ring $R$, denote
by $\widecheck P$ the corresponding contraherent cosheaf on $\Spec R$.
 Let $f\:X\rarrow\Spec R$ be a morphism of schemes and $\F$ be
a quasi-coherent sheaf on $X$.
 Then whenever $\F$ is a very flat quasi-coherent sheaf and $f$~is
a very flat morphism, there is a natural isomorphism of contraherent
cosheaves $\Cohom_R(\F,P)\simeq\Cohom_X(\F,f^!\widecheck P)$ on $X$.
 Indeed, for any affine open subscheme $U\sub X$ one has
$$
 \Hom_R(\F(U),P)\simeq\Hom_{\O_X(U)}(\F(U),\Hom_R(\O_X(U),P))
 \simeq\Hom_R(\F(U),(f^!\widecheck P)[U]).
$$

 The same isomorphism holds whenever $\F$ is a flat quasi-coherent sheaf,
$f$~is a flat morphism, and $P$ is a cotorsion $R$\+module.
 Finally, for any quasi-coherent sheaf $\M$ on $X$, any morphism
$f\:X\rarrow\Spec R$, and any injective $R$\+module $J$ there is
a natural isomorphism of locally cotorsion contraherent cosheaves
$\Cohom_R(\M,J)\simeq\Cohom_X(\M,f^!\widecheck J)$ on~$X$.

\subsection{Contraherent cosheaves of $\fHom$ between quasi-coherent
sheaves}  \label{fHom-subsection}
 A quasi-coherent sheaf $\cP$ on a scheme $X$ is said to be
\emph{cotorsion}~\cite{EE} if $\Ext_X^1(\F,\cP)=0$ for any flat
quasi-coherent sheaf $\F$ on~$X$.
 Here $\Ext_X$ denotes the $\Ext$ groups in the abelian category of
quasi-coherent sheaves on~$X$.
 A quasi-coherent sheaf $\cP$ on $X$ is called \emph{contraadjusted}
if one has $\Ext^1_X(\F,\cP)=0$ for any very flat quasi-coherent
sheaf $\F$ on~$X$ (see Section~\ref{very-flat-morphisms-subsect}
for the definition of the latter notion).

 Clearly, the two classes of quasi-coherent sheaves on $X$ so defined
are closed under extensions, so they form full exact subcategories in
the abelian category of quasi-coherent sheaves.
 Also, these exact subcategories are closed under the passage to
direct summands of objects.

 For any affine morphism of schemes $f\:Y\rarrow X$, any flat
quasi-coherent sheaf $\F$ on $X$, and any quasi-coherent sheaf $\cP$
on $Y$ there is a natural isomorphism of the extension groups
$\Ext_Y^1(f^*\F,\cP)\simeq\Ext_X^1(\F,f_*\cP)$ (see, e.~g.,
\cite[Lemma~1.7(e)]{Pal}).
 Hence the classes of contraadjusted and cotorsion quasi-coherent
sheaves on schemes are preserved by the direct images with respect to
affine morphisms.

 In fact, more generally, for any quasi-compact quasi-separated
morphism $f\:Y\rarrow X$, any flat quasi-coherent sheaf $\F$ on $X$,
and any quasi-coherent sheaf $\cP$ on $Y$, there is a natural
injective map of abelian groups $\Ext_X^1(\F,f_*\cP)\rarrow
\Ext_Y^1(f^*\F,\cP)$ \,\cite[Lemma~1.7(c)]{Pal}.
 Here the conditions on the morphism~$f$ are imposed in order to make
sure that the direct images of quasi-coherent sheaves are quasi-coherent.
 Thus the classes of contraadjusted and cotorsion quasi-coherent
sheaves on schemes are even preserved by direct images with respect to
quasi-compact quasi-separated morphisms (see also
Corollary~\ref{cta-cot-direct} below).

 Let $\F$ be a quasi-coherent sheaf on a scheme~$X$.
 Suppose that an associative ring $R$ acts on $\F$ from the right
by quasi-coherent sheaf endomorphisms.
 Let $M$ be a left $R$\+module.
 Define a contravariant functor $\F\ot_R M$ from the category of
affine open subschemes $U\sub X$ to the category of abelian groups
by the rule $(\F\ot_R M)(U) = \F(U)\ot_R M$.
 The natural $\O_X(U)$\+module structures on the groups $(\F\ot_R M)(U)$
are compatible with the restriction maps $(\F\ot_R M)(U)\rarrow
(\F\ot_R M)(V)$ for embedded affine open subschemes $V\sub U\sub X$,
and the quasi-coherence condition $$(\F\ot_R M)(V) \simeq
\O_X(V)\ot_{\O_X(U)}(\F\ot_R M)(U)$$ holds
(see Remark~\ref{quasi-coherence-rem}).
 Therefore, the functor $\F\ot_R M$ extends uniquely to
a quasi-coherent sheaf on $X$, which we will denote also by $\F\ot_R M$.

 Let $\cP$ be a quasi-coherent sheaf on~$X$.
 Then the abelian group $\Hom_X(\F,\cP)$ of morphisms in the category of
quasi-coherent sheaves on $X$ has a natural left $R$\+module structure.
 One can easily construct a natural isomorphism of abelian groups
$\Hom_X(\F\ot_R M\;\cP)\simeq\Hom_R(M,\Hom_X(\F,\cP))$.

\begin{lem}  \label{ext-tensor-sheaf-mod}
 Suppose that\/ $\Ext_X^i(\F,\cP)=0$ for\/ $0 < i\le i_0$ and either
\par
\textup{(a)} $M$ is a flat left $R$\+module, or \par
\textup{(b)} the right $R$\+modules\/ $\F(U)$ are flat for all
affine open subschemes $U\sub X$. \par
 Then there is a natural isomorphism of abelian groups\/
$\Ext_X^i(\F\ot_R M\;\cP)\simeq\Ext_R^i(M,\Hom_X(\F,\cP))$
for all\/ $0\le i\le i_0$.
\end{lem}

\begin{proof}
 Replace $M$ by its projective $R$\+module resolution~$L_\bu$.
 Then $\Ext_X^i(\F\ot_R L_j\;\cP)=0$ for all $0 <i \le i_0$ and
all~$j$.
 Due to the flatness condition (a) or~(b), the complex of quasi-coherent
sheaves $\F\ot_R L_\bu$ is a resolution of the sheaf $\F\ot_R M$.
 Hence the complex of abelian groups $\Hom_X(\F\ot_R L_\bu\;\cP)$
computes $\Ext_X^i(\F\ot_R M\;\cP)$ for $0\le i\le i_0$.
 On the other hand, this complex is isomorphic to the complex
$\Hom_R(L_\bu,\Hom_R(\F,\cP))$, which computes
$\Ext_R^i(M,\Hom_X(\F,\cP))$.
\end{proof}

 Let $\F$ be a quasi-coherent sheaf with a right action of a ring~$R$
on a scheme $X$, and let $f\:Y\rarrow X$ be a morphism of schemes.
 Then $f^*\F$ is a quasi-coherent sheaf on $Y$ with a right action
of~$R$, and for any left $R$\+module $M$ there is a natural isomorphism
of quasi-coherent sheaves $f^*(\F\ot_R M)\simeq f^*\F\ot_R M$.
 Analogously, if $\G$ is a quasi-coherent sheaf on $Y$ with a right
action of~$R$ and $f$ is a quasi-compact quasi-separated morphism, then
$f_*\G$ is a quasi-coherent sheaf on $X$ with a right action of~$R$,
and for any left $R$\+module $M$ there is a natural morphism of
quasi-coherent sheaves $f_*\G\ot_R M\rarrow f_*(\G\ot_R M)$ on~$X$.
 If the morphism~$f$ is affine or the $R$\+module $M$ is flat, then
this map is an isomorphism of quasi-coherent sheaves on~$X$.

 Let $\F$ be a very flat quasi-coherent sheaf on a semi-separated
scheme $X$, and let $\cP$ be a contraadjusted quasi-coherent sheaf
on~$X$.
 Define a contraherent cosheaf $\fHom_X(\F,\cP)$ by the rule
$U\mpsto\Hom_X(j_*j^*\F,\cP)$ for any affine open subscheme $U\sub X$,
where $j\:U\rarrow X$ denotes the identity open embedding.
 Given two embedded affine open subschemes $V\sub U\sub X$ with
the identity embeddings $j\:U\rarrow X$ and $k\:V\rarrow X$,
the adjunction provides a natural map of quasi-coherent sheaves
$j_*j^*\F\rarrow k_*k^*\F$.
 There is also a natural action of the ring $\O_X(U)$ on
the quasi-coherent sheaf $j_*j^*\F$.
 Thus our rule defines a covariant functor with an $\O_X$\+module
structure on the category of affine open subschemes in~$X$.

 Let us check that the contraadjustedness and contraherence conditions
are satisfied.
 For a very flat $\O_X(U)$\+module $G$, we have
\begin{multline*}
 \Ext^1_{\O_X(U)}(G,\Hom_X(j_*j^*\F,\cP)) \\ \simeq
 \Ext^1_X((j_*j^*\F)\ot_{\O_X(U)}G\;\cP) \simeq
 \Ext^1_X(j_*(j^*\F\ot_{\O_X(U)}G)\;\cP) = 0,
\end{multline*}
by Lemma~\ref{ext-tensor-sheaf-mod} for $R=\O_X(U)$ and $i_0=1$, and
because $j_*(j^*\F\ot_{\O_X(U)}G)$ is a very flat quasi-coherent sheaf
on~$X$.
 The assumption that $j$~is an affine morphism (i.~e., $X$ is
semi-separated) is used here.
 For a pair of embedded affine open subschemes $V\sub U\sub X$,
we have
\begin{multline*}
 \Hom_{\O_X(U)}(\O_X(V),\Hom_X(j_*j^*\F,\cP))\simeq
 \Hom_X((j_*j^*\F)\ot_{\O_X(U)}\O_X(V)\;\cP) \\ \simeq
 \Hom_X(j_*(j^*\F\ot_{\O_X(U)}\O_X(V))\;\cP)\simeq
 \Hom_X(k_*k^*\F,\cP).
\end{multline*}

 Similarly one defines a locally cotorsion contraherent cosheaf
$\fHom_X(\F,\cP)$ for a flat quasi-coherent sheaf $\F$ and
a cotorsion quasi-coherent sheaf $\cP$ on~$X$.
 When $\F$ is a flat quasi-coherent sheaf and $\J$ is an injective
quasi-coherent sheaf on $X$, the contraherent cosheaf $\fHom_X(\F,\J)$
is locally injective.

 Finally, let $\M$ be a quasi-coherent sheaf on a quasi-separated
scheme $X$, and let $\J$ be an injective quasi-coherent sheaf on~$X$.
 Then a locally cotorsion contraherent cosheaf $\fHom_X(\M,\J)$
is defined by the very same rule.
 The proof of the cotorsion and contraherence conditions is the same
as above.

\begin{lem}  \label{fHom-cosections}
 Let $Y\sub X$ be a quasi-compact open subscheme in a semi-separated
scheme such that the identity open embedding $j\:Y\rarrow X$
is an affine morphism.  Then \par
\textup{(a)} for any very flat quasi-coherent sheaf\/ $\F$ and
contraadjusted quasi-coherent sheaf\/ $\cP$ on $X$, there is
a natural isomorphism of\/ $\O(Y)$\+modules\/ $\fHom_X(\F,\cP)[Y]
\simeq\Hom_X(j_*j^*\F,\cP)$; \par
\textup{(b)} for any flat quasi-coherent sheaf\/ $\F$ and cotorsion
quasi-coherent sheaf\/ $\cP$ on $X$, there is a natural isomorphism
of\/ $\O(Y)$\+modules\/ $\fHom_X(\F,\cP)[Y]\simeq
\Hom_X(j_*j^*\F,\cP)$. \par \smallskip
 Now let $Y\sub X$ be any quasi-compact open subscheme in
a quasi-separated scheme; let $j\:Y\rarrow X$ denote the identity
open embedding.  Then \par
\textup{(c)} for any quasi-coherent sheaf\/ $\M$ and injective
quasi-coherent sheaf\/ $\J$ on $X$, there is a natural isomorphism of\/
$\O(Y)$\+modules\/ $\fHom_X(\M,\J)[Y]\simeq\Hom_X(j_*j^*\M,\J)$.
\end{lem}

\begin{proof}
 Let $Y=\bigcup_{\alpha=1}^N U_\alpha$ be a finite affine open
covering of a quasi-separated scheme and $\G$ be a quasi-coherent
sheaf on~$Y$.
 Denote by $k_{\alpha_1,\dotsc,\alpha_i}$ the open embeddings
$U_{\alpha_1}\cap\dotsb\cap U_{\alpha_i}\rarrow Y$.
 Then there is a finite \v Cech exact sequence
\begin{multline} \label{cech-quasi} \textstyle
 0\lrarrow\G\lrarrow\bigoplus_\alpha k_\alpha{}_*k_\alpha^*\G
 \lrarrow\bigoplus_{\alpha<\beta}k_{\alpha,\beta}{}_*k_{\alpha,\beta}^*\G
 \\ \lrarrow\dotsb\lrarrow k_{1,\dotsc,N}{}_*k_{1,\dotsc,N}^*\G\lrarrow 0
\end{multline}
of quasi-coherent sheaves on~$Y$ (to check the exactness, it suffices
to consider the restrictions of this sequence to the open subschemes
$U_\alpha$, over each of which it is contractible).
 Set $\G=j^*\F$ or $j^*\M$.

 When the embedding morphism $j\:Y\rarrow X$ is affine (parts~(a\+b)),
the functor~$j_*$ preserves exactness of sequences of quasi-coherent
sheaves.
 So we obtain a finite exact sequence of quasi-coherent
sheaves on~$X$
\begin{multline} \label{cech-quasi-direct-image} \textstyle
 0\lrarrow j_*j^*\F\lrarrow\bigoplus_\alpha h_\alpha{}_*h_\alpha^*\F
 \lrarrow\bigoplus_{\alpha<\beta}h_{\alpha,\beta}{}_*h_{\alpha,\beta}^*\F
 \\ \lrarrow\dotsb\lrarrow h_{1,\dotsc,N}{}_*h_{1,\dotsc,N}^*\F
 \lrarrow 0,
\end{multline}
where $h_{\alpha_1,\dots,\alpha_i}$ denote the open embedings
$U_{\alpha_1}\cap\dotsb\cap U_{\alpha_i}\rarrow X$.

 Now \eqref{cech-quasi-direct-image}~is a sequence of very flat
quasi-coherent sheaves in the case~(a) and a sequence of flat
quasi-coherent sheaves in the case~(b).
 The functor $\Hom_X({-},\cP)$ transforms it into an exact sequence
of $\O(Y)$\+modules ending in
\begin{multline} \label{cech-quasi-direct-image-hom-from} \textstyle
 \bigoplus_{\alpha<\beta}\fHom_X(\F,\cP)[U_\alpha\cap U_\beta]
 \lrarrow \bigoplus_\alpha\fHom_X(\F,\cP)[U_\alpha] \\ \lrarrow
 \Hom_X(j_*j^*\F,\cP)\lrarrow0,
\end{multline}
and it remains to compare~\eqref{cech-quasi-direct-image-hom-from} with
the construction~\eqref{cosheaf-alt-recover} of the $\O(Y)$\+module
$\fHom_X(\F,\cP)[Y]$ in terms of the modules
$\fHom_X(\F,\cP)[U_\alpha]$ and $\fHom_X(\F,\cP)[U_\alpha\cap U_\beta]$.

 In the context of part~(c), the functor~$j_*$ is left exact, so we
have a left exact sequence of quasi-coherent sheaves on~$X$
$$ \textstyle
 0\lrarrow j_*j^*\M\lrarrow\bigoplus_\alpha h_\alpha{}_*h_\alpha^*\M
 \lrarrow\bigoplus_{\alpha<\beta}
 h_{\alpha,\beta}{}_*h_{\alpha,\beta}^*\M.
$$
 Applying the exact functor $\Hom_X({-},\J)$, we obtain a right exact
sequence of $\O(Y)$\+modules
\begin{multline*} \textstyle
 \bigoplus_{\alpha<\beta}
 \Hom_X(h_{\alpha,\beta}{}_*h_{\alpha,\beta}^*\M,\J)
 \lrarrow \bigoplus_\alpha\fHom_X(\M,\J)[U_\alpha] \\ \lrarrow
 \Hom_X(j_*j^*\M,\J)\lrarrow0.
\end{multline*}
 Let $U_\alpha\cap U_\beta=\bigcup_\gamma W_{\alpha\beta\gamma}$ be
a finite affine open covering of the quasi-compact open subscheme
$U_\alpha\cap U_\beta\sub X$.
 Then the similar right exact sequence for $U_\alpha\cap U_\beta
\sub X$ instead of $Y$ implies, in particular, that the natural map
$\bigoplus_\gamma\fHom_X(\M,\J)[W_{\alpha\beta\gamma}]\rarrow
\Hom_X(h_{\alpha,\beta}{}_*h_{\alpha,\beta}^*\M,\J)$ is surjective
for every pair of indices $\alpha$ and~$\beta$.
 Hence we arrive to a right exact sequence
\begin{multline*} \textstyle
 \bigoplus_{\alpha<\beta}^{\gamma}
 \fHom_X(\M,\J)[W_{\alpha\beta\gamma}]
 \lrarrow \bigoplus_\alpha\fHom_X(\M,\J)[U_\alpha] \\ \lrarrow
 \Hom_X(j_*j^*\M,\J)\lrarrow0
\end{multline*}
and it remains to refer once again to~\eqref{cosheaf-alt-recover}.
\end{proof}

 For any affine morphism of schemes $f\:Y\rarrow X$ and any
quasi-coherent sheaves $\M$ on $X$ and $\N$ on $Y$ there is
a natural isomorphism
\begin{equation}  \label{qcoh-projection}
 f_*(f^*\M\ot_{\O_Y}\N) \simeq \M\ot_{\O_X} f_*\N
\end{equation}
of quasi-coherent sheaves on~$X$ (``the projection formula'').
 In particular, for any quasi-coherent sheaves $\M$ and $\K$ on $X$
there is a natural isomorphism
\begin{equation} \label{push-pull-tensor-sh}
 f_*f^*(\M\ot_{\O_X}\K) \simeq \M\ot_{\O_X} f_*f^*\K
\end{equation}
of quasi-coherent sheaves on~$X$.
 Assuming that the quasi-coherent sheaf $\M$ on $X$ is flat, the same
isomorphisms hold for any quasi-compact quasi-separated morphism of
schemes $f\:Y\rarrow X$.

 For any embedding $j\:U\rarrow X$ of an affine open subscheme into
a semi-separated scheme $X$ and any quasi-coherent sheaves $\K$ and
$\M$ on $X$ there is a natural isomorphism
\begin{equation} \label{push-pull-tensor-mod}
 j_*j^*(\K\ot_{\O_X}\M) \simeq j_*j^*\K\ot_{\O_X(U)}\M(U)
\end{equation}
of quasi-coherent sheaves on~$X$.
 Assuming that the $\O_X(U)$\+module $\M(U)$ is flat, the same
isomorphism holds for a quasi-separated scheme~$X$.

 Recall that the \emph{quasi-coherent internal Hom} sheaf
$\qHom_{X\qc}(\M,\cP)$ for quasi-coherent sheaves $\M$ and $\cP$
on a scheme $X$ is defined as the quasi-coherent sheaf for
which there is a natural isomorphism of abelian groups
$\Hom_X(\K,\qHom_{X\qc}(\M,\cP))\simeq\Hom_X(\M\ot_{\O_X}\K\;\cP)$
for any quasi-coherent sheaf $\K$ on~$X$.
 The sheaf $\qHom_{X\qc}(\M,\cP)$ can be constructed by applying
the coherator functor~\cite[Sections~B.12\+-B.14]{TT} to
the sheaf of $\O_X$\+modules $\qHom_{\O_X}(\M,\cP)$.

\begin{lem} \label{ext-qhom-qc}
 Let $X$ be a scheme.  Then \par
\textup{(a)} for any very flat quasi-coherent sheaf\/ $\F$ and
contraadjusted quasi-coherent sheaf\/ $\cP$ on $X$,
the quasi-coherent sheaf\/ $\qHom_{X\qc}(\F,\cP)$ on $X$ is
contraadjusted; \par
\textup{(b)} for any flat quasi-coherent sheaf\/ $\F$ and
cotorsion quasi-coherent sheaf\/ $\cP$ on $X$,
the quasi-coherent sheaf\/ $\qHom_{X\qc}(\F,\cP)$ on $X$ is
cotorsion; \par
\textup{(c)} for any quasi-coherent sheaf\/ $\M$ and any
injective quasi-coherent sheaf\/ $\J$ on $X$,
the quasi-coherent sheaf\/ $\qHom_{X\qc}(\M,\J)$ on $X$ is
cotorsion; \par
\textup{(d)} for any flat quasi-coherent sheaf\/ $\F$ and any
injective quasi-coherent sheaf\/ $\J$ on $X$,
the quasi-coherent sheaf\/ $\qHom_{X\qc}(\F,\J)$ is injective.
\end{lem}

\begin{proof}
 We will prove part~(a); the proofs of the other parts are similar.
 Let $\G$ be a very flat quasi-coherent sheaf on~$X$.
 We will show that the functor $\Hom_X({-},\qHom_{X\qc}(\F,\cP))$
transforms any short exact sequence of quasi-coherent
sheaves $0\rarrow\K\rarrow\L\rarrow\G\rarrow0$ into a short
exact sequence of abelian groups.
 Indeed, the sequence of quasi-coherent sheaves
$0\rarrow \F\ot_{\O_X}\K\rarrow\F\ot_{\O_X}\L\rarrow\F\ot_{\O_X}\G
\rarrow0$ is exact, because $\F$ is flat (or because $\G$ is flat).
 Since $\F\ot_{\O_X}\G$ is very flat by
Lemma~\ref{very-tensor-hom}(a) and $\cP$ is contraadjusted,
the functor $\Hom_X({-},\cP)$ transforms the latter sequence of
sheaves into a short exact sequence of abelian groups.
\end{proof}

 It follows from the isomorphism~\eqref{push-pull-tensor-sh} that
for any very flat quasi-coherent sheaves $\F$ and $\G$ on
a semi-separated scheme $X$ and any contraadjusted quasi-coherent
sheaf $\cP$ on $X$ there is a natural isomorphism of contraherent
cosheaves
\begin{equation} \label{flat-flat-fhom-qhom}
 \fHom_X(\F\ot_{\O_X}\G\;\cP)\simeq\fHom_X(\G,\qHom_{X\qc}(\F,\cP)).
\end{equation}
 Similarly, for any flat quasi-coherent sheaves $\F$ and $\G$ 
and a cotorsion quasi-coherent sheaf $\cP$ on $X$ there is
a natural isomorphism~\eqref{flat-flat-fhom-qhom} of locally
cotorsion contraherent cosheaves.
 Finally, for any flat quasi-coherent sheaf $\F$, quasi-coherent
sheaf $\M$, and injective quasi-coherent sheaf $\J$ on $X$
there are natural isomorphisms of locally cotorsion contraherent
cosheaves
\begin{equation} \label{flat-inj-fhom-qhom}
 \fHom_X(\M\ot_{\O_X}\F\;\J)\simeq\fHom_X(\M,\qHom_{X\qc}(\F,\J))
 \simeq\fHom_X(\F,\qHom_{X\qc}(\M,\J)).
\end{equation}
 The leftmost isomorphism in~\eqref{flat-inj-fhom-qhom}
holds over any quasi-separated scheme~$X$.

 It follows from the isomorphism~\eqref{push-pull-tensor-mod} that
for any very flat quasi-coherent sheaves $\F$ and $\G$ and
a contraadjusted quasi-coherent sheaf $\cP$ on a semi-separated
scheme $X$ there is a natural isomorphism of contraherent cosheaves
\begin{equation} \label{flat-flat-fhom-cohom}
 \fHom_X(\G\ot_{\O_X}\F\;\cP)\simeq\Cohom_X(\F,\fHom_X(\G,\cP)).
\end{equation}
 Similarly, for any flat quasi-coherent sheaves $\F$ and $\G$
and a cotorsion quasi-coherent sheaf $\cP$ on $X$ there is
a natural isomorphism~\eqref{flat-flat-fhom-cohom} of locally
cotorsion contraherent cosheaves.
 Finally, for any flat quasi-coherent sheaf $\F$, quasi-coherent
sheaf $\K$, and injective quasi-coherent sheaf $\J$ on $X$
there are natural isomorphisms of locally cotorsion contraherent
cosheaves
\begin{equation} \label{flat-fhom-cohom-inj}
 \fHom_X(\K\ot_{\O_X}\F\;\J)\simeq\Cohom_X(\F,\fHom_X(\K,\J))
 \simeq\Cohom_X(\K,\fHom_X(\F,\J)).
\end{equation}
 The leftmost isomorphism in~\eqref{flat-fhom-cohom-inj}
holds over any quasi-separated scheme~$X$.

\begin{rem}  \label{cotorsion-sheaf-ambiguity}
 One can slightly generalize the constructions and results of this
section by weakening the definitions of contraadjusted and cotorsion
quasi-coherent sheaves.
 Namely, a quasi-coherent sheaf $\cP$ on $X$ may be called
weakly cotorsion if the functor $\Hom_X({-},\cP)$ transforms
short exact sequences of flat quasi-coherent sheaves on $X$ into
short exact sequences of abelian groups.
 The weakly contraadjusted quasi-coherent sheaves are defined similarly
(with the flat quasi-coherent sheaves replaced by very flat ones).
 Appropriate versions of Lemmas~\ref{ext-tensor-sheaf-mod}
and~\ref{ext-qhom-qc} can be proved in this setting, and
the contraherent cosheaves $\fHom$ can be defined.

 On a quasi-compact semi-separated scheme $X$ (or more generally, on
a scheme where there are enough flat or very flat quasi-coherent
sheaves), there is no difference between the weak and ordinary
cotorsion/contraadjusted quasi-coherent sheaves
(see Section~\ref{quasi-compact-quasi-coherent} below;
cf.~\cite[Sections~5.1.4 and~5.3]{Psemi}).
 One reason why we chose to use the stronger versions of these
conditions here rather that the weaker ones is that it is not
immediately clear whether the classes of weakly
cotorsion/contraadjusted quasi-coherent sheaves are closed under
extensions, or how the exact categories of such sheaves
should be defined.
\end{rem}

\subsection{Contratensor product of sheaves and cosheaves}
\label{contratensor-subsect}
 Let $X$ be a quasi-separated scheme and $\bB$ be an (initially fixed)
base of open subsets of $X$ consisting of some affine open subschemes.
 Let $\M$ be a quasi-coherent sheaf on $X$ and $\P$ be a cosheaf
of $\O_X$\+modules.

 The \emph{contratensor product} $\M\ocn_X\P$ (computed on
the base~$\bB$) is a quasi-coherent sheaf on $X$ defined as
the (nonfiltered) inductive limit of the following diagram of
quasi-coherent sheaves on $X$ indexed by affine open subschemes
$U\in\bB$ (cf.~\cite[Section~0.3.2]{Groth} and
Section~\ref{cosheaf-mod-subsect} above).

 To any affine open subscheme $U\in\bB$ with the identity open
embedding $j\:U\rarrow X$ we assign the quasi-coherent sheaf
$j_*j^*\M\ot_{\O_X(U)}\P[U]$ on $X$.
 For any pair of embedded affine open subschemes $V\sub U$, \ 
$V$,~$U\in\bB$ with the embedding maps $j\:U\rarrow X$ and
$k\:V\rarrow X$ there is the morphism of quasi-coherent sheaves
$$
 k_*k^*\M\ot_{\O_X(V)}\P[V]\lrarrow
 j_*j^*\M\ot_{\O_X(U)}\P[U]
$$
defined in terms of the natural isomorphism $k_*k^*\M\simeq
j_*j^*\M\ot_{\O_X(U)}\O_X(V)$ of quasi-coherent sheaves on~$X$
and the $\O_X(U)$\+module morphism $\P[V]\rarrow\P[U]$.

 Let $\M$ and $\J$ be quasi-coherent sheaves on a quasi-separated
scheme $X$ for which the contraherent cosheaf $\fHom_X(\M,\J)$ is
defined (i.~e., one of the sufficient conditions given in
Section~\ref{fHom-subsection} for the construction of $\fHom$
to make sense is satisfied).
 Then for any cosheaf of $\O_X$\+modules $\P$ there is a natural
isomorphism of abelian groups 
\begin{equation} \label{fHom-contratensor-adjunction}
 \Hom_X(\M\ocn_X\P\;\J)\simeq\Hom^{\O_X}(\P,\fHom_X(\M,\J)).
\end{equation}
 In other words, the functor $\M\ocn_X{-}$ is left adjoint to
the functor $\fHom_X(\M,{-})$ ``wherever the latter is defined''.

 Indeed, both groups of homomorphisms consist of all the compatible
collections of morphisms of quasi-coherent sheaves
$$
 j_*j^*\M\ot_{\O_X(U)}\P[U]\lrarrow \J
$$
on $X$, or equivalently, all the compatible collections of
morphisms of $\O_X(U)$\+modules
$$
 \P[U]\lrarrow\Hom_X(j_*j^*\M,\J)
$$
defined for all the identity embeddings $j\:U\rarrow X$ of affine open
subschemes $U\in\bB$.
 The compatibility is with respect to the identity embeddings of affine
open subschemes $h\:V\rarrow U$, \ $V$,~$U\in\bB$, into one another.

 In particular, the adjunction
isomorphism~\eqref{fHom-contratensor-adjunction} holds for any
quasi-coherent sheaf $\M$, cosheaf of $\O_X$\+modules $\P$, and
injective quasi-coherent sheaf~$\J$.
 Since there are enough injective quasi-coherent sheaves, it follows
that the quasi-coherent sheaf of contratensor product $\M\ocn_X\P$
does not depend on the base of open affines $\bB$ that was used
to construct it.

 More generally, let $\bD$ be a partially ordered set endowed with
an order-preserving map into the set of all affine open subschemes
of $X$, which we will denote by $a\mpsto U(a)$, i.~e., one has
$U(b)\sub U(a)$ whenever $b\le a\in\bD$.
 Suppose that $X=\bigcup_{a\in\bD}U(a)$ and for any $a$, $b\in\bD$
the intersection $U(a)\cap U(b)\sub X$ is equal to the union
$\bigcup_{c\le a,b} U(c)$.
 Then the inductive limit of the diagram $j_a{}_*j_a^*\M\ot_{\O_X(U_a)}
\P[U_a]$ indexed by $a\in\bD$, where $j_a$~denotes the open embedding
$U_a\rarrow X$, is naturally isomorphic to the contratensor product
$\M\ocn_X\P$.
 The point of this generalization is that a topology base $\bB$ has to
contain arbitrarily small affine open subschemes in $X$, but a diagram
$\bD$ does not need to.

 Indeed, given a cosheaf of $\O_X$\+modules $\P$ and a contraherent
cosheaf $\Q$ on $X$, an arbitrary collection of morphisms of
$\O_X(U_a)$\+modules $\P[U_a]\rarrow\Q[U_a]$ compatible with
the corestriction maps for $b\le a$ uniquely determines a morphism
of cosheaves of $\O_X$\+modules $\P\rarrow\Q$
(see Lemma~\ref{hom-into-contraherent}).
 In particular, this applies to the case of a contraherent
cosheaf $\Q=\fHom_X(\M,\J)$.

 The isomorphism $j_*j^*(\M\ot_{\O_X}\K) \simeq \M\ot_{\O_X} j_*j^*\K$
for an embedding of affine open subscheme $j\:U\rarrow X$ and
quasi-coherent sheaves $\M$ and $\K$ on $X$
(see~\eqref{push-pull-tensor-sh}) allows to construct a natural
isomorphism of quasi-coherent sheaves
\begin{equation}  \label{tensor-contratensor-assoc}
 \M\ot_{\O_X}(\K\ocn_X\P)\simeq(\M\ot_{\O_X}\K)\ocn_X\P
\end{equation}
for any quasi-coherent sheaves $\M$ and $\K$ and any cosheaf of
$\O_X$\+modules $\P$ on a semi-separated scheme~$X$.
 The same isomorphism holds over a quasi-separated scheme $X$,
assuming that the quasi-coherent sheaf $\M$ is flat.

\Section{Locally Contraherent Cosheaves} \label{loc-contra-sect}

\subsection{Exact category of locally contraherent cosheaves}
\label{locally-contraherent}
 A cosheaf of $\O_X$\+mod\-ules $\P$ on a scheme $X$ is called
\emph{locally contraherent} if every point $x\in X$ has an open
neighborhood $x\in W\sub X$ such that the cosheaf of $\O_W$\+modules
$\P|_W$ is contraherent.

 Given an open covering $\bW=\{W\}$ of a scheme $X$, a cosheaf of
$\O_X$\+modules $\P$ is called \emph{$\bW$\+locally contraherent}
if for every open subscheme $W\sub X$ belonging to $\bW$ the cosheaf
of $\O_W$\+modules $\P|_W$ is contraherent on~$W$.
 Obviously, a cosheaf of $\O_X$\+modules $\P$ is locally contraherent
if and only if there exists an open covering $\bW$ of the scheme $X$
such that $\P$ is $\bW$\+locally contraherent.

 Let us call an open subscheme of a scheme $X$ \emph{subordinate}
to an open covering $\bW$ if it is contained in one of the open
subsets of $X$ belonging to~$\bW$.
 An open covering $X=\bigcup_\alpha U_\alpha$ is said to be
subordinate to $\bW$ if every affine open subscheme $U_\alpha$ is
subordinate to~$\bW$.
 Notice that, by the definition of a contraherent cosheaf,
the property of a cosheaf of $\O_X$\+modules to be $\bW$\+locally
contraherent only depends on the collection of all affine open
subschemes $U\sub X$ subordinate to~$\bW$.

\begin{thm}
 Let\/ $\bW$ be an open covering of a scheme~$X$.
 Then the restriction of cosheaves of\/ $\O_X$\+modules to 
the base of open subsets of $X$ consisting of all the affine open
subschemes subordinate to\/ $\bW$ induces an equivalence between
the category of\/ $\bW$\+locally contraherent cosheaves on $X$
and the category of covariant functors with\/ $\O_X$\+module
structures on the category of affine open subschemes of $X$
subordinate to\/ $\bW$, satisfying the contraadjustedness and
contraherence conditions~\textup{(i\+ii)} of
Section~\textup{\ref{contraherent-definition}} for all
affine open subschemes $V\sub U\sub X$ subordinate to\/~$\bW$.
\end{thm}

\begin{proof}
 The same as in Theorem~\ref{contraherent-base}, except that the base
of affine open subschemes of $X$ subordinate to $\bW$ is considered
throughout.
\end{proof}

 Let $X$ be a scheme and $\bW$ be its open covering.
 By Theorem~\ref{cosheaf-base-thm}, the category of cosheaves of
$\O_X$\+modules is a full subcategory of the category of covariant
functors with $\O_X$\+module structures on the category of affine
open subschemes of $X$ subordinate to~$\bW$.
 The category of such functors with $\O_X$\+module structures
is clearly abelian, has exact functors of infinite direct sum and
infinite product, and the functors of cosections over a particular
affine open subscheme subordinate to $\bW$ are exact on it and
preserve infinite direct sums and products.

 The full subcategory of cosheaves of $\O_X$\+modules in this
abelian category is closed under extensions, cokernels, and
infinite direct sums.
 Indeed, cosheaves of $\O_X$\+modules are characterized among arbitrary
covariant functors with $\O_X$\+module structures by the right
exactness condition~\eqref{base-cosheaf}, which is clearly preserved
by extensions, cokernels, and infinite direct sums.
 For the quasi-compactness reasons explained in
Remark~\ref{scheme-topology}, this full subcategory is also closed
under infinite products.

 Therefore, the category of cosheaves of $\O_X$\+modules acquires
the induced exact category structure with exact functors of
infinite direct sum and product, and exact functors of cosections
on affine open subschemes subordinate to~$\bW$.
 Let us denote the category of cosheaves of $\O_X$\+modules endowed
with this exact category structure (which, of course, depends on
the choice of a covering~$\bW$) by $\O_X\cosh_\bW$.
 Along the way we have shown that infinite products
exist in the additive category of cosheaves of $\O_X$\+modules on
a scheme $X$, and the functors of cosections over quasi-compact
quasi-separated open subschemes of $X$ preserve them.

 The full subcategory of $\bW$\+locally contraherent cosheaves is
closed under extensions, cokernels of admissible monomorphisms,
and infinite products in the exact category $\O_X\cosh_\bW$.
 Indeed, the contraadjustedness condition is preserved by these
operations, because so is the full subcategory of contraadjusted modules
in the abelian category of modules over a commutative ring~$R$.
 Short exact sequences and infinite products of contraadjusted
$R$\+modules are also preserved by the $\Hom_R$ functors from very flat
$R$\+modules, such as the $\O(U)$\+modules $\O(V)$ for affine open
subschemes $V$ in affine schemes~$U$; and it remains to observe that
isomorphisms of modules are preserved by extensions, cokernels of
monomorphisms, and the infinite products in order to conclude that
the contraherence condition is preserved.
 Thus the category of $\bW$\+locally contraherent cosheaves
has the induced exact category structure with exact functors of
infinite product, and exact functors of cosections over affine open
subschemes subordinate to~$\bW$.
 We denote this exact category of $\bW$\+locally contraherent cosheaves
on a scheme $X$ by $X\lcth_\bW$.

 More explicitly, a short sequence of $\bW$\+locally contraherent
cosheaves $0\rarrow\P\rarrow\Q\rarrow\R\rarrow0$ is exact in
$X\lcth_\bW$ if the sequence of cosection modules $0\rarrow\P[U]
\rarrow\Q[U]\rarrow\R[U]\rarrow0$ is exact for every affine open
subscheme $U\sub X$ subordinate to~$\bW$.
 Passing to the inductive limit with respect to refinements of
the coverings $\bW$, we obtain the exact category structure on
the category of locally contraherent cosheaves $X\lcth$
on the scheme~$X$.

 A $\bW$\+locally contraherent cosheaf $\P$ on $X$ is said to be
\emph{locally cotorsion} if for any affine open subscheme $U\sub X$
subordinate to $\bW$ the $\O_X(U)$\+module $\P[U]$ is cotorsion.
 By Lemma~\ref{cotors-inj-covering}(a), this definition can be
equivalently rephrased by saying that a locally contraherent cosheaf
$\P$ on $X$ is locally cotorsion if and only if for every affine open
subscheme $U\sub X$ such that the cosheaf $\P|_U$ is contraherent
on the scheme $U$ the $\O(U)$\+module $\P[U]$ is cotorsion.
 It suffices to check the local cotorsion condition for affine open
subschemes $U$ belonging to any chosen affine open covering of
the scheme $X$ subordinate to~$\bW$.
 
 A $\bW$\+locally contraherent cosheaf $\gJ$ on $X$ is called
\emph{locally injective} if for any affine open subscheme $U\sub X$
subordinate to $\bW$ the $\O_X(U)$\+module $\gJ[U]$ is injective.
 By Lemma~\ref{cotors-inj-covering}(b), a locally contraherent
cosheaf $\gJ$ on $X$ is locally injective if and only if for every
affine open subscheme $U\sub X$ such that the cosheaf $\gJ|_U$ is
contraherent on the scheme $U$ the $\O(U)$\+module $\gJ[U]$
is injective.
 It suffices to check the local injectivity condition for affine open
subschemes $U$ belonging to any chosen affine open covering of
the scheme $X$ subordinate to~$\bW$.

 One defines the exact categories $X\lcth_\bW^\lct$ and
$X\lcth_\bW^\lin$ of locally cotorsion and locally injective
$\bW$\+locally contraherent cosheaves on $X$ in the same way as above.
 These are full subcategories closed under extensions, infinite
products, and cokernels of admissible monomorphisms in $X\lcth_\bW$,
with the induced exact category structures.
 The argument for the proof of the closure properties is the same
as above.
 Passing to the inductive limit with respect to refinements, we obtain
the exact categories $X\lcth^\lct$ and $X\lcth^\lin$ of locally
cotorsion and locally injective locally contraherent cosheaves on~$X$.

 To emphasize that we are speaking about some locally contraherent
cosheaves that need not be locally cotorsion or locally injective,
we sometimes call arbitrary locally contraherent cosheaves
``locally contraadjusted''.

 The following lemma may help the reader to feel more comfortable.

\begin{lem} \label{acyclicity-in-lcth-criterion}
 Let $X$ be a scheme with an open covering\/~$\bW$. \par
\textup{(a)} A complex\/ $\P^\bu$ in the exact category of\/
$\bW$\+locally contraherent cosheaves\/ $X\lcth_\bW$ is acyclic
if and only if, for every affine open subscheme $U\sub X$
subordinate to\/ $\bW$, the complex of\/ $\O_X(U)$\+modules\/
$\P^\bu[U]$ is acyclic.
 It suffices to check this condition for affine open subschemes
belonging to any chosen affine open covering of $X$ subordinate
to\/~$\bW$. \par
\textup{(b)} A complex\/ $\P^\bu$ in the exact category of
locally cotorsion\/ $\bW$\+locally contraherent cosehaves\/
$X\lcth_\bW^\lct$ is acyclic if and only if, for every affine
open subscheme $U\sub X$ subordinate to\/ $\bW$, the complex of\/
$\O_X(U)$\+modules\/ $\P^\bu[U]$ is acyclic.
 It suffices to check this condition for affine open subschemes
belonging to any chosen affine open covering of $X$ subordinate
to\/~$\bW$.
\end{lem}

\begin{proof}
 Part~(a): notice that the modules of cocycles in any acyclic complex
of contraadjusted modules are contraadjusted, too, since the class of
contraadjusted modules is closed under quotients.
 Now let $U\sub X$ be an affine open subscheme subordinate to $\bW$,
and let $\gK^n[U]$, \,$n\in\boZ$, denote the modules of cocycles in
the acyclic complex of $\O_X(U)$\+modules $\P^\bu[U]$.
 Let $V$ be an affine open subscheme in~$U$.
 Applying the left exact functor $\Hom_{\O(U)}(\O(V),{-})$ to
the left exact sequence of $\O(U)$\+modules $0\rarrow\gK^n[U]
\rarrow\P^n[U]\rarrow\P^{n+1}[U]$ and using the structure
isomorphisms $\P^i[V]\simeq\Hom_{\O(U)}(\O(V),\P^i[U])$ of
the $\bW$\+locally contraherent cosheaves $\P^i$, we obtain a natural
isomorphism of $\O_X(V)$\+modules $\gK^n[V]\simeq
\Hom_{\O(U)}(\O(V),\gK^n[U])$.
 So we can construct $\bW$\+locally contraherent cosheaves $\gK^n$
on $X$, and it is clear that the original complex of $\bW$\+locally
contraherent cosheaves $\P^\bu$ can be obtained by splicing short
exact sequences $0\rarrow\gK^n\rarrow\P^n\rarrow\gK^{n+1}\rarrow0$.
 This proves the first assertion.
 To deduce the second one, one can use
Lemma~\ref{cta-cot-acyclicity-colocal}(a).

 The proof of the first assertion of part~(b) is similar and based on
Theorem~\ref{cotorsion-periodicity}.
 The second assertion is provable using
Lemma~\ref{cta-cot-acyclicity-colocal}(b).
\end{proof}

\subsection{Contraherent and locally contraherent cosheaves}
\label{counterex-subsect}
 By Lemma~\ref{very-exact-local}(a), a short sequence
of $\bW$\+locally contraherent cosheaves on $X$ is exact
in $X\lcth$ (i.~e., after some refinement of the covering)
if and only if it is exact in $X\lcth_\bW$.
 By Lemma~\ref{very-exact-local}(b), a morphism of $\bW$\+locally
contraherent cosheaves is an admissible epimorphism in
$X\lcth$ if and only if it is an admissible epimorphism in
$X\lcth_\bW$.

 Analogously, by Lemma~\ref{cotors-exact-local}(a), a short sequence
of locally cotorsion $\bW$\+locally contraherent cosheaves on $X$
is exact in $X\lcth^\lct$ if and only if it is exact in
$X\lcth_\bW^\lct$.
 By Lemma~\ref{cotors-exact-local}(b), a morphism of locally
cotorsion $\bW$\+locally contraherent cosheaves on $X$ is an admissible
epimorphism in $X\lcth^\lct$ if and only if it is an admissible
epimorphism in $X\lcth^\lct_\bW$.
 The similar assertions hold for locally injective locally contraherent
cosheaves, and they are provable in the same way using
Lemma~\ref{injective-exact-local}.
 
 On the other hand, a morphism in $X\lcth_\bW$, \ $X\lcth_\bW^\lct$, or
$X\lcth_\bW^\lin$ is an admissible monomorphism if and only if it acts
injectively on the modules of cosections over all the affine open
subschemes $U\sub X$ subordinate to~$\bW$.
 The following counterexample shows that this condition \emph{does}
change when the covering~$\bW$ is refined.

 In other words, the full subcategory $X\lcth_\bW\sub X\lcth$ is closed
under the passage to the kernels of admissible epimorphisms, but
\emph{not} to the cokernels of admissible monomorphisms in $X\lcth$.
 Once we show that, it will also follow that there \emph{do} exist
locally contraherent cosheaves that are not contraherent.
 The locally cotorsion and locally injective contraherent cosheaves
have all the same problems.

\begin{ex}  \label{loc-contraherent-counterexample}
 Let $R$ be a commutative ring and $f$, $g\in R$ be two elements
generating the unit ideal.
 Let $M$ be a nonzero $R$\+module containing no $f$\+divisible
or $g$\+divisible elements, i.~e., $\Hom_R(R[f^{-1}],M)=0=
\Hom_R(R[g^{-1}],M)$.

 Let $M\rarrow P$ be an embedding of $M$ into a contraadjusted
$R$\+module $P$, and let $Q$ be the cokernel of this embedding.
 Then $Q$ is also a contraadjusted $R$\+module.
 One can take $R$ to be a Dedekind domain, so that it has homological
dimension~$1$; then whenever $P$ is a cotorsion or injective $R$\+module,
$Q$ has the same property.

 Consider the morphism of contraherent cosheaves $\widecheck P\rarrow
\widecheck Q$ on $\Spec R$ related to the surjective morphism of
contraadjusted (cotorsion, or injective) $R$\+modules $P\rarrow Q$.
 In restriction to the covering of $\Spec R$ by the two principal
affine open subsets $\Spec R[f^{-1}]$ and $\Spec R[g^{-1}]$, we
obtain two morphisms of contraherent cosheaves related to the two
morphisms of contraadjusted modules $\Hom_R(R[f^{-1}],P)\rarrow
\Hom_R(R[f^{-1}],Q)$ and $\Hom_R(R[g^{-1}],P)\rarrow\Hom_R(R[g^{-1}],Q)$
over the rings $\Spec R[f^{-1}]$ and $\Spec R[g^{-1}]$.

 Due to the condition imposed on $M$, the latter two morphisms of
contraadjusted modules are injective.
 On the other hand, the morphism of contraadjusted $R$\+modules
$P\rarrow Q$ is not.
 It follows that the cokernel $\R$ of the morphism of contraherent
cosheaves $\widecheck P\rarrow\widecheck Q$ taken in the category of all
cosheaves of $\O_{\Spec R}$\+modules (or equivalently, in the category of
copresheaves of $\O_{\Spec R}$\+modules) is contraherent in restriction
to $\Spec R[f^{-1}]$ and $\Spec R[g^{-1}]$, but not over $\Spec R$.
 In fact, one has $\R[\Spec R]=0$ (since the morphism $P\rarrow Q$
is surjective).

 Put $X=\Spec R$, and let $\bW$ be the covering of $X$ by two affine
open subschemes $\Spec R[f^{-1}]$ and $\Spec R[g^{-1}]$.
 Then the short sequence $0\rarrow\widecheck P\rarrow\widecheck Q
\rarrow\R\rarrow0$ is admissible exact in the exact category
$X\lcth_\bW$ of $\bW$\+locally contraherent cosheaves on $\Spec R$
(even though the map $\widecheck P[X]\rarrow\widecheck Q[X]$ is not
injective).
\end{ex}

 Let us point out that for any cosheaf of $\O_X$\+modules $\P$ on
a scheme $X$ such that the $\O_X(U)$\+modules $\P[U]$ are
contraadjusted for all affine open subschemes $U\sub X$ subordinate
to a particular open covering $\bW$, the $\O_X(U)$\+modules
$\P[U]$ are contraadjusted for \emph{all} affine open subschemes
$U\sub X$.
 This is so simply because the class of contraadjusted modules is
closed under finite direct sums, restrictions of scalars, and
cokernels.
 So the contraadjustedness condition~(ii) of
Section~\textup{\ref{contraherent-definition}} is, in fact, local;
it is the contraherence condition~(i) that isn't.

 In the rest of the section we will explain how to distinguish
the contraherent cosheaves among all the locally contraherent ones.
 Let $X$ be a semi-separated scheme, $\bW$ be its open covering, and
$\{U_\alpha\}$ be an affine open covering subordinate to~$\bW$
(i.~e., consisting of affine open subschemes subordinate to~$\bW$).

 Let $\P$ be a $\bW$\+locally contraherent cosheaf on~$X$.
 Consider the homological \v Cech complex of abelian groups
(or $\O(X)$\+modules) $C_\bu(\{U_\alpha\},\P)$ of the form
\begin{equation} \label{cech-homol} \textstyle
 \dotsb\lrarrow
 \bigoplus_{\alpha<\beta<\gamma} \P[U_\alpha\cap U_\beta\cap U_\gamma]
 \lrarrow \bigoplus_{\alpha<\beta} \P[U_\alpha\cap U_\beta]
 \lrarrow \bigoplus_\alpha \P[U_\alpha].
\end{equation}
 Here (as in the sequel) our notation presumes the indices~$\alpha$
to be linearly ordered.
 More generally, the complex~\eqref{cech-homol} can be considered for
any open covering $U_\alpha$ of a topological space~$X$ and any
cosheaf of abelian groups $\P$ on~$X$.
 Let $\Delta(X,\P)=\P[X]$ denote the functor of global cosections of
(locally contraherent) cosheaves on~$X$; then, by the definition,
we have $\Delta(X,\P)\simeq H_0C_\bu(\{U_\alpha\},\P)$.

\begin{lem}  \label{global-contraherence-criterion}
 Let $U$ be an affine scheme with an open covering\/ $\bW$ and
a finite affine open covering\/ $\{U_\alpha\}$ subordinate to\/~$\bW$.
 Then a\/ $\bW$\+locally contraherent cosheaf\/ $\P$ on $U$ is
contraherent if and only if $H_{>0}C_\bu(\{U_\alpha\},\P)=0$.
\end{lem}

\begin{proof}
 The ``only if'' part is provided by Lemma~\ref{very-open-covering}(b).
 Let us prove ``if''.
 If the \v Cech complex $C_\bu(\{U_\alpha\},\P)$ has no higher
homology, then it is a finite resolution of the $\O(U)$\+module $\P[U]$
by contraadjusted $\O(U)$\+modules.
 As we have explained above, the $\O(U)$\+module $\P[U]$ is
contraadjusted, too.

 For any affine open subscheme $V\sub U$, consider the \v Cech
complex $C_\bu(\{V\cap U_\alpha\}\;\P|_V)$ related to the restrictions
of our cosheaf $\P$ and our covering $U_\alpha$ to the open
subscheme~$V$.
 The complex $C_\bu(\{V\cap U_\alpha\}\;\P|_V)$ can be obtained from
the complex $C_\bu(\{U_\alpha\},\P)$ by applying the functor
$\Hom_{\O(U)}(\O(V),{-})$.
 We have
$$
 H_0C_\bu(\{U_\alpha\},\P)\simeq\P[U] \quad
 \text{and} \quad
 H_0C_\bu(\{V\cap U_\alpha\}\;\P|_V)\simeq\P[V].
$$

 Since the functor $\Hom_{\O(U)}(\O(V),{-})$ preserves exactness
of short sequences of contraadjusted $\O(U)$\+modules, we conclude
that $\P[V]\simeq\Hom_{\O(U)}(\O(V),\P[U])$.
 Both the contraadjustedness and contraherence conditions now
have been verified.
\end{proof}

\begin{cor}  \label{affine-lctrh-extensions}
 If a\/ $\bW$\+locally contraherent cosheaf\/ $\Q$ on an affine scheme
$U$ is an extension of two contraherent cosheaves\/ $\P$ and\/ $\R$ in
the exact category $U\lcth_\bW$ (or\/ $\O_U\cosh_\bW$), then\/ $\Q$ is
also a contraherent cosheaf on~$U$.
\end{cor}

\begin{proof}
 Pick a finite affine open covering $\{U_\alpha\}$ of the affine scheme
$U$ subordinate to the covering~$\bW$.
 Then the complex of abelian groups $C_\bu(\{U_\alpha\},\Q)$ is
an extension of the complexes of abelian groups
$C_\bu(\{U_\alpha\},\P)$ and $C_\bu(\{U_\alpha\},\R)$.
 Hence whenever the latter two complexes have no higher homology,
neither does the former one.
\end{proof}

\begin{cor} \label{nonaffine-lctrh-extensions}
 For any scheme $X$ and any its open covering\/ $\bW$, the full
subcategory of\/ $\bW$\+locally contraherent cosheaves on $X$ is
closed under extensions in the exact category of locally contraherent
cosheaves on~$X$.
 In particular, the full subcategory of contraherent cosheaves
on $X$ is closed under extensions in the exact category of locally
contraherent (or\/ $\bW$\+locally contraherent) cosheaves on~$X$.
\end{cor}

\begin{proof}
 Follows easily from Corollary~\ref{affine-lctrh-extensions}.
\end{proof}

 In the terminology of~\cite[Section~4]{Kel}, one can say that
$X\lcth_\bW$ is a fully exact subcategory of $X\lcth$.
 The exact category structure on $X\lcth_\bW$ is inherited from
$X\lcth$ (as mentioned in the beginning of this section).

\subsection{Direct and inverse images of locally contraherent cosheaves}
\label{direct-inverse-loc-contra}
 Let $\bW$ be an open covering of a scheme $X$ and $\bT$ be an open
covering of a scheme~$Y$.
 A morphism of schemes $f\:Y\rarrow X$ is called
\emph{$(\bW,\bT)$\+affine} if for any affine open subscheme $U\sub X$
subordinate to $\bW$ the open subscheme $f^{-1}(U)\sub Y$ is affine
and subordinate to~$\bT$.
 Any $(\bW,\bT)$\+affine morphism is affine.

 Let $f\:Y\rarrow X$ be a $(\bW,\bT)$\+affine morphism of schemes and
$\Q$ be a $\bT$\+locally contraherent cosheaf on~$Y$.
 Then the cosheaf of $\O_X$\+modules $f_!\Q$ on $X$ is $\bW$\+locally
contraherent.
 The proof of this assertion is similar to that of its global version
in Section~\ref{contra-direct-inverse}.
 We have constructed an exact functor of direct image
$f_!\:Y\lcth_\bT\rarrow X\lcth_\bW$ between the exact categories of
$\bT$\+locally contraherent cosheaves on $Y$ and $\bW$\+locally
contraherent cosheaves on~$X$.

 A morphism of schemes $f\:Y\rarrow X$ is called
\emph{$(\bW,\bT)$\+coaffine} if for any affine open subscheme $V\sub Y$
subordinate to $\bT$ there exists an affine open subscheme $U\sub X$
subordinate to $\bW$ such that $f(V)\sub U$.
 Notice that for any fixed open covering $\bW$ of a scheme $X$ and any
morphism of schemes $f\:Y\rarrow X$ the covering $\bT$ of the scheme $Y$
consisting of all the full preimages $f^{-1}(U)$ of affine open
subschemes $U\sub X$ has the property that the morphism $f\:Y\rarrow X$
is $(\bW,\bT)$\+coaffine.

 If the morphism~$f$ is affine, it is also $(\bW,\bT)$\+affine with
respect to the covering $\bT$ constructed in this way.
 A morphism $f\:Y\rarrow X$ is simultaneously $(\bW,\bT)$\+affine and
$(\bW,\bT)$\+coaffine if and only if it is affine and the set
of all affine open subschemes $V\sub Y$ subordinate to $\bT$ consists
precisely of all affine open subschemes $V$ for which there exists
an affine open subscheme $U\sub X$ subordinate to $\bW$ such that
$f(V)\sub U$.

 Let $f\:Y\rarrow X$ be a very flat $(\bW,\bT)$\+coaffine morphism,
and let $\P$ be a $\bW$\+locally contraherent cosheaf on~$X$.
 Define a $\bT$\+locally contraherent cosheaf $f^!\P$ on $Y$ in
the following way.
 Let $V\sub Y$ be an affine open subscheme subordinate to~$\bT$.
 For any affine open subscheme $U\sub X$ subordinate to $\bW$ 
and such that $f(V)\sub U$, we set $(f^!\P)[V]_U =
\Hom_{\O_X(U)}(\O_Y(V),\P[U])$.

 The $\O_Y(V)$\+module $(f^!\P)[V]_U$ is contraadjusted by
Lemma~\ref{very-scalars-veryflat-case}(a).
 The contraherence isomorphism $(f^!\P)[V']_U\simeq
\Hom_{\O_Y(V)}(\O_Y(V'),(f^!\P)[V]_U)$ clearly holds for
any affine open subscheme $V'\sub V$.
 The ($\bW$\+local) contraherence condition on $\P$ implies a natural
isomorphism  of $\O_Y(V)$\+modules $(f^!\P)[V]_{U'}\simeq
(f^!\P)[V]_{U''}$ for any embedded affine open subschemes
$U'\sub U''$ in $X$ subordinate to $\bW$ and containing $f(V)$
(cf.\ Section~\ref{contra-direct-inverse}).
 It remains to construct such an isomorphism for any two (not
necessarily embedded) affine open subschemes $U'$, $U''\sub X$.

 The case of a semi-separated scheme $X$ is clear.
 In the general case, let $U'\cap U''=\bigcup_\alpha U_\alpha$ be
an affine open covering of the intersection.
 Since $V$ is quasi-compact, the image $f(V)$ is covered by
a finite subset of the affine open schemes~$U_\alpha$.
 Let $V_\alpha$ denote the preimages of $U_\alpha$ with respect to
the morphism $V\rarrow U'\cap U''$; then $V=\bigcup_\alpha V_\alpha$
is an affine open covering of the affine scheme~$V$.

 The restrictions $f_{U'}\:V\rarrow U'$ and $f_{U''}\:V\rarrow U''$
of the morphism~$f$ are very flat morphisms of affine schemes,
while the restrictions $\P|_{U'}$ and $\P|_{U''}$ of the cosheaf $\P$
are contraherent cosheaves on $U'$ and~$U''$.
 Consider the contraherent cosheaves $f_{U'}^!\P|_{U'}$ and
$f_{U''}^!\P|_{U''}$ on~$V$ (as defined
in Section~\ref{contra-direct-inverse}).
 Their cosection modules $(f_{U'}^!\P|_{U'})[V_\alpha]$ and
$(f_{U''}^!\P|_{U''})[V_\alpha]$ are naturally isomorphic for
all~$\alpha$, since $f(V_\alpha)\sub U_\alpha\sub U'\cap U''$.
 Similarly, there are natural isomorphisms
$(f_{U'}^!\P|_{U'})[V_\alpha\cap V_\beta]\simeq
(f_{U''}^!\P|_{U''})[V_\alpha\cap V_\beta]$ forming commutative
diagrams with the corestrictions from
$V_\alpha\cap V_\beta$ to $V_\alpha$ and $V_\beta$, since
$f(V_\alpha\cap V_\beta)\sub U_\alpha\cap U_\beta$ and
the intersections $U_\alpha\cap U_\beta$ are affine schemes.
 
 Now the cosheaf axiom~\eqref{cosheaf-definition} for contraherent
cosheaves $f_{U'}^!\P|_{U'}$ and $f_{U''}^!\P|_{U''}$ and
the covering $V=\bigcup_\alpha V_\alpha$ provides the desired
isomorphism between the $\O_Y(V)$\+modules $(f^!\P)[V]_{U'} =
(f_{U'}^!\P|_{U'})[V]$ and $(f^!\P)[V]_{U''} =
(f_{U''}^!\P|_{U''})[V]$.
 One can easily see that such isomorphisms form a commutative
diagram for any three affine open subschemes $U'$, $U''$, $U'''
\sub X$ containing $f(V)$.

 The $\bT$\+locally contraherent cosheaf $f^!\P$ on $Y$ is
constructed.
 We have obtained an exact functor of inverse image
$f^!\:X\lcth_\bW\lrarrow Y\lcth_\bT$.

 Let $f\:Y\rarrow X$ be a flat $(\bW,\bT)$\+coaffine morphism of
schemes, and let $\P$ be a locally cotorsion $\bW$\+locally
contraherent cosheaf on~$X$.
 Then the same procedure as above defines a locally cotorsion
$\bT$\+locally contraherent cosheaf $f^!\P$ on~$Y$.
 So we obtain an exact functor $f^!\:X\lcth_\bW^\lct
\rarrow Y\lcth_\bT^\lct$.
 Finally, for any $(\bW,\bT)$\+coaffine morphism of schemes
$f\:Y\rarrow X$ and any locally injective $\bW$\+locally
contraherent cosheaf $\gJ$ on $X$ the same rule defines a locally
injective $\bT$\+locally contraherent cosheaf $f^!\gJ$ on~$Y$.
 We obtain an exact functor of inverse image $f^!\:X\lcth_\bW^\lin
\rarrow Y\lcth_\bT^\lin$.

 Passing to the inductive limits of exact categories with respect to
the refinements of coverings and taking into account the above remark
about $(\bW,\bT)$\+coaffine morphisms, we obtain an exact functor of
inverse image $f^!\:X\lcth\rarrow Y\lcth$ for any very flat morphism
of schemes $f\:Y\rarrow X$.

 For an open embedding $j\:Y\rarrow X$, the inverse image~$j^!$
coincides with the restriction functor $\P\mpsto\P|_Y$ on
the locally contraherent cosheaves $\P$ on~$X$.
 For a flat morphism~$f$, we have an exact functor $f^!\:X\lcth^\lct
\rarrow Y\lcth^\lct$, and for an arbitrary morphism of schemes
$f\:Y\rarrow Y$ there is an exact functor of inverse image of
locally injective locally contraherent cosheaves
$f^!\:X\lcth^\lin\rarrow Y\lcth^\lin$.

 If $f\:Y\rarrow X$ is a $(\bW,\bT)$\+affine morphism and $\Q$ is
a locally cotorsion $\bT$\+locally contraherent cosheaf on $Y$,
then $f_!\Q$ is a locally cotorsion $\bW$\+locally contraherent
cosheaf on~$X$.
 So the direct image functor~$f_!$ restricts to an exact functor
$f_!\:Y\lcth_\bT^\lct\rarrow X\lcth_\bW^\lct$.
 If $f$ is a flat $(\bW,\bT)$\+affine morphism and $\gI$ is
a locally injective $\bT$\+locally contraherent cosheaf on $Y$,
then $f_!\gI$ is a locally injective $\bW$\+locally contraherent
cosheaf on~$X$.
 Hence in this case the direct image also restricts to an exact
functor $f_!\:Y\lcth_\bT^\lin\rarrow X\lcth_\bW^\lin$.

 Let $f\:Y\rarrow X$ be a $(\bW,\bT)$\+affine $(\bW,\bT)$\+coaffine
morphism.
 Then for any $\bT$\+locally contraherent cosheaf $\Q$ on $Y$ and
locally injective $\bW$\+locally contraherent cosheaf $\gJ$ on $X$
there is an adjunction isomorphism
\begin{equation} \label{direct-inverse-lin-adjunction}
 \Hom^X(f_!\Q,\gJ)\simeq\Hom^Y(\Q,f^!\gJ),
\end{equation}
where $\Hom^X$ and $\Hom^Y$ denote the abelian groups of morphisms
in the categories of locally contraherent cosheaves on $X$ and~$Y$.

 If the morphism~$f$ is, in addition, flat, then the isomorphism
\begin{equation} \label{direct-inverse-lct-adjunction}
 \Hom^X(f_!\Q,\P)\simeq\Hom^Y(\Q,f^!\P)
\end{equation}
holds for any $\bT$\+locally contraherent cosheaf $\Q$ on $Y$ and
locally cotorsion $\bW$\+locally contraherent cosheaf $\P$ on~$X$.
 In particular, $f_!$ and $f^!$ form an adjoint pair of functors
between the exact categories $Y\lcth_\bT^\lct$ and $X\lcth_\bW^\lct$.
 When the morphism~$f$ is very flat, the functor $f^!\:
X\lcth_\bW\rarrow Y\lcth_\bT$ is right adjoint to the functor
$f_!\:Y\lcth_\bT\rarrow X\lcth_\bW$.

 Most generally, there is an adjunction isomorphism
\begin{equation} \label{direct-inverse-cosheaf-lin-adjunction}
 \Hom^{\O_X}(f_!\Q,\gJ)\simeq\Hom^{\O_Y}(\Q,f^!\gJ)
\end{equation}
for any morphism of schemes~$f$, a cosheaf of $\O_Y$\+modules $\Q$,
and a locally injective locally contraherent cosheaf $\gJ$ on~$X$.
 Similarly, there is an isomorphism
\begin{equation} \label{direct-inverse-cosheaf-lct-adjunction}
 \Hom^{\O_X}(f_!\Q,\P)\simeq\Hom^{\O_Y}(\Q,f^!\P)
\end{equation}
for any flat morphism~$f$, a cosheaf of $\O_Y$\+modules $\Q$, and
a locally cotorsion locally contraherent cosheaf $\P$ on $X$, and also
for a very flat morphism~$f$, a cosheaf of $\O_Y$\+modules $\Q$, and
a locally contraherent cosheaf $\P$ on~$X$.

 In all the mentioned cases, both the abelian groups $\Hom^X(f_!\Q,\P)$
or $\Hom^{\O_X}(f_!\Q,\allowbreak\P)$ (etc.)\ and $\Hom^Y(\Q,f^!\P)$ or
$\Hom^{\O_Y}(\Q,f^!\P)$ (etc.)\ are identified with
the group of all the compatible collections of homomorphisms of
$\O_X(U)$\+modules $\Q[V]\rarrow\P[U]$ defined for all affine open
subschemes $U\sub X$ and $V\sub Y$ subordinate to, respectively,
$\bW$ and $\bT$, and such that $f(V)\sub U$.
 In other words, the functor~$f^!$ is right adjoint to
the functor~$f_!$ ``wherever the former functor is defined''.

 All the functors between exact categories of locally contraherent
cosheaves constructed above are exact.
 The functors acting between exact categories of locally contraherent
cosheaves with respect to fixed open coverings $\bW$ or $\bT$
specifying the extension of the locality preserve infinite products.
 The functor of direct image of cosheaves of $\O$\+modules~$f_!$
preserves infinite products whenever a morphism of schemes
$f\:Y\rarrow X$ is quasi-compact and quasi-separated.

\medskip

 Let $f\:Y\rarrow X$ be a morphism of schemes and $j\:U\rarrow X$ be
an open embedding.
 Set $V=U\times_XY$, and denote by $j'\:V\rarrow Y$ and $f'\:V\rarrow U$
the natural morphisms.
 Then for any cosheaf of $\O_Y$\+modules $\Q$, there is a natural
isomorphism of cosheaves of $\O_U$\+modules
$(f_!\Q)|_U\simeq f'_!(\Q|_V)$.

 In particular, suppose $f$~is a $(\bW,\bT)$\+affine morphism.
 Define the open coverings $\bW|_U$ and $\bT|_V$ as the collections
of all intersections of the open subsets $W\in\bW$ and $T\in\bT$ with
$U$ and $V$, respectively.
 Then $f'$~is a $(\bW|_U,\bT|_V$)\+affine morphism.
 For any $\bT$\+locally contraherent cosheaf $\Q$ on $Y$
there is a natural isomorphism of $\bW|_U$\+locally contraherent
cosheaves $j^!f_!\Q\simeq f'_!j'{}^!\Q$ on~$U$.

 More generally, let $f\:Y\rarrow X$ and $g\:x\rarrow X$ be morphisms
of schemes.
 Set $y=x\times_X Y$; let $f'\:y\rarrow x$ and $g'\:y\rarrow Y$ be
the natural morphisms.
 Let $\bW$, $\bT$, and $\bw$ be open coverings of, respectively,
$X$, $Y$, and~$x$ such that the morphism~$f$ is $(\bW,\bT)$\+affine,
while the morphism $g$ is $(\bW,\bw)$\+coaffine.

 Define two coverings $\bt'$ and~$\bt''$ of the scheme~$y$ by the rules
that $\bt'$~consists of all the full preimages of affine open subschemes
in~$x$ subordinate to~$\bw$, while $\bt''$~is the collection of all
the full preimages of affine open subschemes in $Y$ subordinate to~$\bT$.
 One can easily see that the covering~$\bt'$ is subordinate to~$\bt''$.
 Let $\bt$ be any covering of~$y$ such that $\bt'$ is subordinate
to~$\bt$ and $\bt$ is subordinate to~$\bt''$.
 Then the former condition guarantees that the morphism~$f'$ is
$(\bw,\bt)$\+affine, while under the latter condition the morphism~$g'$
is $(\bT,\bt)$\+coaffine.

 Assume that the morphisms~$g$ and~$g'$ are very flat.
 Then for any $\bT$\+locally contraherent cosheaf $\P$ on $Y$ there is
a natural isomorphism $g^!f_!\P\simeq f'_!g'{}^!\P$ of $\bw$\+locally
contraherent cosheaves on~$x$.

 Alternatively, assume that the morphism~$g$ is flat.
 Then for any locally cotorsion $\bT$\+locally contraherent cosheaf
$\P$ on $Y$ there is a natural isomorphism $g^!f_!\P\simeq f'_!g'{}^!\P$
of locally cotorsion $\bw$\+locally contraherent cosheaves on~$x$.
 
 As a third alternative, assume that the morphism~$f$ is flat
(while $g$~may be arbitrary).
 Then for any locally injective $\bT$\+locally contraherent cosheaf $\gJ$
on $Y$ there is a natural isomorphism $g^!f_!\gJ\simeq f'_!g'{}^!\gJ$ of
locally injective $\bw$\+locally contraherent cosheaves on~$x$.

 All these isomorphisms of locally contraherent cosheaves are constructed
using the natural isomorphism of $r$\+modules
$\Hom_R(r,P)\simeq\Hom_S(S\ot_R r\;P)$ for any commutative ring
homomorphisms $R\rarrow S$ and $R\rarrow r$, and any $S$\+module~$P$.
 More specifically, given an affine open subscheme $u\sub x$
subordinate to~$\bw$, pick an affine open subscheme $U\sub X$
subordinate to $\bW$ such that $g(u)\sub U$.
 Then $V=f^{-1}(U)$ is an affine open subscheme in $Y$ subordinate
to~$\bT$, and $v=u\times_UV=f'{}^{-1}(u)$ is an affine open subscheme
in~$y$ subordinate to~$\bt$.
 Take $r=\O_x(u)$, \ $R=\O_X(U)$, and $S=\O_Y(V)$; then
$S\ot_Rr=\O_y(v)$; and consider $P=\P[V]$.

 In other words, the direct images of $\bT$\+locally contraherent
cosheaves under $(\bW,\bT)$\+affine morphisms commute with the inverse
images in those base change situations when all the functors involved
are defined.

\medskip

 The following particular case will be important for us.
 Let $\bW$ be an open covering of a scheme $X$ and $j\:Y\rarrow X$ be
an affine open embedding subordinate to $\bW$ (i.~e., $Y$ is contained
in one of the open subsets of $X$ belonging to~$\bW$).
 Then one can endow the scheme $Y$ with the open covering
$\bT=\{Y\}$ consisting of the only open subset $Y\sub Y$.
 This makes the embedding~$j$ both $(\bW,\bT)$\+affine and
$(\bW,\bT)$\+coaffine.
 Also, the morphism~$j$, being an open embedding, is very flat.

 Therefore, the inverse and direct images $j^!$ and $j_!$ form a pair
of adjoint exact functors between the exact category
$X\lcth_\bW$ of $\bW$\+locally contraherent cosheaves on $X$
and the exact category $Y\ctrh$ of contraherent cosheaves on~$Y$.
 Moreover, the image of the functor~$j_!$ is contained in
the full exact subcategory of (globally) contraherent cosheaves
$X\ctrh\sub X\lcth_\bW$.
 Both the functors preserve the subcategories of locally cotorsion and
locally injective cosheaves.

 Now let $\bW$ be an open covering of a quasi-compact semi-separated
scheme $X$ and let $X=\bigcup_{\alpha=1}^N U_\alpha$ be a finite affine
covering of $X$ subordinate to~$\bW$.
 Denote by $j_{\alpha_1,\dotsc,\alpha_k}$ the open embeddings
$U_{\alpha_1}\cap\dotsb\cap U_{\alpha_k}\rarrow X$.
 Then for any $\bW$\+locally contraherent cosheaf $\P$ on $X$
the cosheaf \v Cech sequence (cf.~\eqref{cech-quasi})
\begin{multline} \label{contraherent-cech} \textstyle
 0 \lrarrow j_{1,\dotsc,N}{}_!j_{1,\dotsc,N}^!\P \lrarrow \dotsb \\
 \textstyle\lrarrow \bigoplus_{\alpha<\beta}j_{\alpha,\beta}{}_!
 j_{\alpha,\beta}^!\P \lrarrow\bigoplus_\alpha j_\alpha{}_!j_\alpha^!\P
 \lrarrow \P\lrarrow0
\end{multline}
is exact in the exact category of $\bW$\+locally contraherent
cosheaves on~$X$.
 Indeed, the corresponding sequence of cosections over every
affine open subscheme $U\sub X$ subordinate to $\bW$ is exact
by Lemma~\ref{very-open-covering}(b).
 We have constructed a finite resolution of a $\bW$\+locally
contraherent cosheaf $\P$ by contraherent cosheaves.
 When $\P$ is a locally cotorsion or locally injective $\bW$\+locally
contraherent cosheaf, the sequence~\eqref{contraherent-cech} is
exact in the category $X\lcth_\bW^\lct$ or $X\lcth_\bW^\lin$,
respectively.

\subsection{Coflasque contraherent cosheaves}  \label{coflasque}
 Let $X$ be a topological space and $\gF$ be a cosheaf of abelian
groups on~$X$.
 A cosheaf $\gF$ is called \emph{coflasque} if for all open subsets
$V\sub U\sub X$ the corestriction map $\gF[V]\rarrow\gF[U]$ is
injective.
 Obviously, the class of coflasque cosheaves of abelian groups is
preserved by the restrictions to open subsets and the direct images
with respect to continuous maps.

 In this section and below, we denote by $H^*(X,{-})$ the functor
of cohomology of sheaves of abelian groups on~$X$.
 When $X$ is a ringed space, the functor $H^*(X,{-})$ agrees with
the right derived functor of global sections computed in
the category of sheaves of $\O_X$\+modules (because injective
sheaves of $\O_X$\+modules are flasque).

\begin{lem}  \label{coflasque-cosheaves}
 Let $X=\bigcup_\alpha U_\alpha$ be an open covering.  Then \par
\textup{(a)} a cosheaf\/ $\gF$ on $X$ is coflasque if and only if
its restriction $\gF|_{U_\alpha}$ to each open subset
$U_\alpha$ is coflasque; \par
\textup{(b)} if the cosheaf\/ $\gF$ is coflasque, then the \v Cech
complex~\textup{\eqref{cech-homol}} has no higher homology
groups, $H_{>0}C_\bu(\{U_\alpha\},\gF)=0$.
\end{lem}

\begin{proof}
 One can either check these assertions directly or deduce them from
the similar results for flasque sheaves of abelian groups using
the construction of the sheaf $U\mpsto\Hom_\boZ(\gF[U],I)$ from
the proof of Theorem~\ref{cosheaf-base-thm}.
 Here $I$ is an injective abelian group; clearly, the sheaf so
obtained is flasque for all~$I$ (or specifically for $I=\boQ/\boZ$)
if and only if the original cosheaf $\gF$ was coflasque.
 The sheaf (dual) versions of assertions~(a\+b) are
well-known~\cite[Section~II.3.1 and Th\'eor\`eme~II.5.2.3(a)]{God}.
\end{proof}

\begin{cor}  \label{coflasque-contraherent}
 Let\/ $X$ be a scheme with an open covering\/ $\bW$ and\/ $\gF$
be a\/ $\bW$\+locally contraherent cosheaf on~$X$.
 Suppose that the cosheaf\/ $\gF$ is coflasque.
 Then\/ $\gF$ is a (globally) contraherent cosheaf on~$X$.
\end{cor}

\begin{proof}
 Follows from Lemmas~\ref{global-contraherence-criterion}
and~\ref{coflasque-cosheaves}(b).
\end{proof}

 Assume that the topological space $X$ has a base of the topology
consisting of quasi-compact open subsets.
 Then one has $\gF[Y]\simeq\varinjlim_{U\sub Y}\gF[U]$ for any
cosheaf of abelian groups $\gF$ on $X$ and any open subset $Y\sub X$,
where the filtered inductive limit is taken over all the quasi-compact
open subsets $U\sub Y$.
 It follows easily that a cosheaf $\gF$ is coflasque if and only if
the corestriction map $\gF[V]\rarrow\gF[U]$ is injective for any
pair of embedded quasi-compact open subsets $V\sub U\sub X$.
 Notice that the dual version of this argument does \emph{not} seem
to work for sheaves, because of nonexactness of filtered projective
limits.

 Moreover, if $X$ is a scheme then it follows from
Lemma~\ref{coflasque-cosheaves}(a) that a cosheaf of abelian groups
$\gF$ on $X$ is coflasque if and only if the corestriction map
$\gF[V]\rarrow\gF[U]$ is injective for any affine open subscheme
$U\sub X$ and quasi-compact open subscheme $V\sub U$.
 Therefore, an infinite product of a family of coflasque
contraherent cosheaves on $X$ is coflasque.
 The following counterexample shows, however, that coflasqueness of
contraherent cosheaves on schemes cannot be checked on the pairs of
embedded affine open subschemes.

\begin{ex}
 Let $X$ be a Noetherian scheme.
 It is well-known that any injective quasi-coherent sheaf $\J$ on $X$
is a flasque sheaf of abelian groups~\cite[Proposition~III.3.4
and Corollary~III.3.6]{HarAG}; indeed, any injective quasi-coherent
sheaf on a locally Noetherian scheme is injective as a sheaf of
$\O$\+modules~\cite[Theorem~II.7.18]{Har}, which clearly implies
flasqueness.
 Let $\F$ be a quasi-coherent sheaf on $X$ and $\F\rarrow\J$ be
an injective morphism.
 Then one has $(\J/\F)(U)\simeq\J(U)/\F(U)$ for any affine open
subscheme $U\sub X$, so the map $(\J/\F)(U)\rarrow(\J/\F)(V)$ is
surjective for any pair of embedded affine open subschemes
$V\sub U$.
 On the other hand, if the quotient sheaf $\J/\F$ were flasque, one
would have $H^{i+1}(X,\F)\simeq H^i(X,\J/\F)=0$ for $i\ge1$, which
is clearly not the case in general.

 Now let $X$ be a scheme over an affine scheme $\Spec R$.
 Let $\M$ be a quasi-coherent sheaf on $X$ and $J$ be an injective
$R$\+module.
 Then for any quasi-compact quasi-separated open subscheme $Y\sub X$
one has $\Cohom_R(\M,J)[Y]\simeq\Hom_R(\M(Y),J)$
(see Section~\ref{cohom-subsection}).
 For a locally Noetherian scheme $X$, it follows that the cosheaf
$\Cohom_R(\M,J)$ is coflasque for all $J$ if and only if the sheaf
$\M$ is flasque (cf.\ Lemma~\ref{co-flasque-preservation}(b) below).
\end{ex}

\begin{cor}  \label{coflasque-acyclic}
 Let $X$ be a scheme with an open covering\/ $\bW$ and\/
$0\rarrow\P\rarrow\Q\rarrow\R\rarrow0$ be a short exact sequence in\/
$\O_X\cosh_\bW$ (e.~g., a short exact sequence of\/ $\bW$\+locally
contraherent cosheaves on $X$).
 Then \par
\textup{(a)} the cosheaf\/ $\Q$ is coflasque whenever both
the cosheaves\/ $\P$ and $\R$ are; \par
\textup{(b)} the cosheaf\/ $\P$ is cosflasque whenever both
the cosheaves\/ $\Q$ and $\R$ are; \par
\textup{(c)} if the cosheaf\/ $\R$ is coflasque, then the short
sequence $0\rarrow\P[Y]\rarrow\Q[Y]\rarrow \R[Y]\rarrow0$
is exact for any open subscheme $Y\sub X$.
\end{cor}

\begin{proof}
 The assertion actually holds for any short exact sequence in
the exact category of cosheaves of abelian groups on
a scheme $X$ with the exact category structure related to the base
of affine open subschemes subordinate to a covering~$\bW$.

 Part~(c): let us first consider the case of a semi-separated open
subscheme~$Y$.
 Pick an affine open covering $Y=\bigcup_\alpha U_\alpha$ subordinate
to~$\bW$.
 The intersection of any nonempty finite subset of $U_\alpha$ being
also an affine open subscheme in $X$, the short sequence of complexes
$0\rarrow C_\bu(\{U_\alpha\},\P|_Y)\rarrow C_\bu(\{U_\alpha\},\Q|_Y)
\rarrow C_\bu(\{U_\alpha\},\R|_Y)\rarrow0$ is exact.
 Recall that $\gF[Y]\simeq H_0C_\bu(\{U_\alpha\},\gF|_Y)$ for
any cosheaf of abelian groups $\gF$ on~$X$.
 By Lemma~\ref{coflasque-cosheaves}, one has
$H_{>0}C_\bu(\{U_\alpha\},\R|_Y)=0$, and it follows that the short
sequence $0\rarrow\P[Y]\rarrow\Q[Y]\rarrow\R[Y]\rarrow0$ is exact.

 Now for any open subscheme $Y\sub X$, pick an open covering
$Y=\bigcup_\alpha U_\alpha$ of $Y$ by semi-separated open
subschemes~$U_\alpha$.
 The intersection of any nonempty finite subset of $U_\alpha$ being
semi-separated, the short sequence of \v Cech complexes for
the covering $U_\alpha$ and the short exact sequence of cosheaves
$0\rarrow\P\rarrow\Q\rarrow\R\rarrow0$ is exact according to
the above, and the same argument concludes the proof.

 Parts~(a\+b): let $U\sub V\sub X$ be any pair of embedded open
subschemes.
 Assuming that the cosheaf $\R$ is coflasque, according to part~(c)
we have short exact sequences of the modules of cosections
$0\rarrow\P[V]\rarrow\Q[V]\rarrow\R[V]\rarrow0$ and
$0\rarrow\P[U]\rarrow\Q[U]\rarrow\R[U]\rarrow0$.
 If the morphism from the former sequence to the latter is injective
on the rightmost terms, then it is injective on the middle terms if
and only if it is injective on the leftmost terms.
\end{proof}

\begin{rem}
 Let $U$ be an affine scheme.
 Then a cosheaf of abelian groups $\gF$ on $U$ is coflasque if and only
if the following three conditions hold:
\begin{enumerate} \renewcommand{\theenumi}{\roman{enumi}}
\item for any affine open subscheme $V\sub U$, the corestriction map
$\gF[V]\rarrow\gF[U]$ is injective;
\item for any two affine open subschemes $V'$, $V''\sub U$,
the image of the map $\gF[V'\cap V'']\rarrow\gF[U]$ is equal to
the intersection of the images of the maps $\gF[V']\rarrow\gF[U]$
and $\gF[V'']\rarrow\gF[U]$;
\item for any finite collection of affine open subschemes
$V_\alpha\sub U$, the images of the maps $\gF[V_\alpha]\rarrow\gF[U]$
generate a distributive sublattice in the lattice of all subgroups
of the abelian group $\gF[U]$.
\end{enumerate}
 Indeed, the condition~(i) being assumed, in view of the exact sequence
$\gF[V'\cap V'']\rarrow\gF[V']\oplus\gF[V'']\rarrow\gF[V'\cup V'']
\rarrow0$ the condition~(ii) is equivalent to injectivity of
the corestriction map $\gF[V'\cup V'']\rarrow\gF[U]$.

 Furthermore, set $W=\bigcup_\alpha V_\alpha$.
 The condition~(iii) for any proper subcollection of $V_\alpha$ and
the conditions~(i\+ii) being assumed, the condition~(iii)
for the collection $\{V_\alpha\}$ becomes equivalent to vanishing
of the higher homology of the \v Cech complex
$C_\bu(\{V_\alpha\},\gF|_W)$ together with
injectivity of the map $\gF[W]\rarrow\gF[U]$.
 Indeed, the complex $C_\bu(\{V_\alpha\},\gF|_W)\rarrow\gF[U]$
can be identified with (the abelian group version of)
the ``cobar complex'' $B^\bu(\gF[U];\>\gF[V_\alpha])$
from~\cite[Proposition~7.2(c*) of Chapter~1]{PP}
(see also~\cite[Lemma~11.4.3.1(c*)]{Psemi}
or~\cite[Lemma~2.18(c*)]{Prel}).
\end{rem}

 Given a scheme $X$, a cosheaf of $\O_X$\+modules $\gG$ is said to be
\emph{flat} if the $\O_X(U)$\+mod\-ule $\gG[U]$ is flat for every
affine open subscheme $U\sub X$.
 A more detailed discussion of this definition can be found in
Section~\ref{contraherent-tensor} below.

\begin{lem}  \label{co-flasque-preservation}
\textup{(a)} Let\/ $\M$ be a flasque quasi-coherent sheaf and\/ $\G$
be a flat quasi-coherent sheaf on a locally Noetherian scheme~$X$.
 Then the quasi-coherent sheaf\/ $\M\ot_{\O_X}\G$ on $X$ is flasque.
\par
\textup{(b)} Let\/ $\M$ be a flasque quasi-coherent sheaf and\/
$\gJ$ be a locally injective contraherent cosheaf on a scheme~$X$.
 Then the contraherent cosheaf\/ $\Cohom_X(\M,\gJ)$ on $X$ is
coflasque. \par
\textup{(c)} Let\/ $\M$ be a flasque quasi-coherent sheaf and\/ $\J$ be
an injective quasi-coherent sheaf on an affine scheme~$U$.
 Then the contraherent cosheaf\/ $\fHom_U(\M,\J)$ on $U$ is coflasque.
\par
\textup{(d)} Let\/ $\M$ be a flasque quasi-coherent sheaf and\/
$\gG$ be a flat cosheaf of\/ $\O_U$\+modules on a Noetherian
affine scheme~$U$.
 Then the quasi-coherent sheaf\/ $\M\ocn_U\gG$ on $U$ is flasque.
\end{lem}

\begin{proof}
 The (co)flasqueness of (co)sheaves being a local property, it suffices
to consider the case of an affine scheme $X=U$ in all assertions~(a\+d).
 Then (c) becomes a restatement of~(b) and (d) a restatement of~(a).
 It also follows that one can extend the assertion~(b) to locally
injective locally contraherent cosheaves~$\gJ$
(see Section~\ref{cohom-loc-der-contrahereable} for the definition
of $\Cohom_X$ in this context).

 To prove part~(a) in the affine case, one notices the isomorphism
$(\M\ot_{\O_U}\G)(V)\simeq\M(V)\ot_{\O(U)}\G(U)$ holding for any
quasi-compact open subscheme $V$ in an affine scheme $U$,
quasi-coherent sheaf $\M$, and a flat quasi-coherent sheaf $\G$ on~$U$.
 More generally, one has $(\M\ot_R G)(Y)\simeq\M(Y)\ot_R G$ for any
quasi-coherent sheaf $\M$ with a right $R$\+module structure on
a quasi-compact quasi-separated scheme $Y$ and a flat left
$R$\+module~$G$.
 Part~(b) follows from the similar isomorphism
$\Cohom_U(\M,\gJ)[V]\simeq\Hom_{\O(U)}(\M(V),\gJ[U])$ holding for
any quasi-coherent sheaf $\M$ and a (locally) injective contraherent
cosheaf $\gJ$ on~$U$.
 More generally, $\Cohom_R(\M,J)[Y]\simeq\Hom_R(\M(Y),J)$ for any
quasi-compact quasi-separated scheme $Y$ over $\Spec R$, quasi-coherent
sheaf $\M$ on $Y$, and an injective $R$\+module~$J$.

 It is helpful to keep in mind that it suffices to check
the coflasqueness, but \emph{not} the flasqueness, on quasi-compact
open subschemes $V$ in affine schemes~$U$.
 That is why we assume the Noetherianity in (a) and~(d).
\end{proof}

\begin{lem}  \label{grothendieck-vanishing}
\textup{(a)} Let $X$ be a Noetherian topological space of
finite Krull dimension~$\le d+1$ (where $d\ge-1$ is an integer),
and let\/ $0\rarrow\E\rarrow\F^0\rarrow\dotsb\rarrow\F^d\rarrow\F
\rarrow0$ be an exact sequence of sheaves of abelian groups on~$X$.
 Suppose that the sheaves\/ $\F^i$ are flasque.
 Then the sheaf\/ $\F$ is flasque. \par
\textup{(b)} Let $X$ be a scheme with an open covering\/ $\bW$ and\/
$0\rarrow\gF\rarrow\gF_d\rarrow\dotsb\rarrow\gF_0\rarrow\gE\rarrow0$
be an exact sequence in\/ $\O_X\cosh_\bW$.
 Suppose that the cosheaves\/ $\gF_i$ are coflasque and the underlying
topological space of the scheme $X$ is Noetherian of finite Krull
dimension\/~$\le d+1$.
 Then the cosheaf\/ $\gF$ is coflasque.
\end{lem}

\begin{proof}
 Part~(a) is a version of Grothendieck's vanishing
theorem~\cite[Th\'eor\`eme~3.6.5]{Toh}; it can be either proved directly
along the lines of Grothendieck's proof, or, making a slightly stronger
assumption that the dimension of $X$ does not exceed~$d$, deduced
from the assersion of Grothendieck's theorem.
 Indeed, let $\G$ denote the image of the morphism of sheaves
$\F^{d-1}\rarrow\F^d$; by Grothendieck's theorem, one has
$H^1(V,\G|_V)\simeq H^{d+1}(V,\E|_V)=0$ for any open subset $V\sub X$.
 Hence the map $\F^d(V)\rarrow\F(V)$ is surjective, and it follows
that the map $\F(X)\rarrow\F(V)$ is surjective, too.

 For a direct proof of part~(a) as stated in the lemma, let us use
the term ``flasque dimension'' for the coresolution dimension of
a sheaf of abelian groups with respect to the coresolving subcategory
of flasque sheaves of abelian groups in the abelian category of
sheaves of abelian groups on a topological space.
 (We refer to Section~\ref{finite-resolutions-subsect} for a general
discussion of the coresolution dimension.)
 We need to show that the flasque dimension of any sheaf of abelian
groups $\E$ on a Noetherian topological space of Krull
dimension~$\le d+1$ does not exceed~$d+1$.
 We will use the properties of the coresolution dimension dual to
the ones listed in Lemma~\ref{fdim-properties}(a,c).

 Notice that the direct images of sheaves with respect to embeddings
of closed subspaces do not increase the flasque dimension, and neither
do the direct limits of sheaves on a Noetherian topological space.
 The argument from the proof of~\cite[Th\'eor\`eme~3.6.5]{Toh} reduces
the question to the case when $X$ is irreducible.
 Then the argument from the proof~\cite[Proposition~3.6.1]{Toh}
reduces the question further to the case when the sheaf of abelian
groups $\E$ is generated by a single section over some (nonempty) open
subset $U\sub X$.
 Replacing $U$ by a smaller nonempty open subset $V\sub U$, we can
assume that the sheaf $\E|_V\simeq A_V$ is a constant sheaf on $V$
corresponding to some (cyclic) abelian group~$A$.

 Denoting by $j\:V\rarrow X$ and $i\:Z=(X\setminus V)\rarrow X$
the open and closed embedding maps, we have a short exact sequence
$0\rarrow j_!j^*\E\rarrow\E\rarrow i_*i^*\E\rarrow0$ of sheaves of
abelian groups on~$X$.
 There is also a similar short exact sequence $0\rarrow j_!A_V\rarrow
A_X \rarrow i_*A_Z\rarrow0$ for the constant sheaf $A_X$ on~$X$.
 By the assumption of induction on~$d$, the flasque dimensions of
$i^*\E$ and $A_Z$ on $Z$ do not exceed~$d$.
 The constant sheaf $A_X$ on $X$ is flasque (since $X$ is irreducible).
 The reference to the dual version of Lemma~\ref{fdim-properties}(a,c)
finishes the proof of part~(a).

 Part~(b): given an injective abelian group~$I$, the sequence
$0\rarrow\Hom_\boZ(\gE[U],I)\rarrow\Hom_\boZ(\gF_0[U],I)\rarrow
\dotsb\rarrow\Hom_\boZ(\gF_d[U],I)\rarrow\Hom_\boZ(\gF[U],I)\rarrow0$
is exact for any affine open subscheme $U\sub X$ subordinate
to~$\bW$.
 Hence the construction assigning to a cosheaf of $\O_X$\+modules $\P$
on $X$ the sheaf $V\mpsto\Hom_\boZ(\P[V],I)$, where $V\sub X$ are
arbitrary open subschemes, transforms our sequence of cosheaves into
an exact sequence of sheaves of $\O_X$\+modules $0\rarrow
\Hom_\boZ(\gE,I)\rarrow\Hom_\boZ(\gF_0,I)\rarrow\dotsb\rarrow
\Hom_\boZ(\gF_d,I)\rarrow\Hom_\boZ(\gF,I)\rarrow0$.
 Now the sheaves $\Hom_\boZ(\gF_i,I)$ are flasque, and by part~(a)
it follows that the sheaf $\Hom_\boZ(\gF,I)$ is.
 Therefore, the cosheaf $\gF$ is coflasque.
 (Cf.\ Section~\ref{finite-resolutions-subsect}.)
\end{proof}

 We have shown, in particular, that coflasque contraherent cosheaves
on a scheme $X$ form a full subcategory of $X\ctrh$ closed under
extensions, kernels of admissible epimorphisms, and infinite products.
 Hence this full subcategory acquires an induced exact category
structure, which we will denote by $X\ctrh_\cfq$.

 Similarly, coflasque locally cotorsion contraherent cosheaves on $X$
form a full subcategory of $X\ctrh^\lct$ closed under extensions,
kernels of admissible epimorphisms, and infinite products.
 We denote the induced exact category structure on this full
subcategory by $X\ctrh^\lct_\cfq$.

 The following corollary addresses a problem pointed out in
Remark~\ref{noncontraherent-direct-image-remark}
and Example~\ref{noncontraherent-direct-image-example}.
 For another approach, see Corollary~\ref{clp-direct} below.

\begin{cor} \label{coflasque-direct}
 Let $f\:Y\rarrow X$ be a quasi-compact quasi-separated morphism of
schemes.
 Then \par
\textup{(a)} the functor of direct image of cosheaves of\/
$\O$\+modules $f_!$ takes coflasque contraherent cosheaves on $Y$
to coflasque contraherent cosheaves on $X$, and induces an exact
functor $f_!\:Y\ctrh_\cfq\rarrow X\ctrh_\cfq$ between these exact
categories; \par
\textup{(b)} the functor of direct image of cosheaves of\/
$\O$\+modules $f_!$ takes coflasque locally cotorsion contraherent
cosheaves on $Y$ to coflasque locally cotorsion contraherent
cosheaves on $X$, and induces an exact functor $f_!\:Y\ctrh^\lct_\cfq
\rarrow X\ctrh^\lct_\cfq$ between these exact categories.
\end{cor}

\begin{proof}
 Since the contraadjustedness/contraherence conditions, local
cotorsion, and exactness of short sequences in $X\ctrh$ or
$X\ctrh^\lct$ only depend on the restrictions to affine open
subschemes, while the coflasqueness is preserved by restrictions to
any open subschemes, we can assume that the scheme $X$ is affine.
 Then the scheme $Y$ is quasi-compact and quasi-separated.

 Let $Y=\bigcup_\alpha V_\alpha$ be a finite affine open covering
and $\gF$ be a coflasque contraherent cosheaf on~$Y$.
 Then, by Lemma~\ref{coflasque-cosheaves}(b), the \v Cech complex
$C_\bu(\{V_\alpha\},\gF)$ is a finite resolution of
the $\O(X)$\+module $\gF[Y]$.
 Let us first consider the case when $Y$ is semi-separated, so
the intersection of any nonempty subset of $V_\alpha$ is affine.

 By Lemma~\ref{very-scalars-always}(a), our resolution consists of
contraadjusted $\O(X)$\+modules.
 When $\gF$ is a locally cotorsion contraherent cosheaf, by
Lemma~\ref{cotors-restrict}(a) this resolution even consists
of cotorsion $\O(X)$\+modules.
 It follows that the $\O(X)$\+module $(f_!\gF)[X]=\gF[Y]$ is
contraadjusted in the former case and cotorsion in the latter one.

 Now let $U\sub X$ be an affine open subscheme.
 Then one has
$$
 \gF[f^{-1}(U)\cap V]\simeq\Hom_{\O_Y(V)}(\O_Y(f^{-1}(U)\cap V)\;\gF[V])
 \simeq\Hom_{\O(X)}(\O_X(U),\gF[V])
$$
for any affine open subscheme $V\sub Y$.
 So the complex $C_\bu(\{f^{-1}(U)\cap V_\alpha\}\;\gF|_{f^{-1}(U)})$
can be obtained by applying the functor $\Hom_{\O(X)}(\O_X(U),{-})$
to the complex $C_\bu(\{V_\alpha\},\gF)$.
 It follows that $(f_!\gF)[U]\simeq\gF[f^{-1}(U)]\simeq
\Hom_{\O(X)}(\O_X(U),\gF[Y])$ and the contraherence condition
holds for $f_!\gF$.

 Finally, let us turn to the general case.
 According to the above, for any semi-separated quasi-compact open
subscheme $V\sub Y$ the $\O(X)$\+module $\gF[V]$ is contraadjusted
(and even cotorsion if $\gF$ is locally cotorsion), and for any
affine open subscheme $U\sub X$ one has 
$\gF[f^{-1}(U)\cap V]\simeq\Hom_{\O(X)}(\O_X(U),\gF[V])$.
 Given that the intersection of any nonempty subset of $V_\alpha$
is separated (being quasi-affine) and quasi-compact, the same
argument as above goes through.
\end{proof}

 The notation $\Delta(X,\P)=\P[X]$ for the $\O(X)$\+module of
global cosections of a (locally contraherent) cosheaf $\P$ on
a scheme $X$ was introduced in Section~\ref{counterex-subsect}.

\begin{cor}  \label{coflasque-complexes-direct-correct}
 Let $f\:Y\rarrow X$ be a quasi-compact quasi-separated morphism of
schemes and\/ $\bT$ be an open covering of~$Y$.
 Then \par
\textup{(a)} for any complex\/ $\F^\bu$ of flasque quasi-coherent
sheaves on $Y$ that is acyclic as a complex over $Y\qcoh$,
the complex $f_*\F^\bu$ of flasque quasi-coherent sheaves on $X$
is acyclic as a complex over $X\qcoh$; \par
\textup{(b)} for any complex\/ $\gF^\bu$ of coflasque contraherent
cosheaves on $Y$ that is acyclic as a complex over $Y\lcth_\bT$,
the complex $f_!\gF^\bu$ of coflasque contraherent cosheaves on $X$
is acyclic as a complex over $X\ctrh$; \par
\textup{(c)} for any complex\/ $\gF^\bu$ of coflasque locally
cotorsion contraherent cosheaves on $Y$ that is acyclic as
a complex over $Y\lcth_\bT^\lct$, the complex $f_!\gF^\bu$ of
coflasque locally cotorsion contraherent cosheaves on $X$ is
acyclic as a complex over $X\ctrh^\lct$.
\end{cor}

\begin{proof}
 Let us prove part~(c), parts~(a\+b) being analogous.
 In view of Lemma~\ref{acyclicity-in-lcth-criterion},
it suffices to consider the case of an affine scheme $X$ and
a quasi-compact quasi-separated scheme~$Y$.
 The proof of part~(b) of that lemma was based on the cotorsion
periodicity theorem (Theorem~\ref{cotorsion-periodicity}), but
the following argument does not even use the cotorsion periodicity.
 We have to show that the complex of cotorsion $\O(X)$\+modules
$\Delta(Y,\gF^\bu)$ is acyclic over $\O(X)\modl^\cot$.
 Considering the case of a semi-separated scheme $Y$ first and
the general case second, one can assume that $Y$ has a finite
open covering $Y=\bigcup_{\alpha=1}^N U_\alpha$ by quasi-compact
quasi-separated schemes $U_\alpha$ subordinate to $\bT$ such that
for any intersection $V$ of a nonempty subset of $U_\alpha$
the complex $\Delta(V,\gF^\bu)$ is acyclic over $\O(X)\modl^\cot$.

 Consider the \v Cech bicomplex $C_\bu(\{U_\alpha\},\gF^\bu)$
of the complex of cosheaves $\gF^\bu$ on the scheme $Y$ with
the open covering~$U_\alpha$.
 There is a natural morphism of bicomplexes
$C_\bu(\{U_\alpha\},\gF^\bu)\rarrow\Delta(Y,\gF^\bu)$; and
for every degree~$i$ the complex $C_\bu(\{U_\alpha\},\gF^i)
\rarrow\Delta(Y,\gF^i)$ is a finite (and uniformly bounded)
acyclic complex of cotorsion $\O(X)$\+modules, essentially by
Lemma~\ref{coflasque-cosheaves}(b) (since the cosheaf $\gF^i$
is coflasque).
 It follows that the induced morphism of total complexes is
a quasi-isomorphism (in fact, a morphism with an absolutely
acyclic cone, in the sense of Section~\ref{derived-second-kind})
of complexes over $\O(X)\modl^\cot$.
 On the other hand, for every~$k$ the complex
$C_k(\{U_\alpha\},\gF^\bu)$ is by assumption acyclic over
$\O(X)\modl^\cot$.
 It follows that the total complex of $C_\bu(\{U_\alpha\},\gF^\bu)$,
and hence also the complex $\Delta(Y,\gF^\bu)$, is acyclic
over $\O(X)\modl^\cot$.
\end{proof}

\subsection{Contrahereable cosheaves and the contraherator}
\label{contrahereable-subsect}
 A cosheaf of $\O_X$\+mod\-ules $\P$ on a scheme $X$ is
said to be \emph{derived contrahereable} if
\begin{enumerate}
\renewcommand{\theenumi}{\roman{enumi}$^\circ$}
\item for any affine open subscheme $U\sub X$ and its finite affine
open covering $U=\bigcup_{\alpha=1}^N U_\alpha$ the homological \v Cech
sequence (cf.~\eqref{cech-homol})
\begin{multline} \label{exactness-condition} \textstyle
 0\lrarrow \P[U_1\cap\dotsb\cap U_N] \lrarrow\dotsb \\ \textstyle
 \lrarrow \bigoplus_{\alpha<\beta} \P[U_\alpha\cap U_\beta]
 \lrarrow \bigoplus_\alpha \P[U_\alpha] \lrarrow \P[U]\lrarrow 0
\end{multline}
is exact; and
\renewcommand{\theenumi}{\roman{enumi}}
\item for any affine open subscheme $U\sub X$, the $\O_X(U)$\+module
$\P[U]$ is contraadjusted.
\end{enumerate}
 We will call~(i$^\circ$) the \emph{exactness condition} and (ii)
the \emph{contraadjustedness condition}.

 Notice that the present contraadjustedness condition~(ii) is
an equivalent restatement of the contraadjustedness condition~(ii) of
Section~\ref{contraherent-definition}, while the exactness
condition~(i$^\circ$) is weaker than the contraherence
condition~(i) of Section~\ref{contraherent-definition}
(provided that the condition~(ii) is assumed).
 Naturally occuring examples of derived contrahereable, but not
(locally) contraherent cosheaves will be considered below in
Section~\ref{contraherent-tensor}.
 The condition~(i$^\circ$) is also weaker than the coflasqueness
condition on a cosheaf of $\O_X$\+modules discussed in
Section~\ref{coflasque}.

 By Remark~\ref{scheme-topology}, the exactness condition~(i$^\circ$)
can be thought of as a strengthening of the cosheaf
property~\eqref{base-cosheaf} of a covariant functor with
an $\O_X$\+module structure on the category of affine open subschemes
of~$X$.
 Any such functor satisfying~(i$^\circ$) can be extended to
a cosheaf of $\O_X$\+modules in a unique way.

 Let $\bW$ be an open covering of a scheme~$X$.
 A cosheaf of $\O_X$\+modules $\P$ on $X$ is called \emph{$\bW$\+locally
derived contrahereable} if its restrictions $\P|_W$ to all the open
subschemes $W\in\bW$ are derived contrahereable on~$W$.
 In other words, this means that the conditions (i$^\circ$) and~(ii)
must hold for all the affine open subschemes $U\sub X$ subordinate
to~$\bW$.
 A cosheaf of $\O_X$\+modules is called locally derived contrahereable
if it is $\bW$\+locally derived contrahereable for some open
covering~$\bW$.
 Any contraherent cosheaf is derived contrahereable, and any
$\bW$\+locally contraherent cosheaf is $\bW$\+locally derived
contrahereable.
 Conversely, according to Lemma~\ref{global-contraherence-criterion},
if a $\bW$\+locally contraherent cosheaf is derived contrahereable,
then it is contraherent.

 A $\bW$\+locally derived contrahereable cosheaf $\P$ on $X$ is called
\emph{locally cotorsion} (respectively, \emph{locally injective}) if
the $\O_X(U)$\+module $\P[U]$ is cotorsion (resp., injective) for any
affine open subscheme $U\sub X$ subordinate to~$\bW$.
 A locally derived contrahereable cosheaf $\P$ is locally
cotorsion (resp., locally injective) if and only if
the $\O_X(U)$\+module $\P[U]$ is cotorsion (resp., injective)
for every affine open subscheme $U\sub X$ such that the cosheaf
$\P|_U$ is derived contrahereable.

 In the exact category of cosheaves of $\O_X$\+modules $\O_X\cosh_\bW$
defined in Section~\ref{locally-contraherent}, the $\bW$\+locally derived
contrahereable cosheaves form an exact subcategory closed under
extensions, infinite products, and cokernels of admissible monomorphisms.
 Hence there is the induced exact category structure on the category of
$\bW$\+locally derived contrahereable cosheaves on~$X$.

\medskip

 Let $U$ be an affine scheme and $\Q$ be a derived contrahereable
cosheaf on~$U$.
 The \emph{contraherator} $\Cr\Q$ of the cosheaf $\Q$ is defined
in this simplest case as the contraherent cosheaf on $U$ corresponding
to the contraadjusted $\O(U)$\+module $\Q[U]$, that is
$\Cr\Q=\widecheck{\Q[U]}$.
 There is a natural morphism of derived contrahereable cosheaves
$\Q\rarrow\Cr\Q$ on~$U$ (see Lemma~\ref{hom-into-contraherent}).
 For any affine open subscheme $V\sub U$ there is a natural
morphism of contraherent cosheaves $\Cr(\Q|_V)\rarrow(\Cr\Q)|_V$ on~$V$.
 Our next goal is to extend this construction to an appropriate class
of cosheaves of $\O_X$\+modules on quasi-compact semi-separated
schemes~$X$ (cf.~\cite[Appendix~B]{TT}).

 Let $X$ be such a scheme with an open covering $\bW$ and a finite
affine open covering $X=\bigcup_{\alpha=1}^N U_\alpha$ subordinate
to~$\bW$.
 The \emph{contraherator complex} $\Cr_\bu(\{U_\alpha\},\P)$ of
a $\bW$\+locally derived contrahereable cosheaf $\P$ on $X$ is a finite
\v Cech complex of (globally) contraherent cosheaves on $X$ of the form
\begin{multline}  \label{contraherator-complex}  \textstyle
 0 \lrarrow j_{1,\dotsc,N}{}_!\Cr (\P|_{U_1\cap\dotsb\cap U_N})
 \lrarrow \dotsb \\ \textstyle
 \lrarrow \bigoplus_{\alpha<\beta}j_{\alpha,\beta}{}_!
 \Cr(\P|_{U_\alpha\cap U_\beta})\lrarrow
 \bigoplus_\alpha j_\alpha{}_!\Cr(\P|_{U_\alpha}),
\end{multline}
where $j_{\alpha_1,\dotsc,\alpha_k}$ is the open embedding
$U_{\alpha_1}\cap\dotsb\cap U_{\alpha_k}\rarrow X$ and
the notation $\Cr(\Q)$ was explained above.
 The differentials in this complex are constructed using the adjunction
of the direct and inverse images of contraherent cosheaves and
the above morphisms $\Cr(\Q|_V)\rarrow(\Cr\Q)|_V$.

\begin{lem}  \label{contraherator-well-defined}
 Let\/ $\P$ be a\/ $\bW$\+locally derived contrahereable cosheaf
on a quasi-compact semi-separated scheme~$X$.
 Then the object of the bounded derived category\/ $\sD^\b(X\ctrh)$ of
the exact category of contraherent cosheaves on $X$ represented by
the complex\/ $\Cr_\bu(\{U_\alpha\},\P)$ does not depend on
the choice of a finite affine open covering\/ $\{U_\alpha\}$ of $X$
subordinate to the covering\/~$\bW$.
\end{lem}

\begin{proof}
 Let us adjoin another affine open subscheme $V\sub X$, subordinate
to~$\bW$, to the covering~$\{U_\alpha\}$.
 Then the complex $\Cr_\bu(\{U_\alpha\},\P)$ embeds into the complex
$\Cr_\bu(\{V,U_\alpha\},\P)$ by a termwise split morphism of complexes
with the cokernel isomorphic to the complex
$k_!\Cr_\bu(\{V\cap U_\alpha\},\P|_V)\rarrow k_!\Cr(\P|_V)$,
where $k\:V\rarrow X$ is the open embedding morphism.

 The complex of contraherent cosheaves
$\Cr_\bu(\{V\cap U_\alpha\},\P|_V)\rarrow\Cr(\P|_V)$ on $V$
corresponds to the acyclic complex of contraadjusted
$\O(V)$\+modules~\eqref{exactness-condition} for the covering
of an affine open subscheme $V\sub X$ by the affine open subschemes
$V\cap U_\alpha$.
 Hence the cokernel of the admissible monomorphism of complexes
$\Cr_\bu(\{U_\alpha\},\P)\rarrow \Cr_\bu(\{V,U_\alpha\},\P)$
is an acyclic complex of contraherent cosheaves on~$X$.

 Now, given two finite affine open coverings $X=\bigcup_\alpha U_\alpha
=\bigcup_\beta V_\beta$ subordinate to $\bW$, one compares both
the complexes $\Cr_\bu(\{U_\alpha\},\P)$ and $\Cr_\bu(\{V_\beta\},\P)$
with the complex $\Cr_\bu(\{U_\alpha,V_\beta\},\P)$ corresponding
to the union of the two coverings $\{U_\alpha,V_\beta\}$.
\end{proof}

 A $\bW$\+locally derived contrahereable cosheaf $\P$ on
a quasi-compact semi-separated scheme $X$ is called \emph{$\bW$\+locally
contrahereable} if the complex $\Cr_\bu(\{U_\alpha\},\P)$ for
some particular (or equivalently, for any) finite affine open
covering $X=\bigcup_\alpha U_\alpha$ is isomorphic in
$\sD^\b(X\lcth_\bW)$ to a $\bW$\+locally contraherent cosheaf on $X$
(viewed as a complex concentrated in homological degree~$0$).
 The $\bW$\+locally contraherent cosheaf that appears here is called
the (\emph{$\bW$\+local}) \emph{contraherator} of $\P$ and denoted
by $\Cr\P$.
 A derived contrahereable cosheaf $\P$ on $X$ is called
\emph{contrahereable} if it is locally contrahereable with respect
to the covering $\{X\}$.

 Any derived contrahereable cosheaf $\Q$ on an affine scheme $U$ is
contrahereable, because the contraherator complex
$\Cr_\bu(\{U\},\Q)$ is concentrated in homological degree~$0$.
 The contraherator cosheaf $\Cr\Q$ constructed in this way coincides
with the one defined above specifically in the affine scheme case;
so our notation and terminology is consistent.
 Any $\bW$\+locally contraherent cosheaf $\P$ on a quasi-compact
semi-separated scheme $X$ is $\bW$\+locally contrahereable;
the corresponding contraherator complex $\Cr_\bu(\{U_\alpha\},\P)$
is the contraherent \v Cech resolution~\eqref{contraherent-cech}
of a $\bW$\+locally contraherent cosheaf $\P=\Cr\P$.

 The contraherator complex construction $\Cr_\bu(\{U_\alpha\},{-})$
is an exact functor from the exact category of $\bW$\+locally derived
contrahereable cosheaves to the exact category of finite complexes
of contraherent cosheaves on $X$, or to the bounded derived category
$\sD^\b(X\ctrh)$.
 The full subcategory of $\bW$\+locally contrahereable cosheaves in
the exact category of $\bW$\+locally derived contrahereable cosheaves
is closed under extensions and infinite products.
 Hence it acquires the induced exact category structure.
 The contraherator $\Cr$ is an exact functor from the exact category
of $\bW$\+locally contrahereable cosheaves to that of $\bW$\+locally
contraherent ones.

 All of the above applies to locally cotorsion and locally injective
$\bW$\+locally derived contrahereable cosheaves as well.
 These form full exact subcategories closed under extensions,
infinite products, and cokernels of admissible monomorphisms in
the exact category of $\bW$\+locally derived contrahereable cosheaves.
 The contraherator complex construction $\Cr_\bu(\{U_\alpha\},{-})$
takes locally cotorsion (resp., locally injective) $\bW$\+locally
derived contrahereable cosheaves to finite complexes of locally
cotorsion (resp., locally injective) contraherent cosheaves on~$X$.

 A locally cotorsion (resp., locally injective) $\bW$\+locally
derived contrahereable cosheaf is called $\bW$\+locally contrahereable
if it is $\bW$\+locally contrahereable as a locally contraadjusted
$\bW$\+locally derived contrahereable cosheaf.
 The contraherator functor $\Cr$ takes locally cotorsion (resp.,
locally injective) $\bW$\+locally contrahereable cosheaves to
locally cotorsion (resp., locally injective) $\bW$\+locally
contraherent cosheaves.

\medskip

 Let $f\:Y\rarrow X$ be a $(\bW,\bT)$\+affine morphism of schemes.
 Then the direct image functor~$f_!$ takes $\bT$\+locally
derived contrahereable cosheaves on $Y$ to $\bW$\+locally derived
contrahereable cosheaves on~$X$.
 For a $(\bW,\bT)$\+affine morphism~$f$ of quasi-compact
semi-separated schemes, the functor~$f_!$ also commutes with
the contraherator complex construction, as there is 
a natural isomorphism of complexes of contraherent cosheaves
$$
 f_!\Cr_\bu(\{f^{-1}(U_\alpha)\},\P)\simeq
 \Cr_\bu(\{U_\alpha\},f_!\P)
$$
for any finite affine open covering $U_\alpha$ of $X$ subordinate
to~$\bW$.
 It follows that the functor~$f_!$ takes $\bT$\+locally
contrahereable cosheaves to $\bW$\+locally contrahereable
cosheaves and commutes with the functor~$\Cr$.

 For any $\bW$\+locally contrahereable cosheaf $\P$ and
$\bW$\+locally contraherent cosheaf $\Q$ on a quasi-compact
semi-separated scheme $X$, there is a natural isomorphism
of the groups of morphisms 
\begin{equation}  \label{contraherator-adjunction}
 \Hom^{\O_X}(\P,\Q)\simeq\Hom^X(\Cr\P,\Q).
\end{equation}
 In other words, the functor $\Cr$ is left adjoint to the identity
embedding functor of the category of $\bW$\+locally contrahereable
cosheaves into the category of $\bW$\+locally contraherent ones.
 Indeed, applying the contraherator complex construction
$\Cr_\bu(\{U_\alpha\},{-})$ to a morphism $\P\rarrow\Q$ and passing to
the degree-zero homology, we obtain the corresponding morphism
$\Cr\P\rarrow\Q$.
 (Notice that the homology of a complex in an exact category is not
defined in general, but one can speak about the degree-zero homology
functor on the category of complexes quasi-isomorphic to complexes
concentrated in degree~$0$.)
 Conversely, any cosheaf of $\O_X$\+modules $\P$ is the cokernel of
the rightmost arrow of the complex
\begin{multline}  \label{cech-cosheaf} \textstyle
 0 \lrarrow j_{1,\dotsc,N}{}_! (\P|_{U_1\cap\dotsb\cap U_N})
 \lrarrow \dotsb \\ \textstyle
 \lrarrow \bigoplus_{\alpha<\beta}j_{\alpha,\beta}{}_!
 (\P|_{U_\alpha\cap U_\beta})\lrarrow
 \bigoplus_\alpha j_\alpha{}_!(\P|_{U_\alpha})
\end{multline}
in the additive category of cosheaves of $\O_X$\+modules.
 A cosheaf $\P$ satisfying the exactness condition~(i$^\circ$) for
affine open subschemes $U\sub X$ subordinate to $\bW$ is also
quasi-isomorphic to the whole complex~\eqref{cech-cosheaf} in
the exact category $\O_X\cosh_\bW$.
 Passing to the degree-zero homology of the natural morphism between
the complexes~\eqref{cech-cosheaf} and~\eqref{contraherator-complex},
we produce the desired adjunction morphism $\P\rarrow\Cr\P$.

 For any $\bW$\+locally contrahereable cosheaf $\P$ and any
quasi-coherent sheaf $\M$ on $X$ the natural morphism of
cosheaves of $\O_X$\+modules $\P\rarrow\Cr\P$ induces an isomorphism
of the contratensor products
\begin{equation}  \label{contraherator-contratensor}
 \M\ocn_X\P\simeq\M\ocn_X\Cr(\P).
\end{equation}
 Indeed, for any injective quasi-coherent sheaf $\J$ on $X$ one has
\begin{multline*}
 \Hom_X(\M\ocn_X\P\;\J)\simeq\Hom^{\O_X}(\P,\fHom_X(\M,\J)) \\ \simeq
 \Hom^X(\Cr\P,\fHom_X(\M,\J))\simeq\Hom_X(\M\ocn_X\Cr\P\;\J)
\end{multline*}
in view of the isomorphism~\eqref{fHom-contratensor-adjunction}.

\subsection{$\Cohom$ into a locally derived contrahereable cosheaf}
\label{cohom-loc-der-contrahereable}
 We start with discussing the $\Cohom$ from a quasi-coherent sheaf
to a locally contraherent cosheaf.

 Let $\bW$ be an open covering of a scheme~$X$.
 Let $\F$ be a very flat quasi-coherent sheaf on $X$, and let $\P$
be a $\bW$\+locally contraherent cosheaf on~$X$.
 The $\bW$\+locally contraherent cosheaf $\Cohom_X(\F,\P)$ on
the scheme~$X$ is defined by the rule $U\mpsto
\Hom_{\O_X(U)}(\F(U),\P[U])$ for any affine open subscheme $U\sub X$
subordinate to~$\bW$.
 The contraadjustedness and $\bW$\+local contraherence conditions
can be verified in the same way as it was done in
Section~\ref{cohom-subsection}.

 Similarly, if $\F$ is a flat quasi-coherent sheaf and $\P$ is
a locally cotorsion $\bW$\+locally contraherent cosheaf on $X$,
then the locally cotorsion $\bW$\+locally contraherent cosheaf
$\Cohom_X(\F,\P)$ is defined by the same rule $U\mpsto
\Hom_{\O_X(U)}(\F(U),\P[U])$ for any affine open subscheme $U\sub X$
subordinate to~$\bW$.

 Finally, if $\M$ is a quasi-coherent sheaf and $\gJ$ is
a locally injective $\bW$\+locally contraherent cosheaf on $X$,
then the locally cotorsion $\bW$\+locally contraherent cosheaf
$\Cohom_X(\M,\gJ)$ is defined by the same rule as above.
 If $\F$ is a flat quasi-coherent sheaf and $\gJ$ is a locally
injective $\bW$\+locally contraherent cosheaf on $X$, then
the $\bW$\+locally contraherent cosheaf $\Cohom_X(\F,\gJ)$ is
locally injective.

 For any two very flat quasi-coherent sheaves $\F$ and $\G$ on
a scheme $X$ and any $\bW$\+locally contraherent cosheaf $\P$ on $X$
there is a natural isomorphism of $\bW$\+locally contraherent
cosheaves
\begin{equation} \label{flat-flat-cohom-cohom}
 \Cohom_X(\F\ot_{\O_X}\G\;\P)\simeq\Cohom_X(\G,\Cohom_X(\F,\P)).
\end{equation}
 Similarly, for any two flat quasi-coherent sheaves $\F$ and $\G$
and a locally cotorsion $\bW$\+locally contraherent cosheaf $\P$
on $X$ there is a natural isomorphism~\eqref{flat-flat-cohom-cohom}
of locally cotorsion $\bW$\+locally contraherent cosheaves.
 Finally, for any flat quasi-coherent sheaf $\F$, quasi-coherent
sheaf $\M$, and locally injective $\bW$\+locally contraherent
cosheaf $\gJ$ on $X$ there are natural isomorphisms of locally
cotorsion $\bW$\+locally contraherent cosheaves
\begin{multline}  \label{flat-lin-cohom-cohom}
 \Cohom_X(\M\ot_{\O_X}\F\;\gJ) \\ \simeq\Cohom_X(\M,\Cohom_X(\F,\gJ))
 \simeq\Cohom_X(\F,\Cohom_X(\M,\gJ)).
\end{multline}

\medskip

 More generally, let $\F$ be a very flat quasi-coherent sheaf on $X$,
and let $\P$ be a $\bW$\+locally derived contrahereable cosheaf
on~$X$.
 The $\bW$\+locally derived contrahereable cosheaf $\Cohom_X(\F,\P)$
on the scheme $X$ is defined by the rule $U\mpsto
\Hom_{\O_X(U)}(\F(U),\P[U])$ for any affine open subscheme $U\sub X$
subordinate to~$\bW$.
 For any pair of embedded affine open subschemes $V\sub U$
the restriction and corestriction morphisms $\F(U)\rarrow\F(V)$
and $\P[V]\rarrow\P[U]$ induce a morphism of $\O_X(U)$\+modules
$$
 \Hom_{\O_X(V)}(\F(V),\P[V])\rarrow\Hom_{\O_X(U)}(\F(U),\P[U]),
$$
so our rule defines a covariant functor with $\O_X$\+module
structure on the category of affine open subschemes $U\sub X$
subordinate to~$\bW$.
 The contraadjustedness condition clearly holds; and to check
the exactness condition~(i$^\circ$) of
Section~\ref{contrahereable-subsect} for this covariant functor,
it suffices to apply the functor $\Hom_{\O_X(U)}(\F(U),{-})$ to
the exact sequence of contraadjusted
$\O_X(U)$\+modules~\eqref{exactness-condition} for the cosheaf~$\P$.

 Similarly, if $\F$ is a flat quasi-coherent sheaf and $\P$ is
a locally cotorsion $\bW$\+locally derived contrahereable cosheaf
on $X$, then the locally cotorsion $\bW$\+locally derived
contrahereable cosheaf $\Cohom_X(\F,\P)$ is defined by the same rule
$U\mpsto \Hom_{\O_X(U)}\allowbreak(\F(U),\P[U])$ for any affine open
subscheme $U\sub X$ subordinate to~$\bW$.
 Finally, if $\M$ is a quasi-coherent sheaf and $\gJ$ is a locally
injective $\bW$\+locally derived contrahereable cosheaf on $X$,
then the locally cotorsion $\bW$\+locally derived contrahereable
cosheaf $\Cohom_X(\M,\gJ)$ is defined by the very same rule.
 If $\F$ is a flat quasi-coherent sheaf and $\gJ$ is a locally
injective $\bW$\+locally derived contrahereable cosheaf on $X$,
then the $\bW$\+locally derived contrahereable cosheaf
$\Cohom_X(\F,\gJ)$ is locally injective.

 The associativity isomorphisms~(\ref{flat-flat-cohom-cohom}\+-%
\ref{flat-lin-cohom-cohom}) hold for the $\Cohom$ into 
$\bW$\+locally derived contrahereable cosheaves under
the assumptions similar to the ones made above in the locally
contraherent case.

\subsection{Contraherent tensor product}
\label{contraherent-tensor}
 Let $\P$ be a cosheaf of $\O_X$\+modules on a scheme $X$ and $\M$
be a quasi-coherent sheaf on~$X$.
 Define a covariant functor with an $\O_X$\+module structure
on the category of affine open subschemes of $X$ by the rule
$U\mpsto\M(U)\ot_{\O_X(U)}\P[U]$.
 To a pair of embedded affine open subschemes $V\sub U\sub X$ this
functor assigns the $\O_X(U)$\+module homomorphism
$$
 \M(V)\ot_{\O_X(V)}\P[V] \simeq \M(U)\ot_{\O_X(U)}\P[V]
 \lrarrow \M(U)\ot_{\O_X(U)}\P[U].
$$
 Obviously, this functor satisfies the condition~\eqref{base-cosheaf}
of Theorem~\ref{cosheaf-base-thm} (since the restriction of
the cosheaf $\P$ to affine open subschemes of $X$ does).
 Hence the functor $U\mpsto\M(U)\ot_{\O_X(U)}\P[U]$ extends uniquely
to a cosheaf of $\O_X$\+modules on $X$, which we will denote by
$\M\ot_X\P$.

 For any quasi-coherent sheaf $\M$, cosheaf of $\O_X$\+modules $\P$,
and locally injective $\bW$\+locally derived contrahereable cosheaf
$\gJ$ on a scheme $X$ there is a natural isomorphism of abelian
groups
\begin{equation} \label{cosheaf-tensor-product-adjunction}
 \Hom^{\O_X}(\M\ot_X\P\;\gJ)\simeq\Hom^{\O_X}(\P,\Cohom_X(\M,\gJ)).
\end{equation}
 The analogous adjunction isomorphism holds in the other cases
mentioned in Section~\ref{cohom-loc-der-contrahereable} when
the functor $\Cohom$ from a quasi-coherent sheaf to a locally
derived contrahereable cosheaf is defined.
 In other words, the functors $\Cohom_X(\M,{-})$ and $\M\ot_X{-}$
between subcategories of the category of cosheaves of $\O_X$\+modules
are adjoint ``wherever the former functor is defined''.

 For a locally free sheaf of finite rank $\E$ and a $\bW$\+locally
contraherent (resp., $\bW$\+locally derived contrahereable) cosheaf
$\P$ on a scheme $X$, the cosheaf of $\O_X$\+mod\-ules $\E\ot_X\P$
is $\bW$\+locally contraherent (resp., $\bW$\+locally derived
contrahereable).
 There is a natural isomorphism of $\bW$\+locally contraherent
(resp., $\bW$\+locally derived contrahereable) cosheaves
$\qHom_{\O_X}(\E,\O_X)\ot_X\P\simeq \Cohom_X(\E,\P)$ on~$X$.

 The isomorphism $j_*j^*(\K\ot_{\O_X}\M)\simeq j_*j^*\K
\ot_{\O_X(U)}\M(U)$ \eqref{push-pull-tensor-mod} for quasi-coherent
sheaves $\K$ and $\M$ and the embedding of an affine open subscheme
$j\:U\rarrow X$ allows to construct a natural isomorphism of
quasi-coherent sheaves
\begin{equation}  \label{sheaf-cosheaf-tensor-assoc}
 (\K\ot_{\O_X}\M)\ocn_X\P\simeq \K\ocn_X(\M\ot_X\P)
\end{equation}
for any quasi-coherent sheaves $\K$ and $\M$ and a cosheaf
of $\O_X$\+modules $\P$ on a semi-separated scheme~$X$.

\medskip

 Let $\bW$ be an open covering of a scheme~$X$.
 We will call a cosheaf of $\O_X$\+modules~$\gF$ \ \emph{$\bW$\+flat}
if the $\O_X(U)$\+module $\gF[U]$ is flat for every affine open
subscheme $U\sub X$ subordinate to~$\bW$.
 A cosheaf of $\O_X$\+modules $\gF$ is said to be \emph{flat} if it is
$\{X\}$\+flat.
 Clearly, the direct image of a $\bT$\+flat cosheaf of $\O_Y$\+modules
with respect to a flat $(\bW,\bT)$\+affine morphism of schemes
$f\:Y\rarrow X$ is $\bW$\+flat.

 One can easily see that whenever a cosheaf $\gF$ is $\bW$\+flat and
satisfies the ``exactness condition''~(i$^\circ$) of
Section~\ref{contrahereable-subsect} for finite affine open coverings
of affine open subschemes $U\sub X$ subordinate to $\bW$,
the cosheaf $\M\ot_X\gF$ also satisfies the condition~(i$^\circ$)
for such open affines $U\sub X$.
 Similarly, whenever a quasi-coherent sheaf $\F$ is flat and
a cosheaf of $\O_X$\+modules $\P$ satisfies the condition~(i$^\circ$),
so does the cosheaf $\F\ot_X\P$.

 The full subcategory of $\bW$\+flat $\bW$\+locally contraherent 
cosheaves is closed under extensions and kernels of admissible
epimorphisms in the exact category of $\bW$\+locally contraherent
cosheaves $X\lcth_\bW$ on~$X$.
 Hence it acquires the induced exact category structure, which we will
denote by $X\lcth_\bW^\fl$.
 The category $X\lcth_{\{X\}}^\fl$ will be denoted by $X\ctrh^\fl$.
 Similarly, there is the exact category structure on the category
of $\bW$\+flat $\bW$\+locally derived contrahereable cosheaves on~$X$
inherited from the exact category of $\bW$\+locally derived
contrahereable cosheaves.

 A quasi-coherent sheaf $\K$ on a scheme $X$ is called \emph{coadjusted}
if the $\O_X(U)$\+module $\K(U)$ is coadjusted (see
Section~\ref{coherent-coadjusted}) for every affine open subscheme
$U\sub X$.
 By Lemma~\ref{coadjusted-local}, the coadjustedness of a quasi-coherent
sheaf is a local property.
 By the definition, if a cosheaf of $\O_X$\+modules $\P$ on $X$ 
satisfies the contraadjustedness condition~(ii) of
Section~\ref{contraherent-definition} or~\ref{contrahereable-subsect}
and a quasi-coherent sheaf $\K$ on $X$ is coadjusted, then
the cosheaf of $\O_X$\+modules $\K\ot_X\P$ also satisfies
the condition~(ii).

 The full subcategory of coadjusted quasi-coherent sheaves is closed
under extensions and the passage to quotient objects in the abelian
category of quasi-coherent sheaves $X\qcoh$ on~$X$.
 Hence it acquires the induced exact category structure, which we will
denote by $X\qcoh^\coa$.

 Let $\K$ be a coadjusted quasi-coherent sheaf and $\gF$ be a
$\bW$\+flat $\bW$\+locally contraherent (or more generally, $\bW$\+flat
$\bW$\+locally derived contrahereable) cosheaf on a scheme~$X$.
 Then the tensor product $\K\ot_X\gF$ satisfies both
conditions (i$^\circ$) and~(ii) for affine open subschemes
$U\sub X$ subordinate to~$\bW$, i.~e., it is $\bW$\+locally derived
contrahereable.
 (Of course, the cosheaf $\K\ot_X\gF$ is \emph{not} in general locally
contraherent, even if the cosheaf $\gF$ was $\bW$\+locally
contraherent.)

 Assuming that the scheme $X$ is quasi-compact and semi-separated,
the contraherator complex construction now allows to assign to
this cosheaf of $\O_X$\+modules a complex of contraherent cosheaves
$\Cr_\bu(\{U_\alpha\}\;\K\ot_X\gF)$ on~$X$.
 To a short exact sequence of coadjusted quasi-coherent sheaves or
$\bW$\+flat $\bW$\+locally contraherent (or $\bW$\+flat $\bW$\+locally
derived contrahereable) cosheaves on~$X$, the functor
$\Cr_\bu(\{U_\alpha\}\;{-}\ot_X{-})$ assigns a short exact sequence
of complexes of contraherent cosheaves.
 
 By Lemma~\ref{contraherator-well-defined}, the corresponding
object of the bounded derived category $\sD^\b(X\ctrh)$ does not
depend on the choice of a finite affine open covering $\{U_\alpha\}$.
 We will denote it by $\K\ot_{X\ct}^\boL\gF$ and call
the \emph{derived contraherent tensor product} of a coadjusted
quasi-coherent sheaf $\K$ and a $\bW$\+flat $\bW$\+locally
contraherent cosheaf $\gF$ on a quasi-compact semi-separated
scheme~$X$.

 When the derived category object $\K\ot_{X\ct}^\boL\gF$, viewed
as an object of the derived category $\sD^\b(X\lcth_\bW)$
via the embedding of exact categories $X\ctrh\rarrow X\lcth_\bW$,
turns out to be isomorphic to an object of the exact category
$X\lcth_\bW$, we say that the (underived) \emph{contraherent
tensor product} of $\K$ and $\gF$ is defined, and denote
the corresponding object by $\K\ot_{X\ct}\gF\.\in X\lcth_\bW$.
 In other words, for a coadjusted quasi-coherent sheaf $\K$ and
a $\bW$\+flat $\bW$\+locally contraherent cosheaf $\gF$ on
a quasi-compact semi-separated scheme $X$ one sets
$\K\ot_{X\ct}\gF = \Cr(\K\ot_X\gF)$ whenever the right-hand side
is defined (where $\Cr$ denotes the $\bW$\+local contraherator).

\medskip

 A scheme $X$ is called \emph{locally coherent}~\cite{CEI,EG} if
for every affine open subscheme $U\sub X$ the ring $\O_X(U)$ is
coherent.
 By Lemma~\ref{coherence-of-rings-local}, it suffices to check
this property for affine open subschemes belonging to any chosen
affine open covering of~$X$.
 Following~\cite[Section~9.1]{PS5}, we say that a scheme $X$ is
\emph{coherent} if $X$ is locally coherent, quasi-compact, and
quasi-separated.

 Let $X$ be a locally coherent scheme.
 Then, by Corollary~\ref{coherent-flat-local}(a), a contraherent
cosheaf $\gF$ on an affine open subscheme $U\sub X$ is flat if and
only if the contraadjusted $\O(U)$\+module $\gF[U]$ is flat.
 Besides, the full subcategory of $\bW$\+flat $\bW$\+locally
contraherent cosheaves in $X\lcth_\bW$ is closed under
infinite products.

 Now assume that the scheme $X$ is locally Noetherian.
 Then any (locally) coherent sheaf on $X$ is coadjusted.
 Any injective quasi-coherent sheaf on $X$ and any quasi-coherent
quotient sheaf of an injective one are coadjusted, too.

 For any injective quasi-coherent sheaf $\J$ and any $\bW$\+flat
$\bW$\+locally derived contrahereable cosheaf $\gF$ on $X$,
the tensor product $\J\ot_X\gF$ is a locally injective
$\bW$\+locally derived contrahereable cosheaf on~$X$.
 For any flat quasi-coherent sheaf $\F$ and any locally injective
$\bW$\+locally derived contrahereable cosheaf $\gJ$ on $X$,
the tensor product $\F\ot_X\gJ$ is a locally injective
$\bW$\+locally derived contrahereable cosheaf.
 These assertions hold since the tensor product of a flat module and
an injective module over a Noetherian ring is injective.

 Let $\M$ be a coherent sheaf and $\gF$ be a $\bW$\+flat
$\bW$\+locally contraherent cosheaf on~$X$.
 Then the cosheaf of $\O_X$\+modules $\M\ot_X\gF$ is $\bW$\+locally
contraherent.
 Indeed, by Corollary~\ref{coherent-very}(a),
the $\O_X(U)$\+module $\M(U)\ot_{\O_X(U)}\gF[U]$ is contraadjusted
for any affine open subscheme $U\sub X$.
 For a pair of embedded affine open subschemes $V\sub U\sub X$
subordinate to the covering~$\bW$, one has
\begin{multline*}
 \M(V)\ot_{\O_X(V)}\gF[V] \simeq \M(U)\ot_{\O_X(U)}\gF[V] \\
 \simeq\M(U)\ot_{\O_X(U)}\Hom_{\O_X(U)}(\O_X(V),\gF[U]) \\
 \simeq \Hom_{\O_X(U)}(\O_X(V)\;\M(U)\ot_{\O_X(U)}\gF[U])
\end{multline*}
according to Corollary~\ref{coherent-very}(c).
 If the scheme $X$ is semi-separated and Noetherian, the contraherent
tensor product $\M\ot_{X\ct}\gF$ is defined and isomorphic to
the tensor product $\M\ot_X\gF$.

 Similarly, it follows from Corollary~\ref{coherent-cotors} that
for any coherent sheaf $\M$ and $\bW$\+flat locally cotorsion
$\bW$\+locally contraherent cosheaf $\gF$ on $X$ the tensor product
$\M\ot_X\gF$ is a locally cotorsion $\bW$\+locally contraherent
cosheaf on~$X$.
 (See Sections~\ref{lct-projective-loc-noetherian}\+-%
\ref{flat-contra-subsection} and
Lemma~\ref{noetherian-contraherent-tensor} for further discussion.)

\subsection{Compatibility of direct and inverse images with
the tensor operations} \label{compatibility-subsect}
 Let $\bW$ be an open covering of a scheme $X$ and $\bT$ be an
open covering of a scheme~$Y$.
 Let $f\:Y\rarrow X$ be a $(\bW,\bT)$\+coaffine morphism.

 Let $\F$ be a flat quasi-coherent sheaf and $\gJ$ be a locally
injective $\bW$\+locally contraherent cosheaf on the scheme~$X$.
 Then there is a natural isomorphism of locally injective 
$\bT$\+locally contraherent cosheaves
\begin{equation}
 f^!\Cohom_X(\F,\gJ)\simeq\Cohom_Y(f^*\F,f^!\gJ)
\end{equation}
on the scheme~$Y$.
 Indeed, for any affine open subscheme $U\sub X$ subordinate to $\bW$
and any affine open subscheme $V\sub Y$ subordinate to $\bT$ such that
$f(V)\sub U$, one has
\begin{multline*}
 f^!\Cohom_X(\F,\gJ)[V]\simeq
 \Hom_{\O_X(U)}(\O_Y(V),\Hom_{\O_X(U)}(\F(U),\gJ[U])) \\ \simeq
 \Hom_{\O_Y(V)}(\O_Y(V)\ot_{\O_X(U)}\F(U)\;
 \Hom_{\O_X(U)}(\O_Y(V),\gJ[U])) \\ \simeq
 \Cohom_Y(f^*\F,f^!\gJ)[V].
\end{multline*}

 Assume additionally that $f$~is a flat morphism.
 Let $\M$ be a quasi-coherent sheaf and $\gJ$ be a locally injective
$\bW$\+locally contraherent cosheaf on~$X$.
 Then there is a natural isomorphism of locally cotorsion
$\bT$\+locally contraherent cosheaves
\begin{equation}  \label{flat-lin-inverse-cohom}
 f^!\Cohom_X(\M,\gJ)\simeq\Cohom_Y(f^*\M,f^!\gJ)
\end{equation}
on~$Y$.
 Analogously, if $\F$ is a flat quasi-coherent sheaf and $\P$ is
a locally cotorsion $\bW$\+locally contraherent cosheaf on $X$,
then there is a natural isomorphism of locally cotorsion
$\bT$\+locally contraherent cosheaves
\begin{equation}  \label{flat-cotors-inverse-cohom}
 f^!\Cohom_X(\F,\P)\simeq\Cohom_Y(f^*\F,f^!\P)
\end{equation}
on the scheme~$Y$.

 Assume that, moreover, $f$~is a very flat morphism.
 Let $\F$ be a very flat quasi-coherent sheaf and $\P$ be
a $\bW$\+locally contraherent cosheaf on~$X$.
 Then there is the natural
isomorphism~\eqref{flat-cotors-inverse-cohom}
of $\bT$\+locally contraherent cosheaves on~$Y$.

 Let $f\:Y\rarrow X$ be a $(\bW,\bT)$\+affine $(\bW,\bT)$\+coaffine
morphism.
 Let $\N$ be a quasi-coherent sheaf on $Y$ and $\gJ$ be a locally
injective $\bW$\+locally contraherent cosheaf on~$X$.
 Then there is a natural isomorphism of locally cotorsion
$\bW$\+locally contraherent cosheaves
\begin{equation}  \label{lin-projection-cohom}
 \Cohom_X(f_*\N,\gJ)\simeq f_!\Cohom_Y(\N,f^!\gJ)
\end{equation}
on the scheme~$X$.
 Indeed, for any affine open subscheme $U\sub X$ subordinate to $\bW$
and its preimage $V=f^{-1}(U)\sub Y$, one has
\begin{multline*}
 \Cohom_X(f_*\N,\gJ)[U]\simeq
 \Hom_{\O_X(U)}(\N(V),\gJ[U]) \\
 \simeq\Hom_{\O_Y(V)}(\N(V),\Hom_{\O_X(U)}(\O_Y(V),\gJ[U])\simeq
 f_!\Cohom_Y(\N,f^!\gJ)[U].
\end{multline*}
 This is one version of the projection formula for the $\Cohom$ from
a quasi-coherent sheaf to a contraherent cosheaf.

 Assume additionally that $f$~is a flat morphism.
 Let $\G$ be a flat quasi-coherent sheaf on $Y$ and $\P$ be
a locally cotorsion $\bW$\+locally contraherent cosheaf on~$X$.
 Then there is a natural isomorphism of locally cotorsion
$\bW$\+locally contraherent cosheaves
\begin{equation}  \label{flat-cotors-projection-cohom}
 \Cohom_X(f_*\G,\P)\simeq f_!\Cohom_Y(\G,f^!\P)
\end{equation}
on the scheme~$X$.

 Assume that, moreover, $f$~is a very flat morphism.
 Let $\G$ be a very flat quasi-coherent sheaf on $Y$ and $\P$ be
a $\bW$\+locally contraherent cosheaf on~$X$.
 Then there is the natural
isomorphism~\eqref{flat-cotors-projection-cohom} of $\bW$\+locally
contraherent cosheaves on~$X$.

 Let $f\:Y\rarrow X$ be a $(\bW,\bT)$\+affine morphism.
 Let $\F$ be a very flat quasi-coherent sheaf on $X$ and $\Q$ be
a $\bT$\+locally contraherent cosheaf on~$Y$.
 Then there is a natural isomorphism of $\bW$\+locally contraherent
cosheaves
\begin{equation}
 \Cohom_X(\F,f_!\Q)\simeq f_!\Cohom_Y(f^*\F,\Q)
\end{equation}
on the scheme~$X$.
 Indeed, for any affine open subscheme $U\sub X$ subordinate to $\bW$
and its preimage $V=f^{-1}(U)\sub Y$, one has
\begin{multline*}
 \Cohom_X(\F,f_!\Q)[U]\simeq
 \Hom_{\O_X(U)}(\F(U),\Q[V]) \\ \simeq
 \Hom_{\O_Y(V)}(\O_Y(V)\ot_{\O_X(U)}\F(U)\;\Q[V]) \simeq
 f_!\Cohom_Y(f^*\F,\Q)[U].
\end{multline*}
 The similar isomorphism of locally cotorsion $\bW$\+locally
contraherent cosheaves on $X$ holds for any flat quasi-coherent
sheaf $\F$ on $X$ and locally cotorsion $\bT$\+locally contraherent
cosheaf $\Q$ on~$Y$.
 This is another version of the projection formula for~$\Cohom$.

 Assume additionally that $f$ is a flat morphism.
 Let $\M$ be a quasi-coherent sheaf on $X$ and $\gI$ be
a locally injective $\bT$\+locally contraherent cosheaf on~$Y$.
 Then there is a natural isomorphism of locally cotorsion
$\bW$\+locally contraherent cosheaves
\begin{equation} \label{lin-projection-II-cohom}
 \Cohom_X(\M,f_!\gI)\simeq f_!\Cohom_Y(f^*\M,\gI)
\end{equation}
on the scheme~$X$.

\medskip

 Let $f\:Y\rarrow X$ be a quasi-compact morphism of semi-separated
schemes.
 Let $\F$ be a very flat quasi-coherent sheaf on $X$ and $\cQ$ be
a contraadjusted quasi-coherent sheaf on~$Y$.
 Then there is a natural isomorphism of contraherent cosheaves
\begin{equation}  \label{flat-cta-fHom-projection}
 \fHom_X(\F,f_*\cQ)\simeq f_!\fHom_Y(f^*\F,\cQ)
\end{equation}
on~$X$.
 Here the quasi-coherent sheaf $f_*\cQ$ on $X$ is contraadjusted
according to Section~\ref{fHom-subsection}.
 The right-hand side is, by construction, a contraherent cosheaf
if the morphism~$f$ is affine, and a cosheaf of $\O_X$\+modules
otherwise (see Section~\ref{contra-direct-inverse}).
 Both sides are, in fact, contraherent in the general case, because
the isomorphism holds and the left-hand side is.
 This is a version of the projection formula for~$\fHom$.

 Indeed, let $j\:U\rarrow X$ be an embedding of an affine open
subscheme; set $V=U\times_X Y$.
 Let $j'\:V\rarrow Y$ and $f'\:V\rarrow U$ be the natural morphisms.
 Then one has
\begin{multline*}
 \fHom_X(\F,f_*\cQ)[U] \simeq
 \Hom_X(j_*j^*\F,f_*\cQ)\simeq\Hom_Y(f^*\!\.j_*j^*\F,\cQ) \\
 \simeq \Hom_Y(j'_*f'{}^*\!\.j^*\F,\cQ)\simeq
 \Hom_Y(j'_*j'{}^*\!\.f^*\F,\cQ) \simeq \fHom_Y(f^*\F,\cQ)[V].
\end{multline*}
 Here we are using the fact that the direct images of quasi-coherent
sheaves with respect to affine morphisms of schemes commute with
the inverse images in the base change situations.
 Notice that, when the morphism $f$~is not affine, neither is
the scheme~$V$; however, the scheme $V$ is quasi-compact and the open
embedding morphism $j'\:V\rarrow Y$ is affine, so
Lemma~\ref{fHom-cosections}(a) applies.
 The similar isomorphism of locally cotorsion contraherent cosheaves
on $X$ holds for any flat quasi-coherent sheaf $\F$ on $X$ and
any cotorsion quasi-coherent sheaf $\cQ$ on~$Y$
(by Lemma~\ref{fHom-cosections}(b)).

 Now let $f$~be a flat quasi-compact morphism of quasi-separated schemes.
 Let $\M$ be a quasi-coherent sheaf on $X$ and $\I$ be an injective
quasi-coherent sheaf on~$Y$.
 Then there is a natural isomorphism of locally cotorsion
contraherent cosheaves
\begin{equation}  \label{inj-fHom-projection}
 \fHom_X(\M,f_*\I)\simeq f_!\fHom_Y(f^*\M,\I)
\end{equation}
on the scheme~$Y$.
 The proof is similar to the above and uses
Lemma~\ref{fHom-cosections}(c).
 Here it is helpful to keep in mind that the direct images of injective
quasi-coherent sheaves under flat quasi-compact quasi-separated
morphisms of schemes are injective quasi-coherent sheaves again.
 It is also important for this argument that the direct images of
quasi-coherent sheaves with respect to quasi-compact quasi-separated 
morphisms commute with the inverse images with respect to flat morphisms
of schemes in the base change situations.

 Let $f\:Y\rarrow X$ be an affine morphism of semi-separated schemes.
 Let $\M$ be a quasi-coherent sheaf on $X$ and $\Q$ be a cosheaf of
$\O_Y$\+modules.
 Then there is a natural isomorphism of quasi-coherent sheaves
\begin{equation} \label{contratensor-projection}
 \M\ocn_X f_!\Q \simeq f_*(f^*\M\ocn_Y\Q)
\end{equation}
on the scheme~$X$.
 This is a version of the projection formula for the contratensor
product of quasi-coherent sheaves and cosheaves of $\O$\+modules.
 (See Section~\ref{direct-contratensor-nonaffine-subsect} below
for partial generalizations to nonaffine morphisms~$f$.)

 Indeed, in the notation above, for any affine open subscheme
$U\sub X$ we have
\begin{multline*}
 j_*j^*\M\ot_{\O_X(U)}(f_!\Q)[U] \simeq
 j_*j^*\M\ot_{\O_X(U)}\Q[V] \\ \simeq
 (j_*j^*\M\ot_{\O_X(U)}\O_Y(V))\ot_{\O_Y(V)}\Q[V] \simeq
 j_*(j^*\M\ot_{\O_X(U)}\O_Y(V))\ot_{\O_Y(V)}\Q[V] \\ \simeq
 j_*f'_*f'{}^*\!\.j^*\M\ot_{\O_Y(V)}\Q[V] \simeq
 f_*j'_*j'{}^*\!\.f^*\M\ot_{\O_Y(V)}\Q[V] \\ \simeq
 f_*(j'_*j'{}^*\!\.f^*\M\ot_{\O_Y(V)}\Q[V]).
\end{multline*}
 As it was explained Section~\ref{contratensor-subsect}, the contratensor
product $f^*\M\ocn_Y\Q$ can be computed as the inductive limit over
the diagram $\bD$ formed by the affine open subschemes $V\sub Y$ of
the form $V=U\times_XY$, where $U$ are affine open subschemes in~$X$.
 It remains to use the fact that the direct image of quasi-coherent
sheaves with respect to an affine morphism is an exact functor
preserving infinite direct sums.
 The same isomorphism~\eqref{contratensor-projection} holds for a flat
affine morphism~$f$ of quasi-separated schemes (the flatness
assumption on~$f$ is needed here so that the functor~$j_*$ commutes
with the functor ${-}\ot_{\O_X(U)}\O_Y(V)$).

 Furthermore, there is a natural morphism from the left-hand side to
the right-hand side of~\eqref{contratensor-projection} for any
quasi-compact morphism of quasi-separated schemes $f\:Y\rarrow X$.
 It is constructed as the composition
\begin{multline*} \textstyle
 \M\ocn_Xf_!\Q\.=\.\varinjlim_U j_*j^*\M\ot_{\O_X(U)} \Q[f^{-1}(U)]
 \.\simeq\.\varinjlim_U j_*j^*\M\ot_{\O_X(U)}
 (\varinjlim_{V\sub f^{-1}(U)}\Q[V]) \\ \textstyle \simeq\.
 \varinjlim_U\varinjlim_{f(V)\sub U} j_*j^*\M\ot_{\O_X(U)} \Q[V]
 \.\simeq\. \varinjlim_{f(V)\sub U} j_*j^*\M\ot_{\O_X(U)} \Q[V]
 \lrarrow \\ \textstyle
 \varinjlim_{f(V)\sub U} j_*(j^*\M\ot_{\O_X(U)}\O_Y(V))
 \ot_{\O_Y(V)}\Q[V] \.\simeq\.
 \varinjlim_{f(V)\sub U} f_*j'_*j'{}^*\!\.f^*\M\ot_{\O_Y(V)}\Q[V]
 \\ \textstyle \lrarrow
 f_*\varinjlim_{f(V)\sub U} j'_*j'{}^*\!\.f^*\M\ot_{\O_Y(V)}\Q[V]
 \.\simeq\. f_*(f^*\M\ocn_Y\Q).
\end{multline*}
 Here the inductive limit is taken firstly over affine open subschemes
$U\sub X$, then over affine open subschemes $V\sub Y$ such that
$f(V)\sub U$, and eventually over pairs of affine open subschemes
$U\sub X$ and $V\sub Y$ such that $f(V)\sub U$.
 The open embeddings $U\rarrow X$ and $V\rarrow Y$ are denoted by $j$
and~$j'$, while the morphism $V\rarrow U$ is denoted by~$f'$.
 The final isomorphism holds, since the contratensor product
$f^*\M\ocn_Y\Q$ can be computed over the diagram $\bD$ formed by all
the pairs of affine open subschemes $(U,V)$ such that $f(V)\sub U$.

 In particular, the isomorphism
\begin{equation} \label{flat-contratensor-projection}
 \M\ocn_X h_!\gF \simeq h_*(h^*\M\ocn_Y\gF)
\end{equation}
holds for any open embedding $h\:Y\rarrow X$ of an affine scheme $Y$
into a quasi-separated scheme $X$, any quasi-coherent sheaf $\M$
on $X$, and any flat cosheaf of $\O_Y$\+modules $\gF$ on~$Y$.
 Indeed, in this case one has
\begin{multline*} \textstyle
 \varinjlim_{h(V)\sub U} h_*j'_*j'{}^*h^*\M\ot_{\O_Y(V)}\gF[V]
 \simeq\varinjlim_V h_*j'_*j'{}^*h^*\M\ot_{\O_Y(V)}\gF[V] \\
 \simeq h_*h^*\M\ot_{\O(Y)}\gF[Y]\simeq h_*(h^*\M\ot_{\O(Y)}\gF[Y]),
\end{multline*}
where the $\varinjlim_V$ is taken over all affine open subschemes
$V\sub Y$, which is clearly equivalent to considering $V=Y$ only.

 Once again, we refer to
Section~\ref{direct-contratensor-nonaffine-subsect} below
for a discussion of the projection formula isomorphisms for
the direct images of contratensor products under more general classes
of morphisms of schemes than in~(\ref{contratensor-projection}\+-%
\ref{flat-contratensor-projection}).

\medskip

 Let $f\:Y\rarrow X$ be an affine morphism of schemes.
 Then for any quasi-coherent sheaf $\M$ on $X$ and any cosheaf
of $\O_Y$\+modules $\Q$ there is a natural isomorphism of
cosheaves of $\O_X$\+modules
\begin{equation}  \label{cosheaf-tensor-projection}
 f_!(f^*\M\ot_Y\Q) \simeq \M\ot_X f_!\Q.
\end{equation}
 This is a version of the projection formula for the tensor product
of quasi-coherent sheaves and cosheaves of $\O_X$\+modules.
 Indeed, for any affine open subscheme $U\sub X$ and its preimage
$V=f^{-1}(U)\sub Y$, one has
\begin{multline*}
 f_!(f^*\M\ot_Y\Q)[U]\simeq
 (\O_Y(V)\ot_{\O_X(U)}\M(U))\ot_{\O_Y(V)}\Q[V] \\ \simeq
 \M(U)\ot_{\O_X(U)}\Q[V]\simeq
 (\M\ot_X f_!\Q)[U].
\end{multline*}

\Section{Derived Categories on Quasi-Compact Semi-Separated Schemes}
\label{derived-on-quasi-compact-sect}

\subsection{Contraadjusted and cotorsion quasi-coherent sheaves}
\label{quasi-compact-quasi-coherent}
 Recall that the definition of a very flat quasi-coherent sheaf was
given in Section~\ref{very-flat-morphisms-subsect} and the definition
of a contraadjusted quasi-coherent sheaf in Section~\ref{fHom-subsection}
(cf.\ Remark~\ref{cotorsion-sheaf-ambiguity}).

 In particular, a quasi-coherent sheaf $\cP$ over an affine scheme $U$
is very flat (respectively, contraadjusted) if and only if
the $\O(U)$\+module $\cP(U)$ is very flat (respectively, contraadjusted).
 The class of very flat quasi-coherent sheaves is preserved by inverse
images with respect to arbitrary morphisms of schemes and direct images
with respect to very flat affine morphisms (which includes affine
open embeddings).
 The class of contraadjusted quasi-coherent sheaves is preserved by
direct images with respect to affine (and even arbitrary quasi-compact
quasi-separated) morphisms of schemes.

 The class of very flat quasi-coherent sheaves on any scheme $X$
is closed under kernels of surjective morphisms.
 Both the full subcategories of very flat and contraadjusted
quasi-coherent sheaves are closed under extensions in the abelian
category of quasi-coherent sheaves.
 Hence they acquire the induced exact category structures, which
we denote by $X\qcoh_\vfl$ and $X\qcoh^\cta$, respectively.

 Let us introduce one bit of categorical terminology.
 Given an exact category $\sE$ and a class of objects $\sC\sub\sE$,
we say that an object $X\in\sE$ is a \emph{finitely iterated 
extension} of objects from $\sC$ if there exists a nonnegative
integer $N$ and a sequence of admissible monomorphisms
$0=X_0\rarrow X_1\rarrow \dotsb\rarrow X_{N-1}\rarrow X_N=X$
in $\sE$ such that the cokernels of all the morphisms
$X_{i-1}\rarrow X_i$ belong to~$\sC$ (cf.\
Section~\ref{very-eklof-trlifaj-subsect}).

 Let $X$ be a quasi-compact semi-separated scheme.

\begin{lem}  \label{quasi-very-flat-cover}
 Any quasi-coherent sheaf\/ $\M$ on $X$ can be included in a short
exact sequence\/ $0\rarrow\cP\rarrow\F\rarrow\M\rarrow0$, where\/ $\F$
is a very flat quasi-coherent sheaf and\/ $\cP$ is a finitely iterated
extension of the direct images of contraadjusted quasi-coherent sheaves
from affine open subschemes in~$X$.
\end{lem}

\begin{proof}
 This is a particular case of
Proposition~\ref{quasi-coherent-gluing-precover}; see
Remark~\ref{quasi-compact-quasi-coherent-remark} below.
 The proof is based on the construction from~\cite[Section~A.1]{EP}
and Theorem~\ref{eklof-trlifaj-very}(b).
 We argue by a kind of induction in the number of affine open
subschemes covering~$X$.
 Assume that for some open subscheme $h\:V\rarrow X$ there is
a short exact sequence $0\rarrow\cQ\rarrow\K\rarrow\M\rarrow0$
of quasi-coherent sheaves on $X$ such that the restriction $h^*\K$
of the sheaf $\K$ to the open subscheme $V$ is very flat, while
the sheaf $\cQ$ is a finitely iterated extension of the direct
images of contraadjusted quasi-coherent sheaves from affine open
subschemes in~$X$.
 Let $j\:U\rarrow X$ be an affine open subscheme; we will construct
a short exact sequence $0\rarrow\cP\rarrow\F\rarrow\M\rarrow0$
having the same properties with respect to the open subscheme
$U\cup V\sub X$.

 Pick a short exact sequence $0\rarrow\cR\rarrow\G\rarrow j^*\K
\rarrow0$ of quasi-coherent sheaves on the affine scheme $U$ such that
the sheaf $\G$ is very flat and the sheaf $\cR$ is contraadjusted.
 Consider its direct image $0\rarrow j_*\cR\rarrow j_*\G\rarrow
j_*j^*\K\rarrow 0$ with respect to the affine open embedding~$j$,
and take its pull-back with respect to the adjunction morphism
$\K\rarrow j_*j^*\K$.
 Let $\F$ denote the middle term of the resulting short exact sequence
of quasi-coherent sheaves on~$X$.
 
 By Lemma~\ref{very-open-covering}(a), it suffices to show that
the restrictions of $\F$ to $U$ and $V$ are very flat in order
to conclude that the restriction to $U\cup V$ is.
 We have $j^*\F\simeq\G$, which is very flat by the construction.
 On the other hand, the sheaf $j^*\K$ is very flat over $V\cap U$,
hence so is the sheaf~$\cR$, as the kernel of a surjective map
$\G\rarrow j^*\K$.
 The embedding $U\cap V\rarrow V$ is a very flat affine morphism,
so the sheaf $j_*\cR$ is very flat over~$V$.
 Now it is clear from the short exact sequence $0\rarrow j_*\cR
\rarrow\F\rarrow\K\rarrow 0$ that the sheaf $\F$
is very flat over~$V$.

 Finally, the kernel $\cP$ of the composition of surjective morphisms
$\F\rarrow\K\rarrow\M$ is an extension of the sheaves $\cQ$
and $j_*\cR$, the latter of which is the direct image of
a contraadjusted quasi-coherent sheaf from an affine open subscheme
of $X$, and the former is a finitely iterated extension of such.
\end{proof}

\begin{cor}  \label{quasi-cta-characterizations}
\textup{(a)} A quasi-coherent sheaf\/ $\cP$ on $X$ is contraadjusted if
and only if the functor\/ $\Hom_X({-},\cP)$ takes short exact sequences
of very flat quasi-coherent sheaves on $X$ to short exact
sequences of abelian groups. \par
\textup{(b)} A quasi-coherent sheaf\/ $\cP$ on $X$ is contraadjusted if
and only if\/ $\Ext_X^{>0}(\F,\cP)=0$ for any very flat quasi-coherent
sheaf\/ $\F$ on~$X$. \par
\textup{(c)} The class of contraadjusted quasi-coherent sheaves on $X$
is closed with respect to the passage to the cokernels of injective
morphisms.
\end{cor}

\begin{proof}
 While the condition in part~(a) is \emph{a priori} weaker and
the condition in part~(b) is \emph{a priori} stronger than our
definition of a contraadjusted quasi-coherent sheaf $\cP$ by
the condition $\Ext_X^1(\F,\cP)=0$ for all very flat quasi-coherent
sheaves $\F$, all the three conditions are easily seen to be equivalent
provided that every quasi-coherent sheaf on $X$ is the quotient sheaf
of a very flat one (cf.~\cite[Lemma~1.3]{BHP}, \cite[Lemma~7.1]{PS6},
\cite[Lemma~1.1(a)]{Pal}).
 That much we know from Lemma~\ref{quasi-very-flat-cover}.
 The condition in~(b) clearly has the property~(c).
\end{proof}

\begin{lem}  \label{quasi-cta-envelope}
 Any quasi-coherent sheaf\/ $\M$ on $X$ can be included in a short
exact sequence\/ $0\rarrow\M\rarrow\cP\rarrow\F\rarrow0$, where\/
$\F$ is a very flat quasi-coherent sheaf and\/ $\cP$ is a finitely
iterated extension of the direct images of contraadjusted
quasi-coherent sheaves from affine open subschemes in~$X$.
\end{lem}

\begin{proof}
 Any quasi-coherent sheaf on a quasi-compact quasi-separated scheme
can be embedded into a finite direct sum of direct images of
injective quasi-coherent sheaves from affine open subschemes
constituting a finite covering.
 So an embedding $\M\rarrow\J$ of a sheaf $\M$ into a sheaf $\J$ with
the desired (and even stronger) properties exists, and it remains
to make sure that the quotient sheaf has the desired properties.

 One does this using Lemma~\ref{quasi-very-flat-cover} and
(the construction from) the Salce lemma, see Lemma~\ref{salce-lemma}(b).
 Present the quotient sheaf $\J/\M$ as the quotient sheaf of
a very flat sheaf $\F$ by a subsheaf $\cQ$ representable as
a finitely iterated extension of the desired kind.
 Set $\cP$ to be the fibered product of $\J$ and $\F$ over $\J/\M$;
then $\cP$ is an extension of $\J$ and $\cQ$, and there is
a natural injective morphism $\M\rarrow\cP$ with the cokernel~$\F$.
\end{proof}

\begin{cor}  \label{quasi-very-cta-cor}
\textup{(a)} Any quasi-coherent sheaf on $X$ admits a surjective
map onto it from a very flat quasi-coherent sheaf such that
the kernel is contraadjusted. \par
\textup{(b)} Any quasi-coherent sheaf on $X$ can be embedded into
a contraadjusted quasi-coherent sheaf in such a way that
the cokernel is very flat. \par
\textup{(c)} A quasi-coherent sheaf on $X$ is contraadjusted if and
only if it is a direct summand of a finitely iterated extension of
the direct images of contraadjusted quasi-coherent sheaves from
affine open subschemes of~$X$.
\end{cor}

\begin{proof}
 Parts (a) and~(b) follow from Lemmas~\ref{quasi-very-flat-cover}
and~\ref{quasi-cta-envelope}, respectively.
 The proof of part~(c) uses the argument from
Lemma~\ref{cotorsion-pair-direct-summands-lemma}.
 Given a contraadjusted quasi-coherent sheaf $\cP$, use
Lemma~\ref{quasi-cta-envelope} to embed it into a finitely iterated
extension $\cQ$ of the desired kind in such a way that the cokernel
$\F$ is a very flat quasi-coherent sheaf.
 Since $\Ext^1_X(\F,\cP)=0$ by the definition, we can conclude
that $\cP$ is a direct summand of~$\cQ$.
\end{proof}

\begin{lem} \label{very-flat-contraadjusted-quasi}
 A quasi-coherent sheaf on $X$ is very flat and contraadjusted if and
only if it is a direct summand of a finite direct sum of the direct
images of very flat contraadjusted quasi-coherent sheaves from affine
open subschemes of~$X$.
\end{lem}

\begin{proof}
 The ``if'' assertion is clear.
 To prove ``only if'', notice that the very flat contraadjusted
quasi-coherent sheaves are the injective objects of the exact category
of very flat quasi-coherent sheaves (cf.\
Section~\ref{contraadjusted-exact-cat}).
 So it remains to show that there are enough injectives of the kind
described in the formulation of Lemma in the exact category
$X\qcoh_\vfl$.

 Indeed, let $\F$ be a very flat quasi-coherent sheaf on~$X$ and
$X=\bigcup_\alpha U_\alpha$ be a finite affine open covering.
 Denote by~$j_\alpha$ the identity open embeddings $U_\alpha\rarrow X$.
 For each~$\alpha$, pick an injective morphism $j_\alpha^*\F\rarrow
\G_\alpha$ from a very flat quasi-coherent sheaf $j_\alpha^*\F$ to
a very flat contraadjusted quasi-coherent sheaf $\G_\alpha$ on $U_\alpha$
such that the cokernel $\G_\alpha/j_\alpha^*\F_\alpha$ is a very flat
(using Corollary~\ref{very-rel-proj-inj}(a)).
 Then $\bigoplus_\alpha j_\alpha{}_*\G_\alpha$ is a very flat
contraadjusted quasi-coherent sheaf on $X$ and the cokernel of
the natural morphism $\F\rarrow\bigoplus_\alpha j_\alpha{}_*\G_\alpha$
is very flat (since its restriction to each $U_\alpha$ is).
\end{proof}

\begin{lem} \label{flat-contraadjusted-quasi}
 A quasi-coherent sheaf on $X$ is flat and contraadjusted if and only
if it is a direct summand of a finitely iterated extension of
the direct images of flat contraadjusted quasi-coherent sheaves from
affine open subschemes of~$X$.
\end{lem}

\begin{proof}
 Given a flat quasi-coherent sheaf $\E$, we apply the construction
of Lemma~\ref{very-flat-contraadjusted-quasi} (based on
Theorem~\ref{eklof-trlifaj-very}(a)) in order to embed $\E$
into a direct sum $\bigoplus_\alpha j_\alpha{}_*\G_\alpha$ of the direct
images of flat contraadjusted quasi-coherent sheaves $\G_\alpha$ on
$U_\alpha$ in such a way that the cokernel $\M$ of the morphism
$\E\rarrow\bigoplus_\alpha j_\alpha{}_*\G_\alpha$ is flat.
 Alternatively, one can use the similar construction from the proof
of Lemma~\ref{flat-cotorsion-quasi} below, making $\G_\alpha$ even
flat cotorsion quasi-coherent sheaves on~$U_\alpha$.

 Then we apply the construction of
Lemma~\ref{quasi-very-flat-cover} to the flat quasi-coherent
sheaf~$\M$.
 One can verify step by step that, given a flat quasi-coherent sheaf
$\M$ as the input, the whole construction of
Lemma~\ref{quasi-very-flat-cover} is performed entirely inside
the class/exact category of flat quasi-coherent sheaves, so
the quasi-coherent sheaf $\cP$ it produces is a finitely iterated 
extension of the direct images of flat contraadjusted quasi-coherent
sheaves from affine open subschemes of~$X$.

 Finally, we perform the construction of
Lemma~\ref{quasi-cta-envelope} using the injective morphism
$\E\rarrow\bigoplus_\alpha j_\alpha{}_*\G_\alpha$ with the cokernel $\M$
and the surjective morphism $\cP\rarrow\M$ as the inputs.
 This produces a short exact sequence of quasi-coherent sheaves
$0\rarrow\E\rarrow\cQ\rarrow\F\rarrow0$ with a very flat quasi-coherent
sheaf $\F$ and a quasi-coherent sheaf $\cQ$ that is again a finitely iterated extension of the direct images of flat contraadjusted
quasi-coherent sheaves from~$U_\alpha$.

 Now if the sheaf $\E$ was also contraadjusted, then the short exact
sequence $0\rarrow\E\rarrow\cQ\rarrow\F\rarrow0$ splits by
the definition of contraadjustedness for quasi-coherent sheaves
(cf.\ Corollary~\ref{quasi-cta-characterizations}(b)),
providing the desired result.
\end{proof}

 The following corollary provides equivalent definitions of
contraadjusted and very flat quasi-coherent sheaves on a quasi-compact
semi-separated scheme resembling the corresponding definitions for
modules over a ring in Section~\ref{very-eklof-trlifaj-subsect}.

\begin{cor} \label{quasi-cta-vfl-equivalent-definitions}
\textup{(a)} A quasi-coherent sheaf\/ $\cP$ on $X$ is contraadjusted if
and only if\/ $\Ext_X^{>0}(j_*j^*\O_X,\cP)=0$ for any affine open
embedding of schemes $j\:Y\rarrow X$. \par
\textup{(b)} A quasi-coherent sheaf\/ $\F$ on $X$ is very flat if and
only if\/ $\Ext_X^1(\F,\cP)=0$ for any contraadjusted quasi-coherent
sheaf\/ $\cP$ on~$X$.
\end{cor}

\begin{proof}
 Part~(a): the ``only if'' assertion follows from
Corollary~\ref{quasi-cta-characterizations}(b).
 To prove the ``if'', notice that any very flat quasi-coherent sheaf
$\F$ on $X$ has a finite \v Cech coresolution by finite direct sums of
sheaves of the form $j_*j^*\F$, where $j\:U\rarrow X$ are embeddings of
affine open subschemes.
 Hence the condition $\Ext_X^{>0}(j_*j^*\F,\cP)=0$ for all such~$j$
implies $\Ext_X^{>0}(\F,\cP)=0$.

 Furthermore, a very flat quasi-coherent sheaf $j^*\F$ on $U$
is a direct summand of a transfinitely iterated extension of
the direct images of the structure sheaves of principal affine open
subschemes $V\sub U$ (by Corollary~\ref{very-flat-transfinite}).
 Since the direct images with respect to affine morphisms preserve
transfinitely iterated extensions, it remains to use
the quasi-coherent sheaf version of the result that
$\Ext^1$\+orthogonality is preserved by transfinitely iterated
extensions in the first argument (see~\cite[Lemma~1]{ET}
or Lemma~\ref{eklof-lemma-general}).

 Part~(b): ``only if'' holds by the definition of contraadjusted
sheaves, and ``if'' can be deduced
from Corollary~\ref{quasi-very-cta-cor}(a) by an argument similar to
the proof of Corollary~\ref{very-flat-transfinite}
(and dual to that of Corollary~\ref{quasi-very-cta-cor}(c)).
\end{proof}

 Now we proceed to formulate the analogues of the above assertions
for cotorsion quasi-coherent sheaves.
 The definition of these was given in Section~\ref{fHom-subsection}.
 A quasi-coherent sheaf $\cP$ on an affine scheme $U$ is cotorsion
if and only if the $\O(U)$\+module $\cP(U)$ is cotorsion.
 The class of cotorsion quasi-coherent sheaves is closed under extensions
in the abelian category of quasi-coherent sheaves on a scheme and under
the direct images with respect to affine (and even quasi-compact
quasi-separated) morphisms of schemes.
 We denote the induced exact category structure on the category of
cotorsion quasi-coherent sheaves on a scheme $X$ by $X\qcoh^\cot$.

 The exact category of flat quasi-coherent sheaves on a scheme $X$ is
denoted by $X\qcoh_\fl$.
 Let us also introduce the notation
\begin{align*}
 X\qcoh_\vfl^\cta &= X\qcoh_\vfl\cap X\qcoh^\cta, \\
 X\qcoh_\fl^\cot  &= X\qcoh_\fl\cap X\qcoh^\cot, \\
 X\qcoh_\fl^\cta  &= X\qcoh_\fl\cap X\qcoh^\cta
\end{align*}
for the intersections of classes.
 The exact category structure on $X\qcoh_\fl^\cta$ is inherited from
$X\qcoh$.
 The additive categories $X\qcoh_\vfl^\cta$ and $X\qcoh_\fl^\cot$ are
endowed with trivial exact category structures.
 The additive category of injective quasi-coherent sheaves on~$X$
(also endowed with the trivial exact structure) is denoted by
$X\qcoh^\inj$.

 As above, in the sequel $X$ denotes a quasi-compact semi-separated
scheme.

\begin{lem}  \label{quasi-flat-cover}
 Any quasi-coherent sheaf\/ $\M$ on $X$ can be included in a short
exact sequence\/ $0\rarrow\cP\rarrow\F\rarrow\M\rarrow0$, where\/
$\F$ is a flat quasi-coherent sheaf and\/ $\cP$ is a finitely
iterated extension of the direct images of cotorsion quasi-coherent
sheaves from affine open subschemes in~$X$.
\end{lem}

\begin{proof}
 This is also a special case of
Proposition~\ref{quasi-coherent-gluing-precover}; see
Remark~\ref{quasi-compact-quasi-coherent-remark} below.
 The proof is similar to that of Lemma~\ref{quasi-very-flat-cover},
except that Theorem~\ref{flat-cover-thm}(b) is being used in place of
Theorem~\ref{eklof-trlifaj-very}(b).
\end{proof}

\begin{cor}  \label{quasi-cotors-characterizations}
\textup{(a)} A quasi-coherent sheaf\/ $\cP$ on $X$ is cotorsion if
and only if the functor\/ $\Hom_X({-},\cP)$ takes short exact sequences
of flat quasi-coherent sheaves on $X$ to short exact
sequences of abelian groups. \par
\textup{(b)} A quasi-coherent sheaf\/ $\cP$ on $X$ is cotorsion if
and only if\/ $\Ext_X^{>0}(\F,\cP)=0$ for any flat quasi-coherent
sheaf\/ $\F$ on~$X$. \par
\textup{(c)} The class of cotorsion quasi-coherent sheaves on $X$
is closed with respect to the passage to the cokernels of injective
morphisms.
\end{cor}

\begin{proof}
 Similar to that of Corollary~\ref{quasi-cta-characterizations}.
\end{proof}

\begin{lem}  \label{quasi-cotors-envelope}
 Any quasi-coherent sheaf\/ $\M$ on $X$ can be included in a short
exact sequence\/ $0\rarrow\M\rarrow\cP\rarrow\F\rarrow0$, where\/
$\F$ is a flat quasi-coherent sheaf and\/ $\cP$ is a finitely
iterated extension of the direct images of cotorsion quasi-coherent
sheaves from affine open subschemes in~$X$.
\end{lem}

\begin{proof}
 Similar to that of Lemma~\ref{quasi-cta-envelope}.
\end{proof}

\begin{cor}  \label{quasi-cotors-cor}
\textup{(a)} Any quasi-coherent sheaf on $X$ admits a surjective
map onto it from a flat quasi-coherent sheaf such that the kernel is
cotorsion. \par
\textup{(b)} Any quasi-coherent sheaf on $X$ can be embedded into
a cotorsion quasi-coherent sheaf in such a way that the cokernel is flat.
\par
\textup{(c)} A quasi-coherent sheaf on $X$ is cotorsion if and
only if it is a direct summand of a finitely iterated extension of
the direct images of cotorsion quasi-coherent sheaves from
affine open subschemes of~$X$.
\end{cor}

\begin{proof}
 Similar to that of Corollary~\ref{quasi-very-cta-cor}.
\end{proof}

\begin{lem} \label{flat-cotorsion-quasi}
 A quasi-coherent sheaf on $X$ is flat and cotorsion if and only if it
is a direct summand of a finite direct sum of the direct images of flat
cotorsion quasi-coherent sheaves from affine open subschemes of~$X$.
\end{lem}

\begin{proof}
 Similar to that of Lemma~\ref{very-flat-contraadjusted-quasi}
and based on Theorem~\ref{flat-cover-thm}(a).
\end{proof}

\begin{rem} \label{quasi-compact-quasi-coherent-remark}
 Most of the results above in this section are particular cases of
the discussion in Section~\ref{gluing-cotorsion-in-qcoh-subsect}.
 Specifically, let $\sR$ be the local class of all commutative
rings~$R$ (in the terminology of Section~\ref{local-classes-subsect}).

 For Lemmas and Corollaries~\ref{quasi-very-flat-cover}\+-%
\ref{very-flat-contraadjusted-quasi}
and~\ref{quasi-cta-vfl-equivalent-definitions}(b),
take $\sE_R=\sK_R=R\modl$ to be the abelian category of all
$R$\+modules, $\sF_R=R\modl_\vfl$ be the class of very flat
$R$\+modules, and $\sC(R)=R\modl^\cta$ to be the class of contraadjusted
$R$\+modules (as in Example~\ref{cotorsion-pairs-examples}(2)).
 For Lemmas and Corollaries~\ref{quasi-flat-cover}\+-%
\ref{flat-cotorsion-quasi}, take $\sE_R=\sK_R=R\modl$ to be
the abelian category of all $R$\+modules, $\sF_R=R\modl_\fl$ to be
the class of flat $R$\+modules, and $\sC(R)=R\modl^\cot$ to be
the class of cotorsion $R$\+modules (as in
Example~\ref{cotorsion-pairs-examples}(1)).
 Notice that the classes of modules $\sE_R=R\modl$, \
$\sF_R=R\modl_\vfl$, and $\sF_R=R\modl_\fl$ are very local,
as mentioned in Examples~\ref{local-classes-examples}.

 In this context, Lemmas~\ref{quasi-very-flat-cover}
and~\ref{quasi-flat-cover} are special cases of
Proposition~\ref{quasi-coherent-gluing-precover}.
 Lemmas~\ref{quasi-cta-envelope} and~\ref{quasi-cotors-envelope}
become particular cases of
Proposition~\ref{quasi-coherent-gluing-preenvelope}.
 The results of Corollaries~\ref{quasi-cta-characterizations}(b\+c),
\ref{quasi-very-cta-cor}, \ref{quasi-cta-vfl-equivalent-definitions}(b),
\ref{quasi-cotors-characterizations}(b\+c), \ref{quasi-cotors-cor}
and Lemmas~\ref{very-flat-contraadjusted-quasi},
\ref{flat-cotorsion-quasi} are covered by the assertions of
Theorem~\ref{quasi-coherent-gluing-theorem}.

 Finally, to deduce Lemma~\ref{flat-contraadjusted-quasi} from
Theorem~\ref{quasi-coherent-gluing-theorem}, one needs to consider
the exact subcategory of flat modules $\sE_R=R\modl_\fl$ in
the abelian category of modules $\sK_R=R\modl$, and the hereditary
complete cotorsion pair in $\sE_R$ formed by the classes of
very flat modules $\sF_R=R\modl_\vfl$ and flat contraadjusted modules
$\sC(R)=R\modl_\fl^\cta$, as in
Example~\ref{cotorsion-pairs-examples}(3).
\end{rem} 

 The following result shows that contraadjusted (and in particular,
cotorsion) quasi-coherent sheaves are adjusted to direct images
with respect to nonaffine morphisms of quasi-compact semi-separated
schemes (cf.\ Corollary~\ref{clp-direct} below).
 The first assertions of both parts~(a) and~(b) were explained
already in greater generality in Section~\ref{fHom-subsection};
and the second assertions of both parts follow easily from the results
of the next Section~\ref{dilute-subsect}, such as
Lemma~\ref{direct-sums-injective-contraadjusted-dilute}(b) and
Corollary~\ref{dilute-direct}.
 But we prefer to include also an independent direct proof.

\begin{cor} \label{cta-cot-direct} \hbadness=1700
 Let $f\:Y\rarrow X$ be a morphism of quasi-compact semi-separated
schemes.  Then \par
\textup{(a)} the functor $f_*\:Y\qcoh\rarrow X\qcoh$ takes the full
exact subcategory $Y\qcoh^\cta\sub Y\qcoh$ into the full exact
subcategory $X\qcoh^\cta\sub X\qcoh$, and induces an exact functor
between these exact categories; \par
\textup{(b)} the functor $f_*\:Y\qcoh\rarrow X\qcoh$ takes the full
exact subcategory $Y\qcoh^\cot\sub Y\qcoh$ into the full exact
subcategory $X\qcoh^\cot\sub X\qcoh$, and induces an exact functor
between these exact categories.
\end{cor}

\begin{proof}
 For any flat affine morphism $g\:V\rarrow Y$ into a quasi-compact
semi-separated scheme $Y$, the inverse image functor~$g^*$ takes
quasi-coherent sheaves that can be represented as finitely iterated
extensions of the direct images of quasi-coherent sheaves from
affine open subschemes in $Y$ to quasi-coherent sheaves of
the similar type on~$V$.
 This follows easily from the fact that direct images of quasi-coherent
sheaves with respect to affine morphisms of schemes commute with
inverse images in the base change situations.
 In particular, it follows from Corollary~\ref{quasi-very-cta-cor}(c)
that the functor~$g^*$  takes contraadjusted quasi-coherent sheaves
on $Y$ to quasi-coherent sheaves that are direct summands of
finitely iterated extensions of the direct images of quasi-coherent
sheaves from affine open subschemes $W\sub V$.

 The quasi-coherent sheaves on $V$ that can be represented as such
iterated extensions form a full exact subcategory in
the abelian category of quasi-coherent sheaves.
 The functor of global sections $\Gamma(V,{-})$ is exact on this
exact category.
 Indeed, there is a natural isomorphism of the Ext groups
$\Ext_V^*(\F,h_*\G)\simeq\Ext_W^*(h^*\F,\G)$ for any quasi-coherent
sheaves $\F$ on $V$ and $\G$ on $W$, and a flat affine morphism
$h\:W\rarrow V$ (see, e.~g., \cite[Lemma~1.7(a)]{Pal}
or~\cite[Lemma~6.1]{PS6}).
 Applying this isomorphism in the case when $h$~is the embedding of
an affine open subscheme and $\F=\O_V$, one concludes that
$\Ext_V^{>0}(\O_V,\G)=0$ for all quasi-coherent sheaves $\G$ from
the exact category in question.
 Alternatively, one can refer to
Corollary~\ref{finitely-iterated-extension-dilute} below.

 Specializing to the case of the open subschemes $V=U\times_XY\sub Y$,
where $U$ are affine open subschemes in $X$, we deduce the assertion
that the functor $f_*\:Y\qcoh^\cta\rarrow X\qcoh$ is exact.
 It remains to recall that the direct images of contraadjusted
quasi-coherent sheaves with respect to affine morphisms of schemes
are contraadjusted in order to show that $f_*$~takes
$Y\qcoh^\cta$ to $X\qcoh^\cta$.
 Since the direct images of cotorsion quasi-coherent sheaves with
respect to affine morphisms of schemes are cotorsion, it similarly
follows that $f_*$~takes $Y\qcoh^\cot$ to $X\qcoh^\cot$.
\end{proof}

\subsection{Dilute quasi-coherent sheaves} \label{dilute-subsect}
 Contraadjusted quasi-coherent sheaves are convenient for many
purposes, but they have a problem that an infinite direct sum of
contraadjusted quasi-coherent sheaves need not be contraadjusted.
 To fix this issue, we suggest using the class of \emph{dilute
quasi-coherent sheaves} introduced by Murfet in~\cite[Section~4.1]{M-n}.

 As in Section~\ref{coflasque}, we denote by $H^*(X,{-})$ the functor
of cohomology of sheaves of abelian groups on a topological space~$X$.
 When $X$ is a ringed space, the same cohomology groups can be obtained
as the right derived functor of global sections computed in
the category of sheaves of $\O_X$\+modules.

 Given a quasi-compact semi-separated scheme $X$, we denote by
$\boR^*\Gamma(X,{-})$ the right derived functor of global sections
computed in the abelian category of quasi-coherent sheaves on $X$
(using injective resolutions in $X\qcoh$).
 In fact, the functors $H^*(X,{-})$ and $\boR^*\Gamma(X,{-})$ agree
with each other on any quasi-compact semi-separated scheme~$X$; so
one has $\boR^i\Gamma(X,\M)\simeq H^i(X,\M)$ for all $\M\in X\qcoh$
and $i\ge0$.
 We refer to~\cite[Proposition~B.8]{TT}, \cite[Section~6]{M-n0},
or~\cite[Lemma Tag~01XD or Lemma Tag~0BDY]{SP} for a proof of
this result.

 We will not need to consider codilute cosheaves in this book; but in
order to facilitate an extension of the results of this section to
cosheaves, we avoid using the sheaf cohomology functor $H^*(X,{-})$
and the related results cited above in this section.

 Let $X$ be a quasi-compact semi-separated scheme.
 We will be interested in open subschemes $W\sub X$ for which
the open embedding morphism $j\:W\rarrow X$ is affine.
 Notice that if $X=\bigcup_\alpha U_\alpha$ is an affine open
covering of $X$, then $W=\bigcup_\alpha(W\cap U_\alpha)$ is
an affine open covering of~$W$; so $W$ is quasi-compact.
 If $X$ is affine, then so is~$W$.

 A quasi-coherent sheaf $\M$ on $X$ is said to be
\emph{dilute}~\cite[Section~4.1]{M-n} if, for every open subscheme
$W\sub X$ with affine open embedding morphism $j\:W\rarrow X$,
one has $\boR^{>0}\Gamma(W,\M|_W)=0$.
 Obviously, if this is the case for all $W$, then the quasi-coherent
sheaf $\M|_W=j^*\M$ on $W$ is also dilute.
 Furthermore, any quasi-coherent sheaf on an affine scheme $U$ is
dilute (as the functor $\Gamma(U,{-})$ is exact on $U\qcoh$).

 It is clear that the full subcategory of dilute quasi-coherent sheaves
is closed under extensions and cokernels of injective morphisms in
$X\qcoh$.
 So it inherits an exact category structure; we will denote this full
subcategory with the inherited exact category structure by
$X\qcoh^\dil\sub X\qcoh$.

 The flatness and affineness assumptions in the following lemma can be
dropped, as we will see below in this section.

\begin{lem} \label{flat-affine-direct-image-dilute}
 Let $f\:Y\rarrow X$ be a flat affine morphism of quasi-compact
semi-separated schemes.
 Then the direct image functor $f_*\:Y\qcoh\rarrow X\qcoh$ takes
the full subcategory $Y\qcoh^\dil$ into the full subcategory
$X\qcoh^\dil$.
\end{lem}

\begin{proof}
 Let $\N$ be a quasi-coherent sheaf on $Y$ and $W\sub X$ be an open
subscheme with an affine open embedding morphism $j\:W\rarrow X$.
 Put $T=W\times_XY$; then the open embedding morphism $j'\:T\rarrow Y$
is also affine.
 Denote by $f'\:T\rarrow W$ the natural morphism.
 Then we have $j^*f_*\N\simeq f'_*j'{}^*\N$.
 Replacing $X$ by $W$ and $Y$ by $T$, the question reduces to showing
that $\boR^{>0}\Gamma(Y,\N)=0$ implies $\boR^{>0}\Gamma(X,f_*\N)=0$ for
a flat affine morphism~$f$.
 
 Let $\N\rarrow\J^\bu$ be an injective coresolution of $\N$ in $Y\qcoh$.
 Then the complex $0\rarrow f_*(\N)\rarrow f_*(\J^\bu)$ is exact in
$X\qcoh$ (since $f$~is affine) and $f_*(\J^\bu)$ is a complex of
injective quasi-coherent sheaves on~$X$ (since $f$~is flat).
 The complex of abelian groups $\Gamma(X,f_*\J^\bu)$ is isomorphic to
the complex $\Gamma(Y,\J^\bu)$; so acyclicity of the latter implies
acyclicity of the former.
\end{proof}

\begin{cor} \label{finitely-iterated-extension-dilute}
 Let $X$ be a quasi-compact semi-separated scheme with a finite affine
open covering $X=\bigcup_\alpha U_\alpha$.
 Then any finitely iterated extension of the direct images of
quasi-coherent sheaves on $U_\alpha$ is a dilute quasi-coherent
sheaves on~$X$.
\end{cor}

\begin{proof}
 Follows from Lemma~\ref{flat-affine-direct-image-dilute} applied to
the open embedding morphisms $j_\alpha\:U_\alpha\rarrow X$ in
the role of~$f$.
\end{proof}

\begin{lem} \label{direct-sums-injective-contraadjusted-dilute}
\textup{(a)} Any injective quasi-coherent sheaf on $X$ is dilute.
 Moreover, any infinite direct sum of injective quasi-coherent sheaves
on $X$ is dilute. \par
\textup{(b)} More generally, any contraadjusted quasi-coherent sheaf
on $X$ is dilute.
 Moreover, any infinite direct sum of contraadjusted quasi-coherent
sheaves on $X$ is dilute.
\end{lem}

\begin{proof}
 All the assertions are particular cases of
Corollary~\ref{finitely-iterated-extension-dilute}.
 Specifically, part~(a) holds because all injective quasi-coherent
sheaves on $X$ are direct summands of finite direct sums of the direct
images of (injective) quasi-coherent sheaves from~$U_\alpha$.
 Similarly, part~(b) follows from Corollary~\ref{quasi-very-cta-cor}(c)
(taking into account that it suffices to consider affine open subschemes
$U=U_\alpha$ belonging to a fixed affine open covering $X=\bigcup_\alpha
U_\alpha$ of the scheme $X$ in the latter corollary, and the length of
the finitely iterated extensions in bounded in it).
\end{proof}

 It follows from the first assertion of
Lemma~\ref{direct-sums-injective-contraadjusted-dilute}(a) and
the discussion above that the full subcategory of dilute quasi-coherent
sheaves $X\qcoh^\dil$ is coresolving in $X\qcoh$ (in the sense of
the definition in Section~\ref{infinite-resolutions-subsect}).
 One easily concludes that the cohomology of quasi-coherent sheaves
on a quasi-compact semi-separated scheme $X$ can be computed using
dilute coresolutions.

\begin{cor} \label{direct-sums-dilute}
 The class of dilute quasi-coherent sheaves $X\qcoh^\dil$ is closed
under infinite direct sums in $X\qcoh$.
\end{cor}

\begin{proof}
 Let $\M_i$ be a family of dilute quasi-coherent sheaves on $X$ and
$W\sub X$ be an open subscheme with an affine open embedding morphism
$j\:W\rarrow X$.
 For every index~$i$, choose an injective coresolution
$j^*\M_i\rarrow\J_i^\bu$ of the quasi-coherent sheaf $j^*\M_i$ on~$W$.
 Then the second assertion of
Lemma~\ref{direct-sums-injective-contraadjusted-dilute}(a) tells us
that $\bigoplus_i j^*\M_i\rarrow\bigoplus_i\J_i^\bu$ is a dilute
coresolution of the quasi-coherent sheaf $\bigoplus_i j^*\M_i$ on~$W$.
 Recall that the restriction to $W$ of any dilute quasi-coherent sheaf
on $X$ is dilute.
 Accordingly, the complex $\Gamma(W,\bigoplus_i\J_i^\bu)\simeq
\bigoplus_i\Gamma(W,\J_i^\bu)$ computes the cohomology of
the quasi-coherent sheaf $\bigoplus_i j^*\M_i$ on~$W$, and therefore
$\boR^{>0}\Gamma(W,\bigoplus_i j^*\M_i)\simeq
\bigoplus_i\boR^{>0}\Gamma(W,j^*\M_i)=0$.
\end{proof}

\begin{lem} \label{affine-direct-image-dilute}
 Let $f\:Y\rarrow X$ be an affine morphism of quasi-compact
semi-separated schemes.
 Then the direct image functor $f_*\:Y\qcoh\rarrow X\qcoh$ takes
dilute quasi-coherent sheaves to dilute quasi-coherent sheaves.
\end{lem}

\begin{proof}
 Arguing as in the first paragraph of the proof of
Lemma~\ref{flat-affine-direct-image-dilute}, the question reduces to
showing that $\boR^{>0}\Gamma(Y,\N)=0$ implies
$\boR^{>0}\Gamma(X,f_*\N)=0$ for an affine morphism~$f$.

 Let $X=\bigcup_\alpha U_\alpha$ be a finite affine open covering
of~$X$.
 For any quasi-coherent sheaf $\M$ on $X$, the \v Cech
coresolution~\eqref{cech-quasi} of $\M$ on $X$ is a coresolution
by dilute quasi-coherent sheaves
(by Lemma~\ref{flat-affine-direct-image-dilute}
or Corollary~\ref{finitely-iterated-extension-dilute}).
 So one can compute the derived functor $\boR^*\Gamma(X,{-})$ using
the \v Cech coresolutions (cf.~\cite[Lemma Tag~01XD]{SP}).

 Put $V_\alpha=f^{-1}(U_\alpha)\sub Y$; then $Y=\bigcup_\alpha V_\alpha$
is a finite affine open covering of $Y$, and one can compute
the derived functor $\boR^*\Gamma(Y,{-})$ using the \v Cech
coresolutions for this covering.
 It remains to observe that the resulting \v Cech complexes
$C^\bu(\{U_\alpha\},f_*\N)$ and $C^\bu(\{V_\alpha\},\N)$ computing
$\boR^*\Gamma(X,f_*\N)$ and $\boR^*\Gamma(Y,\N)$ using the affine
coverings $X=\bigcup_\alpha U_\alpha$ and $Y=\bigcup_\alpha V_\alpha$
are isomorphic to each other, $C^\bu(\{U_\alpha\},f_*\N)\simeq
C^\bu(\{V_\alpha\},\N)$.
 Hence $\boR^*\Gamma(X,f_*\N)\simeq\boR^*\Gamma(Y,\N)$
for any quasi-coherent sheaf $\N$ on~$Y$.
\end{proof}

\begin{lem} \label{direct-image-of-injective-is-dilute}
 Let $f\:Y\rarrow X$ be a morphism of quasi-compact semi-separated
schemes.
 Then the direct image functor $f_*\:Y\qcoh\rarrow X\qcoh$ takes
injective quasi-coherent sheaves on $Y$ to dilute quasi-coherent
sheaves on~$X$.
\end{lem}

\begin{proof}
 Let $Y=\bigcup_\beta V_\beta$ be a finite affine open covering of
the scheme~$Y$.
 Then any injective quasi-coherent sheaf on $Y$ is a direct summand
of a finite direct sum of the direct images of injective quasi-coherent
sheaves from~$V_\beta$.
 It remains to observe that the compositions $V_\beta\rarrow Y
\rarrow X$ are affine morphisms of schemes (since $X$ is
semi-separated), and refer to Lemma~\ref{affine-direct-image-dilute}.
\end{proof}

\begin{cor} \label{dilute-direct}
 Let $f\:Y\rarrow X$ be a morphism of quasi-compact semi-separated
schemes.
 Then the direct image functor $f_*\:Y\qcoh\rarrow X\qcoh$ takes
the full exact subcategory $Y\qcoh^\dil\sub Y\qcoh$ into the full exact
subcategory $X\qcoh^\dil\sub X\qcoh$, and induces an exact functor
between these exact categories.
\end{cor}

\begin{proof}
 First let us prove that the direct image functor~$f_*$ takes short
exact sequences in $Y\qcoh^\dil$ to short exact sequences in $X\qcoh$.
 Indeed, let $U$ be an affine open subscheme in~$X$.
 Then $V=f^{-1}(U)$ is an open subscheme in $Y$ with an affine open
embedding morphism $j'\:V\rarrow Y$ (since the open embedding morphism
$j\:U\rarrow X$ is affine, $X$ being semi-separated).
 Now for any quasi-coherent sheaf $\N$ on $Y$ we have
$(f_*\N)(U)=\N(V)$.
 As the functor $\N\longmapsto\Gamma(V,j'{}^*\N)$ is exact on
$Y\qcoh^\dil$ by the definition, it follows that the functor
$\N\longmapsto(f_*\N)(U)$ is exact on $Y\qcoh^\dil$ as well.

 It remains to show that the functor~$f_*$ takes $Y\qcoh^\dil$ to
$X\qcoh^\dil$.
 Arguing as in the first paragraph of the proof of
Lemma~\ref{flat-affine-direct-image-dilute}, the question reduces to
showing that $\boR^{>0}\Gamma(X,f_*\N)=0$ for any dilute quasi-coherent
sheaf $\N$ on~$Y$.
 Let $\N\rarrow\J^\bu$ be an injective coresolution of $\N$ in $Y\qcoh$.
 Then the complex $0\rarrow\N\rarrow\J^\bu$ is exact in the exact
category $Y\qcoh^\dil$ (because the cokernels of injective morphisms
of dilute quasi-coherent sheaves are dilute, and injective
quasi-coherent sheaves are dilute by
Lemma~\ref{direct-sums-injective-contraadjusted-dilute}(a)).
 According to the previous paragraph, it follows that the complex
$0\rarrow f_*\N\rarrow f_*\J^\bu$ is exact in $X\qcoh$.
 By Lemma~\ref{direct-image-of-injective-is-dilute}, \,$f_*\J^\bu$ is
a complex of dilute quasi-coherent sheaves on~$X$.
 So $f_*\J^\bu$ is a dilute coresolution of $f_*\N$ in $X\qcoh$.
 Thus the complex $\Gamma(X,f_*\J^\bu)$ computes
$\boR^*\Gamma(X,f_*\N)$.
 Finally, the natural isomorphism of complexes of abelian groups
$\Gamma(X,f_*\J^\bu)\simeq\Gamma(Y,\J^\bu)$ induces an isomorphism of
the cohomology groups $\boR^{>0}\Gamma(X,f_*\N)\simeq
\boR^{>0}\Gamma(Y,\N)=0$.
\end{proof}

\begin{exs}
 The converse assertion to
Corollary~\ref{finitely-iterated-extension-dilute} is \emph{not} true.
 Let us present two counterexamples to this effect.
 
 (1)~Let $X=\Spec R$ be an affine scheme and $X=\bigcup_\alpha U_\alpha$
be its finite affine open covering.
 Then an arbitrary quasi-coherent sheaf on $X$ \emph{need not} be
a direct summand of finitely iterated extensions of the direct images
of quasi-coherent sheaves from~$U_\alpha$.
 For example, if $X=\Spec\boZ$ is covered by $\Spec\boZ[1/2]$ and
$\Spec\boZ[1/3]$, then the abelian group/$\boZ$\+module $\boZ$ is
\emph{not} a direct summand if a finitely iterated extension of
$\boZ[1/2]$\+modules and $\boZ[1/3]$\+modules~\cite[Example~4.1]{Pal}.

 This example shows that the converse assertion to
Corollary~\ref{finitely-iterated-extension-dilute} does not hold when
the affine open covering $X=\bigcup_\alpha U_\alpha$ is fixed.
 The next counterexample demonstrates that the converse assertion is
still not true when the affine open covering is allowed to vary.

 (2)~Let $X=\bigcup_\alpha U_\alpha$ be a finite affine open covering of
a quasi-compact semi-separated scheme, and let $j_\alpha\:U_\alpha
\rarrow X$ be the open embedding morphisms.
 Then the class of all (direct summands of) finitely iterated extensions
of the direct images of quasi-coherent sheaves from $U_\alpha$ is
preserved by the tensor product functor $\F\ot_{\O_X}{-}\,\nobreak\:
\allowbreak X\qcoh\rarrow X\qcoh$ for any flat quasi-coherent sheaf
$\F$ on~$X$.
 Indeed, the functor $\F\ot_{\O_X}{-}$ is exact, so it preserves
extensions.
 On the other hand, by the projection formula for the affine
morphism~$j_\alpha$, the essential image of the direct image functor
$j_\alpha{}_*\:U_\alpha\qcoh\rarrow X\qcoh$ is preserved by
the tensor product functor $\G\ot_{\O_X}{-}\,\nobreak\:\allowbreak
X\qcoh\rarrow X\qcoh$ for any quasi-coherent sheaf $\G$ on~$X$.

 Now let $X=\mathbb P^1_k$ be the projective line over a field~$k$.
 Then any open subscheme in $X$ is either an affine scheme, or coincides
with the whole of~$X$.
 Therefore, a quasi-coherent sheaf $\M$ on $X$ is dilute if and only if
$\boR^1\Gamma(X,\M)=0$ (notice that one always has
$\boR^i\Gamma(X,\M)=0$ for $i\ge2$).
 In particular, the invertible quasi-coherent sheaf $\O(n)$ on $X$ is
dilute \emph{if and only if} $n\ge-1$.
 Thus the class of dilute quasi-coherent sheaves on $X$ is not even
preserved by the twists (tensor products) with invertible sheaves.
\end{exs}

\subsection{Antilocal contraherent cosheaves}
\label{clp-subsection}
 The terminology of ``antilocality'' and ``antilocal classes''
was introduced in the paper~\cite{Pal}.
 The contraherent or $\bW$\+locally contraherent cosheaves discussed
in the present section were called ``colocally projective''
in~\cite[Section~4.2]{Pcosh}.
 The terminological discussion in~\cite[Section~5.2]{Pphil} suggested
the terminology change to ``antilocal'', and we follow this suggestion
here.

 Let $X$ be a scheme and $\bW$ be its open covering.
 A $\bW$\+locally contraherent cosheaf $\P$ on $X$ is called
\emph{antilocal} if for any short exact sequence
$0\rarrow\gI\rarrow\gJ\rarrow\gK\rarrow0$ of locally injective
$\bW$\+locally contraherent cosheaves on $X$ the short sequence of
abelian groups $0\rarrow\Hom^X(\P,\gI)\rarrow\Hom^X(\P,\gJ)\rarrow
\Hom^X(\P,\gK)\rarrow0$ is exact.

 Obviously, the class of antilocal $\bW$\+locally contraherent cosheaves
on $X$ is closed under direct summands.
 It follows from the adjunction
isomorphism~\eqref{direct-inverse-lin-adjunction} of
Section~\ref{direct-inverse-loc-contra} that the functor of direct
image of $\bT$\+locally contraherent cosheaves~$f_!$ with
respect to any $(\bW,\bT)$\+affine $(\bW,\bT)$\+coaffine morphism
of schemes $f\:Y\rarrow X$ takes antilocal $\bT$\+locally contraherent
cosheaves on $Y$ to antilocal $\bW$\+locally contraherent cosheaves
on~$X$.
 It is also clear that \emph{any} contraherent cosheaf on an affine
scheme $U$ with the covering $\{U\}$ is antilocal.

\begin{lem}  \label{coflasque-clp}
 On any scheme $X$ with an open covering\/~$\bW$, any coflasque
contraherent cosheaf is antilocal.
\end{lem}

\begin{proof}
 First of all, we recall that any coflasque locally contraherent
cosheaf is (globally) contraherent by
Corollary~\ref{coflasque-contraherent}.
 We will prove a somewhat stronger assertion than stated in the lemma:
any short exact sequence $0\rarrow\gI\rarrow\Q\rarrow\gF\rarrow0$
in $X\lcth_\bW$ with $\gF\in X\ctrh_\cfq$ and $\gI\in X\lcth_\bW^\lin$
splits.
 It will follow easily that the functor $\Hom^X(\gF,{-})$ takes
any short exact sequence in $X\lcth_\bW^\lin$ to a short exact
sequence of abelian groups.

 We proceed by applying Zorn's lemma to the partially ordered set
of sections $\phi_Y\:\gF|_Y\rarrow\Q|_Y$ of the morphism of
cosheaves $\Q\rarrow\gF$ defined over open subsets $Y\sub X$.
 Since it suffices to define a morphism of cosheaves on the modules
of cosections over affine open subschemes, which are quasi-compact,
a compatible system of sections $\phi_{Y_i}$ defined over a linearly
ordered family of open subsets $Y_i\sub X$ extends uniquely
to a section over the union $\bigcup_i Y_i$.

 Now let $\phi_Y$ be a section over~$Y$ and $U\sub X$ be an affine
open subscheme subordinate to~$\bW$. 
 Set $V=Y\cap U$; by assumption, the $\O(U)$\+module homomorphism
$\gF[V]\rarrow\gF[U]$ is injective and the $\O(U)$\+module $\gI[U]$
is injective.
 The short exact sequence of $\O(U)$\+modules $0\rarrow\gI[U]\rarrow
\Q[U]\rarrow\gF[U]\rarrow0$ splits, and the difference between two
such splittings $\gF[U]\birarrow\Q[U]$ is an arbitrary $\O(U)$\+module
morphism $\gF[U]\rarrow\gI[U]$.
 The composition $\gF[V]\rarrow\Q[U]$ of the morphism $\phi_Y[V]\:\gF[V]
\rarrow\Q[V]$ and the corestriction morphism $\Q[V]\rarrow\Q[U]$
can therefore be extended to an $\O(U)$\+linear section
$\gF[U]\rarrow\Q[U]$ of the surjection $\Q[U]\rarrow\gF[U]$.
 
 We have constructed a morphism of contraherent cosheaves
$\phi_U\:\gF|_U\rarrow\Q|_U$ whose restriction to~$V$ coincides
with the restriction of the morphism $\phi_Y\:\gF|_Y\rarrow\Q|_Y$.
 Indeed, for any affine open subscheme $W\subset V$, the two
sections $\phi_Y[W]$ and $\phi_U[W]\:\gF[W]\rarrow\Q[W]$ are equal to
each other, because they both form commutative square diagrams with
one and the same section $\gF[U]\rarrow\Q[U]$.
 The point is that two $\O(W)$\+module morphisms $\gF[W]\rarrow\Q[W]$
are equal to each other whenever they have equal compositions with
the corestriction map $\Q[W]\rarrow\Q[U]$ (in other words,
the isomorphism $\Q[W]\simeq\Hom_{\O(U)}(\O(W),\Q[U])$ implies that
the $\O(W)$\+module $\Q[W]$ is cogenerated by~$\Q[U]$).
 Set $Z=Y\cup U$; the pair of morphisms of cosheaves $\phi_Y$ and
$\phi_U$ extends uniquely to a morphism of cosheaves $\phi_Z\:
\gF|_Z\rarrow\Q|_Z$.
 Since the morphisms $\phi_Y$ and $\phi_U$ were some sections of
the surjection $\Q\rarrow\gF$ over $Y$ and $U$, the morphism $\phi_Z$
is a section of this surjection over~$Z$.
\end{proof}

 Generally speaking, according to the above definition, the antilocality
property of a locally contraherent cosheaf $\P$ on a scheme $X$ may
depend not only on the cosheaf $\P$ itself, but also on
the covering~$\bW$.
 No such dependence occurs on quasi-compact semi-separated schemes.
 Indeed, we will see below in this section that on such a scheme
any antilocal $\bW$\+locally contraherent cosheaf is
(globally) contraherent.
 Moreover, the class of antilocal $\bW$\+locally contraherent cosheaves
coincides with the class of antilocal contraherent cosheaves and does
not depend on the covering~$\bW$.

 Let $X$ be a quasi-compact semi-separated scheme and $\bW$ be
its open covering.

\begin{lem}  \label{lin-envelope}
 Let $X=\bigcup_{\alpha=1}^N U_\alpha$ be a finite affine open
covering subordinate to\/~$\bW$.
 Then \par
\textup{(a)} any\/ $\bW$\+locally contraherent cosheaf\/ $\gM$ on $X$
can be included in an (admissible) short exact sequence\/
$0\rarrow\gM\rarrow\gJ\rarrow\P\rarrow0$, where\/ $\gJ$ is a locally
injective\/ $\bW$\+locally contraherent cosheaf on $X$ and\/ $\P$ is
a finitely iterated extension of the direct images of contraherent
cosheaves from the affine open subschemes $U_\alpha\sub X$; \par
\textup{(b)} any locally cotorsion\/ $\bW$\+locally contraherent
cosheaf\/ $\gM$ on $X$ can be included in a short exact sequence\/
$0\rarrow\gM\rarrow\gJ \rarrow\P\rarrow0$, where\/ $\gJ$ is
a locally injective\/ $\bW$\+locally contraherent cosheaf on $X$
and\/ $\P$ is a finitely iterated extension of the direct images of
locally cotorsion contraherent cosheaves from the affine open
subschemes $U_\alpha\sub X$.
\end{lem}

\begin{proof}
 The argument is a dual version of the proofs of
Lemmas~\ref{quasi-very-flat-cover} and~\ref{quasi-flat-cover}.
 Both parts~(a) and~(b) are particular cases of
Proposition~\ref{loc-contraherent-gluing-preenvelope};
see Remark~\ref{to-clp-subsection-remark}.
 Let us prove part~(a); the proof of part~(b) is completely analogous.

 Arguing by induction in $1\le\beta\le N$, we consider the open
subscheme $V=\bigcup_{\alpha<\beta}U_\alpha$ with the induced covering
$\bW|_V=\{V\cap W\mid W\in\bW\}$ and the identity embedding
$h\:V\rarrow X$.
 Assume that we have constructed a short exact sequence $0\rarrow\gM
\rarrow\gK\rarrow\Q\rarrow0$ of $\bW$\+locally contraherent cosheaves
on $X$ such that the restriction $h^!\gK$ of the $\bW$\+locally
contraherent cosheaf $\gK$ to the open subscheme $V\sub X$ is locally
injective, while the cosheaf $\Q$ on $X$ is a finitely iterated
extension of the direct images of contraherent cosheaves from the affine
open subschemes $U_\alpha\sub X$, \ $\alpha<\beta$.
 When $\beta=1$, it suffices to take $\gK=\gM$ and $\Q=0$ for
the induction base.
 Set $U=U_\beta$ and denote by $j\:U\rarrow X$ the identity open
embedding.

 Let $0\rarrow j^!\gK\rarrow\gI\rarrow\R\rarrow0$ be a short exact
sequence of contraherent cosheaves on the affine scheme $U$ such that
the contraherent cosheaf $\gI$ is (locally) injective.
 Consider its direct image $0\rarrow j_!j^!\gK\rarrow j_!\gI\rarrow
j_!\R\rarrow0$ with respect to the affine open embedding~$j$, and
take its push-out with respect to the adjunction morphism
$j_!j^!\gK\rarrow\gK$.
 Let us show that in the resulting short exact sequence $0\rarrow\gK
\rarrow\gJ\rarrow j_!\R\rarrow0$ the $\bW$\+locally contraherent cosheaf
$\gJ$ is locally injective in the restriction to $U\cup V$.
 By Lemma~\ref{cotors-inj-covering}(b), it suffices to show that
the restrictions of $\gJ$ to $U$ and $V$ are locally injective.

 Indeed, in the restriction to $U$ we have $j^!j_!j^!\gK\simeq j^!\gK$,
hence $j^!\gJ\simeq j^!j_!\gI\simeq\gI$ is a (locally) injective
contraherent cosheaf.
 On the other hand, if $j'\:U\cap V\rarrow V$ and $h'\:U\cap V\rarrow U$
denote the embeddings of $U\cap V$, then $h^!j_!\R\simeq j'_!h'{}^!\R$
(as explained in the end of Section~\ref{direct-inverse-loc-contra}).
 Notice that the contraherent cosheaf $h'{}^!j^!\gK\simeq
j'{}^!h^!\gK$ is locally injective (since the $\bW|_V$\+locally
contraherent cosheaf $h^!\gK$ is locally injective by the induction
assumption), hence the contraherent cosheaf $h'{}^!\R$ is locally
injective as the cokernel of the admissible monomorphism of locally
injective contraherent cosheaves $h'{}^!j^!\gK\rarrow h'{}^!\gI$.
 Since the local injectivity of $\bT$\+locally contraherent cosheaves
is preserved by the direct images with respect to flat
$(\bW,\bT)$\+affine morphisms, the contraherent cosheaf $j'_!h'^!\R$
is locally injective, too.
 Now in the short exact sequence $0\rarrow h^!\gK\rarrow h^!\gJ\rarrow
h^!j_!\R\rarrow0$ of $\bW|_V$\+locally contraherent cosheaves on
$V$ the middle term is locally injective, because so are the other
two terms.

 Finally, the composition of admissible monomorphisms of 
$\bW$\+locally contraherent cosheaves $\gM\rarrow\gK\rarrow\gJ$ 
on $X$ is an admissible monomorphism with the cokernel isomorphic
to an extension of the contraherent cosheaves $j_!\R$ and $\Q$,
hence also a finitely iterated extension of the direct images of
contraherent cosheaves from the affine open subschemes
$U_\alpha\sub X$, \ $\alpha\le\beta$.
 The induction step is finished, and the whole lemma is proved.
\end{proof}

 We denote by $\Ext^{X,*}({-},{-})$ the $\Ext$ groups in the exact
category of $\bW$\+locally contraherent cosheaves on~$X$.
 Notice that these do not depend on the covering $\bW$ and coincide
with the $\Ext$ groups in the whole category of locally contraherent
cosheaves $X\lcth$.
 Indeed, the full exact subcategory $X\lcth_\bW$ is closed under
extensions and the passage to kernels of admissible epimorphisms
in $X\lcth$ (see Section~\ref{counterex-subsect}), and for any object
$\P\in X\lcth$ there exists an admissible epimorphism onto $\P$ from
an object of $X\ctrh\sub X\lcth_\bW$ (see
the resolution~\eqref{contraherent-cech} in
Section~\ref{direct-inverse-loc-contra}).
 So the full subcategory $X\lcth_\bW$ is resolving in $X\lcth$, 
and the result of~\cite[Theorem~12.1(b)]{Kel}
or~\cite[Proposition~13.2.2(i)]{KS} can be applied.

 For the same reasons (up to duality), the $\Ext$ groups computed in
the exact subcategories of locally cotorsion and locally injective
$\bW$\+locally contraherent cosheaves $X\lcth_\bW^\lct$ and
$X\lcth_\bW^\lin$ agree with those in $X\lcth_\bW$ (and also in
$X\lcth^\lct$ and $X\lcth^\lin$).
 Indeed, the full exact subcategories $X\lcth_\bW^\lct$ and
$X\lcth_\bW^\lin$ are closed under extensions and the passage to
cokernels of admissible monomorphisms in $X\lcth_\bW$
(see Section~\ref{locally-contraherent}), and we have just constructed
in Lemma~\ref{lin-envelope} an admissible monomorphism from
any $\bW$\+locally contraherent cosheaf to a locally injective one.
 In other words, the full subcategories $X\lcth_\bW^\lct$ and
$X\lcth_\bW^\lin$ are coresolving in $X\lcth_\bW$.
 We refer to Proposition~\ref{fully-faithful-prop}
or~\ref{infinite-resolutions}(a) for further details.

\begin{cor}  \label{clp-characterizations}
\textup{(a)} A\/ $\bW$\+locally contraherent cosheaf\/ $\P$ on $X$
is antilocal if and only if\/ $\Ext^{X,1}(\P,\gJ)=0$
and if and only if\/ $\Ext^{X,>0}(\P,\gJ)=0$ for all locally
injective\/ $\bW$\+locally contraherent cosheaves\/ $\gJ$ on $X$. \par
\textup{(b)} The class of antilocal\/ $\bW$\+locally contraherent
cosheaves on $X$ is closed under extensions and the passage to kernels
of admissible epimorphisms in the exact category\/ $X\lcth_\bW$.
\end{cor}

\begin{proof}
 Part~(a) follows from the existence of an admissible monomorphism from
any $\bW$\+locally contraherent cosheaf on $X$ into a locally injective
$\bW$\+locally contraherent cosheaf (a weak form of
Lemma~\ref{lin-envelope}(a)), cf.~\cite[Lemma~1.3]{BHP},
\cite[Lemma~7.1]{PS6}, or~\cite[Lemma~1.1(b)]{Pal}.
 Part~(b) follows from part~(a).
\end{proof}

\begin{lem}  \label{clp-cover}
 Let $X=\bigcup_\alpha U_\alpha$ be a finite affine open
covering subordinate to\/~$\bW$.
 Then \par
\textup{(a)} any\/ $\bW$\+locally contraherent cosheaf\/ $\gM$ on $X$
can be included in a short exact sequence\/ $0\rarrow\gJ\rarrow\P
\rarrow\gM\rarrow0$, where\/ $\gJ$ is a locally injective\/
$\bW$\+locally contraherent cosheaf on $X$ and\/ $\P$ is a finitely
iterated extension of the direct images of contraherent cosheaves
from the affine open subschemes $U_\alpha\sub X$; \par
\textup{(b)} any locally cotorsion\/ $\bW$\+locally contraherent
cosheaf\/ $\gM$ on $X$ can be included in a short exact sequence\/
$0\rarrow\gJ\rarrow\P\rarrow\gM\rarrow0$, where\/ $\gJ$ is a locally
injective\/ $\bW$\+locally contraherent cosheaf on $X$ and\/ $\P$ is
a finitely iterated extension of the direct images of locally 
cotorsion contraherent cosheaves from the affine open subschemes
$U_\alpha\sub X$.
\end{lem}

\begin{proof}
 There is an admissible epimorphism
$\bigoplus_\alpha j_\alpha{}_!j_\alpha^!\gM\rarrow\gM$
(see~\eqref{contraherent-cech} for the notation and explanation)
onto any $\bW$\+locally contraherent cosheaf $\gM$ from
a finite direct sum of the direct images of contraherent cosheaves
from the affine open subschemes~$U_\alpha$.
 When $\gM$ is a locally cotorsion $\bW$\+locally contraherent
cosheaf, this is an admissible epimorphism in the category of
locally cotorsion $\bW$\+locally contraherent cosheaves, and
$j_\alpha^!\gM$ are locally cotorsion contraherent cosheaves
on~$U_\alpha$.

 Given that, the desired short exact sequences in Lemma can be obtained
from those of Lemma~\ref{lin-envelope} by (the construction from)
the Salce lemma; see Lemma~\ref{salce-lemma}(a) (cf.\ the proofs of
Lemmas~\ref{eklof-trlifaj-cta}, \ref{quasi-cta-envelope},
and~\ref{quasi-cotors-envelope}).
\end{proof}

\begin{cor} \label{clp-cor}
\textup{(a)} For any\/ $\bW$\+locally contraherent cosheaf\/ $\gM$ on
$X$ there exists an admissible monomorphism from\/ $\gM$ into a locally
injective\/ $\bW$\+locally contraherent cosheaf\/ $\gJ$ on $X$ such
that the cokernel\/ $\P$ is an antilocal\/ $\bW$\+locally
contraherent cosheaf. \par
\textup{(b)} For any\/ $\bW$\+locally contraherent cosheaf\/ $\gM$ on
$X$ there exists an admissible epimorphism onto\/ $\gM$ from
an antilocal\/ $\bW$\+locally contraherent cosheaf\/ $\P$ on $X$ such
that the kernel\/ $\gJ$ is a locally injective\/ $\bW$\+locally
contraherent cosheaf. \par
\textup{(c)} Let $X=\bigcup_\alpha U_\alpha$ be a finite affine open
covering subordinate to\/~$\bW$.
 Then a $\bW$\+locally contraherent cosheaf on $X$ is antilocal
if and only if it is (a contraherent cosheaf and) a direct summand
of a finitely iterated extension of the direct images of contraherent
cosheaves from the affine open subschemes $U_\alpha\sub X$.
\end{cor}

\begin{proof}
 The ``if'' assertion in part~(c) follows from
Corollary~\ref{clp-characterizations}(b) together with our
preliminary remarks in the beginning of this section.
 This having been shown, part~(a) follows from
Lemma~\ref{lin-envelope}(a) and part~(b) from Lemma~\ref{clp-cover}(a).

 The ``only if'' assertion in~(c) follows from
Corollary~\ref{clp-characterizations}(a) and Lemma~\ref{clp-cover}(a)
by the argument from the proof of Corollary~\ref{very-flat-transfinite}
(cf.\ Corollaries~\ref{quasi-very-cta-cor}(c)
and~\ref{quasi-cotors-cor}(c)); see
Lemma~\ref{cotorsion-pair-direct-summands-lemma} for the general
formulation.
 Notice that the functors of direct image with respect to the open
embeddings $U_\alpha\rarrow X$ take contraherent cosheaves to
contraherent cosheaves, and the full subcategory of contraherent
cosheaves $X\ctrh\sub X\lcth$ is closed under extensions.
\end{proof}

 By an antilocal locally cotorsion $\bW$\+locally contraherent
cosheaf we will mean a $\bW$\+locally contraherent cosheaf that is
simultaneously antilocal and locally cotorsion.

\begin{cor} \label{clp-lct-cor}
\textup{(a)} For any locally cotorsion\/ $\bW$\+locally contraherent
cosheaf\/ $\gM$ on $X$ there exists an admissible monomorphism from\/
$\gM$ into a locally injective\/ $\bW$\+locally contraherent cosheaf\/
$\gJ$ on $X$ such that the cokernel\/ $\P$ is an antilocal
locally cotorsion\/ $\bW$\+locally contraherent cosheaf. \par
\textup{(b)} For any\/ locally cotorsion $\bW$\+locally contraherent
cosheaf\/ $\gM$ on $X$ there exists an admissible epimorphism onto\/
$\gM$ from an antilocal locally cotorsion\/ $\bW$\+locally
contraherent cosheaf\/ $\P$ on $X$ such that the kernel\/ $\gJ$ is
a locally injective\/ $\bW$\+locally contraherent cosheaf. \par
\textup{(c)} Let $X=\bigcup_\alpha U_\alpha$ be a finite affine open
covering subordinate to\/~$\bW$.
 Then a locally cotorsion\/ $\bW$\+locally contraherent cosheaf on $X$
is antilocal if and only if it is (a contraherent cosheaf and) a direct
summand of a finitely iterated extension of the direct images of
locally cotorsion contraherent cosheaves from the affine open
subschemes $U_\alpha\sub X$.
\end{cor}

\begin{proof}
 Same as Corollary~\ref{clp-cor}, except that parts~(b) of
Lemmas~\ref{lin-envelope} and~\ref{clp-cover} need to be used.
 Parts~(a\+b) can be also easily deduced from
Corollary~\ref{clp-cor}(a\+b).
\end{proof}

\begin{cor} \label{clp-independence}
 The full subcategory of antilocal\/ $\bW$\+locally contraherent
cosheaves in the exact category of all locally contraherent cosheaves
on $X$ does not depend on the choice of the open covering\/~$\bW$.
\end{cor}

\begin{proof}
 Given two open coverings $\bW'$ and $\bW''$ of the scheme~$X$, pick
a finite affine open covering $X=\bigcup_{\alpha=1}^N U_\alpha$
subordinate to both $\bW'$ and $\bW''$, and apply
Corollary~\ref{clp-cor}(c).
\end{proof}

 As a full subcategory closed under extensions and kernels of admissible
epimorphisms in $X\ctrh$, the category of antilocal contraherent
cosheaves on $X$ acquires the induced exact category structure.
 We denote this exact category by $X\ctrh_\al$.
 The (similarly constructed) exact category of antilocal locally
cotorsion contraherent cosheaves on $X$ is denoted by $X\ctrh_\al^\lct$.

 The full subcategory of contraherent cosheaves that are simultaneously
antilocal and locally injective will be denoted by
$X\ctrh_\al^\lin$.
 Clearly, any extension of two objects from this subcategory is
trivial in $X\ctrh$, so the category of antilocal locally injective
contraherent cosheaves is naturally viewed as an additive category
endowed with the trivial exact category structure.

 It follows from Corollary~\ref{clp-cor}(a\+b) that the additive
category $X\ctrh_\al^\lin$ is simultaneously the category of
projective objects in $X\lcth_\bW^\lin$ and the category of
injective objects in $X\ctrh_\al$, and that it has enough of both
such projectives and injectives.

\begin{cor}  \label{clp-lin}
 Let $X=\bigcup_\alpha U_\alpha$ be a finite affine open covering.
 Then a contraherent cosheaf on $X$ is antilocal and locally
injective if and only if it is isomorphic to a direct summand of
a finite direct sum of the direct images of (locally) injective
contraherent cosheaves from the open embeddings $U_\alpha\rarrow X$.
\end{cor}

\begin{proof}
 For any locally injective $\bW$\+locally contraherent cosheaf $\gJ$
on $X$, the map $\bigoplus_\alpha j_\alpha{}_! j_\alpha^!\gJ
\rarrow\gJ$ is an admissible epimorphism in the category of
locally injective $\bW$\+locally contraherent cosheaves.
 Now if $\gJ$ is also antilocal, then the extension splits.
\end{proof}

\begin{rem} \label{to-clp-subsection-remark}
 Most of the results above in this section are particular cases of
the discussion in Section~\ref{gluing-cotorsion-in-lcth-subsect}.
 Specifically, as in Remark~\ref{quasi-compact-quasi-coherent-remark},
let $\sR$ be the local class of all commutative rings~$R$.

 For Lemmas and Corollaries~\ref{lin-envelope}(a),
\ref{clp-characterizations}(b), \ref{clp-cover}(a), \ref{clp-cor},
\ref{clp-independence}, and~\ref{clp-lin}, take $\sE^R=\sK^R=
R\modl^\cta$ to be the exact category of contraadjusted $R$\+modules,
$\sF(R)=\sE^R=R\modl^\cta$, and $\sC^R=R\modl^\inj$ to be the class
of injective $R$\+modules; so $(\sF(R),\sC^R)$ is a trivial hereditary
complete cotorsion pair in~$\sE^R$.
 For Lemmas and Corollaries~\ref{lin-envelope}(b), \ref{clp-cover}(b),
and~\ref{clp-lct-cor}, take $\sE^R=R\modl^\cot\sub R\modl^\cta=\sK^R$
to be the exact  category of cotorsion $R$\+modules, $\sF(R)=\sE^R=
R\modl^\cot$, and $\sC^R=R\modl^\inj$ to be the class of injective
$R$\+modules; once again, $(\sF(R),\sC^R)$ is a trivial hereditary
complete cotorsion pair in~$\sE^R$.
 Notice that the classes of modules $\sE^R=R\modl^\cta$, \
$\sE^R=R\modl^\cot$, and $\sC^R=R\modl^\inj$ are very colocal, as
mentioned in Examples~\ref{colocal-classes-examples}.

 In this context, Lemma~\ref{lin-envelope}(a\+b) presents two particular
cases of Proposition~\ref{loc-contraherent-gluing-preenvelope}.
 Lemma~\ref{clp-cover}(a\+b) gives two particular cases of
Proposition~\ref{loc-contraherent-gluing-precover}.
 The results of Corollaries~\ref{clp-characterizations}(b)
and~\ref{clp-cor}\+-\ref{clp-lin} are covered by the assertions
of Theorem~\ref{loc-contraherent-gluing-theorem} (in view of
the first assertion of Corollary~\ref{clp-characterizations}(a)).
\end{rem}

\begin{cor}  \label{clp-products}
 The three full subcategories of antilocal cosheaves
$X\ctrh_\al$, \ $X\ctrh^\lct_\al$, and $X\ctrh^\lin_\al$ are closed
under infinite products in the category $X\ctrh$.
\end{cor}

\begin{proof}
 The assertions are easily deduced from the descriptions of the full
subcategories of antilocal cosheaves given in
Corollaries~\ref{clp-cor}(c), \ref{clp-lct-cor}(c), and~\ref{clp-lin}
together with the fact that the functor of direct image of
contraherent cosheaves with respect to an affine morphism of
schemes preserves infinite products.
\end{proof}

\begin{cor} \label{clp-inverse}
 Let $f\:Y\rarrow X$ be an affine morphism of quasi-compact
semi-separated schemes.  Then \par
\textup{(a)} the functor of inverse image of locally injective
locally contraherent cosheaves $f^!\:X\lcth^\lin\rarrow
Y\lcth^\lin$ takes the full subcategory $X\ctrh^\lin_\al$ into
$Y\ctrh^\lin_\al$; \par
\textup{(b)} assuming that the morphism~$f$ is also flat, the functor
of inverse image of locally cotorsion locally contraherent
cosheaves $f^!\:X\lcth^\lct\rarrow Y\lcth^\lct$ takes the full
subcategory $X\ctrh^\lct_\al$ into $Y\ctrh^\lct_\al$; \par
\textup{(c)} assuming that the morphism~$f$ is also very flat,
the functor of inverse image of locally contraherent cosheaves
$f^!\:X\lcth \rarrow Y\lcth$ takes the full subcategory $X\ctrh_\al$
into $Y\ctrh_\al$. \par
\end{cor}

\begin{proof}
 Parts~(a\+c) follow from Corollaries~\ref{clp-cor}(c),
\ref{clp-lct-cor}(c), and~\ref{clp-lin}, respectively, together
with the base change results from the second half of
Section~\ref{direct-inverse-loc-contra}.
 The point is that, in all the three contexts, the inverse image
functors~$f^!$ are exact and commute with the direct images with respect
to flat affine morphisms (in particular, with the direct image
functors~$j_!$ with respect to the embeddings of affine open subschemes
$j\:U\rarrow X$).
 Furthermore, the open subscheme $V=U\times_XY\sub Y$ is affine
for every affine open subscheme $U\sub X$.
\end{proof}

\subsection{Antilocally flat contraherent cosheaves}
\label{clf-subsection}
 Let $X$ be a scheme and $\bW$ be its open covering.
 A $\bW$\+locally contraherent cosheaf $\gF$ on $X$ is called
\emph{antilocally flat} if for any short exact sequence $0\rarrow\P
\rarrow\Q\rarrow\R\rarrow0$ of locally cotorsion $\bW$\+locally
contraherent cosheaves on $X$ the short sequence of abelian groups
$0\rarrow\Hom^X(\gF,\P)\rarrow\Hom^X(\gF,\Q)\rarrow\Hom^X(\gF,\R)
\rarrow0$ is exact.

 Clearly, any antilocally flat $\bW$\+locally contraherent cosheaf
is antilocal.
 It follows from the adjunction
isomorphism~\eqref{direct-inverse-lct-adjunction} that the functor
of direct image of $\bT$\+locally contraherent cosheaves~$f_!$
with respect a flat $(\bW,\bT)$\+affine $(\bW,\bT)$\+coaffine
morphism of schemes $f\:Y\rarrow X$ takes antilocally flat
$\bT$\+locally contraherent cosheaves on $Y$ to antilocally flat
$\bW$\+locally contraherent cosheaves on~$X$.
 A contraherent cosheaf $\gF$ on an affine scheme $U$ with
the covering $\{U\}$ is antilocally flat if and only if
the contraadjusted $\O(U)$\+module $\gF[U]$ is flat
(the ``if'' is obvious, and ``only if'' can be obtained from
Theorem~\ref{flat-cover-thm}(b) or
Example~\ref{cotorsion-pairs-examples}(4)
(cf.\ Corollary~\ref{clf-cor}(c) below).

 Let $X$ be a quasi-compact semi-separated scheme.
 It follows from the results of Section~\ref{clp-subsection} that
any antilocally flat $\bW$\+locally contraherent cosheaf on $X$
is contraherent.
 We will see below in this section that the class of antilocally flat
$\bW$\+locally contraherent cosheaves on $X$ coincides with
the class of antilocally flat contraherent cosheaves and does not
depend on the covering~$\bW$.

\begin{lem}  \label{lct-envelope}
 Let $X=\bigcup_\alpha U_\alpha$ be a finite affine open covering
subordinate to\/~$\bW$.
 Then any\/ $\bW$\+locally contraherent cosheaf\/ $\gM$ on $X$ can be
included in a short exact sequence $0\rarrow\gM\rarrow\P\rarrow\gF
\rarrow0$, where\/ $\P$ is a locally cotorsion\/ $\bW$\+locally
contraherent cosheaf on $X$ and\/ $\gF$ is a finitely iterated extension
of the direct images of contraherent cosheaves on $U_\alpha$
corresponding to flat contraadjusted\/ $\O(U_\alpha)$\+modules.
\end{lem}

\begin{proof}
 This is also a particular case of
Proposition~\ref{loc-contraherent-gluing-preenvelope};
see Remark~\ref{to-clf-subsection-remark} below.
 The proof is similar to that Lemma~\ref{lin-envelope}, except that
Theorem~\ref{flat-cover-thm}(a) needs to be used to resolve
a contraherent cosheaf on an affine open subscheme $U\sub X$
(cf.\ the proof of Lemma~\ref{quasi-very-flat-cover}).
 Besides, one has to use Lemma~\ref{cotors-inj-covering}(a) and
the fact that the class of locally cotorsion contraherent cosheaves
is preserved by direct images with respect to affine morphisms.
\end{proof}

\begin{cor}  \label{clf-characterizations}
\textup{(a)} A\/ $\bW$\+locally contraherent cosheaf\/ $\gF$ on $X$
is antilocally flat if and only if\/ $\Ext^{X,1}(\gF,\P)=0$ and
if and only if\/ $\Ext^{X,>0}(\gF,\P)=0$ for all locally cotorsion\/
$\bW$\+locally contraherent cosheaves\/ $\P$ on~$X$. \par
\textup{(b)} The class of antilocally flat\/ $\bW$\+locally contraherent
cosheaves on $X$ is closed under extensions and the passage to
kernels of admissible epimorphisms in the exact category
$X\lcth_\bW$.
\end{cor}

\begin{proof}
 Similar to the proof of Corollary~\ref{clp-characterizations}.
\end{proof}

\begin{lem}  \label{clf-cover}
 Let $X=\bigcup_\alpha U_\alpha$ be a finite affine open covering
subordinate to\/~$\bW$.
 Then any\/ $\bW$\+locally contraherent cosheaf\/ $\gM$ on $X$ can be
included in a short exact sequence $0\rarrow\P\rarrow\gF\rarrow\gM
\rarrow0$, where $\P$ is a locally cotorsion\/ $\bW$\+locally
contraherent cosheaf on $X$ and\/ $\gF$ is a finitely iterated extension
of the direct images of contraherent cosheaves on $U_\alpha$
corresponding to flat contraadjusted\/ $\O(U_\alpha)$\+modules.
\end{lem}

\begin{proof}
 The proof is similar to that of Lemma~\ref{clp-cover} and based
on Lemma~\ref{lct-envelope}.
 The key is to show that there is an admissible epimorphism in
the exact category $X\lcth_\bW$ onto any $\bW$\+locally contraherent
cosheaf $\gM$ from a finitely iterated extension (in fact, even
a finite direct sum) of the direct images of contraherent cosheaves
on $U_\alpha$ corresponding to flat contraadjusted\/
$\O(U_\alpha)$\+modules.

 Here it suffices to pick admissible epimorphisms from such
contraherent cosheaves $\gF_\alpha$ on $U_\alpha$ onto
the restrictions $j_\alpha^!\gM$ of $\gM$ to $U_\alpha$ and consider
the corresponding morphism $\bigoplus_\alpha j_\alpha{}_!\gF_\alpha
\rarrow\gM$ of $\bW$\+locally contraherent cosheaves on~$X$.
 To check that this is an admissible epimorphism, one can, e.~g.,
notice that it is so in the restriction to each $U_\alpha$ and
recall that being an admissible epimorphism of $\bW$\+locally
contraherent cosheaves is a local property
(see Section~\ref{counterex-subsect}).

 Notice that a morphism $\gG\rarrow\gN$ of contraherent cosheaves on
an affine scheme $U$ being an admissible epimorphism means that
the map of $\O(U)$\+modules $\gG[U]\rarrow\gN[U]$ is surjective with
a contraadjusted kernel.
 So one needs to use either Theorem~\ref{eklof-trlifaj-very}(b) or
Theorem~\ref{flat-cover-thm}(b) in order to establish the existence
of the required admissible epimorphisms $\gF_\alpha\rarrow
j_\alpha^!\gM$.
\end{proof}

\begin{cor} \label{clf-cor}
\textup{(a)} For any\/ $\bW$\+locally contraherent cosheaf\/ $\gM$ on $X$
there exists an admissible monomorphism from\/ $\gM$ into a locally
cotorsion\/ $\bW$\+locally contraherent cosheaf\/ $\P$ on $X$ such that
the cokernel\/ $\gF$ is an antilocally flat\/ $\bW$\+locally contraherent
cosheaf. \par
\textup{(b)} For any\/ $\bW$\+locally contraherent cosheaf\/ $\gM$ on $X$
there exists an admissible epimorphism onto\/ $\gM$ from an antilocally
flat\/ $\bW$\+locally contraherent cosheaf\/ $\gF$ on $X$ such that
the kernel\/ $\P$ is a locally cotorsion\/ $\bW$\+locally contraherent
cosheaf. \par
\textup{(c)} Let $X=\bigcup_\alpha U_\alpha$ be a finite affine open
covering subordinate to\/~$\bW$.
 Then a\/ $\bW$\+locally contraherent cosheaf on $X$ is antilocally flat
if and only if it is a direct summand of a finitely iterated extension
of the direct images of contraherent cosheaves on $U_\alpha$
corresponding to flat contraadjusted\/ $\O(U_\alpha)$\+modules.
\end{cor}

\begin{proof}
 Similar to the proof of Corollary~\ref{clp-cor} and based on
Lemmas~\ref{lct-envelope}, \ref{clf-cover} and
Corollary~\ref{clf-characterizations}.
\end{proof}

\begin{cor} \label{clf-independence}
 The full subcategory of antilocally flat\/ $\bW$\+locally contraherent
cosheaves in the exact category $X\lcth$ does not depend on
the open covering\/~$\bW$.
\end{cor}

\begin{proof}
 Similar to the proof of Corollary~\ref{clp-independence} and
based on Corollary~\ref{clf-cor}(c).
\end{proof}

\begin{rem} \label{to-clf-subsection-remark}
 The results above in this section are also particular cases of
the discussion in Section~\ref{gluing-cotorsion-in-lcth-subsect}.
 As in Remark~\ref{to-clp-subsection-remark}, let $\sR$ be the local
class of all commutative rings~$R$ (in the sense of
Section~\ref{local-classes-subsect}).

 Take $\sE^R=\sK^R=R\modl^\cta$ to be the exact category of 
contraadjusted $R$\+modules, $\sF(R)=R\modl_\fl^\cta$ to be the class
of flat contraadjusted $R$\+modules, and $\sC^R=R\modl^\cot$ to be
the class of cotorsion $R$\+modules (as in
Example~\ref{cotorsion-pairs-examples}(4)).
 Notice that both the classes of $R$\+modules $\sE^R=R\modl^\cta$
and $\sC^R=R\modl^\cot$ are very colocal, as mentioned in
Examples~\ref{colocal-classes-examples}.

 In this context, Lemma~\ref{lct-envelope} is a particular case
of Proposition~\ref{loc-contraherent-gluing-preenvelope}.
 Lemma~\ref{clf-cover} is a particular case of
Proposition~\ref{loc-contraherent-gluing-precover}.
 The results of Corollaries~\ref{clf-characterizations}(b),
\ref{clf-cor}, and~\ref{clf-independence} are covered by
the assertions of Theorem~\ref{loc-contraherent-gluing-theorem}
(in view of the first assertion of
Corollary~\ref{clf-characterizations}(a)).
\end{rem}

 As a full subcategory closed under extensions and kernels of
admissible epimorphisms in $X\ctrh$, the category of antilocally flat
contraherent cosheaves on $X$ acquires the induced exact category
structure, which we denote by $X\ctrh_\alf$.

\medskip

 We refer to Section~\ref{contraherent-tensor} for the definition of
a coherent scheme.
 The following corollary is to be compared with
Corollaries~\ref{finite-krull-flat-contraherent}(a)
and~\ref{finite-krull-flat-clf-cor}(a).

\begin{cor}  \label{clf-noetherian-flat}
 Any antilocally flat contraherent cosheaf over a semi-separated
coherent scheme is flat.
\end{cor}

\begin{proof}
 Follows from Corollary~\ref{clf-cor}(c) together with the remarks
about flat contraherent cosheaves over coherent affine schemes and
the direct images of flat cosheaves of $\O$\+modules in
Section~\ref{contraherent-tensor}.
 Specifically, the contraherent cosheaf $\widecheck F$ over
a coherent affine scheme $U$ is flat whenever the corresponding
contraadjusted $\O(U)$\+module $F$ is flat; and the class of flat
cosheaves of $\O$\+modules on schemes is preserved by extensions
and direct images with respect to flat affine morphisms.
\end{proof}

\begin{cor}  \label{clf-products}
 Over a semi-separated coherent scheme $X$, the full subcategory
$X\ctrh_\alf$ is closed with respect to infinite products in $X\ctrh$.
\end{cor}

\begin{proof}
 In addition to what has been said in the proof of
Corollary~\ref{clp-products}, it is also important here that
infinite products of flat modules over a coherent ring are flat.
\end{proof}

\begin{cor}  \label{clf-restriction}
 Let $X$ be a semi-separated coherent scheme and $j\:Y\rarrow X$
be an affine open embedding.
 Then the inverse image functor~$j^!$ takes antilocally flat
contraherent cosheaves to antilocally flat contraherent cosheaves.
\end{cor}

\begin{proof}
 Similar to the proof of Corollary~\ref{clp-inverse} and based on
Corollary~\ref{clf-cor}(c).
 The only difference is that one has to use also
Corollary~\ref{coherent-flat-local}(a) to establish the desired
assertion in the case of affine schemes $X$ and $Y$ (from which
the general case is then deduced using Corollary~\ref{clf-cor}(c)).
\end{proof}

\subsection{Projective contraherent cosheaves}
\label{projective-contraherent}
 Let $X$ be a quasi-compact semi-sep\-a\-rated scheme and $\bW$ be its
affine open covering.

\begin{lem}  \label{loc-contra-proj}
\textup{(a)} The exact category of\/ $\bW$\+locally contraherent
cosheaves on $X$ has enough projective objects. \par
\textup{(b)} Let $X=\bigcup_\alpha U_\alpha$ be a finite affine open
covering subordinate to\/~$\bW$.
 Then a\/ $\bW$\+locally contraherent cosheaf on $X$ is projective
if and only if it is a direct summand of a direct sum over\/~$\alpha$
of the direct images of contraherent cosheaves on $U_\alpha$
corresponding to very flat contraadjusted\/ $\O(U_\alpha)$\+modules.
\end{lem}

\begin{proof}
 The assertion ``if'' in part~(b) follows from the adjunction of
the direct and inverse image functors for the embeddings
$U_\alpha\rarrow X$ together with the fact that the very flat
contraadjusted modules are the projective objects of the exact
categories of contraadjusted modules over $\O(U_\alpha)$
(see Section~\ref{contraadjusted-exact-cat}).

 It remains to show that there exists an admissible epimorphism onto
any $\bW$\+locally contraherent cosheaf $\gM$ on $X$ from
a direct sum of the direct images of contraherent cosheaves on
$U_\alpha$ corresponding to very flat contraadjusted modules.
 The construction is similar to the one used in the proof of
Lemma~\ref{clf-cover} and based on Theorem~\ref{eklof-trlifaj-very}(b).
 One picks admissible epimorphisms from contraherent cosheaves
$\gF_\alpha$ of the desired kind on the affine schemes $U_\alpha$
onto the restrictions $j_\alpha^!\gM$ of the cosheaf $\gM$ and
considers the corresponding morphism
$\bigoplus_\alpha j_\alpha{}_!\gF_\alpha\rarrow\gM$.
\end{proof}

\begin{cor}  \label{ctrh-lcth-proj}
\textup{(a)} There are enough projective objects in the exact
category $X\lcth$ of locally contraherent cosheaves on $X$,
and all these projective objects belong to the full subcategory
of contraherent cosheaves $X\ctrh\sub X\lcth$. \par
\textup{(b)} The full subcategories of projective objects in
the three exact categories $X\ctrh\sub X\lcth_\bW\sub X\lcth$
coincide. \qed
\end{cor}

\begin{lem}  \label{loc-lct-proj}
\textup{(a)} The exact category of\/ locally cotorsion $\bW$\+locally
contraherent cosheaves on $X$ has enough projective objects. \par
\textup{(b)} Let $X=\bigcup_\alpha U_\alpha$ be a finite affine open
covering subordinate to\/~$\bW$.
 Then a locally cotorsion\/ $\bW$\+locally contraherent cosheaf on $X$
is projective if and only if it is a direct summand of a direct sum
over\/~$\alpha$ of the direct images of locally cotorsion contraherent
cosheaves on $U_\alpha$ corresponding to flat cotorsion\/
$\O(U_\alpha)$\+modules.
\end{lem}

\begin{proof}
 Similar to the proof of Lemma~\ref{loc-contra-proj}
and based on Theorem~\ref{flat-cover-thm}(b).
\end{proof}

\begin{cor}  \label{ctrh-lcth-lct-proj}
\textup{(a)} There are enough projective objects in the exact
category $X\lcth^\lct$ of locally cotorsion locally contraherent
cosheaves on $X$, and all these projective objects belong to
the full subcategory of locally cotorsion contraherent cosheaves
$X\ctrh^\lct\sub X\lcth^\lct$. \par
\textup{(b)} The full subcategories of projective objects in
the three exact categories $X\ctrh^\lct\sub X\lcth_\bW^\lct\sub
X\lcth^\lct$ coincide. \qed
\end{cor}

 We denote the additive category of projective (objects in the category
of) contraherent cosheaves on $X$ by $X\ctrh_\prj$, and the additive
category of projective (objects in the category of) locally
cotorsion contraherent cosheaves $X$ by $X\ctrh^\lct_\prj$.

 Let us issue a \emph{warning} that both the terminology and notation are
misleading here: a projective locally cotorsion contraherent cosheaf on
$X$ does \emph{not} have to be a projective contraherent cosheaf.
 Indeed, a flat cotorsion module over a commutative ring would not
be in general very flat.

 On the other hand we notice that, by the definition, both additive
categories $X\ctrh_\prj$ and $X\ctrh^\lct_\prj$ are contained in
the exact category of antilocally flat contraherent cosheaves
$X\ctrh_\alf$ (and consequently also in the exact category of
antilocal contraherent cosheaves $X\ctrh_\al$).
 Moreover, by the definition one clearly has $X\ctrh^\lct_\prj =
X\ctrh^\lct\cap X\ctrh_\alf$.
 Finally, we notice that $X\ctrh_\prj$ is the category of
\emph{projective} objects in $X\ctrh_\alf$, while $X\ctrh^\lct_\prj$
is the category of \emph{injective} objects in $X\ctrh_\alf$ (and
there are enough of both).

\medskip

 A version of part~(a) of the following corollary that is valid
in a different generality will be obtained in
Section~\ref{flat-contra-subsection}.

\begin{cor}  \label{proj-flat}
 Let $X$ be a semi-separated coherent scheme.  Then \par
\textup{(a)} any cosheaf from $X\ctrh_\prj$ is flat; \par
\textup{(b)} any cosheaf from $X\ctrh^\lct_\prj$ is flat.
\end{cor}

\begin{proof}
 Follows from Corollary~\ref{clf-noetherian-flat}.
\end{proof}

 Alternative versions of the next corollary and of part~(b) of
the previous one will be obtained in
Section~\ref{lct-projective-loc-noetherian}.
 In both cases, we will see that the semi-separatedness and
quasi-compactness assumptions can be dropped if one assumes
local Noetherianity.

\begin{cor}  \label{proj-products}
 Over a semi-separated coherent scheme $X$, the full subcategory
$X\ctrh^\lct_\prj$ of projective locally cotorsion contraherent
cosheaves is closed under infinite products in $X\ctrh$.
\end{cor}

\begin{proof}
 Follows from Corollary~\ref{clf-products}.
\end{proof}

\begin{ex} \label{X-ctrh-prj-not-closed-under-products-ex}
 Notice that the full subcategory $X\ctrh_\prj$ of projective
contraherent cosheaves in \emph{not} closed under infinite products
in $X\ctrh$ even for Noetherian affine schemes~$X$.
 In fact, it suffices to consider the case of the scheme $X=\Spec\boZ$.
 Let $0\rarrow\boZ\rarrow C\rarrow V\rarrow0$ be a short exact
sequence of abelian groups with a very flat contraadjusted abelian
group $C$ and a very flat abelian group $V$, as per
Corollary~\ref{very-rel-proj-inj}(a).
 Then the group $V$ is also contraadjusted as a quotient group of
a contraadjusted group~$C$.
 Consider the short exact sequence $0\rarrow\prod_{n=0}^\infty\boZ
\rarrow\prod_{n=0}^\infty C\rarrow\prod_{n=0}^\infty V\rarrow0$.
 If both the groups $\prod_{n=0}^\infty C$ and $\prod_{n=0}^\infty V$
were very flat, it would follow that the group
$\prod_{n=0}^\infty\boZ$ is very flat (as the kernel of a surjective
morphism of very flat abelian groups), contradicting
Example~\ref{baer-specker-not-very-flat-ex}.
 Thus the class of all very flat contraadjusted abelian groups
$\boZ\modl_\vfl^\cta=X\ctrh_\prj$ is not closed under infinite products
in $\boZ\modl^\cta=X\ctrh$.
\end{ex}

 A more general version of parts~(a\+b) of the following Corollary
will be proved in the next Section~\ref{homology-subsection}, while
alternative versions of parts~(b\+c) will be also obtained in
Section~\ref{lct-projective-loc-noetherian}.
 In both cases, we will see that the affineness assumption on
the morphism~$f$ is unnecessary (but local Noetherianity is assumed
instead in Section~\ref{lct-projective-loc-noetherian}).

\begin{cor} \label{proj-direct-inverse}
 Let $f\:Y\rarrow X$ be an affine morphism of quasi-compact
semi-separated schemes.  Then \par
\textup{(a)} if the morphism~$f$ is very flat, then the direct image
functor~$f_!$ takes projective contraherent cosheaves to projective
contraherent cosheaves; \par
\textup{(b)} if the morphism~$f$ is flat, then the direct image
functor~$f_!$ takes projective locally cotorsion contraherent cosheaves
to projective locally cotorsion contraherent cosheaves; \par
\textup{(c)} if the scheme $X$ is coherent and the morphism~$f$ is
an open embedding, then the inverse image functor~$f^!$ takes projective
locally cotorsion contraherent cosheaves to projective locally cotorsion
contraherent cosheaves.
\end{cor}

\begin{proof}
 Part~(a) holds, since in its assumptions the functor $f_!\:Y\ctrh
\rarrow X\ctrh$ is ``partially left adjoint'' to the exact functor
$f^!\:X\lcth\rarrow Y\lcth$.
 Here the partial adjunction means a natural isomorphism
$\Hom^X(f_!\Q,\P)\simeq\Hom^Y(\Q,f^!\P)$
\,\eqref{direct-inverse-lct-adjunction}
for all $\Q\in Y\ctrh$ and $\P\in X\lcth$.
 The proof of part~(b) is similar (alternatively, it can be deduced
from the facts that the functor~$f_!$ takes $Y\ctrh^\lct$ to
$X\ctrh^\lct$ and $Y\ctrh_\alf$ to $X\ctrh_\alf$).
 Part~(c) follows from Corollary~\ref{clf-restriction}.
\end{proof}

 For a discussion of homotopy projective complexes of contraherent
cosheaves, homotopy projective complexes of locally cotorsion
contraherent cosheaves, and antilocaly flat complexes of contraherent
cosheaves, see Section~\ref{qc-ss-homotopy-projective-subsect} below.

\subsection{Homology of locally contraherent cosheaves}
\label{homology-subsection}
 The functor $\Delta(X,{-})$ of global cosections of locally
contraherent cosheaves on a scheme $X$, which assigns to a cosheaf $\gE$
the abelian group (or even the $\O(X)$\+module) $\gE[X]$, is right exact
as a functor on the exact category of locally contraherent cosheaves
$X\lcth$ on~$X$.
 In other words, if $0\rarrow\gK\rarrow\gL\rarrow\gM\rarrow0$ is a short
exact sequence of locally contraherent cosheaves on $X$, then
the sequence of abelian groups
$$
 \Delta(X,\gK)\lrarrow\Delta(X,\gL)\lrarrow\Delta(X,\gM)\lrarrow0
$$
is exact.
 Indeed, the procedure recovering the groups of cosections of cosheaves
$\gE$ on $X$ from their groups of cosections over affine open
subschemes $U\sub X$ subordinate to a particular covering~$\bW$ and
the corestriction maps between such groups uses the operations of
the infinite direct sum and the cokernel of a morphism (or in other
words, the nonfiltered inductive limit) only
(see~\eqref{cosheaf-definition}, \eqref{cosheaf-recover},
or~\eqref{cech-homol}).

 Recall that for any $(\bW,\bT)$\+affine morphism of schemes
$f\:Y\rarrow X$ the functor of direct image~$f_!$ takes $\bT$\+locally
contraherent cosheaves on $Y$ to $\bW$\+locally contraherent cosheaves
on~$X$.
 By the definition, there is a natural isomorphism of $\O(X)$\+modules
$\gE[Y]\simeq(f_!\gE)[X]$ for any cosheaf of $\O_Y$\+modules~$\gE$.

 Now let $X$ be a quasi-compact semi-separated scheme.
 Then the left derived functor of the functor of global cosections
of locally contraherent cosheaves on $X$ can be defined in
the conventional way using projective resolutions in the exact
category $X\lcth$ (see Lemma~\ref{loc-contra-proj} and
Corollary~\ref{ctrh-lcth-proj}).
 Notice that the derived functors of $\Delta(X,{-})$ (and in fact,
any left derived functors) computed in the exact category $X\lcth_\bW$
for a particular open covering~$\bW$ and in the whole category
$X\lcth$ agree.
 We denote this derived functor by $\boL_*\Delta(X,{-})$.
 The groups $\boL_i\Delta(X,\gE)$ are called the \emph{homology groups}
of a locally contraherent cosheaf $\gE$ on the scheme~$X$.

 Let us point out that the functor $\Delta(U,{-})$ of global cosections
of contraherent cosheaves on an affine scheme $U$ is exact, so
the groups $\boL_{>0}\Delta(U,\gE)$ vanish when $U$ is affine and
$\gE$ is contraherent.

 By Corollary~\ref{proj-direct-inverse}(a), for any very flat
$(\bW,\bT)$\+affine morphism of quasi-compact semi-separated schemes
$f\:Y\rarrow X$ the exact functor $f_!\:Y\lcth_\bT\rarrow X\lcth_\bW$
takes projective contraherent cosheaves on $Y$ to projective
contraherent cosheaves on~$X$.
 It also makes a commutative diagram with the restrictions of
the functors $\Delta(X,{-})$ and $\Delta(Y,{-})$ to the categories
$X\lcth_\bW$ and $Y\lcth_\bT$.
 Hence one has $\boL_*\Delta(Y,\gE)\simeq\boL_*\Delta(X,f_!\gE)$
for any $\bT$\+locally contraherent cosheaf $\gE$ on~$Y$.

 In particular, the latter assertion applies to the embeddings
of affine open subschemes $j\:U\rarrow X$, so $\boL_{>0}
\Delta(X,j_!\gE)=0$ for all contraherent cosheaves $\gE$ on~$U$.
 Since the derived functor $\boL_*\Delta$ takes short exact
sequences of locally contraherent cosheaves to long exact sequences
of abelian groups, it follows from Corollary~\ref{clp-cor}(c)
that $\boL_{>0}\Delta(X,\P)=0$ for any antilocal
contraherent cosheaf $\P$ on~$X$.

 Therefore, the derived functor $\boL_*\Delta$ can be computed
using antilocal resolutions.
 Now we also see that the derived functors $\boL_*\Delta$ defined
in the theories of arbitrary (i.~e., locally contraadjusted)
contraherent cosheaves and of locally cotorsion contraherent
cosheaves agree.

 Let $\gE$ be a $\bW$\+locally contraherent cosheaf on~$X$, and let
$X=\bigcup_\alpha U_\alpha$ be a finite affine open covering of $X$
subordinate to~$\bW$.
 Then the contraherent \v Cech resolution~\eqref{contraherent-cech}
for $\gE$ is an antilocal resolution of a locally contraherent cosheaf
$\gE$, and one can use it to compute
the derived functor $\boL_*(X,\gE)$.
 In other words, the homology of a $\bW$\+locally contraherent
cosheaf $\gE$ on a quasi-compact semi-separated scheme $X$ are
computed by the homological \v Cech complex $C_*(\{U_\alpha\},\gE)$
(see~\eqref{cech-homol}) related to any finite affine open
covering $X=\bigcup_\alpha U_\alpha$ subordinate to~$\bW$.

 The following result is to be compared with the cohomological
criterion of affineness of schemes~\cite[Th\'eor\`eme~5.2.1]{Groth1}.

\begin{cor} \label{contraherence-homological-criterion}
 A locally contraherent cosheaf\/ $\gE$ on an affine scheme $U$ is
contraherent if and only if its higher homology\/
$\boL_{>0}\Delta(X,\gE)$ vanish.
\end{cor}

\begin{proof}
 See Lemma~\ref{global-contraherence-criterion}.
\end{proof}

\begin{cor} \label{contraherence-nilpotent-reduction}
 Let $R$ be a commutative ring, $I\sub R$ be a nilpotent ideal,
and $S=R/I$ be the quotient ring.
 Let $i\:\Spec S\rarrow\Spec R$ denote the corresponding
homeomorphic closed embedding of affine schemes.
 Then a locally injective locally contraherent cosheaf\/ $\gJ$ on
$U=\Spec R$ is contraherent if and only if its inverse image
$i^!\gJ$ on $V=\Spec S$ is contraherent.
\end{cor}

\begin{proof}
 Since any morphism into an affine scheme is coaffine,
the  ``only if'' assertion is obvious.
 To prove the ``if'', we apply again
Lemma~\ref{global-contraherence-criterion}.
 Assuming that $\gJ$ is $\bW$\+locally contraherent on $U$
and $U=\bigcup_\alpha U_\alpha$ is a finite affine open covering
of $U$ subordinate to $\bW$, the homological \v Cech complex
$C_*(\{U_\alpha\},\gJ)$ is a finite complex of injective $R$\+modules.
 Set $V_\alpha=i^{-1}(U_\alpha)\sub V$; then the complex
$C_*(\{V_\alpha\},i^!\gJ)$ is the maximal subcomplex of $R$\+modules
in $C_*(\{U_\alpha\},\gJ)$ annihilated by the action of~$I$.

 Since the maximal submodule ${}_IM\sub M$ annihilated by $I$
is nonzero for any nonzero $R$\+module $M$, one easily proves
by induction that a finite complex $K^\bu$ of injective $R$\+modules
is acyclic at all its terms except perhaps the rightmost one
whenever so is the complex of injective $R/I$\+modules ${}_IK^\bu$.
 Notice that this argument does not seem to apply to the maximal
reduced closed subscheme $\Spec R/J$ of an arbitrary affine scheme
$\Spec R$ in general, as the maximal submodule ${}_JK\sub K$
annihilated by the nilradical $J\sub R$ may well be zero even
for a nonzero injective $R$\+module $K$ (set $K=\Hom_k(F,k)$ in
the example from Remark~\ref{nilradical-reduction-remark}).
\end{proof}

 The following result is to be compared with
Corollaries~\ref{coflasque-direct} and~\ref{cta-cot-direct} (for
another comparison, see Corollaries~\ref{clp-inverse}
and~\ref{clf-restriction}).
 See Remark~\ref{noncontraherent-direct-image-remark}
and Example~\ref{noncontraherent-direct-image-example}
for a discussion of the general problem of noncontraherence
of the direct images.

\begin{cor} \label{clp-direct}
 Let $f\:Y\rarrow X$ be a morphism of quasi-compact semi-separated
schemes.  Then \par
\textup{(a)} the functor of direct image of cosheaves of\/
$\O$\+modules $f_!$ takes antilocal contraherent cosheaves on $Y$ to
antilocal contraherent cosheaves on $X$, and induces an exact functor
$f_!\:Y\ctrh_\al\rarrow X\ctrh_\al$ between these exact categories; \par
\textup{(b)} the functor of direct image of cosheaves of\/
$\O$\+modules $f_!$ takes antilocal locally cotorsion contraherent
cosheaves on $Y$ to antilocal locally cotorsion contraherent cosheaves
on $X$, and induces an exact functor $f_!\:Y\ctrh_\al^\lct\rarrow
X\ctrh_\al^\lct$ between these exact categories; \par
\textup{(c)} if the morphism~$f$ is flat, then the functor of direct
image of cosheaves of\/ $\O$\+modules $f_!$ takes antilocally flat
contraherent cosheaves on $Y$ to antilocally flat contraherent cosheaves
on $X$, and induces an exact functor $f_!\:Y\ctrh_\alf\rarrow
X\ctrh_\alf$ between these exact categories; \par
\textup{(d)} if the morphism~$f$ is flat, then the functor of
direct image of cosheaves of\/ $\O$\+modules $f_!$ takes antilocal
locally injective contraherent cosheaves on $Y$ to antilocal
locally injective contraherent cosheaves on $X$.
\end{cor}

\begin{proof}
 Part~(a): by Corollary~\ref{clp-inverse}(a), the inverse image of
an antilocal contraherent cosheaf on $Y$ with respect to an affine
open embedding $j\:V\rarrow Y$ is antilocal.
 As we have seen above, the global cosections of antilocal
contraherent cosheaves is an exact functor.
 It follows that the functor~$f_!$ takes short exact sequences in
$Y\ctrh_\al$ to short exact sequences in the exact category of
cosheaves of $\O_X$\+modules (with the exact category structure
$\O_X\cosh_{\{X\}}$ related to the covering $\{X\}$ of the scheme~$X$;
see Section~\ref{locally-contraherent}).

 Since $X\ctrh_\al$ is a full exact subcategory closed under
extensions in $X\ctrh$, and the latter exact category is such
an exact subcategory in $\O_X\cosh_{\{X\}}$, in view of
Corollary~\ref{clp-cor}(c) it remains to recall that the direct images
of antilocal contraherent cosheaves with respect to affine morphisms
of schemes are antilocal (see the remarks in the beginning of
Section~\ref{clp-subsection}).

 Part~(b) is similar.
 The proof of part~(c) is also similar and based on part~(a)
together with the remarks about direct images in the beginning of
Section~\ref{clf-subsection} and Corollary~\ref{clf-cor}(c).
 To prove part~(d), one only needs to recall that the direct images
of locally injective contraherent cosheaves with respect to flat
affine morphisms of schemes are locally injective and
use Corollary~\ref{clp-lin}.
\end{proof}

 Let $f\:Y\rarrow X$ be a morphism of quasi-compact semi-separated
schemes.
 By the result of Section~\ref{direct-inverse-loc-contra}
(see~\eqref{direct-inverse-cosheaf-lin-adjunction}),
the adjunction isomorphism~\eqref{direct-inverse-lin-adjunction}
holds, in particular, for any antilocal contraherent cosheaf $\Q$ on $Y$
and any locally injective locally contraherent cosheaf $\gJ$ on~$X$.

 If the morphism~$f$ is flat then, according
to~\eqref{direct-inverse-cosheaf-lct-adjunction}, the adjunction
isomorphism
\begin{equation} \label{direct-inverse-clp-adjunction}
 \Hom^X(f_!\Q,\gM)\simeq\Hom^Y(\Q,f^!\gM)
\end{equation}
holds for any antilocal contraherent cosheaf $\Q$ on $Y$ and any
locally cotorsion locally contraherent cosheaf $\gM$ on~$X$.
 When the morphism~$f$ is also affine, the restrictions of
$f_!$ and~$f^!$ form an adjoint pair of functors between the exact
categories $Y\ctrh_\al^\lct$ and $X\ctrh_\al^\lct$.
 In addition, these functors take the additive categories
$Y\ctrh_\al^\lin$ and $X\ctrh_\al^\lin$ into one another.

 If the morphism~$f$ is very flat, the same adjunction
isomorphism~\eqref{direct-inverse-clp-adjunction} holds for any
antilocal contraherent cosheaf $\Q$ on $Y$ and any
locally contraherent cosheaf $\gM$ on~$X$.
 When the morphism~$f$ is also affine, the restrictions of $f_!$
and~$f^!$ form an adjoint pair of functors between the exact categories
$Y\ctrh_\al$ and $X\ctrh_\al$.

\begin{cor} \label{proj-direct-gen}
 Let $f\:Y\rarrow X$ be a morphism of quasi-compact semi-separated
schemes.  Then \par
\textup{(a)} if the morphism~$f$ is very flat, then the direct image
functor $f_!\:Y\ctrh_\al\rarrow X\ctrh_\al$ takes projective
contraherent cosheaves to projective contraherent cosheaves; \par
\textup{(b)} if the morphism~$f$ is flat, then the direct image
functor $f_!\:Y\ctrh_\al^\lct\rarrow X\ctrh_\al^\lct$ takes projective
locally cotorsion contraherent cosheaves to projective locally
cotorsion contraherent cosheaves.
\end{cor}

\begin{proof}
 Follows from the above partial
adjunctions~\eqref{direct-inverse-clp-adjunction} between
the exact functors $f_!$ and~$f^!$.
 Part~(b) can be also deduced from Corollary~\ref{clp-direct}(b\+c).
\end{proof}

\subsection{Derived categories of sheaves and cosheaves}
\label{derived-of-sheaves-and-cosheaves-subsect}
 We refer to Sections~\ref{derived-second-kind} and~\ref{becker-subsect}
for the definitions of the exotic derived categories mentioned below.
 Our treatment of the coderived and contraderived categories in the rest
of Chapter~\ref{derived-on-quasi-compact-sect} is rather
cursory; the detailed discussion is postponed to subsequent chapters.

 Given an exact category $\sE$, let $\Hot^\st(\sE)$ denote the homotopy
category $\Hot(\sE)$ if $\bst=\abs$, $\co$, $\ctr$, $\bco$, $\bctr$,
or~$\empt$; the category $\Hot^+(\sE)$ if $\bst=\abs+$ or~$+$;
the category $\Hot^-(\sE)$ if $\bst=\abs-$ or~$-$; and the category
$\Hot^\b(\sE)$ if $\bst=\b$.

 We also refer to the definition of the resolution dimension
$\rsd_{\sF/\sE} E$ of an object $E$ of an exact category $\sE$
with respect to a resolving subcategory $\sF\sub\sE$ given in
Section~\ref{finite-resolutions-subsect}.
 The coresolution dimension with respect to a coresolving
subcategory\/~$\sC$ (or the \emph{$\sC$\+coresolution dimension})
$\crd_{\sC/\sE} E$ is defined in the dual way.

 Let $X$ be a quasi-compact semi-separated scheme and $\bW$ be its
open covering.

\begin{lem}  \label{dil-cta-clp-finite-dim}
\textup{(a)} If $X=\bigcup_{\alpha=1}^N U_\alpha$ is a finite affine
open covering, then the coresolution dimension of any quasi-coherent
sheaf on $X$ with respect to the coresolving subcategory of dilute
quasi-coherent sheaves $X\qcoh^\dil\sub X\qcoh$ does not exceed~$N-1$.
\par
\textup{(b)} If $X=\bigcup_{\alpha=1}^N U_\alpha$ is a finite affine
open covering, then the coresolution dimension of any quasi-coherent 
sheaf on $X$ with respect to the coresolving subcategory of
contraadjusted quasi-coherent sheaves $X\qcoh^\cta\sub X\qcoh$ does
not exceed~$N$. \par
\textup{(c)} If $X=\bigcup_{\alpha=1}^N U_\alpha$ is a finite affine
open covering subordinate to\/ $\bW$, then the resolution dimension of
any\/ $\bW$\+locally contraherent cosheaf on $X$ with respect to
the resolving subcategory of antilocal contraherent cosheaves
$X\ctrh_\al\sub X\lcth_\bW$ does not exceed~$N-1$.
 Consequently, the same bound holds for the resolution dimension
of any object of $X\lcth_\bW$ with respect to the resolving
subcategory $X\ctrh$.
\end{lem}

\begin{proof}
 Part~(a): first of all, the full subcategory $X\qcoh^\dil$ is
coresolving in $X\qcoh$ according to the discussion in
Section~\ref{dilute-subsect}.
 Now any quasi-coherent sheaf $\F$ on $X$ has a \v Cech
coresolution~\eqref{cech-quasi} of length $N-1$ by quasi-coherent
sheaves that are dilute in view of
Lemma~\ref{flat-affine-direct-image-dilute}
or Corollary~\ref{finitely-iterated-extension-dilute}.

 Part~(b): the full subcategory $X\qcoh^\cta$ is
coresolving in $X\qcoh$ by
Corollaries~\ref{quasi-cta-characterizations}(c)
and~\ref{quasi-very-cta-cor}(b), and the discussion in the beginning
of Section~\ref{quasi-compact-quasi-coherent}.
 Furthermore, the coresolution dimension of any module over
a commutative ring $R$ with respect to the coresolving subcategory of
contraadjusted $R$\+modules $R\modl^\cta\sub R\modl$ does not exceed~$1$.
 It follows easily that the coresolution dimension of any quasi-coherent
sheaf of the form $j_*\G$, where $j\:U\rarrow X$ is an affine open
subscheme and $\G\in U\qcoh$, with respect to the coresolving
subcategory $X\qcoh^\cta\sub X\qcoh$ does not exceed~$1$, either.
 Once again, any quasi-coherent sheaf $\F$ on $X$ has a \v Cech
coresolution~\eqref{cech-quasi} of length~$N-1$ by finite direct sums
of quasi-coherent sheaves of the form~$j_*\G$.
 It remains to use the dual version of
Corollary~\ref{fdim-resolution}(a).

 Part~(c): the full subcategory $X\ctrh_\al$ is resolving in
$X\lcth_\bW$ by Corollaries~\ref{clp-characterizations}(b)
and~\ref{clp-cor}(b).
 The full subcategory $X\ctrh$ is resolving in $X\lcth_\bW$ for
the reasons explained in Section~\ref{counterex-subsect}.
 It remains to recall the resolution~\eqref{contraherent-cech}.
\end{proof}

\begin{lem}  \label{lct-lin-clp-finite-dim}
 Let $X = \bigcup_{\alpha=1}^N U_\alpha$ be a finite affine open covering
subordinate to\/~$\bW$.  Then \par
\textup{(a)} the resolution dimension of any locally cotorsion\/
$\bW$\+locally contraherent cosheaf on $X$ with respect to the resolving
subcategory of antilocal locally cotorsion contraherent cosheaves
$X\ctrh_\al^\lct\sub X\lcth_\bW^\lct$ does not exceed $N-1$.
 Consequently, the same bound holds for the resolution dimension
of any object of $X\lcth_\bW^\lct$ with respect to the resolving
subcategory $X\ctrh^\lct$; \par
\textup{(b)} the resolution dimension of any locally injective\/
$\bW$\+locally contraherent cosheaf on $X$ with respect to the resolving
subcategory of antilocal locally injective contraherent cosheaves
$X\ctrh_\al^\lin\sub X\lcth_\bW^\lin$ does not exceed $N-1$.
 Consequently, the same bound holds for the resolution dimension
of any object of $X\lcth_\bW^\lin$ with respect to the resolving
subcategory $X\ctrh^\lin$; \par
\textup{(c)} the homological dimension of the exact category\/
$X\lcth_\bW^\lin$ does not exceed $N-1$.
\end{lem}

\begin{proof}
 The proofs of parts~(a\+b) are similar to that of
Lemma~\ref{dil-cta-clp-finite-dim}(c).
 Part~(c) follows from~(b), as $X\ctrh_\al^\lin$ is the full
subcategory of projective objects in the exact category
$X\lcth_\bW^\lin$, and there are enough such projective objects
(see Section~\ref{clp-subsection}).

 Alternatively, for part~(c) one can prove the seemingly stronger
(but actually equivalent) assertion that $\Ext^{X,>N-1}(\P,\gJ)=0$
for all $\bW$\+locally contraherent cosheaves $\P$ and locally
injective $\bW$\+locally contraherent cosheaves $\gJ$ on~$X$
(cf.\ the dual version of
Lemma~\ref{finite-homol-dim-finite-coresol-dim}).
 This assertion follows from Lemma~\ref{dil-cta-clp-finite-dim}(c).
 The same equivalent version of part~(c) is provided by a \v Cech
resolution argument dual-analogous to~\cite[Theorem~6.3(b)]{PS6}
(cf.\ Lemma~\ref{vfl-cta-finite-dim}(c) below).
\end{proof}

\begin{cor} \label{qcoh-dil-cta-derived-equiv-cor}
\textup{(a)} For any symbol\/ $\bst=\b$, $+$, $-$, $\empt$, $\abs+$,
$\abs-$, $\bco$, $\co$, or\/~$\abs$, the triangulated functor\/
$\sD^\st(X\qcoh^\dil)\rarrow\sD^\st(X\qcoh)$ induced by the embedding
of exact/abelian categories $X\qcoh^\dil\rarrow X\qcoh$ is
an equivalence of triangulated categories\/
$\sD^\st(X\qcoh^\dil)\simeq\sD^\st(X\qcoh)$. \par
\textup{(b)} For any symbol\/ $\bst=\b$, $+$, $-$, $\empt$, $\abs+$,
$\abs-$, $\bco$, or\/~$\abs$, the triangulated functor\/
$\sD^\st(X\qcoh^\cta)\rarrow\sD^\st(X\qcoh)$ induced by the embedding
of exact/abelian categories $X\qcoh^\cta\rarrow X\qcoh$ is
an equivalence of triangulated categories\/
$\sD^\st(X\qcoh^\cta)\simeq\sD^\st(X\qcoh)$.
\end{cor}

\begin{proof}
 All the assertions follow from the dual versions of
Propositions~\ref{finite-resolutions}
and~\ref{becker-contraderived-finite-resolutions} together
with Lemma~\ref{dil-cta-clp-finite-dim}(a\+b).
 Notice that the full subcategory $X\qcoh^\dil$ is closed under
infinite direct sums in $X\qcoh$ by
Corollary~\ref{direct-sums-dilute} (so the coderived category
symbol $\bst=\co$ is allowed in part~(a), but not in part~(b)).
\end{proof}

 For a version of the assertion of
Corollary~\ref{qcoh-dil-cta-derived-equiv-cor} for $\bst=\bco$
and the category $X\qcoh^\cot$ instead of $X\qcoh^\cta$,
see Corollary~\ref{dil-cta-cot-into-all-equiv-on-bco}(c) below.

\begin{cor}  \label{ctrh-lcth-cor}
\textup{(a)} For any symbol\/ $\bst=\b$, $+$, $-$, $\empt$, $\abs+$,
$\abs-$, $\bctr$, $\ctr$, or\/~$\abs$, the triangulated functors\/
$\sD^\st(X\ctrh_\al)\rarrow\sD^\st(X\ctrh)\rarrow\sD^\st(X\lcth_\bW)$
induced by the embeddings of exact categories $X\ctrh_\al\rarrow
X\ctrh\rarrow X\lcth_\bW$ are equivalences of triangulated categories
$$
 \sD^\st(X\ctrh_\al)\simeq\sD^\st(X\ctrh)\simeq\sD^\st(X\lcth_\bW).
$$ \par
\textup{(b)} For any symbol\/ $\bst=\b$ or\/~$-$, the triangulated
functor\/ $\sD^\st(X\lcth_\bW)\rarrow\sD^\st(X\lcth)$ induced by
the embedding of exact categories $X\lcth_\bW\rarrow X\lcth$ is
an equivalence of triangulated categories.
\end{cor}

 The reason why most unbounded derived categories aren't mentioned in
part~(b) is because one needs a uniform restriction on the extension
of locality of locally contraherent cosheaves in order to work
simultaneously with infinite collections of these.
 In particular, infinite products exist in $X\lcth_\bW$, but not
necessarily in $X\lcth$, so the contraderived category of the latter
exact category is not well-defined.

\begin{proof}[Proof of Corollary~\ref{ctrh-lcth-cor}]
 Part~(a) follows from Propositions~\ref{finite-resolutions}
and~\ref{becker-contraderived-finite-resolutions} together
with Lemma~\ref{dil-cta-clp-finite-dim}(c).
 Part~(b) in the case $\bst=\b$ is obtained from part~(a) by passing
to the inductive limit over refinements of coverings, while in
the case $\bst=-$ it is provided by
Proposition~\ref{infinite-resolutions}(a).
\end{proof}

 The following two corollaries are similar to the previous one.
 The only difference in the proofs is that
Lemma~\ref{lct-lin-clp-finite-dim}(a\+b) is being used in place
of Lemma~\ref{dil-cta-clp-finite-dim}(c).

\begin{cor}  \label{lct-ctrh-lcth-cor}
\textup{(a)} For any symbol\/ $\bst=\b$, $+$, $-$, $\empt$, $\abs+$,
$\abs-$, $\bctr$, $\ctr$, or\/~$\abs$, the triangulated functors\/
$\sD^\st(X\ctrh_\al^\lct)\rarrow\sD^\st(X\ctrh^\lct)\rarrow
\sD^\st(X\lcth_\bW^\lct)$ induced by the embeddings of exact categories
$X\ctrh_\al^\lct\rarrow X\ctrh^\lct\rarrow X\lcth_\bW^\lct$ are 
equivalences of triangulated categories
$$
 \sD^\st(X\ctrh_\al^\lct)\simeq\sD^\st(X\ctrh^\lct)\simeq
 \sD^\st(X\lcth_\bW^\lct).
$$ \par
\textup{(b)} For any symbol\/ $\bst=\b$ or\/~$-$, the triangulated
functor\/ $\sD^\st(X\lcth_\bW^\lct)\rarrow\sD^\st(X\lcth^\lct)$ induced
by the embedding of exact categories $X\lcth_\bW^\lct\rarrow X\lcth^\lct$
is an equivalence of triangulated categories. \qed
\end{cor}

\begin{cor}  \label{lin-ctrh-lcth-cor}
\textup{(a)} For any symbol\/ $\bst=\b$, $+$, $-$, $\empt$, $\abs+$,
$\abs-$, $\bctr$, $\ctr$, or\/~$\abs$, the triangulated functors\/
$\Hot^\st(X\ctrh_\al^\lin)\rarrow\sD^\st(X\ctrh^\lin)\rarrow
\sD^\st(X\lcth_\bW^\lin)$ induced by the embeddings of additive/exact
categories $X\ctrh_\al^\lin\rarrow X\ctrh^\lin\rarrow X\lcth_\bW^\lin$
are equivalences of triangulated categories
$$
 \Hot^\st(X\ctrh_\al^\lin)\simeq\sD^\st(X\ctrh^\lin)\simeq
 \sD^\st(X\lcth_\bW^\lin).
$$
 In particular, the classes of acyclic, contraacyclic,
Becker-contraacyclic, and absolutely acyclic complexes in
the exact category $X\lcth_\bW^\lin$ coincide,
$$
 \Acycl(X\lcth_\bW^\lin)=\Acycl^\bctr(X\lcth_\bW^\lin)=
 \Acycl^\ctr(X\lcth_\bW^\lin)=\Acycl^\abs(X\lcth_\bW^\lin).
$$ \par
\textup{(b)} For any symbol\/ $\bst=\b$ or\/~$-$, the triangulated
functor\/ $\sD^\st(X\lcth_\bW^\lin)\rarrow\sD^\st(X\lcth^\lin)$
induced by the embedding of exact categories
$X\lcth_\bW^\lin\rarrow X\lcth^\lin$ is an equivalence of
triangulated categories.
\end{cor}

\begin{proof}
 The second assertion of part~(a) follows from the first one.
 It can be also deduced from Lemma~\ref{lct-lin-clp-finite-dim}(c)
using Lemma~\ref{psemi-remark21} and/or
Theorem~\ref{finite-homol-dim-becker-co-contra-derived}(b).
 For a slight/partial generalization of the second assertion of
part~(a), see
Corollary~\ref{vfl-lin-finite-dim-all-derived-coincide}(b) below.
\end{proof}

 In the rest of this section, our aim is to prove
Theorem~\ref{derived-loc-cta-loc-cot-theorem}.
 The argument is based on the cotorsion periodicity theorem
(Theorem~\ref{cotorsion-periodicity}) and the construction of
gluing of cotorsion pairs in the categories of locally contraherent
cosheaves (Proposition~\ref{loc-contraherent-gluing-preenvelope}).
 In the next two lemmas we use the terminology and notation of
Appendix~\ref{cotorsion-pairs-appx}.

\begin{lem} \label{acyclic-complexes-of-cta-flats-cotorsion-pair}
 Let\/ $\sR$ be the local class of all commutative rings $R$ and\/
$\sE^R=\sK^R=\Com(R\modl^\cta)$ be the exact category of complexes
of contraadjusted $R$\+modules.
 Let\/ $\sF(R)\sub\sE^R$ be the class of all acyclic complexes of
flat contraadjusted $R$\+modules with flat (and contraadjusted)
modules of cocycles, and let\/ $\sC^R=\Com(R\modl^\cot)\sub\sE^R$
be the class of all complexes of cotorsion $R$\+modules.
 Then\/ $\sE$ and\/ $\sC$ are very colocal classes, and the pair of
classes $(\sF(R),\sC^R)$ is a hereditary complete cotorsion pair in\/
$\sE^R$ for every ring $R\in\sR$.
\end{lem}

\begin{proof}
 The classes $\sE$ and $\sC$ are very colocal because the classes
of contraadjusted $R$\+modules and cotorsion $R$\+modules are very
colocal; see Examples~\ref{colocal-classes-examples}.
 In any acyclic complex of contraadjusted $R$\+modules, the modules
of cocycles are contraadjusted since all quotient modules of
contraadjusted modules are contraadjusted.
 The full subcategory $\sC^R$ is closed under cokernels of admissible
monomorphisms in $\sE^R$ because the cokernels of injective maps of
cotorsion modules are cotorsion.
 It remains to explain why $(\sF(R),\sC^R)$ is a complete cotorsion
pair in~$\sE^R$.
 For this, we refer to~\cite[Example~7.12]{Pal}.
\end{proof}

\begin{lem} \label{complexes-loc-cot-preenvelope}
 Any complex of\/ $\bW$\+locally contraherent cosheaves\/ $\gM^\bu$
on $X$ can be included in a (termwise admissible) short exact sequence
of complexes\/ $0\rarrow\gM^\bu\rarrow\P^\bu\rarrow\gF^\bu\rarrow0$,
where\/ $\P^\bu$ is a complex of locally cotorsion\/ $\bW$\+locally
contraherent cosheaves on $X$ and\/ $\gF^\bu$ is an acyclic complex
in the exact category of antilocally flat contraherent cosheaves on~$X$.
\end{lem}

\begin{proof}
 Let $X=\bigcup_\alpha U_\alpha$ be a finite affine open covering
subordinate to~$\bW$.
 We apply Proposition~\ref{loc-contraherent-gluing-preenvelope} to
the datum of classes $\sR$, \,$\sE^R$, \,$\sC^R$, and $\sF(R)$ from
Lemma~\ref{acyclic-complexes-of-cta-flats-cotorsion-pair}.
 This produces, for a given complex of $\bW$\+locally contraherent
cosheaves $\gM^\bu$, a short exact sequence of complexes
$0\rarrow\gM^\bu\rarrow\P^\bu\rarrow\gF^\bu\rarrow0$, where $\P^\bu$
is a complex of locally cotorsion contraherent cosheaves and
$\gF^\bu$ is a finitely iterated extension in $\Com(X\lcth_\bW)$ of
direct images under the open embedding morphisms $j_\alpha\:U_\alpha
\rarrow X$ of acyclic complexes $\gF_\alpha^\bu$ in the exact
categories of $U_\alpha\ctrh_\alf$ of antilocally flat contraherent
cosheaves on~$U_\alpha$.
 It remains to recall that the direct image functor
$j_\alpha{}_!\:U_\alpha\ctrh\rarrow X\ctrh$ is exact and takes
$U_\alpha\ctrh_\alf$ into $X\ctrh_\alf$ (see the beginning of
Section~\ref{clf-subsection}, cf.\ Corollary~\ref{clp-direct}(c)).
\end{proof}

 In fact, we will see in Chapter~\ref{becker-on-qcomp-qsep-sect}
that any acyclic complex in the exact category $X\ctrh_\alf$ is
a direct summand of a finitely iterated extension of direct images
of acyclic complexes in the exact categories $U_\alpha\ctrh_\alf$.
 This is the result of
Corollary~\ref{acyclic-complexes-of-alf-antilocal};
cf.\ Corollary~\ref{acyclic-in-alf-all-of-lct-loc-contraherent-pair}.

\begin{thm} \label{derived-loc-cta-loc-cot-theorem}
 Let $X$ be a quasi-compact semi-separated scheme with an open
covering\/~$\bW$.
 Then, for any symbol\/ $\bst=+$ or\/~$\empt$, the triangulated
functor\/ $\sD^\st(X\lcth_\bW^\lct)\rarrow\sD^\st(X\lcth_\bW)$ induced
by the embedding of exact categories $X\lcth_\bW^\lct\rarrow X\lcth_\bW$
is an equivalence of triangulated categories.
\end{thm}

\begin{proof}
 In the easy case of $\bst=+$ it suffices to refer to the dual version
of Proposition~\ref{infinite-resolutions}(a); the point is that
the full subcategory $X\lcth_\bW^\lct$ is coresolving in
$X\lcth_\bW$, as explained in the discussion in
Section~\ref{clp-subsection}.
 The case of fully unbounded complexes, $\bst=\empt$, is interesting.
 Notice that the dual version of Proposition~\ref{finite-resolutions}
is \emph{not} applicable, as the objects of $X\lcth_\bW$ need not 
have finite coresolution dimensions with respect to $X\lcth_\bW^\lct$
(see~\cite[Section~7]{Pphil} for a discussion).

 Here is the argument.
 Lemma~\ref{complexes-loc-cot-preenvelope} tells us, in particular, that
for every complex of $\bW$\+locally contraherent cosheaves $\gM^\bu$
on $X$ there exists a complex of locally cotorsion $\bW$\+locally
contraherent cosheaves $\P^\bu$ on $X$ together with a quasi-isomorphism
$\gM^\bu\rarrow\P^\bu$ of complexes in the exact category $X\lcth_\bW$.
 (Indeed, any acyclic complex in $X\ctrh_\alf$ is also acyclic in
$X\lcth_\bW$.)
 In view of Lemma~\ref{pkoszul-lemma16}(b), it remains to show that any
complex in $X\lcth^\lct_\bW$ that is acyclic in $X\lcth_\bW$ is also
acyclic in $X\lcth^\lct_\bW$.
 This follows immediately from the fact that in any acyclic complex of
cotorsion $R$\+modules (say, with contraadjusted $R$\+modules of
cocycles), the $R$\+modules of cocycles are actually cotorsion;
see Theorem~\ref{cotorsion-periodicity}.

 Alternatively, one can use the result of
Lemma~\ref{complexes-loc-hot-inj-of-inj-preenvelope} below in
this proof instead of Lemma~\ref{complexes-loc-cot-preenvelope}.
 (Indeed, any locally injective $\bW$\+locally contraherent cosheaf
is locally cotorsion.)
\end{proof}

 Comparing Theorem~\ref{derived-loc-cta-loc-cot-theorem} with
Corollaries~\ref{ctrh-lcth-cor}(a) and~\ref{lct-ctrh-lcth-cor}(a),
we come to a commutative square of triangulated equivalences of
conventional derived categories of the exact categories of locally
contraherent cosheaves, with the derived category symbols
$\bst=+$ or~$\empt$:
\begin{equation}
\begin{gathered}
 \xymatrix{
  \sD^\st(X\ctrh^\lct) \ar@<2pt>[r] \ar@<-2pt>@{-}[r]
  \ar@<-2pt>[d] \ar@<2pt>@{-}[d]
  & \sD^\st(X\lcth^\lct_\bW) \ar@<2pt>[d] \ar@<-2pt>@{-}[d] \\
  \sD^\st(X\ctrh) \ar@<-2pt>[r] \ar@<2pt>@{-}[r]
  & \sD^\st(X\lcth_\bW)
 }
\end{gathered}
\end{equation}

\subsection{The ``na\"\i ve'' co-contra correspondence}
\label{naive-co-contra-subsect}
 Let $X$ be a quasi-compact semi-separated scheme with an open
covering~$\bW$.
 The following theorem is the main result of
Chapter~\ref{derived-on-quasi-compact-sect}.

\begin{thm}  \label{naive-co-contra-thm}
 For any symbol\/ $\bst=\b$, $+$, $-$, $\empt$, $\abs+$, $\abs-$,
or\/~$\abs$ there is a natural equivalence of triangulated categories\/
$\sD^\st(X\qcoh)\simeq \sD^\st(X\ctrh)$.
 These equivalences of derived categories form commutative
diagrams with the natural functors\/ $\sD^\b\rarrow\sD^\pm\rarrow\sD$,
\ $\sD^\b\rarrow\sD^{\abs\pm}\rarrow\sD^\abs$, \
$\sD^{\abs\pm}\rarrow\sD^\pm$, \ $\sD^\abs\rarrow\sD$
between different versions of derived categories
of the same exact category.
\end{thm}

 Notice that Theorem~\ref{naive-co-contra-thm} does not say anything
about the coderived and contraderived categories $\sD^\co$ and $\sD^\ctr$
of quasi-coherent sheaves and contraherent cosheaves (neither does
Corollary~\ref{ctrh-lcth-cor} mention the coderived categories).
 The reason is that infinite products are not exact in the abelian
category of quasi-coherent sheaves and infinite direct sums may not
exist in the exact category of contraherent cosheaves.
 So only the coderived category $\sD^\co(X\qcoh)$ and the contraderived
category $\sD^\ctr(X\ctrh)$ are well-defined.
 Comparing these two requires a different approach; the entire
Chapter~\ref{dualizing-complex-sect} will be devoted to that.

\begin{proof}[Proof of Theorem~\ref{naive-co-contra-thm}]
 By Corollaries~\ref{qcoh-dil-cta-derived-equiv-cor}(b)
and~\ref{ctrh-lcth-cor}(a), the triangulated functors
$\sD^\st(X\qcoh^\cta)\rarrow\sD^\st(X\qcoh)$ and $\sD^\st(X\ctrh_\al)
\rarrow\sD^\st(X\ctrh)$ induced by the respective embedings of
exact categories are all equivalences of triangulated categories.
 Hence it suffices to construct a natural equivalence of exact
categories $X\qcoh^\cta\simeq X\ctrh_\al$ in order to prove all
assertions of Theorem.

 According to Sections~\ref{fHom-subsection}
and~\ref{contratensor-subsect}, there are natural functors
$$
 \fHom_X(\O_X,{-})\: X\qcoh^\cta\lrarrow X\ctrh
$$
and
$$
 \O_X\ocn_X{-}\: X\lcth\lrarrow X\qcoh
$$
related by the adjunction
isomorphism~\eqref{fHom-contratensor-adjunction}, which holds
for those objects for which the former functor is defined.
 So it remains to prove the following lemma.
\end{proof}

\begin{lem} \label{cta-clp-equivalence} \hbadness=2600
 On a quasi-compact semi-separated scheme $X$, the functor\/
$\fHom_X(\O_X,{-})$ takes $X\qcoh^\cta$ to $X\ctrh_\al$,
the functor\/ $\O_X\ocn_X{-}$ takes $X\ctrh_\al$ to $X\qcoh^\cta$,
and the restrictions of these functors to these subcategories
are mutually inverse equivalences of exact categories.
\end{lem}

\begin{proof}
 Obviously, on an affine scheme $U$ the functor $\fHom_U(\O_U,{-})$
takes a contraadjusted quasi-coherent sheaf $\cQ$ with
the contraadjusted $\O(U)$\+module of global sections $\cQ(U)$
to the contraherent cosheaf $\Q$ with the contraadjusted
$\O(U)$\+module of global cosections $\Q[U]=\cQ(U)$.
 Furthermore, if $j\:U\rarrow X$ is the embedding of an affine open
subscheme, then by the formula~\eqref{flat-cta-fHom-projection}
of Section~\ref{compatibility-subsect} there is a natural isomorphism
$\fHom_X(\O_X,j_*\cQ)\simeq j_!\Q$ of contraherent cosheaves on~$X$.

 Analogously, the functor $\O_U\ocn_U{-}$ takes a contraherent
cosheaf $\Q$ with the contraadjusted $\O(U)$\+module of global
cosections $\Q[U]$ to the contraadjusted quasi-coherent sheaf
$\cQ$ with the $\O(U)$\+module of global sections $\cQ(U)$ on~$U$.
 If an embedding of affine open subscheme $j\:U\rarrow X$ is given,
then by the formula~\eqref{contratensor-projection} there is
a natural isomorphism $\O_X\ocn_X j_!\Q\simeq j_*\cQ$ of
quasi-coherent sheaves on~$X$.

 By Corollary~\ref{quasi-very-cta-cor}(c), any sheaf from
$X\qcoh^\cta$ is a direct summand of a finitely iterated extension
of the direct images of contraadjusted quasi-coherent sheaves
from affine open subschemes of~$X$.
 It is clear from the definition of the functor $\fHom_X(\O_X,{-})$
that it preserves exactness of short sequences of contraadjusted
quasi-coherent sheaves; hence it preserves, in particular, such
iterated extensions.

 By Corollary~\ref{clp-cor}(c), any cosheaf from $X\ctrh_\al$ is
a direct summand of a finitely iterated extension of the direct
images of contraherent cosheaves from affine open subschemes of~$X$.
 Let us show that the functor $\O_X\ocn_X{-}$ preserves exactness
of short sequences of antilocal contraherent cosheaves on $X$, and
therefore, in particular, preserves such extensions.
 Indeed, the adjunction isomorphism
$$
 \Hom_X(\O_X\ocn_X\P\;\F)\simeq
 \Hom^X(\P,\fHom_X(\O_X,\F))
$$
holds for any contraherent cosheaf $\P$ and contraadjusted
quasi-coherent sheaf~$\F$.
 Besides, the contraherent cosheaf $\fHom_X(\O_X,\J)$ is locally
injective for any injective quasi-coherent sheaf $\J$ on~$X$.
 By Corollary~\ref{clp-characterizations}(a), it follows that the functor
$\P\mpsto \Hom_X(\O_X\ocn_X\P\;\J)$ preserves exactness of short
sequences of antilocal contraherent cosheaves on~$X$, and consequently
so does the functor $\P\mpsto \O_X\ocn_X\P$.

 Now one can easily deduce that the adjunction morphisms
$$
 \P\rarrow\fHom_X(\O_X\;\O_X\ocn_X\P)
 \quad \text{and} \quad
\O_X\ocn_X\fHom_X(\O_X,\F)\rarrow\F
$$
are isomorphisms for all antilocal contraherent cosheaves $\P$ and
contraadjusted quasi-coherent sheaves $\F$, as a morphism of finitely
filtered objects inducing an isomorphism of the associated graded
objects is also itself an isomorphism.
 The proof of Lemma, and hence also of
Theorem~\ref{naive-co-contra-thm}, is finished.
\end{proof}

\begin{cor} \label{cta-clp-direct-image-cor}
 Let $f\:Y\rarrow X$ be a morphism of quasi-compact semi-separated
schemes.
 Then the equivalences of categories $X\qcoh^\cta\simeq X\ctrh_\al$
and $Y\qcoh^\cta\simeq Y\ctrh_\al$ from
Lemma~\textup{\ref{cta-clp-equivalence}} transform the functor
$f_*\:Y\qcoh^\cta\rarrow X\qcoh^\cta$ from
Corollary~\textup{\ref{cta-cot-direct}(a)} into the functor
$f_!\:Y\ctrh_\al\rarrow X\ctrh_\al$ from
Corollary~\textup{\ref{clp-direct}(a)}.
 In other words, the following diagram of exact functors and
exact category equivalences is commutative:
\begin{equation} \label{cta-clp-direct-image-diagram}
\begin{gathered}
 \xymatrix{
  Y\qcoh^\cta \ar@{=}[r] \ar[d]_{f_*} & Y\ctrh_\al \ar[d]^{f_!} \\
  X\qcoh^\cta \ar@{=}[r] & X\ctrh_\al
 }
\end{gathered}
\end{equation}
\end{cor}

\begin{proof}
 The assertion follows from the formula~\eqref{flat-cta-fHom-projection}
of Section~\ref{compatibility-subsect}.
 For an affine morphism~$f$,
the formula~\eqref{contratensor-projection} is also applicable,
leading to the same conclusion.
 More generally, for any morphism $f\:Y\rarrow X$ one can use
the formula~\eqref{flat-al-contratensor-projection-eqn} from
Lemma~\ref{flat-al-contratensor-projection-lemma} below.
\end{proof}

\begin{lem} \label{cta-clp-restricts-to-cot-inj}
\textup{(a)} The equivalence of exact categories
$X\qcoh^\cta\simeq X\ctrh_\al$ from
Lemma~\textup{\ref{cta-clp-equivalence}} restricts to an equivalence of
full exact subcategories $X\qcoh^\cot\simeq X\ctrh_\al^\lct$. \par
\textup{(b)} The equivalence of exact categories
$X\qcoh^\cta\simeq X\ctrh_\al$ from
Lemma~\textup{\ref{cta-clp-equivalence}} restricts to an equivalence of
full additive subcategories $X\qcoh^\inj\simeq X\ctrh_\al^\lin$.
\end{lem}

\begin{proof}
 Part~(a): the assertion that the equivalence of exact categories
$X\qcoh^\cta\simeq X\ctrh_\al$ identifies the full subcategory
$X\qcoh^\cot\sub X\qcoh^\cta$ with the full subcategory
$X\ctrh_\al^\lct\sub X\ctrh_\al$ follows from
Corollaries~\ref{quasi-cotors-cor}(c) and~\ref{clp-lct-cor}(c).

 Part~(b): the assertion that the equivalence of exact categories
$X\qcoh^\cta\simeq X\ctrh_\al$ identifies the full subcategory
$X\qcoh^\inj\sub X\qcoh^\cta$ with the full subcategory
$X\ctrh_\al^\lin\sub X\ctrh_\al$ is easily obtained from
Corollary~\ref{clp-lin} together with the fact that any injective
quasi-coherent sheaf on $X$ is a direct summand of a finite direct
sum of the direct images of injective quasi-coherent sheaves from
open embeddings $U_\alpha\rarrow X$ forming a covering.
\end{proof}

\begin{lem} \label{cta-clp-restricts-to-prj-clf}
\textup{(a)} The equivalence of exact categories
$X\qcoh^\cta\simeq X\ctrh_\al$ from
Lemma~\textup{\ref{cta-clp-equivalence}} restricts to an equivalence of
full additive subcategories $X\qcoh^\cta_\vfl\simeq X\ctrh_\prj$. \par
\textup{(b)} The equivalence of exact categories
$X\qcoh^\cta\simeq X\ctrh_\al$ from
Lemma~\textup{\ref{cta-clp-equivalence}} restricts to an equivalence of
full additive subcategories $X\qcoh^\cot_\fl\simeq X\ctrh^\lct_\prj$.
\par
\textup{(c)} The equivalence of exact categories
$X\qcoh^\cta\simeq X\ctrh_\al$ from
Lemma~\textup{\ref{cta-clp-equivalence}} restricts to an equivalence of
full exact subcategories $X\qcoh^\cta_\fl\simeq X\ctrh_\alf$.
\end{lem}

\begin{proof}
 This is similar to Lemma~\ref{cta-clp-restricts-to-cot-inj}.
 Part~(a) follows from Lemmas~\ref{very-flat-contraadjusted-quasi}
and~\ref{loc-contra-proj}(b).
 Part~(b) follows from Lemmas~\ref{flat-cotorsion-quasi}
and~\ref{loc-lct-proj}(b).
 Part~(c) follows from Lemma~\ref{flat-contraadjusted-quasi} and
Corollary~\ref{clf-cor}(c).
\end{proof}


\begin{cor}  \label{inj-co-contra-cor}
\textup{(a)} For any symbol\/ $\bst=\b$, $+$, $-$, $\empt$, $\abs+$,
$\abs-$, $\bctr$, $\ctr$, or\/~$\abs$, there is a natural equivalence
of triangulated categories\/ $\sD^\st(X\qcoh^\cot)\simeq
\sD^\st(X\ctrh^\lct)$. \par
\textup{(b)} For any symbol\/ $\bst=\b$, $+$, $-$, $\empt$, $\abs+$,
$\abs-$, $\bctr$, $\ctr$, or\/~$\abs$, there is a natural equivalence
of triangulated categories\/ $\Hot^\st(X\qcoh^\inj)\simeq
\sD^\st(X\ctrh^\lin)$.
\end{cor}

\begin{proof} \hbadness=2400
 By Corollaries~\ref{lct-ctrh-lcth-cor}(a)
and~\ref{lin-ctrh-lcth-cor}(a), the functors
$\sD^\st(X\ctrh_\al^\lct)\rarrow\sD^\st(X\lcth_\bW^\lct)$ and
$\Hot^\st(X\ctrh_\al^\lin)\rarrow\sD^\st(X\lcth_\bW^\lin)$ induced
by the respective embeddings of exact categories are equivalences
of triangulated categories.
 It remains to refer to Lemma~\ref{cta-clp-restricts-to-cot-inj}.
\end{proof}

\begin{cor} \label{naive-co-contra-quadrality}
 Let $X$ be a quasi-compact semi-separated scheme with an open
covering\/~$\bW$.
 Then, for any symbol\/ $\bst=+$ or\/~$\empt$, there is a commutative
square diagram of triangulated equivalences of conventional derived
categories of quasi-coherent sheaves and locally contraherent cosheaves,
\begin{equation}
\begin{gathered}
 \xymatrix{
  \sD^\st(X\qcoh^\cot) \ar@{=}[r]
  \ar@<-2pt>[d] \ar@<2pt>@{-}[d]
  & \sD^\st(X\lcth^\lct_\bW) \ar@<2pt>[d] \ar@<-2pt>@{-}[d] \\
  \sD^\st(X\qcoh) \ar@{=}[r] & \sD^\st(X\lcth_\bW)
 }
\end{gathered}
\end{equation}
 Here the lower horizontal equivalence in provided by
Theorem~\textup{\ref{naive-co-contra-thm}} with
Corollary~\textup{\ref{ctrh-lcth-cor}(a)}, while the upper horizontal
equivalence is Corollary~\textup{\ref{inj-co-contra-cor}(a)}.
 The vertical equivalences are induced by the embeddings of
exact categories $X\qcoh^\cot\rarrow X\qcoh$ and
$X\lcth_\bW^\lct\rarrow X\lcth_\bW$.
\end{cor}

\begin{proof}
 It is clear from the constructions that the diagram is commutative.
 The rightmost vertical arrow is a triangulated equivalence by
Theorem~\ref{derived-loc-cta-loc-cot-theorem}.

 Concerning the leftmost vertical arrow, in the easy case $\bst=+$ it
is a triangulated equivalence by~\cite[Proposition~13.2.2(i)]{KS} or
the dual version of Proposition~\ref{infinite-resolutions}(a).
 This simple argument, based on the observations that any quasi-coherent
sheaf can be embedded into an injective one and all injective
quasi-coherent sheaves are cotorsion, applies to an arbitrary
(not necessarily quasi-compact or semi-separated) scheme~$X$.

 In the interesting case $\bst=\empt$, the assertion that the functor
$\sD(X\qcoh^\cot)\rarrow\sD(X\qcoh)$ is a triangulated equivalence
can be found in~\cite[Corollary~10.7]{PS6}.
 The key part of the argument is the cotorsion periodicity theorem
for quasi-coherent sheaves~\cite[Theorem~10.2 or Corollary~10.4]{PS6}.
 For a discussion of another approach to the latter theorem, due to
Estrada et al., see~\cite[Remark~10.3]{PS6}.
 A (co)local approach to quasi-coherent cotorsion periodicity based on
the contraherent cosheaf theory and using
Lemmas~\ref{cta-clp-equivalence}
and~\ref{cta-clp-restricts-to-cot-inj}(a) above is outlined
in~\cite[Section~8]{Pphil}.
\end{proof}

\subsection{Homotopy locally injective complexes}
\label{homotopy-lin-subsect}
 Let $X$ be a quasi-compact semi-separated scheme and $\bW$ be its
open covering.
 The goal of this section is to construct a full subcategory in
the homotopy category $\Hot(X\lcth_\bW^\lin)$ that would be equivalent
to the unbounded derived category $\sD(X\lcth_\bW)$.
 The significance of this construction is best illustrated using
the duality-analogy between the contraherent cosheaves and
the quasi-coherent sheaves.

 As usually, the notation $\sD(X\qcoh_\fl)$ refers to the unbounded
derived category of the exact category of flat quasi-coherent sheaves
on~$X$.
 The full triangulated subcategory $\sD(X\qcoh_\fl)_\fl\sub
\sD(X\qcoh_\fl)$ of \emph{homotopy flat complexes} of flat quasi-coherent
sheaves on $X$ is defined as the minimal triangulated subcategory in
$\sD(X\qcoh_\fl)$ containing the objects of $X\qcoh_\fl$ and closed
under infinite direct sums (cf.\ Section~\ref{homotopy-adjusted}).
 The definition of homotopy flat complexes of flat modules over
a commutative (or even arbitary associative) ring is similar
(cf.\ Example~\ref{homotopy-adjusted-complexes-of-modules-examples}(3)).
 The following result is essentially due to Spaltenstein~\cite{Spal}.

\begin{thm}  \label{homotopy-flat-thm}
\textup{(a)} The composition of natural triangulated functors\/ 
$\sD(X\qcoh_\fl)_\fl\allowbreak\rarrow\sD(X\qcoh_\fl)\rarrow\sD(X\qcoh)$
is an equivalence of triangulated categories. \par
\textup{(b)} A complex of flat quasi-coherent sheaves\/ $\F^\bu$ on $X$
belongs to\/ $\sD(X\qcoh_\fl)_\fl$ if and only if its tensor product\/
$\F^\bu\ot_{\O_X}\M^\bu$ with any acyclic complex of quasi-coherent
sheaves\/ $\M^\bu$ on $X$ is also an acyclic complex of
quasi-coherent sheaves. \par
\textup{(c)} A complex of flat quasi-coherent sheaves\/ $\F^\bu$ on $X$
belongs to\/ $\sD(X\qcoh_\fl)_\fl$ if and only if the complex of flat\/
$\O_X(U)$\+modules\/ $\F^\bu(U)$ is homotopy flat for every affine open
subscheme $U\sub X$.
 It suffices to check this condition for affine open subschemes
forming any chosen affine open covering of~$X$.
\end{thm}

\begin{proof}
 Part~(a) is a particular case of Proposition~\ref{spal-for-exact}.
 To prove part~(b), let us first show that the tensor product of
any complex of quasi-coherent sheaves and a complex of sheaves
acyclic with respect to the exact category $X\qcoh_\fl$ is acyclic.
 Indeed, any complex in an abelian category is a locally
stabilizing inductive limit of finite complexes; so it suffices to
notice that the tensor product of any quasi-coherent sheaf with
a complex acyclic with respect to $X\qcoh_\fl$ is acyclic.
 Hence the class of all complexes of flat quasi-coherent sheaves
satisfying the condition in part~(b) can be viewed as a strictly
full triangulated subcategory in $\sD(X\qcoh_\fl)$.

 Now the ``only if'' assertion easily follows from the facts that
the tensor products of quasi-coherent sheaves preserve infinite
direct sums and the tensor product with a flat quasi-coherent
sheaf is an exact functor.
 In view of (the proof of) part~(a), it suffices to show that any
complex $\F^\bu$ over $X\qcoh_\fl$ satisfying the tensor product
condition of part~(b) and acyclic with respect to $X\qcoh$ is also
acyclic with respect to $X\qcoh_\fl$ in order to prove the ``if''.

 Notice that the tensor product of a bounded above complex of
flat quasi-coherent sheaves and an acyclic complex of quasi-coherent
sheaves is an acyclic complex.
 Since any quasi-coherent sheaf $\M$ over $X$ has a flat resolution,
it follows that the complex of quasi-coherent sheaves
$\M\ot_{\O_X}\F^\bu$ is acyclic.
 One easily concludes that the complex $\F^\bu$ is acyclic with
respect to $X\qcoh_\fl$.

 Part~(c) is a corollary of~(b), as the tensor product condition
in~(b) is clearly local.
 It is helpful to keep in mind that any $\O(U)$\+module can be
extended to a quasi-coherent sheaf on~$X$.
\end{proof}

 The full triangulated subcategory $\sD(X\lcth_\bW^\lin)^\lin$ of
\emph{homotopy locally injective complexes} of locally injective
$\bW$\+locally contraherent cosheaves on $X$ is defined as
the minimal full triangulated subcategory in $\sD(X\lcth_\bW^\lin)$
containing the objects of $X\lcth_\bW^\lin$ and closed under
infinite products (cf.\ Section~\ref{homotopy-adjusted}).
 The definition of homotopy injective complexes of injective modules
over a ring is similar (cf.\
Example~\ref{homotopy-adjusted-complexes-of-modules-examples}(2)).
 A terminological discussion can be found in~\cite[Remark~6.4]{PS4}.

 Given a complex of quasi-coherent sheaves $\M^\bu$ and a complex
of $\bW$\+locally contraherent cosheaves $\P^\bu$ on $X$ such that
the $\bW$\+locally contraherent cosheaf $\Cohom_X(\M^i,\P^j)$ is
defined for all $i$, $j\in\boZ$ (see
Sections~\ref{cohom-subsection}
and~\ref{cohom-loc-der-contrahereable}), we define the complex
$\Cohom_X(\M^\bu,\P^\bu)$ as the total complex of the bicomplex
$\Cohom_X(\M^i,\P^j)$ constructed by taking infinite products of
$\bW$\+locally contraherent cosheaves along the diagonals.

\begin{lem} \label{homotopy-local-injectivity-characterizations}
 Let\/ $\gJ^\bu$ be a complex of locally injective\/ $\bW$\+locally
contraherent cosheaves on~$X$.
 Then the following three conditions are equivalent:
\begin{enumerate}
\item the complex\/ $\gJ^\bu$ is homotopy locally injective, i.~e.,
belongs to\/ $\sD(X\lcth_\bW^\lin)^\lin$;
\item for every acyclic complex of quasi-coherent sheaves\/ $\M^\bu$
on $X$, the complex\/ $\Cohom_X(\M^\bu,\gJ^\bu)$ is acyclic in
the exact category $X\lcth_\bW$ (or, at one's choice, in the exact 
category $X\lcth_\bW^\lct$);
\item for every affine open subscheme $U\sub X$ subordinate to\/ $\bW$,
the complex of injective\/ $\O_X(U)$\+modules\/ $\gJ^\bu[U]$ is
homotopy injective.
 It suffices to check this condition for affine open subschemes
forming any chosen affine open covering of $X$ subordinate to\/~$\bW$.
\end{enumerate}
\end{lem}

 Let us introduce a temporary terminology ``locally homotopy injective
complexes'' for complexes of locally injective $\bW$\+locally 
contraherent cosheaves satisfying the local homotopy injectivity
condition~(3) of
Lemma~\ref{homotopy-local-injectivity-characterizations}.
 Notice the distinction between ``homotopy locally injective'' and
``locally homotopy injective''.
 The lemma claims that there is no difference.
 
\begin{proof}
 (1)~$\Longrightarrow$~(2) Start from the observation that $\Cohom_X$
from any complex of quasi-coherent sheaves into a complex of locally
contraherent cosheaves acyclic with respect to $X\lcth_\bW^\lin$ is
acyclic with respect to $X\lcth_\bW^\lct$.
 Indeed, any complex of quasi-coherent sheaves is a locally stabilizing
inductive limit of a sequence of finite complexes.
 So it remains to recall that $\Cohom_X$ from a quasi-coherent sheaf
into a complex acyclic with respect to $X\lcth_\bW^\lin$ is a complex
acyclic with respect to $X\lcth_\bW^\lct$, and use Lemma~\ref{telescope}.

 Therefore, the class of all complexes over $X\lcth_\bW^\lin$ satisfying
either one of the two $\Cohom$ conditions in~(2) can be viewed as
a strictly full triangulated subcategory in $\sD(X\lcth_\bW^\lin)$.
 Now the assertion follows from the preservation of infinite products
in the second argument by the functor $\Cohom_X$ and its exactness as
a functor $\Cohom_X({-},\gJ)\:(X\qcoh)^\op\rarrow X\lcth_\bW^\lct$ for
any fixed locally injective $\bW$\+locally  contraherent cosheaf $\gJ$
in the second argument.
 (For a dual-analogous/complementary result, see
Corollary~\ref{conventional-derived-Cohom-well-defined} below.)

 (2)~$\Longrightarrow$~(3) Follows from the definitions and
Example~\ref{homotopy-adjusted-complexes-of-modules-examples}(2).
 Once again, it is helpful to keep in mind that any $\O(U)$\+module
can be extended to a quasi-coherent sheaf on~$X$.

 (3)~$\Longrightarrow$~(1) Let $\gJ^\bu$ be a locally homotopy injective
complex in $X\lcth_\bW^\lin$.
 Choose a finite affine open covering $X=\bigcup_{\alpha=1}^N U_\alpha$
subordinate to~$\bW$, and consider the \v Cech
resolution~\eqref{contraherent-cech}
\begin{multline} \label{cech-resolution-of-lin-complex} \textstyle
 0 \lrarrow j_{1,\dotsc,N}{}_!j_{1,\dotsc,N}^!\gJ^\bu \lrarrow \dotsb \\
 \textstyle\lrarrow \bigoplus_{\alpha<\beta}j_{\alpha,\beta}{}_!
 j_{\alpha,\beta}^!\gJ^\bu \lrarrow
 \bigoplus_\alpha j_\alpha{}_!j_\alpha^!\gJ^\bu\lrarrow\gJ^\bu\lrarrow0.
\end{multline}
 Notice that all the terms of~\eqref{cech-resolution-of-lin-complex},
except perhaps the rightmost one, are homotopy locally injective
complexes (indeed, the class of homotopy locally injective complexes
of locally injective contraherent cosheaves is obviously preserved by
direct images with respect to flat affine morphisms).
 By the definition, it follows that the total complex of
the truncated resolution~\eqref{cech-resolution-of-lin-complex} (i.~e.,
of \eqref{cech-resolution-of-lin-complex} without the rightmost term)
is homotopy locally injective.
 But this total complex is isomorphic to $\gJ^\bu$ in
$\sD(X\lcth_\bW^\lin)$, and even in $\sD^\abs(X\lcth_\bW^\lin)$.
\end{proof}

 The following pair of lemmas constitutes another application of
Proposition~\ref{loc-contraherent-gluing-preenvelope}.
 We use the terminology and notation of
Appendix~\ref{cotorsion-pairs-appx}.

\begin{lem} \label{acyclic-complexes-of-cta-cotorsion-pair}
 Let\/ $\sR$ be the local class of all commutative rings $R$ and\/
$\sE^R=\sK^R=\Com(R\modl^\cta)$ be the exact category of complexes of
contraadjusted $R$\+modules.
 Let\/ $\sF(R)\sub\sE^R$ be the class of all acyclic complexes of
contraadjusted $R$\+modules, and let $\sC^R\sub\sE^R$ be the class of
all homotopy injective complexes of injective $R$\+modules.
 Then\/ $\sE$ and\/ $\sC$ are very colocal classes, and the pair of
classes $(\sF(R),\sC^R)$ is a hereditary complete cotorsion pair
in\/~$\sE^R$.
\end{lem}

\begin{proof}
 The class $\sE^R$ is very colocal because the class of contraadjusted
$R$\+modules is very colocal; and the class $\sC^R$ is very colocal
by Examples~\ref{colocal-classes-examples}.
 The full subcategory $\sF(R)$ is obviously closed under the kernels
of admissible epimorphisms in~$\sE^R$.
 It remains to mention that $(\sF(R),\sC^R)$ is a complete cotorsion
pair in $\sE^R$ according to~\cite[Example~7.8]{Pal}.
\end{proof}

\begin{lem} \label{complexes-loc-hot-inj-of-inj-preenvelope}
 Any complex of\/ $\bW$\+locally contraherent cosheaves\/ $\gM^\bu$ on
$X$ can be included in a (termwise admissible) short exact sequence of
complexes\/ $0\rarrow\gM^\bu\rarrow\gJ^\bu\rarrow\Q^\bu\rarrow0$,
where\/ $\gJ^\bu$ is a locally homotopy injective complex of locally
injective\/ $\bW$\+locally contraherent cosheaves and\/ $\Q^\bu$ is
an acyclic complex in the exact category of antilocal contraherent
cosheaves on~$X$.
\end{lem}

\begin{proof}
 This is similar to Lemma~\ref{complexes-loc-cot-preenvelope}.
 Let $X=\bigcup_\alpha U_\alpha$ be an affine open covering
subordinate to~$\bW$.
 We apply Proposition~\ref{loc-contraherent-gluing-preenvelope} to
the datum of classes $\sR$, \,$\sE^R$, \,$\sC^R$, and $\sF(R)$ from
Lemma~\ref{acyclic-complexes-of-cta-cotorsion-pair}.
 This produces, for a given complex $\gM^\bu$, a short exact sequence of
complexes $0\rarrow\gM^\bu\rarrow\gJ^\bu\rarrow\Q^\bu\rarrow0$, where
$\gJ^\bu$ is a locally homotopy injective complex of locally injective
$\bW$\+locally contraherent cosheaves and $\Q^\bu$ is a finitely
iterated extension in $\Com(X\lcth_\bW)$ of direct images of acyclic
complexes $\Q_\alpha^\bu$ in the exact categories $U_\alpha\ctrh$
of contraherent cosheaves on~$U_\alpha$.
 It remains to recall that the direct image functor
$j_\alpha{}_!\:U_\alpha\ctrh\rarrow X\ctrh_\al\sub X\ctrh$ is exact
(see the beginning of Section~\ref{clp-subsection};
cf.\ Corollary~\ref{clp-direct}(a)).
\end{proof}

 In fact, we will see in Chapter~\ref{becker-on-qcomp-qsep-sect}
that any acyclic complex in the exact category $X\ctrh_\al$ is
a direct summand of a finitely iterated extension of direct images
of acyclic complexes in the exact categories $U_\alpha\ctrh$.
 This is the result of
Corollary~\ref{acyclic-complexes-of-lcta-lct-antilocal}(a).

\begin{thm} \label{homotopy-lin-thm}
 The composition of natural triangulated functors\/
$\sD(X\lcth_\bW^\lin)^\lin\allowbreak\rarrow\sD(X\lcth_\bW^\lin)\rarrow
\sD(X\lcth_\bW)$ is an equivalence of triangulated categories.
\end{thm}

\begin{proof}
 For any complex of $\bW$\+locally contraherent cosheaves $\gM^\bu$,
Lemma~\ref{complexes-loc-hot-inj-of-inj-preenvelope} provides
a quasi-isomorphism $\gM^\bu\rarrow\gJ^\bu$ of complexes in
$X\lcth_\bW$ from $\gM^\bu$ into a locally homotopy injective complex of
locally injective $\bW$\+locally contraherent cosheaves~$\gJ^\bu$.
 (Indeed, any acyclic complex in $X\ctrh_\al$ is also acyclic in
$X\lcth_\bW$.)
 Furthermore, Lemma~\ref{homotopy-local-injectivity-characterizations}%
\,(3)\,$\Rightarrow$\,(1) tells us that all locally homotopy injective
complexes are homotopy locally injective.

 From this point on, the argument proceeds similarly (or
dual-analogously) to the proof of Proposition~\ref{spal-for-exact}.
 By Lemma~\ref{pkoszul-lemma16}(b), it follows from the previous
paragraph, in particular, that $\sD(X\lcth_\bW)$ is equivalent to
the quotient category of $\sD(X\lcth_\bW^\lin)$ by the thick subcategory
of complexes over $X\lcth_\bW^\lin$ acyclic in $X\lcth_\bW$.
 By the dual version of Corollary~\ref{hom-vanishing-cor}, the latter
thick subcategory is semiorthogonal to $\sD(X\lcth_\bW^\lin)^\lin$
in $\sD(X\lcth_\bW^\lin)$.
 Now the same construction of a quasi-isomorphism in $X\lcth_\bW$,
specialized to complexes with the terms in $X\lcth_\bW^\lin$, also
implies that these two full subcategories form a semiorthogonal
decomposition of $\sD(X\lcth_\bW^\lin)$.
\end{proof}

\begin{lem} \label{clp-linlin-orthogonality}
 Let\/ $\P^\bu$ be a complex over the exact category $X\ctrh_\al$
and\/ $\gJ^\bu$ be a complex over $X\lcth_\bW^\lin$ belonging to\/
$\sD(X\lcth_\bW^\lin)^\lin$.
 Then the natural morphism of graded abelian groups
$H^*\Hom^X(\P^\bu,\gJ^\bu)\simeq\Hom_{\Hot(X\lcth_\bW)}(\P^\bu,\gJ^\bu)
\rarrow\Hom_{\sD(X\lcth_\bW)}(\P^\bu,\gJ^\bu)$ is an isomorphism
(in other words, the complex\/ $\Hom^X(\P^\bu,\gJ^\bu)$ computes
the groups of morphisms in the derived category\/ $\sD(X\lcth_\bW)$).
\end{lem}

\begin{proof}
 It is convenient to represent morphisms $\P^\bu\rarrow\gJ^\bu$
in $\sD(X\lcth_\bW)$ by fractions $\P^\bu\larrow\Q^\bu\rarrow\gJ^\bu$
of morphisms in $\Hot(X\lcth_\bW)$.
 Since we know that any complex over $X\lcth_\bW$ admits
a quasi-isomorphism into it from a complex over $X\ctrh_\al$ and any
complex over $X\ctrh_\al$ acyclic over $X\lcth_\bW$ is also acyclic over
$X\ctrh_\al$ (by Corollary~\ref{ctrh-lcth-cor}(a)), it suffices to show
that the complex of abelian groups $\Hom^X(\P^\bu,\gJ^\bu)$ is acyclic
for any complex $\P^\bu$ acyclic with respect to $X\ctrh_\al$ and any
complex $\gJ^\bu\in\sD(X\lcth_\bW^\lin)^\lin$.
 For a complex $\gJ^\bu$ obtained from objects of $X\lcth_\bW^\lin$ by
iterating the operations of cone and infinite product the latter
assertion is obvious (see Corollary~\ref{clp-characterizations}(a)),
so it remains to consider the case of a complex $\gJ^\bu$ acyclic
with respect to $X\lcth_\bW^\lin$.

 In this case we will show that the complex $\Hom^X(\P^\bu,\gJ^\bu)$
is acyclic for any complex $\P^\bu$ over $X\ctrh_\al$.
 Indeed, by Corollary~\ref{lin-ctrh-lcth-cor}(a) any acyclic complex
over $X\lcth_\bW^\lin$ is absolutely acyclic.
 Then it remains to recall that the functor $\Hom^X$ from
an antilocal contraherent cosheaf preserves exactness of short
exact sequences in $X\lcth_\bW^\lin$, by the definition.
\end{proof}

\subsection{Derived functors of direct and inverse image}
\label{derived-direct-inverse}
 For the rest of Chapter~\ref{derived-on-quasi-compact-sect}, the upper
index~$\bst$ in the notation for derived and homotopy categories
stands for one of the symbols $\b$, $+$, $-$, $\empt$, $\abs+$,
$\abs-$, $\bco$, $\bctr$, $\co$, $\ctr$, or~$\abs$.
 The discussion of the Becker co/contraderived category symbols
$\bco$ and $\bctr$ is actually postponed to
Sections~\ref{direct-images-in-becker-subsect}\+-%
\ref{derived-tensor-operations-subsect}, but we include these symbols
in the list for consistency.

 Let $f\:Y\rarrow X$ be a morphism of quasi-compact semi-separated
schemes.
 Then for any symbol $\bst\ne\co$, $\ctr$, $\bco$, $\bctr$ the right
derived functor of direct image
\begin{equation}  \label{qcoh-direct}
 \boR f_*\:\sD^\st(Y\qcoh)\lrarrow\sD^\st(X\qcoh)
\end{equation}
is constructed in the following way.
 By Corollary~\ref{qcoh-dil-cta-derived-equiv-cor}(b), the natural
functor $\sD^\st(Y\qcoh^\cta)\rarrow\sD^\st(Y\qcoh)$ is an equivalence
of triangulated categories (as is the similar functor for sheaves
over~$X$).
 By Corollary~\ref{cta-cot-direct}(a), the restriction of the functor
of direct image $f_*\:Y\qcoh\rarrow X\qcoh$ provides an exact functor
$Y\qcoh^\cta\rarrow X\qcoh^\cta$.
 Now the derived functor $\boR f_*$ is defined by restricting
the functor of direct image $f_*\:\Hot^\st(Y\qcoh)\rarrow
\Hot^\st(X\qcoh)$ to the full subcategory of complexes of
contraadjusted quasi-coherent sheaves on~$Y$.

 More generally, for any symbol $\bst\ne\ctr$, $\bco$, $\bctr$
the right derived functor~$\boR f_*$ \eqref{qcoh-direct} can be
constructed using dilute coresolutions.
 By Corollary~\ref{qcoh-dil-cta-derived-equiv-cor}(a), the natural
functor $\sD^\st(Y\qcoh^\dil)\rarrow\sD^\st(Y\qcoh)$ is an equivalence
of triangulated categories (as is the similar functor for sheaves
over~$X$).
 By Corollary~\ref{dilute-direct}, the restriction of the direct image
functor $f_*\:Y\qcoh\rarrow X\qcoh$ is an exact functor
$f_*\:Y\qcoh^\dil\rarrow X\qcoh^\dil$.
 Now the derived functor $\boR f_*$ \eqref{qcoh-direct} is defined as
the induced functor $\sD^\st(Y\qcoh^\dil)\rarrow\sD^\st(X\qcoh^\dil)$.
 For a construction of the functor $\boR f_*$ in the case
$\bst=\bco$, see Section~\ref{direct-images-in-becker-subsect}.

 For any symbol $\bst\ne\co$, $\bco$, $\bctr$, the left derived functor
of direct image
\begin{equation}   \label{ctrh-direct}
 \boL f_!\:\sD^\st(Y\ctrh)\lrarrow\sD^\st(X\ctrh)
\end{equation}
is constructed in the following way.
 By Corollary~\ref{ctrh-lcth-cor}(a), the natural functor
$\sD^\st(Y\ctrh_\al)\rarrow\sD^\st(Y\ctrh)$ is an equivalence of
triangulated categories (as is the similar functor for cosheaves
over~$X$).
 By Corollary~\ref{clp-direct}(a), there is an exact functor of
direct image $f_!\:Y\ctrh_\al\rarrow X\ctrh_\al$.
 The derived functor $\boL f_!$ is defined as the induced functor
$\sD^\st(Y\ctrh_\al)\rarrow\sD^\st(X\ctrh_\al)$.
 For a construction of the functor $\boL f_!$ in the case $\bst=\bctr$,
see Section~\ref{direct-images-in-becker-subsect}.
 
 Similarly one defines the left derived functor of direct image
\begin{equation}   \label{ctrh-lct-direct}
 \boL f_!\:\sD^\st(Y\ctrh^\lct)\rarrow\sD^\st(X\ctrh^\lct)
\end{equation}
for any symbol $\bst\ne\co$, $\bco$, $\bctr$.
 Once again, for the case $\bst=\bctr$ we refer to
Section~\ref{direct-images-in-becker-subsect}.

 By construction, the triangulated functors~\eqref{ctrh-direct}
and~\eqref{ctrh-lct-direct} form a commutative diagram with
the triangulated functors $\sD^\st(Y\ctrh^\lct)\rarrow\sD^\st(Y\ctrh)$
and $\sD^\st(X\ctrh^\lct)\rarrow\sD^\st(X\ctrh)$ induced by
the respective embeddings of exact categories.
 In particular, for $\bst=+$ or~$\empt$, we obtain a commutative
diagram of triangulated functors and triangulated equivalences
provided by Theorem~\ref{derived-loc-cta-loc-cot-theorem},
\begin{equation} \label{derived-lct-lcta-direct-images-compatible}
\begin{gathered}
 \xymatrix{
  \sD^\st(Y\ctrh^\lct) \ar@<2pt>[r] \ar@<-2pt>@{-}[r] \ar[d]^{\boL f_!}
  & \sD^\st(Y\ctrh) \ar[d]^{\boL f_!} \\
  \sD^\st(X\ctrh^\lct) \ar@<2pt>[r] \ar@<-2pt>@{-}[r] & \sD^\st(X\ctrh)
 }
\end{gathered}
\end{equation}

\begin{thm} \label{direct-images-identified}
 Let $f\:Y\rarrow X$ be a morphism of quasi-compact semi-separated 
schemes.
 Then, for any symbol\/ $\bst\ne\co$, $\ctr$, $\bco$, $\bctr$,
the equivalences of categories\/ $\sD^\st(Y\qcoh)\simeq\sD^\st(Y\ctrh)$
and\/ $\sD^\st(X\qcoh)\simeq \sD^\st(X\ctrh)$ from
Theorem~\textup{\ref{naive-co-contra-thm}} transform the right derived
functor\/ $\boR f_*$ \textup{\eqref{qcoh-direct}} into the left derived
functor\/ $\boL f_!$~\textup{\eqref{ctrh-direct}}.
 In other words, the following square diagram of triangulated functors
and triangulated equivalences is commutative:
\begin{equation} \label{derived-cta-clp-direct-image-diagram}
\begin{gathered}
 \xymatrix{
  \sD^\st(Y\qcoh) \ar@{=}[r] \ar[d]_{\boR f_*}
  & \sD^\st(Y\ctrh) \ar[d]^{\boL f_!} \\
  \sD^\st(X\qcoh) \ar@{=}[r] & \sD^\st(X\ctrh)
 }
\end{gathered}
\end{equation}
\end{thm}

\begin{proof}
 It suffices to show that the equivalences of exact categories
$Y\qcoh^\cta\simeq Y\ctrh_\al$ and $X\qcoh^\cta\simeq X\ctrh_\al$
from Lemma~\ref{cta-clp-equivalence} transform the functor~$f_*$
into the functor~$f_!$.
 This is the result of
Corollary~\ref{cta-clp-direct-image-cor}.
\end{proof}

 Let $f\:Y\rarrow X$ be a morphism of schemes.
 Let $\bW$ and $\bT$ be open coverings of the schemes $X$ and $Y$,
respectively, for which the morphism~$f$ is $(\bW,\bT)$\+coaffine.
 According to Section~\ref{direct-inverse-loc-contra}, there is an exact
functor of inverse image $f^!\:X\lcth_\bW^\lin\rarrow Y\lcth_\bT^\lin$;
for a flat morphism~$f$, there is also an exact functor
$f^!\:X\lcth_\bW^\lct\rarrow Y\lcth_\bT^\lct$, and for a very flat
morphism~$f$, an exact functor $f^!\:X\lcth_\bW\rarrow Y\lcth_\bT$.

 For a quasi-compact semi-separated scheme~$X$, it follows from
Corollary~\ref{quasi-cotors-cor}(a) and
Proposition~\ref{infinite-resolutions}(a) that the natural functor
$\sD^-(X\qcoh_\fl)\rarrow\sD^-(X\qcoh)$ is an equivalence of categories.
 Similarly, it follows from Corollary~\ref{clp-cor}(a) and
the dual version of Proposition~\ref{infinite-resolutions}(a) that
the natural functor $\sD^+(X\lcth_\bW^\lin)\rarrow\sD^+(X\lcth_\bW)$
is an equivalence of triangulated categories.
 This allows to define, for any morphism $f\:Y\rarrow X$ into
a quasi-compact semi-separated scheme $X$ and coverings $\bW$, $\bT$
as above, the derived functors of inverse image
\begin{equation}  \label{qcoh-minus-inverse}
 \boL f^*\:\sD^-(X\qcoh)\lrarrow\sD^-(Y\qcoh)
\end{equation}
and
$$  
 \boR f^!\:\sD^+(X\lcth_\bW)\lrarrow\sD^+(Y\lcth_\bT)
$$  
by applying the functors $f^*$ and~$f^!$ to (appropriately bounded)
complexes of flat sheaves and locally injective cosheaves.
 When both schemes are quasi-compact and semi-separated, one can
take into account the equivalences of categories from
Corollary~\ref{ctrh-lcth-cor}(a) in order to produce the right derived
functor
\begin{equation}  \label{ctrh-plus-inverse}
 \boR f^!\:\sD^+(X\ctrh)\lrarrow\sD^+(Y\ctrh),
\end{equation}
which clearly does not depend on the choice of the coverings
$\bW$ and~$\bT$.

 According to Section~\ref{homotopy-lin-subsect}, the natural functors
$\sD(X\qcoh_\fl)_\fl\rarrow\sD(X\qcoh)$ and $\sD(X\lcth_\bW^\lin)^\lin
\rarrow\sD(X\lcth_\bW)$ are equivalences of categories for any
quasi-compact semi-separated scheme $X$ with an open covering~$\bW$.
 This allows to define, for any morphism $f\:Y\rarrow X$ into
a quasi-compact semi-separated scheme $X$ and coverings $\bW$, $\bT$
as above, the derived functors of inverse image
\begin{equation}  \label{qcoh-inverse}
 \boL f^*\:\sD(X\qcoh)\lrarrow\sD(Y\qcoh)
\end{equation}
and
$$  
 \boR f^!\:\sD(X\lcth_\bW)\lrarrow\sD(Y\lcth_\bT)
$$  
by applying the functors $f^*\:\Hot(X\qcoh)\rarrow\Hot(Y\qcoh)$
and $f^!\:\Hot(X\lcth_\bW^\lin)\allowbreak\rarrow\Hot(Y\lcth_\bT^\lin)$
to homotopy flat complexes of flat quasi-coherent sheaves and homotopy
locally injective complexes of locally injective $\bW$\+locally
contraherent cosheaves, respectively.
 Of course, this construction is well-known for quasi-coherent
sheaves~\cite{Spal,N-bb}; we discuss here the sheaf and cosheaf
situations together in order to emphasize the duality-analogy
between them.

 When both schemes are quasi-compact and semi-separated, one can
use the equivalences of categories from Corollary~\ref{ctrh-lcth-cor}(a)
in order to obtain the right derived functor
\begin{equation}  \label{ctrh-inverse}
 \boR f^!\:\sD(X\ctrh)\lrarrow\sD(Y\ctrh),
\end{equation}
which does not depend on the choice of the coverings $\bW$ and~$\bT$.
 Notice also that the restriction of the functor~$f^*$ takes
$\Hot(X\qcoh_\fl)_\fl$ into $\Hot(Y\qcoh_\fl)_\fl$ and the restriction
of the functor~$f^!$ takes $\Hot(X\lcth_\bW^\lin)^\lin$ into
$\Hot(Y\lcth_\bT^\lin)^\lin$.

 It is easy to see that for any morphism of quasi-compact semi-separated
schemes $f\:Y\rarrow X$ the functor $\boL f^*$ \eqref{qcoh-minus-inverse}
is left adjoint to the functor $\boR f_*\:\sD^-(Y\qcoh)\rarrow
\sD^-(X\qcoh)$ \eqref{qcoh-direct} and the functor $\boR f^!$
\eqref{ctrh-plus-inverse} is right adjoint to the functor
$\boL f_!\:\sD^+(Y\ctrh)\allowbreak\rarrow \sD^+(X\ctrh)$
\eqref{ctrh-direct}.
 Essentially, one uses the partial adjunctions on the level of exact
categories together with the fact that the derived functor constructions
are indeed those of the ``left'' and ``right'' derived functors,
as stated (cf.~\cite[Lemma~8.3]{Psemi}).  {\hbadness=1400\par}

 The point is that the isomorphisms with (co)resolution complexes in
the derived categories used for the constructions of triangulated
category equivalences producing the derived functors are not just
fractions of morphisms in the homotopy categories, but actual
quasi-isomorphisms acting in a particular direction.
 The direction in which the quasi-isomorphism with the (co)resolution
complex acts distinguishes the left derived functors from
the right ones.

 Similarly, one concludes from the construction in the proof of
Theorem~\ref{homotopy-flat-thm} that the functor $\boL f^*$
\eqref{qcoh-inverse} is left adjoint to the functor
$\boR f_*\:\sD(Y\qcoh)\rarrow\sD(X\qcoh)$ \eqref{qcoh-direct}.
 And in order to show that the functor $\boR f^!$ \eqref{ctrh-inverse}
is right adjoint to the functor $\boL f_!\:\sD(Y\ctrh)\rarrow
\sD(X\ctrh)$ \eqref{ctrh-direct}, one can look into the construction
in the proof of Theorem~\ref{homotopy-lin-thm} or use
Lemma~\ref{clp-linlin-orthogonality}.
 So we have obtained a new proof of the following classical
result~\cite{Del}, \cite[Example~4.2]{N-bb}, \cite[Chapter~4]{Lip}.

\begin{cor}  \label{neemans-functor-construction}
 For any morphism of quasi-compact semi-separated schemes
$f\:Y\rarrow X$, the derived direct image functor\/ $\boR f_*\:
\sD(Y\qcoh)\rarrow\sD(X\qcoh)$ has a right adjoint functor
$f^!\:\sD(X\qcoh)\rarrow\sD(Y\qcoh)$.
\end{cor}

\begin{proof}
 We have explicitly constructed the functor~$f^!$ as
the right derived functor $\boR f^!\:\sD(X\ctrh)\rarrow\sD(Y\ctrh)$
\eqref{ctrh-inverse} of the exact functor $f^!\:X\lcth^\lin\rarrow
Y\lcth^\lin$ between exact subcategories of the exact categories of
locally contraherent cosheaves on $X$ and~$Y$.
 The above construction of the functor~$f^!$ for bounded below
complexes~\eqref{ctrh-plus-inverse} is particularly simple and
elementary.
 In either case, the construction is based on the identification of
the functor $\boR f_*$ of derived direct image of quasi-coherent
sheaves with the functor $\boL f_!$ of derived direct image of
contraherent cosheaves, which is provided by
Theorems~\ref{naive-co-contra-thm} and~\ref{direct-images-identified}.
\end{proof}

\subsection{Morphisms of finite flat dimension}
\label{finite-dim-morphisms-I}
 A morphism of schemes $f\:Y\rarrow X$ is said to have \emph{flat
dimension not exceeding~$D$} if for all affine open subschemes $U\sub X$
and $V\sub Y$ such that $f(V)\sub U$ the $\O_X(U)$\+module $\O_Y(V)$
has flat dimension not exceeding~$D$.
 The morphism~$f$ has \emph{very flat dimension not exceeding~$D$} if
the similar bound holds for the very flat dimension of
the $\O_X(U)$\+modules $\O_Y(V)$.

 For any morphism $f\:Y\rarrow X$ of finite flat dimension~$\le D$ into
a quasi-compact semi-separated scheme $X$ and any symbol $\bst\ne\ctr$,
$\bco$, $\bctr$, the left derived functor
\begin{equation}  \label{qcoh-inverse-ffd-morphism}
 \boL f^*\:\sD^\st(X\qcoh)\lrarrow\sD^\st(Y\qcoh)
\end{equation}
is constructed in the following way.
 Let us call a quasi-coherent sheaf $\F$ on $X$ \emph{adjusted to~$f$}
if for all affine open subschemes $U\sub X$ and $V\sub Y$ such that
$f(V)\sub U$ one has $\Tor^{\O_X(U)}_{>0}(\O_Y(V),\F(U))=0$.
 Quasi-coherent sheaves $\F$ on $X$ adjusted to~$f$ form a resolving
subcategory $X\qcoh_\fadj\sub X\qcoh$ closed under infinite direct sums
and such that the resolution dimension of any quasi-coherent sheaf
on $X$ with respect to $X\qcoh_\fadj$ does not exceed~$D$.
 By Proposition~\ref{finite-resolutions}, it follows that
the natural functor $\sD^\st(X\qcoh_\fadj)\rarrow\sD^\st(X\qcoh)$ is
an equivalence of triangulated categories.

 The right exact functor $f^*\:X\qcoh\rarrow Y\qcoh$ restricts to
an exact functor $f^*\:X\qcoh_\fadj\rarrow Y\qcoh$.
 In view of the above equivalence of categories, the induced
functor on the derived categories $f^*\:\sD^\st(X\qcoh_\fadj)\rarrow
\sD^\st(Y\qcoh)$ provides the desired derived functor $\boL f^*$.
 (For a construction of the functor $\boL f^*$ in the case $\bst=\bco$,
we refer to Section~\ref{finite-dim-morphisms-derived-inverse-subsect}.)
 For any morphism of finite flat dimension $f\:Y\rarrow X$ between
quasi-compact semi-separated schemes $Y$ and $X$, the functor
$\boL f^*$ is left adjoint to the functor $\boR f_*$ \eqref{qcoh-direct}
from Section~\ref{derived-direct-inverse}
(cf.~\cite[Section~1.9]{EP}).

 For any morphism $f\:Y\rarrow X$ of finite very flat dimension~$\le D$
into a quasi-compact semi-separated scheme $X$, any open coverings $\bW$
and $\bT$ of the schemes $X$ and $Y$ for which the morphism~$f$ is
$(\bW,\bT)$\+coaffine, and any symbol $\bst\ne\co$, $\bco$, $\bctr$,
the right derived functor
\begin{equation}  \label{lcth-inverse-ffd-morphism}
 \boR f^!\:\sD^\st(X\lcth_\bW)\lrarrow\sD^\st(Y\lcth_\bT)
\end{equation}
is constructed in the following way.
 Let us call a $\bW$\+locally contraherent cosheaf $\P$ on $X$
\emph{adjusted to~$f$} if for all affine open subschemes $U\sub X$
and $V\sub Y$ such that $U$ is subordinate to $\bW$ and $f(V)\sub U$
one has $\Ext_{\O_X(U)}^{>0}(\O_Y(V),\P[U])=0$.
 One can easily see that the adjustedness condition does not change
when restricted to open subschemes $V$ subordinate to $\bT$, nor
it is changed by a refinement of the covering~$\bW$ (the argument
is based on the \v Cech (co)resolutions~\eqref{cech-quasi}
and~\eqref{contraherent-cech}).

 Locally contraherent cosheaves on $X$ adjusted to~$f$ form
a coresolving subcategory $X\lcth^\fadj\sub X\lcth$; the subcategory
$X\lcth_\bW^\fadj = X\lcth^\fadj\cap X\lcth_\bW$ is also closed
under infinite products in $X\lcth_\bW$.
 The coresolution dimension of any $\bW$\+locally contraherent cosheaf
on $X$ with respect to $X\lcth_\bW^\fadj$ does not exceed~$D$.
 By the dual version of Proposition~\ref{finite-resolutions},
the natural functor $\sD^\st(X\lcth_\bW^\fadj)\rarrow
\sD^\st(X\lcth_\bW)$ is an equivalence of triangulated categories.
 The construction of the exact functor $f^!\:X\lcth_\bW^\lin\rarrow
Y\lcth_\bT^\lin$ from Section~\ref{direct-inverse-loc-contra} extends
without any changes to the case of cosheaves from $X\lcth_\bW^\fadj$,
defining an exact functor
$$
 f^!\:X\lcth_\bW^\fadj\lrarrow Y\lcth_\bT.
$$
 Instead of Lemma~\ref{very-scalars-veryflat-case}(a) used in
Sections~\ref{contra-direct-inverse}
and~\ref{direct-inverse-loc-contra},
one can use Lemma~\ref{adjusted-scalars-contra}(a) in order to check
that the contraadjustedness condition is preserved here.

 In view of the above equivalence of triangulated categories,
the induced functor $f^!\:\sD^\st(X\lcth_\bW^\fadj)\rarrow
\sD^\st(Y\lcth_\bT)$ provides the desired functor $\boR f^!$
\eqref{lcth-inverse-ffd-morphism}.
 (For a construction of the functor $\boR f^!$ in the case $\bst=\bctr$,
we refer to Section~\ref{finite-dim-morphisms-derived-inverse-subsect}.)
 When both schemes $X$ and $Y$ are quasi-compact and semi-separated,
one can use the equivalences of categories from
Corollary~\ref{ctrh-lcth-cor}(a) in order to obtain the right
derived functor
\begin{equation}  \label{ctrh-inverse-ffd-morphism}
 \boR f^!\:\sD^\st(X\ctrh)\lrarrow\sD^\st(Y\ctrh),
\end{equation}
which is right adjoint to the functor $\boL f_!$ \eqref{ctrh-direct}
from Section~\ref{derived-direct-inverse}.

 For a morphism $f\:Y\rarrow X$ of flat dimension~$\le D$ into
a quasi-compact semi-separated scheme $X$, any open coverings $\bW$
and $\bT$ of the schemes $X$ and $Y$ for which the morphism~$f$
is $(\bW,\bT)$\+coaffine, and any symbol $\bst\ne\co$, $\bco$, $\bctr$,
one can similarly construct the right derived functor
\begin{equation}  \label{lcth-lct-inverse-ffd-morphism}
 \boR f^!\:\sD^\st(X\lcth_\bW^\lct)\lrarrow\sD^\st(Y\lcth_\bT^\lct).
\end{equation}
 More precisely, a locally cotorsion $\bW$\+locally contraherent cosheaf
$\P$ on $X$ is said to be \emph{adjusted to~$f$} if for all affine
open subschemes $U\sub X$ and $V\sub Y$ such that $U$ is subordinate
to~$\bW$ and $f(V)\sub U$, and for all flat $\O_Y(V)$\+modules $G$,
one has $\Ext_{\O_X(U)}^{>0}(G,\P[U])=0$.
 This condition does not change when restricted to open subschemes
$V\sub Y$ subordinate to $\bT$, nor it is changed by any refinement
of the covering~$\bW$ of the scheme~$X$.

 As above, locally cotorsion $\bW$\+locally contraherent cosheaves
on $X$ adjusted to~$f$ form a coresolving subcategory
$X\lcth_\bW^{\lct,\,\fadj}\sub X\lcth_\bW^\lct$ closed under
infinite products.
 It follows from Lemma~\ref{affine-d+D-lemma}(a) (for $d=0$) that
the coresolution dimension of any locally cotorsion $\bW$\+locally
contraherent cosheaf on $X$ with respect to $X\lcth_\bW^{\lct,\,\fadj}$
does not exceed~$D$.
 Hence the natural functor $\sD^\st(X\lcth_\bW^{\lct,\,\fadj})\rarrow
\sD^\st(X\lcth_\bW^\lct)$ is an equivalence of triangulated categories.
 
 Now Lemma~\ref{adjusted-scalars-contra}(b) allows to construct
an exact functor
$$
 f^!\:X\lcth_\bW^{\lct,\,\fadj}\lrarrow Y\lcth_\bT^\lct,
$$
and the induced functor $f^!\:\sD^\st(X\lcth_\bW^{\lct,\,\fadj})
\rarrow\sD^\st(Y\lcth_\bT^\lct)$ provides the desired derived
functor~\eqref{lcth-lct-inverse-ffd-morphism}.
 (Once again, for the case $\bst=\bctr$ we refer to
Section~\ref{finite-dim-morphisms-derived-inverse-subsect}.)
 When both schemes $X$ and $Y$ are quasi-compact and semi-separated,
one can use Corollary~\ref{lct-ctrh-lcth-cor}(a) in order to obtain
the right derived functor
\begin{equation}  \label{ctrh-lct-inverse-ffd-morphism}
 \boR f^!\:\sD^\st(X\ctrh^\lct)\lrarrow\sD^\st(Y\ctrh^\lct),
\end{equation}
which is right adjoint to the functor $\boL f_!$ \eqref{ctrh-lct-direct}.

 By construction, for a morphism~$f$ of finite very flat dimension,
the triangulated functors~\eqref{ctrh-inverse-ffd-morphism}
and~\eqref{ctrh-lct-inverse-ffd-morphism} form a commutative diagram
with the triangulated functors $\sD^\st(X\ctrh^\lct)\rarrow
\sD^\st(X\ctrh)$ and $\sD^\st(Y\ctrh^\lct)\rarrow\sD^\st(Y\ctrh)$
induced by the respective embeddings of exact categories.
 In particular, for quasi-compact semi-separated schemes $X$ and $Y$,
a morphism $f\:Y\rarrow X$ of finite very flat dimension, and
$\bst=+$ or~$\empt$, we obtain a commutative diagram of triangulated
functors and triangulated equivalences provided by
Theorem~\ref{derived-loc-cta-loc-cot-theorem},
\begin{equation} \label{derived-lct-lcta-inverse-images-compatible}
\begin{gathered}
 \xymatrix{
  \sD^\st(Y\ctrh^\lct) \ar@<2pt>[r] \ar@<-2pt>@{-}[r]
   & \sD^\st(Y\ctrh) \\
  \sD^\st(X\ctrh^\lct) \ar@<2pt>[r] \ar@<-2pt>@{-}[r] \ar[u]_{\boR f^!}
  & \sD^\st(X\ctrh) \ar[u]_{\boR f^!}
 }
\end{gathered}
\end{equation}

\Section{Becker's Contraderived Categories on Quasi-Compact
Semi-Separated Schemes} \label{becker-on-qcomp-qsep-sect}

\subsection{Antilocality of complexes of contraadjusted and
cotorsion sheaves} \label{quasi-cta-cot-antilocality-subsect}
 The exposition in Chaper~\ref{becker-on-qcomp-qsep-sect} is largely
based on Appendix~\ref{cotorsion-pairs-appx}.
 In particular, we use the terminology and notation of
Section~\ref{cotorsion-prelim-subsect} in connection with cotorsion
pairs in exact categories, and of Sections~\ref{local-classes-subsect}
and~\ref{colocal-classes-subsect} in connection with local/colocal
classes of commutative rings, modules, and complexes.

 For any additive category $\sE$, denote by $\Com(\sE)$
the additive category of (unbounded) complexes in~$\sE$.
 When $\sE$ is an exact category, we endow $\Com(\sE)$ with the termwise
exact structure (see Section~\ref{cotorsion-prelim-subsect}).
 For any exact category $\sE$, denote by $\Acycl(\sE)\sub\Com(\sE)$
the class of all acyclic complexes in~$\sE$ (as in
Section~\ref{small-object-argument-subsect}).

 Let $X$ be a quasi-compact semi-separated scheme.
 The result of the following lemma is obtained by combining
the small object argument in the category of complexes as
in~\cite[Lemma~4.9]{Gil3} with a suitable contraadjusted/cotorsion
periodicity theorem.

\begin{lem} \label{acyclic-of-vfl-flat-arbitrary-of-cta-cot-pairs}
\textup{(a)} The pair of classes of acyclic complexes of very flat
quasi-coherent sheaves with very flat sheaves of cocycles\/
$\sF=\Acycl(X\qcoh_\vfl)$ and arbitrary complexes of contraadjusted
quasi-coherent sheaves\/ $\sC=\Com(X\qcoh^\cta)$ is a hereditary
complete cotorsion pair in the abelian category\/ $\Com(X\qcoh)$. \par
\textup{(b)} The pair of classes of acyclic complexes of flat
quasi-coherent sheaves with flat sheaves of cocycles\/
$\sF=\Acycl(X\qcoh_\fl)$ and arbitrary complexes of cotorsion
quasi-coherent sheaves\/ $\sC=\Com(X\qcoh^\cot)$ is a hereditary
complete cotorsion pair in the abelian category\/ $\Com(X\qcoh)$.
\end{lem}

\begin{proof}
 Part~(a): the class $R\modl_\vfl$ of very flat $R$\+modules is
deconstructible in $R\modl$ (see
Examples~\ref{deconstructible-examples}).
 Hence Lemma~\ref{locally-L-sheaves-deconstructible} tells us that
the class $X\qcoh_\vfl$ of very flat quasi-coherent sheaves is
deconstructible in $X\qcoh$ for any scheme~$X$.
 By~\cite[Proposition~4.4]{Sto-hill}, it follows that the class
$\Acycl(X\qcoh_\vfl)$ of all acyclic complexes of very flat
quasi-coherent sheaves with very flat quasi-coherent sheaves of
cocycles is deconstructible in $\Com(X\qcoh)$.
 When $X$ is quasi-compact and semi-separated,
Lemma~\ref{quasi-very-flat-cover} tells us that the class $X\qcoh_\vfl$
is generating in $X\qcoh$, and it follows easily that the class
$\Acycl(X\qcoh_\vfl)$ is generating in $\Com(X\qcoh)$ (in fact,
any complex of quasi-coherent sheaves on $X$ is a quotient complex
of a contractible complex of very flat quasi-coherent sheaves).
 Obviously, the class $\Acycl(X\qcoh_\vfl)$ is closed under direct
summands in $\Com(X\qcoh)$.
 Applying Corollary~\ref{eklof-trlifaj-corollary}, we conclude that
there exists a complete cotorsion pair $(\sF,\sC)$ in $\Com(X\qcoh)$
with $\sF=\Acycl(X\qcoh_\vfl)$.
 The class $\Acycl(X\qcoh_\vfl)$ is closed under kernels of epimorphisms
in $\Com(X\qcoh)$, so the cotorsion pair $(\sF,\sC)$ is hereditary.

 It remains to show that $\sC$ is the class of all complexes of
contraadjusted quasi-coherent sheaves on~$X$.
 Indeed, the class $\sF$ contains all the contractible two-term
complexes of very flat sheaves $\dotsb\rarrow0\rarrow\F\overset=\rarrow
\F\rarrow0\rarrow\dotsb$, where $\F\in X\qcoh_\vfl$.
 By Lemma~\ref{G-plus-G-minus-complexes-Ext-adjunction} (for $i=1$),
it follows that $\sC\sub\Com(X\qcoh^\cta)$.
 Conversely, let $\cP^\bu\in\Com(X\qcoh^\cta)$ be a complex of
contraadjusted quasi-coherent sheaves.
 We argue similarly to~\cite[proof of Theorem~5.3]{BCE}.
 Consider a special preenvelope sequence $0\rarrow\cP^\bu\rarrow\cQ^\bu
\rarrow\G^\bu\rarrow0$ \,\eqref{special-preenvelope-sequence} with
$\cQ^\bu\in\sC$ and $\G^\bu\in\sF$.
 Then both $\cP^\bu$ and $\cQ^\bu$ are complexes of contraadjusted
sheaves, hence (by Corollary~\ref{quasi-cta-characterizations}(c))
so is the complex $\G^\bu$.
 On the other hand, $\G^\bu$ is a acyclic complex (with flat sheaves
of cocycles, but we do not need to use this condition).
 Since any quasi-coherent sheaf on $X$ has finite coresolution dimension
with respect to $X\qcoh^\cta$ (by Lemma~\ref{dil-cta-clp-finite-dim}(b)),
it follows by virtue of the dual version of
Corollary~\ref{fdim-acyclic-cor} that all the sheaves of cocycles of
the complex $\G^\bu$ are also contraadjusted.
 Now Lemma~\ref{F-tilde-orthogonal-to-C-tilde} tells us that
$\G^\bu\in\sF^{\perp_1}$.
 Finally, acyclicity of the complexes $\Hom_X(\F^\bu,\cQ^\bu)$ and
$\Hom_X(\F^\bu,\G^\bu)$ for all $\F^\bu\in\sF$ implies acyclicity
of the complex $\Hom_X(\F^\bu,\cP^\bu)$, and it remains to
recall Lemma~\ref{Ext-1-as-homotopy-Hom} in order to conclude
that $\cP^\bu\in\sF^{\perp_1}=\sC$.

 Part~(b) was first stated in~\cite[Theorem~3.3]{CET} (but there was
an error in the proof, fixed in a subsequent paper;
see~\cite[Remark~10.3]{PS6} for a discussion).
 The proof is similar to part~(a), except that instead of the references
to Lemma~\ref{dil-cta-clp-finite-dim}(b) and
Corollary~\ref{fdim-acyclic-cor} one needs to use the cotorsion
periodicity theorem for quasi-coherent sheaves.
 This theorem claims that, in any acyclic complex of cotorsion
quasi-coherent sheaves on a quasi-compact semi-separated scheme,
the sheaves of cocycles are also cotorsion (\cite[Corollary~10.4]{PS6}
or~\cite[Theorem~8.3]{Pphil}).
 One also has to use the fact that the class of all flat $R$\+modules
is deconstructible (see Examples~\ref{deconstructible-examples}).
\end{proof}

\begin{cor} \label{complexes-of-quasi-cta-cot-antilocal}
 Let $X=\bigcup_\alpha U_\alpha$ be a finite affine open covering
of a quasi-compact semi-separated scheme~$X$.  Then \par
\textup{(a)} any complex of contraadjusted quasi-coherent sheaves
on $X$ is a direct summand of a finitely iterated extension of
the direct images of complexes of contraadjusted quasi-coherent
sheaves from~$U_\alpha$; \par
\textup{(b)} any complex of cotorsion quasi-coherent sheaves on $X$
is a direct summand of a finitely iterated extension of the direct
images of complexes of cotorsion quasi-coherent sheaves from~$U_\alpha$.
\end{cor}

\begin{proof}
 Part~(a): in the notation and terminology of
Sections~\ref{local-classes-subsect}\+-%
\ref{gluing-cotorsion-in-qcoh-subsect}, let $\sR$ be the local class of
all commutative rings $R$ and $\sE_R=\Com(R\modl)$ be the abelian
category of complexes of $R$\+modules.
 Let $\sF_R=\Acycl(R\modl_\vfl)\sub\sE_R$ be the class of all acyclic
complexes of very flat $R$\+modules with very flat $R$\+modules of
cocycles, and let $\sC(R)=\Com(R\modl^\cta)$ be the class of all
complexes of contraadjusted $R$\+modules.
 Then the classes $\sE$ and $\sF$ are very local by
Examples~\ref{local-classes-examples} and~\cite[Example~6.11]{Pal}.
 The pair of classes $(\sF_R,\sC(R))$ is a hereditary complete
cotorsion pair in $\sE_R$ by
Lemma~\ref{acyclic-of-vfl-flat-arbitrary-of-cta-cot-pairs}(a)
(for the affine scheme $\Spec R$).

 Applying Theorem~\ref{quasi-coherent-gluing-theorem} to the datum
of the classes $\sR$, \,$\sE_R$, \,$\sF_R$, and $\sC(R)$, we obtain
a cotorsion pair $(\sF_X,\sC(X))$ in the abelian category
$\sE_X=\Com(X\qcoh)$.
 Here $\sF_X$ is the class of all locally\+$\sF$ complexes, while
$\sC(X)$ is the class of all direct summands of finitely iterated
extensions of direct images of complexes from~$\sC(U_\alpha)$
in~$\sE_X$.
 Now it is easy to see that $\sF_X=\Acycl(X\qcoh_\vfl)$,
and the comparison with
Lemma~\ref{acyclic-of-vfl-flat-arbitrary-of-cta-cot-pairs}(a)
(for the scheme~$X$) shows that $\Com(X\qcoh^\cta)=\sC(X)$, as desired.

 The proof of part~(b) is similar and uses
Lemma~\ref{acyclic-of-vfl-flat-arbitrary-of-cta-cot-pairs}(b)
together with~\cite[Example~6.12]{Pal}.
\end{proof}

\begin{cor} \label{complexes-of-quasi-vfl-flat-cta-cot-antilocal}
 Let $X=\bigcup_\alpha U_\alpha$ be a finite affine open covering
of a quasi-compact semi-separated scheme~$X$.  Then \par
\textup{(a)} any complex of very flat contraadjusted quasi-coherent
sheaves on $X$ is a direct summand of a finitely iterated extension of
the direct images of complexes of very flat contraadjusted
quasi-coherent sheaves from~$U_\alpha$; \par
\textup{(b)} any complex of flat cotorsion quasi-coherent sheaves on $X$
is a direct summand of a finitely iterated extension of the direct
images of complexes of flat cotorsion quasi-coherent sheaves
from~$U_\alpha$; \par
\textup{(c)} any complex of flat contraadjusted quasi-coherent
sheaves on $X$ is a direct summand of a finitely iterated extension of
the direct images of complexes of flat contraadjusted quasi-coherent
sheaves from~$U_\alpha$.
\end{cor}

\begin{proof}
 First of all, consider the hereditary complete cotorsion pair
$(\sF,\sC)$ from
Lemma~\ref{acyclic-of-vfl-flat-arbitrary-of-cta-cot-pairs}(a)
in the abelian category $\sK=\Com(X\qcoh)$, and consider the exact
subcategory of complexes of very flat quasi-coherent sheaves
$\sE=\Com(X\qcoh_\vfl)\sub\sK$.
 By Lemmas~\ref{restricting-hereditary-cotorsion}
and~\ref{restricting-cotorsion-pairs-lemma}(a), the cotorsion pair
$(\sF,\sC)$ restricts to a hereditary complete cotorsion pair in~$\sE$.
 So the pair of classes $\sF=\Acycl(X\qcoh_\vfl)$ and
$\sE\cap\sC=\Com(X\qcoh_\vfl^\cta)$ is a hereditary complete cotorsion
pair in~$\sE$.
 In the rest of this proof, we redenote $\sE\cap\sC$ by~$\sC$.

 Let $\sR$ be the local class of all commutative rings $R$ and
$\sE_R=\Com(R\modl_\vfl)\sub\sK_R=\Com(R\modl)$ be the exact category
of complexes of very flat $R$\+modules.
 Put $\sF_R=\Acycl(R\modl_\vfl)$ and $\sC(R)=\Com(R\modl_\vfl^\cta)$,
where the notation is $R\modl_\vfl^\cta=R\modl_\vfl\cap R\modl^\cta$,
as usual.
 Then the classes $\sE$ and $\sF$ are very local by
Examples~\ref{local-classes-examples} and~\cite[Example~6.11]{Pal}.
 The pair of classes $(\sF_R,\sC(R))$ is a hereditary complete
cotorsion pair in $\sE_R$, as we have just seen.

 Applying Theorem~\ref{quasi-coherent-gluing-theorem} to the datum
of the classes $\sR$, \,$\sE_R$, $\sF_R$, and $\sC(R)$, we obtain
a cotorsion pair $(\sF_X,\sC(X))$ in the abelian category
$\sE_X$ of locally\+$\sE$ complexes of quasi-coherent sheaves on~$X$.
 Clearly, $\sE_X=\Com(X\qcoh_\vfl)$ and $\sF_X=\Acycl(X\qcoh_\vfl)$.
 The theorem tells us that $\sC(X)$ is the class of all direct summands
of finitely iterated extensions of direct images of complexes of
very flat contraadjusted quasi-coherent sheaves from~$U_\alpha$.
 Comparing the two cotorsion pairs in $\Com(X\qcoh_\vfl)$ that we
have constructed, we conclude that $\Com(X\qcoh_\vfl^\cta)=\sC(X)$,
as desired.

 The proof of part~(b) is similar; one needs to restrict the cotorsion
pair from Lemma~\ref{acyclic-of-vfl-flat-arbitrary-of-cta-cot-pairs}(b)
to the exact subcategory $\Com(X\qcoh_\fl)\sub\Com(X\qcoh)$.
 The proof of part~(c) is also similar; one needs to restrict
the cotorsion pair from
Lemma~\ref{acyclic-of-vfl-flat-arbitrary-of-cta-cot-pairs}(a)
to the exact subcategory $\Com(X\qcoh_\fl)\sub\Com(X\qcoh)$.
\end{proof}

\begin{cor} \label{complexes-of-quasi-flat-cta-cocycles-antilocal}
 Let $X=\bigcup_\alpha U_\alpha$ be a finite affine open covering
of a quasi-compact semi-separated scheme~$X$.  Then any acyclic
complex of flat contraadjusted quasi-coherent sheaves on $X$ with flat
(and contraadjusted) sheaves of cocycles is a direct summand of
a finitely iterated extension of the direct images from $U_\alpha$
of acyclic complexes of flat contraadjusted quasi-coherent sheaves
with flat sheaves of cocycles.
\end{cor}

\begin{proof}
 Any acyclic complex of contraadjusted quasi-coherent sheaves has
contraadjusted quasi-coherent sheaves of cocycles
by Lemma~\ref{dil-cta-clp-finite-dim}(b) and the dual version of
Corollary~\ref{fdim-acyclic-cor}.
 To prove the corollary, consider the hereditary complete cotorsion
pair $(\sF,\sC)$ from
Lemma~\ref{acyclic-of-vfl-flat-arbitrary-of-cta-cot-pairs}(a) in
the abelian category $\sK=\Com(X\qcoh)$, and consider the exact
subcategory of acyclic complexes of flat quasi-coherent sheaves with
flat sheaves of cocycles $\sE=\Acycl(X\qcoh_\fl)\sub\sK$.
 By Lemmas~\ref{restricting-hereditary-cotorsion}
and~\ref{restricting-cotorsion-pairs-lemma}(a), the cotorsion pair
$(\sF,\sC)$ restricts to a hereditary complete cotorsion pair in~$\sE$.
 So the pair of classes $\sF=\Acycl(X\qcoh_\vfl)$ and
$\sE\cap\sC=\Acycl(X\qcoh_\fl)\cap\Com(X\qcoh^\cta)=
\Acycl(X\qcoh_\fl^\cta)$ is a hereditary complete cotorsion pair
in~$\sE$.
 The rest of the argument is similar to the proof of
Corollary~\ref{complexes-of-quasi-vfl-flat-cta-cot-antilocal}.
\end{proof}

 The definition of the full subcategory of homotopy flat complexes
of flat quasi-coherent sheaves $\sD(X\qcoh_\fl)_\fl\sub\sD(X\qcoh_\fl)$
in the derived category $\sD(X\qcoh_\fl)$ was given in
Section~\ref{homotopy-lin-subsect}.
 Similarly, we define the full subcategory of \emph{homotopy very flat
complexes} of very flat quasi-coherent sheaves
$\sD(X\qcoh_\vfl)_\vfl\sub\sD(X\qcoh_\vfl)$ as the minimal triangulated
subcategory in the derived category $\sD(X\qcoh_\vfl)$ containing
the objects of $X\qcoh_\vfl$ and closed under infinite direct sums
(cf.\ Section~\ref{homotopy-adjusted}).
 The definition of homotopy very flat complexes of very flat modules
over a commutative ring is similar.

\begin{lem} \label{homotopy-very-flatness-local}
 Let $X$ be a quasi-compact semi-separated scheme.
 Then a complex of very flat quasi-coherent sheaves $\F^\bu$ on $X$
belongs to $\sD(X\qcoh_\vfl)_\vfl$ if and only if the complex of
very flat\/ $\O_X(U)$\+modules\/ $\F^\bu(U)$ is homotopy very flat
for every affine open subscheme $U\sub X$.
 It suffices to check this condition for affine open subschemes
forming any chosen affine open covering of~$X$.
\end{lem}

\begin{proof}
 This is dual-analogous to the proof of
Lemma~\ref{homotopy-local-injectivity-characterizations}%
\,(1)\,$\Leftrightarrow$\,(3).
 The ``only if'' assertion is obvious.
 To prove the ``if'', choose a finite affine open covering
$X=\bigcup_{\alpha=1}^N U_\alpha$ of the scheme $X$, and consider
the \v Cech coresolution~\eqref{cech-quasi}
\begin{multline} \label{cech-coresolution-of-vfl-complex} \textstyle
 0\lrarrow\F^\bu\lrarrow\bigoplus_\alpha k_\alpha{}_*k_\alpha^*\F^\bu
 \lrarrow\bigoplus_{\alpha<\beta}k_{\alpha,\beta}
 {}_*k_{\alpha,\beta}^*\F^\bu \\ \lrarrow\dotsb\lrarrow
 k_{1,\dotsc,N}{}_*k_{1,\dotsc,N}^*\F^\bu\lrarrow 0.
\end{multline}
 Notice that all the terms of~\eqref{cech-coresolution-of-vfl-complex},
except perhaps the leftmost one, are homotopy very flat complexes
(indeed, the class of homotopy very flat complexes of very flat
quasi-coherent sheaves is obviously preserved by inverse images with
respect to all morphisms of schemes and direct images with respect to
very flat affine morphisms).
 By the definition, it follows that the total complex of the truncated
coresolution~\eqref{cech-coresolution-of-vfl-complex} (i.~e., of
\eqref{cech-coresolution-of-vfl-complex} without the leftmost term)
is homotopy very flat.
 But this total complex is isomorphic to $\F^\bu$ in $\sD(X\qcoh_\vfl)$,
and even in $\sD^\abs(X\qcoh_\vfl)$.
\end{proof}

 Denote by $\Com(X\qcoh_\vfl)_\vfl\sub\Com(X\qcoh_\vfl)$ the class of
homotopy very flat complexes of very flat quasi-coherent sheaves on $X$
viewed as objects of the abelian category of complexes $\Com(X\qcoh)$.
 Similarly, let $\Com(X\qcoh_\fl)_\fl\sub\Com(X\qcoh_\fl)$ denote
the class of homotopy flat complexes of flat quasi-coherent sheaves
on $X$ viewed as objects of the abelian category $\Com(X\qcoh)$.
 Given a commutative ring $R$, the notation $\Com(R\modl_\vfl)_\vfl
\sub\Com(R\modl_\vfl)$ and $\Com(R\modl_\fl)_\fl\sub\Com(R\modl_\fl)$
is understood similarly.

 In the following lemma and corollary, as before in this section, it is
helpful to keep in mind that all acyclic complexes of contraadjusted
quasi-coherent sheaves have contraadjusted sheaves of cocycles
by Lemma~\ref{dil-cta-clp-finite-dim}(b) and the dual version of
Corollary~\ref{fdim-acyclic-cor}; and all acyclic complexes of cotorsion
quasi-coherent sheaves have cotorsion sheaves of cocycles
by~\cite[Corollary~10.4]{PS6} or~\cite[Theorem~8.3]{Pphil}.

\begin{lem} \label{homotopy-vfl-flat-acyclic-of-cta-cot-pairs}
\textup{(a)} The pair of classes of homotopy very flat complexes of
very flat quasi-coherent sheaves\/ $\sF'=\Com(X\qcoh_\vfl)_\vfl$ and
acyclic complexes of contraadjusted quasi-coherent sheaves (with
contraadjusted sheaves of cocycles)\/ $\sC'=\Acycl(X\qcoh^\cta)$ is
a hereditary complete cotorsion pair in the abelian category\/
$\Com(X\qcoh)$. {\hbadness=1750\par}
\textup{(b)} The pair of classes of homotopy flat complexes of flat
quasi-coherent sheaves\/ $\sF'=\Com(X\qcoh_\fl)_\fl$ and acyclic
complexes of cotorsion quasi-coherent sheaves (with cotorsion sheaves
of cocycles)\/ $\sC'=\Acycl(X\qcoh^\cot)$ is a hereditary complete
cotorsion pair in the abelian category\/ $\Com(X\qcoh)$.
\end{lem}

\begin{proof}
 Both parts~(a) and~(b) are special cases of
Proposition~\ref{dg-F-acycl-C-cotorsion-pair} applied to
the cotorsion pairs $(\sF,\sC)=(X\qcoh_\vfl\;X\qcoh^\cta)$ and
$(\sF,\sC)=(X\qcoh_\fl\;X\qcoh^\cot)$, respectively, in
the Grothendieck abelian category $\sA=X\qcoh$.
 In particular, in part~(a) the pair of classes
$(\sF,\sC)=(X\qcoh_\vfl\;X\qcoh^\cta)$ is a complete cotorsion
pair in $X\qcoh$ by Corollary~\ref{quasi-very-cta-cor}(a\+b) with
Lemma~\ref{cotorsion-pair-direct-summands-lemma}; this cotorsion
pair is hereditary since the class $X\qcoh_\vfl$ is closed under
kernels of epimorphisms in $X\qcoh$ (or by
Corollary~\ref{quasi-cta-characterizations}(c)); and the class
$X\qcoh_\vfl$ is deconstructible in $X\qcoh$, as explained in
the beginning of the proof of
Lemma~\ref{acyclic-of-vfl-flat-arbitrary-of-cta-cot-pairs}.
 The proof of part~(b) is similar, and based on
Corollaries~\ref{quasi-cotors-cor}(a\+b)
and~\ref{quasi-cotors-characterizations}(c)
(see also~\cite[Theorem~6.7]{Gil1}).
\end{proof}

\begin{cor} \label{acyclic-complexes-of-quasi-cta-cot-antilocal}
 Let $X=\bigcup_\alpha U_\alpha$ be a finite affine open covering
of a quasi-compact semi-separated scheme~$X$.  Then \par
\textup{(a)} any acyclic complex of contraadjusted quasi-coherent
sheaves on $X$ is a direct summand of a finitely iterated extension of
the direct images of acyclic complexes of contraadjusted quasi-coherent
sheaves from~$U_\alpha$; \par
\textup{(b)} any acyclic complex of cotorsion quasi-coherent sheaves on
$X$ is a direct summand of a finitely iterated extension of the direct
images of acyclic complexes of cotorsion quasi-coherent sheaves
from~$U_\alpha$.
\end{cor}

\begin{proof}
 This is similar to
Corollary~\ref{complexes-of-quasi-cta-cot-antilocal},
except that Lemma~\ref{homotopy-vfl-flat-acyclic-of-cta-cot-pairs}
is used in place of
Lemma~\ref{acyclic-of-vfl-flat-arbitrary-of-cta-cot-pairs}.
 Part~(a): let $\sR$ be the local class of all commutative rings $R$
and $\sE_R=\Com(R\modl)$ be the abelian category of complexes of
$R$\+modules.
 Let $\sF_R=\Com(R\modl_\vfl)_\vfl$ be the class of all homotopy
very flat complexes of very flat $R$\+modules, and let $\sC(R)=
\Acycl(R\modl^\cta)$ be the class of acyclic complexes of
contraadjusted $R$\+modules.
 Then the class $\sE$ is obviously very local, while the class $\sF$
is very local by Lemma~\ref{homotopy-very-flatness-local} and its proof.
 The pair of classes $(\sF_R,\sC(R))$ is a hereditary complete
cotorsion pair in $\sE_R$ by
Lemma~\ref{homotopy-vfl-flat-acyclic-of-cta-cot-pairs}
(for the affine scheme $\Spec R$).

 Applying Theorem~\ref{quasi-coherent-gluing-theorem} to the datum of
the classes $\sR$, \,$\sE_R$, \,$\sF_R$, and $\sC(R)$, we obtain
a cotorsion pair $(\sF_X,\sC(X))$ in the abelian category
$\sE_X=\Com(X\qcoh)$.
 Now Examples~\ref{local-classes-examples} and
Lemma~\ref{homotopy-very-flatness-local} tell us that the class
of all locally\+$\sF$ complexes $\sF_X$ coincides with
$\Com(X\qcoh_\vfl)_\vfl$.
 Finally, the comparison with
Lemma~\ref{homotopy-vfl-flat-acyclic-of-cta-cot-pairs}(a)
(for the scheme~$X$) shows that the class $\Acycl(X\qcoh^\cta)$
coindices with the class $\sC(X)$ of all direct summands of finitely
iterated extensions of direct images of complexes from $\sC(U_\alpha)$,
as desired.

 The proof of part~(b) is similar and uses
Lemma~\ref{homotopy-vfl-flat-acyclic-of-cta-cot-pairs}(b)
together with Theorem~\ref{homotopy-flat-thm}(c).
\end{proof}

 In the following corollary we consider acyclic complexes of flat or
very flat quasi-coherent sheaves \emph{with not necesarily flat or
very flat sheaves of cocycles}.

\begin{cor} \label{complexes-of-quasi-vfl-flat-cta-cot-acycl-antilocal}
 Let $X=\bigcup_\alpha U_\alpha$ be a finite affine open covering
of a quasi-compact semi-separated scheme~$X$.  Then \par
\textup{(a)} any acyclic complex of very flat contraadjusted
quasi-coherent sheaves on $X$ is a direct summand of a finitely
iterated extension of the direct images of acyclic complexes of very
flat contraadjusted quasi-coherent sheaves from~$U_\alpha$; \par
\textup{(b)} any acyclic complex of flat cotorsion quasi-coherent
sheaves on $X$ is a direct summand of a finitely iterated extension of
the direct images of acyclic complexes of flat cotorsion quasi-coherent
sheaves from~$U_\alpha$; \par
\textup{(c)} any acyclic complex of flat contraadjusted quasi-coherent
sheaves on $X$ is a direct summand of a finitely iterated extension of
the direct images of acyclic complexes of flat contraadjusted
quasi-coherent sheaves from~$U_\alpha$.
\end{cor}

\begin{proof}
 Part~(a): the argument is based on
Lemma~\ref{acyclic-of-vfl-flat-arbitrary-of-cta-cot-pairs}(a).
 It is similar to the proof of
Corollary~\ref{complexes-of-quasi-vfl-flat-cta-cot-antilocal}(a),
except that one has to consider the exact subcategory of acyclic
complexes of very flat quasi-coherent sheaves $\sE=\Acycl(X\qcoh)\cap
\Com(X\qcoh_\vfl)\sub\sK=\Com(X\qcoh)$.
 Both the acyclicity of complexes and their termwise very flatness
are very local properties for Zariski open coverings of quasi-compact
semi-separated (in particular, affine) schemes, in the obvious sense.
 Hence their conjunction (the intersection of the two classes) is
also very local.
 The proofs of parts~(b\+c) are similar to the proofs of the respective
parts of Corollary~\ref{complexes-of-quasi-vfl-flat-cta-cot-antilocal},
with the similar replacement of $\sE=\Com(X\qcoh_\fl)$ by
$\sE=\Acycl(X\qcoh)\cap\Com(X\qcoh_\fl)$.

 Alternatively, the proofs of all the three parts~(a\+c) can be based
on Lemma~\ref{homotopy-vfl-flat-acyclic-of-cta-cot-pairs}.
 To prove part~(a), one needs to restrict the cotorsion pair from
Lemma~\ref{homotopy-vfl-flat-acyclic-of-cta-cot-pairs}(a) to the exact
subcategory $\sE=\Com(X\qcoh_\vfl)\sub\sK=\Com(X\qcoh)$ of all
complexes of very flat quasi-coherent sheaves.
 To prove part~(b), one restricts the cotorsion pair from
Lemma~\ref{homotopy-vfl-flat-acyclic-of-cta-cot-pairs}(b) to the exact
subcategory $\sE=\Com(X\qcoh_\fl)\sub\sK=\Com(X\qcoh)$ of all
complexes of flat quasi-coherent sheaves.
 To establish part~(c), one restricts the cotorsion pair from
Lemma~\ref{homotopy-vfl-flat-acyclic-of-cta-cot-pairs}(a) to the exact
subcategory $\sE=\Com(X\qcoh_\fl)$ of all complexes of flat
quasi-coherent sheaves.
\end{proof}

 We will use the definition of \emph{Becker-contraacyclic complexes}
$\Acycl^\bctr(\sE)$ in an exact category $\sE$, as presented in
Section~\ref{becker-subsect}.
 Notice that the class of very flat contraadjusted quasi-coherent
sheaves, $X\qcoh_\vfl^\cta$, on a quasi-compact semi-separated scheme
$X$, is the class of all projective objects in
the exact category $X\qcoh^\cta$.
 The same class $X\qcoh_\vfl^\cta$ is the class of all projective
objects in the exact category of flat contraadjusted quasi-coherent
sheaves $X\qcoh_\fl^\cta$.
 Indeed, the complete cotorsion pair ($X\qcoh_\vfl$, $X\qcoh^\cta$) in
$X\qcoh$ restricts to a complete cotorsion pair in the full exact
subcategory $X\qcoh_\fl\sub X\qcoh$.
 So a complex in $X\qcoh_\fl^\cta$ is Becker-contraacyclic if and only
if it is Becker-contraacyclic in $X\qcoh^\cta$.
 Similarly, the class of flat cotorsion quasi-coherent sheaves,
$X\qcoh_\fl^\cot$, is the class of all projective objects in the exact
category $X\qcoh^\cot$.

 The cotorsion pairs described in the next lemma appeared
in~\cite[Proposition~3.2, Theorem~5.5, and Section~5.3]{Gil2}
and~\cite[Lemma~4.9]{Gil3}.

\begin{lem} \label{all-vfl-flat-contraacyclic-of-cta-cot-pairs}
\textup{(a)} The pair of classes of all complexes of very flat
quasi-coherent sheaves\/ $\sF''=\Com(X\qcoh_\vfl)$ and
Becker-contraacyclic complexes of contraadjusted quasi-coherent
sheaves\/ $\sC''=\Acycl^\bctr(X\qcoh^\cta)$ is a hereditary complete
cotorsion pair in the abelian category\/ $\Com(X\qcoh)$. \par
\textup{(b)} The pair of classes of all complexes of flat
quasi-coherent sheaves\/ $\sF''=\Com(X\qcoh_\fl)$ and
Becker-contraacyclic complexes of cotorsion quasi-coherent sheaves\/
$\sC''=\Acycl^\bctr(X\qcoh^\cot)$ is a hereditary complete cotorsion
pair in the abelian category\/ $\Com(X\qcoh)$. \hbadness=2600
\end{lem}

\begin{proof}
 Part~(a): as mentioned in the beginning of the proof of
Lemma~\ref{acyclic-of-vfl-flat-arbitrary-of-cta-cot-pairs},
the class $X\qcoh_\vfl$ is deconstructible in $X\qcoh$.
 By~\cite[Proposition~4.3]{Sto-hill}, it follows that the class
$\Com(X\qcoh_\vfl)$ of all complexes of very flat quasi-coherent sheaves
is deconstructible in $\Com(X\qcoh)$.
 The class $\Com(X\qcoh_\vfl)$ is also closed under direct summands
and, assuming $X$ to be quasi-compact and semi-separated, this class
is generating in $\Com(X\qcoh)$.
 Applying Corollary~\ref{eklof-trlifaj-corollary}, we conclude
that there exists a complete cotorsion pair $(\sF'',\sC'')$ in
$\Com(X\qcoh)$ with $\sF''=\Com(X\qcoh_\vfl)$.
 The class $\Com(X\qcoh_\vfl)$ is closed under kernels of epimorphisms
in $\Com(X\qcoh)$, so the cotorsion pair $(\sF'',\sC'')$ is hereditary.

 It remains to show that $\sC''$ is the class of all
Becker-contraacyclic complexes in $X\qcoh^\cta$.
 Indeed, by Lemma~\ref{G-plus-G-minus-complexes-Ext-adjunction}
(for $i=1$), we have $\sC''\sub\Com(X\qcoh^\cta)$.
 By Lemmas~\ref{restricting-hereditary-cotorsion}
and~\ref{restricting-cotorsion-pairs-lemma}(b), the cotorsion pair
$(\sF'',\sC'')$ restricts to a hereditary complete cotorsion pair
$(\sF,\sC)$ in the exact category $\Com(X\qcoh^\cta)$.
 So we have a cotorsion pair $\sF=\Com(X\qcoh_\vfl^\cta)$ and
$\sC=\sC''$ in $\Com(X\qcoh^\cta)$.
 Now it is clear from Lemma~\ref{Ext-1-as-homotopy-Hom} that $\sC''$
is the class of all Becker-contraacyclic complexes in $X\qcoh^\cta$.
 The proof of part~(b) is similar.
\end{proof}

\begin{cor} \label{quasi-cta-cot-contraacyclic-direct-images}
 Let $f\:Y\rarrow X$ be a morphism of quasi-compact semi-separated
schemes.  Then \par
\textup{(a)} the direct image functor~$f_*$ takes Becker-contraacyclic
complexes in $Y\qcoh^\cta$ to Becker-contraacyclic complexes in
$X\qcoh^\cta$; \par
\textup{(b)} the direct image functor~$f_*$ takes Becker-contraacyclic
complexes in $Y\qcoh^\cot$ to Becker-contraacyclic complexes in
$X\qcoh^\cot$.
\end{cor}

\begin{proof}
 In both the contexts of parts~(a) and~(b), the assertion follows from
the description of the class of Becker-contraacyclic complexes
as $\sC''=\sF''{}^{\perp_1}$
provided by Lemma~\ref{all-vfl-flat-contraacyclic-of-cta-cot-pairs}.
 One needs to use~\cite[Lemma~1.7(c)]{Pal}, applying it to the pair of
adjoint functors $f^*\:\Com(X\qcoh)\rarrow\Com(Y\qcoh)$ and
$f_*\:\Com(Y\qcoh)\rarrow\Com(X\qcoh)$ between the abelian categories
of complexes of quasi-coherent sheaves, for a suitable version of
$\Ext^1$\+adjunction.
 The point is that the functor~$f^*$ takes complexes of (very) flat
quasi-coherent sheaves to complexes of (very) flat quasi-coherent
sheaves, and the complexes of flat quasi-coherent sheaves are
adjusted to~$f^*$.
\end{proof}

\begin{cor} \label{contraacyclic-complexes-of-quasi-cta-cot-antilocal}
 Let $X=\bigcup_\alpha U_\alpha$ be a finite affine open covering of
a quasi-compact semi-separated scheme~$X$.  Then \par
\textup{(a)} any Becker-contraacyclic complex of contraadjusted
quasi-coherent sheaves on $X$ is a direct summand of a finitely
iterated extension of the direct images of Becker-contraacyclic
complexes of contraadjusted quasi-coherent sheaves from~$U_\alpha$; \par
\textup{(b)} any Becker-contraacyclic complex of cotorsion
quasi-coherent sheaves on $X$ is a direct summand of a finitely
iterated extension of the direct images of Becker-contraacyclic
complexes of cotorsion quasi-coherent sheaves from~$U_\alpha$.
\end{cor}

\begin{proof}
 This is similar to
Corollaries~\ref{complexes-of-quasi-cta-cot-antilocal}
and~\ref{acyclic-complexes-of-quasi-cta-cot-antilocal}, and based on
Lemma~\ref{all-vfl-flat-contraacyclic-of-cta-cot-pairs}.
 In part~(a), the point is that the classes $\sE_R=\Com(R\modl)$ and
$\sF_R=\Com(R\modl_\vfl)$ are very local for all commutative rings~$R$
(see Examples~\ref{local-classes-examples}), and moreover, the class
of all locally\+$\sF$ complexes of quasi-coherent sehaves on $X$
coincides with $\Com(X\qcoh_\vfl)$.
 Comparing the cotorsion pair in $\Com(X\qcoh)$ described in
Lemma~\ref{all-vfl-flat-contraacyclic-of-cta-cot-pairs}(a) with the one
provided by Theorem~\ref{quasi-coherent-gluing-theorem}, one comes to
the assertion of the corollary.
 Part~(b) is similar.
\end{proof}

\begin{cor} \label{complexes-of-quasi-flat-cta-contraacyclic-antilocal}
 Let $X=\bigcup_\alpha U_\alpha$ be a finite affine open covering of
a quasi-compact semi-separated scheme~$X$.  Then any
Becker-contraacyclic complex of flat contraadjusted quasi-coherent
sheaves on $X$ is a direct summand of a finitely iterated extension of
the direct images of Becker-contraacyclic complexes of flat
contraadjusted quasi-coherent sheaves from~$U_\alpha$.
\end{cor}

\begin{proof}
 The argument is similar to the proofs of
Corollaries~\ref{complexes-of-quasi-vfl-flat-cta-cot-antilocal}(c),
\ref{complexes-of-quasi-flat-cta-cocycles-antilocal},
and~\ref{complexes-of-quasi-vfl-flat-cta-cot-acycl-antilocal}(c).
 Consider the hereditary complete cotorsion pair $(\sF'',\sC'')$ from
Lemma~\ref{all-vfl-flat-contraacyclic-of-cta-cot-pairs}(a) in
the abelian category $\sK=\Com(X\qcoh)$, and consider the exact
subcategory of complexes of flat quasi-coherent sheaves
$\sE=\Com(X\qcoh_\fl)\sub\sK$.
 By Lemmas~\ref{restricting-hereditary-cotorsion}
and~\ref{restricting-cotorsion-pairs-lemma}(a), the cotorsion
pair $(\sF'',\sC'')$ restricts to a hereditary complete cotorsion
pair in~$\sE$.
 So the pair of classes $\sF=\sF''=\Com(X\qcoh_\vfl)$ and
$\sC=\sE\cap\sC''=\Acycl^\bctr(X\qcoh_\fl^\cta)$ is
a hereditary complete cotorsion pair in~$\sE$.

 Let $\sR$ be the local class of all commutative rings $R$ and
$\sE_R=\Com(R\modl_\fl)\sub\sK_R=\Com(R\modl)$ be the exact 
category of complexes of flat $R$\+modules.
 Put $\sF_R=\Com(R\modl_\vfl)$ and
$\sC(R)=\Acycl^\bctr(R\modl_\fl^\cta)$.
 Then the classes $\sE$ and $\sF$ are very local by
Examples~\ref{local-classes-examples}.
 The pair of classes $(\sF_R,\sC(R))$ is a hereditary complete
cotorsion pair in $\sE_R$, as we have just seen.

 Applying Theorem~\ref{quasi-coherent-gluing-theorem} to the datum
of the classes $\sR$, \,$\sE_R$, \,$\sF_R$, and $\sC(R)$, we obtain
a cotorsion pair $(\sF_X,\sC(X))$ in the abelian category
$\sE_X=\Com(X\qcoh_\fl)$, with $\sF_X=\Com(X\qcoh_\vfl)$.
 The theorem tells us that $\sC(X)$ is the class of all direct summands
of finitely iterated extensions of direct images of
Becker-contraacyclic complexes in the categories of flat
contraadjusted quasi-coherent sheaves on~$U_\alpha$.
 Comparing the two cotorsion pairs in $\Com(X\qcoh_\fl)$ that we
have constructed, we conclude that $\Acycl^\bctr(X\qcoh_\fl^\cta)
=\sC(X)$, as desired.
\end{proof}

\begin{cor} \label{acyclic=bctraacyclic-in-flat-cta}
 A complex of flat contraadjusted quasi-coherent sheaves on
a quasi-compact semi-separated scheme $X$ is Becker-contraacyclic
in $X\qcoh_\fl^\cta$ if and only if it is acyclic with flat (and
contraadjusted) sheaves of cocycles.
 In other words,
$$
 \Acycl^\bctr(X\qcoh_\fl^\cta)=\Acycl(X\qcoh_\fl^\cta).
$$
\end{cor}

\begin{proof}
 Corollaries~\ref{complexes-of-quasi-flat-cta-cocycles-antilocal}
and~\ref{complexes-of-quasi-flat-cta-contraacyclic-antilocal} are
stated as implications in one (nontrivial) direction, but the converse
implications can be easily seen to hold (and the proofs of the two
corollaries actually prove the ``if and only if'').
 So comparing the two corollaries reduces the question to the case
of an affine scheme $U=\Spec R$.
 In this case, Proposition~\ref{becker-contraacyclicity-for-cta-cot}(a)
tells us that a complex of contraadjusted modules is
Becker-contraacyclic in $R\modl^\cta$ if and only if it is
Becker-contraacyclic in $R\modl$.
 Thus the question reduces further to showing that a complex of flat
modules is Becker-contraacyclic in $R\modl$ if and only if it is
acyclic with flat modules of cocycles.
 This is Theorem~\ref{flat-projective-periodicity}(b) together with
Proposition~\ref{flat-projective-periodicity-complements}(a).
\end{proof}

 Our discussion of the result of
Corollary~\ref{acyclic=bctraacyclic-in-flat-cta} will continue
in Section~\ref{Becker-contraderived-subsect}; see
Corollary~\ref{becker-co-contraderived-of-fl-cta}.

\subsection{Becker's coderived category}
\label{Becker-coderived-subsect}
 We refer to Section~\ref{becker-subsect} for the definitions of
\emph{Becker-coacyclic complexes} and
\emph{Becker's coderived category}.

 Let us start with several elementary observations collected in
the following corollary.
 Notice first of all that, for any scheme $X$, the category $X\qcoh$
of quasi-coherent sheaves on $X$ is a Grothendieck abelian category.
 In particular, the assertion that $X\qcoh$ has a set of generators
is due to Gabber~\cite[Proposition Tag~077P]{SP}.

\begin{cor} \label{homotopy-derived-coderived-cor}
 For any scheme $X$, the triangulated functor\/ $\Hot^+(X\qcoh^\inj)
\rarrow\sD^+(X\qcoh)$ induced by the embedding of additive/abelian
categories $X\qcoh^\inj\rarrow X\qcoh$ is an equivalence of triangulated
categories.
 For any symbol\/ $\bst=\b$, $\abs+$, $\abs-$, $\bco$, $\co$, or $\abs$,
the triangulated functor\/ $\Hot^\st(X\qcoh^\inj)\rarrow\sD^\st(X\qcoh)$
is fully faithful.
\end{cor}

\begin{proof}
 The first assertion is classical; it goes back, at least,
to~\cite[Propositions~I.4.7 and~II.1.1]{Har} and follows from the fact
that there are enough injective objects in $X\qcoh$.
 It is a particular case of the dual version of
Proposition~\ref{infinite-resolutions}(a).
 The second assertion, suitably stated, holds for any exact category
(or for any exact category with exact functors of infinite direct sum,
as appropriate) by Lemma~\ref{homotopy-inj-proj-fully-faithful}(a).
 In the case of the Becker coderived category, it follows immediately
from the definitions.
\end{proof}

 The following theorem is a stronger version of
Corollary~\ref{homotopy-derived-coderived-cor} for $\bst=\bco$.

\begin{thm} \label{quasi-coherent-becker-coderived}
 For any scheme $X$, the pair of classes of Becker-coacyclic
complexes of quasi-coherent sheaves\/ $\sF=\Acycl^\bco(X\qcoh)$
and all complexes of injective quasi-coherent sheaves\/
$\sC=\Com(X\qcoh^\inj)$ is a hereditary complete cotorsion pair
in the abelian category\/ $\Com(X\qcoh)$.
 Consequently, the composition of triangulated functors
$$
 \Hot(X\qcoh^\inj)\rarrow\Hot(X\qcoh)\rarrow\sD^\bco(X\qcoh)
$$
is a triangulated equivalence\/ $\Hot(X\qcoh^\inj)\simeq
\sD^\bco(X\qcoh)$.
\end{thm}

\begin{proof}
 This is a special case of
Theorem~\ref{coderived-of-grothendieck-contraderived-of-lpacepo}(a)
and its proof (cf.\ the proof of
Corollary~\ref{becker-coderived-of-cta-cot} below).
\end{proof}

 The following lemma should be compared with~\cite[Lemma~A.7]{Psemten}.

\begin{lem} \label{qcoh-flat-inverse-image-becker-coacyclic}
 For any flat quasi-compact quasi-separated morphism of schemes
$f\:Y\rarrow X$, the inverse image functor $f^*\:X\qcoh\rarrow
Y\qcoh$ takes Becker-coacyclic complexes of quasi-coherent sheaves
on $X$ to Becker-coacyclic complexes of quasi-coherent sheaves on~$Y$.
\end{lem}

\begin{proof}
 This is a particular case of
Lemma~\ref{exact-with-adjoint-preservation-lemma}(a), as
the functor~$f^*$ is exact and has a right adjoint functor
$f_*\:Y\qcoh\rarrow X\qcoh$.
 In fact, one can drop the quasi-compactness and quasi-separatedness
conditions on~$f$ by noticing that the functor~$f^*$ always preserves
infinite direct sums and is exact whenever $f$~is flat; then it remains
to refer to~\cite[Lemma~A.5]{Psemten}.
 The point is that the suitable ``quasi-coherent direct image
functor'' (the direct image of sheaves of $\O$\+modules restricted
to quasi-coherent sheaves on $Y$ and composed with
the coherator~\cite[Section~B.12]{TT}, \cite[Proposition Tag~077P]{SP}
on~$X$) is right adjoint to~$f^*$ for any morphism of schemes~$f$.
\end{proof}

 The next lemma should be compared with~\cite[Lemma~A.8]{Psemten}.

\begin{lem} \label{qcoh-affine-direct-image-becker-coacyclic}
 For any affine morphism of schemes $f\:Y\rarrow X$, the direct image
functor $f_*\:Y\qcoh\rarrow X\qcoh$ takes Becker-coacyclic complexes
of quasi-coherent sheaves on $Y$ to Becker-coacyclic complexes of
quasi-coherent sheaves on~$X$.
\end{lem}

\begin{proof}
 This is also a particular case of
Lemma~\ref{exact-with-adjoint-preservation-lemma}(a)
or~\cite[Lemma~A.5]{Psemten}.
 The point is that the direct image functor~$f_*$ is exact and
preserves infinite direct sums, so it has a right adjoint by
the Special Adjoint Functor Theorem.
\end{proof}

\begin{cor} \label{qcoh-becker-coacyclicity-is-local}
 Let $X=\bigcup_\alpha U_\alpha$ be a finite affine open covering of
a quasi-compact semi-separated scheme.
 Then a complex of quasi-coherent sheaves\/ $\cA^\bu$ on $X$ is
Becker-coacyclic if and only if, for every index~$\alpha$, the complex
of\/ $\O(U_\alpha)$\+modules $\cA^\bu(U_\alpha)$ is Becker-coacyclic.
\end{cor}

\begin{proof}
 This is~\cite[Proposition~A.10]{Psemten}.
 The argument is similar to the one in~\cite[Remark~1.3]{EP} and
uses the \v Cech coresolution together with
Lemmas~\ref{qcoh-flat-inverse-image-becker-coacyclic},
\ref{qcoh-affine-direct-image-becker-coacyclic},
and~\ref{Positselski-trivial-are-Becker-trivial}(a).
\end{proof}

\begin{cor} \label{complexes-of-quasi-injective-antilocal}
 Let $X=\bigcup_\alpha U_\alpha$ be a finite affine open covering
of a quasi-compact semi-separated scheme~$X$.
 Then any complex of injective quasi-coherent sheaves on $X$ is
a direct summand of a finitely iterated extension of the direct images
of complexes of injective quasi-coherent sheaves from~$U_\alpha$.
\end{cor}

\begin{proof}
 The proof is similar to that of
Corollary~\ref{complexes-of-quasi-cta-cot-antilocal} and based on
Theorem~\ref{quasi-coherent-becker-coderived}.
 Let $\sR$ be the local class of all commutative rings $R$ and
$\sE_R=\Com(R\modl)$ be the abelian category of complexes of
$R$\+modules.
 Let $\sF_R=\Acycl^\bco(R\modl)$ be the class of all Becker-coacyclic
complexes of $R$\+modules, and let $\sC(R)=\Com(R\modl^\inj)$ be
the class of all complexes of injective $R$\+modules.
 Then the class $\sE$ is obviously very local, while the class $\sF$
is very local by Lemma~\ref{qcoh-affine-direct-image-becker-coacyclic}
and Corollary~\ref{qcoh-becker-coacyclicity-is-local} (applied to
affine schemes).
 The pair of classes $(\sF_R,\sC(R))$ is a hereditary complete cotorsion
pair in $\sE_R$ by Theorem~\ref{quasi-coherent-becker-coderived}
(for the affine scheme $\Spec R$).

 Applying Theorem~\ref{quasi-coherent-gluing-theorem} to the datum of
the classes $\sR$, \,$\sE_R$, \,$\sF_R$, and $\sC(R)$, we obtain
a cotorsion pair $(\sF_X,\sC(X))$ in the abelian category
$\sE_X=\Com(X\qcoh)$.
 Here $\sF_X$ is the class of all locally\+$\sF$ complexes, so
Corollary~\ref{qcoh-becker-coacyclicity-is-local} (for the scheme~$X$)
tells us that $\sF_X=\Acycl^\bco(X\qcoh)$.
 By Theorem~\ref{quasi-coherent-becker-coderived}, it follows that
$\sC(X)=\Com(X\qcoh^\inj)$.
 On the other hand, Theorem~\ref{quasi-coherent-gluing-theorem} provides
a description of $\sC(X)$ as the class of all direct summands of
finitely iterated extensions of direct images of complexes from
$\sC(U_\alpha)$ in~$\sE_X$.
 Comparing the two descriptions of the class $\sC(X)$, we arrive to
the desired conclusion.
\end{proof}

 Assume that $X$ is a quasi-compact semi-separated scheme.
 In this case, the full subcategories $X\qcoh^\dil$, $X\qcoh^\cta$,
and $X\qcoh^\cot$ are coresolving in $X\qcoh$, in the sense of
Section~\ref{infinite-resolutions-subsect}
(by Corollaries~\ref{quasi-cta-characterizations}(c)
and~\ref{quasi-cotors-characterizations}(c), and
Section~\ref{dilute-subsect}), so it follows that the classes of
injective objects in the exact categories $X\qcoh^\dil$, $X\qcoh^\cta$,
and $X\qcoh^\cot$ coincide with the class $X\qcoh^\inj$ of injective
objects in $X\qcoh$.
 Therefore, a complex in $X\qcoh^\cta$ is Becker-coacyclic in
$X\qcoh^\cta$ if and only if it is Becker-coacyclic in $X\qcoh$;
and similarly for complexes in $X\qcoh^\dil$ or $X\qcoh^\cot$.

 There are enough injective objects in $X\qcoh^\dil$, $X\qcoh^\cta$,
and $X\qcoh^\cot$.
 Consequently, all the triangulated functors $\Hot^+(X\qcoh^\inj)
\rarrow\sD^+(X\qcoh^\cot)\rarrow\sD^+(X\qcoh^\cta)\rarrow
\sD^+(X\qcoh^\dil)\rarrow \sD^+(X\qcoh)$ are triangulated equivalences
(by the dual version of Proposition~\ref{infinite-resolutions}(a)).

\begin{cor} \label{becker-coderived-of-cta-cot}
\textup{(a)} The hereditary complete cotorsion pair
\textup{(}$\Acycl^\bco(X\qcoh)$, $\Com(X\qcoh^\inj)$\textup{)} in
the abelian category\/ $\Com(X\qcoh)$ restricts to a hereditary
complete cotorsion pair in the exact subcategory\/
$\Com(X\qcoh^\cta)\sub\Com(X\qcoh)$.
 Consequently, the composition of triangulated functors
$$
 \Hot(X\qcoh^\inj)\rarrow\Hot(X\qcoh^\cta)\rarrow\sD^\bco(X\qcoh^\cta)
$$
is a triangulated equivalence\/ $\Hot(X\qcoh^\inj)\simeq
\sD^\bco(X\qcoh^\cta)$. \par
\textup{(b)} The hereditary complete cotorsion pair
\textup{(}$\Acycl^\bco(X\qcoh)$, $\Com(X\qcoh^\inj)$\textup{)} in
the abelian category\/ $\Com(X\qcoh)$ restricts to a hereditary
complete cotorsion pair in the exact subcategory\/
$\Com(X\qcoh^\cot)\sub\Com(X\qcoh)$.
 Consequently, the composition of triangulated functors
$$
 \Hot(X\qcoh^\inj)\rarrow\Hot(X\qcoh^\cot)\rarrow\sD^\bco(X\qcoh^\cot)
$$
is a triangulated equivalence\/ $\Hot(X\qcoh^\inj)\simeq
\sD^\bco(X\qcoh^\cot)$.
\end{cor}

\begin{proof}
 In both parts~(a) and~(b), the first assertion follows from the first
assertion of Theorem~\ref{quasi-coherent-becker-coderived} in view of
Lemmas~\ref{restricting-hereditary-cotorsion}
and~\ref{restricting-cotorsion-pairs-lemma}(b).
 The second assertion is then provable similarly
to~\cite[proof of Corollary~9.5]{PS4}.

 Let us spell out part~(b).
 It is clear from the definitions that the triangulated functor
$\Hot(X\qcoh^\inj)\rarrow\sD^\bco(X\qcoh^\cot)$ is fully faithful.
 In order to prove that it is a triangulated equivalence, it suffices
to find, for any complex of cotorsion quasi-coherent sheaves
$\cP^\bu$ on~$X$, a complex of injective quasi-coherent sheaves
$\J^\bu$ together with a morphism of complexes $\cP^\bu\rarrow\J^\bu$
with a Becker-coacyclic cone.
 Indeed, let $0\rarrow\cP^\bu\rarrow\J^\bu\rarrow\cA^\bu\rarrow0$ be
a special preenvelope sequence~\eqref{special-preenvelope-sequence}
for the complex $\cP^\bu$ in the complete cotorsion pair of
the first assertion of~(b).
 So $\J^\bu$ is a complex of injective quasi-coherent sheaves on $X$,
while $\cA^\bu$ is a Becker-coacyclic complex of cotorsion
quasi-coherent sheaves.
 The total complex of the short exact sequence of complexes
$0\rarrow\cP^\bu\rarrow\J^\bu\rarrow\cA^\bu\rarrow0$ is absolutely
acyclic, hence Becker-coacyclic in $X\qcoh^\cot$ by
Lemma~\ref{Positselski-trivial-are-Becker-trivial}(a).
 Since $\cA^\bu$ is also Becker-coacyclic, it follows that the cone
of the morphism $\cP^\bu\rarrow\J^\bu$ is Becker-coacyclic in
$X\qcoh^\cot$ as well.

 Alternatively, for the second assertion of part~(a) one can refer to
Theorem~\ref{quasi-coherent-becker-coderived} together with
Corollary~\ref{qcoh-dil-cta-derived-equiv-cor}(b).
\end{proof}

\begin{cor} \label{dil-cta-cot-into-all-equiv-on-bco}
\textup{(a)} The inclusion of exact/abelian categories
$X\qcoh^\dil\rarrow X\qcoh$ induces a triangulated equivalence of
the Becker coderived categories
$$
 \sD^\bco(X\qcoh^\dil)\simeq\sD^\bco(X\qcoh).
$$ \par
\textup{(b)} The inclusion of exact/abelian categories
$X\qcoh^\cta\rarrow X\qcoh$ induces a triangulated equivalence of
the Becker coderived categories
$$
 \sD^\bco(X\qcoh^\cta)\simeq\sD^\bco(X\qcoh).
$$ \par
\textup{(c)} The inclusion of exact/abelian categories
$X\qcoh^\cot\rarrow X\qcoh$ induces a triangulated equivalence of
the Becker coderived categories
$$
 \sD^\bco(X\qcoh^\cot)\simeq\sD^\bco(X\qcoh).
$$
\end{cor}

\begin{proof}
 Once again, a complex in $X\qcoh^\dil$, $X\qcoh^\cta$, or $X\qcoh^\cot$
is Becker-coacyclic if and only if it is Becker-coacyclic in $X\qcoh$;
so the induced triangulated functor on the Becker coderived
categories is well-defined in the three cases.
 Part~(a) is a part of
Corollary~\ref{qcoh-dil-cta-derived-equiv-cor}(a).
 Part~(b) is a part of
Corollary~\ref{qcoh-dil-cta-derived-equiv-cor}(b).
 Part~(c) can be obtained by comparing
Theorem~\ref{quasi-coherent-becker-coderived} with
Corollary~\ref{becker-coderived-of-cta-cot}(b).

 Alternatively, one can prove part~(c) directly as follows:
in view of Lemma~\ref{pkoszul-lemma16}(b), it suffices to find for any
complex $\M^\bu$ in $X\qcoh$ a complex $\cP^\bu$ in $X\qcoh^\cot$
together with a morphism of complexes $\M^\bu\rarrow\cP^\bu$ with
a Becker-coacyclic cone.
 But this is a part of Theorem~\ref{quasi-coherent-becker-coderived},
which tells us that one can even find a complex $\cP^\bu=\J^\bu\in
X\qcoh^\inj$ with this property.
\end{proof}

\begin{cor} \label{becker-coacyclic-cta-cot-antilocal}
 Let $X=\bigcup_\alpha U_\alpha$ be an affine open covering
of a quasi-compact semi-separated scheme.  Then \par
\textup{(a)} any Becker-coacyclic complex in $X\qcoh^\cta$ is
a direct summand of a finitely iterated extension of direct images
of Becker-coacyclic complexes in $U_\alpha\qcoh^\cta$; \par
\textup{(b)} any Becker-coacyclic complex in $X\qcoh^\cot$ is
a direct summand of a finitely iterated extension of direct images
of Becker-coacyclic complexes in $U_\alpha\qcoh^\cot$.
\end{cor}

\begin{proof}
 The proof is similar to that of
Corollaries~\ref{complexes-of-quasi-vfl-flat-cta-cot-antilocal}(a\+b)
and~\ref{complexes-of-quasi-vfl-flat-cta-cot-acycl-antilocal}(a\+b), and
based on Lemma~\ref{acyclic-of-vfl-flat-arbitrary-of-cta-cot-pairs}.
 Part~(a): consider the hereditary complete cotorsion pair
$(\sF,\sC)$~$=$ ($\Acycl(X\qcoh_\vfl)$, $\Com(X\qcoh^\cta)$)
from Lemma~\ref{acyclic-of-vfl-flat-arbitrary-of-cta-cot-pairs}(a)
in the abelian category $\sK=\Com(X\qcoh)$, and consider the full
subcategory of Becker-coacyclic complexes of quasi-coherent sheaves
$\sE=\Acycl^\bco(X\qcoh)\subset\sK$.
 Notice that the full subcategory $\sE$ is closed under extensions in
$\sK$ by Lemma~\ref{Positselski-trivial-are-Becker-trivial}(a).
 Furthermore, for any complex $\F^\bu\in\sF$ and any complex of
contraadjusted quasi-coherent sheaves $\cP^\bu\in\sC$, the complex
$\Hom_X(\F^\bu,\cP^\bu)$ is acyclic by
Lemma~\ref{Ext-1-as-homotopy-Hom}.
 In particular, this holds for any complex of injective quasi-coherent
sheaves $\cP^\bu=\J^\bu$; hence we have $\sF\sub\sE$.
 By Lemmas~\ref{restricting-hereditary-cotorsion}
and~\ref{restricting-cotorsion-pairs-lemma}(a), the cotorsion pair
$(\sF,\sC)$ restricts to a hereditary complete cotorsion pair in~$\sE$.
 So the pair of classes $\sF=\Acycl(X\qcoh_\vfl)$ and
$\sE\cap\sC=\Acycl^\bco(X\qcoh^\cta)$ is a hereditary complete
cotorsion pair in~$\sE$.
 In the rest of this proof, we redenote $\sE\cap\sC$ by~$\sC$.

 Let $\sR$ be the local class of all commutative rings $R$ and
$\sE_R=\Acycl^\bco(R\modl)\sub\sK_R=\Com(R\modl)$ be the exact
category of Becker-coacyclic complexes of $R$\+modules.
 Put $\sF_R=\Acycl(R\modl_\vfl)$ and $\sC(R)=\Acycl^\bco(R\modl^\cta)$.
 Then the classes $\sE$ and $\sF$ are very local
by~\cite[Examples~6.7 and~6.11]{Pal}.
 The pair of classes $(\sF_R,\sC(R))$ is a hereditary complete
cotorsion pair in $\sE_R$, as we have just seen.

 Applying Theorem~\ref{quasi-coherent-gluing-theorem} to the datum of
the classes $\sR$, \,$\sE_R$, \,$\sF_R$, and $\sC(R)$, we obtain
a cotorsion pair $(\sF_X,\sC(X))$ in the abelian category $\sE_X$
of locally\+$\sE$ complexes of quasi-coherent sheaves on~$X$.
 By Corollary~\ref{qcoh-becker-coacyclicity-is-local}, we have
$\sE_X=\Acycl^\bco(X\qcoh)$.
 Clearly, $\sF_X=\Acycl(X\qcoh_\vfl)$.
 The theorem tells us that $\sC(X)$ is the class of all direct summands 
of finitely iterated extensions of direct images of Becker-coacyclic
complexes of contraadjusted quasi-coherent sheaves from~$U_\alpha$.
 Comparing the two cotorsion pairs in $\Acycl^\bco(X\qcoh)$ that we
have constructed, we conclude that $\Acycl^\bco(X\qcoh^\cta)=\sC(X)$,
as desired.

 The proof of part~(b) is similar and uses
Lemma~\ref{acyclic-of-vfl-flat-arbitrary-of-cta-cot-pairs}(b)
together with~\cite[Example~6.12]{Pal}.
 Alternatively, both the assertions of the corollary can be obtained
from the proofs of the respective assertions of
Lemma~\ref{coacyclic-of-antilocal-all-lin-cotorsion-pairs} below
(by moving the question from the world of quasi-coherent sheaves to
the realm of contraherent cosheaves).
\end{proof}

\begin{cor} \label{dil-cta-cot-direct-image-becker-coacyclic}
 Let $f\:Y\rarrow X$ be a morphism of quasi-compact semi-separated
schemes.
 Then \par
\textup{(a)} the direct image functor $f_*\:Y\qcoh^\dil\rarrow
X\qcoh^\dil$ takes Becker-coacyclic complexes in $Y\qcoh^\dil$ to
Becker-coacyclic complexes in $X\qcoh^\dil$; \par
\textup{(b)} the direct image functor $f_*\:Y\qcoh^\cta\rarrow
X\qcoh^\cta$ takes Becker-coacyclic complexes in $Y\qcoh^\cta$ to
Becker-coacyclic complexes in $X\qcoh^\cta$; \par
\textup{(c)} the direct image functor $f_*\:Y\qcoh^\cot\rarrow
X\qcoh^\cot$ takes Becker-coacyclic complexes in $Y\qcoh^\cot$ to
Becker-coacyclic complexes in $X\qcoh^\cot$.
\end{cor}

\begin{proof}
 The direct image functor $f_*\:Y\qcoh\rarrow X\qcoh$ takes
contraadjusted quasi-coherent sheaves to contraadjusted ones and
cotorsion quasi-coherent sheaves to cotorsion ones by
Corollary~\ref{cta-cot-direct}.
 Since a complex in $X\qcoh^\cot$ is Becker-coacyclic in $X\qcoh^\cot$
if and only if it is Becker-coacyclic in $X\qcoh$, and the same applies
to Becker-coacyclic complexes in $X\qcoh^\cta$, part~(c) follows
from part~(b).
 The assertion of part~(b) follows from
Corollary~\ref{becker-coacyclic-cta-cot-antilocal}(a)
and Lemma~\ref{qcoh-affine-direct-image-becker-coacyclic}.
 The assertion of part~(c) also follows from
Corollary~\ref{becker-coacyclic-cta-cot-antilocal}(b)
and Lemma~\ref{qcoh-affine-direct-image-becker-coacyclic}.

 The direct image functor~$f_*$ takes dilute quasi-coherent sheaves
to dilute quasi-coherent sheaves by Corollary~\ref{dilute-direct}.
 To deduce part~(a) from part~(b), consider a Becker-coacyclic
complex $\N^\bu$ in $Y\qcoh^\dil$.
 Arguing as in the proof of
Lemma~\ref{second-kind-complex-resolution} (with the arrows inverted),
we construct an acyclic complex of complexes $0\rarrow\N^\bu\rarrow
\cP^{0,\bu}\rarrow\cP^{1,\bu}\rarrow\cP^{2,\bu}\rarrow\dotsb$ in
$Y\qcoh$ such that $\cP^{n,\bu}$ is a contractible complex in
$Y\qcoh^\cta$ for every $n\ge0$.
 By Lemma~\ref{dil-cta-clp-finite-dim}(b) and the dual version of
Corollary~\ref{fdim-cor}, it follows that there is a finite acyclic
complex of complexes $0\rarrow\N^\bu\rarrow\cP^{0,\bu}\rarrow\cP^{1,\bu}
\rarrow\dotsb\rarrow\cP^{N-1,\bu}\rarrow\cQ^\bu\rarrow0$ in $Y\qcoh$
such that $\cP^{n,\bu}$ is a contractible complex in $Y\qcoh^\cta$ for
every $0\le n\le N-1$, while $\cQ^\bu$ is a complex in $Y\qcoh^\cta$.
 Clearly, our finite complex of complexes $0\rarrow\N^\bu\rarrow
\cP^{\bu,\bu}\rarrow\cQ^\bu\rarrow0$ is also acyclic as a complex of
complexes in $Y\qcoh^\dil$.

 The total complex of the bicomplex $0\rarrow\N^\bu\rarrow\cP^{\bu,\bu}
\rarrow\cQ^\bu\rarrow0$ is Becker-coacyclic in $Y\qcoh^\dil$ by
Lemma~\ref{Positselski-trivial-are-Becker-trivial}(a).
 Since the complex $\N^\bu$ is Becker-coacyclic in $Y\qcoh^\dil$, while
the complexes $\cP^{n,\bu}$ are contractible, it follows that
the complex $\cQ^\bu$ is Becker-coacyclic in $Y\qcoh^\dil$, hence also
in $Y\qcoh^\cta$.
 By part~(b), we can conclude that the complex $f_*\cQ^\bu$ is
Becker-coacyclic in $X\qcoh^\cta$, hence also in $X\qcoh^\dil$.
 On the other hand, Corollary~\ref{dilute-direct} tells us that
$0\rarrow f_*\N^\bu\rarrow f_*\cP^{\bu,\bu}\rarrow f_*\cQ^\bu\rarrow0$
is an acyclic finite complex of complexes in $X\qcoh^\dil$.
 Once again, the total complex of the bicomplex
$0\rarrow f_*\N^\bu\rarrow f_*\cP^{\bu,\bu}\rarrow f_*\cQ^\bu\rarrow0$
is Becker-coacyclic in $X\qcoh^\dil$ by
Lemma~\ref{Positselski-trivial-are-Becker-trivial}(a).
 Since the complex $f_*\cQ^\bu$ is Becker-coacyclic in
$X\qcoh^\dil$, while the complexeses $f_*\cP^{n,\bu}$ are contractible,
we arrive to the final conclusion that the complex $f_*\N$ is
Becker-coacyclic in $X\qcoh^\dil$.
\end{proof}

\begin{cor} \label{bco-qcoh-lin-derived-equivalence}
 Let $X$ be a quasi-compact semi-separated scheme with an open
covering\/~$\bW$.
 Then there are natural equivalences of triangulated categories
$$
 \sD^+(X\qcoh)\simeq\sD^+(X\lcth_\bW^\lin)=\sD^{\abs+}(X\lcth_\bW^\lin)
$$
and
\begin{multline*}
 \sD^\bco(X\qcoh)\simeq\sD(X\lcth_\bW^\lin)=\sD^\bctr(X\lcth_\bW^\lin)
 \\ =\sD^\ctr(X\lcth_\bW^\lin)=\sD^\abs(X\lcth_\bW^\lin).
\end{multline*}
\end{cor}

\begin{proof}
 For the first assertion, compare
Corollary~\ref{homotopy-derived-coderived-cor} with
Corollary~\ref{inj-co-contra-cor}(b).
 For the second assertion, compare
Theorem~\ref{quasi-coherent-becker-coderived}
with Corollary~\ref{inj-co-contra-cor}(b).
\end{proof}

\subsection{Antilocality of complexes of antilocally flat cosheaves}
\label{antilocality-of-complexes-of-cosheaves-subsect}
 In this section we use the results of
Section~\ref{naive-co-contra-subsect} in order to transfer the results
of Section~\ref{quasi-cta-cot-antilocality-subsect} from the world of
quasi-coherent sheaves into the realm of contraherent cosheaves.

\begin{cor} \label{complexes-of-al-antilocal}
 Let $X=\bigcup_\alpha U_\alpha$ be a finite affine open covering of
a quasi-compact semi-separated scheme.  Then \par
\textup{(a)} any complex of antilocal contraherent cosheaves on $X$ is
a direct summand of a finitely iterated extension of the direct images
of complexes of contraherent cosheaves from~$U_\alpha$; \par
\textup{(b)} any complex of antilocal locally cotorsion contraherent
cosheaves on $X$ is a direct summand of a finitely iterated extension
of the direct images of complexes of locally cotorsion contraherent
cosheaves from~$U_\alpha$.
\end{cor}

\begin{proof}
 It suffices to apply the equivalences of exact categories from
Lemma~\ref{cta-clp-equivalence} (for part~(a)) and
Lemma~\ref{cta-clp-restricts-to-cot-inj}(a) (for part~(b))
to the results of 
Corollary~\ref{complexes-of-quasi-cta-cot-antilocal}(a\+b).
 One also needs to take into account the fact that the two mentioned
equivalences of exact categories tranform the direct images of
quasi-coherent sheaves~$f_*$ into the direct images of contraherent
cosheaves~$f_!$ for any affine morphism of semi-separated schemes~$f$
(see Corollary~\ref{cta-clp-direct-image-cor}).
 This indirect argument is based on
Theorem~\ref{quasi-coherent-gluing-theorem}.

 Alternatively, one can devise a more direct argument based on
Theorem~\ref{loc-contraherent-gluing-theorem}.
 Let us sketch this argument for part~(a).
 Let $\sR$ be the local class of all commutative rings $R$ and
$\sE^R=\Com(R\modl^\cta)$ be the exact category of complexes of
contraadjusted $R$\+modules.
 Let $\sF(R)=\sE^R$ be the class of all complexes of contraadjusted
$R$\+modules and $\sC^R$ be the class of contractible complexes of
injective $R$\+modules.
 In other words, $\sC^R$ is the class of all injective objects in
the exact category~$\sE^R$.
 There are enough of them, so $(\sF(R),\sC^R)$ is a hereditary
complete cotorsion pair in~$\sE^R$.
 Both the classes $\sE$ and $\sC$ are very colocal
(see Examples~\ref{colocal-classes-examples} and~\cite[proof of
Corollary~7.4]{Pal}).

 Let $\bW$ be any open covering of $X$ to which the covering
$X=\bigcup_\alpha U_\alpha$ is subordinate.
 Applying Theorem~\ref{loc-contraherent-gluing-theorem}, we obtain
a hereditary complete cotorsion pair $(\sF(X),\sC_\bW^X)$ in
the exact category $\sE_\bW^X=\Com(X\lcth_\bW)$.
 Here $\sF(X)$ is the class of all direct summands of finitely iterated
extensions of direct images of complexes of contraherent cosheaves
from $U_\alpha$, while $\sC_\bW^X$ is the class of all locally\+$\sC$
complexes of $\bW$\+locally contraherent cosheaves on~$X$.
 The latter one is the class of all locally contractible complexes of
locally injective $\bW$\+locally contraherent cosheaves on~$X$.

 Let $\P^\bu$ be a complex of antilocal contraherent cosheaves and
$\gJ^\bu\in\sC_\bW^X$ be a locally contractible complex of locally
injective $\bW$\+locally contraherent cosheaves on~$X$.
 It remains to show that the complex of abelian groups
$\Hom^X(\P^\bu,\gJ^\bu)$ is acyclic; then
(by Lemma~\ref{Ext-1-as-homotopy-Hom}) we will be able to conclude
that $\P^\bu\in\sF(X)$.
 In fact, for a complex of locally injective $\bW$\+locally
contraherent cosheaves, being locally contractible is the same thing
as being acyclic with respect to the exact category $X\lcth_\bW^\lin$.
 Now we can refer to the second paragraph of the proof of
Lemma~\ref{clp-linlin-orthogonality}.
\end{proof}

\begin{cor} \label{complexes-of-prj-alf-antilocal}
 Let $X=\bigcup_\alpha U_\alpha$ be a finite affine open covering of
a quasi-compact semi-separated scheme.  Then \par
\textup{(a)} any complex of projective contraherent cosheaves on $X$
is a direct summand of a finitely iterated extension of the direct
images of complexes of projective contraherent cosheaves on~$U_\alpha$;
\par
\textup{(b)} any complex of projective locally cotorsion contraherent
cosheaves on $X$ is a direct summand of a finitely iterated extension
of the direct images of complexes of projective locally cotorsion
contraherent cosheaves on~$U_\alpha$; \par
\textup{(c)} any complex of antilocally flat contraherent cosheaves
on $X$ is a direct summand of a finitely iterated extension of
the direct images of complexes of antilocally flat contraherent cosheaves
on~$U_\alpha$.
\end{cor}

\begin{proof}
 Here we are not aware of any direct argument comparable to
the alternative proof of Corollary~\ref{complexes-of-al-antilocal}
above, and our only recourse is to
Corollary~\ref{complexes-of-quasi-vfl-flat-cta-cot-antilocal}
in conjunction with Lemma~\ref{cta-clp-restricts-to-prj-clf}.
 Corollary~\ref{cta-clp-direct-image-cor}, establishing
the compatibility of the direct images in the quasi-coherent sheaf
and contraherent cosheaf realms, plays an important role in this
argument.
\end{proof}

\begin{cor} \label{acyclic-complexes-of-alf-antilocal}
 Let $X=\bigcup_\alpha U_\alpha$ be a finite affine open covering of
a quasi-compact semi-separated scheme.
 Then any acyclic complex in the exact category of antilocally flat
contraherent cosheaves on $X$ is a direct summand of a finitely
iterated extension of the direct images of acyclic complexes in
the exact categories of antilocally flat contraherent cosheaves
on~$U_\alpha$.
\end{cor}

\begin{proof}
 Follows from
Corollary~\ref{complexes-of-quasi-flat-cta-cocycles-antilocal}
and Lemma~\ref{cta-clp-restricts-to-prj-clf}(c) (similarly to
the previous proofs in this section).
\end{proof}

\begin{cor} \label{acyclic-complexes-of-lcta-lct-antilocal}
 Let $X=\bigcup_\alpha U_\alpha$ be a finite affine open covering of
a quasi-compact semi-separated scheme.  Then \par
\textup{(a)} any acyclic complex in the exact category of antilocal
contraherent cosheaves on $X$ is a direct summand of a finitely
iterated extension of the direct images of acyclic complexes in
the exact categories of contraherent cosheaves on~$U_\alpha$;
\par
\textup{(b)} any acyclic complex in the exact category of antilocal
locally cotorsion contraherent cosheaves on $X$ is a direct summand of
a finitely iterated extension of the direct images of acyclic complexes
in the exact categories of locally cotorsion contraherent
cosheaves on~$U_\alpha$.
\end{cor}

\begin{proof}
 This is similar to Corollary~\ref{complexes-of-al-antilocal}.
 Part~(a) follows from
Corollary~\ref{acyclic-complexes-of-quasi-cta-cot-antilocal}(a)
and Lemma~\ref{cta-clp-equivalence}.
 Part~(b) follows from
Corollary~\ref{acyclic-complexes-of-quasi-cta-cot-antilocal}(b)
and Lemma~\ref{cta-clp-restricts-to-cot-inj}(a).
 This indirect argument is based on
Theorem~\ref{quasi-coherent-gluing-theorem}.

 Alternatively, there is a more direct argument based on
Theorem~\ref{loc-contraherent-gluing-theorem}.
 Here is a sketch for part~(a).
 Let $\sR$ be the local class of all commutative rings $R$ and
$\sE^R=\Com(R\modl^\cta)$ be the exact category of complexes of
contraadjusted $R$\+modules.
 Let $\sF(R)=\Acycl(R\modl^\cta)\sub\sE^R$ be the class of all acyclic
complexes of contraadjusted $R$\+modules and $\sC^R\sub\sE^R$ be
the class of all homotopy injective complexes of injective $R$\+modules.
 This is the setting of
Lemmas~\ref{acyclic-complexes-of-cta-cotorsion-pair}
and~\ref{complexes-loc-hot-inj-of-inj-preenvelope}.

 Let $\bW$ be any open covering of $X$ to which the covering
$X=\bigcup_\alpha U_\alpha$ is subordinate.
 Applying Theorem~\ref{loc-contraherent-gluing-theorem}, we obtain
a hereditary complete cotorsion pair $(\sF(X),\sC_\bW^X)$ in
the exact category $\sE_\bW^X=\Com(X\lcth_\bW)$.
 Here $\sF(X)$ is the class of all direct summands of finitely
iterated extensions of direct images of acyclic complexes in the exact
categories of contraherent cosheaves on $U_\alpha$, while $\sC_\bW^X$
is the class of all locally\+$\sC$ complexes of $\bW$\+locally
contraherent cosheaves on~$X$.
 The latter one is the class of all locally homotopy injective
complexes of locally injective $\bW$\+locally contraherent cosheaves
on~$X$ (in the terminology of Section~\ref{homotopy-lin-subsect}).

 Let $\P^\bu$ be an acyclic complex in the exact category of antilocal
contraherent cosheaves on $X$ and $\gJ^\bu\in\sC_\bW^X$ be a locally
homotopy injective complex of locally injective $\bW$\+locally
contraherent cosheaves on~$X$.
 It remains to show that the complex of abelian groups
$\Hom^X(\P^\bu,\gJ^\bu)$ is acyclic; then
(in view of Lemma~\ref{Ext-1-as-homotopy-Hom}) we will be able to
conclude that $\P^\bu\in\sF(X)$.
 In fact, any locally homotopy injective complex is homotopy locally
injective by Lemma~\ref{homotopy-local-injectivity-characterizations}\,%
(1)\,$\Leftrightarrow$\,(3), and we can refer to the proof of
Lemma~\ref{clp-linlin-orthogonality}.
\end{proof}

\begin{cor} \label{complexes-of-ctrh-prj-alf-acycl-antilocal}
 Let $X=\bigcup_\alpha U_\alpha$ be a finite affine open covering of
a quasi-compact semi-separated scheme.  Then \par
\textup{(a)} any complex of projective contraherent cosheaves on $X$
that is acyclic in the exact category of antilocal contraherent
cosheaves $X\ctrh_\al$ is a direct summand of a finitely iterated
extension of the direct images of complexes of projective contraherent
cosheaves on $U_\alpha$ acyclic in the exact categories
$U_\alpha\ctrh_\al$; \par
\textup{(a)} any complex of projective locally cotorsion contraherent
cosheaves on $X$ that is acyclic in the exact category of antilocal
locally cotorsion contraherent cosheaves $X\ctrh^\lct_\al$ is a direct
summand of a finitely iterated extension of the direct images of
complexes of projective locally cotorsion contraherent cosheaves on
$U_\alpha$ acyclic in the exact categories $U_\alpha\ctrh^\lct_\al$;
\par
\textup{(c)} any complex of antilocally flat contraherent cosheaves
on $X$ that is acyclic in the exact category of antilocal contraherent
cosheaves $X\ctrh_\al$ is a direct summand of a finitely iterated
extension of the direct images of complexes of antilocally flat
contraherent cosheaves on $U_\alpha$ acyclic in the exact categories
$U_\alpha\ctrh_\al$.
\end{cor}

\begin{proof}
 This is similar to Corollary~\ref{complexes-of-prj-alf-antilocal}.
 The assertions follow
from Corollary~\ref{complexes-of-quasi-vfl-flat-cta-cot-acycl-antilocal}
together with Lemma~\ref{cta-clp-restricts-to-prj-clf}.
 Lemmas~\ref{cta-clp-equivalence}
and~\ref{cta-clp-restricts-to-cot-inj}(a) are also helpful.
\end{proof}

\begin{cor} \label{acyclic-in-alf-all-of-lct-loc-contraherent-pair}
 Let\/ $X$ be a quasi-compact semi-separated scheme with an open
covering\/~$\bW$.
 Then the pair of classes of acyclic complexes in the exact category of
antilocally flat contraherent cosheaves\/ $\sF=\Acycl(X\ctrh_\alf)$ and
all complexes of locally cotorsion\/ $\bW$\+locally contraherent
cosheaves\/ $\sC=\Com(X\lcth_\bW^\lct)$ is a hereditary complete
cotorsion pair in the exact category\/ $\Com(X\lcth_\bW)$.
\end{cor}

\begin{proof}
 The argument is based on Theorem~\ref{loc-contraherent-gluing-theorem}
and Corollary~\ref{acyclic-complexes-of-alf-antilocal}.
 Let $\sR$ be the local class of all commutative rings $R$ and
$\sE^R=\Com(R\modl^\cta)$ be the exact category of complexes of
contraadjusted $R$\+modules.
 Let\/ $\sF(R)=\Acycl(R\modl_\fl^\cta)$ be the class of all acyclic
complexes of flat contraadjusted $R$\+modules with flat modules of
cocycles, and let $\sC^R\sub\sE^R$ be the class of all complexes of
cotorsion $R$\+modules.
 This is the datum of the classes $\sR$, \,$\sE^R$, \,$\sC^R$, and
$\sF(R)$ described in
Lemma~\ref{acyclic-complexes-of-cta-flats-cotorsion-pair}.

 Theorem~\ref{loc-contraherent-gluing-theorem} produces a hereditary
complete cotorsion pair $(\sF(X),\sC^X_\bW)$ in $\sE^X_\bW=
\Com(X\lcth_\bW)$, where $\sC^X_\bW=\Com(X\lcth_\bW^\lct)$, while
$\sF(X)$ is the class of all direct summands of finitely iterated
extensions of direct images of acyclic complexes in the exact
categories of antilocally flat contraherent cosheaves on~$U_\alpha$
(for any finite affine open covering $X=\bigcup_\alpha U_\alpha$
subordinate to~$\bW$).
 By Corollary~\ref{acyclic-complexes-of-alf-antilocal}, we have
$\sF(X)=\Acycl(X\ctrh_\alf)$.
\end{proof}

 The following lemma is dual-analogous to
Lemma~\ref{all-vfl-flat-contraacyclic-of-cta-cot-pairs}.

\begin{lem} \label{coacyclic-of-antilocal-all-lin-cotorsion-pairs}
 Let\/ $X$ be a quasi-compact semi-separated scheme with an open
covering\/~$\bW$.  Then \par
\textup{(a)} the pair of classes of Becker-coacyclic complexes
in the exact category of antilocal contraherent cosheaves\/
$\sF'''=\Acycl^\bco(X\ctrh_\al)$ and all complexes of locally
injective\/ $\bW$\+locally contraherent cosheaves\/
$\sC'''=\Com(X\lcth_\bW^\lin)$ is a hereditary complete cotorsion pair
in the exact category\/ $\Com(X\lcth_\bW)$; \par
\textup{(b)} the pair of classes of Becker-coacyclic complexes in
the exact category of antilocal locally cotorsion contraherent
cosheaves\/ $\sF'''=\Acycl^\bco(X\ctrh^\lct_\al)$ and all complexes of
locally injective\/ $\bW$\+locally contraherent cosheaves\/
$\sC'''=\Com(X\lcth_\bW^\lin)$ is a hereditary complete cotorsion pair
in the exact category\/ $\Com(X\lcth^\lct_\bW)$.
\end{lem}

\begin{proof}
 Part~(a): the argument is based on
Theorem~\ref{loc-contraherent-gluing-theorem}.
 Let $\sR$ be the local class of all commutative rings $R$ and
$\sE^R=\Com(R\modl^\cta)$ be the exact category of complexes of
contraadjusted $R$\+modules.
 Let $\sF(R)=\Acycl^\bco(R\modl^\cta)$ be the class of
Becker-coacyclic complexes of contraadjusted $R$\+modules, and
let $\sC^R=\Com(R\modl^\inj)$ be the class of all complexes of
injective $R$\+modules.
 By Corollary~\ref{becker-coderived-of-cta-cot}(a) (for the affine
scheme $\Spec R$), the pair of classes $(\sF(R),\sC^R)$ is
a hereditary complete cotorsion pair in~$\sE^R$.
 Both the classes $\sE$ and $\sC$ are very colocal
by Examples~\ref{colocal-classes-examples}.

 Applying Theorem~\ref{loc-contraherent-gluing-theorem}, we obtain
a hereditary complete cotorsion pair $(\sF''',\sC''')\allowbreak=
(\sF(X),\sC_\bW^X)$ in the exact category $\sE_\bW^X=\Com(X\lcth_\bW)$,
with the right-hand class $\sC_\bW^X=\Com(X\lcth^\lin_\bW)$.
 The theorem also tells us that the class $\sF(X)$ consists of some
complexes of antilocal contraherent cosheaves, $\sF(X)\sub
\Com(X\ctrh_\al)$.
 By Lemmas~\ref{restricting-hereditary-cotorsion}
and~\ref{restricting-cotorsion-pairs-lemma}(a), the cotorsion
pair $(\sF(X),\sC^X)$ restricts to a hereditary complete cotorsion
pair $(\sF,\sC)$ in the exact category $\Com(X\ctrh_\al)$.
 So we have a cotorsion pair $\sF=\sF'''=\sF(X)$ and
$\sC=\Com(X\ctrh_\al^\lin)$ in $\Com(X\ctrh_\al)$.
 Now it is clear from Lemma~\ref{Ext-1-as-homotopy-Hom} that
$\sF'''$ is the class of all Becker-coacyclic complexes in
the exact category $X\ctrh_\al$.
 The cotorsion pair in part~(b) can be constructed similarly or
obtained by restricting the cotorsion pair of part~(a) to
the exact subcategory $\Com(X\lcth^\lct_\bW)\sub\Com(X\lcth_\bW)$.
\end{proof}

\subsection{Locality of Becker-contraacyclicity}
\label{locality-of-Becker-contra-subsect}
 Once again, we refer to Section~\ref{becker-subsect} for
the definition of the full subcategory of
\emph{Becker-contraacyclic complexes} $\Acycl^\bctr(\sE)\sub\Hot(\sE)$
for an exact category~$\sE$.
 In particular, for any scheme $X$ with an open covering $\bW$,
we have the full subcategory of \emph{Becker-contraacyclic\/
$\bW$\+locally contraherent cosheaves} $\Acycl^\bctr(X\lcth_\bW)
\sub\Hot(X\lcth_\bW)$ and the full subcategory of
\emph{Becker-contraacyclic locally cotorsion\/ $\bW$\+locally
contraherent cosheaves} $\Acycl^\bctr(X\lcth_\bW^\lct)\sub
\Hot(X\lcth_\bW^\lct)$. {\hbadness=2350\par}

 Recall that, for a quasi-compact semi-separated scheme $X$,
all projective contraherent cosheaves and all projective locally
cotorsion contraherent cosheaves are (globally) contraherent,
so the classes of such cosheaves do not depend on the open
covering~$\bW$ (see Corollaries~\ref{ctrh-lcth-proj}(b)
and~\ref{ctrh-lcth-lct-proj}(b)).
 Consequently, the property of a complex of (locally contraadjusted
or locally cotorsion) $\bW$\+locally contraherent cosheaves on $X$
to be Becker-contraacyclic does not change when one refines
the covering~$\bW$ to another open covering~$\bW'$.

 The aim of this section is to show that the Becker-contraacyclicity
can be checked locally on quasi-compact semi-separated schemes.
 Furthermore, we will see that a complex in $X\lcth_\bW^\lct$ is
Becker-contraacyclic if and only if it is Becker-contraacyclic as
a complex in $X\lcth_\bW$.

\begin{lem} \label{inverse-images-preserve-contraacyclicity}
 Let $f\:Y\rarrow X$ be either an affine morphism of schemes, or
a morphism of quasi-compact semi-separated schemes.
 Let\/ $\bW$ be an open covering of $X$ and\/ $\bT$ be an open covering
of $Y$ such that the morphism~$f$ is $(\bW,\bT)$\+coaffine (and
$(\bW,\bT)$\+affine if we are in the setting when $f$~is affine).
 Then \par
\textup{(a)} if the morphism~$f$ is very flat, then the inverse image
functor $f^!\:X\lcth_\bW\rarrow Y\lcth_\bT$ takes Becker-contraacyclic
complexes of\/ $\bW$\+locally contraherent cosheaves to
Becker-contraacyclic complexes of\/ $\bT$\+locally contraherent 
cosheaves; \par
\textup{(b)} if the morphism~$f$ is flat, then the inverse image
functor $f^!\:X\lcth_\bW^\lct\rarrow Y\lcth_\bT^\lct$ takes
Becker-contraacyclic complexes of locally cotorsion\/ $\bW$\+locally
contraherent cosheaves to Becker-contraacyclic complexes of locally
cotorsion\/ $\bT$\+locally contraherent cosheaves.
\end{lem}

\begin{proof}
 Part~(a): if the morphism~$f$ is $(\bW,\bT)$\+affine and
$(\bW,\bT)$\+coaffine, then we have a pair of adjoint functors
$f_!\:Y\lcth_\bT\rarrow X\lcth_\bW$ and $f^!\:X\lcth_\bW\rarrow
Y\lcth_\bT$ (see formula~\eqref{direct-inverse-lct-adjunction}
in Section~\ref{direct-inverse-loc-contra}).
 The functor~$f^!$ is exact (as well as the functor~$f_!$, but we do not
use this fact), so Lemma~\ref{exact-with-adjoint-preservation-lemma}(b)
is applicable.

 When the schemes $X$ and $Y$ are quasi-compact and semi-separated,
a similar argument goes through using the partial
adjunction~\eqref{direct-inverse-clp-adjunction} from
Section~\ref{homology-subsection}.
 Notice that the functor $f_!\:Y\ctrh_\al\rarrow X\ctrh_\al$ takes
projective contraherent cosheaves to projective contraherent
cosheaves by Corollary~\ref{proj-direct-gen}(a).

 The proof of part~(b) is similar.
\end{proof}

\begin{thm} \label{Becker-contraacyclicity-local-on-qcomp-qsep}
 Let $X$ be a quasi-compact semi-separated scheme with an open
covering\/ $\bW$, and let $X=\bigcup_\alpha U_\alpha$ be a finite
affine open covering subordinate to\/~$\bW$.  Then \par
\textup{(a)} a complex\/ $\gB^\bu$ in $X\lcth_\bW$ is
Becker-contraacyclic if and only if, for every~$\alpha$,
the complex\/ $\gB^\bu[U_\alpha]$ is Becker-contraacyclic
in\/ $\O(U_\alpha)\modl^\cta$; \par
\textup{(b)} a complex\/ $\gB^\bu$ in $X\lcth_\bW^\lct$ is
Becker-contraacyclic if and only if, for every~$\alpha$,
the complex\/ $\gB^\bu[U_\alpha]$ is Becker-contraacyclic
in\/ $\O(U_\alpha)\modl^\cot$.
\end{thm}

\begin{proof}
 Part~(a): the ``only if'' implication follows immediately from
Lemma~\ref{inverse-images-preserve-contraacyclicity}(a) (for
the very flat affine open embeddings of quasi-compact semi-separated
schemes $j_\alpha\:U_\alpha\rarrow X$).
 Our proof of the ``if'' assertion is a version
of~\cite[Theorem~5.1(a)]{Pal}.
 The argument is based on
Corollary~\ref{complexes-of-prj-alf-antilocal}(a).
 Indeed, let $\P^\bu$ be a complex in $X\ctrh_\prj$.
 By Lemma~\ref{Ext-1-as-homotopy-Hom}, instead of showing that
the complex $\Hom_{X\lcth_\bW}(\P^\bu,\gB^\bu)$ is acyclic, we can,
equivalently, show that $\Ext^1_{\Com(X\lcth_\bW)}(\P^\bu,\gB^\bu)=0$.
 The latter condition is obviously preserved by finitely iterated
extensions (and the passages to direct summands) of
the complexes~$\P^\bu$.
 Alternatively, one can say that extensions of complexes of projective
objects are termwise split, so they correspond to distinguished
triangles in the homotopy category $\Hot(X\lcth_\bW)$.
 So Corollary~\ref{complexes-of-prj-alf-antilocal}(a) allows to assume
that $\P^\bu=j_!\Q^\bu$, where $j\:U\rarrow X$ is the embedding of
an affine open subscheme $U=U_\alpha$ for some index~$\alpha$ and
$\Q^\bu$ is a complex in $U\ctrh_\prj$.
 It remains to use the adjunction $\Hom_{X\lcth_\bW}(j_!\Q^\bu,\gB^\bu)
\simeq\Hom_{U\ctrh}(\Q^\bu,j^!\gB^\bu)$
(see~\eqref{direct-inverse-lct-adjunction}
or~\eqref{direct-inverse-clp-adjunction}).

 The proof of part~(b) is similar and based on
Lemma~\ref{inverse-images-preserve-contraacyclicity}(b)
and Corollary~\ref{complexes-of-prj-alf-antilocal}(b).
\end{proof}

\begin{cor} \label{Becker-contraacyclic-are-acyclic}
 Let $X$ be a quasi-compact semi-separated scheme with an open
covering\/~$\bW$.  Then \par
\textup{(a)} all Becker-contraacyclic complexes in the exact
category $X\lcth_\bW$ are acyclic in $X\lcth_\bW$; \par
\textup{(b)} all Becker-contraacyclic complexes in the exact
category $X\lcth_\bW^\lct$ are acyclic in $X\lcth_\bW^\lct$.
\end{cor}

\begin{proof}
 One generally expects contraacyclic complexes to be acyclic in
reasonable contexts, but we do \emph{not} have a general result of this
kind applicable in the situation at hand; cf.\
Lemmas~\ref{bounded-Becker-trivial-are-acyclic}(b)\+-%
\ref{with-co-kernels-Becker-trivial-are-acyclic}(b) and
Remark~\ref{Becker-co-contra-acyclic-are-acyclic-remark}.
 Instead, we notice that the Becker-contraacyclicity in $X\lcth_\bW$
or $X\lcth_\bW^\lct$ is a local property by
Theorem~\ref{Becker-contraacyclicity-local-on-qcomp-qsep}, and so
is acyclicity of complexes in these exact categories
(by Lemma~\ref{acyclicity-in-lcth-criterion}).
 This reduces the question to the case of an affine scheme $U$ with
the open covering~$\{U\}$; in other words, we need to prove
the acyclicity of Becker-contraacyclic complexes in the exact
categories $R\modl^\cta$ and $R\modl^\cot$ for a commutative ring~$R$.
 Here we recall that a complex in $R\modl^\cta$ is acyclic in
$R\modl^\cta$ if and only if it is acyclic in $R\modl$; and a complex
in $R\modl^\cot$ is acyclic in $R\modl^\cot$ if and only if it is
acyclic in $R\modl$ (by Theorem~\ref{cotorsion-periodicity}).
 So, by Proposition~\ref{becker-contraacyclicity-for-cta-cot},
it remains to prove that any Becker-contraacyclic complex in $R\modl$
is acyclic.
 Finally, Lemma~\ref{with-co-kernels-Becker-trivial-are-acyclic}(b)
is applicable now.
\end{proof}

\begin{cor} \label{lct-Becker-contraacyclic-iff-as-lcta}
 Let $X$ be a quasi-compact semi-separated scheme with an open
covering\/~$\bW$.  Then a complex in $X\lcth_\bW^\lct$ is
Becker-contraacyclic in $X\lcth_\bW^\lct$ if and only if it is
Becker-contraacyclic in $X\lcth_\bW$,
$$
 \Acycl^\bctr(X\lcth_\bW^\lct)=\Com(X\lcth_\bW^\lct)\cap
 \Acycl^\bctr(X\lcth_\bW).
$$
\end{cor}

\begin{proof}
 Once again, Theorem~\ref{Becker-contraacyclicity-local-on-qcomp-qsep}
reduces the question to the case of an affine scheme $U$ with the open
covering~$\{U\}$.
 So we need to prove that, for any commutative ring $R$, a complex
in $R\modl^\cot$ is Becker-contraacyclic in $R\modl^\cot$ if and only
if it is Becker-contraacyclic in $R\modl^\cta$.
 For this purpose, it suffices to compare the two parts of
Proposition~\ref{becker-contraacyclicity-for-cta-cot}.
\end{proof}

\begin{lem} \label{morphism-of-affines-direct-image-Becker-contra}
 Let $f\:V\rarrow U$ be a morphism of affine schemes.  Then \par
\textup{(a)} the direct image functor $f_!\:V\ctrh\rarrow U\ctrh$
takes Becker-contraacyclic complexes to Becker-contraacyclic complexes;
\par
\textup{(b)} the direct image functor $f_!\:V\ctrh^\lct\rarrow
U\ctrh^\lct$ takes Becker-contraacyclic complexes to
Becker-contraacyclic complexes.
\end{lem}

\begin{proof}
 Given a commutative ring homomorphism $R\rarrow S$, we have to show
that any Becker-contraacyclic complex in $S\modl^\cta$ is
Becker-contraacyclic as a complex in $R\modl^\cta$, and any
Becker-contraacyclic complex in $S\modl^\cot$ is Becker-contraacyclic
as a complex in $R\modl^\cot$.
 In view of Proposition~\ref{becker-contraacyclicity-for-cta-cot},
it suffices to show that any Becker-contraacyclic complex in $S\modl$
is Becker-contraacyclic in $R\modl$.
 Now the exact functor of restriction of scalars $S\modl\rarrow R\modl$
has a left adjoint functor of extension of scalars $S\ot_R{-}\,\:
R\modl\rarrow S\modl$, and we can refer to
Lemma~\ref{exact-with-adjoint-preservation-lemma}(b).

 Alternatively, the assertions of the lemma are the particular case
of Corollary~\ref{quasi-cta-cot-contraacyclic-direct-images}
for affine schemes.
\end{proof}

\begin{cor} \label{affine-morphism-direct-image-Becker-contra}
 Let $f\:Y\rarrow X$ be an affine morphism of quasi-compact
semi-separated schemes.
 Let\/ $\bW$ be an open covering of $X$ and\/ $\bT$ be an open covering
of $Y$ such that the morphism~$f$ is $(\bW,\bT)$\+affine.  Then \par
\textup{(a)} the direct image functor $f_!\:Y\lcth_\bT\rarrow
X\lcth_\bW$ takes Becker-contraacyclic complexes to Becker-contraacyclic
complexes; \par
\textup{(b)} the direct image functor $f_!\:Y\lcth_\bT^\lct\rarrow
X\lcth_\bW^\lct$ takes Becker-contraacyclic complexes to
Becker-contraacyclic complexes.
\end{cor}

\begin{proof}
 In view of Corollary~\ref{lct-Becker-contraacyclic-iff-as-lcta},
it suffices to prove part~(a).
 Let $U\sub X$ be an affine open subscheme subordinate to~$\bW$,
and let $j\:U\rarrow X$ be the open embedding morphism.
 Put $V=U\times_XY$, and denote by $j'\:V\rarrow Y$ and
$f'\:V\rarrow U$ the related morphisms.
 Let $\gB^\bu$ be a Becker-contraacyclic complex in $Y\lcth_\bT$.
 We have the base change isomorphism $j^!f_!\gB^\bu\simeq f'_!j'{}^!
\gB^\bu$ of complexes in $U\ctrh$.
 By Lemma~\ref{inverse-images-preserve-contraacyclicity}(a),
the complex $j'{}^!\gB^\bu$ is Becker-contraacyclic in $V\ctrh$.
 By Lemma~\ref{morphism-of-affines-direct-image-Becker-contra}(a),
it follows that the complex $f'_!j'{}^!\gB^\bu$ is Becker-contraacyclic
in $U\ctrh$.
 We have shown that the complex $j^!f_!\gB^\bu$ is Becker-contraacyclic
in $U\ctrh$, and it remains to refer to
Theorem~\ref{Becker-contraacyclicity-local-on-qcomp-qsep}(a) to deduce
the desired conclusion that the complex $f_!\gB^\bu$ is
Becker-contraacyclic in $X\lcth_\bW$.
\end{proof}

\begin{rem} \label{locality-of-Becker-contraacyclicity-remark}
 The usual proof of locality of Positselski-coacyclicity for complexes
of quasi-coherent sheaves (or matrix factorizations)
\cite[Remark~1.3]{EP} is much simpler than our very indirect proof of
Theorem~\ref{Becker-contraacyclicity-local-on-qcomp-qsep}.
 The same usual argument, based on the \v Cech coresolution, also
proves locality of Becker-coacyclicity for complexes of quasi-coherent
sheaves on quasi-compact semi-separated
schemes~\cite[Proposition~A.10]{Psemten} (see
Corollary~\ref{qcoh-becker-coacyclicity-is-local} above).
 The same (or rather, dual-analogous) argument can be also used for
proving locality of Positselski-contraacyclicity for complexes of
$\bW$\+locally contraherent cosheaves.
 (See also
Corollary~\ref{locality-of-co-contra-acyclicity-on-loc-Notherian}
below for the Noetherian case, where a somewhat different approach
applies.)

 But for Becker-contraacyclicity, this argument appears to break down on
a subtle point that there is no direct proof of preservation of
contraacyclicity by the direct images (even with respect to embeddings
of affine open subschemes into quasi-compact semi-separated schemes).
 Indeed, our proof of
Corollary~\ref{affine-morphism-direct-image-Becker-contra} above was
based on Theorem~\ref{Becker-contraacyclicity-local-on-qcomp-qsep},
and not the other way around.
 The problem is that the direct image functor $j_!\:U\ctrh\rarrow
X\lcth_\bW$ has \emph{no} left adjoint, so
Lemma~\ref{exact-with-adjoint-preservation-lemma}(b) is not applicable
to it.
 Even in the fully affine case, in the context of
Lemma~\ref{morphism-of-affines-direct-image-Becker-contra}, we had
to leave the realm of contraherent cosheaves and move the argument to
the world of arbitrary modules over commutative rings in order to
apply Lemma~\ref{exact-with-adjoint-preservation-lemma}(b).

 An alternative approach is to deduce the preservation of
Becker-contraacyclicity of complexes of contraherent cosheaves by
the direct images with respect to the embeddings of affine open
subschemes from
Corollary~\ref{quasi-cta-cot-contraacyclic-direct-images}
(as mentioned in the proof of
Corollary~\ref{morphism-qcomp-qsep-direct-image-of-bcontra-al} below).
\end{rem}

 Let $X$ be a quasi-compact semi-separated scheme with an open
covering~$\bW$.
 Now we apply the same definition of a Becker-contraacyclic complex
from Section~\ref{becker-subsect} to the exact category of antilocal
contraherent cosheaves $X\ctrh_\al$ and the exact category of
antilocal locally cotorsion contraherent cosheaves $X\ctrh^\lct_\al$.
 Notice that the projective objects in $X\ctrh_\al$ coincide with
those in $X\lcth_\bW$; so a complex in $X\ctrh_\al$ is
Becker-contraacyclic if and only if it is Becker-contraacyclic in
$X\lcth_\bW$.
 Similarly, the projective objects in $X\ctrh^\lct_\al$ coincide with
those in $X\lcth_\bW^\lct$; so a complex in $X\ctrh^\lct_\al$ is
Becker-contraacyclic if and only if it is Becker-contraacyclic in
$X\lcth_\bW^\lct$.

 We start with a corollary establishing the antilocality of
Becker-contraacyclicity of complexes of antilocal contraherent
cosheaves (in the spirit of
Section~\ref{antilocality-of-complexes-of-cosheaves-subsect})
before continuing the discussion of preservation of
Becker-contraacyclicity by the direct images.

\begin{cor} \label{contraacyclic-complexes-of-lcta-lct-antilocal}
 Let $X=\bigcup_\alpha U_\alpha$ be a finite affine open covering of
a quasi-compact semi-separated scheme.  Then \par
\textup{(a)} any Becker-contraacyclic complex in the exact category of
antilocal contraherent cosheaves on $X$ is a direct summand of
a finitely iterated extension of the direct images of
Becker-contraacyclic complexes in the exact categories of contraherent
cosheaves on~$U_\alpha$; \par
\textup{(b)} any Becker-contraacyclic complex in the exact category of
antilocal locally cotorsion contraherent cosheaves on $X$ is a direct
summand of a finitely iterated extension of the direct images of
Becker-contraacyclic complexes in the exact categories of locally
cotorsion contraherent cosheaves on~$U_\alpha$.
\end{cor}

\begin{proof}
 This is similar to Corollaries~\ref{complexes-of-al-antilocal}
and~\ref{acyclic-complexes-of-lcta-lct-antilocal}.
 Part~(a) follows from
Corollary~\ref{contraacyclic-complexes-of-quasi-cta-cot-antilocal}(a)
and Lemma~\ref{cta-clp-equivalence}.
 Part~(b) follows from
Corollary~\ref{contraacyclic-complexes-of-quasi-cta-cot-antilocal}(b)
and Lemma~\ref{cta-clp-restricts-to-cot-inj}(a).
 Corollary~\ref{cta-clp-direct-image-cor} also needs to be used.
\end{proof}

\begin{lem} \label{into-affine-direct-images-of-bcontra-al}
 Let $f\:V\rarrow U$ be a morphism from a quasi-compact semi-separated
scheme $V$ to an affine scheme~$U$.  Then \par
\textup{(a)} the direct image functor $f_!\:V\ctrh_\al\rarrow U\ctrh$
takes Becker-contraacyclic complexes to Becker-contraacyclic complexes;
\par
\textup{(b)} the direct image functor $f_!\:V\ctrh^\lct_\al\rarrow
U\ctrh^\lct$ takes Becker-contraacyclic complexes to
Becker-contraacyclic complexes.
\end{lem}

\begin{proof}
 In view of Corollary~\ref{lct-Becker-contraacyclic-iff-as-lcta},
it suffices to prove part~(a).
 Let $V=\bigcup_{\alpha=1}^N V_\alpha$ be a finite affine open covering
of the scheme~$V$, and let $\gB^\bu$ be a Becker-contraacyclic
complex of antilocal contraherent cosheaves on~$V$.
 Consider the \v Cech resolution~\eqref{contraherent-cech}
\begin{multline} \label{cech-resolution-of-bcontra-complex} \textstyle
 0 \lrarrow k_{1,\dotsc,N}{}_!k_{1,\dotsc,N}^!\gB^\bu \lrarrow \dotsb \\
 \textstyle\lrarrow \bigoplus_{\alpha<\beta}k_{\alpha,\beta}{}_!
 k_{\alpha,\beta}^!\gB^\bu \lrarrow
 \bigoplus_\alpha k_\alpha{}_!k_\alpha^!\gB^\bu\lrarrow\gB^\bu\lrarrow0,
\end{multline}
where $k_{\alpha_1,\dotsc,\alpha_i}$ denote the open embeddings
$V_{\alpha_1}\cap\dotsb\cap V_{\alpha_i}\rarrow V$.
 By Lemma~\ref{inverse-images-preserve-contraacyclicity}(a), all
the complexes $k_{\alpha_1,\dotsc,\alpha_i}^!\gB^\bu$ are
Becker-contraacyclic.
 All the terms of the complex~\eqref{cech-resolution-of-bcontra-complex}
are obviously complexes of antilocal contraherent cosheaves on~$X$.

 The finite complex~\eqref{cech-resolution-of-bcontra-complex} is
acyclic as a complex in $\Com(X\ctrh_\al)$, since the full subcategory
$X\ctrh_\al$ is closed under kernels of epimorphisms in $X\ctrh$
(see Corollary~\ref{clp-characterizations}(b)).
 By Corollary~\ref{clp-direct}(a), the functor~$f_!$
takes~\eqref{cech-resolution-of-bcontra-complex} to an acyclic finite
complex in $\Com(U\ctrh)$.
 By Lemma~\ref{morphism-of-affines-direct-image-Becker-contra}(a), all
the terms of the complex obtained by applying~$f_!$
to~\eqref{cech-resolution-of-bcontra-complex}, except perhaps
the rightmost one, are Becker-contraacyclic in $U\ctrh$.
 The total complex of the bicomplex obtained by applying~$f_!$
to~\eqref{cech-resolution-of-bcontra-complex} is absolutely acyclic,
hence Becker-contraacyclic in $U\ctrh$
(see Lemma~\ref{Positselski-trivial-are-Becker-trivial}(a)).
 It follows that the complex $f_!\gB^\bu$ is also Becker-contraacyclic
in $U\ctrh$.

 Alternatively, the assertions of the lemma can be deduced from
Lemma~\ref{morphism-of-affines-direct-image-Becker-contra} together with
Corollary~\ref{contraacyclic-complexes-of-lcta-lct-antilocal} and
Corollary~\ref{clp-direct}(a\+b).
\end{proof}

\begin{cor} \label{morphism-qcomp-qsep-direct-image-of-bcontra-al}
 Let $f\:Y\rarrow X$ be a morphism of quasi-compact semi-separated
schemes.  Then \par
\textup{(a)} the direct image functor $f_!\:Y\ctrh_\al\rarrow
X\ctrh_\al$ takes Becker-contraacyclic complexes to Becker-contraacyclic
complexes; \par
\textup{(b)} the direct image functor $f_!\:Y\ctrh_\al^\lct\rarrow
X\ctrh_\al^\lct$ takes Becker-contraacyclic complexes to
Becker-contraacyclic complexes.
\end{cor}

\begin{proof}
 Let us prove part~(a).
 The existence of the direct image functor $f_!\:Y\ctrh_\al\allowbreak
\rarrow X\ctrh_\al$ was established in Corollary~\ref{clp-direct}(a).
 Let $U\sub X$ be an affine open subscheme with the open embedding
morphism $j\:U\rarrow X$.
 Put $V=U\times_XY$, and denote by $j'\:V\rarrow Y$ and
$f'\:V\rarrow U$ the related morphisms.
 Let $\gB^\bu$ be a Becker-contraacyclic complex in $Y\ctrh_\al$.
 As in the proof of
Corollary~\ref{affine-morphism-direct-image-Becker-contra}, we have
the base change isomorphism $j^!f_!\gB^\bu\simeq f'_!j'{}^!
\gB^\bu$ of complexes in $U\ctrh$ (see
Section~\ref{direct-inverse-loc-contra}).
 By Lemma~\ref{inverse-images-preserve-contraacyclicity}(a) (for
the very flat affine open embedding of quasi-compact semi-separated
schemes $j'\:V\rarrow Y$), the complex $j'{}^!\gB^\bu$ is
Becker-contraacyclic in $V\ctrh$.
 By Corollary~\ref{clp-inverse}(c), the terms of the complex
$j'{}^!\gB^\bu$ are antilocal contraherent cosheaves on~$V$.
 So Lemma~\ref{into-affine-direct-images-of-bcontra-al}(a) tells us
that the complex $f'_!j'{}^!\gB^\bu$ is Becker-contraacyclic in
$U\ctrh$, and it remains to invoke
Theorem~\ref{Becker-contraacyclicity-local-on-qcomp-qsep}(a).

 Alternatively, one can pass from contraherent cosheaves to
quasi-coherent sheaves via Lemmas~\ref{cta-clp-equivalence}
and~\ref{cta-clp-restricts-to-cot-inj}(a) with
Corollary~\ref{cta-clp-direct-image-cor}, and deduce the desired
assertions of the corollary from
Corollary~\ref{quasi-cta-cot-contraacyclic-direct-images}.
\end{proof}

\subsection{Becker's contraderived categories}
\label{Becker-contraderived-subsect}
 As in Section~\ref{Becker-coderived-subsect}, we start with
a collection of elementary observations.

\begin{cor} \label{homotopy-derived-contraderived-cor}  \hbadness=1700
\textup{(a)} For any quasi-compact semi-separated scheme $X$,
the natural functors\/
$\Hot^-(X\ctrh_\prj)\rarrow\sD^-(X\ctrh)\rarrow\sD^-(X\lcth_\bW)
\rarrow\sD^-(X\lcth)$ are equivalences of triangulated categories.
 For any scheme $X$ with an open covering\/ $\bW$ and any symbol\/
$\bst=\b$, $\abs+$, $\abs-$, $\bctr$, $\ctr$, or\/~$\abs$, the natural
triangulated functor\/ $\Hot^\st(X\ctrh_\prj)\rarrow\sD^\st(X\lcth_\bW)$
is fully faithful. \par
\textup{(b)} For any quasi-compact semi-separated scheme $X$,
the natural functors\/ $\Hot^-(X\ctrh^\lct_\prj)\rarrow
\sD^-(X\ctrh^\lct)\rarrow\sD^-(X\lcth_\bW^\lct)\rarrow
\sD^-(X\lcth^\lct)$ are equivalences of triangulated categories.
 For any scheme $X$ with an open covering\/ $\bW$ and any symbol\/
$\bst=\b$, $\abs+$, $\abs-$, $\bctr$, $\ctr$, or\/~$\abs$, the natural
triangulated functor\/ $\Hot^\st(X\ctrh^\lct_\prj)\rarrow
\sD^\st(X\lcth^\lct_\bW)$ is fully faithful.
\end{cor}

\begin{proof}
 This is dual-analogous to
Corollary~\ref{homotopy-derived-coderived-cor}.
 In each part~(a\+b), the first assertion follows from there being
enough projective objects in the respective exact category,
together with Proposition~\ref{infinite-resolutions}(a).
 In the second assertions of both parts, the notation
$X\ctrh_\prj$ or $X\ctrh^\lct_\prj$ stands for the full additive
subcategories of projective objects in the exact categories
$X\ctrh$ or $X\ctrh^\lct$ on an arbitrary scheme~$X$.
 Irrespectively of there being enough such projectives,
these kind of assertions hold in any exact category (or in any exact
category with exact functors of infinite product, as appropriate)
by Lemma~\ref{homotopy-inj-proj-fully-faithful}(b).
 In the case of the Becker contraderived category, the assertion holds
essentially by the definition.
\end{proof}

 Let $X$ be a quasi-compact semi-separated scheme with an open
covering~$\bW$.
 The following theorem and corollary tell us that the Becker
contraderived categories of $\bW$\+locally contraherent cosheaves
and locally cotorsion $\bW$\+locally contraherent cosheaves on $X$
are well-behaved.

\begin{thm} \label{becker-contraderived-cotorsion-pair}
\textup{(a)} The pair of classes of all complexes of projective
contraherent cosheaves\/ $\sF=\Com(X\ctrh_\prj)$ and
Becker-contraacyclic complexes of\/ $\bW$\+locally contraherent
cosheaves\/ $\sC=\Acycl^\bctr(X\lcth_\bW)$ is a hereditary complete
cotorsion pair in the exact category\/ $\Com(X\lcth_\bW)$. \par
\textup{(b)} The pair of classes of all complexes of projective locally
cotorsion contraherent cosheaves\/ $\sF=\Com(X\ctrh^\lct_\prj)$ and
Becker-contraacyclic complexes of locally cotorsion\/ $\bW$\+locally
contraherent cosheaves\/ $\sC=\Acycl^\bctr(X\lcth_\bW^\lct)$ is
a hereditary complete cotorsion pair in the exact category\/
$\Com(X\lcth_\bW^\lct)$.
\end{thm}

\begin{proof}
 The argument is based on
Theorem~\ref{loc-contraherent-gluing-theorem}.
 To prove part~(a), let $\sR$ be the local class of all commutative
rings $R$ and $\sE^R=\Com(R\modl^\cta)$ be the exact category of
complexes of contraadjusted $R$\+modules.
 Let $\sF(R)=\Com(R\modl_\vfl^\cta)\sub\sE^R$ be the class of
all complexes of very flat contraadjusted $R$\+modules, and let
$\sC^R=\Acycl^\bctr(R\modl^\cta)$ be the class of Becker-contraacyclic
complexes of contraadjusted $R$\+modules.
 Then $\sE$ is a very colocal class by
Examples~\ref{colocal-classes-examples}, while $\sC$ is a very colocal
class by Theorem~\ref{Becker-contraacyclicity-local-on-qcomp-qsep}(a)
(for the affine scheme $\Spec R$ with the open covering
$\bW_{\Spec R}=\{\Spec R\}$) and
Lemma~\ref{morphism-of-affines-direct-image-Becker-contra}(a).
 Furthermore, the pair of classes $(\sF(R),\sC^R)$ is a hereditary
complete cotorsion pair in $\Com(R\modl^\cta)$ by the proof of
Lemma~\ref{all-vfl-flat-contraacyclic-of-cta-cot-pairs}(a)
(for the affine scheme $\Spec R$).

 Applying Theorem~\ref{loc-contraherent-gluing-theorem}, we obtain
a hereditary complete cotorsion pair $(\sF(X),\sC^X_\bW)$ in
the exact category $\sE^X_\bW=\Com(X\lcth_\bW)$.
 Here $\sF(X)$ is the class of all direct summands of finitely iterated
extensions of direct images of complexes of projective contraherent
cosheaves from $U_\alpha$, while $\sC^X_\bW$ is the class of all
locally\+$\sC$ complexes of $\bW$\+locally contraherent cosheaves on~$X$.
 By Corollary~\ref{complexes-of-prj-alf-antilocal}(a)
(for the scheme~$X$), \,$\sF(X)$ is the class of all complexes of
projective contraherent cosheaves on~$X$.
 The class $\sC^X_\bW$ can be called the class of all ``locally
Becker-contraacyclic complexes of $\bW$\+locally contraherent cosheaves
on~$X$\,''.
 By Theorem~\ref{Becker-contraacyclicity-local-on-qcomp-qsep}(a)
(for the scheme~$X$), we have $\sC^X_\bW=\Acycl^\bctr(X\lcth_\bW)$,
as desired. {\hbadness=1650\par}

 The proof of part~(b) is similar and applies
Theorem~\ref{loc-contraherent-gluing-theorem} to the exact categories
$\sE^R=\Com(R\modl^\cot)\sub\sK^R=\Com(R\modl^\cta)$.
 Theorem~\ref{Becker-contraacyclicity-local-on-qcomp-qsep}(b),
Lemma~\ref{morphism-of-affines-direct-image-Becker-contra}(b),
Lemma~\ref{all-vfl-flat-contraacyclic-of-cta-cot-pairs}(b), and
Corollary~\ref{complexes-of-prj-alf-antilocal}(b) play the key roles.
\end{proof}

\begin{cor} \label{becker-contraderived-of-lcta-lct-well-behaved}
\textup{(a)} The composition of triangulated functors
$$
 \Hot(X\ctrh_\prj)\lrarrow\Hot(X\lcth_\bW)\lrarrow
 \sD^\bctr(X\lcth_\bW)
$$
is a triangulated equivalence\/ $\Hot(X\ctrh_\prj)\simeq
\sD^\bctr(X\lcth_\bW)$. \par
\textup{(b)} The composition of triangulated functors
$$
 \Hot(X\ctrh^\lct_\prj)\lrarrow\Hot(X\lcth^\lct_\bW)\lrarrow
 \sD^\bctr(X\lcth^\lct_\bW)
$$
is a triangulated equivalence\/ $\Hot(X\ctrh_\prj^\lct)\simeq
\sD^\bctr(X\lcth^\lct_\bW)$.
\end{cor}

\begin{proof}
 This is a standard argument similar to~\cite[proof of
Corollary~7.4]{PS4} and dual-analogous to the argument in
Corollary~\ref{becker-coderived-of-cta-cot}.
 Let us spell out part~(a).
 It is clear from the definitions that the triangulated functor
$\Hot(X\ctrh_\prj)\rarrow\sD^\bctr(X\lcth_\bW)$ is fully faithful.
 In order to prove that it is a triangulated equivalence, it suffices
to find, for any complex of $\bW$\+locally contraherent  cosheaves
$\gM^\bu$ on $X$, a complex of projective contraherent cosheaves
$\P^\bu$ together with a morphism of complexes $\P^\bu\rarrow\gM^\bu$
whose cone is Becker-contraacyclic in $X\lcth_\bW$.
 Indeed, let $0\rarrow\gB^\bu\rarrow\P^\bu\rarrow\gM^\bu\rarrow0$
be a special precover sequence~\eqref{special-precover-sequence}
for the complex $\gM^\bu$ in the complete cotorsion pair of 
Theorem~\ref{becker-contraderived-cotorsion-pair}(a).
 So $\P^\bu$ is a complex of projective contraherent cosheaves
and $\gB^\bu$ is a Becker-contraacyclic complex of $\bW$\+locally
contraherent cosheaves on~$X$.
 The total complex of the short exact sequence of complexes
$0\rarrow\gB^\bu\rarrow\P^\bu\rarrow\gM^\bu\rarrow0$ is absolutely
acyclic, hence Becker-contraacyclic in $X\lcth_\bW$ by
Lemma~\ref{Positselski-trivial-are-Becker-trivial}(a).
 Since $\gB^\bu$ is also Becker-contraacyclic, it follows that
the cone of the morphism $\P^\bu\rarrow\gM^\bu$ is Becker-contraacyclic
in $X\lcth_\bW$ as well.
\end{proof}

 In particular, one can see from
Corollary~\ref{becker-contraderived-of-lcta-lct-well-behaved}
that the Becker contraderived categories $\sD^\bctr(X\lcth_\bW)$
and $\sD^\bctr(X\lcth^\lct_\bW)$ do not depend on the open covering
$\bW$ of a quasi-compact semi-separated scheme $X$ (i.~e., they
do not change when the open covering is refined).
 Actually, this fact was established already in
Corollaries~\ref{ctrh-lcth-cor}(a) and~\ref{lct-ctrh-lcth-cor}(a).

\begin{cor} \label{becker-contraderived-of-al-lct-alf-well-behaved}
\textup{(a)} The hereditary complete cotorsion pair
\textup{(}$\Com(X\ctrh_\prj)$, $\Acycl^\bctr(X\lcth_\bW)$\textup{)}
in the exact category\/ $\Com(X\lcth_\bW)$ restricts to a hereditary
complete cotorsion pair \textup{(}$\Com(X\ctrh_\prj)$,
$\Acycl^\bctr(X\ctrh_\al)$\textup{)} in the exact subcategory\/
$\Com(X\ctrh_\al)\sub\Com(X\lcth_\bW)$.
 Consequently, the composition of triangulated functors
$$
 \Hot(X\ctrh_\prj)\lrarrow\Hot(X\ctrh_\al)\lrarrow
 \sD^\bctr(X\ctrh_\al)
$$
is a triangulated equivalence\/ $\Hot(X\ctrh_\prj)\simeq
\sD^\bctr(X\ctrh_\al)$. \par
\textup{(b)} The hereditary complete cotorsion pair
\textup{(}$\Com(X\ctrh^\lct_\prj)$,
$\Acycl^\bctr(X\lcth^\lct_\bW)$\textup{)} in the exact category\/
$\Com(X\lcth^\lct_\bW)$ restricts to a hereditary complete cotorsion
pair \textup{(}$\Com(X\ctrh^\lct_\prj)$,
$\Acycl^\bctr(X\ctrh^\lct_\al)$\textup{)}in the exact subcategory\/
$\Com(X\ctrh^\lct_\al)\sub\Com(X\lcth^\lct_\bW)$.
 Consequently, the composition of triangulated functors
$$
 \Hot(X\ctrh^\lct_\prj)\lrarrow\Hot(X\ctrh^\lct_\al)\lrarrow
 \sD^\bctr(X\ctrh^\lct_\al)
$$
is a triangulated equivalence\/ $\Hot(X\ctrh^\lct_\prj)\simeq
\sD^\bctr(X\ctrh^\lct_\al)$. \par
\textup{(c)} The hereditary complete cotorsion pair
\textup{(}$\Com(X\ctrh_\prj)$, $\Acycl^\bctr(X\lcth_\bW)$\textup{)}
in the exact category\/ $\Com(X\lcth_\bW)$ restricts to a hereditary
complete cotorsion pair \textup{(}$\Com(X\ctrh_\prj)$,
$\Acycl^\bctr(X\ctrh_\alf)$\textup{)} in the exact subcategory\/
$\Com(X\ctrh_\alf)\sub\Com(X\lcth_\bW)$.
 Consequently, the composition of triangulated functors
$$
 \Hot(X\ctrh_\prj)\lrarrow\Hot(X\ctrh_\alf)\lrarrow
 \sD^\bctr(X\ctrh_\alf)
$$
is a triangulated equivalence\/ $\Hot(X\ctrh_\prj)\simeq
\sD^\bctr(X\ctrh_\alf)$.
\end{cor}

\begin{proof}
 In part~(a), the first assertion follows from
Theorem~\ref{becker-contraderived-cotorsion-pair}(a) in view of
Lemmas~\ref{restricting-hereditary-cotorsion}
and~\ref{restricting-cotorsion-pairs-lemma}(a).
 One should keep in mind that $X\ctrh_\al$ is a resolving subcategory
in $X\lcth_\bW$ by Corollary~\ref{clp-characterizations}(b).
 The second assertion is then provable similarly to the proof of
Corollary~\ref{becker-contraderived-of-lcta-lct-well-behaved}(a).

 The proof of part~(b) is similar and based on
Theorem~\ref{becker-contraderived-cotorsion-pair}(b).
 Alternatively, for the second assertions of parts~(a\+b) one can refer
to Corollary~\ref{becker-contraderived-of-lcta-lct-well-behaved}
together with Corollaries~\ref{ctrh-lcth-cor}(a)
and~\ref{lct-ctrh-lcth-cor}(a).
 The proof of part~(c) is also similar and based on
Theorem~\ref{becker-contraderived-cotorsion-pair}(a).
 One should keep in mind that $X\ctrh_\alf$ is a resolving subcategory
in $X\lcth_\bW$ by Corollary~\ref{clf-characterizations}(b).
\end{proof}

\begin{lem} \label{acyclic-in-alf-all-of-al-lct-cotorsion-pair}
 The pair of classes of acyclic complexes in the exact category of
antilocally flat contraherent cosheaves\/ $\sF=\Acycl(X\ctrh_\alf)$ and
all complexes of antilocal locally cotorsion contraherent cosheaves\/
$\sC=\Com(X\ctrh_\al^\lct)$ is a hereditary complete cotorsion pair in
the exact category\/ $\Com(X\ctrh_\al)$.
\end{lem}

\begin{proof}
 It suffices to restrict the cotorsion pair
of Corollary~\ref{acyclic-in-alf-all-of-lct-loc-contraherent-pair}
to the exact subcategory $\Com(X\ctrh_\al)\sub\Com(X\lcth_\bW)$,
using Lemmas~\ref{restricting-hereditary-cotorsion}
and~\ref{restricting-cotorsion-pairs-lemma}(a).

 Alternatively, one can restrict the cotorsion pair of
Lemma~\ref{acyclic-of-vfl-flat-arbitrary-of-cta-cot-pairs}(b)
to the exact subcategory $\Com(X\qcoh^\cta)\sub\Com(X\qcoh)$,
using Lemmas~\ref{restricting-hereditary-cotorsion}
and~\ref{restricting-cotorsion-pairs-lemma}(b).
 This produces a hereditary complete cotorsion pair
$\sF=\Acycl(X\qcoh_\fl^\cta)$ and $\sC=\Com(X\qcoh^\cot)$
in the exact category $\sE=\Com(X\qcoh^\cta)$.
 Then it remains to transfer the latter cotorsion pair from the world
of quasi-coherent sheaves into the realm of contraherent cosheaves
using the equivalence of exact categories from
Lemma~\ref{cta-clp-equivalence} together with
Lemmas~\ref{cta-clp-restricts-to-cot-inj}(a)
and~\ref{cta-clp-restricts-to-prj-clf}(c).
\end{proof}

\begin{cor} \label{acyclic-in-alf-all-of-lct-prj-cotorsion-pair}
 The pair of classes of acyclic complexes in the exact category of
antilocally flat contraherent cosheaves\/ $\sF=\Acycl(X\ctrh_\alf)$ and
all complexes of projective locally cotorsion contraherent cosheaves\/
$\sC=\Com(X\ctrh^\lct_\prj)$ is a hereditary complete cotorsion pair
in the exact category\/ $\Com(X\ctrh_\alf)$.
\end{cor}

\begin{proof}
 Restrict the cotorsion pair of
Lemma~\ref{acyclic-in-alf-all-of-al-lct-cotorsion-pair} to the exact
subcategory $\Com(X\ctrh_\alf)\sub\Com(X\ctrh_\al)$ using
Lemmas~\ref{restricting-hereditary-cotorsion}
and~\ref{restricting-cotorsion-pairs-lemma}(a), and recall that
$X\ctrh^\lct_\prj=X\ctrh_\alf\cap X\ctrh^\lct_\al=
X\ctrh_\alf\cap X\lcth^\lct_\bW$
(see Section~\ref{projective-contraherent}). \hbadness=1375
\end{proof}

\begin{cor} \label{acycl=bctracycl=bcoacycl-in-alf}
 In the exact category $X\ctrh_\alf$ of antilocally flat contraherent
cosheaves on $X$, the three classes of acyclic, Becker-coacyclic, and
Becker-contraacyclic complexes coincide.
\end{cor}

\begin{proof}
 This is a globalized version of~\cite[first assertions of
Theorems~7.14 and~7.18]{Pphil}.
 By Corollary~\ref{acyclic=bctraacyclic-in-flat-cta}, we have
$\Acycl(X\qcoh_\fl^\cta)=\Acycl^\bctr(X\qcoh_\fl^\cta)$, which in view
of the equivalence of exact categories from
Lemma~\ref{cta-clp-restricts-to-prj-clf}(c) translates into
$\Acycl(X\ctrh_\alf)=\Acycl^\bctr(X\ctrh_\alf)$.

 On the other hand, the injective objects of the exact category
$X\ctrh_\alf$ are precisely the projective locally cotorsion
contraherent cosheaves (see Section~\ref{projective-contraherent}).
 So the equality $\Acycl(X\ctrh_\alf)=\Acycl^\bco(X\ctrh_\alf)$
follows immediately from
Corollary~\ref{acyclic-in-alf-all-of-lct-prj-cotorsion-pair}
in view of Lemma~\ref{Ext-1-as-homotopy-Hom}.
\end{proof}

\begin{cor} \label{becker-coderived-of-alf-well-behaved}
 The composition of triangulated functors
$$
 \Hot(X\ctrh^\lct_\prj)\lrarrow\Hot(X\ctrh_\alf)\lrarrow
 \sD^\bco(X\ctrh_\alf)
$$
is a triangulated equivalence\/
$$
 \Hot(X\ctrh^\lct_\prj)\simeq\sD^\bco(X\ctrh_\alf)=
 \sD(X\ctrh_\alf)=\sD^\bctr(X\ctrh_\alf).
$$
\end{cor}

\begin{proof}
 Follows from
Corollary~\ref{acyclic-in-alf-all-of-lct-prj-cotorsion-pair} in view of
Corollary~\ref{acycl=bctracycl=bcoacycl-in-alf} by the standard
argument similar to the proofof~\cite[Corollary~9.5]{PS4} or
Corollary~\ref{becker-coderived-of-cta-cot}, and dual-analogous to
the proof of
Corollary~\ref{becker-contraderived-of-lcta-lct-well-behaved}.
\end{proof}

 The following lemma will be useful in
Sections~\ref{qc-ss-homotopy-projective-subsect}
and~\ref{antilocal-semicontraderived-subsect}.

\begin{lem} \label{all-alf-becker-contraacyclic-lct-cotorsion-pair}
 The pair of classes of all complexes of antilocally flat contraherent
cosheaves\/ $\sF=\Com(X\ctrh_\alf)$ and Becker-contraacyclic complexes
of locally cotorsion\/ $\bW$\+locally contraherent cosheaves\/
$\sC=\Acycl^\bctr(X\lcth_\bW^\lct)$ is a hereditary complete
cotorsion pair in the exact category\/ $\Com(X\lcth_\bW)$.
\end{lem}

\begin{proof}
 The argument is similar to the proof of
Theorem~\ref{becker-contraderived-cotorsion-pair} and
based on Theorem~\ref{loc-contraherent-gluing-theorem}.
 Let $\sR$ be the local class of all commutative rings $R$ and
$\sE^R=\Com(R\modl^\cta)$ be the exact category of complexes of
contraadjusted $R$\+modules.
 Let $\sF(R)=\Com(R\modl^\cta_\fl)\sub\sE^R$ be the class of all
complexes of flat contraadjusted $R$\+modules, and let
$\sC^R=\Acycl^\bctr(R\modl^\cot)$ be the class of Becker-contraacyclic
complexes (in the exact category) of cotorsion $R$\+modules.
 Then $\sE$ is a very colocal class by
Examples~\ref{colocal-classes-examples}, while $\sC$ is a very colocal
class by Theorem~\ref{Becker-contraacyclicity-local-on-qcomp-qsep}(b)
(for the affine scheme $\Spec R$) and
Lemma~\ref{morphism-of-affines-direct-image-Becker-contra}(b).

 In order to show that the pair of classes $(\sF(R),\sC^R)$ is
a hereditary complete cotorsion pair in $\Com(R\modl^\cta)$, put
$U=\Spec R$ and consider the hereditary complete cotorsion pair
$(\sF'',\sC'')$ from
Lemma~\ref{all-vfl-flat-contraacyclic-of-cta-cot-pairs}(b)
in the exact category $\Com(U\qcoh)=\Com(R\modl)$.
 So $\sF''=\Com(R\modl_\fl)$ and $\sC''=\Acycl^\bctr(X\qcoh^\cot)$.
 By Lemmas~\ref{restricting-hereditary-cotorsion}
and~\ref{restricting-cotorsion-pairs-lemma}(b), the cotorsion pair
$(\sF'',\sC'')$ restricts to a hereditary complete cotorsion pair in
the exact subcategory $\Com(R\modl^\cta)\sub\Com(R\modl)$,
providing the desired cotorsion pair $(\sF(R),\sC^R)$ in~$\sE^R$.

 Applying Theorem~\ref{loc-contraherent-gluing-theorem}, we obtain
a hereditary complete cotorsion pair $(\sF(X),\sC^X_\bW)$ in
the exact category $\sE^X_\bW=\Com(X\lcth_\bW)$.
 Here $\sF(X)$ is the class of all direct summands of finitely
iterated extensions of direct images of complexes of antilocally
flat contraherent cosheaves from $U_\alpha$, while $\sC^X_\bW$ is
the class of all locally\+$\sC$ complexes of $\bW$\+locally
contraherent cosheaves on~$X$.
 By Corollary~\ref{complexes-of-prj-alf-antilocal}(c), \,$\sF(X)$ is
the class of all complexes of antilocally flat contraherent cosheaves
on~$X$.
 By Theorem~\ref{Becker-contraacyclicity-local-on-qcomp-qsep}(b),
we have $\sC^X_\bW=\Acycl^\bctr(X\lcth_\bW^\lct)$. \hbadness=1625
\end{proof}

\begin{thm} \label{becker-contraderived-lcta-lct-equivalent}
 For any quasi-compact semi-separated scheme $X$ with an open covering\/
$\bW$, the embedding of exact categories $X\lcth^\lct_\bW\rarrow
X\lcth_\bW$ induces an equivalence of the Becker contraderived
categories
$$
 \sD^\bctr(X\lcth^\lct_\bW)\simeq\sD^\bctr(X\lcth_\bW).
$$
\end{thm}

\begin{proof}
 This is a globalized version of~\cite[Corollaries~7.21
and~7.30]{Pphil}.
 Notice first of all that all Becker-contraacyclic complexes in
$X\lcth^\lct_\bW$ are also Becker-contraacyclic in $X\lcth_\bW$
by Corollary~\ref{lct-Becker-contraacyclic-iff-as-lcta}.
 Therefore, the triangulated functor
$$
 \sD^\bctr(X\lcth^\lct_\bW)\lrarrow\sD^\bctr(X\lcth_\bW)
$$
induced by the embedding of exact categories $X\lcth^\lct_\bW\rarrow
X\lcth_\bW$ is well-defined.

 The assertion that this functor is a category equivalence is now
obtained by comparing the results of
Corollaries~\ref{becker-contraderived-of-lcta-lct-well-behaved}(a\+b),
\ref{becker-contraderived-of-al-lct-alf-well-behaved}(c),
and~\ref{becker-coderived-of-alf-well-behaved}.
 To spell out this argument, consider the commutative diagram of
triangulated functors induced by embeddings of additive/exact categories
\begin{equation} \label{becker-contrader-lcta-lct-equivalent-diagram}
\begin{gathered}
 \xymatrix{
  \Hot(X\ctrh^\lct_\prj) \ar@<2pt>[r] \ar@<-2pt>@{-}[r]
  \ar@<-2pt>[dd] \ar@<2pt>@{-}[dd]
  & \sD^{\bctr=\empt=\bco}(X\ctrh_\alf)
  \ar@<2pt>[rdd] \ar@<-2pt>@{--}[rdd]
  & \Hot(X\ctrh_\prj) \ar@<-2pt>[l] \ar@<2pt>@{-}[l]
  \ar@<2pt>[dd] \ar@<-2pt>@{-}[dd] \\ \\
  \sD^\bctr(X\lcth^\lct_\bW) \ar@<-2pt>[rr] \ar@<2pt>@{--}[rr]
  && \sD^\bctr(X\lcth_\bW)
 }
\end{gathered}
\end{equation}
 The two upper horizontal arrows are triangulated equivalences by
Corollaries~\ref{becker-coderived-of-alf-well-behaved}
and~\ref{becker-contraderived-of-al-lct-alf-well-behaved}(c), while
the two vertical arrows are triangulated equivalences by
Corollary~\ref{becker-contraderived-of-lcta-lct-well-behaved}(a\+b).
 It follows that the diagonal arrow and the lower horizonal arrow
also triangulated equivalences.
\end{proof}

\begin{cor} \label{becker-co-contraderived-of-fl-cta}
 In the exact category $X\qcoh_\fl^\cta$ of flat contraadjusted
quasi-coherent sheaves on $X$, the three classes of acyclic,
Becker-coacyclic, and Becker-contraacyclic complexes coincide,
$$
 \Acycl^\bco(X\qcoh_\fl^\cta)=\Acycl(X\qcoh_\fl^\cta)=
 \Acycl^\bctr(X\qcoh_\fl^\cta).
$$
 The composition of triangulated functors
$$
 \Hot(X\qcoh_\vfl^\cta)\lrarrow
 \Hot(X\qcoh_\fl^\cta)\lrarrow
 \sD^\bctr(X\qcoh_\fl^\cta)
$$
is a triangulated equivalence
$$
 \Hot(X\qcoh_\vfl^\cta)\simeq\sD^\bctr(X\qcoh_\fl^\cta)=
 \sD(X\qcoh_\fl^\cta)=\sD^\bco(X\qcoh_\fl^\cta).
$$
 The composition of triangulated functors
$$
 \Hot(X\qcoh_\fl^\cot)\lrarrow
 \Hot(X\qcoh_\fl^\cta)\lrarrow
 \sD^\bco(X\qcoh_\fl^\cta)
$$
is also a triangulated equivalence
$$
 \Hot(X\qcoh_\fl^\cot)\simeq\sD^\bco(X\qcoh_\fl^\cta)=
 \sD(X\qcoh_\fl^\cta)=\sD^\bctr(X\qcoh_\fl^\cta).
$$
\end{cor}

\begin{proof}
 The equality $\Acycl(X\qcoh_\fl^\cta)=\Acycl^\bctr(X\qcoh_\fl^\cta)$
is the result of Corollary~\ref{acyclic=bctraacyclic-in-flat-cta}.
 All the assertions of the corollary are equivalent restatements of
Corollaries~\ref{becker-contraderived-of-al-lct-alf-well-behaved}(c),
\ref{acycl=bctracycl=bcoacycl-in-alf},
and~\ref{becker-coderived-of-alf-well-behaved}, in view
of Lemma~\ref{cta-clp-restricts-to-prj-clf}.
\end{proof}

\begin{lem}  \label{vfl-cta-finite-dim}
 Let $X=\bigcup_{\alpha=1}^N U_\alpha$ be a finite affine open covering.
 Then \par
\textup{(a)} the coresolution dimension of any very flat quasi-coherent
sheaf on $X$ with respect to the coresolving subcategory of
contraadjusted very flat quasi-coherent sheaves $X\qcoh_\vfl^\cta\sub
X\qcoh_\vfl$ does not exceed~$N$; \par
\textup{(b)} the coresolution dimension of any flat quasi-coherent sheaf
on $X$ with respect to the coresolving subcategory of contraadjusted
flat quasi-coherent sheaves $X\qcoh_\fl^\cta\sub X\qcoh_\fl$ does not
exceed~$N$; \par
\textup{(c)} the homological dimension of the exact category
$X\qcoh_\vfl$ does not exceed~$N$.
\end{lem}

\begin{proof}
 In parts~(a\+b), the two full subcategories in question are
coresolving in their respective ambient exact categories due to
Corollary~\ref{quasi-very-cta-cor}(b), so it remains to apply
Lemma~\ref{dil-cta-clp-finite-dim}(b) and the dual version of
Corollary~\ref{fdim-subcategory-cor}.
 Part~(c) follows immediately from part~(a), because $X\qcoh_\vfl^\cta$
is the full subcategory of injective objects in the exact category
$X\qcoh_\vfl$.

 Alternatively, for part~(c) one can prove the seemingly stronger
(but actually equivalent) assertion that $\Ext_X^{>N}(\F,\M)=0$ for all
very flat quasi-coherent sheaves $\F$ and all quasi-coherent sheaves
$\M$ on~$X$ (cf.\ Lemma~\ref{finite-homol-dim-finite-coresol-dim}).
 This assertion follows from Lemma~\ref{dil-cta-clp-finite-dim}(b).
 The same equivalent version of part~(c) is also a particular case
of~\cite[Theorem~6.3(b)]{PS6}.
 (Cf.\ Lemma~\ref{lct-lin-clp-finite-dim}(c).)
\end{proof}

 In the following theorem and corollary, we use the notation
$\Hot^\st(\sE)$ from
Section~\ref{derived-of-sheaves-and-cosheaves-subsect}.

\begin{thm} \label{vlf-cta-fl-cot-derived-vlf-fl-equivalences}
\textup{(a)} For any symbol\/ $\bst=\b$, $+$, $-$, $\empt$, $\abs+$,
$\abs-$, $\bco$, $\co$, or\/~$\abs$, the triangulated functor\/
$\Hot^\st(X\qcoh_\vfl^\cta)\rarrow\sD^\st(X\qcoh_\vfl)$ induced by
the embedding of additive/exact categories $X\qcoh_\vfl^\cta\rarrow
X\qcoh_\vfl$ is an equivalence of triangulated categories,
$$
 \Hot^\st(X\qcoh_\vfl^\cta)\simeq\sD^\st(X\qcoh_\vfl).
$$
 In particular, the classes of acyclic, coacyclic, Becker-coacyclic,
and absolutely acyclic complexes coincide in the exact category
$X\qcoh_\vfl$. \par
\textup{(b)} For any symbol\/ $\bst=\b$, $+$, $-$, $\empt$, $\abs+$,
$\abs-$, $\bco$, or\/~$\abs$, the triangulated functor\/
$\sD^\st(X\qcoh_\fl^\cta)\rarrow\sD^\st(X\qcoh_\fl)$ induced by
the embedding of additive/exact categories $X\qcoh_\fl^\cta\rarrow
X\qcoh_\fl$ is an equivalence of triangulated categories,
$$
 \sD^\st(X\qcoh_\fl^\cta)\simeq\sD^\st(X\qcoh_\fl).
$$ \par
\textup{(c)} For any symbol\/ $\bst=+$, $\bco$, or\/~$\empt$,
the triangulated functor\/ $\Hot^\st(X\qcoh_\fl^\cot)\rarrow
\sD^\st(X\qcoh_\fl)$ induced by the embedding of additive/exact
categories $X\qcoh_\fl^\cot\rarrow X\qcoh_\fl$ is an equivalence
of triangulated categories,
$$
 \Hot^\st(X\qcoh_\fl^\cot)\simeq\sD^\st(X\qcoh_\fl).
$$
 In particular, the classes of acyclic and Becker-coacyclic complexes
coincide in the exact category $X\qcoh_\fl$,
$$
 \Acycl^\bco(X\qcoh_\fl)=\Acycl(X\qcoh_\fl).
$$
\end{thm}

\begin{proof}
 Part~(a): assuming $\bst\ne\co$, the assertion follows from
Lemma~\ref{vfl-cta-finite-dim}(a) in view of the dual versions of
Propositions~\ref{finite-resolutions}
and~\ref{becker-contraderived-finite-resolutions}.
 In particular, we have proved that $\sD(X\qcoh_\vfl)=
\sD^\abs(X\qcoh_\vfl)$, so (a)~holds for $\bst=\co$ as well.
 Alternatively, one can refer to Lemma~\ref{vfl-cta-finite-dim}(c),
Theorem~\ref{finite-homol-dim-becker-co-contra-derived}(a),
and Lemma~\ref{psemi-remark21} for a proof of all assertions
of part~(a) except the one about $\bst=\co$.
 The proof of part~(b) is similar and based on
Lemma~\ref{vfl-cta-finite-dim}(b).
 Part~(c) for $\bst=+$ follows from
Corollaries~\ref{quasi-cotors-characterizations}(c)
and~\ref{quasi-cotors-cor}(b) in view of the dual version of
Proposition~\ref{infinite-resolutions}(a).

 For a slight/partial generalization of the second assertion of
part~(a), see
Corollary~\ref{vfl-lin-finite-dim-all-derived-coincide}(a) below.

 The proof of part~(c) for $\bst=\bco$ and~$\empt$ is based on
Lemma~\ref{acyclic-of-vfl-flat-arbitrary-of-cta-cot-pairs}(b).
 Recall that the flat cotorsion quasi-coherent sheaves are
the injective objects of the category $X\qcoh_\fl$.
 So, first of all, it follows from
Lemma~\ref{acyclic-of-vfl-flat-arbitrary-of-cta-cot-pairs}(b)
together with Lemma~\ref{Ext-1-as-homotopy-Hom} that all acyclic
complexes in the exact category $X\qcoh_\fl$ are Becker-coacyclic.
 Now it remains to construct, for any complex of flat
quasi-coherent sheaves $\F^\bu$ on $X$, a complex of flat
cotorsion quasi-coherent sheaves $\cP^\bu$ together with
a morphism of complexes $\F^\bu\rarrow\cP^\bu$ with the cone
in $\Acycl(X\qcoh_\fl)$.

 Indeed, let $0\rarrow\F^\bu\rarrow\cP^\bu\rarrow\G^\bu\rarrow0$
be a special preenvelope short exact sequence in the cotorsion pair of
Lemma~\ref{acyclic-of-vfl-flat-arbitrary-of-cta-cot-pairs}(b).
 So $\cP^\bu$ is a complex in $X\qcoh^\cot$ and $\G^\bu$ is
an acyclic complex in the exact category $X\qcoh_\fl$.
 Then both $\F^\bu$ and $\G^\bu$ are complexes of flat quasi-coherent
sheaves, hence so is~$\cP^\bu$.
 It remains to observe that the total complex of the short exact
sequence $0\rarrow\F^\bu\rarrow\cP^\bu\rarrow\G^\bu\rarrow0$ is
absolutely acyclic (hence acyclic) in $X\qcoh_\fl$.
 As $\G^\bu$ is acyclic in $X\qcoh_\fl$, it follows that the cone of
of the morphism $\F^\bu\rarrow\cP^\bu$ is acyclic in $X\qcoh_\fl$.
\end{proof}

\begin{cor}  \label{vfl-co-contra-cor-expanded}
\textup{(a)} For any symbol\/ $\bst=\b$, $+$, $-$, $\empt$, $\abs+$,
$\abs-$, $\bco$, $\co$, or\/~$\abs$, there is a natural equivalence of
triangulated categories\/ $\sD^\st(X\qcoh_\vfl)\simeq
\Hot^\st(X\ctrh_\prj)$. \par
\textup{(b)} For any symbol\/ $\bst=\b$, $+$, $-$, $\empt$, $\abs+$,
$\abs-$, $\bco$, or\/~$\abs$, there is a natural equivalence of
triangulated categories\/ $\sD^\st(X\qcoh_\fl)\simeq
\sD^\st(X\ctrh_\alf)$. \par
\textup{(c)} For any symbol\/ $\bst=+$, $\bco$, or\/~$\empt$, there is
a natural equivalence of triangulated categories\/
$\sD^\st(X\qcoh_\fl)\simeq\Hot^\st(X\ctrh^\lct_\prj)$.
\end{cor}

\begin{proof}
 Compare Theorem~\ref{vlf-cta-fl-cot-derived-vlf-fl-equivalences}
with Lemma~\ref{cta-clp-restricts-to-prj-clf}.
\end{proof}

 The following lemma comes from~\cite[Proposition~4.4]{ES}, where it is
attributed to \v St\!'ov\'\i\v cek.

\begin{lem} \label{veryflat-flat-periodicity}
 Let $R$ be a commutative ring.
 If all the modules of cocycles in an acyclic complex of very flat
$R$\+modules $F^\bu$ are flat, then the modules of cocycles are
actually very flat.
\end{lem}

\begin{proof}
 The claim is that the triangulated functor $\sD(R\modl_\vfl)\rarrow
\sD(R\modl_\fl)$ induced by the inclusion of exact categories
$R\modl_\vfl\rarrow R\modl_\fl$ has zero kernel (i.~e., takes nonzero
objects to nonzero objects).
 Indeed, the similar functor $\sD(R\modl_\prj)\rarrow\sD(R\modl_\vfl)$
is a triangulated equivalence by
Proposition~\ref{finite-resolutions}.
 Hence it suffices to prove that the triangulated functor
$\sD(R\modl_\prj)\rarrow\sD(R\modl_\fl)$ has zero kernel.
 This is the result of Theorem~\ref{flat-projective-periodicity}(a).
 (In fact, Theorem~\ref{flat-projective-periodicity}(b)
with Proposition~\ref{flat-projective-periodicity-complements}(b)
tell us that the functor $\sD(R\modl_\prj)\rarrow\sD(R\modl_\fl)$ is
a triangulated equivalence for any associative ring~$R$.)
\end{proof}

 The nontrivial ($\bst=\empt$) case of the following theorem is due to
Estrada and Sl\'avik~\cite[Corollary~6.2]{ES}.

\begin{thm} \label{derived-vfl-fl-equivalence}
 Let $X$ be a quasi-compact semi-separated scheme.
 Then, for any symbol\/ $\bst=-$ or\/~$\empt$, the triangulated
functor\/ $\sD^\st(X\qcoh_\vfl)\rarrow\sD^\st(X\qcoh_\fl)$ induced by
the embedding of exact categories $X\qcoh_\vfl\rarrow X\qcoh_\fl$ is
an equivalence of triangulated categories,
$$
 \sD^\st(X\qcoh_\vfl)\simeq\sD^\st(X\qcoh_\fl).
$$
\end{thm}

\begin{proof}
 The easy case $\bst=-$ follows from
Lemma~\ref{quasi-very-flat-cover} and
Proposition~\ref{infinite-resolutions}(a).
 The interesting case is $\bst=\empt$.
 The following argument is a differently worded version of the proof
in~\cite{ES}.
 Notice that, by Lemma~\ref{veryflat-flat-periodicity}, any acyclic
complex in the exact category $X\qcoh_\fl$ with the terms belonging
to $X\qcoh_\vfl$ is acyclic in the exact category $X\qcoh_\vfl$.
 In view of Lemma~\ref{pkoszul-lemma16}(a), it remains to construct,
for any given complex of flat quasi-coherent sheaves $\G^\bu$,
a complex of very flat quasi-coherent sheaves $\F^\bu$ on $X$ together
with a morphism of complex $\F^\bu\rarrow\G^\bu$ whose cone is acyclic
in the exact category $X\qcoh_\fl$.

 For this purpose, one can use a special precover short exact sequence
$0\rarrow\cH^\bu\rarrow\F^\bu\rarrow\G^\bu\rarrow0$ of the complex
$\G^\bu$ in the cotorsion pair of
Lemma~\ref{all-vfl-flat-contraacyclic-of-cta-cot-pairs}(a).
 So $\F^\bu$ is a complex of very flat quasi-coherent sheaves and
$\cH^\bu$ is a Becker-contraacyclic complex of contraadjusted
quasi-coherent sheaves.
 Since both $\G^\bu$ and $\F^\bu$ are complexes of flat quasi-coherent
sheaves, so is~$\cH^\bu$.
 Now $\cH^\bu$ is a Becker-contraacyclic complex of flat contraadjusted
quasi-coherent sheaves, and
Corollary~\ref{acyclic=bctraacyclic-in-flat-cta} tells us that it is
a complex of flat quasi-coherent sheaves with flat sheaves of cocycles.
 As the totalization of the short exact sequence
$0\rarrow\cH^\bu\rarrow\F^\bu\rarrow\G^\bu\rarrow0$ is absolutely
acyclic in $X\qcoh_\fl$ and the complex $\cH^\bu$ is acyclic in
$X\qcoh_\fl$, we can conclude that the cone of the morphism
$\F^\bu\rarrow\G^\bu$ is acyclic in $X\qcoh_\fl$.
\end{proof}

 Taken together, the results of
Corollary~\ref{becker-co-contraderived-of-fl-cta},
Theorem~\ref{vlf-cta-fl-cot-derived-vlf-fl-equivalences}(a,c)
(for $\bst=\empt$), and Theorem~\ref{derived-vfl-fl-equivalence} can be
described by the following commutative diagram of triangulated
equivalences induced by embeddings of exact/additive categories:
\begin{equation} \label{flat-vfl-derived-equivalences-diagram}
\begin{gathered}
 \xymatrix{
  \Hot(X\qcoh_\fl^\cot) \ar@<2pt>[r] \ar@<-2pt>@{-}[r]
  \ar@<-2pt>[dd] \ar@<2pt>@{-}[dd]
  & \sD^{\bctr=\empt=\bco}(X\qcoh_\fl^\cta)
  \ar@<-2pt>[ldd] \ar@<2pt>@{-}[ldd]
  & \Hot(X\qcoh_\vfl^\cta) \ar@<-2pt>[l] \ar@<2pt>@{-}[l]
  \ar@<2pt>[dd] \ar@<-2pt>@{-}[dd] \\ \\
  \sD(X\qcoh_\fl) && \sD(X\qcoh_\vfl)
  \ar@<2pt>[ll] \ar@<-2pt>@{-}[ll]
 }
\end{gathered}
\end{equation}

\begin{cor} \label{bctr-lcth-vfl-fl-derived-equivalences}
 For any quasi-compact semi-separated scheme $X$ with an open covering\/
$\bW$, there is a commutative square diagram of triangulated
equivalences
\begin{equation} \label{derived-vfl-fl-qcoh-becker-contrader-lcth-equiv}
\begin{gathered}
 \xymatrix{
  \sD(X\qcoh_\fl) \ar@{=}[d] & \sD(X\qcoh_\vfl) \ar@{=}[d]
  \ar@<-2pt>[l] \ar@<2pt>@{-}[l] \\
  \sD^\bctr(X\lcth_\bW^\lct) \ar@<-2pt>[r] \ar@<2pt>@{-}[r]
  & \sD^\bctr(X\lcth_\bW)
 }
\end{gathered}
\end{equation}
 Here the horizontal arrows are induced by the respective embeddings
of exact categories.
 The vertical equivalences are obtained by comparing
the commutative
diagrams~\eqref{becker-contrader-lcta-lct-equivalent-diagram}
and~\eqref{flat-vfl-derived-equivalences-diagram} and identifying
the upper lines of these two diagrams using the exact/additive category
equivalences of Lemma~\textup{\ref{cta-clp-restricts-to-prj-clf}}.
\qed
\end{cor}

\subsection{Homotopy projective complexes}
\label{qc-ss-homotopy-projective-subsect}
 In this section we apply the techniques developed above in this
Chapter~\ref{becker-on-qcomp-qsep-sect} in order to study
the conventional derived category of contraherent cosheaves
(which is the subject of the previous
Chapter~\ref{derived-on-quasi-compact-sect}).

 Let $X$ be a scheme with an open covering~$\bW$.
 A complex $\J^\bu$ of quasi-coherent sheaves on $X$ is called
\emph{homotopy injective} if the complex of abelian groups
$\Hom_X(\M^\bu,\J^\bu)$ is acyclic for any acyclic complex of
quasi-coherent sheaves $\M^\bu$ on~$X$.
 The full subcategory of homotopy injective complexes in $\Hot(X\qcoh)$
is denoted by $\Hot(X\qcoh)^\inj$ and the full subcategory of
complexes of injective quasi-coherent sheaves that are also homotopy
injective is denoted by $\Hot(X\qcoh^\inj)^\inj\sub\Hot(X\qcoh^\inj)$.
{\hbadness=1575\par}

 Similarly, a complex $\gF^\bu$ of $\bW$\+locally contraherent cosheaves
on $X$ is called \emph{homotopy projective} if the complex of abelian
groups $\Hom^X(\gF^\bu,\gM^\bu)$ is acyclic for any acyclic complex
$\gM^\bu$ over the exact category $X\lcth_\bW$.
 The full subcategory of homotopy projective complexes in
$\Hot(X\lcth_\bW)$ is denoted by $\Hot(X\lcth_\bW)_\prj$.
 We will see below in this section that the property of a complex of
$\bW$\+locally contraherent cosheaves on a quasi-compact semi-separated
scheme $X$ to be homotopy projective does not change when the covering
$\bW$ is replaced by its refinement.

 Finally, a complex $\P^\bu$ of locally cotorsion $\bW$\+locally
contraherent cosheaves on $X$ is called \emph{homotopy projective} if
the complex of abelian groups $\Hom^X(\P^\bu,\gM^\bu)$ is acyclic for
any acyclic complex $\gM^\bu$ over the exact category $X\lcth_\bW^\lct$.
 The full subcategory of homotopy projective complexes in
$\Hot(X\lcth_\bW^\lct)$ is denoted by $\Hot(X\lcth_\bW^\lct)_\prj$.
 Let us issue a \emph{warning} that our terminology is misleading:
a homotopy projective complex of locally cotorsion $\bW$\+locally
contraherent cosheaves need not be homotopy projective as a complex
of $\bW$\+locally contraherent cosheaves.
 It will be shown below that the property of a complex of locally
cotorsion $\bW$\+locally contraherent cosheaves on a quasi-compact
semi-separated scheme $X$ to be homotopy projective does not change
when the covering $\bW$ is replaced by its refinement.

\begin{lem}  \label{qcomp-ssep-homotopy-proj-independence}
 Let $X$ be a quasi-compact semi-separated scheme with
an open covering\/~$\bW$.  Then \par
\textup{(a)} a complex\/ $\gF^\bu$ over $X\ctrh_\prj$ belongs to\/
$\Hot(X\lcth_\bW)_\prj$ if and only if the complex\/
$\Hom^X(\gF^\bu,\gE^\bu)$ is acyclic for any complex\/ $\gE^\bu$ over
$X\ctrh_\prj$ acyclic with respect to $X\ctrh_\al$; \par
\textup{(b)} a complex\/ $\P^\bu$ over $X\ctrh^\lct_\prj$ belongs to\/
$\Hot(X\lcth_\bW^\lct)_\prj$ if and only if the complex\/
$\Hom^X(\P^\bu,\gE^\bu)$ is acyclic for any complex\/ $\gE^\bu$ over
$X\ctrh^\lct_\prj$ acyclic with respect to $X\ctrh^\lct_\al$.
\end{lem}

\begin{proof}
 We will prove part~(a), part~(b) being similar.
 The ``only if'' assertion holds by the definition.
 To check the ``if'', consider a complex $\gM^\bu$ over $X\lcth_\bW$.
 By (the proof of)
Corollary~\ref{becker-contraderived-of-lcta-lct-well-behaved}(a),
there exists a complex $\gE^\bu$ over $X\ctrh_\prj$ together with
a morphism of complexes of locally contraherent cosheaves
$\gE^\bu\rarrow\gM^\bu$ with a cone Becker-contraacyclic
in $X\lcth_\bW$.
 Moreover, the complex $\Hom^X$ from any complex of projective 
contraherent cosheaves to a Becker-contraacyclic complex in $X\lcth_\bW$
is acyclic by the definition.
 Hence the morphism $\Hom^X(\gF^\bu,\gE^\bu)\rarrow\Hom^X(\gF^\bu,
\gM^\bu)$ is a quasi-isomorphism.
 Finally, if the complex $\gM^\bu$ is acyclic in $X\lcth_\bW$,
then so is the complex $\gE^\bu$ (by
Corollary~\ref{Becker-contraacyclic-are-acyclic}(a)).
 By Lemma~\ref{dil-cta-clp-finite-dim}(c) and
Corollary~\ref{fdim-acyclic-cor}, it follows that the complex
$\gE^\bu$ is also acyclic with respect to $X\ctrh_\al$.
\end{proof}

 According to Lemma~\ref{qcomp-ssep-homotopy-proj-independence},
the property of a complex over $X\ctrh_\prj$ (respectively, over
$X\ctrh^\lct_\prj$) to belong to $\Hot(X\lcth_\bW)_\prj$ (resp.,
$\Hot(X\lcth_\bW^\lct)_\prj$) does not depend on the covering~$\bW$
(for a quasi-compact semi-separated scheme).
 We will denote the full subcategory in $\Hot(X\ctrh_\prj)$ (resp.,
$\Hot(X\ctrh^\lct_\prj)$) consisting of the homotopy projective
complexes by $\Hot(X\ctrh_\prj)_\prj$ (resp.,
$\Hot(X\ctrh^\lct_\prj)_\prj$).
 It is a standard fact (provable using silly truncations and
Lemma~\ref{telescope} in the category of abelian groups) that
bounded above complexes of projectives are homotopy projective,
$\Hot^-(X\ctrh_\prj)\sub\Hot(X\ctrh_\prj)_\prj$ and
$\Hot^-(X\ctrh^\lct_\prj)\sub\Hot(X\ctrh^\lct_\prj)_\prj$.

 We will denote by $\Com(X\ctrh_\prj)_\prj\sub\Com(X\ctrh)$
and $\Com(X\ctrh^\lct_\prj)_\prj\sub\Com(X\ctrh^\lct)$ the classes of
all homotopy projective complexes of projective (locally contraadjusted
or locally cotorsion) contraherent cosheaves, viewed as full
subcategories in the category of complexes $\Com(X\ctrh)$.

\medskip

 We start with a reminder of the classical case of quasi-coherent
sheaves.
 The following result is well-known.

\begin{thm}  \label{quasi-homotopy-injective}
 For any scheme $X$, the natural functors\/ $\Hot(X\qcoh^\inj)^\inj
\rarrow\Hot(X\qcoh)^\inj\rarrow\sD(X\qcoh)$ are equivalences of
triangulated categories.
\end{thm}

\begin{proof}
 It is clear that both functors are fully faithful.
 This assertion that they are triangulated equivalences holds for any
Grothendieck abelian category $\sA$ in the role of $X\qcoh$.
 Denote by $\Hot(\sA)^\inj\sub\Hot(\sA)$ the full subcategory of
homotopy injective complexes in $\sA$ and by $\Hot(\sA^\inj)^\inj
\sub\Hot(\sA^\inj)$ the full subcategory of homotopy injective
complexes of injective objects.
 Then the functor $\Hot(\sA)^\inj\rarrow\sD(\sA)$ is an equilalence
of categories by~\cite[Theorem~5.4]{AJS} (for an even more
general statement, see~\cite[Theorem~6]{K-app}).
 The stronger result that the functor $\Hot(\sA^\inj)^\inj\rarrow
\sD(\sA)$ is an equivalence of categories can be found
in~\cite[Theorem~3.13 and Lemma~3.7(ii)]{Ser}, \cite[Example~3.2]{Hov},
\cite[Corollary~7.1]{Gil1}, or~\cite[Corollary~8.5]{PS4}.

 In order to show that any homotopy injective complex in $\sA$ is
homotopy equivalent to a homotopy injective complex of injective
objects, one can also argue as follows.
 By (the proof of)
Theorem~\ref{coderived-of-grothendieck-contraderived-of-lpacepo}(a)
(Theorem~\ref{quasi-coherent-becker-coderived} for quasi-coherent
sheaves), for any complex $M^\bu$ over $\sA$ there exists a closed
morphism $M^\bu\rarrow J^\bu$ into a complex $J^\bu$ over $A^\inj$
with a Becker-coacyclic cone~$E^\bu$.
 If the complex $M^\bu$ was homotopy injective, then the morphism
$E^\bu\rarrow M^\bu[1]$ is homotopic to zero (since $E^\bu$ is acyclic
by Lemma~\ref{with-co-kernels-Becker-trivial-are-acyclic}(a)), hence
the complex $E^\bu$ is a direct summand of $J^\bu$ in $\Hot(\sA)$.
 Any morphism $E^\bu\rarrow J^\bu$ being also homotopic to zero,
it follows that the complex $E^\bu$ is contractible and the morphism
$M^\bu\rarrow J^\bu$ is a homotopy equivalence.
\end{proof}

 Let $X$ be a quasi-compact semi-separated scheme with an open
covering~$\bW$.

\begin{thm} \label{lcth-homotopy-projective-cotorsion-pairs}
\textup{(a)} The pair of classes of all homotopy projective complexes
of projective contraherent cosheaves\/ $\sF=\Com(X\ctrh_\prj)_\prj$
and acyclic complexes of\/ $\bW$\+locally contraherent cosheaves\/
$\sC=\Acycl(X\lcth_\bW)$ is a hereditary complete cotorsion pair in
the exact category\/ $\Com(X\lcth_\bW)$. \par
\textup{(b)} The pair of classes of all homotopy projective complexes
of projective locally cotorsion contraherent cosheaves\/
$\sF=\Com(X\ctrh^\lct_\prj)_\prj$ and acyclic complexes of\/
locally cotorsion\/ $\bW$\+locally contraherent cosheaves\/
$\sC=\Acycl(X\lcth_\bW^\lct)$ is a hereditary complete cotorsion pair
in the exact category\/ $\Com(X\lcth_\bW^\lct)$.
\end{thm}

\begin{proof}
 Part~(a): let us start with the case of an affine scheme $U=\Spec R$
with the open covering~$\{U\}$.
 Lemma~\ref{homotopy-vfl-flat-acyclic-of-cta-cot-pairs}(a) provides
a hereditary complete cotorsion pair $(\sF',\sC')$ in the abelian
category $\Com(U\qcoh)=\Com(R\modl)$ such that $\sC'=
\Acycl(U\qcoh^\cta)=\Acycl(R\modl^\cta)$ is the class of acyclic
complexes of contraadjusted $R$\+modules.
 By Lemmas~\ref{restricting-hereditary-cotorsion}
and~\ref{restricting-cotorsion-pairs-lemma}(b), the cotorsion pair
$(\sF',\sC')$ restricts to a hereditary complete cotorsion pair
$(\sF,\sC)$ in the exact subcategory $\sE=\Com(R\modl^\cta)\sub
\Com(R\modl)$.
 So we have $\sF=\Com(R\modl^\cta)\cap\sF'=\Com(R\modl^\cta)\cap
\Com(R\modl_\vfl)_\vfl$ and $\sC=\sC'=\Acycl(R\modl^\cta)$.
 Now we identify the exact category $\Com(R\modl^\cta)$ with
$\Com(U\ctrh)$, and notice that
$\sF\sub\Com(R\modl^\cta_\vfl)=\Com(U\ctrh_\prj)$ (this follows also
from the equality $\sC=\Acycl(U\ctrh)$ by means of
Lemma~\ref{G-plus-G-minus-complexes-Ext-adjunction}) and moreover,
$\sF=\Com(U\ctrh_\prj)_\prj$ (by Lemma~\ref{Ext-1-as-homotopy-Hom}).
{\hbadness=1075\par}

 In the general case of a quasi-compact semi-separated scheme $X$
with an open covering $\bW$, the argument is based on
Theorem~\ref{loc-contraherent-gluing-theorem}.
 Let $\sR$ be the local class of all commutative rings $R$ and
$\sE^R=\Com(R\modl^\cta)$ be the exact category of complexes of
contraadjusted $R$\+modules.
 Put $\sF(R)=\Com(R\modl^\cta)\cap\Com(R\modl_\vfl)_\vfl$ and
$\sC^R=\Acycl(R\modl^\cta)\sub\sE^R$, as above.
 Then $\sE$ is a very colocal class by
Examples~\ref{colocal-classes-examples}, while $\sC$ is a very
colocal class by Lemma~\ref{acyclicity-in-lcth-criterion}(a).

 Applying Theorem~\ref{loc-contraherent-gluing-theorem}, we obtain
a hereditary complete cotorsion pair $(\sF(X),\sC_\bW^X)$ in
the exact category $\sE_\bW^X=\Com(X\lcth_\bW)$.
 Here $\sF(X)$ is the class of all direct summands of finitely
iterated extensions of direct images of homotopy projective complexes
of projective contraherent cosheaves from $U_\alpha$, while
$\sC_\bW^X$ is the class of all locally\+$\sC$ complexes of
$\bW$\+locally contraherent cosheaves on~$X$.
 By Lemma~\ref{acyclicity-in-lcth-criterion}(a), we have
$\sC_\bW^X=\Acycl(X\lcth_\bW)$.
 Once again, it remains to notice that
$\sF(X)\sub\Com(X\ctrh_\prj)$ (e.~g., by
Lemma~\ref{G-plus-G-minus-complexes-Ext-adjunction}, or by
Corollary~\ref{proj-direct-inverse}(a)) and moreover,
$\sF(X)=\Com(X\ctrh_\prj)_\prj$ (by Lemma~\ref{Ext-1-as-homotopy-Hom}).
{\hbadness=1625\par}

 The proof of part~(b) is similar.
 One needs to restrict the hereditary complete cotorsion pair
$(\sF',\sC')$ of
Lemma~\ref{homotopy-vfl-flat-acyclic-of-cta-cot-pairs}(b) for
the affine scheme $U=\Spec R$ to the exact subcategory
$\sE=\Com(R\modl^\cot)\sub\Com(R\modl)$ and use
Lemma~\ref{acyclicity-in-lcth-criterion}(b)
(cf.~\cite[Example~7.14]{Pal}).
 Corollary~\ref{proj-direct-inverse}(b) may be also used.
\end{proof}

\begin{cor} \label{qcomp-ssep-homotopy-projective}
\textup{(a)} The natural functors\/ $\Hot(X\ctrh_\prj)_\prj\rarrow
\Hot(X\lcth_\bW)_\prj\allowbreak\rarrow\sD(X\lcth_\bW)$ are equivalences
of triangulated categories. \par
\textup{(b)} The natural functors\/ $\Hot(X\ctrh^\lct_\prj)_\prj
\rarrow\Hot(X\lcth_\bW^\lct)_\prj\rarrow\sD(X\lcth^\lct_\bW)$
are equivalences of triangulated categories.
\end{cor}

\begin{proof}
 This is similar to~\cite[Corollary~6.7]{PS4}.
 We will prove part~(a).
 As in Theorem~\ref{quasi-homotopy-injective}, it is clear that both
functors are fully faithful.
 To prove the essential surjectivity, suppose given a complex $\gM^\bu$
in $X\lcth_\bW$.
 We need to find a complex $\gF^\bu$ in $\Hot(X\ctrh_\prj)_\prj$
together with a quasi-isomorphism $\gF^\bu\rarrow\gM^\bu$ of complexes
in $X\lcth_\bW$.
 For this purpose, consider a special precover short exact sequence
$0\rarrow\gE^\bu\rarrow\gF^\bu\rarrow\gM^\bu$ in the cotorsion pair
of Theorem~\ref{lcth-homotopy-projective-cotorsion-pairs}(a).
 So we have $\gF^\bu\in\Com(X\ctrh_\prj)_\prj$ and $\gE^\bu\in
\Acycl(X\lcth_\bW)$.
 Since the complex $\gE^\bu$ is acyclic in $X\lcth_\bW$, it follows
that the morphism of complexes $\gF^\bu\rarrow\gM^\bu$ is
a quasi-isomorphism.
\end{proof}

 Comparing Lemma~\ref{qcomp-ssep-homotopy-proj-independence}(a) with
Corollary~\ref{qcomp-ssep-homotopy-projective}(a), we now come to
the promised conclusion that the property of a complex of $\bW$\+locally
contraherent cosheaves on $X$ to be homotopy projective is not affected
by refinements of the covering~$\bW$.
 Similarly, Lemma~\ref{qcomp-ssep-homotopy-proj-independence}(b) with
Corollary~\ref{qcomp-ssep-homotopy-projective}(b) imply that
the property of a complex of locally cotorsion $\bW$\+locally
contraherent cosheaves on $X$ to be homotopy projective is not changed
by refinements of the covering~$\bW$.

\medskip

 Let $X$ be a scheme.
 A complex $\gF^\bu$ of $\bW$\+locally contraherent cosheaves on $X$ is
called \emph{homotopy antilocally flat} if the complex of abelian groups
$\Hom^X(\gF^\bu,\gM^\bu)$ is acyclic for any complex of locally
cotorsion $\bW$\+locally contraherent cosheaves $\gM^\bu$ acyclic over
the exact category $X\lcth_\bW^\lct$.
 By Theorem~\ref{cotorsion-periodicity} or
Lemma~\ref{acyclicity-in-lcth-criterion}, the latter condition is
equivalent to the acyclicity over $X\lcth_\bW$.
 The full subcategory of homotopy antilocally flat complexes in
$\Hot(X\lcth_\bW)$ is denoted by $\Hot(X\lcth_\bW)_\alf$.
 Notice that the property of a complex of $\bW$\+locally contraherent
cosheaves to be homotopy antilocally flat may possibly change when
the covering $\bW$ is replaced by its refinement.

\begin{lem}  \label{qcomp-ssep-homotopy-flat-independence}
 Let $X$ be a quasi-compact semi-separated scheme with an open
covering\/~$\bW$.
 Then a complex of antilocally flat contraherent cosheaves\/ $\gF^\bu$
on $X$ is homotopy antilocally flat if and only if the complex\/
$\Hom^X(\gF^\bu,\gE^\bu)$ is acyclic for any complex\/ $\gE^\bu$ over
$X\ctrh^\lct_\prj$ acyclic with respect to $X\ctrh_\al$.
\end{lem}

\begin{proof}
 The argument is similar to the proof of
Lemma~\ref{qcomp-ssep-homotopy-proj-independence} and based on
Corollary~\ref{becker-contraderived-of-lcta-lct-well-behaved}(b),
Corollary~\ref{Becker-contraacyclic-are-acyclic}(b), and
Lemma~\ref{lct-lin-clp-finite-dim}(a).
 The only difference is that one has to use
Lemma~\ref{all-alf-becker-contraacyclic-lct-cotorsion-pair}
together with Lemma~\ref{Ext-1-as-homotopy-Hom} in order to show that
the complex $\Hom^X$ from any complex of antilocally flat contraherent
cosheaves on $X$ to a Becker-contraacyclic complex over
$X\lcth_\bW^\lct$ is acyclic.
\end{proof}

 According to Lemma~\ref{qcomp-ssep-homotopy-flat-independence},
the property of a complex over $X\ctrh_\alf$ to belong to
$\Hot(X\lcth_\bW)_\alf$ does not depend on the covering~$\bW$
(for a quasi-compact semi-separated scheme~$X$).
 We denote the full subcategory in $\Hot(X\ctrh_\alf)$ consisting
of the homotopy antilocally flat complexes by $\Hot(X\ctrh_\alf)_\alf$
and the corresponding full subcategory in $\Com(X\ctrh_\alf)$ by
$\Com(X\ctrh_\alf)_\alf$.
 One can easily check that bounded above complexes of antilocally
flat contraherent cosheaves are homotopy antilocally flat,
$\Hot^-(X\ctrh_\alf)\sub\Hot(X\ctrh_\alf)_\alf$.

\begin{thm}  \label{qcomp-ssep-homotopy-flat}
 Let $X$ be a quasi-compact semi-separated scheme.
 Then the quotient category of the homotopy category of homotopy
antilocally flat complexes of antilocally flat contraherent cosheaves\/
$\Hot(X\ctrh_\alf)_\alf$ on $X$ by its thick subcategory of acyclic
complexes over the exact category $X\ctrh_\alf$ is equivalent to
the derived category\/ $\sD(X\ctrh)$.
\end{thm}

\begin{proof}
 It follows from
Corollary~\ref{acyclic-in-alf-all-of-lct-loc-contraherent-pair}
and Lemma~\ref{Ext-1-as-homotopy-Hom} that $\Acycl(X\ctrh_\alf)
\sub\Com(X\ctrh_\alf)_\alf$.
 According to (the proof of)
Corollary~\ref{qcomp-ssep-homotopy-projective}(b), there is
a quasi-isomorphism into any complex over $X\ctrh$ from a complex
belonging to $\Hot(X\ctrh_\prj)_\prj\sub\Hot(X\ctrh_\alf)_\alf$.
 In view of Lemma~\ref{pkoszul-lemma16}(a), it remains to show that
any homotopy antilocally flat complex of antilocally flat contraherent
cosheaves that is acyclic over $X\ctrh$ is also acyclic over
$X\ctrh_\alf$.

 According to (the proof of)
Corollary~\ref{becker-coderived-of-alf-well-behaved},
any complex $\gF^\bu$ over $X\ctrh_\alf$ admits a morphism into
a complex $\P^\bu$ over $X\ctrh^\lct_\prj$ with a cone $\gE^\bu$
acyclic with respect to $X\ctrh_\alf$.
 If the complex $\gF^\bu$ was homotopy antilocally flat, it follows
that the complex $\P^\bu$ is homotopy antilocally flat, too.
 This means that $\P^\bu$ is a homotopy projective complex of locally
cotorsion contraherent cosheaves on~$X$.
 If the complex $\gF^\bu$ was also acyclic over $X\ctrh$, so is
the complex $\P^\bu$.
 By Theorem~\ref{cotorsion-periodicity} or
Lemma~\ref{acyclicity-in-lcth-criterion}, it follows that $\P^\bu$ is
acyclic over $X\ctrh^\lct$.
 Any homotopy projective complex of locally cotorsion contraherent
cosheaves that is acyclic in $X\ctrh^\lct$ is obviously
contractible.
 Thus the complex $\gF^\bu$ is homotopy equivalent to $\gE^\bu[-1]$,
hence also acyclic over $X\ctrh_\alf$.
\end{proof}

\subsection{Derived direct images in Becker's co/contraderived
categories} \label{direct-images-in-becker-subsect}
 In the rest of Chapter~\ref{becker-on-qcomp-qsep-sect}, as in
Sections~\ref{derived-direct-inverse}\+-\ref{finite-dim-morphisms-I},
the upper index~$\bst$ in the notation for derived and homotopy
categories stands for one of the symbols $\b$, $+$, $-$, $\empt$,
$\abs+$, $\abs-$, $\bco$, $\bctr$, $\co$, $\ctr$, or~$\abs$.
 This section complements the discussion in
Section~\ref{derived-direct-inverse}.

 Let $f\:Y\rarrow X$ be a morphism of quasi-compact semi-separated
schemes.
 Then the right derived functor of direct image
\begin{equation}  \label{bco-qcoh-direct}
 \boR f_*\:\sD^\bco(Y\qcoh)\lrarrow\sD^\bco(X\qcoh)
\end{equation}
is constructed in the following way.
 By Corollary~\ref{qcoh-dil-cta-derived-equiv-cor}(b)
or~\ref{dil-cta-cot-into-all-equiv-on-bco}(b), the natural functor
$\sD^\bco(Y\qcoh^\cta)\rarrow\sD^\bco(Y\qcoh)$ is an equivalence of
triangulated categories (as is the similar functor for sheaves
over~$X$).
 By Corollary~\ref{dil-cta-cot-direct-image-becker-coacyclic}(b),
the direct image functor~$f_*$ takes Becker-coacyclic complexes
in $Y\qcoh^\cta$ to Becker-coacyclic complexes in $X\qcoh^\cta$.
 Now the derived functor~$\boR f_*$ is defined by restricting
the functor of direct image $f_*\:\Hot(Y\qcoh)\rarrow\Hot(X\qcoh)$
to the full subcategory of complexes of contraadjusted quasi-coherent
sheaves on~$Y$.

 Alternatively, Corollary~\ref{qcoh-dil-cta-derived-equiv-cor}(a)
or~\ref{dil-cta-cot-into-all-equiv-on-bco}(a) tells us that the natural 
functor $\sD^\bco(Y\qcoh^\dil)\rarrow\sD^\bco(Y\qcoh)$ is an equivalence
of triangulated categories (as is the similar functor for sheaves
over~$X$).
 By Corollary~\ref{dil-cta-cot-direct-image-becker-coacyclic}(a),
the direct image functor~$f_*$ takes Becker-coacyclic complexes
in $Y\qcoh^\dil$ to Becker-coacyclic complexes in $X\qcoh^\dil$.
 The derived functor~$\boR f_*$ \eqref{bco-qcoh-direct} can be defined
by restricting the functor of direct image $f_*\:\Hot(Y\qcoh)\rarrow
\Hot(X\qcoh)$ to the full subcategory of complexes of dilute
quasi-coherent sheaves on~$Y$.

 The left derived functor of direct image
\begin{equation}   \label{bctr-ctrh-direct}
 \boL f_!\:\sD^\bctr(Y\ctrh)\lrarrow\sD^\bctr(X\ctrh)
\end{equation}
is constructed in the following way.
 By Corollary~\ref{ctrh-lcth-cor}(a), the natural functor
$\sD^\bctr(Y\ctrh_\al)\rarrow\sD^\bctr(Y\ctrh)$ is an equivalence of
triangulated categories (as is the similar functor for cosheaves
over~$X$).
 By Corollary~\ref{morphism-qcomp-qsep-direct-image-of-bcontra-al}(a),
the direct image functor~$f_!$ takes Becker-contraacyclic complexes
in $Y\ctrh_\al$ to Becker-contraacyclic complexes in $X\ctrh_\al$.
 The derived functor $\boL f_!$ is defined as the induced functor
$\sD^\st(Y\ctrh_\al)\rarrow\sD^\st(X\ctrh_\al)$.

 Similarly one defines the left derived functor of direct image
\begin{equation}   \label{bctr-ctrh-lct-direct}
 \boL f_!\:\sD^\bctr(Y\ctrh^\lct)\rarrow\sD^\bctr(X\ctrh^\lct).
\end{equation}
 Recall that, by Corollary~\ref{lct-Becker-contraacyclic-iff-as-lcta},
a complex in $\sD^\bctr(Y\ctrh^\lct)$ is Becker-contraacyclic if and
only if it is Becker-contraacyclic in $\sD^\bctr(Y\ctrh)$ (and
similarly for cosheaves over~$X$).
 Hence the functors~\eqref{bctr-ctrh-direct}
and~\eqref{bctr-ctrh-lct-direct} agree with each other, and we obtain
a commutative diagram of triangulated functors and triangulated
equivalences provided by
Theorem~\ref{becker-contraderived-lcta-lct-equivalent},
\begin{equation} \label{bcontra-lct-lcta-direct-images-compatible}
\begin{gathered}
 \xymatrix{
  \sD^\bctr(Y\ctrh^\lct) \ar@<2pt>[r] \ar@<-2pt>@{-}[r]
  \ar[d]^{\boL f_!}
  & \sD^\bctr(Y\ctrh) \ar[d]^{\boL f_!} \\
  \sD^\bctr(X\ctrh^\lct) \ar@<2pt>[r] \ar@<-2pt>@{-}[r]
  & \sD^\bctr(X\ctrh)
 }
\end{gathered}
\end{equation}

\begin{thm} \label{becker-co-contra-derived-direct-inverse-adjunction}
\textup{(a)} Let\/ $\bW$ be an open covering of $X$ and\/ $\bT$ be
an open covering of $Y$ such that the morphism $f\:Y\rarrow X$ is
$(\bW,\bT)$\+coaffine.
 Then equivalences of categories\/
$\sD^\bco(Y\qcoh)\simeq\sD(Y\lcth_\bT^\lin)$ and\/
$\sD^\bco(X\qcoh)\simeq\sD(X\lcth_\bW^\lin)$ from
Corollary~\textup{\ref{bco-qcoh-lin-derived-equivalence}} transform
the derived functor\/ $\boR f_*$ \eqref{bco-qcoh-direct} into
a left adjoint functor to the triangulated functor
$f^!\:\sD(X\lcth_\bW^\lin)\rarrow\sD(Y\lcth_\bT^\lin)$ induced
by the exact functor of inverse image
$f^!\:X\lcth_\bW^\lin\rarrow Y\lcth_\bT^\lin$,
\begin{equation} \label{bco-qcoh-direct-adjoint-to-lin-inverse}
\begin{gathered}
 \xymatrix{
  \sD^\bco(Y\qcoh) \ar@{=}[r] \ar[d]_{\boR f_*}
  & \sD(Y\lcth_\bT^\lin) \\
  \sD^\bco(X\qcoh) \ar@{=}[r]
  & \sD(X\lcth_\bW^\lin) \ar[u]_{f^!}
 }
\end{gathered}
\end{equation} \par
\textup{(b)} The equivalences of categories\/
$\sD^\bctr(Y\ctrh)\simeq\sD(Y\qcoh_\vfl)$ and\/
$\sD^\bctr(X\ctrh)\allowbreak\simeq\sD(X\qcoh_\vfl)$ from
Corollary~\textup{\ref{bctr-lcth-vfl-fl-derived-equivalences}}
transform the derived functor\/ $\boL f_!$ \eqref{bctr-ctrh-direct}
into a right adjoint functor to the triangulated functor
$f^*\:\sD(X\qcoh_\vfl)\rarrow\sD(Y\qcoh_\vfl)$ induced by the exact
functor of inverse image $f^*\:X\qcoh_\vfl\rarrow Y\qcoh_\vfl$,
\begin{equation} \label{bctr-lcth-direct-adjoint-to-vfl-inverse}
\begin{gathered}
 \xymatrix{
  \sD(Y\qcoh_\vfl) \ar@{=}[r]
  & \sD^\bctr(Y\ctrh) \ar[d]^{\boL f_!} \\
  \sD(X\qcoh_\vfl) \ar@{=}[r] \ar[u]^{f^*}
  & \sD^\bctr(X\ctrh)
 }
\end{gathered}
\end{equation} \par \hfuzz=10pt
\textup{(c)} The equivalences of categories\/
$\sD^\bctr(Y\ctrh^\lct)\simeq\sD(Y\qcoh_\fl)$ and\/
$\sD^\bctr(X\ctrh^\lct)\allowbreak\simeq\sD(X\qcoh_\fl)$ from
Corollary~\textup{\ref{bctr-lcth-vfl-fl-derived-equivalences}}
transform the derived functor\/ $\boL f_!$ \eqref{bctr-ctrh-lct-direct}
into a right adjoint functor to the triangulated functor
$f^*\:\sD(X\qcoh_\fl)\rarrow\sD(Y\qcoh_\fl)$ induced by the exact
functor of inverse image $f^*\:X\qcoh_\fl\rarrow Y\qcoh_\fl$,
\begin{equation} \label{bctr-lcth-lct-direct-adjoint-to-fl-inverse}
\begin{gathered}
 \xymatrix{
  \sD(Y\qcoh_\fl) \ar@{=}[r]
  & \sD^\bctr(Y\ctrh^\lct) \ar[d]^{\boL f_!} \\
  \sD(X\qcoh_\fl) \ar@{=}[r] \ar[u]^{f^*}
  & \sD^\bctr(X\ctrh^\lct)
 }
\end{gathered}
\end{equation} \par
 In all the three
diagrams~\textup{(\ref{bco-qcoh-direct-adjoint-to-lin-inverse}\+-%
\ref{bctr-lcth-lct-direct-adjoint-to-fl-inverse})}, the left vertical
functor is left adjoint to the right vertical functor.
\end{thm}

\begin{proof}[First proof]
 Part~(a): let us interpret both the functors in question as acting
between the derived categories of locally injective locally
contraherent cosheaves.
 Moved to the world of locally contraherent cosheaves, the construction
of the functor $\boR f_*$ \eqref{bco-qcoh-direct} does the following.
 Given a complex $\Q^\bu$ in $Y\lcth_\bT^\lin$, one finds a complex
$\P^\bu\in Y\ctrh_\al^\lin$ together with a morphism of complexes
$\P^\bu\rarrow\Q^\bu$ with the cone absolutely acyclic in
$Y\lcth_\bT^\lin$.
 Applying the direct image functor~$f_!$, one obtains a complex
$f_!\P^\bu$ in the exact category $X\ctrh_\al$.
 Notice that the terms of the complex $f_!\P^\bu$ are still antilocal
(by Corollary~\ref{clp-direct}(a)), but they need not be locally
injective in general.
 Let $0\rarrow f_!\P^\bu\rarrow\gK^\bu\rarrow\gA^\bu\rarrow0$ be
a special preenvelope sequence~\eqref{special-preenvelope-sequence} in
the cotorsion pair $(\sF,\sC)$ from the proof
of Lemma~\ref{coacyclic-of-antilocal-all-lin-cotorsion-pairs}(a).
 So $\gK^\bu\in\Com(X\ctrh_\al^\lin)$ is a complex of antilocal
locally injective contraherent cosheaves, while $\gA^\bu\in
\Acycl^\bco(X\ctrh_\al)$ is a Becker-coacyclic complexes of antilocal
contraherent cosheaves.
 Then the functor $\boR f_*$, translated from quasi-coherent sheaves
into locally contraherent cosheaves, takes the complex $\Q^\bu$ to
the complex~$\gK^\bu$.

 Now let $\gJ^\bu\in\Com(X\lcth^\lin_\bW)$ be a complex of locally
injective $\bW$\+locally contraherent cosheaves on~$X$.
 Then we have a short exact sequence of complexes of abelian groups
$0\rarrow\Hom^X(\gA^\bu,\gJ^\bu)\rarrow\Hom^X(\gK^\bu,\gJ^\bu)\rarrow
\Hom^X(f_!\P^\bu,\gJ^\bu)\rarrow0$, and the complex $\Hom^X(\gA^\bu,
\gJ^\bu)$ is acyclic by Lemma~\ref{Ext-1-as-homotopy-Hom} applied
in the context of the cotorsion pair $(\sF''',\sC''')$ from
Lemma~\ref{coacyclic-of-antilocal-all-lin-cotorsion-pairs}(a).
 Hence the map $\Hom^X(\gK^\bu,\gJ^\bu)\allowbreak\rarrow
\Hom^X(f_!\P^\bu,\gJ^\bu)$ is a quasi-isomorphism.
 Recall the isomorphism of complexes $\Hom^X(f_!\P^\bu,\gJ^\bu)\simeq
\Hom^Y(\P^\bu,f^!\gJ^\bu)$ provided by
the formula~\eqref{direct-inverse-cosheaf-lin-adjunction}.
 One can finish the argument by passing to the inductive limit with
respect to morphisms $\gJ^\bu\rarrow\gJ'{}^\bu$ in
$\Hot(X\lcth^\lin_\bW)$ with the cones acyclic in $X\lcth^\lin_\bW$,
or just notice that the complexes $\Hom^X(\gK^\bu,\gJ^\bu)$ and
$\Hom^Y(\P^\bu,f^!\gJ^\bu)$ actually compute $\Hom$ in the respective
derived categories $\sD(X\lcth_\bW^\lin)=\sD^\bctr(X\lcth_\bW^\lin)$
and $\sD(Y\lcth_\bT^\lin)=\sD^\bctr(X\lcth_\bW^\lin)$, because
the terms of the complexes $\gK^\bu$ and $\P^\bu$ are antilocal.

 Part~(c): we interpret both the functors in question as acting between
the derived categories of flat quasi-coherent sheaves.
 Moved to the world of quasi-coherent sheaves, the construction of
the functor $\boL f_!$ \eqref{bctr-ctrh-lct-direct} does the following.
 Given a complex $\G^\bu$ in $Y\qcoh_\fl$, one finds a complex
$\cQ^\bu\in Y\qcoh_\fl^\cot$ together with a morphism of complexes
$\G^\bu\rarrow\cQ^\bu$ with the cone acyclic in $X\qcoh_\fl$.
 Applying the direct image functor~$f_*$, one obtains a complex
$f_*\cQ^\bu$ in the exact category $X\qcoh^\cot$.
 Notice that the terms of the complex $f_*\cQ^\bu$ are still cotorsion
(by Corollary~\ref{cta-cot-direct}(b)), but they need not be flat
in general.
 Let $0\rarrow\cB^\bu\rarrow\cH^\bu\rarrow f_*\cQ^\bu\rarrow0$ be
a special precover sequence~\eqref{special-precover-sequence} in
the cotorsion pair $(\sF,\sC)$ from the proof of
Lemma~\ref{all-vfl-flat-contraacyclic-of-cta-cot-pairs}(b).
 So $\cH^\bu\in\Com(X\qcoh^\cot_\fl)$ is a complex of flat cotorsion
quasi-coherent sheaves, while $\cB^\bu\in\Acycl^\bctr(X\qcoh^\cot)$
is a Becker-contraacyclic complex of cotorsion quasi-coherent sheaves.
 Then the functor $\boL f_!$, translated from contraherent cosheaves
to quasi-coherent sheaves, takes the complex $\G^\bu$ to
the complex~$\cH^\bu$.

 Now let $\F^\bu\in\Com(X\qcoh_\fl)$ be a complex of flat quasi-coherent
sheaves on~$X$.
 Then we have a short exact sequence of complexes of abelian groups
$0\rarrow\Hom_X(\F^\bu,\cB^\bu)\rarrow\Hom_X(\F^\bu,\cH^\bu)\rarrow
\Hom_X(\F^\bu,f_*\cQ^\bu)\rarrow0$, and the complex $\Hom_X(\F^\bu,
\cB^\bu)$ is acyclic by Lemma~\ref{Ext-1-as-homotopy-Hom} applied
in the context of the cotorsion pair $(\sF'',\sC'')$ from
Lemma~\ref{all-vfl-flat-contraacyclic-of-cta-cot-pairs}(b).
 Hence the map $\Hom_X(\F^\bu,\cH^\bu)\allowbreak\rarrow
\Hom_X(\F^\bu,f_*\cQ^\bu)$ is a quasi-isomorphism.
 The argument finishes dual-analogously to the proof of part~(a).
 The proof of part~(b) is similar to that of part~(c).
\end{proof}

\begin{proof}[Second proof]
 Part~(a): notice first of all that $\sD(X\lcth_\bW^\lin)=
\sD^\abs(X\lcth_\bW^\lin)$ by Corollary~\ref{lin-ctrh-lcth-cor}(a),
and similarly for $Y\lcth_\bT^\lin$.
 Let $\I^\bu$ be a complex over $Y\qcoh^\inj$ representing
a given object of $\sD^\bco(Y\qcoh)$ (see
Theorem~\ref{quasi-coherent-becker-coderived})
and $\gK^\bu$ be a complex over $X\ctrh_\al^\lin$ representing
a given object of $\sD^\abs(X\lcth_\bW^\lin)$.
 Let $f_*\I^\bu\rarrow\J^\bu$ be a morphism from the complex
$f_*\I^\bu$ over $X\qcoh^\cta$ (see Corollary~\ref{cta-cot-direct}(a))
to a complex $\J^\bu$ over $X\qcoh^\inj$ with a cone Becker-coacyclic
with respect to $X\qcoh$.
 We have to construct a natural isomorphism of abelian groups
$\Hom_{\sD^\abs(X\lcth^\lin_\bW)}(\fHom_X(\O_X,\J^\bu),\gK^\bu)\simeq
\Hom_{\sD^\abs(Y\lcth^\lin_\bT)}(\fHom_Y(\O_Y,\I^\bu),f^!\gK^\bu)$.

 Both $\fHom_X(\O_X,\J^\bu)$ and $\gK^\bu$ are complexes over
$X\ctrh_\al^\lin$ (by Lemma~\ref{cta-clp-restricts-to-cot-inj}(b)).
 Hence one has
\begin{multline*}
 \Hom_{\sD^\abs(X\lcth^\lin_\bW)}(\fHom_X(\O_X,\J^\bu),\gK^\bu)\simeq
 \Hom_{\Hot(X\ctrh^\lin_\al)}(\fHom_X(\O_X,\J^\bu),\gK^\bu) \\
 \simeq \Hom_{\Hot(X\qcoh^\inj)}(\J^\bu\;\O_X\ocn_X\gK^\bu) \\
 \simeq\Hom_{\Hot(X\qcoh^\cta)}(f_*\I^\bu\;\O_X\ocn_X\gK^\bu)
 \simeq\Hom_{\Hot(X\ctrh_\al)}(\fHom_X(\O_X,f_*\I^\bu),\gK^\bu).
\end{multline*}
 Here the first isomorphism holds because $\sD^\abs(X\ctrh^\lin)
\simeq\Hot(X\ctrh^\lin_\al)$ (see
Corollary~\ref{lin-ctrh-lcth-cor}(a)) and the second isomorphism is
provided by the equivalence of categories $X\ctrh^\lin_\al\simeq
X\qcoh^\inj$ (see Lemma~\ref{cta-clp-restricts-to-cot-inj}(b)).
 The third isomorphism holds because $\O_X\ocn_X\gK^\bu$ is
a complex over $X\qcoh^\inj$, and the fourth one comes from
the equivalence of categories $X\qcoh^\cta\simeq X\ctrh_\al$
(see Lemma~\ref{cta-clp-equivalence}).

 Furthermore, by~\eqref{flat-cta-fHom-projection}
and~\eqref{direct-inverse-cosheaf-lin-adjunction} one has
\begin{multline*}
 \Hom_{\Hot(X\ctrh_\al)}(\fHom_X(\O_X,f_*\I^\bu),\gK^\bu)\simeq
 \Hom_{\Hot(X\ctrh_\al)}(f_!\fHom_Y(\O_Y,\I^\bu),\gK^\bu) \\
 \simeq\Hom_{\Hot(Y\lcth_\bT^\lin)}(\fHom_Y(\O_Y,\I^\bu),f^!\gK^\bu)
 \simeq\Hom_{\sD^\abs(Y\lcth^\lin_\bT)}
 (\fHom_Y(\O_Y,\I^\bu),f^!\gK^\bu),
\end{multline*}
where the last isomorphism follows from the proof of
Lemma~\ref{homotopy-inj-proj-fully-faithful}(b), as
$\fHom_Y(\O_Y,\I^\bu)$ is a complex over $Y\ctrh_\al^\lin$ and
the objects of $Y\ctrh_\al^\lin$ are projective in the exact category
$Y\lcth_\bT^\lin$. {\hbadness=1500\par}

 Part~(c): notice first of all that $\sD(X\qcoh_\fl)=
\sD^\bco(X\qcoh_\fl)$
by Theorem~\ref{vlf-cta-fl-cot-derived-vlf-fl-equivalences}(c),
and similarly for $Y\qcoh_\fl$.
 Let $\L^\bu$ be a complex over $X\qcoh^\cot_\fl$ representing
a given object of $\sD^\bco(X\qcoh_\fl)$ and $\gG^\bu$ be a complex
over $Y\ctrh^\lct_\prj$ representing a given object of
$\sD^\bctr(X\ctrh^\lct)$ (see
Corollary~\ref{becker-contraderived-of-lcta-lct-well-behaved}(b)).
 Let $\gF^\bu\rarrow f_!\gG^\bu$ be a morphism into the complex
$f_!\gG^\bu$ over $X\ctrh^\lct_\al$ (see Corollary~\ref{clp-direct}(b))
from a complex $\gF^\bu$ over $X\ctrh^\lct_\prj$ with a cone
Becker-contraacyclic with respect to $X\ctrh^\lct$.
 We have to construct a natural isomorphism of abelian groups
$\Hom_{\sD^\bco(X\qcoh_\fl)}(\L^\bu\;\O_X\ocn_X\gF^\bu)\simeq
\Hom_{\sD^\bco(Y\qcoh_\fl)}(f^*\L^\bu\;\O_Y\ocn_Y\gG^\bu)$.

 Both $\L^\bu$ and $\O_X\ocn_X\gF^\bu$ are complexes over
$X\qcoh^\cot_\fl$ (by Lemma~\ref{cta-clp-restricts-to-prj-clf}(b)).
 Hence one has
\begin{multline*}
 \Hom_{\sD^\bco(X\qcoh_\fl)}(\L^\bu\;\O_X\ocn_X\gF^\bu)\simeq
 \Hom_{\Hot(X\qcoh^\cot_\fl)}(\L^\bu\;\O_X\ocn_X\gF^\bu) \\
 \simeq\Hom_{\Hot(X\ctrh^\lct_\prj)}(\fHom_X(\O_X,\L^\bu),\gF^\bu) \\
 \simeq\Hom_{\Hot(X\ctrh^\lct_\al)}(\fHom_X(\O_X,\L^\bu),f_!\gG^\bu)
 \simeq\Hom_{\Hot(X\qcoh^\cot)}(\L^\bu\;\O_X\ocn_X f_!\gG^\bu).
\end{multline*}
 Here the first isomorphism holds because $\sD^\bco(X\qcoh_\fl)
\simeq\Hot(X\qcoh^\cot_\fl)$ (see
Theorem~\ref{vlf-cta-fl-cot-derived-vlf-fl-equivalences}(c)) and
the second isomorphism is provided by the equivalence of categories
$X\qcoh^\cot_\fl\simeq X\ctrh^\lct_\prj$
(see Lemma~\ref{cta-clp-restricts-to-prj-clf}(b)).
 The third isomorphism holds because $\fHom_X(\O_X,\L^\bu)$ is
a complex over $X\ctrh^\lct_\prj$, and the fourth one comes
from the equivalence of categories $X\ctrh^\lct_\al\simeq X\qcoh^\cot$
(see Lemma~\ref{cta-clp-restricts-to-cot-inj}(a)).

 Furthermore, one has
\begin{multline*}
 \Hom_{\Hot(X\qcoh^\cot)}(\L^\bu\;\O_X\ocn_X f_!\gG^\bu)\simeq
 \Hom_{\Hot(X\qcoh^\cot)}(\L^\bu\;f_*(\O_Y\ocn_Y\gG^\bu)) \\
 \simeq\Hom_{\Hot(Y\qcoh_\fl)}(f^*\L^\bu\;\O_Y\ocn_Y\gG^\bu)
 \simeq\Hom_{\sD^\bco(Y\qcoh_\fl)}(f^*\L^\bu\;\O_Y\ocn_Y\gG^\bu).
\end{multline*}
 Here the first isomorphism holds by
Corollary~\ref{cta-clp-direct-image-cor} (since $\gG^\bu$ is
a complex of antilocal contraherent cosheaves on $Y$), and
the last one is true because $\O_Y\ocn_Y\gG^\bu$ is a complex over
$Y\qcoh^\cot_\fl$ and the objects of $Y\qcoh^\cot_\fl$ are injective
in the exact category $Y\qcoh_\fl$.
\end{proof}

\subsection{Derived inverse images for morphisms of finite flat
dimension} \label{finite-dim-morphisms-derived-inverse-subsect}
 The definitions of the \emph{flat dimension} and the \emph{very
flat dimension} of a morphism of schemes $f\:Y\rarrow X$ were
given in Section~\ref{finite-dim-morphisms-I}.

\begin{lem} \label{Tor-Ext-perp-cotorsion-pairs-lemma}
 Let $R$ be a commutative ring and $M$ be an $R$\+module.  Then \par
\textup{(a)} there exists a hereditary complete cotorsion pair
$(\sF,\sC)$ in the abelian category $R\modl$ such that\/ $\sF$
is the class of all modules $F\in R\modl$ satisfying the condition\/
$\Tor^R_{>0}(M,F)=0$; \par
\textup{(b)} there exists a hereditary complete cotorsion pair
$(\sF,\sC)$ in the exact category $R\modl^\cta$ such that\/ $\sC$
is the class of all modules $C\in R\modl^\cta$ satisfying
the condition $\Ext_R^{>0}(M,C)=0$; \par
\textup{(c)} there exists a hereditary complete cotorsion pair
$(\sF,\sC)$ in the exact category $R\modl^\cot$ such that\/ $\sC$
is the class of all modules $C\in R\modl^\cot$ satisfying
the condition $\Ext_R^{>0}(M,C)=0$.
\end{lem}

\begin{proof}
 Part~(a): this assertion actually holds over any associative ring
$R$ (for a right $R$\+module~$M$).
 Let $\dotsb\rarrow G_2\rarrow G_1\rarrow G_0\rarrow M\rarrow0$ be
a flat resolution of the $R$\+module~$M$; denote by $Z_n$ the kernels
of the morphisms $G_n\rarrow G_{n-1}$.
 Then the desired cotorsion pair $(\sF,\sC)$ is cogenerated by
the set of pure-injective modules $\Hom_\boZ(M,\boQ/\boZ)$ and
$\Hom_\boZ(Z_n,\boQ/\boZ)$, \ $n\ge0$.
 This proves that $(\sF,\sC)$ is a complete cotorsion pair
by~\cite[Theorem~6.19]{GT}.
 The class $\sF$ is clearly closed under the kernels of epimorphisms,
so the cotorsion pair $(\sF,\sC)$ is hereditary.

 Part~(b): let $\dotsb\rarrow P_2\rarrow P_1\rarrow P_0\rarrow M
\rarrow0$ be a projective resolution of the $R$\+module~$M$;
denote by $Z_n$ the kernels of the morphisms $P_n\rarrow P_{n-1}$.
 Consider the cotorsion pair $(\sF',\sC')$ in $R\modl$ generated by
the $R$\+modules $R[r^{-1}]$, \,$r\in R$, \ $M$, and $Z_n$, \,$n\ge0$.
 Then $(\sF',\sC')$ is a complete cotorsion pair in $R\modl$
by~\cite[Theorems~2 and~10]{ET} or
Theorem~\ref{eklof-trlifaj-general}(a).
 The class $\sC=\sC'$ is clearly closed under the cokernels of
monomorphisms, so the cotorsion pair $(\sF',\sC')$ is hereditary.
 Finally, by Lemmas~\ref{restricting-hereditary-cotorsion}
and~\ref{restricting-cotorsion-pairs-lemma}(b), the cotorsion pair
$(\sF',\sC')$ restricts to the desired hereditary complete cotorsion
pair $(\sF,\sC)$ in the exact subcategory $R\modl^\cta\sub R\modl$.

 The assertion of part~(c) also holds over any associative ring $R$
(for a left $R$\+mod\-ule~$M$).
 The proof is similar to part~(b), except that one needs to replace
the set of $R$\+modules $R[r^{-1}]$, \,$r\in R$, by a set of
left $R$\+modules $\sS$ such that $R\modl_\fl=\Fil(\sS)$
(see Examples~\ref{deconstructible-examples}).
\end{proof}

 Let $X$ be a quasi-compact semi-separated scheme and $f\:Y\rarrow X$
be a morphism of schemes.
 Continuing the discussion in Section~\ref{finite-dim-morphisms-I},
we will say that a quasi-coherent sheaf $\F$ on an open subscheme
$X'\sub X$ is \emph{adjusted to~$f$} if for all affine open subschemes
$U\sub X'$ and $V\sub Y$ such that $f(V)\sub U$ one has
$\Tor^{\O_X(U)}_{>0}(\O_Y(V),\F(U))=0$.
 By the change-of-scalars lemma~\cite[Proposition~VI.4.1.2]{CaE}
or~\cite[Lemma~4.1(b)]{PSl1}, for an affine open subscheme $X'\sub X$
it suffices to check this condition for $U=X'$.

 Let $\bW$ be an open covering of the scheme~$X$.
 We will say that a $\bW|_{X'}$\+locally contraherent cosheaf $\P$ on
an open subscheme $X'\sub X$ is \emph{adjusted to~$f$} if for all
affine open subschemes $U\sub X'$ subordinate to $\bW$ and $V\sub Y$
such that $f(V)\sub U$ one has $\Ext_{\O_X(U)}^{>0}(\O_Y(V),\P[U])=0$.
 By the change-of-scalars lemma~\cite[Proposition~VI.4.1.4]{CaE}
or~\cite[Lemma~3.4(b)]{Pal}, for an affine open subscheme $X'\sub X$
subordinate to $\bW$ it suffices to check this condition for $U=X'$.

 Furthermore, we will say that a locally cotorsion $\bW|_{X'}$\+locally
contraherent cosheaf $\P$ on an open subscheme $X'\sub X$ is
\emph{adjusted to~$f$} (as a locally cotorsion locally contraherent
cosheaf) if for all affine open subschemes $U\sub X'$ subordinate
to $\bW$ and $V\sub Y$ such that $f(V)\sub U$, and all flat
$\O_Y(V)$\+modules $G$, one has $\Ext_{\O_X(U)}^{>0}(G,\P[U])=0$.
 Once again, for an affine open subscheme $X'\sub X$ subordinate to
$\bW$, it suffices to check this condition for $U=X'$.

 Using \v Cech (co)resolutions~\eqref{cech-quasi}
and~\eqref{contraherent-cech}, one shows that the property of
a quasi-coherent sheaf, a locally contraherent cosheaf, or a locally
cotorsion locally contraherent cosheaf to be adjusted to~$f$ is
local with respect to open coverings of open subschemes of~$X$.
 (For cosheaves, this was mentioned already in
Section~\ref{finite-dim-morphisms-I}.)

\begin{lem} \label{f-adjusted-direct-image-lemma}
 Let $X''\sub X'\sub X$ be a pair of embedded open subschemes.
 Put\/ $\bW'=\bW|_{X'}$ and\/ $\bW''=\bW|_{X''}$, and
assume that the open embedding morphism $h\:X''\rarrow X'$ is affine.
 Then \par
\textup{(a)} the direct image functor $h_*\:X''\qcoh\rarrow X'\qcoh$
takes quasi-coherent sheaves adjusted to~$f$ to quasi-coherent sheaves
adjusted to~$f$; \par
\textup{(b)} the direct image functor $h_!\:X''\lcth_{\bW''}\rarrow
X'\lcth_{\bW'}$ takes locally contraherent cosheaves adjusted to~$f$ to
locally contraherent cosheaves adjusted to~$f$; \par
\textup{(c)} the direct image functor $h_!\:X''\lcth_{\bW''}^\lct
\rarrow X'\lcth_{\bW'}^\lct$ takes locally cotorsion locally
contraherent cosheaves adjusted to~$f$ to locally cotorsion
locally contraherent cosheaves adjusted to~$f$.
\end{lem}

\begin{proof}
 Since the property is local, the question reduces to the case of
affine open subschemes $X'$ and $X''$ in~$X$ (subordinate to $\bW$
in the context of parts~(b\+c)), when once again it suffices to
invoke the change-of-scalars lemma~\cite[Propositions~VI.4.1.1
and~VI.4.1.3]{CaE} or~\cite[Lemma~4.1(a\+b)]{PSl1}.
\end{proof}

\begin{cor} \label{f-adjusted-cotorsion-pairs}
 Let $X$ be a quasi-compact semi-separated scheme with an open
covering\/~$\bW$, and let $f\:Y\rarrow X$ be a morphism of schemes.
 Then \par
\textup{(a)} there exists a hereditary complete cotorsion pair
$(\sF,\sC)$ in the abelian category $X\qcoh$ such that\/
$\sF=X\qcoh_\fadj$ is the class of all quasi-coherent sheaves on $X$
adjusted to~$f$; \par
\textup{(b)} there exists a hereditary complete cotorsion pair
$(\sF,\sC)$ in the exact category $X\lcth_\bW$ such that\/
$\sC=X\lcth_\bW^\fadj$ is the class of all\/ $\bW$\+locally
contraherent cosheaves on $X$ adjusted to~$f$; \par
\textup{(c)} there exists a hereditary complete cotorsion pair
$(\sF,\sC)$ in the exact category $X\lcth_\bW^\lct$ such that\/
$\sC=X\lcth_\bW^{\lct,\fadj}$ is the class of all locally cotorsion\/
$\bW$\+locally contraherent cosheaves on $X$ adjusted to~$f$.
\end{cor}

\begin{proof}
 Part~(a) is provable using Theorem~\ref{eklof-trlifaj-general} in
the Grothendieck abelian category $X\qcoh$ and a suitable version of
Lemma~\ref{locally-L-sheaves-deconstructible}.
 The assumption that $X$ is quasi-compact and semi-separated still
needs to be used in order to check that the class $X\qcoh_\fadj$ is
generating in $X\qcoh$.
 Alternatively, one can build an argument based on
Theorem~\ref{quasi-coherent-gluing-theorem} with
Remark~\ref{quasi-coherent-gluing-nonuniversal-remark}.
 Lemmas~\ref{Tor-Ext-perp-cotorsion-pairs-lemma}(a)
and~\ref{f-adjusted-direct-image-lemma}(a) need to be used.
 The proofs of parts~(b\+c) are based on
Theorem~\ref{loc-contraherent-gluing-theorem}
with Remark~\ref{loc-contraherent-gluing-nonuniversal-remark}.

 Let us explain the proof of part~(c).
 For every affine open subscheme $U\sub X$ subordinate to $\bW$,
put $\sE^U=U\qcoh^\cot\sub\sK^{\O(U)}=\O(U)\modl^\cta$.
 Put $\sC^U=U\ctrh^{\lct,\fadj}=U\lcth_{\{U\}}^{\lct,\fadj}\sub\sE^U$.
 Put $R=\O(U)$, and consider the following $R$\+module~$M$.
 For every affine open subscheme $V\sub Y$ such that $f(V)\sub U$,
pick a set of flat $\O(V)$\+modules $\sS_V$ such that
$\O(V)\modl_\fl=\Fil(\sS_V)$.
 Let $M$ be the direct sum of all $\O(V)$\+modules from $\sS_V$,
taken over all $V$ as above and viewed as an $R$\+module.
 Lemma~\ref{Tor-Ext-perp-cotorsion-pairs-lemma}(c) applied to
the $R$\+module $M$ tells us that there is a hereditary complete
cotorsion pair $(\sF(U),\sC^U)$ in the exact category~$\sE^U$.
 Lemma~\ref{f-adjusted-direct-image-lemma}(c) and the preceding
discussion imply that the system of classes $\sC^U$ is colocal and
preserved by direct images with respect to the identity embeddings
of affine open subschemes subordinate to $\bW$ in~$X$.
 So Remark~\ref{loc-contraherent-gluing-nonuniversal-remark}
is applicable.
\end{proof}

\begin{thm} \label{f-adjusted-becker-co-contra-derived-equivalences}
 Let $X$ be a quasi-compact semi-separated scheme with an open
covering\/~$\bW$, and let $f\:Y\rarrow X$ be a scheme morphism of
finite flat dimension\/~$\le D$.  In this context: \par
\textup{(a)} A complex in $X\qcoh_\fadj$ is Becker-coacyclic in
$X\qcoh_\fadj$ if and only if it is Becker-coacyclic in $X\qcoh$,
$$
 \Acycl^\bco(X\qcoh_\fadj)=\Com(X\qcoh_\fadj)\cap\Acycl^\bco(X\qcoh).
$$
 The inclusion of exact/abelian categories $X\qcoh_\fadj\rarrow X\qcoh$
induces an equivalence of the Becker coderived categories
$$
 \sD^\bco(X\qcoh_\fadj)\simeq\sD^\bco(X\qcoh).
$$ \par
\textup{(b)} Assume that the morphism~$f$ has finite very flat
dimension\/~$\le D$.
 Then a complex in $X\lcth_\bW^\fadj$ is Becker-contraacyclic in
$X\lcth_\bW^\fadj$ if and only if it is Becker-contraacyclic in
$X\lcth_\bW$,
$$
 \Acycl^\bctr(X\lcth_\bW^\fadj)=
 \Com(X\lcth_\bW^\fadj)\cap\Acycl^\bctr(X\lcth_\bW).
$$
 The inclusion of exact categories $X\lcth_\bW^\fadj\rarrow X\lcth_\bW$ 
induces an equivalence of the Becker contraderived categories
$$
 \sD^\bctr(X\lcth_\bW^\fadj)\simeq\sD^\bctr(X\lcth_\bW).
$$ \par
\textup{(c)} A complex in $X\lcth_\bW^{\lct,\fadj}$ is
Becker-contraacyclic in $X\lcth_\bW^{\lct,\fadj}$ if and only if it is
Becker-contraacyclic in $X\lcth_\bW^\lct$,
$$
 \Acycl^\bctr(X\lcth_\bW^{\lct,\fadj})=
 \Com(X\lcth_\bW^{\lct,\fadj})\cap\Acycl^\bctr(X\lcth_\bW^\lct).
$$
 The inclusion of exact categories $X\lcth_\bW^{\lct,\fadj}\rarrow
X\lcth_\bW^\lct$  induces an equivalence of the Becker contraderived
categories
$$
 \sD^\bctr(X\lcth_\bW^{\lct,\fadj})\simeq\sD^\bctr(X\lcth_\bW^\lct).
$$
\end{thm}

\begin{proof}
 Part~(a): in the assumptions of the theorem, every quasi-coherent
sheaf on $X$ has a finite resolution dimension~$\le D$ with respect
to the resolving subcategory $\sF=X\qcoh_\fadj$.
 A hereditary complete cotorsion pair $(\sF,\sC)$ in $X\qcoh$ exists
by Corollary~\ref{f-adjusted-cotorsion-pairs}(a).
 Now the dual version of
Proposition~\ref{becker-contraderived-finite-coresolutions} provides
both the desired assertions.
 The proofs of parts~(b\+c) are similar and based on
Corollary~\ref{f-adjusted-cotorsion-pairs}(b\+c) together with
Proposition~\ref{becker-contraderived-finite-coresolutions}.
\end{proof}

\begin{lem} \label{f-adjusted-inverse-image-preserves-lemma}
 Let $X$ and $Y$ be quasi-compact semi-separated schemes, and let
$f\:Y\rarrow X$ be a morphism of finite flat dimension\/~$\le D$.
 Let\/ $\bW$ be an open covering of $X$ and\/ $\bT$ be an open
covering of $Y$ such that the morphism~$f$ is $(\bW,\bT)$\+coaffine.
 In this context: \par
\textup{(a)} The inverse image functor $f^*\:X\qcoh\rarrow Y\qcoh$
takes Becker-coacyclic complexes of $f$\+adjusted quasi-coherent
sheaves on $X$ to Becker-coacyclic complexes of quasi-coherent
sheaves on~$Y$. \par
\textup{(b)} Assume that the morphism~$f$ has finite very flat
dimension\/~$\le D$.
 Then the inverse image functor $f^!\:X\lcth_\bW^\fadj\rarrow
Y\lcth_\bT$ takes Becker-contraacyclic complexes of $f$\+adjusted\/
$\bW$\+locally contraherent cosheaves on $X$ to Becker-contraacyclic
complexes of\/ $\bT$\+locally contraherent cosheaves on~$Y$. \par
\textup{(c)} The inverse image functor $f^!\:X\lcth_\bW^{\lct,\fadj}
\rarrow Y\lcth_\bT^\lct$ takes Becker-con\-tra\-acyclic complexes of
$f$\+adjusted locally cotorsion\/ $\bW$\+locally contraherent cosheaves
on $X$ to Becker-contraacyclic complexes of locally cotorsion\/
$\bT$\+locally contraherent cosheaves on~$Y$.
\end{lem}

\begin{proof}
 The difficulty here is that, although the functor
$f^*\:X\qcoh_\fadj\rarrow Y\qcoh$ in part~(a) and the functors~$f^!$
in parts~(b\+c) are exact,
Lemma~\ref{exact-with-adjoint-preservation-lemma} is still
not applicable, because the required adjoint functors taking values
in exact categories do not exist.
 So a direct argument is needed.

 Part~(a): by
Theorem~\ref{f-adjusted-becker-co-contra-derived-equivalences}(a)
and Corollary~\ref{qcoh-becker-coacyclicity-is-local}, the question
is local in $X$ and $Y$, so it reduces to affine schemes.
 Given an associative ring homomorphism $f\:R\rarrow S$, let us say,
in the context of this proof of part~(a), that a left $R$\+module $A$
is \emph{adjusted to~$f$} if $\Tor^R_{>0}(S,A)=0$.
 We need to show that, if the flat dimension of the right $R$\+module
$S$ does not exceed~$D$, then the functor of extension of scalars
$S\ot_R\nobreak{-}\,\:R\modl\rarrow S\modl$ takes Becker-coacyclic
complexes of $f$\+adjusted $R$\+modules to Becker-coacyclic complexes
of $S$\+modules.

 Indeed, let $A^\bu$ be a complex of $f$\+adjusted left $R$\+modules
that is Becker-coacyclic in $R\modl$, and let $K^\bu$ be a complex
of injective left $S$\+modules.
 We have to show that the complex of abelian groups
$\Hom_S(S\ot_RA^\bu\;K^\bu)\simeq\Hom_R(A^\bu,K^\bu)$ is acyclic.
 By Lemma~\ref{affine-d+D-lemma}(b) (for $d=0$), \,$K^\bu$ is
a complex of left $R$\+modules of injective dimension~$\le D$.
 Furthermore, by~\cite[Proposition~VI.4.1.3]{CaE}
or~\cite[Lemma~4.1(a)]{PSl1}, we have $\Ext_R^{>0}(A^i,K^j)\simeq
\Ext_S^{>0}(S\ot_R\nobreak A^i,K^j)=0$ for all $i$, $j\in\boZ$.

 Now let $0\rarrow K^\bu\rarrow J^{0,\bu}\rarrow\dotsb\rarrow
J^{D,\bu}\rarrow0$ be a finite exact complex of complexes of
left $R$\+modules with injective $R$\+modules~$J^{k,j}$.
 Let $J^\bu$ denote the total complex of the bicomplex $J^{\bu,\bu}$
and $C^\bu$ denote the total complex of the bicomplex
$K^\bu\rarrow J^{\bu,\bu}$.
 Then the complex $\Hom_R(A^\bu,J^\bu)$ is acyclic since $A^\bu$
is Becker-coacyclic in $R\modl$, and the complex $\Hom_R(A^\bu,C^\bu)$
is acyclic because the exact sequence of complexes
$0\rarrow K^\bu\rarrow J^{0,\bu}\rarrow\dotsb\rarrow J^{D,\bu}\rarrow0$
remains exact after applying $\Hom_R(A^\bu,{-})$.
 Thus the complex $\Hom_R(A^\bu,K^\bu)$ is acyclic.

 Part~(c): by
Theorems~\ref{f-adjusted-becker-co-contra-derived-equivalences}(c)
and~\ref{Becker-contraacyclicity-local-on-qcomp-qsep}(b),
the question is local in $X$ and $Y$, so it reduces to affine schemes.
 Given an associative ring homomorphism $f\:R\rarrow S$, let us say,
in the context of this proof of part~(c), that a left $R$\+module $B$
is \emph{adjusted to~$f$} if $\Ext_R^{>0}(G,B)=0$ for all flat left
$S$\+modules~$G$.
 We need to show that, if the flat dimension of the left $R$\+module
$S$ does not exceed~$D$, then the functor of coextension of scalars
$\Hom_R(S,{-})\:R\modl\rarrow S\modl$ takes Becker-contraacyclic
complexes of $f$\+adjusted cotorsion $R$\+modules to
Becker-contraacyclic complexes of cotorsion $S$\+modules.

 Indeed, let $B^\bu$ be a complex of $f$\+adjusted cotorsion left
$R$\+modules that is Becker-contraacyclic in $R\modl^\cot$, and
let $Q^\bu$ be a complex of flat cotorsion $S$\+modules.
 We have to show that the complex of abelian groups 
$\Hom_S(Q^\bu,\Hom_R(S,B^\bu))\simeq\Hom_R(Q^\bu,B^\bu)$ is acyclic.
 By Lemma~\ref{affine-d+D-lemma}(a) (for $d=0$), \,$Q^\bu$ is
a complex of left $R$\+modules of flat dimension~$\le D$.
 Furthermore, by~\cite[Proposition~VI.4.1.4]{CaE}
or~\cite[Lemma~3.4(b)]{Pal}, we have $\Ext_R^{>0}(Q^j,B^i)\simeq
\Ext_S^{>0}(Q^j,\Hom_R(S,B^i))=0$ for all $i$, $j\in\boZ$, since
the left $S$\+modules $\Hom_R(S,B^i)$ are cotorsion by
Lemma~\ref{adjusted-scalars-contra}(b).

 Now let $0\rarrow F_D^\bu\rarrow F_{D-1}^\bu\rarrow\dotsb\rarrow
F_0^\bu\rarrow Q^\bu\rarrow0$ be a finite exact complex of
complexes of left $R$\+modules with flat $R$\+modules~$F_k^j$.
 Let $F^\bu$ denote the total complex of the bicomplex $F_\bu^\bu$
and $C^\bu$ denote the total complex of the bicomplex
$F_\bu^\bu\rarrow Q^\bu$.
 Then the complex $\Hom_R(F^\bu,B^\bu)$ is acyclic since $F^\bu$
is a complex of flat $R$\+modules and the complex $B^\bu$ is
Becker-contraacyclic in $R\modl^\cot$.
 The point is that ($\Com(R\modl_\fl)$, $\Acycl^\bctr(R\modl^\cot)$)
is a cotorsion pair in $\Com(R\modl)$; see~\cite[Proposition~3.2 and Theorem~5.5]{Gil2} or~\cite[Lemma~4.9]{Gil3}, and also
Lemma~\ref{all-vfl-flat-contraacyclic-of-cta-cot-pairs}(b)
above for a quasi-coherent sheaf version; then it remains to recall
Lemma~\ref{Ext-1-as-homotopy-Hom}.
 On the other hand, the complex $\Hom_R(C^\bu,B^\bu)$ is acyclic
because the exact sequence of complexes $0\rarrow F_D^\bu\rarrow
F_{D-1}^\bu\rarrow\dotsb\rarrow F_0^\bu\rarrow Q^\bu\rarrow0$
remains exact after applying $\Hom_R({-},B^\bu)$.
 Thus the complex $\Hom_R(Q^\bu,B^\bu)$ is acyclic.

 The proof of part~(b) is similar to part~(c).
\end{proof}

 Let $X$, $Y$ be quasi-compact semi-separated schemes and
$f\:Y\rarrow X$ be a morphism of finite flat dimension.
 Then the left derived functor
\begin{equation}  \label{bco-qcoh-inverse-ffd-morphism}
 \boL f^*\:\sD^\bco(X\qcoh)\lrarrow\sD^\bco(Y\qcoh)
\end{equation}
is constructed as the functor $f^*\:\sD^\bco(X\qcoh_\fadj)\rarrow
\sD^\bco(Y\qcoh)$ induced by the functor $f^*\:X\qcoh_\fadj\rarrow
Y\qcoh$.
 Theorem~\ref{f-adjusted-becker-co-contra-derived-equivalences}(a) and
Lemma~\ref{f-adjusted-inverse-image-preserves-lemma}(a) are used here.
 The functor $\boL f^*$ is left adjoint to the functor $\boR f_*$ 
\eqref{bco-qcoh-direct}
from Section~\ref{direct-images-in-becker-subsect}.

 Let $\bW$ and $\bT$ be open coverings of the schemes $X$ and $Y$
for which the morphism~$f$ is $(\bW,\bT)$\+coaffine.
 Assuming that the morphism~$f$ has finite very flat dimension,
the right derived functor
\begin{equation}  \label{bctr-lcth-inverse-ffd-morphism}
 \boR f^!\:\sD^\bctr(X\lcth_\bW)\lrarrow\sD^\bctr(Y\lcth_\bT)
\end{equation}
is constructed as the functor $f^!\:\sD^\bctr(X\lcth_\bW^\fadj)
\rarrow\sD^\bctr(Y\lcth_\bT)$ induced by the functor
$f^!\:X\lcth_\bW^\fadj\rarrow Y\lcth_\bT$.
 Lemma~\ref{adjusted-scalars-contra}(a),
Theorem~\ref{f-adjusted-becker-co-contra-derived-equivalences}(b), and
Lemma~\ref{f-adjusted-inverse-image-preserves-lemma}(b) are used here.
 Then one can use the equivalences of categories from
Corollary~\ref{ctrh-lcth-cor}(a) in order to obtain the right
derived functor
\begin{equation}  \label{bctr-ctrh-inverse-ffd-morphism}
 \boR f^!\:\sD^\bctr(X\ctrh)\lrarrow\sD^\bctr(Y\ctrh),
\end{equation}
which is right adjoint to the functor $\boL f_!$
\eqref{bctr-ctrh-direct}
from Section~\ref{direct-images-in-becker-subsect}.

 For quasi-compact semi-separated schemes $X$, $Y$, a morphism
$f\:Y\rarrow X$ of finite flat dimension, and any open coverings
$\bW$ and $\bT$ of the schemes $X$ and $Y$ for which the morphism~$f$
is $(\bW,\bT)$\+coaffine, the right derived functor
\begin{equation}  \label{bctr-lcth-lct-inverse-ffd-morphism}
 \boR f^!\:\sD^\bctr(X\lcth_\bW^\lct)\lrarrow\sD^\bctr(Y\lcth_\bT^\lct)
\end{equation}
is constructed as the functor $f^!\:\sD^\bctr(X\lcth_\bW^{\lct,\fadj})
\rarrow\sD^\bctr(Y\lcth_\bT^\lct)$ induced by the functor
$f^!\:X\lcth_\bW^{\lct,\fadj}\rarrow Y\lcth_\bT^\lct$.
 Lemma~\ref{adjusted-scalars-contra}(b),
Theorem~\ref{f-adjusted-becker-co-contra-derived-equivalences}(c), and
Lemma~\ref{f-adjusted-inverse-image-preserves-lemma}(c) are used here.
 Then one can use the equivalences of categories from
Corollary~\ref{lct-ctrh-lcth-cor}(a) in order to obtain the right
derived functor
\begin{equation}  \label{bctr-ctrh-lct-inverse-ffd-morphism}
 \boR f^!\:\sD^\bctr(X\ctrh^\lct)\lrarrow\sD^\bctr(Y\ctrh^\lct),
\end{equation}
which is right adjoint to the functor
$\boL f_!$~\eqref{bctr-ctrh-lct-direct}.

 For a morphism~$f$ of finite very flat dimension,
the functors~\eqref{bctr-ctrh-inverse-ffd-morphism}
and~\eqref{bctr-ctrh-lct-inverse-ffd-morphism} agree with each other,
i.~e., there is a commutative diagram of triangulated functors and
triangulated equivalences provided by
Theorem~\ref{becker-contraderived-lcta-lct-equivalent},
\begin{equation} \label{bcontra-lct-lcta-inverse-images-compatible}
\begin{gathered}
 \xymatrix{
  \sD^\bctr(Y\ctrh^\lct) \ar@<2pt>[r] \ar@<-2pt>@{-}[r]
  & \sD^\bctr(Y\ctrh) \\
  \sD^\bctr(X\ctrh^\lct) \ar@<2pt>[r] \ar@<-2pt>@{-}[r]
  \ar[u]_{\boR f^!} & \sD^\bctr(X\ctrh) \ar[u]_{\boR f^!}
 }
\end{gathered}
\end{equation}

\begin{rem}
 Let $f\:Y\rarrow X$ be a morphism of quasi-compact semi-separated
schemes, and let $\bW$ and $\bT$ be open coverings of $X$ and $Y$
for which the morphism~$f$ is $(\bW,\bT)$\+coaffine.
 Then
Theorem~\ref{becker-co-contra-derived-direct-inverse-adjunction}(a)
tells us that the functor $f^!\:\sD(X\lcth_\bW^\lin)\rarrow
\sD(Y\lcth_\bT^\lin)$ has a left adjoint.
 When the morphism~$f$ has finite flat dimension,
formula~\eqref{bco-qcoh-inverse-ffd-morphism} tells us that
the functor left adjoint to~$f^!$, in turn, has a left adjoint.
 Both the adjoint functors are explicitly constructed in terms of
the Becker coderived categories of quasi-coherent sheaves.

 Dual-analogously, for any morphism $f\:Y\rarrow X$ of quasi-compact
semi-separated schemes,
Theorem~\ref{becker-co-contra-derived-direct-inverse-adjunction}(c)
tells us that the functor $f^*\:\sD(X\qcoh_\fl)\rarrow\sD(Y\qcoh_\fl)$
has a right adjoint.
 When the morphism~$f$ has finite flat dimension,
formula~\eqref{bctr-ctrh-lct-inverse-ffd-morphism} tells us that
the functor right adjoint to~$f^*$, in turn, has a right adjoint.
 Both the adjoint functors are explicitly constructed in terms of
the Becker contraderived categories of contraherent cosheaves.

 Is the finite flat dimension condition on~$f$ necessary for
the existence of the second adjoints?
 Let us point out that it certainly \emph{cannot} be simply dropped.
 Consider the case of a morphism of affine schemes $\Spec S\rarrow
\Spec R$ corresponding to a ring homomorphism $f\:R\rarrow S$.
 More generally, let $f\:R\rarrow S$ be a homomorphism of associative
rings.
 We are interested in the functor $f^*=S\ot_R\nobreak{-}\,\:\allowbreak
\sD(R\modl_\fl)\rarrow\sD(S\modl_\fl)$.
 The equivalences of categories $\Hot(R\modl_\prj)\simeq
\sD(R\modl_\fl)$ and $\Hot(S\modl_\prj)\simeq\sD(S\modl_\fl)$
provided by Theorem~\ref{flat-projective-periodicity}(b)
and Proposition~\ref{flat-projective-periodicity-complements}(b)
identify the functor~$f^*$ with the functor $S\ot_R{-}\,\:
\Hot(R\modl_\prj)\rarrow\Hot(S\modl_\prj)$.

 Assume that the rings $R$ and $S$ are right coherent.
 Then the categories $\Hot(R\modl_\prj)$ and $\Hot(S\modl_\prj)$
are compactly generated~\cite[Proposition~7.14]{N-f}
(see~\cite[Theorem~4.10]{M-th} and
Theorem~\ref{contraderived-compactly-generated}(c\+e) below for
a quasi-coherent sheaf and a contraherent cosheaf versions,
and~\cite[Theorem~11.2]{PS8} for a generalization to contramodules
over topological rings).
 By~\cite[Theorem~5.1]{N-bb}, a direct sum-preserving functor $F$
between compactly generated triangulated categories has
a second right adjoint if and only if $F$ takes compact objects to
compact objects.
 The full subcategory of compact objects in $\Hot(R\modl_\prj)$ is
anti-equivalent to the bounded derived category of finitely presented
right $R$\+modules $\sD^\b(\modrfp R)$, and similarly for
the compact objects in $\Hot(S\modl_\prj)$
\,\cite[Theorem~3.2]{Jorg}, \cite[Proposition~7.12]{N-f}.

 Looking into the construction of this equivalence, one comes to
the conclusion that the functor $S\ot_R{-}\,\:
\Hot(R\modl_\prj)\rarrow\Hot(S\modl_\prj)$ takes compact objects to
compact objects (i.~e., has a second right adjoint) if and only
if the derived functor ${-}\ot_R^\boL S\:\sD^-(\modrfp R)\rarrow
\sD^-(\modrfp S)$ takes the full subcategory of complexes with
bounded cohomology $\sD^\b(\modrfp R)\sub\sD^-(\modrfp R)$ into
$\sD^\b(\modrfp S)\sub\sD^-(\modrfp S)$.
 Of course, this is not true in general, and the condition on the left
$R$\+module $S$ to have finite flat dimension is relevant here.
\end{rem}

\subsection{Sheaves and cosheaves of finite flat and locally injective
dimension} \label{finite-flat-lin-dim-subsect}
 A quasi-coherent sheaf $\F$ on a scheme $X$ is said to have \emph{flat
dimension not exceeding~$d$} if the flat dimension of
the $\O_X(U)$\+module $\F(U)$ does not exceed~$d$ for every affine open
subscheme $U\sub X$.
 If a quasi-coherent sheaf $\F$ on $X$ admits a resolution by flat
quasi-coherent sheaves (e.~g., $X$ is quasi-compact and
semi-separated), then the flat dimension of $\F$ is equal to the minimal
length of such resolution.

 The property of a quasi-coherent sheaf to have flat dimension not
exceeding~$d$ is local, since so is the property of a quasi-coherent
sheaf to be flat
(cf.\ Lemma~\ref{flat-veryflat-cotors-inj-dim-local}(a)).
 Quasi-coherent sheaves of finite flat dimension form a full subcategory
$X\qcoh_\fd\sub X\qcoh$ closed under extensions and kernels of
surjective morphisms; the full subcategory $X\qcoh_\ffdd\sub X\qcoh_\fd$
of quasi-coherent sheaves of flat dimension not exceeding~$d$ is
closed under the same operations, and also under infinite direct sums.

 Let us say that a quasi-coherent sheaf $\F$ on a scheme $X$ has
\emph{very flat dimension not exceeding~$d$} if the very flat dimension
of the $\O_X(U)$\+module $\F(U)$ does not exceed~$d$ for every affine
open subscheme $U\sub X$ (see Section~\ref{veryflat-cotors-dim-subsect}
for the definition).
 Over a quasi-compact semi-separated scheme $X$, a quasi-coherent sheaf
has very flat dimension~$\le d$ if and only if it admits a very flat
resolution of length~$\le\nobreak d$.

 Since the property of a quasi-coherent sheaf to be very flat is local,
so is its property to have very flat dimension not exceeding~$d$
(cf.\ Lemma~\ref{flat-veryflat-cotors-inj-dim-local}(b)).
 Quasi-coherent sheaves of very flat dimension~$\le\nobreak d$ form
a full subcategory $X\qcoh_\fvfdd\sub X\qcoh$ closed under extensions,
kernels of surjective morphisms, and infinite direct sums.
 We denote the inductive limit of the exact categories $X\qcoh_\fvfdd$
as $d\to\infty$ by $X\qcoh_\vfd$.

 A $\bW$\+locally contraherent cosheaf $\P$ on a scheme $X$ is said to
have \emph{locally injective dimension not exceeding~$d$} if
the injective dimension of the $\O_X(U)$\+module $\P[U]$ does not
exceed~$d$ for every affine open subscheme $U\sub X$ subordinate
to~$\bW$.
 Over a quasi-compact semi-separated scheme $X$, a $\bW$\+locally
contraherent cosheaf has locally injective dimension~$\le d$ if and only
if it admits a locally injective coresolution of length~$\le d$
in the exact category $X\lcth_\bW$.

 The property of a $\bW$\+locally contraherent cosheaf to have locally
injective dimension not exceeding~$d$ is local and refinements of
the covering~$\bW$ do not change it
(see Lemma~\ref{flat-veryflat-cotors-inj-dim-local}(d)).
 $\bW$\+locally contraherent cosheaves of finite locally injective
dimension form a full subcategory $X\lcth_\bW^\lid$ closed under
extensions and cokernels of admissible monomorphisms; the full
subcategory $X\lcth_\bW^\flidd\sub X\lcth_\bW^\lid$ of quasi-coherent
sheaves of locally injective dimension not exceeding~$d$ is
closed under the same operations, and also under infinite products.
 We set $X\ctrh^\lid=X\lcth_{\{X\}}^\lid$ and $X\ctrh^\flidd=
X\lcth_{\{X\}}^\flidd$.

 For the rest of the section, let $X$ be a quasi-compact semi-separated
scheme with an open covering~$\bW$.

\begin{cor} \label{finite-vfl-lin-dim-cor}
\textup{(a)} For any symbol\/ $\bst\ne\ctr$, $\bco$, $\bctr$ and any
(finite) integer $d\ge0$, the triangulated functor\/
$\sD^\st(X\qcoh_\fl)\rarrow\sD^\st(X\qcoh_\ffdd)$ induced by
the embedding of exact categories $X\qcoh_\fl\rarrow X\qcoh_\ffdd$
is an equivalence of triangulated categories. \par
\textup{(b)} For any symbol\/ $\bst\ne\ctr$, $\bco$, $\bctr$ and any
(finite) integer $d\ge0$, the triangulated functor\/
$\sD^\st(X\qcoh_\vfl)\rarrow\sD^\st(X\qcoh_\fvfdd)$ induced by
the embedding of exact categories $X\qcoh_\vfl\rarrow X\qcoh_\fvfdd$
is an equivalence of triangulated categories. \par
\textup{(c)} For any symbol\/ $\bst\ne\co$, $\bco$, $\bctr$ and any
(finite) integer $d\ge0$, the triangulated functor\/
$\sD^\st(X\lcth_\bW^\lin)\rarrow\sD^\st(X\lcth_\bW^\flidd)$ induced by
the embedding of exact categories $X\lcth_\bW^\lin\rarrow
X\lcth_\bW^\flidd$ is an equivalence of triangulated categories.
\end{cor}

\begin{proof}
 Parts~(a\+b) follow from Proposition~\ref{finite-resolutions},
while part~(c) follows from the dual version of the same.
\end{proof}

\begin{cor}  \label{limit-vfl-lin-dim-cor}
\textup{(a)} For any symbol\/ $\bst=\b$ or\/~$-$, the triangulated
functor\/ $\sD^\st(X\qcoh_\fl)\rarrow\sD^\st(X\qcoh_\fd)$ induced
by the embedding of exact categories $X\qcoh_\fl\rarrow X\qcoh_\fd$ is
an equivalence of triangulated categories. \par
\textup{(b)} For any symbol\/ $\bst=\b$ or\/~$-$, the triangulated
functor\/ $\sD^\st(X\qcoh_\vfl)\rarrow\sD^\st(X\qcoh_\vfd)$ induced by
the embedding of exact categories $X\qcoh_\vfl\rarrow X\qcoh_\vfd$ is
an equivalence of triangulated categories. \par
\textup{(c)} For any symbol\/ $\bst=\b$ or\/~$+$, the triangulated
functor\/ $\sD^\st(X\lcth_\bW^\lin)\rarrow\sD^\st(X\lcth_\bW^\lid)$
induced by the embedding of exact categories $X\lcth_\bW^\lin\rarrow
X\lcth_\bW^\lid$ is an equivalence of triangulated categories.
\end{cor}

\begin{proof}
 The assertions concerning the case $\bst=\b$ follow from
the respective assertions of Corollary~\ref{finite-vfl-lin-dim-cor}
by passage to the inductive limit as $d\to\infty$.
 The assertions concerning the case $\bst=-$ in parts~(a\+b)
follow from Proposition~\ref{infinite-resolutions}(a),
while the assertion about $\bst=+$ in part~(c) follows from
the dual version of it.
\end{proof}

\begin{lem} \label{flid-ctrh-dimension}
 If\/ $X=\bigcup_{\alpha=1}^N U_\alpha$ is a finite affine open covering
subordinate to\/ $\bW$, then the resolution dimension of any object of
the exact category $X\lcth_\bW^\flidd$ with respect to the resolving
subcategory $X\ctrh^\flidd$ does not exceed $N-1$.
\end{lem}

\begin{proof}
 In view of Lemma~\ref{dil-cta-clp-finite-dim}(c) and
Corollary~\ref{fdim-subcategory-cor}, it suffices to show that every
object of $X\lcth_\bW^\flidd$ admits an admissible epimorphism with
respect to the exact category $X\lcth_\bW^\flidd$ from an object of
$X\ctrh^\flidd$.
 We will do more and show that the exact
sequence~\eqref{contraherent-cech} is a resolution of an object
$\P\in X\lcth_\bW^\flidd$ by objects of $X\ctrh^\flidd$.
 Indeed, the functor of inverse image with respect to a very
flat $(\bW,\bT)$\+coaffine morphism $f\:Y\rarrow X$ takes
$X\lcth_\bW^\flidd$ into $Y\lcth_\bT^\flidd$, while the functor of
direct image with respect to a flat $(\bW,\bT)$\+affine morphism~$f$
takes $Y\lcth_\bT^\flidd$ into $X\lcth_\bW^\flidd$.
 The sequence~\eqref{contraherent-cech} is exact over
$X\lcth_\bW^\flidd$, since it is exact over $X\lcth_\bW$ and
$X\lcth_\bW^\flidd$ is closed under cokernels of admissible
monomorphisms in $X\lcth_\bW$.
\end{proof}

\begin{lem} \label{vfl-lin-finite-homol-dim}
\textup{(a)} If $X=\bigcup_{\alpha=1}^N U_\alpha$ is a finite affine
open covering, then the homological dimension of the exact category 
$X\qcoh_\fvfdd$ does not exceed $N+d$. \par
\textup{(b)} If $X=\bigcup_{\alpha=1}^N U_\alpha$ is a finite affine
open covering subordinate to\/ $\bW$, then the homological dimension of
the exact category $X\lcth_\bW^\flidd$ does not exceed $N-1+d$.
\end{lem}

\begin{proof}
 Part~(a): actually, one proves the (seemingly stronger) assertion that
$\Ext^{>N+d}_X(\F,\M)=0$ for any quasi-coherent sheaf $\M$ and any
quasi-coherent sheaf of very flat dimension~$\le d$ over~$X$
(also, the $\Ext$ groups in the exact category $X\qcoh_\fvfdd$ agree
with those in the abelian category $X\qcoh$).
 Since any object of $X\qcoh_\fvfdd$ has a finite resolution
of length~$\le d$ by objects $X\qcoh_\vfl$, it suffices to consider
the case of $\F\in X\qcoh_\vfl$, which is covered by
Lemma~\ref{vfl-cta-finite-dim}(c).
 The proof of part~(b) is similar and based on
Lemma~\ref{lct-lin-clp-finite-dim}(c). \hbadness=1100
\end{proof}

\begin{cor} \label{vfl-lin-finite-dim-all-derived-coincide} \hfuzz=3pt
 \textup{(a)} The natural triangulated functors\/
$\sD^\abs(X\qcoh_\fvfdd)\rarrow\sD^\co(X\qcoh_\fvfdd)\rarrow
\sD(X\qcoh_\fvfdd)$ and\/ $\sD^{\abs\pm}(X\qcoh_\fvfdd)\rarrow
\sD^\pm(X\qcoh_\fvfdd)$ are equivalences of triangulated categories.
 In particular, such functors between the derived categories of
the exact category $X\qcoh_\vfl$ are equivalences of categories. \par
\textup{(b)} The natural triangulated functors\/
$\sD^\abs(X\lcth_\bW^\flidd)\rarrow\sD^\ctr(X\lcth_\bW^\flidd)\rarrow
\sD(X\lcth_\bW^\flidd)$ and\/ $\sD^{\abs\pm}(X\lcth_\bW^\flidd)\rarrow
\sD^\pm(X\lcth_\bW^\flidd)$ are equivalences of triangulated categories.
 In particular, such functors between the derived categories of
the exact category $X\lcth_\bW^\lin$ are equivalences of categories.
\end{cor}

\begin{proof}
 Follows from the respective parts of
Lemma~\ref{vfl-lin-finite-homol-dim} together with
Lemma~\ref{psemi-remark21}.
 The ``in particular'' clauses are also covered by
Theorem~\ref{vlf-cta-fl-cot-derived-vlf-fl-equivalences}(a)
and Corollary~\ref{lin-ctrh-lcth-cor}(a).
\end{proof}

 As a matter of notational convenience, set the triangulated category
$\sD^\st(X\qcoh_\fd)$ to be the inductive limit of (the equivalences of
categories of) $\sD^\st(X\qcoh_\ffdd)$ as $d\to\infty$ for any symbol
$\bst\ne\ctr$, $\bco$, $\bctr$.
 Furthermore, set $\sD^\st(X\qcoh_\vfd)$ to be the inductive limit of
(the equivalences of categories of) $\sD^\st(X\qcoh_\fvfdd)$ as
$d\to\infty$.
 For any morphism $f\:Y\rarrow X$ into a quasi-compact
semi-separated scheme $X$ one constructs the left derived functor
\begin{equation}  \label{qcoh-inverse-ffd-sheaves}
 \boL f^*\:\sD^\st(X\qcoh_\fd)\lrarrow\sD^\st(Y\qcoh_\fd)
\end{equation}
as the functor on the derived categories induced by the exact functor
$f^*\:X\qcoh_\fl\rarrow Y\qcoh_\fl$. 
 So, by the definition, the functor $\boL f^*$
\eqref{qcoh-inverse-ffd-sheaves} is essentially the functor~$f^*$ from
Theorem~\ref{becker-co-contra-derived-direct-inverse-adjunction}(c).

 The left derived functor
\begin{equation}  \label{qcoh-inverse-fvfd}
 \boL f^*\:\sD^\st(X\qcoh_\vfd)\lrarrow\sD^\st(Y\qcoh_\vfd)
\end{equation}
is constructed in the similar way.
 So, by the definition, the functor $\boL f^*$
\eqref{qcoh-inverse-fvfd} is essentially the functor~$f^*$ from
Theorem~\ref{becker-co-contra-derived-direct-inverse-adjunction}(b).
 By construction, the triangulated
functors~\eqref{qcoh-inverse-ffd-sheaves}
and~\eqref{qcoh-inverse-fvfd} form a commutative diagram
with the triangulated functors
$\sD^\st(X\qcoh_\vfd)\rarrow\sD^\st(X\qcoh_\fd)$ and
$\sD^\st(Y\qcoh_\vfd)\rarrow\sD^\st(Y\qcoh_\fd)$
induced by the respective embeddings of exact categories.
{\hbadness=1275\par}

 In particular, for $\bst=-$ or~$\empt$, we obtain a commutative
diagram of triangulated functors and triangulated equivalences
provided by Theorem~\ref{derived-vfl-fl-equivalence},
\begin{equation} \label{derived-vfl-fl-inverse-images-compatible}
\begin{gathered}
 \xymatrix{
  \sD^\st(Y\qcoh_\vfd) \ar@<2pt>[r] \ar@<-2pt>@{-}[r]
  & \sD^\st(Y\qcoh_\fd)  \\
  \sD^\st(X\qcoh_\vfd) \ar@<2pt>[r] \ar@<-2pt>@{-}[r] \ar[u]^{\boL f^*}
  & \sD^\st(X\qcoh_\fd) \ar[u]^{\boL f^*}
 }
\end{gathered}
\end{equation}

 Analogously, set $\sD^\st(X\lcth_\bW^\lid)$ to be the inductive limit
of (the equivalences of categories) $\sD^\st(X\lcth_\bW^\flidd)$ as
$d\to\infty$ for any symbol $\bst\ne\co$, $\bco$, $\bctr$.
 For any morphism $f\:Y\rarrow X$ into a quasi-compact semi-separated
scheme $X$ and any open coverings $\bW$ and $\bT$ of the schemes $Y$ and
$X$ for which the morphism~$f$ is $(\bW,\bT)$\+coaffine, the right
derived functor
$$
 \boR f^!\:\sD^\st(X\lcth_\bW^\lid)\lrarrow\sD^\st(Y\lcth_\bT^\lid)
$$
is constructed as the functor on the derived categories induced by
the exact functor $f^!\:X\lcth_\bW^\lin\rarrow Y\lcth_\bT^\lin$.
 So, by the definition, the functor $\boR f^!$ is essentially
the functor~$f^!$ from
Theorem~\ref{becker-co-contra-derived-direct-inverse-adjunction}(a).

 As usually, we set $\sD^\st(X\ctrh^\lid)=\sD^\st(X\lcth_{\{X\}}^\lid)$.
 Now Lemma~\ref{flid-ctrh-dimension} together with
Proposition~\ref{finite-resolutions} provide a natural equivalence
of triangulated categories $\sD^\st(X\ctrh^\lid)\simeq
\sD^\st(X\lcth_\bW^\lid)$.
 For a morphism $f\:Y\rarrow X$ of quasi-compact semi-separated
schemes, such equivalences allow one to define the derived functor
\begin{equation} \label{ctrh-inverse-flid-cosheaves}
 \boR f^!\:\sD^\st(X\ctrh^\lid)\lrarrow\sD^\st(Y\ctrh^\lid),
\end{equation}
which clearly does not depend on the choice of the coverings $\bW$
and~$\bT$.

\subsection{Morphisms of finite flat dimension~III}
\label{finite-dim-morphisms-III}
 Let $X$ be a quasi-compact semi-separated scheme with an open
covering~$\bW$.
 Let us introduce notation for intersections of full subcategories
\begin{align*}
 X\qcoh^\cta_\ffdd &= X\qcoh^\cta\cap X\qcoh_\ffdd, \\
 X\qcoh^\dil_\ffdd &= X\qcoh^\dil\cap X\qcoh_\ffdd, \\
 X\qcoh^\cta_\fvfdd &= X\qcoh^\cta\cap X\qcoh_\fvfdd, \\
 X\qcoh^\dil_\fvfdd &= X\qcoh^\dil\cap X\qcoh_\fvfdd, \\
 X\ctrh_\al^\flidd &= X\ctrh_\al\cap X\lcth_\bW^\flidd,
\end{align*}
and notice that the latter full subcategory does not depend on
an open covering~$\bW$ (just as the notation suggests).

\begin{lem}  \label{dil-cta-clp-finite-dim-ffd-flid}
 Let $X=\bigcup_{\alpha=1}^N U_\alpha$ be a finite affine open covering.
 Then \par
\textup{(a)} the coresolution dimension of any quasi-coherent sheaf of
flat dimension\/~$\le\nobreak d$ on $X$ with respect to the coresolving 
subcategory $X\qcoh^\dil_\ffdd\sub X\qcoh_\ffdd$ does not exceed~$N-1$;
\par
\textup{(b)} the coresolution dimension of any quasi-coherent sheaf of
flat dimension\/~$\le\nobreak d$ on $X$ with respect to the coresolving
subcategory $X\qcoh^\cta_\ffdd\sub X\qcoh_\ffdd$ does not exceed~$N$;
\par
\textup{(c)} the coresolution dimension of any quasi-coherent sheaf of
very flat dimension\/~$\le\nobreak d$ on $X$ with respect to
the coresolving subcategory $X\qcoh^\dil_\fvfdd\sub X\qcoh_\fvfdd$
does not exceed~$N-1$; \par
\textup{(d)} the coresolution dimension of any quasi-coherent sheaf of
very flat dimension\/~$\le\nobreak d$ on $X$ with respect to
the coresolving subcategory $X\qcoh^\cta_\fvfdd\sub X\qcoh_\fvfdd$
does not exceed~$N$; \par
\textup{(e)} assuming that the affine open covering
$X=\bigcup_{\alpha=1}^N U_\alpha$ is subordinate to\/ $\bW$,
the resolution dimension of any\/ $\bW$\+locally contraherent cosheaf
of locally injective dimension\/~$\le\nobreak d$ on $X$ with respect to
the resolving subcategory $X\ctrh_\al^\flidd\sub X\lcth_\bW^\flidd$
does not exceed $N-1$.
\end{lem}

\begin{proof}
 Parts~(a\+b): in view of Lemma~\ref{dil-cta-clp-finite-dim}(a\+b) and
the dual version of Corollary~\ref{fdim-subcategory-cor}, it suffices
to show that there exists an injective morphism from any given
quasi-coherent sheaf belonging to $X\qcoh_\ffdd$ into a quasi-coherent
sheaf belonging to $X\qcoh^\cta_\ffdd$ with the cokernel belonging
to $X\qcoh_\ffdd$.
 This follows from Corollary~\ref{quasi-very-cta-cor}(b)
or~\ref{quasi-cotors-cor}(b).
 The proofs of parts~(c\+d) are similar.

 The proof of part~(e) is analogous up to duality, and based on
Lemma~\ref{dil-cta-clp-finite-dim}(c) and Corollary~\ref{clp-cor}(b)
(alternatively, the argument from the proof of
Lemma~\ref{flid-ctrh-dimension} is sufficient in this case).
\end{proof}

\begin{lem}  \label{cta-clp-lem-ffd-flid}
 Let $X=\bigcup_\alpha U_\alpha$ be a finite affine open covering.
Then \par
\textup{(a)} a quasi-coherent sheaf on $X$ belongs to
$X\qcoh^\cta_\ffdd$ if and only if it is a direct summand of a finitely
iterated extension of the direct images of quasi-coherent sheaves from
$U_\alpha\qcoh^\cta_\ffdd$; \par
\textup{(b)} a quasi-coherent sheaf on $X$ belongs to
$X\qcoh^\cot_\ffdd$ if and only if it is a direct summand of a finitely
iterated extension of the direct images of quasi-coherent sheaves from
$U_\alpha\qcoh^\cot_\ffdd$; \par
\textup{(c)} a quasi-coherent sheaf on $X$ belongs to
$X\qcoh^\cta_\fvfdd$ if and only if it is a direct summand of a finitely
iterated extension of the direct images of quasi-coherent sheaves from
$U_\alpha\qcoh^\cta_\fvfdd$; \par
\textup{(d)} a contraherent cosheaf on $X$ belongs to
$X\ctrh_\al^\flidd$ if and only if it is a direct summand of a finitely
iterated extension of the direct images of contraherent cosheaves
from $U_\alpha\ctrh^\flidd$.
\end{lem}

\begin{proof}
 One can repeat the arguments in
Sections~\ref{quasi-compact-quasi-coherent}
and~\ref{clp-subsection} working with, respectively, quasi-coherent
sheaves of (very) flat dimension~$\le d$ only or locally contraherent
cosheaves of locally injective dimension~$\le d$ only throughout
(cf.\ Lemma~\ref{flat-contraadjusted-quasi}).
 Alternatively, one can refer to
Theorem~\ref{quasi-coherent-gluing-theorem} for parts~(a\+c) and to
Theorem~\ref{loc-contraherent-gluing-theorem} for part~(d).
 Let us spell out this argument for parts~(a) and~(d).

 Part~(a): consider the hereditary complete cotorsion pair
$(\sF,\sC)=(X\qcoh_\vfl$, $X\qcoh^\cta$) from
Corollary~\ref{quasi-very-cta-cor}(a\+b) (see also
Lemma~\ref{cotorsion-pair-direct-summands-lemma}) in the abelian
category $\sK=X\qcoh$, and consider the exact subcategory
$\sE=X\qcoh_\ffdd\sub\sK$ of quasi-coherent sheaves of flat
dimension~$\le d$.
 By Lemmas~\ref{restricting-hereditary-cotorsion}
and~\ref{restricting-cotorsion-pairs-lemma}(a), the cotorsion pair
$(\sF,\sC)$ restricts to a hereditary complete cotorsion pair in~$\sE$.
 So the pair of classes $\sF=X\qcoh_\vfl$ and $\sE\cap\sC=
X\qcoh^\cta_\ffdd$ is a hereditary complete cotorsion pair in~$\sE$.
 In the rest of the argument, we redenote $\sE\cap\sC$ by~$\sC$.

 Let $\sR$ be the local class of all commutative rings $R$ and
$\sE_R=R\modl_\ffdd\sub\sK_R=R\modl$ be the exact subcategory of
$R$\+modules of flat dimension~$\le d$.
 Put $\sF_R=R\modl_\vfl$ and $\sC(R)=R\modl^\cta_\ffdd=R\modl^\cta
\cap R\modl_\ffdd$.
 Then the class $\sF$ is very local by
Examples~\ref{local-classes-examples}, and it is easy to check
that the class $\sE$ is very local, too (see
Lemma~\ref{flat-veryflat-cotors-inj-dim-local}(a)).
 The pair of classes $(\sF_R,\sC(R))$ is a hereditary complete
cotorsion pair in~$\sE_R$, as we have just seen.

 Applying Theorem~\ref{quasi-coherent-gluing-theorem} to the datum
of the classes $\sR$, \,$\sE_R$, \,$\sF_R$, and $\sC(R)$, we obtain
a cotorsion pair $(\sF_X,\sC(X))$ in the abelian category $\sE_X$ of
locally\+$\sE$ quasi-coherent sheaves on~$X$.
 By the definitions, one has $\sE_X=X\qcoh_\ffdd$ and
$\sF_X=X\qcoh_\vfl$.
 The theorem tells us that $\sC(X)$ is the class of all direct summands
of finitely iterated extensions of direct images of quasi-coherent
sheaves from $U_\alpha\qcoh^\cta_\ffdd$.
 Comparing the two cotorsion pairs in $X\qcoh_\ffdd$ that we have
constructed, we conclude that $X\qcoh^\cta_\ffdd=\sC(X)$, as desired.

 Part~(d): let $\bW$ be any open covering of $X$ to which the covering
$X=\bigcup_\alpha U_\alpha$ is subordinate.
 Consider the hereditary complete cotorsion pair $(\sF,\sC)=
(X\ctrh_\al$, $X\lcth_\bW^\lin$) from Corollary~\ref{clp-cor}(a\+b)
(see also Lemma~\ref{cotorsion-pair-direct-summands-lemma}) in
the exact category $\sK=X\lcth_\bW$, and consider the exact
subcategory $\sE=X\lcth_\bW^\flidd\sub\sK$ of $\bW$\+locally
contraherent cosheaves of locally injective dimension~$\le d$.
 By Lemmas~\ref{restricting-hereditary-cotorsion}
and~\ref{restricting-cotorsion-pairs-lemma}(b), the cotorsion pair
$(\sF,\sC)$ restricts to a hereditary complete cotorsion pair in~$\sE$.
 So the pair of classes $\sE\cap\sF=X\ctrh_\al^\flidd$ and
$\sC=X\lcth_\bW^\lin$ is a hereditary complete cotorsion pair in~$\sE$.
 In the rest of the argument, we redenote $\sE\cap\sF$ by~$\sF$.
 
 Let $\sR$ be the local class of all commutative rings $R$ and
$\sE^R=R\modl^\fidd\cap R\modl^\cta\sub\sK^R=R\modl^\cta$ be
the exact subcategory of contraadjusted $R$\+modules whose injective
dimension does not exceed~$d$.
 Here the notation $R\modl^\fidd\sub R\modl$ stands for the full
subcategory of all $R$\+modules of injective dimension~$\le d$.
 Put $\sF(R)=\sE^R$ and $\sC^R=R\modl^\inj$.
 There are enough injective objects in the exact category $\sE^R$,
so $(\sF(R),\sC^R)$ is a hereditary complete cotorsion pair in~$\sE^R$.
 The class $\sC^R$ is very colocal by
Examples~\ref{colocal-classes-examples}, and one can check that
the class $\sE^R$ is very colocal, too (see
Lemma~\ref{flat-veryflat-cotors-inj-dim-local}(d)).

 Applying Theorem~\ref{loc-contraherent-gluing-theorem} to the datum
of the classes $\sR$, \,$\sE^R$, \,$\sC^R$, and $\sF(R)$, we obtain
a hereditary complete cotorsion pair $(\sF(X),\sC^X_\bW)$ in
the exact category $\sE^X_\bW=X\lcth_\bW^\flidd$.
 Here $\sF(X)$ is the class of all direct summands of finitely
iterated extensions of direct images of contraherent cosheaves
from $U_\alpha\ctrh^\flidd$, while $\sC^X_\bW=X\lcth_\bW^\lin$ is
the class of locally injective $\bW$\+locally contraherent cosheaves.
 Comparing the two cotorsion pairs in $X\lcth_\bW^\flidd$ that we have
constructed, we conclude that $X\ctrh_\al^\flidd=\sF(X)$, as desired.
\end{proof}

 Let $f\:Y\rarrow X$ be a morphism of finite flat dimension~$\le D$
between quasi-compact semi-separated schemes $X$ and~$Y$.

\begin{cor}  \label{dil-cta-clp-direct-ffd-flid}
\textup{(a)} The exact functor $f_*\:Y\qcoh^\dil\rarrow X\qcoh^\dil$
takes objects of $Y\qcoh^\dil_\ffdd$ to objects of
$X\qcoh^\dil_\ffdD$. \par
\textup{(b)} The exact functor $f_*\:Y\qcoh^\cta\rarrow X\qcoh^\cta$
takes objects of $Y\qcoh^\cta_\ffdd$ to objects of
$X\qcoh^\cta_\ffdD$. \par
\textup{(c)} If the morphism~$f$ has very flat dimension not
exceeding~$D$, then the exact functor $f_*\:Y\qcoh^\dil\rarrow
X\qcoh^\dil$ takes objects of $Y\qcoh^\dil_\fvfdd$ to
objects of $X\qcoh^\dil_\fvfdD$. \par
\textup{(d)} If the morphism~$f$ has very flat dimension not
exceeding~$D$, then the exact functor $f_*\:Y\qcoh^\cta\rarrow
X\qcoh^\cta$ takes objects of $Y\qcoh^\cta_\fvfdd$ to
objects of $X\qcoh^\cta_\fvfdD$. \par
\textup{(e)} The exact functor $f_!\:Y\ctrh_\al\rarrow X\ctrh_\al$
takes objects of $Y\ctrh_\al^\flidd$ to objects of
$X\ctrh_\al^\flidD$.
\end{cor}

\begin{proof}
 Part~(b) follows from Lemma~\ref{cta-clp-lem-ffd-flid}(a) together
with the fact that the direct image with respect to an affine
morphism of flat dimension~$\le D$ takes quasi-coherent sheaves of
flat dimension~$\le d$ to quasi-coherent sheaves of
flat dimension~$\le d+D$.
 The latter is provided by Lemma~\ref{affine-d+D-lemma}(a).
 The proof of part~(d) is similar and based on
Lemmas~\ref{cta-clp-lem-ffd-flid}(c) and~\ref{affine-very-d+D-lemma}(b).

 Let us explain how to deduce part~(a) from part~(b).
 Let $\N\in Y\qcoh^\dil_\ffdd$ be a dilute quasi-coherent sheaf of
flat dimension~$\le d$ on~$Y$.
 By Lemma~\ref{dil-cta-clp-finite-dim-ffd-flid}(b), there exists
an exact sequence $0\rarrow\N\rarrow\cP^0\rarrow\cP^1\rarrow
\dotsb\rarrow\cP^N\rarrow0$ in $Y\qcoh_\ffdd$ with $\cP^n\in
X\qcoh^\cta_\ffdd$ for all $0\le n\le N$.
 Clearly, the sequence $0\rarrow\N\rarrow\cP^0\rarrow\cP^1\rarrow
\dotsb\rarrow\cP^N\rarrow0$ is also exact in $Y\qcoh^\dil_\ffdd$.
 By Corollary~\ref{dilute-direct} and part~(b), the sequence
of direct images $0\rarrow f_*\N\rarrow f_*\cP^0\rarrow f_*\cP^1\rarrow
\dotsb\rarrow f_*\cP^N\rarrow0$ is exact in $X\qcoh^\dil$ with
$f_*\cP^n\in X\qcoh^\cta_\ffdd$.
 Applying Corollary~\ref{fdim-resolution}(b), we conclude that
$f_*\N\in X\qcoh_\ffdd$.
 Similarly one uses
Lemma~\ref{dil-cta-clp-finite-dim-ffd-flid}(d) in order to deduce
part~(c) from part~(d).

 Finally, part~(e) follows from Lemma~\ref{cta-clp-lem-ffd-flid}(d)
together with the fact that the direct image with respect to
a $(\bW,\bT)$\+affine morphism of flat dimension~$\le D$ takes
$\bT$\+locally contraherent cosheaves of locally injective
dimension~$\le d$ to $\bW$\+locally contraherent cosheaves of locally
injective dimension~$\le d+D$.
 The latter is provided by Lemma~\ref{affine-d+D-lemma}(b).
\end{proof}

 According to Lemma~\ref{dil-cta-clp-finite-dim-ffd-flid}(b) and
the dual version of Proposition~\ref{finite-resolutions}, for any
symbol $\bst\ne\co$, $\ctr$, $\bco$, $\bctr$ the natural functor
$\sD^\st(Y\qcoh^\cta_\ffdd)\rarrow\sD^\st(Y\qcoh_\ffdd)$
is an equivalence of triangulated categories (as is the similar
functor for sheaves over~$X$).
 So one can construct the right derived functor
\begin{equation} \label{qcoh-direct-ffdd-ffdD}
 \boR f_*\: \sD^\st(Y\qcoh_\ffdd)\lrarrow\sD^\st(X\qcoh_\ffdD)
\end{equation}
as the functor on the derived categories induced by the exact functor
$f_*\:Y\qcoh^\cta_\ffdd\allowbreak\rarrow X\qcoh^\cta_\ffdD$
from Corollary~\ref{dil-cta-clp-direct-ffd-flid}(b).

 More generally, Lemma~\ref{dil-cta-clp-finite-dim-ffd-flid}(a) and
the dual version of Proposition~\ref{finite-resolutions} tell us that
for any symbol $\bst\ne\ctr$, $\bco$, $\bctr$ the natural functor
$\sD^\st(Y\qcoh^\dil_\ffdd)\rarrow\sD^\st(Y\qcoh_\ffdd)$
is an equivalence of triangulated categories (as is the similar
functor for sheaves over~$X$).
 This allows us to construct the right derived functor~$\boR f_*$
\eqref{qcoh-direct-ffdd-ffdD} as the functor on the derived categories
induced by the exact functor $f_*\:Y\qcoh^\dil_\ffdd\allowbreak\rarrow
X\qcoh^\dil_\ffdD$ from Corollary~\ref{dil-cta-clp-direct-ffd-flid}(a).
 Passing to the inductive limits as $d\to\nobreak\infty$, we obtain
the right derived functor
\begin{equation}  \label{qcoh-direct-ffd}
 \boR f_*\:\sD^\st(Y\qcoh_\fd)\lrarrow\sD^\st(X\qcoh_\fd),
\end{equation}
which is right adjoint to the functor $\boL f^*$
\eqref{qcoh-inverse-ffd-sheaves}.

 Comparing this adjunction for $\bst=\empt$ with the adjunction of
Theorem~\ref{becker-co-contra-derived-direct-inverse-adjunction}(c),
we arrive to a commutative diagram of triangulated functors and
triangulated equivalences from
Corollary~\ref{bctr-lcth-vfl-fl-derived-equivalences},
\begin{equation} \label{ffd-bctr-lcth-lct-direct-images-agree}
\begin{gathered}
 \xymatrix{
  \sD(Y\qcoh_\fd) \ar@{=}[r] \ar[d]_{\boR f_*}
  & \sD^\bctr(Y\ctrh^\lct) \ar[d]^{\boL f_!} \\
  \sD(X\qcoh_\fd) \ar@{=}[r]
  & \sD^\bctr(X\ctrh^\lct)
 }
\end{gathered}
\end{equation}
for any morphism $f\:Y\rarrow X$ of finite flat dimension.

 For a morphism $f$ of finite very flat dimension, the right derived
functor
\begin{equation}  \label{qcoh-direct-fvfd}
 \boR f_*\:\sD^\st(Y\qcoh_\vfd)\lrarrow\sD^\st(X\qcoh_\vfd)
\end{equation}
right adjoint to the functor $\boL f^*$ \eqref{qcoh-inverse-fvfd}
is constructed (for any symbol $\bst\ne\ctr$, $\bco$, $\bctr$)
in the similar way.
 Comparing this adjunction for $\bst=\empt$ with the adjunction of
Theorem~\ref{becker-co-contra-derived-direct-inverse-adjunction}(b),
we arrive to a commutative diagram of triangulated functors and
triangulated equivalences from
Corollary~\ref{bctr-lcth-vfl-fl-derived-equivalences},
\begin{equation} \label{fvfd-bctr-lcth-direct-images-agree}
\begin{gathered}
 \xymatrix{
  \sD(Y\qcoh_\vfd) \ar@{=}[r] \ar[d]_{\boR f_*}
  & \sD^\bctr(Y\ctrh) \ar[d]^{\boL f_!} \\
  \sD(X\qcoh_\vfd) \ar@{=}[r]
  & \sD^\bctr(X\ctrh)
 }
\end{gathered}
\end{equation}
for any morphism $f\:Y\rarrow X$ of finite very flat dimension.

 By construction, the triangulated
functors~\eqref{qcoh-direct-ffd} and~\eqref{qcoh-direct-fvfd} form
a commutative diagram with the triangulated functors
$\sD^\st(Y\qcoh_\vfd)\rarrow\sD^\st(Y\qcoh_\fd)$ and
$\sD^\st(X\qcoh_\vfd)\rarrow\sD^\st(X\qcoh_\fd)$
induced by the respective embeddings of exact categories.
 In particular, for $\bst=-$ or~$\empt$, we obtain a commutative
diagram of triangulated functors and triangulated equivalences
provided by Theorem~\ref{derived-vfl-fl-equivalence},
\begin{equation} \label{derived-vfl-fl-direct-images-compatible}
\begin{gathered}
 \xymatrix{
  \sD^\st(Y\qcoh_\vfd) \ar@<2pt>[r] \ar@<-2pt>@{-}[r] \ar[d]_{\boR f_*}
  & \sD^\st(Y\qcoh_\fd) \ar[d]_{\boR f_*} \\
  \sD^\st(X\qcoh_\vfd) \ar@<2pt>[r] \ar@<-2pt>@{-}[r]
  & \sD^\st(X\qcoh_\fd) 
 }
\end{gathered}
\end{equation}

 Analogously, according to
Lemma~\ref{dil-cta-clp-finite-dim-ffd-flid}(e)
and Proposition~\ref{finite-resolutions}, for any symbol $\bst\ne\co$,
$\bco$, $\bctr$ the natural functor $\sD^\st(Y\ctrh_\al^\flidd)\rarrow
\sD^\st(Y\ctrh^\flidd)$ is an equivalence of triangulated categories
(as is the similar functor for cosheaves over~$X$).
 Thus one can construct the left derived functor
$$
 \boL f_!\: \sD^\st(Y\ctrh^\flidd)\lrarrow\sD^\st(X\ctrh^\flidD)
$$
for a morphism $f\:Y\rarrow X$ of finite flat dimension~$\le D$
as the functor on the derived categories induced by the exact functor
$f_!\:Y\ctrh_\al^\flidd\allowbreak\rarrow X\ctrh_\al^\flidD$
from Corollary~\ref{dil-cta-clp-direct-ffd-flid}(e).
 Passing to the inductive limits as $d\to\nobreak\infty$, we obtain
the left derived functor
\begin{equation}  \label{ctrh-direct-flid}
 \boL f_!\:\sD^\st(Y\ctrh^\lid)\lrarrow\sD^\st(X\ctrh^\lid),
\end{equation}
which is left adjoint to the functor $\boR f^!$
\eqref{ctrh-inverse-flid-cosheaves}.

 Comparing this adjunction for $\bst=\empt$ with the adjunction of
Theorem~\ref{becker-co-contra-derived-direct-inverse-adjunction}(a),
we arrive to a commutative diagram of triangulated functors and
triangulated equivalences from
Corollary~\ref{bco-qcoh-lin-derived-equivalence},
\begin{equation} \label{bco-qcoh-lin-direct-images-agree}
\begin{gathered}
 \xymatrix{
  \sD^\bco(Y\qcoh) \ar@{=}[r] \ar[d]_{\boR f_*}
  & \sD(Y\ctrh^\lid)  \ar[d]^{\boL f_!} \\
  \sD^\bco(X\qcoh) \ar@{=}[r] & \sD(X\ctrh^\lid)
 }
\end{gathered}
\end{equation}
for any morphism $f\:Y\rarrow X$ of finite flat dimension.

\subsection{Finite injective and projective dimension}
\label{finite-inj-proj-dim-subsect}
 For any scheme $X$, we denote the full subcategory of objects of
injective dimension~$\le d$ in the abelian category $X\qcoh$ by
$X\qcoh^\fidd$.
 For a quasi-compact semi-separated scheme $X$, the full subcategory
of objects of projective dimension~$\le d$ in the exact category
$X\lcth_\bW$ is denoted by $X\lcth_{\bW,\,\fpdd}$ and the full
subcategory of objects of projective dimension~$\le d$ in the exact
category $X\lcth_\bW^\lct$ by $X\lcth^\lct_{\bW,\,\fpdd}$.
 Set $X\ctrh_\fpdd=X\lcth_{\{X\},\,\fpdd}$ and $X\ctrh^\lct_\fpdd=
X\lcth^\lct_{\{X\},\,\fpdd}$.

 Furthermore, the \emph{antilocally flat dimension} of a $\bW$\+locally
contraherent cosheaf on a quasi-compact semi-separated scheme $X$ is
defined as its resolution dimension with respect to the resolving
subcategory $X\ctrh_\alf\sub X\lcth_\bW$
(see Section~\ref{finite-resolutions-subsect}).
 The full subcategory $X\ctrh_\alf$ is resolving in $X\lcth_\bW$ by
Corollaries~\ref{clf-characterizations}(b) and~\ref{clf-cor}(b).
 A $\bW$\+locally contraherent cosheaf $\gF$ on $X$ has antilocally
flat dimension~$\le d$ if and only if $\Ext^{>d}(\gF,\P)=0$ for all
locally cotorsion $\bW$\+locally contraherent cosheaves $\P$ on~$X$.
 The full subcategory of objects of antilocally flat dimension~$\le d$
in $X\lcth_\bW$ is denoted by $X\lcth_{\bW,\,\alfdd}$.
 We set $X\ctrh_\alfdd=X\lcth_{\{X\},\,\alfdd}$.

 Since the full subcategories of projective objects in $X\lcth_\bW$ and
$X\lcth_\bW^\lct$ do not depend on the covering $\bW$, and
the full subcategories $X\lcth_\bW\sub X\lcth$ and $X\lcth_\bW^\lct
\sub X\lcth^\lct$ are closed under kernels of admissible epimorphisms,
the projective dimension of an object of $X\lcth_\bW$ or
$X\lcth_\bW^\lct$ does not change when the open covering $\bW$ is
replaced by its refinement.
 Similarly, the antilocally flat dimension of a $\bW$\+locally
contraherent cosheaf on $X$ does not depend on the covering~$\bW$.

 One can easily see that the full subcategory $X\qcoh^\fidd\sub X\qcoh$
is closed under extensions and cokernels of admissible monomorphisms,
while the full subcategories $X\lcth^\lct_{\bW,\,\fpdd}\sub
X\lcth_\bW^\lct$, \ $X\lcth_{\bW,\,\fpdd}\sub X\lcth_\bW$,
and $X\lcth_{\bW,\,\alfdd}\sub X\lcth_\bW$ are closed
under extensions and kernels of admissible epimorphisms.

\begin{cor} \label{derived-fid-fpd-cor}  \hbadness=1800 \hfuzz=13pt
\textup{(a)} For any scheme $X$, the natural triangulated functors\/
$\Hot(X\qcoh^\inj)\rarrow\sD^{\abs=\bco}(X\qcoh^\fidd)\rarrow
\sD(X\qcoh^\fidd)$, \ $\Hot^\pm(X\qcoh^\inj)\rarrow
\sD^{\abs\pm}(X\qcoh^\fidd)\rarrow\sD^\pm(X\qcoh^\fidd)$, and\/
$\Hot^\b(X\qcoh^\inj)\rarrow \sD^\b(X\qcoh^\fidd)$ are equivalences
of categories. \par
\textup{(b)} For any quasi-compact semi-separated scheme $X$,
the natural triangulated functors\/ $\Hot(X\ctrh_\prj)\rarrow
\sD^{\abs=\bctr}(X\lcth_{\bW,\,\fpdd})\rarrow
\sD(X\lcth_{\bW,\,\fpdd})$, \ $\Hot^\pm(X\ctrh_\prj)\rarrow
\sD^{\abs\pm}(X\lcth_{\bW,\,\fpdd})\rarrow
\sD^\pm(X\lcth_{\bW,\,\fpdd})$, and\/ $\Hot^\b(X\ctrh_\prj)\allowbreak
\rarrow\sD^\b(X\lcth_{\bW,\,\fpdd})$ are equivalences of categories.
\par
\textup{(c)} For any quasi-compact semi-separated scheme $X$ and any
symbol\/ $\bst\ne\co$, $\ctr$, $\bco$, the natural triangulated
functors\/ $\sD^\st(X\ctrh_\alf)\rarrow\sD^\st(X\lcth_{\bW,\,\alfdd})$
are equivalences of categories. \par
\textup{(d)} For any quasi-compact semi-separated scheme $X$,
the natural triangulated functors\/ $\Hot(X\ctrh^\lct_\prj)\rarrow
\sD^{\abs=\bctr}(X\lcth^\lct_{\bW,\,\fpdd})\rarrow
\sD(X\lcth^\lct_{\bW,\,\fpdd})$, \ $\Hot^\pm(X\ctrh^\lct_\prj)\rarrow
\sD^{\abs\pm}(X\lcth^\lct_{\bW,\,\fpdd})\rarrow
\sD^\pm(X\lcth^\lct_{\bW,\,\fpdd})$, and\/ $\Hot^\b(X\ctrh^\lct_\prj)
\allowbreak\rarrow\sD^\b(X\lcth^\lct_{\bW,\,\fpdd})$ are equivalences
of categories.
\end{cor}

\begin{proof}
 All the assertions of parts~(b\+d) follow from
Propositions~\ref{finite-resolutions}
and~\ref{becker-contraderived-finite-resolutions}, while part~(a)
follows from the dual versions of these.
 See also Theorem~\ref{finite-homol-dim-becker-co-contra-derived}.
\end{proof}

 Let $X$ be a quasi-compact semi-separated scheme.

\begin{lem}  \label{cta-clp-finite-flat-inj-dim-identified}
 The equivalence of exact categories $X\qcoh^\cta\simeq X\ctrh_\al$
from Lemma~\textup{\ref{cta-clp-equivalence}} identifies
the exact subcategories \par
\textup{(a)} $X\qcoh^\cta\cap X\qcoh^\fidd\sub X\qcoh^\cta$ with
$X\ctrh_\al^\flidd\sub X\ctrh_\al$, \par
\textup{(b)} $X\qcoh^\cta_\fvfdd\sub X\qcoh^\cta$ with
$X\ctrh_\al\cap X\ctrh_\fpdd\sub X\ctrh_\al$, \par
\textup{(c)} $X\qcoh^\cta_\ffdd\sub X\qcoh^\cta$ with
$X\ctrh_\al\cap X\ctrh_\alfdd \sub X\ctrh_\al$, \par
\textup{(d)} $X\qcoh^\cot_\ffdd\sub X\qcoh^\cta$ with
$X\ctrh_\al\cap X\ctrh^\lct_\fpdd\sub X\ctrh_\al$.
\end{lem}

\begin{proof}
 Part~(a): since the functor $\O_X\ocn_X{-}$ takes short exact sequences
in $X\ctrh_\al$ to short exact sequences in $X\qcoh^\cta$ and commutes
with the direct images of contraherent cosheaves and quasi-coherent
sheaves from the affine open subschemes of $X$, it follows from
Lemma~\ref{cta-clp-lem-ffd-flid}(d) that this functor takes
$X\ctrh_\al^\flidd$ into $X\qcoh^\cta\cap X\qcoh^\fidd$.
 To prove the converse, consider a contraadjusted quasi-coherent
sheaf $\cP$ of injective dimension~$d$ on $X$, and let
$0\rarrow\cP\rarrow\J^0\rarrow\dotsb\rarrow\J^d\rarrow0$ be its
injective coresolution of length~$d$ in $X\qcoh$.

 This coresolution is an exact sequence in the exact category
$X\qcoh^\cta$, so it is transformed to an exact sequence in the exact
category $X\ctrh_\al$ by the functor $\fHom_X(\O_X,{-})$.
 The contraherent cosheaves $\fHom_X(\O_X,\J^i)$ being locally
injective according to Lemma~\ref{cta-clp-restricts-to-cot-inj}(b),
it follows that the contraherent cosheaf $\fHom_X(\O_X,\cP)$ admits
a coresolution of length~$d$ by locally injective contraherent
cosheaves, i.~e., has locally injective dimension~$\le\nobreak d$.

 Part~(d): according to Lemma~\ref{cta-clp-restricts-to-cot-inj}(a),
the equivalence of exact categories $X\qcoh^\cta\simeq X\ctrh_\al$
identifies $X\qcoh^\cot$ with $X\ctrh^\lct_\al$.
 Furthermore, it follows from Lemma~\ref{cta-clp-lem-ffd-flid}(b)
that the functor $\fHom_X(\O_X,{-})$ takes $X\qcoh^\cot_\ffdd$ into
$X\ctrh_\al\cap X\ctrh^\lct_\fpdd$, since
a cotorsion module of flat dimension~$\le d$ over a commutative
ring~$R$ corresponds to a (locally) cotorsion contraherent cosheaf
of projective dimension~$\le d$ over $\Spec R$.
 Conversely, a projective resolution of length~$d$ of an object of
$X\ctrh^\lct_\al$ is transformed by the functor $\O_X\ocn_X{-}$
into a flat resolution of length~$d$ of the corresponding cotorsion
quasi-coherent sheaf (see Lemma~\ref{cta-clp-restricts-to-prj-clf}(b)).

 The proof of part~(b) is similar and based on
Lemmas~\ref{cta-clp-lem-ffd-flid}(c)
and~\ref{cta-clp-restricts-to-prj-clf}(a), while 
the proof of part~(c) is based on Lemmas~\ref{cta-clp-lem-ffd-flid}(a)
and~\ref{cta-clp-restricts-to-prj-clf}(c).
 In both cases, it is important that every projective (or, respectively,
antilocally flat) contraherent cosheaf is antilocal.
\end{proof}

\begin{cor}  \label{cta-clp-cor-fid-fpd}
 Let $X=\bigcup_\alpha U_\alpha$ be a finite affine open covering.
 Then \par
\textup{(a)} a quasi-coherent sheaf on $X$ belongs to
$X\qcoh^\cta\cap X\qcoh^\fidd$ if and only if it is a direct summand
of a finitely iterated extension of the direct images of
quasi-coherent sheaves from
$U_\alpha\qcoh^\cta\cap U_\alpha\qcoh^\fidd$; \par
\textup{(b)} a contraherent cosheaf on $X$ belongs to
$X\ctrh_\al\cap X\ctrh_\fpdd$ if and only if it is a direct summand
of a finitely iterated extension of the direct images of
contraherent cosheaves from
$U_\alpha\ctrh_\al\cap U_\alpha\ctrh_\fpdd$; \par
\textup{(c)} a contraherent cosheaf on $X$ belongs to
$X\ctrh_\al\cap X\ctrh_\alfdd$ if and only if it is a direct summand
of a finitely iterated extension of the direct images of
contraherent cosheaves from
$U_\alpha\ctrh_\al\cap U_\alpha\ctrh_\alfdd$; \par
\textup{(d)} a contraherent cosheaf on $X$ belongs to
$X\ctrh_\al\cap X\ctrh^\lct_\fpdd$ if and only if it is a direct 
summand of a finitely iterated extension of the direct images of
contraherent cosheaves from
$U_\alpha\ctrh_\al\cap U_\alpha\ctrh^\lct_\fpdd$.
\end{cor}

\begin{proof}
 Follows from Lemmas~\ref{cta-clp-lem-ffd-flid}
and~\ref{cta-clp-finite-flat-inj-dim-identified}.
\end{proof}

 Let $f\:Y\rarrow X$ be a morphism of quasi-compact semi-separated
schemes.

\begin{lem}  \label{finite-flat-dim-inj-proj-direct}
\textup{(a)} Whenever the flat dimension of the morphism~$f$ does not
exceed~$D$, the functor of direct image~$f_*$ takes injective
quasi-coherent sheaves on $Y$ to quasi-coherent sheaves of
injective dimension~$\le D$ on~$X$. \par
\textup{(b)} Whenever the very flat dimension of the morphism~$f$ does
not exceed~$D$, the functor of direct image~$f_!$ takes projective
contraherent cosheaves on $Y$ to contraherent cosheaves of
projective dimension\/~$\le D$ on~$X$. \par
\textup{(c)} Whenever the flat dimension of the morphism~$f$ does
not exceed~$D$, the functor of direct image~$f_!$ takes antilocally
flat contraherent cosheaves on $Y$ to contraherent cosheaves of
antilocally flat dimension\/~$\le D$ on~$X$. \par
\textup{(d)} Whenever the flat dimension of the morphism~$f$ does
not exceed~$D$, the functor of direct image~$f_!$ takes projective
locally cotorsion contraherent cosheaves on $Y$ to locally cotorsion
contraherent cosheaves of projective dimension\/~$\le D$ on~$X$.
\end{lem}

\begin{proof}
 In part~(a) it actually suffices to assume that the scheme $Y$
is quasi-compact and quasi-separated (while the scheme $X$ has to be
quasi-compact and semi-separated).
 Let us prove parts~(b) and~(c), parts~(a) and (d) being analogous.

 Part~(b): the functor~$f_!$ takes $Y\ctrh_\prj$ to $X\ctrh$ by
Corollary~\ref{clp-direct}(a), so it remains to check that
$\Ext^{X,>D}(f_!\gF,\Q)=0$ for all cosheaves $\gF\in Y\ctrh_\prj$
and $\Q\in X\ctrh$.
 According to the adjunction of derived functors
$\boL f_!$~\eqref{ctrh-direct} and
$\boR f^!$~\eqref{ctrh-inverse-ffd-morphism} for $\bst=\b$,
one has $\Ext^{X,*}(f_!\gF,\Q)\simeq
\Hom_{\sD^\b(Y\lcth)}(\gF\;\boR f^!(\Q)[*])$.
 Hence it suffices to show that the object $\boR f^!(\Q)\in
\sD^\b(Y\lcth)$ can be represented by a finite complex over
$Y\lcth$ concentrated in the cohomological degrees~$\le D$.
 The latter is true because any $\bW$\+locally contraherent cosheaf
on $X$ admits a finite coresolution of length~$\le D$ by $f$\+adjusted
$\bW$\+locally contraherent cosheaves.

 Part~(c): the functor~$f_!$ takes $Y\ctrh_\alf$ to $X\ctrh$ by
Corollary~\ref{clp-direct}(a); hence it remains to check that
$\Ext^{X,>D}(f_!\gF,\P)=0$ for all cosheaves $\gF\in Y\ctrh_\alf$
and $\P\in X\ctrh^\lct$.
 For this purpose, we need to use a partial adjunction of the derived
functors $\boL f_!$~\eqref{ctrh-direct} and
$\boR f^!$~\eqref{ctrh-lct-inverse-ffd-morphism} for $\bst=\b$.
 The existence of such partial adjunction is established in the way
similar to the other (full) adjunctions mentioned in
Section~\ref{finite-dim-morphisms-I}.
 It is helpful to keep in mind that $X\lcth_\bW^\lct$ is a coresolving
subcategory in $X\lcth_\bW$.
 So one has $\Ext^{X,*}(f_!\gF,\P)\simeq
\Hom_{\sD^\b(Y\lcth)}(\gF\;\boR f^!(\P)[*])$.

 Hence it suffices to show that the object $\boR f^!(\P)\in
\sD^\b(Y\lcth^\lct)$ can be represented by a finite complex over
$Y\lcth^\lct$ concentrated in the cohomological degrees~$\le D$.
 The latter is true because any locally cotorsion $\bW$\+locally
contraherent cosheaf on $X$ admits a finite coresolution of
length~$\le D$ by $f$\+adjusted locally cotorsion $\bW$\+locally
contraherent cosheaves.
\end{proof}

 Set the triangulated category $\sD^\st(X\qcoh^\fid)$ to be
the inductive limit of (the equivalences of categories of)
$\sD^\st(X\qcoh^\fidd)$ as $d\to\infty$ for any $\bst\ne\co$, $\ctr$,
$\bctr$.
 For any morphism of finite flat dimension $f\:Y\rarrow X$ between
quasi-compact semi-separated schemes, one constructs the right
derived functor
\begin{equation} \label{qcoh-direct-fid}
 \boR f_*\:\sD^\st(Y\qcoh^\fid)\lrarrow\sD^\st(X\qcoh^\fid)
\end{equation}
as the functor on the homotopy/derived categories induced by
the additive functor $f_*\:Y\qcoh^\inj\rarrow X\qcoh^\fiD$
from Lemma~\ref{finite-flat-dim-inj-proj-direct}(a).
 By construction, the triangulated
functors~\eqref{qcoh-direct-fid} for $\bst=\empt$ or~$\bco$
and~\eqref{bco-qcoh-direct} agree with
each other, and we obtain a commutative diagram of triangulated
functors and triangulated equivalences provided by
Theorem~\ref{quasi-coherent-becker-coderived},
\begin{equation} \label{bco-fid-direct-images-compatible}
\begin{gathered}
 \xymatrix{
  \sD(Y\qcoh^\fid) \ar@<2pt>[r] \ar@<-2pt>@{-}[r]
  \ar[d]_{\boR f_*}
  & \sD^\bco(Y\qcoh) \ar[d]_{\boR f_*} \\
  \sD(X\qcoh^\fid) \ar@<2pt>[r] \ar@<-2pt>@{-}[r]
  & \sD^\bco(X\qcoh)
 }
\end{gathered}
\end{equation}

 Analogously, set $\sD^\st(X\ctrh_\fpd)$ to be the inductive limit
of (the equivalences of categories of) $\sD^\st(X\ctrh_\fpdd)$ as
$d\to\infty$ for any $\bst\ne\co$, $\ctr$, $\bco$.
 Besides, set $\sD^\st(X\ctrh_\alfd)$ to be the inductive limit of
(the equivalences of categories of) $\sD^\st(X\ctrh_\alfdd)$ as
$d\to\infty$, and $\sD^\st(X\ctrh^\lct_\fpd)$ to be the inductive
limit of (the equivalences of categories of)
$\sD^\st(X\ctrh^\lct_\fpdd)$ as $d\to\infty$.
 For any morphism of finite flat dimension $f\:Y\rarrow X$
between quasi-compact semi-separated schemes, one constructs
the left derived functor
\begin{equation}  \label{ctrh-direct-lct-fpd}
 \boL f_!\:\sD^\st(Y\ctrh^\lct_\fpd)\lrarrow\sD^\st(X\ctrh^\lct_\fpd)
\end{equation}
as the functor on the homotopy/derived categories induced by
the additive functor $f_!\:Y\ctrh^\lct_\prj\rarrow X\ctrh^\lct_\fpD$
from Lemma~\ref{finite-flat-dim-inj-proj-direct}(d).
 The left derived functor
\begin{equation}  \label{ctrh-direct-clfd}
 \boL f_!\:\sD^\st(Y\ctrh_\alfd)\lrarrow\sD^\st(X\ctrh_\alfd)
\end{equation}
is constructed in the similar way.
 Finally, for a morphism~$f$ of finite very flat dimension, one
can similarly construct the left derived functor
\begin{equation}  \label{ctrh-direct-fpd}
 \boL f_!\:\sD^\st(Y\ctrh_\fpd)\lrarrow\sD^\st(X\ctrh_\fpd).
\end{equation}

 By construction, for a morphism~$f$ of finite flat dimension,
the triangulated functors~\eqref{ctrh-direct-lct-fpd}
and~\eqref{ctrh-direct-clfd} form a commutative diagram with
the triangulated functors $\sD^\st(Y\ctrh^\lct_\fpd)\rarrow
\sD^\st(Y\ctrh_\alfd)$ and $\sD^\st(X\ctrh^\lct_\fpd)\rarrow
\sD^\st(X\ctrh_\alfd)$ induced by the respective embeddings of
exact categories.
 For a morphism~$f$ of finite very flat dimension,
the triangulated functors~\eqref{ctrh-direct-clfd}
and~\eqref{ctrh-direct-fpd} also form a commutative diagram with
the triangulated functors $\sD^\st(Y\ctrh_\fpd)\rarrow
\sD^\st(Y\ctrh_\alfd)$ and $\sD^\st(X\ctrh_\fpd)\rarrow
\sD^\st(X\ctrh_\alfd)$.

 In particular, for $\bst=\empt$, we obtain a commutative diagram
of triangulated functors and triangulated equivalences from
the upper line of the commutative
diagram~\eqref{becker-contrader-lcta-lct-equivalent-diagram},
\begin{equation} \label{derived-fpd-lct-alfd-direct-images-compatible}
\begin{gathered}
 \xymatrix{
  \sD(Y\ctrh^\lct_\fpd) \ar@<2pt>[r] \ar@<-2pt>@{-}[r]
  \ar[d]^{\boL f_!}
  & \sD(Y\ctrh_\alfd) \ar[d]^{\boL f_!}
  & \sD(Y\ctrh_\fpd) \ar@<-2pt>[l] \ar@<2pt>@{-}[l]
  \ar[d]^{\boL f_!} \\
  \sD(X\ctrh^\lct_\fpd) \ar@<2pt>[r] \ar@<-2pt>@{-}[r]
  & \sD(X\ctrh_\alfd)
  & \sD(X\ctrh_\fpd) \ar@<-2pt>[l] \ar@<2pt>@{-}[l]  
 }
\end{gathered}
\end{equation}
 Here the left-hand part of
the diagram~\eqref{derived-fpd-lct-alfd-direct-images-compatible}
is well-defined and commutative for a morphism~$f$ of finite flat
dimension, while the right-hand part of the diagram is well-defined
and commutative for a morphism~$f$ of finite very flat dimension.

 The triangulated functor~\eqref{ctrh-direct-lct-fpd} for
$\bst=\empt$ or~$\bctr$ agrees with the triangulated
functor~\eqref{bctr-ctrh-lct-direct}, while the triangulated
functor~\eqref{ctrh-direct-fpd} for $\bst=\empt$ or~$\bctr$ agrees
with the triangulated functor~\eqref{bctr-ctrh-direct}.
 So we obtain commutative diagrams of triangulated functors and
triangulated equivalences provided by
Corollary~\ref{becker-contraderived-of-lcta-lct-well-behaved},
\begin{equation} \label{bctr-fpd-lct-lcta-direct-images-compatible}
\begin{gathered}
 \xymatrix{
  \sD(Y\ctrh^\lct_\fpd) \ar@<2pt>[r] \ar@<-2pt>@{-}[r]
  \ar[d]^{\boL f_!}
  & \sD^\bctr(Y\ctrh^\lct) \ar[d]^{\boL f_!} \\
  \sD(X\ctrh^\lct_\fpd) \ar@<2pt>[r] \ar@<-2pt>@{-}[r]
  & \sD^\bctr(X\ctrh^\lct)
 }
 \qquad
  \xymatrix{
  \sD(Y\ctrh_\fpd) \ar@<2pt>[r] \ar@<-2pt>@{-}[r]
  \ar[d]^{\boL f_!}
  & \sD^\bctr(Y\ctrh) \ar[d]^{\boL f_!} \\
  \sD(X\ctrh_\fpd) \ar@<2pt>[r] \ar@<-2pt>@{-}[r]
  & \sD^\bctr(X\ctrh)
 }
\end{gathered}
\end{equation}
 Here the left-hand diagram
in~\eqref{bctr-fpd-lct-lcta-direct-images-compatible} is well-defined
and commutative for a morphism~$f$ of finite flat dimension, while
the right-hand diagram is well-defined and commutative for
a morphism~$f$ of finite very flat dimension.

 The following corollary provides three restricted versions of
Theorem~\ref{direct-images-identified}.
 The case of $\bst=\empt$ can be interpreted as a special case of
Theorem~\ref{becker-co-contra-derived-direct-inverse-adjunction}
(where the morphism~$f$ is \emph{not} assumed to have finite flat
dimension).

\begin{cor}  \label{inj-vfl-direct-images-identified}
\textup{(a)} Assume that the morphism~$f$ has finite flat dimension.
 Then for any symbol\/ $\bst\ne\co$, $\ctr$, $\bco$, $\bctr$
the equivalences of categories\/
$\Hot^\st(Y\qcoh^\inj)\simeq\sD^\st(Y\ctrh^\lin)$ and\/
$\Hot^\st(X\qcoh^\inj)\simeq\sD^\st(X\ctrh^\lin)$ from
Corollary~\textup{\ref{inj-co-contra-cor}(b)} transform the right
derived functor\/ $\boR f_*$~\textup{\eqref{qcoh-direct-fid}} into
the left derived functor\/ $\boL f_!$~\textup{\eqref{ctrh-direct-flid}},
\begin{equation}
\begin{gathered}
 \xymatrix{
  \Hot^\st(Y\qcoh^\inj) \ar@<2pt>[r] \ar@<-2pt>@{-}[r]
  & \sD^\st(Y\qcoh^\fid) \ar@{=}[r] \ar[d]_{\boR f_*}
  & \sD^\st(Y\ctrh^\lid) \ar[d]^{\boL f_!}
  & \sD^\st(Y\ctrh^\lin) \ar@<-2pt>[l] \ar@<2pt>@{-}[l] \\
  \Hot^\st(Y\qcoh^\inj) \ar@<-2pt>[r] \ar@<2pt>@{-}[r]
  & \sD^\st(X\qcoh^\fid) \ar@{=}[r] & \sD^\st(X\ctrh^\lid)
  & \sD^\st(X\ctrh^\lin) \ar@<2pt>[l] \ar@<-2pt>@{-}[l]
 }
\end{gathered}
\end{equation} \par
\textup{(b)} Assume that the morphism~$f$ has finite very flat
dimension.
 Then for any symbol\/ $\bst\ne\co$, $\ctr$, $\bco$, $\bctr$
the equivalences of categories\/
$\sD^\st(Y\qcoh_\vfl)\simeq\Hot^\st(Y\ctrh_\prj)$ and\/
$\sD^\st(X\qcoh_\vfl)\simeq\Hot^\st(X\ctrh_\prj)$ from
Corollary~\textup{\ref{vfl-co-contra-cor-expanded}(a)}
transform the right derived functor\/
$\boR f_*$~\textup{\eqref{qcoh-direct-fvfd}} into the left derived
functor\/ $\boL f_!$~\textup{\eqref{ctrh-direct-fpd}},
\begin{equation}
\begin{gathered}
 \xymatrix{
  \sD^\st(Y\qcoh_\vfl) \ar@<2pt>[r] \ar@<-2pt>@{-}[r]
  & \sD^\st(Y\qcoh_\vfd) \ar@{=}[r] \ar[d]_{\boR f_*}
  & \sD^\st(Y\ctrh_\fpd) \ar[d]^{\boL f_!}
  & \Hot^\st(Y\ctrh_\prj) \ar@<-2pt>[l] \ar@<2pt>@{-}[l] \\
  \sD^\st(Y\qcoh_\vfl) \ar@<-2pt>[r] \ar@<2pt>@{-}[r]
  & \sD^\st(X\qcoh_\vfd) \ar@{=}[r] & \sD^\st(X\ctrh_\fpd)
  & \Hot^\st(X\ctrh_\prj) \ar@<2pt>[l] \ar@<-2pt>@{-}[l]
 }
\end{gathered}
\end{equation} \par
\textup{(c)} Assume that the morphism~$f$ has finite flat dimension.
 Then for any symbol\/ $\bst\ne\co$, $\ctr$, $\bco$, $\bctr$
the equivalences of categories\/
$\sD^\st(Y\qcoh_\fl)\simeq\sD^\st(Y\ctrh_\alf)$ and\/
$\sD^\st(X\qcoh_\fl)\simeq\sD^\st(X\ctrh_\alf)$ from
Corollary~\textup{\ref{vfl-co-contra-cor-expanded}(b)}
transform the right derived functor\/
$\boR f_*$~\textup{\eqref{qcoh-direct-ffd}} into the left derived
functor\/ $\boL f_!$~\textup{\eqref{ctrh-direct-clfd}},
\begin{equation}
\begin{gathered}
 \xymatrix{
  \sD^\st(Y\qcoh_\fl) \ar@<2pt>[r] \ar@<-2pt>@{-}[r]
  & \sD^\st(Y\qcoh_\fd) \ar@{=}[r] \ar[d]_{\boR f_*}
  & \sD^\st(Y\ctrh_\alfd) \ar[d]^{\boL f_!}
  & \sD^\st(Y\ctrh_\alf) \ar@<-2pt>[l] \ar@<2pt>@{-}[l] \\
  \sD^\st(Y\qcoh_\fl) \ar@<-2pt>[r] \ar@<2pt>@{-}[r]
  & \sD^\st(X\qcoh_\fd) \ar@{=}[r] & \sD^\st(X\ctrh_\alfd)
  & \sD^\st(X\ctrh_\alf) \ar@<2pt>[l] \ar@<-2pt>@{-}[l]  
 }
\end{gathered}
\end{equation} \par
\textup{(d)} Assume that the morphism~$f$ has finite flat dimension.
 Then for any symbol\/ $\bst=+$ or\/~$\empt$,
the equivalences of categories\/
$\sD^\st(Y\qcoh_\fl)\simeq\Hot^\st(Y\ctrh^\lct_\prj)$ and\/
$\sD^\st(X\qcoh_\fl)\simeq\Hot^\st(X\ctrh^\lct_\prj)$ from
Corollary~\textup{\ref{vfl-co-contra-cor-expanded}(c)}
transform the right derived functor\/
$\boR f_*$~\textup{\eqref{qcoh-direct-ffd}} into the left derived
functor\/ $\boL f_!$~\textup{\eqref{ctrh-direct-lct-fpd}},
\begin{equation}
\begin{gathered}
 \xymatrix{
  \sD^\st(Y\qcoh_\fl) \ar@<2pt>[r] \ar@<-2pt>@{-}[r]
  & \sD^\st(Y\qcoh_\fd) \ar@{=}[r] \ar[d]_{\boR f_*}
  & \sD^\st(Y\ctrh^\lct_\fpd) \ar[d]^{\boL f_!}
  & \Hot^\st(Y\ctrh^\lct_\prj) \ar@<-2pt>[l] \ar@<2pt>@{-}[l] \\
  \sD^\st(Y\qcoh_\fl) \ar@<-2pt>[r] \ar@<2pt>@{-}[r]
  & \sD^\st(X\qcoh_\fd) \ar@{=}[r] & \sD^\st(X\ctrh^\lct_\fpd)
  & \Hot^\st(X\ctrh^\lct_\prj) \ar@<2pt>[l] \ar@<-2pt>@{-}[l]
 }
\end{gathered}
\end{equation}
\end{cor}

\begin{proof}
 Part~(a): the composition of equivalences of triangulated categories
$\sD^\st(X\allowbreak\qcoh^\fiD)\simeq\Hot(X\qcoh^\inj)\simeq
\sD^\st(X\ctrh^\lin)\simeq\sD^\st(X\ctrh^\fliD)$ can be constructed
directly in terms of the equivalence of exact categories
$X\qcoh^\cta\cap X\qcoh^\fiD\simeq X\ctrh_\al^\fliD$
from Lemma~\ref{cta-clp-finite-flat-inj-dim-identified}(a).
 Now this equivalence of exact categories together with the equivalence
of additive categories $Y\qcoh^\inj\simeq Y\ctrh^\lin_\al$
from Lemma~\ref{cta-clp-restricts-to-cot-inj}(b) transform
the functor~$f_*$ into the functor~$f_!$, which implies the desired
assertion.
 Parts~(b\+d) are similarly proved using
Lemma~\ref{cta-clp-finite-flat-inj-dim-identified}(b\+d).
\end{proof}

\subsection{Derived tensor operations}
\label{derived-tensor-operations-subsect}
 Let $X$ be a quasi-compact semi-separated scheme.
 The functor of tensor product
\begin{equation}  \label{flat-coderived-tensor}
 \ot_{\O_X}\:\sD^\st(X\qcoh_\fl)\times\sD^\st(X\qcoh_\fl)
 \lrarrow\sD^\st(X\qcoh_\fl)
\end{equation}
is well-defined for any symbol $\bst=\co$, $\abs+$, $\abs-$, or~$\abs$,
since the tensor product of a complex coacyclic (respectively,
absolutely acyclic) over the exact category $X\qcoh_\fl$ and any
complex over $X\qcoh_\fl$ is coacyclic (resp., absolutely acyclic)
over $X\qcoh_\fl$.
 The functor of tensor product
\begin{equation}  \label{veryflat-coderived-tensor}
 \ot_{\O_X}\:\sD^\st(X\qcoh_\vfl)\times\sD^\st(X\qcoh_\vfl)
 \lrarrow\sD^\st(X\qcoh_\vfl)
\end{equation}
is well-defined for the similar reasons (though in this case
the situation is actually simpler; see below).

 Furthermore, the functor of tensor product
\begin{equation}   \label{flat-arbitrary-coderived-tensor}
 \ot_{\O_X}\:\sD^\st(X\qcoh_\fl)\times\sD^\st(X\qcoh)
 \lrarrow\sD^\st(X\qcoh)
\end{equation}
is well-defined for $\bst=\co$, $\abs+$, $\abs-$, or~$\abs$, since
the tensor product of a complex coacyclic (resp., absolutely acyclic)
over the exact category $X\qcoh_\fl$ and any complex over $X\qcoh$
is coacyclic (resp., absolutely acyclic) over $X\qcoh$, as is
the tensor product of any complex over $X\qcoh_\fl$ and a complex
coacyclic (resp., absolutely acyclic) over $X\qcoh$.
 In view of Lemma~\ref{co-contra-bounded-fully-faithful}(a),
the functors~(\ref{flat-coderived-tensor}\+-%
\ref{flat-arbitrary-coderived-tensor}) are also well-defined
for $\bst=+$; and it is a standard fact that they are well-defined
for $\bst=\b$ or~$-$.
 The point is that any bounded below acyclic complex is coacyclic,
while any bounded above complex belongs to the minimal triangulated
subcategory in the homotopy category containing the terms of
the complex and closed under countable direct sums (in particular,
any bounded above complex in $X\qcoh_\fl$ is homotopy flat; see below).

 In addition, the following lemma tells us that the tensor product
functors~(\ref{flat-coderived-tensor}\+-%
\ref{flat-arbitrary-coderived-tensor}) are also well-defined for
$\bst=\bco$. {\hbadness=1850\par}

\begin{lem} \label{becker-coderived-tensor-products-well-defined}
\textup{(a)} The tensor product of a Becker-coacyclic complex in
$X\qcoh_\vfl$ with an arbitrary complex in $X\qcoh_\vfl$ is
a Becker-coacyclic complex in $X\qcoh_\vfl$. \par
\textup{(b)} The tensor product of a Becker-coacyclic complex in
$X\qcoh_\fl$ with an arbitrary complex in $X\qcoh_\fl$ is
a Becker-coacyclic complex in $X\qcoh_\fl$. \par
\textup{(c)} The tensor product of an arbitrary complex in
$X\qcoh_\fl$ with a Becker-coacyclic complex in $X\qcoh$ is
a Becker-coacyclic complex in $X\qcoh$. \par
\textup{(d)} The tensor product of a Becker-coacyclic complex in
$X\qcoh_\fl$ with an arbitrary complex in $X\qcoh$ is a Becker-coacyclic
complex in $X\qcoh$.
\end{lem}

\begin{proof}
 This lemma (particularly, part~(d)) is an improved version
of~\cite[Lemma~5.1]{Psemten}.
 In part~(a), we have $\Acycl^\bco(X\qcoh_\vfl)=
\Acycl^\co(X\qcoh_\vfl)=\Acycl^\abs(X\qcoh_\vfl)$ by
Theorem~\ref{vlf-cta-fl-cot-derived-vlf-fl-equivalences}(a),
so the assertion follows from the preceding discussion.
 The following ``global'' argument proving parts~(b\+-d) is based
on the construction of the quasi-coherent internal $\qHom$ functor
$\qHom_{X\qc}$ (see Section~\ref{fHom-subsection}).

 Given two complexes $\M^\bu$ and $\cP^\bu$ over $X\qcoh$, we set
$\qHom_{X\qc}(\M^\bu,\cP^\bu)$ to be the total complex of
the bicomplex $\qHom_{X\qc}(\M^i,\cP^j)$ over $X\qcoh$ constructed
by taking infinite products along the diagonals of the bicomplex.
 Here the infinite products are taken in the Grothendieck abelian
category $X\qcoh$; they can be constructed by applying
the coherator functor to the infinite products taken in
the category of sheaves of $\O_X$\+modules.
 For any three complexes $\M^\bu$, $\N^\bu$, and $\cP^\bu$ over
$X\qcoh$, one has a natural isomorphism of complexes of abelian groups
\begin{equation} \label{tensor-qHom-adjunction-for-complexes}
 \Hom_X(\M^\bu\ot_{\O_X}\N^\bu\;\cP^\bu)\simeq
 \Hom_X(\N^\bu,\qHom_{X\qc}(\M^\bu,\cP^\bu)).
\end{equation}
 Notice that the class of cotorsion quasi-coherent sheaves $X\qcoh^\cot$
is closed under infinite products in $X\qcoh$
by~\cite[Corollary 8.3]{CoFu} or~\cite[Corollary A.2]{CoSt}
(the same applies to the class of contraadjusted quasi-coherent sheaves
$X\qcoh^\cta$ and the class of injective quasi-coherent sheaves
$X\qcoh^\inj$).

 To prove part~(b), let $\F^\bu$ be a Becker-coacyclic complex in
$X\qcoh_\fl$ and $\G^\bu$ be an arbitrary complex in $X\qcoh_\fl$.
 We need to show that $\F^\bu\ot_{\O_X}\G^\bu$ is a Becker-coacyclic
complex in $X\qcoh_\fl$; this means that the complex
$\Hom_X(\F^\bu\ot_{\O_X}\G^\bu\;\cP^\bu)$ is acyclic for any complex
of flat cotorsion quasi-coherent sheaves~$\cP^\bu$.
 By Lemma~\ref{ext-qhom-qc}(b) and the observation above,
$\qHom_{X\qc}(\G^\bu,\cP^\bu)$ is a complex of cotorsion quasi-coherent
sheaves.
 By Theorem~\ref{vlf-cta-fl-cot-derived-vlf-fl-equivalences}(c),
the complex $\F^\bu$ is acyclic in $X\qcoh_\fl$; so
Lemma~\ref{acyclic-of-vfl-flat-arbitrary-of-cta-cot-pairs}(b) with
Lemma~\ref{Ext-1-as-homotopy-Hom} tell us that the complex
$\Hom_X(\F^\bu,\qHom_{X\qc}(\G^\bu,\cP^\bu))$ is acyclic.

 To prove part~(c), let $\F^\bu$ be a complex in $X\qcoh_\fl$ and
$\M^\bu$ be a Becker-coacyclic complex in $X\qcoh$.
 We need to show that $\F^\bu\ot_{\O_X}\M^\bu$ is a Becker-coacyclic
complex in $X\qcoh$; this means that the complex
$\Hom_X(\F^\bu\ot_{\O_X}\M^\bu\;\J^\bu)$ is acyclic for any complex
of injective quasi-coherent sheaves~$\J^\bu$.
 By Lemma~\ref{ext-qhom-qc}(d) and the observation above,
$\qHom_{X\qc}(\F^\bu,\J^\bu)$ is a complex of injective quasi-coherent
sheaves.
 Hence the complex $\Hom_X(\M^\bu,\qHom_{X\qc}(\F^\bu,\J^\bu))$
is acyclic.

 For part~(d), let $\F^\bu$ be a Becker-coacyclic complex in
$X\qcoh_\fl$ and $\M^\bu$ be a complex in $X\qcoh$.
 We need to show that the complex
$\Hom_X(\F^\bu\ot_{\O_X}\M^\bu\;\J^\bu)$ is acyclic for any complex
of injective quasi-coherent sheaves~$\J^\bu$.
 By Lemma~\ref{ext-qhom-qc}(c) and the observation above,
$\qHom_{X\qc}(\M^\bu,\J^\bu)$ is a complex of cotorsion quasi-coherent
sheaves.
 Similarly to the proof of part~(b), it follows that the complex
$\Hom_X(\F^\bu,\qHom_{X\qc}(\M^\bu,\J^\bu))$ is acyclic.

 Alternatively, part~(d) follows from the facts that all
acyclic complexes in $X\qcoh_\fl$ are direct limits of absolutely
acyclic complexes~\cite[Theorems~5.2 and~6.3]{PS6},
\cite[Example~2.5, Remarks~2.6 and~5.1, and Theorem~6.5]{Pflcc}
and the class of all Becker-coacyclic complexes in $X\qcoh$ is
closed under direct limits~\cite[Corollary~7.17]{PS5}.
 Another alternative approach to all parts~(a\+d) is to use
the locality of Becker-coacyclicity
(Corollary~\ref{qcoh-becker-coacyclicity-is-local}) in order to
reduce the question to affine schemes in the spirit of the proof
of Lemma~\ref{becker-co-contra-derived-cohom-well-defined} below.
\end{proof}

 Thus, in view of
Lemma~\ref{becker-coderived-tensor-products-well-defined} and
the preceding discussion, $\sD^\st(X\qcoh_\fl)$ and
$\sD^\st(X\qcoh_\vfl)$ become tensor triangulated categories for
any $\bst\ne\empt$, $\ctr$, $\bctr$; and $\sD^\st(X\qcoh)$ is
a triangulated module category over $\sD^\st(X\qcoh_\fl)$.
 Let us mention once again that, according to
Theorem~\ref{vlf-cta-fl-cot-derived-vlf-fl-equivalences}(a,c), one has
$\Acycl(X\qcoh_\vfl)=\Acycl^\abs(X\qcoh_\vfl)=\Acycl^\bco(X\qcoh_\vfl)$ 
and $\Acycl(X\qcoh_\fl)\allowbreak=\Acycl^\bco(X\qcoh_\fl)$; so
the conventional unbounded derived categories of the exact categories
$X\qcoh_\vfl$ and $X\qcoh_\fl$ coincide with their respective Becker
coderived categories.
 Therefore, $\sD^\st(X\qcoh_\fl)$ and $\sD^\st(X\qcoh_\vfl)$ are also
tensor triangulated categories for $\bst=\empt$ (equivalent to each
other by Theorem~\ref{derived-vfl-fl-equivalence}); but
$\sD^\st(X\qcoh)$ is \emph{not} a triangulated module category over
this tensor triangulated category for $\bst=\empt$. {\hfuzz=3pt\par}

\medskip

 The definition of a homotopy flat complex of flat quasi-coherent
sheaves was given in Section~\ref{homotopy-lin-subsect}.
 The left derived functor of tensor product
\begin{equation}  \label{derived-tensor}
 \ot_{\O_X}^\boL\:\sD^\st(X\qcoh)\times\sD^\st(X\qcoh)\lrarrow
 \sD^\st(X\qcoh)
\end{equation}
is constructed for $\bst=\empt$ by applying the functor $\ot_{\O_X}$
of tensor product of complexes of quasi-coherent sheaves to homotopy
flat complexes of flat quasi-coherent sheaves in one of the arguments
and arbitrary complexes of quasi-coherent sheaves in the other one.
 This derived functor is well-defined by
Theorem~\ref{homotopy-flat-thm}(a\+b).
 In the case of $\bst=-$, the derived functor~\eqref{derived-tensor}
can be constructed by applying the functor $\ot_{\O_X}$ to bounded
above complexes of flat quasi-coherent sheaves in one of the arguments.
 So $\sD^\st(X\qcoh)$ is a tensor triangulated category for any
symbol $\bst=-$ or~$\empt$.

\medskip

 The full triangulated subcategory of \emph{homotopy very flat
complexes} of very flat quasi-coherent sheaves $\sD(X\qcoh_\vfl)_\vfl$
in the unbounded derived category of the exact category of very flat
quasi-coherent sheaves $\sD(X\qcoh_\vfl)$ on $X$ was defined in
Section~\ref{quasi-cta-cot-antilocality-subsect}.
 By Corollary~\ref{quasi-very-cta-cor}(a) and
Proposition~\ref{spal-for-exact}, the composition of natural
triangulated functors $\sD(X\qcoh_\vfl)_\vfl\rarrow\sD(X\qcoh_\vfl)
\rarrow\sD(X\qcoh)$ is an equivalence of triangulated categories. 

 By Theorem~\ref{vlf-cta-fl-cot-derived-vlf-fl-equivalences}(a), any
acyclic complex over $X\qcoh_\vfl$ is absolutely acyclic
over $X\qcoh_\vfl$.
 Hence the composition of functors $\sD(X\qcoh_\vfl)_\vfl\rarrow
\sD(X\qcoh_\vfl)=\sD^\abs(X\qcoh_\vfl)\rarrow\sD^\co(X\qcoh_\fl)$
provides a natural functor
\begin{equation}  \label{quasi-derived-to-flat-coderived}
 \sD(X\qcoh)\rarrow\sD^\co(X\qcoh_\fl),
\end{equation}
which is clearly a tensor triangulated functor between these tensor
triangulated categories.
 In fact, this functor is also fully faithful, and left adjoint to
the composition of Verdier localization functors $\sD^\co(X\qcoh_\fl)
\rarrow\sD(X\qcoh_\fl)\rarrow\sD(X\qcoh)$
(see Corollary~\ref{homotopy-second-category-left-adjoint-cor}).
 Composing the functor~\eqref{quasi-derived-to-flat-coderived} with
the tensor action functor~\eqref{flat-arbitrary-coderived-tensor}
for $\bst=\co$, we obtain the left derived functor
\begin{equation}  \label{quasi-tensor-action}
 \ot_{\O_X}^{\boL'}\: \sD(X\qcoh)\times\sD^\co(X\qcoh)
 \lrarrow\sD^\co(X\qcoh)
\end{equation}
making $\sD^\co(X\qcoh)$ a triangulated module category over
the triangulated tensor category $\sD(X\qcoh)$.

 Similarly, the functor $\sD(X\qcoh_\fl)_\fl\rarrow
\sD(X\qcoh_\fl)=\sD^\bco(X\qcoh_\fl)$ can be interpreted as a functor
\begin{equation}  \label{quasi-derived-to-flat-becker-coderived}
 \sD(X\qcoh)\rarrow\sD^\bco(X\qcoh_\fl),
\end{equation}
which is also clearly a fully faithful tensor triangulated functor
left adjoint to the Verdier localization functor
$\sD^\bco(X\qcoh_\fl)=\sD(X\qcoh_\fl)\rarrow\sD(X\qcoh)$
(see Corollary~\ref{induced-by-inclusion-left-adjoint}).
 Composing the functor~\eqref{quasi-derived-to-flat-becker-coderived}
with the tensor action functor~\eqref{flat-arbitrary-coderived-tensor}
for $\bst=\bco$, we obtain the left derived functor
\begin{equation}  \label{quasi-tensor-action-becker}
 \ot_{\O_X}^{\boL'}\: \sD(X\qcoh)\times\sD^\bco(X\qcoh)
 \lrarrow\sD^\bco(X\qcoh)
\end{equation}
making $\sD^\bco(X\qcoh)$ a triangulated module category over
the triangulated tensor category $\sD(X\qcoh)$ as well
(cf.~\cite[Section~1.4]{Gai}).

\medskip

 Given a symbol $\bst=\b$, $+$, $-$, $\empt$, $\abs+$, $\abs-$, $\co$,
$\ctr$, $\bco$, $\bctr$, or~$\abs$, we set $\bst'$ to be
the ``dual''symbol $\b$, $-$, $+$, $\empt$, $\abs-$, $\abs+$, $\ctr$,
$\co$, $\bctr$, $\bco$, or~$\abs$, respectively.
 Furthermore, let $\bW$ be an open covering of the scheme~$X$.
 Given two complexes $\F^\bu$ over $X\qcoh$ and $\gJ^\bu$ over
$X\lcth_\bW$ such that the $\bW$\+locally contraherent cosheaf
$\Cohom_X(\F^i,\gJ^j)$ is well-defined by the constructions of
Section~\ref{cohom-loc-der-contrahereable} for every pair
$(i,j)\in\boZ^2$, we set $\Cohom_X(\F^\bu,\gJ^\bu)$ to be the total
complex of the bicomplex $\Cohom_X(\F^i,\gJ^j)$ over $X\lcth_\bW$
constructed by taking infinite products of $\bW$\+locally contraherent
cosheaves along the diagonals of the bicomplex (as in
Section~\ref{homotopy-lin-subsect}).

 The functor of cohomomorphisms
\begin{equation}  \label{flat-lct-second-kind-tensor}
 \Cohom_X\:\sD^\st(X\qcoh_\fl)^\op\times\sD^{\st'}(X\lcth_\bW^\lct)
 \lrarrow\sD^{\st'}(X\lcth_\bW^\lct)
\end{equation}
is well-defined for any symbol $\bst=\co$, $\abs+$, $\abs-$, or~$\abs$,
since the $\Cohom$ from a complex coacyclic (respectively, absolutely
acyclic) over $X\qcoh_\fl$ into any complex over $X\lcth_\bW^\lct$ is
a contraacyclic (resp., absolutely acyclic) complex over
$X\lcth_\bW^\lct$, as is the $\Cohom$ from any complex over $X\qcoh_\fl$
into a complex contraacyclic (resp., absolutely acyclic) over
$X\lcth_\bW^\lct$.
 The functor of cohomomorphisms
\begin{equation}  \label{vfl-contra-second-kind-tensor}
  \Cohom_X\:\sD^\st(X\qcoh_\vfl)^\op\times\sD^{\st'}(X\lcth_\bW)
  \lrarrow\sD^{\st'}(X\lcth_\bW)
\end{equation}
is well-defined for the similar reasons.
 In view of Lemma~\ref{co-contra-bounded-fully-faithful},
the functors~(\ref{flat-lct-second-kind-tensor}\+-%
\ref{vfl-contra-second-kind-tensor}) are also well-defined
for $\bst=+$; and one can straightforwardly check that they are
well-defined for $\bst=\b$ or~$-$.
 Once again, the point is that any bounded above acyclic complex is
contraacyclic, while any bounded below complex belongs to the minimal
triangulated subcategory in the homotopy category containing
the terms of the complex and closed under countable products (in
particular, any bounded below complex in $X\lcth_\bW^\lin$ is
homotopy locally injective; see below).

 Furthermore, the following lemma tells us that the $\Cohom$
functors~(\ref{flat-lct-second-kind-tensor}\+-%
\ref{vfl-contra-second-kind-tensor}) are also well-defined
for $\bst=\bco$.

\begin{lem} \label{becker-co-contra-derived-cohom-well-defined}
\hfuzz=14.5pt
\textup{(a)} For any complex\/ $\cA^\bu\in\Acycl^\bco(X\qcoh_\vfl)$
and any complex\/ $\P^\bu\in\Com(X\lcth_\bW)$, the complex\/
$\Cohom_X(\cA^\bu,\P^\bu)$ belongs to\/ $\Acycl^\bctr(X\lcth_\bW)$. \par
\textup{(b)} For any complex\/ $\F^\bu\in\Com(X\qcoh_\vfl)$ and
any complex\/ $\gB^\bu\in\Acycl^\bctr(X\lcth_\bW)$, the complex\/
$\Cohom_X(\F^\bu,\gB^\bu)$ belongs to\/ $\Acycl^\bctr(X\lcth_\bW)$. \par
\textup{(c)} For any complex\/ $\cA^\bu\in\Acycl^\bco(X\qcoh_\fl)$
and any complex\/ $\P^\bu\in\Com(X\lcth_\bW^\lct)$, the complex\/
$\Cohom_X(\cA^\bu,\P^\bu)$ belongs to\/ $\Acycl^\bctr(X\lcth_\bW^\lct)$.
\par
\textup{(d)} For any complex\/ $\F^\bu\in\Com(X\qcoh_\fl)$ and
any complex\/ $\gB^\bu\in\Acycl^\bctr(X\lcth_\bW^\lct)$, the complex\/
$\Cohom_X(\F^\bu,\gB^\bu)$ belongs to\/ $\Acycl^\bctr(X\lcth_\bW^\lct)$.
\end{lem}

\begin{proof}
 The property of a complex in $X\qcoh_\vfl$ or $X\qcoh_\fl$ to be
acyclic is obviously local for open coverings, while the classes
of Becker-coacyclic complexes in these exact categories coincide with
the respective classes of acyclic complexes by
Theorem~\ref{vlf-cta-fl-cot-derived-vlf-fl-equivalences}(a,c).
 The property of a complex in $X\lcth_\bW$ or $X\lcth_\bW^\lct$ to
be Becker-contraacyclic is also local by
Theorem~\ref{Becker-contraacyclicity-local-on-qcomp-qsep}.
 Hence the question reduces to the case of an affine scheme.

 Similarly to the proof of
Lemma~\ref{becker-coderived-tensor-products-well-defined}, part~(a)
follows from the preceding discussion, because any absolutely
acyclic complex in $X\lcth_\bW$ is Becker-contraacyclic by
Lemma~\ref{Positselski-trivial-are-Becker-trivial}(a).
 Now let $R$ be a commutative ring.
 In part~(b), we need to prove that, for any complex
$F^\bu\in\Com(R\modl_\vfl)$ and any complex
$B^\bu\in\Acycl^\bctr(R\modl^\cta)$, the complex $\Hom_R(F^\bu,B^\bu)$
belongs to $\Acycl^\bctr(R\modl^\cta)$.
 Given a complex $G^\bu\in\Com(R\modl_\vfl^\cta)$, we need to check
that the complex $\Hom_R(G^\bu,\Hom_R(F^\bu,B^\bu))$ is acyclic.
 Here it suffices to observe that $\Hom_R(G^\bu,\Hom_R(F^\bu,P^\bu))
\simeq\Hom_R(F^\bu\ot_R G^\bu\;B^\bu)$ and $F^\bu\ot_RG^\bu$ is
a complex in $R\qcoh_\vfl$, and use
Lemma~\ref{all-vfl-flat-contraacyclic-of-cta-cot-pairs}(a) (for
the affine scheme $\Spec R$) with Lemma~\ref{Ext-1-as-homotopy-Hom}.
 The proof of part~(d) is similar and uses
Lemma~\ref{all-vfl-flat-contraacyclic-of-cta-cot-pairs}(b).
 
 In part~(c), we need to show that, for any complex
$A^\bu\in\Acycl^\bco(R\modl_\fl)$ and any complex
$P^\bu\in\Com(R\modl^\cot)$, the complex $\Hom_R(A^\bu,P^\bu)$
belongs to $\Acycl^\bctr(R\modl^\cot)$.
 Given a complex $G^\bu\in\Com(R\modl_\fl^\cot)$, we need to check
that the complex $\Hom_R(G^\bu,\Hom_R(A^\bu,P^\bu))$ is acyclic.
 Here we observe that $\Hom_R(G^\bu,\Hom_R(A^\bu,P^\bu))
\simeq\Hom_R(A^\bu,\Hom_R(G^\bu,P^\bu))$ and $\Hom_R(G^\bu,P^\bu)$
is a complex in $R\modl^\cot$, and use
Lemma~\ref{acyclic-of-vfl-flat-arbitrary-of-cta-cot-pairs}(b) (for
$\Spec R$) with Lemma~\ref{Ext-1-as-homotopy-Hom}.
\end{proof}

 Thus, in view of
Lemma~\ref{becker-co-contra-derived-cohom-well-defined} and
the preceding discussion, the category opposite to
$\sD^{\st'}(X\lcth_\bW^\lct)$ becomes a triangulated module category
over the tensor triangulated category $\sD^\st(X\qcoh_\fl)$, and
the category opposite to $\sD^{\st'}(X\lcth_\bW)$ is a triangulated
module category over the tensor triangulated category
$\sD^\st(X\qcoh_\vfl)$, for any symbol $\bst\ne\empt$, $\ctr$,~$\bctr$.

\medskip

 The definition of a homotopy locally injective complex of
locally injective $\bW$\+locally contraherent cosheaves was given in
Section~\ref{homotopy-lin-subsect}.
 The following corollary is a dual-analogous version of
(and a complement to)
Lemma~\ref{homotopy-local-injectivity-characterizations}\,%
(1)\,$\Rightarrow$\,(2).

\begin{cor} \label{conventional-derived-Cohom-well-defined}
\textup{(a)} For any homotopy very flat complex of very flat
quasi-coherent sheaves\/ $\F^\bu$ and any acyclic complex of\/
$\bW$\+locally contraherent cosheaves\/ $\P^\bu$ on $X$, the complex
of\/ $\bW$\+locally contraherent cosheaves\/ $\Cohom_X(\F^\bu,\P^\bu)$
on $X$ is acyclic. \par
\textup{(b)} For any homotopy flat complex of flat quasi-coherent
sheaves\/ $\F^\bu$ and any acyclic complex of locally cotorsion\/
$\bW$\+locally contraherent cosheaves\/ $\P^\bu$ on $X$, the complex
of locally cotorsion\/ $\bW$\+locally contraherent cosheaves\/
$\Cohom_X(\F^\bu,\P^\bu)$ on $X$ is acyclic.
\end{cor}

\begin{proof}
 Let us prove part~(b).
 Start from the observation that $\Cohom_X$ from any complex of
quasi-coherent sheaves acyclic with respect to $X\qcoh_\fl$ into
any complex of locally cotorsion $\bW$\+locally contraherent cosheaves
is acyclic in $X\lcth_\bW^\lct$.
 Indeed, any acyclic complex is Becker-coacyclic in $X\qcoh_\fl$ by
Theorem~\ref{vlf-cta-fl-cot-derived-vlf-fl-equivalences}(c) and
any Becker-contraacyclic complex is acyclic in $X\lcth_\bW^\lct$ by
Corollary~\ref{Becker-contraacyclic-are-acyclic}(b), so it remains to
refer to Lemma~\ref{becker-co-contra-derived-cohom-well-defined}(c).

 Therefore, the class of all complexes of flat quasi-coherent sheaves
$\F^\bu$ satisfying the conclusion of part~(b) can be viewed as
a strictly full triangulated subcategory in $\sD(X\qcoh_\fl)$.
 The class of all homotopy flat complexes of flat quasi-coherent
sheaves is a strictly full subcategory in $\sD(X\qcoh_\fl)$ by
the definition (see Section~\ref{homotopy-lin-subsect}).
 It remains to observe that the functor $\Cohom_X$ takes infinite
direct sums of complexes in the first argument to infinite products,
and is exact as a functor $\Cohom_X(\F,{-})\:X\lcth_\bW^\lct\rarrow
X\lcth_\bW^\lct$ for any fixed flat quasi-coherent sheaf $\F$ in
the first argument.
\end{proof}

 The right derived functor of cohomomorphisms
\begin{equation}  \label{derived-cohom}
 \boR\Cohom_X\:\sD^\st(X\qcoh)^\op\times\sD^{\st'}(X\lcth_\bW)
 \lrarrow\sD^{\st'}(X\lcth_\bW)
\end{equation}
is constructed for $\bst=\empt$ by applying the functor $\Cohom_X$
to homotopy very flat complexes of very flat quasi-coherent sheaves
in the first argument and arbitrary complexes of $\bW$\+locally
contraherent cosheaves in the second argument, or alternatively,
to arbitrary complexes of quasi-coherent sheaves in the first
argument and homotopy locally injective complexes of locally injective
$\bW$\+locally contraherent cosheaves in the second argument.
 This derived functor is well-defined by Theorem~\ref{homotopy-lin-thm},
Lemma~\ref{homotopy-local-injectivity-characterizations}\,%
(1)\,$\Rightarrow$\,(2), Proposition~\ref{spal-for-exact}, and
Corollary~\ref{conventional-derived-Cohom-well-defined}(a)
(cf.~\cite[Lemma~2.7]{Psemi}).
 In the case of $\st=-$, the derived functor~\eqref{derived-cohom}
can be constructed by applying the functor $\Cohom_X$ to bounded above
complexes of very flat quasi-coherent sheaves in the first argument or
bounded below complexes of locally injective $\bW$\+locally contraherent
cosheaves in the second argument.
 So $\sD^{\st'}(X\lcth_\bW)^\op$ is a triangulated module category over
the tensor triangulated category $\sD^\st(X\qcoh)$ for any symbol
$\bst=-$ or~$\empt$.

 If one prefers the locally cotorsion contraherent cosheaves, then
the right derived functor of cohomomorphisms
\begin{equation}  \label{derived-lct-cohom}
 \boR\Cohom_X\:\sD^\st(X\qcoh)^\op\times\sD^{\st'}(X\lcth_\bW^\lct)
 \lrarrow\sD^{\st'}(X\lcth_\bW^\lct)
\end{equation}
can be constructed for $\bst=\empt$ by applying the functor $\Cohom_X$
to homotopy flat complexes of flat quasi-coherent sheaves in the first
argument and arbitrary complexes of locally cotorsion $\bW$\+locally
contraherent cosheaves in the second argument, or alternatively,
to arbitrary complexes of quasi-coherent sheaves in the first
argument and homotopy locally injective complexes of locally injective
$\bW$\+locally contraherent cosheaves in the second argument.
 This derived functor is well-defined by Theorem~\ref{homotopy-lin-thm},
Lemma~\ref{homotopy-local-injectivity-characterizations}\,%
(1)\,$\Rightarrow$\,(2), Theorem~\ref{homotopy-flat-thm}(a), and
Corollary~\ref{conventional-derived-Cohom-well-defined}(b).
 In the case of $\st=-$, the derived functor~\eqref{derived-lct-cohom}
can be constructed by applying the functor $\Cohom_X$ to bounded above
complexes of flat quasi-coherent sheaves in the first argument or
bounded below complexes of locally injective $\bW$\+locally contraherent
cosheaves in the second argument.
 Thus $\sD^{\st'}(X\lcth_\bW^\lct)^\op$ is also a triangulated module
category over the tensor triangulated category $\sD^\st(X\qcoh)$
for any symbol $\bst=-$ or~$\empt$.

\medskip

 Composing the functor~\eqref{quasi-derived-to-flat-coderived}
with the tensor action functor~\eqref{flat-lct-second-kind-tensor}
for $\bst=\co$, one obtains the right derived functor
\begin{equation}  \label{lct-tensor-action}
 \boR'\Cohom_X\: \sD(X\qcoh)^\op\times\sD^\ctr(X\lcth_\bW^\lct)
 \lrarrow\sD^\ctr(X\lcth_\bW^\lct)
\end{equation}
making $\sD^\ctr(X\lcth_\bW^\lct)^\op$ a triangulated module category
over the triangulated tensor category $\sD(X\qcoh)$.
 Similarly, the composition $\sD(X\qcoh_\vfl)_\vfl\rarrow
\sD(X\qcoh_\vfl)=\sD^\co(X\qcoh_\vfl)$ provides a natural
functor
\begin{equation}  \label{quasi-derived-to-very-flat-coderived}
 \sD(X\qcoh)\rarrow\sD^\co(X\qcoh_\vfl),
\end{equation}
which is a tensor triangulated functor between these tensor
triangulated categories.
 It is also fully faithful, and left adjoint to the composition of
Verdier localization functors $\sD^\co(X\qcoh_\vfl)\rarrow
\sD(X\qcoh_\vfl)\rarrow\sD(X\qcoh)$.
 Composing the functor~\eqref{quasi-derived-to-very-flat-coderived}
with the tensor action functor~\eqref{vfl-contra-second-kind-tensor}
for $\bst=\co$, one obtains the right derived functor
\begin{equation}  \label{contra-tensor-action}
 \boR'\Cohom_X\: \sD(X\qcoh)^\op\times\sD^\ctr(X\lcth_\bW)
 \lrarrow\sD^\ctr(X\lcth_\bW)
\end{equation}
making $\sD^\ctr(X\lcth_\bW)^\op$ a triangulated module category
over the triangulated tensor category $\sD(X\qcoh)$.

 Similarly, composing
the functor~\eqref{quasi-derived-to-flat-becker-coderived}
with the tensor action functor~\eqref{flat-lct-second-kind-tensor}
for $\bst=\bco$, one obtains the right derived functor
\begin{equation}  \label{becker-lct-tensor-action}
 \boR'\Cohom_X\: \sD(X\qcoh)^\op\times\sD^\bctr(X\lcth_\bW^\lct)
 \lrarrow\sD^\bctr(X\lcth_\bW^\lct)
\end{equation}
making $\sD^\bctr(X\lcth_\bW^\lct)^\op$ a triangulated module category
over the triangulated tensor category $\sD(X\qcoh)$.
 Composing the functor $\sD(X\qcoh)\rarrow\sD^\co(X\qcoh_\vfl)=
\sD^\bco(X\qcoh_\vfl)$ \eqref{quasi-derived-to-very-flat-coderived}
with the tensor action functor~\eqref{vfl-contra-second-kind-tensor}
for $\bst=\bco$, one obtains the right derived functor
\begin{equation}  \label{becker-contra-tensor-action}
 \boR'\Cohom_X\: \sD(X\qcoh)^\op\times\sD^\bctr(X\lcth_\bW)
 \lrarrow\sD^\bctr(X\lcth_\bW)
\end{equation}
making $\sD^\bctr(X\lcth_\bW)^\op$ a triangulated module category
over the triangulated tensor category $\sD(X\qcoh)$.

\Section{Noetherian Schemes}
\label{dualizing-complex-sect}

\subsection{Projective locally cotorsion contraherent cosheaves}
\label{lct-projective-loc-noetherian}
 Let $X$ be a (not necessarily semi-separated) locally Noetherian scheme
and $\bW$ be its open covering.
 The following theorem, providing a description of projective locally
cotorsion contraherent cosheaves on $X$, is to be compared to
Hartshorne's classification of injective quasi-coherent sheaves on~$X$
\,\cite[Proposition~II.7.17]{Har} based on Matlis' classification of
injective modules over a Noetherian commutative ring~\cite{Mat}.
 Our theorem is based on Enochs' classification of flat cotorsion
modules over a Noetherian commutative ring~\cite{En} (stated above
as Theorem~\ref{flat-cotorsion-classification}).

 Theorem~\ref{proj-lct-classification} can be also compared to
the results of~\cite[Section~13]{Pcta}, which describe arbitrary
locally cotorsion contraherent cosheaves on (affine) Noetherian schemes
of Krull dimension~$1$ (see specifically~\cite[Corollaries~13.12
and~13.13]{Pcta}, as well as the discussions in the introduction
to~\cite{Pcta} and~\cite[Remark~12.13]{Pcta}).

\begin{thm}  \label{proj-lct-classification}
\textup{(a)} There are enough projective objects in the exact categories
of locally cotorsion locally contraherent cosheaves $X\lcth_\bW^\lct$
and $X\lcth^\lct$, and all these projective objects belong to the full
subcategory of locally cotorsion contraherent cosheaves $X\ctrh^\lct$.
 The full subcategories of projective objects in the three exact
categories $X\ctrh^\lct\sub X\lcth_\bW^\lct\sub X\lcth^\lct$ coincide.
\par
\textup{(b)} For any scheme point $x\in X$, denote by\/
$\widehat\O_{x,X}$ the completion of the local ring\/ $\O_{x,X}$ and
by $\iota_x\:\Spec\O_{x,X}\rarrow X$ the natural morphism.
 Then a locally cotorsion contraherent cosheaf\/ $\gF$ on $X$ is
projective if and only if it is isomorphic to an infinite product\/
$\prod_{x\in X}\iota_x{}_!\widecheck F_x$ of the direct images
$\iota_x{}_!\widecheck F_x$ of the contraherent cosheaves
$\widecheck F_x$ on\/ $\Spec\O_{x,X}$ corresponding to some free
contramodules\/ $F_x$ over the complete Noetherian local rings\/
$\widehat\O_{x,X}$ (viewed as\/ $\O_{x,X}$\+modules via
the restriction of scalars).
\end{thm}

\begin{proof}
 For a quasi-compact semi-separated scheme~$X$, the assertions of
part~(a) were proved in Section~\ref{projective-contraherent}. 
 In the general case, we will prove parts (a) and~(b) simultaneously.
 The argument is based on Theorem~\ref{flat-cotorsion-classification}.

 First of all, let us show that the cosheaf of $\O_X$\+modules
$\iota_x{}_!\widecheck P_x$ is a locally cotorsion contraherent cosheaf
on $X$ for any contramodule $P_x$ over $\widehat\O_{x,X}$.
 It suffices to check that the restriction of $\iota_x{}_!
\widecheck P_x$ to every affine open subscheme $U\sub X$ is
a (locally) cotorsion contraherent cosheaf.
 If the point~$x$ belongs to $U$, then the morphism $\iota_x$
factorizes through $U$; denoting the morphism $\Spec\O_{x,X}
\rarrow U$ by $\kappa_x$, one has $(\iota_x{}_!\widecheck P_x)|_U
\simeq \kappa_x{}_!\widecheck P_x$.
 The morphism~$\kappa_x$ being affine and flat, and the contraherent
cosheaf $\widecheck P_x$ over $\Spec\O_{x,X}$ being (locally)
cotorsion by Proposition~\ref{contramodules-cotorsion}(a),
$\kappa_x{}_!\widecheck P_x$ is a (locally) cotorsion contraherent
cosheaf over~$U$.

 If the point~$x$ does not belong to $U$, we will show that
the $\O_X(U)$\+module $(\iota_x{}_!\widecheck P_x)[U]\simeq
\widecheck P_x[\iota_x^{-1}(U)]$ vanishes.
 So it will follow that the restriction of $\iota_x{}_!\widecheck P_x$
to the open complement of the closure of the point $x$ in $X$ is 
a zero cosheaf.
 Indeed, it suffices to check that the module of cosections
$\widecheck P_x[V]$ vanishes for any principal affine open subscheme
$V\sub \Spec\O_{x,X}$ that does not contain the closed point.
 In other words, it has to be shown that $\Hom_{\O_{x,X}}
(\O_{x,X}[s^{-1}],P_x)=0$ for any element~$s$ from
the maximal ideal of the local ring~$\O_{x,X}$.
 This holds for any $\widehat\O_{x,X}$\+contramodule~$P_x$;
see~\cite[Theorem~B.1.1(2c)]{Pweak} (cf.\ the definition of
an $I$\+contramodule $R$\+module in~\cite[Section~2]{Pmgm}
and~\cite[Section~9]{Pcta}, the latter one based on the discussion
in~\cite[Section~2 and Theorem~5.1]{Pcta}).

 It follows that cosheaves of the form
$\prod_x\iota_x{}_!\widecheck P_x$ on $X$, where $P_x$ are some
$\widehat\O_{x,X}$\+contramodules, belong to $X\ctrh^\lct$.
 Let us show that cosheaves $\gF=\prod_x\iota_x{}_!\widecheck F_x$,
where $F_x$ denote some free $\widehat\O_{x,X}$\+contramodules,
are projective objects in the exact category $X\lcth_\bW^\lct$.
 Pick an affine open covering $U_\alpha$ of the scheme $X$ 
subordinate to the covering $\bW$, and choose a well-ordering
of the set of indices~$\{\alpha\}$.

 Given an index~$\alpha$, denote by $\gF_\alpha$ the product of
the cosheaves $\iota_z{}_!\widecheck F_z$ on $X$ taken over all
points $z\in S_\alpha=U_\alpha\setminus\bigcup_{\gamma<\alpha}U_\gamma$.
 Clearly, one has $\gF\simeq\prod_\alpha\gF_\alpha$.
 Furthermore, since $X$ is locally Noetherian, for any affine
(or even quasi-compact) open subscheme $U\sub X$ one has
$U\cap S_\alpha=\empt$ for all but a finite number of indices~$\alpha$.
 Let us emphasize that this claim is based on Noetherianity rather than
just quasi-compactness: the point is that any ascending chain of
open subschemes in a Noetherian scheme terminates.
 Consequently, $\gF_\alpha[U]=0$ for all but a finite number of
indices~$\alpha$.
 We conclude that the cosheaf $\gF$ is also the direct sum
of the cosheaves $\gF_\alpha$ (taken in the category of cosheaves
of $\O_X$\+modules).
 Hence it suffices to check that each $\gF_\alpha$ is a projective
object in $X\lcth_\bW^\lct$.

 Denoting by~$j_\alpha$ the open embedding $U_\alpha\rarrow X$
and by~$\kappa_z$ the natural morphisms $\Spec\O_{z,X}\rarrow U_\alpha$,
we notice that $\gF_\alpha=j_\alpha{}_!\gG_\alpha$, where $\gG_\alpha =
\prod_{z\in S_\alpha}\kappa_z{}_!\widecheck F_z$.
 The observation that the direct images with respect to quasi-compact
quasi-separated morphisms of schemes preserve infinite products of
cosheaves of $\O$\+modules (see Section~\ref{locally-contraherent})
is used here.
 By Theorem~\ref{flat-cotorsion-classification}, \,$\gG_\alpha$
is a projective object in $U_\alpha\ctrh^\lct$; and 
by the adjunction~\eqref{direct-inverse-cosheaf-lct-adjunction} from
Section~\ref{direct-inverse-loc-contra} it follows that $\gF_\alpha$
is a projective object in $X\lcth_\bW^\lct$.

 Now let us construct for any locally cotorsion $\bW$\+locally
contraherent cosheaf $\Q$ on $X$ an admissible epimorphism onto
$\Q$ from an object $\gF=\prod_{x\in X}\iota_x{}_!\widecheck F_x$
in the exact category $X\lcth_\bW^\lct$.
 Choose an affine open covering $U_\alpha$ as above, and proceed
by transfinite induction in~$\alpha$, constructing contraherent
cosheaves $\gF_\alpha=\prod_{z\in S_\alpha}\iota_z{}_!\widecheck F_z$
and morphisms of locally contraherent cosheaves $\gF_\alpha\rarrow\Q$.

 Suppose that such cosheaves and morphisms have been constructed for
all $\alpha<\beta$; then, as it was explained above, there is
the induced morphism of locally contraherent cosheaves
$\prod_{\alpha<\beta}\gF_\alpha\rarrow\Q$.
 Assume that the related morphism of the modules of cosections
$\prod_{\alpha<\beta}\gF_\alpha[U]=\bigoplus_{\alpha<\beta}
\gF_\alpha[U]\rarrow\Q[U]$ is an admissible epimorphism of
cotorsion $\O_X(U)$\+modules for any affine open subscheme
$U\sub\bigcup_{\alpha<\beta} U_\alpha\sub X$ subordinate to~$\bW$.
 We are going to construct a contraherent cosheaf $\gF_\beta=
\prod_{z\in S_\beta}\iota_z{}_!\widecheck F_z$ and a morphism
of locally contraherent cosheaves $\gF_\beta\rarrow\Q$ such that
the induced morphism $\prod_{\alpha\le\beta}\gF_\alpha\rarrow\Q$
satisfies the above condition for any affine open subscheme
$U\sub\bigcup_{\alpha\le\beta}U_\alpha\sub X$ subordinate to~$\bW$.

 Pick an admissible epimorphism in $U_\beta\ctrh^\lct$ onto
the (locally) cotorsion contraherent cosheaf $j_\beta^!\Q$ from
a projective (locally) cotorsion contraherent cosheaf $\gG$ on~$U_\beta$.
 By Theorem~\ref{flat-cotorsion-classification}, the cosheaf
$\gG$ decomposes into a direct product $\gG\simeq\prod_{x\in U_\beta}
\lambda_x{}_!\widecheck G_x$, where $G_x$ are some free
$\widehat\O_{x,X}$\+contramodules and $\lambda_x\:\Spec\O_{x,X}\rarrow
U_\beta$ are the natural morphisms.
 Let us rewrite this product as $\gG=
\prod_{z\in S_\beta}\lambda_z{}_!\widecheck G_z\oplus
\prod_{y\in U_\beta\setminus S_\beta}\lambda_y{}_!\widecheck G_y$;
denote the former direct summand by $\gG_\beta$ and the latter one
by~$\gE$.
 Set $\gF=\prod_{\alpha<\beta}\gF_\alpha$ and $\gF_\beta =
j_\beta{}_!\gG_\beta$.

 The property of a morphism of locally cotorsion $\bT$\+locally
contraherent cosheaves to be an admissible epimorphism being local,
we only need to show that the natural morphism $\gF\oplus\gF_\beta
\rarrow\Q$ becomes an admissible epimorphism in $U_\beta\ctrh^\lct$
when restricted to~$U_\beta$.
 In other words, this means that the morphism $j_\beta^!\gF\oplus
\gG_\beta\rarrow j_\beta^!\Q$ should be an admissible epimorphism
in $U_\beta\ctrh^\lct$.
 It suffices to check that the admissible epimorphism $\gG\rarrow
j_\beta^!\Q$ factorizes through the morphism in question, or that
the morphism $\gE\rarrow j_\beta^!\Q$ factorizes through
the morphism $j_\beta^!\gF\rarrow j_\beta^!\Q$.

 Denote by $j$ the open embedding $U_\beta\setminus S_\beta\rarrow
U_\beta$; then one has $\gE\simeq j_!\gL$, where $\gL$ is
a projective object in $(U_\beta\setminus S_\beta)\ctrh^\lct$
(as we have proved above).
 Then the morphism $j_!\gL\rarrow j_\beta^!\Q$ factorizes through
the morphism $j_\beta^!\gF\rarrow j_\beta^!\Q$, since the morphism
$\gL\rarrow j^!j_\beta^!\Q$ factorizes through an admissible
epimorphism of locally cotorsion contraherent cosheaves
$j^!j_\beta^!\gF\rarrow j^!j_\beta^!\Q$.
 Here we are using the adjunction
isomorphism~\eqref{open-cosheaf-direct-inverse-adjunction}
from Section~\ref{contra-direct-inverse}, as well as the assumption
that the morphism $\gF\rarrow\Q$ is an admissible epimorphism
over $\bigcup_{\alpha<\beta}U_\alpha\supset U_\beta\setminus S_\beta$.

 We have proved part~(a); and to finish the proof of part~(b) it
remains to show that the class of contraherent cosheaves of
the form $\prod_x\iota_x{}_!\widecheck F_x$ on $X$ is closed 
under direct summands.
 The argument is based on the following two lemmas.

\begin{lem} \label{proj-lct-orthogonality}
 \textup{(a)} Suppose that the set of all scheme points of $X$ is
presented as a union of two nonintersecting subsets $X=S\sqcup T$
such that for any points $z\in S$ and $y\in T$ the closure of~$z$
in $X$ does not contain~$y$.
 Then for any cosheaves of\/ $\O$\+modules\/ $\P_y$ over\/
$\Spec \O_{y,X}$ and any contramodules $P_z$ over $\widehat\O_{z,X}$
one has\/ $\Hom^{\O_X}(\prod_{y\in T}\iota_y{}_!\P_y\;
\prod_{z\in S}\iota_z{}_!\widecheck P_z)=0$. \par
 \textup{(b)} For any scheme point~$x\in X$, the functor assigning
to an $\widehat\O_{x,X}$\+contramodule $P_x$ the locally cotorsion
contraherent cosheaf $\iota_x{}_!\widecheck P_x$ on $X$ is
fully faithful.
\end{lem}

\begin{proof}
 Part~(a): by the definition of the infinite product, 
it suffices to show that $\Hom^{\O_X}(\prod_{y\in T}\iota_y{}_!\P_y\;
\iota_z{}_!\widecheck P_z)=0$ for any $z\in S$.
 Let $Z$ be the closure of $z$ in $X$ and $Y=X\setminus Z$ be its
complement; then one has $T\sub Y$ and $Y$ is an open subscheme in~$X$.
 Let $j$~denote the open embedding $Y\rarrow X$.
 Given $y\in T$, denote the natural morphism $\Spec\O_{y,X}
\rarrow Y$ by $\kappa_y$, so $\iota_y=j\circ\kappa_y$.

 Now we have $\prod_{y\in T}\iota_y{}_!\P_y\simeq
j_!\prod_{y\in T}\kappa_y{}_!\P_y$ and, according to
the adjunction~\eqref{direct-inverse-cosheaf-lct-adjunction},
$$ \textstyle
\Hom^{\O_X}(j_!\prod_{y\in T}\kappa_y{}_!\P_y\;
\iota_z{}_!\widecheck P_z)\simeq\Hom^{\O_Y}(\prod_{y\in T}
\kappa_y{}_!\P_y\;j^!\iota_z{}_!\widecheck P_z).
$$
 It was shown above that
$j^!\iota_z{}_!\widecheck P_z=0$, so we are done.

 Part~(b): it was explained in the paragraph before
Proposition~\ref{contramodules-cotorsion}
in Section~\ref{cotorsion-modules} that the functor assigning to
an $\widehat\O_{x,X}$\+contramodule $P_x$ the locally cotorsion
contraherent cosheaf $\widecheck P_x$ on $\Spec\O_{x,X}$ is fully
faithful.
 The morphism~$\iota_x$ being flat and coaffine,
the adjunction~\eqref{direct-inverse-cosheaf-lct-adjunction} applies
and it suffices to show that the adjunction morphism $\widecheck P_x
\rarrow\iota_x^!\iota_x{}_!\widecheck P_x$ is an isomorphism in
$\Spec\O_{x,X}\ctrh^\lct$.
 One can replace the scheme $X$ by any affine open
subscheme $U\sub X$ containing~$x$; and it remains to use 
the isomorphism $\O_{x,X}\otimes_{\O(U)}\O_{x,X}\simeq\O_{x,X}$,
which implies the isomorphism $\Hom_{\O(U)}(\O_{x,X},P_x)\simeq P_x$
for any $\O_{x,X}$\+module~$P_x$.
\end{proof}

 The next lemma will be applied first to a certain category of
cosheaves on a Noetherian scheme $X$, and then to the \emph{opposite}
category to a certain category of cosheaves on a locally
Noetherian scheme~$X$.

\begin{lem} \label{triangular-idempotent-image-lemma}
 Let\/ $\sA$ be an additive category and $F_\alpha$ be a family of
objects of\/ $\sA$ indexed by a well-ordered set of
indices~$\{\alpha\}$.
 Assume that, for every index~$\beta$, the products
$\prod_{\alpha<\beta}F_\alpha$ and $\prod_{\alpha\ge\beta}F_\alpha$
exist in\/ $\sA$, and\/ $\Hom_\sA(\prod_{\alpha\ge\beta}F_\alpha\;
\prod_{\alpha<\beta}F_\alpha)=0$.
 Assume further that the images of all idempotent endomorphisms of
the objects $F_\alpha$ exist in\/~$\sA$.
 Then the images of all idempotent endomorphisms of the object\/
$\prod_\alpha F_\alpha$ exist in\/ $\sA$ as well; moreover,
they are isomorphic to products over~$\alpha$ of images of
idempotent endomorphisms of~$F_\alpha$.
\end{lem}

\begin{proof}
 It follows from the assumptions of the lemma that any endomorphism
$f\:\prod_\alpha F_\alpha\rarrow\prod_\alpha F_\alpha$ induces
endomorphisms $f_{<\beta}\:\prod_{\alpha<\beta}F_\alpha\rarrow
\prod_{\alpha<\beta}F_\alpha$ such that the morphisms $f_{<\beta}$
and $f_{<\gamma}$ form commutative diagrams with the direct summand
projections $\prod_{\alpha<\gamma}F_\alpha\rarrow\prod_{\alpha<\beta}
F_\alpha$ for all pairs of indices $\beta<\gamma$.
 There also naturally induced endomorphisms $f_\alpha\:F_\alpha
\rarrow F_\alpha$ forming obvious commutative diagrams.

 Now let $e\:\prod_\alpha F_\alpha\rarrow\prod_\alpha F_\alpha$ be
an idempotent endomorphism.
 Proceeding by transfinite induction on~$\beta$, one proves that
the morphism~$e_{<\beta}$ has an image in $\sA$ isomorphic to
the product of the images of~$e_\alpha$ over the indices $\alpha<\beta$.
 Moreover, the isomorphisms in question can be chosen in such a way
that they form commutative diagrams with the subproduct projections
for pairs of indices $\beta<\gamma$ as above.
 Since the passages to the images of idempotent endomorphisms preserve
projective limits, the question reduces to the case of a successor
ordinal, where the following observation is sufficient.
 Let $0\rarrow A\rarrow B\rarrow C\rarrow0$ be a split short exact
sequence in $\sA$, and let $g\:(0\to A\to B\to C\to0)\rarrow(0\to A\to
B\to C\to0$) be its idempotent endomorphism.
 Assume that the components $g_A\:A\rarrow A$ and $g_C\:C\rarrow C$
of~$g$ have images $A'$ and $C'$ in~$\sA$.
 Then the component $g_B\:B\rarrow B$ also has an image $B'$ in $\sA$,
and the induced short sequence $0\rarrow A'\rarrow B'\rarrow C'
\rarrow0$ is split exact.
\end{proof}

 Now we can finish the proof of Theorem.
 Firstly, consider the case of a Noetherian (i.~e., quasi-compact)
scheme~$X$.
 Then there is \emph{no} infinite sequence of points $x_1$, $x_2$,
$x_3$,~\dots~$\in X$ such that $x_{n+1}\ne x_n$ and $x_{n+1}$ belongs
to the closure of~$x_n$ for every $n\ge1$.
 Consequently, the specialization order on the set of all points of $X$
can be refined to a well-ordering.
 So there exists a well-ordering relation~$<$ on $X$ such that
$x<y$ whenever $x$~belongs to the closure of~$y$.
 Now we apply Lemma~\ref{triangular-idempotent-image-lemma}
to the full additive subcategory $\sA$ of objects of the form
$\prod_{x\in X}\iota_x{}_!\widecheck F_x$ in $X\ctrh$, and claim that
the images of all idempotent endomorphisms exist in~$\sA$.
 Indeed, the first assumption of 
Lemma~\ref{triangular-idempotent-image-lemma} holds by
Lemma~\ref{proj-lct-orthogonality}(a).
 To check the second assumption, use
Lemma~\ref{proj-lct-orthogonality}(b) together with the fact that
the image of any idempotent endomorphism of a free
$\widehat\O_{x,X}$\+contramodule $F_x$ is also a free
$\widehat\O_{x,X}$\+contramodule (\cite[Lemma~1.3.2]{Pweak}
or~\cite[Corollary~10.7]{Pcta}).

 For an arbitrary locally Noetherian scheme $X$, we consider
the \emph{direct sum} decomposition $\gF=\prod_\alpha\gF_\alpha=
\bigoplus_\alpha\gF_\alpha$ of the object
$\gF=\prod_{x\in X}\iota_x{}_!\widecheck F_x$ in the category
of cosheaves of $\O_X$\+modules, as per the explanation above
in this proof.
 Here $X=\bigcup_\alpha U_\alpha$ is an affine open covering of $X$
indexed by a well-ordered set of indices~$\alpha$.
 We apply Lemma~\ref{triangular-idempotent-image-lemma} to
the \emph{opposite category} $\sA$ to the full additive subcategory of
objects of the form $\prod_{x\in X}\iota_x{}_!\widecheck F_x$ in
the category of cosheaves of $\O_X$\+modules, and claim that the images
of all idempotent endomorphisms exist in~$\sA$.
 Once again, the first assumption of 
Lemma~\ref{triangular-idempotent-image-lemma} holds by
Lemma~\ref{proj-lct-orthogonality}(a).
 To check the second assumption, recall that
$\gF_\alpha=j_\alpha{}_!\gG_\alpha$, where
$\gG_\alpha=\prod_{z\in S_\alpha}\kappa_z{}_!\widecheck F_z$.
 The functor $j_\alpha{}_!$ is fully faithful in view of
the adjunction~\eqref{open-cosheaf-direct-inverse-adjunction}
or~\eqref{direct-inverse-cosheaf-lct-adjunction}; so the images of
all idempotent endomorphisms of the cosheaves $\gF_\alpha$ have
the desired form by the argument in the previous paragraph applied to
the Noetherian affine scheme $U_\alpha$ in the role of~$X$.
\end{proof}

 As in Section~\ref{projective-contraherent}, we denote the category
of projective locally cotorsion contraherent cosheaves on $X$ by
$X\ctrh^\lct_\prj$.
 The following corollary says that being a projective locally cotorsion
contraherent cosheaf on a locally Noetherian scheme is a local
property.
 Besides, all such cosheaves are coflasque (see Section~\ref{coflasque}).
 The similar properties of injective quasi-coherent sheaves
are usually deduced from~\cite[Theorem~II.7.18]{Har}.

\begin{cor}  \label{lct-proj-local}
\textup{(a)} Let $Y\sub X$ be an open subscheme.
 Then for any cosheaf\/ $\gF\in X\ctrh^\lct_\prj$ the cosheaf\/
$\gF|_Y$ belongs to $Y\ctrh^\lct_\prj$. \par
\textup{(b)} In the situation of~\textup{(a)}, the corestriction map\/
$\gF[Y]\rarrow\gF[X]$ is injective.  If the scheme $X$ is affine,
then $0\rarrow\gF[Y]\rarrow\gF[X]\rarrow\gF[X]/\gF[Y]\rarrow0$ is a split
short exact sequence of flat cotorsion\/ $\O(X)$\+modules. \par
\textup{(c)} Let $X=\bigcup_\alpha Y_\alpha$ be an open covering.
Then a locally contraherent cosheaf on $X$ belongs
to $X\ctrh^\lct_\prj$ if and only if its restrictions to $Y_\alpha$
belong to $Y_\alpha\ctrh^\lct_\prj$ for all\/~$\alpha$.
\end{cor}

\begin{proof}
 Part~(a): let $j$~denote the open embedding $Y\rarrow X$.
 Then by Theorem~\ref{proj-lct-classification} one has
$j^!\gF\simeq j^!(\prod_{x\in X}\iota_x{}_!\widecheck F_x)\simeq
\prod_{x\in X}j^!\iota_x{}_!\widecheck F_x$.
 Furthermore, it was explained in the proof of the same Theorem
that $j^!\iota_z{}_!\widecheck F_z=0$ for any $z\in X\setminus Y$.
 Denoting by~$\kappa_y$ the natural morphism $\Spec\O_{y,X}
\rarrow Y$ for a point $y\in Y$, one clearly has
$j^!\iota_y{}_!\widecheck F_y\simeq\kappa_y{}_!\widecheck F_y$.
 Hence the isomorphism $j^!\gF\simeq\prod_{y\in Y}
\kappa_y{}_!\widecheck F_y$, proving that $j^!\gF\in Y\ctrh^\lct_\prj$.

 Part~(b): by the definition, one has
$(\iota_x{}_!\widecheck F_x)[X]\simeq F_x$.
 Since the cosections over a quasi-compact quasi-separated scheme
commute with infinite products of contraherent cosheaves, we conclude
from the computation above that $\gF[U]\simeq\prod_{x\in U}F_x$ for
any quasi-compact open subscheme $U\sub X$.
 The cosheaf axiom~\eqref{cosheaf-definition} now implies
the isomorphism $\gF[X]\simeq\varinjlim_{U\sub X}\gF[U]$, the filtered
inductive limit being taken over all quasi-compact open subschemes
$U\sub X$, and similarly for $\gF[Y]$.
 So the corestriction map $\gF[Y]\rarrow\gF[X]$ is the embedding of
a direct summand.

 Part~(c): the ``only if'' assertion follows from part~(a); let us
prove the ``if''.
 A transfinite induction in (any well-ordering of) the set of
indices~$\alpha$ reduces the question to the following assertion.
 Suppose $X$ is presented as the union of two open subschemes $W\cup Y$.
 Furthermore, the restriction of a locally contraherent cosheaf $\gF$ on
$X$ onto $W$ is identified with the direct product
$\prod_{w\in W}\kappa_w{}_!\widecheck F_w$, where $F_w$ are some
free $\widehat\O_{w,X}$\+contramodules, while $\kappa_w$~are
the natural morphisms $\Spec\O_{w,X}\rarrow W$.
 Assume also that the restriction of $\gF$ onto $Y$ is a projective
locally cotorsion contraherent cosheaf.
 Then there exist some free $\widehat\O_{z,X}$\+contramodules $F_z$
defined for all $z\in X\setminus W$ and an isomorphism of locally
contraherent cosheaves $\gF\simeq\prod_{x\in X}\iota_x{}_!\widecheck F_x$
whose restriction to $W$ coincides with the given isomorphism
$\gF|_W\simeq\prod_{w\in W}\kappa_w{}_!\widecheck F_w$.

 Indeed, by Theorem~\ref{proj-lct-classification} we have
$\gF|_Y\simeq\prod_{y\in Y}\iota'_y{}_!\widecheck F'_y$, where
$\iota'_y$~denote the natural morphisms $\Spec\O_{y,X}\rarrow Y$
and $F'_y$ are some free contramodules over $\widehat\O_{y,X}$.
 Restricting to the intersection $V=W\cap Y$, we obtain
an isomorphism of contraherent cosheaves
$\prod_{v\in V}\kappa'_v{}_!\widecheck F_v\simeq
\prod_{v\in V}\kappa'_v{}_!\widecheck F'_v$ on $V$, where
$\kappa'_v$~denotes the natural morphisms $\Spec\O_{v,X}\rarrow V$.
 It is clear from the arguments in the second half of the proof
of Theorem~\ref{proj-lct-classification}(b) that such an isomorphism
of infinite products induces an isomorphism of
$\widehat\O_{v,X}$\+contramodules $F_v\simeq F'_v$.
 Let us identify $F'_v$ with $F_v$ using this isomorphism, and set
$F_y=F'_y$ for $y\in X\setminus W\sub Y$.
 Then our isomorphism of infinite products can be viewed as
an automorphism~$\phi$ of the contraherent cosheaf
$\prod_{v\in V}\kappa'_v{}_!\widecheck F_v$ on~$V$.
 We would like to show that the cosheaf $\gF$ is isomorphic to
$\prod_{x\in X}\iota_x{}_!\widecheck F_x$.

 Since a cosheaf of $\O_X$\+modules is determined by its restrictions
to $W$ and $Y$ together with the induced isomorphism between
the restrictions of these to the intersection $V=W\cap Y$, it
suffices to check that the automorphism~$\phi$ can be lifted to
an automorphism of the contraherent cosheaf
$\prod_{v\in V}\iota'_v\widecheck F_v$ on~$Y$.
 In fact, the rings of endomorphisms of the two contraherent cosheaves
are isomorphic.
 Indeed, the functor of direct image of cosheaves of
$\O$\+modules~$j_!$ with respect to the open embedding $j\:V\rarrow Y$
is fully faithful by~\eqref{open-cosheaf-direct-inverse-adjunction};
the morphism~$j$ being quasi-compact and quasi-separated, this functor
also preserves infinite products.
\end{proof}

 We refer to Section~\ref{contraherent-tensor} for the definitions of
flat and $\bW$\+flat cosheaves of $\O_X$\+modules, and to
Section~\ref{clf-subsection} for the definition of antilocally flat
$\bW$\+locally contraherent cosheaves.

\begin{cor}  \label{lct-prj-flat}
 The classes of projective locally cotorsion contraherent cosheaves
and\/ $\bW$\+flat locally cotorsion\/ $\bW$\+locally contraherent
cosheaves on $X$ coincide.
 In particular, any\/ $\bW$\+flat locally cotorsion\/ $\bW$\+locally
contraherent cosheaf on $X$ is flat, contraherent, and antilocally flat.
\end{cor}

\begin{proof}
 All projective locally cotorsion contraherent cosheaves on $X$ are
flat by Corollary~\ref{proj-flat}(b) (for affine schemes)
and Corollary~\ref{lct-proj-local}(a).
 To prove that any $\bW$\+flat locally cotorsion $\bW$\+locally
contraherent cosheaf $\gF$ belongs to $X\ctrh^\lct_\prj$, pick an affine
open covering $U_\alpha$ of the scheme $X$ subordinate to~$\bW$.
 Then the $\O_X(U_\alpha)$\+modules $\gF[U_\alpha]$ are,
by the definition, flat and cotorsion.
 Hence the locally cotorsion contraherent cosheaves $\gF|_{U_\alpha}$
are projective, and it remains to apply
Corollary~\ref{lct-proj-local}(c).
\end{proof}

 At this point we can prove an assertion that was promised in
Section~\ref{coherent-coadjusted}.

\begin{proof}[Proof of
Corollary~\ref{coherent-flat-local}(b), the cotorsion case]
 Let $R$ be a Noetherian commutative ring and $R\rarrow S_\alpha$ be
a collection of homomorphisms of commutative rings for which
the corresponding collection of morphisms of affine schemes
$\Spec S_\alpha\rarrow\Spec R$ is an open covering.
 Let $P$ be a cotorsion $R$\+module for which all
the $S_\alpha$\+modules $\Hom_R(S_\alpha,P)$ are flat.
 We need to show that the $R$\+module $P$ is flat.

 Put $X=\Spec R$, and let $\bW$ be the open covering of $X$ formed
by all the affine open subschemes $W_\alpha=\Spec S_\alpha\subset X$.
 Consider the locally cotorsion contraherent cosheaf
$\gF=\widecheck P$ on~$X$.
 Then Corollary~\ref{coherent-flat-local}(a) implies that
the contraherent cosheaf $\gF$ on $X$ is $\bW$\+flat.
 By Corollary~\ref{lct-prj-flat}, it follows that $\gF$ is a flat
contraherent cosheaf on~$X$.
 Thus the $R$\+module $P=\gF[X]$ is flat.
\end{proof}

\begin{cor}  \label{loc-noetherian-proj-products}
 The full subcategory $X\ctrh^\lct_\prj$ is closed under infinite 
products in $X\ctrh^\lct$ or $X\ctrh$.
\end{cor}

\begin{proof}
 Easily deduced either from Corollary~\ref{lct-prj-flat}, or
directly from Theorem~\ref{proj-lct-classification}(b)
(cf.~\cite[Lemma~1.3.7]{Pweak}).
\end{proof}

 Let us point out once again that the full subcategory $X\ctrh_\prj$ is
\emph{not} closed under infinite products in $X\ctrh$ even for
Noetherian affine schemes~$X$.
 See Example~\ref{X-ctrh-prj-not-closed-under-products-ex}.

\begin{cor}  \label{proj-lct-direct}
 Let $f\:Y\rarrow X$ be a quasi-compact morphism of locally Noetherian
schemes.  Then \par
\textup{(a)} the functor of direct image of cosheaves of\/ $\O$\+modules
$f_!$ takes projective locally cotorsion contraherent cosheaves on $Y$
to locally cotorsion contraherent cosheaves on~$X$; \par
\textup{(b)} if the morphism~$f$ is flat, then the functor of direct
image of cosheaves of\/ $\O$\+modules $f_!$ takes projective locally
cotorsion contraherent cosheaves on $Y$ to projective locally cotorsion
contraherent cosheaves on~$X$.
\end{cor}

\begin{proof}
 Part~(a): let us show that the functor~$f_!$ takes locally cotorsion
contraherent cosheaves of the form $\P=\prod_{y\in Y}\kappa_y{}_!
\widecheck P_y$ on $Y$, where $\kappa_y$~are the natural morphisms
$\Spec\O_{y,Y}\rarrow Y$ and $P_y$ are some
$\widehat\O_{y,Y}$\+contramodules, to locally cotorsion contraherent
cosheaves of the same form on~$X$.
 More precisely, one has $f_!\P\simeq\prod_{x\in X}\iota_x{}_!
\widecheck P_x$, where $\iota_x$~are the natural morphisms
$\Spec\O_{x,X}\rarrow X$ and $P_x=\prod_{f(y)=x}P_y$,
the $\widehat\O_{x,X}$\+contramodule structures on $P_y$ being obtained
by the (contra)restriction of scalars with respect to the homomorphisms
of complete Noetherian rings $\widehat\O_{x,X}\rarrow\widehat\O_{y,Y}$
(see~\cite[Section~1.8]{Pweak} or~\cite[Section~2.9]{Pproperf}).

 Indeed, the morphism~$f$ being quasi-compact and quasi-separated,
the functor~$f_!$ preserves infinite products of cosheaves of
$\O$\+modules.
 So it suffices to consider the case of a locally cotorsion contraherent
cosheaf $\kappa_y{}_!\widecheck P_y$ on~$Y$.
 Now the composition of morphisms of schemes $\Spec\O_{y,Y}\rarrow Y
\rarrow X$ is equal to the composition $\Spec\O_{y,Y}\rarrow
\Spec\O_{x,X}\rarrow X$, and it remains to use the compatibility of
the direct images of cosheaves of $\O$\+modules with the compositions
of morphisms of schemes.
 Alternatively, part~(a) is a particular case of
Corollary~\ref{coflasque-direct}(b).

 To deduce part~(b), one can apply the assertion that
an $\widehat\O_{x,X}$\+contramodule is projective if and only if it is
a flat $\O_{x,X}$\+module~\cite[Corollary~B.8.2(c)]{Pweak}.
 Alternatively, use
the adjunction~\eqref{direct-inverse-cosheaf-lct-adjunction} together
with exactness of the inverse image of locally cotorsion locally
contraherent cosheaves with respect to a flat morphism of schemes.
\end{proof}

\subsection{Flat contraherent cosheaves}  \label{flat-contra-subsection}
 In this section we complete our treatment of flat contraherent
cosheaves on locally Noetherian schemes, which was started in
Sections~\ref{coherent-coadjusted} and~\ref{contraherent-tensor} and
continued in Sections~\ref{clf-subsection}
and~\ref{lct-projective-loc-noetherian}.

 A $\bW$\+locally contraherent cosheaf $\gM$ on a scheme $X$ is said to
have \emph{locally cotorsion dimension not exceeding~$d$} if
the cotorsion dimension of the $\O_X(U)$\+module $\gM[U]$ does not
exceed~$d$ for any affine open subscheme $U\sub X$ subordinate to~$\bW$
(cf.\ Sections~\ref{veryflat-cotors-dim-subsect}
and~\ref{finite-flat-lin-dim-subsect}).
 Clearly, the locally cotorsion dimension of a $\bW$\+locally
contraherent cosheaf does not change when the covering $\bW$ is replaced
by its refinement (see
Lemma~\ref{flat-veryflat-cotors-inj-dim-local}(c)).
 If a $\bW$\+locally contraherent cosheaf $\gM$ has a coresolution by
locally cotorsion $\bW$\+locally contraherent cosheaves in $X\lcth_\bW$,
then its locally cotorsion dimension is equal to the minimal length
of such coresolution.

\begin{lem}  \label{loc-cotors-dim}
 Let $X$ be a semi-separated Noetherian scheme of Krull dimension~$D$
with an open covering\/~$\bW$.
 Then \par
\textup{(a)} the coresolution dimension of any\/ $\bW$\+locally
contraherent cosheaf on $X$ with respect to the coresolving subcategory
$X\lcth_\bW^\lct\sub X\lcth_\bW$ does not exceed~$D$; \par
\textup{(b)} the coresolution dimension of any\/ $\bW$\+flat
$\bW$\+locally contraherent cosheaf on $X$ with respect to
the coresolving subcategory of projective locally cotorsion contraherent
cosheaves $X\ctrh^\lct_\prj\sub X\lcth_\bW^\fl$ does not exceed~$D$.
\end{lem}

\begin{proof}
 The locally cotorsion dimension of any locally contraherent cosheaf
on a locally Noetherian scheme of Krull dimension $D$ does not exceed $D$
by Corollary~\ref{raynaud-gruson-cotors-cor}.
 The full subcategory $X\lcth_\bW^\lct$ is coresolving in $X\lcth_\bW$ 
due to Corollary~\ref{clp-cor}(a) or~\ref{clf-cor}(a) and the results
of Section~\ref{locally-contraherent}, hence the coresolution dimension
with respect to $X\lcth_\bW^\lct$ is well-defined and obviously equal
to the locally cotorsion dimension.
 By Corollary~\ref{lct-prj-flat}, we have
$X\ctrh^\lct_\prj=X\lcth_\bW^\lct\cap X\lcth_\bW^\fl$.
 The full subcategory $X\ctrh^\lct_\prj$ is coresolving in
$X\lcth_\bW^\fl$ by Corollaries~\ref{clf-cor}(a),
\ref{clf-independence}, and~\ref{clf-noetherian-flat};
so part~(b) follows from~(a) in view of the dual version of
Corollary~\ref{fdim-subcategory-cor}.
 For the generalization of part~(b) to non-semi-separated Noetherian
schemes, see Corollary~\ref{lct-prj-envelope}(b) below.
\end{proof}

 Part~(a) of the following corollary should be compared with
Corollary~\ref{finite-krull-flat-clf-cor} below.
 See Remark~\ref{flat-antilocally-flat-remark} for a discussion.

\begin{cor}  \label{finite-krull-flat-contraherent}
\textup{(a)} On a semi-separated Noetherian scheme $X$ of finite Krull
dimension, the classes of\/ $\bW$\+flat\/ $\bW$\+locally contraherent
cosheaves and antilocally flat contraherent cosheaves coincide. \par
\textup{(b)} On a locally Noetherian scheme $X$ of finite Krull
dimension, any\/ $\bW$\+flat\/ $\bW$\+locally contraherent cosheaf is
flat and contraherent, so the categories $X\lcth_\bW^\fl$ and
$X\ctrh^\fl$ coincide. \par
\textup{(c)} On a locally Noetherian scheme $X$ of finite Krull
dimension, any flat contraherent cosheaf\/ $\gF$ is coflasque, so
one has $X\ctrh^\fl\sub X\ctrh_\cfq$.
 If\/ $X$ is affine and $Y\sub X$ is an open subscheme, then
$0\rarrow\gF[Y]\rarrow\gF[X]\rarrow\gF[X]/\gF[Y]\rarrow0$ is
a short exact sequence of flat contraadjusted\/ $\O(X)$\+modules.
\end{cor}

\begin{proof}
 Part~(a): any antilocally flat contraherent cosheaf on $X$ is flat
by Corollary~\ref{clf-noetherian-flat} (see also
Corollary~\ref{clf-independence}).
 On the other hand, by Lemma~\ref{loc-cotors-dim}(b),
any $\bW$\+flat $\bW$\+locally contraherent cosheaf $\gF$ on $X$ has
a finite coresolution by cosheaves from $X\ctrh^\lct_\prj$,
which are antilocally flat by definition.
 By Corollary~\ref{clf-characterizations}(b), it follows that $\gF$
is antilocally flat.

 Part~(b): given an affine open subscheme $U\sub X$, denote by
$\bW|_U$ the collection of all open subsets $U\cap W$ with $W\in\bW$.
 Then the restriction $\gF|_U$ of any $\bW$\+flat $\bW$\+locally
contraherent cosheaf $\gF$ on $X$ onto $U$ is, by definition,
$\bW|_U$\+flat and $\bW|_U$\+locally contraherent.
 Applying part~(a), we conclude that the cosheaf $\gF|_U$ is
contraherent and flat.
 This being true for any affine open subscheme $U\sub X$ means
precisely that $\gF$ is contraherent and flat on~$X$.

 Part~(c): coflasqueness being a local property by
Lemma~\ref{coflasque-cosheaves}(a), it suffices to consider
the case of an affine scheme~$X$.
 Then we have seen above that $\gF$ has a finite coresolution
by cosheaves from $X\ctrh^\lct_\prj$, which have the properties
listed in part~(c) by Corollary~\ref{lct-proj-local}(b).
 By Corollary~\ref{coflasque-acyclic}(b), our coresolution is
an exact sequence in $X\ctrh_\cfq$ and the cosheaf $\gF$ is
coflasque.

 Furthermore, by Corollary~\ref{coflasque-acyclic}(c), the related
sequence of cosections over any open subscheme $Y\sub X$ is exact.
 Passing to the corestriction maps related to the pair of embedded
open subschemes $Y\sub X$, we obtain an injective morphism of
exact sequences of $\O(X)$\+modules; so the related sequence of
cokernels is also exact.
 All of its terms except perhaps the leftmost one being flat, it
follows that the leftmost term $\gF[X]/\gF[Y]$ is a flat
$\O(X)$\+module, too.
 The $\O(X)$\+module $\gF[Y]$ is contraadjusted by
Corollary~\ref{coflasque-direct}(a) applied to the embedding
morphism $Y\rarrow X$.
\end{proof}

 Now we can prove another assertion promised in
Section~\ref{coherent-coadjusted}.

\begin{proof}[Proof of
Corollary~\ref{coherent-flat-local}(b), the finite Krull dimension case]
 Let $R$ be a Noetherian commutative ring of finite Krull dimension
and $R\rarrow S_\alpha$ be a collection of homomorphisms of commutative
rings for which the corresponding collection of morphisms of affine
schemes $\Spec S_\alpha\rarrow\Spec R$ is an open covering.
 Let $P$ be a contraadjusted $R$\+module for which all
the $S_\alpha$\+modules $\Hom_R(S_\alpha,P)$ are flat.
 We have to show that the $R$\+module $P$ is flat.

 Put $X=\Spec R$, and let $\bW$ be the open covering of $X$ formed
by all the affine open subschemes $W_\alpha=\Spec S_\alpha\subset X$.
 Consider the contraherent cosheaf $\gF=\widecheck P$ on~$X$.
 Then Corollary~\ref{coherent-flat-local}(a) tells us that
the contraherent cosheaf $\gF$ on $X$ is $\bW$\+flat.
 By Corollary~\ref{finite-krull-flat-contraherent}(b), it follows
that $\gF$ is a flat contraherent cosheaf on~$X$.
 Hence the $R$\+module $P=\gF[X]$ is flat.
\end{proof}

\begin{lem}  \label{flat-direct-lemma}
 Let $f\:Y\rarrow X$ be a flat quasi-compact morphism from
a semi-separated locally Noetherian scheme $Y$ of finite Krull
dimension to a locally Noetherian scheme~$X$.
 Then the functor of direct image of cosheaves of\/ $\O$\+modules
$f_!$ takes flat contraherent cosheaves on $Y$ to flat
contraherent cosheaves on~$X$, and induces an exact functor
$Y\ctrh^\fl\rarrow X\ctrh^\fl$.
\end{lem}

\begin{proof}
 Clearly, it suffices to consider the case of a Noetherian affine
scheme~$X$.
 Then the scheme $Y$ is Noetherian and semi-separated.
 By Corollary~\ref{finite-krull-flat-contraherent}(a), any flat
contraherent cosheaf on $Y$ is antilocally flat.
 By Corollary~\ref{clp-direct}(c), the functor~$f_!$ restricts
to an exact functor $Y\ctrh_\alf\rarrow X\ctrh_\alf$.
 By Corollary~\ref{clf-noetherian-flat}, any antilocally flat
contraherent cosheaf on~$X$ is flat.
 (A proof working in a greater generality will be given below in
Corollary~\ref{finite-krull-flat-direct}(b).)
\end{proof}

 Let $X$ be a Noetherian scheme of finite Krull dimension~$D$ with
an open covering $\bW$ and a finite affine open covering $U_\alpha$
subordinate to~$\bW$.

\begin{cor}  \label{finite-krull-contrah-projective}
\textup{(a)} There are enough projective objects in the exact
categories of locally contraherent cosheaves $X\lcth_\bW$ and $X\lcth$
on $X$, and all these projective objects belong to the full subcategory
of contraherent cosheaves $X\ctrh$.
 The full subcategories of projective objects in the three exact
categories $X\ctrh\sub X\lcth_\bW\sub X\lcth$ coincide. \par
\textup{(b)} A contraherent cosheaf on $X$ is projective if and only
if it is isomorphic to a direct summand of a finite direct sum of
the direct images of projective contraherent cosheaves from~$U_\alpha$.
 In particular, any projective contraherent cosheaf on $X$ is flat
(and consequently, coflasque).
\end{cor}

\begin{proof}
 The argument is similar to the proofs of Lemmas~\ref{clf-cover},
\ref{loc-contra-proj}, and~\ref{loc-lct-proj}; the only difference
is that, the scheme $X$ being not necessarily semi-separated,
one has to also use Lemma~\ref{flat-direct-lemma}.
 Let $j_\alpha\:U_\alpha\rarrow X$ denote the open embedding
morphisms.
 According to Corollary~\ref{proj-flat}(a) and
Lemma~\ref{flat-direct-lemma}, the cosheaf of $\O_X$\+modules
$j_\alpha{}_!\gF_\alpha$ is flat and contraherent for any
projective contraherent cosheaf $\gF_\alpha$ on~$U_\alpha$.
 The adjunction \eqref{open-cosheaf-direct-inverse-adjunction}
or~\eqref{direct-inverse-cosheaf-lct-adjunction} shows that it is
also a projective object in $X\lcth$.

 It remains to show that there are enough projective objects of
the form $\bigoplus_\alpha j_\alpha{}_!\gF_\alpha$ in $X\lcth_\bW$.
 Let $\Q$ be a $\bW$\+locally contraherent cosheaf on~$X$.
 For every~$\alpha$, pick an admissible epimorphism onto
the contraherent cosheaf $j_\alpha^!\Q$ from a projective
contraherent cosheaf $\gF_\alpha$ on~$U_\alpha$.
 Then the same adjunction provides a natural morphism
$\bigoplus_\alpha j_\alpha{}_!\gF_\alpha\rarrow\Q$.
 This morphism is an admissible epimorphism of locally contraherent
cosheaves, because it is so in restriction to each open
subscheme $U_\alpha\sub X$.

 Finally, any flat contraherent cosheaf on $X$ is coflasque by
Corollary~\ref{finite-krull-flat-contraherent}(c).
\end{proof}

 As in Section~\ref{projective-contraherent}, we denote the category
of projective contraherent cosheaves on $X$ by $X\ctrh_\prj$.
 Clearly, the objects of $X\ctrh_\prj$ are also the projective objects
of the exact category $X\ctrh^\fl$ (and there are enough of them).
 We will see below in this section (in the proof of
Corollary~\ref{lct-prj-envelope}) that the objects of
$X\ctrh^\lct_\prj$ are the \emph{injective} objects of $X\ctrh^\fl$
(and there are enough of them).

\begin{lem}  \label{coflasque-lct-envelope}
\textup{(a)} Any coflasque contraherent cosheaf\/ $\gE$ on $X$ can be
included in a short exact sequence\/ $0\rarrow\gE\rarrow\P\rarrow\gF
\rarrow0$ in $X\ctrh_\cfq$, where\/ $\P$ is a coflasque locally
cotorsion contraherent cosheaf on $X$ and\/ $\gF$ is a finitely
iterated extension of the direct images of flat contraherent cosheaves
from~$U_\alpha$. \par
\textup{(b)} The coresolution dimension of any coflasque
contraherent cosheaf on $X$ with respect to the coresolving subcategory
of coflasque locally cotorsion contraherent cosheaves
$X\ctrh^\lct_\cfq\sub X\ctrh_\cfq$ does not exceed~$D$.
\end{lem}

\begin{proof}
 Part~(a): the proof is similar to that of Lemma~\ref{lct-envelope};
the only difference is that, the scheme $X$ being not necessarily
semi-separated, one has to also use Corollary~\ref{coflasque-direct}.
 Notice that flat contraherent cosheaves on $U_\alpha$ are coflasque
by Lemma~\ref{finite-krull-flat-contraherent}(c) and
finitely iterated extensions of coflasque contraherent cosheaves in
$X\ctrh$ are coflasque by Lemma~\ref{coflasque-acyclic}(a).

 One proceeds by induction in a linear ordering of
the indices~$\alpha$, considering the open subscheme
$V=\bigcup_{\alpha<\beta} U_\alpha\sub X$.
 Assume that we have constructed a short exact sequence $0\rarrow
\gE\rarrow\gK\rarrow\gL\rarrow0$ of coflasque contraherent cosheaves
on $X$ such that the contraherent cosheaf $\gK|_V$ is locally cotorsion,
while the cosheaf $\gL$ on $X$ is a finitely iterated extension of
the direct images of flat contraherent cosheaves from the affine
open subschemes $U_\alpha$ with $\alpha<\beta$.
 Let $j\:U=U_\beta\rarrow X$ be the identity open embedding.

 Pick a short exact sequence $0\rarrow j^!\gK\rarrow\Q\rarrow\gG
\rarrow0$ of contraherent cosheaves on $U$ such that $\Q$ is a locally
cotorsion contraherent cosheaf and $\gG$ is a flat contraherent cosheaf
(see Theorem~\ref{flat-cover-thm}(a) and
Corollary~\ref{coherent-flat-local}(a)); then the cosheaf $\gG$ is
coflasque by Lemma~\ref{finite-krull-flat-contraherent}(c) and
the cosheaf $\Q$ is coflasque by Lemma~\ref{coflasque-acyclic}(a).
 By Corollary~\ref{coflasque-direct}(a), the related sequence of direct
images $0\rarrow j_!j^!\gK\rarrow j_!\Q\rarrow j_!\gG\rarrow0$ is
a short exact sequence of coflasque contraherent cosheaves on $X$;
by part~(b) of the same corollary, the cosheaf $j_!\Q$ belongs to
$X\ctrh^\lct_\cfq$.

 Let $0\rarrow\gK\rarrow\R\rarrow j_!\gG\rarrow0$ denote
the push-out of the short exact sequence $0\rarrow j_!j^!\gK
\rarrow j_!\Q\rarrow j_!\gG\rarrow0$ with respect to the natural
(adjunction) morphism $j_!j^!\gK\rarrow\gK$.
 We will show that the coflasque contraherent cosheaf $\R$ on $X$ is
locally cotorsion in restriction to $U\cup V$.

 Indeed, in the restriction to $U$ one has $j^!\R\simeq\Q$.
 On the other hand, denoting by~$j'$ the embedding
$U\cap V\rarrow V$, one has $(j_!\gG)|_V\simeq j'_!(\gG|_{U\cap V})$.
 The contraherent cosheaf $\gK|_{U\cap V}$ being locally cotorsion,
so is the cokernel $\gG|_{U\cap V}$ of the admissible monomorphism of
locally cotorsion contraherent cosheaves $\gK|_{U\cap V}\rarrow
\Q|_{U\cap V}$.

 By Corollary~\ref{coflasque-direct}(b), $j'_!(\gG|_{U\cap V})$ is
a coflasque locally cotorsion contraherent cosheaf on~$V$.
 Now in the short exact sequence $0\rarrow\gK|_V\rarrow\R|_V\rarrow
(j_!\gG)|_V\rarrow0$ the middle term is locally cotorsion, since
the two other terms are.

 Finally, the composition $\gE\rarrow\gK\rarrow\R$ of admissible
monomorphisms in $X\ctrh_\cfq$ is again an admissible monomorphism
with the cokernel isomorphic to an extension of the flat
contraherent cosheaves $j_!\gG$ and $\gL$.

 Part~(b): the full subcategory $X\ctrh^\lct_\cfq$ is coresolving
in $X\ctrh_\cfq$ by part~(a), and the coresolution dimension does not
exceed $D$ by Corollary~\ref{raynaud-gruson-cotors-cor}
(cf.\ the proof of Lemma~\ref{loc-cotors-dim}(a)).
\end{proof}

\begin{cor} \label{lct-prj-envelope}
\textup{(a)} Any flat contraherent cosheaf\/ $\gG$ on $X$ can be
included in a short exact sequence\/ $0\rarrow\gG\rarrow\P\rarrow\gF
\rarrow0$ in $X\ctrh^\fl$, where\/ $\P$ is a projective locally
cotorsion contraherent cosheaf on $X$ and\/ $\gF$ is a finitely
iterated extension of the direct images of flat contraherent cosheaves
from~$U_\alpha$. \par
\textup{(b)} The coresolution dimension of any flat contraherent
cosheaf on $X$ with respect to the coresolving subcategory of projective
locally cotorsion contraherent cosheaves $X\ctrh^\lct_\prj\sub
X\ctrh^\fl$ does not exceed~$D$; the homological dimension of
the exact category $X\ctrh^\fl$ does not exceed~$D$; and
the resolution dimension of any flat contraherent cosheaf on $X$ with
respect to the resolving category of projective contraherent cosheaves
$X\ctrh_\prj\sub X\ctrh^\fl$ does not exceed~$D$.
\end{cor}

\begin{proof}
 Part~(a) follows from Lemma~\ref{coflasque-lct-envelope}(a) together
with Corollaries~\ref{finite-krull-flat-contraherent}(c)
and~\ref{lct-prj-flat}.
 Part~(b): the full subcategory $X\ctrh^\lct_\prj$ is coresolving
in $X\ctrh^\fl$ by part~(a); and the coresolution dimension does not
exceed~$D$ by Lemma~\ref{coflasque-lct-envelope}(b) and the dual version
of Corollary~\ref{fdim-subcategory-cor}.
 This proves the first assertion.

 By the dual version of Proposition~\ref{infinite-resolutions}(a)
or~\ref{finite-resolutions}, it follows that the natural functor
$\Hot^+(X\ctrh^\lct_\prj)\rarrow\sD^+(X\ctrh^\fl)$ is fully
faithful (the exact category structure on $X\ctrh^\lct_\prj$
being trivial).
 Applying the first assertion of part~(b) again, we conclude that
the homological dimension of the exact category $X\ctrh^\fl$ does not
exceed $D$, and consequently that the resolution dimension of
any object of $X\ctrh^\fl$ with respect to its subcategory of
projective objects $X\ctrh_\prj$ does not exceed~$D$.

 Perhaps the following way to phrase the proof of the second two
assertions of part~(b) is more illuminating.
 The full subcategory $X\ctrh^\lct_\prj$ is coresolving and closed
under direct summands in $X\ctrh^\fl$, and the induced exact
category structure on $X\ctrh^\lct_\prj$ is trivial (any short
exact sequence in $X\ctrh^\lct_\prj$ is split).
 It follows from these observations already that $X\ctrh^\lct_\prj$
is the full subcategory of injective objects in $X\ctrh^\fl$, and
there are enough such injective objects.
 Therefore, the supremum of the coresolution dimensions of
the objects of $X\ctrh^\fl$ with respect to $X\ctrh^\lct_\prj$
is the homological dimension of the exact category $X\ctrh^\fl$.
\end{proof}

\begin{lem}  \label{flat-iterated-extension-covering}
\textup{(a)} Any coflasque contraherent cosheaf\/ $\gE$ on $X$ can be
included in a short exact sequence\/ $0\rarrow\P\rarrow\gF\rarrow\gE
\rarrow0$ in $X\ctrh_\cfq$, where\/ $\P$ is a coflasque locally
cotorsion contraherent cosheaf on $X$ and\/ $\gF$ is a finitely
iterated extension of the direct images of flat contraherent cosheaves
from~$U_\alpha$. \par
\textup{(b)} Any flat contraherent cosheaf\/ $\gG$ on $X$ can be
included in a short exact sequence\/ $0\rarrow\P\rarrow\gF\rarrow\gG
\rarrow0$ in $X\ctrh^\fl$, where $\P$ is a projective locally cotorsion
contraherent cosheaf on $X$ and\/ $\gF$ is a finitely iterated
extension of the direct images of flat contraherent cosheaves
from~$U_\alpha$.
\end{lem}

\begin{proof}
 According to Corollary~\ref{finite-krull-contrah-projective},
there is an admissible epimorphism $\Q\rarrow\gG$ in the exact
category $X\ctrh$ onto any given contraherent cosheaf $\gG$ from
a finite direct sum $\Q$ of the direct images of flat (and even
projective) contraherent cosheaves from~$U_\alpha$.
 Furthermore, any admissible epimorphism in $X\ctrh$ between objects
from $X\ctrh_\cfq$ is also an admissible epimorphism in $X\ctrh_\cfq$
(by Lemma~\ref{coflasque-acyclic}(b)), and any admissible epimorphism
in $X\ctrh$ between objects from $X\ctrh^\fl$ is an admissible
epimorphism in $X\ctrh^\fl$.
 The rest of the argument in both parts~(a) and~(b) is similar to
the proofs of Lemmas~\ref{clf-cover} and~\ref{clp-cover},
and based on Corollaries~\ref{coflasque-lct-envelope}(a) and
\ref{lct-prj-envelope}(a), respectively (applied to the kernel of
the morphism $\Q\rarrow\gG$).
 See Lemma~\ref{salce-lemma}(a) for a generalization.
 One can also deduce part~(b) as a particular case of part~(a)
using Corollaries~\ref{finite-krull-flat-contraherent}(c)
and~\ref{lct-prj-flat}.
\end{proof}

 Let us emphasize that the full subcategory $X\ctrh^\lct$ is \emph{not}
known to be coresolving in $X\ctrh$ when $X$ is not semi-separated
(cf.~\cite[Theorem~2.2]{SlSt}).
 The following corollary is to be compared with
Corollaries~\ref{homotopy-derived-contraderived-cor}
and~\ref{acycl=bctracycl=bcoacycl-in-alf} above;
see also Corollaries~\ref{finite-krull-ctrh-lcth-derived},
\ref{derived-contra-lct-cor}, and
Theorem~\ref{derived-inj-proj-resolutions} below.

\begin{cor}  \label{finite-krull-derived-equivalences} \hbadness=1300
\textup{(a)} For any Noetherian scheme $X$ of finite Krull dimension,
the natural functors\/ $\Hot(X\ctrh_\prj)$, $\Hot(X\ctrh^\lct_\prj)
\rarrow\sD^\abs(X\ctrh^\fl)\rarrow\sD^{\ctr=\bctr=\bco}(X\ctrh^\fl)
\rarrow\sD(X\ctrh^\fl)$ are equivalences of triangulated categories, as
are the natural functors $\Hot^\pm(X\ctrh_\prj)$,
$\Hot^\pm(X\ctrh^\lct_\prj)\rarrow\sD^{\abs\pm}(X\ctrh^\fl)\rarrow
\sD^\pm(X\ctrh^\fl)$ and\/ $\Hot^\b(X\ctrh_\prj)$,
$\Hot^\b(X\ctrh^\lct_\prj)\rarrow\sD^\b(X\ctrh^\fl)$ . \par
\textup{(b)} For any locally Noetherian scheme $X$ with an open
covering\/ $\bW$, the natural functors\/ $\Hot^-(X\ctrh^\lct_\prj)
\rarrow\sD^-(X\ctrh^\lct)\rarrow\sD^-(X\lcth_\bW^\lct)\rarrow
\sD^-(X\lcth^\lct)$ are equivalences of triangulated categories. \par
\textup{(c)} For any Noetherian scheme $X$ of finite Krull dimension
with an open covering\/ $\bW$, the natural functors\/
$\Hot^-(X\ctrh_\prj)\rarrow\sD^-(X\ctrh^\fl)\rarrow\sD^-(X\ctrh)
\rarrow\sD^-(X\lcth_\bW)\rarrow\sD^-(X\lcth)$ are equivalences of
triangulated categories. \par
\textup{(d)} For any Noetherian scheme $X$ of finite Krull dimension
with an open covering\/ $\bW$, the natural functors\/
$\sD^-(X\lcth_\bW^\lct)\rarrow\sD^-(X\lcth_\bW)$ and\/
$\sD^-(X\lcth^\lct)\rarrow\sD^-(X\lcth)$ are equivalences of
triangulated categories.
\end{cor}

\begin{proof}
 In view of Corollary~\ref{lct-prj-envelope}(b), all assertions of
part~(a) follow from Proposition~\ref{finite-resolutions} (together
with its dual version) and Lemma~\ref{psemi-remark21}.
 For the assertions concerning the Becker contraderived and coderived
categories, see Theorem~\ref{finite-homol-dim-becker-co-contra-derived}.
 Part~(b) is provided by Proposition~\ref{infinite-resolutions}(a)
together with Theorem~\ref{proj-lct-classification}(a), while
part~(c) follows from the same Proposition together with
Corollary~\ref{finite-krull-contrah-projective}.
 Finally, part~(d) is obtained by comparing parts~(a\+c).
\end{proof}

 It follows from Corollary~\ref{finite-krull-derived-equivalences}
that the $\Ext$ groups computed in the exact categories $X\lcth_\bW$
and $X\lcth_\bW^\lct$ agree with each other and with the $\Ext$
groups computed in the exact categories $X\lcth$ and $X\lcth^\lct$.
 Besides, these also agree with the $\Ext$ groups computed in
the exact categories $X\ctrh^\fl$ and $X\ctrh^\lct_\prj$ (the latter
being endowed with the trivial exact category structure).
 As in Section~\ref{clp-subsection}, we denote these $\Ext$ groups by
$\Ext^{X,*}({-},{-})$.

\begin{cor}  \label{finite-krull-flat-clf-cor}
\textup{(a)} Let $X$ be a Noetherian scheme of finite Krull dimension
with an open covering\/~$\bW$.
 Then one has\/ $\Ext^{X,>0}(\gF,\Q)=0$ for any flat contraherent
cosheaf\/ $\gF$ and locally cotorsion\/ $\bW$\+locally contraherent
cosheaf\/ $\Q$ on~$X$.
 Consequently, any flat contraherent cosheaf on $X$ is antilocally
flat. \par
\textup{(b)} Let $X$ be a Noetherian scheme of finite Krull dimension
with a finite affine open covering $X=\bigcup_\alpha U_\alpha$.
 Then a contraherent cosheaf on $X$ is flat if and only if it is
a direct summand of a finitely iterated extension of the direct images
of flat contraherent cosheaves from~$U_\alpha$.
\end{cor}

\begin{proof}
 According to Corollary~\ref{lct-prj-envelope}(b), any flat contraherent
cosheaf on $X$ has a finite coresolution by projective locally
cotorsion contraherent cosheaves.
 Since the $\Ext$ groups in the exact categories $X\lcth_\bW$
and $X\lcth_\bW^\lct$ agree, the assertion~(a) follows.
 Now the ``only if'' assertion in part~(b) is easily deduced from
Lemma~\ref{flat-iterated-extension-covering}(b) together with part~(a),
while the ``if'' is provided by Lemma~\ref{flat-direct-lemma}.
\end{proof}

\begin{rem} \label{flat-antilocally-flat-remark}
 For the sake of clarity and avoiding confusion, let us recall and
summarize our results concerning flat and antilocally flat
contraherent cosheaves.

 On a quasi-compact semi-separated scheme, the class of antilocally
flat $\bW$\+locally contraherent cosheaves does not depend on the choice
of an open covering~$\bW$ (by Corollary~\ref{clf-independence}).
 On a semi-separated coherent scheme, any antilocally flat contraherent cosheaf is flat (Corollary~\ref{clf-noetherian-flat}).

 On a locally Noetherian scheme $X$ of finite Krull dimension, any
$\bW$\+flat $\bW$\+locally contraherent cosheaf is flat and contraherent
(Corollary~\ref{finite-krull-flat-contraherent}(b)).
 So the class of $\bW$\+flat $\bW$\+locally contraherent cosheaves on
$X$ does not depend on an open covering~$\bW$.
 On a Noetherian scheme of finite Krull dimension, any flat contraherent
cosheaf is antilocally flat (also as a $\bW$\+locally contraherent
cosheaf for any open covering~$\bW$) by
Corollary~\ref{finite-krull-flat-clf-cor}.

 In particular, on a semi-separated Noetherian scheme of finite Krull
dimension, the classes of flat and antilocally flat contraherent
cosheaves coincide (Corollary~\ref{finite-krull-flat-contraherent}(a)).
\end{rem}

 In the next corollary we use the terminology of
Section~\ref{cotorsion-prelim-subsect}.

\begin{cor} \label{coflasque-flat-loc-cotorsion-cotorsion-pair}
 Let $X$ be a Noetherian scheme of finite Krull dimension.
 Then the pair of classes of objects $(X\ctrh^\fl$, $X\ctrh^\lct_\cfq)$
is a hereditary complete cotorsion pair in the exact category of
coflasque contraherent cosheaves $X\ctrh_\cfq$.
\end{cor}

\begin{proof}
 The $\Ext^1$ (and moreover, $\Ext^{>0}$) orthogonality holds by
Corollary~\ref{finite-krull-flat-clf-cor}(a),
the special preenvelope sequences are provided by
Lemma~\ref{coflasque-lct-envelope}(a) with
Corollary~\ref{finite-krull-flat-clf-cor}(b), and
the special precover sequences are provided by
Lemma~\ref{flat-iterated-extension-covering}(a).
 It remains to notice that both the classes $X\ctrh^\fl$ and
$X\ctrh^\lct_\cfq$ are obviously closed under direct summands in
$X\ctrh$, and refer to Lemma~\ref{cotorsion-pair-direct-summands-lemma}.
\end{proof}

 The following corollary is to be compared with
Corollaries~\ref{coflasque-direct}, \ref{proj-direct-inverse},
\ref{clp-direct}, \ref{proj-direct-gen}, and~\ref{proj-lct-direct}.

\begin{cor}  \label{finite-krull-flat-direct}
 Let $f\:Y\rarrow X$ be a quasi-compact morphism of locally Noetherian
schemes such that the scheme $Y$ has finite Krull dimension.
 Then \par
\textup{(a)} the functor of direct image of cosheaves of\/
$\O$\+modules $f_!$ takes flat contraherent cosheaves on $Y$ to
contraherent cosheaves on $X$, and induces an exact functor
$f_!\:Y\ctrh^\fl\rarrow X\ctrh$ between these exact categories; \par
\textup{(b)} if the morphism~$f$ is flat, then the functor of
direct image of cosheaves of\/ $\O$\+modules $f_!$ takes flat
contraherent cosheaves on $Y$ to flat contraherent cosheaves on $X$,
and induces an exact functor $f_!\:Y\ctrh^\fl\rarrow X\ctrh^\fl$
between these exact categories; \par
\textup{(c)} if the scheme $Y$ is Noetherian and the morphism~$f$
is very flat, then the functor of direct image of cosheaves of\/
$\O$\+modules $f_!$ takes projective contraherent cosheaves on $Y$
to projective contraherent cosheaves on~$X$.
\end{cor}

\begin{proof}
 Part~(a) is a particular case of Corollary~\ref{coflasque-direct}(a)
in view of Corollary~\ref{finite-krull-flat-contraherent}(c).
 In part~(b), one can assume that $X$ is an affine scheme, so
the scheme $Y$ is Noetherian of finite Krull dimension.
 It suffices to show that the $\O(X)$\+module $(f_!\gF)[X]=
\gF[Y]$ is flat for any flat contraherent cosheaf $\gF$ on~$Y$.
 For this purpose, consider a finite coresolution of the cosheaf $\gF$
by objects from $Y\ctrh^\lct_\prj$ in $Y\ctrh^\fl$, and apply
the functor $\Delta(Y,{-})$ to it.
 It follows from Corollaries~\ref{finite-krull-flat-contraherent}(c)
and~\ref{coflasque-acyclic}(c) that the sequence will remain exact.
 By Corollary~\ref{proj-lct-direct}(b), we obtain a finite 
coresolution of the module $\gF[Y]$ by flat $\O(X)$\+modules, implying
that the $\O(X)$\+module $\gF[Y]$ is also flat.
 Part~(c) follows from part~(a) or~(b) together with
Corollary~\ref{finite-krull-contrah-projective}(b) and
the adjunction~\eqref{direct-inverse-cosheaf-lct-adjunction}.
\end{proof}

\subsection{Homology of locally cotorsion locally
contraherent cosheaves}  \label{lct-homology-subsection}
 Let $X$ be a locally Noetherian scheme.
 Then the left derived functor of the functor of global cosections
$\Delta(X,{-})$ of locally cotorsion locally contraherent cosheaves
on $X$ is defined using projective resolutions in the exact category
$X\lcth^\lct$ (see Theorem~\ref{proj-lct-classification}(a)).

 Notice that the derived functors $\boL_*\Delta(X,{-})$ computed in
the exact category $X\lcth_\bW^\lct$ for a particular open covering
$\bW$ and in the whole category $X\lcth^\lct$ agree.
 The groups $\boL_i\Delta(X,\gE)$ are called the \emph{homology groups}
of a locally cotorsion locally contraherent cosheaf $\gE$ on
the scheme~$X$.

 Let us show that $\boL_{>0}\Delta(X,\gF)=0$ for any coflasque locally
cotorsion contraherent cosheaf~$\gF$.
 By Corollary~\ref{lct-proj-local}(b), any projective locally cotorsion
contraherent cosheaf on $X$ is coflasque.
 In view of Corollary~\ref{coflasque-acyclic}(b), any resolution of
an object of $X\ctrh^\lct_\cfq$ by objects of $X\ctrh^\lct_\prj$ in
the category $X\lcth^\lct$ is exact with respect to the exact
category $X\ctrh^\lct_\cfq$.
 By part~(c) of the same Corollary, the functor $\Delta(X,{-})$
preserves exactness of such sequences.
 Hence one can compute the derived functor $\boL_*\Delta(X,{-})$
using coflasque locally cotorsion contraherent resolutions.

 Similarly, let $X$ be a Noetherian scheme of finite Krull dimension.
 Then the left derived functor of the functor of global cosections
$\Delta(X,{-})$ of locally contraherent cosheaves on $X$ is defined
using projective resolutions in the exact category $X\lcth$
(see Corollary~\ref{finite-krull-contrah-projective}).
 The derived functors $\boL_*\Delta(X,{-})$ computed in the exact
category $X\lcth_\bW$ for any particular open covering $\bW$ and in
the whole category $X\lcth$ agree.
 Furthermore, one can compute the derived functor $\boL_*\Delta(X,{-})$
using coflasque resolutions.
 The groups $\boL_i\Delta(X,\gE)$ are called the \emph{homology groups}
of a locally contraherent cosheaf $\gE$ on the scheme~$X$.

 It is clear from the above that the two definitions agree when they
are both applicable; they also agree with the definitions given in
Section~\ref{homology-subsection} (cf.\ Lemma~\ref{coflasque-clp}).

 Using injective resolutions, one can similarly define the cohomology
of quasi-coherent sheaves $\boR^*\Gamma(X,{-})$ on a locally Noetherian
scheme~$X$ (cf.\ the discussions in
Sections~\ref{coflasque} and~\ref{dilute-subsect}).
 Injective quasi-coherent sheaves on $X$ being flasque, this definition
agrees with the conventional sheaf-theoretical one; so
$\boR^i\Gamma(X,\F)\simeq H^i(X,\F)$ for all $\F\in X\qcoh$ and $i\ge0$
when $X$ is locally Noetherian.
 The full subcategory of flasque quasi-coherent sheaves in $X\qcoh$ is
closed under extensions, cokernels of injective morphisms and infinite
direct sums; we denote the induced exact category structure on it
by $X\qcoh^\fq$.

 The following lemma is to be compared with
Lemmas~\ref{dil-cta-clp-finite-dim} and~\ref{lct-lin-clp-finite-dim}.

\begin{lem}  \label{finite-krull-ctrh-lcth-finite-dim}
\textup{(a)} Let $X$ be a locally Noetherian scheme of Krull
dimension~$D$.
 Then the coresolution dimension of any quasi-coherent sheaf on $X$
with respect to the coresolving subcategory of flasque quasi-coherent
sheaves $X\qcoh^\fq\sub X\qcoh$ does not exceed~$D$. \par
\textup{(b)} Let $X$ be a locally Noetherian scheme of Krull
dimension~$D$.
 Then the resolution dimension of any locally cotorsion locally
contraherent cosheaf on $X$ with respect to the resolving subcategory of
coflasque locally cotorsion contraherent cosheaves $X\ctrh_\cfq^\lct
\sub X\lcth^\lct$ does not exceed~$D$.
 Consequently, the same bound holds for the resolution dimension
of any object of $X\lcth^\lct$ with respect to the resolving subcategory
$X\ctrh^\lct$. \par
\textup{(c)} Let $X$ be a Noetherian scheme of Krull dimension~$D$.
 Then the resolution dimension of any locally contraherent
cosheaf on $X$ with respect to the resolving subcategory of coflasque
contraherent cosheaves $X\ctrh_\cfq\sub X\lcth$ does not exceed~$D$.
 Consequently, the same bound holds for the resolution dimension
of any object of $X\lcth$ with respect to the resolving subcategory
$X\ctrh$.
\end{lem}

\begin{proof}
 In part~(a), the full subcategory $X\qcoh^\fq$ is coresoving in
$X\qcoh$ because there are enough injective quasi-coherent cosheaves
on $X$ and all of them are flasque.
 In part~(b), the full subcategory $X\ctrh_\cfq^\lct$ is resolving
in $X\lcth^\lct$ by Theorem~\ref{proj-lct-classification}(a)
and Corollary~\ref{lct-proj-local}(b).
 In part~(c), the full subcategory $X\ctrh_\cfq$ is resolving in
$X\lcth$ by Corollary~\ref{finite-krull-contrah-projective}
(see also Corollary~\ref{coflasque-acyclic}(a\+b)).
 The remaining assertions follow from Lemma~\ref{grothendieck-vanishing}
and Corollary~\ref{coflasque-contraherent}.
\end{proof}

\begin{cor}  \label{coflasque-resolutions-finite}
\textup{(a)} Let $X$ be a locally Noetherian scheme of finite
Krull dimension.
 Then for any symbol\/ $\bst=\b$, $+$, $-$, $\empt$, $\abs+$, $\abs-$,
$\bco$, $\co$, or\/ $\abs$ the triangulated functor\/
$\sD^\st(X\qcoh^\fq)\rarrow\sD^\st(X\qcoh)$ induced by the embedding
of exact categories $X\qcoh^\fq\rarrow X\qcoh$ is an equivalence of
triangulated categories. \par
\textup{(b)} Let $X$ be a locally Noetherian scheme of finite
Krull dimension.
 Then for any symbol\/ $\bst=\b$, $+$, $-$, $\empt$, $\abs+$, $\abs-$,
$\bctr$, $\ctr$, or\/ $\abs$ the triangulated functor\/
$\sD^\st(X\ctrh^\lct_\cfq)\rarrow\sD^\st(X\lcth_\bW^\lct)$
induced by the embedding of exact categories $X\ctrh^\lct_\cfq\rarrow
X\lcth_\bW^\lct$ is an equivalence of triangulated categories. \par
 For any symbol\/ $\bst=\b$, $+$, $-$, $\empt$, $\abs+$, $\abs-$,
$\bctr$, or\/ $\abs$, the triangulated functor\/
$\sD^\st(X\ctrh^\lct_\cfq)\rarrow\sD^\st(X\lcth^\lct)$ is an equivalence
of categories. \par
\textup{(c)} Let $X$ be a Noetherian scheme of finite Krull dimension.
 Then for any symbol\/ $\bst=\b$, $+$, $-$, $\empt$, $\abs+$, $\abs-$,
$\bctr$, $\ctr$, or\/ $\abs$ the triangulated functor\/
$\sD^\st(X\ctrh_\cfq)\rarrow\sD^\st(X\lcth_\bW)$ induced by
the embedding of exact categories $X\ctrh_\cfq\rarrow X\lcth_\bW$
is an equivalence of triangulated categories. \par
 For any symbol\/ $\bst=\b$, $+$, $-$, $\empt$, $\abs+$, $\abs-$,
$\bctr$, or\/ $\abs$, the triangulated functor\/ $\sD^\st(X\ctrh_\cfq)
\rarrow\sD^\st(X\lcth)$ is an equivalence of categories.
\end{cor}

\begin{proof}
 Follows from Lemma~\ref{finite-krull-ctrh-lcth-finite-dim}
together with Propositions~\ref{finite-resolutions}
and~\ref{becker-contraderived-finite-resolutions}.
 For more general assertions in the cases $\bst=\empt$, $\bco$, $\co$,
$\bctr$, or~$\ctr$ in the context of parts~(a\+b), see
Corollary~\ref{coflasque-resolutions-infinite} below.
\end{proof}

\begin{cor}  \label{finite-krull-ctrh-lcth-derived}
\textup{(a)} Let $X$ be a locally Noetherian scheme of finite
Krull dimension.
 Then for any symbol\/ $\bst=\b$, $+$, $-$, $\empt$, $\abs+$, $\abs-$,
$\bctr$, $\ctr$, or\/ $\abs$, the triangulated functor\/
$\sD^\st(X\ctrh^\lct)\rarrow\sD^\st(X\lcth_\bW^\lct)$ induced by
the embedding of exact categories $X\ctrh^\lct\rarrow
X\lcth_\bW^\lct$ is an equivalence of categories. \par
 For any symbol\/ $\bst=\b$, $+$, $-$, $\empt$, $\abs+$, $\abs-$,
$\bctr$, or\/ $\abs$, the triangulated functor\/
$\sD^\st(X\lcth_\bW^\lct) \rarrow\sD^\st(X\lcth^\lct)$ is
an equivalence of categories. \par
\textup{(b)} Let $X$ be a Noetherian scheme of finite Krull dimension.
 Then for any symbol\/ $\bst=\b$, $+$, $-$, $\empt$, $\abs+$, $\abs-$,
$\bctr$, $\ctr$, or\/ $\abs$, the triangulated functor\/
$\sD^\st(X\ctrh)\rarrow\sD^\st(X\lcth_\bW)$ induced by
the embedding of exact categories $X\ctrh\rarrow X\lcth_\bW$
is an equivalence of triangulated categories. \par
 For any symbol\/ $\bst=\b$, $+$, $-$, $\empt$, $\abs+$, $\abs-$,
$\bctr$, or\/ $\abs$, the triangulated functor\/ $\sD^\st(X\lcth_\bW)
\rarrow\sD^\st(X\lcth)$ is an equivalence of categories.
\end{cor}

\begin{proof}
 Follows from Lemma~\ref{finite-krull-ctrh-lcth-finite-dim}(b\+c)
(see Corollaries~\ref{ctrh-lcth-cor}\+-\ref{lin-ctrh-lcth-cor}
for comparison).
\end{proof}

 The results below in this section purport to replace the above
homological dimension estimates based on the Krull dimension with
the ones based on the numbers of covering open affines.

\begin{lem}  \label{ldelta-rgamma-dimension}
 Let $X$ be a Noetherian scheme and $X=\bigcup_{\alpha=1}^N U_\alpha$
be its finite affine open covering.
 For each subset $1\le \alpha_1<\dotsb<\alpha_k\le N$ of
indices~$\{\alpha\}$, let $U_{\alpha_1}\cap\dotsb\cap U_{\alpha_k}=
\bigcup_{\beta=1}^n V_\beta$, where $n=n_{\alpha_1,\dotsc,\alpha_k}$,
be a finite affine open covering of the intersection.
 Let $M$ denote the supremum of the expressions $k-1+
n_{\alpha_1,\dotsc,\alpha_k}$ taken over all the nonempty subsets
of indices $\alpha_1$,~\dots,~$\alpha_k$.
 Then one has \par
\textup{(a)} $\boR^{\ge M}\Gamma(X,\E)=0$ for any
quasi-coherent sheaf\/ $\E$ on~$X$; \par
\textup{(b)} $\boL_{\ge M}\Delta(X,\gE)=0$ for any locally cotorsion\/
$\bW$\+locally contraherent cosheaf\/ $\gE$ on $X$, provided that
the affine open covering $\{U_\alpha\}$ is subordinate to\/~$\bW$.
 Assuming additionally that the Krull dimension of $X$ is finite,
the same bound holds for any\/ $\bW$\+locally contraherent cosheaf\/
$\gE$ on~$X$.
\end{lem}

\begin{proof}
 We will prove part~(b).
 The first assertion: let $\gF_\bu\rarrow\gE$ be a projective
resolution of an object $\gE\in X\lcth_\bW^\lct$.
 Consider the \v Cech resolution~\eqref{contraherent-cech} 
for each cosheaf $\gF_i\in X\ctrh^\lct_\prj$.
 By Corollaries~\ref{lct-proj-local}(a) and~\ref{proj-lct-direct}(b),
this is a sequence of projective locally cotorsion contraherent
cosheaves.
 Its restriction to each open subset $U_\alpha\sub X$ being naturally
contractible, this finite sequence is exact in $X\ctrh^\lct$,
and consequently also exact (i.~e., even contractible)
in $X\ctrh^\lct_\prj$.
 
 Denote the bicomplex of cosheaves we have obtained (without
the rightmost term that is being resolved) by
$\Cr_\bu(\{U_\alpha\},\gF_\bu)$ and the corresponding bicomplex of
the groups of global cosections by $C_\bu(\{U_\alpha\},\gF_\bu)$.
 Now the total complex of the bicomplex $C_\bu(\{U_\alpha\},\gF_\bu)$
is quasi-isomorphic (in fact, in this case even homotopy equivalent)
to the complex $\Delta(X,\gF_\bu)$ computing the homology groups
$\boL_*\Delta(X,\gE)$.

 On the other hand, for each $1\le k\le N$, the complex
$C_k(\{U_\alpha\},\gF_\bu)$ computes the direct sum of the homology
of the cosheaves $j_{\alpha_1,\dotsc,\alpha_k}^!\gE$ on
$U_{\alpha_1}\cap\dotsb\cap U_{\alpha_k}$ over all
$1\le\alpha_1<\dotsb<\alpha_k\le N$.
 The schemes $U_{\alpha_1}\cap\dotsb\cap U_{\alpha_k}$ being
quasi-compact and separated, and the cosheaves
$j_{\alpha_1,\dotsc,\alpha_k}^!\gE$ being contraherent,
the latter homology can be also computed by the \v Cech complexes
$C_\bu(\{V_\beta\},j_{\alpha_1,\dotsc,\alpha_k}^!\gE)$
(see Section~\ref{homology-subsection}) and consequently vanish
in the homological degrees $\ge n_{\alpha_1,\dotsc,\alpha_k}$.

 To prove the second assertion, one uses a flat (or more generally,
coflasque) resolution $\gF_\bu$ of a contraherent cosheaf $\gE$ and
argues as above using Corollary~\ref{finite-krull-flat-direct}(b)
(or Corollary~\ref{coflasque-direct}(a), respectively).
\end{proof}

 Let us say that a locally cotorsion locally contraherent
cosheaf $\P$ on a locally Noetherian scheme $X$ is \emph{acyclic} if
$\boL_{>0}\Delta(X,\P)=0$.
 Acyclic locally cotorsion $\bW$\+locally contraherent cosheaves
form a full subcategory in $X\lcth_\bW^\lct$ closed under extensions
and kernels of admissible epimorphisms; when $X$ is Noetherian (i.~e.,
quasi-compact), this subcategory is also closed under infinite products.
 Hence it acquires the induced exact category structure, which we
denote by $X\lcth_{\bW,\,\ac}^\lct$.

 Similarly, we say that a locally contraherent cosheaf $\P$ on
a Noetherian scheme $X$ of finite Krull dimension is \emph{acyclic} if
$\boL_{>0}\Delta(X,\P)=0$.
 Acyclic $\bW$\+locally contraherent cosheaves form a full subcategory
in $X\lcth_\bW$ closed under extensions, kernels of admissible
epimorphisms, and infinite products.
 The induced exact category structure on this subcategory is denoted
by $X\lcth_{\bW,\,\ac}$.
 Clearly, one has $X\ctrh_\cfq^\lct\sub X\lcth_{\bW,\,\ac}^\lct$ and
$X\ctrh_\cfq\sub X\lcth_{\bW,\,\ac}$ (under the respective assumptions).
 
 Finally, a quasi-coherent sheaf $\cP$ on a locally Noetherian scheme
$X$ is said to be \emph{acyclic} if $\boR^{>0}\Gamma(X,\cP)=0$.
 Acyclic quasi-coherent sheaves form a full subcategory in
$X\qcoh$ closed under extensions and cokernels of injective morphisms;
when $X$ is Noetherian, this subcategory is also closed
under infinite direct sums.
 The induced exact category structure on it is denoted by
$X\qcoh^\ac$.

\begin{cor} \label{acyclic-finite-dimension}
 Let $X$ be a Noetherian scheme with an open covering\/~$\bW$
and $M$ be the minimal possible value of the nonnegative integer
defined in Lemma~\textup{\ref{ldelta-rgamma-dimension}}
(depending on an affine covering $X=\bigcup_\alpha U_\alpha$
subordinate to\/~$\bW$ and affine coverings of its intersections
$U_{\alpha_1}\cap\dotsb\cap U_{\alpha_k}$).
 Then \par
\textup{(a)} the coresolution dimension of any quasi-coherent
sheaf on $X$ with respect to the coresolving subcategory $X\qcoh^\ac
\sub X\qcoh$ does not exceed~$M-1$; \par
\textup{(b)} the resolution dimension of any locally cotorsion\/
$\bW$\+locally contraherent cosheaf on $X$ with respect to the resolving
subcategory $X\lcth_{\bW,\,\ac}^\lct\sub X\lcth_\bW^\lct$ does not
exceed~$M-1$; \par
\textup{(c)} assuming $X$ has finite Krull dimension, the resolution
dimension of any\/ $\bW$\+locally contraherent cosheaf on $X$ with
respect to the resolving subcategory $X\lcth_{\bW,\,\ac}\sub X\lcth_\bW$
does not exceed~$M-1$.
\end{cor}

\begin{proof}
 Part~(b): since there are enough projectives in $X\lcth_\bW^\lct$
and these belong to $X\lcth_{\bW,\,\ac}^\lct$, the full subcategory
$X\lcth_{\bW,\,\ac}^\lct$ is resolving in $X\lcth_\bW^\lct$.
 Now if $0\rarrow\Q\rarrow\P_{M-2}\rarrow\dotsb\rarrow\P_0\rarrow\gE
\rarrow0$ is an exact sequence in $X\lcth_\bW^\lct$ with
$\P_i\in X\lcth_{\bW,\,\ac}^\lct$, then it is clear from
Lemma~\ref{ldelta-rgamma-dimension}(b) that
$\Q\in X\lcth_{\bW,\,\ac}^\lct$.
 The proofs of parts~(a) and~(c) are similar.
\end{proof}

 Let $f\:Y\rarrow X$ be a quasi-compact morphism of locally
Noetherian schemes.
 By Corollary~\ref{finite-krull-derived-equivalences}(b),
the natural functor $\Hot^-(Y\ctrh^\lct_\prj)\rarrow\sD^-(Y\lcth^\lct)$
is an equivalence of triangulated categories.
 The derived functor $\boL f_!\:\sD^-(Y\lcth^\lct)\rarrow
\sD^-(X\lcth^\lct)$ is constructed by applying the functor
$f_!\:Y\ctrh^\lct_\prj\rarrow X\ctrh^\lct$ from
Corollary~\ref{proj-lct-direct}(a) termwise to bounded above
complexes of projective locally cotorsion contraherent cosheaves.
 By Corollary~\ref{coflasque-direct}(b), one can compute
the derived functor $\boL f_!$ using resolutions of (bounded above)
complexes in $Y\lcth^\lct$ by complexes of coflasque locally cotorsion
contraherent cosheaves.

 Similarly, if the scheme $Y$ is Noetherian of finite Krull dimension,
by Corollary~\ref{finite-krull-derived-equivalences}(c),
the natural functor $\Hot^-(Y\ctrh_\prj)\rarrow\sD^-(Y\lcth)$
is an equivalence of triangulated categories.
 The derived functor $\boL f_!\:\sD^-(Y\lcth)\rarrow\sD^-(X\lcth)$ is
constructed by applying the functor $f_!\:Y\ctrh_\prj\rarrow X\ctrh$
from Corollary~\ref{finite-krull-flat-direct}(a) termwise to bounded
above complexes of projective contraherent cosheaves.
 By Corollary~\ref{coflasque-direct}(a), one can compute the derived
functor $\boL f_!$ using resolutions by coflasque contraherent
cosheaves.

 Let $\bW$ and $\bT$ be open coverings of the schemes $X$ and~$Y$.
 We will call a locally cotorsion $\bT$\+locally contraherent cosheaf
$\Q$ on $Y$ \emph{acyclic with respect to~$f$ over\/~$\bW$} (or
\emph{$f/\bW$\+acyclic}) if the object
$\boL f_!(\Q)\in\sD^-(X\lcth^\lct)$ belongs to the full subcategory
$X\lcth_\bW^\lct\sub X\lcth^\lct\sub\sD^-(X\lcth^\lct)$.
 In other words, the complex $\boL f_!(\Q)$ should have resolution
dimension not exceeding~$0$ with respect to the resolving subcategory
$X\lcth_\bW^\lct\sub X\lcth^\lct$ (in the sense of
Section~\ref{finite-resolutions-subsect}).

 Similarly, if the scheme $Y$ is Noetherian of finite Krull dimension,
we will call a $\bT$\+locally contraherent cosheaf $\Q$ on $Y$
\emph{acyclic with respect to~$f$ over\/~$\bW$} (or
\emph{$f/\bW$\+acyclic}) if the object
$\boL f_!(\Q)\in\sD^-(X\lcth)$ belongs to the full subcategory
$X\lcth_\bW\sub X\lcth\sub\sD^-(X\lcth)$.
 In other words, the complex $\boL f_!(\Q)$ must have resolution
dimension at most~$0$ with respect to the resolving subcategory
$X\lcth_\bW\sub X\lcth$.
 According to Corollary~\ref{fdim-subcategory-cor}, an object of
$Y\lcth_\bT^\lct$ is $f/\bW$\+acyclic if and only if it is
$f/\bW$\+acyclic as an object of $Y\lcth_\bT$.
 Any coflasque (locally cotorsion) contraherent cosheaf on $Y$ is
$f/\bW$\+acyclic.

 It is easy to see, using Lemma~\ref{fdim-triangle}(a\+b), that
the full subcategory of $f/\bW$\+acyclic locally cotorsion
$\bT$\+locally contraherent cosheaves in $Y\lcth_\bT^\lct$ is closed
under extensions, kernels of admissible epimorphisms, and infinite
products.
 The full subcategory of $f/\bW$\+acyclic $\bT$\+locally contraherent
cosheaves in $Y\lcth_\bT$ (defined under the appropriate assumptions
above) has the same properties.
 We denote these subcategories with their induced exact category
structures by $Y\lcth^\lct_{\bT,\,\fWac}$ and
$Y\lcth_{\bT,\,\fWac}$, respectively.

 Finally, the natural functor $\Hot^+(Y\qcoh^\inj)\rarrow
\sD^+(X\qcoh)$ is an equivalence of triangulated categories
by Corollary~\ref{homotopy-derived-coderived-cor}; and the derived
functor $\boR f_*\:\sD^+(Y\qcoh)\rarrow\sD^+(X\qcoh)$ is constructed
by applying the functor~$f_*$ termwise to bounded below complexes
of injective quasi-coherent sheaves.
 One can also compute the derived functor $\boR f_*$ using
coresolutions by flasque quasi-coherent sheaves.

 A quasi-coherent sheaf $\cQ$ on $Y$ is called \emph{acyclic
with respect to~$f$} (or \emph{$f$\+acyclic}) if the object
$\boR f_*(\cQ)\in\sD^+(X\qcoh)$ belongs to the full subcategory
$X\qcoh\sub\sD^+(X\qcoh)$.
 The full subcategory of $f$\+acyclic quasi-coherent sheaves in
$Y\qcoh$ is closed under extensions, cokernels of injective
morphisms, and infinite direct sums.
 We denote this exact subcategory by $Y\qcoh^\fac$.

\begin{lem}
\textup{(a)} Let\/ $\Q$ be an $f/\bW$\+acyclic locally cotorsion\/
$\bT$\+locally contraherent cosheaf on~$Y$.
 Then the cosheaf of\/ $\O_X$\+modules $f_!\Q$ is locally
cotorsion\/ $\bW$\+locally contraherent, and the object represented
by it in the derived category\/ $\sD^-(X\lcth^\lct)$ is naturally
isomorphic to\/~$\boL f_!\Q$. \par
\textup{(b)} Assuming that the scheme $Y$ is Noetherian of finite
Krull dimension, let\/ $\Q$ be an $f/\bW$\+acyclic\/ $\bT$\+locally
contraherent cosheaf on~$Y$.
 Then the cosheaf of\/ $\O_X$\+modules $f_!\Q$ is\/ $\bW$\+locally
contraherent, and the object represented by it in the derived category\/
$\sD^-(X\lcth)$ is naturally isomorphic to\/~$\boL f_!\Q$.
\end{lem}

\begin{proof}
 We will prove part~(a); the proof of part~(b) is similar.
 A complex $\dotsb\rarrow\P_2\rarrow\P_1\rarrow\P_0$ over
$X\lcth^\lct_\bW$ being isomorphic to an object
$\P\in X\lcth^\lct_\bW$ in $\sD^-(X\lcth^\lct)$ means that for each
affine open subscheme $U\sub X$ subordinate to $\bW$ the complex
of cotorsion $\O(U)$\+modules $\dotsb\rarrow\P_2[U]\rarrow
\P_1[U]\rarrow\P_0[U]$ is acyclic except at the rightmost term,
its $\O(U)$\+modules of cocycles are cotorsion, and the cokernel of
the morphism $\P_1\rarrow\P_0$ taken in the category of cosheaves
of $\O_X$\+modules, that is the cosheaf $U\mpsto\coker(\P_1[U]
\to\P_0[U])$, is identified with~$\P$.

 Now let $\dotsb\rarrow\gF_2\rarrow\gF_1\rarrow\gF_0$ be a projective
resolution of the object $\Q$ in the exact category $Y\lcth_\bT^\lct$;
then the cosheaves $f_!\gF_i$ on $X$ belong to
$X\ctrh^\lct\sub X\lcth_\bW^\lct$.
 It remains to notice that the functor~$f_!$ preserves cokernels
taken in the categories of cosheaves of $\O_Y$\+ and $\O_X$\+modules.
\end{proof}

\begin{lem} \label{noetherian-morphism-lct-acyclic-finite-dim}
 Let $f\:Y\rarrow X$ be a morphism of Noetherian schemes
with open coverings\/~$\bT$ and\/~$\bW$.
 Assume that either \par
\textup{(a)} $\bW$ is a finite affine open covering of~$X$, or \par
\textup{(b)} one of the schemes $X$ or $Y$ has finite
Krull dimension. \par
 Then any locally cotorsion\/ $\bT$\+locally contraherent cosheaf
on $Y$ has finite resolution dimension with respect to the resolving
subcategory $Y\lcth^\lct_{\bT,\,\fWac}\sub Y\lcth_\bT^\lct$.
\end{lem}

\begin{proof}
 Since $Y\ctrh^\lct_\prj\sub Y\lcth^\lct_{\bT,\,\fWac}$ and by
Lemma~\ref{fdim-triangle}(a\+b), the full subcategory
$Y\lcth^\lct_{\bT,\,\fWac}$ is resolving in $Y\lcth_\bT^\lct$ for any
quasi-compact morphism of locally Noetherian schemes $f\:Y\rarrow X$.
 In view of Corollaries~\ref{proj-lct-direct}(a) and~\ref{fdim-cor},
the resolution dimension of a cosheaf $\gE\in Y\lcth_\bT^\lct$
with respect to $Y\lcth^\lct_{\bT,\,\fWac}$ does not exceed (in
fact, is equal to) the resolution dimension of the complex
$\boL f_!(\gE)\in\sD^-(X\lcth^\lct)$ with respect to $X\lcth_\bW^\lct$.

 To prove part~(a), denote by $M(Z,\bt)$ the minimal value of
the nonnegative integer $M$ defined in
Lemma~\ref{ldelta-rgamma-dimension} and
Corollary~\ref{acyclic-finite-dimension} for a given Noetherian
scheme~$Z$ with an open covering~$\bt$.
 Let $M$ be the maximal value of $M(f^{-1}(W),\bT|_{f^{-1}(W)})$
taken over all affine open subschemes $W\in\bW$.
 We will show that the resolution dimension of any object of
$Y\lcth_\bT^\lct$ with respect to $Y\lcth^\lct_{\bT,\,\fWac}$
does not exceed $M-1$. 

 Set $Z=f^{-1}(W)$ and $\bt=\bT|_Z$.
 It suffices to check that for any Noetherian affine scheme $W$,
a Noetherian scheme~$Z$ with an open covering~$\bt$, a morphism
of schemes $g\:Z\rarrow W$, and a cosheaf $\gE\in Z\lcth_\bt^\lct$,
the complex of cotorsion $\O(W)$\+modules $(\boL g_!\gE)[W]$
is isomorphic to a complex of cotorsion $\O(W)$\+modules
concentrated in the homological degrees~$\le M-1$ in
the derived category $\sD^-(\O(W)\modl^\cot)$.
 The argument below follows the proof of
Lemma~\ref{ldelta-rgamma-dimension}.

 Let $\gF_\bu\rarrow\gE$ be a projective resolution of
the locally cotorsion locally contraherent cosheaf $\gE$ and
$U_\alpha$ be a finite affine open covering of the scheme $Z$
subordinate to~$\bt$.
 Then the total complex of the bicomplex $C_\bu(\{U_\alpha\},\gF_\bu)$
is homotopy equivalent to the complex $(g_!\gF_\bu)[W]$ computing
the object $(\boL g_!\gE)[W]\in\sD^-(\O(W)\modl^\cot)$.
 On the other hand, for each $1\le k\le N$ the complex
$C_k(\{U_\alpha\},\gF_\bu)$ is the direct sum of the complexes
$(g_{\alpha_1,\dotsc,\alpha_k!}j^!_{\alpha_1,\dotsc,\alpha_k}
\gF_\bu)[W]$, where $g_{\alpha_1,\dotsc,\alpha_k}$ denotes
the composition $g\circ\nobreak j_{\alpha_1,\dotsc,\alpha_k}\:
U_{\alpha_1}\cap\dotsb\cap U_{\alpha_k}\rarrow W$.

 Let $U_{\alpha_1}\cap\dotsb\cap U_{\alpha_k}=\bigcup_\beta V_\beta$
be a finite affine open covering of the intersection
$U_{\alpha_1}\cap\dotsb\cap U_{\alpha_k}$.
 Then the complex $(g_{\alpha_1,\dotsc,\alpha_k!}
j^!_{\alpha_1,\dotsc,\alpha_k}\gF_\bu)[W]$ is homotopy equivalent
to the total complex of the \v Cech bicomplex $C_\bu(\{V_\beta\},
j_{\alpha_1,\dotsc,\alpha_k}^!\gF_\bu)$.
 For each $1\le l\le n_{\alpha_1,\dotsc,\alpha_k}$, the complex
$C_l(\{V_\beta\},j^!_{\alpha_1,\dotsc,\alpha_k}\gF_\bu)$ is
the direct sum of the complexes $(h_{\beta_1,\dotsc,\beta_l!}
e^!_{\beta_1,\dotsc,\beta_l}\gF_\bu)[W]$, where $e_{\beta_1,\dotsc,
\beta_l}$ are the embeddings $V_{\beta_1}\cap\dotsb\cap V_{\beta_l}
\rarrow Z$ and $h_{\beta_1,\dotsc,\beta_l}$ are the compositions
$g\circ e_{\beta_1,\dotsc,\beta_l}\: V_{\beta_1}\cap\dotsb\cap
V_{\beta_l}\rarrow W$.

 Finally, the schemes $V_{\beta_1}\cap\dotsb\cap V_{\beta_l}$ being
affine and $e_{\beta_1,\dotsc,\beta_l}^!\gF_\bu$ being a projective
resolution of a locally cotorsion contraherent cosheaf
$e_{\beta_1,\dotsc,\beta_l}^!\gE$ on $V_{\beta_1}\cap\dotsb
\cap V_{\beta_l}$, the complex of cotorsion $\O(W)$\+modules
$(h_{\beta_1,\dotsc,\beta_l!}e^!_{\beta_1,\dotsc,\beta_l}\gF_\bu)[W]$
is isomorphic to the cotorsion $\O(W)$\+module
$\gE[V_{\beta_1}\cap\dotsb\cap V_{\beta_l}] \simeq
(h_{\beta_1,\dotsc,\beta_l!}e^!_{\beta_1,\dotsc,\beta_l}\gE)[W]$
in $\sD^-(\O(W)\modl^\cot)$.

 Part~(a) is proved.
 Similarly one can show that the resolution dimension of any object
of $Y\lcth_\bT$ with respect to $Y\lcth_{\bT,\,\fWac}$ does not exceed
$M-1$ (assuming that the Krull dimension of $Y$ is finite).

  Now if the scheme $X$ has finite Krull dimension $D$, pick a finite
affine open covering $X=\bigcup_{\alpha=1}^N U_\alpha$ subordinate
to~$\bW$.
 Then, by Lemma~\ref{finite-krull-ctrh-lcth-finite-dim}(b),
the resolution dimension of any object of $X\lcth_{\{U_\alpha\}}^\lct$
with respect to $X\lcth_\bW^\lct\subset X\lcth_{\{U_\alpha\}}^\lct$
does not exceed~$D$.
 It follows that the resolution dimension of any complex from
$\sD^-(X\lcth^\lct)$ with respect to $X\lcth_\bW^\lct$ does not exceed
its resolution dimension with respect to $X\lcth_{\{U_\alpha\}}^\lct$
plus~$D$.
 Hence the resolution dimension of any locally cotorsion locally
contraherent cosheaf from $Y\lcth_\bT^\lct$ with respect to
$Y\lcth^\lct_{\bT,\,\fWac}$ does not exceed its resolution dimension
with respect to $Y\lcth^\lct_{\bT,\,f/\{U_\alpha\}\dashac}$ plus~$D$.
 Here the former summand is finite by part~(a).

 If the scheme $Y$ has finite Krull dimension $D$, then
the resolution dimension in question does not exceed~$D$ by
Lemma~\ref{finite-krull-ctrh-lcth-finite-dim}(b), since
$Y\ctrh^\lct_\cfq\sub Y\lcth^\lct_{\bT,\,\fWac}$.
\end{proof}

\begin{lem} \label{noetherian-morphism-qcoh-acyclic-finite-dim}
 Let $f\:Y\rarrow X$ be a morphism of Noetherian schemes.
 Then any quasi-coherent sheaf on $Y$ has finite coresolution
dimension with respect to the exact subcategory
$Y\qcoh^\fac\sub Y\qcoh$.
\end{lem}

\begin{proof}
 Similar to (and simpler than)
Lemma~\ref{noetherian-morphism-lct-acyclic-finite-dim}.
\end{proof}

\begin{cor}  \label{f-acyclic-resolutions} 
 Let $f\:Y\rarrow X$ be a morphism of Noetherian schemes.  Then \par
\textup{(a)} for any symbol\/ $\bst=\b$, $+$, $-$, $\empt$, $\abs+$,
$\abs-$, $\bco$, $\co$, or\/ $\abs$, the triangulated functor\/
$\sD^\st(Y\qcoh^\fac)\rarrow\sD^\st(Y\qcoh)$ induced by the embedding
of exact categories $Y\qcoh^\fac\rarrow Y\qcoh$ is an equivalence
of categories; \par
\textup{(b)} for any symbol\/ $\bst=-$, $\bctr$, or\/ $\ctr$,
the triangulated functor\/ $\sD^\st(Y\lcth^\lct_{\bT,\,\fWac})
\allowbreak\rarrow\sD^\st(Y\lcth_\bT^\lct)$ induced by the embedding
of exact categories $Y\lcth^\lct_{\bT,\,\fWac}\rarrow Y\lcth_\bT^\lct$
is an equivalence of categories; \par
\textup{(c)} assuming that one of the conditions of
Lemma~\textup{\ref{noetherian-morphism-lct-acyclic-finite-dim}} holds,
for any symbol\/ $\bst=\b$, $+$, $-$, $\empt$, $\abs+$, $\abs-$,
$\bctr$, $\ctr$, or\/ $\abs$ the triangulated functor\/
$\sD^\st(Y\lcth^\lct_{\bT,\,\fWac})\rarrow \sD^\st(Y\lcth_\bT^\lct)$
induced by the embedding of exact categories $Y\lcth^\lct_{\bT,\,\fWac}
\rarrow Y\lcth_\bT^\lct$ is an equivalence of categories.
\end{cor}

\begin{proof}
 Part~(a) follows from
Lemma~\ref{noetherian-morphism-qcoh-acyclic-finite-dim} together
with the dual versions of Propositions~\ref{finite-resolutions}
and~\ref{becker-contraderived-finite-resolutions}.
 Part~(b) is provided by Propositions~\ref{infinite-resolutions}
and~\ref{becker-contraderived-infinite-resolutions}.
 Part~(c) follows from
Lemma~\ref{noetherian-morphism-lct-acyclic-finite-dim} together
with Propositions~\ref{finite-resolutions}
and~\ref{becker-contraderived-finite-resolutions}.
\end{proof}

\subsection{Background equivalences of triangulated categories}
\label{more-background}
 The results of this section complement those of
Sections~\ref{derived-of-sheaves-and-cosheaves-subsect}\+-%
\ref{naive-co-contra-subsect}, \ref{Becker-contraderived-subsect},
and~\ref{finite-flat-lin-dim-subsect}\+-%
\ref{finite-inj-proj-dim-subsect}.

 Let $X$ be a semi-separated Noetherian scheme.
 The \emph{cotorsion dimension} of a quasi-coherent sheaf on $X$ is
defined as its coresolution dimension with respect to the coresolving
subcategory $X\qcoh^\cot\sub X\qcoh$, i.~e., the minimal length of
a coresolution by cotorsion quasi-coherent sheaves.
 For the definition of the \emph{very flat dimension} of
a quasi-coherent sheaf on $X$, we refer to
Section~\ref{finite-flat-lin-dim-subsect}.

 Recall also the notation $X\qcoh_\fl^\cot=X\qcoh_\fl\cap
X\qcoh^\cot$ from Section~\ref{quasi-compact-quasi-coherent}.

\begin{lem}  \label{flat-cotors-homol-dim}
 Let $X=\bigcup_{\alpha=1}^N U_\alpha$ be a finite affine open covering,
and let $D$ denote the Krull dimension of the scheme~$X$.
 Then \par
\textup{(a)} the very flat dimension of any flat quasi-coherent sheaf
on $X$ does not exceed~$D$; \par
\textup{(b)} the homological dimension of the exact category of flat
quasi-coherent sheaves on $X$ does not exceed $N-1+D$; \par
\textup{(c)} the cotorsion dimension of any quasi-coherent sheaf on
$X$ does not exceed $N-1+D$; \par
\textup{(d)} the coresolution dimension of any flat quasi-coherent
sheaf on $X$ with respect to the exact subcategory of flat cotorsion
quasi-coherent sheaves $X\qcoh_\fl^\cot\sub X\qcoh_\fl$ does not
exceed $N-1+D$.
\end{lem}

\begin{proof}
 All the assertions follow more or less directly from
Theorem~\ref{raynaud-gruson-flat-thm}.
 In view of that result from~\cite{RG}, part~(a) is obvious.
 For parts~(b\+d), recall that the pair of full subcategories
($X\qcoh_\fl$, $X\qcoh^\cot$) is a hereditary complete cotorsion
pair in the abelian category $X\qcoh$ (for any quasi-compact
semi-separated scheme $X$) by
Corollaries~\ref{quasi-cotors-characterizations}(b\+c)
and~\ref{quasi-cotors-cor}(a\+b) together with
Lemma~\ref{cotorsion-pair-direct-summands-lemma}.
 Consequently, for a semi-separated Noetherian scheme,
parts~(b\+d) are equivalent restatements of one another, corresponding
to the equivalent conditions~(4), (1), and~(2) of
Lemma~\ref{finite-homol-dim-finite-coresol-dim}.
 To prove these properties, it suffices to check the condition of
Lemma~\ref{finite-homol-dim-finite-coresol-dim}(3), which follows
from Theorem~\ref{raynaud-gruson-flat-thm} in view
of~\cite[Theorem~6.3(b)]{PS6} (cf.\ Lemma~\ref{vfl-cta-finite-dim}).
\end{proof}

 The notation $X\qcoh_\fd$ and $X\qcoh_\ffdd$ for exact categories of
quasi-coherent sheaves of finite flat dimension was introduced in
Section~\ref{finite-flat-lin-dim-subsect}.

\begin{cor}  \label{quasi-finite-flat-dim-all-derived-coincide}
 Let $X$ be a semi-separated Noetherian scheme of finite Krull dimension
and $d\ge0$ be any (finite) integer. 
 Then the natural triangulated functors\/ $\sD^\abs(X\qcoh_\ffdd)\rarrow
\sD^\co(X\qcoh_\ffdd)\rarrow\sD(X\qcoh_\ffdd)$ and\/
$\sD^{\abs\pm}(X\qcoh_\ffdd)\allowbreak\rarrow\sD^\pm(X\qcoh_\ffdd)$ are
equivalences of triangulated categories.
 In particular, such functors between the derived categories of
the exact category $X\qcoh_\fl$ are equivalences of categories.
 Moreover, the natural triangulated functors\/ $\sD^\abs(X\qcoh_\fl)
\rarrow\sD^\co(X\qcoh_\fl)\rarrow\sD^\bco(X\qcoh_\fl)\rarrow
\sD(X\qcoh_\fl)$ (are well-defined and) are equivalences of triangulated
categories.
\end{cor}

\begin{proof}
 Concerning Becker's coderived category, one has
$\Acycl^\bco(X\qcoh_\fl)=\Acycl(X\qcoh_\fl)$, hence
$\sD^\bco(X\qcoh_\fl)=\sD(X\qcoh_\fl)$ for any quasi-compact
semi-separated scheme $X$ by
Theorem~\ref{vlf-cta-fl-cot-derived-vlf-fl-equivalences}(c).
 The remaining assertions follow from
Lemma~\ref{flat-cotors-homol-dim}(b) together with
Lemma~\ref{psemi-remark21}.  \hbadness=1400
\end{proof}

\begin{cor}  \label{derived-fl-vfl-cor}
 Let $X$ be a semi-separated Noetherian scheme of finite Krull
dimension.
 Then for any symbol\/ $\bst=\b$, $+$, $-$, $\empt$, $\abs+$,
$\abs-$, $\bco$, $\co$, or\/~$\abs$, the triangulated functor\/
$\sD^\st(X\qcoh_\vfl)\rarrow\sD^\st(X\qcoh_\fl)$ induced by
the embedding of exact categories $X\qcoh_\vfl\rarrow X\qcoh_\fl$
is an equivalence of triangulated categories.
\end{cor}

\begin{proof}
 Concerning Becker's coderived categories, one has
$\Acycl^\bco(X\qcoh_\vfl)=\Acycl(X\qcoh_\vfl)$
and $\Acycl^\bco(X\qcoh_\fl)=\Acycl(X\qcoh_\fl)$ for any
quasi-compact semi-separated scheme $X$ by
Theorem~\ref{vlf-cta-fl-cot-derived-vlf-fl-equivalences}(a,c).
 It follows that the induced triangulated functor on Becker's
coderived categories is well-defined.
 The assertions that the induced triangulated functors are
triangulated equivalences for $\bst=-$, $\empt$, and~$\bco$ hold
for any quasi-compact semi-separated scheme~$X$; see
Theorem~\ref{derived-vfl-fl-equivalence}.
 The remaining cases follow from Lemma~\ref{flat-cotors-homol-dim}(a)
together with Proposition~\ref{finite-resolutions}.
\end{proof}

\begin{cor}  \label{derived-contra-lct-cor}
 Let $X$ be a Noetherian scheme of finite Krull dimension.  Then \par
\textup{(a)} for any symbol\/ $\bst=\b$, $+$, $-$, $\empt$, $\abs+$,
$\abs-$, $\bctr$, $\ctr$, or\/~$\abs$, the triangulated functor\/
$\sD^\st(X\ctrh^\lct_\cfq)\rarrow \sD^\st(X\ctrh_\cfq)$ induced by
the embedding of exact categories $X\ctrh^\lct_\cfq\rarrow X\ctrh_\cfq$
is an equivalence of triangulated categories; \par
\textup{(b)} for any symbol\/ $\bst=\b$, $+$, $-$, $\empt$, $\abs+$,
$\abs-$, $\bctr$, $\ctr$, or\/~$\abs$, the triangulated functor\/
$\sD^\st(X\lcth_\bW^\lct)\rarrow \sD^\st(X\lcth_\bW)$ induced by
the embedding of exact categories $X\lcth_\bW^\lct\rarrow X\lcth_\bW$
is an equivalence of triangulated categories; \par
\textup{(c)} for any symbol\/ $\bst=\b$, $+$, $-$, $\empt$, $\abs+$,
$\abs-$, $\bctr$, or\/~$\abs$, the triangulated functor\/
$\sD^\st(X\lcth^\lct)\rarrow\sD^\st(X\lcth)$ induced by
the embedding of exact categories $X\lcth^\lct\rarrow X\lcth$
is an equivalence of triangulated categories.
\end{cor}

\begin{proof}
 Part~(a) for all the symbols except $\bctr$ is provided by
Lemma~\ref{coflasque-lct-envelope}(b) together
with the dual version of Proposition~\ref{finite-resolutions}.
 In the case $\bst=\bctr$, one can apply
Proposition~\ref{becker-contraderived-finite-coresolutions},
whose assumptions are satisfied by
Corollaries~\ref{finite-krull-contrah-projective},
\ref{coflasque-flat-loc-cotorsion-cotorsion-pair},
and Lemma~\ref{coflasque-lct-envelope}(b).

 Parts~(b\+c) follow from part~(a) and
Corollary~\ref{coflasque-resolutions-finite}(b\+c).
 Let us just explain why the induced triangulated functor in
the case $\bst=\bctr$ is well-defined.
 Let $\gB^\bu$ be a Becker-contraacyclic complex in
the exact category $X\lcth^\lct$; we need to show that $\gB^\bu$
is Becker-contraacyclic in $X\lcth$.
 Using Lemma~\ref{finite-krull-ctrh-lcth-finite-dim}(b) and
the construction from the beginning of the proof of
Proposition~\ref{finite-resolutions}, we produce a finite acyclic
complex $0\rarrow\gF_D^\bu\rarrow\gF_{D-1}^\bu\rarrow\dotsb
\rarrow\gF_0^\bu\rarrow\gB^\bu\rarrow0$ of complexes in
$X\lcth^\lct$ such that $\gF_i^\bu$ are complexes of coflasque
locally cotorsion contraherent cosheaves for all $0\le i\le D$.
 Denote by $\gF^\bu$ the total complex of the bicomplex
$\gF_D^\bu\rarrow\gF_{D-1}^\bu\rarrow\dotsb\rarrow\gF_0^\bu$.
 Then the cone of the natural morphism of complexes $\gF^\bu\rarrow
\gB^\bu$ is absolutely acyclic in $X\lcth^\lct$, hence also absolutely
acyclic in $X\lcth$.
 By Lemma~\ref{Positselski-trivial-are-Becker-trivial}(a),
it follows that this cone is Becker-contraacyclic in $X\lcth^\lct$
and in $X\lcth$.
 Since $\gB^\bu$ is Becker-contraacyclic in $X\lcth^\lct$, we can
conclude that $\gF^\bu$ is Becker-contraacyclic in $X\lcth^\lct$.
 As $\gF^\bu$ is a complex in $X\ctrh^\lct_\cfq$ and the classes
of projective objects in $X\lcth^\lct$ and $X\ctrh^\lct_\cfq$
coincide by Corollary~\ref{lct-proj-local}(b), this means that
the complex $\gF^\bu$ is Becker-contraacyclic in $X\ctrh^\lct_\cfq$.
 By part~(a), it follows that $\gF^\bu$ is Becker-contraacyclic
in $X\ctrh_\cfq$.
 Since the classes of projective objects in $X\lcth$ and $X\ctrh_\cfq$
coincide by Corollary~\ref{finite-krull-contrah-projective}, this
means that the complex $\gF^\bu$ is Becker-contraacyclic in
$X\lcth$.
 Finally, as we know that the cone of the morphism $\gF^\bu\rarrow
\gB^\bu$ is Becker-contraacyclic in $X\lcth$, it follows that
the complex $\gB^\bu$ is Becker-contraacyclic in $X\lcth$, as desired.

 Alternatively, in part~(b) one can reduce the case $\bst=\bctr$ to
the case $\bst=\ctr$ by referring to
Theorem~\ref{derived-inj-proj-resolutions}(b,d) below.
 In the case of a semi-separated Noetherian scheme $X$
of finite Krull dimension, the assertions~(b\+c) can be obtained
directly from Lemma~\ref{loc-cotors-dim}(a) using
Propositions~\ref{finite-resolutions}
and~\ref{becker-contraderived-finite-coresolutions}.
\end{proof}

 The following corollary is another restricted version of
Theorem~\ref{naive-co-contra-thm}; it is to be compared with
Corollary~\ref{vfl-co-contra-cor-expanded}.

\begin{cor}  \label{flat-noetherian-co-contra-cor}
 Let $X$ be a semi-separated Noetherian scheme of finite Krull
dimension.
 Then for any symbol\/ $\bst=\b$, $+$, $-$, $\empt$, $\abs+$,
$\abs-$, $\bco$, $\co$, or\/~$\abs$ there is a natural equivalence of
triangulated categories\/ $\sD^\st(X\qcoh_\fl)\simeq
\Hot^\st(X\ctrh^\lct_\prj)$.
\end{cor}

\begin{proof}
 The cases $\bst=+$, $\empt$, or~$\bco$ are covered by
Corollary~\ref{vfl-co-contra-cor-expanded}(c) (which holds for any
quasi-compact semi-separated scheme~$X$).
 Assuming $\bst\ne\bco$, $\co$, by Lemma~\ref{flat-cotors-homol-dim}(d)
together with the dual version of Proposition~\ref{finite-resolutions}
the triangulated functor $\Hot^\st(X\qcoh_\fl^\cot)\rarrow
\sD^\st(X\qcoh_\fl)$ is an equivalence of categories.
 In view of Corollary~\ref{quasi-finite-flat-dim-all-derived-coincide},
the same assertion holds for $\bst=\bco$ or~$\co$.
 Hence it remains to recall that the equivalence of categories from
Lemma~\ref{cta-clp-equivalence} identifies $X\qcoh_\fl^\cot$ with
$X\ctrh^\lct_\prj$ (see Lemma~\ref{cta-clp-restricts-to-prj-clf}(b)).
\end{proof}

\begin{cor}  \label{fl-lct-prj-direct-images-identified}
 Let $f\:Y\rarrow X$ be a morphism of finite flat dimension between
semi-separated Noetherian schemes of finite Krull dimension.
 Then for any symbol\/ $\bst=\b$, $+$, $-$, $\empt$, $\abs+$, $\abs-$,
or\/~$\abs$ the equivalences of triangulated categories\/
$\sD^\st(Y\qcoh_\fl)\simeq\Hot^\st(Y\ctrh^\lct_\prj)$
and\/ $\sD^\st(X\qcoh_\fl)\simeq \Hot^\st(X\ctrh^\lct_\prj)$ from
Corollary~\textup{\ref{flat-noetherian-co-contra-cor}} transform
the right derived functor\/ $\boR f_*$~\textup{\eqref{qcoh-direct-ffd}}
into the left derived functor\/
$\boL f_!$~\textup{\eqref{ctrh-direct-lct-fpd}}.
 In other words, the following diagram of triangulated functors
and triangulated equivalences is commutative:
\begin{equation} \label{fl-lct-prj-direct-image-diagram}
\begin{gathered}
 \xymatrix{
  \sD^\st(Y\qcoh_\fl) \ar@<2pt>[r] \ar@<-2pt>@{-}[r]
  & \sD^\st(Y\qcoh_\fd) \ar@{=}[r] \ar[d]_{\boR f_*}
  & \sD^\st(Y\ctrh^\lct_\fpd) \ar[d]^{\boL f_!}
  & \Hot^\st(Y\ctrh^\lct_\prj) \ar@<-2pt>[l] \ar@<2pt>@{-}[l] \\
  \sD^\st(Y\qcoh_\fl) \ar@<-2pt>[r] \ar@<2pt>@{-}[r]
  & \sD^\st(X\qcoh_\fd) \ar@{=}[r] & \sD^\st(X\ctrh^\lct_\fpd)
  & \Hot^\st(X\ctrh^\lct_\prj) \ar@<2pt>[l] \ar@<-2pt>@{-}[l]  
 }
\end{gathered}
\end{equation}
\end{cor}

\begin{proof}
 The cases $\bst=+$ or~$\empt$ are covered by
Corollary~\ref{inj-vfl-direct-images-identified}(d) (which holds
for any quasi-compact semi-separated schem~$X$).
 All the cases can be either deduced from
Corollary~\ref{inj-vfl-direct-images-identified}(c), or proved
directly in the similar way using
Lemma~\ref{cta-clp-finite-flat-inj-dim-identified}(d).
\end{proof}

 Let $X$ be a locally Noetherian scheme with an open covering~$\bW$.
 As in Section~\ref{finite-inj-proj-dim-subsect}, we denote by
$X\qcoh^\fidd$ the full subcategory of objects of injective
dimension~$\le d$ in $X\qcoh$ and by $X\lcth^\lct_{\bW,\,\fpdd}$
the full subcategory of objects of projective dimension~$\le d$
in $X\lcth_\bW^\lct$.
 For a Noetherian scheme $X$ of finite Krull dimension, let
$X\lcth_{\bW,\,\fpdd}$ denote the full subcategory of objects of
projective dimension~$\le d$ in $X\lcth_\bW$.
 We set $X\ctrh^\lct_\fpdd=X\lcth_{\{X\},\fpdd}^\lct$ and
$X\ctrh_\fpdd=X\lcth_{\{X\},\fpdd}$.
 Clearly, the projective dimension of an object of $X\lcth_\bW^\lct$
or $X\lcth_\bW$ does not change when the open covering $\bW$ is
replaced by its refinement.

 The full subcategory $X\qcoh^\fidd\sub X\qcoh$ is closed under
extensions, cokernels of injective morphisms, and infinite direct sums.
 The full subcategory $X\lcth^\lct_{\bW,\,\fpdd}\sub X\lcth_\bW^\lct$
is closed under extensions, kernels of admissible epimorphisms, and
infinite products (see Corollary~\ref{loc-noetherian-proj-products}).
 The full subcategory $X\lcth_{\bW,\,\fpdd}\sub X\lcth_\bW$ is closed
under extensions and kernels of admissible epimorphisms.

\begin{cor} \label{loc-noetherian-derived-fid-fpd} \hfuzz=4pt
\textup{(a)} Let $X$ be a locally Noetherian scheme.
Then the natural triangulated functors\/
$\Hot(X\qcoh^\inj)\rarrow\sD^\abs(X\qcoh^\fidd)\rarrow
\sD^{\co=\bco}(X\qcoh^\fidd)\rarrow\sD(X\qcoh^\fidd)$,
\ $\Hot^\pm(X\qcoh^\inj)\rarrow\sD^{\abs\pm}(X\qcoh^\fidd)\rarrow
\sD^\pm(X\qcoh^\fidd)$, and\/ $\Hot^\b(X\qcoh^\inj)\rarrow
\sD^\b(X\qcoh^\fidd)$ are equivalences of categories. \par
\textup{(b)} Let $X$ be a locally Noetherian scheme.
Then the natural triangulated functors\/ $\Hot(X\ctrh^\lct_\prj)
\rarrow\sD^\abs(X\lcth^\lct_{\bW,\,\fpdd})\rarrow
\sD^{\ctr=\bctr}(X\lcth^\lct_{\bW,\,\fpdd})\rarrow
\sD(X\lcth^\lct_{\bW,\,\fpdd})$, \ $\Hot^\pm(X\ctrh^\lct_\prj)\rarrow
\sD^{\abs\pm}(X\lcth^\lct_{\bW,\,\fpdd})\rarrow
\sD^\pm(X\lcth^\lct_{\bW,\,\fpdd})$, and\/ $\Hot^\b(X\ctrh^\lct_\prj)
\rarrow\sD^\b(X\lcth^\lct_{\bW,\,\fpdd})$ are equivalences
of categories. \par
\textup{(c)} Let $X$ be a Noetherian scheme of finite Krull dimension.
Then the natural triangulated functors\/ $\Hot(X\ctrh_\prj)
\rarrow\sD^{\abs=\bctr}(X\lcth_{\bW,\,\fpdd})\rarrow
\sD(X\lcth_{\bW,\,\fpdd})$, \ $\Hot^\pm(X\ctrh_\prj)\rarrow
\sD^{\abs\pm}(X\lcth_{\bW,\,\fpdd})\rarrow
\sD^\pm(X\lcth_{\bW,\,\fpdd})$, and\/ $\Hot^\b(X\ctrh_\prj)
\allowbreak\rarrow\sD^\b(X\lcth_{\bW,\,\fpdd})$ are equivalences
of categories. \hfuzz=11pt
\end{cor}

\begin{proof}
 Part~(a) follows from Corollary~\ref{derived-fid-fpd-cor}(a) and
the fact that the functors of infinite direct sum are exact in
$X\qcoh^\fidd$ for a locally Noetherian scheme~$X$.
 Parts~(b\+c) follow from Proposition~\ref{finite-resolutions}
and Lemma~\ref{psemi-remark21}.
 For the contraderived category in part~(b), it is important that
the functors of infinite product are exact in
$X\lcth^\lct_{\bW,\,\fpdd}$.
 In connection with Becker's coderived and contraderived categories,
see Theorem~\ref{finite-homol-dim-becker-co-contra-derived}.
\end{proof}

 A cosheaf of $\O_X$\+modules $\gG$ on a scheme $X$ is said to have
\emph{$\bW$\+flat dimension not exceeding~$d$} if the flat dimension of
the $\O_X(U)$\+module $\gG[U]$ does not exceed~$d$ for any affine open
subscheme $U\sub X$ subordinate to~$\bW$.
 The flat dimension of a cosheaf of $\O_X$\+modules is defined as its
$\{X\}$\+flat dimension.
 $\bW$\+locally contraherent cosheaves of $\bW$\+flat dimension not
exceeding~$d$ on a locally Noetherian scheme $X$ form a full subcategory
$X\lcth_\bW^\ffdd\sub X\lcth_\bW$ closed under extensions, kernels of
admissible epimorphisms, and infinite products.
 We set $X\ctrh^\ffdd=X\lcth_{\{X\}}^\ffdd$.

 The flat dimension of a contraherent cosheaf $\gG$ on a coherent
affine scheme $U$ is equal to the flat dimension of
the $\O_X(U)$\+module~$\gG[U]$ (see Section~\ref{contraherent-tensor}).
 Over a semi-separated coherent scheme $X$, a $\bW$\+locally
contraherent cosheaf has $\bW$\+flat dimension~$\le d$ if and only if it
admits a resolution of length~$\le d$ by $\bW$\+flat $\bW$\+locally
contraherent cosheaves (see Corollary~\ref{proj-flat}(a)).

 Hence it follows from Corollary~\ref{finite-krull-flat-contraherent}(b)
(applied to affine open subschemes $U\sub X$) that the $\bW$\+flat
dimension of a $\bW$\+locally contraherent cosheaf on a locally
Noetherian scheme $X$ of finite Krull dimension does not change when
the covering $\bW$ is replaced by its refinement.
 According to part~(a) of the same Corollary, on a semi-separated
Noetherian scheme of finite Krull dimension the $\bW$\+flat dimension
of a $\bW$\+locally contraherent cosheaf is equal to its
antilocally flat dimension; so $X\lcth_\bW^\ffdd=X\lcth_{\bW,\,\alfdd}$.
 By Corollary~\ref{lct-prj-flat}, the $\bW$\+flat dimension of
a locally cotorsion $\bW$\+locally contraherent cosheaf on a locally
Noetherian scheme $X$ coincides with its projective dimension in
$X\lcth_\bW^\lct$ (and also does not depend on~$\bW$).
 So one has $X\lcth_\bW^\ffdd\cap X\lcth_\bW^\lct = 
X\lcth^\lct_{\bW,\,\fpdd}$.

\begin{lem}  \label{flat-contra-finite-proj-dimension}
 Let $X$ be a Noetherian scheme of finite Krull dimension~$D$.
 Then a\/ $\bW$\+locally contraherent cosheaf on $X$ has finite
projective dimension in the exact category $X\lcth_\bW$ if and only if
it has finite\/ $\bW$\+flat dimension.
 More precisely, the inclusions of full subcategories
$X\lcth_{\bW,\,\fpdd} \sub X\lcth_\bW^\ffdd\sub
X\lcth_{\bW,\,\fpdD}$ hold in the category $X\lcth_\bW$.
\end{lem}

\begin{proof}
 The inclusion $X\lcth_{\bW,\,\fpdd}\sub X\lcth_\bW^\ffdd$ holds due
to Corollary~\ref{finite-krull-contrah-projective}.
 Conversely, by the same Corollary any $\bW$\+locally contraherent
cosheaf $\gM$ on $X$ has a resolution by flat contraherent cosheaves,
so the $\bW$\+flat dimension of $\gM$ is equal to its
resolution dimension with respect to $X\ctrh^\fl\sub X\lcth_\bW$
(see Corollary~\ref{finite-krull-flat-contraherent}(b)).
 It remains to apply the last assertion of
Corollary~\ref{lct-prj-envelope}(b).
\end{proof}

\begin{cor} \label{derived-contra-ffd-cor}
 For any Noetherian scheme $X$ of finite Krull dimension
and any (finite) integer $d\ge0$, the natural triangulated functors\/
$\Hot(X\ctrh_\prj)\rarrow \sD^\abs(X\lcth_\bW^\ffdd)\rarrow
\sD^{\ctr=\bctr}(X\lcth_\bW^\ffdd)\rarrow\sD(X\lcth_\bW^\ffdd)$, \
$\Hot^\pm(X\ctrh_\prj)\rarrow\sD^{\abs\pm}(X\lcth_\bW^\ffdd)\rarrow
\sD^\pm(X\lcth_\bW^\ffdd)$, and\/ $\Hot^\b(X\ctrh_\prj)\rarrow
\sD^\b(X\lcth_\bW^\ffdd)$ are equivalences of triangulated categories.
\end{cor}

\begin{proof}
 It is clear from Lemma~\ref{flat-contra-finite-proj-dimension} that
the homological dimension of the exact category $X\lcth_\bW^\ffdd$ is
finite, so it remains to apply Lemma~\ref{psemi-remark21} (to obtain
the equivalences between various derived categories of this exact
category) and Proposition~\ref{finite-resolutions} (to identify
the absolute derived categories with the homotopy categories of
projective objects).
 Alternatively, one can use~\cite[Theorem~3.6 and Remark~3.6]{Pkoszul}.
 Concerning the contraderived category, it is important that infinite
products are exact in $X\lcth_\bW^\ffdd$.
 For the Becker contraderived category, see
see Theorem~\ref{finite-homol-dim-becker-co-contra-derived}.
\end{proof}

 The following theorem is the main result of this section.
 It should be compared with
Theorem~\ref{quasi-coherent-becker-coderived}
and Corollary~\ref{becker-contraderived-of-lcta-lct-well-behaved}.

\begin{thm} \label{derived-inj-proj-resolutions}
\textup{(a)} For any locally Noetherian scheme $X$, the classes of
Positselski-coacyclic and Becker-coacyclic complexes in the abelian
category $X\qcoh$ coincide.
 The natural functor\/ $\Hot(X\qcoh^\inj)\rarrow\sD^{\co=\bco}(X\qcoh)$
is an equivalence of triangulated categories. {\hbadness=1200\par}
\textup{(b)} For any locally Noetherian scheme $X$, the classes of
Positselski-contraacyclic and Becker-contraacyclic complexes in
the exact category $X\lcth_\bW^\lct$ coincide.
 The natural functor\/ $\Hot(X\ctrh_\prj^\lct)\rarrow\sD^{\ctr=\bctr}
(X\lcth_\bW^\lct)$ is an equivalence of triangulated categories. \par
\textup{(c)} For any semi-separated coherent scheme $X$ and any
symbol\/ $\bst=-$, $\ctr$, or\/~$\bctr$, the natural functor\/
$\sD^\st(X\ctrh^\fl)\rarrow\sD^\st(X\lcth_\bW)$ is an equivalence of
triangulated categories. \par
\textup{(d)} For any Noetherian scheme $X$ of finite Krull dimension,
the classes of Positselski-contraacyclic and Becker-contraacyclic
complexes in the exact category $X\lcth_\bW$ coincide.
 The natural functors\/ $\Hot(X\ctrh_\prj)\rarrow\sD^\abs(X\ctrh^\fl)
\rarrow\sD^{\ctr=\bctr}(X\lcth_\bW)$ are equivalences of triangulated
categories.
\end{thm}

\begin{proof}
 Part~(a) is a standard result (see, e.~g., \cite[Lemma~1.7(b)]{EP})
which is a particular case of
Proposition~\ref{positselski-becker-co-contra-derived}(a)
and can be also obtained from the dual version of
Proposition~\ref{infinite-resolutions}(b).
 The key observation is that there are enough injectives in $X\qcoh$
and the full subcategory $X\qcoh^\inj$ they form is closed under
infinite direct sums.
 Similarly, part~(b) is a particular case of
Proposition~\ref{positselski-becker-co-contra-derived}(b)
and can be also obtained from Proposition~\ref{infinite-resolutions}(b).
 In any case, the argument is based on
Theorem~\ref{proj-lct-classification}(a)
and Corollary~\ref{loc-noetherian-proj-products}.
 Part~(c) holds by Propositions~\ref{infinite-resolutions}
and~\ref{becker-contraderived-infinite-resolutions} together with
Lemma~\ref{clf-cover} or Lemmas~\ref{loc-contra-proj}(a)
and~\ref{proj-flat}(a).

 Finally, in part~(d) the functors $\Hot(X\ctrh_\prj)\rarrow
\sD^\abs(X\ctrh^\fl)\rarrow\sD^{\ctr=\bctr}\allowbreak(X\ctrh^\fl)$ are
equivalences of categories by Corollary~\ref{derived-contra-ffd-cor},
while the functors $\sD^\ctr(X\ctrh^\fl)\rarrow\sD^\ctr(X\lcth_\bW)$
and $\sD^\bctr(X\ctrh^\fl)\rarrow\sD^\bctr(X\lcth_\bW)$
are equivalences of categories by
Propositions~\ref{infinite-resolutions}(b)
and~\ref{becker-contraderived-infinite-resolutions} together with
Corollary~\ref{finite-krull-contrah-projective}.
 A direct proof of the equivalences $\Hot(X\ctrh_\prj)\rarrow
\sD^\ctr(X\lcth_\bW)\rarrow\sD^\bctr(X\lcth_\bW)$ is also possible;
it proceeds along the following lines.

 One has to use the fully general form of
Proposition~\ref{positselski-becker-co-contra-derived}(b) with
the assumption of finite projective dimension of countable products
of projective objects.
 Alternatively, one can apply the even more general
Corollary~\ref{finite-homol-dim-equivalence-cor}.
 Let $X=\bigcup_\alpha U_\alpha$ be a finite affine open covering;
then it follows from Corollary~\ref{finite-krull-contrah-projective}(b)
that an infinite product of projective contraherent cosheaves on $X$
is a direct summand of a direct sum over~$\alpha$ of the direct
images of contraherent cosheaves on $U_\alpha$ correspoding to infinite
products of very flat contraadjusted $\O(U_\alpha)$\+modules.

 Infinite products of such modules may not be very flat, but they
are certainly flat and contraadjusted. 
 By the last assertion of Corollary~\ref{lct-prj-envelope}(b),
one can conclude that the projective dimensions of infinite products
of projective objects in $X\ctrh$ do not exceed the Krull dimension
$D$ of the scheme~$X$.
 So Proposition~\ref{positselski-becker-co-contra-derived}(b) is
applicable to the exact category $X\lcth_\bW$, or in other words,
the assumption of Corollary~\ref{finite-homol-dim-equivalence-cor} is
satisfied by the pair of exact categories $X\ctrh_\prj\sub X\lcth_\bW$.
\end{proof}

 The following corollary is to be compared with
Corollaries~\ref{finite-krull-derived-equivalences}(b)
and~\ref{finite-krull-ctrh-lcth-derived}(b).

\begin{cor}  \label{loc-noetherian-ctrh-lcth-contraderived}
 For any locally Noetherian scheme~$X$ with an open covering\/ $\bW$,
the triangulated functor\/ $\sD^{\ctr=\bctr}(X\ctrh^\lct)\rarrow
\sD^{\ctr=\bctr}(X\lcth_\bW^\lct)$ induced by the embedding of exact
categories $X\ctrh^\lct\rarrow X\lcth_\bW^\lct$ is an equivalence of
triangulated categories.
\end{cor}

\begin{proof}
 Follows from Theorem~\ref{derived-inj-proj-resolutions}(b) applied
to the coverings $\{X\}$ and $\bW$ of the scheme~$X$.
 Alternatively, for $\bst=\ctr$ or $\bst=\bctr$ separately one can
apply directly Proposition~\ref{infinite-resolutions}(b)
or~\ref{becker-contraderived-infinite-resolutions} together with
Theorem~\ref{proj-lct-classification}(a).
\end{proof}

 The next corollary is to be compared with
Corollary~\ref{qcoh-becker-coacyclicity-is-local},
Theorem~\ref{Becker-contraacyclicity-local-on-qcomp-qsep},
and Remark~\ref{locality-of-Becker-contraacyclicity-remark}.

\begin{cor} \label{locality-of-co-contra-acyclicity-on-loc-Notherian}
 Let $X$ be a locally Noetherian scheme, $\bW$ be its open covering,
and $X=\bigcup_\alpha U_\alpha$ be an open covering of $X$ subordinate
to~$\bW$.
 Denote by $j_\alpha\:U_\alpha\rarrow X$ the open embedding morphisms.
 In this setting: \par
\textup{(a)} A complex of quasi-coherent sheaves\/ $\cA^\bu$ on $X$
is coacyclic in $X\qcoh$ if and only if the complexes
$j_\alpha^*\cA^\bu$ are coacyclic in $U_\alpha\qcoh$ for all~$\alpha$.
\par
\textup{(b)} A complex of locally cotorsion\/ $\bW$\+locally
contraherent cosheaves\/ $\gB^\bu$ on $X$ is contraacyclic in
$X\lcth_\bW^\lct$ if and only if the complexes $j_\alpha^!\gB^\bu$ are
contraacyclic in $U_\alpha\ctrh^\lct$ for all~$\alpha$. \par
\textup{(c)} Assume that the scheme $X$ is Noetherian of finite
Krull dimension.
 Then a complex of\/ $\bW$\+locally contraherent cosheaves\/ $\gB^\bu$
on $X$ is contraacyclic in $X\lcth_\bW$ if and only if the complexes
$j_\alpha^!\gB^\bu$ are contraacyclic in $U_\alpha\ctrh$
for all~$\alpha$.
\end{cor}

\begin{proof}
 Notice first of all that, in the context of part~(a), the Becker
and Positselski coacyclicity properties are equivalent by
Theorem~\ref{derived-inj-proj-resolutions}(a).
 Similarly, in the contexts of parts~(b) and~(c), the Becker and
Positselski contraacyclicity properties are equivalent by
Theorem~\ref{derived-inj-proj-resolutions}(b,d).

 The argument from~\cite[Lemma~A.9]{Psemten} proves part~(a).
 One needs to notice that the functors~$j_\alpha^*$ take
Positselski-coacyclic complexes to Positselski-coacyclic complexes
(since these functors are exact and preserve infinite direct sums),
and use Theorem~\ref{derived-inj-proj-resolutions}(a) together with
the fact that injectivity of quasi-coherent sheaves is a local
property on locally Noetherian schemes.

 The proof of part~(b) is dual-analogous to~(a).
 The functors~$j_\alpha^!$ take Positselski-contraacyclic complexes
to Positselski-contraacyclic complexes, since these functors are
exact and preserve infinite products.
 Conversely, given a complex $\gB^\bu$ in $X\lcth_\bW^\lct$, one can
find a complex $\P^\bu$ in $X\ctrh^\lct_\prj$ together with a morphism
of complexes $\P^\bu\rarrow\gB^\bu$ with a cone contraacyclic
in $X\lcth_\bW^\lct$ (by Theorem~\ref{derived-inj-proj-resolutions}(b)).
 Then contraacyclicity of the complexes $j_\alpha^!\gB^\bu$ implies
contraacyclicity of the complexes $j_\alpha^!\P^\bu$ in
$U_\alpha\ctrh^\lct$.
 By Corollary~\ref{lct-proj-local}(a), the complexes $j_\alpha^!\P^\bu$
are complexes of projective objects in $U_\alpha\ctrh^\lct$; so it
follows that these complexes are contractible.
 In other words, this means that the complexes $j_\alpha^!\P^\bu$
are acyclic in $U_\alpha\ctrh^\lct$ with the objects of cocycles
belonging $U_\alpha\ctrh^\lct_\prj$.
 By Lemma~\ref{acyclicity-in-lcth-criterion}(b), it follows that
the complex $\P^\bu$ is acyclic in $X\ctrh^\lct$; and
Corollary~\ref{lct-proj-local}(c) implies that the objects of
cocycles of $\P^\bu$ belong to $X\ctrh^\lct_\prj$.
 Thus the complex $\P^\bu$ is contractible, and we can conclude that
the complex $\gB^\bu$ is contraacyclic in $X\lcth_\bW^\lct$.

 In part~(c), the functors~$j_\alpha^!$ take Positselski-contraacyclic
complexes to Positselski-contraacyclic complexes for the same reason
as in part~(b).
 Conversely, given a complex $\gB^\bu$ in $X\lcth_\bW$, let us find
a complex $\gF^\bu$ in $X\ctrh^\fl$ together with a morphism of
complexes $\gF^\bu\rarrow\gB^\bu$ with a cone contraacyclic in
$X\lcth_\bW$ (using Theorem~\ref{derived-inj-proj-resolutions}(d)).
 Then contraacyclicity of the complexes $j_\alpha^!\gB^\bu$ implies
contraacyclicity of the complexes $j_\alpha^!\gF^\bu$
in $U_\alpha\ctrh$.
 Obviously, $j_\alpha^!\gF^\bu$ are complexes in $U_\alpha\ctrh^\fl$;
and Theorem~\ref{derived-inj-proj-resolutions}(d) tells us
that the contraacyclicity of these complexes in $U_\alpha\ctrh$ implies
their (absolute) acyclicity in $U_\alpha\ctrh^\fl$.
 So the complexes $j_\alpha^!\gF^\bu$ are acyclic in $U_\alpha\ctrh$
with the objects of cocycles belonging to $U_\alpha\ctrh^\fl$.
 By Lemma~\ref{acyclicity-in-lcth-criterion}(a), it follows that
the complex $\gF^\bu$ is acyclic in $X\ctrh$; and
Corollary~\ref{finite-krull-flat-contraherent}(b) tells us that
the objects of cocycles of $\gF^\bu$ belong to $X\ctrh^\fl$.
 So the complex $\gF^\bu$ is acyclic in $X\ctrh^\fl$; 
by Corollary~\ref{finite-krull-derived-equivalences}(a), it follows
that $\gF^\bu$ is contraacyclic (in fact, absolutely acyclic)
in $X\ctrh^\fl$, hence also in $X\ctrh$.
 Now we can conclude that the complex $\gB^\bu$ is contraacyclic
in $X\lcth_\bW$.
\end{proof}

\subsection{Co-contra correspondence over a regular scheme}
 Let $X$ be a regular semi-separated Noetherian scheme of
finite Krull dimension.

\begin{thm}
\textup{(a)} The triangulated functor\/ $\sD^\co(X\qcoh_\fl)\rarrow
\sD^\co(X\qcoh)$ induced by the embedding of exact categories
$X\qcoh_\fl\rarrow X\qcoh$ is an equivalence of triangulated
categories. \par
\textup{(b)} The triangulated functor\/ $\sD^\ctr(X\ctrh^\lin)\rarrow
\sD^\ctr(X\ctrh)$ induced by the embedding of exact categories
$X\ctrh^\lin\rarrow X\ctrh$ is an equivalence of triangulated
categories. \par
\textup{(c)} There is a natural equivalence of triangulated
categories\/ $\sD^\co(X\qcoh)\simeq\sD^\ctr(X\ctrh)$ provided by
the derived functors\/ $\boR\fHom_X(\O_X,{-})$ and\/
$\O_X\ocn_X^\boL{-}$.
\end{thm}

\begin{proof}
 Part~(a) actually holds for any symbol $\bst\ne\ctr$, $\bctr$
in the upper indices of the derived category signs.
 For $\bst\ne\bco$, $\ctr$, $\bctr$, it is a particular
case of Corollary~\ref{finite-vfl-lin-dim-cor}(a).
 Indeed, one has $X\qcoh=X\qcoh_\ffdd$ provided that $d$~is greater or
equal to the Krull dimension of~$X$.
 For $\bst=\bco$, we recall that $\sD^\co(X\qcoh_\fl)=
\sD^\bco(X\qcoh_\fl)$ by
Corollary~\ref{quasi-finite-flat-dim-all-derived-coincide} and
$\sD^\co(X\qcoh)=\sD^\bco(X\qcoh)$ by
Theorem~\ref{derived-inj-proj-resolutions}(a).

 Similarly, part~(b) actually holds for any symbol $\bst\ne\co$, $\bco$
in the upper indices.
 For $\bst\ne\bctr$, $\co$, $\bco$, it is a particular case of
Corollary~\ref{finite-vfl-lin-dim-cor}(c).
 Indeed, one has $X\lcth_\bW = X\lcth_\bW^\flidd$ provided that $d$~is
greater or equal to the Krull dimension of~$X$.
 For $\bst=\bctr$, we recall that $\sD^\ctr(X\lcth_\bW^\lin)=
\sD^\bctr(X\lcth_\bW^\lin)$ by Corollary~\ref{lin-ctrh-lcth-cor}(a)
and $\sD^\ctr(X\lcth_\bW)=\sD^\bctr(X\lcth_\bW)$ by
Theorem~\ref{derived-inj-proj-resolutions}(d).

 To prove part~(c), notice that all the triangulated functors
$\sD^\abs(X\qcoh)\rarrow\sD^\co(X\qcoh)\rarrow\sD(X\qcoh)$
and $\sD^\abs(X\lcth_\bW)\rarrow\sD^\ctr(X\lcth_\bW)\rarrow
\sD(X\lcth_\bW)$ are equivalences of categories by
Corollary~\ref{vfl-lin-finite-dim-all-derived-coincide}
(since one also has $X\qcoh=X\qcoh_\fvfdd$ provided that $d$~is greater
or equal to the Krull dimension of~$X$).
 So it remains to apply Theorem~\ref{naive-co-contra-thm}. \hfuzz=7.5pt
\end{proof}

\subsection{Co-contra correspondence over a Gorenstein scheme}
 Let $X$ be a Gorenstein semi-separated Noetherian scheme of
finite Krull dimension.
 We will use the following formulation of the Gorenstein condition:
for any affine open subscheme $U\sub X$, the classes of
$\O_X(U)$\+modules of finite flat dimension, of finite projective
dimension, and of finite injective dimension coincide.

 Notice that neither of these dimensions can exceed the Krull dimension
$D$ of the scheme~$X$.
 Accordingly, the class of $\O_X(U)$\+modules defined by the above
finite homological dimension conditions is closed under both
infinite direct sums and infinite products.
 It is also closed under extensions and the passages to the cokernels
of embeddings and the kernels of surjections.

 Moreover, since the injectivity of a quasi-coherent sheaf on
a Noetherian scheme is a local property, the full subcategories of
quasi-coherent sheaves of finite flat dimension and of finite
injective dimension coincide in $X\qcoh$.
 Similarly, the full subcategories of locally contraherent cosheaves
of finite flat dimension and of finite locally injective dimension
coincide in $X\lcth$.
 Neither of these dimensions can exceed~$D$.

\begin{thm}
\textup{(a)} The triangulated functors\/ $\sD^\co(X\qcoh_\fl)\rarrow
\sD^\co(X\qcoh_\fd)\allowbreak\rarrow\sD^\co(X\qcoh)$ induced by
the embeddings of exact categories $X\qcoh_\fl\rarrow X\qcoh_\fd\rarrow
X\qcoh$ are equivalences of triangulated categories. \par
\textup{(b)} The triangulated functors\/ $\sD^\ctr(X\ctrh^\lin)\rarrow
\sD^\ctr(X\ctrh^\lid)\rarrow\sD^\ctr(X\ctrh)$ induced by the embeddings
of exact categories $X\ctrh^\lin\rarrow X\ctrh^\lid\rarrow X\ctrh$
are equivalences of triangulated categories. \par
\textup{(c)} There is a natural equivalence of triangulated
categories\/ $\sD^\co(X\qcoh_\fd)\simeq\sD^\ctr(X\ctrh^\lid)$ provided by
the derived functors\/ $\boR\fHom_X(\O_X,{-})$ and\/
$\O_X\ocn_X^\boL{-}$.
\end{thm}

\begin{proof}
 Parts~(a\+b): by Corollary~\ref{finite-vfl-lin-dim-cor}(a,c),
the functors $\sD^\st(X\qcoh_\fl)\rarrow\sD^\st(X\allowbreak\qcoh_\fd)$
are equivalences of categories for any symbol $\bst\ne\bco$, $\ctr$,
$\bctr$ and the functors $\sD^\st(X\lcth_\bW^\lin)\rarrow
\sD^\st(X\lcth_\bW^\lid)$ are equivalences of categories for any symbol
$\bst\ne\bctr$, $\co$, $\bco$.
 One also has $\sD^\co(X\qcoh_\fl)=\sD^\bco(X\qcoh_\fl)$ by
Corollary~\ref{quasi-finite-flat-dim-all-derived-coincide} and
$\sD^\co(X\qcoh)=\sD^\bco(X\qcoh)$ by
Theorem~\ref{derived-inj-proj-resolutions}(a), as well as
$\sD^\ctr(X\lcth_\bW^\lin)=\sD^\bctr(X\lcth_\bW^\lin)$ by
Corollary~\ref{lin-ctrh-lcth-cor}(a) and
$\sD^\ctr(X\lcth_\bW)=\sD^\bctr(X\lcth_\bW)$ by
Theorem~\ref{derived-inj-proj-resolutions}(d).

 To prove that the functor $\sD^\co(X\qcoh_\fd)\rarrow\sD^\co(X\qcoh)$
is an equivalence of categories, notice that one has $X\qcoh^\inj\sub
X\qcoh^\fiD=X\qcoh_\fd$ and the functor $\Hot(X\qcoh^\inj)\rarrow
\sD^\co(X\qcoh^\fiD)$ is an equivalence of categories by
Corollary~\ref{loc-noetherian-derived-fid-fpd}(a), while the composition
$\Hot(X\qcoh^\inj)\rarrow\sD^\co(X\qcoh^\fiD)\rarrow\sD^\co(X\qcoh)$
is an equivalence of categories by
Theorem~\ref{derived-inj-proj-resolutions}(a).

 Similarly, to prove that the functor $\sD^\ctr(X\lcth_\bW^\lid)\rarrow
\sD^\ctr(X\lcth_\bW)$ is an equialence of categories, notice that one
has $X\ctrh_\prj\sub X\lcth_\bW^\ffD = X\lcth_\bW^\lid$ and the functor
$\Hot(X\ctrh_\prj)\rarrow\sD^\ctr(X\lcth_\bW^\ffD)$ is an equivalence of
categories by Corolary~\ref{derived-contra-ffd-cor}, while
the composition $\Hot(X\ctrh_\prj)\rarrow\sD^\ctr(X\lcth_\bW^\ffD)
\rarrow\sD^\ctr(X\lcth_\bW)$ is an equivalence of categories by
Theorem~\ref{derived-inj-proj-resolutions}(d).

 To prove part~(c), notice that the functors $\sD^\abs(X\qcoh_\fd)
\rarrow\sD^\co(X\qcoh_\fd)\rarrow\sD(X\qcoh_\fd)$ are equivalences of
categories by Corollary~\ref{quasi-finite-flat-dim-all-derived-coincide},
while the functors $\sD^\abs(X\lcth_\bW^\lid)\rarrow
\sD^\ctr(X\lcth_\bW^\lid)\rarrow\sD(X\lcth_\bW^\lid)$ are equivalences of
categories by Corollary~\ref{vfl-lin-finite-dim-all-derived-coincide}(b).

 Furthermore, consider the intersections $X\qcoh^\cta_\fd=
X\qcoh^\cta\cap X\qcoh_\fd$ and $X\ctrh_\al^\lid=
X\ctrh_\al\cap X\lcth_\bW^\lid$.
 As was explained in Section~\ref{finite-dim-morphisms-III},
the functor $\sD^\st(X\qcoh^\cta_\fd)\allowbreak\rarrow
\sD^\st(X\qcoh_\fd)$ is an equivalence of triangulated categories
for any $\bst\ne\co$, $\ctr$, $\bco$, $\bctr$, while the functor
$\sD^\st(X\ctrh_\al^\lid)\rarrow\sD^\st(X\lcth_\bW^\lid)$ is
an equivalence of triangulated categories for any $\bst\ne\co$,
$\bco$, $\bctr$.

 Finally, it is clear from Lemma~\ref{cta-clp-lem-ffd-flid}(a,d)
(see also Lemma~\ref{cta-clp-finite-flat-inj-dim-identified}) that
the equivalence of exact categories $X\qcoh^\cta\simeq X\ctrh_\al$ of
Lemma~\ref{cta-clp-equivalence} identifies their full exact
subcategories $X\qcoh^\cta_\fd$ and $X\ctrh_\al^\lid$ under
our assumptions.
 So the induced equivalence of the derived categories $\sD^\abs$
or~$\sD$,
$$
 \sD^\abs(X\qcoh^\cta_\fd)\simeq\sD^\abs(X\ctrh_\al^\lid)
 \quad\text{or}\quad
 \sD(X\qcoh^\cta_\fd)\simeq\sD(X\ctrh_\al^\lid)
$$
provides the desired equivalence of triangulated categories
in part~(c).
\end{proof}

\subsection{Co-contra correspondence over a scheme with
a dualizing complex}
 Let $X$ be a semi-separated Noetherian scheme with a dualizing
complex~$\D_X^\bu$ \,\cite[Chapter~V]{Har},
\cite[Section Tag~0A7A]{SP}, which we will view as a finite complex of
injective quasi-coherent sheaves on~$X$.
 Notice that any Noetherian scheme with a dualizing complex has finite
Krull dimension~\cite[Corollary~V.7.2]{Har}, \cite[Lemma Tag~0A80]{SP}.
 The following result complements the covariant Serre--Grothendieck
duality theory as developed in the papers, the thesis, and
the book~\cite{IK,N-f,M-th,EP,Pfp,Psemten}.

\begin{thm} \label{co-contra-dualizing}
 There are natural equivalences between the four triangulated
categories\/ $\sD^{\abs=\co=\empt}(X\qcoh_\fl)$, \ $\sD^\co(X\qcoh)$, \
$\sD^\ctr(X\ctrh)$, and\/ $\sD^{\abs=\ctr=\empt}(X\ctrh^\lin)$.
 (Here the notation $\abs=\co=\empt$ and $\abs=\ctr=\empt$ presumes
the assertions that the corresponding derived categories of the first
and second kind coincide for the exact category in question.)
 Among these, the equivalences\/ $\sD^\abs(X\qcoh_\fl)\simeq\sD^\ctr
(X\ctrh)$ and\/ $\sD^\co(X\qcoh)\simeq\sD^\abs(X\ctrh^\lin)$ do not
require a dualizing complex and do not depend on it; all the remaining
equivalences do and do.

 More specifically, there is a commutative diagram of triangulated
equivalences provided by derived and underived functors as on
the diagram
\begin{equation} \label{co-contra-quadrality-diagram}
\begin{gathered}
 \xymatrix{
  \sD^\abs(X\qcoh_\fl) \ar@<3pt>[rrrrrr]^{\D_X^\subbu\ot_{\O_X}{-}}   
  \ar@<3pt>[ddd]^{\boR\fHom_X(\O_X,{-})}
  &&&&&& \sD^\co(X\qcoh)
  \ar@<3pt>[llllll]^{\boR\qHom_{X\qc}(\D_X^\subbu,{-})}
  \ar@<-3pt>[llllllddd]_{\boR\fHom_X(\D_X^\subbu,{-})\quad}
  \ar@<3pt>[ddd]^{\boR\fHom_X(\O_X,{-})} \\ \\ \\
  \sD^\ctr(X\ctrh) \ar@<3pt>[rrrrrr]^{\D_X^\subbu\ot_{X\ct}^\boL{-}}
  \ar@<3pt>[uuu]^{\O_X\ocn_X^\boL{-}}
  \ar@<-3pt>[rrrrrruuu]_{\quad\D_X^\subbu\ocn_X^\boL{-}}
  &&&&&& \sD^\abs(X\ctrh^\lin)
  \ar@<3pt>[llllll]^{\Cohom_X(\D_X^\subbu,{-})}
  \ar@<3pt>[uuu]^{\O_X\ocn_X^\boL{-}}
 }
\end{gathered}
\end{equation}
\end{thm}

\begin{proof}
 For any quasi-compact semi-separated scheme $X$ with an open covering
$\bW$, one has $\sD^\abs(X\lcth_\bW^\lin)=\sD^\ctr(X\lcth_\bW^\lin)
=\sD(X\lcth_\bW^\lin)$ by
Corollary~\ref{vfl-lin-finite-dim-all-derived-coincide}(b).
 For any semi-separated Noetherian scheme $X$ of finite Krull dimension, 
one has $\sD^\abs(X\qcoh_\fl)=\sD^\co(X\qcoh_\fl)=\sD(X\qcoh_\fl)$
by Corollary~\ref{quasi-finite-flat-dim-all-derived-coincide}.

 Notice also that for any quasi-compact semi-separated scheme $X$ with
an open covering $\bW$, one has 
$\sD^\ctr(X\lcth_\bW^\lin)=\sD^\bctr(X\lcth_\bW^\lin)
=\sD(X\lcth_\bW^\lin)$ by Corollary~\ref{lin-ctrh-lcth-cor}(a).
 For any semi-separated Noetherian scheme $X$ of finite Krull dimension,
one has $\sD^\co(X\qcoh_\fl)=\sD^\bco(X\qcoh_\fl)=\sD(X\qcoh_\fl)$
by Corollary~\ref{quasi-finite-flat-dim-all-derived-coincide}.
 For any Noetherian scheme $X$ of finite Krull dimension with an open
covering $\bW$, one has $\sD^\ctr(X\lcth_\bW)=\sD^\bctr(X\lcth_\bW)$
by Theorem~\ref{derived-inj-proj-resolutions}(d).
 For any locally Noetherian scheme $X$, one has
$\sD^\co(X\qcoh)=\sD^\bco(X\qcoh)$ by
Theorem~\ref{derived-inj-proj-resolutions}(a).
 So there is really no difference between the Positselski and
the Becker co/contraderived categories in the context of
Theorem~\ref{co-contra-dualizing}.

 For any semi-separated Noetherian scheme $X$, one has $\sD^\co(X\qcoh)
\simeq\Hot(X\allowbreak\qcoh^\inj)$ by
Theorem~\ref{derived-inj-proj-resolutions}(a) and $\Hot(X\qcoh^\inj)
\simeq\sD^\abs(X\lcth_\bW^\lin)$ by Corollary~\ref{inj-co-contra-cor}(b).
 Hence the desired equivalence $\sD^\co(X\qcoh)\simeq
\sD^\abs(X\lcth_\bW^\lin)$, which is provided by the derived functors
$$
\boR\fHom_X(\O_X,{-})\:\sD^\co(X\qcoh)\lrarrow\sD^\abs(X\lcth_\bW^\lin)
$$
and
$$
\O_X\ocn_X^\boL{-}\:\sD^\abs(X\lcth_\bW^\lin)\lrarrow\sD^\co(X\qcoh).
$$
 Alternatively, one can refer to
Corollary~\ref{bco-qcoh-lin-derived-equivalence}.

 For any semi-separated Noetherian scheme $X$ of finite Krull dimension,
one has $\sD^\abs(X\qcoh_\fl)\simeq\Hot(X\ctrh^\lct_\prj)$ by
Corollary~\ref{flat-noetherian-co-contra-cor},
$\Hot(X\ctrh^\lct_\prj)\simeq\sD^\ctr(X\lcth_\bW^\lct)$ by
Theorem~\ref{derived-inj-proj-resolutions}(b), and
$\sD^\ctr(X\lcth_\bW^\lct)\simeq\sD^\ctr(X\lcth_\bW)$
by Corollary~\ref{derived-contra-lct-cor}(b).
 Alternatively, one can refer to the equivalence
$\sD^\abs(X\qcoh_\fl)\simeq\sD^\abs(X\qcoh_\vfl)$ holding by
Corollary~\ref{derived-fl-vfl-cor},
$\sD^\abs(X\qcoh_\vfl)\simeq\Hot(X\ctrh_\prj)$ by
Corollary~\ref{vfl-co-contra-cor-expanded}(a), and
$\Hot(X\ctrh_\prj)\simeq\sD^\ctr(X\lcth_\bW)$ by
Theorem~\ref{derived-inj-proj-resolutions}(d).
 Either way, one gets the same desired equivalence $\sD^\abs(X\qcoh_\fl)
\simeq\sD^\ctr(X\lcth_\bW)$, which is provided by the derived functors
$$
\boR\fHom_X(\O_X,{-})\:\sD^\abs(X\qcoh_\fl)\lrarrow\sD^\ctr(X\lcth_\bW)
$$
and
$$
\O_X\ocn_X^\boL{-}\:\sD^\ctr(X\lcth_\bW)\lrarrow\sD^\abs(X\qcoh_\fl).
$$
 Alternatively, one can refer to
Corollary~\ref{bctr-lcth-vfl-fl-derived-equivalences}.

 Now we are going to construct a commutative diagram of equivalences
of triangulated categories
\begin{equation} \label{co-contra-quadrality-diagram-symbol-st}
\begin{gathered}
 \xymatrix{
  \sD^\st(X\qcoh_\fl) \ar@<3pt>[rrrrrr]^{\D_X^\subbu\ot_{\O_X}{-}}   
  \ar@<3pt>[ddd]^{\boR\fHom_X(\O_X,{-})}
  &&&&&& \Hot^\st(X\qcoh^\inj)
  \ar@<3pt>[llllll]^{\qHom_{X\qc}(\D_X^\subbu,{-})}
  \ar@<-3pt>[llllllddd]_{\fHom_X(\D_X^\subbu,{-})\quad}
  \ar@<3pt>[ddd]^{\fHom_X(\O_X,{-})} \\ \\ \\
  \Hot^\st(X\ctrh^\lct_\prj)
  \ar@<3pt>[rrrrrr]^{\D_X^\subbu\ot_{X\ct}{-}}
  \ar@<3pt>[uuu]^{\O_X\ocn_X{-}}
  \ar@<-3pt>[rrrrrruuu]_{\quad\D_X^\subbu\ocn_X{-}}
  &&&&&& \sD^\st(X\lcth_\bW^\lin)
  \ar@<3pt>[llllll]^{\Cohom_X(\D_X^\subbu,{-})}
  \ar@<3pt>[uuu]^{\O_X\ocn_X^\boL{-}}
 }
\end{gathered}
\end{equation}
for any symbol $\bst=\b$, $\abs+$, $\abs-$, or~$\abs$.

 The exterior vertical functors are constructed by applying the additive
functors $\O_X\ocn_X{-}$ and $\fHom_X(\O_X,{-})$ to the given complexes
termwise.
 The interior (derived) vertical functors have been defined in
Corollaries~\ref{inj-co-contra-cor}(b)
and~\ref{flat-noetherian-co-contra-cor}.
 All the functors invoking the dualizing complex $\D_X^\bu$ are
constructed by applying the respective exact functors of two
arguments to $\D_X^\bu$ and the given unbounded complex termwise and
totalizing the bicomplexes so obtained.

 First of all, one notices that the functors in the interior upper
triangle are right adjoint to the ones in the exterior.
 This follows from the adjunction~\eqref{fHom-contratensor-adjunction}
together with the adjunction of the tensor product of quasi-coherent
sheaves and the quasi-coherent internal Hom.

 The upper horizontal functors $\D_X^\bu\ot_{\O_X}{-}$ and
$\qHom_{X\qc}(\D_X^\bu,{-})$ are mutually inverse for the reasons
explained in~\cite[Theorem~8.4 and Proposition~8.9]{M-th}
and~\cite[Theorem~2.5]{EP}.
 The argument in~\cite{EP} is based on the observations that
the morphism of finite complexes of flat quasi-coherent sheaves
$$
 \F\lrarrow\qHom_{X\qc}(\D_X^\bu\;\D_X^\bu\ot_{\O_X}\F)
$$
is a quasi-isomorphism for any sheaf $\F\in X\qcoh_\fl$ and the morphism
of finite complexes of injective quasi-coherent sheaves
$$
 \D_X^\bu\ot_{\O_X}\qHom_{X\qc}(\D_X^\bu,\J)\lrarrow\J
$$
is a quasi-isomorphism for any sheaf $\J\in X\qcoh^\inj$.

 Let us additionally point out that, according to
Lemma~\ref{ext-qhom-qc}(c) and~\cite[Lemma~8.7]{M-th} (see
also~\cite[Lemma~2.5]{EP} or~\cite[Lemma~4.10]{Psemten}),
$\qHom_{X\qc}(\D_X^\bu,\J^\bu)$ is a complex of flat cotorsion
quasi-coherent sheaves for any complex $\J^\bu$ over $X\qcoh^\inj$.
 So the functor $\qHom_{X\qc}(\D_X^\bu,{-})$ actually lands in
$\Hot^\st(X\qcoh_\fl^\cot)$ (as does the functor $\O_X\ocn_X{-}$ on
the left-hand side of
the diagram~\eqref{co-contra-quadrality-diagram-symbol-st},
according to Lemma~\ref{cta-clp-restricts-to-prj-clf}(b) or
the proof of Corollary~\ref{flat-noetherian-co-contra-cor}).
 The interior upper triangle is commutative due to
the natural isomorphism~\eqref{flat-inj-fhom-qhom}.
 The exterior upper triangle is commutative due to
the natural isomorphism~\eqref{tensor-contratensor-assoc}.

 In order to discuss the equivalence of categories in the lower
horizontal line, we will need the following lemma.
 It is based on the definitions of the $\Cohom$ functor in
Section~\ref{cohom-loc-der-contrahereable} and the contraherent
tensor product functor $\ot_{X\ct}$ in Section~\ref{contraherent-tensor}.

\begin{lem}  \label{noetherian-contraherent-tensor}
 Let\/ $\J$ be an injective quasi-coherent sheaf on a semi-separated
Noetherian scheme~$X$ with an open covering\/~$\bW$.
 Then there are two well-defined exact functors
$$
 \Cohom_X(\J,{-})\: X\lcth_\bW^\lin\lrarrow X\ctrh_\alf
$$
and
$$
 \J\ot_{X\ct}{-}\: X\ctrh_\alf\lrarrow X\lcth_\bW^\lin
$$
between the exact categories $X\lcth_\bW^\lin$ and $X\ctrh_\alf$
of locally injective\/ $\bW$\+locally contraherent cosheaves and
antilocally flat contraherent cosheaves on~$X$.
 The functor\/ $\J\ot_{X\ct}{-}$ is left adjoint to the functor\/
$\Cohom_X(\J,{-})$.
 Besides, the functor\/ $\Cohom_X(\J,{-})$ takes values in
the additive subcategory $X\ctrh^\lct_\prj\sub X\ctrh_\alf$,
while the functor\/ $\J\ot_{X\ct}{-}$ takes values in the additive
subcategory $X\ctrh^\lin_\al\sub X\lcth_\bW^\lin$.
 For any quasi-coherent sheaf\/ $\M$ and any antilocally flat
contraherent cosheaf\/ $\gF$ on $X$ there is a natural isomorphism
\begin{equation}  \label{cotensor-contraherent-tensor-assoc}
 \M\ocn_X(\J\ot_{X\ct}\gF)\simeq (\M\ot_{\O_X}\J)\ocn_X\gF
\end{equation}
of quasi-coherent sheaves on~$X$.
\end{lem}

\begin{proof}
 Notice that when $X$ has finite Krull dimension (as in
Theorem~\ref{co-contra-dualizing}), the classes of flat and antilocally
flat contraherent cosheaves on $X$ coincide by
Corollary~\ref{finite-krull-flat-contraherent}(a).
 Otherwise, we only know that any antilocally flat contraherent cosheaf
on $X$ is flat (by Corollary~\ref{clf-noetherian-flat}).
 See Remark~\ref{flat-antilocally-flat-remark} for a discussion.

 Let us show that the locally cotorsion $\bW$\+locally contraherent
cosheaf $\Cohom_X(\J,\gK)$ is projective for any locally injective
$\bW$\+locally contraherent cosheaf $\gK$ on~$X$.
 Indeed, $\J$ is a direct summand of a finite direct sum of the direct
images of injective quasi-coherent sheaves $\I$ from the embeddings of
affine open subschemes $j\:U\rarrow X$ subordinate to~$\bW$.
 So it suffices to consider the case of $\J=j_*\I$.  {\hfuzz=1.5pt\par}

 According to~\eqref{lin-projection-cohom}, there is a natural
isomorphism of locally cotorsion ($\bW$\+lo\-cally) contraherent
cosheaves $\Cohom_X(j_*\I,\gK)\simeq j_!\Cohom_U(\I,j^!\gK)$ on~$X$.
 The $\O(U)$\+modules $\I(U)$ and $\gK[U]$ are injective, so
$\Hom_{\O(U)}(\I(U),\gK[U])$ is a flat cotorsion $\O(U)$\+module.
 In other words, the locally cotorsion contraherent cosheaf
$\Cohom_U(\I,j^!\gK)$ is projective on $U$, and therefore
its direct image with respect to~$j$ is projective on $X$ (see
Lemma~\ref{loc-lct-proj}(b) or
Corollary~\ref{proj-direct-inverse}(b)).

 Now let $\gF$ be an antilocally flat contraherent cosheaf on~$X$.
 Then, in particular, $\gF$ is a flat contraherent cosheaf
(Corollary~\ref{clf-noetherian-flat}), so the tensor product
$\J\ot_X\gF$ is a locally injective derived contrahereable cosheaf
on~$X$ (see Section~\ref{contraherent-tensor}).

 Moreover, by Corollary~\ref{clf-cor}(c), $\gF$ is a direct summand of
a finitely iterated extension of the direct images of flat contraherent
cosheaves from affine open subschemes of~$X$.
 It was explained in Section~\ref{contrahereable-subsect} that
derived contrahereable cosheaves on affine schemes are contrahereable
and the direct images of cosheaves with respect to affine morphisms
preserve contrahereability.
 Besides, the full subcategory of contrahereable cosheaves on $X$ is
closed under extensions in the exact category of derived contrahereable
cosheaves, and the functor $\J\ot_X{-}$ takes short exact sequences of
flat contraherent cosheaves to short exact sequences of derived
contrahereable cosheaves on~$X$ (see Section~\ref{contraherent-tensor}).
 So it follows from the isomorphism~\eqref{cosheaf-tensor-projection}
that $\J\ot_X\gF$ is a locally injective contrahereable cosheaf.

 Its contraherator $\J\ot_{X\ct}\gF=\Cr(\J\ot_X\gF)$ is consequently
a locally injective contraherent cosheaf on~$X$.
 Furthermore, according to Section~\ref{contrahereable-subsect}
the (global) contraherator construction is an exact functor commuting
with the direct images with respect to affine morphisms.
 Hence the contraherent cosheaf $\J\ot_{X\ct}\gF$ is a direct summand
of a finitely iterated extension of the direct images of (locally)
injective contraherent cosheaves from affine open subschemes of $X$,
i.~e., $\J\ot_{X\ct}\gF$ is an antilocal locally injective
contraherent cosheaf.

 We have constructed the desired exact functors.
 A combination of the adjunction isomorphisms
\eqref{cosheaf-tensor-product-adjunction}
and~\eqref{contraherator-adjunction} makes them adjoint to each
other.
 Finally, for any $\M\in X\qcoh$ and $\gF\in X\ctrh_\alf$ one has
$$
 \M\ocn_X (\J\ot_{X\ct}\gF) = \M\ocn_X\Cr(\J\ot_X\gF)\simeq
 \M\ocn_X(\J\ot_X\gF) \simeq (\M\ot_{\O_X}\J)\ocn_X\gF
$$
according to the isomorphisms~\eqref{contraherator-contratensor}
and~\eqref{sheaf-cosheaf-tensor-assoc}.
\end{proof}

 Now we can return to the proof of Theorem~\ref{co-contra-dualizing}.
 The functors in the interior lower triangle are left adjoint to
the ones in the exterior, as it follows from
the adjunction~\eqref{fHom-contratensor-adjunction} and
Lemma~\ref{noetherian-contraherent-tensor}.
 Let us show that the lower horizontal functors are mutually inverse.
 
 According to Corollary~\ref{lin-ctrh-lcth-cor}(a), the functor
$\Hot^\st(X\ctrh_\al^\lin)\rarrow\sD^\st(X\lcth_\bW^\lin)$
induced by the embedding $X\ctrh_\al^\lin\rarrow X\lcth_\bW^\lin$
is an equivalence of triangulated categories.
 Therefore, it suffices to show that for any cosheaf $\gJ\in
X\ctrh_\al^\lin$ the morphism of finite complexes of
contraherent cosheaves
$$
 \D_X^\bu\ot_{X\ct}\Cohom_X(\D_X^\bu,\gJ)\lrarrow\gJ
$$
is a homotopy equivalence (or just a quasi-isomorphism in $X\ctrh$), 
and for any cosheaf $\P\in X\ctrh^\lct_\prj$ the morphism of finite
complexes of contraherent cosheaves
$$
 \P\lrarrow\Cohom_X(\D_X^\bu\;\D_X^\bu\ot_{X\ct}\P)
$$
is a homotopy equivalence (or just a quasi-isomorphism in $X\ctrh$).

 According to Corollary~\ref{clp-lin} and Lemma~\ref{loc-lct-proj},
any object of $X\ctrh_\al^\lin$ or $X\ctrh^\lct_\prj$ is a direct
summand of a finite direct sum of direct images of objects in
the similar categories on affine open subschemes of~$X$.
 According to \eqref{lin-projection-II-cohom}
and~\eqref{cosheaf-tensor-projection} together with the results
of Section~\ref{contrahereable-subsect}, both functors
$\Cohom_X(\D_X^\bu,{-})$ and $\D_X^\bu\ot_{X\ct}{-}$ commute with
such direct images.
 So the question reduces to the case of an affine scheme $U$, for
which the distinction between quasi-coherent sheaves and
contraherent cosheaves mostly loses its significance, as both are
identified with (appropriate classes of) $\O(U)$\+modules.
 For this reason, the desired quasi-isomorphisms follow from
the similar quasi-isomorphisms for quasi-coherent sheaves obtained
in~\cite[proof of Theorem~2.5]{EP} (as quoted above).

 Alternatively, one can argue in the way similar to the proof
in~\cite{EP}.
 Essentially, this means using an ``inverse image localization''
procedure (known also as locality and colocality~\cite{Pal})
in place of the ``direct image localization'' (known otherwise
as ``antilocality''~\cite{Pal}) above.
 The argument proceeds as follows.

 Let ${}'\D_X^\bu\rarrow\D_X^\bu$ be a quasi-isomorphism between
a finite complex ${}'\D_X^\bu$ of coherent sheaves over $X$ and
the complex of injective quasi-coherent sheaves~$\D_X^\bu$.
 Then the tensor product ${}'\D_X^\bu\ot_X\gF$ is a finite complex of
contraherent cosheaves for any flat contraherent cosheaf $\gF$ on~$X$
(see the end of Section~\ref{contraherent-tensor}).
 Furthermore, the morphism ${}'\D_X^\bu\rarrow\D_X^\bu$ is
a quasi-isomorphism of finite complexes over the exact category of
coadjusted quasi-coherent sheaves $X\qcoh^\coa$ on $X$, hence
the induced morphism ${}'\D_X^\bu\ot_X\gF\rarrow\D_X^\bu\ot_X\gF$ is
a quasi-isomorphism of finite complexes over the exact category of
derived contrahereable cosheaves on~$X$.
 It follows that the morphism ${}'\D_X^\bu\ot_X\gF\simeq{}'\D_X^\bu
\ot_{X\ct}\gF\rarrow\D_X^\bu\ot_{X\ct}\gF$ is a quasi-isomorphism of
finite complexes of contraherent cosheaves on~$X$ for any flat
contraherent cosheaf~$\gF$.

 Let $\gJ$ be a cosheaf from $X\lcth_\bW^\lin$.
 In order to show that the morphism of finite complexes
$\D_X^\bu\ot_{X\ct}\Cohom_X(\D_X^\bu,\gJ)\rarrow\gJ$ is
a quasi-isomorphism over $X\lcth_\bW^\lin$, it suffices to check
that the composition ${}'\D_X^\bu\ot_X\Cohom_X(\D_X^\bu,\gJ)
\rarrow\D_X^\bu\ot_{X\ct}\Cohom_X(\D_X^\bu,\gJ)\rarrow\gJ$
is a quasi-isomorphism over $X\lcth_\bW$.
 The latter assertion can be checked locally, i.~e., it simply means
that for any affine open subscheme $U\sub X$ subordinate to $\bW$
the morphism ${}'\D_X^\bu(U)\ot_{\O_X(U)}\Hom_{\O_X(U)}
(\D_X^\bu(U),\gJ[U])\rarrow\gJ[U]$ is a quasi-isomorphism of
complexes of $\O_X(U)$\+modules.
 This can be deduced from the condition that the morphism
$\O_X(U)\rarrow\Hom_{\O_X(U)}({}'\D_X^\bu(U),\D_X^\bu(U))$ is
a quasi-isomorphism, as explained in the proof in~\cite{EP}
(see also~\cite[Lemma~4.2(b)]{Pfp}.

 Let $\gF$ be a flat contraherent cosheaf on~$X$.
 Pick a bounded above complex of very flat quasi-coherent sheaves
${}''\D_X^\bu$ over $X$ together with a quasi-isomorphism
${}''\D_X^\bu\rarrow{}'\D_X^\bu$.
 Then the bounded below complex of contraherent cosheaves
$\Cohom_X({}''\D_X^\bu\;{}'\D_X^\bu\ot_X\gF)$ is well-defined.
 The morphisms of bounded below complexes
$\Cohom_X({}''\D_X^\bu\;{}'\D_X^\bu\ot_X\gF)\rarrow
\Cohom_X({}''\D_X^\bu\;\D_X^\bu\ot_{X\ct}\gF)$ and
$\Cohom_X(\D_X^\bu\;\allowbreak\D_X^\bu\ot_{X\ct}\gF)\rarrow
\Cohom_X({}''\D_X^\bu\;{}\D_X^\bu\ot_{X\ct}\gF)$
are quasi-isomorphisms over $X\ctrh$.
 Thus in order to show that the morphism $\gF\rarrow\Cohom_X
(\D_X^\bu\;\D_X^\bu\ot_{X\ct}\gF)$ is a quasi-isomorphism in
$X\ctrh^\fl$, it suffices to check that the morphism $\gF\rarrow
\Cohom_X({}''\D_X^\bu\;\allowbreak{}'\D_X^\bu\ot_X\gF)$ is
a quasi-isomorphism of bounded below complexes over $X\ctrh$.

 The latter is again a local assertion, meaning simply that
the morphism $\gF[U]\rarrow\Hom_{\O_X(U)}({}''\D_X^\bu(U)\;
{}'\D_X^\bu(U)\ot_{\O_X(U)}\gF[U])$ is a quasi-isomorphism
of complexes of $\O_X(U)$\+modules for any affine open subscheme
$U\sub X$.
 One proves it by replacing ${}''\D_X^\bu(U)$ by a quasi-isomorphic
bounded above complex ${}'''\D_X^\bu(U)$ of finitely generated
projective $\O_X(U)$\+modules, and reducing again to the condition
that the morphism $\O_X(U)\rarrow\Hom_{\O_X(U)}({}'''\D_X^\bu(U),
{}'\D_X^\bu(U))$ is a quasi-isomorphism (cf.~\cite[proof of
Theorem~2.5]{EP} or~\cite[Lemmas~4.1(a) and~4.2(a)]{Pfp}).

\medskip

 According to Lemma~\ref{cta-clp-restricts-to-cot-inj}(b) or the proof
of Corollary~\ref{inj-co-contra-cor}(b), the functor
$\fHom_X(\O_X,{-})$ on the right-hand side of
the diagram~\eqref{co-contra-quadrality-diagram-symbol-st}
actually lands in $\Hot^\st(X\ctrh_\al^\lin)$
(as does the functor $\D_X^\bu\ot_{X\ct}{-}$, according to
Lemma~\ref{noetherian-contraherent-tensor}).
 The exterior lower triangle is commutative due to the natural
isomorphism~\eqref{flat-fhom-cohom-inj}.
 The interior lower triangle is commutative due to the natural
isomorphism~\eqref{cotensor-contraherent-tensor-assoc}.

 The assertion that the two diagonal functors on the diagram are
mutually inverse follows from the above.
 It can be also proved directly in the manner of the former of
the above two proofs of the assertion that the two lower horizontal
functors are mutually inverse.
 One needs to use the natural isomorphisms \eqref{inj-fHom-projection}
and~\eqref{contratensor-projection} for commutation with
the direct images.
 A more general assertion is provided by (the proof of)
Theorem~\ref{non-semi-separated-co-contra} below.
\end{proof}

\subsection{Co-contra correspondence over a non-semi-separated
scheme}
 The goal of this section is to obtain partial generalizations of
Theorems~\ref{naive-co-contra-thm} and~\ref{co-contra-dualizing}
to the case of a non-semi-separated Noetherian scheme.

\begin{thm}  \label{non-semi-separated-naive-co-contra}
 Let $X$ be a Noetherian scheme of finite Krull dimension.
 Then for any symbol\/ $\bst=\b$, $+$, $-$, or\/~$\empt$ there is
a natural equivalence of triangulated categories\/
$\sD^\st(X\qcoh)\simeq\sD^\st(X\ctrh)$.
\end{thm}

\begin{proof}
 According to Corollaries~\ref{finite-krull-ctrh-lcth-derived}
and~\ref{derived-contra-lct-cor}, one has
$\sD^\st(X\ctrh)\simeq\sD^\st(X\lcth_\bW)\allowbreak\simeq\sD^\st(X\lcth)
\simeq\sD^\st(X\lcth^\lct)\simeq\sD^\st(X\lcth_\bW^\lct)\simeq
\sD^\st(X\ctrh^\lct)$ for any open covering $\bW$ of the scheme~$X$.
 We will construct an equivalence of triangulated categories
$\sD(X\qcoh)\simeq\sD(X\lcth^\lct)$ and then show that it takes
the full subcategories $\sD^\st(X\qcoh)\sub\sD(X\qcoh)$ into
the full subcategories $\sD^\st(X\lcth^\lct)\sub\sD(X\lcth^\lct)$
and back for all symbols $\bst=\b$, $+$, or~$-$.

 By Lemma~\ref{grothendieck-vanishing}(a), the sheaf $\O_X$
has a finite coresolution by flasque quasi-coherent sheaves.
 We fix such a coresolution $\E^\bu$ for the time being.
 Given a complex $\M^\bu$ over $X\qcoh$, we pick a complex $\J^\bu$
over $X\qcoh^\inj$ quasi-isomorphic to $\M^\bu$ over $X\qcoh$ (see
Theorem~\ref{derived-inj-proj-resolutions}(a), cf.\
Theorem~\ref{quasi-homotopy-injective}) and assign to $\M^\bu$
the total complex of the bicomplex $\fHom_X(\E^\bu,\J^\bu)$
over $X\ctrh^\lct$.
 Given a complex $\P^\bu$ over $X\lcth^\lct$, we pick a complex
$\gF^\bu$ over $X\ctrh^\lct_\prj$ quasi-isomorphic to $\P^\bu$ over
$X\lcth^\lct$ (see Theorem~\ref{derived-inj-proj-resolutions}(b), cf.\
Theorem~\ref{finite-krull-homotopy-projective}(a) below)
and assign to $\P^\bu$ the total complex of the bicomplex
$\E^\bu\ocn_X\gF^\bu$ over $X\qcoh$.

 Let us first show that the complex $\fHom_X(\E^\bu,\J^\bu)$ over
$X\ctrh^\lct$ is acyclic whenever a complex $\J^\bu$ over
$X\qcoh^\inj$ is acyclic over $X\qcoh$.
 For any scheme point $x\in X$, let $\m_{x,X}$ denote the maximal
ideal of the local ring $\O_{x,X}$.
 By~\cite[Proposition~II.7.17]{Har}, any injective quasi-coherent
sheaf $\I$ on $X$ can be presented as an infinite direct sum
$\I=\bigoplus_{x\in X}\iota_x{}_*\widetilde I_x$, where $\iota_x\:
\Spec\O_{x,X}\rarrow X$ are the natural morphisms and
$\widetilde I_x$ are the quasi-coherent sheaves on $\Spec\O_{x,X}$
corresponding to infinite direct sums of copies of the injective
envelopes of the $\O_{x,X}$\+modules $\O_{x,X}/\m_{x,X}$.
 Let $X=\bigcup_{\alpha=1}^N U_\alpha$ be a finite affine open
covering.
 Set $S_\beta\sub X$ to be the set-theoretic complement to
$\bigcup_{\alpha<\beta}U_\alpha$ in $U_\beta$, and consider
the direct sum decomposition $\I=\bigoplus_{\alpha=1}^N\I_\alpha$
with $\I_\alpha=\bigoplus_{z\in S_\alpha}\iota_z{}_*\widetilde I_z$.

 The associated decreasing filtration $\I_{\ge\alpha}=
\bigoplus_{\beta\ge\alpha}\I_\beta$ is preserved by all morphisms
of injective quasi-coherent sheaves $\I$ on $X$ (cf.\
Theorem~\ref{proj-lct-classification} and
Lemma~\ref{proj-lct-orthogonality}).
 We obtain a termwise split filtration $\J^\bu_{\ge\alpha}$ on
the complex $\J^\bu$ with the associated quotient complexes
$\J^\bu_\alpha$ isomorphic to the direct images
$j_\alpha{}_*\K_\alpha^\bu$ of complexes of injective quasi-coherent
sheaves $\K_\alpha^\bu$ from the open embeddings
$j_\alpha\:U_\alpha\rarrow X$.
 Moreover, for $\alpha=1$ the complex of quasi-coherent sheaves
$\K_1^\bu\simeq j_1^*\J^\bu$ is acyclic, since the complex $\J^\bu$ is;
and the complex $j_1{}_*\K_1^\bu$ is acyclic by
Corollary~\ref{coflasque-complexes-direct-correct}(a) or
Lemma~\ref{grothendieck-vanishing}(a).
 It follows by induction that all the complexes $\K_\alpha^\bu$
over $U_\alpha\qcoh^\inj$ are acyclic over $U_\alpha\qcoh$.

 Now one has $\fHom_X(\E^\bu,j_\alpha{}_*\K_\alpha^\bu)\simeq
j_\alpha{}_!\fHom_{U_\alpha}(j_\alpha^*\E^\bu,\K_\alpha^\bu)$
by~\eqref{inj-fHom-projection}.
 The complex $\fHom_{U_\alpha}(j_\alpha^*\E^\bu,\K_\alpha^\bu)$
over $U_\alpha\ctrh^\lct$ is quasi-isomorphic to
$\fHom_{U_\alpha}(\O_{U_\alpha},\K_\alpha^\bu)$, since
$\K_\alpha^\bu$ is a complex over $U_\alpha\qcoh^\inj$, while
the complex $\fHom_{U_\alpha}(\O_{U_\alpha},\K_\alpha^\bu)$ is
acyclic over $U_\alpha\ctrh^\lct$, since the complex $\K_\alpha^\bu$
is acyclic over $U_\alpha\qcoh^\cot$ (see
Corollary~\ref{raynaud-gruson-cotors-cor} or
Lemma~\ref{flat-cotors-homol-dim}(c), or more generally,
Theorem~\ref{cotorsion-periodicity} or~\cite[Corollary~10.4]{PS6}).
 So the complex $\fHom_{U_\alpha}(j_\alpha^*\E^\bu,\K_\alpha^\bu)$
is acyclic over $U_\alpha\ctrh^\lct$; by
Lemma~\ref{co-flasque-preservation}(c), it is also a complex
of coflasque contraherent cosheaves.
 By Corollary~\ref{coflasque-complexes-direct-correct}(c), or
alternatively by Lemma~\ref{grothendieck-vanishing}(b) together with
Corollary~\ref{coflasque-direct}(b), it follows that the complex
$\fHom_X(\E^\bu,j_\alpha{}_*\K_\alpha^\bu)$ is acyclic over
$X\ctrh^\lct$.

 Therefore, so is the complex $\fHom_X(\E^\bu,\J^\bu)$.
 Similarly one proves that the complex $\E^\bu\ocn_X\gF^\bu$ is
acyclic over $X\qcoh$ whenever a complex $\gF^\bu$ over
$X\ctrh^\lct_\prj$ is acyclic over $X\ctrh^\lct$.
 One has to use Theorem~\ref{proj-lct-classification}
and Lemma~\ref{proj-lct-orthogonality} (see the proof of
Theorem~\ref{contraderived-compactly-generated}(c) below),
the isomorphism~\eqref{flat-contratensor-projection}, and
Lemma~\ref{co-flasque-preservation}(d).

 We have shown that the derived functors $\boR\fHom_X(\E^\bu,{-})$
and $\E^\bu\ocn_X^\boL{-}$ are well defined by the above rules
$\M^\bu\mpsto\fHom_X(\E^\bu,\J^\bu)$ and
$\P^\bu\mpsto\E^\bu\ocn_X\gF^\bu$.
 It is a standard fact that
the adjunction~\eqref{fHom-contratensor-adjunction} makes such
two triangulated functors adjoint to each other
(cf.~\cite[Lemma~8.3]{Psemi}).
 Let us check that the adjunction morphism $\E^\bu\ocn_X^\boL
\fHom_X(\E^\bu,\J^\bu)\rarrow\J^\bu$ is an isomorphism in
$\sD(X\qcoh)$ for any complex $\J^\bu$ over $X\qcoh^\inj$.
 For the reasons explained above, one can assume $\J^\bu =
j_\alpha{}_*\K_\alpha^\bu$ for some complex $\K_\alpha^\bu$ over
$U_\alpha\qcoh^\inj$.
 Then $\fHom_X(\E^\bu,\J^\bu)\simeq j_\alpha{}_!\fHom_{U_\alpha}
(j_\alpha^*\E^\bu,\K_\alpha^\bu)$.

 Let $\gG_\alpha^\bu$ be a complex over $U_\alpha\ctrh^\lct_\prj$
endowed with a quasi-isomorphism $\gG_\alpha^\bu\rarrow
\fHom_{U_\alpha}(j_\alpha^*\E^\bu,\K_\alpha^\bu)$ over
$U_\alpha\ctrh^\lct$.
 Then $j_\alpha{}_!\gG_\alpha^\bu$ is a complex over
$X\ctrh^\lct_\prj$, and the morphism $j_\alpha{}_!\gG_\alpha^\bu
\rarrow j_\alpha{}_!\fHom_{U_\alpha}(j_\alpha^*\E^\bu,\K_\alpha^\bu)$
is a quasi-isomorphism over $X\ctrh^\lct$ (as both $\gG_\alpha^\bu$
and $\fHom_{U_\alpha}(j_\alpha^*\E^\bu,\K_\alpha^\bu)$ are complexes
of coflasque contraherent cosheaves).
 So one has $\E^\bu\ocn_X^\boL\fHom_X(\E^\bu,\J^\bu) \simeq
\E^\bu\ocn_X j_\alpha{}_!\gG_\alpha^\bu\simeq
j_\alpha{}_*(j_\alpha^*\E^\bu\ocn_{U_\alpha}\gG_\alpha^\bu)$
by~\eqref{flat-contratensor-projection}.
 Both $\K_\alpha^\bu$ and $j_\alpha^*\E^\bu\ocn_{U_\alpha}\gG_\alpha^\bu$
being complexes of flasque quasi-coherent sheaves on $U_\alpha$ by
Lemma~\ref{co-flasque-preservation}(d), it remains to show that
the natural morphism  $j_\alpha^*\E^\bu\ocn_{U_\alpha}\gG_\alpha^\bu
\rarrow\K_\alpha^\bu$ is a quasi-isomorphism over $U_\alpha\qcoh$. 
 Now the morphisms $\O_{U_\alpha}\ocn_{U_\alpha}\gG_\alpha^\bu
\rarrow j_\alpha^*\E^\bu\ocn_{U_\alpha}\gG_\alpha^\bu$ and
$\O_{U_\alpha}\ocn_{U_\alpha}\gG_\alpha^\bu\rarrow\O_{U_\alpha}
\ocn_{U_\alpha}\fHom_{U_\alpha}(j_\alpha^*\E^\bu,\K_\alpha^\bu)
\rarrow\O_{U_\alpha}\ocn_{U_\alpha}\fHom_{U_\alpha}(\O_{U_\alpha},
\K_\alpha^\bu)\rarrow\K_\alpha^\bu$ are quasi-isomorphisms,
and the desired assertion follows.

 Similarly one shows that the adjunction morphism $\boR\fHom_X(\E^\bu\;
\E^\bu\ocn_X\gF^\bu)\rarrow\gF^\bu$ is an isomorphism in
$\sD(X\ctrh^\lct)$ for any complex $\gF^\bu$ over $X\ctrh^\lct_\prj$.
 This finishes the construction of the equivalence of categories
$\sD(X\qcoh)\simeq\sD(X\lcth^\lct)$.
 To show that it does not depend on the choice of a flasque resolution
$\E^\bu$ of the sheaf $\O_X$, consider an acyclic finite complex
$\L^\bu$ of flasque quasi-coherent sheaves on~$X$.
 Then for any injective quasi-coherent sheaf $\J$ on $X$ the complex
$\fHom_X(\L^\bu,\J)$ over $X\ctrh^\lct$ is acyclic by construction.
 Here it is helpful to notice that the complex of (flasque)
quasi-coherent sheaves $j_*j^*\L^\bu$ on $X$ is acyclic for any
open embedding morphism $j\:U\rarrow X$.
 To show that the complex $\L^\bu\ocn_X\gF$ is acyclic over $X\qcoh$
for any cosheaf $\gF\in X\ctrh^\lct_\prj$, one reduces the question
to the case of an affine scheme~$X$ using
Theorem~\ref{proj-lct-classification}(b) and
Lemma~\ref{co-flasque-preservation}(d).

 Finally, it remains to show that the equivalence of categories
$\sD(X\qcoh)\simeq\sD(X\lcth^\lct)$ that we have constructed takes
bounded above (resp., below) complexes to bounded above (resp.,
below) complexes and vice versa (up to quasi-isomorphism).
 If a complex $\M^\bu$ over $X\qcoh$ is bounded below, it has
a bounded below injective resolution $\J^\bu$ and the complex
$\fHom_X(\E^\bu,\J^\bu)$ over $X\ctrh^\lct$ is also bounded below.
 Now assume that a complex $\J^\bu$ over $X\qcoh^\inj$
has bounded above cohomology.

 Arguing as above, consider its decreasing filtration
$\J_{\ge\alpha}^\bu$ with the associated quotient complexes
$\J_\alpha^\bu\simeq j_\alpha{}_*\K_\alpha^\bu$.
 Using Lemma~\ref{grothendieck-vanishing}(a), one shows that
the cohomology sheaves of the complexes $\K_\alpha^\bu$ over
$U_\alpha\qcoh^\inj$ are also bounded above.
 By Corollary~\ref{raynaud-gruson-cotors-cor}, the coresolution
dimension of any bounded complex in $U_\alpha\qcoh$ with respect to
the coresolving subcategory $U_\alpha\qcoh^\cot\sub U_\alpha\qcoh$
is finite, and it follows that the complex
$\fHom_{U_\alpha}(\O_X,\K_\alpha^\bu)$ is quasi-isomorphic to
a bounded above complex over $U_\alpha\ctrh^\lct$.
 The complex $\fHom_{U_\alpha}(j_\alpha^*\E^\bu,\K_\alpha^\bu)$
over $U_\alpha\ctrh^\lct$ is quasi-isomorphic to
$\fHom_{U_\alpha}(\O_{U_\alpha},\K_\alpha^\bu)$.
 Finally, the complex $\fHom_X(\E^\bu,\J_\alpha^\bu)\allowbreak\simeq
j_\alpha{}_!\fHom_{U_\alpha}(j_\alpha^*\E^\bu,\K_\alpha^\bu)$ 
is quasi-isomorphic to a bounded above complex over $X\ctrh^\lct$
by Lemma~\ref{co-flasque-preservation}(c),
Corollary~\ref{coflasque-complexes-direct-correct}(c), and
Corollary~\ref{coflasque-resolutions-finite}(b).

 Similarly one can show that, for any complex $\gF^\bu$ over
$X\ctrh^\lct_\prj$ quasi-isomorphic to a bounded below complex
over $X\ctrh^\lct$, the complex $\E^\bu\ocn_X\gF^\bu$ over $X\qcoh$
has bounded below cohomology sheaves.
 A further discussion can be found in
Corollary~\ref{naive-co-contra-works-with-cplxes-of-flat-cosheaves}
and the proof of Theorem~\ref{roos-axiom-for-finite-Krull-dim} below.
\end{proof}

 Now let $X$ be a Noetherian scheme with a dualizing
complex~$\D_X^\bu$ \,\cite[Chapter~V]{Har}, \cite[Section Tag~0A7A]{SP}.
 As above, we will consider $\D_X^\bu$ as a finite complex of injective
quasi-coherent sheaves on~$X$.
 The following partial version of the covariant Serre--Grothendieck
duality holds without the semi-separatedness assumption on~$X$.

\begin{thm}  \label{non-semi-separated-co-contra}
 The choice of a dualizing complex\/ $\D_X^\bu$ induces a natural
equivalence of triangulated categories\/ $\sD^\co(X\qcoh)\simeq
\sD^\ctr(X\ctrh)$.
\end{thm}

\begin{proof}
 According to Corollary~\ref{derived-contra-lct-cor}(b),
one has $\sD^\ctr(X\lcth_\bW^\lct)\simeq\sD^\ctr(X\lcth_\bW)$ and
$\sD^\bctr(X\lcth_\bW^\lct)\simeq\sD^\bctr(X\lcth_\bW)$ for
any open covering $\bW$ of the scheme~$X$.
 By Theorem~\ref{derived-inj-proj-resolutions}(a\+b), one has
$\Hot(X\qcoh^\inj)\simeq\sD^{\co=\bco}(X\qcoh)$ and
$\Hot(X\allowbreak\ctrh^\lct_\prj)\simeq
\sD^{\ctr=\bctr}(X\lcth_\bW^\lct)$.
 We will show that the functors $\fHom_X(\D_X^\bu,{-})$ and
$\D_X^\bu\ocn_X{-}$ induce an equivalence of the homotopy categories
$\Hot^\st(X\qcoh^\inj)\simeq\Hot^\st(X\ctrh^\lct_\prj)$ for any
symbol $\bst=\b$, $+$, $-$, or~$\empt$.

 Let $\I$ be a quasi-coherent sheaf on $X$ and $j\:U\rarrow X$ be
the embedding of an affine open subscheme.
 Then the results of Section~\ref{compatibility-subsect} provide
a natural isomorphism of contraherent cosheaves $\fHom_X(\I,j_*\J)
\simeq j_!\fHom_U(j^*\I,\J)$ on $X$ for any injective quasi-coherent
sheaf $\J$ on $U$ and a natural isomorphism of quasi-coherent sheaves
$\I\ocn_X j_!\gG\simeq j_*(j^*\I\ocn_X \gG)$ on $X$ for any flat
cosheaf of $\O_U$\+modules $\gG$ on~$U$.

 Notice that the functor~$j_*$ takes injective quasi-coherent sheaves
to injective quasi-coherent sheaves and the functor~$j_!$ takes
projective locally cotorsion contraherent cosheaves to projective 
locally cotorsion contraherent cosheaves 
(Corollary~\ref{proj-lct-direct}(b)).
 Furthermore, let $X=\bigcup_\alpha U_\alpha$ be a finite affine
open covering.
 It is clear from the classification theorems (see
Theorem~\ref{proj-lct-classification}(b)) that any injective
quasi-coherent sheaf or projective locally cotorsion contraherent 
cosheaf on $X$ is a finite direct sum of the direct images of
similar (co)sheaves from~$U_\alpha$.

 Now assume that $\I$ is an injective quasi-coherent sheaf on~$X$.
 Then it follows from the discussion above together with
Lemmas~\ref{cotors-hom}(b) and~\ref{coherent-tensor-hom-lemma} that
the functors $\fHom_X(\I,{-})$ and $\I\ocn_X{-}$ take injective
quasi-coherent sheaves to projective locally cotorsion contraherent
cosheaves on~$X$ and back.
 By~\eqref{fHom-contratensor-adjunction}, these are two adjoint 
functors between the additive categories $X\qcoh^\inj$ and
$X\ctrh^\lct_\prj$.
 Substituting $\D_X^\bu$ in place of $\I$ and totalizing the finite
complexes of complexes of (co)sheaves, we obtain two adjoint functors
$\fHom_X(\D_X^\bu,{-})$ and $\D_X^\bu\ocn_X{-}$ between the homotopy
categories $\Hot^\st(X\qcoh^\inj)$ and $\Hot^\st(X\ctrh^\lct_\prj)$.

 In order to show that these are mutually inverse equivalences, it
suffices to check that the adjunction morphisms
$\D_X^\bu\ocn_X\fHom_X(\D_X^\bu,\J)\rarrow\J$ and
$\P\rarrow\fHom_X(\D_X^\bu\;\allowbreak\D_X^\bu\ocn_X\P)$
are quasi-isomorphisms/homotopy equivalences of finite complexes for
any $\J\in X\qcoh^\inj$ and $\P\in X\ctrh^\lct_\prj$.
 Presenting $\J$ and $\P$ as finite direct sums of the direct images 
of similar (co)sheaves from affine open subschemes of $X$
and taking again into account the isomorphisms
\eqref{inj-fHom-projection}, \eqref{flat-contratensor-projection}
reduces the question to the case of an affine scheme, where
the assertion is already known.

 Alternatively, one can work directly in the greater generality of
arbitrary (not necessarily locally cotorsion) and flat contraherent
cosheaves.
 According to Theorem~\ref{derived-inj-proj-resolutions}(d), one has
$\sD^\abs(X\ctrh^\fl)\simeq\sD^\ctr(X\lcth_\bW)$.
 Let us show that the functors $\fHom_X(\D_X^\bu,{-})$ and
$\D_X^\bu\ocn_X{-}$ induce an equivalence of triangulated categories
$\Hot^\st(X\qcoh^\inj)\simeq\sD^\st(X\ctrh^\fl)$ for any symbol
$\bst=\b$, $\abs+$, $\abs-$, or~$\abs$.

 Given an injective quasi-coherent sheaf $\I$ on $X$, let us first
check that the functor $\gF\mpsto \I\ocn_X\gF$ takes short exact
sequences of flat contraherent cosheaves to short exact sequences of
quasi-coherent sheaves on~$X$.
 By the adjunction isomorphism~\eqref{fHom-contratensor-adjunction},
for any injective quasi-coherent sheaf $\J$ on $X$ one has
$\Hom_X(\I\ocn_X\gF\;\J)\simeq \Hom^X(\gF,\fHom_X(\I,\J))$.
 The contraherent cosheaf $\Q=\fHom_X(\I,\J)$ being locally cotorsion,
the functor $\gF\mpsto\Hom^X(\gF,\Q)$ is exact on $X\ctrh^\fl$ by
Corollary~\ref{finite-krull-flat-clf-cor}(a).

 Furthermore, by part~(b) of the same Corollary any flat contraherent
cosheaf $\gF$ on $X$ is a direct summand of a finitely iterated
extension of the direct images $j_!\gG$ of flat contraherent cosheaves
$\gG$ on affine open subschemes $U\sub X$.
 Using the isomorphism~\eqref{flat-contratensor-projection}, we
conclude that the quasi-coherent sheaf $\I\ocn_X\gF$ is injective.
 It follows, in particular, that the complex of quasi-coherent
sheaves $\I\ocn_X\gF^\bu$ is contractible for any acyclic complex
$\gF^\bu$ over the exact category $X\ctrh^\fl$.
 Therefore, the same applies to the complex $\D_X^\bu\ocn_X\gF^\bu$
over $X\qcoh^\inj$.

 Finally, to prove that the map $\gF\rarrow\fHom_X(\D_X^\bu\;
\D_X^\bu\ocn_X\gF)$ is a quasi-isomorphism for any flat contraherent
cosheaf $\gF$, it suffices again to consider the case $\gF=j_!\gG$,
when the assertion follows from the isomorphisms
\eqref{inj-fHom-projection},~\eqref{flat-contratensor-projection}.
 Hence the morphism $\gF^\bu\rarrow\fHom_X(\D_X^\bu\;
\D_X^\bu\ocn_X\gF^\bu)$ has a cone absolutely acyclic with respect
to $X\ctrh^\fl$ for any complex $\gF^\bu$ over $X\ctrh^\fl$.
\end{proof}

\subsection{Acyclic complexes of flat contraherent cosheaves}
\label{acyclic-complexes-of-flat-cosheaves-subsect}
 The aim of this section is to construct the derived functor
$\E^\bu\ocn_X^\boL{-}\,\:\sD(X\ctrh)\rarrow\sD(X\qcoh)$ providing
the triangulated equivalence of
Theorem~\ref{non-semi-separated-naive-co-contra} in terms of complexes
of locally contraadjusted (rather than locally cotorsion) contraherent
cosheaves on $X$ and their resolutions by complexes of projective or
flat contraherent cosheaves on~$X$.
 For this purpose, we need to describe the complexes of flat
contraherent cosheaves on $X$ that are \emph{acyclic in $X\ctrh$}
(rather than in $X\ctrh^\fl$).

 Let $X$ be a Noetherian scheme of finite Krull dimension with an open
covering~$\bW$.
 By Theorem~\ref{derived-inj-proj-resolutions}(b), the classes of
Becker-contraacyclic and Positselski-contraacyclic complexes in
the exact category $X\lcth_\bW^\lct$ coincide.
 By Theorem~\ref{derived-inj-proj-resolutions}(d), the same applies
to the exact category $X\lcth_\bW$.
 By Corollary~\ref{derived-contra-lct-cor}(b), a complex in
$X\lcth_\bW^\lct$ is contraacyclic in $X\lcth_\bW^\lct$ if and only if
it is contraacyclic in $X\lcth_\bW$.
 Corollary~\ref{locality-of-co-contra-acyclicity-on-loc-Notherian}(b\+c)
tells us that contraacyclicity of complexes in $X\lcth_\bW^\lct$ or
$X\lcth_\bW$ is a local property.
 The class of contraacyclic complexes also remains unchanged when
the covering $\bW$ is replaced by its refinement.

 Moreover, by Corollary~\ref{coflasque-resolutions-finite}(b),
a complex in the exact category of coflasque locally cotorsion
contraherent cosheaves $X\ctrh^\lct_\cfq$ is (Becker or Positselski)
contraacyclic in $X\ctrh^\lct_\cfq$ if and only if it is contraacyclic
in $X\lcth_\bW^\lct$.
 Similarly, by Corollary~\ref{coflasque-resolutions-finite}(c),
a complex in the exact category of coflasque contraherent cosheaves
$X\ctrh_\cfq$ is (Becker or Positselski) contraacyclic in $X\ctrh_\cfq$
if and only if it is contraacyclic in $X\lcth_\bW$.
 So there is no ambiguity in the formulation of the following theorem.

\begin{thm} \label{flat-contraacyclic-lct-coflasque-cotorsion-pair}
\textup{(a)} The pair of classes of arbitrary complexes of flat
contraherent cosheaves and contraacyclic complexes of coflasque
locally cotorsion contraherent cosheaves $(\Com(X\ctrh^\fl)$,
$\Acycl^\ctr(X\ctrh_\cfq^\lct))$ is a hereditary complete cotorsion
pair in the exact category\/ $\Com(X\ctrh_\cfq)$ of complexes of\/
coflasque contraherent cosheaves on~$X$. \par
\textup{(b)} Let $X=\bigcup_\alpha U_\alpha$ be a finite affine open
covering of the scheme~$X$.
 Then any complex of flat contraherent cosheaves on $X$ is a direct
summand of a finitely iterated extension of the direct images of
complexes of flat contraherent cosheaves from~$U_\alpha$.
\end{thm}

\begin{proof}
 This is a version of
Lemma~\ref{all-alf-becker-contraacyclic-lct-cotorsion-pair}
for non-semi-separated Noetherian schemes of finite Krull dimension.
 Let us first prove that $\Ext^1_{\Com(X\lcth_\bW)}(\gF^\bu,\Q^\bu)=0$
for all complexes $\gF^\bu\in\Com(X\ctrh^\fl)$ and
$\Q^\bu\in\Acycl^\ctr(X\lcth_\bW^\lct)$.
 In view of Corollary~\ref{finite-krull-flat-clf-cor}(a) and
Lemma~\ref{Ext-1-as-homotopy-Hom}, it suffices to show that any
morphism of complexes $\gF^\bu\rarrow\Q^\bu$ is homotopic to zero.
 This is clear from the definition of a Positselski-contraacylic
complex in $X\lcth_\bW^\lct$ together with the same
Corollary~\ref{finite-krull-flat-clf-cor}(a).

 Alternatively, by Corollary~\ref{lct-prj-envelope}(b) and (the proof
of) the dual version of Proposition~\ref{finite-resolutions}, there
exists a complex of projective locally cotorsion contraherent cosheaves
$\gG^\bu\in\Com(X\ctrh^\lct_\prj)$ together with a morphism of complexes
$\gF^\bu\rarrow \gG^\bu$ whose cone $\gA^\bu$ is absolutely acyclic
in $X\ctrh^\fl$.
 Now any morphism of complexes $\gG^\bu\rarrow\Q^\bu$ is homotopic to
zero by the definition of Becker-contraacyclicity in $X\lcth_\bW^\lct$,
while it is clear from Corollary~\ref{finite-krull-flat-clf-cor}(a)
that any morphism of complexes $\gA^\bu\rarrow\Q^\bu$ is homotopic
to zero.

 Given this initial observation, both the assertions of the theorem
follow from the next two lemmas
(Lemmas~\ref{coflasque-contraacyclic-lct-preenvelope}
and~\ref{coflasque-termwise-flat-precover}) by virtue of
Lemma~\ref{cotorsion-pair-direct-summands-lemma}.
 Notice that the class of flat contraherent cosheaves on Noetherian
schemes of finite Krull dimension is preserved by (extensions and)
direct images with respect to flat morphisms of schemes by
Corollary~\ref{finite-krull-flat-direct}(b).
\end{proof}

\begin{lem} \label{coflasque-contraacyclic-lct-preenvelope}
 Let $X=\bigcup_\alpha U_\alpha$ be a finite affine open covering of
a Noetherian scheme $X$ of finite Krull dimension.
 Then any complex of coflasque contraherent cosheaves\/ $\gE^\bu$ on $X$
can be included in a short exact sequence\/ $0\rarrow\gE^\bu\rarrow
\Q^\bu\rarrow\gF^\bu\rarrow0$ in\/ $\Com(X\ctrh_\cfq)$ such that\/
$\Q^\bu\in\Acycl^\ctr(X\ctrh_\cfq^\lct)$, while\/ $\gF^\bu$ is
a finitely iterated extension of the direct images of complexes of
flat contraherent cosheaves from~$U_\alpha$.
\end{lem}

\begin{proof}
 The argument is similar to the proof of
Lemma~\ref{coflasque-lct-envelope}(a); cf.\ the proof of
Lemma~\ref{coflasque-acyclic-preenvelope} below.
 Let us discuss this argument as a version of
Theorem~\ref{loc-contraherent-gluing-theorem} (or more specifically,
Proposition~\ref{loc-contraherent-gluing-preenvelope}) for complexes
of coflasque contraherent cosheaves on a non-semi-separated
Noetherian scheme.
 As the scheme is not assumed to be semi-separated, it is important
that all the constructions proceed within the realm of (complexes of)
coflasque contraherent cosheaves.
 The fact that the direct images with respect to quasi-compact
semi-separated (or quasi-separated) morphisms of schemes preserve
coflasqueness and restrict to exact functors between the exact
categories of coflasque contraherent cosheaves
(Corollary~\ref{coflasque-direct}) is used here.

 Given a commutative ring $R$, denote by $R\modl_\cfq^\cot\sub
R\modl_\cfq^\cta\sub R\modl$ the full subcategories in $R\modl$
formed by the modules corresponding to coflasque (locally cotorsion
or locally contraadjusted) contraherent cosheaves on $\Spec R$.
 So, for a Noetherian commutative ring $R$ of finite Krull dimension,
we have $R\modl_\fl^\cta\sub R\modl_\cfq^\cta\sub R\modl^\cta$
(by Corollaries~\ref{coherent-flat-local}(a)
and~\ref{finite-krull-flat-contraherent}(c)) and
$R\modl_\fl^\cot\sub R\modl_\cfq^\cot\sub R\modl^\cot$.

 Let us explain how to adopt the approach of
Sections~\ref{colocal-classes-subsect}\+-%
\ref{gluing-cotorsion-in-lcth-subsect} to the situation at hand.
 Let $\sR$ be the local class of all Noetherian commutative rings $R$
of finite Krull dimension.
 Put $\sE^R=\sK^R=\Com(R\modl_\cfq^\cta)$, and consider the pair of
classes $\sF(R)$ and $\sC^R\sub\sE^R$ with $\sF(R)=
\Com(R\modl_\fl^\cta)$ and $\sC^R=\Acycl^\ctr(R\modl_\cfq^\cot)$.
 The pair of classes $\sF(R)$ and $\sC^R$ is a hereditary complete
cotorsion pair in~$\sE^R$; it can be obtained by restricting
the hereditary complete cotorsion pair of
Lemma~\ref{all-alf-becker-contraacyclic-lct-cotorsion-pair}
(for the affine scheme $\Spec R$) to the exact subcategory
$\Com(R\modl_\cfq^\cta)\sub\Com(R\modl^\cta)$ and using
Lemmas~\ref{restricting-hereditary-cotorsion}
and~\ref{restricting-cotorsion-pairs-lemma}(a) together with
Corollary~\ref{coflasque-acyclic}(a\+b).
 The class $\sE^R$ is colocal by Lemma~\ref{coflasque-cosheaves}(a).
 In view of Theorem~\ref{Becker-contraacyclicity-local-on-qcomp-qsep}(b)
or Corollary~\ref{locality-of-co-contra-acyclicity-on-loc-Notherian}(b),
it follows that the class $\sC^R$ is colocal, too.

 As the scheme $X$ is not assumed to be semi-separated, it is important
for applicability of the construction that the direct images with
respect to nonaffine open embedding morphisms take locally\+$\sE$
complexes to locally\+$\sE$ complexes and locally\+$\sC$ complexes to
locally $\sC$\+complexes.
 This holds by Corollary~\ref{coflasque-direct}; notice that any
exact direct product-preserving functor between exact categories takes
Positselski-contraacyclic complexes to Positselski-contraacyclic
complexes.

 With these observations in mind, the construction from the proof of
Lemma~\ref{coflasque-lct-envelope}(a) or from the proof of
Proposition~\ref{loc-contraherent-gluing-preenvelope} applies, proving
the lemma.
\end{proof}

\begin{lem} \label{coflasque-termwise-flat-precover}
 Let $X=\bigcup_\alpha U_\alpha$ be a finite affine open covering of
a Noetherian scheme $X$ of finite Krull dimension.
 Then any complex of coflasque contraherent cosheaves\/ $\gE^\bu$ on $X$
can be included in a short exact sequence\/ $0\rarrow\Q^\bu\rarrow
\gF^\bu\rarrow\gE^\bu\rarrow0$ in\/ $\Com(X\ctrh_\cfq)$ such that\/
$\Q^\bu\in\Acycl^\ctr(X\ctrh_\cfq^\lct)$, while\/ $\gF^\bu$ is
a finitely iterated extension of the direct images of complexes of
flat contraherent cosheaves from~$U_\alpha$.
\end{lem}

\begin{proof}
 Any complex in $X\ctrh_\cfq$ is the target of an admissible epimorphism
in $\Com(X\ctrh_\cfq)$ from a contractible complex of projective
contraherent cosheaves on~$X$.
 All projective contraherent cosheaves on $X$ are flat, and all
contractible complexes of projective contraherent cosheaves can be also
obtained as direct summands of finite direct sums of the direct images
of contractible complexes of projective contraherent cosheaves on
$U_\alpha$ (essentially by
Corollary~\ref{finite-krull-contrah-projective}).
 Consequently, the assertion of the lemma follows from
Lemma~\ref{coflasque-contraacyclic-lct-preenvelope} by virtue of
Lemma~\ref{salce-lemma}(a).
\end{proof}

\begin{thm} \label{acyclic-flat-contraacyclic-lct-coflasque-pair}
\textup{(a)} The pair $(\sF,\sC)$ of classes of complexes of flat
contraherent cosheaves acyclic as complexes of contraherent cosheaves,
$\sF=\Com(X\ctrh^\fl)\cap\Acycl(X\ctrh)$, and contraacyclic complexes
of coflasque locally cotorsion contraherent cosheaves,
$\sC=\Acycl^\ctr(X\ctrh_\cfq^\lct)$, is a hereditary complete cotorsion
pair in the exact category\/ $\Acycl(X\ctrh_\cfq)$ of acyclic
complexes of\/ coflasque contraherent cosheaves on~$X$. \par
\textup{(b)} Let $X=\bigcup_\alpha U_\alpha$ be a finite affine open
covering of the scheme~$X$.
 Then any complex of flat contraherent cosheaves on $X$ acyclic in
$X\ctrh$ is a direct summand of a finitely iterated extension of
the direct images of complexes of flat contraherent cosheaves
on $U_\alpha$ acyclic in $U_\alpha\ctrh$.
\end{thm}

\begin{proof}
 Notice first of all that a complex in $X\ctrh_\cfq$ is acyclic in
$X\ctrh_\cfq$ if and only if it is acyclic in $X\lcth_\bW$ (by
Lemma~\ref{finite-krull-ctrh-lcth-finite-dim}(c) and
Corollary~\ref{fdim-acyclic-cor}).
 The hereditary complete cotorsion pair in part~(a) can be obtained
from the one of
Theorem~\ref{flat-contraacyclic-lct-coflasque-cotorsion-pair}(a)
by restricting to the exact subcategory $\Acycl(X\ctrh_\cfq)\sub
\Com(X\ctrh_\cfq)$ (use Lemmas~\ref{restricting-hereditary-cotorsion}
and~\ref{restricting-cotorsion-pairs-lemma}(b)).
 In order to prove part~(b), we need to go through the rounds of
the coflasque version of the technique of
Sections~\ref{colocal-classes-subsect}\+-%
\ref{gluing-cotorsion-in-lcth-subsect} again.
 The assertion of part~(b) follows from part~(a) and
Lemma~\ref{coflasque-acyclic-termwise-flat-precover} below by
the argument of Lemma~\ref{cotorsion-pair-direct-summands-lemma}.
\end{proof}

\begin{lem} \label{acyclic-coflasque-contraacyclic-lct-preenvelope}
 Let $X=\bigcup_\alpha U_\alpha$ be a finite affine open covering of
a Noetherian scheme $X$ of finite Krull dimension.
 Then any acyclic complex of coflasque contraherent cosheaves\/
$\gE^\bu$ on $X$ can be included in a short exact sequence\/ $0\rarrow
\gE^\bu\rarrow\Q^\bu\rarrow\gF^\bu\rarrow0$ in\/ $\Acycl(X\ctrh_\cfq)$
such that\/ $\Q^\bu\in\Acycl^\ctr(X\ctrh_\cfq^\lct)$, while\/ $\gF^\bu$
is a finitely iterated extension of the direct images of complexes of
flat contraherent cosheaves on $U_\alpha$ acyclic in $U_\alpha\ctrh$.
\end{lem}

\begin{proof}
 This is similar to Lemma~\ref{coflasque-contraacyclic-lct-preenvelope}.
 Let $\sR$ be the local class of all Noetherian commutative rings of
finite Krull dimension.
 Put $\sE^R=\Acycl(R\modl_\cfq^\cta)\sub\Com(R\modl_\cfq^\cta)=\sK^R$,
and consider the pair of classes $\sF(R)=\Com(R\modl_\fl^\cta)\cap
\Acycl(R\modl^\cta)$ and $\sC^R=\Acycl^\ctr(R\modl_\cfq^\cot)\sub\sE^R$.
 The pair of classes $\sF(R)$ and $\sC^R$ is a hereditary complete
cotorsion pair in $\sE^R$ by
Theorem~\ref{acyclic-flat-contraacyclic-lct-coflasque-pair}(a),
which we have already proved.
 The class $\sE^R$ is colocal by
Lemmas~\ref{acyclicity-in-lcth-criterion}(a)
and~\ref{coflasque-cosheaves}(a).
 The class $\sC^R$ is colocal according to the proof of
Lemma~\ref{coflasque-contraacyclic-lct-preenvelope}.

 The direct images with respect to nonaffine open embedding morphisms
take locally\+$\sE$ complexes (i.~e., acyclic complexes of
coflasque contraherent cosheaves) to locally\+$\sE$ complexes by
Corollary~\ref{coflasque-direct}(a).
 With these observations in mind, one can apply the construction from
the proof of Lemma~\ref{coflasque-lct-envelope}(a) or
Proposition~\ref{loc-contraherent-gluing-preenvelope}.
\end{proof}

\begin{lem} \label{coflasque-acyclic-termwise-flat-precover}
 Let $X=\bigcup_\alpha U_\alpha$ be a finite affine open covering of
a Noetherian scheme $X$ of finite Krull dimension.
 Then any acyclic complex of coflasque contraherent cosheaves\/
$\gE^\bu$ on $X$ can be included in a short exact sequence\/
$0\rarrow\Q^\bu\rarrow\gF^\bu\rarrow\gE^\bu\rarrow0$ in\/
$\Com(X\ctrh_\cfq)$ such that\/
$\Q^\bu\in\Acycl^\ctr(X\ctrh_\cfq^\lct)$, while\/ $\gF^\bu$ is
a finitely iterated extension of the direct images of complexes of
flat contraherent cosheaves on $U_\alpha$ acyclic in $U_\alpha\ctrh$.
\end{lem}

\begin{proof}
 Similar to the proof of
Lemma~\ref{coflasque-termwise-flat-precover} and based on
Lemma~\ref{acyclic-coflasque-contraacyclic-lct-preenvelope}.
\end{proof}

 Now that both
Theorems~\ref{flat-contraacyclic-lct-coflasque-cotorsion-pair}
and~\ref{acyclic-flat-contraacyclic-lct-coflasque-pair} are proved,
we can deduce the corollary that we have been aiming at.

\begin{cor} \label{naive-co-contra-works-with-cplxes-of-flat-cosheaves}
 Let $X$ be a Noetherian scheme of finite Krull dimension, and let\/
$\E^\bu$ be a finite coresolution of the structure sheaf\/ $\O_X$ by
flasque quasi-coherent sheaves on~$X$.
 Then, for any complex\/ $\gF^\bu$ in $X\ctrh^\fl$ that is acyclic in
$X\ctrh$, the complex of flasque quasi-coherent sheaves\/
$\E^\bu\ocn_X\gF^\bu$ on $X$ is acyclic.
\end{cor}

\begin{proof}
 First of all, for any quasi-coherent sheaf $\M$ on $X$,
the contratensor product functor
$\M\ocn_X{-}\,\:X\ctrh^\fl\rarrow X\qcoh$ is exact,
as it was essentially explained in the proof of
Theorem~\ref{non-semi-separated-co-contra}.
 This argument is based on
the adjunction isomorphism~\eqref{fHom-contratensor-adjunction}
and Corollary~\ref{finite-krull-flat-clf-cor}(a).
 Using Corollary~\ref{finite-krull-flat-clf-cor}(b),
Lemma~\ref{co-flasque-preservation}(d), and
the isomorphism~\eqref{flat-contratensor-projection}, one can show
that the quasi-coherent sheaf $\M\ocn_X\gF$ is flasque for any flasque
quasi-coherent sheaf $\M$ and any flat contraherent cosheaf $\gF$
on~$X$ (cf.\ the proof of
Lemma~\ref{open-embed-fq-flat-ctrtensor-projection-lemma} below).

 Similarly, by
Theorem~\ref{acyclic-flat-contraacyclic-lct-coflasque-pair}(b), any
complex $\gF^\bu$ in $X\ctrh^\fl$ that is acyclic in $X\ctrh$ is
a direct summand of a finitely iterated extension of the direct
images of some complexes $\gF_\alpha^\bu$ in $U_\alpha\ctrh^\fl$
that are acyclic in $U_\alpha\ctrh$.
 Denoting by $j_\alpha\:U_\alpha\rarrow X$ the open embedding
morphisms, we have $\E^\bu\ocn_X j_\alpha{}_!\gF_\alpha^\bu\simeq
j_\alpha{}_*(j_\alpha^*\E^\bu\ocn_{U_\alpha}\gF_\alpha^\bu)$
by formula~\eqref{flat-contratensor-projection}.

 Since $j_\alpha^*\E^\bu$ is a finite coresolution of $\O_{U_\alpha}$
and $\gF_\alpha^\bu[U_\alpha]$ is a complex of flat
$\O(U_\alpha)$\+modules acyclic in $\O(U_\alpha)\modl$, one can
easily see that $j_\alpha^*\E^\bu\ocn_{U_\alpha}\gF_\alpha^\bu$
is an acyclic complex of quasi-coherent sheaves on~$U_\alpha$.
 By Lemma~\ref{co-flasque-preservation}(d), it is also a complex of
flasque quasi-coherent sheaves.
 According to Corollary~\ref{coflasque-complexes-direct-correct}(a),
it follows that
$j_\alpha{}_*(j_\alpha^*\E^\bu\ocn_{U_\alpha}\gF_\alpha^\bu)$ is
an acyclic complex of (flasque) quasi-coherent sheaves on~$X$.
 It remains to mention that a finitely iterated extension of
acyclic complexes is again an acyclic complex in $X\qcoh$.
\end{proof}

\subsection{Roos axiom for Noetherian schemes of finite Krull
dimension}
 In this section we work out an application of
Theorem~\ref{non-semi-separated-naive-co-contra} to
category-theoretic properties of the Grothendieck abelian categories
of quasi-coherent sheaves $\sA=X\qcoh$ on Noetherian schemes $X$ of
finite Krull dimension.

 Let $\sA$ be an abelian category with infinite products and enough
injective objects, and let $n\ge0$ be an integer.
 Then one says that $\sA$ satisfies the \emph{Roos axiom
$\mathrm{AB}4^*$\+$n$} \cite[Definition~1.1]{Roos} if, for any set
$\Lambda$, the derived functor of $\Lambda$\+indexed products
$\prod_{\lambda\in\Lambda}^{(*)}\:\sA^\Lambda\rarrow\sA$ has
homological dimension~$\le n$, i.~e.,
$\prod_{\lambda\in\Lambda}^{(i)}A_\lambda=0$ for every
$\Lambda$\+indexed family of objects
$(A_\lambda\in\sA)_{\lambda\in\Lambda}$ and all $i>n$.
 This simply means that if $J^\bu_\lambda$ are injective coresolutions
of the objects $A_\lambda$ in $\sA$, then
$H^i(\prod_{\lambda\in\Lambda}J^\bu_\lambda)=0$ for $i>n$.
 The Roos axiom for the categories of quasi-coherent sheaves on
quasi-compact semi-separated schemes and stacks is
discussed in the preprints~\cite[Theorem~1.1 and Remark~3.3]{HX},
\cite[Theorem~8.27]{HPSV}, and~\cite[Theorem~1.4,
Corollary~3.11, and Remark~3.13]{PRoos}.
 For an abstract category-theoretic discussion of the approach
worked out in the following theorem, see~\cite[Corollary~4.8]{PRoos}
(and for an alternative elementary approach,
\cite[Theorem~2.12]{PRoos}).
 The results of the preprint~\cite[Theorems~2.9 and~4.3]{HPSV}
illustrate the importance of the Roos axiom.

\begin{thm} \label{roos-axiom-for-finite-Krull-dim}
 Let $X$ be a Noetherian scheme of finite Krull dimension~$\le D$.
 Then the abelian category of quasi-coherent sheaves $X\qcoh$
satisfies the Roos axiom $\mathrm{AB}4^*$\+$n$ with the parameter
$n=2D+1$.
\end{thm}

\begin{proof}
 The argument is based on
Theorem~\ref{non-semi-separated-naive-co-contra}.
 The point is that the direct product functors are not exact in
the abelian category of quasi-coherent sheaves $X\qcoh$, but they are
exact in the exact category of contraherent cosheaves $X\ctrh$.
 So a derived equivalence $\sD(X\qcoh)\simeq\sD(X\ctrh)$ allows to
compute the derived direct products in $X\qcoh$.
 In particular, this approach is applicable if one wants to prove
that the derived direct products have finite homological dimension
in $X\qcoh$, and one knows that the derived equivalence
$\sD(X\qcoh)\simeq\sD(X\ctrh)$ identifies the bounded derived category
$\sD^\b(X\qcoh)\sub\sD(X\qcoh)$ with the bounded derived category
$\sD^\b(X\ctrh)\sub\sD(X\ctrh)$.

 If one wishes to extract a specific bound on the constant~$n$ in
the Roos axiom ($n=2D+1$ in the assertion of the theorem), then one
has to look into the proof of
Theorem~\ref{non-semi-separated-naive-co-contra}.
 In this context, the result of
Corollary~\ref{naive-co-contra-works-with-cplxes-of-flat-cosheaves}
helps to obtain a better bound.
 Let us spell out the details.
 We start with a basic lemma.

\begin{lem} \label{flasque-adjunction-kernel-lemma}
 Let $X$ be a topological space, $Y\sub X$ be an open subspace, and\/
$\L$ be a flasque sheaf of abelian groups on~$X$.
 Denote by $j\:Y\rarrow X$ the open embedding map.
 Then the adjunction morphism of sheaves\/ $\L\rarrow j_*j^*\L$ is
surjective, and its kernel is a flasque sheaf of abelian groups on~$X$.
\end{lem}

\begin{proof}
 For any open subset $U\sub X$, we have $(j_*j^*\L)(U)=\L(Y\cap U)$.
 So the assumption of flasqueness of $\L$ implies surjectivity of
the natural map $\L(U)\rarrow (j_*j^*\L)(U)$.
 This is an even stronger property than surjectivity of the morphism of
sheaves $\L\rarrow j_*j^*\L$.
 Now let $V\sub U$ be a smaller open subset.
 We need to show that the natural map
$$
 \ker[\L(U)\to\L(Y\cap U)]\lrarrow\ker[\L(V)\to\L(Y\cap\nobreak V)]
$$
is surjective.
 For this purpose, put $W=V\cup(Y\cap U)\sub U$.
 It remains to observe that the abelian group $\L(W)$ is the pullback
of the pair of restriction maps $\L(V)\rarrow\L(Y\cap V)$ and
$\L(Y\cap U)\rarrow\L(Y\cap V)$ (by the sheaf axiom), and recall that
the restriction map $\L(U)\rarrow\L(W)$ is surjective (by flasqueness).
\end{proof}

 Following the proof of 
Theorem~\ref{non-semi-separated-naive-co-contra}, we choose a finite
coresolution $\E^\bu$ of the structure sheaf $\O_X$ by flasque
quasi-coherent sheaves $\E^i$ on~$X$.
 By Lemma~\ref{grothendieck-vanishing}(a), we can assume that
the complex $\E^\bu$ is concentrated in the cohomological degrees
$0\le i\le D$.
 Let $X=\bigcup_{\alpha=1}^N U_\alpha$ be a finite affine open covering
of the scheme~$X$, and $j_\alpha\:U_\alpha\rarrow X$ be the open
embedding morphisms.

 Let us start with a discussion of the construction of the functor
$\boR\fHom_X(\E^\bu,\M)$ for a single quasi-coherent sheaf $\M$ on~$X$.
 Following the proof of
Theorem~\ref{non-semi-separated-naive-co-contra}, we choose
an injective coresolution $\J^\bu$ of the object $\M$ in $X\qcoh$, 
and consider the termwise split finite decreasing filtration
$\J_{\ge\alpha}^\bu$ on the complex~$\J^\bu$.

 Let us show by induction that, for every $1\le\alpha\le N$,
the complex of injective quasi-coherent sheaves $\J_{\ge\alpha}^\bu$
on $X$ is quasi-isomorphic to a complex of flasque quasi-coherent
sheaves $\L_{\ge\alpha}^\bu$ concentrated in the cohomological degrees
$0\le i\le D$.
 For the complex $\J_{\ge1}^\bu=\J^\bu$, which is a coresolution of
the sheaf $\M$, this assertion holds by
Lemma~\ref{grothendieck-vanishing}(a); so we have the induction base.

 By construction, we have $\K_1^\bu=j_1^*\J^\bu$ and
$\J^\bu/\J_{\ge2}^\bu\simeq j_1{}_*\K_1^\bu$.
 The quasi-isomorphism $\L_{\ge1}^\bu\rarrow\J_{\ge1}^\bu$ induces
a quasi-isomorphism $j_1^*\L_{\ge1}^\bu\rarrow j_1^*\J_{\ge1}^\bu$,
which induces a quasi-isomorphism
$j_1{}_*j_1^*\L_{\ge1}^\bu\rarrow j_1{}_*j_1^*\J_{\ge1}^\bu$
(by Corollary~\ref{coflasque-complexes-direct-correct}(a)).
 Both the morphisms of complexes $\L_{\ge1}^\bu\rarrow
j_1{}_*j_1^*\L_{\ge1}^\bu$ and $\J_{\ge1}^\bu\rarrow
j_1{}_*j_1^*\J_{\ge1}^\bu$ are termwise surjective
(by Lemma~\ref{flasque-adjunction-kernel-lemma}), the so induced
morphism of the kernels is a quasi-isomorphism of complexes, too.
 Now the kernel of the morphism of complexes $\J_{\ge1}^\bu\rarrow
j_1{}_*j_1^*\J_{\ge1}^\bu$ is the complex~$\J_{\ge2}^\bu$.
 It remains to denote by $\L_{\ge2}^\bu$ the kernel of the morphism
of complexes $\L_{\ge1}^\bu\rarrow j_1{}_*j_1^*\L_{\ge1}^\bu$,
and refer to Lemma~\ref{flasque-adjunction-kernel-lemma} for
the assertion that $\L_{\ge2}^\bu$ is a complex of flasque
quasi-coherent sheaves on~$X$.

 On the next step, we have $\K_2^\bu=j_2^*\J_{\ge2}^\bu$ and
$\J_{\ge2}^\bu/\J_{\ge3}^\bu\simeq j_2{}_*\K_2^\bu$.
 Similarly to the previous paragraph, we denote by $\L_{\ge3}^\bu$
the kernel of the termwise surjective morphism of complexes
$\L_{\ge2}^\bu\rarrow j_2{}_*j_2^*\L_{\ge2}^\bu$, etc.
 Proceeding in this way, we construct the desired complexes
$\L_{\ge\alpha}^\bu$ for all $1\le\alpha\le N$.
 It follows, in particular, that the complex $\K_\alpha^\bu=
j_\alpha^*\J_{\ge\alpha}^\bu$ of injective quasi-coherent sheaves
on $U_\alpha$ is quasi-isomorphic to a complex of (flasque)
quasi-coherent sheaves $j_\alpha^*\L_{\ge\alpha}^\bu$ concentrated
in the cohomological degrees $0\le i\le D$.
 That is what we wanted to know.

 As the scheme $U_\alpha$ is affine, it follows that the complex of
injective quasi-coherent sheaves $\K_\alpha^\bu$ on $U_\alpha$ is
quasi-isomorphic to a suitable complex of contraadjusted quasi-coherent
sheaves $\cP_\alpha^\bu$ concentrated in the cohomological degrees
$0\le i\le D+1$.
 Now the complex of locally cotorsion contraherent cosheaves
$\fHom_{U_\alpha}(j_\alpha^*\E^\bu,\K_\alpha^\bu)$ on $U_\alpha$ is
quasi-isomorphic to the complex of locally injective contraherent
cosheaves $\fHom_{U_\alpha}(\O_{U_\alpha},\K_\alpha^\bu)$, which
is in turn quasi-isomorphic to the complex of locally contraadjusted
contraherent cosheaves $\fHom_{U_\alpha}(\O_{U_\alpha},\cP_\alpha^\bu)$
as a complex in $U_\alpha\ctrh$.

 Recall that $\fHom_{U_\alpha}(j^*\E^\bu,\K_\alpha^\bu)$ is
a complex of coflasque (locally cotorsion) contraherent cosheaves on
$U_\alpha$ by Lemma~\ref{co-flasque-preservation}(c).
 The full subcategory of coflasque contraherent cosheaves $X\ctrh_\cfq$
is resolving in $X\ctrh$ (see the proof of
Lemma~\ref{finite-krull-ctrh-lcth-finite-dim}(c)), and in fact
the resolution dimension does not exceed~$D$.
 Consequently, the complex of contraherent cosheaves
$\fHom_{U_\alpha}(\O_{U_\alpha},\cP_\alpha^\bu)$ on $U_\alpha$ is
quasi-isomorphic to a complex of coflasque contraherent cosheaves
concentrated in the cohomological degrees $-D\le i\le D+1$.
 Applying Corollary~\ref{coflasque-complexes-direct-correct}(b),
we can conclude that the complex of locally cotorsion contraherent
cosheaves $j_\alpha{}_!\fHom_{U_\alpha}(j_\alpha^*\E^\bu,\K_\alpha^\bu)
\simeq\fHom_X(\E^\bu,j_\alpha{}_*\K_\alpha^\bu)$ on $X$ is
quasi-isomorphic to a complex of locally contraadjusted contraherent
cosheaves concentrated in the cohomological degrees $-D\le i\le D+1$.

 Passing to a finitely iterated extension, we deduce the assertion that
the complex of locally cotorsion contraherent cosheaves
$\fHom_X(\E^\bu,\J^\bu)$ on $X$ is quasi-isomorphic to a complex of
locally contraadjusted contraherent cosheaves concentrated in
the cohomological degrees $i\le D+1$ (or more specifically,
$-D\le i\le D+1$).
 We refer to~\cite[Section~4.1]{PRoos} for a general discussion of
cohomologically bounded complexes in exact categories.

 Now let us return to the point that the direct product functors are
not exact in $X\qcoh$, but they are exact in $X\ctrh$.
 The direct products in the derived category $\sD(X\qcoh)$ can be
computed by taking the products of homotopy injective complexes in
the homotopy category $\Hot(X\qcoh)$.
 In particular, given a family of quasi-coherent sheaves
$(\M_\lambda)_{\lambda\in\Lambda}$ on $X$, choose an injective
coresolution $\J_\lambda^\bu$ in $X\qcoh$ for each sheaf~$\M_\lambda$;
then the direct product of the complexes $\prod_{\lambda\in\Lambda}
\J_\lambda^\bu$, taken termwise, represents the direct product of
the objects $\M_\lambda$ in the triangulated category $\sD(X\qcoh)$.
 On the other hand, the direct products in $\sD(X\ctrh)$ can be
computed by taking the direct products of arbitrary complexes of
contraherent cosheaves termwise (see~\cite[Lemma~1.5]{BN}
or~\cite[Proposition~1.2.1 and Lemma~3.2.10]{N-tr}).

 As any equivalence of categories, the triangulated equivalence
$\sD(X\qcoh)\simeq\sD(X\ctrh)$ preserves direct products.
 The same conclusion can be made from the explicit construction of
this triangulated equivalence.
 Notice that, for any quasi-coherent sheaf $\E$ on $X$,
the functor $\fHom_X(\E,{-})\:X\qcoh^\inj\rarrow X\ctrh^\lct$
preserves infinite products (by construction).
 Hence, in the situation at hand, the functor $\boR\fHom(\E^\bu,{-})$
takes termwise direct products of arbitrary complexes of injective
quasi-coherent sheaves to direct products in $\sD(X\ctrh)$.
 It is important for us, however, that the direct products preserve
quasi-isomorphisms of complexes in the exact category of contraherent
cosheaves.

 Returning to a family of quasi-coherent sheaves
$(\M_\lambda)_{\lambda\in\Lambda}$ on $X$ and the injective
coresolutions $\J_\lambda^\bu$ of $\M_\lambda$, we recall
that each complex of locally cotorsion contraherent cosheaves
$\fHom_X(\E^\bu,\J_\lambda^\bu)$ on $X$ is quasi-isomorphic to
a complex of locally contraadjusted contraherent cosheaves
concentrated in the cohomological degrees $-D\le i\le D+1$
(as it was shown above).
 Consequently, the complex of locally cotorsion contraherent
cosheaves $\fHom_X(\E^\bu\;\prod_{\lambda\in\Lambda}\J_\lambda^\bu)
\simeq\prod_{\lambda\in\Lambda}\fHom_X(\E^\bu,\J_\lambda^\bu)$
is also quasi-isomorphic to a complex of locally contraadjusted
contraherent cosheaves concentrated in the cohomological
degrees $-D\le i\le D+1$.

 Finally, let $\gF^\bu$ be a complex of flat (e.~g., projective)
contraherent cosheaves quasi-isomorphic to
$\fHom_X(\E^\bu\;\prod_{\lambda\in\Lambda}\J_\lambda^\bu)$ as
a complex in $X\ctrh$.
 We can choose the complex $\gF^\bu$ so that it is concentrated
in the cohomological degrees $i\le D+1$.
 In the proof of Theorem~\ref{non-semi-separated-naive-co-contra},
the derived functor $\E^\bu\ocn_X^\boL{-}$ was constructed by
applying the functor $\E^\bu\ocn_X{-}$ to complexes of projective
locally cotorsion contraherent cosheaves; but
Corollary~\ref{naive-co-contra-works-with-cplxes-of-flat-cosheaves}
implies that the same derived functor can be obtained by applying
the functor $\E^\bu\ocn_X{-}$ to arbitrary complexes of
flat contraherent cosheaves on~$X$.
 As the functors $\boR\fHom_X(\E^\bu,{-})$ and $\E^\bu\ocn_X^\boL{-}$
in Theorem~\ref{non-semi-separated-naive-co-contra} are inverse
to each other, we can conclude that the complex of quasi-coherent
sheaves $\E^\bu\ocn_X\gF^\bu$ is quasi-isomorphic to the complex
$\prod_{\lambda\in\Lambda}\J_\lambda^\bu$ in $X\qcoh$.

 Since the complex $\E^\bu$ is concentrated in the cohomological
degrees $0\le i\le D$ and the complex $\gF^\bu$ is concentrated
in the cohomological degrees $i\le D+1$, the contratensor product
complex $\E^\bu\ocn_X\gF^\bu$ is concentrated in the cohomological
degrees $i\le 2D+1$.
 We have shown that the complex of quasi-coherent sheaves
$\prod_{\lambda\in\Lambda}\J_\lambda^\bu$ is quasi-isomorphic
to a complex of quasi-coherent sheaves concentrated in
the cohomological degrees $i\le 2D+1$.
 Therefore, the quasi-coherent sheaves of cohomology of the complex
$\prod_{\lambda\in\Lambda}\J_\lambda^\bu$ vanish in
the cohomological degrees $i>2D+1$, as desired.
\end{proof}

\Section{Morphisms of Finite Type Between Noetherian Schemes}
\label{finite-type-morphisms-sect}

\subsection{Compact generators} \label{compact-generators-subsect}
 Let $\sD$ be a triangulated category where arbitrary infinite
direct sums exist.
 We recall that object $C\in\sD$ is called \emph{compact} if
the functor $\Hom_\sD(C,{-})$ takes infinite direct sums in $\sD$
to infinite direct sums of abelian groups~\cite{N-c,N-bb,Rou}.
 A set of compact objects $\sC\sub\sD$ is said to \emph{generate}
$\sD$ if any object $X\in\sD$ such that $\Hom_\sD(C,X[*])=0$ for
all $C\in\sC$ vanishes in~$\sD$.

 Equivalently, this means that any full triangulated subcategory of
$\sD$ containing $\sC$ and closed under infinite direct sums
coincides with~$\sD$ \,\cite[Lemma~3.2]{N-bb},
\cite[Theorem~4.22(2)]{Rou}.
 If $\sC$ is a set of compact generators for $\sD$, then an object
of $\sD$ is compact if and only if it belongs to the minimal thick
subcategory of $\sD$ containing~$\sC$ \,\cite[Lemma~2.2]{N-c},
\cite[Theorem~4.22(3)]{Rou}.

 Let $X$ be a scheme with an open covering~$\bW$.
 When $X$ is Noetherian, we denote by $X\coh$ the abelian category
of coherent sheaves on~$X$.

\begin{thm}  \label{contraderived-compactly-generated}
\textup{(a)} For any scheme $X$, the coderived category\/
$\sD^\co(X\qcoh)$ admits arbitrary infinite direct sums, while
the contraderived categories\/ $\sD^\ctr(X\lcth_\bW^\lct)$ and\/
$\sD^\ctr(X\lcth_\bW)$ admit infinite products. 
 Becker's coderived and contraderived categories have similar
properties: the triangulated category $\sD^\bco(X\qcoh)$ admits
arbitrary infinite direct sums, while the triangulated categories\/
$\sD^\bctr(X\lcth_\bW^\lct)$ and\/ $\sD^\bctr(X\lcth_\bW)$ admit
infinite products. \par
\textup{(b)} For any Noetherian scheme $X$, the coderived category\/
$\sD^{\co=\bco}(X\qcoh)$ is compactly generated.
 The triangulated functor\/ $\sD^\b(X\coh)\rarrow\sD^\co(X\qcoh)$
induced by the embedding of abelian categories $X\coh\rarrow X\qcoh$
is fully faithful, and its image is the full subcategory of compact
objects in\/ $\sD^\co(X\qcoh)$. \par
\textup{(c)} For any Noetherian scheme $X$, the contraderived
category\/ $\sD^{\ctr=\bctr}(X\lcth_\bW^\lct)$ admits arbitrary
infinite direct sums and is compactly generated. \par
\textup{(d)} For any Noetherian scheme $X$ of finite Krull dimension,
the contraderived category\/ $\sD^{\ctr=\bctr}(X\lcth_\bW)$ admits
arbitrary infinite direct sums and is compactly generated. \par
\textup{(e)} For any semi-separated Noetherian scheme $X$,
the Becker contraderived category\/ $\sD^\bctr(X\lcth_\bW)$ admits
arbitrary infinite direct sums and is compactly generated.
\end{thm}

\begin{proof}
 Part~(a) holds, because the abelian category $X\qcoh$ admits arbitrary
infinite direct sums and the full subcategory of coacyclic complexes
in $\Hot(X\qcoh)$ is closed under infinite direct sums
(see~\cite[Lemma~1.5]{BN} or~\cite[Proposition~1.2.1 and
Lemma~3.2.10]{N-tr}).
 Analogously, the exact categories $X\lcth_\bW^\lct$ and $X\lcth_\bW$
admit arbitrary infinite products and the full subcategories of
contraacyclic complexes in $\Hot(X\lcth_\bW^\lct)$ and
$\Hot(X\lcth_\bW)$ are closed under infinite products.
 The argument for Becker's coderived and contraderived categories is
the same; one only needs to observe that, for any exact category $\sE$
with infinite direct sums, the class $\Acycl^\bco(\sE)$ is closed under
infinite direct sums in $\Hot(\sE)$.

 Under the assumption of part~(b), the abelian category $X\qcoh$ is
a locally Noetherian Grothendieck category, so the assertions
hold by Theorem~\ref{derived-inj-proj-resolutions}(a)
and~\cite[Proposition~2.3]{K-st} (see also
Lemma~\ref{co-contra-bounded-fully-faithful}).
 A more generally applicable assertion/argument can be found
in~\cite[Proposition~1.5(d)]{EP}
and/or~\cite[Section~3.11]{Pkoszul}.
 For an even more general version, applicable to coherent schemes,
see~\cite[Proposition~9.38 and Theorem~9.39]{Pedg}
or~\cite[Theorem~8.19 and Corollary~9.3]{PS5}.

 Under the assumptions of part~(d), all the triangulated categories
$\sD^\ctr(X\lcth_\bW^\lct)$ and $\sD^\ctr(X\lcth_\bW)$ are equivalent
to each other by Corollaries~\ref{finite-krull-ctrh-lcth-derived}
and~\ref{derived-contra-lct-cor}(b).
 The same applies to Becker's contraderived categories; moreover,
one has $\sD^\ctr(X\lcth_\bW^\lct)=\sD^\bctr(X\lcth_\bW^\lct)$ and
$\sD^\ctr(X\lcth_\bW)=\sD^\bctr(X\lcth_\bW)$ by
Theorem~\ref{derived-inj-proj-resolutions}(b,d).
 So part~(d) follows from part~(c).
 Furthermore, if the scheme $X$ admits a dualizing complex,
then the assertions of parts~(c\+d) follow from
Theorem~\ref{non-semi-separated-co-contra} and part~(b).

 Similarly, in the context of part~(e), one has
$\sD^\bctr(X\lcth_\bW)\simeq\sD^\bctr(X\lcth_\bW^\lct)$ by
Theorem~\ref{becker-contraderived-lcta-lct-equivalent}, so part~(e)
also follows from part~(c).
 In fact, in this case one has $\sD^\bctr(X\lcth_\bW^\lct)\simeq
\sD(X\qcoh_\fl)$ by
Corollary~\ref{bctr-lcth-vfl-fl-derived-equivalences}, so the assertion
of part~(e) is covered by the result of~\cite[Theorem~4.10]{M-th}.

 Part~(c): one has $\sD^\ctr(X\lcth_\bW^\lct)=
\sD^\bctr(X\lcth_\bW^\lct)$ by
Theorem~\ref{derived-inj-proj-resolutions}(b).
 The triangulated category $\sD^{\ctr=\bctr}(X\lcth_\bW^\lct)$
does not depend on $\bW$ by
Corollary~\ref{loc-noetherian-ctrh-lcth-contraderived}.
 The following somewhat complicated argument based on the results of
the paper~\cite{Rou} and the dissertation~\cite{M-th} allows to prove
part~(c) in the stated generality.

 By Theorem~\ref{derived-inj-proj-resolutions}(b), the triangulated
category in question is equivalent to the homotopy category
$\Hot(X\ctrh^\lct_\prj)$.
 Our argument proceeds in the typical fashion: first one covers
a semi-separated Noetherian scheme by affine open subschemes with
affine intersections; then one covers an arbitrary Noetherian scheme
by semi-separated open subschemes with semi-separated intersections.
 Let us first consider the case when the scheme $X$ is semi-separated.
 Then Corollary~\ref{vfl-co-contra-cor-expanded}(c) identifies our
triangulated category with $\sD^\bco(X\qcoh_\fl)=\sD(X\qcoh_\fl)$.
 It follows immediately that this triangulated category admits
arbitrary infinite direct sums.

 In the case of an affine Noetherian scheme $U$ of finite Krull
dimension, another application of Proposition~\ref{finite-resolutions}
allows to identify $\sD^\co(U\qcoh_\fl)$ with the homotopy category
of complexes of projective $\O(U)$\+modules, which is compactly
generated by~\cite[Theorem~2.4]{Jorg}.
 More generally, the category $\sD(U\qcoh_\fl)$ is equivalent to
the homotopy category of projective $\O(U)$\+modules for any affine
scheme $U$ by Theorem~\ref{flat-projective-periodicity}(b)
and Proposition~\ref{flat-projective-periodicity-complements}(b),
and $\sD(U\qcoh_\fl)$ is compactly generated for any affine Noetherian
scheme $U$ by~\cite[Proposition~7.14]{N-f} (see also~\cite{N-a}).
 Besides, the homotopy category $\Hot(U\ctrh^\lct_\prj)$ is equivalent
to $\sD(U\qcoh_\fl)$ for any affine scheme $U$
by~\cite[Corollary~5.8]{Sto}, \cite[Theorem~7.18]{Pphil},
or Corollary~\ref{becker-coderived-of-alf-well-behaved} above.
 Finally, for any semi-separated Noetherian scheme $X$
the triangulated category $\sD(X\qcoh_\fl)$ is compactly generated
by~\cite[Theorem~4.10]{M-th}.

 Now let us turn to the general case.
 First we have to show that the category $\Hot(X\ctrh^\lct_\prj)$
admits arbitrary infinite direct sums.
 Let $X=\bigcup_{\alpha=1}^N U_\alpha$ be a finite affine open covering,
and let $S_\beta\sub X$ denote the set-theoretic complement to
$\bigcup_{\alpha<\beta}U_\alpha$ in~$U_\beta$.
 Let $j_\alpha\:U_\alpha\rarrow X$ denote the open embedding morphisms;
then the direct image functor $j_\alpha{}_!\:
\Hot(U_\alpha\ctrh^\lct_\prj)\rarrow\Hot(X\ctrh^\lct_\prj)$ is left
adjoint to the inverse image functor $j_\alpha^!\:
\Hot(X\ctrh^\lct_\prj)\rarrow\Hot(U_\alpha\ctrh^\lct_\prj)$
(see Corollaries~\ref{lct-proj-local}(a) and~\ref{proj-lct-direct}(b),
and the adjunction~\eqref{direct-inverse-cosheaf-lct-adjunction}).
 Hence the functor~$j_\alpha{}_!$ preserves infinite direct sums.

 As explained in the proof of Theorem~\ref{proj-lct-classification},
any projective locally cotorsion contraherent cosheaf $\gF$ on $X$
decomposes into a direct sum $\gF=\bigoplus_{\alpha=1}^N\gF_\alpha$,
where each direct summand $\gF_\alpha$ is an infinite product over
the points $z\in S_\alpha$ of the direct images of contraherent
cosheaves on $\Spec\O_{z,X}$ corresponding to free contramodules
over $\widehat\O_{z,X}$.
 According to Lemma~\ref{proj-lct-orthogonality}, the associated
increasing filtration $\gF_{\le\alpha} = \bigoplus_{\beta\le\alpha}
\gF_\beta$ on $\gF$ is preserved by all morphisms of cosheaves
$\gF\in X\ctrh^\lct_\prj$.

 Given a family ${}^{(i)}\gF^\bu$ of complexes over $X\ctrh^\lct_\prj$,
we now see that every complex ${}^{(i)}\gF^\bu$ is endowed with
a finite termwise split filtration ${}^{(i)}\gF^\bu_{\le\alpha}$ such
that the family of associated quotient complexes
${}^{(i)}\gF^\bu_\alpha$ can be obtained by applying the direct image
functor~$j_\alpha{}_!$ to a family of complexes over
$U_\alpha\ctrh^\lct_\prj$.
 It follows that the object $\bigoplus_i {}^{(i)}\gF^\bu_\alpha$
exists in $\Hot(X\ctrh^\lct_\prj)$, and it remains to apply
the following lemma (which is slightly stronger
than~\cite[Proposition~1.2.1]{N-tr}).

\begin{lem}
 Let $A_i\rarrow B_i\rarrow C_i\rarrow A_i[1]$ be a family of
distinguished triangles in a triangulated category\/~$\sD$.
 Suppose that the infinite direct sums\/ $\bigoplus_i A_i$ and\/
$\bigoplus_i B_i$ exist in\/~$\sD$.
 Then a cone $C$ of the natural morphism\/ $\bigoplus_i A_i\rarrow
\bigoplus_i B_i$ is the infinite direct sum of the family of
objects $C_i$ in\/~$\sD$.
\end{lem}

\begin{proof}
 Set $A=\bigoplus_i A_i$ and $B=\bigoplus_i B_i$.
 By one of the triangulated category axioms, there exist morphisms
of distinguished triangles $(A_i\to B_i\to C_i\to A_i[1])\rarrow
(A\to B\to C\to A[1])$ whose components $A_i\rarrow A$ and
$B_i\rarrow B$ are the natural embeddings.
 For any object $E\in\sD$, apply the functor $\Hom_\sD({-},E)$
to this family of morphisms of triangles and pass to the infinite
product (of abelian groups) over~$i$.
 The resulting morphism from the long exact sequence
$\dotsb\rarrow\Hom_\sD(A[1],E)\rarrow\Hom_\sD(C,E)\rarrow
\Hom_\sD(B,E)\rarrow\Hom_\sD(A,E)\rarrow\dotsb$
to the long exact sequence
$\dotsb\rarrow\prod_i\Hom_\sD(A_i[1],E)\rarrow
\prod_i\Hom_\sD(C_i,E)\rarrow\prod_i\Hom_\sD(B_i,E)\rarrow
\prod_i\Hom_\sD(A_i,E)\rarrow\dotsb$
is an isomorphism at the two thirds of all the terms, and
consequently an isomorphism at the remaining terms, too.
\end{proof}

 Denote temporarily the homotopy category $\Hot(X\ctrh^\lct_\prj)$
by $\sD(X)$.
 To show that the category $\sD(X)$ is compactly generated, we will
use the result of~\cite[Theorem~5.15]{Rou}.
 Let $Y\sub X$ be an open subscheme such that the category $\sD(Y)$ is
compactly generated (e.~g., we already know this to hold when $Y$ is
semi-separated).
 Let $j\:Y\rarrow X$ denote the open embedding morphism.

 The composition $j^!j_!$ of the direct image and inverse image
functors $j_!\:\sD(Y)\rarrow\sD(X)$ and $j^!\:\sD(X)\rarrow\sD(Y)$
is isomorphic to the identity endofunctor of $\sD(Y)$, so
the functor~$j_!$ is fully faithful and the functor~$j^!$
is a Verdier localization functor.
 Applying again Lemma~\ref{proj-lct-orthogonality}, we conclude
that the kernel of~$j^!$ is the homotopy category of projective
locally cotorsion contraherent cosheaves on $X$ with vanishing
restrictions to~$Y$.
 Denote this homotopy category by $\sD(Z,X)$, where $Z=X\setminus Y$,
and its identity embedding functor by $i_!\:\sD(Z,X)\rarrow\sD(X)$.

 The functor~$j_!$ is known to preserve infinite products, and
the triangulated category $\sD(Y)$ is assumed to be compactly
generated; so it follows that there exists a triangulated functor
$j^*\:\sD(X)\rarrow\sD(Y)$ left adjoint to~$j_!$
(see~\cite[Remark~6.4.5 and Theorem~8.6.1]{N-tr}
and~\cite[Proposition~3.3(2)]{K-st}).
 The existence of the functor~$j_!$ left adjoint to~$j^!$ implies
existence of a functor $i^*\:\sD(X)\rarrow\sD(Z,X)$ left adjoint
to~$i_!$; and the existence of the functor~$j^*$ left adjoint to~$j_!$
implies existence of a functor $i_+\:\sD(Z,X)\rarrow\sD(X)$ left
adjoint to~$i^*$.

 The functors $j^*$ and~$i_+$ have double right adjoints (i.~e.,
the right adjoints and the right adjoints to the right adjoints),
hence they not only preserve infinite direct sums, but also take
compact objects to compact objects.
 Furthermore, for any open subscheme $W\sub X$ with the embedding
morphism $h\:W\rarrow X$ one has the base change isomorphism
$h^!j_!\simeq j'_!h'{}^!$, where $j'$ and~$h'$ denote the open
embeddings $W\cap Y\rarrow W$ and $W\cap Y\rarrow Y$.
 If the triangulated category $\sD(W\cap Y)$ is compactly generated,
one can pass to the left adjoint functors, obtaining an isomorphism
of triangulated functors $j^*h_!\simeq h'_!\.j'{}^*$.

 Let $X=\bigcup_\alpha U_\alpha$ be a finite affine (or, more
generally, semi-separated) open covering, $Z_\alpha =
X\setminus U_\alpha$ be the corresponding closed complements,
and $i_\alpha{}_+\:\sD(Z_\alpha,X)\rarrow\sD(X)$ be
the related fully faithful triangulated functors.
 It follows from the above that the images of the functors
$i_\alpha{}_+$ form a collection of Bousfield subcategories in
$\sD(X)$ pairwise \emph{intersecting properly} in the sense
of~\cite[Lemma~5.7(2)]{Rou} or~\cite[Lemma~2.12(iii)]{M-th}.
 Indeed, in the notation of the previous paragraph one has
$j^*h_!h^*\P^\bu\simeq h'_!\,j'{}^*h^*\P^\bu\simeq
h'_!h'{}^*j^*\P^\bu=0$ for any object $\P^\bu\in\sD(X)$ such that
$j^*\P^\bu=0$ in $\sD(Y)$.
 Furthermore, the category $\sD(X)$ being generated by the images
of the functors $j_\alpha{}_!$, the intersection of the kernels of
the localization functors $j_\alpha^*\:\sD(X)\rarrow\sD(U_\alpha)$
is zero.
 These kernels coincide with the images of the functors~$i_\alpha{}_+$.
 Thus the triangulated subcategories $i_\alpha{}_+\sD(Z_\alpha,X)
\subset\sD(X)$ form a \emph{cocovering} (in the sense
of~\cite[Section~5.3.3]{Rou}).

 In order to satisfy the assumptions of~\cite[Theorem~5.15]{Rou}
or~\cite[Theorem~2.13]{M-th}, it remains to check that intersections
of the images of $i_\beta{}_+\sD(Z_\beta,X)\sub\sD(X)$ under
the localization functor $\sD(X)\rarrow\sD(U_\alpha)$ are compactly
generated in $\sD(U_\alpha)$ (in the sense
of~\cite[Section~3.3.1]{Rou}).
 Now for any subset of indices $\beta_1$,~\dots, $\beta_k$,
the intersection of the kernels of the localization functors
$j_{\beta_s}^*\:\sD(X)\rarrow\sD(U_{\beta_s})$, taken over all
$s=1$,~\dots, $k$, is equal to the kernel of the localization functor
$j^*\:\sD(X)\rarrow\sD(U_{\beta_1}\cup\dotsb\cup U_{\beta_k})$.
 The latter kernel coincides with the image of the fully faithful
triangulated functor $i_+\:\sD(Z,X)\rarrow\sD(X)$, where
$Z=X\setminus(U_{\beta_1}\cup\dotsb\cup U_{\beta_k})$.

 Furthermore, for any index~$\beta$, the image of
$i_\beta{}_+\sD(Z_\beta,X)\sub\sD(X)$ under the localization functor
$\sD(X)\rarrow\sD(U_\alpha)$ coincides with the kernel of
the localization functor $\sD(U_\alpha)\rarrow
\sD(U_\alpha\cap U_\beta)$.
 Thus the intersection of the images of
$i_{\beta_s}{}_+\sD(Z_{\beta_s},X)$ under the localization functor
$\sD(X)\rarrow\sD(U_\alpha)$ coincides with the kernel of
the localization functor $\sD(U_\alpha)\rarrow
\sD(U_\alpha\cap(U_{\beta_1}\cup\dotsb\cup U_{\beta_k}))$.

 Let $Y$ be a semi-separated Noetherian scheme of finite Krull
dimension and $V\sub Y$ be an open subscheme with the closed
complement $Z=Y\setminus V$.
 Let $j\:V\rarrow Y$ and $i\:Z\rarrow Y$ be the related morphisms.
 We will show that the image of the fully faithful triangulated
functor $i_+\:\sD(Z,Y)\rarrow\sD(Y)$ is compactly generated in
$\sD(Y)$; this is clearly sufficient.

 The result of Corollary~\ref{vfl-co-contra-cor-expanded}(c)
identifies $\sD(Y)=\Hot(Y\ctrh^\lct_\prj)$ with $\sD(Y\qcoh_\fl)$
and $\sD(V)=\Hot(V\ctrh^\lct_\prj)$ with $\sD(V\qcoh_\fl)$.
 According to Corollary~\ref{inj-vfl-direct-images-identified}(d),
this identification transforms the functor $j_!\:\sD(V)\rarrow\sD(Y)$
into the derived functor $\boR j_*\:\sD(V\qcoh_\fl)\rarrow
\sD(Y\qcoh_\fl)$ constructed in~\eqref{qcoh-direct-ffd}.
 The latter functor is right adjoint to the functor
$j^*\:\sD(Y\qcoh_\fl)\rarrow\sD(V\qcoh_\fl)$, which is therefore
identified with the functor $j^*\:\sD(Y)\rarrow\sD(V)$.

 Finally, we refer to~\cite[Proposition~4.5 and Theorem~4.10]{M-th}
for the assertion that the kernel of the functor
$j^*\:\sD(Y\qcoh_\fl)\rarrow\sD(V\qcoh_\fl)$ is compactly
generated in $\sD(Y\qcoh_\fl)$.
\end{proof}

\begin{thm}  \label{derived-compactly-generated} \hbadness=1500
\textup{(a)} For any scheme $X$, the derived category\/ $\sD(X\qcoh)$
admits infinite direct sums, while the derived categories\/
$\sD(X\lcth_\bW^\lct)$ and\/ $\sD(X\lcth_\bW)$ admit infinite products.
\par
\textup{(b)} For any quasi-compact semi-separated scheme $X$,
the derived category\/ $\sD(X\qcoh)$ is compactly generated.
 The full triangulated subcategory of perfect complexes in\/
$\sD^\b(X\qcoh_\vfl)\sub\sD^\b(X\qcoh_\fl)\sub\sD^\b(X\qcoh)\sub
\sD(X\qcoh)$ is the full subcategory of compact objects in\/
$\sD(X\qcoh)$. \par
\textup{(c)} For any Noetherian scheme $X$, the derived category\/
$\sD(X\qcoh)$ is compactly generated.
 The full triangulated subcategory of perfect complexes in\/
$\sD^\b(X\coh)\sub\sD^\b(X\qcoh)\sub\sD(X\qcoh)$ is the full
subcategory of compact objects in\/ $\sD(X\qcoh)$. \par
\textup{(d)} For any quasi-compact semi-separated scheme $X$,
the derived category\/ $\sD(X\lcth_\bW)$ admits infinite direct sums
and is compactly generated. \par
\textup{(e)} For any Noetherian scheme $X$ of finite Krull dimension,
the derived categories\/ $\sD(X\lcth_\bW^\lct)$ and\/ $\sD(X\lcth_\bW)$
admit infinite direct sums and are compactly generated.
\end{thm}

\begin{proof}
 The proof of part~(a) is similar to that of
Theorem~\ref{contraderived-compactly-generated}(a):
the assersions hold, since the class of acyclic complexes over
$X\qcoh$ is closed under infinite direct sums, and the classes of
acyclic complexes over $X\lcth_\bW^\lct$ and $X\lcth_\bW$ are
closed under infinite products.

 Parts~(b) and~(c) are particular cases of~\cite[Theorem~6.8]{Rou},
according to which the derived category $\sD(X)$ of complexes of
sheaves of $\O_X$\+modules with quasi-coherent cohomology sheaves is
compactly generated for any quasi-compact quasi-separated scheme~$X$.
 Here one needs to know that the natural functor $\sD(X\qcoh)\rarrow
\sD(X)$ is an equivalence of categories when $X$ is either
quasi-compact and semi-separated, or else Noetherian
(cf.~\cite[Appendix~B]{TT}).
 In the semi-separated case, this was proved
in~\cite[Sections~5\+-6]{BN}.
 The proof in the Noetherian case is similar.
 Both proofs can be found in~\cite[Propositions Tags~08DB and~09T4]{SP}.

 Alternatively, one can prove parts~(b) and~(c) directly in the way
analogous to the argument in~\cite{Rou}.
 In either approach, one needs to know that the functor $\boR j_*$
of derived direct image of complexes over $Y\qcoh$ with respect to
an open embedding $j\:Y\rarrow X$ of schemes from the class under
consideration is well-behaved.
 E.~g., it needs to be local in the base, or form a commutative
square with the derived functor of direct image of complexes of
$\O_Y$\+modules, etc.

 In the semi-separated case, one can establish such properties using
contraadjusted coresolutions and (the proof of)
Corollary~\ref{cta-cot-direct}(a) (see the construction of the functor
$\boR f_*$ in Section~\ref{derived-direct-inverse} above).
 Better yet, dilute coresolutions and
Corollary~\ref{dilute-direct} can be used.
 The exposition in~\cite[Theorems~31 and~42]{M-n} can be recommended.

 In the Noetherian case, one needs to use flasque coresolutions and
Corollary~\ref{coflasque-complexes-direct-correct}(a)
(see the construction of the functor $\boR f_*$ in
Section~\ref{derived-direct-special} below).
 The point is that, even though the notion of a homotopy injective
complex of (injective) quasi-coherent sheaves is not local on
Noetherian schemes, one can use the local notion of a complex of
flasque quasi-coherent sheaves instead.
 Even though the flasque coresolution dimension of a quasi-coherent
sheaf on $X$ may be infinite if the Krull dimension of $X$ is infinite,
all complexes of flasque quasi-coherent sheaves on Noetherian schemes
are still adjusted to direct images by
Corollary~\ref{coflasque-complexes-direct-correct}(a).
 See also the second assertion of
Corollary~\ref{coflasque-resolutions-infinite}(a) below.

 Part~(d) follows from part~(b) together with
Theorem~\ref{naive-co-contra-thm}.
 Part~(e) follows from part~(c) together with
Theorem~\ref{non-semi-separated-naive-co-contra}
and Corollary~\ref{derived-contra-lct-cor}(b).
\end{proof}

\subsection{Homotopy projective complexes on Noetherian schemes}
 Let $X$ be a scheme with an open covering~$\bW$.
 We refer to Section~\ref{qc-ss-homotopy-projective-subsect} for
the definitions of homotopy injective complexes of quasi-coherent
sheaves on $X$ and homotopy projective complexes of (locally cotorsion)
$\bW$\+locally contraherent cosheaves on~$X$.
 The related notation $\Hot(X\lcth_\bW)_\prj\sub\Hot(X\lcth_\bW)$
and $\Hot(X\lcth_\bW^\lct)_\prj\sub\Hot(X\lcth_\bW^\lct)$ was
also defined in Section~\ref{qc-ss-homotopy-projective-subsect}.

 We will see below in this section that the property of a complex of
$\bW$\+locally contraherent cosheaves on a Noetherian scheme $X$ of
finite Krull dimension to be homotopy projective does not change when
the covering $\bW$ is replaced by its refinement. 
 Similarly, the property of a complex of locally cotorsion
$\bW$\+locally contraherent cosheaves on a Noetherian scheme $X$ of
finite Krull dimension to be homotopy projective does not change
when the covering $\bW$ is changed.

 The following lemma is a non-semi-separated (but locally Noetherian)
version of Lemma~\ref{qcomp-ssep-homotopy-proj-independence}.

\begin{lem}  \label{homotopy-proj-independence}
\textup{(a)} Let $X$ be a locally Noetherian scheme with an open
covering\/~$\bW$.
 Then a complex\/ $\P^\bu$ over $X\ctrh^\lct_\prj$ belongs to\/
$\Hot(X\lcth_\bW^\lct)_\prj$ if and only if the complex\/
$\Hom^X(\P^\bu,\gE^\bu)$ is acyclic for any complex\/ $\gE^\bu$ over
$X\ctrh^\lct_\prj$ acyclic with respect to $X\ctrh^\lct$.
 If the scheme $X$ has finite Krull dimension, then the complex\/
$\gE^\bu$ can be assumed to be acyclic with respect to
$X\ctrh^\lct_\cfq$. \par
\textup{(b)} Let $X$ be a Noetherian scheme of finite Krull dimension
with an open covering\/~$\bW$.
 Then a complex\/ $\gF^\bu$ over $X\ctrh_\prj$ belongs to\/
$\Hot(X\lcth_\bW)_\prj$ if and only if the complex\/
$\Hom^X(\gF^\bu,\gE^\bu)$ is acyclic for any complex\/ $\gE^\bu$ over
$X\ctrh_\prj$ acyclic with respect to $X\ctrh_\cfq$.
\end{lem}

\begin{proof}
 We will prove part~(a), part~(b) being similar.
 The ``only if'' assertion holds by the definition.
 To check the ``if'', consider a complex $\gM^\bu$ over
$X\lcth_\bW^\lct$.
 By (the proof of) Theorem~\ref{derived-inj-proj-resolutions}(b),
there exists a complex $\gE^\bu$ over $X\ctrh^\lct_\prj$ together
with a morphism of complexes of locally contraherent cosheaves
$\gE^\bu\rarrow\gM^\bu$ with a cone contraacyclic with respect
to $X\lcth_\bW^\lct$.
 Moreover, the complex $\Hom^X$ from any complex of projective locally
cotorsion contraherent cosheaves to a contraacyclic complex over
$X\lcth_\bW^\lct$ is acyclic.
 Hence the morphism $\Hom^X(\P^\bu,\gE^\bu)\rarrow\Hom^X(\P^\bu,
\gM^\bu)$ is a quasi-isomorphism.
 Finally, if the complex $\gM^\bu$ is acyclic over $X\lcth_\bW^\lct$,
then so is the complex $\gE^\bu$, and by
Lemma~\ref{acyclicity-in-lcth-criterion}(b) it follows that
the complex $\gE^\bu$ is also acyclic with respect to $X\ctrh^\lct$.
 If the scheme $X$ has finite Krull dimension, then
Corollary~\ref{lct-proj-local}(b) and
Lemma~\ref{finite-krull-ctrh-lcth-finite-dim}(b) with
Corollary~\ref{fdim-acyclic-cor} tell us that the complex $\gE^\bu$ is
acyclic with respect to $X\ctrh^\lct_\cfq$.
\end{proof}

 According to Lemma~\ref{homotopy-proj-independence}, the property
of a complex over $X\ctrh^\lct_\prj$ (respectively, over
$X\ctrh_\prj$) to belong to $\Hot(X\lcth_\bW^\lct)_\prj$
(resp., $\Hot(X\lcth_\bW)_\prj$) does not depend on the covering~$\bW$
(in the assumptions of the respective part of the lemma).
 We will denote the full subcategory in $\Hot(X\ctrh^\lct_\prj)$
(resp., $\Hot(X\ctrh_\prj)$) consisting of the homotopy projective
complexes by $\Hot(X\ctrh^\lct_\prj)_\prj$ (resp.,
$\Hot(X\ctrh_\prj)_\prj$).
 It is a standard fact that bounded above complexes of projectives are
homotopy projective, $\Hot^-(X\ctrh^\lct_\prj)\subset
\Hot(X\ctrh^\lct_\prj)_\prj$ and $\Hot^-(X\ctrh_\prj)\subset
\Hot(X\ctrh_\prj)_\prj$.

 We refer to Section~\ref{qcomp-ssep-homotopy-proj-independence}
for the notation $\Com(X\ctrh_\prj)_\prj\sub\Com(X\ctrh)$
and $\Com(X\ctrh^\lct_\prj)_\prj\sub\Com(X\ctrh^\lct)$.
 The full subcategories $\Com(X\ctrh_\prj)_\prj$ and
$\Com(X\ctrh^\lct_\prj)_\prj$ are closed under extensions and kernels
of admissible epimorphisms in the ambient exact categories
$\Com(X\ctrh)$ and $\Com(X\ctrh^\lct)$, respectively.
 (Notice that any such extension or admissible epimorphism of complexes
belonging to these full subcategories is termwise split.)

 Furthermore, the class of homotopy projective complexes of projective
contraherent cosheaves is preserved by the direct images with respect to
very flat morphisms of Noetherian schemes of finite Krull dimension.
 Similarly, the class of homotopy projective complexes of projective
locally cotorsion contraherent cosheaves is preserved by the direct
images with respect to flat quasi-compact morphisms of locally
Noetherian schemes.
 These assertions follow from
Corollaries~\ref{finite-krull-flat-direct}(c)
and~\ref{proj-lct-direct}(b),
and the adjunction~\eqref{direct-inverse-cosheaf-lct-adjunction}.

 The following theorem is a non-semi-separated (but Noetherian)
version of Corollary~\ref{qcomp-ssep-homotopy-projective}.

\begin{thm}  \label{finite-krull-homotopy-projective}
 Let $X$ be a Noetherian scheme with an open covering\/~$\bW$.
 In this setting: \par
\textup{(a)} the natural functors\/
$\Hot(X\ctrh^\lct_\prj)_\prj\rarrow\Hot(X\lcth_\bW^\lct)_\prj\rarrow
\sD(X\lcth^\lct_\bW)$ are equivalences of triangulated categories; \par
\textup{(b)} if $X$ has finite Krull dimension, then the natural
functors\/ $\Hot(X\ctrh_\prj)_\prj\rarrow\Hot(X\lcth_\bW)_\prj\rarrow
\sD(X\lcth_\bW)$ are equivalences of triangulated categories.
\end{thm}

\begin{proof}[First proof]
 We will prove part~(b), part~(a) being similar.
 Since both functors are clearly fully faithful, it suffices to show
that the composition $\Hot(X\ctrh_\prj)_\prj\rarrow\sD(X\lcth_\bW)$
is an equivalence of categories.
 This is equivalent to saying that the localization functor
$\Hot(X\lcth_\bW)\rarrow\sD(X\lcth_\bW)$ has a left adjoint
whose image is essentially contained in $\Hot(X\ctrh_\prj)$.

 The functor in question factorizes into the composition of two
localization functors $\Hot(X\lcth_\bW)\rarrow\sD^\ctr(X\lcth_\bW)
\rarrow\sD(X\lcth_\bW)$.
 The functor $\Hot(X\lcth_\bW)\allowbreak\rarrow\sD^\ctr(X\lcth_\bW)$
has a left adjoint functor $\sD^\ctr(X\lcth_\bW)\simeq\Hot(X\ctrh_\prj)
\hookrightarrow\Hot(X\lcth_\bW)$ provided by
Theorem~\ref{derived-inj-proj-resolutions}(d); so it remains to show
that the functor $\sD^\ctr(X\lcth_\bW)\rarrow\sD(X\lcth_\bW)$ has
a left adjoint.

 Since the latter functor preserves infinite products, the assertion
follows from Theorem~\ref{contraderived-compactly-generated}(d)
and the covariant Brown representability
theorem~\cite[Proposition~3.3(2)]{K-st}.
 Here one also needs to know that the derived category
$\sD(X\lcth_\bW^\lct)$ or $\sD(X\lcth_\bW)$ ``exists'' (i.~e., morphisms
between any given two objects form a set rather than a proper class).
 This observation was, of course, already presumed in
Theorem~\ref{derived-compactly-generated}(e); but it holds in
greater generality.

 The fact that the Hom sets are small can be established by noticing
that the classes of quasi-isomorphisms are locally small in
$\Hot(X\lcth_\bW^\lct)$ and $\Hot(X\lcth_\bW)$
(see~\cite[Section~10.3.6 and Proposition~10.4.4]{Wei}).
 Specifically, for any scheme $X$ with an open covering $\bW$,
the result of~\cite[Corollary~10.4]{PS4} is applicable, proving
``existence'' of the derived categories $\sD(X\lcth_\bW^\lct)$
and $\sD(X\lcth_\bW)$.
 Alternatively, one can use
Theorem~\ref{non-semi-separated-naive-co-contra} in order to
reduce the question to ``existence'' of the derived category
$\sD(X\qcoh)$, which follows, e.~g., from
Theorem~\ref{quasi-homotopy-injective}.

 Similarly one can prove Theorem~\ref{quasi-homotopy-injective} for
a Noetherian scheme $X$ using
Theorems~\ref{derived-inj-proj-resolutions}(a)
and~\ref{contraderived-compactly-generated}(b); the only difference
is that this time one needs a right adjoint functor,
so the contravariant Brown
representability~\cite[Proposition~3.3(1)]{K-st} has to be applied.
 One again, with this approach one needs to know that the derived
category $\sD(X\qcoh)$ ``exists''.
 (Cf.~\cite[Section~5.5]{Pkoszul}).
\end{proof}

\begin{proof}[Second proof]
 Part~(b): it suffices to construct, for every complex $\gM^\bu$ over
$X\lcth_\bW$, a complex $\P^\bu\in\Hot(X\ctrh_\prj)_\prj$ together
with a quasi-isomorphism $\P^\bu\rarrow\gM^\bu$ of complexes in
$X\lcth_\bW$.
 To begin with, by (the proof of)
Corollary~\ref{coflasque-resolutions-finite}(c), there exists
a complex of coflasque contraherent cosheaves $\gE^\bu$ on $X$
together with a quasi-isomorphism $\gE^\bu\rarrow\gM^\bu$ of
complexes in $X\lcth_\bW$.
 Now the desired complex $\P^\bu\in\Hot(X\ctrh_\prj)_\prj$ together
with a quasi-isomorphism $\P^\bu\rarrow\gE^\bu$ of complexes in
$X\ctrh_\cfq$ is provided by
Lemma~\ref{coflasque-homotopy-projective-precover}(b) below.

 The proof of part~(a) is similar.
 It suffices to construct, for every complex $\gM^\bu$ over
$X\lcth^\lct_\bW$, a complex $\P^\bu\in\Hot(X\ctrh^\lct_\prj)_\prj$
together with a quasi-isomorphism $\P^\bu\rarrow\gM^\bu$ of complexes
in $X\lcth^\lct_\bW$.
 By Theorem~\ref{derived-inj-proj-resolutions}(b), there exists
a complex of projective locally cotorsion contraherent cosheaves
$\gE^\bu$ on $X$ together with a morphism of complexes $\gE^\bu
\rarrow\gM^\bu$ with a cone contraacyclic in $X\lcth^\lct_\bW$.
 All Positselski-contraacyclic complexes are obviously acyclic;
so $\gE^\bu\rarrow\gM^\bu$ is a quasi-isomorphism in $X\lcth^\lct_\bW$.
 By Corollary~\ref{lct-proj-local}(b), \,$\gE^\bu$ is a complex of
coflasque locally cotorsion contraherent cosheaves on~$X$.
  Now the desired complex $\P^\bu\in\Hot(X\ctrh^\lct_\prj)_\prj$
together with a quasi-isomorphism $\P^\bu\rarrow\gE^\bu$ of complexes
in $X\ctrh^\lct$ is provided by
Lemma~\ref{coflasque-homotopy-projective-precover}(a) below.
\end{proof}

\begin{lem} \label{coflasque-acyclic-preenvelope}
 Let $X$ be a Noetherian scheme.  In this setting: \par
\textup{(a)} any complex\/ $\gE^\bu$ in the exact category
$X\ctrh^\lct_\cfq$ can be included in a short exact sequence\/
$0\rarrow\gE^\bu\rarrow\gA^\bu\rarrow\P^\bu\rarrow0$, where\/ $\gA^\bu$
is a complex in $X\ctrh^\lct_\cfq$ that is acyclic in $X\ctrh^\lct$,
and\/ $\P^\bu\in\Com(X\ctrh^\lct_\prj)_\prj$ is a homotopy projective
complex of projective locally cotorsion contraherent cosheaves; \par
\textup{(b)} if $X$ has finite Krull dimension, then any complex
$\gE^\bu$ in the exact category $X\ctrh_\cfq$ can be included in
a short exact sequence\/ $0\rarrow\gE^\bu\rarrow\gA^\bu\rarrow\P^\bu
\rarrow0$, where\/ $\gA^\bu$ is an acyclic complex in $X\ctrh_\cfq$
and\/ $\P^\bu\in\Com(X\ctrh_\prj)_\prj$ is a homotopy projective
complex of projective contraherent cosheaves.
\end{lem}

\begin{proof}
 The argument is similar to the proof of
Lemma~\ref{coflasque-lct-envelope}(a);
see also a more general discussion in the proofs in
Section~\ref{acyclic-complexes-of-flat-cosheaves-subsect}.
 Let us prove part~(b).
 Let $X=\bigcup_\alpha U_\alpha$ be a finite affine open covering of
the scheme~$X$.
 We proceed by induction in a linear ordering of the indices~$\alpha$,
considering the open subscheme $V=\bigcup_{\alpha<\beta}U_\alpha\sub X$.
 Assume that we have constructed a short exact sequence $0\rarrow
\gE^\bu\rarrow\gB^\bu\rarrow\gF^\bu\rarrow0$ of complexes of coflasque
contraherent cosheaves on $X$ such that the complex $\gB^\bu|_V$ is
acyclic in $V\ctrh_\cfq$, while $\gF^\bu\in\Com(X\ctrh_\prj)_\prj$.
 Let $j\:U=U_\beta\rarrow X$ be the identity open embedding.

 Pick a short exact sequence $0\rarrow j^!\gB^\bu\rarrow\gC^\bu
\rarrow\Q^\bu\rarrow0$ of complexes of contraherent cosheaves on $U$
such that the complex $\gC^\bu$ is acyclic in $U\ctrh$ and
$\Q^\bu\in\Com(U\ctrh_\prj)_\prj$
(see Theorem~\ref{lcth-homotopy-projective-cotorsion-pairs}(a)).
 Then $\Q^\bu$ is a complex of coflasque contraherent cosheaves on $U$
by Corollary~\ref{finite-krull-contrah-projective}(b) and $\gC^\bu$
is a complex of coflasque contraherent cosheaves on $U$ by
Corollary~\ref{coflasque-acyclic}(a).
 By Lemma~\ref{finite-krull-ctrh-lcth-finite-dim}(b) with
Corollary~\ref{fdim-acyclic-cor}, the complex $\gC^\bu$ is acyclic
in $U\ctrh_\cfq$.

 By Corollary~\ref{coflasque-direct}(a), the related sequence of direct
images $0\rarrow j_!j^!\gB^\bu\rarrow j_!\gC^\bu\rarrow j_!\Q^\bu
\rarrow0$ is a short exact sequence of complexes of coflasque
contraherent cosheaves on~$X$; by the same corollary, the complex
$j_!\gC^\bu$ is acyclic in $X\ctrh_\cfq$.
 According to the discussion above in this section, we have
$j_!\Q^\bu\in\Com(X\ctrh_\prj)_\prj$.

 Let $0\rarrow\gB^\bu\rarrow\gA^\bu\rarrow j_!\Q^\bu\rarrow0$
denote the push-out of the short exact sequence $0\rarrow
j_!j^!\gB^\bu\rarrow j_!\gC^\bu\rarrow j_!\Q^\bu\rarrow0$ with respect
to the natural (adjunction) morphism $j_!j^!\gB^\bu\rarrow\gB^\bu$.
 We will show that the complex of coflasque contraherent cosheaves
$\gA^\bu$ on $X$ is acyclic in restriction to $U\cup V$.

 Indeed, in the restriction to $U$ one has $j^!\gA^\bu\simeq\gC^\bu$.
 On the other hand, denoting by~$j'$ the open embedding
$U\cap V\rarrow V$, one has $(j_!\Q^\bu)|_V\simeq
j'_!(\Q^\bu|_{U\cap V})$.
 The complex of coflasque contraherent cosheaves $\gB^\bu|_{U\cap V}$
being acyclic, so is the cokernel $\Q^\bu|_{U\cap V}$ of the admissible 
monomorphism of acyclic complexes of coflasque contraherent cosheaves
$\gB^\bu|_{U\cap V}\rarrow\gC^\bu|_{U\cap V}$.

 By Corollary~\ref{coflasque-direct}(a), \,$j'_!(\Q^\bu|_{U\cap V})$ is
an acyclic complex of coflasque contraherent cosheaves on~$V$.
 Now in the short exact sequence of complexes of coflasque
contraherent cosheaves $0\rarrow\gB^\bu|_V\rarrow\gA^\bu|_V\rarrow
(j_!\Q^\bu)|_V\rarrow0$ the middle term is acyclic, since
the two other terms are.

 Finally, the composition $\gE^\bu\rarrow\gB^\bu\rarrow\gA^\bu$ of
admissible monomorphisms in $\Com(X\ctrh_\cfq)$ is again an admissible
monomorphism with the cokernel isomorphic to an extension of two
complexes $j_!\Q^\bu$ and $\gF^\bu$ belonging to
$\Com(X\ctrh_\prj)_\prj$.

 The proof of part~(a) is similar, except that
Lemma~\ref{finite-krull-ctrh-lcth-finite-dim} is not applicable,
since the Krull dimension of $X$ is not assumed to be finite.
 Consequently, we cannot obtain acyclic complexes in $X\ctrh^\lct_\cfq$,
and have to do with complexes in $X\ctrh^\lct_\cfq$ acyclic in
$X\ctrh^\lct$ instead.
 For the assertion that the class of such complexes is preserved by
the direct images with respect to quasi-compact quasi-separated
morphisms, we need to refer to
Corollary~\ref{coflasque-complexes-direct-correct}(c).
\end{proof}

\begin{lem} \label{coflasque-homotopy-projective-precover}
 Let $X$ be a Noetherian scheme.  In this setting: \par
\textup{(a)} any complex\/ $\gE^\bu$ in the exact category
$X\ctrh^\lct_\cfq$ can be included in a short exact sequence\/
$0\rarrow\gA^\bu\rarrow\P^\bu\rarrow\gE^\bu\rarrow0$, where\/ $\gA^\bu$
is a complex in $X\ctrh^\lct_\cfq$ that is acyclic in $X\ctrh^\lct$,
and\/ $\P^\bu\in\Com(X\ctrh^\lct_\prj)_\prj$ is a homotopy projective
complex of projective locally cotorsion contraherent cosheaves; \par
\textup{(b)} if $X$ has finite Krull dimension, then any complex
$\gE^\bu$ in the exact category $X\ctrh_\cfq$ can be included in
a short exact sequence\/ $0\rarrow\gA^\bu\rarrow\P^\bu\rarrow\gE^\bu
\rarrow0$, where\/ $\gA^\bu$ is an acyclic complex in $X\ctrh_\cfq$
and\/ $\P^\bu\in\Com(X\ctrh_\prj)_\prj$ is a homotopy projective
complex of projective contraherent cosheaves.
\end{lem}

\begin{proof}
 The proof is similar to that of
Lemma~\ref{flat-iterated-extension-covering}(a) and
based on Lemma~\ref{coflasque-acyclic-preenvelope}.
 Part~(b): following the proof of
Lemma~\ref{flat-iterated-extension-covering}(a), there exists
an admissible epimorphism in $X\ctrh_\cfq$ onto any object of
$X\ctrh_\cfq$ from a projective contraherent cosheaf.
 It follows easily that there exists an admissible epimorphism in
$\Com(X\ctrh_\cfq)$ onto an object of $\Com(X\ctrh_\cfq)$ from
a contractible (hence homotopy projective) complex of projective
contraherent cosheaves.
 Now the construction of Lemma~\ref{salce-lemma}(a) is applicable.

 Part~(a): given an object $\gE^\bu$ in $\Com(X\ctrh^\lct_\cfq)$,
we need to find an admissible epimorphism in
$\Com(X\ctrh^\lct_\cfq)$ from a contractible complex of projective
locally cotorsion contraherent cosheaves onto~$\gE^\bu$.
 This can be done similarly to the argument for part~(b) using
Theorem~\ref{proj-lct-classification}(a),
Corollary~\ref{lct-proj-local}(b),
and Corollary~\ref{coflasque-acyclic}(b).
\end{proof}

 In line with the general context of this section, we will use
the terminology \emph{homotopy flat complexes of\/ $\bW$\+locally
contraherent cosheaves} for what was called ``homotopy antilocally
flat complexes of $\bW$\+locally contraherent cosheaves'' in
Section~\ref{qc-ss-homotopy-projective-subsect}.
 So, we say that a complex $\gF^\bu$ of $\bW$\+locally contraherent
cosheaves on $X$ is homotopy flat if, for any acyclic complex of
locally cotorsion $\bW$\+locally contraherent cosheaves $\gM^\bu$ on
$X$, any morphism of complexes $\gF^\bu\rarrow\gM^\bu$ is homotopic
to zero.
 The full subcategory of homotopy flat complexes in $\Hot(X\lcth_\bW)$
is denoted by $\Hot(X\lcth_\bW)^\fl$.
 The following lemma is a version of
Lemma~\ref{qcomp-ssep-homotopy-flat-independence}.

\begin{lem}  \label{homotopy-flat-independence}
 Let $X$ be a Noetherian scheme of finite Krull dimension with an open
covering\/~$\bW$.
 Then a complex of flat contraherent cosheaves\/ $\gF^\bu$ on $X$ is
homotopy flat if and only if the complex\/ $\Hom^X(\gF^\bu,\gE^\bu)$
is acyclic for any complex\/ $\gE^\bu$ over $X\ctrh^\lct_\prj$
acyclic with respect to $X\ctrh_\cfq$.
\end{lem}

\begin{proof}
 The argument is similar to the proof of
Lemma~\ref{homotopy-proj-independence} and based on
Theorem~\ref{derived-inj-proj-resolutions}(b) together with
Corollary~\ref{lct-proj-local}(b),
Lemma~\ref{finite-krull-ctrh-lcth-finite-dim}(b), and
Corollary~\ref{fdim-acyclic-cor}.
 The only difference is that one has to use
Corollary~\ref{finite-krull-flat-clf-cor}(a) in order to show that
the complex $\Hom^X$ from any complex of flat contraherent cosheaves
on $X$ to a contraacyclic complex over $X\lcth_\bW^\lct$ is acyclic.
\end{proof}

 According to Lemma~\ref{homotopy-flat-independence}, the property
of a complex over $X\ctrh^\fl$ to belong to $\Hot(X\lcth_\bW)^\fl$
does not depend on the covering~$\bW$ (on a Noetherian scheme~$X$
of finite Krull dimension).
 We denote the full subcategory in $\Hot(X\ctrh^\fl)$ consisting
of the homotopy flat complexes by $\Hot(X\ctrh^\fl)^\fl$.
 One can easily check that bounded above complexes of flat contraherent
cosheaves are homotopy flat, $\Hot^-(X\ctrh^\fl)\sub
\Hot(X\ctrh^\fl)^\fl$.

\begin{thm}  \label{finite-krull-homotopy-flat}
 Let $X$ be a Noetherian scheme of finite Krull dimension.
Then the quotient category of the homotopy category of
homotopy flat complexes of flat contraherent cosheaves\/
$\Hot(X\ctrh^\fl)^\fl$ on $X$ by its thick subcategory of acyclic
complexes over the exact category $X\ctrh^\fl$ is equivalent to
the derived category\/ $\sD(X\ctrh)$.
\end{thm}

\begin{proof}
 By Corollary~\ref{lct-prj-envelope}(b) and Lemma~\ref{psemi-remark21},
any acyclic complex over $X\ctrh^\fl$ is absolutely acyclic; and
one can see from Corollary~\ref{finite-krull-flat-clf-cor}(a) that
any absolutely acyclic complex over $X\ctrh^\fl$ is homotopy flat.
 According to (the proof of)
Theorem~\ref{finite-krull-homotopy-projective}(b), there is
a quasi-isomorphism into any complex over $X\ctrh$ from a complex
belonging to $\Hot(X\ctrh_\prj)_\prj\sub\Hot(X\ctrh^\fl)^\fl$.
 In view of Lemma~\ref{pkoszul-lemma16}(a), it remains to show that
any homotopy flat complex of flat contraherent cosheaves that is acyclic
over $X\ctrh$ is also acyclic over $X\ctrh^\fl$.
 
 According again to Corollary~\ref{lct-prj-envelope}(b) and
the dual version of (the proof of) Proposition~\ref{finite-resolutions},
any complex $\gF^\bu$ over $X\ctrh^\fl$ admits a morphism into
a complex $\P^\bu$ over $X\ctrh^\lct_\prj$ with a cone absolutely
acyclic with respect to $X\ctrh^\fl$.
 If the complex $\gF^\bu$ was homotopy flat, it follows that
the complex $\P^\bu$ is homotopy flat, too.
 This means that $\P^\bu$ is a homotopy projective complex of
(projective) locally cotorsion contraherent cosheaves on~$X$.
 If the complex $\gF^\bu$ was also acyclic over $X\ctrh$, so is
the complex $\P^\bu$.
 It follows that $\P^\bu$ is acyclic over $X\ctrh^\lct$
(see Lemma~\ref{acyclicity-in-lcth-criterion},
Theorem~\ref{cotorsion-periodicity}, or
Corollary~\ref{raynaud-gruson-cotors-cor}), and therefore
contractible.
 We have proved that the complex $\gF^\bu$ is absolutely acyclic
over $X\ctrh^\fl$.
\end{proof}

\subsection{Special inverse image of contraherent cosheaves}
\label{special-inverse-subsection}
 Recall that an affine morphism of schemes $f\:Y\rarrow X$ is called
\emph{finite} if for any affine open subscheme $U\sub X$
the ring $\O_Y(f^{-1}(U))$ is a finitely generated module over
the ring~$\O_X(U)$ \,\cite[Sections Tags~0517 and~01WG]{SP}.
 One can easily see that this condition on a morphism~$f$ is local
in~$X$ \,\cite[Lemma Tag~01WI]{SP}.

 Let $f\:Y\rarrow X$ be a finite morphism of locally
Noetherian schemes.
 Given a quasi-coherent sheaf $\M$ on $X$, one defines
the quasi-coherent sheaf $f^!\M$ on $Y$ by the rule
$$
(f^!\M)(V)=\O_Y(V)\ot_{\O_Y(f^{-1}(U))}\Hom_{\O_X(U)}
(\O_Y(f^{-1}(U))\;\M(U))
$$
for any affine open subschemes $V\sub Y$ and $U\sub X$ such that
$f(V)\sub U$ \cite[Section~III.6]{Har}.
 The construction is well-defined, since for any pair of embedded
affine open subschemes $U'\sub U\sub X$ one has
\begin{multline*}
 \Hom_{\O_X(U')}(\O_Y(f^{-1}(U')\;\M(U')) \\ \simeq
 \Hom_{\O_X(U')}(\O_X(U')\ot_{\O_X(U)}\O_Y(f^{-1}(U))\;
 \O_X(U')\ot_{\O_X(U)}\M(U)) \\ \simeq
 \Hom_{\O_X(U)}(\O_Y(f^{-1}(U))\;\O_X(U')\ot_{\O_X(U)}\M(U)) \\
 \simeq\O_X(U')\ot_{\O_X(U)}\Hom_{\O_X(U)}(\O_Y(f^{-1}(U))\;\M(U)).
\end{multline*}
 Indeed, one has $\Hom_R(L\;F\ot_R M)\simeq F\ot_R\Hom_R(L,M)$
for any module $M$, finitely presented module $L$, and flat
module $F$ over a commutative ring~$R$.
 (See Section~\ref{direct-inverse-loc-contra} for a treatment of 
the non-semi-separatedness issue.)

 The functor $f^!\:X\qcoh\rarrow Y\qcoh$ is right adjoint to
the exact functor $f_*\:Y\qcoh\rarrow X\qcoh$.
 Indeed, it suffices to define a morphism of quasi-coherent sheaves
on $Y$ on the modules of sections over the affine open subschemes
$f^{-1}(U)\sub Y$.
 So given quasi-coherent sheaves $\M$ on $X$ and $\N$ on $Y$, both
groups of morphisms $\Hom_X(f_*\N,\M)$ and $\Hom_Y(\N,f^!\M)$ are
identified with the group of all compatible collections of morphisms
of $\O_X(U)$\+modules $\N(f^{-1}(U))\rarrow\M(U)$, or equivalently,
compatible collections of morphisms of $\O_Y(f^{-1}(U))$\+modules
$\N(f^{-1}(U))\allowbreak\rarrow\Hom_{\O_X(U)}(\O_Y(f^{-1}(U))\;\M(U))$.
 It follows that the functor~$f^!$ takes injective quasi-coherent
sheaves on $X$ to injective quasi-coherent sheaves on~$Y$.

 Let $i\:Z\rarrow X$ be a closed embedding of locally
Noetherian schemes.
 Let $\bW$ be an open covering of $X$ and $\bT$ be an open
covering of $Z$ such that $i$~is a $(\bW,\bT)$\+coaffine morphism.
 Given a $\bW$\+flat $\bW$\+locally contraherent cosheaf $\gF$ on $X$,
one defines a $\bT$\+flat $\bT$\+locally contraherent cosheaf
$i^*\gF$ on $Z$ by the rule
$$
 (i^*\gF)[i^{-1}(U)]=\O_Z(i^{-1}(U))\ot_{\O_X(U)}\gF[U]
$$
for any affine open subscheme $U\sub X$ subordinate to~$\bW$.
 Clearly, affine open subschemes of the form $i^{-1}(U)$ constitute
a base of the topology of~$Z$.
 The construction is well-defined, since for any pair of embedded
affine open subschemes $U'\sub U\sub X$ subordinate to $\bW$ one has
\begin{multline*}
 \O_Z(i^{-1}(U'))\ot_{\O_X(U')}\gF[U'] \\ \simeq
 (\O_Z(i^{-1}(U))\ot_{\O_X(U)}\O_X(U'))\ot_{\O_X(U')}
 \Hom_{\O_X(U)}(\O_X(U'),\gF[U]) \\
 \simeq\O_Z(i^{-1}(U))\ot_{\O_X(U)}\Hom_{\O_X(U)}(\O_X(U'),\gF[U])
 \displaybreak[0] \\
 \simeq\Hom_{\O_X(U)}(\O_X(U')\;\O_Z(i^{-1}(U))\ot_{\O_X(U)}\gF[U]) \\
 \simeq\Hom_{\O_Z(i^{-1}(U))}(\O_Z(i^{-1}(U))\ot_{\O_X(U)}\O_X(U')\;
 \O_Z(i^{-1}(U))\ot_{\O_X(U)}\gF[U]) \\ \simeq
 \Hom_{\O_Z(i^{-1}(U))}(\O_Z(i^{-1}(U'))\;\O_Z(i^{-1}(U))
 \ot_{\O_X(U)}\gF[U]),
\end{multline*}
where the third isomorphism (or the combination of the third and fourth
isomorphisms) holds by
Lemma~\ref{quotient-scalars-contraadjusted}(c).
 The $\O_Z(i^{-1}(U))$\+module $\O_Z(i^{-1}(U))\ot_{\O_X(U)}\gF[U]$
is contraadjusted by Lemma~\ref{quotient-scalars-contraadjusted}(b).

 Assuming that the morphism~$i$ is $(\bW,\bT)$\+affine and
$(\bW,\bT)$\+coaffine, the exact functor $i^*\:X\lcth_\bW^\fl
\rarrow Z\lcth_\bT^\fl$ is ``partially left adjoint'' to the exact
functor $i_!\:Z\lcth_\bT\rarrow X\lcth_\bW$.
 In other words, for any $\bW$\+flat $\bW$\+locally contraherent
cosheaf $\gF$ on $X$ and any $\bT$\+locally contraherent cosheaf $\Q$
on $Z$ there is a natural adjunction isomorphism
\begin{equation}  \label{special-closed-adjunction}
 \Hom^X(\gF,i_!\Q)\simeq \Hom^Z(i^*\gF,\Q).
\end{equation}
 Indeed, it suffices to define a morphism of $\bT$\+locally
contraherent cosheaves on $Z$ on the modules of cosections over
the open subschemes $i^{-1}(U)\sub Z$ for all affine open
subschemes $U\sub X$ subordinate to~$\bW$.
 So both groups of morphisms in question are identified with
the group of all compatible collections of morphisms of
$\O_X(U)$\+modules $\gF[U]\rarrow\Q[i^{-1}(U)]$, or equivalently,
compatible collections of morphisms of $\O_Z((i^{-1}(U))$\+modules
$\O_Z((i^{-1}(U))\ot_{\O_X(U)}\gF[U]\rarrow\Q[i^{-1}(U)]$.

 Notice that for any open covering $\bT$ of a closed subscheme
$Z\sub X$ there exists an open covering $\bW$ of the scheme $X$
for which the embedding morphism~$i$ is $(\bW,\bT)$\+affine and
$(\bW,\bT)$\+coaffine.
 For a locally Noetherian scheme $X$ of finite Krull dimension and
its closed subscheme $Z$, one has $X\lcth_\bW^\fl=X\ctrh^\fl$ and
$Z\lcth_\bT^\fl=Z\ctrh^\fl$ for any open coverings $\bW$ and $\bT$
(see Corollary~\ref{finite-krull-flat-contraherent}(b)).
 In this case, the adjunction
isomorphism~\eqref{special-closed-adjunction} holds for any
flat contraherent cosheaf $\gF$ on $X$ and locally contraherent
cosheaf $\Q$ on~$Z$.
 Most generally, the isomorphism
$$
 \Hom^{\O_X}(\gF,i_!\Q)\simeq \Hom^{\O_Z}(i^*\gF,\Q)
$$
holds for any $\bW$\+flat $\bW$\+locally contraherent cosheaf $\gF$
on a locally Noetherian scheme $X$ and any cosheaf of $\O_Z$\+modules
$\Q$ on a closed subscheme $Z\sub X$.

 Let $f\:Y\rarrow X$ be a finite morphism of locally Noetherian schemes.
 Given a projective locally cotorsion (locally) contraherent cosheaf
$\P$ on $X$, one defines a projective locally cotorsion (locally) 
contraherent cosheaf $f^*\P$ on $Y$ by the rule
$$
 (f^*\P)[V] = \Hom_{\O_Y(f^{-1}(U))}(\O_Y(V)\;
 \O_Y(f^{-1}(U))\ot_{\O_X(U)}\P[U])
$$
for any affine open subschemes $V\sub Y$ and $U\sub X$ such that
$f(V)\sub U$.
 The construction is well-defined, since for any affine open subschemes
$U'\sub U\sub X$ there is a natural isomorphism of
$\O_Y(f^{-1}(U'))$\+modules
$$
 \O_Y(f^{-1}(U'))\ot_{\O_X(U')}\P[U'] \simeq
 \Hom_{\O_Y(f^{-1}(U))}(\O_Y(f^{-1}(U'))\;\O_Y(f^{-1}(U))
 \ot_{\O_X(U)}\P[U])
$$ 
obtained in the way similar to the above computation for the closed
embedding case, except that Lemma~\ref{quotient-scalars-cotorsion}(c)
is being applied.
 The $\O_Y(f^{-1}(U))$\+module $\O_Y(f^{-1}(U))\allowbreak
\ot_{\O_X(U)}\P[U]$ is flat and cotorsion by
Lemma~\ref{quotient-scalars-cotorsion}(a);
hence the $\O_Y(V)$\+module $(f^*\P)[V]$ is flat and cotorsion
by Corollary~\ref{coherent-flat-local}(a)
and Lemma~\ref{cotors-coexten}(a).
 A contraherent cosheaf that is locally flat and cotorsion on
a locally Noetherian scheme $Y$ belongs to $Y\ctrh^\lct_\prj$
by Corollary~\ref{lct-prj-flat}.

 Assuming that the morphism~$f$ is $(\bW,\bT)$\+affine, the functor
$f^*\:X\ctrh^\lct_\prj\rarrow Y\ctrh^\lct_\prj$ is ``partially left
adjoint'' to the exact functor $f_!\:Y\lcth_\bT\rarrow X\lcth_\bW$.
 In other words, for any $\bW$\+flat locally cotorsion $\bW$\+locally
contraherent cosheaf $\P$ on $X$ and any contraherent cosheaf $\Q$
on $Y$ there is a natural adjunction isomorphism
\begin{equation}
 \Hom^X(\P,f_!\Q)\simeq\Hom^Y(f^*\P,\Q).
\end{equation}
 Indeed, it suffices to define a morphism of $\bT$\+locally
contraherent cosheaves on $Y$ on the modules of cosections over
the open subschemes $f^{-1}(U)\sub Y$ for all affine open subschemes
$U\sub X$ subordinate to $\bW$, and the construction proceeds
exactly in the same way as in the above case of a closed embedding~$i$.
 Generally, the isomorphism
$$
 \Hom^{\O_X}(\P,f_!\Q)\simeq \Hom^{\O_Y}(f^*\P,\Q)
$$
holds for any projective locally cotorsion contraherent cosheaf $\P$ on
$X$ and any cosheaf of $\O_Y$\+modules $\Q$ on~$Y$.

 The functor $f^!\:X\qcoh\rarrow Y\qcoh$ for a finite morphism of
locally Noetherian schemes $f\:Y\rarrow X$ preserves infinite
direct sums of quasi-coherent sheaves.
 The functor $i^*\:X\lcth_\bW^\fl\rarrow Z\lcth_\bT^\fl$ for 
a $(\bW,\bT)$\+affine closed embedding of locally Noetherian
schemes $i\:Z\rarrow X$ preserves infinite products of flat
locally contraherent cosheaves.
 The functor $f^*\:X\ctrh^\lct_\prj\rarrow Y\ctrh^\lct_\prj$ for
a finite morphism of locally Noetherian schemes $f\:Y\rarrow X$
preserves infinite products of projective locally cotorsion
contraherent cosheaves.

\subsection{Derived functors of direct and special inverse image}
\label{derived-direct-special}
 For the rest of Chapter~\ref{finite-type-morphisms-sect}, the upper
index~$\bst$ in the notation for derived and homotopy categories
stands for one of the symbols $\b$, $+$, $-$, $\empt$, $\abs+$,
$\abs-$, $\bco$, $\bctr$, $\co$, $\ctr$, or~$\abs$.

 Let $X$ be a locally Noetherian scheme with an open covering~$\bW$.
 The following corollary is to be compared with
Corollary~\ref{coflasque-resolutions-finite}.

\begin{cor} \label{coflasque-resolutions-infinite}
\textup{(a)} The triangulated functor\/ $\sD^{\co=\bco}(X\qcoh^\fq)
\rarrow\sD^{\co=\bco}(X\qcoh)$ is an equivalence of categories.
 The category\/ $\sD(X\qcoh)$ is equivalent to the quotient category
of the homotopy category\/ $\Hot(X\qcoh^\fq)$ by the thick subcategory
of complexes over $X\qcoh^\fq$ that are acyclic over $X\qcoh$.
{\hfuzz=12pt\par}
\textup{(b)} The triangulated functor\/
$\sD^{\ctr=\bctr}(X\ctrh^\lct_\cfq)\rarrow\sD^{\ctr=\bctr}
(X\lcth_\bW^\lct)$ is an equivalence of categories.
 The category\/ $\sD(X\lcth_\bW^\lct)$ is equivalent to the quotient
category of the homotopy category\/ $\Hot(X\ctrh^\lct_\cfq)$ by
the thick subcategory of complexes over $X\ctrh^\lct_\cfq$ that are
acyclic over $X\lcth_\bW^\lct$.
\end{cor}

\begin{proof}
 First of all, one has $\sD^\co(X\qcoh)=\sD^\bco(X\qcoh)$ and
$\sD^\ctr(X\lcth_\bW^\lct)=\sD^\bctr(X\lcth_\bW^\lct)$ by
Theorem~\ref{derived-inj-proj-resolutions}(a\+b).
 The first assertions in parts~(a) and~(b) hold, because there are
enough injective quasi-coherent sheaves, which are flasque, and enough
projective locally cotorsion contraherent cosheaves, which are
coflasque.
 Since the class of flasque quasi-coherent sheaves on a locally
Noetherian scheme is also closed under infinite direct sums, while
the class of coflasque contraherent cosheaves on any scheme is closed
under infinite products (see Section~\ref{coflasque}), the desired
assertions are provided by Propositions~\ref{infinite-resolutions}(b)
and~\ref{becker-contraderived-infinite-resolutions} (for part~(b))
and the dual versions of these propositions (for part~(a)).
 Now we know that there is a quasi-isomorphism from any complex over
$X\qcoh$ into a complex over $X\qcoh^\fq$ and onto any complex over
$X\lcth_\bW^\lct$ from a complex over $X\ctrh^\lct_\cfq$, so the second
assertions in parts~(a) and~(b) follow by Lemma~\ref{pkoszul-lemma16}.
\end{proof}

 Let $f\:Y\rarrow X$ be a quasi-compact morphism of locally
Noetherian schemes.
 As it was mentioned in Section~\ref{lct-homology-subsection},
the functor
\begin{equation}  \label{loc-noetherian-quasi-direct-plus}
 \boR f_*\:\sD^+(Y\qcoh)\lrarrow\sD^+(X\qcoh)
\end{equation}
can be constructed using injective or flasque resolutions.
 More generally, the right derived functor
\begin{equation} \label{loc-noetherian-quasi-direct}
 \boR f_*\:\sD^\st(Y\qcoh)\lrarrow\sD^\st(X\qcoh)
\end{equation}
can be constructed for $\bst=\co$ or~$\bco$ using injective
(see Theorem~\ref{derived-inj-proj-resolutions}(a)) or flasque
(see Corollary~\ref{coflasque-resolutions-infinite}(a)) resolutions,
and for $\bst=\empt$ using homotopy injective
(see Theorem~\ref{quasi-homotopy-injective}) or flasque (cf.\
Corollary~\ref{coflasque-complexes-direct-correct}(a)) resolutions.

 When the scheme $Y$ has finite Krull dimension, the functor $\boR f_*$
can be constructed for any symbol $\bst\ne\ctr$, $\bctr$ using flasque
resolutions (see Corollary~\ref{coflasque-resolutions-finite}(a)).
 Finally, when both schemes $X$ and $Y$ are Noetherian,
the functor~\eqref{loc-noetherian-quasi-direct} can be constructed
for any symbol $\bst\ne\ctr$, $\bctr$ using $f$\+acyclic resolutions
(see Corollary~\ref{f-acyclic-resolutions}(a)).

 The functor
\begin{equation}  \label{loc-noetherian-lct-direct-minus}
 \boL f_!\:\sD^-(Y\lcth^\lct)\lrarrow\sD^-(X\lcth^\lct)
\end{equation}
was constructed in Section~\ref{lct-homology-subsection} using
projective or coflasque resolutions.
 More generally, the left derived functor
\begin{equation} \label{loc-noetherian-lct-direct}
 \boL f_!\:\sD^\st(Y\ctrh^\lct)\lrarrow\sD^\st(X\ctrh^\lct)
\end{equation}
can be constructed for $\bst=\ctr$ or~$\bctr$ using projective
(see Theorem~\ref{derived-inj-proj-resolutions}(b)) or coflasque
(see Corollary~\ref{coflasque-resolutions-infinite}(b)) resolutions;
and for $\bst=\empt$ using homotopy projective (when $X$ is Noetherian,
see Theorem~\ref{finite-krull-homotopy-projective}(a))
or coflasque (in the general case, cf.\
Corollary~\ref{coflasque-complexes-direct-correct}(c)) resolutions.

 When the scheme $Y$ has finite Krull dimension,
the functor~\eqref{loc-noetherian-lct-direct} can be constructed
for any symbol $\bst\ne\co$, $\bco$ using coflasque resolutions (see
Corollary~\ref{coflasque-resolutions-finite}(b)).
 Finally, when both schemes $X$ and $Y$ are Noetherian and one of
the conditions of Lemma~\ref{noetherian-morphism-lct-acyclic-finite-dim}
is satisfied, one can construct the left derived functor
$$
 \boL f_!\:\sD^\st(Y\lcth_\bT^\lct)\lrarrow\sD^\st(X\lcth_\bW^\lct)
$$
for any symbol $\bst\ne\co$, $\bco$ using $f/\bW$\+acyclic resolutions
(see Corollary~\ref{f-acyclic-resolutions}(c)).

 Now assume that the scheme $Y$ is Noetherian of finite Krull dimension.
 The functor
\begin{equation}  \label{loc-noetherian-contra-direct-minus}
 \boL f_!\:\sD^-(Y\lcth)\lrarrow\sD^-(X\lcth)
\end{equation}
was constructed in Section~\ref{lct-homology-subsection} using
projective or coflasque resolutions.
 More generally, the left derived functor
\begin{equation}  \label{loc-noetherian-contra-direct}
 \boL f_!\:\sD^\st(Y\ctrh)\lrarrow\sD^\st(X\ctrh)
\end{equation}
can be constructed for any symbol $\bst\ne\co$, $\bco$ using coflasque
resolutions (see Corollary~\ref{coflasque-resolutions-finite}(c);
cf.\ Corollary~\ref{finite-krull-ctrh-lcth-derived}).

 Assuming that the scheme $Y$ is Noetherian of finite Krull dimension,
for any symbol $\bst\ne\co$, $\bco$, the triangulated
functors~\eqref{loc-noetherian-lct-direct}
and~\eqref{loc-noetherian-contra-direct} form a commutative diagram
with the triangulated functors $\sD^\st(Y\ctrh^\lct)\rarrow
\sD^\st(Y\ctrh)$ and $\sD^\st(X\ctrh^\lct)\rarrow\sD^\st(X\ctrh)$
induced by the respective embeddings of exact categories.
 In particular, when both the schemes $X$ and $Y$ are Noetherian
of finite Krull dimension, we obtain a commutative diagram of
triangulated functors and triangulated equivalences provided by
Corollary~\ref{derived-contra-lct-cor}(b),
\begin{equation} \label{noeth-derived-lct-lcta-direct-images-compatible}
\begin{gathered}
 \xymatrix{
  \sD^\st(Y\ctrh^\lct) \ar@<2pt>[r] \ar@<-2pt>@{-}[r] \ar[d]^{\boL f_!}
  & \sD^\st(Y\ctrh) \ar[d]^{\boL f_!} \\
  \sD^\st(X\ctrh^\lct) \ar@<2pt>[r] \ar@<-2pt>@{-}[r] & \sD^\st(X\ctrh)
 }
\end{gathered}
\end{equation}

 For any quasi-compact morphism~$f$ of locally Noetherian schemes and
any symbol $\bst=\empt$, $\bco$, or~$\co$, the functor $\boR f_*$
\eqref{loc-noetherian-quasi-direct} preserves infinite direct sums,
as it is clear from its construction in terms of complexes of
flasque quasi-coherent sheaves (cf.~\cite[Lemma~1.4]{N-bb}).
 By Theorems~\ref{contraderived-compactly-generated}(b),
\ref{derived-compactly-generated}(c),
and~\cite[Proposition~3.3(1)]{K-st}, it follows that whenever
the scheme $Y$ is Noetherian there exists a triangulated functor
\begin{equation}  \label{loc-noetherian-quasi-adjoint}
 f^!\:\sD^\st(X\qcoh)\lrarrow\sD^\st(Y\qcoh)
\end{equation}
right adjoint to~$\boR f_*$.

 For any symbol $\bst=\empt$, $\bctr$, or~$\ctr$, the functor $\boL f_!$
\eqref{loc-noetherian-lct-direct} preserves infinite products,
as it is clear from its construction in terms of complexes of
coflasque locally cotorsion contraherent cosheaves.
 By Theorems~\ref{contraderived-compactly-generated}(c),
\ref{derived-compactly-generated}(e),
and~\cite[Proposition~3.3(2)]{K-st}, it follows that whenever
the scheme $Y$ is Noetherian of finite Krull dimension there exists
a triangulated functor
\begin{equation}  \label{loc-noetherian-lct-adjoint}
 f^*\:\sD^\st(X\ctrh^\lct)\lrarrow\sD^\st(Y\ctrh^\lct)
\end{equation}
left adjoint to~$\boL f_!$.
 (We refer to the first proof of
Theorem~\ref{finite-krull-homotopy-projective} for the discussion of
``existence'' of the derived category $\sD^\st(X\ctrh^\lct)$, which
is needed for applicability of~\cite[Proposition~3.3(2)]{K-st}.)
 Assuming again that the scheme $Y$ is Noetherian of finite Krull
dimension, the functor $\boL f_!$ \eqref{loc-noetherian-contra-direct}
preserves infinite products, and it follows that there exists
a triangulated functor
\begin{equation}   \label{loc-noetherian-contra-adjoint}
 f^*\:\sD^\st(X\ctrh)\lrarrow\sD^\st(Y\ctrh)
\end{equation}
left adjoint to~$\boL f_!$.

 Assuming that both the schemes $X$ and $Y$ are Noetherian of finite
Krull dimension, and passing to the adjoint functors in the commutative
diagram~\eqref{noeth-derived-lct-lcta-direct-images-compatible}, we
obtain a commutative diagram of triangulated functors and triangulated
equivalences
\begin{equation} \label{noeth-lct-lcta-special-inverse-images-compat}
\begin{gathered}
 \xymatrix{
  \sD^\st(Y\ctrh^\lct) \ar@<2pt>[r] \ar@<-2pt>@{-}[r]
  & \sD^\st(Y\ctrh) \\
  \sD^\st(X\ctrh^\lct) \ar@<2pt>[r] \ar@<-2pt>@{-}[r]
   \ar[u]_{f^*} & \sD^\st(X\ctrh) \ar[u]_{f^*}
 }
\end{gathered}
\end{equation}
for any symbol $\bst=\empt$, $\bctr$, or~$\ctr$.

 Now let $f\:Y\rarrow X$ be a finite morphism of locally Noetherian
schemes.
 Notice that any finite morphism is affine, so the functor
$f_*\:Y\qcoh\rarrow X\qcoh$ is exact; hence we have the induced functor
\begin{equation}  \label{quasi-direct-induced}
 f_*\:\sD^\st(Y\qcoh)\lrarrow\sD^\st(X\qcoh)
\end{equation}
defined for any symbol $\bst\ne\ctr$, $\bctr$
(recall that $\sD^\bco(Y\qcoh)=\sD^\co(Y\qcoh)$ by
Theorem~\ref{derived-inj-proj-resolutions}(a)).
 The right derived functor
\begin{equation}   \label{quasi-derived-special-inverse}
 \boR f^!\:\sD^\st(X\qcoh)\lrarrow\sD^\st(Y\qcoh)
\end{equation}
of the special inverse image functor $f^!\:X\qcoh\rarrow Y\qcoh$
from Section~\ref{special-inverse-subsection} is constructed for
$\bst=+$, $\bco$, or $\co$ in terms of injective resolutions
(see Theorem~\ref{derived-inj-proj-resolutions}(a)), and for
$\bst=\empt$ in terms of homotopy injective resolutions
(see Theorem~\ref{quasi-homotopy-injective}).
 The right derived functor $\boR f^!$
\eqref{quasi-derived-special-inverse} is right adjoint to
the induced functor $f_*$~\eqref{quasi-direct-induced}.

 Similarly, the functor $f_!\:Y\ctrh^\lct\rarrow X\ctrh^\lct$ is
well-defined and exact, as is the functor
$f_!\:Y\ctrh\rarrow X\ctrh$; hence we have the induced functors
\begin{equation}  \label{lct-direct-induced}
 f_!\:\sD^\st(Y\ctrh^\lct)\lrarrow\sD^\st(X\ctrh^\lct)
\end{equation}
defined for any symbol $\bst\ne\co$, $\bco$, and 
\begin{equation}  \label{contra-direct-induced}
 f_!\:\sD^\st(Y\ctrh)\lrarrow\sD^\st(X\ctrh)
\end{equation}
defined for any symbol $\bst\ne\co$, $\bco$, $\bctr$.
 (Recall that $\sD^\bctr(Y\ctrh^\lct)=\sD^\ctr(Y\ctrh^\lct)$ by
Theorem~\ref{derived-inj-proj-resolutions}(a); when $Y$ is Noetherian
of finite Krull dimension, one also has $\sD^\bctr(Y\ctrh)=
\sD^\ctr(Y\ctrh)$ by Theorem~\ref{derived-inj-proj-resolutions}(d).)

 The left derived functor 
\begin{equation}   \label{lct-derived-special-inverse}
 \boL f^*\:\sD^\st(X\ctrh^\lct)\lrarrow\sD^\st(Y\ctrh^\lct)
\end{equation}
of the special inverse image functor $f^*\:X\ctrh^\lct_\prj\rarrow
Y\ctrh^\lct_\prj$ from Section~\ref{special-inverse-subsection}
is constructed for $\bst=-$, $\bctr$, or $\ctr$ in terms of projective
(locally cotorsion) resolutions (see
Theorem~\ref{derived-inj-proj-resolutions}(b)).
 When the scheme $X$ is Noetherian,
the functor~\eqref{lct-derived-special-inverse} is constructed
for $\bst=\empt$ in terms of homotopy projective resolutions
(see Theorem~\ref{finite-krull-homotopy-projective}(a)).
 The left derived functor $\boL f^*$ \eqref{lct-derived-special-inverse}
is left adjoint to the induced functor $f_!$~\eqref{lct-direct-induced}.
 
 Finally, let $i\:Z\rarrow X$ be a closed embedding of Noetherian schemes
of finite Krull dimension.
 Then the left derived functor
\begin{equation}   \label{contra-derived-special-inverse}
 \boL i^*\:\sD^\st(X\ctrh)\lrarrow\sD^\st(Z\ctrh)
\end{equation}
of the special inverse image functor $i^*\:X\ctrh^\fl\rarrow
Z\ctrh^\fl$ from Section~\ref{special-inverse-subsection}
is constructed for $\bst=-$, $\bctr$, or $\ctr$ in terms of flat
resolutions (see Theorem~\ref{derived-inj-proj-resolutions}(d)),
and for $\bst=\empt$ in terms of homotopy projective 
(see Theorem~\ref{finite-krull-homotopy-projective}(b))
or flat and homotopy flat (see Theorem~\ref{finite-krull-homotopy-flat})
resolutions.
 The left derived functor $\boL i^*$
\eqref{contra-derived-special-inverse} is left adjoint to the induced
functor $i_!$~\eqref{contra-direct-induced}.

 Clearly, the constructions of the derived functors
\eqref{lct-derived-special-inverse}
and~\eqref{contra-derived-special-inverse} agree
wherever both are defined.
 In other words, for any closed embedding $i\:Z\rarrow X$ of Noetherian
schemes of finite Krull dimension, and any symbol $\bst=-$, $\empt$,
$\bctr$, or~$\ctr$, we have a commutative diagram of triangulated
functors and triangulated equivalences
\begin{equation} \label{closed-emb-lct-lcta-special-inverse-compat}
\begin{gathered}
 \xymatrix{
  \sD^\st(Z\ctrh^\lct) \ar@<2pt>[r] \ar@<-2pt>@{-}[r]
  & \sD^\st(Z\ctrh) \\
  \sD^\st(X\ctrh^\lct) \ar@<2pt>[r] \ar@<-2pt>@{-}[r]
   \ar[u]_{\boL i^*} & \sD^\st(X\ctrh) \ar[u]_{\boL i^*}
 }
\end{gathered}
\end{equation}
which can be also viewed as a particular case of the commutative
diagram~\eqref{noeth-lct-lcta-special-inverse-images-compat}
for $\bst=\empt$, $\bctr$, or~$\ctr$.

 The following theorem generalizes
Corollaries~\ref{inj-vfl-direct-images-identified}
and~\ref{fl-lct-prj-direct-images-identified} to morphisms~$f$
of not necessarily finite flat dimension (between
Noetherian schemes).

\begin{thm}  \label{finite-krull-inverse-images-identified}
 \textup{(a)} Let $f\:Y\rarrow X$ be a morphism of semi-separated
Noetherian schemes.
 Then the equivalences of triangulated categories\/ $\sD^\co(X\qcoh)
\simeq\sD^\abs(X\ctrh^\lin)$ and\/ $\sD^\co(Y\qcoh)\simeq
\sD^\abs(Y\ctrh^\lin)$ from Theorem~\textup{\ref{co-contra-dualizing}}
transform  the triangulated functor $f^!\:\sD^\co(X\qcoh)\rarrow
\sD^\co(Y\qcoh)$ \textup{\eqref{loc-noetherian-quasi-adjoint}} into
the triangulated functor $f^!\:\sD^\abs(X\ctrh^\lin)\rarrow
\sD^\abs(Y\ctrh^\lin)$ \textup{\eqref{ctrh-inverse-flid-cosheaves}},
\begin{equation}
\begin{gathered}
 \xymatrix{
  \sD^\co(Y\qcoh) \ar@{=}[r] & \sD^\abs(Y\ctrh^\lin) \\
  \sD^\co(X\qcoh) \ar@{=}[r] \ar[u]^{f^!}
  & \sD^\abs(X\ctrh^\lin) \ar[u]_{f^!}
 }
\end{gathered}
\end{equation} \par
 \textup{(b)} Let $f\:Y\rarrow X$ be a morphism of semi-separated
Noetherian schemes of finite Krull dimension.
 Then the equivalences of triangulated categories\/
$\sD^\abs(X\qcoh_\fl)\simeq\sD^\ctr(X\ctrh)$ and\/
$\sD^\abs(Y\qcoh_\fl)\simeq\sD^\ctr(Y\ctrh)$ from
Theorem~\textup{\ref{co-contra-dualizing}}
transform the triangulated functor $f^*\:\sD^\abs(X\qcoh_\fl)\rarrow
\sD^\abs(Y\qcoh_\fl)$ \textup{\eqref{qcoh-inverse-ffd-sheaves}} into
the triangulated functor $f^*\:\sD^\ctr(X\ctrh)\rarrow\sD^\ctr(Y\ctrh)$
\textup{\eqref{loc-noetherian-contra-adjoint}},
\begin{equation}
\begin{gathered}
 \xymatrix{
  \sD^\abs(Y\qcoh_\fl) \ar@{=}[r] & \sD^\ctr(Y\ctrh) \\
  \sD^\abs(X\qcoh_\fl) \ar@{=}[r] \ar[u]^{f^*}
  & \sD^\ctr(X\ctrh) \ar[u]_{f^*}
 }
\end{gathered}
\end{equation}
\end{thm}

\begin{proof}
 This is the special case of
Theorem~\ref{becker-co-contra-derived-direct-inverse-adjunction}
for Noetherian schemes.
 Recall that one has $\sD^\abs(X\ctrh^\lin)=\sD(X\ctrh^\lin)$ for
any quasi-compact semi-separated scheme $X$ by
Corollary~\ref{lin-ctrh-lcth-cor}(a),
$\sD^\co(X\qcoh)=\sD^\bco(X\qcoh)$ for any locally Noetherian
scheme $X$ by Theorem~\ref{derived-inj-proj-resolutions}(a),
$\sD^\abs(X\qcoh_\fl)=\sD(X\qcoh_\fl)$ for any semi-separated
Noetherian scheme of finite Krull dimension by
Corollary~\ref{quasi-finite-flat-dim-all-derived-coincide},
and $\sD^\ctr(X\ctrh)=\sD^\bctr(X\ctrh)\simeq\sD^\bctr(X\ctrh^\lct)
=\sD^\ctr(X\ctrh^\lct)$ for any Noetherian scheme of finite Krull
dimension by Theorem~\ref{derived-inj-proj-resolutions}(b,d)
and Corollary~\ref{derived-contra-lct-cor}(b).
\end{proof}

\subsection{Adjoint functors and bounded complexes}
 The following theorem is a version of
Theorem~\ref{direct-images-identified} for non-semi-separated Noetherian
schemes.
 One would like to have it for an arbitrary morphism of Noetherian
schemes of finite Krull dimension, but we are only able to present
a proof in the case of a flat morphism.

\begin{thm}  \label{non-semi-separated-direct-images-identified}
 Let $f\:Y\rarrow X$ be a flat morphism of Noetherian schemes of finite
Krull dimension.
 Then the equivalences of triangulated categories\/ $\sD(Y\qcoh)\simeq
\sD(Y\ctrh)$ and $\sD(X\qcoh)\simeq\sD(X\ctrh)$ from
Theorem~\textup{\ref{non-semi-separated-naive-co-contra}} transform
the right derived functor\/ $\boR f_*\:\sD(Y\qcoh)\rarrow\sD(X\qcoh)$
\textup{\eqref{loc-noetherian-quasi-direct}} into the left derived
functor\/ $\boL f_!\:\sD(Y\ctrh)\rarrow\sD(X\ctrh)$
\textup{\eqref{loc-noetherian-contra-direct}}.
 In other words, the following square diagram of triangulated functors
and triangulated equivalences is commutative:
\begin{equation} \label{non-ssep-derived-cta-clp-direct-image-diagram}
\begin{gathered}
 \xymatrix{
  \sD^\st(Y\qcoh) \ar@{=}[r] \ar[d]_{\boR f_*}
  & \sD^\st(Y\ctrh) \ar[d]^{\boL f_!} \\
  \sD^\st(X\qcoh) \ar@{=}[r] & \sD^\st(X\ctrh)
 }
\end{gathered}
\end{equation}
\end{thm}

\begin{proof}
 In view of the commutative
diagram~\eqref{noeth-derived-lct-lcta-direct-images-compatible},
it suffices to check that the equivalences of triangulated
categories $\sD(Y\qcoh)\simeq\sD(Y\ctrh^\lct)$ and
$\sD(X\qcoh)\simeq\sD(X\ctrh^\lct)$ constructed in the proof of
Theorem~\ref{non-semi-separated-naive-co-contra} transform
the functor~\eqref{loc-noetherian-quasi-direct} into
the functor~\eqref{loc-noetherian-lct-direct}.

 Let $\O_X\rarrow\E_X^\bu$ be a finite resolution of the sheaf $\O_X$
by flasque quasi-coherent sheaves on~$X$.
 Then the morphism $\O_Y\rarrow f^*\E_X^\bu$ is a quasi-isomorphism
in $Y\qcoh$.
 Pick a finite resolution $f^*\E_X^\bu\rarrow\E_Y^\bu$ of the complex
$f^*\E_X^\bu$ by flasque quasi-coherent sheaves on~$Y$.
 Then the composition $\O_Y\rarrow\E_Y^\bu$ is also
a quasi-isomorphism.

 For any complex $\J^\bu$ over $Y\qcoh^\inj$ there is a natural
isomorphism of complexes of locally cotorsion
contraherent cosheaves $f_!\fHom_Y(f^*\E_X,\J^\bu)\simeq
\fHom_X(\E_X^\bu,f_*\J^\bu)$ on~$X$ \,\eqref{inj-fHom-projection}.
 Composing this isomorphism with the morphism
$f_!\fHom_Y(\E_Y^\bu,\J^\bu)\rarrow f_!\fHom_Y(f^*\E_X^\bu,\J^\bu)$
induced by the quasi-isomorphism $f^*\E_X^\bu\rarrow\E_Y^\bu$,
we obtain a natural morphism of complexes of cosheaves
of $\O_X$\+modules $f_!\fHom_Y(\E_Y^\bu,\J^\bu)\rarrow
\fHom_X(\E_X^\bu,f_*\J^\bu)$.
 Finally, pick a quasi-isomorphism $\gF^\bu\rarrow\fHom_Y(\E_Y^\bu,
\J^\bu)$ of complexes over the exact category $Y\ctrh^\lct$ acting
from a complex $\gF^\bu$ over $Y\ctrh^\lct_\cfq$ to the complex
$\fHom_Y(\E_Y^\bu,\J^\bu)$.
 Applying the functor~$f_!$ and composing again, we obtain a morphism
$\boL f_!\fHom_Y(\E_Y^\bu,\J^\bu)=f_!\gF^\bu\rarrow
\fHom_X(\E_X^\bu,f_*\J^\bu)$ of complexes over $X\ctrh^\lct$.
 Notice that $f_*\J^\bu$ is a complex of injective quasi-coherent
sheaves on~$X$ (since the morphism~$f$ is flat).
 We have constructed a natural transformation
$$
 \boL f_!\,\.\boR\fHom_Y(\E_Y^\bu,{-})\lrarrow
 \boR\fHom_X(\E_X^\bu,\boR f_*({-}))
$$
of functors $\sD(Y\qcoh)\rarrow\sD(X\ctrh^\lct)$.
 One can easily check that such natural transformations are
compatible with the compositions of flat morphisms~$f$.

 Now one can cover the scheme $X$ with affine open subschemes $U_\alpha$
and the scheme $Y$ with affine open subschemes $V_\beta$ so that for
every~$\beta$ there exists~$\alpha$ for which $f(V_\beta)\sub U_\alpha$.
 The triangulated category $\sD(Y\qcoh)$ is generated by the derived
direct images of objects from $\sD(V_\alpha\qcoh)$ (see, e.~g.,
the proof of Theorem~\ref{non-semi-separated-naive-co-contra}).
 So it suffices to show that our natural transformation is
an isomorphism whenever either both schemes $X$ and $Y$ are
semi-separated, or the morphism~$f$ is the embedding of an affine open
subscheme.
 The former case is covered by Theorem~\ref{direct-images-identified},
while in the latter situation one can use
the isomorphism~\eqref{inj-fHom-projection} together with
Lemma~\ref{co-flasque-preservation}(c).
 The point is that, when $f$~is an open embedding, the restriction
functor~$f^*$ preserves flasqueness of quasi-coherent sheaves, and
one can take $\E_Y^\bu=f^*\E_X^\bu$ (while one needs $Y$ to be affine
for applicability of Lemma~\ref{co-flasque-preservation}(c)).

 Alternatively, according to Section~\ref{compatibility-subsect} for
any complex $\P^\bu$ over $Y\ctrh^\lct_\prj$ there is a natural
morphism $\E_X^\bu\ocn_X f_!\P^\bu\rarrow f_*(f^*\E_X^\bu\ocn_Y\P^\bu)$
of complexes of quasi-coherent sheaves on~$X$.
 Composing it with the morphism $f_*(f^*\E_X^\bu\ocn_Y\P^\bu)\rarrow
f_*(\E_Y^\bu\ocn_Y\P^\bu)$ induced by the quasi-isomorphism
$f^*\E_X^\bu\rarrow\E_Y^\bu$, we obtain a natural morphism of complexes
of quasi-coherent sheaves $\E_X^\bu\ocn_X f_!\P^\bu\rarrow
f_*(\E_Y^\bu\ocn_Y\P^\bu)$ on~$X$.
 Finally, pick a quasi-isomorphism $\E_Y^\bu\ocn_Y\P^\bu\rarrow\F^\bu$
of complexes over $Y\qcoh$ acting from the complex
$\E_Y^\bu\ocn_Y\P^\bu$ to a complex $\F^\bu$ over $Y\qcoh^\fq$.
 Applying the functor~$f_*$ and composing again, we obtain a morphism
$\E_X^\bu\ocn_X f_!\P^\bu\rarrow f_*\F^\bu=
\boR f_*(\E_Y^\bu\ocn_Y\P^\bu)$ of complexes over $X\qcoh$.
 Notice that $f_!\P^\bu$ is a complex of projective locally cotorsion
contraherent cosheaves on $X$ (by Corollary~\ref{proj-lct-direct}(b)).
 We have constructed a natural transformation
$$
 \E_X^\bu\ocn_X^\boL\boL f_!({-})\lrarrow
 \boR f_*(\E_Y^\bu\ocn_Y^\boL{-})
$$
of functors $\sD(Y\ctrh^\lct)\rarrow\sD(X\qcoh)$.
 To finish the proof, one continues to argue as above, using
the isomorphism~\eqref{flat-contratensor-projection} and
Lemma~\ref{co-flasque-preservation}(d).
\end{proof}

\begin{rem}
 Let $f\:Y\rarrow X$ be a morphism of Noetherian schemes of finite
Krull dimension.
 If the conclusion of
Theorem~\ref{non-semi-separated-direct-images-identified} holds for~$f$,
it follows by adjunction that the equivalences of categories
$\sD(X\qcoh)\simeq\sD(X\ctrh)$ and $\sD(Y\qcoh)\simeq\sD(Y\ctrh)$
from Theorem~\ref{non-semi-separated-naive-co-contra} transform
the functor $f^!\:\sD(X\qcoh)\rarrow\sD(Y\qcoh)$ 
\eqref{loc-noetherian-quasi-adjoint} into a functor
right adjoint to $\boL f_!\:\sD(Y\ctrh)\rarrow\sD(X\ctrh)$
\eqref{loc-noetherian-contra-direct} and the functor
$f^*\:\sD(X\ctrh)\rarrow\sD(Y\ctrh)$
\eqref{loc-noetherian-contra-adjoint} into a functor left adjoint to
$\boR f_*\:\sD(Y\qcoh)\rarrow\sD(X\qcoh)$
\eqref{loc-noetherian-quasi-direct}.

 Of course, these are supposed to be the conventional derived functors
of inverse image of quasi-coherent sheaves and contraherent cosheaves.
 One can notice, however, that the conventional constructions of such
functors involve some difficulties when one is working outside of
the situations covered by Theorems~\ref{direct-images-identified}
and~\ref{non-semi-separated-direct-images-identified}.
 The case of semi-separated schemes $X$ and $Y$ is covered by our
exposition in Chapter~\ref{derived-on-quasi-compact-sect}
(see~(\ref{qcoh-inverse}\+-\ref{ctrh-inverse})), and in the case
of a flat morphism~$f$ the underived inverse image would do
(if one is working with locally cotorsion contraherent cosheaves).

 In the general case, it is not clear if there exist enough flat
quasi-coherent sheaves or locally injective contraherent cosheaves
to make the derived functor constructions work.
 One can construct the derived functor
$\boL f^*\:\sD(X\qcoh)\rarrow\sD(Y\qcoh)$ using complexes of
sheaves of $\O$\+modules with quasi-coherent cohomology sheaves.
 We do \emph{not} know how to define a derived functor
$\boR f^!\:\sD(X\ctrh)\rarrow\sD(Y\ctrh)$ or
$\boR f^!\:\sD(X\ctrh^\lct)\rarrow\sD(Y\ctrh^\lct)$ for an arbitrary
morphism~$f$ of Noetherian schemes of finite Krull dimension.
\end{rem}

\begin{lem}  \label{derived-bounded-restrict} \hbadness=2350\hfuzz=3pt
 \textup{(a)} For any morphism $f\:Y\rarrow X$ from a Noetherian
scheme $Y$ to a locally Noetherian scheme $X$ such that either
the scheme $X$ is Noetherian or the scheme $Y$ has finite Krull
dimension, the triangulated functor $f^!\:\sD(X\qcoh)\rarrow
\sD(Y\qcoh)$ \textup{\eqref{loc-noetherian-quasi-adjoint}}
takes\/ $\sD^+(X\qcoh)$ into\/ $\sD^+(Y\qcoh)$ and induces
a triangulated functor $f^!\:\sD^+(X\qcoh)\rarrow\sD^+(Y\qcoh)$
right adjoint to the right derived functor\/
$\boR f_*$~\textup{\eqref{loc-noetherian-quasi-direct-plus}}. \par
 \textup{(b)} For any morphism of Noetherian schemes of finite
Krull dimension $f\:Y\rarrow X$, the triangulated functor
$f^*\:\sD(X\ctrh^\lct)\rarrow\sD(Y\ctrh^\lct)$
\textup{\eqref{loc-noetherian-lct-adjoint}} takes\/
$\sD^-(X\ctrh^\lct)$ into\/ $\sD^-(Y\ctrh^\lct)$ and induces
a triangulated functor $f^*\:\sD^-(X\ctrh^\lct)\allowbreak\rarrow
\sD^-(Y\ctrh^\lct)$ left adjoint to the left derived functor\/
$\boL f_!$~\textup{\eqref{loc-noetherian-lct-direct-minus}}. \par
 \textup{(c)} For any morphism of Noetherian schemes of finite
Krull dimension $f\:Y\rarrow X$, the triangulated functor
$f^*\:\sD(X\ctrh)\rarrow\sD(Y\ctrh)$
\textup{\eqref{loc-noetherian-contra-adjoint}} takes\/
$\sD^-(X\ctrh)$ into\/ $\sD^-(Y\ctrh)$ and induces a triangulated
functor $f^*\:\sD^-(X\ctrh)\rarrow\sD^-(Y\ctrh)$ left adjoint to
the left derived functor\/
$\boL f_!$~\textup{\eqref{loc-noetherian-contra-direct-minus}}.
\end{lem}

\begin{proof}
 In each part~(a\+c), the second assertion follows immediately
from the first one (since the derived functors of direct image
of bounded and unbounded complexes agree).
 In view of Corollary~\ref{derived-contra-lct-cor}(b), part~(c)
is also equivalent to part~(b).

 To prove part~(a), notice that a complex of quasi-coherent sheaves
$\N^\bu$ over $Y$ has its cohomology sheaves concentrated in
the cohomological degrees~$\ge -N$ if and only if one has
$\Hom_{\sD(Y\qcoh)}(\L^\bu,\N^\bu)=0$ for any complex $\L^\bu$
over $Y\qcoh$ concentrated in the cohomological degrees~$<-N$.
 This is true for any abelian category in place of $Y\qcoh$.
 Now given a complex $\M^\bu$ over $X\qcoh$, one has
$\Hom_{\sD(Y\qcoh)}(\L^\bu,f^!\M^\bu)\simeq
\Hom_{\sD(X\qcoh)}(\boR f_*\L^\bu,\M^\bu)=0$ whenever $\M^\bu$
is concentrated in the cohomological degrees~$\ge -N+M$, where
$M$~is a certain fixed constant.
 Indeed, the functor $\boR f_*$ raises the cohomological degrees
by at most the Krull dimension of the scheme~$Y$ (if it has
finite Krull dimension), or by at most the constant from
Lemmas~\ref{noetherian-morphism-lct-acyclic-finite-dim}\+-%
\ref{noetherian-morphism-qcoh-acyclic-finite-dim} (if both
schemes are Noetherian).

 In order to deduce parts~(b\+c), we will use
Theorem~\ref{non-semi-separated-direct-images-identified}.
 The equivalence of triangulated categories $\sD(Y\qcoh)\simeq
\sD(Y\ctrh)$ from Theorem~\ref{non-semi-separated-naive-co-contra}
identifies $\sD^-(Y\qcoh)$ with $\sD^-(Y\ctrh)$.
 Given a complex $\gM^\bu$ from $\sD^-(X\ctrh)$, we would like to
show that the complex $\N^\bu$ over $Y\qcoh$ corresponding
to the complex $\gN^\bu=f^*\gM^\bu$ over $Y\ctrh$ belongs to
$\sD^-(Y\qcoh)$.
 This is equivalent to saying that the complex $j'{}^*\N^\bu$
over $V\qcoh$ belongs to $\sD^-(V\qcoh)$ for the embedding of any
small enough affine open subscheme $j'\:V\rarrow Y$.

 By Theorem~\ref{non-semi-separated-direct-images-identified}
applied to the morphism~$j'$, the complex $j'{}^*\N^\bu$ corresponds
to the complex $j'{}^*\gN^\bu$ over $V\ctrh$.
 We can assume that the composition $f\circ j'\:V\rarrow X$ factorizes
through the embedding of an affine open subscheme $j\:U\rarrow X$.
 Let $f'$~denote the related morphism $V\rarrow U$.
 Then one has $j'{}^*\gN^\bu\simeq f'{}^*\!\.j^*\gM^\bu$ in
$\sD(V\ctrh)$.
 Let $\M^\bu$ denote the complex over $X\qcoh$ corresponding
to~$\gM^\bu$; by Theorem~\ref{non-semi-separated-naive-co-contra},
one has $\M^\bu\in\sD^-(X\qcoh)$.
 Applying again Theorem~\ref{non-semi-separated-direct-images-identified}
to the morphism~$j$ and Theorem~\ref{direct-images-identified} to
the morphism~$f'$, we conclude that the complex
$f'{}^*\!\.j^*\gM$ over $\sD(V\ctrh)$ corresponds to the complex
$\boL f'{}^*j^*\M$ over $\sD(V\qcoh)$.
 The latter is clearly bounded above, and the desired assertion
is proved.
\end{proof}

\begin{cor}  \label{co-contra-derived-bounded-restrict}
\textup{(a)} For any morphism $f\:Y\rarrow X$ from a Noetherian scheme
$Y$ to a locally Noetherian scheme $X$ such that either the scheme
$X$ is Noetherian or the scheme $Y$ has finite Krull dimension,
the triangulated functor $f^!\:\sD^\co(X\qcoh)\rarrow\sD^\co(Y\qcoh)$
\textup{\eqref{loc-noetherian-quasi-adjoint}} takes\/ $\sD^+(X\qcoh)$
into\/ $\sD^+(Y\qcoh)$, and the induced triangulated functor
$f^!\:\sD^+(X\qcoh)\rarrow\sD^+(Y\qcoh)$ coincides with the one
obtained in Lemma~\textup{\ref{derived-bounded-restrict}}. \par
 \textup{(b)} For any morphism of Noetherian schemes of finite
Krull dimension $f\:Y\rarrow X$, the triangulated functor
$f^*\:\sD^\ctr(X\ctrh^\lct)\rarrow\sD^\ctr(Y\ctrh^\lct)$
\textup{\eqref{loc-noetherian-lct-adjoint}} takes\/
$\sD^-(X\ctrh^\lct)$ into\/ $\sD^-(Y\ctrh^\lct)$, and the induced
triangulated functor $f^*\:\allowbreak\sD^-(X\ctrh^\lct)\rarrow
\sD^-(Y\ctrh^\lct)$ coincides with the one obtained in
Lemma~\textup{\ref{derived-bounded-restrict}}. {\hbadness=2500\par}
 \textup{(c)} For any morphism of Noetherian schemes of finite
Krull dimension $f\:Y\rarrow X$, the triangulated functor
$f^*\:\sD^\ctr(X\ctrh)\rarrow\sD^\ctr(Y\ctrh)$
\textup{\eqref{loc-noetherian-contra-adjoint}} takes\/
$\sD^-(X\ctrh)$ into\/ $\sD^-(Y\ctrh)$, and the induced
triangulated functor $f^*\:\sD^-(X\ctrh^\lct)\allowbreak\rarrow
\sD^-(Y\ctrh^\lct)$ coincides with the one obtained in
Lemma~\textup{\ref{derived-bounded-restrict}}. \hbadness=2100
\end{cor}

\begin{proof}
 Notice first of all that there are natural fully faithful functors
$\sD^+(X\qcoh)\rarrow\sD^\co(X\qcoh)$, \ $\sD^-(X\ctrh^\lct)\rarrow
\sD^\ctr(X\ctrh^\lct)$, \ $\sD^-(X\ctrh)\rarrow\sD^\ctr(X\ctrh)$
(and similarly for~$Y$) by Lemma~\ref{co-contra-bounded-fully-faithful}.
 Furthermore, one can prove the first assertion of part~(a) in the way
similar to the proof of Lemma~\ref{derived-bounded-restrict}(a).

 A complex of quasi-coherent sheaves $\N^\bu$ over $Y$ is isomorphic
in $\sD^\co(Y\qcoh)$ to a complex whose terms are concentrated in
the cohomological degrees~$\ge -N$ if and only if one has
$\Hom_{\sD^\co(Y\qcoh)}(\L^\bu,\N^\bu)=0$ for any complex $\L^\bu$
over $Y\qcoh$ whose terms are concentrated in the cohomological
degrees~$<N$.
 This is true for any abelian category with exact functors of infinite
direct sum in place of $Y\qcoh$, since complexes with the terms
concentrated in the degrees~$\le0$ and~$\ge0$ form a t\+structure
on the coderived category~\cite[Remark~4.1]{Psemi}.
 For a more explicit argument, see~\cite[Lemma~2.2]{K-st}.
 The rest of the proof of the first assertion is similar to that of
Lemma~\ref{derived-bounded-restrict}(a); and both assertions of
part~(a) can be proved in the way dual-analogous to the following
proof of part~(c).

 It is clear from the constructions of the functors
$\boL f_!$~\eqref{loc-noetherian-contra-direct} in terms of coflasque
resolutions that the functors $\boL f_!\:\sD^\ctr(Y\ctrh)\rarrow
\sD^\ctr(X\ctrh)$ and $\boL f_!\:\sD(Y\ctrh)\rarrow\sD(X\ctrh)$ form
a commutative diagram with the Verdier localization functors
$\sD^\ctr(Y\ctrh)\rarrow\sD(Y\ctrh)$ and $\sD^\ctr(X\ctrh)\rarrow
\sD(X\ctrh)$.
 According to~\eqref{loc-noetherian-contra-adjoint} and (the proof of)
Theorem~\ref{finite-krull-homotopy-projective}(b), all the functors
in this commutative square have left adjoints, which therefore also
form a commutative square. {\hbadness=1450\par}

 The functor $\sD(X\ctrh)\rarrow\sD^\ctr(X\ctrh)$ left adjoint to
the localization functor can be constructed as the functor assigning
to a complex over $X\ctrh$ its homotopy projective resolution, viewed
as an object of $\sD^\ctr(X\ctrh)$ (and similarly for~$Y$).
 Since the homotopy projective resolution $\P^\bu$ of a bounded above
complex $\gM^\bu$ over $X\ctrh$ can be chosen to be also bounded above,
and the cone of the morphism $\P^\bu\rarrow\gM^\bu$ is contraacyclic,
part~(c) follows from Lemma~\ref{derived-bounded-restrict}(c).
\end{proof}

 We refer to~\cite[D\'efinition~5.4.1]{Groth1}
or~\cite[Section Tag~01W0]{SP} for the definition of
a \emph{proper morphism} of schemes.

\begin{prop}  \label{proper-finite-flat-dim-extraordinary-inverse}
 Let $f\:Y\rarrow X$ be a proper morphism (of finite type and) of
finite flat dimension between semi-separated Noetherian schemes.
 Assume that the object $f^!\O_X$ is compact in\/ $\sD(Y\qcoh)$
(i.~e., it is a perfect complex on~$Y$).
 Then the functor $f^!\:\sD^\co(X\qcoh)\rarrow\sD^\co(Y\qcoh)$
\textup{\eqref{loc-noetherian-quasi-adjoint}} is naturally
isomorphic to the functor $f^!\O_X\ot_{\O_Y}^{\boL'}\boL f^*$,
where $f^!\O_X\in\sD(Y\qcoh)$, while\/ $\boL f^*\:\sD^\co(X\qcoh)
\rarrow\sD^\co(Y\qcoh)$ is the functor constructed
in~\textup{(\ref{qcoh-inverse-ffd-morphism},~%
\ref{bco-qcoh-inverse-ffd-morphism})} and\/
$\ot_{\O_Y}^{\boL'}$ is the tensor action
functor~\textup{\eqref{quasi-tensor-action}}.
\end{prop}

\begin{proof}
 Instead of proving the desired assertion from scratch, we will
deduce it step by step from the related results
of~\cite[Section~5]{N-bb}, using the results above in this section
and Section~\ref{compact-generators-subsect} to bridge the gap between
the derived and coderived categories.

\begin{lem}
 Let $f\:Y\rarrow X$ be a morphism of finite flat dimension between
semi-separated Noetherian schemes.
 Then for any objects\/ $\L^\bu\in\sD(Y\qcoh)$ and\/
$\M^\bu\in\sD^\co(X\qcoh)$ there is a natural isomorphism\/
$\boR f_*\L^\bu\ot_{\O_X}^{\boL'}\M^\bu\simeq
\boR f_*(\L^\bu\ot_{\O_Y}^{\boL'}\boL f^*\M^\bu)$ in\/
$\sD^\co(X\qcoh)$, where\/ $\boR f_*$ denotes the derived direct
image functors~\eqref{qcoh-direct}.
\end{lem}

\begin{proof}
 For any objects $\K^\bu\in\sD(X\qcoh)$ and $\M^\bu\in\sD^\co(X\qcoh)$
one easily constructs a natural isomorphism $\boL f^*(\K^\bu\ot_{\O_X}
^{\boL'}\M^\bu)\simeq\boL f^*\K^\bu\ot_{\O_Y}^{\boL'}\boL f^*\M^\bu$
in $\sD^\co(Y\qcoh)$ (where $\boL f^*$ denotes the derived functors
(\ref{qcoh-inverse},~\ref{qcoh-inverse-ffd-morphism})).
 Substituting $\K^\bu=\boR f_*\L^\bu$ with $\L^\bu\in\sD(Y\qcoh)$,
one can consider the composition
$\boL f^*(\boR f_*\L^\bu\ot_{\O_X}^{\boL'}\M^\bu)\simeq
\boL f^*\.\boR f_*\L^\bu\ot_{\O_Y}^{\boL'}\boL f^*\M^\bu\rarrow
\L^\bu\ot_{\O_Y}^{\boL'}\boL f^*\M^\bu$ of morphisms in
$\sD^\co(Y\qcoh)$.
 By adjunction, we obtain the natural transformation
$\boR f_*\L^\bu\ot_{\O_X}^{\boL'}\M^\bu\rarrow
\boR f_*(\L^\bu\ot_{\O_Y}^{\boL'}\boL f^*\M^\bu)$ of functors
$\sD(Y\qcoh)\times\sD^\co(X\qcoh)\rarrow\sD^\co(X\qcoh)$.

 Since all the functors involved preserve infinite direct sums, it
suffices to check that our morphism is an isomorphism for compact
generators of the categories $\sD(Y\qcoh)$ and $\sD^\co(X\qcoh)$,
that is one can assume $\L^\bu$ to be a perfect complex on $X$ and
$\M^\bu$ to be a finite complex of coherent sheaves on~$Y$
(see Theorems~\ref{derived-compactly-generated}(b)
and~\ref{contraderived-compactly-generated}(b)).
 In this case all the complexes involved are bounded below and,
in view of Lemma~\ref{co-contra-bounded-fully-faithful},
the question reduces to the similar assertion for the conventional
derived categories, which is known due to~\cite[Proposition~5.3]{N-bb}.
\end{proof}

 In particular, for any object $\M^\bu\in\sD^\co(X\qcoh)$ we have
a natural isomorphism $\boR f_*f^!\O_X\ot_{\O_X}^{\boL'}\M^\bu
\simeq\boR f_*(f^!\O_X\ot_{\O_Y}^{\boL'}\boL f^*\M^\bu)$.
 Composing it with the morphism induced by the adjunction morphism
$\boR f_*f^!\O_X\rarrow\O_X$, we obtain a natural morphism
$\boR f_*(f^!\O_X\ot_{\O_Y}^{\boL'}\boL f^*\M^\bu)\rarrow
\M^\bu$ in $\sD^\co(X\qcoh)$.
 We have constructed a natural transformation
$f^!\O_X\ot_{\O_Y}^{\boL'}\boL f^*\M^\bu\rarrow f^!\M^\bu$ of
functors $\sD^\co(X\qcoh)\rarrow\sD^\co(Y\qcoh)$.

 Since the morphism~$f$ is proper, the functor $\boR f_*\:
\sD^\co(Y\qcoh)\rarrow\sD^\co(X\qcoh)$ takes compact objects
to compact objects~\cite[Th\'eor\`eme~3.2.1]{Groth2},
\cite[Proposition Tag~02O5]{SP}.
 By~\cite[Theorem~5.1]{N-bb}, it follows that the functor
$f^!\:\sD^\co(X\qcoh)\rarrow\sD^\co(Y\qcoh)$ preserves infinite
direct sums.
 So does the functor in the left-hand side of our morphism;
therefore, it suffices to check that this morphism is an isomorphism
when $\M^\bu$ is a finite complex of coherent sheaves on~$X$.
 Since we assume $f^!\O_X$ to be a perfect complex, in this case all
the complexes involved are bounded below, and in view of
Corollary~\ref{co-contra-derived-bounded-restrict}(a), the question
again reduces to the similar assertion for the conventional derived
categories, which is provided
by~\cite[Example~5.2 and Theorem~5.4]{N-bb}
(cf.~\cite[Proposition~4.7.1]{Lip}).
\end{proof}

\begin{rem}
 The condition that $f^!\O_X$ is a perfect complex, which was not
needed in~\cite[Section~5]{N-bb}, cannot be dropped in the above
Proposition, as one can see already in the case when $X$ is
the spectrum of a field and $Y\rarrow X$ is a finite morphism.
 The problem arises because of the difference between the left
and right adjoint functors to the Verdier localization functor
$\sD^\co(Y\qcoh)\rarrow\sD(Y\qcoh)$.

 To construct a specific counterexample, let $k$~be a field and
$R$ be the quotient ring $R=k[x,y]/(x^2,xy,y^2)$, so that
the images of the elements $1$, $x$, and $y$ form a basis in $R$
over~$k$.
 Set $X=\Spec k$ and $Y=\Spec R$.
 Then the injective coherent sheaf $\widetilde{R^*}$ on $Y$ 
corresponding to the $R$\+module $R^*=\Hom_k(R,k)$ represents
the object $f^!\O_X\in\sD(Y\qcoh)$ as well as the object
$f^!\O_X\in\sD^\co(Y\qcoh)$.
 Looking into the construction of the tensor action functor
$\ot_{\O_Y}^{\boL'}$ \eqref{quasi-tensor-action}, in order to show
that $f^!\O_X\not\simeq f^!\O_X\ot_{\O_Y}^{\boL'}\O_Y$, we have to
prove that the $R$\+module $R^*$ is not isomorphic in $\sD^\co(R\modl)$
to its projective $R$\+module resolution~$P_\bu$.

 Indeed, the complex $\Hom_R(P_\bu,M^\bu)$ is acyclic for any acyclic
complex of $R$\+modules $M^\bu$, while the complex $\Hom_R(A^\bu,J^\bu)$
is acyclic for any coacyclic complex of $R$\+modules $A^\bu$ and any
complex of injective $R$\+modules~$J^\bu$.
 So, if $R^*$ were isomorphic to $P_\bu$ in $\sD^\co(R\modl)$, then
it would follow that the complex $\Hom_R(R^*,J^\bu)$ is acyclic for
any acyclic complex of injective $R$\+modules~$J^\bu$.
 Since there is a short exact sequence of $R$\+modules
$0\rarrow k^{\oplus 3}\rarrow R^{\oplus 2}\rarrow R^*\rarrow0$,
the complex $\Hom_R(k,J^\bu)$ would then be also acyclic, implying that
the complex of injective $R$\+modules $J^\bu$ is contractible
(see, e.~g., \cite[Lemma~1.4.3]{Pweak}).
 This would mean that the coderived category of $R$\+modules
coincides with their derived category, which cannot be true,
as their subcategories of compact objects are clearly different.
 For example, the $R$\+module~$k$ is compact in $\sD^\co(R\modl)$,
but not in $\sD(R\modl)$.

 On the other hand, one can get rid of the semi-separatedness
assumptions
in Proposition~\ref{proper-finite-flat-dim-extraordinary-inverse}
by constructing the functor $\boL f^*\:\sD(X\qcoh)\rarrow
\sD(Y\qcoh)$ in terms of complexes of sheaves of $\O$\+modules
with quasi-coherent cohomology sheaves (as it is done in~\cite{N-bb};
cf.\ the proof of Theorem~\ref{derived-compactly-generated}(b\+c))
and obtaining the functor $\boL f^*\:\sD^\co(X\qcoh)\rarrow
\sD^\co(Y\qcoh)$ from it using the techniques of~\cite{Gai}.
\end{rem}

\subsection{Compatibilities for a smooth morphism}
 We refer to~\cite[Section Tag~01V4]{SP} for a general discussion of
smooth morphisms of schemes.
 From our perspective in this section, it is important that smooth
morphisms are of finite type/finite presentation, flat, and with regular
(or at least, Gorenstein) fibers of bounded Krull dimension (cf.\
the discussion of ``weakly smooth morphisms''
in~\cite[Chapter~10]{Psemten}).

 Let $f\:Y\rarrow X$ be a smooth morphism of Noetherian schemes.
 Let $\D_X^\bu$ be a dualizing complex for $X$; then $f^*\D_X^\bu$
is a dualizing complex for~$Y$ \,\cite[Theorem~V.8.3]{Har},
\cite[Lemma Tag~0E4D]{SP}, \cite[Lemma~10.20]{Psemten}.
 The complex $f^*\D_X^\bu$ being not necessarily a complex of
injectives, let us pick a finite complex over $Y\qcoh$ quasi-isomorphic
to $f^*\D_X^\bu$ and denote it temporarily by~$\D_Y^\bu$.

\begin{cor}  \label{tensor-cohom-dualizing-compatibility} \hfuzz=4pt
 Assume the schemes $X$ and $Y$ to be semi-separated. \par
\textup{(a)} The equivalences of triangulated categories\/
$\sD^\abs(X\qcoh_\fl)\simeq\sD^\co(X\qcoh)$ and\/ $\sD^\abs(Y\qcoh_\fl)
\simeq\sD^\co(Y\qcoh)$ from Theorem~\textup{\ref{co-contra-dualizing}}
related to the choice of the dualizing complexes\/ $\D_X^\bu$ and\/
$\D_Y^\bu$ on $X$ and $Y$ transform the inverse image functor
$f^*\:\sD^\abs(X\qcoh_\fl)\rarrow\sD^\abs(Y\qcoh_\fl)$
\textup{\eqref{qcoh-inverse-ffd-sheaves}} into the (underived, as
the morphism~$f$ is flat) inverse image functor
$f^*\:\sD^\co(X\qcoh)\rarrow\sD^\co(Y\qcoh)$
\textup{(\ref{qcoh-inverse-ffd-morphism},
\ref{bco-qcoh-inverse-ffd-morphism})},
\begin{equation}
\begin{gathered}
 \xymatrix{
  \sD^\abs(Y\qcoh_\fl) \ar@{=}[r] & \sD^\co(Y\qcoh) \\
  \sD^\abs(X\qcoh_\fl) \ar@{=}[r] \ar[u]^{f^*}
  & \sD^\co(X\qcoh) \ar[u]^{f^*}
 }
\end{gathered}
\end{equation} \par
\textup{(b)} The equivalences of triangulated categories\/
$\sD^\abs(X\ctrh^\lin)\simeq\sD^\ctr(X\ctrh)$ and\/ $\sD^\abs
(Y\ctrh^\lin)\simeq\sD^\ctr(Y\ctrh)$ from
Theorem~\textup{\ref{co-contra-dualizing}} related to the choice of
the dualizing complexes\/ $\D_X^\bu$ and\/ $\D_Y^\bu$ on $X$ and $Y$
transform the inverse image functor
$f^!\:\sD^\abs(X\ctrh^\lin)\rarrow\sD^\abs(Y\ctrh^\lin)$
\textup{\eqref{ctrh-inverse-flid-cosheaves}} into the (underived, as
the morphism~$f$ is very flat) inverse image functor
$f^!\:\sD^\ctr(X\ctrh)\rarrow\sD^\ctr(Y\ctrh)$
\textup{(\ref{ctrh-inverse-ffd-morphism},
\ref{bctr-ctrh-inverse-ffd-morphism})},
\begin{equation}
\begin{gathered}
 \xymatrix{
  \sD^\ctr(Y\ctrh) \ar@{=}[r] & \sD^\abs(Y\ctrh^\lin) \\
  \sD^\ctr(X\ctrh) \ar@{=}[r] \ar[u]_{f^!}
  & \sD^\abs(X\ctrh^\lin) \ar[u]_{f^!}
 }
\end{gathered}
\end{equation}
\end{cor}

\begin{rem}
 By Theorem~\ref{very-flat-theorem}(a), any flat morphism of finite type
between Noetherian schemes is very flat (as is, more generally,
any flat morphism of finite presentation between arbitrary schemes).
 Without using the very flat conjecture/theorem, the morphism~$f$
being at least flat, there is an underived inverse image functor
$f^!\:\sD^\ctr(X\ctrh^\lct)\rarrow\sD^\ctr(Y\ctrh^\lct)$
\,(\ref{ctrh-lct-inverse-ffd-morphism},
\ref{bctr-ctrh-lct-inverse-ffd-morphism}).
\end{rem}

\begin{proof}[Proof of
Corollary~\ref{tensor-cohom-dualizing-compatibility}]
 Part~(a) (cf.~\cite[Section~3.8]{EP}): notice that for any
finite complex $\E^\bu$ of quasi-coherent sheaves on $Y$ the functor
$\sD^\abs(Y\qcoh_\fl)\rarrow\sD^\co(Y\qcoh)$ constructed by tensoring
complexes over $Y\qcoh_\fl$ with the complex $\E^\bu$ over $\O_Y$
is well-defined, and replacing $\E^\bu$ by a quasi-isomorphic finite
complex leads to an isomorphic functor.
 Indeed, the tensor product of any absolutely acyclic complex in
$Y\qcoh_\fl$ with any complex in $Y\qcoh$ is an absolutely acyclic
complex in $Y\qcoh$, as is the tensor product of any complex in
$Y\qcoh_\fl$ with an absolutely acyclic complex in $Y\qcoh$.
 So it remains to use the isomorphism $f^*(\D_X^\bu\ot_{\O_X}\F)\simeq
f^*\D_X^\bu\ot_{\O_X}f^*\F$ holding for any flat (or even arbitrary)
quasi-coherent sheaf $\F$ on~$X$.

 Part~(b): let $\bW$ and $\bT$ be open coverings of the schemes $X$
and $Y$ such that the morphism~$f$ is $(\bW,\bT)$\+coaffine.
 For any finite complex $\E^\bu$ of quasi-coherent sheaves on $Y$,
the functor $\sD^\abs(Y\lcth_\bT^\lin)\rarrow\sD^\ctr(Y\lcth_\bT^\lct)
\rarrow\sD^\ctr(Y\lcth_\bT)$ constructed by taking $\Cohom_Y$ from
$\E^\bu$ to complexes over $Y\lcth_\bT^\lin$ is well-defined, and
replacing $\E^\bu$ by a quasi-isomorphic complex leads to an isomorphic
functor.
 Indeed, for any complex $\E^\bu$ in $Y\qcoh$ and any absolutely
acyclic complex $\gJ^\bu$ in $Y\lcth_\bT^\lin$, the complex
$\Cohom_Y(\E^\bu,\gJ^\bu)$ (constructed by totalizing the bicomplex of
$\Cohom$ by taking infinite products along the diagonals) is absolutely
acyclic in $Y\lcth_\bT^\lct$, and the same holds for any absolutely
acyclic complex $\E^\bu$ in $X\qcoh$ and any complex $\gJ^\bu$ in
$Y\lcth_\bT^\lin$.
 So it remains to use the isomorphism~\eqref{flat-lin-inverse-cohom}.
\end{proof}

 The following theorem is to be compared with
Theorems~\ref{direct-images-identified}
and~\ref{non-semi-separated-direct-images-identified}, and
Corollaries~\ref{inj-vfl-direct-images-identified}
and~\ref{fl-lct-prj-direct-images-identified}.

\begin{thm}  \label{smooth-direct-images-identified}
 The equivalences of triangulated categories\/
$\sD^\co(Y\qcoh)\simeq\sD^\ctr(Y\ctrh)$ and\/
$\sD^\co(X\qcoh)\simeq\sD^\ctr(X\ctrh)$ from
Theorem~\textup{\ref{co-contra-dualizing}}
or~\textup{\ref{non-semi-separated-co-contra}} related to the choice of
the dualizing complexes\/ $\D_Y^\bu$ and\/ $\D_X^\bu$ on $Y$ and $X$
transform the right derived functor\/ $\boR f_*\:\sD^\co(Y\qcoh)\rarrow
\sD^\co(X\qcoh)$ \textup{\eqref{loc-noetherian-quasi-direct}}
into the left derived functor\/ $\boL f_!\:\sD^\ctr(Y\ctrh)\rarrow
\sD^\ctr(X\ctrh)$ \textup{\eqref{loc-noetherian-contra-direct}}.
 In other words, the following square diagram of triangulated functors
and triangulated equivalences is commutative: \hbadness=1475
\begin{equation} \label{smooth-direct-images-diagram}
\begin{gathered}
 \xymatrix{
  \sD^\co(Y\qcoh) \ar@{=}[r] \ar[d]_{\boR f_*}
  & \sD^\ctr(Y\ctrh) \ar[d]^{\boL f_!} \\
  \sD^\co(X\qcoh) \ar@{=}[r] & \sD^\ctr(X\ctrh)
 }
\end{gathered}
\end{equation}
\end{thm}

\begin{proof}
 In view of the commutative
diagram~\eqref{noeth-derived-lct-lcta-direct-images-compatible},
it suffices to show that the equivalences of triangulated categories
$\sD^\co(Y\qcoh)\simeq\sD^\ctr(Y\ctrh^\lct)$ and
$\sD^\co(X\qcoh)\simeq\sD^\ctr(X\ctrh^\lct)$ from the proof of
Theorem~\ref{non-semi-separated-co-contra} transform
the functor~\eqref{loc-noetherian-quasi-direct} into
the functor~\eqref{loc-noetherian-lct-direct}.
 The latter equivalences were constructed on the level of injective
and projective resolutions as the equivalences of homotopy categories
$\Hot(Y\qcoh^\inj)\simeq\Hot(Y\ctrh^\lct_\prj)$ induced by
the functors $\fHom_Y(\D_Y^\bu,{-})$ and $\D_Y^\bu\ocn_Y{-}$,
and similarly for~$X$.
 In particular, it was shown that these functors take complexes over
$Y\qcoh^\inj$ to complexes over $Y\ctrh^\lct_\prj$ and back.

 Furthermore, we notice that such complexes are adjusted to
the derived functors $\boR f_*$ and $\boL f_!$ acting between
the co/contraderived categories.
 The morphism~$f$ being flat, the direct image functors $f_*$
and~$f_!$ take $Y\qcoh^\inj$ into $X\qcoh^\inj$ and $Y\ctrh^\lct_\prj$
into $X\ctrh^\lct_\prj$ (see Corollary~\ref{proj-lct-direct}(b)).

 For any complex $\J^\bu$ over $Y\qcoh^\inj$ there is a natural 
isomorphism of complexes of locally cotorsion contraherent cosheaves
$f_!\fHom_Y(f^*\D_X^\bu,\J^\bu)\simeq\fHom_X(\D_X^\bu,f_*\J^\bu)$
on~$X$ \,\eqref{inj-fHom-projection}.
 Composing this isomorphism with the morphism
$f_!\fHom_Y(\D_Y^\bu,\J^\bu)\rarrow f_!\fHom_Y(f^*\D_X^\bu,\J^\bu)$
induced by the quasi-isomorphism $f^*\D_X^\bu\rarrow\D_Y^\bu$,
we obtain a natural morphism
\begin{equation}  \label{fHom-direct-transformation}
 f_!\fHom_Y(\D_Y^\bu,\J^\bu)\lrarrow\fHom_X(\D_X^\bu,f_*\J^\bu)
\end{equation}
of complexes over $X\ctrh^\lct_\prj$.
 Similarly, according to Section~\ref{compatibility-subsect} for any
complex $\P^\bu$ over $Y\ctrh^\lct_\prj$ there is a natural morphism
$\D_X^\bu\ocn_X f_!\P^\bu\rarrow f_*(f^*\D_X^\bu\ocn_Y\P^\bu)$
of complexes over $X\qcoh$.
 Composing it with the morphism $f_*(f^*\D_X^\bu\ocn_Y\P^\bu)\rarrow
f_*(\D_Y^\bu\ocn_Y\P^\bu)$ induced by the quasi-isomorphism
$f^*\D_X^\bu\rarrow\D_Y^\bu$, we obtain a natural morphism
\begin{equation}  \label{contratensor-direct-transformation}
 \D_X^\bu\ocn_X f_!\P^\bu\lrarrow f_*(\D_Y^\bu\ocn_Y\P^\bu)
\end{equation}
of complexes over $X\qcoh^\inj$.
 The natural morphisms~(\ref{fHom-direct-transformation}\+-%
\ref{contratensor-direct-transformation}) are compatible with
the adjunction~\eqref{fHom-contratensor-adjunction} and with
the compositions of the morphisms of schemes~$f$.

 It is explained in Remark~\ref{adjunction-compatibility-remark} below
what the compatibility with the adjunction means.
 The point is that it suffices to show that
the morphism~\eqref{fHom-direct-transformation} is a homotopy
equivalence for all $\J^\bu\in\Hot(Y\qcoh^\inj)$; then it follows
that the morphism~\eqref{contratensor-direct-transformation}
is a homotopy equivalence for all $\P^\bu\in\Hot(X\ctrh^\lct_\prj)$.
 Conversely, if it is shown that
the morphism~\eqref{contratensor-direct-transformation} is a homotopy
equivalence for all $\P^\bu\in\Hot(X\ctrh^\lct_\prj)$, then it would
follow that the morphism~\eqref{fHom-direct-transformation} is
a homotopy equivalence for all $\J^\bu\in\Hot(Y\qcoh^\inj)$.

 Furthermore, it suffices to show that
the morphism~\eqref{fHom-direct-transformation} is a homotopy
equivalence (or just a quasi-isomorphism over $X\ctrh$) for any
one-term complex $\J^\bu=\J$ over $Y\qcoh^\inj$, or equivalently,
the morphism~\eqref{contratensor-direct-transformation} is a homotopy
equivalence (or just a quasi-isomorphism over $X\qcoh$) for any one-term
complex $\P^\bu=\P$ over $Y\ctrh^\lct_\prj$.
 For this purpose, we notice that the functors $f_!\:Y\ctrh^\lct\rarrow
X\ctrh^\lct$ and $f_*\:Y\qcoh\rarrow X\qcoh$ are exact for any affine
morphism~$f$.
 Besides, the functor $\fHom_Y\:Y\qcoh^\op\times Y\qcoh^\inj
\rarrow Y\ctrh^\lct$ is exact for any semi-separated scheme $Y$, while
the functor $\ocn_Y\:Y\qcoh\times Y\ctrh^\lct_\prj\rarrow Y\qcoh$
is exact for any affine scheme~$Y$.
 On the other hand, the morphism $f^*\D_X^\bu\rarrow\D_Y^\bu$ is
a homotopy equivalence whenever the morphism~$f$ is an open embedding.

 Hence, in view of~\eqref{inj-fHom-projection},
the morphism~\eqref{fHom-direct-transformation} is actually
an isomorphism in the homotopy category whenever either $Y$ is
semi-separated and $f$~is affine, or $f$~is an open embedding.
 In view of~(\ref{contratensor-projection}\+-%
\ref{flat-contratensor-projection}),
the morphism~\eqref{contratensor-direct-transformation} is
an isomorphism in the homotopy category whenever the scheme $Y$
is affine and either the scheme $X$ is semi-separated, or
the morphism~$f$ is an open embedding.

 Moreover, the functor $\ocn_Y\:Y\qcoh\times Y\ctrh_\alf\rarrow Y\qcoh$
is actually exact for any quasi-compact semi-separated scheme $Y$,
as one can see from the natural
isomorphism~\eqref{fHom-contratensor-adjunction} together with the fact
that $\fHom_X(\M,\J)$ is a locally cotorsion contraherent cosheaf
for any quasi-coherent sheaf $\M$ and any injective quasi-coherent sheaf
$\J$ on~$Y$.
 Using formulas~\eqref{contratensor-projection} above
and~\eqref{open-embed-fq-flat-ctrtensor-projection-eqn} below, one can
see that the morphism~\eqref{contratensor-direct-transformation} is
an isomorphism in the homotopy category whenever either $X$ and $Y$ are
semi-separated and $f$~is affine, or $f$~is an open embedding.

 The argument finishes similarly to the proof of
Theorem~\ref{non-semi-separated-direct-images-identified}.
 Cover the scheme $X$ with affine open subschemes $U_\alpha$ and
the scheme $Y$ with affine open subschemes $V_\beta$ so that
for every~$\beta$ there exists~$\alpha$ for which $f(V_\beta)\sub
U_\alpha$.
 Decompose an object $\J\in Y\qcoh^\inj$ into a direct sum of
the direct images of injective quasi-coherent sheaves from~$V_\beta$.
 Since we know~\eqref{fHom-direct-transformation} to be a homotopy
equivalence for the morphisms $V_\beta\rarrow Y$, \ $U_\alpha\rarrow X$,
and $V_\beta\rarrow U_\alpha$, it follows that this map is a homotopy
equivalence of complexes over $X\ctrh^\lct_\prj$ for the morphism~$f$.
 Alternatively, decompose an object $\P\in Y\ctrh^\lct_\prj$ into
a direct sum of the direct images of projective locally cotorsion
contraherent cosheaves from $V_\beta$, etc.
\end{proof}

\begin{rem} \label{adjunction-compatibility-remark}
 Suppose given a pair of adjoint functors $R\:\sC\rarrow\sD$ and
$L\:\sD\rarrow\sC$ between two categories $\sC$ and~$\sD$ (with
the functor $L$ left adjoint to~$R$).
 Suppose further given a similar pair of adjoint functors $R'\:\sC'
\rarrow\sD'$ and $L'\:\sD'\rarrow\sC'$ between two categories
$\sC'$ and~$\sD'$.
 Finally, suppose given two functors $F\:\sC'\rarrow\sC$ and
$G\:\sD'\rarrow\sD$, and two natural transformations
$GR'\rarrow RF$ and $LG\rarrow FL'$.

 In the context of the proof of
Theorem~\ref{smooth-direct-images-identified}, we consider
the categories of complexes $\sC=\Com(X\qcoh^\inj)$ and
$\sC'=\Com(Y\qcoh^\inj)$, while $\sD=\Com(X\ctrh^\lct_\prj)$ and
$\sD'=\Com(Y\ctrh^\lct_\prj)$.
 The adjoint pair $(L,R)$ is formed by the functors
$L=\D_X^\bu\ocn_X{-}$ and $R=\fHom_X(\D_X^\bu,{-})$ and similarly,
$L'=\D_Y^\bu\ocn_Y{-}$ and $R'=\fHom_Y(\D_Y^\bu,{-})$.
 The functors $F$ and $G$ are $F=f_*$ and $G=f_!$.
 The natural transformation $GR'\rarrow RF$
is~\eqref{fHom-direct-transformation}, and
the natural transformation $LG\rarrow FL'$
is~\eqref{contratensor-direct-transformation}.

 What does it mean that the two natural transformations are
compatible with the two adjunctions?
 It means that the latter natural transformation can be obtained from
the former one as the composition $LG\rarrow LGR'L'\rarrow LRFL'
\rarrow FL'$.
 Equivalently, the former natural transformation can be obtained from
the latter one as the composition $GR'\rarrow RLGR'\rarrow RFL'R'
\rarrow RF$.

 Now if this is the case, and if inverting suitable classes of
morphisms in $\sC$, $\sD$, $\sC'$, $\sD'$ turns the adjunctions into
equivalences, and if moreover the functors $F$ and $G$ descend to
well-defined functors between the localized categories, then
the natural transformation $GR'\rarrow RF$ becomes an isomorphism
in the localized category if and only if so becomes the natural
transformation $LG\rarrow FL'$.
\end{rem}

\subsection{Compatibilities for finite and proper morphisms}
 Let $X$ be a Noetherian scheme with a dualizing complex $\D_X^\bu$,
which we will view as a finite complex over $X\qcoh^\inj$.
 Let $f\:Y\rarrow X$ be a finite morphism of schemes and
$f^!\:X\qcoh\rarrow Y\qcoh$ be the special inverse image functor from
Section~\ref{special-inverse-subsection}.
 Then $\D_Y^\bu=f^!\D_X^\bu$ is a dualizing complex on~$Y$
\,\cite[Proposition~V.2.4]{Har}, \cite[Lemma Tag~0AX0]{SP}.

\begin{thm}  \label{finite-morphism-special-inverse-images-identified}
 The equivalences of triangulated categories\/
$\sD^\co(X\qcoh)\simeq\sD^\ctr(X\ctrh^\lct)$ and\/
$\sD^\co(Y\qcoh)\simeq\sD^\ctr(Y\ctrh^\lct)$ from
Theorem~\textup{\ref{co-contra-dualizing}}
or~\textup{\ref{non-semi-separated-co-contra}} related to the choice of
the dualizing complexes\/ $\D_X^\bu$ and\/ $\D_Y^\bu$ on $X$ and $Y$
transform the right derived functor\/ $\boR f^!\:\sD^\co(X\qcoh)
\rarrow\sD^\co(Y\qcoh)$ \textup{\eqref{quasi-derived-special-inverse}}
into the left derived functor\/ $\boL f^*\:\sD^\ctr(X\ctrh^\lct)\rarrow
\sD^\ctr(Y\ctrh^\lct)$ \textup{\eqref{lct-derived-special-inverse}}.
 In other words, the following square diagram of triangulated functors
and triangulated equivalences is commutative: \hbadness=1325
\begin{equation}
\begin{gathered}
 \xymatrix{
  \sD^\co(Y\qcoh) \ar@{=}[r] & \sD^\ctr(Y\ctrh^\lct) \\
  \sD^\co(X\qcoh) \ar@{=}[r] \ar[u]^{\boR f^!}
  & \sD^\ctr(X\ctrh^\lct) \ar[u]_{\boL f^*}
 }
\end{gathered}
\end{equation}
\end{thm}

\begin{proof}
 We will show that the equivalences of homotopy categories
$\Hot(X\qcoh^\inj)\simeq\Hot(X\ctrh^\lct_\prj)$ and
$\Hot(Y\qcoh^\inj)\simeq\Hot(Y\ctrh^\lct_\prj)$ from the proof of
Theorem~\ref{non-semi-separated-co-contra} transform
the functor~$f^!$ into the functor~$f^*$.
 Notice that complexes over $X\qcoh^\inj$ and $X\ctrh^\lct_\prj$
are adjusted to the derived functors $\boR f^!$ and $\boL f^*$
acting between the co/contraderived categories.
 The special inverse image functors $f^!$ and~$f^*$ take $X\qcoh^\inj$
to $Y\qcoh^\inj$ and $X\ctrh^\lct_\prj$ to $Y\ctrh^\lct_\prj$,
as explained in Section~\ref{special-inverse-subsection}.

 First of all, we will need the following base change lemma.

\begin{lem}  \label{special-base-change}
 Let $g\:x\rarrow X$ be a morphism of locally Noetherian schemes
and $f\:Y\rarrow X$ be a finite morphism.
 Set\/ $y=x\times_X Y$, and denote the natural morphisms by
$f'\:y\rarrow x$ and\/ $g'\:y\rarrow Y$.
 Then \par
\textup{(a)} for any quasi-coherent sheaf\/ $\M$ on $X$ there is
a natural morphism of quasi-coherent sheaves $g'{}^*\!\.f^!\M\rarrow
f'{}^!g^*\M$ on~$y$; this map is an isomorphism for any\/ $\M$
whenever the morphism~$g$ is flat; \par
\textup{(b)} whenever the morphism~$g$ is an open embedding, for any
projective locally cotorsion contraherent cosheaf\/ $\P$ on $X$
there is a natural isomorphism of projective locally cotorsion
contraherent cosheaves $g'{}^!f^*\P\simeq f'{}^*g^!\.\P$ on~$y$; \par
\textup{(c)} whenever the morphism~$g$ is quasi-compact, for any
quasi-coherent sheaf\/~$\mm$ on\/~$x$ there is a natural isomorphism
of quasi-coherent sheaves $g'_*f'{}^!\.\mm\simeq f^!g_*\mm$ on~$Y$; \par
\textup{(d)} whenever the morphism~$g$ is flat and quasi-compact,
for any projective locally cotorsion contraherent cosheaf\/~$\p$
on\/~$x$ there is a natural isomorphism of projective locally
cotorsion contraherent cosheaves $g'_!f'{}^*\p\simeq f^*g_!\.\p$ on~$Y$.
\end{lem}

\begin{proof}
 Notice that the functors being composed in parts~(b) and~(d) preserve
the class of projective locally cotorsion contraherent cosheaves by
Corollaries~\ref{lct-proj-local}(a) and~\ref{proj-lct-direct}(b)
together with the discussion in the paragraph preceding this lemma.

 Furthermore, parts~(a) and~(b) were essentially proved in
the computations in Section~\ref{special-inverse-subsection}.
 Indeed, the assertions~(a\+b) are local in $X$ and~$x$, so one can
assume all the four schemes to be affine.
 Then part~(a) holds due to the natural isomorphism of
$(r\ot_R\nobreak S)$\+modules $S\ot_R\Hom_R(r,M)\simeq
\Hom_R(r\;S\ot_RM)$ valid for any commutative ring homomorphisms
$R\rarrow S$ and $R\rarrow r$ making $S$ a flat $R$\+module and $r$
a finitely presented $R$\+module, and for any $R$\+module~$M$.
 Part~(b) holds due to Lemma~\ref{quotient-scalars-cotorsion}(c).

 The assertions of parts~(c\+d) are local in $X$, so one can assume
the schemes $X$ and $Y$ to be affine.
 Then if $x$ is also affine, the assertion~(c) reduces to the obvious
natural isomorphism of $S$\+modules $\Hom_r(r\ot_RS\;m)\simeq
\Hom_R(S,m)$ for any commutative ring homomorphisms $R\rarrow S$
and $R\rarrow r$, and any $r$\+module~$m$.
 The general case of a Noetherian scheme~$x$ is handled by computing
the global sections of the sheaf~$\mm$ in terms of a finite affine
covering~$u_\alpha$ of~$x$ and finite affine coverings of
the intersections $u_\alpha\cap u_\beta$.
 The point is that the global sections of~$\mm$ over~$x$ are computed
as a certain kernel, and the functor $\Hom_R(S,{-})$ preserves kernels.

 The proof of part~(d) is similar to that of part~(c).
 If $x$ is affine, the assertion~(d) reduces to the obvious natural
isomorphism of $S$\+modules $(r\ot_RS)\ot_rp\simeq S\ot_Rp$ for
any $r$\+module~$p$.
 In the general case, the cosections of~$\p$ over~$x$ are computed
as a certain cokernel, and the functor $S\ot_R{-}$ preserves cokernels.
\end{proof}

 For an affine open subscheme $U\sub X$ with the open embedding
morphism $j\:U\rarrow X$, denote by $V$ the full preimage
$U\times_XY\sub Y$ and by $j'\:V\rarrow Y$ its open embedding morphism.
 Let $\I$ and $\J$ be injective quasi-coherent sheaves on~$X$. 
 {Then the composition
$$
 \Hom_X(j_*j^*\I,\J)\lrarrow\Hom_Y(f^!j_*j^*\I,f^!\J)
 \.\simeq\.\Hom_Y(j'_*j'{}^*\!\.f^!\I,f^!\J)
$$ \emergencystretch=3em \hbadness=3350
of the map induced by the functor~$f^!$ with the isomorphism
provided by Lemma~\ref{special-base-change}(a,c) induces
a morphism of $\O(V)$\+modules
\begin{multline*}
 (f^*\!\.\fHom_X(\I,\J))[V] \.=\.
 \O(V)\ot_{\O(U)}\Hom_X(j_*j^*\I,\J) \\ \lrarrow
 \Hom_Y(j'_*j'{}^*\!\.f^!\I,f^!\J) \.=\.
 \fHom_Y(f^!\I,f^!\J)[V],
\end{multline*}
defining a morphism of projective locally cotorsion contraherent
cosheaves $f^*\!\.\fHom_X(\I,\J)\allowbreak
\rarrow\fHom_Y(f^!\I,f^!\J)$ on~$Y$. \par}

 To show that this morphism is an isomorphism, decompose the sheaf $\J$
into a finite direct sum of the direct images of injective
quasi-coherent sheaves $\K$ from embeddings $k\:W\rarrow X$
of affine open subschemes $W\sub X$.
 Then the isomorphisms of Lemma~\ref{special-base-change}(a,c\+d)
together with the isomorphism~\eqref{inj-fHom-projection} reduce
the question to the case of an affine scheme~$X$.
 Indeed, denoting the natural morphisms by $f'\:Y\times_XW\rarrow W$
and $k'\:Y\times_XW\rarrow Y$, we have natural isomorphisms
$$
 f^*\!\.\fHom_X(\I,k_*\K)\simeq f^*k_!\fHom_W(k^*\I,\K)\simeq
 k'_!f'{}^*\!\.\fHom_W(k^*\I,\K)
$$
and
\begin{multline*}
 \fHom_Y(f^!\I,f^!k_*\K)\simeq\fHom_Y(f^!\I,k'_*f'{}^!\K) \\
 \simeq k'_!\fHom_{Y\times_XW}(k'{}^*\!\.f^!\I,f'{}^!\K)\simeq
 k'_!\fHom_{Y\times_XW}(f'{}^!k^*\I,f'{}^!\K).
\end{multline*}
 In the case of an affine scheme $X$, we have
\begin{multline*}
 \O(Y)\ot_{\O(X)}\Hom_{\O(X)}(\I(X),\J(X))\simeq
 \Hom_{\O(X)}(\Hom_{\O(X)}(\O(Y),\I(X)),\.\J(X))\\ \simeq
 \Hom_{\O(Y)}(\Hom_{\O(X)}(\O(Y),\I(X))\;\Hom_{\O(X)}(\O(Y),\J(X))),
\end{multline*}
since $\O(Y)$ is a finitely presented $\O(X)$\+module and $\J(X)$ is
an injective $\O(X)$\+module.

 Alternatively, let $\M$ be a quasi-coherent sheaf and $\P$ be
a projective locally cotorsion contraherent cosheaf on~$X$.
 We keep the notation above for the open embeddings $j\:U\rarrow X$
and $j'\:V\rarrow Y$.
 Then the composition
\begin{multline*}
 j'_*j'{}^*\!\.f^!\M\ot_{\O(V)}(f^*\P)[V] \.\simeq\.
 j'_*j'{}^*\!\.f^!\M\ot_{\O(U)}\P[U] \\ \simeq\.
 f^!j_*j^*\M\ot_{\O(U)}\P[U] \lrarrow
 f^!(j_*j^*\M\ot_{\O(U)}\P[U])
\end{multline*}
provides a natural morphism $j'_*j'{}^*\!\.f^!\M\ot_{\O(V)}(f^*\P)[V]
\rarrow f^!(j_*j^*\M\ot_{\O(U)}\P[U])$ of quasi-coherent sheaves on~$Y$.
 Passing to the inductive limit over $U$ and noticing that
the contratensor product $f^!\M\ocn_Y f^*\P$ can be computed on
the diagram of affine open subschemes $V=U\times_X Y\sub Y$
(see Section~\ref{contratensor-subsect}), we obtain a morphism of
quasi-coherent sheaves $f^!\M\ocn_Y f^*\P\rarrow f^!(\M\ocn_X\P)$
on~$Y$.
 When the quasi-coherent sheaf $\M=\I$ is injective, this is a morphism
of injective quasi-coherent sheaves.

 To show that this morphism is an isomorphism, decompose the cosheaf
$\P$ into a finite direct sum of the direct images of projective
locally cotorsion contraherent cosheaves $\Q$ from embeddings
$k\:W\rarrow X$ of affine open subschemes $W\sub X$.
 Then the isomorphisms of Lemma~\ref{special-base-change}(a,c\+d)
together with the isomorphism~\eqref{flat-contratensor-projection}
reduce the question to the case of an affine scheme~$X$.
 Indeed, in the above notation for the morphisms of schemes, we have
\begin{multline*}
 f^!\M\ocn_Y f^*k_!\.\Q\simeq f^!\M\ocn_Y k'_!f'{}^*\Q \\
 \simeq k'_*(k'{}^*\!\.f^!\M\ocn_{Y\times_XW}f'{}^*\Q)\simeq
 k'_*(f'{}^!k^*\M\ocn_{Y\times_XW}f'{}^*\Q)
\end{multline*}
and
$$
 f^!(\M\ocn_X k_!\.\Q)\simeq f^!k_*(k^*\M\ocn_W\Q)\simeq
 k'_*f'{}^!(k^*\M\ocn_W\Q).
$$
 In the case of an affine scheme $X$, we have
\begin{multline*}
 \Hom_{\O(X)}(\O(Y),\M(X))\ot_{\O(Y)}(\O(Y)\ot_{\O(X)}\P[X]) \\
 \simeq\Hom_{\O(X)}(\O(Y),\M(X))\ot_{\O(X)}\P[X] \simeq
 \Hom_{\O(X)}(\O(Y)\;\M(X)\ot_{\O(X)}\P[X]),
\end{multline*}
since the $\O(X)$\+module $\O(Y)$ is finitely presented and
the $\O(X)$\+module $\P[X]$ is flat.
\end{proof}

 Now let $X$ be a semi-separated Noetherian scheme with a dualizing
complex~$\D_X^\bu$.
 Let $f\:Y\rarrow X$ be a proper morphism (of finite type);
see the references before
Proposition~\ref{proper-finite-flat-dim-extraordinary-inverse}.
 Notice that, the morphism~$f$ being separated, the scheme~$Y$ is
consequently semi-separated (and Noetherian), too.

 Set $\D_Y^\bu$ to be a finite complex over $Y\qcoh^\inj$
quasi-isomorphic to $f^!\D_X^\bu$, where $f^!$ denotes
the triangulated functor~\eqref{loc-noetherian-quasi-adjoint}
right adjoint to the derived direct image functor~$\boR f_*$.
 Then $\D_Y^\bu$ is a dualizing complex on~$Y$
\,\cite[Remark before Proposition~V.8.5]{Har}
or~\cite[Section Tag~0A9Y and Lemma Tag~0AA3]{SP}.
 Notice that the triangulated functors $f^!\:\sD(X\qcoh)\rarrow
\sD(Y\qcoh)$ and $f^!\:\sD^\co(X\qcoh)\rarrow\sD^\co(Y\qcoh)$
agree on bounded below complexes by
Corollary~\ref{co-contra-derived-bounded-restrict}(a).

 The following theorem is supposed to be applied together with
Theorem~\ref{finite-krull-inverse-images-identified}
and/or Corollary~\ref{tensor-cohom-dualizing-compatibility}.

\begin{thm}  \label{proper-inverse-images-identified}
 The equivalences of triangulated categories\/ $\sD^\abs(X\qcoh_\fl)
\simeq\sD^\co(X\qcoh)$ and\/ $\sD^\abs(Y\qcoh_\fl)\simeq
\sD^\co(Y\qcoh)$ from Theorem~\textup{\ref{co-contra-dualizing}}
related to the choice of the dualizing complexes\/ $\D_X^\bu$ and\/
$\D_Y^\bu$ on $X$ and $Y$ transform the triangulated functor
$f^*\:\sD^\abs(X\qcoh_\fl)\rarrow\sD^\abs(Y\qcoh_\fl)$
\textup{\eqref{qcoh-inverse-ffd-sheaves}} into the triangulated functor
$f^!\:\sD^\co(X\qcoh)\rarrow\sD^\co(Y\qcoh)$
\textup{\eqref{loc-noetherian-quasi-adjoint}}.
 In other words, the following square diagram of triangulated functors
and triangulated equivalences is commutative:
\begin{equation}
\begin{gathered}
 \xymatrix{
  \sD^\abs(Y\qcoh_\fl) \ar@{=}[r] & \sD^\co(Y\qcoh) \\
  \sD^\abs(X\qcoh_\fl) \ar@{=}[r] \ar[u]^{f^*}
  & \sD^\co(X\qcoh) \ar[u]^{f^!}
 }
\end{gathered}
\end{equation}
\end{thm}

\begin{proof}
 For any complex of flat quasi-coherent sheaves $\F^\bu$ on $X$
there is a natural isomorphism $f_*(\D_Y^\bu\ot_{\O_Y} f^*\F^\bu)\simeq
f_*(\D_Y^\bu)\ot_{\O_X}\F^\bu$ of complexes over $X\qcoh$
(see~\eqref{qcoh-projection}).
 Hence the adjunction morphism $f_*\D_Y^\bu = \boR f_* f^!\D_X^\bu
\rarrow \D_X^\bu$ in $\sD^+(X\qcoh)\sub\sD^\co(X\qcoh)$ induces
a natural morphism $\boR f_*(\D_Y^\bu\ot_{\O_Y}f^*\F^\bu)=
f_*(\D_Y^\bu\ot_{\O_Y}f^*\F^\bu)\rarrow
\D_X^\bu\ot_{\O_X}\F^\bu$ in $\sD^\co(X\qcoh)$
(see~\eqref{flat-arbitrary-coderived-tensor}).
 Notice that $\D_Y^\bu\ot_{\O_Y}f^*\F^\bu$ is a complex of injective
quasi-coherent sheaves on $Y$, so it is adjusted to $\boR f_*\:
\sD^\co(Y\qcoh)\rarrow\sD^\co(X\qcoh)$
as per~\eqref{loc-noetherian-quasi-direct}.
 We have constructed a natural transformation $\D_Y^\bu\ot_{\O_Y}
f^*\F^\bu \rarrow f^!(\D_X^\bu\ot_{\O_X}\F^\bu)$ of functors
$\sD^\abs(X\qcoh_\fl)\rarrow\sD^\co(Y\qcoh)$.

 Since finite complexes of coherent sheaves $\N^\bu$ on $Y$ form
a set of compact generators of $\sD^\co(Y\qcoh)$
(see Theorem~\ref{contraderived-compactly-generated}(b)), in order
to prove that our morphism is an isomorphism it suffices to show
that the induced morphism of abelian groups
$\Hom_{\sD^\co(Y\qcoh)}(\N^\bu\;\D_Y^\bu\ot_{\O_Y}f^*\F^\bu)
\rarrow\Hom_{\sD^\co(Y\qcoh)}(\N^\bu\;f^!(\D_X^\bu\ot_{\O_X}\F^\bu)
\simeq\Hom_{\sD^\co(X\qcoh)}(\boR f_*\N^\bu\;\D_X^\bu\ot_{\O_X}\F^\bu)$
is an isomorphism for any~$\N^\bu$.
 Both $\D_Y^\bu\ot_{\O_Y}f^*\F^\bu$ and $\D_X^\bu\ot_{\O_X}\F^\bu$
being complexes of injective quasi-coherent sheaves, the Hom into
them in the coderived categories coincides with the one taken in
the homotopy categories of complexes of quasi-coherent sheaves.

 The complexes $\D_Y^\bu$, \ $\D_X^\bu$ and $\N^\bu$ being finite,
the latter Hom only depends on a finite fragment of
the complex~$\F^\bu$.
 This reduces the question to the case of a single coherent sheaf
$\N$ on $Y$ and a single flat quasi-coherent sheaf $\F$ on~$X$;
one has to show that the natural morphism of complexes of abelian
groups $\Hom_Y(\N\;\D_Y^\bu\ot_{\O_Y}f^*\F)\rarrow
\Hom_X(\boR f_*\N\;\D_X^\bu\ot_{\O_X}\F)$ is a quasi-isomorphism.
 Notice that in the case $\F=\O_X$ this is the definition of
$\D_Y^\bu\simeq f^!\D_X^\bu$.

 In the case when there are enough vector bundles (locally free
sheaves of finite rank) on $X$, one can pick a resolution $\L_\bu$
of the flat quasi-coherent sheaf $\F$ by infinite direct sums of
vector bundles and argue as above, reducing the question further to
the case when $\F$ is an infinite direct sum of vector bundles, when
the assertion follows by compactness.
 In the general case, using the \v Cech resolution~\eqref{cech-quasi}
of a flat quasi-coherent sheaf $\F$, one can assume $\F$ to be
the direct image of a flat quasi-coherent sheaf $\G$ from an affine
open subscheme $U\sub X$.

 By the same projection formula~\eqref{qcoh-projection} applied to
the open embeddings $j\:U\rarrow X$ and $j'\:V=U\times_X Y\rarrow Y$
(which are affine morphisms), one has $\D_X^\bu\ot_{\O_X}j_*\G \simeq
j_*(j^*\D_X^\bu\ot_{\O_U}\G)$ and $\D_Y^\bu\ot_{\O_Y}f^*\!\.j_*\G
\simeq\D_Y^\bu\ot_{\O_Y}j'_*f'^*\G\simeq
j'_*(j'^*\D_Y^\bu\ot_{\O_V}f'^*\G)$.
 Here we are also using the base change isomorphism
$f^*\!\.j_*\G\simeq j'_*f'^*\G$, and denote by~$f'$ the morphism
$V\rarrow U$.
 Using the adjunction of the direct and inverse image functors
$j_*$ and $j^*$ together with the isomorphism
$j^*\.\boR f_*\N\simeq\boR f'_*\,j'{}^*\N$ in $\sD^\b(U\coh)$,
one can replace the morphism $f\:Y\rarrow X$ by the morphism
$f'\:V\rarrow U$ into an affine scheme $U$ in the desired assertion.
 Writing it down explicitly, we have
\begin{multline*}
 \Hom_Y(\N\;\D_Y^\bu\ot_{\O_Y}f^*j_*\G)\simeq
 \Hom_Y(\N\;j'_*(j'^*\D_Y^\bu\ot_{\O_V}f'^*\G)) \\
 \simeq\Hom_V(j'{}^*\N\;j'^*\D_Y^\bu\ot_{\O_V}f'^*\G)
\end{multline*}
and
\begin{multline*}
 \Hom_X(\boR f_*\N\;\D_X^\bu\ot_{\O_X}j_*\G)\simeq
 \Hom_X(\boR f_*\N\;j_*(j^*\D_X^\bu\ot_{\O_U}\G)) \\
 \simeq\Hom_U(j^*\boR f_*\N\;j^*\D_X^\bu\ot_{\O_U}\G)\simeq
 \Hom_U(\boR f'_*\,j'{}^*\N\;j^*\D_X^\bu\ot_{\O_U}\G).
\end{multline*}
 In this argument one also needs to know that the natural morphism
$j'^*\D_Y^\bu\simeq j'^*\!\.f^!\D_X^\bu\rarrow f'{}^!j^*\D_X^\bu$
is an isomorphism in $\sD^+(V\qcoh)$.
 This is a result of Deligne~\cite{Del}; see~\cite[Lemma Tag~0ATX]{SP}
or Theorem~\ref{coderived-deligne-extraordinary-inverse} below.

 It remains to point out that there are obviously enough vector
bundles on~$U$.
\end{proof}

\subsection{The extraordinary inverse image}
 Let $f\:Y\rarrow X$ be a separated morphism of finite type between
Noetherian schemes.
 (What we call) the \emph{Hartshorne--Deligne extaordinary inverse image
functor} $f^+\:\sD^\co(X\qcoh)\rarrow\sD^\co(Y\qcoh)$ is a triangulated
functor defined by the following rules:
\begin{enumerate}
\renewcommand{\theenumi}{\roman{enumi}}
\item whenever $f$~is an open embedding, $f^+=f^*$ is the conventional
inverse image functor (\ref{qcoh-inverse-ffd-morphism},~%
\ref{bco-qcoh-inverse-ffd-morphism}) (induced by the exact functor
$f^*\:Y\qcoh\rarrow X\qcoh$ and left adjoint to the derived direct
image functor $\boR f_*\:\sD^\co(Y\qcoh)\rarrow\sD^\co(X\qcoh)$
(\ref{qcoh-direct},~\ref{bco-qcoh-direct},~%
\ref{loc-noetherian-quasi-direct}));
\item whenever $f$~is a proper morphism, $f^+=f^!$ is 
(what we call) the \emph{Lipman--Neeman extraordinary inverse image
functor}~\eqref{loc-noetherian-quasi-adjoint}, i.~e., the functor right
adjoint to the derived direct image functor~$\boR f_*$
(cf.\ Corollary~\ref{neemans-functor-construction});
\item given two composable separated morphisms of finite type
$f\:Y\rarrow X$ and $g\:Z\rarrow Y$ between Noetherian schemes,
one has a natural isomorphism of triangulated functors
$(fg)^+\simeq g^+f^+$.
\end{enumerate}

 Deligne constructs his extraordinary inverse image functor (which
we denote by) $f^+\:\sD^+(X\qcoh)\rarrow\sD^+(Y\qcoh)$ in~\cite{Del}
(cf.~\cite[Theorem~III.8.7]{Har} and~\cite[Section Tag~0A9Y]{SP}).
 Notice that, by Lemma~\ref{derived-bounded-restrict}(a) and
Corollary~\ref{co-contra-derived-bounded-restrict}(a),
the two versions of the Lipman--Neeman extraordinary inverse image
$f^!\:\sD(X\qcoh)\rarrow\sD(Y\qcoh)$ and $f^!\:\sD^\co(X\qcoh)
\rarrow\sD^\co(Y\qcoh)$ have isomorphic restrictions to
$\sD^+(X\qcoh)\sub\sD(X\qcoh)$, $\sD^\co(X\qcoh)$, both acting
from $\sD^+(X\qcoh)$ to $\sD^+(Y\qcoh)$ and being right adjoint
to $\boR f_*\:\sD^+(Y\qcoh)\rarrow\sD^+(X\qcoh)$
\eqref{loc-noetherian-quasi-direct-plus}.
 So there is no ambiguity here.

 According to a counterexample of Neeman's~\cite[Example~6.5]{N-bb},
there \emph{cannot} exist a triangulated functor $f^+\:\sD(X\qcoh)
\rarrow\sD(Y\qcoh)$ defined for all, say, locally closed embeddings
of affine schemes of finite type over a fixed field and satisfying
(i)~for open embeddings, (ii)~for closed embeddings, and
(iii)~for compositions.
 Our aim is to show that there is \emph{no} similar inconsistency in
the rules~(i\+iii) to prevent existence of a functor~$f^+$ acting
between the coderived categories of (unbounded complexes of)
quasi-coherent sheaves.
 Given that, and using the fact of compactifiability of separated
morphisms of finite type between Noetherian schemes~\cite{Conr},
\cite[Theorem Tag~0F41]{SP} (a version of Nagata's theorem,
cf.~\cite{Del}), one can to proceed with the construction of
the functor $f^+\:\sD^\co(X\qcoh)\rarrow\sD^\co(Y\qcoh)$ as
in~\cite[Section Tag~0A9Y]{SP} (cf.~\cite[Sections~5 and~6]{Gai}).

\begin{thm}  \label{coderived-deligne-extraordinary-inverse}
 Let $g\:Y\rarrow X$ and\/ $g'\:V\rarrow U$ be proper morphisms (of
finite type) between Noetherian schemes, forming a commutative diagram
with open embeddings $h\:U\rarrow X$ and\/ $h'\:V\rarrow Y$.
 Then there is a natural isomorphism of triangulated functors
$h'{}^*g^!\simeq g'{}^!h^*\:\sD^\co(X\qcoh)\rarrow\sD^\co(V\qcoh)$.
\end{thm}

\begin{proof}
 We follow the approach in~\cite{Del} (with occasional technical
points borrowed from~\cite{Gai}).
 The argument is based on the results of~\cite[Chapitre~3]{Groth2}.

 Set $V'=U\times_X Y$.
 Then the natural morphism $V\rarrow V'$ is both an open embedding
and a proper morphism (by~\cite[Lemma Tag~01W6]{SP}), that is, $V'$
is a disconnected union of $V$ and $V'\setminus V$.
 Hence one can assume that $V=V'$, i.~e., the square is Cartesian.

 The construction of the derived functor
$\boR f_*$~\eqref{loc-noetherian-quasi-direct} being local in
the base (since the flasqueness/injectivity properties of
quasi-coherent sheaves are local), there is a natural isomorphism
of triangulated functors $h^*\boR g_*\simeq \boR g'_*\,h'^*$.
 Given an object $\M^\bu\in\sD^\co(X\qcoh)$, we now have a natural
morphism $\boR g'_*\, h'{}^*g^!\M^\bu\simeq h^*\boR g_*\.g^!\M^\bu
\rarrow h^*\M^\bu$ in $\sD^\co(U\qcoh)$, inducing a natural morphism
$h'{}^*g^!\M^\bu\rarrow g'{}^!h^*\M^\bu$ in $\sD^\co(V\qcoh)$.
 We have constructed a morphism of functors $h'{}^*g^!\rarrow
g'{}^!h^*$.

 In order to prove that this morphism is an isomorphism in
$\sD^\co(V\qcoh)$, we will show that the induced morphism
$\Hom_{\sD^\co(V\qcoh)}(\L^\bu\;h'{}^*g^!\M^\bu)\rarrow
\Hom_{\sD^\co(V\qcoh)}(\L^\bu\;\allowbreak g'{}^!h^*\M^\bu)$
is an isomorphism of abelian groups for any finite complex of
coherent sheaves $\L^\bu$ on~$V$.
 This would be sufficient in view of the fact that finite complexes
of coherents form a set of compact generators of $\sD^\co(V\qcoh)$
(see Theorem~\ref{contraderived-compactly-generated}(b)).

 For any category $\sC$, we denote by $\Pro\sC$ the category of
pro-objects in $\sC$ (see, e.~g., \cite[Chapter~6]{KS}).
 For our purposes, it suffices to consider projective systems
indexed by the nonnegative integers.
 To any open embedding of Noetherian schemes $j\:W\rarrow Z$ one
assigns an exact functor (between abelian categories)
$j_!\:W\coh\rarrow\Pro(Z\coh)$ defined by the following rule.

 Let $\I\sub\O_Z$ denote the sheaf of ideals corresponding to
some closed subscheme structure on the complement $Z\setminus W$.
 Given a coherent sheaf $\F$ on $W$, pick coherent sheaf $\F'$
on $Z$ with $\F'|_W\simeq\F$.
 The functor~$j_!$ takes the sheaf $W$ to the projective system
formed by the coherent sheaves $\I^n\F'\sub\F'$, \ $n\ge0$ on $Z$
(and their natural embeddings).
 One can check that this pro-object does not depend (up to a natural
isomorphism) on the choice of a coherent extension~$\F'$.
 In fact, one has
\begin{equation} \label{abelian-Deligne-adjunction}
 \Hom_W(\F,j^*\G)\simeq\Hom_Z(j_!\F,\G)
\end{equation}
for any $\G\in Z\qcoh$ \cite[Proposition~4]{Del}.

 Passing to the bounded derived categories, one obtains a triangulated
functor $j_!\:\sD^\b(W\coh)\rarrow\sD^\b(\Pro Z\coh)$.
 Let $\pro\sD^\b(Z\coh)\sub\Pro\sD^\b(Z\coh)$ denote the full
subcategory formed by projective systems of complexes with uniformly
bounded cohomology sheaves.
 The system of cohomology functors $\pro H^i\:\pro\sD^\b(Z\coh)\rarrow
\Pro(Z\coh)$ is conservative (a result applicable to any abelian
category~\cite[Proposition~3]{Del}).
 Furthermore, there is a natural functor $\sD^\b(\Pro Z\coh)\rarrow
\pro\sD^\b(Z\coh)$.
 Composing it with the above functor~$j_!$ and passing to pro-objects
in $\sD^\b(W\coh)$ one constructs the functor $j_!\:\pro\sD^\b(W\coh)
\rarrow\pro\sD^\b(Z\coh)$.

 On the other hand, the functor $h^*\:\sD^\b(Z\coh)\rarrow\sD^\b(W\coh)$
induces a natural functor $\pro h^*\:\pro\sD^\b(Z\coh)\rarrow
\pro\sD^\b(W\coh)$.

\begin{lem}  \label{derived-extbyzero-adjunction}
\textup{(a)} For any objects\/ $\L^\bu\in\pro\sD^\b(W\coh)$ and\/
$\M^\bu\in\sD^\co(Z\qcoh)$ there is a natural isomorphism of abelian
groups\/ $\Hom_{\Pro\sD^\co(W\qcoh)}(\L^\bu,h^*\M^\bu)\simeq
\Hom_{\Pro\sD^\co(Z\qcoh)}(h_!\L^\bu,\M^\bu)$. \par
\textup{(b)} For any objects\/ $\F^\bu\in\pro\sD^\b(W\coh)$ and\/
$\G^\bu\in\pro\sD^\b(Z\coh)$ there is a natural isomorphism of abelian
groups\/ $\Hom_{\pro\sD^\b(W\coh)}(\F^\bu,\pro h^*\G^\bu)\simeq
\Hom_{\pro\sD^\b(Z\coh)}(h_!\F^\bu,\G^\bu)$. \hbadness=2675
\end{lem}

\begin{proof}
 Notice that the assertion of part~(a) is \emph{not} true for
the conventional unbounded derived categories.
 To prove it for the coderived categories, one assumes $\M^\bu$ to be
a complex over $Z\qcoh^\inj$, so that $\Hom$ into $\M^\bu$ in
the coderived category is isomorphic to the similar $\Hom$ in
the homotopy category of complexes.
 It is important here that $h^*\M^\bu$ is a complex over $W\qcoh^\inj$
in these assumptions (while the restriction of a homotopy injective
complex of injective quasi-coherent sheaves on a Noetherian scheme to
an open subscheme \emph{need not} be homotopy injective~\cite{Bel},
\cite[Example~2.9]{Pal}).
 The desired adjunction on the level of (pro)derived categories then
follows from the above adjunction~\eqref{abelian-Deligne-adjunction}
on the level of abelian categories.
 To deduce part~(b), one can use the fact that the natural functor
$\sD^\b(Z\coh)\rarrow\sD^\co(Z\qcoh)$ is fully faithful
(and similarly for~$W$).
\end{proof}

 Given a proper morphism $f\:T\rarrow Z$ (of finite type) between
Noetherian schemes, there is the direct image functor $\boR f_*\:
\sD^\b(T\coh)\rarrow\sD^\b(Z\coh)$ \cite[Th\'eor\`eme~3.2.1]{Groth2},
\cite[Proposition Tag~02O5]{SP}.
 Passing to pro-objects, one obtains the induced functor
$\pro\boR f_*\:\pro\sD^\b(T\coh)\rarrow\pro\sD^\b(Z\coh)$.
{\hbadness=1400\par}

\begin{lem}  \label{derived-extbyzero-commutation}
 Let $V\rarrow U$, $Y\rarrow X$ be a Cartesian square formed by
proper morphisms of Noetherian schemes $g\:Y\rarrow X$ and\/
$g'\:V\rarrow U$, and open embeddings $h\:U\rarrow X$ and\/
$h'\:V\rarrow Y$ (as above).
 Then there is a natural isomorphism of functors\/
$\pro\boR g_*\circ h'_!\simeq h_!\circ\pro\boR g'_*\:
\pro\sD^\b(V\coh)\rarrow\pro\sD^\b(X\coh)$.
\end{lem}

\begin{proof}
 Given an object $\L^\bu\in\pro\sD^b(V\coh)$, one has a natural
morphism $\pro\boR g'_*\L^\bu\allowbreak\rarrow
\pro\boR g'_*\pro h'{}^*\,h'_!\L^\bu\simeq
\pro h^*\pro\boR g_*\.h'_!\L^\bu$ in $\pro\sD^\b(U\coh)$, inducing
a natural morphism $h_!\pro\boR g'_*\L^\bu\rarrow
\pro\boR g_*\.h'_!\L^\bu$ in $\pro\sD^\b(X\coh)$
by Lemma~\ref{derived-extbyzero-adjunction}(b).
 We have constructed a morphism of functors $h_!\circ\pro\boR g'_*
\rarrow\pro\boR g_*\circ h'_!$.
 In order to check that this morphism is an isomorphism, it suffices
to consider the case of $\L^\bu\in\sD^\b(V\coh)$ and show that
the morphism in question becomes an isomorphism after applying
the cohomology functors $\pro H^i$ taking values in $\Pro(X\coh)$.

 Furthermore, one can restrict oneself to the case of a single
coherent sheaf $\L\in V\coh$.
 So we have to show that the morphisms $h_!\.\boR^ig'_*(\L)\rarrow
\Pro\boR^ig_*(h'_!\L)$ are isomorphisms in $\Pro(X\coh)$ for any
coherent sheaf $\L$ on~$V$ and all $i\ge0$.
 Finally, one can assume $X$ to be an affine scheme (as all
the constructions of the functors involved are local in the base,
and the property of a morphism in $\Pro(X\coh)$ to be
an isomorphism is local in~$X$).

 Let $\L'$ be a quasi-coherent sheaf on $Y$ extending~$\L$.
 Set $R=\O(X)$ and $I=\O_X(\I)$, where $\I$ is a sheaf of ideals
in $\O_X$ corresponding to a closed subscheme structure on
$X\setminus U$.
 The question reduces to showing that the natural morphism between
the pro-objects represented by the projective systems $I^nH^i(Y,\L')$
and $H^i(Y,I^n\L')$ is an isomorphism in $\Pro(R\modl)$.
 According to~\cite[Corollaire~3.3.2]{Groth2}
or~\cite[Lemma Tag~02O8]{SP}, the $(\bigoplus_{n=0}^\infty I^n)$\+module
$\bigoplus_{n=0}^\infty H^i(Y,I^n\L')$ is finitely generated;
the desired assertion is deduced straightforwardly from this result.
\end{proof}

 Now we can finish the proof of the theorem.
 Given a finite complex $\L^\bu$ over $V\coh$ and a complex $\M^\bu$
over $X\qcoh$, one has
\begin{multline*}
\Hom_{\sD^\co(V\qcoh)}(\L^\bu\;h'{}^*g^!\M^\bu)\simeq
\Hom_{\Pro\sD^\co(Y\qcoh)}(h'_!\L^\bu\;g^!\M^\bu) \\
\simeq\Hom_{\Pro\sD^\co(X\qcoh)}(\pro\boR g_*\.h'_!\L^\bu\;\M^\bu)
\end{multline*}
and
\begin{multline*}
\Hom_{\sD^\co(V\qcoh)}(\L^\bu\;g'{}^!h^*\M^\bu)\simeq
\Hom_{\sD^\co(U\qcoh)}(\boR g'_*\L^\bu\;h^*\M^\bu) \\
\simeq\Hom_{\Pro\sD^\co(X\qcoh)}(h_!\.\boR g'_*\L^\bu\;\M^\bu)
\end{multline*}
by Lemma~\ref{derived-extbyzero-adjunction}(a), so it remains to apply
Lemma~\ref{derived-extbyzero-commutation}.
\end{proof}

 The following theorem, which is the main result of this section,
follows from the several previous results.
 Let $f\:Y\rarrow X$ be a morphism of finite type between
Noetherian schemes.
 Let $\D_X^\bu$ be a dualizing complex on~$X$; set $\D_Y^\bu$ to be
a finite complex over $Y\qcoh^\inj$ quasi-isomorphic to $f^+\D_X^\bu$
\cite[Remark before Proposition~V.8.5]{Har}, \cite[Lemma Tag~0AA3]{SP}.

\begin{thm}
 The equivalences of categories\/ $\sD^\co(X\qcoh)\simeq
\sD^\ctr(X\ctrh)$ and\/ $\sD^\co(Y\qcoh)\simeq\sD^\ctr(Y\ctrh)$
from Theorem~\textup{\ref{co-contra-dualizing}}
or~\textup{\ref{non-semi-separated-co-contra}}
related to the choice of the dualizing complexes\/ $\D_X^\bu$ and\/
$\D_Y^\bu$ on $X$ and $Y$ transform the Hartshorne--Deligne
extraordinary inverse image functor $f^+\:\sD^\co(X\qcoh)\rarrow
\sD^\co(Y\qcoh)$ into the functor~$f^*\:\sD^\ctr(X\ctrh)\rarrow
\sD^\ctr(Y\ctrh)$ \textup{\eqref{loc-noetherian-contra-adjoint}}
left adjoint to the derived functor\/
$\boL f_!$~\textup{\eqref{loc-noetherian-contra-direct}},
at least, in either of the following two situations:
\begin{enumerate}
\renewcommand{\theenumi}{\alph{enumi}}
\item $f$ is a separated morphism of semi-separated Noetherian
schemes that can be factorized into an open embedding followed
by a proper morphism;
\item $f$ is morphism of Noetherian schemes that can be factorized
into a finite morphism followed by an open embedding followed by
a smooth proper morphism.
\end{enumerate}
 In other words, the following square diagram of triangulated functors
and triangulated equivalences is commutative:
\begin{equation}
\begin{gathered}
 \xymatrix{
  \sD^\co(Y\qcoh) \ar@{=}[r] & \sD^\ctr(Y\ctrh) \\
  \sD^\co(X\qcoh) \ar@{=}[r] \ar[u]^{f^+}
  & \sD^\ctr(X\ctrh) \ar[u]_{f^*}
 }
\end{gathered}
\end{equation}
\end{thm}

\begin{proof}
 Notice that all separated morphisms of (semi-separated) Noetherian
schemes satisfy~(a) by~\cite{Conr} or~\cite[Theorem Tag~0F41]{SP}.
 According to
Theorem~\ref{finite-morphism-special-inverse-images-identified},
the desired assertion is true for finite morphisms of Noetherian
schemes with dualizing complexes $f\:Y\rarrow X$.
 Comparing Theorem~\ref{finite-krull-inverse-images-identified}(b)
with Theorem~\ref{proper-inverse-images-identified}, one comes to
the same conclusion in the case of a proper morphism
$f\:Y\rarrow X$ (of finite type) between semi-separated
Noetherian schemes with dualizing complexes.

 On the other hand, let us consider the case of a smooth
morphism~$f$.
 Then it is essentially known (\cite[Corollary~VII.4.3]{Har},
\cite[Theorem~5.4]{N-bb}, \cite[Lemma Tag~0B6U, part~(e)]{SP},
Proposition~\ref{proper-finite-flat-dim-extraordinary-inverse}
above) that the Hartshorne--Deligne extraordinary inverse image
functor $f^+\:\sD^\co(X\qcoh)\rarrow\sD^\co(Y\qcoh)$
only differs from the conventional inverse image functor
$f^*\:\sD^\co(X\qcoh)\rarrow\sD^\co(Y\qcoh)$ by a shift and a twist.
 Namely, for any complex $\M^\bu$ over $X\qcoh$
one has $f^!\M^\bu\simeq \omega_{Y/X}[d]\ot_{\O_Y}f^*\M^\bu$,
where $d$~is the relative dimension and $\omega_{Y/X}$ is
the line bundle of relative top forms on~$Y$.
 In particular, in our present notation one would have
$\D_Y^\bu\simeq\omega_{Y/X}[n]\ot_{\O_Y}f^*\D_X^\bu$
(up to a quasi-isomorphism of finite complexes over $Y\qcoh$).

 All the references in the previous paragraph, with the exception of
our Proposition~\ref{proper-finite-flat-dim-extraordinary-inverse},
treat either the case of bounded below derived categories $\sD^+$,
or the conventional unbounded derived categories~$\sD$.
 Proposition~\ref{proper-finite-flat-dim-extraordinary-inverse}
applies to the coderived categories (which we are interested in),
but assumes the morphism~$f$ to be proper in addition to smooth.
 In the case of an open embedding~$f$, one has $f^+=f^*$ by
the definition.
 Hence our assumptions in part~(b).

 By Theorem~\ref{smooth-direct-images-identified},
the equivalences of triangulated categories
$\sD^\co(Y\qcoh)\simeq\sD^\ctr(Y\ctrh)$ and
$\sD^\co(X\qcoh)\simeq\sD^\ctr(X\ctrh)$ related to the choice of
the dualizing complexes $f^*\D_X^\bu$ and $\D_X^\bu$ on $Y$ and $X$
transform the functor $\boR f_*\:\sD^\co(Y\qcoh)\rarrow
\sD^\co(X\qcoh)$ \eqref{loc-noetherian-quasi-direct} into the functor
$\boL f_!\:\sD^\ctr(Y\ctrh)\rarrow\sD^\ctr(X\ctrh)$
\eqref{loc-noetherian-contra-direct}.
 Passing to the left adjoint functors, we conclude that the same
equivalences transform the functor $f^*\:\sD^\co(X\qcoh)\rarrow
\sD^\co(Y\qcoh)$ (induced by the exact functor $f^*\:X\qcoh\allowbreak
\rarrow Y\qcoh$, cf.~(\ref{qcoh-inverse-ffd-morphism},
\ref{bco-qcoh-inverse-ffd-morphism})) into
the functor $f^*\:\sD^\ctr(X\ctrh)\rarrow\sD^\ctr(Y\ctrh)$
\eqref{loc-noetherian-contra-adjoint}.
 
 It remains to take the twist and the shift into account in order
to deduce the desired assertion for the smooth morphism~$f$.
 As a particular case, the above argument also covers the situation
when $f$~is an open embedding.
\end{proof}

\Section{Semiderived Categories}  \label{semiderived-sect}

\subsection{Antilocal definition of semicoderived category}
\label{antilocal-semicoderived-subsect}
 The \emph{semiderived category} is a generic term for two dual
constructions: the \emph{semicoderived} category, and
the \emph{semicontraderived} one.
 In the context of this book, the construction of the semicoderived
category is applied to quasi-coherent sheaves, and the construction
of the semicontraderived category to contraherent cosheaves.

 Accordingly, we will sometimes drop the co/contra prefix, and speak
simply of the semiderived category when it is clear from the context
whether we are dealing with quasi-coherent sheaves or with
contraherent cosheaves.
 Conversely, we will occasionally use the term ``semicoderived'' to
indicate that we are working with quasi-coherent sheaves (as in
the title of this section), and the term ``semicontraderived'' to mean
that contraherent cosheaves are involved.

 The constructions of the semiderived categories can be based either
on the constructions of the co/contraderived categories in
the sense of Positselski, or on the co/contraderived categories in
the sense of Becker.
 Following~\cite[Appendix~A]{Psemten} (as opposed to the exposition
in the main body of the book~\cite{Psemten}), we use
the co/contraderived categories in the sense of Becker in this
chapter.

 There are two equivalent constructions of the semiderived categories
discussed in this chapter: the local and the antilocal one.
 The local definition of the semi(co)derived category of quasi-coherent
sheaves on $Y$ for a morphism of schemes $f\:Y\rarrow X$ with
a quasi-compact semi-separated scheme $X$ was worked out
in~\cite[Section~A.3]{Psemten}.
 Let briefly recall it before proceeding to dissuss the antilocal
definition.

 We start with a special case described in the following lemma.

\begin{lem} \label{affine-locality-semicoacyclicity}
 Let $\pi\:Y\rarrow X$ be a morphism of affine schemes, and let
$Y=\bigcup_\alpha V_\alpha$ be an affine open covering of
the affine scheme~$Y$.
 Denote by $j_\alpha\:V_\alpha\rarrow Y$ the open embedding morphisms.
 Let $\N^\bu$ be a complex of quasi-coherent sheaves on~$Y$.
 Then the complex of quasi-coherent sheaves $\pi_*\N^\bu$ on $X$ is
Becker-coacyclic in $X\qcoh$ if and only if, for every index~$\alpha$,
the complex of quasi-coherent sheaves $\pi_*j_\alpha{}_*j_\alpha^*
\N^\bu$ on $X$ Becker-coacyclic in $X\qcoh$.
\end{lem}

\begin{proof}
 This is~\cite[Lemma~A.15]{Psemten}.
 For a proof of the contraherent cosheaf version, see
Lemmas~\ref{colocalization-becker-semicontraacyclic}\+-%
\ref{affine-covered-by-affines-semicontraacyclic} below.
\end{proof}

 Let $X$ be a quasi-compact semi-separated scheme and
$\pi\:Y\rarrow X$ be a morphism of schemes.
 Choose an open covering $Y=\bigcup_\alpha V_\alpha$ of the scheme
$Y$ such that the compositions $V_\alpha\rarrow Y\rarrow X$ are
affine morphisms of schemes (e.~g., one can take any affine open
covering of~$Y$).
 Let $j_\alpha\:V_\alpha\rarrow Y$ denote the open embedding morphism.
 We say that a complex $\N^\bu$ of quasi-coherent sheaves on $Y$
is \emph{semi}(\emph{co})\emph{acyclic} (over~$X$) if, for every
index~$\alpha$, the complex $\pi_*j_\alpha{}_*j_\alpha^*\N^\bu$
of quasi-coherent sheaves on $X$ is Becker-coacyclic.
 According
to~\cite[Proposition~A.18\,(a)\,$\Leftrightarrow$\,(f)]{Psemten},
this property of a complex $\N^\bu$ on $Y$ does not depend on
the choice of an open covering $Y=\bigcup_\alpha V_\alpha$.
 This assertion is more general than
Lemma~\ref{affine-locality-semicoacyclicity}; a proof of
the contraherent cosheaf version can be found in
Corollary~\ref{morphism-into-qcom-ssep-semicontraacyclic} below.

 Let us issue a \emph{warning} that our terminology is highly
misleading: the semiacyclicity is a stronger property than
the acyclicity.
 In fact, any Becker-coacyclic complex in $Y\qcoh$ is semi(co)acyclic
over $X$, and any semi(co)acyclic complex is acyclic in $Y\qcoh$.
 The converse implications do not hold in general.
 For example, if $X=Y$ and $\pi=\id$, then the semiacyclic complexes
are precisely the coacyclic ones
(by Corollary~\ref{qcoh-becker-coacyclicity-is-local}); while if
$X$ is a regular Noetherian scheme of finite Krull dimension,
then the semiacyclic complexes are precisely the acyclic ones
(by Lemma~\ref{psemi-remark21}).
 We refer to~\cite[Remark~A.19]{Psemten} for further details and
to Remark~\ref{semicontraacyclicity-remark} below for a contraherent
cosheaf version.

 We will denote the full subcategory of semiacyclic complexes of
quasi-coherent sheaves on $Y$ (relative to~$X$) by
$\Acycl^\si_X(Y\qcoh)$.
 So we have
$$
 \Acycl^\bco(Y\qcoh)\sub\Acycl^\si_X(Y\qcoh)\sub\Acycl(Y\qcoh).
$$

 The full subcategory of semiacyclic complexes $\Acycl^\si_X(Y\qcoh)$
is a thick subcategory closed under infinite direct sums in
the homotopy category $\Hot(Y\qcoh)$.
 Viewed as a full subcategory in the abelian category of complexes
$\Com(Y\qcoh)$, the full subcategory of semiacyclic complexes is
closed under the kernels of epimorphisms, the cokernels of
monomorphisms, extensions, and infinite direct sums.
 (Recall that all absolutely acyclic complexes in $Y\qcoh$
are Becker-coacyclic by
Lemma~\ref{Positselski-trivial-are-Becker-trivial}(a), hence also
semiacyclic.
 Furthermore, the class of Becker-coacyclic complexes in $X\qcoh$
is closed under infinite direct sums.)

 The triangulated Verdier quotient category
$$
 \sD^\si_X(Y\qcoh)=\Hot(Y\qcoh)/\Acycl^\si_X(Y\qcoh)
$$
is called the \emph{semi}(\emph{co})\emph{derived category} of
quasi-coherent sheaves on $Y$ (relative to~$X$).

\begin{lem} \label{semicoacyclic-direct-image}
 Let $X$ be a quasi-compact semi-separated scheme, $\pi\:Y\rarrow X$ be
a morphism of schemes, and $g\:Z\rarrow Y$ be an affine morphism of
schemes.
 Then the direct image functor $g_*\:Z\qcoh\rarrow Y\qcoh$ takes
semiacyclic complexes of quasi-coherent sheaves on~$Z$ (relative to
the morphism $\pi g\:Z\rarrow X$) to semiacyclic complexes of
quasi-coherent sheaves on~$Y$ (relative to the morphism
$\pi\:Y\rarrow X$).
 In other words, $g_*\Acycl^\si_X(Z\qcoh)\sub\Acycl^\si_X(Y\qcoh)$.
 Conversely, given a complex\/ $\L^\bu$ in $Z\qcoh$, if the complex
$g_*\L^\bu$ on $Y$ is semiacyclic over $X$, then the complex\/ $\L^\bu$
on $Z$ is also semiacyclic over~$X$.
\end{lem}

\begin{proof}
 Let $Y=\bigcup_\alpha V_\alpha$ be an open covering of $Y$ with
open embedding morphisms $j_\alpha\:V_\alpha\rarrow Y$ such that
the compositions $\pi j_\alpha\:V_\alpha\rarrow X$ are affine
morphisms.
 Put $W_\alpha=V_\alpha\times_YZ$, and denote the natural morphisms
by $j_\alpha'\:W_\alpha\rarrow Z$ and $g_\alpha'\:W_\alpha\rarrow
V_\alpha$.
 Then the compositions $\pi g j_\alpha'=\pi j_\alpha g_\alpha'\:
W_\alpha\rarrow X$ are affine morphisms (because
the morphisms~$g_\alpha'$ are affine).
 It remains to point out the natural isomorphisms
$\pi_*j_\alpha{}_*j_\alpha^*g_*\L^\bu\simeq
\pi_*j_\alpha{}_*g'_\alpha{}_*j^{\prime*}_\alpha\L^\bu\simeq
\pi_*g_*j'_\alpha{}_*j^{\prime*}_\alpha\L^\bu$ of complexes of
quasi-coherent sheaves on $X$ holding for any complex of
quasi-coherent sheaves $\L^\bu$ on~$Z$.
\end{proof}

 The following lemma and corollary explain the antilocal point of view
on the semicoderived category.

\begin{lem} \label{semicoacyclic-of-injectives-antilocal}
 Let $\pi\:Y\rarrow X$ be a morphism of quasi-compact semi-separated
schemes, and let $Y=\bigcup_\alpha V_\alpha$ be a finite affine open
covering of~$Y$.
 Let\/ $\J^\bu$ be a complex of injective quasi-coherent sheaves on~$Y$.
 Then the complex\/ $\J^\bu$ is semiacyclic over $X$ if and only if
it is a direct summand of a finitely iterated extension of the direct
images of complexes of injective quasi-coherent sheaves on $V_\alpha$
that are semiacyclic over~$X$.
\end{lem}

\begin{proof}
 The argument is based on
Theorem~\ref{quasi-coherent-gluing-theorem} with
Remark~\ref{quasi-coherent-gluing-nonuniversal-remark}.
 For any affine open subscheme $U\sub Y$, put
$\sE_U=\Acycl^\si_X(U\qcoh)\sub\sK_{\O(U)}=\Com(\O(U)\modl)$.
 Put $\sF_U=\Acycl^\bco(U\qcoh)\sub\sE_U$.
 Then the system of classes $\sE_U$ is local essentially by
the definition, and Lemma~\ref{semicoacyclic-direct-image} implies
that it is satisfies the direct image condition; so it is very local
(with respect to identity open embeddings of affine open subschemes
of~$Y$).
 The system of classes $\sF_U$ is very local by
Lemma~\ref{qcoh-affine-direct-image-becker-coacyclic} and
Corollary~\ref{qcoh-becker-coacyclicity-is-local}.

 By Theorem~\ref{quasi-coherent-becker-coderived},
the pair of classes ($\Acycl^\bco(U\qcoh)$, $\Com(U\qcoh^\inj)$)
is a hereditary complete cotorsion pair in the abelian category
$\Com(U\qcoh)$.
 By Lemmas~\ref{restricting-hereditary-cotorsion}
and~\ref{restricting-cotorsion-pairs-lemma}(a), this cotorsion
pair restricts to the exact subcategory $\sE_U=\Acycl^\si_X(U\qcoh)\sub
\Com(U\qcoh)$, providing a hereditary complete cotorsion pair
$(\sF_U,\sC(U))$ in $\sE_U$ with the right class
$\sC(U)=\Com(U\qcoh^\inj)\cap\Acycl^\si_X(U\qcoh)$.
 Similarly we obtain a hereditary complete cotorsion pair $(\sF,\sC)$
in the exact category $\Acycl^\si_X(Y\qcoh)$ with the left class
$\sF=\Acycl^\bco(Y\qcoh)$ and the right class
$\sC=\Com(Y\qcoh^\inj)\cap\Acycl^\si_X(Y\qcoh)$.

 Now Remark~\ref{quasi-coherent-gluing-nonuniversal-remark} is
applicable.
 It produces a hereditary complete cotorsion pair $(\sF_Y,\sC(Y))$
in the exact category $\sE_Y=\Acycl^\si_X(Y\qcoh)$ of localy\+$\sE$
complexes of quasi-coherent sheaves on~$Y$.
 The left class $\sF_Y$ is the class of locally\+$\sF$ complexes of
quasi-coherent sheaves on $Y$, so $\sF_Y=\Acycl^\bco(Y\qcoh)=\sF$
by Corollary~\ref{qcoh-becker-coacyclicity-is-local}.
 Hence the right class is 
$\sC(Y)=\sC=\Com(Y\qcoh^\inj)\cap\Acycl^\si_X(Y\qcoh)$.

 On the other hand, Theorem~\ref{quasi-coherent-gluing-theorem} with
Remark~\ref{quasi-coherent-gluing-nonuniversal-remark} provide
a description of the class $\sC(Y)$ as the class of all direct
summands of finitely iterated extensions of direct images of
complexes from $\sC(V_\alpha)$.
 Comparing these two descriptions of the class $\sC(Y)$, we arrive
to the assertion of the lemma.
\end{proof}

\begin{cor} \label{antilocal-description-of-semicoacyclic}
 Let $\pi\:Y\rarrow X$ be a morphism of quasi-compact semi-separated
schemes, and let $Y=\bigcup_\alpha V_\alpha$ be a finite affine open
covering of~$Y$.
 Denote by $j_\alpha\:V_\alpha\rarrow Y$ the open embedding morphisms.
 Then the full subcategory of semiacyclic complexes\/
$\Acycl^\si_X(Y\qcoh)$ coincides with the minimal thick subcategory
of\/ $\Hot(Y\qcoh)$ containing the subcategories
$j_\alpha{}_*\Acycl^\si_X(V_\alpha\qcoh)$ for all the indices~$\alpha$
and the subcategory of absolutely acyclic complexes\/
$\Acycl^\abs(Y\qcoh)$.
\end{cor}

\begin{proof}
 Denote temporarily the minimal thick subcategory of $\Hot(Y\qcoh)$
containing $\Acycl^\abs(Y\qcoh)$ and
$j_\alpha{}_*\Acycl^\si_X(V_\alpha\qcoh)$ for all~$\alpha$ by~$\sT$.
 Then $\sT\sub\Acycl^\si_X(Y\qcoh)$ by
Lemma~\ref{semicoacyclic-direct-image}.
 To prove the converse inclusion, consider a complex $\N^\bu
\in\Acycl^\si_X(Y\qcoh)$.
 By Theorem~\ref{quasi-coherent-becker-coderived}, there exists
a complex of injective quasi-coherent sheaves $\J^\bu$ on $Y$ together
with a morphism of complexes of quasi-coherent sheaves $\N^\bu\rarrow
\J^\bu$ with a Becker-coacyclic cone.
 As any Becker-coacyclic complex in $Y\qcoh$ is semiacyclic, it
follows that the complex $\J^\bu$ is semiacyclic.
 By Lemma~\ref{semicoacyclic-of-injectives-antilocal}, we have
$\J^\bu\in\sT$ (notice that the direct images of injective
quasi-coherent sheaves under open embeddings of schemes are injective,
and any extension of complexes of injective quasi-coherent sheaves
is termwise split, so it corresponds to a distinguished triangle
in $\Hot(Y\qcoh)$).
 It remains to show that $\Acycl^\bco(Y\qcoh)\sub\sT$.

 Consider a complex $\cB^\bu\in\Acycl^\bco(Y\qcoh)$.
 By Corollary~\ref{qcoh-dil-cta-derived-equiv-cor}(b) (for $\bst=\abs$),
there exists a complex of contraadjusted quasi-coherent sheaves
$\C^\bu$ on $Y$ together with a morphism of complexes $\cB^\bu\rarrow
\C^\bu$ with an absolutely acyclic cone.
 Since the absolutely acyclic complexes are Becker-coacyclic,
the complex $\C^\bu$ is also Becker-coacyclic in $Y\qcoh$,
or equivalently, in $Y\qcoh^\cta$.
 By Corollary~\ref{becker-coacyclic-cta-cot-antilocal}(a), it follows
that the complex $\C^\bu$ (hence also the complex~$\cB^\bu$) belongs
to the minimal thick subcategory of $\Hot(Y\qcoh)$ containing
$\Acycl^\abs(Y\qcoh)$ and $j_\alpha{}_*\Acycl^\bco(V_\alpha\qcoh)$
for all indices~$\alpha$.
\end{proof}

 By a \emph{semi}(\emph{co})\emph{acyclic complex of dilute
quasi-coherent sheaves} on $Y$ we will mean a complex in $Y\qcoh^\dil$
that is semiacyclic as a complex in $Y\qcoh$.
 The same applies to complexes of contraadjusted, cotorsion, or
injective quasi-coherent sheaves on~$Y$.
 We will denote the full subcategories of semiacyclic complexes of
dilute, contraadjusted, cotorsion, and injective quasi-coherent
sheaves on~$Y$ (relative to~$X$) by $\Acycl^\si_X(Y\qcoh^\dil)$,
\ $\Acycl^\si_X(Y\qcoh^\cta)$, \ $\Acycl^\si_X(Y\qcoh^\cot)$, and
$\Acycl^\si_X(Y\qcoh^\inj)$.

\begin{cor} \label{semicoacyclic-dil-cta-cot-direct-image-to-X}
 Let $\pi\:Y\rarrow X$ be a morphism of quasi-compact semi-separated
schemes.  Then \par
\textup{(a)} the direct image functor $\pi_*\:Y\qcoh\rarrow X\qcoh$
takes semiacyclic complexes of dilute quasi-coherent sheaves on~$Y$
(relative to~$X$) to Becker-coacyclic complexes of dilute quasi-coherent
sheaves on~$X$; \par
\textup{(b)} the direct image functor $\pi_*\:Y\qcoh\rarrow X\qcoh$
takes semiacyclic complexes of contraadjusted quasi-coherent sheaves
on~$Y$ (relative to~$X$) to Becker-coacyclic complexes of contraadjusted
quasi-coherent sheaves on~$X$; \par
\textup{(c)} the direct image functor $\pi_*\:Y\qcoh\rarrow X\qcoh$
takes semiacyclic complexes of cotorsion quasi-coherent sheaves on~$Y$
(relative to~$X$) to Becker-coacyclic complexes of cotorsion
quasi-coherent sheaves on~$X$.
\end{cor}

\begin{proof}
 The direct image functor~$\pi_*$ preserves the diluteness,
contraadjustedness, and cotorsion properties of quasi-coherent sheaves
by Corollaries~\ref{cta-cot-direct} and~\ref{dilute-direct}.
 Consequently, part~(c) follows from part~(b).

 To prove part~(b), we recall that it was essentially shown in
the proof of Corollary~\ref{antilocal-description-of-semicoacyclic}
that $\Acycl^\si_X(Y\qcoh^\cta)$ is the thick subcategory of
$\Hot(Y\qcoh^\cta)$ generated by
$j_\alpha{}_*\Acycl^\si_X(V_\alpha\qcoh^\cta)$ and
$\Acycl^\abs(Y\qcoh^\cta)$.
 By the definition of semicoacyclicity for the affine morphism
$\pi j_\alpha\:V_\alpha\rarrow X$, the functor~$\pi_*$ takes
$j_\alpha{}_*\Acycl^\si_X(V_\alpha\qcoh)$ into
$\Acycl^\bco(X\qcoh)$ (hence it also takes
$j_\alpha{}_*\Acycl^\si_X(V_\alpha\qcoh^\cta)$ into
$\Acycl^\bco(X\qcoh^\cta)$).
 Since the functor $\pi_*\:Y\qcoh^\cta\rarrow X\qcoh^\cta$ is exact by
Corollary~\ref{cta-cot-direct}(a), it takes $\Acycl^\abs(Y\qcoh^\cta)$ into $\Acycl^\abs(X\qcoh^\cta)$.

 Part~(a) is deduced from part~(b) in the way similar to the proof
of Corollary~\ref{dil-cta-cot-direct-image-becker-coacyclic}(a).
 Consider a semiacyclic complex $\N^\bu$ in $Y\qcoh^\dil$.
 There exists a finite acyclic complex of complexes
$0\rarrow\N^\bu\rarrow\cP^{0,\bu}\rarrow\cP^{1,\bu}
\rarrow\dotsb\rarrow\cP^{N-1,\bu}\rarrow\cQ^\bu\rarrow0$ in $Y\qcoh$
such that $\cP^{n,\bu}$ is a contractible complex in $Y\qcoh^\cta$ for
every $0\le n\le N-1$, while $\cQ^\bu$ is a complex in $Y\qcoh^\cta$.
 Now it follows that the complex $\cQ^\bu$ is semiacyclic; and by
part~(b) we can conclude that the complex $\pi_*\cQ^\bu$ is
Becker-coacyclic in $X\qcoh^\cta$.
 Therefore, the complex $\pi_*\N^\bu$ is Becker-coacyclic in
$X\qcoh^\dil$ (because $0\rarrow f_*\N^\bu\rarrow f_*\cP^{\bu,\bu}
\rarrow f_*\cQ^\bu\rarrow0$ is an acyclic finite complex of complexes
in $X\qcoh^\dil$ by Corollary~\ref{dilute-direct}).
\end{proof}

 The right derived functor of direct image
\begin{equation} \label{semicoderived-to-X-derived-direct-image}
 \boR\pi_*\:\sD^\si_X(Y\qcoh)\rarrow\sD^\bco(X\qcoh)
\end{equation}
is constructed in the following way.
 In view of (the proof of)
Corollary~\ref{qcoh-dil-cta-derived-equiv-cor}(a) and
Lemma~\ref{pkoszul-lemma16}(b), the natural functor
$$
 \Hot(Y\qcoh^\dil)/\Acycl^\si_X(Y\qcoh^\dil)\lrarrow
 \sD^\si_X(Y\qcoh)
$$
is an equivalence of triangulated categories.
 By Corollary~\ref{semicoacyclic-dil-cta-cot-direct-image-to-X}(a),
the direct image functor~$\pi_*$ takes semiacyclic complexes in
$Y\qcoh^\dil$ to Becker-coacyclic complexes in $X\qcoh^\dil$.
 Now the derived functor $\boR\pi_*$ is defined by restricting
the functor of direct image $\pi_*\:\Hot(Y\qcoh)\rarrow\Hot(X\qcoh)$
to the full subcategory of complexes of dilute quasi-coherent sheaves
on~$Y$.
 Alternatively, one can use contraadjusted, cotorsion, or injective
coresolutions instead of the dilute ones.

 When the morphism~$\pi$ is affine, the direct image functor
$\pi_*\:Y\qcoh\rarrow X\qcoh$ takes semiacyclic complexes to
Becker-coacyclic ones, so it need not be derived.
 In this case, the functor~$\pi_*$ induces a well-defined
triangulated functor
\begin{equation} \label{semicoderived-to-X-underived-direct-image}
 \pi_*\:\sD^\si_X(Y\qcoh)\rarrow\sD^\bco(X\qcoh),
\end{equation}
which agrees with~\eqref{semicoderived-to-X-derived-direct-image}.

\subsection{Local definition of semicontraderived category}
\label{local-semicontraderived-subsect}
 The exposition in this section is a dual-analogous (contraherent
cosheaf) version of~\cite[Section~A.3]{Psemten}.

 First of all, let us mention an error in the formulation
of~\cite[Lemma~A.12]{Psemten}: the condition that countable direct
limits must be exact in $\sA$ is missing.
 Otherwise, one would have to assume a dual Mittag-Leffler condition
for the direct system of complexes $(M^\bu_n)_{n\ge0}$, or ensure
exactness of the telescope short sequence in some other way.
 This does not affect the validity of the rest of the exposition
in~\cite[Appendix~A]{Psemten}, but forces us to formulate the next
lemma for modules only.

\begin{lem} \label{colocalization-becker-semicontraacyclic}
 Let $R\rarrow S$ be a morphism of commutative rings. \par
\textup{(a)} Let $Q^\bu$ be a complex of contraadjusted $S$\+modules
and $s\in S$ be an element.
 Assume that $Q^\bu$ is Becker-contraacyclic as a complex of
contraadjusted $R$\+modules.
 Then the complex\/ $\Hom_S(S[s^{-1}],Q^\bu)$ is also
Becker-contraacyclic as a complex of contraadjusted $R$\+modules. \par
\textup{(b)} Let $Q^\bu$ be a complex of cotorsion $S$\+modules
and $s\in S$ be an element.
 Assume that $Q^\bu$ is Becker-contraacyclic as a complex of
cotorsion $R$\+modules.
 Then the complex\/ $\Hom_S(S[s^{-1}],Q^\bu)$ is also
Becker-contraacyclic as a complex of cotorsion $R$\+modules.
\end{lem}

\begin{proof}
 Let us recall, first of all, that a complex in $R\modl^\cta$ is
Becker-contraacyclic in $R\modl^\cta$ if and only if it is
Becker-contraacyclic in $R\modl$, while a complex in $R\modl^\cot$
is Becker-contraacyclic in $R\modl^\cot$ if and only if it is
Becker-contraacyclic (in $R\modl^\cta$ or) in $R\modl$
(see Proposition~\ref{becker-contraacyclicity-for-cta-cot}).
 Furthermore, the complex $\Hom_S(S[s^{-1}],Q^\bu)$ is a complex
of contraadjusted $R$\+modules in the assumptions of part~(a) by
Lemmas~\ref{very-tensor-hom}(b) and~\ref{very-scalars-always}(a),
and this complex is a complex of cotorsion $R$\+modules in
the assumptions of part~(b) by Lemmas~\ref{cotors-hom}(a)
and~\ref{cotors-restrict}(a).
 Therefore, part~(a) implies part~(b).
 One can also prove part~(b) independently in the same way as
the proof of part~(a) below.

 Let $0\rarrow P_1\rarrow P_0\rarrow S[s^{-1}]\rarrow0$ be a projective
resolution of the $S$\+module $S[s^{-1}]$ (cf.\
Section~\ref{very-eklof-trlifaj-subsect}).
 Since $Q^\bu$ is a complex of contraadjusted $S$\+modules, the short
sequence of complexes $0\rarrow\Hom_S(S[s^{-1}],Q^\bu)\rarrow
\Hom_S(P_0,Q^\bu)\rarrow\Hom_S(P_1,Q^\bu)\rarrow0$ is exact.
 Moreover, it is a short exact sequence of complexes in $R\modl^\cta$.
 Now the complexes of contraadjusted $R$\+modules $\Hom_S(P_0,Q^\bu)$
and $\Hom_S(P_1,Q^\bu)$ are Becker-contraacyclic, because
the Becker-contraacyclicity is preserved by infinite products.
 In view of Lemma~\ref{Positselski-trivial-are-Becker-trivial}(a),
it follows that the complex $\Hom_S(S[s^{-1}],Q^\bu)$ is
Becker-contraacyclic in $R\modl^\cta$, too.
\end{proof}

\begin{lem} \label{affine-covered-by-principal-semicontraacyclic}
 Let $\pi\:Y\rarrow X$ be a morphism of affine schemes and
$Y=\bigcup_\alpha V_\alpha$ be a covering of the scheme $Y$
by principal affine open subschemes $V_\alpha\sub Y$.
 Denote by $j_\alpha\:V_\alpha\rarrow Y$ the open embedding morphisms.
 Let\/ $\Q^\bu$ be a complex of (globally) contraherent cosheaves on~$Y$.
 Then the complex $\pi_!\Q^\bu$ of contraherent cosheaves on $X$ is
Becker-contraacyclic if and only if, for every index~$\alpha$,
the complex $\pi_!j_\alpha{}_!j_\alpha^!\Q^\bu$ of contraherent
cosheaves on $X$ is Becker-contraacyclic.
\end{lem}

\begin{proof}
 This is a dual-analogous version of~\cite[Lemma~A.14]{Psemten}.
 If the complex $\pi_*\Q^\bu$ is Becker-contraacyclic in $X\ctrh$,
then Lemma~\ref{colocalization-becker-semicontraacyclic} tells us
that the complex $\pi_!j_\alpha{}_!j_\alpha^!\Q^\bu$ is
Becker-contraacyclic in $X\ctrh$ for every~$\alpha$.

 Conversely, we can assume the set of indices~$\alpha$ to be finite.
 Similarly to the proof of~\cite[Lemma~A.14]{Psemten}, we choose
a linear order on the set of indices, and for any $\alpha_1<\dotsb
<\alpha_k$ denote by $j_{\alpha_1,\dotsc,\alpha_k}\:
\bigcap_{s=1}^k V_{\alpha_s}\rarrow Y$ the natural open embedding.
 Assume that the complex $\pi_!j_\alpha{}_!j_\alpha^!\Q^\bu$ is
Becker-contraacyclic in $X\ctrh$ for every~$\alpha$.
 Then, by Lemma~\ref{colocalization-becker-semicontraacyclic},
the complex $\pi_!j_{\alpha_1,\dotsc,\alpha_k}{}_!
j_{\alpha_1,\dotsc,\alpha_k}^!\Q^\bu$ is Becker-contraacyclic in
$X\ctrh$ for all $\alpha_1<\dotsb<\alpha_k$ and $k>0$.

 Consider the \v Cech resolution~\eqref{cech-contra}
or~\eqref{contraherent-cech}
\begin{equation} \label{cech-resolution-of-complex} \textstyle
 0 \lrarrow\dotsb\lrarrow \bigoplus_{\alpha<\beta}j_{\alpha,\beta}{}_!
 j_{\alpha,\beta}^!\Q^\bu \lrarrow\bigoplus_\alpha
 j_\alpha{}_!j_\alpha^!\Q^\bu\lrarrow \Q^\bu\lrarrow0
\end{equation}
of the complex of contraherent cosheaves $\Q^\bu$ on~$Y$.
 The direct image functor $\pi_!\:Y\ctrh\rarrow X\ctrh$ for
a morphism of affine schemes $\pi\:Y\rarrow X$ is exact, so
it takes an exact complex~\eqref{cech-resolution-of-complex}
in $Y\ctrh$ to an exact complex in $X\ctrh$:
\begin{equation} \label{direct-image-of-cech-resolution-of-complex}
 \textstyle
 0\lrarrow\dotsb\lrarrow
 \bigoplus_{\alpha<\beta}\pi_!j_{\alpha,\beta}{}_!
 j_{\alpha,\beta}^!\Q^\bu \lrarrow\bigoplus_\alpha
 \pi_!j_\alpha{}_!j_\alpha^!\Q^\bu\lrarrow\pi_!\Q^\bu\lrarrow0.
\end{equation}

 Finally, the total complex of
the bicomplex~\eqref{direct-image-of-cech-resolution-of-complex}
is Becker-contraacyclic in $X\ctrh$ by
Lemma~\ref{Positselski-trivial-are-Becker-trivial}(a).
 We have seen that all the terms of this complex of complexes,
except perhaps the rightmost one, are Becker-contraacyclic.
 It follows that the rightmost term $\pi_!\Q^\bu$ is
Becker-contraacyclic in $X\ctrh$, too.
\end{proof}

 The next corollary is the contraherent cosheaf version of
Lemma~\ref{affine-locality-semicoacyclicity}.

\begin{cor} \label{affine-covered-by-affines-semicontraacyclic}
 Let $\pi\:Y\rarrow X$ be a morphism of affine schemes and
$Y=\bigcup_\alpha V_\alpha$ be a covering of the scheme $Y$ by
(not necessarily principal) affine open subschemes $V_\alpha\sub Y$.
 Denote by $j_\alpha\:V_\alpha\rarrow Y$ the open embedding morphisms.
 Let\/ $\Q^\bu$ be a complex of (globally) contraherent cosheaves on~$Y$.
 Then the complex $\pi_!\Q^\bu$ of contraherent cosheaves on $X$ is
Becker-contraacyclic if and only if, for every index~$\alpha$,
the complex $\pi_!j_\alpha{}_!j_\alpha^!\Q^\bu$ of contraherent
cosheaves on $X$ is Becker-contraacyclic.
\end{cor}

\begin{proof}
 This is a formal consequence of
Lemma~\ref{affine-covered-by-principal-semicontraacyclic}.
 The argument from~\cite[proof of Lemma~A.15]{Psemten} applies.
\end{proof}

\begin{cor} \label{morphism-into-affine-semicontraacyclic}
 Let $X$ be an affine scheme and $\pi\:Y\rarrow X$ be a morphism of
schemes.
 Let\/ $\bT$ be an open covering of $Y$ and let\/ $\Q^\bu$ be
a complex of\/ $\bT$\+locally contraherent cosheaves on~$Y$.
 Let\/ $\bigcup_\beta V_\beta=Y=\bigcup_\gamma W_\gamma$ be two affine
open coverings of~$Y$, both subordinate to\/~$\bT$.
 Denote by $j_\beta\:V_\beta\rarrow Y$ and
$k_\gamma\:W_\gamma\rarrow Y$ the open embedding morphisms.
 Then the complexes of contraherent cosheaves
$\pi_!j_\beta{}_!j_\beta^!\Q^\bu$ on $X$ are Becker-contraacyclic
for all~$\beta$ if and only if the complexes of contraherent cosheaves
$\pi_!k_\gamma{}_!k_\gamma^!\Q^\bu$ on $X$ are Becker-contraacyclic
for all~$\gamma$.
\end{cor}

\begin{proof}
 This is a formal consequence of
Corollary~\ref{affine-covered-by-affines-semicontraacyclic}.
 The argument from~\cite[proof of Lemma~A.16]{Psemten} applies.
\end{proof}

\begin{lem} \label{affine-morphism-semicontraacyclic}
 Let $\pi\:Y\rarrow X$ be an affine morphism of quasi-compact
semi-separated schemes.
 Let\/ $\bW$ and\/ $\bT$ be open coverings of the schemes $X$ and
$Y$ such that the morphism~$\pi$ is $(\bW,\bT)$\+affine.
 Let $X=\bigcup_\alpha U_\alpha$ be an affine open covering of $X$
subordinate to~$\bW$.
 Put $V_\alpha=U_\alpha\times_XY$ (so $Y=\bigcup_\alpha V_\alpha$ is
an affine open covering of~$Y$).
 Denote the natural morphisms by $j_\alpha\:U_\alpha\rarrow X$, \
$j'_\alpha\:V_\alpha\rarrow Y$, and
$\pi_\alpha\:V_\alpha\rarrow U_\alpha$.
 Let\/ $\Q^\bu$ be a complex of\/ $\bT$\+locally contraherent
cosheaves on~$Y$.
 Then the following conditions are equivalent:
\begin{enumerate}
\item the complex $\pi_!\Q^\bu$ is Becker-contraacyclic in
$X\lcth_\bW$;
\item the complexes $\pi_!j'_\alpha{}_!j^{\prime!}_\alpha\Q^\bu$ are
Becker-contraacyclic in $X\lcth_\bW$ for all~$\alpha$;
\item the complexes $\pi_\alpha{}_!j^{\prime!}_\alpha\Q^\bu$ are
Becker-contraacyclic in $U_\alpha\ctrh$ for all~$\alpha$.
\end{enumerate}
\end{lem}

\begin{proof}
 This is the contraherent cosheaf version of~\cite[Lemma~A.17]{Psemten}.

 (1)\,$\Longleftrightarrow$\,(3)
 In view of the natural isomorphism $j_\alpha^!\pi_!\Q^\bu\simeq
\pi_\alpha{}_!j^{\prime!}_\alpha\Q^\bu$ of complexes of contraherent
cosheaves on $U_\alpha$, the assertion follows from the locality of
Becker-contraacyclicity of complexes of locally contraherent
cosheaves on quasi-compact semi-separated schemes
(Theorem~\ref{Becker-contraacyclicity-local-on-qcomp-qsep}(a)).

 (2)\,$\Longleftrightarrow$\,(3)
 We have $\pi_!j'_\alpha{}_!j^{\prime!}_\alpha\Q^\bu\simeq
j_\alpha{}_!\pi_\alpha{}_!j^{\prime!}_\alpha\Q^\bu$, since
$\pi j'_\alpha=j_\alpha\pi_\alpha$.
 It remains to observe that, for any quasi-compact semi-separated
scheme $X$ with an open covering $\bW$ and any affine open subscheme
$U\sub X$ subordinate to $\bW$, with the open embedding morphism
$j\:U\rarrow X$, a complex of contraherent cosheaves $\P^\bu$ on $U$
is Becker-contraacyclic if and only if the complex of $\bW$\+locally
contraherent cosheaves $j_!\P^\bu$ on $X$ is Becker-contraacyclic.
 This follows from
Lemma~\ref{inverse-images-preserve-contraacyclicity}(a) and
Corollary~\ref{affine-morphism-direct-image-Becker-contra}(a).
\end{proof}

\begin{cor} \label{morphism-into-qcom-ssep-semicontraacyclic}
 Let $X$ be a quasi-compact semi-separated scheme and $\pi\:Y\rarrow X$
be a morphism of schemes.
 Let\/ $\bW$ be an open covering of $X$ and\/ $\bT$ be an open covering
of~$Y$.
 Let $X=\bigcup_\alpha U_\alpha$ be an affine open covering of $X$
subordinate to\/ $\bW$, and let\/ $\bigcup_\beta V_\beta=Y=
\bigcup_\gamma v_\gamma$ be two open coverings of $Y$ subordinate
to\/~$\bT$.
 Denote by\/ $\bT|_{V_\beta}$ and\/ $\bT|_{v_\gamma}$ the restrictions
of the open covering\/ $\bT$ to the open subschemes $V_\beta$ and
$v_\gamma\sub Y$.
 Assume that the compositions $V_\beta\rarrow Y\rarrow X$ and
$v_\gamma\rarrow Y\rarrow X$ are affine morphisms; moreover, assume
that the morphisms $V_\beta\rarrow X$ are
$(\bW,\bT|_{V_\beta})$\+affine, while the morphisms
$v_\gamma\rarrow X$ are $(\bW,\bT|_{v_\gamma})$\+affine.

 Denote by $j_\alpha\:U_\alpha\rarrow X$, \ $k_\beta\:V_\beta\rarrow
Y$, and $l_\gamma\:v_\gamma\rarrow Y$ the open embedding morphisms.
 Put $S_{\alpha,\beta}=U_\alpha\times_XV_\beta$ and $s_{\alpha,\gamma}
=U_\alpha\times_Xv_\gamma$, and denote by $g_{\alpha,\beta}\:
S_{\alpha,\beta}\rarrow Y$ and $h_{\alpha,\gamma}\:s_{\alpha,\gamma}
\rarrow Y$ the natural open embeddings.
 Denote by $\pi'_{\alpha,\beta}\:S_{\alpha,\beta}\rarrow U_\alpha$
and $\pi''_{\alpha,\gamma}\:s_{\alpha,\gamma}\rarrow U_\alpha$
the natural morphisms of affine schemes.

 Let\/ $\Q^\bu$ be a complex of\/ $\bT$\+locally contraherent cosheaves
on~$Y$.
 Then the following conditions are equivalent:
\begin{enumerate}
\item the complexes $\pi_!k_\beta{}_!k_\beta^!\Q^\bu$ are
Becker-contraacyclic in $X\lcth_\bW$ for all~$\beta$;
\item the complexes $\pi_!g_{\alpha,\beta}{}_!g_{\alpha,\beta}^!\Q^\bu$
are Becker-contraacyclic in $X\lcth_\bW$ for all~$\alpha$ and~$\beta$;
\item the complexes $\pi'_{\alpha,\beta}{}_!g_{\alpha,\beta}^!\Q^\bu$
are Becker-contraacyclic in $U_\alpha\ctrh$
for all~$\alpha$ and~$\beta$;
\item the complexes $\pi''_{\alpha,\gamma}{}_!h_{\alpha,\gamma}^!\Q^\bu$
are Becker-contraacyclic in $U_\alpha\ctrh$
for all~$\alpha$ and~$\gamma$;
\item the complexes
$\pi_!h_{\alpha,\gamma}{}_!h_{\alpha,\gamma}^!\Q^\bu$ are
Becker-contraacyclic in $X\lcth_\bW$ for all~$\alpha$ and~$\gamma$;
\item the complexes $\pi_!l_\gamma{}_!l_\gamma^!\Q^\bu$ are
Becker-contraacyclic in $X\lcth_\bW$ for all~$\gamma$.
\end{enumerate}
\end{cor}

\begin{proof}
 This is a formal consequence of
Corollary~\ref{morphism-into-affine-semicontraacyclic} and
Lemma~\ref{affine-morphism-semicontraacyclic}.
 The argument from~\cite[proof of Proposition~A.18]{Psemten} applies.

 Specifically, the equivalences (1)~$\Longleftrightarrow$ (2)
$\Longleftrightarrow$~(3) hold by
Lemma~\ref{affine-morphism-semicontraacyclic} applied to
the affine morphism of schemes $\pi k_\beta\:V_\beta\rarrow X$
and the complex of $\bT|_{V_\beta}$\+locally contraherent cosheaves
$k_\beta^!\Q^\bu$ on~$V_\beta$.
 Similarly, the equivalences (4)~$\Longleftrightarrow$ (5)
$\Longleftrightarrow$~(6) hold by
Lemma~\ref{affine-morphism-semicontraacyclic} applied to
the affine morphism of schemes $\pi l_\gamma\:v_\gamma\rarrow X$
and the complex of $\bT|_{v_\gamma}$\+locally contraherent cosheaves
$l_\gamma^!\Q^\bu$ on~$v_\gamma$.
 Finally, the equivalence (3)~$\Longleftrightarrow$~(4) is provided
by Corollary~\ref{morphism-into-affine-semicontraacyclic}.
 Notice that $\bigcup_\beta S_{\alpha,\beta}=U_\alpha\times_X\nobreak Y
=\bigcup_\gamma s_{\alpha,\gamma}$ are two affine open coverings
of the scheme $U_\alpha\times_XY$, both subordinate to the open covering
$\bT|_{U_\alpha\times_XY}$ of $U_\alpha\times_XY$.
\end{proof}

 Now we can formulate the local definition of the semi(contra)derived
category of $\bT$\+locally contraherent cosheaves on~$Y$.
 Let $\pi\:Y\rarrow X$ be a morphism of schemes; assume that the scheme
$X$ is quasi-compact and semi-separated.
 Let $\bW$ be an open covering of $X$ and $\bT$ be an open covering
of~$Y$; and let $\Q^\bu$ be a complex of $\bT$\+locally contraherent
cosheaves on~$Y$.

 Choose a covering of the scheme $Y$ by open subschemes $V_\beta\sub Y$
such that the compositions $V_\beta\rarrow Y\rarrow X$ are affine,
and moreover, $(\bW,\bT|_{V_\beta})$\+affine morphisms of schemes.
 For example, any affine open covering of $Y$ subordinate to $\bT$
satisfies this condition.
 Denote by $k_\beta\:V_\beta\rarrow Y$ the open embedding morphisms.

 We will say that the complex $\Q^\bu$ is
\emph{semi}(\emph{contra})\emph{acyclic} (over~$X$) if, for every
index~$\beta$, the complex $\pi_!k_\beta{}_!k_\beta^!\Q^\bu$ of
$\bW$\+locally contraherent cosheaves on $X$ is Becker-contraacyclic
in $X\lcth_\bW$.
 According to
Corollary~\ref{morphism-into-qcom-ssep-semicontraacyclic}\,%
(1)\,$\Leftrightarrow$\,(6), this property of a complex of
$\bT$\+locally contraherent cosheaves $\Q^\bu$ on $Y$ does not
depend on the choice of an open covering $Y=\bigcup_\beta V_\beta$.
 It follows that the semicontraacyclicity property is not affected
by refinements of the open covering $\bT$ of the scheme~$Y$.
 As the property of a complex in $X\lcth_\bW$ to be
Becker-contraacyclic is not changed by refinements of the open
covering $\bW$ (see Corollary~\ref{ctrh-lcth-cor}(a) or
Theorem~\ref{Becker-contraacyclicity-local-on-qcomp-qsep}(a)),
the semicontraacyclicity property of complexes on $Y$ is also
unaffected by refinements of the open covering $\bW$ of the scheme~$X$.

 Notice also that a complex in $Y\lcth_\bT^\lct$ is semicontraacyclic
``in the world of locally cotorsion locally contraherent cosheaves''
(in the obvious sense) if and only if it is semicontraacyclic as
a complex in $Y\lcth_\bT$.
 This follows from Corollary~\ref{lct-Becker-contraacyclic-iff-as-lcta}.

\begin{rem} \label{semicontraacyclicity-remark}
 Similarly to~\cite[Remark~A.19]{Psemten} (cf.\ the discussion in
Section~\ref{antilocal-semicoderived-subsect}), we should \emph{warn}
the reader that our terminology is misleading.
 The semi(contra)acyclicity of a complex in $Y\lcth_\bT$ is a stronger
property than the acyclicity.
 By design, the semicontraacyclicity is an intermediate property
between the acyclicity and the Becker-contraacyclicity.
{\hbadness=1075\par}

 Indeed, assuming the scheme $Y$ (as well as~$X$) to be quasi-compact
and semi-separated, any Becker-contraacyclic complex in $Y\lcth_\bT$
is semicontraacyclic by
Lemma~\ref{inverse-images-preserve-contraacyclicity}(a) and
Corollary~\ref{affine-morphism-direct-image-Becker-contra}(a).
 On the other hand, any semicontraacyclic complex $\Q^\bu$ in
$Y\lcth_\bT$ is acyclic.
 By Lemma~\ref{acyclicity-in-lcth-criterion}(a), in order to establish
this fact, it suffices to check that the complex of abelian groups
$\Q^\bu[V]$ is acyclic for every affine open subscheme $V\sub Y$
subordinate to~$\bT$.
 Moreover, by the same lemma, without loss of generality we can
assume existence of an affine open subscheme $U\sub X$ such that
$\pi(V)\sub U$.
 Denoting by $\pi'\:V\rarrow U$ the natural morphism of affine schemes
and by $g\:V\rarrow Y$ the open embedding morphism, we see from
Corollary~\ref{morphism-into-qcom-ssep-semicontraacyclic}\,%
(1)\,$\Leftrightarrow$\,(3) that the complex $\pi'_!g^!\Q^\bu$ is
Becker-contraacyclic in $U\ctrh$.
 By Corollary~\ref{Becker-contraacyclic-are-acyclic}(a)
or Proposition~\ref{becker-contraacyclicity-for-cta-cot}(a),
it follows that the complex $\pi'_!g^!\Q^\bu$ is acyclic in
$U\ctrh$.
 In other words, this means that the complex of abelian groups
$\Q^\bu[V]\simeq(\pi'_!g^!\Q^\bu)[U]$ is acyclic, as desired.

 Similarly one proves that any semicontraacyclic complex in
$Y\lcth_\bT^\lct$ is acyclic in $Y\lcth_\bT^\lct$.
 Assuming $Y$ to be quasi-compact and semi-separated, any
Becker-contraacyclic complex in $Y\lcth_\bT^\lct$ is semicontraacyclic.
\end{rem}

 We will denote the full subcategory of semiacyclic complexes of
$\bT$\+locally contraherent cosheaves on~$Y$ (relative to~$X$) by
$\Acycl^\si_X(Y\lcth_\bT)$.
 The full subcategory of semiacyclic complexes of locally cotorsion
$\bT$\+locally contraherent cosheaves on~$Y$ (relative to~$X$) is
denoted by $\Acycl^\si_X(Y\lcth_\bT^\lct)$.
 So we have
$$
 \Acycl^\bctr(Y\lcth_\bT)\sub\Acycl^\si_X(Y\lcth_\bT)\sub
 \Acycl(Y\lcth_\bT)
$$
and similarly
$$
 \Acycl^\bctr(Y\lcth_\bT^\lct)\sub\Acycl^\si_X(Y\lcth_\bT^\lct)
 \sub\Acycl(Y\lcth_\bT^\lct).
$$
 We also put $\Acycl^\si_X(Y\ctrh)=\Acycl^\si_X(Y\lcth_{\{Y\}})$
and $\Acycl^\si_X(Y\ctrh^\lct)=\Acycl^\si_X(Y\lcth_{\{Y\}}^\lct)$.
{\hbadness=1850\par}

 The full subcategory of semiacyclic complexes
$\Acycl^\si_X(Y\lcth_\bT)$ is a thick subcategory closed under
infinite products in the homotopy category $\Hot(Y\lcth_\bT)$.
 Viewed as a full subcategory in the exact category of complexes
$\Com(Y\lcth_\bT)$, the full subcategory of semiacyclic complexes
is closed under the kernels of admissible epimorphisms, the cokernels
of admissible monomorphisms, and infinite products.
 (Recall that all absolutely acyclic complexes in $Y\lcth_\bT$ are
Becker-contraacyclic by
Lemma~\ref{Positselski-trivial-are-Becker-trivial}(a),
hence also semiacyclic.
 Furthermore, the class of Becker-contraacyclic complexes in
$X\lcth_\bW$ is closed under infinite products, while the functors of
inverse image of locally contraherent cosheaves and their direct image
with respect to affine morphisms of schemes preserve infinite products.)

 The triangulated Verdier quotient category
$$
 \sD^\si_X(Y\lcth_\bT)=\Hot(Y\lcth_\bT)/\Acycl^\si_X(Y\lcth_\bT)
$$
is called the \emph{semi}(\emph{contra})\emph{derived category} of
$\bT$\+locally contraherent cosheaves on~$Y$ (relative to~$X$).
 Similarly, the semi(contra)derived category of locally cotorsion
$\bT$\+locally contraherent cosheaves on~$Y$ (relative to~$X$) is
defined as the triangulated Verdier quotient category
$$
 \sD^\si_X(Y\lcth_\bT^\lct)=
 \Hot(Y\lcth_\bT^\lct)/\Acycl^\si_X(Y\lcth_\bT^\lct).
$$
 We also put $\sD^\si_X(Y\ctrh)=\sD^\si_X(Y\lcth_{\{Y\}})$
and $\sD^\si_X(Y\ctrh^\lct)=\sD^\si_X(Y\lcth_{\{Y\}}^\lct)$,
and call these triangulated categories
the \emph{semi}(\emph{contra})\emph{derived category of contraherent
cosheaves} on $Y$ and the \emph{semi}(\emph{contra})\emph{derived category of locally cotorsion contraherent cosheaves} on~$Y$
(relative to~$X$).

\begin{cor} \label{semicontraderived-independent-of-covering}
 Let $\pi\:Y\rarrow X$ be a morphism of quasi-compact semi-separated
schemes, $\bW$ be an open covering of $X$, and\/ $\bT$ be an open
covering of~$Y$.
 Then \par
\textup{(a)} the triangulated functor\/ $\sD^\si_X(Y\ctrh)\rarrow
\sD^\si_X(Y\lcth_\bT)$ induced by the embedding of exact categories
$Y\ctrh\rarrow Y\lcth_\bT$ is an equivalence of triangulated categories;
\par
\textup{(a)} the triangulated functor\/ $\sD^\si_X(Y\ctrh^\lct)\rarrow
\sD^\si_X(Y\lcth_\bT^\lct)$ induced by the embedding of exact categories
$Y\ctrh^\lct\rarrow Y\lcth_\bT^\lct$ is an equivalence of triangulated
categories.
\end{cor}

\begin{proof}
 It was explained in the discussion above that a complex of (locally
contraadjusted or locally cotorsion) contraherent cosheaves on $Y$
is semicontraacyclic if and only if it is semicontraacyclic as
a complex of (locally contraadjusted or locally cotorsion, respectively)
$\bT$\+locally contraherent cosheaves on~$Y$.
 Hence the induced triangulated functors in parts~(a) and~(b) are
well-defined.
 It remains to invoke (the proofs of) Corollaries~\ref{ctrh-lcth-cor}(a)
and~\ref{lct-ctrh-lcth-cor}(a) to the effect that for every complex
of $\bT$\+locally contraherent cosheaves $\gM^\bu$ on $Y$ there exists
a complex of contraherent cosheaves $\P^\bu$ together with a morphism
of complexes $\P^\bu\rarrow\gM^\bu$ with a cone absolutely acyclic in
$Y\lcth_\bT$, and similarly for locally cotorsion contraherent
cosheaves.
 Now all absolutely acyclic complexes are Becker-contraacyclic, all
Becker-contraacyclic complexes are semicontraacyclic, and
Lemma~\ref{pkoszul-lemma16}(a) applies.
\end{proof}

\begin{thm} \label{semicontraderived-lcta-lct-equivalent}
 Let $\pi\:Y\rarrow X$ be a morphism of quasi-compact semi-separated
schemes, $\bW$ be an open covering of $X$, and\/ $\bT$ be an open
covering of~$Y$.
 Then the triangulated functor\/ $\sD^\si_X(Y\lcth_\bT^\lct)\rarrow
\sD^\si_X(Y\lcth_\bT)$ induced by the embedding of exact categories
$Y\lcth_\bT^\lct\rarrow Y\lcth_\bT$ is an equivalence of triangulated
categories.
\end{thm}

\begin{proof}
 The argument is similar to the proof of
Theorem~\ref{becker-contraderived-lcta-lct-equivalent}.
 Introduce notation for the intersections of triangulated subcategories
\begin{align*}
 \Acycl^\si_X(Y\ctrh_\prj) &=
 \Hot(Y\ctrh_\prj)\cap\Acycl^\si_X(Y\lcth_\bT), \\
 \Acycl^\si_X(Y\ctrh^\lct_\prj) &=
 \Hot(Y\ctrh^\lct_\prj)\cap\Acycl^\si_X(Y\lcth_\bT^\lct), \\
 \Acycl^\si_X(Y\ctrh_\alf) &=
 \Hot(Y\ctrh_\alf)\cap\Acycl^\si_X(Y\lcth_\bT)
\end{align*}
and for the triangulated quotient categories
\begin{align*}
 \sD^\si_X(Y\ctrh_\prj) &=
 \Hot(Y\ctrh_\prj)/\Acycl^\si_X(Y\ctrh_\prj), \\
 \sD^\si_X(Y\ctrh^\lct_\prj) &=
 \Hot(Y\ctrh^\lct_\prj)/\Acycl^\si_X(Y\ctrh^\lct_\prj), \\
 \sD^\si_X(Y\ctrh_\alf) &=
 \Hot(Y\ctrh_\alf)/\Acycl^\si_X(Y\ctrh_\alf).
\end{align*}

 Let us show that the triangulated functors
$$
 \sD^\si_X(Y\ctrh_\prj)\rarrow\sD^\si_X(Y\lcth_\bT)
 \quad\text{and}\quad
 \sD^\si_X(Y\ctrh^\lct_\prj)\rarrow\sD^\si_X(Y\lcth_\bT^\lct)
$$
induced by the inclusions of additive/exact categories
$Y\ctrh_\prj\rarrow Y\lcth_\bT$ and
$Y\ctrh^\lct_\prj\rarrow Y\lcth_\bT^\lct$ are equivalences
of triangulated categories.
 Indeed, by (the proof of)
Corollary~\ref{becker-contraderived-of-lcta-lct-well-behaved}(a),
for any complex of $\bT$\+locally contraherent cosheaves $\gM^\bu$
on $Y$ there exists a complex of projective contraherent cosheaves
$\P^\bu$ together with a morphism of complexes $\P^\bu\rarrow\gM^\bu$
with a cone Becker-contraacyclic in $Y\lcth_\bT$.
 Similarly, by (the proof of)
Corollary~\ref{becker-contraderived-of-lcta-lct-well-behaved}(b),
for any complex of locally cotorsion $\bT$\+locally contraherent
cosheaves $\gM^\bu$ on $Y$ there exists a complex of projective
locally cotorsion contraherent cosheaves $\P^\bu$ together with
a morphism of complexes $\P^\bu\rarrow\gM^\bu$ with a cone
Becker-contraacyclic in $Y\lcth_\bT^\lct$.
 All Becker-contraacyclic complexes are semicontraacyclic, so
Lemma~\ref{pkoszul-lemma16}(a) is applicable.

 It remains to show that the functors
$$
 \sD^\si_X(Y\ctrh_\prj)\rarrow\sD^\si_X(Y\ctrh_\alf)
 \quad\text{and}\quad
 \sD^\si_X(Y\ctrh^\lct_\prj)\rarrow\sD^\si_X(Y\ctrh_\alf)
$$
induced by the inclusions of additive/exact categories
$Y\ctrh_\prj\rarrow Y\ctrh_\alf$ and
$Y\ctrh^\lct_\prj\rarrow Y\ctrh_\alf$ are triangulated equivalences.
 Indeed, by (the proof of)
Corollary~\ref{becker-contraderived-of-al-lct-alf-well-behaved}(c),
for any complex of antilocally flat contraherent cosheaves $\gF^\bu$
on $Y$ there exists a complex of projective contraherent cosheaves
$\P^\bu$ together with a morphism of complexes $\P^\bu\rarrow\gF^\bu$
with a cone Becker-contraacyclic in $Y\ctrh_\alf$, or equivalently,
in $Y\lcth_\bT$.
 On the other hand, by (the proof of)
Corollary~\ref{becker-coderived-of-alf-well-behaved}, for any complex
of antilocally flat contraherent cosheaves $\gF^\bu$ there exists
a complex of projective locally cotorsion contraherent cosheaves
$\Q^\bu$ on $Y$ together with a morphism of complexes $\gF^\bu
\rarrow\Q^\bu$ with a cone Becker-coacyclic in $Y\ctrh_\alf$.
 By Corollary~\ref{acycl=bctracycl=bcoacycl-in-alf}, the classes of
Becker-coacyclic and Becker-contraacyclic complexes in
the exact category $Y\ctrh_\alf$ coincide.
 Once again, the argument finishes by applying
Lemma~\ref{pkoszul-lemma16}(a\+b).
\end{proof}

\subsection{Antilocal definition of semicontraderived category}
\label{antilocal-semicontraderived-subsect}
 We start with a contraherent cosheaf version of
Lemma~\ref{semicoacyclic-direct-image}.

\begin{lem} \label{semicontraacyclic-direct-image}
 Let $X$ be a quasi-compact semi-separated scheme, $\pi\:Y\rarrow X$
be a morphism of schemes, and $g\:Z\rarrow Y$ be an affine morphism
of schemes.
 Let\/ $\bT$ be an open covering of $Y$ and\/ $\bS$ be an open covering
of $Z$ such that $g$~is a $(\bT,\bS)$\+affine morphism.
 In this context: \par
\textup{(a)} The direct image functor $g_!\:Z\lcth_\bS\rarrow
Y\lcth_\bT$ takes semiacyclic complexes of\/ $\bS$\+locally contraherent
cosheaves on~$Z$ (relative to the morphism $\pi g\:Z\rarrow X$) to
semiacyclic complexes of\/ $\bT$\+locally contraherent cosheaves on~$Y$
(relative to the morphism $\pi\:Y\rarrow X$).
 In other words, $g_!\Acycl^\si_X(Z\lcth_\bS)\sub
\Acycl^\si_X(Y\lcth_\bT)$.
 Conversely, given a complex\/ $\R^\bu$ in $Z\lcth_\bS$, if the complex
$g_!\R^\bu$ on $Y$ is semiacyclic over $X$, then the complex\/ $\R^\bu$
on $Z$ is also semiacyclic over~$X$. \par
\textup{(b)} The direct image functor $g_!\:Z\lcth_\bS^\lct\rarrow
Y\lcth_\bT^\lct$ takes semiacyclic complexes of locally cotorsion\/
$\bS$\+locally contraherent cosheaves on~$Z$ (relative to the morphism
$\pi g\:Z\rarrow X$) to semiacyclic complexes of locally cotorsion\/
$\bT$\+locally contraherent cosheaves on~$Y$ (relative to the morphism
$\pi\:Y\rarrow X$).
 In other words, $g_!\Acycl^\si_X(Z\lcth_\bS^\lct)\sub
\Acycl^\si_X(Y\lcth_\bT^\lct)$.
 Conversely, given a complex\/ $\R^\bu$ in $Z\lcth_\bS^\lct$, if
the complex $g_!\R^\bu$ on $Y$ is semiacyclic over $X$, then
the complex\/ $\R^\bu$ on $Z$ is also semiacyclic over~$X$.
\end{lem}

\begin{proof}
 Since the notions of semi(contra)acyclicity in the realms of locally
contraadjusted contraherent cosheaves and locally cotorsion contraherent
cosheaves agree (as per the discussion in
Section~\ref{local-semicontraderived-subsect}), part~(b) follows from
part~(a).
 Now let $\bW$ be an open covering of $X$, and let $Y=\bigcup_\beta
V_\beta$ be an open covering of $Y$ with open embedding morphisms
$k_\beta\:V_\beta\rarrow Y$ such that the compositions $\pi k_\beta\:
V_\beta\rarrow X$ are $(\bW,\bT|_{V_\beta})$\+affine morphisms.
 Put $U_\beta=V_\beta\times_YZ$, and denote the natural morphisms by
$k'_\beta\:U_\beta\rarrow Z$ and $g'_\beta\:U_\beta\rarrow V_\beta$.
 Then the compositions $\pi gk'_\beta=\pi k_\beta g'_\beta\:
U_\beta\rarrow X$ are $(\bW,\bS|_{U_\beta})$\+affine morphisms,
because the morphisms~$g'_\beta$ are
$(\bT|_{V_\beta},\bS|_{U_\beta})$\+affine.
 It remains to point out the natural isomorphisms
$\pi_!k_\beta{}_!k_\beta^!g_!\R^\bu\simeq
\pi_!k_\beta{}_!g'_\beta{}_!k^{'!}_\beta\R^\bu
\simeq\pi_!g_!k'_\beta{}_!k^{'!}_\beta\R^\bu$ of complexes
of $\bW$\+locally contraherent cosheaves on $X$, which hold for
any complex of $\bS$\+locally contraherent cosheaves $\R^\bu$ on~$Z$.
\end{proof}

\begin{lem} \label{semicontraacyclic-of-projectives-antilocal}
 Let $\pi\:Y\rarrow X$ be a morphism of quasi-compact semi-separated
schemes, and let $Y=\bigcup_\alpha V_\alpha$ be a finite affine open
covering of~$Y$. \par
\textup{(a)} Let\/ $\P^\bu\in\Com(Y\ctrh_\prj)$ be a complex of
projective contraherent cosheaves on~$Y$.
 Then the complex\/ $\P^\bu$ is semiacyclic over $X$ if and only if it
is a direct summand of a finitely iterated extension of the direct
images of complexes of projective contraherent cosheaves on $V_\alpha$
that are semiacyclic over~$X$. \par
\textup{(b)} Let\/ $\Q^\bu\in\Com(Y\ctrh^\lct_\prj)$ be a complex of
projective locally cotorsion contraherent cosheaves on~$Y$.
 Then the complex\/ $\Q^\bu$ is semiacyclic over $X$ if and only if it
is a direct summand of a finitely iterated extension of the direct
images of complexes of projective locally cotorsion contraherent
cosheaves on $V_\alpha$ that are semiacyclic over~$X$. \par
\textup{(c)} Let\/ $\gF^\bu\in\Com(Y\ctrh_\alf)$ be a complexes of
antilocally flat contraherent cosheaves on~$Y$.
 Then the complex\/ $\gF^\bu$ is semiacyclic over $X$ if and only if it
is a direct summand of a finitely iterated extension of the direct
images of complexes of antilocally flat contraherent cosheaves on
$V_\alpha$ that are semiacyclic over~$X$.
\end{lem}

\begin{proof}
 This is a dual-analogous version of
Lemma~\ref{semicoacyclic-of-injectives-antilocal}.
 The arguments are based on
Theorem~\ref{loc-contraherent-gluing-theorem}
with Remark~\ref{loc-contraherent-gluing-nonuniversal-remark}.

 Part~(a): for any affine open subscheme $U\sub Y$, put
$\sE^U=\Acycl^\si_X(U\ctrh)\sub\sK^{\O(U)}=\Com(\O(U)\modl^\cta)$.
 Put $\sC^U=\Acycl^\bctr(U\ctrh)\sub\sE^U$.
 Then the system of classes $\sE^U$ is colocal essentially by
the definition, and Lemma~\ref{semicontraacyclic-direct-image}(a)
implies that it satisfies the direct image condition; so it is very
colocal (with respect to identity open embeddings of affine open 
subschemes of~$Y$).
 The system of classes $\sC^U$ is very colocal by
Theorem~\ref{Becker-contraacyclicity-local-on-qcomp-qsep}(a)
and Corollary~\ref{affine-morphism-direct-image-Becker-contra}(a).

 By Theorem~\ref{becker-contraderived-cotorsion-pair}(a),
the pair of classes ($\Com(U\ctrh_\prj)$, $\Acycl^\bctr(U\ctrh)$)
is a hereditary complete cotorsion pair in the exact category
$\Com(U\ctrh)$.
 By Lemmas~\ref{restricting-hereditary-cotorsion}
and~\ref{restricting-cotorsion-pairs-lemma}(b), this cotorsion pair
restricts to the exact subcategory $\sE^U=\Acycl^\si_X(U\ctrh)
\sub\Com(U\ctrh)$, providing a hereditary complete cotorsion pair
$(\sF(U),\sC^U)$ in $\sE^U$ with the left class
$\sF(U)=\Com(U\ctrh_\prj)\cap\Acycl^\si_X(U\ctrh)$.
 Similarly we obtain a hereditary complete cotorsion pair $(\sF,\sC)$
in the exact category $\Acycl^\si_X(Y\lcth_\bT)$ with the left class
$\sF=\Com(Y\ctrh_\prj)\cap\Acycl^\si_X(Y\lcth_\bT)$ and the right class
$\sC=\Acycl^\bctr(Y\lcth_\bT)$.

 Now Remark~\ref{loc-contraherent-gluing-nonuniversal-remark} is
applicable.
 It produces a hereditary complete cotorsion pair $(\sF(Y),\sC^Y_\bT)$
in the exact category $\sE^Y_\bT=\Acycl^\si_X(Y\lcth_\bT)$ of
locally\+$\sE$ complexes of $\bT$\+locally contraherent cosheaves
on~$Y$.
 The right class $\sC^Y_\bT$ is the class of locally\+$\sC$
complexes of $\bT$\+locally contraherent cosheaves on $Y$, so
$\sC^Y_\bT=\Acycl^\bctr(Y\lcth_\bT)=\sC$ by
Theorem~\ref{Becker-contraacyclicity-local-on-qcomp-qsep}(a).
 Hence the left class is $\sF(Y)=\sF=
\Com(Y\ctrh_\prj)\cap\Acycl^\si_X(Y\lcth_\bT)$.

 On the other hand,
Theorem~\ref{loc-contraherent-gluing-theorem}
with Remark~\ref{loc-contraherent-gluing-nonuniversal-remark}
provide a description of the class $\sF(Y)$ as the class of all
direct summands of finitely iterated extensions of direct images of
complexes from $\sF(V_\alpha)$.
 Comparing these two descriptions of the class $\sF(Y)$, we arrive
to the assertion of part~(a).

 The proof of part~(b) is similar and based on
Lemma~\ref{semicontraacyclic-direct-image}(b),
Theorem~\ref{Becker-contraacyclicity-local-on-qcomp-qsep}(b),
Corollary~\ref{affine-morphism-direct-image-Becker-contra}(b),
and Theorem~\ref{becker-contraderived-cotorsion-pair}(b).
 One constructs a hereditary complete cotorsion pair $(\sF,\sC)$
in the exact category $\sE^Y_\bT=\Acycl^\si_X(Y\lcth_\bT^\lct)$ with
the left class $\sF=\Com(Y\ctrh^\lct_\prj)\cap
\Acycl^\si_X(Y\lcth_\bT^\lct)$ and the right class
$\sC=\Acycl^\bctr(Y\lcth_\bT^\lct)$.

 The proof of part~(c) is also similar and based on
Lemma~\ref{semicontraacyclic-direct-image}(a),
Theorem~\ref{Becker-contraacyclicity-local-on-qcomp-qsep}(b),
Corollary~\ref{affine-morphism-direct-image-Becker-contra}(b),
and Lemma~\ref{all-alf-becker-contraacyclic-lct-cotorsion-pair}.
 One constructs a hereditary complete cotorsion pair $(\sF,\sC)$
in the exact category $\sE^Y_\bT=\Acycl^\si_X(Y\lcth_\bT)$ with
the left class $\sF=\Com(Y\ctrh_\alf)\cap\Acycl^\si_X(Y\lcth_\bT)$
and the right class $\sC=\Acycl^\bctr(Y\lcth_\bT^\lct)$.
\end{proof}

\begin{cor} \label{antilocal-description-of-semicontraacyclic}
 Let $\pi\:Y\rarrow X$ be a morphism of quasi-compact semi-separated
schemes.
 Let\/ $\bT$ be an open covering of $Y$, and let $Y=\bigcup_\alpha
V_\alpha$ be an affine open covering of $Y$ subordinate to\/~$\bT$.
 Denote by $j_\alpha\:V_\alpha\rarrow Y$ the open embedding morphisms.
 Then \par
\textup{(a)} the full subcategory of semiacyclic complexes\/
$\Acycl^\si_X(Y\lcth_\bT)$ coincides with the minimal thick subcategory
of\/ $\Hot(Y\lcth_\bT)$ containing the subcategories
$j_\alpha{}_!\Acycl^\si_X(V_\alpha\ctrh)$ for all the indices~$\alpha$
and the subcategory of absolutely acyclic complexes\/
$\Acycl^\abs(Y\lcth_\bT)$; \hbadness=1175\par
\textup{(b)} the full subcategory of semiacyclic complexes\/
$\Acycl^\si_X(Y\lcth_\bT^\lct)$ coincides with the minimal thick
subcategory of\/ $\Hot(Y\lcth_\bT^\lct)$ containing the subcategories
$j_\alpha{}_!\Acycl^\si_X(V_\alpha\ctrh^\lct)$ for all
the indices~$\alpha$ and the subcategory of absolutely acyclic
complexes\/ $\Acycl^\abs(Y\lcth_\bT^\lct)$.
\end{cor}

\begin{proof}
 This is the contraherent cosheaf version of
Corollary~\ref{antilocal-description-of-semicoacyclic}.
 We will prove part~(a).
 Denote temporarily the minimal thick subcategory of $\Hot(Y\lcth_\bT)$
containing $\Acycl^\abs(Y\lcth_\bT)$ and
$j_\alpha{}_!\Acycl^\si_X(V_\alpha\ctrh)$ for all~$\alpha$ by~$\sT$.
 Then the inclusion $\sT\sub\Acycl^\si_X(Y\lcth_\bT)$ holds by
Lemma~\ref{semicontraacyclic-direct-image}.
 To prove the converse inclusion, consider a complex
$\Q^\bu\in\Acycl^\si_X(Y\lcth_\bT)$.
 By (the proof of)
Corollary~\ref{becker-contraderived-of-lcta-lct-well-behaved}(a),
there exists a complex of projective contraherent cosheaves $\P^\bu$
on $Y$ together with a morphism of complexes of $\bT$\+locally
contraherent cosheaves $\P^\bu\rarrow\Q^\bu$ a with
a Becker-contraacyclic cone.
 As any Becker-contraacyclic complex in $Y\lcth_\bT$ is 
semiacyclic, it follows that the complex $\P^\bu$ is semiacyclic.
 By Lemma~\ref{semicontraacyclic-of-projectives-antilocal}(a),
we have $\P^\bu\in\sT$ (notice that the functors~$j_\alpha{}_!$
take projective contraherent cosheaves to projective contraherent
cosheaves by Corollary~\ref{proj-direct-inverse}(a), and any
extension of complexes of projective contraherent cosheaves is
termwise split, so it corresponds to a distinguished triangle
in $\Hot(Y\lcth_\bT)$).
 It remains to show that $\Acycl^\bctr(Y\lcth_\bT)\sub\sT$.

 Consider a complex $\gB^\bu\in\Acycl^\bctr(Y\lcth_\bT)$.
 By Corollary~\ref{ctrh-lcth-cor}(a) (for $\bst=\abs$), there exists
a complex of antilocal contraherent cosheaves $\gA^\bu$ on $Y$
together with a morphism of complexes $\gA^\bu\rarrow\gB^\bu$ with
an absolutely acyclic cone.
 Since the absolutely acyclic complexes are Becker-contraacyclic,
the complex $\gA^\bu$ is also Becker-contraacyclic in $Y\lcth_\bT$,
or equivalently, in $Y\ctrh_\al$.
 By Corollary~\ref{contraacyclic-complexes-of-lcta-lct-antilocal}(a),
it follows that the complex~$\gA^\bu$ belongs to the minimal thick
subcategory of $\Hot(Y\ctrh_\al)$ containing
$\Acycl^\abs(Y\ctrh_\al)$ and $j_\alpha{}_!\Acycl^\bctr(V_\alpha\ctrh)$
for all the indices~$\alpha$.
 Hence the complex $\gB^\bu$ belongs to the minimal thick subcategory
of $\Hot(Y\lcth_\bT)$ containing $\Acycl^\abs(Y\lcth_\bT)$ and
$j_\alpha{}_!\Acycl^\bctr(V_\alpha\ctrh)$ for all~$\alpha$.
\end{proof}

 By a \emph{semi}(\emph{contra})\emph{acyclic complex of antilocal
contraherent cosheaves} on $Y$ we will mean a complex in $Y\ctrh_\al$
that is semiacyclic as a complex in $Y\lcth_\bT$.
 The similar definition applies to complexes of antilocal locally
cotorsion contraherent cosheaves on~$Y$.
 We will denote the full subcategories of semiacyclic complexes of
antilocal locally contraadjusted and locally cotorsion contraherent
cosheaves on~$Y$ (relative to~$X$) by $\Acycl^\si_X(Y\ctrh_\al)$ and
$\Acycl^\si_X(Y\ctrh_\al^\lct)$.

\begin{cor} \label{semicontraacyclic-al-lcta-lct-direct-image-to-X}
 Let $\pi\:Y\rarrow X$ be a morphism of quasi-compact semi-separated
schemes.
 Then \par
\textup{(a)} the direct image functor $\pi_!\:Y\ctrh_\al\rarrow
X\ctrh_\al$ takes semiacyclic complexes of antilocal contraherent
cosheaves on~$Y$ (relative to~$X$) to Becker-contraacyclic complexes
of antilocal contraherent cosheaves on~$X$; \par
\textup{(b)} the direct image functor $\pi_!\:Y\ctrh_\al^\lct\rarrow
X\ctrh_\al^\lct$ takes semiacyclic complexes of antilocal locally
cotorsion contraherent cosheaves on~$Y$ (relative to~$X$) to
Becker-contraacyclic complexes of antilocal locally cotorsion
contraherent cosheaves on~$X$.
\end{cor}

\begin{proof}
 The direct image functors $\pi_!\:Y\ctrh_\al\rarrow X\ctrh_\al$ and
$\pi_!\:Y\ctrh_\al^\lct\rarrow X\ctrh_\al^\lct$ are well-defined and
exact by Corollary~\ref{clp-direct}(a\+b).
 Let us prove part~(a).
 It was essentially shown in the proof of
Corollary~\ref{antilocal-description-of-semicontraacyclic} that
$\Acycl^\si_X(Y\ctrh_\al)$ is the thick subcategory of
$\Hot(Y\ctrh_\al)$ generated by
$j_\alpha{}_!\Acycl^\si_X(V_\alpha\ctrh)$ and
$\Acycl^\abs(Y\ctrh_\al)$.
 By the definition of semicontraacyclicity for the affine morphism
$\pi j_\alpha\:V_\alpha\rarrow X$, the functor~$\pi_*$ takes
$j_\alpha{}_!\Acycl^\si_X(V_\alpha\ctrh)$ into $\Acycl^\bctr(X\ctrh)$
(hence also into $\Acycl^\bctr(X\ctrh_\al)$).
 Since the functor $\pi_!\:Y\ctrh_\al\rarrow X\ctrh_\al$ is exact,
it takes $\Acycl^\abs(Y\ctrh_\al)$ into $\Acycl^\abs(X\ctrh_\al$).
\end{proof}

 The left derived functor of direct image
\begin{equation} \label{semicontraderived-to-X-derived-direct-image}
 \boL\pi_!\:\sD^\si_X(Y\ctrh)\rarrow\sD^\bctr(X\ctrh)
\end{equation}
is constructed in the following way.
 In view of (the proof of) Corollary~\ref{ctrh-lcth-cor}(a) and
Lemma~\ref{pkoszul-lemma16}(a), the natural functor
$$
 \Hot(Y\ctrh_\al)/\Acycl^\si_X(Y\ctrh_\al)\lrarrow
 \sD^\si_X(Y\ctrh)
$$
is an equivalence of triangulated categories.
 By Corollary~\ref{semicontraacyclic-al-lcta-lct-direct-image-to-X}(a),
the direct image functor~$\pi_!$ takes semiacyclic complexes in
$Y\ctrh_\al$ to Becker-contraacyclic complexes in $X\ctrh_\al$.
 Now the derived functor $\boL\pi_!$ is defined by restricting
the functor of direct image $\pi_!\:\Hot(Y\ctrh)\rarrow\Hot(X\ctrh)$
to the full subcategory of complexes of antilocal contraherent cosheaves
on~$Y$.

 Similarly one defines the left derived functor of direct image
\begin{equation} \label{semicontraderived-lct-to-X-derived-direct-image}
 \boL\pi_!\:\sD^\si_X(Y\ctrh^\lct)\rarrow\sD^\bctr(X\ctrh^\lct).
\end{equation}
 Recall that, according to the discussion in
Section~\ref{local-semicontraderived-subsect}, the notions of
semiacyclicity for complexes in $Y\ctrh^\lct$ and in $Y\ctrh$ agree.
 Hence the functors~\eqref{semicontraderived-to-X-derived-direct-image}
and~\eqref{semicontraderived-lct-to-X-derived-direct-image} also agree
with each other, and we obtain a commutative diagram of triangulated
functors and triangulated equivalences provided by
Theorems~\ref{semicontraderived-lcta-lct-equivalent}
and~\ref{becker-contraderived-lcta-lct-equivalent},
\begin{equation} \label{semicontra-lct-lcta-direct-images-compatible}
\begin{gathered}
 \xymatrix{
  \sD^\si_X(Y\ctrh^\lct) \ar@<2pt>[r] \ar@<-2pt>@{-}[r]
  \ar[d]^{\boL f_!}
  & \sD^\si_X(Y\ctrh) \ar[d]^{\boL f_!} \\
  \sD^\bctr(X\ctrh^\lct) \ar@<2pt>[r] \ar@<-2pt>@{-}[r]
  & \sD^\bctr(X\ctrh)
 }
\end{gathered}
\end{equation}

 When the morphism~$\pi$ is affine, let $\bW$ be an open covering of
the scheme $X$ and $\bT$ be an open covering of the scheme $Y$
such that the morphism~$\pi$ is $(\bW,\bT)$\+affine.
 Then the direct image functors $\pi_!\:Y\lcth_\bT\rarrow X\lcth_\bW$
and $\pi_!\:Y\lcth_\bT^\lct\rarrow X\lcth_\bW^\lct$ take semiacyclic
complexes to Becker-contraacyclic ones, so they need not be derived.
 In this case, the functors~$\pi_!$ induce well-defined
triangulated functors
\begin{equation} \label{semicontraderived-to-X-underived-direct-image}
 \pi_!\:\sD^\si_X(Y\lcth_\bT)\rarrow\sD^\bctr(X\lcth_\bW)
\end{equation}
and
\begin{equation} \label{semicontrader-lct-to-X-underived-direct-image}
 \pi_!\:\sD^\si_X(Y\lcth_\bT^\lct)\rarrow\sD^\bctr(X\lcth_\bW^\lct)
\end{equation}
which agree with~\eqref{semicontraderived-to-X-derived-direct-image}
and~\eqref{semicontraderived-lct-to-X-derived-direct-image}.

\subsection{Direct images of contratensor products under
nonaffine morphisms} \label{direct-contratensor-nonaffine-subsect}
 The aim of this section is to close a gap in the formalism of
Sections~\ref{contratensor-subsect} and~\ref{compatibility-subsect}
concerning applicability of the projection formula for the contratensor product~(\ref{contratensor-projection}\+-%
\ref{flat-contratensor-projection}) to nonaffine morphisms of
nonaffine schemes. {\hbadness=1200\par}

 Let $f\:Y\rarrow X$ be a quasi-compact morphism of quasi-separated
schemes.
 Let $\M$ be a quasi-coherent sheaf on $X$ and $\Q$ be a cosheaf of
$\O_Y$\+modules on~$Y$.
 Then, as explained in Section~\ref{compatibility-subsect}, there is
a natural morphism
\begin{equation} \label{contratensor-projection-comparison}
 \M\ocn_X f_!\Q\lrarrow f_*(f^*\M\ocn_Y\Q)
\end{equation}
of quasi-coherent sheaves on~$X$.

\begin{lem} \label{flat-al-contratensor-projection-lemma}
 Let $f\:Y\rarrow X$ be a morphism of quasi-compact semi-separated
schemes, $\F$ be a flat quasi-coherent sheaf on $X$, and\/ $\Q$ be
an antilocal contraherent cosheaf on~$Y$.
 Then the natural morphism~\eqref{contratensor-projection-comparison}
is an isomorphism
\begin{equation} \label{flat-al-contratensor-projection-eqn}
 \F\ocn_X f_!\Q\simeq f_*(f^*\F\ocn_Y\Q)
\end{equation}
of quasi-coherent sheaves on~$X$.
\end{lem}

\begin{proof}
 Let us show that both the left-hand and the right-hand sides
of~\eqref{flat-al-contratensor-projection-eqn} are exact functors
$Y\ctrh_\al\rarrow X\qcoh$ of an antilocal contraherent cosheaf $\Q$
on $Y$ for any fixed flat quasi-coherent sheaf $\F$ on~$X$.
 Indeed, the functor $f_!\:Y\ctrh_\al\rarrow X\ctrh_\al$ is exact
by Corollary~\ref{clp-direct}(a).
 In order to show that the functor $\F\ocn_X{-}\,\:X\ctrh_\al\rarrow
X\qcoh$ is exact, it suffices to pick an injective quasi-coherent
sheaf $\J$ on $X$ and recall that $\Hom_X(\F\ocn_X\P\;\J)\simeq
\Hom^X(\P,\fHom_X(\F,\J))$ for any locally contraherent cosheaf $\P$
on $X$ by the formula~\eqref{fHom-contratensor-adjunction}.
 According to Section~\ref{fHom-subsection}, \,$\fHom_X(\F,\J)$ is
a locally injective contraherent cosheaf on~$X$.
 Hence the functor $\P\longmapsto\Hom^X(\P,\fHom_X(\F,\J))$ is exact
for $\P\in X\ctrh_\al$ by Corollary~\ref{clp-characterizations}(a).
 This proves that the left-hand side
of~\eqref{flat-al-contratensor-projection-eqn} is exact as a functor
of~$\Q$.

 Concerning the right-hand side, $f^*\F$ is a flat quasi-coherent
sheaf on $Y$, so the argument from the previous paragraph shows that
the functor $\Q\longmapsto f^*\F\ocn_Y\Q$ is exact on the category
of antilocal contraherent cosheaves $\Q\in Y\ctrh_\al$.
 Let us show that the quasi-coherent sheaf $f^*\F\ocn_Y\Q$ on $Y$
is dilute.
 The class of dilute quasi-coherent sheaves on $Y$ is closed under
extensions in $Y\qcoh$, so in view of Corollary~\ref{clp-cor}(c)
it suffices to check that $f^*\F\ocn_Y j_!\R$ is a dilute
quasi-coherent sheaf on $Y$ for every affine open subscheme $V\sub Y$
with the open embedding morphism $j\:V\rarrow Y$ and any contraherent
cosheaf $\R$ on~$V$.
 By~\eqref{contratensor-projection}, we have $f^*\F\ocn_Y j_!\R
\simeq j_*(j^*f^*\F\ocn_V\R)$, and it remains to recall that the direct
image of any quasi-coherent sheaf from an affine open subscheme of $Y$
is dilute (Lemma~\ref{flat-affine-direct-image-dilute}).
 As the functor of direct image of dilute quasi-coherent sheaves with
respect to~$f$ is exact (Corolary~\ref{dilute-direct}), this proves
that the right-hand side
of~\eqref{flat-al-contratensor-projection-eqn} is exact as a functor
of $\Q$, too.

 Using Corollary~\ref{clp-cor}(c) again, it remains to prove
the assertion of the lemma for $\Q=j_!\R$, where $\R$ is a contraherent
cosheaf on~$V$.
 In this case, the desired assertion follows from
formula~\eqref{contratensor-projection} applied to the two affine
morphisms $j\:V\rarrow Y$ and $fj\:V\rarrow X$.
 One needs to use commutativity of the triangle diagram
$$
 \F\ocn_X f_!j_!\R\lrarrow f_*(f^*\F\ocn_Y j_!\R)\lrarrow
 f_*j_*(j^*f^*\F\ocn_V\R)
$$
of morphisms of quasi-coherent sheaves on~$X$.
\end{proof}

\begin{lem} \label{alf-contratensor-projection-lemma}
 Let $f\:Y\rarrow X$ be a flat morphism of quasi-compact semi-separated
schemes, $\M$ be a quasi-coherent sheaf on $X$, and\/ $\gG$ be
an antilocally flat contraherent cosheaf on~$Y$.
 Then the natural morphism~\eqref{contratensor-projection-comparison}
is an isomorphism
\begin{equation} \label{alf-contratensor-projection-eqn}
 \M\ocn_X f_!\gG\simeq f_*(f^*\M\ocn_Y\gG)
\end{equation}
of quasi-coherent sheaves on~$X$.
\end{lem}

\begin{proof}
 The argument is similar to the proof of
Lemma~\ref{flat-al-contratensor-projection-lemma}.
 We check that both the left-hand and the right-hand sides
of~\eqref{alf-contratensor-projection-eqn} are exact functors
$Y\ctrh_\alf\rarrow X\qcoh$ of an antilocally flat contraherent
cosheaf $\gG$ on $Y$ for any fixed quasi-coherent sheaf $\M$ on~$X$.
 Indeed, the functor $f_!\:Y\ctrh_\alf\rarrow X\ctrh_\alf$ is exact
by Corollary~\ref{clp-direct}(c).
 In order to show that the functor $\M\ocn_X{-}\,\:X\ctrh_\alf\rarrow
X\qcoh$ is exact, we pick an injective quasi-coherent sheaf $\J$ on $X$
and use the isomorphism $\Hom_X(\M\ocn_X\gF\;\J)\simeq
\Hom^X(\gF,\fHom_X(\M,\J))$ for a locally contraherent cosheaf $\gF$
on $X$ (formula~\eqref{fHom-contratensor-adjunction}).
 According to Section~\ref{fHom-subsection}, \,$\fHom_X(\M,\J)$ is
a locally cotorsion contraherent cosheaf on~$X$.
 Hence the functor $\gF\longmapsto\Hom^X(\gF,\fHom_X(\M,\J))$ is exact
for $\gF\in X\ctrh_\alf$ by Corollary~\ref{clf-characterizations}(a).
 This proves that the left-hand side
of~\eqref{alf-contratensor-projection-eqn} is exact as a functor
of~$\gG$.

 Concerning the right-hand side, the argument from the previous
paragraph shows that the functor $\gG\longmapsto f^*\M\ocn_Y\gG$ is
exact on the category of antilocally flat contraherent cosheaves
$\gG\in Y\ctrh_\alf$.
 Similarly to the proof of
Lemma~\ref{flat-al-contratensor-projection-lemma}, one uses
Corollary~\ref{clf-cor}(c) and formula~\eqref{contratensor-projection}
in order to show that the quasi-coherent sheaf $f^*\M\ocn_Y\gG$ on $Y$
is dilute.
 Hence the right-hand side
of~\eqref{alf-contratensor-projection-eqn} is exact as a functor
of $\gG$, too.
 The argument finishes similarly to the proof of
Lemma~\ref{flat-al-contratensor-projection-lemma}, using
Corollary~\ref{clf-cor}(c) and formula~\eqref{contratensor-projection}
again.
\end{proof}

\begin{lem} \label{open-embed-fq-flat-ctrtensor-projection-lemma}
 Let $j\:Y\rarrow X$ be an open embedding of Noetherian schemes of
finite Krull dimension, $\E$ be a flasque quasi-coherent sheaf on $X$,
and\/ $\gG$ be a flat contraherent cosheaf on~$Y$.
 Then the natural morphism~\eqref{contratensor-projection-comparison}
is an isomorphism
\begin{equation} \label{open-embed-fq-flat-ctrtensor-projection-eqn}
 \E\ocn_X j_!\gG\simeq j_*(j^*\E\ocn_Y\gG)
\end{equation}
of quasi-coherent sheaves on~$X$.
\end{lem}

\begin{proof}
 Once again, we show that both the left-hand and the right-hand
sides of~\eqref{open-embed-fq-flat-ctrtensor-projection-eqn} are exact
functors $Y\ctrh^\fl\rarrow X\qcoh$ of a flat contraherent cosheaf
$\gG$ on $Y$ for any fixed flasque quasi-coherent sheaf $\E$ on~$X$.
 Indeed, the functor $j_!\:Y\ctrh^\fl\rarrow X\ctrh^\fl$ is exact
by Corollary~\ref{finite-krull-flat-direct}(b).
 In order to show that the functor $\M\ocn_X\nobreak{-}\,\:\allowbreak
X\ctrh^\fl\rarrow X\qcoh$ is exact for any quasi-coherent sheaf $\M$
on $X$, we pick an injective quasi-coherent sheaf $\J$ on $X$ and use
the isomorphism $\Hom_X(\M\ocn_X\gF\;\J)\simeq
\Hom^X(\gF,\fHom_X(\M,\J))$ for a locally contraherent cosheaf $\gF$
on $X$ (formula~\eqref{fHom-contratensor-adjunction}).
 According to Section~\ref{fHom-subsection}, \,$\fHom_X(\M,\J)$ is
a locally cotorsion contraherent cosheaf on~$X$.
 Hence the functor $\gF\longmapsto\Hom^X(\gF,\fHom_X(\M,\J))$ is exact
for $\gF\in X\ctrh^\fl$ by Corollary~\ref{finite-krull-flat-clf-cor}(a).
 This proves that the left-hand side
of~\eqref{open-embed-fq-flat-ctrtensor-projection-eqn} is exact as
a functor of~$\gG$.

 Concerning the right-hand side, the argument from the previous
paragraph shows that the functor $\gG\longmapsto j^*\E\ocn_Y\gG$ is
exact on the category of flat contraherent cosheaves
$\gG\in Y\ctrh^\fl$.
 Let us show that the quasi-coherent sheaf $f^*\E\ocn_Y\gG$ on $Y$
is flasque.
 The class of flasque quasi-coherent sheaves is closed under extensions
in $Y\qcoh$, so in view of
Corollary~\ref{finite-krull-flat-clf-cor}(b) it suffices to check that
$j^*\E\ocn_Y k_!\gH$ is a flasque quasi-coherent sheaf on $Y$ for every
affine open subscheme $V\sub Y$ with the open embedding morphism
$k\:V\rarrow Y$ and any flat contraherent cosheaf $\gH$ on~$V$.
 By~\eqref{flat-contratensor-projection}, we have $j^*\E\ocn_Y k_!\gH
\simeq k_*(k^*j^*\E\ocn_V\gH)$.
 The quasi-coherent sheaf $k^*j^*\E\ocn_V\gH$ on $V$ is flasque by
Lemma~\ref{co-flasque-preservation}(d), and it remains to point out
that the direct images of flasque sheaves are flasque.
 As the direct images of flasque (quasi-coherent) sheaves are exact
functors, this proves that the right-hand side
of~\eqref{open-embed-fq-flat-ctrtensor-projection-eqn} is exact as
a functor of $\gG$, too.

 Using Corollary~\ref{finite-krull-flat-clf-cor}(b) again, it remains to
prove the assertion of the lemma for $\gG=k_!\gH$, where $\gH$ is a flat
contraherent cosheaf on~$V$.
 In this case, the desired assertion follows from
formula~\eqref{flat-contratensor-projection} applied to the two open
embeddings of an affine scheme $k\:V\rarrow Y$ and $jk\:V\rarrow X$,
similarly to the final paragraph of the proof of
Lemma~\ref{flat-al-contratensor-projection-lemma}.
\end{proof}

\subsection{Semico-semicontra correspondence}
 In this section, we consider a semi-separated Noetherian scheme $X$
with a dualizing complex $\D_X^\bu$, a quasi-compact semi-separated
scheme $Y$, and a flat morphism of schemes $\pi\:Y\rarrow X$.
 Notice that $\sD^\co(X\qcoh)=\sD^\bco(X\qcoh)$ and
$\sD^\ctr(X\lcth_\bW)=\sD^\bctr(X\lcth_\bW)$ by
Theorem~\ref{derived-inj-proj-resolutions}(a,d).
 So the constructions of semiderived categories based on the Becker
co/contraderived categories
(as in Sections~\ref{antilocal-semicoderived-subsect}\+-%
\ref{antilocal-semicontraderived-subsect}) agree with the similar
constructions based on the Positselski co/contraderived categories in
this case.

\begin{thm} \label{semico-semicontra-diagonal}
 The choice of a dualizing complex\/ $\D_X^\bu$ induces a natural
equivalence of triangulated categories\/ $\sD^\si_X(Y\qcoh)\simeq
\sD^\si_X(Y\ctrh)$.
 This triangulated equivalence forms a commutative square diagram with
the triangulated equivalence\/ $\sD^\co(X\qcoh)\simeq
\sD^\ctr(X\ctrh)$ from Theorem~\textup{\ref{co-contra-dualizing}}
or~\textup{\ref{non-semi-separated-co-contra}} and the derived
direct image functors\/
$\boR\pi_*$~\eqref{semicoderived-to-X-derived-direct-image}
and\/ $\boL\pi_!$~\eqref{semicontraderived-to-X-derived-direct-image},
\begin{equation}
\begin{gathered}
 \xymatrix{
  \sD^\si_X(Y\qcoh) \ar@{=}[r] \ar[d]_{\boR\pi_*}
  & \sD^\si_X(Y\ctrh) \ar[d]^{\boL\pi_!} \\
  \sD^\co(X\qcoh) \ar@{=}[r] & \sD^\ctr(X\ctrh)
 }
\end{gathered}
\end{equation}
or, with the notation for the derived functors providing
the triangulated equivalences,
\begin{equation} \label{semico-semicontra-diagonal-direct-images-diagr}
\begin{gathered}
 \xymatrix{
  \sD^\si_X(Y\qcoh) \ar[d]_{\boR\pi_*}
  \ar@<2pt>[rrrr]^{\boR\fHom_Y(\pi^*\D_X^\subbu,{-})}
  &&&& \sD^\si_X(Y\ctrh) \ar[d]^{\boL\pi_!} 
  \ar@<2pt>[llll]^{\pi^*\D_X^\subbu\ocn_Y^\boL{-}} \\
  \sD^\co(X\qcoh)
  \ar@<2pt>[rrrr]^{\boR\fHom_X(\D_X^\subbu,{-})}
  &&&& \sD^\ctr(X\ctrh) \ar@<2pt>[llll]^{\D_X^\subbu\ocn_X^\boL{-}}
 }
\end{gathered}
\end{equation}
\end{thm}

\begin{proof}
 Notice first of all that, according to
Corollary~\ref{semicontraderived-independent-of-covering}
and Theorem~\ref{semicontraderived-lcta-lct-equivalent}, one has
$\sD^\si_X(Y\ctrh)\simeq\sD^\si_X(Y\lcth_\bT)\simeq
\sD^\si_X(Y\lcth_\bT^\lct)\simeq\sD^\si_X(Y\ctrh^\lct)$ for any
open covering $\bT$ of the scheme~$Y$.
 Now let us construct a pair of adjoint functors between
the triangulated categories $\sD^\si_X(Y\qcoh)$ and
$\sD^\si_X(Y\ctrh)\simeq\sD^\si_X(Y\ctrh^\lct)$ and show that
they are mutually inverse equivalences.

 Given a complex $\N^\bu$ in $Y\qcoh$, we pick a complex $\J^\bu$
in $Y\qcoh^\inj$ together with a morphism of complexes
$\N^\bu\rarrow\J^\bu$ with a cone semi(co)acyclic over~$X$.
 Then we assign to $\N^\bu$ the total complex of the bicomplex
$\fHom_Y(\pi^*\D_X^\bu,\J^\bu)$ in the additive/exact category
$Y\ctrh^\lct$.
 Here $\pi^*\D_X^\bu\in\Hot^\b(Y\qcoh)$ is the inverse image of
complex $\D_X^\bu\in\Hot^\b(X\qcoh^\inj)$ with respect to the flat
morphism of schemes $\pi\:Y\rarrow X$.

 Given a complex $\Q^\bu$ in $Y\ctrh$ (or even in $Y\lcth_\bT$),
we pick a complex $\gF^\bu$ in $Y\ctrh_\alf$ together with a morphism
of complexes $\gF^\bu\rarrow\Q^\bu$ with a cone semi(contra)acyclic
over~$X$ (cf.\ the proof of
Theorem~\ref{semicontraderived-lcta-lct-equivalent}).
 Then we assign to $\Q^\bu$ the total complex of the bicomplex
$\pi^*\D_X^\bu\ocn_Y\gF^\bu$ in the abelian category $Y\qcoh$.

 Let $\E^\bu$ be an arbitrary finite complex of quasi-coherent
sheaves on~$X$.
 Let us show that the complex $\fHom_Y(\pi^*\E^\bu,\J^\bu)$ is
semicontraacyclic over $X$ whenever a complex $\J^\bu$ in $Y\qcoh^\inj$
is semicoacyclic over~$X$.
 Indeed, let $Y=\bigcup_\alpha V_\alpha$ be a finite affine open
covering of the scheme~$Y$.
 Denote by $j_\alpha\:V_\alpha\rarrow X$ the open embedding morphisms.
 By Lemma~\ref{semicoacyclic-of-injectives-antilocal},
the complex $\J^\bu$ is a direct summand of a finitely iterated
extension of the direct images of semicoacyclic complexes of
injective quasi-coherent sheaves $\J_\alpha^\bu$ from~$V_\alpha$.
 As such an extension is always termwise split, it remains to check
that the complex $\fHom_Y(\pi^*\E^\bu,j_\alpha{}_*\J_\alpha^\bu)$ is
semicontraacyclic.
 Here the semiacyclicity condition on the complex $\J_\alpha^\bu$
means that the complex $\pi_*j_\alpha{}_*\J_\alpha^\bu$ is
Becker-coacyclic in $X\qcoh$.

 According to~\eqref{inj-fHom-projection}, we have an isomorphism
$$
 \fHom_Y(\pi^*\E^\bu,j_\alpha{}_*\J_\alpha^\bu)\simeq
 j_\alpha{}_!\fHom_{V_\alpha}(j_\alpha^*\pi^*\E^\bu,\J_\alpha^\bu)
$$
of complexes of locally cotorsion contraherent cosheaves on~$Y$.
 As the functor~$j_\alpha{}_!$ takes semiacyclic complexes in
$V_\alpha\ctrh^\lct$ to semiacyclic complexes in $Y\ctrh^\lct$ by
Lemma~\ref{semicontraacyclic-direct-image}(b), it suffices to show that
the complex $\fHom_{V_\alpha}(j_\alpha^*\pi^*\E^\bu,\J_\alpha^\bu)$
is semicontraacyclic over~$X$.
 The latter statement (which we need to prove) means that the complex
$\pi_!j_\alpha{}_!\fHom_{V_\alpha}(j_\alpha^*\pi^*\E^\bu,
\J_\alpha^\bu)$ is Becker-contraacyclic in $X\ctrh^\lct$.
 Applying~\eqref{inj-fHom-projection} again, we have
$$
 \pi_!j_\alpha{}_!\fHom_{V_\alpha}(j_\alpha^*\pi^*\E^\bu,
 \J_\alpha^\bu)\simeq\fHom_X(\E^\bu,\pi_*j_\alpha{}_*\J_\alpha^\bu),
$$
and it remains to observe that the complex
$\pi_*j_\alpha{}_*\J_\alpha^\bu$ is contractible as
a Becker-coacyclic complex of injective quasi-coherent sheaves on~$X$.

 Let us show that the complex $\pi^*\E^\bu\ocn_Y\gF^\bu$ is
semicoacyclic over $X$ whenever a complex $\gF^\bu$ in $Y\ctrh_\alf$
is semicontraacyclic over~$X$.
 Indeed, by Lemma~\ref{semicontraacyclic-of-projectives-antilocal}(c),
the complex $\gF^\bu$ is a direct summand of a finitely iterated
extension of the direct images of semicontraacyclic complexes of
antilocally flat contraherent cosheaves $\gF_\alpha^\bu$
from~$V_\alpha$.
 The direct image functors~$j_\alpha{}_!$ take antilocally flat
contraherent cosheaves on $V_\alpha$ to antilocally flat contraherent 
cosheaves on $Y$, as explained in the beginning of
Section~\ref{clf-subsection}.
 In view of the adjunction~\eqref{fHom-contratensor-adjunction}
and Corollary~\ref{clf-characterizations}(a), the functor
$\M\ocn_Y{-}\:Y\ctrh_\alf\rarrow Y\qcoh$ is exact for any
quasi-coherent sheaf $\M$ on~$Y$.
 Therefore, it remains to check that the complex
$\pi^*\E^\bu\ocn_Y j_\alpha{}_!\gF_\alpha^\bu$ is semicoacyclic.
 Here the semiacyclicity condition on the complex $\gF_\alpha^\bu$
means that the complex $\pi_!j_\alpha{}_!\gF_\alpha^\bu$ is
Becker-contraacyclic in $X\ctrh$.

 According to~\eqref{contratensor-projection}, we have an isomorphism
$$
 \pi^*\E^\bu\ocn_Y j_\alpha{}_!\gF_\alpha^\bu\simeq
 j_\alpha{}_*(j_\alpha^*\pi^*\E^\bu\ocn_{V_\alpha}\gF_\alpha^\bu)
$$
of complexes of quasi-coherent sheaves on~$Y$.
 As the functor~$j_\alpha{}_*$ takes semiacyclic complexes in
$V_\alpha\qcoh$ to semiacyclic complexes in $Y\qcoh$ by
Lemma~\ref{semicoacyclic-direct-image}, it suffices to show that
the complex $j_\alpha^*\pi^*\E^\bu\ocn_{V_\alpha}\gF_\alpha^\bu$
is semicoacyclic over~$X$.
 The latter statement means that the complex
$\pi_*j_\alpha{}_*(j_\alpha^*\pi^*\E^\bu\ocn_{V_\alpha}\gF_\alpha^\bu)$
is Becker-coacyclic in $X\qcoh$.
 Applying~\eqref{contratensor-projection} to
the morphism~$\pi j_\alpha$, we have
$$
 \pi_*j_\alpha{}_*(j_\alpha^*\pi^*\E^\bu\ocn_{V_\alpha}\gF_\alpha^\bu)
 \simeq \E^\bu\ocn_X\pi_!j_\alpha{}_!\gF_\alpha^\bu.
$$

 By Corollary~\ref{clp-direct}(c), the functor~$\pi_!$ also takes
antilocally flat contraherent cosheaves on $Y$ to antilocally flat
contraherent cosheaves on~$X$; so $\pi_!j_\alpha{}_!\gF_\alpha^\bu$
is a Becker-contraacyclic complex of antilocally flat
contraherent cosheaves on~$X$.
 In order to show that the complex
$\E^\bu\ocn_X\pi_!j_\alpha{}_!\gF_\alpha^\bu$ is Becker-coacyclic
in $X\qcoh$, we need to check that the complex of abelian groups
$\Hom_X(\E^\bu\ocn_X\pi_!j_\alpha{}_!\gF_\alpha^\bu\;\I^\bu)$ is
acyclic for any complex $\I^\bu$ in $X\qcoh^\inj$.
 Indeed, the adjunction~\eqref{fHom-contratensor-adjunction} tells us
that
$$
 \Hom_X(\E^\bu\ocn_X\pi_!j_\alpha{}_!\gF_\alpha^\bu\;\I^\bu)\simeq
 \Hom^X(\pi_!j_\alpha{}_!\gF_\alpha^\bu,\fHom_X(\E^\bu,\I^\bu)).
$$
 Now the complex $\pi_!j_\alpha{}_!\gF_\alpha^\bu$ is acyclic in
$X\ctrh_\alf$ by Corollary~\ref{acycl=bctracycl=bcoacycl-in-alf},
the complex $\fHom_X(\E^\bu,\I^\bu)$ is a complex of locally cotorsion
contraherent cosheaves on $X$ according to
Section~\ref{fHom-subsection}, and we can refer to
Corollary~\ref{acyclic-in-alf-all-of-lct-loc-contraherent-pair}
together with Lemma~\ref{Ext-1-as-homotopy-Hom} to the effect
that any morphism of complexes from $\pi_!j_\alpha{}_!\gF_\alpha^\bu$
to $\fHom_X(\E^\bu,\I^\bu)$ is homotopic to zero.

 For any finite complex of quasi-coherent sheaves $\E^\bu$ on $X$,
we have constructed a pair of well-defined derived functors
$$
 \boR\fHom_Y(\pi^*\E^\bu,{-})\:\sD^\si_X(Y\qcoh)\lrarrow
 \sD^\si_X(Y\ctrh)
$$
and
$$
 \pi^*\E^\bu\ocn_Y^\boL{-}\,\:\sD^\si(Y\ctrh)\lrarrow
 \sD^\si(Y\qcoh).
$$
 It follows easily from
the adjunction~\eqref{fHom-contratensor-adjunction} that
the triangulated functor $\pi^*\E^\bu\ocn_Y^\boL{-}$ is left adjoint to
the triangulated functor $\boR\fHom_Y(\pi^*\E^\bu,{-})$.
 It remains to show that, when $\E^\bu=\D_X^\bu$, the adjunction
morphisms for the pair of derived functors
$\boR\fHom_Y(\pi^*\D_X^\bu,{-})$ and $\pi^*\D_X^\bu\ocn_Y^\boL{-}$
are isomorphisms in the semiderived categories $\sD^\si_X(Y\qcoh)$
and $\sD^\si_X(Y\ctrh)$.

 In view of the natural isomorphisms~\eqref{inj-fHom-projection}
and~\eqref{contratensor-projection}, we have commutative square
diagrams of triangulated functors
\begin{equation} \label{embedding-of-affine-semiderived-RfHom-square}
\begin{gathered}
 \xymatrix{
  \sD^\si_X(V_\alpha\qcoh)
  \ar[rrr]^{\boR\fHom_{V_\alpha}(j_\alpha^*\pi^*\E^\subbu,{-})}
  \ar[d]_{j_\alpha{}_*}
  &&& \sD^\si_X(V_\alpha\ctrh) \ar[d]^{j_\alpha{}_!} \\
  \sD^\si_X(Y\qcoh)
  \ar[rrr]^{\boR\fHom_Y(\pi^*\E^\subbu,{-})}
  &&& \sD^\si_X(Y\ctrh)
 }
\end{gathered}
\end{equation}
and
\begin{equation} \label{embedding-of-affine-semiderived-Locn-square}
\begin{gathered}
 \xymatrix{
  \sD^\si_X(V_\alpha\qcoh) \ar[d]_{j_\alpha{}_*}
  &&& \sD^\si_X(V_\alpha\ctrh)
  \ar[lll]_{j_\alpha^*\pi^*\E^\subbu\ocn_{V_\alpha}^\boL{-}}
  \ar[d]^{j_\alpha{}_!} \\
  \sD^\si_X(Y\qcoh) &&& \sD^\si_X(Y\ctrh)
  \ar[lll]_{\pi^*\E^\subbu\ocn_Y^\boL{-}}
 }
\end{gathered}
\end{equation}
 Here the triangulated functors~$j_\alpha{}_*$ and~$j_\alpha{}_!$
are well-defined in view of Lemmas~\ref{semicoacyclic-direct-image}
and~\ref{semicontraacyclic-direct-image}(a).
 Moreover, the commutative
squares~(\ref{embedding-of-affine-semiderived-RfHom-square}\+-%
\ref{embedding-of-affine-semiderived-Locn-square}) are compatible
with the adjunctions between the functors $\ocn^\boL$ and $\boR\fHom$,
in the sense that the direct image functors~$j_\alpha{}_*$
and~$j_\alpha{}_!$ take adjunction morphisms for the adjoint pair
in the upper horizontal lines to the adjunction morphisms for
the adjoint pair in the lower horizontal lines of the two diagrams.

 Furthermore, the same isomorphism~\eqref{inj-fHom-projection}
from Section~\ref{compatibility-subsect} and
the isomorphism~\eqref{alf-contratensor-projection-eqn} from
Lemma~\ref{alf-contratensor-projection-lemma} induce commutative
square diagrams of triangulated functors
\begin{equation} \label{projection-to-X-semiderived-RfHom-square}
\begin{gathered}
 \xymatrix{
  \sD^\si_X(Y\qcoh) \ar[rrr]^{\boR\fHom_Y(\pi^*\E^\subbu,{-})}
  \ar[d]_{\boR \pi_*} &&& \sD^\si_X(Y\ctrh) \ar[d]^{\boL \pi_!} \\
  \sD^\bco(X\qcoh) \ar[rrr]^{\boR\fHom_X(\E^\subbu,{-})}
  &&& \sD^\bctr(X\ctrh)
 }
\end{gathered}
\end{equation}
and
\begin{equation} \label{projection-to-X-semiderived-Locn-square}
\begin{gathered}
 \xymatrix{
  \sD^\si_X(Y\qcoh) \ar[d]_{\boR\pi_*} &&& \sD^\si_X(Y\ctrh)
  \ar[lll]_{\pi^*\E^\subbu\ocn_Y^\boL{-}} \ar[d]^{\boL\pi_!} \\
  \sD^\bco(X\qcoh) &&& \sD^\bctr(X\ctrh)
  \ar[lll]_{\E^\subbu\ocn_X^\boL{-}}
 }
\end{gathered}
\end{equation}
 Here the derived functors $\boR\pi_*$ and $\boL\pi_!$ were constructed
in~\eqref{semicoderived-to-X-derived-direct-image}
and~\eqref{semicontraderived-to-X-derived-direct-image}.
 We recall that the morphism~$\pi$ is flat by assumption; so the functor
$\pi_*\:Y\qcoh\rarrow X\qcoh$ preserves injectivity of quasi-coherent
sheaves, while the functor $\pi_!\:Y\ctrh_\al\rarrow X\ctrh_\al$
preserves antilocal flatness of contraherent cosheaves (according to
Corollary~\ref{clp-direct}(c)).
 The commutative
squares~(\ref{projection-to-X-semiderived-RfHom-square}\+-%
\ref{projection-to-X-semiderived-Locn-square}) are compatible
with the adjunctions between the functors $\ocn^\boL$ and $\boR\fHom$,
in the sense that the direct image functors~$\pi_*$ and~$\pi_!$ take
adjunction morphisms for the adjoint pair in the upper horizontal
lines to the adjunction morphisms for the adjoint pair in the lower
horizontal lines of the two diagrams.

 Now let us prove that the adjunction morphism
$$
 \pi^*\D_X^\bu\ocn_Y^\boL\boR\fHom_Y(\pi^*\D_X^\bu,\N^\bu)\lrarrow\N^\bu
$$
is an isomorphism in $\sD^\si_X(Y\qcoh)$ for every object
$\N^\bu\in\sD^\si_X(Y\qcoh)$.
 It follows from Theorem~\ref{quasi-coherent-becker-coderived}
and Corollary~\ref{complexes-of-quasi-injective-antilocal} that
the Becker coderived category $\sD^\bco(Y\qcoh)$ coincides with
its minimal thick subcategory generated by the subcategories
$j_\alpha{}_*\sD^\bco(V_\alpha\qcoh)$.
 Consequently, the semiderived category $\sD^\si_X(Y\qcoh)$ coincides
with its minimal thick subcategory generated by the subcategories
$j_\alpha{}_*\sD^\si_X(V_\alpha\qcoh)$.
 Therefore, it suffices to check that the adjunction morphism
{\hbadness=2075
$$
 \pi^*\D_X^\bu\ocn_Y^\boL
 \boR\fHom_Y(\pi^*\D_X^\bu,j_\alpha{}_*\N_\alpha^\bu)
 \lrarrow j_\alpha{}_*\N_\alpha^\bu
$$
is} an isomorphism in $\sD^\si_X(Y\qcoh)$ for every object
$\N_\alpha^\bu\in\sD^\si_X(V_\alpha\qcoh)$.
 In view of the commutative
diagrams~(\ref{embedding-of-affine-semiderived-RfHom-square}\+-%
\ref{embedding-of-affine-semiderived-Locn-square}) and their
compatibility with the adjunctions, the question reduces to
showing that the adjunction morphism
$$
 j_\alpha^*\pi^*\D_X^\bu\ocn_{V_\alpha}^\boL
 \boR\fHom_{V_\alpha}(j_\alpha^*\pi^*\D_X^\bu,\N_\alpha^\bu)
 \lrarrow\N_\alpha^\bu
$$
is an isomorphism in $\sD^\si_X(V_\alpha\qcoh)$.

 By the definition of the semicoderived category
$\sD^\si_X(V_\alpha\qcoh)$ for an affine morphism of schemes
$\pi j_\alpha\:V_\alpha\rarrow X$, a morphism~$g$ in
$\sD^\si_X(V_\alpha\qcoh)$ is an isomorphism if and only if
the morphism $(\pi j_\alpha)_*(g)=\boR\pi_* j_\alpha{}_*(g)$
is an isomorphism in $\sD^\bco(X\qcoh)$.
 (In other words, the triangulated functor $(\pi j_\alpha)_*\:
\sD^\si_X(V_\alpha\qcoh)\rarrow\sD^\bco(X\qcoh)$ is conservative.)
 Consequently, in view of the commutative
diagrams~(\ref{embedding-of-affine-semiderived-RfHom-square}\+-%
\ref{projection-to-X-semiderived-Locn-square}) and their
compatibility with the adjunctions, the question reduces to
showing that the adjunction morphism
$$
 \D_X^\bu\ocn^\boL_X
 \boR\fHom_X(\D_X^\bu,\boR\pi_*j_\alpha{}_*\N_\alpha^\bu)
 \lrarrow\boR\pi_*j_\alpha{}_*\N_\alpha^\bu
$$
is an isomorphism in $\sD^\bco(X\qcoh)$.
 The latter fact was established in the proof of
Theorem~\ref{co-contra-dualizing}
or~\ref{non-semi-separated-co-contra}.

 Finally, let us show that the adjunction morphism
$$
 \Q^\bu\lrarrow\boR\fHom_Y(\pi^*\D_X^\bu\;
 \pi^*\D_X^\bu\ocn_Y^\boL\Q^\bu)
$$
is an isomorphism in $\sD^\si_X(Y\ctrh)$ for every object
$\Q^\bu\in\sD_X^\si(Y\ctrh)$.
 It follows from
Corollaries~\ref{becker-contraderived-of-lcta-lct-well-behaved}(a)
and~\ref{complexes-of-prj-alf-antilocal}(a) that the Becker
contraderived category $\sD^\bctr(Y\ctrh)$ coincides with its minimal
thick subcategory generated by the subcategories
$j_\alpha{}_!\sD^\bctr(V_\alpha\ctrh)$.
 Consequently, the semiderived category $\sD^\si_X(Y\ctrh)$ coincides
with its minimal thick subcategory generated by the subcategories
$j_\alpha{}_!\sD^\si_X(V_\alpha\ctrh)$.
 Therefore, it suffices to check that the adjunction morphism
{\hbadness=2075
$$
 j_\alpha{}_!\Q_\alpha^\bu\lrarrow\boR\fHom_Y(\pi^*\D_X^\bu\;
 \pi^*\D_X^\bu\ocn_Y^\boL j_\alpha{}_!\Q_\alpha^\bu)
$$
is} an isomorphism in $\sD^\si_X(Y\ctrh)$ for every object
$\Q_\alpha^\bu\in\sD_X^\si(V_\alpha\ctrh)$.
 In view of the commutative
diagrams~(\ref{embedding-of-affine-semiderived-RfHom-square}\+-%
\ref{embedding-of-affine-semiderived-Locn-square}) and their
compatibility with the adjunctions, the question reduces to
showing that the adjunction morphism
$$
 \Q_\alpha^\bu\lrarrow\boR\fHom_{V_\alpha}(j_\alpha^*\pi^*\D_X^\bu\;
 j_\alpha^*\pi^*\D_X^\bu\ocn_{V_\alpha}^\boL\Q_\alpha^\bu)
$$
is an isomorphism in $\sD^\si_X(V_\alpha\ctrh)$.

 By the definition of the semicontraderived category
$\sD^\si_X(V_\alpha\ctrh)$ for an affine morphism of schemes
$\pi j_\alpha\:V_\alpha\rarrow X$, a morphism~$h$ in
$\sD^\si_X(V_\alpha\ctrh)$ is an isomorphism if and only if
the morphism $(\pi j_\alpha)_!(h)=\boL\pi_!j_\alpha{}_!(h)$ is
an isomorphism in $\sD^\bctr(X\ctrh)$.
 (In other words, the triangulated functor $(\pi j_\alpha)_!\:
\sD^\si_X(V_\alpha\ctrh)\rarrow\sD^\bctr(X\ctrh)$ is conservative.)
  Consequently, in view of the commutative
diagrams~(\ref{embedding-of-affine-semiderived-RfHom-square}\+-%
\ref{projection-to-X-semiderived-Locn-square}) and their
compatibility with the adjunctions, the question reduces to
showing that the adjunction morphism
$$
 \boL\pi_!j_\alpha{}_!\Q_\alpha^\bu\lrarrow
 \boR\fHom_X(\D_X^\bu\;\D_X^\bu\ocn_X^\boL
 \boL\pi_!j_\alpha{}_!\Q_\alpha^\bu)
$$
is an isomorphism in $\sD^\bctr(X\ctrh)$.
 The latter fact was established in the proof of
Theorem~\ref{co-contra-dualizing}
or~\ref{non-semi-separated-co-contra}.
 To facilitate the comparison, let us just mention that the classes
of flat and antilocally flat contraherent cosheaves coincide on
a semi-separated Noetherian scheme $X$ of finite Krull dimension
by Corollary~\ref{finite-krull-flat-contraherent}(a)
(cf.\ Remark~\ref{flat-antilocally-flat-remark}).
\end{proof}

\subsection{Sheaves injective over the base}
\label{sheaves-injective-over-base-subsect}
 The following series of two lemmas and two corollaries is to be
compared with~\cite[Section~A.4]{Psemten}.

\begin{lem} \label{quasi-injectivity-over-base-tensor-product-affine}
 Let $X$ be a Noetherian affine scheme, $Y$ be an affine scheme, and
$\pi\:Y\rarrow X$ be a morphism of schemes.
 Let\/ $\K$ be a quasi-coherent sheaf on $Y$ and\/ $\F$ be a flat
quasi-coherent sheaf on~$Y$.
 Assume that the quasi-coherent sheaf $\pi_*\K$ on $X$ is injective.
 Then the quasi-coherent sheaf $\pi_*(\F\ot_{\O_Y}\K)$ on $X$ is
injective, too.
\end{lem}

\begin{proof}
 The argument from~\cite[Lemma~A.23]{Psemten} applies.
 Essentially, the claim is that for any commutative ring homomorphism
$R\rarrow S$ with a Noetherian ring $R$, any flat $S$\+module $F$,
and any $S$\+module $K$ that is injective over $R$, the $S$\+module
$F\ot_SK$ is also injective over~$R$.
 This is the assertion of
Lemma~\ref{injective-over-base-ring-tensor-product} above.
\end{proof}

\begin{cor} \label{quasi-injectivity-over-base-tensor-product}
 Let $X$ be a locally Noetherian scheme and $\pi\:Y\rarrow X$ be
an affine morphism of schemes.
 Let\/ $\K$ be a quasi-coherent sheaf on $Y$ and\/ $\F$ be a flat
quasi-coherent sheaf on~$Y$.
 Assume that the quasi-coherent sheaf $\pi_*\K$ on $X$ is injective.
 Then the quasi-coherent sheaf $\pi_*(\F\ot_{\O_Y}\K)$ on $X$ is
injective, too.
\end{cor}

\begin{proof}
 The point is that injectivity of quasi-coherent sheaves on a locally
Noetherian scheme $X$ is a local property: a quasi-coherent sheaf on
$X$ is injective if and only if its restrictions to open subschemes
forming some given open covering of $X$ are injective.
 In the situation at hand, it suffices to cover $X$ with affine open
subschemes $U_\alpha\sub X$, put $V_\alpha=U_\alpha\times_XY$, and
apply Lemma~\ref{quasi-injectivity-over-base-tensor-product-affine}
to the morphisms $\pi_\alpha\:V_\alpha\rarrow U_\alpha$.
\end{proof}

\begin{lem} \label{affine-composition-impies-affine}
 Let $f\:Y\rarrow X$ and $g\:V\rarrow Y$ be morphisms of schemes.
 Assume that the scheme $Y$ is semi-separated and the morphism
$fg\:V\rarrow X$ is affine.
 Then the morphism $g\:V\rarrow Y$ is affine, too.
\end{lem}

\begin{proof}
 This is~\cite[Lemma Tag~01SG]{SP}.
 We follow~\cite[proof of Lemma~A.21]{Psemten}.
 Let $W\sub Y$ be an affine open subscheme in~$Y$.
 Then the morphism $W\times_XV\rarrow W$ is affine as a base change
of an affine morphism $V\rarrow X$.
 Hence the scheme $W\times_XV$ is affine.
 Now we have $W\times_YV=(W\times_XV)\times_{Y\times Y}Y$.
 As the morphism $Y\rarrow Y\times Y$ is affine (the scheme $Y$ being
semi-separated by assumption) and the scheme $W\times_XV$ is affine,
it folows that the scheme $W\times_YV$ is affine, too.
\end{proof}

 The following corollary is to be compared
with~\cite[Lemma~A.21]{Psemten}.

\begin{cor} \label{quasi-injectivity-over-base-ascent}
 Let $X$ be a locally Noetherian scheme, $Y$ be a semi-separated scheme,
and $\pi\:Y\rarrow X$ be an affine morphism of schemes.
 Let $V\sub Y$ be an open subscheme in $Y$ with the open
embedding morphism $j\:V\rarrow Y$ such that the composition
$\pi j\:V\rarrow X$ is an affine morphism.
 Let\/ $\K$ be a quasi-coherent sheaf on $Y$ such that the quasi-coherent
sheaf $\pi_*\K$ on $X$ is injective.
 Then the quasi-coherent sheaf $\pi_*j_*j^*\K$ on $X$ is injective, too.
\end{cor}

\begin{proof}
 By Lemma~\ref{affine-composition-impies-affine}, the open embedding
morphism $j\:V\rarrow Y$ is affine.
 Now we have natural isomorphisms $j_*j^*\K\simeq
j_*(\O_V\ot_{\O_V}j^*\K)\simeq j_*\O_V\ot_{\O_Y}\K$ of quasi-coherent
sheaves on~$Y$, and it remains to apply
Corollary~\ref{quasi-injectivity-over-base-tensor-product} to
the quasi-coherent sheaf $\K$ and the flat quasi-coherent sheaf
$\F=j_*\O_V$ on~$Y$.
\end{proof}

 Corollary~\ref{quasi-injectivity-over-base-ascent} provides
the ascent for the property of injectivity of quasi-coherent
sheaves over a Noetherian base scheme.
 The next series of a lemma and two corollaries takes care of
the descent.

\begin{lem} \label{quasi-injectivity-over-base-pure-submodule-affine}
 Let $X$ be a Noetherian affine scheme, $Y$ be an affine scheme,
and $\pi\:Y\rarrow X$ be a morphism of schemes.
 Let\/ $\K$ be a quasi-coherent sheaf on $Y$, and let\/ $0\rarrow\cH
\rarrow\G\rarrow\F\rarrow0$ be a short exact sequence of quasi-coherent
sheaves on $Y$ such that\/ $\F$ is a flat quasi-coherent sheaf.
 Assume that the quasi-coherent sheaf $\pi_*(\G\ot_{\O_Y}\K)$ on $X$
is injective.
 Then the quasi-coherent sheaf $\pi_*(\cH\ot_{\O_Y}\K)$ on $X$ is
injective, too.
\end{lem}

\begin{proof}
 This is a restatement of Lemma~\ref{pure-submodule-is-injective-lemma}.
\end{proof}

\begin{cor} \label{quasi-injectivity-over-base-pure-subsheaf}
 Let $X$ be a locally Noetherian scheme and $\pi\:Y\rarrow X$ be
an affine morphism of schemes.
 Let\/ $\K$ be a quasi-coherent sheaf on $Y$, and let\/ $0\rarrow\cH
\rarrow\G\rarrow\F\rarrow0$ be a short exact sequence of quasi-coherent
sheaves on $Y$ such that\/ $\F$ is a flat quasi-coherent sheaf.
 Assume that the quasi-coherent sheaf $\pi_*(\G\ot_{\O_Y}\K)$ on $X$
is injective.
 Then the quasi-coherent sheaf $\pi_*(\cH\ot_{\O_Y}\K)$ on $X$ is
injective, too.
\end{cor}

\begin{proof}
 The argument from the proof of
Corollary~\ref{quasi-injectivity-over-base-tensor-product}
reduces the question to
Lemma~\ref{quasi-injectivity-over-base-pure-submodule-affine}.
\end{proof}

\begin{cor} \label{quasi-injectivity-over-base-locality}
 Let $X$ be a locally Noetherian scheme, $Y$ be a quasi-compact
semi-separated scheme, and $\pi\:Y\rarrow X$ be an affine morphism of
schemes.
 Let $Y=\bigcup_\alpha V_\alpha$ be an open covering of~$Y$.
 Denote by $j_\alpha\:V_\alpha\rarrow Y$ the open embedding morphisms,
and assume that the compositions $\pi j_\alpha\:V_\alpha\rarrow X$
are affine morphisms.
 Let\/ $\K$ be a quasi-coherent sheaf on~$Y$.
 Then the quasi-coherent sheaf $\pi_*\K$ on $X$ is injective if and
only if the quasi-coherent sheaves $\pi_*j_\alpha{}_*j_\alpha^*\K$
on~$X$ are injective for all~$\alpha$.
\end{cor}

\begin{proof}
 This is a global version of
Corollary~\ref{injectivity-over-noetherian-base-is-local}.
 The ``only if'' implication is provided by
Corollary~\ref{quasi-injectivity-over-base-ascent}.
 To prove the ``if'', recall that the morphisms~$j_\alpha$ are affine
by Lemma~\ref{affine-composition-impies-affine}.
 Without loss of generality we can assume that the set of
indices~$\alpha$ is finite.
 Then the \v Cech coresolution~\eqref{cech-quasi} of the quasi-coherent
sheaf $\O_Y$ on $Y$ with respect to the open covering
$Y=\bigcup_\alpha V_\alpha$ is a finite exact sequence of flat
quasi-coherent sheaves on~$Y$.
 Hence the injective morphism $\O_Y\rarrow\bigoplus_\alpha
j_\alpha{}_*j_\alpha^*\O_Y$ of quasi-coherent sheaves on $Y$ has
a flat cokernel.
 By Corollary~\ref{quasi-injectivity-over-base-pure-subsheaf},
injectivity of the quasi-coherent sheaves
$\pi_*(j_\alpha{}_*j_\alpha^*\O_Y\ot_{\O_Y}\K)$ on $X$ implies
injectivity of the quasi-coherent sheaf $\pi_*\K$, and it remains
to recall the projection formula isomorphisms
$j_\alpha{}_*j_\alpha^*\O_Y\ot_{\O_Y}\K\simeq j_\alpha{}_*j_\alpha^*\K$.
\end{proof}

 Let $X$ be a semi-separated Noetherian scheme, $Y$ be a quasi-compact
semi-sep\-a\-rated scheme, and $\pi\:Y\rarrow X$ be a morphism
of schemes.
 We will say that a quasi-coherent sheaf $\K$ on $Y$ is \emph{injective
over~$X$} (or \emph{$X$\+injective} for brevity) if it satisfies
the following condition.
 Let $Y=\bigcup_\alpha V_\alpha$ be an open covering of $Y$ such that
the compositions $V_\alpha\rarrow Y\rarrow X$ are affine morphisms
of schemes.
 For example, any affine open covering of $Y$ has this property.
 Denote by $j_\alpha\:V_\alpha\rarrow Y$ the open embedding morphisms.
 Then the condition is that, for every index~$\alpha$,
the quasi-coherent sheaf $\pi_*j_\alpha{}_*j_\alpha^*\K$ on $X$ must be
injective.
 It follows easily from
Corollary~\ref{quasi-injectivity-over-base-locality} that this condition
does not depend on the choice of an open covering
$Y=\bigcup_\alpha V_\alpha$.

 The full subcategory of $X$\+injective quasi-coherent sheaves on $Y$
will be denoted by $Y\qcoh^\Xinj\sub Y\qcoh$.
 Clearly, the full subcategory $Y\qcoh^\Xinj$ is closed under extensions,
cokernels of monomorphisms, and infinite direct sums in the abelian
category $Y\qcoh$; so it inherits an exact category structure.
 The result of~\cite[Lemma~A.21]{Psemten} tells us that if
the morphism~$\pi$ is flat, then any injective quasi-coherent sheaf on
$Y$ is $X$\+injective, i.~e., $Y\qcoh^\inj\sub Y\qcoh^\Xinj$.

 The following theorem is our version of~\cite[Theorems~5.1(a)
and~5.2(a)]{Pfp} and~\cite[Proposition~7.9]{Psemten}.
 Notice that our theorem is both less and more general than
the proposition in~\cite{Psemten}, in that a flat affine morphism of
ind-schemes was considered in~\cite{Psemten}, while we consider
a flat nonaffine morphism of schemes.

\begin{thm} \label{quasi-injective-over-base-semiderived}
 Let $X$ be a semi-separated Noetherian scheme, $Y$ be a quasi-compact
semi-separated scheme, and $\pi\:Y\rarrow X$ be a flat morphism of
schemes.
 Then the inclusion of exact/abelian categories $Y\qcoh^\Xinj\rarrow
Y\qcoh$ induces a triangulated equivalence
$$
 \sD(Y\qcoh^\Xinj)\simeq\sD^\si_X(Y\qcoh)
$$
between the conventional derived category of the exact category of
$X$\+injective quasi-coherent sheaves on $Y$ and the semiderived
category of quasi-coherent sheaves on $Y$ relative to~$X$.
\end{thm}

\begin{proof}
 By (the proof of) Theorem~\ref{quasi-coherent-becker-coderived},
for any complex of quasi-coherent sheaves $\M^\bu$ on $Y$ there exists
a complex of injective quasi-coherent sheaves $\J^\bu$ on $Y$ together
with a morphism of complexes $\M^\bu\rarrow\J^\bu$ with
a Becker-coacyclic cone.
 Now any injective quasi-coherent sheaf on $Y$ is $X$\+injective, while
any Becker-coacyclic complex of quasi-coherent sheaves on $Y$ is
semiacyclic.
 So Lemma~\ref{pkoszul-lemma16}(b) is applicable, and it remains to
show that a complex of $X$\+injective quasi-coherent sheaves on $Y$
is semiacyclic if and only if it is acyclic in the exact category
$Y\qcoh^\Xinj$.

 Indeed, let $V\sub Y$ be an open subscheme with the open embedding
morphism $j\:V\rarrow X$ such that the composition $\pi j\:V\rarrow X$
is an affine morphism.
 Let $\cA^\bu$ be an acyclic complex in $Y\qcoh^\Xinj$.
 Then $\pi_*j_*j^*\cA^\bu$ is an acyclic complex in the exact category
$X\qcoh^\inj$, i.~e., a contractible complex of injective
quasi-coherent sheaves on~$X$.
 Any contractible complex in $X\qcoh$ is obviously Becker-coacyclic;
so the complex $\cA^\bu$ is semiacyclic on $Y$ relative to~$X$.

 Conversely, let $\N^\bu$ be a complex in $Y\qcoh^\Xinj$ that is
semiacyclic over $X$ as a complex in $Y\qcoh$.
 Then $\pi_*j_*j^*\N^\bu$ is a complex in $X\qcoh^\inj$ that is
Becker-coacyclic in $X\qcoh$.
 Clearly, this means that $\pi_*j_*j^*\N^\bu$ is a contractible
complex in $X\qcoh^\inj$.
 In other words, the complex $\pi_*j_*j^*\N^\bu$ is acyclic in
$X\qcoh$ with injective sheaves of cocycles.
 As the functor $\pi_*j_*=(\pi j)_*\:V\qcoh\rarrow X\qcoh$ is
exact and faithful, it follows that the complex $j^*\N^\bu$ is
acyclic in $V\qcoh$ and its sheaves of cocycles are taken to
injective quasi-coherent sheaves on $X$ by the functor~$\pi_*j_*$.
 As this holds for any open subscheme $V\sub Y$ affine over $X$,
and such open subschemes form an open covering of $Y$, we can conclude
that the complex $\N^\bu$ is acyclic in $Y\qcoh$ with the sheaves
of cocycles belonging to $Y\qcoh^\Xinj$.
\end{proof}

 Let us introduce notation for the intersections of full subcategories
\begin{align*}
 Y\qcoh^{\cta,\Xinj} &= Y\qcoh^\cta\cap Y\qcoh^\Xinj, \\
 Y\qcoh^{\cot,\Xinj} &= Y\qcoh^\cot\cap Y\qcoh^\Xinj.
\end{align*}

\begin{lem} \label{quasi-X-injective-cta-coresol-dim}
 Let $X$ be a semi-separated Noetherian scheme, $Y$ be a quasi-compact
semi-separated scheme, and $\pi\:Y\rarrow X$ be a flat morphism
of schemes.
 Let $Y=\bigcup_{\alpha=1}^N V_\alpha$ be a finite affine open covering
of~$Y$.
 Then the coresolution dimension of any $X$\+injective quasi-coherent
sheaf on $Y$ with respect to the coresolving subcategory of
contraadjusted $X$\+injective quasi-coherent sheaves
$Y\qcoh^{\cta,\Xinj}\sub Y\qcoh^\Xinj$ does not exceed~$N$.
\end{lem}

\begin{proof}
 Any sheaf from $Y\qcoh^\Xinj$ can be embedded into an injective
quasi-coherent sheaf on $Y$, and all injective quasi-coherent sheaves
on $Y$ are contraadjusted and $X$\+injective
(by~\cite[Lemma~A.21]{Psemten}).
 Therefore, the assertions follow from
Lemma~\ref{dil-cta-clp-finite-dim}(b) by virtue of the dual
version of Corollary~\ref{fdim-subcategory-cor}.
\end{proof}

\begin{cor} \label{quasi-X-injective-cta-derived-equiv}
 Let $X$ be a semi-separated Noetherian scheme, $Y$ be a quasi-compact
semi-separated scheme, and $\pi\:Y\rarrow X$ be a flat morphism
of schemes.
 Then, for any symbol\/ $\bst=\b$, $+$, $-$, $\empt$, $\abs+$, $\abs-$,
$\bco$, or\/ $\abs$, the triangulated functor\/
$\sD^\st(Y\qcoh^{\cta,\Xinj})\rarrow\sD^\st(Y\qcoh^\Xinj)$ induced by
the embedding of exact categories $Y\qcoh^{\cta,\Xinj}\rarrow
Y\qcoh^\Xinj$ is an equivalence of triangulated categories
$$
 \sD^\st(Y\qcoh^{\cta,\Xinj})\simeq\sD^\st(Y\qcoh^\Xinj).
$$
\end{cor}

\begin{proof}
 Follows from Lemma~\ref{quasi-X-injective-cta-coresol-dim}
in view of the opposite versions of
Propositions~\ref{finite-resolutions}
and~\ref{becker-contraderived-finite-resolutions}.
\end{proof}

\begin{thm} \label{quasi-X-injective-cot-derived-equiv}
 Let $X$ be a semi-separated Noetherian scheme, $Y$ be a quasi-compact
semi-separated scheme, and $\pi\:Y\rarrow X$ be a flat morphism
of schemes.
 Then, for any symbol\/ $\bst=+$ or\/~$\empt$, the triangulated
functors\/ $\sD^\st(Y\qcoh^{\cot,\Xinj})\rarrow
\sD^\st(Y\qcoh^{\cta,\Xinj})\rarrow\sD^\st(Y\qcoh^\Xinj)$ induced by
the embeddings of exact categories $Y\qcoh^{\cot,\Xinj}\rarrow
Y\qcoh^{\cta,\Xinj}\rarrow Y\qcoh^\Xinj$ are equivalences of
triangulated categories
$$
 \sD^\st(Y\qcoh^{\cot,\Xinj})\simeq
 \sD^\st(Y\qcoh^{\cta,\Xinj})\simeq\sD^\st(Y\qcoh^\Xinj).
$$
\end{thm}

\begin{proof}
 Taking Corollary~\ref{quasi-X-injective-cta-derived-equiv} into
account, it suffices to show that the functor
$\sD^\st(Y\qcoh^{\cot,\Xinj})\rarrow\sD^\st(Y\qcoh^\Xinj)$ is
a triangulated equivalence.
 Notice first of all that the full subcategory $Y\qcoh^{\cot,\Xinj}$
is coresolving in $Y\qcoh^\Xinj$; this is provable similarly to
the proof of Lemma~\ref{quasi-X-injective-cta-coresol-dim}.
 So the assertion for $\bst=+$ follows by virtue of the dual version of
Proposition~\ref{infinite-resolutions}(a).

 Furthermore, any complex in $Y\qcoh^{\cot,\Xinj}$ that is acyclic in
$Y\qcoh^\Xinj$ is also acyclic in $Y\qcoh^{\cot,\Xinj}$ by the cotorsion
periodicity theorem for quasi-coherent sheaves on $Y$
(see~\cite[Section~10]{PS6}).
 In view of Lemma~\ref{pkoszul-lemma16}(b), in order to prove
the theorem for $\bst=\empt$ it remains to construct, for any complex
$\K^\bu$ in $Y\qcoh^\Xinj$, a complex $\cQ^\bu$ in $Y\qcoh^{\cot,\Xinj}$
together with a morphism of complexes $\K^\bu\rarrow\cQ^\bu$ with
a cone acyclic in $Y\qcoh^\Xinj$.

 By (the proof of) Theorem~\ref{quasi-coherent-becker-coderived},
for any complex $\M^\bu$ in $Y\qcoh$ there exists a complex of injective
quasi-coherent sheaves $\J^\bu$ on $Y$ together with a morphism of
complexes $\M^\bu\rarrow\J^\bu$ whose cone is Becker-coacyclic in
$Y\qcoh$.
 According to the discussion in
Section~\ref{antilocal-semicoderived-subsect}
and~\cite[Remark~A.19]{Psemten}, any Becker-coacyclic complex of
quasi-coherent sheaves on $Y$ is semiacyclic relative to~$X$
(essentially by Lemmas~\ref{qcoh-flat-inverse-image-becker-coacyclic}
and~\ref{qcoh-affine-direct-image-becker-coacyclic}).
 In the case of $\M^\bu=\K^\bu$, the cone of the morphism $\K^\bu
\rarrow\J^\bu$ is a complex of $X$\+injective quasi-coherent
sheaves on $Y$ (since $Y\qcoh^\inj\sub Y\qcoh^\Xinj$).
 By Theorem~\ref{quasi-injective-over-base-semiderived}, any semiacyclic
complex of $X$\+injective quasi-coherent sheaves on $Y$ is
acyclic in $Y\qcoh^\Xinj$.
 It remains to put $\cQ^\bu=\J^\bu$ and recall that injective
quasi-coherent sheaves on $Y$ are both cotorsion and $X$\+injective.
\end{proof}

\begin{rem} \label{sheaves-inj-over-non-noetherian-base-remark}
 The results of Theorem~\ref{quasi-injective-over-base-semiderived},
Lemma~\ref{quasi-X-injective-cta-coresol-dim},
Corollary~\ref{quasi-X-injective-cta-derived-equiv}, and
Theorem~\ref{quasi-X-injective-cot-derived-equiv} remain valid in
a different seting where the morphism~$\pi$ is affine, but
the scheme $X$ need not be Noetherian.
 Let $Y$ be a quasi-compact semi-separated scheme and $\pi\:Y\rarrow X$
be an affine morphism of schemes.
 In this context, one says that a quasi-coherent sheaf $\K$ on $Y$
is \emph{injective over~$X$} if the quasi-coherent sheaf $\pi_*\K$
on $X$ is injective.
 Similarly, one says that a complex of quasi-coherent sheaves $\N^\bu$
on $Y$ is \emph{semi}(\emph{co})\emph{acyclic} (over~$X$) if the complex
of quasi-coherent sheaves $\pi_*\N^\bu$ on $X$ is Becker-coacyclic.
 Keeping the rest of the terminology and notation as above,
the assertions of the lemma, corollary, and two theorems are true for
any quasi-compact semi-separated scheme $Y$ and a flat affine morphism
of schemes $\pi\:Y\rarrow X$.
 The proofs are the same, or actually even simpler.
\end{rem}

\subsection{Sheaves flat over the base}
\label{sheaves-flat-over-base-subsect}
 The notion of a quasi-coherent sheaf flat over a base scheme goes back
to~\cite[Section~2.1]{Groth3}; see also~\cite[Section Tag~01U2]{SP}.
 For an ind-scheme version of the discussion in this section (for
an affine morphism of ind-schemes), see~\cite[Section~7.2]{Psemten}.

\begin{lem} \label{quasi-flatness-over-base-tensor-product-affine}
 Let $\pi\:Y\rarrow X$ be a morphism of affine schemes, $\K$ be
a quasi-coherent sheaf on $Y$, and\/ $\F$ be a flat quasi-coherent
sheaf on~$Y$.
 Assume that the quasi-coherent sheaf $\pi_*\K$ on $X$ is flat.
 Then the quasi-coherent sheaf $\pi_*(\F\ot_{\O_Y}\K)$ on $X$ is flat,
too.
\end{lem}

\begin{proof}
 This is a restatement of
Lemma~\ref{flat-over-base-ring-tensor-product}.
\end{proof}

\begin{cor} \label{quasi-flatness-over-base-tensor-product}
 Let $\pi\:Y\rarrow X$ be an affine morphism of schemes, $\K$ be
a quasi-coherent sheaf on $Y$, and\/ $\F$ be a flat quasi-coherent sheaf
on~$Y$.
 Assume that the quasi-coherent sheaf $\pi_*\K$ on $X$ is flat.
 Then the quasi-coherent sheaf $\pi_*(\F\ot_{\O_Y}\K)$ on $X$ is flat,
too.
\end{cor}

\begin{proof}
 By the definition, flatness of quasi-coherent sheaves on schemes is
a local property; so the argument from the proof of
Corollary~\ref{quasi-injectivity-over-base-tensor-product} deduces
the desired assertion from
Lemma~\ref{quasi-flatness-over-base-tensor-product-affine}.
\end{proof}

\begin{cor} \label{quasi-flatness-over-base-ascent}
 Let $Y$ be a semi-separated scheme and $\pi\:Y\rarrow X$ be an affine
morphism of schemes.
 Let $V\sub Y$ be an open subscheme in $Y$ with the open embedding
morphism $j\:V\rarrow Y$ such that the composition $\pi j\:V\rarrow X$
is an affine morphism.
 Let\/ $\K$ be a quasi-coherent sheaf on $Y$ such that
the quasi-coherent sheaf $\pi_*\K$ on $X$ is flat.
 Then the quasi-coherent sheaf $\pi_*j_*j^*\K$ on $X$ is flat, too.
\end{cor}

\begin{proof}
 The argument from the proof of
Corollary~\ref{quasi-injectivity-over-base-ascent} deduces the desired
assertion from
Corollary~\ref{quasi-flatness-over-base-tensor-product}.
\end{proof}

\begin{cor} \label{quasi-flatness-over-base-locality}
 Let $Y$ be a quasi-compact semi-separated scheme and $\pi\:Y\rarrow X$
be an affine morphism of schemes.
 Let $Y=\bigcup_\alpha V_\alpha$ be an open covering of~$Y$.
 Denote by $j_\alpha\:V_\alpha\rarrow Y$ the open embedding morphisms,
and assume that the compositions $\pi j_\alpha\:V_\alpha\rarrow X$
are affine morphisms.
 Let\/ $\K$ be a quasi-coherent sheaf on~$Y$.
 Then the quasi-coherent sheaf $\pi_*\K$ on $X$ is flat if and only if
the quasi-coherent sheaves $\pi_*j_\alpha{}_*j_\alpha^*\K$ on~$X$ are
flat for all~$\alpha$.
\end{cor}

\begin{proof}
 This is a global version of
Lemma~\ref{flatness-over-base-is-local}.
 The ``only if'' implication is provided by
Corollary~\ref{quasi-flatness-over-base-ascent}.
 To prove the ``if'', we argue as in the proof of
Corollary~\ref{quasi-injectivity-over-base-locality}, considering
the \v Cech coresolution~\eqref{cech-quasi} of the quasi-coherent
sheaf $\K$ on $Y$ with respect to (a finite subcovering of)
the open covering $Y=\bigcup_\alpha V_\alpha$.
 The argument finishes by noticing that, whenever in a finite exact
sequence of quasi-coherent sheaves on $X$, all the sheaves except
possibly the leftmost one are flat, one can conclude that the leftmost
one is flat too.
 Indeed, the class of flat quasi-coherent sheaves on $X$ is closed
under kernels of surjective morphisms.
\end{proof}

 Let $\pi\:Y\rarrow X$ be a morphism of schemes and $\K$ be
a quasi-coherent sheaf on~$Y$.
 The simplest definition of what it means for $\K$ to be \emph{flat
over~$X$} is that, for every pair of affine open subschemes
$U\sub X$ and $V\sub Y$ such that $f(V)\sub U$, the $\O_Y(V)$\+module
$\K(V)$ must be flat \emph{as an\/ $\O_X(U)$\+module}.

 A slightly more fancy definition based on the lemma and corollaries
above presumes that $X$ is a quasi-compact semi-separated scheme.
 Then a quasi-coherent sheaf $\K$ on $Y$ is flat over $X$
(\emph{$X$\+flat} for brevity) if and only if it satisfies the following
condition.
 Let $Y=\bigcup_\alpha V_\alpha$ be an open covering of $Y$ such that
the compositions $V_\alpha\rarrow Y\rarrow X$ are affine morphisms of
schemes.
 Denote by $j_\alpha\:V_\alpha\rarrow Y$ the open embedding morphisms.
 Then the condition is that, for every index~$\alpha$,
the quasi-coherent sheaf $\pi_*j_\alpha{}_*j_\alpha^*\K$ on $X$ must be
flat.
 It follows easily from
Corollary~\ref{quasi-flatness-over-base-locality} that this condition
does not depend on the choice of an open covering
$Y=\bigcup_\alpha V_\alpha$.

 The full subcategory of $X$\+flat quasi-coherent sheaves on $Y$ will
be denoted by $Y\qcoh_\Xfl\sub Y\qcoh$.
 Clearly, the full subcategory $Y\qcoh_\Xfl$ is closed under extensions,
kernels of epimorphisms, and infinite direct sums in the abelian
category $Y\qcoh$; so it inherits an exact category structure.
 If the morphism~$\pi$ is flat, then any flat quasi-coherent sheaf
on $Y$ is $X$\+flat, i.~e., $Y\qcoh_\fl\sub Y\qcoh_\Xfl$.

 The notion of a quasi-coherent sheaf very flat over a base scheme was
defined in Section~\ref{very-flat-morphisms-subsect}, but it has
a difficulty assocated with it, as described in
Example~\ref{very-covering-mod-counterex} and pointed out in
a warning in Section~\ref{very-flat-morphisms-subsect}.
 Because of this problem, for the purposes of this chapter we will use
a weakened version of this notion, defined below for \emph{affine}
morphisms of schemes $\pi\:Y\rarrow X$ only.

 Let $\pi\:Y\rarrow X$ be an affine morphism of schemes.
 We will say that a quasi-coherent sheaf $\K$ on $Y$ is \emph{affinely
very flat over~$X$} (\emph{$X$\+very flat} for brevity) if
the quasi-coherent sheaf $\pi_*(\K)$ on $X$ is very flat.
 So any quasi-coherent sheaf on $Y$ very flat over $X$ is affinely
very flat over $X$, but the converse need not be true.
 The full subcategory of $X$\+very flat quasi-coherent sheaves on $Y$
will be denoted by $Y\qcoh_\Xvfl\sub Y\qcoh$.
 The full subcategory $Y\qcoh_\Xvfl$ is closed under extensions, kernels
of epimorphisms, and infinite direct sums in the abelian category
$Y\qcoh$; so it inherits an exact category structure.
 Obviously, we have $Y\qcoh_\Xvfl\sub Y\qcoh_\Xfl$.
 If the morphism~$\pi$ is very flat, then any very flat quasi-coherent
sheaf on $Y$ is $X$\+very flat, i.~e., $Y\qcoh_\vfl\sub Y\qcoh_\Xvfl$
(see Lemma~\ref{very-scalars-veryflat-case}(b) and
Section~\ref{very-flat-morphisms-subsect}).

 Let us introduce notation for the intersections of full subcategories
\begin{align*}
 Y\qcoh^\cta_\Xfl &= Y\qcoh^\cta\cap Y\qcoh_\Xfl, \\
 Y\qcoh^\cot_\Xfl &= Y\qcoh^\cot\cap Y\qcoh_\Xfl, \\
 Y\qcoh^\cta_\Xvfl &= Y\qcoh^\cta\cap Y\qcoh_\Xvfl.
\end{align*}

\begin{lem} \label{quasi-X-flat-cta-finite-coresol-dim}
 Let $\pi\:Y\rarrow X$ be a flat morphism of quasi-compact
semi-separated schemes.
 Let $Y=\bigcup_{\alpha=1}^N V_\alpha$ be an affine open covering
of~$Y$.
 In this context: \par
\textup{(a)} The coresolution dimension of any $X$\+flat quasi-coherent
sheaf on $Y$ with respect to the coresolving subcategory
$Y\qcoh_\Xfl^\cta\sub Y\qcoh_\Xfl$ does not exceed~$N$. \par
\textup{(b)} Assume that the morphism~$\pi$ is affine and very flat.
 Then the coresolution dimension of any $X$\+very flat quasi-coherent
sheaf on $Y$ with respect to the coresolving subcategory
$Y\qcoh_\Xvfl^\cta\sub Y\qcoh_\Xvfl$ does not exceed~$N$.
\end{lem}

\begin{proof}
 Part~(a): let $\K$ be an $X$\+flat quasi-coherent sheaf on~$Y$.
 By Corollary~\ref{quasi-very-cta-cor}(b) or~\ref{quasi-cotors-cor}(b),
there exists a short exact sequence of quasi-coherent sheaves
$0\rarrow\K\rarrow\cP\rarrow\F\rarrow0$ with a (very) flat
quasi-coherent sheaf $\F$ and a contraadjusted quasi-coherent
sheaf $\cP$ on~$Y$.
 Now both the quasi-coherent sheaves $\K$ and $\F$ are $X$\+flat,
hence so is the quasi-coherent sheaf~$\cP$.
 Thus any object of $Y\qcoh_\Xfl$ can be admissibly embedded into
a contraadjusted $X$\+flat quasi-coherent sheaf on~$Y$.
 It remains to refer to Lemma~\ref{dil-cta-clp-finite-dim}(b) and
the dual version of Corollary~\ref{fdim-subcategory-cor}.
 The proof of part~(b) is similar and based on
Corollary~\ref{quasi-very-cta-cor}(b).
\end{proof}

\begin{cor}  \label{quasi-X-flat-cta-derived-equiv}
 Let $\pi\:Y\rarrow X$ be a flat morphism of quasi-compact
semi-separated schemes.
 In this context: \par
\textup{(a)} For any symbol\/ $\bst=\b$, $+$, $-$, $\empt$, $\abs+$,
$\abs-$, or\/ $\abs$, the triangulated functor\/
$\sD^\st(Y\qcoh^\cta_\Xfl)\rarrow\sD^\st(Y\qcoh_\Xfl)$ induced by
the embedding of exact categories $Y\qcoh^\cta_\Xfl\rarrow
Y\qcoh_\Xfl$ is an equivalence of triangulated categories
$$
 \sD^\st(Y\qcoh^\cta_\Xfl)\simeq\sD^\st(Y\qcoh_\Xfl).
$$ \par
\textup{(b)} Assume that the morphism~$\pi$ is affine and very flat.
 Then, for any symbol\/ $\bst=\b$, $+$, $-$, $\empt$, $\abs+$,
$\abs-$, or\/ $\abs$, the triangulated functor\/
$\sD^\st(Y\qcoh^\cta_\Xvfl)\rarrow\sD^\st(Y\qcoh_\Xvfl)$ induced by
the embedding of exact categories $Y\qcoh^\cta_\Xvfl\rarrow
Y\qcoh_\Xvfl$ is an equivalence of triangulated categories\/
$$
 \sD^\st(Y\qcoh^\cta_\Xvfl)\simeq\sD^\st(Y\qcoh_\Xvfl).
$$
\end{cor}

\begin{proof}
 Follows from Lemma~\ref{quasi-X-flat-cta-finite-coresol-dim} in
view of Proposition~\ref{finite-resolutions}.
\end{proof}

\begin{thm}   \label{quasi-X-flat-cot-derived-equiv}
 Let $\pi\:Y\rarrow X$ be a flat morphism of quasi-compact
semi-separated schemes.
 Then, for any symbol\/ $\bst=+$ or\/~$\empt$, the triangulated
functors\/ $\sD^\st(Y\qcoh^\cot_\Xfl)\rarrow\sD^\st(Y\qcoh^\cta_\Xfl)
\rarrow\sD^\st(Y\qcoh_\Xfl)$ induced by the embeddings of exact
categories $Y\qcoh^\cot_\Xfl\rarrow Y\qcoh^\cta_\Xfl\rarrow
Y\qcoh_\Xfl$ are equivalences of triangulated categories
$$
 \sD^\st(Y\qcoh^\cot_\Xfl)\simeq
 \sD^\st(Y\qcoh^\cta_\Xfl)\simeq\sD^\st(Y\qcoh_\Xfl).
$$
\end{thm}

\begin{proof}
 In view of Corollary~\ref{quasi-X-flat-cta-derived-equiv}(a),
it suffices to show that the functor $\sD^\st(Y\qcoh^\cot_\Xfl)
\allowbreak\rarrow\sD^\st(Y\qcoh_\Xfl)$ is a triangulated equivalence.
 Notice first of all that the full subcategory $Y\qcoh^\cot_\Xfl$
is coresolving in $Y\qcoh_\Xfl$; this is provable similarly to
the proof of Lemma~\ref{quasi-X-flat-cta-finite-coresol-dim}(a)
(using Corollary~\ref{quasi-cotors-cor}(b)).
 So the assertion for $\bst=+$ follows by virtue of the dual version
of Proposition~\ref{infinite-resolutions}(a).
{\hfuzz=2pt\par}

 Furthermore, any complex in $Y\qcoh^\cot_\Xfl$ that is acyclic in
$Y\qcoh_\Xfl$ is also acyclic in $Y\qcoh^\cot_\Xfl$ by the cotorsion
periodicity theorem for quasi-coherent sheaves on~$Y$
(see~\cite[Theorem~10.2, Remark~10.3, and Corollary~10.4]{PS6}).
 In view of Lemma~\ref{pkoszul-lemma16}(b), in order to prove
the theorem for $\bst=\empt$ it remains to construct, for any
complex $\K^\bu$ in $Y\qcoh_\Xfl$, a complex $\cQ^\bu$ in
$Y\qcoh^\cot_\Xfl$ together with a morphism of complexes
$\K^\bu\rarrow\cQ^\bu$ with a cone acyclic in $Y\qcoh_\Xfl$.

 By Lemma~\ref{acyclic-of-vfl-flat-arbitrary-of-cta-cot-pairs}(b),
the pair of classes ($\Acycl(Y\qcoh_\fl)$, $\Com(Y\qcoh^\cot)$) is
a hereditary complete cotorsion pair in the abelian category
$\Com(Y\qcoh)$.
 By Lemmas~\ref{restricting-hereditary-cotorsion}
and~\ref{restricting-cotorsion-pairs-lemma}(a), this cotorsion
pair restricts to the exact subcategory $\Com(Y\qcoh_\Xfl)\sub
\Com(Y\qcoh)$ (since $Y\qcoh_\fl\sub Y\qcoh_\Xfl$).
 So we obtain a hereditary complete cotorsion pair
($\Acycl(Y\qcoh_\fl)$, $\Com(Y\qcoh^\cot_\Xfl)$) in the exact
category $\Com(Y\qcoh_\Xfl)$.
 Now given a complex $\K^\bu$ in $Y\qcoh_\Xfl$, consider a special
preenvelope short exact sequence $0\rarrow\K^\bu\rarrow\cQ^\bu\rarrow
\F^\bu\rarrow0$ with $\cQ^\bu\in\Com(Y\qcoh^\cot_\Xfl)$ and
$\F^\bu\in\Acycl(Y\qcoh_\fl)$.
 Clearly, the cone of the morphism $\K^\bu\rarrow\cQ^\bu$ is acyclic
in $Y\qcoh_\Xfl$.
\end{proof}

 The following lemma is a generalization of
Lemma~\ref{flat-cotors-homol-dim}(a)
and Corollary~\ref{derived-fl-vfl-cor}.

\begin{lem} \label{fl-vlf-over-base-lemma}
 Let $X$ be a semi-separated Noetherian scheme of finite Krull
dimension~$D$ and $\pi\:Y\rarrow X$ be a very flat affine morphism
of schemes.
 Then \par
\textup{(a)} the resolution dimension of any $X$\+flat quasi-coherent
sheaf on $Y$ with respect to the resolving subcategory
$Y\qcoh_\Xvfl\sub Y\qcoh$ does not exceed~$D$; \par
\textup{(b)} for any symbol\/ $\bst=\b$, $+$, $-$, $\empt$, $\abs+$,
$\abs-$, $\co$, or\/~$\abs$, the triangulated functor\/
$\sD^\st(Y\qcoh_\Xvfl)\rarrow\sD^\st(Y\qcoh_\Xfl)$ induced by
the inclusion of exact categories $Y\qcoh_\Xvfl\rarrow Y\qcoh_\Xfl$
is an equivalence of triangulated categories.
\end{lem}

\begin{proof}
 Part~(a): the full subcategory $Y\qcoh_\Xvfl$ is resolving in
$Y\qcoh$ (hence also in $Y\qcoh_\Xfl$) because there are enough
very flat quasi-coherent sheaves on $Y$
(by Lemma~\ref{quasi-very-flat-cover}), and all of them
are $X$\+very flat.
 Given an $X$\+flat quasi-coherent sheaf $\K$ on $Y$, consider
an exact sequence $0\rarrow\L\rarrow\G_{D-1}\rarrow\dotsb\rarrow
\G_0\rarrow\K\rarrow0$ of quasi-coherent sheaves on $Y$ with
$X$\+very flat quasi-coherent sheaves $\G_i$, \ $0\le i\le D-1$.
 Then $0\rarrow\pi_*\L\rarrow\pi_*\G_{D-1}\rarrow\dotsb\rarrow
\pi_*\G_0\rarrow\pi_*\K\rarrow0$ is an exact sequence of quasi-coherent
sheaves on $X$ with a flat quasi-coherent sheaf $\pi_*\K$ and
very flat quasi-coherent sheaves $\pi_*\G_i$.
 By Lemma~\ref{flat-cotors-homol-dim}(a) and
Corollary~\ref{fdim-cor}, it follows that the quasi-coherent sheaf
$\pi_*\L$ on $X$ is very flat; so $\L\in Y\qcoh_\Xvfl$.
 Part~(b) follows from part~(a) in view of
Proposition~\ref{finite-resolutions}.
\end{proof}

 In fact, the triangulated functor
$$
 \sD(Y\qcoh_\Xvfl)\lrarrow\sD(Y\qcoh_\Xfl)
$$
induced by the embedding of exact categories $Y\qcoh_\Xvfl\rarrow
Y\qcoh_\Xfl$ is a triangulated equivalence for any very flat affine
morphism of quasi-compact semi-separated schemes $\pi\:Y\rarrow X$.
 This is a part of Corollary~\ref{flat-side-over-base-equivalences}
below.

\subsection{Cosheaves projective over the base}
\label{cosheaves-projective-over-base-subsect}
 This section is a dual-analogous (contraherent cosheaf) version of
Section~\ref{sheaves-injective-over-base-subsect}.
 Concerning the question of (co)locality of projectivity over the base,
we prove the coascent, but cannot prove the codescent.
 So we only obtain some partial version of the desired dual-analogous
picture, with the main theorem proved for flat \emph{affine}
morphisms $\pi\:Y\rarrow X$.

 For the definitions of a coherent scheme and a flat cosheaf, see
Section~\ref{contraherent-tensor}.

\begin{lem} \label{contrah-flatness-over-base-Cohom-affine}
 Let $X$ be a coherent affine scheme, $Y$ be an affine scheme, and
$\pi\:Y\rarrow X$ be a morphism of schemes.
 Let\/ $\gK$ be a (globally) contraherent cosheaf on $Y$ such that
the contraherent cosheaf\/ $\pi_!\gK$ on $X$ is flat.
 Then \par
\textup{(a)} for any very flat quasi-coherent sheaf\/ $\F$ on $Y$,
the contraherent cosheaf $\pi_!\Cohom_Y(\F,\gK)$ on $X$ is flat; \par
\textup{(b)} if the contraherent cosheaf\/ $\gK$ on $Y$ is (locally)
cotorsion, then for any flat quasi-coherent sheaf\/ $\F$ on $Y$
the locally cotorsion contraherent cosheaf $\pi_!\Cohom_Y(\F,\gK)$
on $X$ is flat.
\end{lem}

\begin{proof}
 This is a restatement of
Lemma~\ref{commutative-coherent-cotorsper-lemma}.
\end{proof}

\begin{cor} \label{contrah-flatness-over-base-Cohom}
 Let $X$ be a locally coherent scheme and $\pi\:Y\rarrow X$ be
an affine morphism of schemes.
 Let\/ $\bW$ be an open covering of $X$ and\/ $\bT$ be an open covering
of $Y$ such that $\pi$~is a $(\bW,\bT)$\+affine morphism.
 Let\/ $\gK$ be a\/ $\bT$\+locally contraherent cosheaf on $Y$ such that
the\/ $\bW$\+locally contraherent cosheaf $\pi_!\gK$ on $X$ is\/
$\bW$\+flat.
 Then \par
\textup{(a)} for any very flat quasi-coherent sheaf\/ $\F$ on $Y$,
the\/ $\bW$\+locally contraherent cosheaf $\pi_!\Cohom_Y(\F,\gK)$ on $X$
is\/ $\bW$\+flat; \par
\textup{(b)} if the locally contraherent cosheaf\/ $\gK$ on $Y$ is
locally cotorsion, then for any flat quasi-coherent sheaf\/ $\F$ on $Y$
the\/ locally cotorsion $\bW$\+locally contraherent cosheaf
$\pi_!\Cohom_Y(\F,\gK)$ on $X$ is\/ $\bW$\+flat.
\end{cor}

\begin{proof}
 Let $U\sub X$ be an affine open subscheme subordinate to~$\bW$.
 Then $V=U\times_XY$ is an affine open subscheme in $Y$ subordinate
to~$\bT$.
 In order to prove that the $\O_X(U)$\+module $\pi_!\Cohom_Y(\F,\gK)[U]$
is flat, it suffices to apply (the respective part of)
Lemma~\ref{contrah-flatness-over-base-Cohom-affine} to
the morphism $\pi_U\:V\rarrow U$.
\end{proof}

\begin{cor} \label{contrah-flatness-over-base-ascent}
 Let $X$ be a locally coherent scheme, $Y$ be a semi-separated scheme,
and $\pi\:Y\rarrow X$ be an affine morphism of schemes.
 Let\/ $\bW$ be an open covering of $X$ and\/ $\bT$ be an open covering
of $Y$ such that $\pi$~is a $(\bW,\bT)$\+affine morphism.
 Let $V\sub Y$ be an open subscheme in $Y$ with the open
embedding morphism $j\:V\rarrow Y$ such that the composition
$\pi j\:V\rarrow X$ is an affine morphism.
 Let\/ $\gK$ be a\/ $\bT$\+locally contraherent cosheaf on $Y$ such that
the\/ $\bW$\+locally contraherent cosheaf $\pi_!\gK$ on $X$ is\/
$\bW$\+flat.
 Then the\/ $\bW$\+locally contraherent cosheaf $\pi_!j_!j^!\gK$ on $X$ 
is\/ $\bW$\+flat, too.
\end{cor}

\begin{proof}
 The open embedding morphism $j\:V\rarrow Y$ is affine by
Lemma~\ref{affine-composition-impies-affine}, hence $j$~is also
$(\bT,\bT|_V)$\+affine and the morphism $\pi j\:V\rarrow X$ is
$(\bW,\bT|_V)$\+affine.
 Besides, the morphism~$j$ is $(\bT,\bT|_V)$\+coaffine.
 Therefore, the cosheaf $\pi_!j_!j^!\gK$ on $X$ is $\bW$\+locally
contraherent.
 According to the projection
formula~\eqref{flat-cotors-projection-cohom}, we have natural
isomorphisms $j_!j^!\gK\simeq j_!\Cohom_V(\O_V,j^!\gK)\simeq
\Cohom_Y(j_*\O_V,\gK)$ of $\bT$\+locally contraherent cosheaves on~$Y$.
 So it remains to apply
Corollary~\ref{contrah-flatness-over-base-Cohom}(a) to
the $\bT$\+locally contraherent cosheaf $\gK$ and the very flat
quasi-coherent sheaf $\F=j_*\O_V$ on $Y$ in order to show that
the $\bW$\+locally contraherent cosheaf $\pi_!j_!j^!\gK$ on $X$
is $\bW$\+flat.
\end{proof}

 We do \emph{not} know how to prove a dual-analogous version of
Lemma~\ref{quasi-injectivity-over-base-pure-submodule-affine} in
the contraherent cosheaf context.
 See Question~\ref{codescent-of-flatness-question} for a discussion.
 Being unable to resolve this problem, we restrict ourselves to
\emph{affine} morphisms $\pi\:Y\rarrow X$ in the exposition below
in this section.

 Let $X$ be a locally Noetherian scheme and $\pi\:Y\rarrow X$ be
an affine morphism of schemes.
 Let $\bW$ be an open covering of $X$ and $\bT$ be an open covering of
$Y$ such that the morphism~$\pi$ is $(\bW,\bT)$\+affine.
 We will say that a locally cotorsion $\bT$\+locally contraherent
cosheaf $\P$ on $Y$ is \emph{projective over~$X$} (or
\emph{$X$\+projective} for brevity) if the locally cotorsion
$\bW$\+locally contraherent cosheaf $\pi_!\P$ on $X$ is projective.
 Notice that any projective locally cotorsion $\bW$\+locally contraherent
cosheaf on $X$ is (globally) contraherent by
Theorem~\ref{proj-lct-classification}(a), and the class of such
cosheaves on $X$ does not depend on~$\bW$; so the property of $\P$
to be projective over~$X$ remains unchanged when the coverings $\bW$
and $\bT$ are refined.

 The full subcategory of $X$\+projective locally cotorsion
$\bT$\+locally contraherent cosheaves on $Y$ will be denoted by
$Y\lcth_{\bT,\Xprj}^\lct\sub Y\lcth_\bT^\lct$.
 Clearly, the full subcategory $Y\lcth_{\bT,\Xprj}^\lct$ is closed
under extensions and kernels of admissible epimorphisms in the exact
category $Y\lcth_\bT^\lct$; so it inherits an exact category structure.
 In view of Corollary~\ref{loc-noetherian-proj-products}, the full
subcategory $Y\lcth_{\bT,\Xprj}^\lct$ is also closed under infinite
products in $Y\lcth_\bT^\lct$.
 We put $Y\ctrh_\Xprj^\lct=Y\lcth_{\{Y\},\Xprj}^\lct$.
 Assuming additionally that the scheme $Y$ is quasi-compact and
semi-separated, and that the morphism~$f$ is flat, the functor~$\pi_!$
takes projective locally cotorsion contraherent cosheaves on $Y$ to
projective locally cotorsion contraherent cosheaves on~$X$ (cf.\
the proof of Corollary~\ref{proj-direct-inverse}(b)); so one has
$Y\ctrh^\lct_\prj\sub Y\ctrh_\Xprj^\lct$.

 Similarly, let $X$ be a Noetherian scheme that is either semi-separated
or has finite Krull dimension.
 We will say that a $\bT$\+locally contraherent cosheaf $\P$
on $Y$ is \emph{projective over~$X$} (\emph{$X$\+projective}) if
the $\bW$\+locally contraherent cosheaf $\pi_!\P$ on $X$ is projective.
 Once again, any projective $\bW$\+locally contraherent cosheaf on $X$
is (globally) contraherent by Corollary~\ref{ctrh-lcth-proj}
or~\ref{finite-krull-contrah-projective}(a), and the class of
such cosheaves on $X$ does not depend on~$\bW$; so the property of $\P$
to be projective over $X$ does not depend on the coverings $\bW$
and~$\bT$.
 Notice that an $X$\+projective locally cotorsion locally contraherent
cosheaf is usually \emph{not} $X$\+projective as a (locally 
contraadjusted) locally contraherent cosheaf.

 The full subcategory of $X$\+projective $\bT$\+locally contraherent
cosheaves on $Y$ will be denoted by $Y\lcth_{\bT,\Xprj}\sub
Y\lcth_\bT$.
 Clearly, the full subcategory $Y\lcth_{\bT,\Xprj}$ is closed under
extensions and kernels of admissible epimorphisms in the exact category
$Y\lcth_\bT$; so it inherits an exact category structure.
 We put $Y\ctrh_\Xprj=Y\lcth_{\{Y\},\Xprj}$.
 Assuming additionally that the scheme $Y$ is (quasi-compact and)
semi-separated and that the morphism~$f$ is very flat,
the functor~$\pi_!$ takes projective contraherent cosheaves on $Y$ to
projective contraherent cosheaves on~$X$ (cf.\ the proof of
Corollary~\ref{proj-direct-inverse}(a)); so one has
$Y\ctrh_\prj\sub Y\ctrh_\Xprj$.

 Finally, we will say that a $\bT$\+locally contraherent cosheaf $\gF$
on $Y$ is \emph{flat over~$X$} (\emph{$X$\+flat}) if the $\bW$\+locally
contraherent cosheaf $\pi_!\gF$ on $X$ is $\bW$\+flat.
 By Corollary~\ref{lct-prj-flat}, over a locally Noetherian base scheme
$X$, a locally cotorsion $\bT$\+locally contraherent cosheaf on $Y$ is
$X$\+projective (as a locally cotorsion locally contraherent cosheaf)
if and only if it is $X$\+flat.
 Assuming that the base scheme $X$ is Noetherian and either
semi-separated or of finite Krull dimension,
Corollary~\ref{proj-flat}(a) or~\ref{finite-krull-contrah-projective}(b)
tells us that any $X$\+projective $\bT$\+locally contraherent cosheaf
is $X$\+flat.
 By Corollary~\ref{finite-krull-flat-contraherent}(b), when
the scheme $X$ is Noetherian of finite Krull dimension, the property of
a $\bT$\+locally contraherent cosheaf $\gF$ to be $X$\+flat does not
depend on the coverings $\bW$ and~$\bT$.

 The full subcategory of $X$\+flat $\bT$\+locally contraherent cosheaves
on $Y$ will be denoted by $Y\lcth_\bT^\Xfl\sub Y\lcth_\bT$.
 The full subcategory $Y\lcth_\bT^\Xfl$ is closed under extensions,
kernels of admissible epimorphisms, and infinite products in the exact
category $Y\lcth_\bT$; so it inherits an exact category structure.
 We put $Y\ctrh^\Xfl=Y\lcth_{\{Y\}}^\Xfl$.
 Assuming that the scheme $X$ is semi-separated and Noetherian and that
the morphism~$f$ is flat, the functor~$\pi_!$ takes antilocally flat
contraherent cosheaves on $Y$ to flat contraherent cosheaves on~$X$
(see Section~\ref{clf-subsection}); so one has
$Y\ctrh_\alf\sub Y\ctrh^\Xfl$.

\begin{thm} \label{contrah-projective-flat-over-base-semiderived}
 Let $X$ be a semi-separated Noetherian scheme and $\pi\:Y\rarrow X$
be a flat affine morphism of schemes.
 Let\/ $\bW$ be an open covering of $X$ and\/ $\bT$ be an open covering
of $Y$ such that the morphism~$\pi$ is $(\bW,\bT)$\+affine.
 In this context: \par
\textup{(a)} The inclusion of exact categories $Y\lcth_{\bT,\Xprj}^\lct
\rarrow Y\lcth_\bT^\lct$ induces a triangulated equivalence
$$
 \sD(Y\lcth_{\bT,\Xprj}^\lct)\simeq\sD^\si_X(Y\lcth_\bT^\lct)
$$
between the conventional derived category of the exact category of
$X$\+projective locally cotorsion\/ $\bT$\+locally contraherent
cosheaves on $Y$ and the semiderived category of locally cotorsion\/
$\bT$\+locally contraherent cosheaves on $Y$ relative to~$X$. \par
\textup{(b)} Assume that the scheme $X$ has finite Krull dimension.
 Then the inclusion of exact categories $Y\lcth_\bT^\Xfl\rarrow
Y\lcth_\bT$ induces a triangulated equivalence
$$
 \sD(Y\lcth_\bT^\Xfl)\simeq\sD^\si_X(Y\lcth_\bT)
$$
between the conventional derived category of the exact category of
$X$\+flat\/ $\bT$\+locally contraherent cosheaves on $Y$ and
the semiderived category of\/ $\bT$\+locally contraherent cosheaves
on $Y$ relative to~$X$. \par
\textup{(c)} Assume that the morphism~$f$ is very flat.
 Then the inclusion of exact categories $Y\lcth_{\bT,\Xprj}\rarrow
Y\lcth_\bT$ induces a triangulated equivalence
$$
 \sD(Y\lcth_{\bT,\Xprj})\simeq\sD^\si_X(Y\lcth_\bT)
$$
between the conventional derived category of the exact category of
$X$\+projective\/ $\bT$\+locally contraherent cosheaves on $Y$ and
the semiderived category of\/ $\bT$\+locally contraherent cosheaves
on $Y$ relative to~$X$. 
\end{thm}

\begin{proof}
 Part~(a): by (the proof of)
Corollary~\ref{becker-contraderived-of-lcta-lct-well-behaved}(b),
for any complex of locally cotorsion $\bT$\+locally contraherent
cosheaves $\gM^\bu$ on $Y$ there exists a complex of projective
locally cotorsion contraherent cosheaves $\P^\bu$ on $Y$ together
with a morphism of complexes $\P^\bu\rarrow\gM^\bu$ with a cone
Becker-contraacyclic in $Y\lcth_\bT^\lct$.
 Now any projective locally cotorsion contraherent cosheaf on $Y$
is $X$\+projective, while any Becker-contraacyclic complex in
$Y\lcth_\bT^\lct$ is semiacyclic.
 So Lemma~\ref{pkoszul-lemma16}(a) is applicable, and it remains to
show that a complex of $X$\+projective locally cotorsion
$\bT$\+locally contraherent cosheaves on $Y$ is semiacyclic if and
only if it is acyclic in the exact category $Y\lcth_{\bT,\Xprj}^\lct$.

 Indeed, let $\gA^\bu$ be an acyclic complex in
$Y\lcth_{\bT,\Xprj}^\lct$.
 Then $\pi_!\gA^\bu$ is an acyclic complex in the exact category
$X\ctrh^\lct_\prj$, i.~e., a contractible complex of projective
locally cotorsion contraherent cosheaves on~$X$.
 Any contractible complex in $X\lcth_\bW^\lct$ is obviously
Becker-contraacyclic; so the complex $\gA^\bu$ is semiacyclic
on $Y$ relative to~$X$.

 Conversely, let $\Q^\bu$ be a complex in $Y\lcth_{\bT,\Xprj}^\lct$
that is semiacyclic over $X$ as a complex in $Y\lcth_\bT^\lct$.
 Then $\pi_!\Q^\bu$ is a complex in $X\ctrh^\lct_\prj$ that is
Becker-contraacyclic in $X\lcth_\bW^\lct$.
 Clearly, this means that $\pi_!\Q^\bu$ is a contractible complex
in $X\ctrh^\lct_\prj$.
 In other words, the complex $\pi_!\Q^\bu$ is acyclic in
$X\lcth_\bW^\lct$ with projective (locally cotorsion contraherent)
cosheaves of cocycles.
 In view of Lemma~\ref{acyclicity-in-lcth-criterion}(b), it follows
that the complex $\Q^\bu$ is acyclic in $Y\lcth_\bT^\lct$ and its
cosheaves of cocycles are taken to projective locally cotorsion
contraherent cosheaves on $X$ by the functor~$\pi_!$.
 So the complex $\Q^\bu$ is acyclic in $Y\lcth_{\bT,\Xprj}^\lct$.

 Part~(b): by (the proof of)
Corollary~\ref{becker-contraderived-of-lcta-lct-well-behaved}(a),
for any complex of $\bT$\+locally contraherent cosheaves $\gM^\bu$
on $Y$ there exists a complex of projective contraherent cosheaves
$\P^\bu$ on $Y$ together with a morphism of complexes $\P^\bu\rarrow
\gM^\bu$ with a cone Becker-contraacyclic in $Y\lcth_\bT$.
 Now any projective contraherent cosheaf on $Y$ is antilocally flat,
and any antilocally flat contraherent cosheaf on $Y$ is $X$\+flat;
while any Becker-contraacyclic complex in $Y\lcth_\bT$ is semiacyclic.
 So Lemma~\ref{pkoszul-lemma16}(a) is applicable, and it remains to
show that a complex of $X$\+flat $\bT$\+locally contraherent cosheaves
on $Y$ is semiacyclic if and only if it is acyclic in $Y\lcth_\bT^\Xfl$.

 Indeed, let $\gA^\bu$ be an acyclic complex in $Y\lcth_\bT^\Xfl$.
 Then $\pi_!\gA^\bu$ is an acyclic complex in the exact category
$X\ctrh^\fl$ (because the direct image functor $\pi_!\:Y\lcth_\bT^\Xfl
\rarrow X\ctrh^\fl$ is exact).
 By Corollary~\ref{finite-krull-derived-equivalences}(a), any
acyclic complex in $X\ctrh^\fl$ is Becker-contraacyclic in
$X\ctrh^\fl$, or equivalently, in $X\ctrh$.
 So the complex $\gA^\bu$ is semiacyclic on $Y$ relative to~$X$.

 Conversely, let $\Q^\bu$ be a complex in $Y\lcth_\bT^\Xfl$ that is
semiacyclic as a complex in $Y\lcth_\bT$.
 Then $\pi_!\Q^\bu$ is a complex in $X\ctrh^\fl$ that is
Becker-contraacyclic in $X\lcth_\bW$.
 Once again, this means that $\pi_!\Q^\bu$ is Becker-contraacyclic in
$X\ctrh^\fl$ (see Corollary~\ref{finite-krull-contrah-projective}(b)),
and by Corollary~\ref{finite-krull-derived-equivalences}(a) it follows
that $\pi_!\Q^\bu$ is acyclic in $X\ctrh^\fl$.
 Using Lemma~\ref{acyclicity-in-lcth-criterion}(b), we can conclude
that the complex $\Q^\bu$ is acyclic in $Y\lcth_\bT^\Xfl$.

 The proof of part~(c) is similar to that of part~(a).
\end{proof}

 Let us introduce notation for the intersections of full subcategories
\begin{align*}
 Y\ctrh_{\al,\Xprj}^\lct &= Y\ctrh_\al^\lct\cap Y\lcth_{\bT,\Xprj}^\lct,
 \\
 Y\ctrh_{\al,\Xprj} &= Y\ctrh_\al\cap Y\lcth_{\bT,\Xprj}, \\
 Y\ctrh_\al^\Xfl &= Y\ctrh_\al\cap Y\lcth_\bT^\Xfl.
\end{align*}

\begin{lem} \label{contrah-X-proj-fl-finite-resol-dim}
 Let $X$ be a Noetherian scheme, $Y$ be a semi-separated scheme, and
$\pi\:Y\rarrow X$ be a flat affine morphism of schemes.
 Let\/ $\bW$ be an open covering of $X$ and\/ $\bT$ be an open covering
of $Y$ such that the morphism~$\pi$ is $(\bW,\bT)$\+affine.
 Let $Y=\bigcup_{\alpha=1}^N V_\alpha$ be a finite affine open covering
of $Y$ subordinate to\/~$\bT$.
 In this context: \par
\textup{(a)} The resolution dimension of any $X$\+projective locally
cotorsion\/ $\bT$\+locally contraherent cosheaf on $X$ with respect
to the resolving subcategory of antilocal $X$\+projective locally
cotorsion contraherent cosheaves $Y\ctrh_{\al,\Xprj}^\lct\sub
Y\lcth_{\bT,\Xprj}^\lct$ does not exceed $N-1$.
 Consequently, the same bound holds for the resolution dimension of
any object of $Y\lcth_{\bT,\Xprj}^\lct$ with respect to the resolving
subcategory $Y\ctrh_\Xprj^\lct$. \par
\textup{(b)} Assume that the scheme $X$ has finite Krull dimension.
 Then the resolution dimension of any $X$\+flat\/ $\bT$\+locally
contraherent cosheaf on $X$ with respect to the resolving subcategory
of antilocal $X$\+flat contraherent cosheaves
$Y\ctrh_\al^\Xfl\sub Y\lcth_\bT^\Xfl$ does not exceed $N-1$.
 Consequently, the same bound holds for the resolution dimension of
any object of $Y\lcth_\bT^\Xfl$ with respect to the resolving
subcategory $Y\ctrh^\Xfl$. \par
\textup{(c)} Assume that the morphism~$f$ is very flat.
 Then the resolution dimension of any $X$\+projective\/ $\bT$\+locally
contraherent cosheaf on $X$ with respect to the resolving subcategory
of antilocal $X$\+projective contraherent cosheaves
$Y\ctrh_{\al,\Xprj}\sub Y\lcth_{\bT,\Xprj}$ does not exceed $N-1$.
 Consequently, the same bound holds for the resolution dimension of
any object of $Y\lcth_{\bT,\Xprj}$ with respect to the resolving
subcategory $Y\ctrh_\Xprj$.
\end{lem}

\begin{proof}
 Part~(a): any object of $Y\lcth_{\bT,\Xprj}^\lct$ is the target of
an admissible epimorphism from a projective locally cotorsion
contraherent cosheaf on $Y$, and all projective locally cotorsion
contraherent cosheaves on $Y$ belong to $Y\ctrh_{\al,\Xprj}^\lct$.
 Therefore, all the assertions follow from
Lemma~\ref{lct-lin-clp-finite-dim}(a) by virtue of
Corollary~\ref{fdim-subcategory-cor}.

 The proofs of parts~(b\+c) are similar and based on
Lemma~\ref{dil-cta-clp-finite-dim}(c).
 In part~(b), one needs to observe that any object of $Y\lcth_\bT^\Xfl$
is the target of an admissible epimorphism from a projective
contraherent cosheaf on $Y$, and all projective contraherent cosheaves
on $Y$ belong to $Y\ctrh_\al^\Xfl$.
 In part~(c), any object of $Y\lcth_{\bT,\Xprj}$ is the target of
an admissible epimorphism from a projective contraherent cosheaf
on $Y$, and all such cosheaves belong to $Y\ctrh_{\al,\Xprj}$.
\end{proof}

\begin{cor} \label{contrah-lcth-projective-flat-over-base}
 Let $X$ be a semi-separated Noetherian scheme and $\pi\:Y\rarrow X$
be a flat affine morphism of schemes.
 Let\/ $\bW$ be an open covering of $X$ and\/ $\bT$ be an open covering
of $Y$ such that the morphism~$\pi$ is $(\bW,\bT)$\+affine.
 In this context: \par
\textup{(a)} For any symbol\/ $\bst=\b$, $+$, $-$, $\empt$, $\abs+$,
$\abs-$, $\bctr$, $\ctr$, or\/ $\abs$, the inclusions of exact
categories $Y\ctrh_{\al,\Xprj}^\lct\rarrow Y\ctrh_\Xprj^\lct
\rarrow Y\lcth_{\bT,\Xprj}^\lct$ induce equivalences of
the derived categories
$$
 \sD^\st(Y\ctrh_{\al,\Xprj}^\lct)\simeq
 \sD^\st(Y\ctrh_\Xprj^\lct)\simeq\sD^\st(Y\lcth_{\bT,\Xprj}^\lct).
$$ \par
\textup{(b)} Assume that the scheme $X$ has finite Krull dimension.
 Then, for any symbol\/ $\bst=\b$, $+$, $-$, $\empt$, $\abs+$,
$\abs-$, $\bctr$, $\ctr$, or\/ $\abs$, the inclusions of exact
categories $Y\ctrh_\al^\Xfl\rarrow Y\ctrh^\Xfl\rarrow Y\lcth_\bT^\Xfl$
induce equivalences of the derived categories
$$
 \sD^\st(Y\ctrh_\al^\Xfl)\simeq
 \sD^\st(Y\ctrh^\Xfl)\simeq\sD^\st(Y\lcth_\bT^\Xfl).
$$ \par
\textup{(c)} Assume that the morphism~$f$ is very flat.
 Then, for any symbol\/ $\bst=\b$, $+$, $-$, $\empt$, $\abs+$,
$\abs-$, $\bctr$, or\/ $\abs$, the inclusions of exact
categories $Y\ctrh_{\al,\Xprj}\rarrow Y\ctrh_\Xprj\rarrow
Y\lcth_{\bT,\Xprj}$ induce equivalences of the derived categories
$$
 \sD^\st(Y\ctrh_{\al,\Xprj})\simeq
 \sD^\st(Y\ctrh_\Xprj)\simeq\sD^\st(Y\lcth_{\bT,\Xprj}).
$$
\end{cor}

\begin{proof}
 All the assertions follow from the respective assertions of
Lemma~\ref{contrah-X-proj-fl-finite-resol-dim} in view of
Propositions~\ref{finite-resolutions}
and~\ref{becker-contraderived-finite-resolutions}.
 Alternatively, in each part~(a\+c), the second equivalence for
$\bst=\empt$ can be obtained as follows.
 Part~(a): compare
Theorem~\ref{contrah-projective-flat-over-base-semiderived}(a)
with Corollary~\ref{semicontraderived-independent-of-covering}(b).
 Part~(b): compare
Theorem~\ref{contrah-projective-flat-over-base-semiderived}(b)
with Corollary~\ref{semicontraderived-independent-of-covering}(a).
 Part~(c): compare
Theorem~\ref{contrah-projective-flat-over-base-semiderived}(c)
with Corollary~\ref{semicontraderived-independent-of-covering}(a).
\end{proof}

\begin{cor} \label{lcth-lcta-lct-projective-flat-over-base}
 Let $X$ be a semi-separated Noetherian scheme of finite Krull
dimension and $\pi\:Y\rarrow X$ be a flat affine morphism of schemes.
 Let\/ $\bW$ be an open covering of $X$ and\/ $\bT$ be an open covering
of $Y$ such that the morphism~$\pi$ is $(\bW,\bT)$\+affine.
 In this context: \par
\textup{(a)} The inclusion of exact categories
$Y\lcth_{\bT,\Xprj}^\lct\rarrow Y\lcth_\bT^\Xfl$ induces an equivalence
of the conventional unbounded derived categories
$$
 \sD(Y\lcth_{\bT,\Xprj}^\lct)\simeq\sD(Y\lcth_\bT^\Xfl).
$$ \par
\textup{(b)} If the morphism~$f$ is very flat, then the inclusion of
exact categories $Y\lcth_{\bT,\Xprj}\allowbreak\rarrow Y\lcth_\bT^\Xfl$
induces an equivalence of the conventional unbounded derived categories
$$
 \sD(Y\lcth_{\bT,\Xprj})\simeq\sD(Y\lcth_\bT^\Xfl).
$$
\end{cor}

\begin{proof}
 To prove part~(a), compare
Theorem~\ref{contrah-projective-flat-over-base-semiderived}(a\+b)
with Theorem~\ref{semicontraderived-lcta-lct-equivalent}.
 Part~(b) follows from
Theorem~\ref{contrah-projective-flat-over-base-semiderived}(b\+c).
\end{proof}

\begin{rem} \label{cosh-proj-over-qcomp-ssep-base-remark}
 We assume (local) Noetherianity of the scheme $X$ in the main
definitions and results of  this section chiefly because the Noetherian
case is relevant for the purposes of our theory of semico-semicontra
correspondence.
 Otherwise, the assumptions of (local) Noetherianity and finite
Krull dimension can be dropped from the above definitions of
$X$\+projective locally contraherent cosheaves on~$Y$; then one would
have to assume $X$ to be quasi-compact and semi-separated from
the outset.
 In this context, the notion of an $X$\+flat locally contraherent
cosheaf $\gF$ would have to be replaced by that of
an \emph{$X$\+antilocally flat} $\bT$\+locally contraherent cosheaf,
defined by the condition that the $\bW$\+locally contraherent cosheaf
$\pi_!\gF$ on $X$ is antilocally flat.
 All the assertions of
Theorem~\ref{contrah-projective-flat-over-base-semiderived},
Lemma~\ref{contrah-X-proj-fl-finite-resol-dim}, and
Corollaries~\ref{contrah-lcth-projective-flat-over-base}
and~\ref{lcth-lcta-lct-projective-flat-over-base} remain valid with
such a change of context.
 The only exceptions are the assertions involving
the Positselski contraderived categories ($\bst=\ctr$) in
Corollary~\ref{contrah-lcth-projective-flat-over-base}(a\+b), which
require the scheme $X$ to be coherent.
 In particular, the references to
Corollary~\ref{finite-krull-derived-equivalences}(a) in the proof of
Theorem~\ref{contrah-projective-flat-over-base-semiderived}(b) would
have to be replaced by references to
Corollary~\ref{acycl=bctracycl=bcoacycl-in-alf} for antilocally
flat contraherent cosheaves on~$X$.
\end{rem}

\subsection{Cosheaves locally injective over the base}
\label{cosheaves-loc-injective-over-base-subsect}
 This section is a partial dual-analogous (contraherent cosheaf)
version of Section~\ref{sheaves-flat-over-base-subsect}.
 We start with a discussion of locally $X$\+contrainjective locally
contraherent cosheaves on $Y$ for a nonaffine morphism of schemes
$\pi\:Y\rarrow X$ (based on Section~\ref{contrainjective-subsect})
before passing to $X$\+locally injective locally contraherent cosheaves
on $Y$ for an affine morphism~$\pi$.

 Let $\pi\:Y\rarrow X$ be a morphism of schemes.
 Let $\bW$ be an open covering of $X$ and $\bT$ be an open covering
of~$Y$.
 We will say that a $\bT$\+locally contraherent cosheaf $\gK$ on $Y$
is \emph{locally contrainjective over~$X$} (or \emph{locally
$X$\+contrainjective} for brevity) if, for every affine open subscheme
$U\sub X$ subordinate to $\bW$ and every affine open subscheme
$V\sub Y$ subordinate to $\bT$ such that $f(V)\sub U$,
the $\O_Y(V)$\+module $\gK[V]$ is injective as an $\O_X(U)$\+module.
 In view of the criterion of Lemma~\ref{contrainjective-lemma}(3),
a $\bT$\+locally contraherent cosheaf $\gK$ on $Y$ is
locally $X$\+contrainjective if and only if, for every $U$ and $V$
as above, the $\O_Y(V)$\+module $\gK[V]$ is
$(\O_Y(V),\O_X(U))$-contrainjective (in the sense of the definition
in Section~\ref{contrainjective-subsect}).

 Lemma~\ref{cotors-inj-covering}(b) tells us that the local
$X$\+contrainjectivity condition does not change when the open covering
$\bW$ of the scheme $X$ is replaced by its refinement.
 By Lemma~\ref{contra-cotors-injective-covering-mod}(a), the condition
also remains unchanged when the open covering $\bT$ of the scheme $Y$
is replaced by its refinement.
 When the morphism~$f$ is $(\bW,\bT)$\+affine, a $\bT$\+locally
contraherent cosheaf $\gK$ on $Y$ is locally $X$\+contrainjective
if and only if, for every open subscheme $U\sub X$ subordinate to $\bW$,
the $\O_Y(f^{-1}(U))$\+module $\gK[f^{-1}(U)]$ is
$(\O_Y(f^{-1}(U)),\O_X(U))$-contrainjective.
 Notice that injectivity of the $\O_X(U)$\+module $\gK[f^{-1}(U)]$ is
\emph{not} sufficient in the latter condition (by
Example~\ref{not-contrainjective-example}).

 The full subcategory of locally $X$\+contrainjective $\bT$\+locally
contraherent cosheaves on $Y$ will be denoted by
$Y\lcth_\bT^\Xlctin\sub Y\lcth_\bT$.
 The full subcategory $Y\lcth_\bT^\Xlctin$ is closed under extensions,
cokernels of admissible monomorphisms, and infinite products in
the exact category $Y\lcth_\bT$; so it inherits an exact category
structure.
 If the morphism~$f$ is flat, then any locally injective $\bT$\+locally
contraherent cosheaf on $Y$ is locally $X$\+contrainjective; so
$X\lcth_\bT^\lin\sub Y\lcth_\bT^\Xlctin$.

 Let $\pi\:Y\rarrow X$ be an affine morphism of schemes.
 Let $\bW$ be an open covering of $X$ and $\bT$ be an open covering of
$Y$ such that the morphism~$\pi$ is $(\bW,\bT)$\+affine.
 We will say that a $\bT$\+locally contraherent cosheaf $\gK$ on $Y$
is (\emph{affinely}) \emph{locally injective over~$X$}
(or \emph{$X$\+locally injective} for brevity) if the $\bW$\+locally
contraherent cosheaf $\pi_!\gK$ on $X$ is locally injective.
 The property of a $\bT$\+locally contraherent cosheaf $\gK$ to be
locally injective over $X$ does not change when the coverings
$\bW$ and $\bT$ are refined.
 Any locally $X$\+contrainjective $\bT$\+locally contraherent cosheaf
on $Y$ is $X$\+locally injective, but the converse need not be true.

 The full subcategory of $X$\+locally injective $\bT$\+locally
contraherent cosheaves on $Y$ will be denoted by $Y\lcth_\bT^\Xlin
\sub Y\lcth_\bT$.
 We will also use the notation $Y\lcth_\bT^{\lct,\Xlin}=
Y\lcth_\bT^\lct\cap Y\lcth_\bT^\Xlin$ for the full subcategory of
$X$\+locally injective locally cotorsion $\bT$\+locally contraherent
cosheaves on~$Y$.
 The full subcategory of $Y\lcth_\bT^\Xlin$ is closed under extensions,
cokernels of admissible monomorphisms, and infinite products in
the exact category $Y\lcth_\bT$; so it inherits an exact category
structure.
 The same applies to the full subcategory
$Y\lcth_\bT^{\lct,\Xlin}\sub Y\lcth_\bT^\lct$.
 If the morphism~$f$ is flat, then any locally injective
$\bT$\+locally contraherent cosheaf on $Y$ is $X$\+locally injective,
i.~e., $Y\lcth_\bT^\lin\sub Y\lcth_\bT^{\lct,\Xlin}
\sub Y\lcth_\bT^\Xlin$.
 We put $Y\ctrh^\Xlin=Y\lcth_{\{Y\}}^\Xlin$ and
$Y\ctrh^{\lct,\Xlin}=Y\lcth_{\{Y\}}^{\lct,\Xlin}$.

\begin{lem}
 Let $Y$ be a quasi-compact semi-separated scheme and $\pi\:Y\rarrow X$
be a flat affine morphism of schemes.
 Let\/ $\bW$ be an open covering of $X$ and\/ $\bT$ be an open covering
of $Y$ such that the morphism~$\pi$ is $(\bW,\bT)$\+affine.
 Then \par
\textup{(a)} the coresolution dimension of any $X$\+locally injective\/
$\bT$\+locally contraherent cosheaf on $X$ with respect to
the coresolving subcategory of locally $X$\+contrainjective\/
$\bT$\+locally contraherent cosheaves $Y\lcth_\bT^\Xlctin\sub
Y\lcth_\bT$ does not exceed~$1$; \par
\textup{(b)} for any symbol\/ $\bst=\b$, $+$, $-$, $\empt$, $\abs+$,
$\abs-$, $\ctr$, or\/ $\abs$, the triangulated functor\/
$\sD^\st(Y\lcth_\bT^\Xlctin)\rarrow\sD^\st(Y\lcth_\bT^\Xlin)$ induced
by the inclusion of exact categories $Y\lcth_\bT^\Xlctin\rarrow
Y\lcth_\bT^\Xlin$ is an equivalence of triangulated categories.
\end{lem}

\begin{proof}
 Part~(a): the full subcategory $Y\lcth_\bT^\Xlctin$ is coresolving in
$Y\lcth_\bT$ (hence also in $Y\lcth_\bT^\Xlin$) because there are
enough locally injective $\bT$\+locally contraherent cosheaves on~$Y$
(by Corollary~\ref{clp-cor}(a)), and all of them are locally
$X$\+contrainjective.
 The coresolution dimension~$1$ assertion follows from
Corollary~\ref{contrainjective-cor}(b) applied to the commutative
ring homomorphisms $\O_X(U)\rarrow\O_Y(f^{-1}(U))$, where $U$ ranges
over the affine open subschemes in $X$ subordinate to~$\bW$.
 Part~(b) follows from part~(a) in view of the opposite version of
Proposition~\ref{finite-resolutions}.
\end{proof}

 Let us introduce notation for the intersections of full subcategories
\begin{align*}
 Y\ctrh_\al^\Xlin &= Y\ctrh_\al\cap Y\lcth_\bT^\Xlin, \\
 Y\ctrh_\al^{\lct,\Xlin} &= Y\ctrh_\al\cap Y\lcth_\bT^{\lct,\Xlin}.
\end{align*}

\begin{lem} \label{contrah-X-loc-injective-finite-resol-dim}
 Let $Y$ be a quasi-compact semi-separated scheme and $\pi\:Y\rarrow X$
be a flat affine morphism of schemes.
 Let\/ $\bW$ be an open covering of $X$ and\/ $\bT$ be an open covering
of $Y$ such that the morphism~$\pi$ is $(\bW,\bT)$\+affine.
 Let $Y=\bigcup_{\alpha=1}^N V_\alpha$ be a finite affine open covering
of $Y$ subordinate to\/~$\bT$.
 In this context: \par
\textup{(a)} The resolution dimension of any $X$\+locally injective\/
$\bT$\+locally contraherent cosheaf on $X$ with respect to the resolving
subcategory of antilocal $X$\+locally injective contraherent cosheaves
$Y\ctrh_\al^\Xlin\sub Y\lcth_\bT^\Xlin$ does not exceed $N-1$.
 Consequently, the same bound holds for the resolution dimension of
any object of $Y\lcth_\bT^\Xlin$ with respect to the resolving
subcategory $Y\ctrh^\Xlin$. \par
\textup{(b)} The resolution dimension of any $X$\+locally injective
locally cotorsion\/ $\bT$\+locally contraherent cosheaf on $X$ with
respect to the resolving subcategory of antilocal $X$\+locally injective 
locally cotorsion contraherent cosheaves $Y\ctrh_\al^{\lct,\Xlin}\sub
Y\lcth_\bT^{\lct,\Xlin}$ does not exceed $N-1$.
 Consequently, the same bound holds for the resolution dimension of
any object of $Y\lcth_\bT^{\lct,\Xlin}$ with respect to the resolving
subcategory $Y\ctrh^{\lct,\Xlin}$.
\end{lem}

\begin{proof}
 Part~(a): let $\gK$ be an $X$\+locally injective $\bT$\+locally
contraherent cosheaf on~$Y$.
 By Corollary~\ref{clp-cor}(b), there exists a short exact sequence of
$\bT$\+locally contraherent cosheaves $0\rarrow\gJ\rarrow\P\rarrow\gK
\rarrow0$ with a locally injective $\bT$\+locally contraherent cosheaf
$\gJ$ and an antilocal contraherent cosheaf $\P$ on~$Y$.
 Then we have $\gJ\in Y\lcth_\bT^\Xlin$, hence
$\P\in Y\ctrh_\al^\Xlin$.
 Thus any object of $Y\lcth_\bT^\Xlin$ is the target of an admissible
epimorphism from an antilocal $X$\+locally injective contraherent
cosheaf on~$Y$.
 It remains to refer to Lemma~\ref{dil-cta-clp-finite-dim}(c) and
Corollary~\ref{fdim-subcategory-cor}.
 The proof of part~(b) is similar and based on
Corollary~\ref{clp-lct-cor}(b) together with
Lemma~\ref{lct-lin-clp-finite-dim}(a).
\end{proof}

\begin{cor} \label{contrah-lcth-loc-injective-over-base}
 Let $Y$ be a quasi-compact semi-separated scheme and $\pi\:Y\rarrow X$
be a flat affine morphism of schemes.
 Let\/ $\bW$ be an open covering of $X$ and\/ $\bT$ be an open covering
of $Y$ such that the morphism~$\pi$ is $(\bW,\bT)$\+affine.
 In this context: \par
\textup{(a)} For any symbol\/ $\bst=\b$, $+$, $-$, $\empt$, $\abs+$,
$\abs-$, $\ctr$, or\/ $\abs$, the inclusions of exact categories
$Y\ctrh_\al^\Xlin\rarrow Y\ctrh^\Xlin\rarrow Y\lcth_\bT^\Xlin$
induce equivalences of the derived categories
$$
 \sD^\st(Y\ctrh_\al^\Xlin)\simeq
 \sD^\st(Y\ctrh^\Xlin)\simeq\sD^\st(Y\lcth_\bT^\Xlin).
$$ \par
\textup{(b)} For any symbol\/ $\bst=\b$, $+$, $-$, $\empt$, $\abs+$,
$\abs-$, $\ctr$, or\/ $\abs$, the inclusions of exact categories
$Y\ctrh_\al^{\lct,\Xlin}\rarrow Y\ctrh^{\lct,\Xlin}\rarrow
Y\lcth_\bT^{\lct,\Xlin}$ induce equivalences of the derived categories
$$
 \sD^\st(Y\ctrh_\al^{\lct,\Xlin})\simeq
 \sD^\st(Y\ctrh^{\lct,\Xlin})\simeq\sD^\st(Y\lcth_\bT^{\lct,\Xlin}).
$$
\end{cor}

\begin{proof}
 All the assertions follow from the respective assertions of
Lemma~\ref{contrah-X-loc-injective-finite-resol-dim} in view of
Proposition~\ref{finite-resolutions}.
\end{proof}

 In rest of this section, our aim is to obtain representations of
the derived categories $\sD(Y\lcth_\bT^\Xlin)$ and
$\sD(Y\lcth_\bT^{\lct,\Xlin})$ as quotient categories of the homotopy
categories of complexes of locally injective (rather than merely
$X$\+locally injective) locally contraherent cosheaves on~$Y$.
 These results will be used in the construction of the derived functor
in Theorem~\ref{semico-semicontra-quadrality-theorem}.

 The following pair of lemmas is to be compared with
Lemmas~\ref{acyclic-complexes-of-cta-cotorsion-pair}\+-%
\ref{complexes-loc-hot-inj-of-inj-preenvelope}.

\begin{lem} \label{becker-coderived-of-cta-cot-affine-case}
\textup{(a)} Let\/ $\sR$ be the local class of all commutative rings $R$
and\/ $\sE^R=\sK^R=\Com(R\modl^\cta)$ be the exact category of
complexes of contraadjusted $R$\+modules.
 Let\/ $\sF(R)\sub\sE^R$ be the class of all Becker-coacyclic complexes
of contraadjusted $R$\+modules, and let\/ $\sC^R=\Com(R\modl^\inj)
\sub\sE^R$ be the class of all complexes of injective $R$\+modules.
 Then\/ $\sE$ and\/ $\sC$ are very colocal classes, and the pair of
classes $(\sF(R),\sC^R)$ is a hereditary complete cotorsion pair
in~$\sE^R$. \par
\textup{(b)} Let\/ $\sR$ be the local class of all commutative rings $R$
and\/ $\sE^R=\Com(R\modl^\cot)\sub\sK^R=\Com(R\modl^\cta)$ be the exact
category of complexes of cotorsion $R$\+modules.
 Let\/ $\sF(R)\sub\sE^R$ be the class of all Becker-coacyclic complexes
of cotorsion $R$\+modules, and let\/ $\sC^R=\Com(R\modl^\inj)\sub
\sE^R$ be the class of all complexes of injective $R$\+modules.
 Then\/ $\sE$ and\/ $\sC$ are very colocal classes, and the pair of
classes $(\sF(R),\sC^R)$ is a hereditary complete cotorsion pair
in\/~$\sE^R$.
\end{lem}

\begin{proof}
 The mentioned classes are very colocal by
Examples~\ref{colocal-classes-examples}.
 In both parts~(a) and~(b), the full subcategory $\sC^R$ is obviously
closed under the cokernels of admissible monomorphisms in~$\sE^R$.
 Discussions of the cotorsion pairs in parts~(a) and~(b) can be found
in~\cite[Examples~7.10 and~7.11]{Pal} (alternatively, these are
the special cases of the respective assertions of
Corollary~\ref{becker-coderived-of-cta-cot} for affine schemes).
\end{proof}

\begin{lem} \label{contrah-glued-from-coderived-cotorsion-pairs}
 Let $Y$ be a quasi-compact semi-separated scheme with an open
covering\/ $\bT$ and a finite affine open covering
$Y=\bigcup_\alpha V_\alpha$ subordinate to\/~$\bT$.
 In this context: \par
\textup{(a)} There is a hereditary complete cotorsion pair $(\sF,\sC)$
in the exact category\/ $\Com(Y\lcth_\bT)$ of complexes of\/
$\bT$\+locally contraherent cosheaves on $Y$ such that\/
$\sC=\Com(Y\lcth_\bT^\lin)$ is the class of all complexes of
locally injective\/ $\bT$\+locally contraherent cosheaves on $Y$,
while\/ $\sF=\sF_Y$ is the class of all direct summands of
finitely iterated extensions in\/ $\Com(Y\ctrh_\bT)$ of the direct
images of Becker-coacyclic complexes\/ $\Q_\alpha^\bu$ in the exact
categories $V_\alpha\ctrh$. \par
\textup{(b)} There is a hereditary complete cotorsion pair $(\sF,\sC)$
in the exact category\/ $\Com(Y\lcth_\bT^\lct)$ of complexes of locally
cotorsion\/ $\bT$\+locally contraherent cosheaves on $Y$ such that\/
$\sC=\Com(Y\lcth_\bT^\lin)$ is the class of all complexes of
locally injective\/ $\bT$\+locally contraherent cosheaves on $Y$,
while\/ $\sF=\sF_Y$ is the class of all direct summands of
finitely iterated extensions in\/ $\Com(Y\ctrh^\lct)$ of the direct
images of Becker-coacyclic complexes\/ $\Q_\alpha^\bu$ in the exact categories $V_\alpha\ctrh^\lct$.
\end{lem}

\begin{proof}
 Part~(a): apply Theorem~\ref{loc-contraherent-gluing-theorem}
to the datum of classes $\sR$, \,$\sE^R$, \,$\sC^R$, and $\sF(R)$ from
Lemma~\ref{becker-coderived-of-cta-cot-affine-case}(a).
 The proof of part~(b) is similarly based on
Lemma~\ref{becker-coderived-of-cta-cot-affine-case}(b).
\end{proof}

\begin{lem} \label{F-Y-F-X-classes-direct-image}
 Let $X$ be a quasi-compact semi-separated scheme and $\pi\:Y\rarrow X$
be an affine morphism of schemes.
 Then \par
\textup{(a)} the direct image functor $\pi_!\:Y\ctrh\rarrow X\ctrh$
takes the full subcategory\/ $\sF_Y$ constructed in
Lemma~\textup{\ref{contrah-glued-from-coderived-cotorsion-pairs}(a)}
into the respective full subcategory\/ $\sF_X\sub\Com(X\ctrh)$; 
\emergencystretch=1em\par
\textup{(b)} the direct image functor $\pi_!\:Y\ctrh^\lct\rarrow
X\ctrh^\lct$ takes the full subcategory\/ $\sF_Y$ constructed in
Lemma~\textup{\ref{contrah-glued-from-coderived-cotorsion-pairs}(b)}
into the respective full subcategory\/ $\sF_X\sub\Com(X\ctrh^\lct)$.
\end{lem}

\begin{proof}
 Part~(a): choose open coverings $\bW$ and $\bT$ of the schemes $X$
and $Y$ such that the morphism~$\pi$ is $(\bW,\bT)$\+affine, and
choose a finite affine open covering $X=\bigcup_\alpha U_\alpha$
of the scheme $X$ subordinate to~$\bW$.
 Let $Y=\bigcup_\alpha V_\alpha$, where $V_\alpha=\pi^{-1}(U_\alpha)
\sub Y$, be the corresponding affine open covering of $Y$
subordinate to~$\bT$.
 Then the description of the full subcategory $\sF_Y$ provided in
Lemma~\ref{contrah-glued-from-coderived-cotorsion-pairs}(a) together
with the fact that the functor $\pi_!\:Y\lcth_\bT\rarrow X\lcth_\bW$
is exact and with Lemma~\ref{qcoh-affine-direct-image-becker-coacyclic}
(for the morphisms of affine schemes $\pi_\alpha\:V_\alpha\rarrow
U_\alpha$) implies the desired assertion.
 Alternatively, one can choose the coverings $\bW$ and $\bT$ so that
the morphism~$\pi$ is both $(\bW,\bT)$\+affine and
$(\bW,\bT)$\+coaffine, and use Lemma~\ref{Ext-1-as-homotopy-Hom}
together with
the partial adjunction~\eqref{direct-inverse-lin-adjunction} and
the existence of a well-defined inverse image functor
$\pi^!\:X\lcth_\bW^\lin\rarrow Y\lcth_\bT^\lin$.
 The proof of part~(b) is similar.
\end{proof}

\begin{thm} \label{loc-injective-descript-of-derived-X-loc-injective}
 Let $X$ be a quasi-compact semi-separated scheme and $\pi\:Y\rarrow X$
be a flat affine morphism of schemes.
 Let\/ $\bW$ be an open covering of $X$ and\/ $\bT$ be an open covering
of $Y$ such that the morphism~$\pi$ is $(\bW,\bT)$\+affine.
 Then \par
\textup{(a)} the triangulated functor
$$
 \frac{\Hot(Y\lcth_\bT^\lin)}
 {\Acycl(Y\lcth_\bT^\Xlin)\cap\Hot(Y\lcth_\bT^\lin)}
 \lrarrow\sD(Y\lcth_\bT^\Xlin)
$$
induced by the embedding of additive categories
$Y\lcth_\bT^\lin\rarrow Y\lcth_\bT^\Xlin$ is an equivalence of
triangulated categories; \par
\textup{(b)} the triangulated functor
$$
 \frac{\Hot(Y\lcth_\bT^\lin)}
 {\Acycl(Y\lcth_\bT^{\lct,\Xlin})\cap\Hot(Y\lcth_\bT^\lin)}
 \lrarrow\sD(Y\lcth_\bT^{\lct,\Xlin})
$$
induced by the embedding of additive categories
$Y\lcth_\bT^\lin\rarrow Y\lcth_\bT^\Xlin$ is an equivalence of
triangulated categories.
\end{thm}

\begin{proof}
 Part~(a): in view of Lemma~\ref{pkoszul-lemma16}(b), it suffices to
construct, for every complex of $X$\+locally injective $\bT$\+locally
contraherent cosheaves $\gK^\bu$ on $Y$, a complex of locally injective
$\bT$\+locally contraherent cosheaves $\gJ^\bu$ together with
a morphism of complexes $\gK^\bu\rarrow\gJ^\bu$ with a cone
acyclic in $Y\lcth_\bT^\Xlin$.
 For this purpose, consider a special preenvelope short exact sequence
$0\rarrow\gK^\bu\rarrow\gJ^\bu\rarrow\Q^\bu\rarrow0$ in the hereditary
complete cotorsion pair of
Lemma~\ref{contrah-glued-from-coderived-cotorsion-pairs}(a).
 Then $\gJ^\bu$ is a complex in $Y\lcth_\bT^\lin$ and $\gK^\bu$
is a complex in $Y\lcth_\bT^\Xlin$, so $\Q^\bu$ is a complex in
$Y\lcth_\bT^\Xlin$.
 Hence $\pi_!\Q^\bu$ is a complex of locally injective
$\bW$\+locally contraherent cosheaves on~$X$.

 On the other hand, $\Q^\bu$ is a complex in the category~$\sF_Y$.
 By Lemma~\ref{F-Y-F-X-classes-direct-image}(a), it follows that
$\pi_!\Q^\bu$ is a complex in the category~$\sF_X$.
 But $\sF_X$ and $\Com(X\lcth_\bW^\lin)$ are $\Ext^1$\+orthogonal
classes
(by Lemma~\ref{contrah-glued-from-coderived-cotorsion-pairs}(a)
applied to the scheme $X$ with its open covering~$\bW$).
 As both the classes are closed under shifts, we can refer to
Lemma~\ref{Ext-1-as-homotopy-Hom} and conclude that the complex
$\pi_!\Q^\bu$ is contractible.
 Hence the complex $\Q^\bu$ is acyclic in $Y\lcth_\bT^\Xlin$
(in view of Lemma~\ref{acyclicity-in-lcth-criterion}(a)), and it follows
that the cone of the morphism $\gK^\bu\rarrow\gJ^\bu$ is also acyclic
in $Y\lcth_\bT^\Xlin$.

 The proof of part~(b) is similar, and based on
Lemmas~\ref{contrah-glued-from-coderived-cotorsion-pairs}(b)
and~\ref{F-Y-F-X-classes-direct-image}(b).
\end{proof}

\begin{cor} \label{contrah-lcta-lct-loc-injective-over-base}
 Let $X$ be a quasi-compact semi-separated scheme and $\pi\:Y\rarrow X$
be a flat affine morphism of schemes.
 Let\/ $\bW$ be an open covering of $X$ and\/ $\bT$ be an open covering
of $Y$ such that the morphism~$\pi$ is $(\bW,\bT)$\+affine.
 Then the triangulated functor
$$
 \sD(Y\lcth_\bT^{\lct,\Xlin})\lrarrow\sD(Y\lcth_\bT^\Xlin)
$$
induced by the embedding of exact categories
$Y\lcth_\bT^{\lct,\Xlin}\rarrow Y\lcth_\bT^\Xlin$ is an equivalence
of triangulated categories.
\end{cor}

\begin{proof}
 It suffices to compare the two parts of
Theorem~\ref{loc-injective-descript-of-derived-X-loc-injective}.
 Recall that any complex in $Y\lcth_\bT^\lct$ that is acyclic in
$Y\lcth_\bT$ is also acyclic in $Y\lcth_\bT^\lct$ (by
Lemma~\ref{acyclicity-in-lcth-criterion}); hence any complex in
$Y\lcth_\bT^\lin$ that is acyclic in $Y\lcth_\bT^\Xlin$ is
also acyclic in $Y\lcth_\bT^{\lct,\Xlin}$.
\end{proof}

\subsection{Equivalence on the $X$-flat side}
 In this section we discuss the definitions from
Section~\ref{sheaves-flat-over-base-subsect} and the ones from
Section~\ref{cosheaves-projective-over-base-subsect} adopted to
the context of quasi-compact semi-separated scheme $X$, as per
Remark~\ref{cosh-proj-over-qcomp-ssep-base-remark}.

 In particular, let us introduce the notation $Y\lcth_{\bT,\Xalf}\sub
Y\lcth_\bT$ for the full subcategory of $X$\+antilocally flat
$\bT$\+locally contraherent cosheaves on $Y$, as defined in
Remark~\ref{cosh-proj-over-qcomp-ssep-base-remark}.
 Put $Y\ctrh_\Xalf=Y\lcth_{\{Y\},\Xalf}$ and
$$
 Y\ctrh_{\al,\Xalf}=Y\ctrh_\al\cap Y\ctrh_\Xalf=
 Y\ctrh_\al\cap Y\lcth_{\bT,\Xalf}.
$$
 Notice that $Y\lcth_{\bT,\Xalf}=Y\lcth_\bT^\Xfl$ when the scheme $X$
is semi-separated and Noetherian of finite Krull dimension,
because $X\ctrh_\alf=X\ctrh^\fl$ in this case by
Corollary~\ref{finite-krull-flat-contraherent}(a)
(cf.\ Remark~\ref{flat-antilocally-flat-remark}).

 The following lemma is a relative version/generalization of
Lemma~\ref{cta-clp-restricts-to-prj-clf}.

\begin{lem} \label{cta-clp-restricts-to-flat-side-over-base}
 Let $\pi\:Y\rarrow X$ be an affine morphism of quasi-compact
semi-separated schemes.
 Then \par
\textup{(a)} the equivalence of exact categories $Y\qcoh^\cot\simeq
Y\ctrh_\al^\lct$ from
Lemma~\textup{\ref{cta-clp-restricts-to-cot-inj}(a)}
resticts to an equivalence of full exact subcategories
$Y\qcoh^\cot_\Xfl\simeq Y\ctrh_{\al,\Xprj}^\lct$; \par
\textup{(b)} the equivalence of exact categories $Y\qcoh^\cta\simeq
Y\ctrh_\al$ from Lemma~\textup{\ref{cta-clp-equivalence}}
restricts to an equivalence of full exact subcategories
$Y\qcoh^\cta_\Xfl\simeq Y\ctrh_{\al,\Xalf}$; \par
\textup{(c)} the equivalence of exact categories $Y\qcoh^\cta\simeq
Y\ctrh_\al$ from Lemma~\textup{\ref{cta-clp-equivalence}}
restricts to an equivalence of full exact subcategories
$Y\qcoh^\cta_\Xvfl\simeq Y\ctrh_{\al,\Xprj}$.
\end{lem}

\begin{proof}
 The proofs of all the three parts are based on
Corollary~\ref{cta-clp-direct-image-cor}.
 In view of that Corollary, part~(a) follows from
Lemma~\ref{cta-clp-restricts-to-prj-clf}(b) applied to
the scheme~$X$ (one has to keep in mind the facts that the direct
image functor~$f_*$ takes $Y\qcoh^\cot$ to $X\qcoh^\cot$ and
the direct image functor~$f_!$ takes $Y\ctrh^\lct$ to $X\ctrh^\lct$).
 Part~(b) follows from Lemma~\ref{cta-clp-restricts-to-prj-clf}(c)
applied to the scheme~$X$, and, similarly, part~(c) follows from
Lemma~\ref{cta-clp-restricts-to-prj-clf}(a) applied to the scheme~$X$.
 In all the three cases, the argument is that the horizontal category
equivalences in commutative diagram~\eqref{cta-clp-direct-image-diagram}
preserve the preimages of full subcategories.
\end{proof}

\begin{cor} \label{flat-side-over-base-equivalences}
 Let $\pi\:Y\rarrow X$ be a flat affine morphism of quasi-compact
semi-separated schemes.
 Let\/ $\bW$ be an open covering of $X$ and\/ $\bT$ be an open covering
of $Y$ such that the morphism~$\pi$ is $(\bW,\bT)$\+affine.
 Then there is the upper half of the following commutative diagram of
triangulated equivalences:
\begin{equation} \label{flat-side-over-base-diagram}
\begin{gathered}
 \xymatrix{
  \sD(Y\qcoh_\Xfl) \ar@<-2pt>@{-}[r] \ar@{=}[d]
  & \sD(Y\qcoh^\cot_\Xfl) \ar@<-2pt>[l] \ar@<-2pt>[d] \ar@{=}[r]
  & \sD(Y\ctrh^\lct_{\al,\Xprj}) \ar@<2pt>[r] \ar@<2pt>[d] 
  & \sD(Y\lcth^\lct_{\bT,\Xprj}) \ar@<2pt>@{-}[l] \ar@<2pt>[d]
  \\
  \sD(Y\qcoh_\Xfl) \ar@<-2pt>@{-}[r] \ar@<2pt>@{--}[d]
  & \sD(Y\qcoh^\cta_\Xfl) \ar@<-2pt>[l]
  \ar@<-2pt>@{-}[u] \ar@<2pt>@{--}[d] \ar@{=}[r]
  & \sD(Y\ctrh_{\al,\Xalf}) \ar@<2pt>[r] \ar@<2pt>@{--}[u]
  \ar@<-2pt>@{--}[d]
  & \sD(Y\lcth_{\bT,\Xalf}) \ar@<2pt>@{-}[l] \ar@<2pt>@{-}[u]
  \ar@<-2pt>@{-}[d]
  \\
  \sD(Y\qcoh_\Xvfl) \ar@<2pt>@{-}[r] \ar@<2pt>[u]
  & \sD(Y\qcoh^\cta_\Xvfl) \ar@<2pt>[l] \ar@{=}[r] \ar@<2pt>[u]
  & \sD(Y\ctrh_{\al,\Xprj}) \ar@<-2pt>[r] \ar@<-2pt>[u]
  & \sD(Y\lcth_{\bT,\Xprj}) \ar@<-2pt>@{-}[l] \ar@<-2pt>[u]
 }
\end{gathered}
\end{equation}
 When the morphism~$\pi$ is very flat, the lower half of the commutative
diagram of triangulated equivalences~\eqref{flat-side-over-base-diagram}
is also valid.
\end{cor}

\begin{proof}
 The middle column of horizontal equality signs is formed by
the equivalences of derived categories induced by the equivalences
of exact categories from
Lemma~\ref{cta-clp-restricts-to-flat-side-over-base}.
 The vertical equality in the upper corner on the left is
the identity functor.
 All the other functors are induced by the obvious embeddings of
exact categories.
 The leftmost column of horizontal equivalences is provided by
Corollary~\ref{quasi-X-flat-cta-derived-equiv}
and Theorem~\ref{quasi-X-flat-cot-derived-equiv}.
 The same applies to the upper vertical equivalence in the middle
left column.
 In the case when $X$ is Noetherian of finite Krull dimension,
the lower vertical equivalence in the leftmost column is provided by
Lemma~\ref{fl-vlf-over-base-lemma}(b); otherwise, it follows from
commutativity of the diagram.
 The rightmost column of horizontal equivalences is provided by
Corollary~\ref{contrah-lcth-projective-flat-over-base} with
Remark~\ref{cosh-proj-over-qcomp-ssep-base-remark}.
 The vertical equivalences in the rightmost column are
Corollary~\ref{lcth-lcta-lct-projective-flat-over-base} with
Remark~\ref{cosh-proj-over-qcomp-ssep-base-remark}.
 All the remaining functors are triangulated equivalences in view of
the commutativity of the diagram.
\end{proof}

\subsection{Equivalence on the $X$-injective side}
 In this section we discuss the definitions from
Section~\ref{cosheaves-loc-injective-over-base-subsect} and the ones
from Section~\ref{sheaves-injective-over-base-subsect} adopted to
the context of a non-Noetherian scheme $X$ and an affine morphism~$\pi$,
as per Remark~\ref{sheaves-inj-over-non-noetherian-base-remark}.

 The following lemma is a relative version/generalization of
Lemma~\ref{cta-clp-restricts-to-cot-inj}(b).

\begin{lem} \label{cta-clp-restricts-to-injective-side-over-base}
 Let $\pi\:Y\rarrow X$ be an affine morphism of quasi-compact
semi-separated schemes.
 Then \par
\textup{(a)} the equivalence of exact categories $Y\qcoh^\cta\simeq
Y\ctrh_\al$ from Lemma~\textup{\ref{cta-clp-equivalence}}
restricts to an equivalence of full exact subcategories
$Y\qcoh^{\cta,\Xinj}\simeq Y\ctrh_\al^\Xlin$; \par
\textup{(b)} the equivalence of exact categories $Y\qcoh^\cot\simeq
Y\ctrh_\al^\lct$ from
Lemma~\textup{\ref{cta-clp-restricts-to-cot-inj}(a)}
resticts to an equivalence of full exact subcategories
$Y\qcoh^{\cot,\Xinj}\simeq Y\ctrh_\al^{\lct,\Xlin}$.
\end{lem}

\begin{proof}
 Both part follow from
Lemma~\ref{cta-clp-restricts-to-cot-inj}(b) (applied to the scheme~$X$)
in view of Corollary~\ref{cta-clp-direct-image-cor}.
 In part~(b), it is helpful to in mind the facts that the direct
image functor~$f_*$ takes $Y\qcoh^\cot$ to $X\qcoh^\cot$ and
the direct image functor~$f_!$ takes $Y\ctrh^\lct$ to $X\ctrh^\lct$.
 In both cases, the argument is that the horizontal category
equivalences in commutative diagram~\eqref{cta-clp-direct-image-diagram}
preserve the preimages of full subcategories.
\end{proof}

\begin{cor} \label{injective-side-over-base-equivalences}
 Let $\pi\:Y\rarrow X$ be a flat affine morphism of quasi-compact
semi-separated schemes.
 Let\/ $\bW$ be an open covering of $X$ and\/ $\bT$ be an open covering
of $Y$ such that the morphism~$\pi$ is $(\bW,\bT)$\+affine.
 Then there is the following commutative diagram of triangulated
equivalences: \hfuzz=6pt
\begin{equation} \label{injective-side-over-base-diagram}
\begin{gathered}
 \!\!\.
 \xymatrix{
  \sD(Y\qcoh^\Xinj) \ar@<-2pt>@{-}[r] \ar@{=}[d]
  & \sD(Y\qcoh^{\cot,\Xinj}) \ar@<-2pt>[l] \ar@<-2pt>[d] \ar@{=}[r]
  & \sD(Y\ctrh_\al^{\lct,\Xlin}) \ar@<2pt>[r] \ar@<2pt>[d] 
  & \sD(Y\lcth_\bT^{\lct,\Xlin}) \ar@<2pt>@{-}[l] \ar@<2pt>[d]
  \\
  \sD(Y\qcoh^\Xinj) \ar@<-2pt>@{-}[r]
  & \sD(Y\qcoh^{\cta,\Xinj}) \ar@<-2pt>[l] \ar@<-2pt>@{-}[u] \ar@{=}[r]
  & \sD(Y\ctrh_\al^\Xlin) \ar@<2pt>[r] \ar@<2pt>@{--}[u]
  & \sD(Y\lcth_\bT^\Xlin) \ar@<2pt>@{-}[l] \ar@<2pt>@{-}[u]
 }
\end{gathered}
\end{equation}
\end{cor}

\begin{proof}
 The middle column of horizontal equality signs is formed by
the equivalences of derived categories induced by the equivalences
of exact categories from
Lemma~\ref{cta-clp-restricts-to-injective-side-over-base}.
 The vertical equality in the leftmost column is the identity functor.
 All the other functors are induced by the obvious embeddings of
exact categories.
 The leftmost commutative square of triangulated equivalences is
provided by Theorem~\ref{quasi-X-injective-cot-derived-equiv} with
Remark~\ref{sheaves-inj-over-non-noetherian-base-remark}.
 The horizontal equivalences in the rightmost commutative square are
Corollary~\ref{contrah-lcth-loc-injective-over-base}.
 The vertical equivalence in the rightmost column is
Corollary~\ref{contrah-lcta-lct-loc-injective-over-base}.
 The vertical functor in the middle right column is a triangulated
equivalence in view of the commutativity of the diagram.
\end{proof}

\subsection{Quadrality diagram}
 In this section, we consider a semi-separated Noetherian scheme $X$
with a dualizing complex $\D_X^\bu$, and a flat affine morphism of
schemes $\pi\:Y\rarrow X$.
 The following semico-semicontra correspondence quadrality theorem
is a generalization of Theorem~\ref{co-contra-dualizing}.

\begin{thm} \label{semico-semicontra-quadrality-theorem}
 There are natural triangulated equivalences between the four
triangulated categories\/ $\sD(Y\qcoh_\Xfl)$, \ $\sD^\si_X(Y\qcoh)$, \
$\sD^\si_X(Y\ctrh)$, and\/ $\sD(Y\ctrh^\Xlin)$, forming a commutative
square diagram
\begin{equation} \label{semico-semicontra-quadrality-diagram}
\begin{gathered}
 \xymatrix{
  \sD(Y\qcoh_\Xfl) \ar@<3pt>[rrrrrr]^{\pi^*\D_X^\subbu\ot_{\O_Y}{-}}   
  \ar@<3pt>[ddd]^{\boR\fHom_Y(\O_Y,{-})}
  &&&&&& \sD^\si_X(Y\qcoh)
  \ar@<3pt>[llllll]^{\boR\qHom_{Y\qc}(\pi^*\D_X^\subbu,{-})}
  \ar@<-3pt>[llllllddd]_{\boR\fHom_Y(\pi^*\D_X^\subbu,{-})\quad}
  \ar@<3pt>[ddd]^{\boR\fHom_Y(\O_Y,{-})} \\ \\ \\
  \sD^\si_X(Y\ctrh) \ar@<3pt>@{-}[rrrrrr]
  \ar@<3pt>[uuu]^{\O_Y\ocn_Y^\boL{-}}
  \ar@<-3pt>[rrrrrruuu]_{\quad\pi^*\D_X^\subbu\ocn_Y^\boL{-}}
  &&&&&& \sD(Y\ctrh^\Xlin)
  \ar@<3pt>[llllll]^{\boR\Cohom_Y(\pi^*\D_X^\subbu,{-})}
  \ar@<3pt>[uuu]^{\O_Y\ocn_Y^\boL{-}}
 }
\end{gathered}
\end{equation}
 Among these equivalences, the equivalences\/
$\sD(Y\qcoh_\Xfl)\simeq\sD^\si_X(Y\ctrh)$ and\/ $\sD^\si_X(Y\qcoh)
\simeq\sD(Y\ctrh^\Xlin)$ do not require a dualizing complex and do not
depend on it; all the remaining equivalences do and do.

 Moreover, the commutative square
diagram~\eqref{semico-semicontra-quadrality-diagram} together with
the commutative square diagram~\eqref{co-contra-quadrality-diagram}
from Theorem~\textup{\ref{co-contra-dualizing}}, the direct image
functors\/ $\pi_*$~\eqref{semicoderived-to-X-underived-direct-image}
and\/ $\pi_!$~\eqref{semicontraderived-to-X-underived-direct-image},
and the obvious underived direct image functors
\begin{equation} \label{X-flat-to-flat-induced-direct-image}
 \pi_*\:\sD(Y\qcoh_\Xfl)\lrarrow\sD(X\qcoh_\fl)
\end{equation}
and
\begin{equation} \label{X-lin-to-lin-induced-direct-image}
 \pi_!\:\sD(Y\ctrh^\Xlin)\lrarrow\sD(X\ctrh^\lin).
\end{equation}
form the diagram of a commutative cube of triangulated functors and
triangulated equivalences.
\end{thm}

\begin{rem} \label{semicosheaves-remark}
 Notice that the commutative diagram of the quadrality of
semico-semicontra
correspondence~\eqref{semico-semicontra-quadrality-diagram} is
\emph{less symmetric} than the similar square diagram of co-contra
correspondence~\eqref{co-contra-quadrality-diagram}.
 In particular, there is \emph{no} arrow pointing rightwards in
the lower horizontal line of
the diagram~\eqref{semico-semicontra-quadrality-diagram}.
 So we offer no explicit construction of the inverse functor
to the triangulated equivalence $\boR\Cohom_Y(\pi^*\D_X^\bu,{-})$
in the lower line.
 The reason is that the dualizing complex $\D_X^\bu$ is a complex
of injective quasi-coherent sheaves on a Noetherian scheme $X$,
and such sheaves are coadjusted (so suitable for the construction
of the contraherent tensor product functor~$\ot_{X\ct}$), as per
Section~\ref{contraherent-tensor}.
 However, the terms of the complex of quasi-coherent sheaves
$\pi^*\D_X^\bu$ on the scheme $Y$ are \emph{not} coadjusted.
 Generally, coadjusted modules over non-Noetherian commutative rings
may be scarce (see Remark~\ref{coadjusted-over-non-Noetherian-remark}).
 But even existence of enough coadjusted quasi-coherent sheaves on $Y$
would not solve the problem yet, as a contraherent tensor product
functor pointing rightwards in the lower line
of~\eqref{semico-semicontra-quadrality-diagram} has to produce a complex
of \emph{$X$\+locally injective} contraherent cosheaves on~$Y$.

 Furthermore, the arrow pointing leftwards in the lower horizontal
line of the diagram~\eqref{co-contra-quadrality-diagram} was
an \emph{underived} functor $\Cohom_X(\D_X^\bu,{-})$, which is
dual-analogous to the underived tensor product functor
$\D_X^\bu\ot_{\O_X}{-}$ pointing rightwards in the upper line
of~\eqref{co-contra-quadrality-diagram}.
 On the diagram~\eqref{semico-semicontra-quadrality-diagram},
the arrow pointing rightwards in the upper line is still an underived
tensor product functor $\pi^*\D_X^\bu\ot_{\O_Y}{-}$, but
the dual-analogous $\Cohom$ functor pointing leftwards in the lower
line of the diagram needs to be \emph{derived}.
 The reason for that is the problem illustrated by
Example~\ref{not-contrainjective-example}:
given a commutative ring homomorphism $R\rarrow S$ (even such that
$S$ is a flat $R$\+module), an $R$\+injective $S$\+contraadjusted
(or even $S$\+cotorsion) $S$\+module $K$, and an $R$\+module $M$,
it does not follow that $\Hom_R(M,K)\simeq\Hom_S(S\ot_RM\;K)$
is a contraadjusted $S$\+module.

 In both these difficulties arising in the construction of
the triangulated equivalence in the lower line
of~\eqref{semico-semicontra-quadrality-diagram}, the fundamental
problem seems to be that we want to produce complexes of contraherent
cosheaves (or at least, locally contraherent cosheaves) on the whole
scheme $Y$, which means that they should assign well-defined modules of
cosections to all (or all small enough) affine open subschemes
$V\sub Y$.
 If we were content to restrict our consideration to the affine open
subschemes in $Y$ of the form $V=\pi^{-1}(U)$, where $U$ ranges over
the affine open subschemes of $X$, the problem would disappear.
 But if one wants to have well-behaved modules of cosections over
smaller affine open subschemes in $Y$, then contraadjustedness along
$Y$ (i.~e., contraadjustedness of $S$\+modules with respect to elements
$s\in S$, rather than just $s\in R$) becomes a requirement that
is harder to satisfy.

 This problem may be solvable, after all: e.~g., to construct
the desired contraherent tensor product functor producing complexes
of $X$\+locally injective contraherent cosheaves on $Y$, one can try
to coresolve $R$\+injective $S$\+modules by $R$\+injective
contraadjusted/cotorsion $S$\+modules, or indeed coresolve whole
$X$\+locally injective ``$X$\+contraherent $X$\+semicosheaves'' on $Y$
(with $\O_Y(V)$\+modules of cosections defined for open subschemes
$V=\pi^{-1}(U)\sub Y$)---coresolve them by actual $X$\+locally
injective contraherent cosheaves on~$Y$.
 Here an ``$X$\+contraherent $X$\+semicosheaf'' on $Y$ is essentially
the same thing as a contraherent cosheaf of modules on $X$ over
the quasi-coherent sheaf of rings $\cA=\pi_*\O_Y$, in the sense
of the discussion in~\cite[Section~3.2]{Pphil}.
 We do not go into these details in this book.

 To put it in other words: There are, basically, two points of view on
affine morphisms of schemes $\pi\:Y\rarrow X$.
 The geometric point of view, adopted in this book, presumes a natural
generalization to nonaffine morphisms~$\pi$.
 The algebraic approach, which is not taken up in this book, describes
affine morphisms~$\pi$ in terms of the related quasi-coherent
commutative algebras $\cA=\pi_*\O_Y$ over the scheme~$X$.
 This leads to a natural generalization to (cosheaves of modules over)
quasi-coherent noncommutative algebras, quasi-coherent quasi-algebras
(as in~\cite[Section~2.3]{Pedg}), etc.
 It appears that the geometric point of view allows a better
understanding of the diagonal in the big square
diagram~\eqref{semico-semicontra-quadrality-diagram}
(in that the diagonal equivalence works for nonaffine morphisms~$\pi$,
as per Theorem~\ref{semico-semicontra-diagonal}), while the algebraic
approach to affine morphisms of schemes sheds more light on the lower
horizontal line in~\eqref{semico-semicontra-quadrality-diagram}.
 The algebraic approach has been developed to some extent in
the very recent preprint~\cite{Pdomc}, but quite in the direction
anticipated in the present remark.
\end{rem}

\begin{proof}[Proof of
Theorem~\ref{semico-semicontra-quadrality-theorem}]
 The triangulated equivalence in the upper horizontal line of
the desired diagram~\eqref{semico-semicontra-quadrality-diagram}
(``on the sheaf side'') is a particular case
of~\cite[Theorem~7.16]{Psemten}.
 The result in~\cite{Psemten} is more general in that it applies to
flat affine morphisms of ind-schemes (rather than just schemes).
 The commutativity of the square diagram formed by the upper
line of~\eqref{semico-semicontra-quadrality-diagram}, the upper line
of~\eqref{co-contra-quadrality-diagram}, and the direct image
functors $\pi_*$~\eqref{X-flat-to-flat-induced-direct-image}
and $\pi_*$~\eqref{semicoderived-to-X-underived-direct-image}
is also explained in the proof of~\cite[Theorem~7.16]{Psemten}.

 The triangulated equivalence in the diagonal of the square
diagram~\eqref{semico-semicontra-quadrality-diagram} is a particular
case of Theorem~\ref{semico-semicontra-diagonal}.
 The result of Theorem~\ref{semico-semicontra-diagonal} is more general
in that it applies to flat \emph{nonaffine} morphisms of schemes.
 The commuativity of the square diagram formed by the diagonal line
of~\eqref{semico-semicontra-quadrality-diagram}, the diagonal line
of~\eqref{co-contra-quadrality-diagram}, and the direct image
functors $\pi_*$~\eqref{semicoderived-to-X-underived-direct-image}
and $\pi_!$~\eqref{semicontraderived-to-X-underived-direct-image}
is also a part of Theorem~\ref{semico-semicontra-diagonal}.

 The vertical triangulated equivalence in the leftmost column
of~\eqref{semico-semicontra-quadrality-diagram}
(``on the $X$\+flat side'') is obtained by combining
Corollary~\ref{flat-side-over-base-equivalences} with
Theorem~\ref{contrah-projective-flat-over-base-semiderived}(b).
 The commutativity of the square diagram formed by the leftmost column
of~\eqref{semico-semicontra-quadrality-diagram}, the leftmost column
of~\eqref{co-contra-quadrality-diagram}, and the direct image
functors $\pi_*$~\eqref{X-flat-to-flat-induced-direct-image}
and $\pi_!$~\eqref{semicontraderived-to-X-underived-direct-image}
follows essentially from Corollary~\ref{cta-clp-direct-image-cor}.

 The vertical triangulated equivalence in the rightmost column
of~\eqref{semico-semicontra-quadrality-diagram}
(``on the $X$\+injective side'') is obtained by combining
Corollary~\ref{injective-side-over-base-equivalences} with
Theorem~\ref{quasi-injective-over-base-semiderived}.
 Once again, the commutativity of the square diagram formed by
the rightmost column
of~\eqref{semico-semicontra-quadrality-diagram}, the rightmost column
of~\eqref{co-contra-quadrality-diagram}, and the direct image functors
$\pi_*$~\eqref{semicoderived-to-X-underived-direct-image}
and $\pi_!$~\eqref{X-lin-to-lin-induced-direct-image} follows
essentially from Corollary~\ref{cta-clp-direct-image-cor}.

 In order to prove the theorem, we mainly need to construct
the derived functor $\boR\Cohom_Y(\pi^*\D_X^\bu,{-})$ in the lower
horizontal line of~\eqref{semico-semicontra-quadrality-diagram},
and show that the diagram of triangulated
functors~\eqref{semico-semicontra-quadrality-diagram} is commutative.
 Other aspects of the picture are also discussed in the arguments below.

 Let $\E^\bu$ be an arbitrary finite complex of injective quasi-coherent
sheaves on~$X$ (we are really interested in the case $\E^\bu=\D_X^\bu$).
 Let us redraw the vertical parts of the big diagram taking into account
the results of Corollary~\ref{flat-side-over-base-equivalences} with
Theorem~\ref{contrah-projective-flat-over-base-semiderived}(b)
in the leftmost column and the results of
Corollary~\ref{injective-side-over-base-equivalences} with
Theorem~\ref{quasi-injective-over-base-semiderived}
in the rightmost column:
\begin{equation} \label{semico-semicontra-big-diagram-partly-redrawn}
\begin{gathered}
 \xymatrix{
  \sD(Y\qcoh_\Xfl^\cta)   
  \ar@<3pt>[ddd]^{\fHom_Y(\O_Y,{-})}
  &&&&&& \sD(Y\qcoh^{\cta,\Xinj})
  \ar@<-3pt>[llllllddd]_{\boR\fHom_Y(\pi^*\E^\subbu,{-})\quad}
  \ar@<3pt>[ddd]^{\fHom_Y(\O_Y,{-})} \\ \\ \\
  \sD(Y\ctrh^\Xfl_\al)
  \ar@<3pt>[uuu]^{\O_Y\ocn_Y{-}}
  \ar@<-3pt>[rrrrrruuu]_{\quad\pi^*\E^\subbu\ocn_Y^\boL{-}}
  &&&&&& \sD(Y\ctrh^\Xlin_\al)
  \ar@<3pt>[uuu]^{\O_Y\ocn_Y{-}}
 }
\end{gathered}
\end{equation}
 Now the vertical functors are underived, and simply constructed by
applying the functors $\O_Y\ocn_Y{-}$ or $\fHom_Y(\O_Y,{-})$ termwise
to complexes from the respective classes.
 The diagonal functors still need to be derived.

 To begin with, let us say a few words about the upper horizontal line
of~\eqref{semico-semicontra-quadrality-diagram}, for the sake of
completeness of the exposition.
 The functor pointing rightwards,
$$
 \pi^*\E^\bu\ot_{\O_Y}{-}\,\:
 \sD(Y\qcoh_\Xfl)\lrarrow\sD(Y\qcoh^\Xinj),
$$
is underived from the outset.
 For any complex of $X$\+flat quasi-coherent sheaves $\G^\bu$ on $Y$,
the tensor product $\pi^*\E^\bu\ot_{\O_Y}\G^\bu$ is a complex of
$X$\+injective quasi-coherent sheaves on $Y$, in view of the projection
formula isomorphism $\pi_*(\pi^*\E^\bu\ot_{\O_Y}\G^\bu)\simeq
\E^\bu\ot_{\O_X}\pi_*\G^\bu$ (recall that the tensor product of a flat
quasi-coherent sheaf and an injective quasi-coherent sheaf on
a Noetherian scheme $X$ is injective).
 If the complex $\G^\bu$ is acyclic in $Y\qcoh_\Xfl$, then the complex
$\pi_*\G^\bu$ is acyclic in $X\qcoh_\fl$, and it follows that
the complex $\E^\bu\ot\O_X\pi_*\G^\bu$ is acyclic (in fact,
contractible) in $X\qcoh^\inj$.
 One easily concludes that the complex $\pi^*\E^\bu\ot_{\O_Y}\G^\bu$
is acyclic in $Y\qcoh^\Xinj$ in this case.

 The right derived functor pointing leftwards,
$$
 \boR\qHom_{Y\qc}(\pi^*\E^\bu,{-})\:
 \sD^\si_X(Y\qcoh)\lrarrow\sD(Y\qcoh_\Xfl),
$$
is constructed by applying the functor
$\qHom_{Y\qc}(\pi^*\E^\bu,{-})$ to complexes of $X$\+injective
quasi-coherent sheaves on~$Y$.
 Let us redraw the upper horizontal line
of~\eqref{semico-semicontra-quadrality-diagram} as follows:
\begin{equation} \label{semico-semicontra-upper-line-redrawn}
\begin{gathered}
 \xymatrix{
  \sD(Y\qcoh_\Xfl) \ar@<3pt>[rrrrrr]^{\pi^*\E^\subbu\ot_{\O_Y}{-}}   
  &&&&&& \sD(Y\qcoh^\Xinj)
  \ar@<3pt>[llllll]^{\qHom_{Y\qc}(\pi^*\E^\subbu,{-})},
 }
\end{gathered}
\end{equation}
where the equivalence $\sD^\si_X(Y\qcoh)\simeq\sD(Y\qcoh^\Xinj)$
is provided by Theorem~\ref{quasi-injective-over-base-semiderived}.
 Now the functor pointing leftwards is underived, too.

 For any complex of $X$\+injective quasi-coherent sheaves $\K^\bu$
on $Y$, the quasi-coherent internal $\qHom$ complex
$\qHom_{Y\qc}(\pi^*\E^\bu,\K^\bu)$ is a complex of $X$\+flat
quasi-coherent sheaves on $Y$, in view of the natural isomorphism
$\pi_*\qHom_{Y\qc}(\pi^*\E^\bu,\K^\bu)\simeq
\qHom_{X\qc}(\E^\bu,\pi_*\K^\bu)$.
 (Recall that the sheaf $\qHom_{X\qc}$ between two injective
quasi-coherent sheaves on a semi-separated Noetherian scheme $X$ is
flat and cotorsion; see the references in the proof of
Theorem~\ref{co-contra-dualizing}.)
 If the complex $\K^\bu$ is acyclic in $Y\qcoh^\Xinj$, then the complex
$\pi_*\K^\bu$ is contractible in $X\qcoh^\inj$, and it follows that
the complex $\qHom_{X\qc}(\E^\bu,\pi_*\K^\bu)$ is acyclic in
$X\qcoh_\fl$ (in fact, contractible in $X\qcoh_\fl^\cot$).
 One concludes that the complex $\qHom_{Y\qc}(\pi^*\E^\bu,\K^\bu)$
is acyclic in $Y\qcoh_\Xfl$ in this case.

 One can easily see that the two functors
in~\eqref{semico-semicontra-upper-line-redrawn} are adjoint to
each other.
 In order to show that this adjunction is an equivalence for
$\E^\bu=\D_X^\bu$, one observes that it agrees with the adjunction
in the upper horizontal line
of~\eqref{co-contra-quadrality-diagram}, in the sense that the direct
image functors $\pi_*$~\eqref{X-flat-to-flat-induced-direct-image}
and $\pi_*$~\eqref{semicoderived-to-X-underived-direct-image} take
the adjunction morphisms for the upper line
of~\eqref{semico-semicontra-quadrality-diagram} to the adjunction
morphisms for the upper line of~\eqref{co-contra-quadrality-diagram}.
 Since the adjunction morphisms for the upper line
of~\eqref{co-contra-quadrality-diagram} are isomorphisms, and
the direct image functors
$\pi_*$~\eqref{X-flat-to-flat-induced-direct-image}
and $\pi_*$~\eqref{semicoderived-to-X-underived-direct-image}
are conservative, it follows that the adjunction morphisms for
the upper line of~\eqref{semico-semicontra-quadrality-diagram} are
isomorphisms, too.
 A similar argument was utilized in the proof of
Theorem~\ref{semico-semicontra-diagonal}.

 Let us explain the construction of the functor
$\boR\Cohom_Y(\pi^*\E^\bu,{-})$ in the lower horizontal line
of~\eqref{semico-semicontra-quadrality-diagram}.
 The right derived functor
$$
 \boR\Cohom_Y(\pi^*\E^\bu,{-})\:\sD(Y\ctrh^\lin)\lrarrow
 \sD^\si_X(Y\ctrh)
$$
is constructed by applying the functor $\Cohom_Y(\pi^*\E^\bu,{-})$ to
complexes of locally injective (rather than merely $X$\+locally
injective) contraherent cosheaves on~$Y$.
 Let us redraw the lower horizontal line as follows:
\begin{equation} \label{semico-semicontra-lower-line-redrawn}
\begin{gathered}
 \xymatrix{
  \sD(Y\ctrh^\Xfl)
  &&&& {\frac{\Hot(Y\ctrh^\lin)}
 {\Acycl(Y\ctrh_\bT^\Xlin)\cap\Hot(Y\ctrh^\lin)}}
  \ar[llll]^-{\Cohom_Y(\pi^*\E^\subbu,{-})}.
 }
\end{gathered}
\end{equation}
 Here we use the equivalence of triangulated categories in
the left-hand side provided by
Theorem~\ref{contrah-projective-flat-over-base-semiderived}(b)
and the equivalence of triangulated categories in
the right-hand side provided by
Theorem~\ref{loc-injective-descript-of-derived-X-loc-injective}(a).
 Now the functor is underived.

 For any complex of locally injective contraherent cosheaves $\gJ^\bu$
on $Y$, the $\Cohom$ complex $\Cohom_Y(\pi^*\E^\bu,\gJ^\bu)$ is
a complex of $X$\+flat (in fact, $X$\+projective locally cotorsion)
contraherent cosheaves on $Y$, in view of the projection formula
isomorphism $\pi_!\Cohom_Y(\pi^*\E^\bu,\gJ^\bu)\simeq
\Cohom_X(\E^\bu,\pi_!\gJ^\bu)$ \eqref{lin-projection-II-cohom}.
 Recall that the $\Cohom$ cosheaf from injective quasi-coherent sheaf
to a locally injective contraherent cosheaf on a semi-separated
Noetherian scheme $X$ is a projective locally cotorsion contraherent
cosheaf, by Lemma~\ref{noetherian-contraherent-tensor} (it is even
easier to see that such a $\Cohom$ cosheaf is flat).
 If the complex $\gJ^\bu$ is acyclic in $Y\ctrh^\Xlin$, then
the complex $\pi_!\gJ^\bu$ is acyclic in $Y\ctrh^\lin$, and it follows
that the complex $\Cohom_X(\E^\bu,\pi_!\gJ^\bu)$ is acyclic in
$X\ctrh^\fl$ (in fact, contractible in $X\ctrh^\lct_\prj$).
 Using Lemma~\ref{acyclicity-in-lcth-criterion}, one concludes that
the complex $\Cohom_Y(\pi^*\E^\bu,\gJ^\bu)$ is acyclic in
$Y\ctrh^\Xfl$ (and even in $Y\ctrh^\lct_\Xprj$) in this case.

 The same $\Cohom$ projection formula~\eqref{lin-projection-II-cohom}
provides a commutative square diagram of triangulated functors
\begin{equation} \label{projection-to-X-Cohom-square}
\begin{gathered}
 \xymatrix{
  \sD^\si_X(Y\ctrh) \ar[d]^{\pi_!} &&& \sD(Y\ctrh^\Xlin)
  \ar[lll]^{\boR\Cohom_Y(\pi^*\E^\subbu,{-})} \ar[d]^{\pi_!} \\
  \sD^\ctr(X\ctrh) &&& \sD(X\ctrh^\lin)
  \ar[lll]^{\Cohom_X(\E^\subbu,{-})}
 }
\end{gathered}
\end{equation}
with the upper line given
by~\eqref{semico-semicontra-lower-line-redrawn}, the lower line
essentially taken from the big square
diagram~\eqref{co-contra-quadrality-diagram}, and the direct image
functors $\pi_!$~\eqref{semicontraderived-to-X-underived-direct-image}
and $\pi_!$~\eqref{X-lin-to-lin-induced-direct-image},
as promised in the formulation of the theorem.

 It remains to prove commutativity of the big square
diagram~\eqref{semico-semicontra-quadrality-diagram} (then it will
follow that the functor in the lower horizontal line is a triangulated
equivalence, since all the other functors are).
 Recall that the diagonal right derived functor
$\boR\fHom_Y(\pi^*\E^\bu,{-})$ was constructed in the proof of
Theorem~\ref{semico-semicontra-diagonal} by applying the functor
$\fHom_Y(\pi^*\E^\bu,{-})$ to complexes of injective (and not merely
$X$\+injective) quasi-coherent sheaves on~$Y$.
 Similarly, the diagonal left derived functor
$\pi^*\E^\bu\ocn_Y^\boL{-}$ was constructed in the same proof by
applying the functor $\pi^*\E^\bu\ocn_Y{-}$ to complexes of
antilocally flat contraherent cosheaves on~$Y$.

 Let us check commutativity of the interior upper triangle
in~\eqref{semico-semicontra-quadrality-diagram}.
 This commutativity holds for any finite complex of injective
quasi-coherent sheaves $\E^\bu$ on $X$ in the role of~$\D_X^\bu$.
 Starting with a complex of injective quasi-coherent sheaves $\J^\bu$
on $Y$, we apply the upper horizontal functor
$\qHom_{Y\qc}(\pi^*\E^\bu,{-})$ and obtain a complex of $X$\+flat
cotorsion quasi-coherent sheaves $\qHom_{Y\qc}(\pi^*\E^\bu,\J^\bu)$
on $Y$ (see Lemma~\ref{ext-qhom-qc}(c) and the discussion above
in this proof).
 Any cotorsion quasi-coherent sheaf is certainly contraadjusted; so
such a complex is adjusted to the derived functor
$\boR\fHom_Y(\O_Y,{-})$ in the leftmost column
of~\eqref{semico-semicontra-quadrality-diagram}, as per
the diagram~\eqref{semico-semicontra-big-diagram-partly-redrawn}.
 So it remains to observe the commutativity of the triangular diagram
of underived functors
\begin{equation} \label{semico-semicontra-interior-upper-triangle}
\begin{gathered}
 \xymatrix{
  \sD(Y\qcoh_\Xfl^\cta)
  \ar[ddd]^{\fHom_Y(\O_Y,{-})}
  &&&&&& \Hot(Y\qcoh^\inj)
  \ar[llllll]^{\qHom_{Y\qc}(\pi^*\E^\subbu,{-})}
  \ar[llllllddd]_{\fHom_Y(\pi^*\E^\subbu,{-})\quad}
  \\ \\ \\
  \sD(Y\ctrh^\Xfl)
 }
\end{gathered}
\end{equation}
which holds due to the natural isomorphism~\eqref{flat-inj-fhom-qhom}.

 Similarly, let us check that the exterior upper triangle
in~\eqref{semico-semicontra-quadrality-diagram} is commutative for
any finite complex of injective quasi-coherent sheaves $\E^\bu$ on $X$
in the role of~$\D_X^\bu$.
 Starting with a complex of antilocally flat contraherent cosheaves
$\gF^\bu$ on $Y$ (notice that antilocally flat contraherent cosheaves
are certainly antilocal), we apply the vertical functor in the leftmost
column and obtain a complex of flat contraadjusted quasi-coherent
sheaves $\O_Y^\bu\ocn_Y\gF^\bu$ on $Y$ (by
Lemma~\ref{cta-clp-restricts-to-prj-clf}(c)).
 The upper horizontal functor $\pi^*\E^\bu\ot_{\O_Y}{-}$ need not be
derived, and can be applied directly to any complex of $X$\+flat
quasi-coherent sheaves on~$Y$.
 It remains to observe the commutativity of the triangular diagram of
underived functors
\begin{equation} \label{semico-semicontra-exterior-upper-triangle}
\begin{gathered}
 \xymatrix{
  \sD(Y\qcoh_\Xfl) \ar[rrrrrr]^{\pi^*\E^\subbu\ot_{\O_Y}{-}}   
  &&&&&& \sD(Y\qcoh^\Xinj)
  \\ \\ \\
  \sD(Y\ctrh_\alf)
  \ar[uuu]^{\O_Y\ocn_Y{-}}
  \ar[rrrrrruuu]_{\quad\pi^*\E^\subbu\ocn_Y{-}}
 }
\end{gathered}
\end{equation}
which holds due to the natural
isomorphism~\eqref{tensor-contratensor-assoc}.

 Finally, it remains to check commutativity of the exterior lower
triangle in~\eqref{semico-semicontra-quadrality-diagram}.
 Once again, this commutativity holds for any finite complex of
injective quasi-coherent sheaves $\E^\bu$ on $X$ in the role
of~$\D_X^\bu$.
 Starting with a complex of injective quasi-coherent sheaves $\J^\bu$
on $Y$ (notice that injective quasi-coherent sheaves are certainly
contraadjusted), we apply the vertical functor in the rightmost column
and obtain a complex of antilocal locally injective contraherent
cosheaves $\fHom_Y(\O_Y,\J^\bu)$ on $Y$ (by
Lemma~\ref{cta-clp-restricts-to-cot-inj}(b)).
 Such a complex is adjusted to the derived functor
$\boR\Cohom_Y(\pi^*\E^\bu,{-})$ in the lower horizontal line
of~\eqref{semico-semicontra-quadrality-diagram}.
 So it remains to observe the commutativity of the triangular
diagram of underived functors
\begin{equation} \label{semico-semicontra-exterior-lower-triangle}
\begin{gathered}
 \xymatrix{
  &&&&&& \Hot(Y\qcoh^\inj)
  \ar[llllllddd]_{\fHom_Y(\pi^*\E^\subbu,{-})\quad}
  \ar[ddd]^{\fHom_Y(\O_Y,{-})} \\ \\ \\
  \sD(Y\ctrh^\Xfl)
  &&&&&& \sD(Y\ctrh^\lin)
  \ar[llllll]^{\Cohom_Y(\pi^*\E^\subbu,{-})}
 }
\end{gathered}
\end{equation}
which holds due to the natural isomorphism~\eqref{flat-fhom-cohom-inj}.
\end{proof}

\appendix

\Section{Derived Categories of Exact Categories and~Resolutions}
\label{exact-derived}

 In this appendix we recall and review some general results about
the derived categories of the first and the second kind for abstract
exact categories and their full subcategories, in presence of finite
or infinite resolutions.
 There is nothing essentially new here.
 Two or three most difficult arguments are omitted or only briefly
sketched with references to the author's previous works containing
elaborated proofs of similar results in different (but more concrete)
settings given in place of the details.

 Note that the present exposition does \emph{not} attain the maximal
natural generality for many results considered here.
 For most assertions concerning derived categories of the second kind,
the full generality is that of exact DG\+categories~\cite[Section~3.2
and Remark~3.5]{Pkoszul}, for which we refer the reader to
the preprint~\cite{Pedg}.

\subsection{Derived categories of the second kind}
\label{derived-second-kind}
 We suggest the survey paper~\cite{Bueh} as the general reference source
on exact categories.
 The papers~\cite[Appendix~A]{Kel0} and~\cite[Appendix~A]{Partin} can be
used as supplementary sources.

 Let $\sE$ be an exact category. 
 The homotopy categories of (finite, bounded above, bounded below,
and unbounded) complexes over $\sE$ will be denoted by
$\Hot^\b(\sE)$, \ $\Hot^-(\sE)$, \ $\Hot^+(\sE)$, and $\Hot(\sE)$,
respectively.
 For the definitions of the conventional derived categories
(of the first kind) $\sD^\b(\sE)$, \ $\sD^-(\sE)$, \ $\sD^+(\sE)$,
and $\sD(\sE)$ we refer to~\cite{N-e,Kel,Bueh}
and~\cite[Section~A.7]{Partin}.

 Let us recall one piece of terminology suggested
in~\cite[Section~2.1]{Psemi} and relevant for non-idempotent-complete
exact categories~$\sE$.
 We call a complex in $\sE$ \emph{exact} if it is composed of
(admissible) short exact sequences.
 A (bounded or unbounded) complex is said to be \emph{acyclic} if it is
homotopy equivalent to an exact complex, or equivalently, if it is
a direct summand of an exact complex.
 A morphism of complexes in an exact category is called
a \emph{quasi-isomorphism} if its cone is acyclic.

 Here are the definitions of the derived categories of the second
kind~\cite{Psemi,Pkoszul,PP2,EP}.
 Notice that all derived categories of the second kind in this appendix
are understood \emph{in the sense of Positselski} (rather than in
the sense of Becker~\cite{Jorg,K-st,N-f,Bec}).
 For definitions of the coderived and contraderived categories in
the sense of Becker, see Section~\ref{becker-subsect}.
 We refer the reader to the paper~\cite[Remark~9.2]{PS4}
and the survey~\cite[Section~7]{Pksurv} for a historical and
philosophical discussion.

 An (unbounded) complex over $\sE$ is said to be \emph{absolutely
acyclic} if it belongs to the minimal thick subcategory of $\Hot(\sE)$
containing all the total complexes of short exact sequences
of complexes over~$\sE$.
 Here a short exact sequence $0\rarrow {}'\.\!K^\bu\rarrow K^\bu\rarrow
{}''\!\.K^\bu\rarrow0$ of complexes over $\sE$ is viewed as a bicomplex
with three rows and totalized as such.
 So a complex is absolutely acyclic if and only if it is a homotopy
direct summand of a complex obtained from totalizations of short exact
sequences of complexes using the operation of passage to the cone
of a closed morphism of complexes repeatedly.
 (Here a ``homotopy direct summand'' means a direct summand in
the homotopy category.)
 The \emph{absolute derived category} $\sD^\abs(\sE)$ of the exact
category $\sE$ is defined as the quotient category of the homotopy
category $\Hot(\sE)$ by the thick subcategory of absolutely
acyclic complexes.

 Similarly, a bounded above (respectively, below) complex over $\sE$ is
called absolutely acyclic if it belongs to the minimal thick subcategory
of $\Hot^-(\sE)$ (resp., $\Hot^+(\sE)$) containing all the total
complexes of short exact sequences of bounded above (resp., below)
complexes over~$\sE$.
 We will see below that a bounded above (resp., below) complex over
$\sE$ is absolutely acyclic if and only if it is absolutely acyclic
as an unbounded complex, so there is no ambiguity in our terminology.
 The bounded above (resp., below) absolute derived category of $\sE$
is defined as the quotient category of $\Hot^-(\sE)$  (resp.,
$\Hot^+(\sE)$) by the thick subcategory of absolutely acyclic complexes
and denoted by $\sD^{\abs-}(\sE)$ (resp., $\sD^{\abs+}(\sE)$).

 We do not define the ``absolute derived category of finite complexes
over~$\sE$'', as it would not be any different from the conventional
bounded derived category $\sD^\b(\sE)$.
 Indeed, any (bounded or unbounded) absolutely acyclic complex is
acyclic; and any finite acyclic complex over an exact category is
absolutely acyclic, since it is homotopy equivalent to a complex
composed of short exact sequences.

 An exact category $\sE$ is said to have \emph{finite homological
dimension}~$\le d$ if $\Ext_\sE^n(X,Y)=0$ for all integers $n>d$
and all objects $X$, $Y\in\sE$.

\begin{lem} \label{psemi-remark21}
 In an exact category of finite homological dimension, every acyclic
complex is absolutely acyclic.
\end{lem}

\begin{proof}
 This is~\cite[Remark~2.1]{Psemi}.
\end{proof}

 For comparison with the results below, we recall that for any exact
category $\sE$ the natural functors $\sD^\b(\sE)\rarrow
\sD^\pm(\sE)\rarrow\sD(\sE)$ are all fully
faithful~\cite[Lemma~11.7]{Kel}, \cite[Corollary~A.11]{Partin}.

\begin{lem}  \label{b-abs-plus-minus-fully-faithful}
 For any exact category\/ $\sE$, the functors\/
$\sD^\b(\sE)\rarrow\sD^{\abs-}(\sE)\rarrow\sD^\abs(\sE)$ and\/
$\sD^\b(\sE)\rarrow\sD^{\abs+}(\sE)\rarrow\sD^\abs(\sE)$ induced
by the natural embeddings of the categories of bounded complexes
into those of unbounded ones are all fully faithful.
\end{lem}

\begin{proof}
 We will show that any morphism in $\Hot(\sE)$ (in an appropriate
direction) between a complex bounded in a particular way and
a complex absolutely acyclic with respect to the class of complexes
unbounded in that particular way factorizes through a complex
absolutely acyclic with respect to the class of correspondingly
bounded complexes.
 For this purpose, it suffices to demonstrate that any absolutely acyclic
complex can be presented as a termwise stabiling filtered inductive (or
projective) limit of complexes absolutely acyclic with respect to
the class of complexes bounded from a particular side.

 Indeed, any short exact sequence of complexes over $\sE$ is
the inductive limit of short exact sequences of their subcomplexes of
silly filtration, which are bounded below.
 One easily concludes that any absolutely acyclic complex is
a homotopy direct summand of a termwise stabilizing inductive limit of
complexes absolutely acyclic with respect to the class of complexes
bounded below, and any absolutely acyclic complex bounded above is
a homotopy direct summand of a termwise stabilizing inductive limit of
finite acyclic complexes.
 This proves that the functors $\sD^\b(\sE)\rarrow\sD^{\abs-}(\sE)$
and $\sD^{\abs+}(\sE)\rarrow\sD^\abs(\sE)$ are fully faithful.

 On the other hand, any absolutely acyclic complex bounded below, being,
by implication, an acyclic complex bounded below, is homotopy equivalent
to an exact complex bounded below.
 The latter is the inductive limit of its subcomplexes of canonical
filtration, which are finite acyclic complexes.
 This shows that the functor $\sD^\b(\sE)\rarrow\sD^{\abs+}(\sE)$ is
fully faithful, too.
 Finally, to prove that the functor $\sD^{\abs-}(\sE)\rarrow
\sD^\abs(\sE)$ is fully faithful, one presents any absolutely
acyclic complex as a homotopy direct summand of a termwise stabilizing
projective limit of complexes absolutely acyclic with respect to
the class of complexes bounded above.
\end{proof}

 Assume that infinite direct sums exist and are exact functors in
the exact category~$\sE$, i.~e., the coproducts of short exact
sequences are short exact sequences.
 Then a complex over $\sE$ is called \emph{coacyclic} if it belongs to
the minimal triangulated subcategory of $\Hot(\sE)$ containing
the total complexes of short exact sequences of complexes over $\sE$
and closed under infinite direct sums.
 Clearly, any coacyclic complex is acyclic
(cf.~\cite[Remark~3.1]{Pphil}).
 The \emph{coderived category} $\sD^\co(\sE)$ of the exact category $\sE$
is defined as the quotient category of the homotopy category
$\Hot(\sE)$ by the thick subcategory of coacyclic complexes.

 Similarly, if the functors of infinite product are everywhere defined
and exact in the exact category $\sE$, one calls a complex over $\sE$
\emph{contraacyclic} if it belongs to the minimal triangulated
subcategory of $\Hot(\sE)$ containing the total complexes of short
exact sequences of complexes over $\sE$ and closed under infinite
products.
 Clearly, any contraacyclic complex is acyclic.
 The \emph{contraderived category} $\sD^\ctr(\sE)$ of the exact
category $\sE$ is the quotient category of $\Hot(\sE)$ by the thick
subcategory of contraacyclic complexes~\cite[Sections~2.1
and~4.1]{Psemi}.

\begin{lem} \label{co-contra-bounded-fully-faithful}
\textup{(a)} For any exact category\/ $\sE$ with exact functors of
infinite direct sum, the full subcategory of bounded below complexes
in\/ $\sD^\co(\sE)$ is equivalent to\/ $\sD^+(\sE)$. \par
\textup{(b)} For any exact category\/ $\sE$ with exact functors of
infinite product, the full subcategory of bounded above complexes
in\/ $\sD^\ctr(\sE)$ is equivalent to\/ $\sD^-(\sE)$.
\end{lem}

\begin{proof}
 By~\cite[Lemmas~2.1 and~4.1]{Psemi}, any bounded below acyclic
complex over $\sE$ is coacyclic and any bounded above acyclic
complex over $\sE$ is contraacyclic.
 Hence there are natural triangulated functors $\sD^+(\sE)\rarrow
\sD^\co(\sE)$ and $\sD^-(\sE)\rarrow\sD^\ctr(\sE)$ (in the respective
assumptions of parts~(a) and~(b)).
 It also follows that the subcomplexes and quotient complexes of
canonical filtration of any co/contraacyclic complex remain
co/contraacyclic.
 Hence any morphism in $\Hot(\sE)$ from a bounded above complex to
a co/contraacyclic complex factorizes through a bounded above
co/contraacyclic complex, and any morphism from a co/contraacyclic
complex to a bounded below complex factorizes through a bounded
below co/contraacyclic complex.
 Therefore, our triangulated functors are fully
faithful~\cite[Remark~4.1]{Psemi}.
\end{proof}

 Denote the full additive subcategory of injective objects in $\sE$
by $\sE^\inj\sub\sE$ and the full additive subcategory of
projective objects by $\sE_\prj\sub\sE$.

\begin{lem}  \label{homotopy-inj-proj-fully-faithful}
\textup{(a)} The triangulated functors\/ $\Hot^\b(\sE^\inj)\rarrow
\sD^\b(\sE)$, \ $\Hot^\pm(\sE^\inj)\rarrow\sD^{\abs\pm}(\sE)$, \
$\Hot(\sE^\inj)\rarrow\sD^\abs(\sE)$, and\/
$\Hot(\sE^\inj)\rarrow\sD^\co(\sE)$ are fully faithful. \par
\textup{(b)} The triangulated functors\/ $\Hot^\b(\sE_\prj)\rarrow
\sD^\b(\sE)$, \ $\Hot^\pm(\sE_\prj)\rarrow\sD^{\abs\pm}(\sE)$, \
$\Hot(\sE_\prj)\rarrow\sD^\abs(\sE)$, and\/
$\Hot(\sE_\prj)\rarrow\sD^\ctr(\sE)$ are fully faithful.
\end{lem}

\begin{proof}
 This is essentially a version of~\cite[Theorem~3.5]{Pkoszul}
and a particular case of~\cite[Remark~3.5]{Pkoszul}
or~\cite[Theorem~5.5]{Pedg}.
 For any total complex $A^\bu$ of a short exact sequence of complexes
over $\sE$ and any complex $J^\bu$ over $\sE^\inj$, the complex of
abelian groups $\Hom_\sE(A^\bu,J^\bu)$ is acyclic.
 Therefore, the same also holds for a complex $A^\bu$ that can be
obtained from such total complexes using the operations of cone
and infinite direct sum (irrespectively even of such operations being
everywhere defined or exact in~$\sE$).
 Similarly, for any total complex $B^\bu$ of a short exact sequence
of complexes over $\sE$ and any complex $P^\bu$ over $\sE_\prj$,
the complex of abelian groups $\Hom_\sE(P^\bu,B^\bu)$ is acyclic
(hence the same also holds for any complexes $B^\bu$ that can be
obtained from such total complexes using the operations of cone
and infinite product).
 This semiorthogonality implies the assertions of Lemma.
\end{proof}

\subsection{Fully faithful functors}  \label{fully-faithful-subsect}
 Let $\sE$ be an exact category and $\sF\sub\sE$ be a full additive
subcategory.
 One says that the full subcategory $\sF$ \emph{inherits an exact
structure} from the ambient exact category $\sE$ if the class of
all short sequences in $\sF$ that are (admissible) exact in $\sE$
is an exact category structure on~$\sE$.

 Any full additive subcategory closed under extensions inherits
an exact category structure; so does any full additive subcategory
closed under the cokernels of admissible monomorphisms and
the kernels of admissible epimorphisms~\cite[Section~A.5(3)]{Partin},
\cite[Section~2]{DS}, \cite[Section~4.5]{Pedg}.
 The inherited exact category structure on $\sF$ is also said to be
\emph{induced} by the exact structure on~$\sE$.
 A criterion characterizing full subcategories inheriting an exact
category structure can be found in~\cite[Theorem~2.6]{DS}
or~\cite[Lemma~4.21]{Pedg}.

 Let $\sE$ be an exact category and $\sF\sub\sE$ be a full subcategory.
 One says that the full subcategory $\sF$ is \emph{self-generating}
in~$\sE$ \,\cite[Section~1]{BHP}, \cite[Section~7]{PS6} if for any
admissible epimorphism $E\rarrow F$ in $\sE$ from an object $E\in\sE$
to an object $F\in\sF$ there exist an object $G\in\sF$, an admissible
epimorphism $G\rarrow F$ in $\sE$, and a morphism $G\rarrow E$ in $\sE$
such that the triangle $G\rarrow E\rarrow F$ is commutative.

 A self-generating full subcategory closed under kernels of admissible
epimorphisms inherits an exact category structure if and only if it is
closed under extensions~\cite[Lemma~7.1]{Pedg}.
 We will say that a full subcategory $\sF\sub\sE$ is
\emph{self-resolving}~\cite[Section~7.1]{Pedg} if it is self-generating,
closed under extensions, and closed under kernels of admissible
epimorphisms.

 Dually, a full subcategory $\sC\sub\sE$ is said to be
\emph{self-cogeneraing} if for any admissible monomorphism
$C\rarrow E$ in $\sE$ from an object $C\in\sC$ to an object $E\in\sE$
there exist an object $D\in\sC$, an admissible monomorphism
$C\rarrow D$ in $\sE$, and a morphism $E\rarrow D$ in $\sE$ such that
the triangle $C\rarrow E\rarrow D$ is commutative.
 A full subcategory $\sC\sub\sE$ is said to be \emph{self-coresolving}
if it is self-cogenerating, closed under extensions, and closed under
cokernels of admissible monomorphisms (cf.~\cite[Section~12]{Kel}).

 For the rest of this section, we consider an exact category $\sE$ and
a self-resolving full subcategory $\sF\sub\sE$, endowed with
the induced exact category structure.

\begin{prop}  \label{fully-faithful-prop}
 For any symbol\/ $\bst=\b$, $-$, $\abs+$, $\abs-$, $\ctr$, or\/~$\abs$,
the triangulated functor\/ $\sD^\st(\sF)\rarrow\sD^\st(\sE)$ induced by
the exact embedding functor\/ $\sF\rarrow\sE$ is fully faithful.
 When\/ $\bst=\ctr$, it is presumed here that the functors of
infinite product are everywhere defined and exact in the exact
category\/ $\sE$ and preserve the full subcategory $\sF\sub\sE$.
\end{prop}

\begin{proof}
 We use the notation $\Hot^\st(\sE)$ for the category $\Hot(\sE)$
if $\bst=\empt$, $\co$, $\ctr$, or~$\abs$, the category $\Hot^-(\sE)$
if $\bst=-$ or~$\abs-$, the category $\Hot^+(\sE)$ if $\bst=+$
or~$\abs+$, and the category $\Hot^\b(\sE)$ if $\bst=\b$.
 Let us call an object of $\Hot^\st(\sE)$ \,\emph{$\bst$\+acyclic}
if it is annihilated by the localization functor $\Hot^\st(\sE)\rarrow
\sD^\st(\sE)$.

 In view of Lemma~\ref{b-abs-plus-minus-fully-faithful}, it suffices
to consider the cases $\bst=-$, $\abs$, and~$\ctr$.
 In each case, we will show that any morphism in $\Hot^\st(\sE)$
from a complex $F^\bu\in\Hot^\st(\sF)$ into a $\bst$\+acyclic
complex $A^\bu\in\Hot^\st(\sE)$ factorizes through a complex
$G^\bu\in\Hot^\st(\sF)$ that is $\bst$\+acyclic as a complex
over~$\sF$.

 We start with the case $\bst=-$.
 Assume for simplicity that both complexes $F^\bu$ and $A^\bu$
are concentrated in the nonpositive cohomological degrees; $F^\bu$
is a complex over $\sF$ and $A^\bu$ is an exact complex over~$\sE$.
 Notice that the morphism $A^{-1}\rarrow A^0$ must be 
an admissible epimorphism in this case.
 Set $G^0=F^0$, and let $B^{-1}\in\sE$ denote the fibered product
of the objects $F^0$ and $A^{-1}$ over~$A^0$.
 Then there exists a unique morphism $A^{-2}\rarrow B^{-1}$ having
a zero composition with the morphism $B^{-1}\rarrow F^0$
and forming a commutative diagram with the morphisms $A^{-2}
\rarrow A^{-1}$ and $B^{-1}\rarrow A^{-1}$.

 One easily checks that the complex $\dotsb \rarrow A^{-3}\rarrow
A^{-2}\rarrow B^{-1}\rarrow F^0\rarrow 0$ is exact.
 Furthemore, there exists a unique morphism $F^{-1}\rarrow B^{-1}$
whose compositions with the morphisms $B^{-1}\rarrow F^0$
and $B^{-1}\rarrow A^{-1}$ are the differential $F^{-1}\rarrow F^0$
and the component $F^{-1}\rarrow A^{-1}$ of the morphism of
complexes $F^\bu\rarrow A^\bu$.
 We have factorized the latter morphism of complexes through
the above exact complex, whose degree-zero term belongs to~$\sF$.

 From this point on we proceed by induction on the homological degree.
 Suppose that our morphism of complexes $F^\bu\rarrow A^\bu$ has been
factorized through an exact complex $\dotsb\rarrow A^{-n-2}\rarrow
A^{-n-1}\rarrow B^{-n}\rarrow G^{-n+1}\rarrow \dotsb\rarrow G^0\rarrow0$,
which coincides with the complex $A^\bu$ in the degrees $-n-1$ and
below, and whose terms belong to $\sF$ in the degrees $-n+1$ and above.
$$
 \xymatrix{
  \dotsb \ar[r] & F^{-n-2} \ar[r]\ar[d] & F^{-n-1} \ar[r]\ar[d]
  & F^{-n} \ar[r]\ar[d] & F^{-n+1} \ar[r]\ar[d] & F^{-n+2} 
  \ar[r]\ar[d] & \dotsb \\
  \dotsb \ar[r] & A^{-n-2} \ar[r]\ar@{=}[d] & A^{-n-1} \ar[r]\ar@{=}[d]
  & B^{-n} \ar[r]\ar[d] & G^{-n+1} \ar[r]\ar[d] & G^{-n+2}
  \ar[r]\ar[d] &\dotsb \\
  \dotsb \ar[r] & A^{-n-2} \ar[r] & A^{-n-1} \ar[r] & A^{-n} \ar[r]
  & A^{-n+1} \ar[r] & A^{-n+2} \ar[r] & \dotsb
 }
$$
 Since the full subcategory $\sF\sub\sE$ is assumed to be closed under
kernels of admissible epimorphisms, the image $Z^{-n}$ of the morphism
$B^{-n}\rarrow G^{-n+1}$ belongs to~$\sF$.
 Let $G^{-n}\rarrow Z^{-n}$ be an admissible epimorphism in $\sF$
factorizable through the admissible epimorphism $B^{-n}\rarrow Z^{-n}$
in~$\sE$.

 Replacing, if necessary, the object $G^{-n}$ by the direct sum
$G^{-n}\oplus F^{-n}$, we can make the morphism $F^{-n}\rarrow B^{-n}$
factorizable through the morphism $G^{-n}\rarrow B^{-n}$.
 Let $H^{-n-1}$ and $Y^{-n-1}$ denote the kernels of the admissible
epimorphisms $G^{-n}\rarrow Z^{-n}$ and $B^{-n}\rarrow Z^{-n}$,
respectively; then there is a natural morphism
$H^{-n-1}\rarrow Y^{-n-1}$.
 Denote by $B^{-n-1}$ the fibered product of the latter morphism with
the admissible epimorphism $A^{-n-1}\rarrow Y^{-n-1}$.
 There is a unique morphism $A^{-n-2}\rarrow B^{-n-1}$ having
a zero composition with the morphism $B^{-n-1}\rarrow H^{-n-1}$ and
forming a commutative diagram with the morphisms $A^{-n-2}\rarrow
A^{-n-1}$ and $B^{-n-1}\rarrow A^{-n-1}$.
$$
 \xymatrix{
  & & & & F^{-n} \ar[d] \\
  & & B^{-n-1} \ar@{..>>}[r] \ar@{..>}[d]
  & H^{-n-1} \ar@{>->}[r] \ar[d] & G^{-n} \ar[d] \ar@{->>}[rd] \\
  A^{-n-3} \ar[r] & A^{-n-2} \ar[r] \ar@{-->}[ru]
  & A^{-n-1} \ar@{->>}[r] & Y^{-n-1} \ar@{>->}[r]
  & B^{-n} \ar@{->>}[r] & Z^{-n} \ar[r] & G^{-n+1}
 }
$$

 We have constructed an exact complex $\dotsb\rarrow A^{-n-3}\rarrow
A^{-n-2}\rarrow B^{-n-1}\rarrow G^{-n}\rarrow G^{-n+1}\rarrow\dotsb
\rarrow G^0\rarrow 0$, coinciding with our previous complex in
the degrees $-n+1$ and above and with the complex $A^\bu$ in
the degrees $-n-2$ and below, and having terms belonging to $\sF$
in the degrees $-n$ and above.
 The morphism from the complex $F^\bu$ into our previous intermediate
complex factorizes through the new one (since the composition $F^{-n-1}
\rarrow F^{-n}\rarrow G^{-n}$ factorizes uniquely through
the admissible monomorphism $H^{-n-1}\rarrow G^{-n}$, and then there
exists a unique morphism $F^{-n-1}\rarrow B^{-n-1}$ whose compositions
with the morphisms $B^{-n-1}\rarrow H^{-n-1}$ and $B^{-n-1}\rarrow
A^{-n-1}$ are equal to the required ones).
$$
 \xymatrix{
  F^{-n-2} \ar[rr] \ar[dd] && F^{-n-1} \ar[rr] \ar@{-->}[d]\ar@{-->}[rd]
  && F^{-n} \ar[rr] \ar[d] && F^{-n+1} \ar[dd] \\
  & &B^{-n-1} \ar[r] \ar[d]
  & H^{-n-1} \ar@{>->}[r] \ar[d] & G^{-n} \ar[d] \ar@{->>}[rd] \\
  A^{-n-2} \ar[rr] \ar[rru] && A^{-n-1} \ar@{->>}[r]
  & Y^{-n-1} \ar@{>->}[r]
  & B^{-n} \ar@{->>}[r] & Z^{-n} \ar[r] & G^{-n+1}
 }
$$

 We have performed one step of the induction procedure, arriving to
a new factorization of the morphism of complexes $F^\bu\rarrow A^\bu$
as on the diagram
$$
 \xymatrix{
  \dotsb \ar[r] & F^{-n-2} \ar[r]\ar[d] & F^{-n-1} \ar[r]\ar[d]
  & F^{-n} \ar[r]\ar[d] & F^{-n+1} \ar[r]\ar[d] & \dotsb \\
  \dotsb \ar[r] & A^{-n-2} \ar[r]\ar@{=}[d] & B^{-n-1} \ar[r]\ar[d]
  & G^{-n} \ar[r]\ar[d] & G^{-n+1} \ar[r]\ar@{=}[d] & \dotsb \\
  \dotsb \ar[r] & A^{-n-2} \ar[r]\ar@{=}[d] & A^{-n-1} \ar[r]\ar@{=}[d]
  & B^{-n} \ar[r]\ar[d] & G^{-n+1} \ar[r]\ar[d] & \dotsb \\
  \dotsb \ar[r] & A^{-n-2} \ar[r] & A^{-n-1} \ar[r] & A^{-n} \ar[r]
  & A^{-n+1} \ar[r] & \dotsb
 }
$$
 Continuing with this procedure ad infinitum provides the desired
exact complex $G^\bu$ over~$\sF$.

\medskip

 The proof in the case $\bst=\abs$ is similar to that
of~\cite[Proposition~1.5]{EP}, and the case $\bst=\ctr$
is provable along the lines of~\cite[Remark~1.5]{EP}
and~\cite[Theorem~4.2.1]{Pweak} (cf.\ the proofs of
Proposition~\ref{finite-homol-dim-fully-faithful} and
Theorem~\ref{co-induced-resolutions}).
 A sketch of this argument, adopted to the situation at hand,
is reproduced below.
 For the full generality, see~\cite[Theorem~7.9 and
Corollary~7.12]{Pedg}.

 One has to show that any morphism from a complex over $\sF$
to a complex absolutely acyclic (contraacyclic) over $\sE$ factorizes
through a complex absolutely acyclic (contraacyclic) over $\sF$ in
the homotopy category $\Hot(\sE)$.
 This is checked by induction in the transformation rules using which
one constructs arbitrary absolutely acyclic (contraacyclic) complexes
over $\sE$ from the total complexes of short exact sequences.

 Let $K^\bu\rarrow L^\bu\rarrow M^\bu\rarrow K^\bu[1]$ be
a distinguished triangle in $\Hot(\sE)$ such that every morphism
from a complex over $\sF$ to $K^\bu$ or $M^\bu$ factorizes through
a complex absolutely acyclic (contraacyclic) over~$\sF$ in $\Hot(\sE)$.
 Let us show that the complex $L^\bu$ has the same property.
 Let $F^\bu$ be a complex over $\sF$ and $F^\bu\rarrow L^\bu$ be
a morphism of complexes.
 Then the composition $F^\bu\rarrow L^\bu\rarrow M^\bu$ factorizes
through a complex absolutely acyclic (contraacyclic) over $\sF$ in
$\Hot(\sE)$.
 It follows that there exists a complex $G^\bu$ over $\sF$ and
a morphism of complexes $G^\bu\rarrow F^\bu$ with a cone absolutely
acyclic (contraacyclic) over $\sF$ such that the composition $G^\bu
\rarrow F^\bu\rarrow L^\bu\rarrow M^\bu$ vanishes in $\Hot(\sE)$.
 It follows that the composition $G^\bu\rarrow F^\bu\rarrow L^\bu$
factorizes through the morphism $K^\bu\rarrow L^\bu$ in $\Hot(\sE)$.
 Hence we obtain a morphism of complexes $G^\bu\rarrow K^\bu$, which
also factorizes through a complex absolutely acyclic (contraacyclic)
over $\sF$ in $\Hot(\sE)$.
 It follows that there exists a complex $H^\bu$ over $\sF$ and
a morphism of complexes $H^\bu\rarrow G^\bu$ with a cone absolutely
acyclic (contraacyclic) over $\sF$ such that the composition
$H^\bu\rarrow G^\bu\rarrow K^\bu$ vanishes in $\Hot(\sE)$.
 Thus the composition $H^\bu\rarrow G^\bu\rarrow F^\bu\rarrow L^\bu$
also vanishes in $\Hot(\sE)$.
 Now the composition $H^\bu\rarrow G^\bu\rarrow F^\bu$ has a cone
absolutely acyclic (contraacyclic) over $\sF$, and the morphism
$F^\bu\rarrow L^\bu$ factorizes through this cone in $\Hot(\sE)$.

 Let $A_\alpha^\bu$ be a family of complexes over $\sE$ such that
every morphism from a complex over $\sF$ to $A_\alpha^\bu$ factorizes
through a complex contraacyclic over $\sF$ in $\Hot(\sE)$.
 Let us show that the complex $\prod_\alpha A_\alpha^\bu$ has the same
property.
 Given a complex $F^\bu$ over $\sF$ and a morphism of complexes
$F^\bu\rarrow\prod_\alpha A_\alpha^\bu$, let $F^\bu\rarrow A_\alpha^\bu$
be the components of this morphism.
 Then there exist complexes $G_\alpha^\bu$ contraacyclic over $\sF$
such that the morphism $F^\bu\rarrow A_\alpha^\bu$ factorizes
through~$G_\alpha^\bu$.
 Hence the morphism $F^\bu\rarrow\prod_\alpha A_\alpha^\bu$ factorizes
through the complex $\prod_\alpha G_\alpha^\bu$, which is contraacyclic
over~$\sF$.
 We have shown that the class of all complexes over $\sE$ we
are interested in is closed under cones (in both the cases $\bst=\abs$
or~$\ctr$) and under infinite products (in the case $\bst=\ctr$).

 Finally, the case of a morphism from a complex over $\sF$ to
the total complex of a short exact sequence of complexes over $\sE$
is treated using the next two lemmas.

\begin{lem}  \label{lifted-contractible}
 Let $U^\bu\rarrow V^\bu\rarrow W^\bu$ be a short exact sequence of
complexes over an exact category\/ $\sE$ and $M^\bu$ be
its total complex.
 Then a morphism $N^\bu\rarrow M^\bu$ of complexes over\/~$\sE$ is
homotopic to zero whenever its component $N^\bu\rarrow W^\bu[-1]$,
which is a morphism of graded objects in\/ $\sE$, can be lifted to
a morphism of graded objects $N^\bu\rarrow V^\bu[-1]$. \qed
\end{lem}

\begin{lem}  \label{induced-liftable}
 Let $U^\bu\rarrow V^\bu\rarrow W^\bu$ be a short exact sequence of
complexes over an exact category\/ $\sE$ and $M^\bu$ be
its total complex.
 Let $Q$ be a graded object in\/ $\sE$ and $\widetilde Q^\bu$ be
the (contractible) complex over\/ $\sE$ freely generated by~$Q$.
 Then a morphism of complexes $\tilde q\:\widetilde Q^\bu\rarrow M^\bu$
has its component $\widetilde Q^\bu\rarrow W^\bu[-1]$ liftable
to a morphism of graded objects $\widetilde Q^\bu\rarrow V^\bu[-1]$
whenever its restriction $q\:Q\rarrow M^\bu$ to the graded subobject
$Q\sub\widetilde Q^\bu$ has the same property. \qed
\end{lem}

 In order to apply the lemmas, one needs to notice that, in our
assumptions on $\sE$ and $\sF$, for any admissible epimorphism
$V\rarrow W$ in $\sE$, any object $F\in\sF$, and any morphism
$F\rarrow W$ in $\sE$ there exist an object $Q\in\sF$, an admissible
epimorphism $Q\rarrow F$ in $\sF$, and a morphism $Q\rarrow V$ in
$\sE$ such that the square $Q\rarrow F$, \ $V\rarrow W$ is commutative.
 Otherwise, the argument is no different from the one
in~\cite{EP,Pweak}.

 Specifically, let $U^\bu\rarrow V^\bu\rarrow W^\bu$ be a short exact
sequence of complexes over $\sE$ and $M^\bu$ be its total complex.
 Given a complex $F^\bu$ over $\sF$ and a morphism of complexes
$F^\bu\rarrow M^\bu$, consider its component $F^\bu\rarrow W^\bu[-1]$.
 It is a morphism of graded objects in~$\sE$.
 Pick a graded object $Q$ in $\sF$ together with an admissible epimorphism
$Q\rarrow F^\bu$ and a morphism $Q\rarrow V^\bu[-1]$ such that
the square $Q\rarrow F^\bu$, \ $V^\bu[-1]\rarrow W^\bu[-1]$
is commutative.

 Let $\widetilde Q^\bu$ be the complex freely generated by $Q$ and
$\widetilde Q^\bu\rarrow F^\bu$ be the morphism of complexes whose
restriction to $Q$ is equal to the admissible epimorphism
$Q\rarrow F^\bu$.
 By Lemma~\ref{induced-liftable}, the component $\widetilde Q^\bu
\rarrow W^\bu[-1]$ of the composition $\widetilde Q^\bu\rarrow F^\bu
\rarrow M^\bu$ is liftable to a morphism of graded objects
$\widetilde Q^\bu\rarrow V^\bu[-1]$.

 Let $R^\bu$ denote the kernel of the admissible epimorphism of
complexes $\widetilde Q^\bu\rarrow F^\bu$ over~$\sF$, and let $G^\bu$
be the cone of the closed admissible monomorpism of complexes
$R^\bu\rarrow \widetilde Q^\bu$.
 Then there is a natural admissible epimorphism of complexes
$G^\bu\rarrow F^\bu$ with a cone absolutely acyclic over~$\sF$.
 It remains to show that the composition $G^\bu\rarrow F^\bu\rarrow
M^\bu$ is contractible.

 According to Lemma~\ref{lifted-contractible}, it suffices to check
that the component $G^\bu\rarrow W^\bu[-1]$ is liftable to
a morphism of graded objects $G^\bu\rarrow V^\bu[-1]$.
 However, the complex $G^\bu$ viewed as a graded object in $\sF$ is
isomorphic to the direct sum $\widetilde Q^\bu\oplus R^\bu[1]$.
 The composition $G^\bu\rarrow F^\bu\rarrow M^\bu$ factorizes
through the projection onto the first direct summand
$\widetilde Q^\bu$, where its component landing in $W^\bu[-1]$ is
liftable to $V^\bu[-1]$, as we have already seen.
 This observation finishes the proof. 
\end{proof}

\subsection{Infinite resolutions}
\label{infinite-resolutions-subsect}
 A class of objects $\sF$ in an exact category $\sE$ is said to be
\emph{generating} if every object of $\sE$ is the image of an admissible
epimorphism from an object belonging to~$\sF$.
 A full subcategory $\sF\sub\sE$ is said to be \emph{resolving} if
$\sF$ is generating, closed under extensions, and closed under kernels
of admissible epimorphisms~\cite[Section~2]{Sto0}.
 Clearly, any generating full subcategory is self-generating
(in the sense of Section~\ref{fully-faithful-subsect}), and any
resolving full subcategory is self-resolving.

 Dually, a class of objects $\sC\sub\sE$ is said to be
\emph{cogenerating} if every object of $\sE$ is an admissible subobject
of an object belonging to~$\sC$.
 A full subcategory $\sC\subset\sE$ is said to be \emph{coresolving} if
$\sC$ is cogenerating, closed under extensions, and closed under 
cokernels of admissible monomorphisms.

 Let $\sE$ be an exact category and $\sF\sub\sE$ be a resolving
full subcategory.
 We endow $\sF$ with the induced structure of an exact category.

\begin{prop} \label{infinite-resolutions}
\textup{(a)} The triangulated functor\/ $\sD^-(\sF)\rarrow\sD^-(\sE)$
induced by the exact embedding functor\/ $\sF\rarrow\sE$ is
an equivalence of triangulated categories. \par
\textup{(b)} If the infinite products are everywhere defined and exact
in the exact category\/ $\sE$, and preserve the full subcategory\/
$\sF\sub\sE$, then the triangulated functor\/
$\sD^\ctr(\sF)\rarrow\sD^\ctr(\sE)$ induced by the embedding\/
$\sF\rarrow\sE$ is an equivalence of categories.
\end{prop}

\begin{proof}
 Part~(a) is well-known; see~\cite[Lemma~I.4.6]{Har}
or~\cite[Proposition~13.2.2(i)]{KS}.
 Our proof of part~(a) is based on part~(b) of the following lemma,
which will be also used in the sequel.

\begin{lem} \label{bounded-complex-resolution}
\textup{(a)} For any finite complex $E^{-d}\rarrow\dotsb\rarrow E^0$
over\/ $\sE$ there exists a finite complex $F^{-d}\rarrow\dotsb
\rarrow F^0$ over\/ $\sF$ together with a morphism of complexes
$F^\bu\rarrow E^\bu$ over\/ $\sE$ such that the morphisms
$F^i\rarrow E^i$ are admissible epimorphisms in\/ $\sE$ and
the cocone (or equivalently, the termwise kernel) of the morphism
$F^\bu\rarrow E^\bu$ is quasi-isomorphic to an object of\/ $\sE$ placed
in the cohomological degree~$-d$. \par
\textup{(b)} For any bounded above complex\/ $\dotsb\rarrow E^{-2}
\rarrow E^{-1}\rarrow E^0$ over\/ $\sE$ there exists a bounded above
complex\/ $\dotsb\rarrow F^{-2}\rarrow F^{-1}\rarrow F^0$ over\/ $\sF$
together with a quasi-isomorphism of complexes $F^\bu\rarrow E^\bu$
over\/ $\sE$ such that the morphisms $F^i\rarrow E^i$ are admissible
epimorphisms in\/~$\sE$.
\end{lem}

\begin{proof}
 Pick an admissible epimorphism $F^0\rarrow E^0$ with $F^0\in\sF$ and
consider the fibered product $G^{-1}$ of the objects $E^{-1}$ and $F^0$
over $E^0$ in~$\sE$.
 Then there exists a unique morphism $E^{-2}\rarrow G^{-1}$ having
a zero composition with the morphism $G^{-1}\rarrow F^0$ and forming
a commutative diagram with the morphisms $E^{-2}\rarrow E^{-1}$ and
$G^{-1}\rarrow E^{-1}$.
 Continuing the construction, pick an admissible epimorphism
$F^{-1}\rarrow G^{-1}$ with $F^{-1}\in\sF$, consider the fibered
product $G^{-2}$ of $E^{-2}$ and $F^{-1}$ over $G^{-1}$, etc.
 In the case~(a), proceed in this way until the object $F^{-d}$ is
constructed; in the case~(b), proceed indefinitely.
 The desired assertions follow from the observation that natural
morphism between the complexes $G^{-d}\rarrow F^{-d+1}\rarrow\dotsb
\rarrow F^0$ and $E^{-d}\rarrow E^{-d+1}\rarrow\dotsb\rarrow E^0$
is a quasi-isomorphism for any $d\ge1$.
\end{proof}

\begin{lem} \label{pkoszul-lemma16}
 Let\/ $\sH$ be a triangulated category and\/ $\sF$, $\sG$, $\sX
\sub\sH$ be (strictly) full triangulated subcategories. \par
\textup{(a)} If for every object $H\in\sH$ there exists an object
$F\in\sF$ together with a morphism $F\rarrow H$ with a cone belonging
to\/ $\sX$, then the natural triangulated functor
$$
 \sF/(\sF\cap\sX)\lrarrow\sH/\sX
$$
is an equivalence of triangulated categories. \par
\textup{(b)} If for every object $H\in\sH$ there exists an object
$G\in\sG$ together with a morphism $H\rarrow G$ with a cone belonging
to\/ $\sX$, then the natural triangulated functor
$$
 \sG/(\sG\cap\sX)\lrarrow\sH/\sX
$$
is an equivalence of triangulated categories.
\end{lem}

\begin{proof}
 This is~\cite[Proposition~10.2.7]{KS} or~\cite[Lemma~1.6]{Pkoszul}.
\end{proof}

 In view of Lemmas~\ref{bounded-complex-resolution}(b)
and~\ref{pkoszul-lemma16}(a), in order to finish the proof of part~(a)
of Proposition it remains to show that any bounded above complex over
$\sF$ that is acyclic over $\sE$ is also acyclic over~$\sF$.
 This follows immediately from the condition that $\sF$ is closed
with respect to the passage to the kernels of admissible
epimorphisms in~$\sE$.

 In the situation of part~(b), the functor in question is fully
faithful by Proposition~\ref{fully-faithful-prop}(b).
 A construction of a morphism with contraacyclic cone onto
a given complex over $\sE$ from an appropriately chosen complex
over $\sF$ is presented below.
 For a generalization, see~\cite[Theorem~7.11 and Corollary~7.12]{Pedg}.

\begin{lem}  \label{second-kind-complex-resolution}
\textup{(a)} For any complex $E^\bu$ over\/ $\sE$, there exists
a complex $P^\bu$ over\/ $\sF$ together with a morphism of complexes\/
$P^\bu\rarrow E^\bu$ such that the morphism $P^i\rarrow E^i$ is
an admissible epimorphism for each $i\in\boZ$. \par
\textup{(b)} For any complex $E^\bu$ over\/ $\sE$, there exists
a bicomplex $P^\bu_\bu$ over\/ $\sF$ together with a morphism of
bicomplexes $P^\bu_\bu\rarrow E^\bu$ over $\sE$ such that the complexes
$P_j^\bu$ vanish for all $j<0$, while for each $i\in\boZ$ the complex\/
$\dotsb\rarrow P_2^i\rarrow P_1^i\rarrow P_0^i\rarrow E^i\rarrow 0$
is acyclic with respect to the exact category\/~$\sE$.
\end{lem}

\begin{proof}
 To prove part~(a), pick admissible epimorphisms $F^i\rarrow E^i$ onto
all the objects $E^i\in\sE$ from some objects $F^i\in\sF$.
 Then the contractible complex $P^\bu$ with the terms $P^i=F^i\oplus
F^{i-1}$ (that is the complex freely generated by the graded object $F^*$
over $\sF$) comes together with a natural morphism of complexes
$P^\bu\rarrow E^\bu$ with the desired property.
 Part~(b) is easily deduced from~(a) by passing to the termwise kernel
of the morphism of complexes $P_0^\bu=P^\bu\rarrow E^\bu$ and iterating
the construction.
\end{proof}

\begin{lem}  \label{telescope}
 Let\/ $\sA$ be an additive category with countable direct products.
 Let\/ $\dotsb\rarrow P^{\bu\bu}(2)\rarrow P^{\bu\bu}(1)\rarrow
P^{\bu\bu}(0)$ be a projective system of bicomplexes over\/~$\sA$.
 Suppose that for every pair of integers $i$, $j\in\boZ$ the projective
system\/ $\dotsb \rarrow P^{ij}(2)\allowbreak\rarrow P^{ij}(1)\rarrow
P^{ij}(0)$ stabilizes, and let $P^{ij}(\infty)$ denote
the corresponding limit.
 Then the total complex of the bicomplex $P^{\bu\bu}(\infty)$ constructed
by taking infinite products along the diagonals is homotopy equivalent to
a complex obtained from the total complexes of the bicomplexes
$P^{\bu\bu}(n)$ (constructed in the same way) by iterated application of
the operations of shift, cone, and countable product.
\end{lem}

\begin{proof}
 Denote by $T(n)$ and $T(\infty)$ the total complexes of, respectively,
the bicomplexes $P^{\bu\bu}(n)$ and $P^{\bu\bu}(\infty)$.
 Then the short sequence of telescope construction
$$\textstyle
 0\rarrow T(\infty)\rarrow\prod_nT(n)\rarrow\prod_nT(n)\rarrow0
$$
is a termwise split short exact sequence of complexes over~$\sA$.
 Indeed, at every term of the complexes the sequence decomposes into
a countable product of sequences corresponding to the projective
systems $P^{ij}(*)$ with fixed indices $i$, $j$ and their limits
$P^{ij}(\infty)$.
 It remains to notice that the telescope sequence of a stabilizing
projective system is split exact, and a product of split exact
sequences is split exact.
\end{proof}

 Returning to part~(b) of Proposition, given a complex $E^\bu$ over
$\sE$, one applies Lemma~\ref{second-kind-complex-resolution}(b)
to obtain a bicomplex $P^\bu_\bu$ over $\sF$ mapping onto~$E^\bu$.
 Let us show that the cone of the morphism onto $E^\bu$ from
the total complex $T^\bu$ constructed by taking infinite products along
the diagonals of the bicomplex $P_\bu^\bu$ is a contraacyclic complex
over~$\sE$.
 For this purpose, augment the bicomplex $P_\bu^\bu$ with the complex
$E^\bu$ and represent the resulting bicomplex as the termwise stabilizing
projective limit of its quotient bicomplexes of canonical filtration
with respect to the lower indices.
 The latter bicomplexes being finite exact sequences of complexes
over $\sE$, the assertion follows from Lemma~\ref{telescope}
(cf.~\cite[Lemma~4.1]{Psemi}).
\end{proof}

\subsection{Homotopy adjusted complexes}  \label{homotopy-adjusted}
 The following simple construction (cf.~\cite{Spal}) will be useful
for us when working with the conventional unbounded derived categories
in Chapter~\ref{derived-on-quasi-compact-sect} (see, specifically,
Section~\ref{homotopy-lin-subsect}).  {\hbadness=6300\par}

 Let $\sE$ be an exact category.
 If the functors of infinite direct sum exist and are exact in $\sE$,
we denote by $\sD(\sE)^\lh\sub\sD(\sE)$ the minimal full triangulated
subcategory in $\sD(\sE)$ containing the objects of $\sE$ and closed
under infinite direct sums.
 Similarly, if the functors of infinite product exist and are exact
in $\sE$, we denote by $\sD(\sE)^\rh$ the minimal full triangulated
subcategory of $\sD(\sE)$ containing the objects of $\sE$ and
closed under infinite products.
 It is not difficult to show (see the proof of
Proposition~\ref{spal-for-exact}) that, in the assumptions
of the respective definitions, $\sD^-(\sE)\sub\sD(\sE)^\lh$ and
$\sD^+(\sE)\sub\sD(\sE)^\rh$. 

 Let $\sF\sub\sE$ be a self-resolving full subcategory in an exact
category, as defined in Section~\ref{fully-faithful-subsect}.

\begin{lem}  \label{acyclic-factorization-lemma}
 Let $B^\bu$ be a bounded above complex over\/ $\sF$ and $C^\bu$ be
an acyclic complex over\/~$\sE$.
 Then any morphism of complexes $B^\bu\rarrow C^\bu$ over\/ $\sE$
factorizes through a bounded above acyclic complex $K^\bu$
over\/~$\sF$.
\end{lem}

\begin{proof}
 The canonical truncation of exact complexes over the exact category
$\sE$ allows to assume the complex $C^\bu$ to be bounded above.
 In this case, the assertion was established in the proof
of Proposition~\ref{fully-faithful-prop}(a).
 (For a different argument leading to a slightly weaker conclusion,
see~\cite[proof of Lemma~2.9]{EP}.)
\end{proof}

\begin{cor}  \label{hom-vanishing-cor}
 Let $B^\bu$ be a bounded above complex over\/ $\sF$ and $C^\bu$ be
a complex over\/ $\sF$ that is acyclic as a complex over\/~$\sE$.
 Then the group\/ $\Hom_{\sD(\sF)}(B^\bu,C^\bu)$ of morphisms in
the derived category\/ $\sD(\sF)$ of the exact category\/ $\sF$
vanishes.
\end{cor}

\begin{proof}
 Indeed, any morphism from $B^\bu$ to $C^\bu$ in $\sD(\sF)$ can be
represented as a fraction $B^\bu\rarrow {}'C^\bu\larrow C^\bu$,
where $B^\bu\rarrow {}'C^\bu$ is a morphism of complexes over~$\sF$
and ${}'C^\bu\rarrow C^\bu$ is a quasi-isomorphism of such complexes.
 Then the complex ${}'C^\bu$ is also acyclic over $\sE$, and it remains
to apply Lemma~\ref{acyclic-factorization-lemma} to the morphism
$B^\bu\rarrow {}'C^\bu$.
\end{proof}

 From now on we assume $\sF$ to be a resolving full subcategory
in $\sE$, as in Section~\ref{infinite-resolutions-subsect}.

\begin{prop}  \label{spal-for-exact}
 Suppose that the exact category\/ $\sE$ is actually abelian, and that
infinite direct sums are everywhere defined and exact in the category\/
$\sE$ and preserve the full exact subcategory\/ $\sF\sub\sE$.
 Then the composition of natural triangulated functors\/ $\sD(\sF)^\lh
\rarrow\sD(\sF)\rarrow\sD(\sE)$ is an equivalence of triangulated
categories.
\end{prop}

\begin{proof}
 We will show that any complex over\/ $\sE$ is the target of
a quasi-isomorphism with the source belonging to $\sD(\sF)^\lh$.
 By Lemma~\ref{pkoszul-lemma16}(a), it will follow, in particular, that
$\sD(\sE)$ is equivalent to the localization of $\sD(\sF)$ by the thick
subcategory of complexes over $\sF$ acyclic over~$\sE$.
 By Corollary~\ref{hom-vanishing-cor}, the latter subcategory is
semiorthogonal to $\sD(\sF)^\lh$, so the same construction of
a quasi-isomorphism with respect to $\sE$, specialized to complexes
with the terms in $\sF$, will also imply that these two full
subcategories form a semiorthogonal decomposition of $\sD(\sF)$.
 This would clearly suffice to prove the desired assertion.

 Let $C^\bu$ be a complex over~$\sE$.
 Consider all of its subcomplexes of canonical truncation, pick
a termwise epic quasi-isomorphism onto each of them from a bounded
above complex over~$\sF$ vanishing starting from the same cohomological
degree (as per Lemma~\ref{bounded-complex-resolution}(b)), and replace
the latter complex with its finite subcomplex of silly filtration with,
say, only two nonzero terms.
 Take the direct sum $B_0^\bu$ of all the obtained complexes over $\sF$
and consider the natural morphism of complexes $B_0^\bu\rarrow C^\bu$.
 This is a termwise epimorphism of complexes which also acts by
epimorphisms on all the objects of coboundaries, cocycles, and
cohomology.
 Next we apply the same construction to the kernel of this morphism of
complexes, etc.

 We have constructed an exact complex of complexes $\dotsb\rarrow B_2^\bu
\rarrow B_1^\bu\rarrow B_0^\bu\rarrow C^\bu\rarrow0$ which remains
exact after replacing all the complexes $B_i^\bu$ and $C^\bu$ with
their cohomology objects (taken in the abelian category~$\sE$).
 All the complexes $B_i^\bu$ belong to $\sD(\sF)^\lh$ by
the construction.
 It remains to show that the totalization of the bicomplex $B_\bu^\bu$
obtained by taking infinite direct sums along the diagonals also
belongs to $\sD(\sF)^\lh$ and maps quasi-isomorphically onto~$C^\bu$.

 The totalization of the bicomplex $B_\bu^\bu$ is a direct limit of
the totalizations of its subbicomplexes of silly filtration in
the lower indices $B_n^\bu\rarrow\dotsb\rarrow B_0^\bu$.
 By the dual version of Lemma~\ref{telescope}, the former assertion
follows.
 To prove the latter one, it suffices to apply the following result due
to Eilenberg and Moore~\cite[Theorem~7.4]{EM0} to the bicomplex
obtained by augmenting $B_\bu^\bu$ with~$C^\bu$.
\end{proof}

\begin{lem}  \label{eilenberg-moore}
 Let $\sA$ be an abelian category with exact functors of countable
direct sum, and let $D^\bu_\bu$ be a bicomplex over $\sA$ such that
the complexes $D_j^\bu$ vanish for all $j<0$, while the complexes
$D_\bu^i$ are acyclic for all $i\in\boZ$, as are the complexes
$H^i(D_\bu^\bu)$.
 Then the total complex of the bicomplex $D_\bu^\bu$ obtained by taking
infinite direct sums along the diagonals is acyclic.
\end{lem}

\begin{proof}
 Denote by $S^\bu(\infty)$ the totalization of the bicomplex $D_\bu^\bu$
and by $S^\bu(n)$ the totalizations of its subbicomplexes of silly
filtration $D_n^\bu\rarrow\dotsb\rarrow D_0^\bu$.
 Consider the telescope short exact sequence
$$ \textstyle
 0\rarrow \bigoplus_n S^\bu(n)\rarrow \bigoplus_n S^\bu(n)
 \rarrow S^\bu(\infty)\rarrow0
$$
and pass to the long exact sequence of cohomology associated with this
short exact sequence of complexes.
 The morphisms $\bigoplus_n H^i(S^\bu(n))\rarrow\bigoplus_n
H^i(S^\bu(n))$ in this long exact sequence are the differentials in
the two-term complexes computing the derived functor of inductive limit
$\varinjlim^*_nH^i(S^\bu(n))$.
 It is clear from the conditions on the bicomplex $D_\bu^\bu$ that
the morphisms of cohomology $H^i(S^\bu(n-1))\rarrow H^i(S^\bu(n))$
induced by the embeddings of complexes $S^\bu(n-1)\rarrow S^\bu(n)$ 
vanish.
 Hence the morphisms $\bigoplus_n H^i(S^\bu(n))\rarrow\bigoplus_n
H^i(S^\bu(n))$ are isomorphisms and $H^*(S^\bu(\infty))=\nobreak0$.
\end{proof}

\begin{exs} \label{homotopy-adjusted-complexes-of-modules-examples}
 Let $R$ be an associative ring and $\sE=R\modl$ be the abelian
category of left $R$\+modules.

 (1)~Let $\sF=R\modl_\prj$ be the full additive subcategory of
projective $R$\+modules (with the induced trivial exact
category structure).
 Then $\sD(\sF)=\Hot(\sF)$.
 The full subcategory $\sF$ is resolving and closed under infinite
direct sums in $\sE$, so Proposition~\ref{spal-for-exact} is
applicable and we obtain a triangulated equivalence
$\sD(\sF)^\lh\simeq\sD(\sE)$.

 A complex of left $R$\+modules $P^\bu$ is said to be \emph{homotopy
projective}~\cite{Spal} if, for every acyclic complex of left
$R$\+modules $B^\bu$, the complex of abelian groups
$\Hom_R(P^\bu,B^\bu)$ is acyclic.
 We refer to~\cite[Remark~6.4]{PS4} for a terminological discussion.

 Let us show that the full subcategory $\sD(\sF)^\lh\sub\sD(\sF)$
coincides with the full subcategory of homotopy projective complexes
of projective $R$\+modules in $\Hot(R\modl_\prj)$.
 Indeed, the class of homotopy projective complexes of projective
$R$\+modules is obviously a strictly full triangulated subcategory
closed under infinite direct sums in $\Hot(\sF)$, and all one-term
complexes corresponding to projective modules are homotopy projective.
 Hence all the objects of $\sD(\sF)^\lh$ are homotopy projective.

 Then the converse inclusion follows from the definition of a homotopy
projective complex and Proposition~\ref{spal-for-exact}.
 See, e.~g., \cite[Theorem~1.4(a)]{Pkoszul} and~\cite[Theorem~6.12]{PS4}
for generalizations.

 (2)~Dually, let $\sC=R\modl^\inj$ be the full subcategory of injective
$R$\+modules in~$\sE$.
 Then $\sD(\sC)=\Hot(\sC)$.
 The full subcategory $\sC$ is coresolving and closed under infinite
products in $\sE$, so the dual version of
Proposition~\ref{spal-for-exact} is applicable and we obtain
a triangulated equivalence $\sD(\sC)^\rh\simeq\sD(\sE)$.

 A complex of left $R$\+modules $J^\bu$ is said to be \emph{homotopy
injective}~\cite{Spal} if, for every acyclic complex of left
$R$\+modules $A^\bu$, the complex of abelian groups
$\Hom_R(A^\bu,J^\bu)$ is acyclic.
 Let us show that the full subcategory $\sD(\sC)^\rh\sub\sD(\sC)$
coincides with the full subcategory of homotopy injective complexes of
injective $R$\+modules in $\Hot(R\modl^\inj)$.
 Indeed, the class of homotopy injective complexes of injective
$R$\+modules is a strictly full triangulated subcategory closed under
infinite products in $\Hot(\sC)$, and all one-term complexes 
corresponding to injective modules are homotopy injective.
 Hence all the objects of $\sD(\sC)^\rh$ are homotopy injective.

 Then the converse inclusion follows from the definition of
a homotopy injective complex and the dual version of
Proposition~\ref{spal-for-exact}.
 See, e.~g., \cite[Theorem~1.5(a)]{Pkoszul} or
even~\cite[Theorem~8.9]{PS4} for generalizations.

 (3)~Now let $\sF=R\modl_\fl$ be the full subcategory of flat
$R$\+modules (with the induced exact category structure).
 Then $\sF$ is a resolving full subcategory closed under infinite
direct sums in $\sE$, and Proposition~\ref{spal-for-exact} provides
a triangulated equivalence $\sD(\sF)^\lh\simeq\sD(\sE)$.

 A complex of left $R$\+modules $F^\bu$ is said to be \emph{homotopy
flat}~\cite{Spal} if, for every acyclic complex of right $R$\+modules
$A^\bu$, the complex of abelian groups $A^\bu\ot_R F^\bu$ is acyclic.
 One can show that a complex of flat $R$\+modules is homotopy flat
if and only if, viewed as an object of the derived category
$\sD(\sF)=\sD(R\modl_\fl)$, it belongs to the full subcategory
$\sD(\sF)^\lh\sub\sD(\sF)$.
 For commutative rings $R$, this a particular case of the quasi-coherent
sheaf version provided by Theorem~\ref{homotopy-flat-thm};
the noncommutative case can be treated similarly.
\end{exs}

\begin{rem}
 In particular, it follows from Proposition~\ref{spal-for-exact} that
$\sD(\sE)=\sD(\sE)^\lh$ for any abelian category $\sE$ with exact
functors of infinite direct sum.
 On the other hand, it is instructive to look into the situation of
Example~\ref{homotopy-adjusted-complexes-of-modules-examples}(1)
with $\sE=R\modl$ and $\sF=R\modl_\prj$.
 In this case, we have $\sD(\sF)=\Hot(\sF)\not\simeq\sD(\sE)$, while
$\sD(\sF)^\lh\subneq\sD(\sF)$ is the full subcategory of 
homotopy projective complexes in $\Hot(R\modl_\prj)$.
 Hence one can see (by considering $\sE=\sF=R\modl_\prj$) that
the assertion of Proposition~\ref{spal-for-exact} is not generally
true when the exact category $\sE$ is not abelian.
\end{rem}

\begin{cor} \label{induced-by-inclusion-left-adjoint}
 In the assumptions of Proposition~\textup{\ref{spal-for-exact}},
the fully faithful functor\/ $\sD(\sE)\simeq\sD(\sF)^\lh\rarrow\sD(\sF)$
is left adjoint to the triangulated functor\/ $\sD(\sF)\rarrow\sD(\sE)$
induced by the embedding of exact categories\/ $\sF\rarrow\sE$.
\end{cor}

\begin{proof}
 Clear from the proof of Proposition~\ref{spal-for-exact}.
\end{proof}

 It follows from Corollary~\ref{induced-by-inclusion-left-adjoint}
that the triangulated functor $\sD(\sF)\rarrow\sD(\sE)$ is
a Verdier quotient/localization functor.

 Keeping the assumptions of Proposition~\ref{spal-for-exact}, assume
additionally that the exact category $\sF$ has finite homological
dimension.
 Then the natural functor $\sD^\co(\sF)\rarrow\sD(\sF)$ is an equivalence
of triangulated categories by Lemma~\ref{psemi-remark21}.
 Consider the composition of triangulated functors $\sD(\sE)\simeq
\sD(\sF)^\lh\rarrow\sD(\sF)\simeq\sD^\co(\sF)\rarrow\sD^\co(\sE)$.
 The following result is a generalization of~\cite[Lemma~2.9]{EP}.

\begin{cor}  \label{homotopy-left-adjoint-cor}
 The functor\/ $\sD(\sE)\rarrow\sD^\co(\sE)$ so constructed is left
adjoint to the Verdier localization functor\/ $\sD^\co(\sE)\rarrow
\sD(\sE)$.
\end{cor}

\begin{proof}
 One has to show that $\Hom_{\sD^\co(\sE)}(B^\bu,C^\bu)=0$ for any
complex $B^\bu\in\sD(\sF)^\lh$ and any acyclic complex $C^\bu$
over~$\sE$.
 This vanishing easily follows from
Lemma~\ref{acyclic-factorization-lemma}.
 Indeed, any acyclic complex over $\sF$ is coacyclic
by Lemma~\ref{psemi-remark21}; and it remains to keep in mind that
any infinite direct sums in $\sD(\sF)\simeq\sD^\co(\sF)$ remain
infinite direct sums in $\sD^\co(\sE)$ (i.~e., the functor
$\sD^\co(\sF)\rarrow\sD^\co(\sE)$ preserves infinite direct sums)
by~\cite[Lemma~1.5]{BN} or~\cite[Lemma~3.2.10]{N-tr}.
\end{proof}

 Finally, let $\sG\sub\sF$ be two full subcategories in an abelian
category $\sE$, each resolving in $\sE$ and satisfying the assumptions
of Proposition~\ref{spal-for-exact}.
 Assume that the exact category $\sG$ has finite homological dimension.
 Consider the composition of triangulated functors
$\sD(\sE)\simeq\sD(\sG)^\lh\rarrow\sD(\sG)\simeq\sD^\co(\sG)\rarrow
\sD^\co(\sF)$.
 The next corollary is a straightforward generalization of
the previous one.

\begin{cor}  \label{homotopy-second-category-left-adjoint-cor}
 The functor\/ $\sD(\sE)\rarrow\sD^\co(\sF)$ constructed above is
left adjoint to the composition of Verdier localization functors\/
$\sD^\co(\sF)\rarrow\sD(\sF)\rarrow\sD(\sE)$.
\end{cor}

\begin{proof}
 The point is that $\Hom_{\sD^\co(\sF)}(B^\bu,C^\bu)=0$ for any
complex $B^\bu\in\sD(\sG)^\lh$ and any complex $C^\bu$ in $\sF$
acyclic over~$\sE$.
\end{proof}

 Our discussion of the full subcategory $\sD(\sF)^\lh\sub\sD(\sF)$
will be continued in Proposition~\ref{dg-F-acycl-C-cotorsion-pair}
in Appendix~\ref{cotorsion-pairs-appx}.

\subsection{Finite resolutions}
\label{finite-resolutions-subsect}
 Let $\sE$ be an exact category and $\sF\sub\sE$ be a resolving
full subcategory, as in Section~\ref{infinite-resolutions-subsect}.
 Assume additionally that the additive category $\sE$ is
``semi-saturated'' in the terminology of~\cite{N-e,Partin}, or
which is the same, ``weakly idempotent-complete'' in the terminology
of~\cite{Bueh,Sto-ICRA} (i.~e., it contains the kernels of its split
epimorphisms, or equivalently, the cokernels of its split monomorphisms).
 Then the additive category $\sF$ has the same property.
 The following results elaborate upon the ideas
of~\cite[Remark~2.1]{EP}.
 For another exposition, see~\cite[Section~2]{Sto0}.

 We will say that an object of the derived category $\sD^-(\sE)$ has
\emph{$\sF$\+resolution dimension not exceeding~$m$} if its isomorphism
class can be represented by a bounded above complex $F^\bu$ over $\sF$
such that $F^i=0$ for $i<-m$.
 By the definition, the full subcategory of objects of finite
$\sF$\+resolution dimension in $\sD^-(\sE)$ is the image of
the fully faithful triangulated functor $\sD^\b(\sF)\rarrow
\sD^-(\sF)\simeq\sD^-(\sE)$.

 Dually, let $\sC\sub\sE$ be a coresolving full subcategory.
 Then an object of the derived category $\sD^+(\sE)$ is said to have
\emph{$\sC$\+coresolution dimension not exceeding~$m$} if its
isomorphism class in the derived category can be represented by
a bounded below complex $C^\bu$ over $\sC$ such that $C^i=0$ for $i>m$.

\begin{lem}  \label{fdim-well-defined}
 If a bounded above complex $G^\bu$ over the resolving subcategory\/
$\sF$, viewed as an object of the derived category\/ $\sD^-(\sE)$,
has $\sF$\+resolution dimension not exceeding~$m$, then
the differential\/ $G^{-m-1}\rarrow G^{-m}$ has a cokernel\/
${}'G^{-m}$ in the additive category\/ $\sF$, and the complex\/
$\dotsb\rarrow G^{-m-1}\rarrow G^{-m}\rarrow {}'G^{-m}\rarrow0$
over\/ $\sF$ is acyclic.
 Consequently, the complex $G^\bu$ over\/ $\sF$ is quasi-isomorphic to
the finite complex\/ $0\rarrow {}'G^{-m}\rarrow G^{-m+1}\rarrow
G^{-m+2}\rarrow\dotsb\rarrow0$.
\end{lem}

\begin{proof}
 In view of the equivalence of categories $\sD^-(\sF)\simeq \sD^-(\sE)$
(see Proposition~\ref{infinite-resolutions}(a)), the assertion really
depends on the exact subcategory $\sF$ only.
 By the definition of the derived category, two complexes representing
isomorphic objects in it are connected by a pair of quasi-isomorphisms.
 Thus it suffices to consider two cases when there is a quasi-isomorphism
acting either in the direction $G^\bu\rarrow F^\bu$, or in the opposite
direction $F^\bu\rarrow G^\bu$ (where $F^\bu$ is a bounded above complex
over $\sF$ such that $F^i=0$ for $i<-m$).

 In the former case, acyclicity of the cone of the morphism
$G^\bu\rarrow F^\bu$ implies the existence of cokernels of
its differentials and the acyclicity of canonical truncations,
which provides the desired conclusion.
 Notice that any bounded above acyclic complex in a weakly
idempotent-complete exact category is exact.
 Furthermore, in the situation at hand the differential $G^{-m-1}
\rarrow G^{-m}$ in the complex $G^\bu$ is also a differential in
the complex $\cone(G^\bu\to F^\bu)$, up to a sign.

 In the latter case, from acyclicity of the cone of the morphism
$F^\bu\rarrow G^\bu$ one can similarly see that the morphism
$G^{-m-2}\rarrow G^{-m-1}\oplus F^{-m}$ with the vanishing component
$G^{-m-2}\rarrow F^{-m}$ has a cokernel, and it follows that
the morphism $G^{-m-2}\rarrow G^{-m-1}$ also does (since the category
$\sF$ is weakly idempotent-complete).
 Denoting the cokernel of the latter morphism by ${}'G^{-m-1}$, one
easily concludes that the complex $\dotsb\rarrow G^{-m-2}\rarrow
G^{-m-1}\rarrow {}'G^{-m-1}\rarrow0$ is acyclic, and it remains to show
that the morphism ${}'G^{-m-1}\rarrow G^{-m}$ is
an admissible monomorphism in the exact category~$\sF$.
 Indeed, the morphism ${}'G^{-m-1}\oplus F^{-m}\rarrow G^{-m}\oplus
F^{-m+1}$ is; and hence so is its composition with the embedding of
a direct summand ${}'G^{-m-1}\rarrow {}'G^{-m-1}\oplus F^{-m}$.
 The latter morphism has a vanishing component
${}'G^{-m-1}\rarrow F^{-m}$, etc.
\end{proof}

\begin{cor} \label{fdim-cor}
 Let $G^\bu$ be a finite complex over the exact category\/ $\sE$ such
that $G^i=0$ for $i<-m$ and $G^i\in\sF$ for $i>-m$.
 Assume that the object represented by $G^\bu$ in $\sD^-(\sE)$ has
$\sF$\+resolution dimension not exceeding~$m$.
 Then the object $G^{-m}$ also belongs to\/~$\sF$.
\end{cor}

\begin{proof}
 Replace the object $G^{-m}$ with its resolution by objects from $\sF$
and apply Lemma~\ref{fdim-well-defined}.
\end{proof}

\begin{lem} \label{fdim-triangle}
 Let $A^\bu\rarrow B^\bu\rarrow C^\bu\rarrow A^\bu[1]$ be
a distinguished triangle in the derived category\/ $\sD^-(\sE)$.
 In this context: \par
\textup{(a)} if the\/ $\sF$\+resolution dimensions of $A^\bu$ and
$C^\bu$ do not exceed~$m$, then the\/ $\sF$\+resolution dimension
of $B^\bu$ does not exceed~$m$; \par
\textup{(b)} if the\/ $\sF$\+resolution dimension of $B^\bu$ does not
exceed~$m$ and the\/ $\sF$\+resolution dimension of $C^\bu$ does not
exceed $m+1$, then the\/ $\sF$\+resolution dimension of $A^\bu$ does not
exceed~$m$; \par
\textup{(c)} if the\/ $\sF$\+resolution dimension of $B^\bu$ does not
exceed~$m$ and the\/ $\sF$\+resolution dimension of $A^\bu$ does not
exceed $m-1$, then the\/ $\sF$\+resolution dimension of $C^\bu$ does not
exceed~$m$.
\end{lem}

\begin{proof}
 It suffices to prove part~(c), as parts~(a) and~(b) can be obtained
by rotating the triangle.
 The morphism $A^\bu\rarrow B^\bu$ in $\sD^-(\sE)\simeq\sD^-(\sF)$
can be represented by a morphism of complexes
${}'\!\.F^\bu\rarrow{}''\!\.F^\bu$ in $\Hot^-(\sF)$.
 By Lemma~\ref{fdim-well-defined}, both complexes can be replaced
by their canonical truncations at the degree~$-m$.
 Obviously, there is the induced morphism between the complexes
truncated in this way, so we can simply assume that ${}'\!\.F^i=0
={}''\!\.F^i$ for $i<-m$.
 Moreover, the complex ${}'\!\.F^\bu$ could be truncated even one step
further, i.~e., the morphism ${}'\!\.F^{-m}\rarrow{}'\!\.F^{-m+1}$
is an admissible monomorphism in the exact category~$\sF$.
 From this one easily concludes that, for the cone $G^\bu$ of
the morphism of complexes ${}'\!\.F^\bu\rarrow{}''\!\.F^\bu$, one
has $G^i=0$ for $i<-m-1$, and the morphism $G^{-m-1}\rarrow G^{-m}$
is an admissible monomorphism in~$\sF$.
\end{proof}

 We say that an object $E\in\sE$ has \emph{$\sF$\+resolution dimension
not exceeding~$m$} \,\cite[Section~2]{Sto0} if the corresponding object
of the derived category $\sD^-(\sE)$ does.
 In other words, $E$ must have a resolution of the length
not exceeding~$m$ by objects of~$\sF$.
 Let us denote the $\sF$\+resolution dimension of an object $E$
by $\rsd_{\sF/\sE}E$.

 Dually, given a coresolving subcategory $\sC\sub\sE$, an object
$E\in\sE$ is said to have \emph{$\sC$\+coresolution dimension
not exceeding~$m$} if $E$ has a coresolution of the length
not exceeding~$m$ by objects of $\sC$ in~$\sE$.
 We denote the $\sC$\+coresolution dimension of an object $E$
by $\crd_{\sC/\sE}E$.

\begin{cor} \label{fdim-acyclic-cor}
 Assume that there exists an integer $m\ge0$ such that all the objects
of\/ $\sE$ have resolution dimensions\/~$\le m$ with respect to\/~$\sF$.
 Then any acyclic complex in\/ $\sE$ with the terms in\/ $\sF$ is
acyclic in\/~$\sF$.
\end{cor}

\begin{proof}
 For any acyclic complex $E^\bu$ in a weakly idempotent-complete
exact category\/ $\sE$, the complex $E^\bu\oplus E^\bu[1]$ is exact.
 Hence, in the situation at hand, it suffices to show that any exact
complex $C^\bu$ in $\sE$ with the terms in $\sF$ is exact in~$\sF$.
 Applying Corollary~\ref{fdim-cor} to the finite quotient complexes of
canonical and silly truncation of the complex $C^\bu$, we conclude
that its objects of cocycles in the exact category $\sE$ belong to
the subcategory $\sF\sub\sE$.
\end{proof}

\begin{cor}  \label{fdim-subcategory-cor}
 Let\/ $\sE'\sub\sE$ be a (strictly) full weakly idempotent-complete
additive subcategory with an induced exact category structure.
 Set\/ $\sF'=\sE'\cap\sF$, and assume that every object of\/ $\sE'$
is the image of an admissible epimorphism in the exact category\/
$\sE'$ acting from an object belonging to\/~$\sF'$.
 Then\/ $\sF'$ is a resolving full subcategory in\/ $\sE'$, and for
any object $E\in\sE'$ one has\/ $\rsd_{\sF/\sE}E =\rsd_{\sF'/\sE'}E$.
 The same applies to bounded above complexes $E^\bu$ in\/~$\sE'$.
\end{cor}

\begin{proof}
 The first assertion follows from the definitions, and the second one
from Corollary~\ref{fdim-cor}.
 To prove the final assertion, use Lemma~\ref{fdim-well-defined}.
\end{proof}

\begin{lem}  \label{fdim-properties}
 Let\/ $0\rarrow E'\rarrow E\rarrow E''\rarrow0$ be an (admissible)
short exact sequence in\/~$\sE$.
 Then \par
\textup{(a)} if\/ $\rsd_{\sF/\sE} E'\le m$ and\/
$\rsd_{\sF/\sE} E''\le m$, then\/ $\rsd_{\sF/\sE} E\le m$; \par
\textup{(b)} if\/ $\rsd_{\sF/\sE} E\le m$ and\/
$\rsd_{\sF/\sE} E''\le m+1$, then\/ $\rsd_{\sF/\sE} E'\le m$; \par
\textup{(c)} if\/ $\rsd_{\sF/\sE} E\le m$ and\/
$\rsd_{\sF/\sE} E'\le m-1$, then\/ $\rsd_{\sF/\sE} E''\le m$.
\end{lem}

\begin{proof}
 This is a particular case of Lemma~\ref{fdim-triangle}.
\end{proof}

\begin{cor} \label{fdim-resolution}
\textup{(a)} Let\/ $0\rarrow E_n\rarrow\dotsb\rarrow E_0\rarrow E
\rarrow 0$ be an exact sequence in\/~$\sE$.
 Then the\/ $\sF$\+resolution dimension\/ $\rsd_{\sF/\sE} E$ does
not exceed the supremum of the expressions\/ $\rsd_{\sF/\sE} E_i+i$
over\/ $0\le i\le n$ (where we set\/ $\rsd_{\sF/\sE}0=-1$). \par
\textup{(b)} Let\/ $0\rarrow E\rarrow E^0\rarrow\dotsb\rarrow E^n
\rarrow 0$ be an exact sequence in\/~$\sE$.
 Then the\/ $\sF$\+resolution dimension\/ $\rsd_{\sF/\sE} E$ does
not exceed the supremum of the expressions\/ $\rsd_{\sF/\sE} E^i-i$
over\/ $0\le i\le n$.
\end{cor}

\begin{proof}
 Part~(a) follows by induction from Lemma~\ref{fdim-properties}(c),
and part~(b) similarly follows from Lemma~\ref{fdim-properties}(b).
\end{proof}

\begin{prop}  \label{finite-resolutions}
 Suppose that the\/ $\sF$\+resolution dimension of all objects
$E\in\sE$ does not exceed a fixed constant~$n$.
 Then the triangulated functor\/ $\sD^\st(\sF)\rarrow\sD^\st(\sE)$
induced by the exact embedding functor\/ $\sF\rarrow\sE$ is
an equivalence of triangulated categories for any symbol\/ $\bst=\b$,
$+$, $-$, $\empt$, $\abs+$, $\abs-$, $\co$, $\ctr$, or\/~$\abs$.

 When\/ $\bst=\co$ (respectively, $\bst=\ctr$), it is presumed here
that the functors of infinite direct sum (resp., infinite product) are
everywhere defined and exact in the category\/ $\sE$ and preserve
the full subcategory\/ $\sF\sub\sE$.
\end{prop}

\begin{proof}
 The cases $\bst=-$ or $\ctr$ were considered in
Proposition~\ref{infinite-resolutions} (and hold in its
weaker assumptions).
 They can be also treated together with the other cases, as it is
explained below.

 Using the construction of Lemma~\ref{second-kind-complex-resolution}
and taking into account Corollary~\ref{fdim-cor}, one produces for any
$\bst$\+bounded complex $E^\bu$ over $\sE$ its finite resolution
$P_\bu^\bu$ of length~$n$ (in the lower indices) by $\bst$\+bounded (in
the upper indices) complexes over~$\sF$.
 The total complex of $P_\bu^\bu$ maps by a morphism with an absolutely
acyclic cone over the exact category $\sE$ onto the complex $E^\bu$.
 By Lemma~\ref{pkoszul-lemma16}(a), it remains to show that any
$\bst$\+bounded complex $C^\bu$ over $\sF$ that is $\bst$\+acyclic as
a complex over $\sE$ is also $\bst$\+acyclic as a complex over~$\sF$.

 In the cases $\bst=\b$, $+$, $-$, or~$\empt$, we can refer to
Corollary~\ref{fdim-acyclic-cor}.
 In view of Lemma~\ref{b-abs-plus-minus-fully-faithful}, it remains
to consider the cases $\bst=\co$, $\ctr$, or~$\abs$.
 The related argument is similar to that in~\cite[Theorem~3.2]{PP2}
and~\cite[Theorem~1.4]{EP}, and goes back to the proof
of~\cite[Theorem~7.2.2]{Psemi}.
 For the full generality, see~\cite[Theorem~6.6 and
Corollary~6.20]{Pedg}.

 For every $0\le d\le n$, denote by $\sF_d\sub\sE$ the full
subcategory formed by all the objects of $\sF$\+resolution
dimension~$\le d$; so $\sF_0=\sF$ and $\sF_n=\sE$.
 By Lemma~\ref{fdim-properties}, the full subcategory $\sF_d$ is
closed under extensions and the kernels of admissible epimorphisms
in $\sE$; so it is an exact subcategory.
 For brevity, we will use the term ``an $\sF_d$\+complex'' instead
of ``a complex over $\sF_d$'' and ``an $\sF_d$\+$\bst$\+acyclic
complex'' instead of ``a $\bst$\+acyclic complex over $\sF_d$''.
 Our aim is to show that every $\sF_n$\+$\bst$\+acyclic
$\sF_0$\+complex is $\sF_0$\+$\bst$\+acyclic.
 For this purpose, we will prove that every $\sF_d$\+$\bst$\+acyclic
$\sF_{d-1}$\+complex is $\sF_{d-1}$\+$\bst$\+acyclic; the desired
assertion will then follow by induction.

 It suffices to construct for any $\sF_d$\+$\bst$\+acyclic complex
$M^\bu$ a short exact sequence of complexes $0\rarrow K^\bu\rarrow L^\bu
\rarrow M^\bu\rarrow0$ in $\sF_d$ such that the complexes
$K^\bu$ and $L^\bu$ are $\sF_{d-1}$\+$\bst$\+acyclic.
 Then if $M^\bu$ is an $\sF_{d-1}$\+complex, it would follow that
both the cone of the morphism $K^\bu\rarrow L^\bu$ and the total
complex of the short exact sequence $K^\bu\rarrow L^\bu\rarrow M^\bu$
are $\sF_{d-1}$\+$\bst$\+acyclic, so $M^\bu$ also is.
 The construction is based on four lemmas similar to those
in~\cite{PP2} and~\cite{EP}.

\begin{lem}
 Let $M^\bu$ be the total complex of a short exact sequence of
complexes ${}'\!\.M^\bu\rarrow {}''\!\.M^\bu\rarrow {}'''\!\.M^\bu$
over\/~$\sF_d$.
 Then there exists an admissible epimorphism onto $M^\bu$ from
a contractible complex $L^\bu$ over\/ $\sF$ with
an\/ $\sF_{d-1}$\+$\bst$\+acyclic kernel~$K^\bu$.  \qed
\end{lem}

\begin{lem}
\textup{(a)} Let ${}'\!\.K^\bu\rarrow{}\!\.'L^\bu\rarrow{}'\!\.N^\bu$
and ${}''\!\.K^\bu\rarrow{}''\!\.L^\bu\rarrow{}''\!\.N^\bu$ be short
exact sequences of complexes in the exact category\/ $\sE$ such that
the complexes ${}'\!\.K^\bu$, ${}'\!\.L^\bu$, ${}''\!\.K^\bu$,
${}''\!\.L^\bu$ are\/ $\sF_{d-1}$\+$\bst$\+acyclic, and let
${}'\!\.N^\bu\rarrow{}''\!\.N^\bu$ be a morphism of complexes
over\/~$\sE$.
 Then there exists a short exact sequence $K^\bu\rarrow L^\bu\rarrow
N^\bu$ of complexes in\/ $\sE$ such that $N^\bu\simeq\cone({}'\!\.N^\bu
\to{}''\!\.N^\bu)$ and the complexes $K^\bu$ and $L^\bu$ are\/
$\sF_{d-1}$\+$\bst$\+acyclic. \par
\textup{(b)} In the situation of\/~\textup{(a)}, assume that
the morphism ${}'\!\.N^\bu\rarrow{}''\!\.N^\bu$ is an admissible
monomorphism of complexes in the exact category\/ $\sE$ whose cokernel
$N_0^\bu$ is a complex over\/~$\sF_d$.
 Then there exists a short exact sequence $K_0^\bu\rarrow L_0^\bu\rarrow
N_0^\bu$ of complexes in\/ $\sF_{d-1}$ such that the complexes $K_0^\bu$
and $L_0^\bu$ are\/ $\sF_{d-1}$\+$\bst$\+acyclic.
\end{lem}

\begin{proof}
 The proof is similar to that in~\cite{PP2}.
\end{proof}

\begin{lem}
 For any contractible complex $M^\bu$ over\/ $\sF_d$ there exists
an admissible epimorphism onto $M^\bu$ from a contractible
complex $L^\bu$ over\/ $\sF$ such that the kernel is a contractible
complex over\/~$\sF_{d-1}$.  \qed
\end{lem}

\begin{lem}
 Let $M^\bu\rarrow{}'\!\.M^\bu$ be a homotopy equivalence of complexes
over\/ $\sF_d$ such that ${}'\!M^\bu$ is the cokernel of
an admissible monomorphism of\/ $\sF_{d-1}$\+$\bst$\+acyclic complexes.
 Then ${}M^\bu$ is also the cokernel of an admissible monomorphism
of\/ $\sF_{d-1}$\+$\bst$\+acyclic complexes.
\end{lem}

\begin{proof}
 The proof is similar to that in~\cite{PP2}.
\end{proof}

 This finishes the proof in the case $\bst=\abs$.
 To deal with the case $\bst=\co$, notice that the full subcategories
$\sF_d\sub\sE$ are closed under infinite direct sums in its
assumptions, and so are the classes of all $\sF_d$\+$\bst$\+acyclic
complexes.
 Finally, the property of a complex over $\sE$ to be presentable
as an admissible monomorphism of $\sF_{d-1}$\+$\bst$\+acyclic complexes
is stable under infinite direct sums.
 The case $\bst=\ctr$ is similar.
\end{proof}

\subsection{Finite homological dimension}
\label{finite-homol-dim-subsect}
 First we return to the setting of a self-resolving full subcategory
$\sF$ in an exact category $\sE$, as in
Section~\ref{fully-faithful-subsect}.
 The following result is a partial generalization of
Lemma~\ref{homotopy-inj-proj-fully-faithful}.

\begin{prop}  \label{finite-homol-dim-fully-faithful}
 Assume that infinite products are everywhere defined and exact in
the exact category\/ $\sE$.
 Suppose also that the exact category\/ $\sF$ has finite homological
dimension.
 Then the triangulated functor\/ $\sD^\abs(\sF)\rarrow\sD^\ctr(\sE)$
induced by the exact embedding functor\/ $\sF\rarrow\sE$ is fully
faithful.
\end{prop}

\begin{proof}
 This is a particular case of~\cite[Theorem~8.15]{Pedg}.
 A combination of the assertions of
Proposition~\ref{fully-faithful-prop}(b) and Lemma~\ref{psemi-remark21}
implies the assertion of Proposition in the case when $\sF$ is 
preserved by the infinite products in~$\sE$.
 The proof in the general case consists in a combination of
the arguments proving the two mentioned results.

 We will show that any morphism in $\Hot(\sE)$ from a complex over $\sF$
to a complex contraacyclic over $\sE$ factorizes through a complex
absolutely acyclic over~$\sF$.
 For this purpose, it suffices to check that the class of complexes
over $\sE$ having this property contains the total complexes of short
exact sequences of complexes over $\sE$ and is closed with respect
to cones and infinite products.

 The former two assersions were proved in
Section~\ref{fully-faithful-subsect}.
 To prove the latter one, suppose that we are given a morphism
of complexes $F^\bu\rarrow\prod_\alpha A_\alpha^\bu$ in $\Hot(\sE)$,
where $F^\bu$ is a complex over~$\sF$.
 Suppose further that each component morphism $F^\bu\rarrow
A_\alpha^\bu$ factorizes through a complex $G_\alpha^\bu$ that is
absolutely acyclic over~$\sF$.
 Clearly, the morphism $F^\bu\rarrow\prod_\alpha A_\alpha^\bu$
factorizes through the complex $\prod_\alpha G_\alpha^\bu$.

 It follows from the proof the proof of Lemma~\ref{psemi-remark21}
(see~\cite[Remark~2.1]{Psemi}; for a more generally applicable argument,
see also~\cite[Section~1.6]{EP} or~\cite[Proposition~8.8]{Pedg})
that there exists an integer~$n$ such that every complex absolutely
acyclic over $\sF$ can be obtained from totalizations of short exact
sequences of complexes over $\sF$ by applying the operation of
the passage to a cone at most $n$~times.
 Therefore, the product $\prod_\alpha G_\alpha^\bu$ is an absolutely
acyclic complex over~$\sE$.
 Hence it follows from what we already know that the morphism
$F^\bu\rarrow\prod_\alpha G_\alpha^\bu$ factorizes through
a complex absolutely acyclic over~$\sF$.
\end{proof}

 Now let us assume that $\sF$ is a resolving full subcategory in
a weakly idempotent-complete exact category $\sE$, as in
Section~\ref{finite-resolutions-subsect}.
 The assumption of the following corollary is essentially
a generalization of the conditions (${*}$)\+-(${*}{*}$)
from~\cite[Sections~3.7\+-3.8]{Pkoszul}.

\begin{cor}  \label{finite-homol-dim-equivalence-cor}
 Assume that infinite products are everywhere defined and exact in
the exact category\/ $\sE$.
 Suppose also that the exact category\/ $\sF$ has finite homological
dimension; and assume additionally that countable products of objects
from\/ $\sF$ taken in the category $\sE$ have finite\/ $\sF$\+resolution
dimensions.
 Then the functor\/ $\sD^\abs(\sF)\rarrow\sD^\ctr(\sE)$ is
an equivalence of triangulated categories.
\end{cor}

\begin{proof}
 In view of Proposition~\ref{finite-homol-dim-fully-faithful},
it suffices to find for any complex over $\sE$ a morphism
into it from a complex over $\sF$ with a cone contraacyclic over~$\sE$.
 The following two-step construction procedure goes back
to~\cite{Pkoszul}.

 Given a complex $E^\bu$ over $\sE$, we proceed as in the above
proof of Proposition~\ref{infinite-resolutions}(b), applying
Lemma~\ref{second-kind-complex-resolution}(b) in order to obtain
a bicomplex $P^\bu_\bu$ over~$\sF$.
 By Lemma~\ref{telescope}, the total complex $T^\bu$ of the bicomplex
$P^\bu_\bu$ constructed by taking infinite products along the diagonals
maps onto $E^\bu$ with a cone contraacyclic with respect to~$\sE$.

 By assumption, the $\sF$\+resolution dimensions of the terms
of the complex $T^\bu$ are finite, and in fact bounded by a fixed
constant~$d$.
 Applying Lemma~\ref{second-kind-complex-resolution}(b) again
together with Corollary~\ref{fdim-cor}, we produce a finite complex
of complexes $Q_\bu^\bu$ of length~$d$ (in the lower indices)
over $\sF$ mapping termwise quasi-isomorphically onto~$T^\bu$.
 Now the total complex of $Q_\bu^\bu$ maps onto $T^\bu$ with
a cone absolutely acyclic over~$\sE$.
\end{proof}

 The following generalization of the last result is straightforward
(cf.~\cite[Theorem~8.16]{Pedg}).
 Let $\sG\sub\sF$ be two resolving full subcategories in an exact
category~$\sE$.
 Suppose that the pair of subcategories $\sF\sub\sE$ satisfies
the assumptions of Proposition~\ref{finite-homol-dim-fully-faithful},
while countable products of objects from\/ $\sG$ taken in
the category $\sE$ have finite $\sF$\+resolution dimension.
 Then the functor\/ $\sD^\abs(\sF)\rarrow\sD^\ctr(\sE)$ is
an equivalence of triangulated categories.

\Section{Cotorsion Pairs}
\label{cotorsion-pairs-appx}

\subsection{Preliminaries} \label{cotorsion-prelim-subsect}
 We refer to the introduction to the paper~\cite{Pctrl} for
an introductory discussion of the concept of a cotorsion pair.
 Discussions of cotorsion pairs in exact categories can be found
in~\cite[Section~5]{Sto-ICRA} and~\cite[Section~1]{Pal}.

 Let $\sE$ be an exact category, and let $\sF$, $\sC\sub\sE$ be two
classes of objects.
 One denotes by $\sF^{\perp_1}\sub\sE$ the class of all objects
$E\in\sE$ such that $\Ext^1_\sE(F,E)=0$ for all $F\in\sF$.
 Similarly, $\sF^{\perp_{\ge1}}\sub\sE$ is the class of all objects
$E\in\sE$ such that $\Ext^n_\sE(F,E)=0$ for all $F\in\sF$ and $n\ge1$.
 Dually, the notation ${}^{\perp_1}\sC\sub\sE$ stands for the class
of all objects $E\in\sE$ such that $\Ext^1_\sE(E,C)=0$ for all
$C\in\sC$, while ${}^{\perp_{\ge1}}\sC\sub\sE$ is the class of
all objects $E\in\sE$ such that $\Ext^n_\sE(E,C)=0$ for all $C\in\sC$
and $n\ge1$.

 A pair of classes of objects $(\sF,\sC)$ in $\sE$ is said to be
a \emph{cotorsion pair} if $\sC=\sF^{\perp_1}$ and $\sF={}^{\perp_1}\sC$.
 One can think of the objects from $\sF$ as of ``somewhat projective
objects'' in $\sE$ and of the objects from $\sC$ as of ``somewhat
injective objects''; then the definition of a cotorsion pair is
interpreted to mean that the two classes of somewhat projectives and
somewhat injectives fit together.

 Recall that a class of objects $\sF\sub\sE$ is said to be
\emph{generating} (cf.\ Section~\ref{infinite-resolutions-subsect}) if
for every object $E\in\sE$ there is an object $F\in\sF$ together with
an admissible epimorphism $F\rarrow E$ in~$\sE$.
 Dually, a class of objects $\sC\sub\sE$ is \emph{cogenerating} if
for every $E\in\sE$ there exists $C\in\sC$ together with an admissible
monomorphism $E\rarrow C$.
 If there are enough projective and injective objects in an exact
category $\sE$, then in any cotorsion pair $(\sF,\sC)$ in $\sE$
the left class $\sF$ is generating and the right class $\sC$ is
cogenerating; but this need not be the case in general.

 Let $(\sF,\sC)$ be a cotorsion pair in~$\sE$; assume that the class
$\sF$ is generating and the class $\sC$ is cogenerating in~$\sE$.
 Then the following four conditions are
equivalent~\cite[Lemma~6.17]{Sto-ICRA}, \cite[Section~1]{Pal},
\cite[Lemma~7.1]{PS6}:
\begin{enumerate}
\renewcommand{\theenumi}{\roman{enumi}}
\item the class $\sF$ is closed under kernels of admissible epimorphisms
in~$\sE$;
\item the class $\sC$ is closed under cokernels of admissible
monomorphisms in~$\sE$;
\item $\Ext^2_\sE(F,C)=0$ for all $F\in\sF$ and $C\in\sC$;
\item $\Ext^n_\sE(F,C)=0$ for all $F\in\sF$, \,$C\in\sC$, and $n\ge1$.
\end{enumerate}
 The condition~(i) can be expressed by saying that $\sF$ is
a resolving subcategory in $\sE$, while the condition~(ii) is
equivalent to the requirement that $\sC$ be a coresolving subcategory
in~$\sE$ (in the sense of Section~\ref{infinite-resolutions-subsect}).
 The condition~(iv) means that $\sF^{\perp_1}=\sF^{\perp_{\ge1}}=\sC$
and ${}^{\perp_1}\sC={}^{\perp_{\ge1}}\sC=\sF$.

 A cotorsion pair $(\sF,\sC)$ in $\sE$ is called \emph{hereditary}
if the class $\sF$ is generating in $\sE$, the class $\sC$ is
cogenerating in $\sE$, and the equivalent conditions~(i\+iv)
are satisfied.

 What does it mean that there are ``enough somewhat projective objects''
and/or ``enough somewhat injective objects'' in the context of
a cotorsion pair $(\sF,\sC)$ in~$\sE$\,?
 Here is the answer: a cotorsion pair $(\sF,\sC)$ is said to be
\emph{complete} if for every object $E\in\sE$ there exist (admissible)
short exact sequences
\begin{gather}
 0\rarrow C'\rarrow F\rarrow E\rarrow0
 \label{special-precover-sequence} \\
 0\rarrow E\rarrow C\rarrow F'\rarrow0
 \label{special-preenvelope-sequence}
\end{gather}
in $\sE$ with objects $F$, $F'\in\sF$ and $C$, $C'\in\sC$.
 A short exact sequence~\eqref{special-precover-sequence} is called
a \emph{special precover sequence}.
 A short exact sequence~\eqref{special-preenvelope-sequence} is called
a \emph{special preenvelope sequence}.
 Collectively, the short exact
sequences~(\ref{special-precover-sequence}\+-%
\ref{special-preenvelope-sequence}) are referred to as
\emph{approximation sequences}.

\begin{lem} \label{salce-lemma}
 Let\/ $\sE$ be an exact category and $(\sF,\sC)$ be a pair of classes
of objects in\/~$\sE$. \par
\textup{(a)} If the class\/ $\sF$ is generating and closed under
extensions in $\sE$, and a short exact
sequence~\eqref{special-preenvelope-sequence} with $C\in\sC$,
\,$F'\in\sF$ exists for every object $E\in\sE$, then a short exact
sequence~\eqref{special-precover-sequence} with $F\in\sF$, \,$C'\in\sC$
exists for every object $E\in\sE$. \par
\textup{(b)} If the class\/ $\sC$ is cogenerating and closed under
extensions in $\sE$, and a short exact
sequence~\eqref{special-precover-sequence} with $F\in\sF$, \,$C'\in\sC$
exists for every object $E\in\sE$, then a short exact
sequence~\eqref{special-preenvelope-sequence} with $C\in\sC$,
\,$F'\in\sF$ exists for every object $E\in\sE$.
\end{lem}

\begin{proof}
 These two dual assertions are known as the \emph{Salce
lemmas}~\cite{Sal}.
 One can find a module-theoretic version of the argument for part~(a)
in~\cite[second half of the proof of Theorem~10]{ET} and in
the proof of Lemma~\ref{eklof-trlifaj-cta}.
 See, e.~g., \cite[Lemma~1.1]{PS4} or~\cite[Lemma~1.3]{Pal} for
category-theoretic versions stated under somewhat more restrictive
assumptions than the formulation above.

 The same argument goes through in the general case.
 Specifically, in part~(a), given an object $E\in\sE$, pick
an admissible epimorphism $F\rarrow E$ in $\sE$ with $F\in\sF$,
and let $0\rarrow K\rarrow F\rarrow E\rarrow0$ be the resulting
short exact sequence.
 Choose a short exact sequence $0\rarrow K\rarrow C\rarrow G\rarrow0$
\,\eqref{special-preenvelope-sequence} for the object $K\in\sE$ with
some objects $C\in\sC$ and $G\in\sF$; and denote by $H$ the push-out of
the two admissible monomorphisms $K\rarrow F$ and $K\rarrow C$.
 Then there is an (admissible) short exact sequence $0\rarrow F\rarrow H
\rarrow G\rarrow0$ in $\sE$, implying that $H\in\sF$; and $0\rarrow C
\rarrow H\rarrow E\rarrow0$ is the desired short exact
sequence~\eqref{special-precover-sequence} for the object $E\in\sE$.
\end{proof}

 Given a class of objects $\sA\subset\sE$, we denote by $\sA^\oplus
\subset\sE$ the class of all direct summands of objects from $\sA$
in~$\sE$.

\begin{lem} \label{cotorsion-pair-direct-summands-lemma}
 Let\/ $\sE$ be an exact category and $(\sF,\sC)$ be a pair of classes
of objects in\/ $\sE$ such that $\Ext^1_\sE(F,C)=0$ for all $F\in\sF$
and $C\in\sC$.
 Assume that the short exact
sequences~(\ref{special-precover-sequence}\+-%
\ref{special-preenvelope-sequence}) with $F$, $F'\in\sF$ and
$C$, $C'\in\sC$ exist for all objects $E\in\sE$.
 Then one has\/ $\sF^{\perp_1}=\sC^\oplus$ and\/
${}^{\perp_1}\sC=\sF^\oplus$.
 In other words, the pair of classes of objects
$(\sF^\oplus,\sC^\oplus)$ is a complete cotorsion pair in\/~$\sE$.
\end{lem}

\begin{proof}
 This is~\cite[Lemma~1.2]{PS4} or~\cite[Lemma~1.4]{Pal};
in the module-theoretic context, a special case of this argument was
used in the proof of Corollary~\ref{very-flat-transfinite}.
 We only need to show that $\sF^{\perp_1}\subset\sC^\oplus$, and
dually, ${}^{\perp_1}\sC\subset\sF^\oplus$.
 Indeed, let $G\in{}^{\perp_1}\sC$ be an object.
 Pick a short exact sequence $0\rarrow C\rarrow F\rarrow G\rarrow0$
\,\eqref{special-precover-sequence} in $\sE$ with $F\in\sF$ and
$C\in\sC$.
 Then $G$ is a direct summand of $F$, since $\Ext^1_\sE(G,C)=0$. 
\end{proof}

 Let $\sK$ be an exact category and $\sE\sub\sK$ be a full subcategory
inheriting an exact category structure (as defined in
Section~\ref{fully-faithful-subsect}).
 Let $(\sF,\sC)$ be a cotorsion pair in~$\sE$.
 We say that $(\sF,\sC)$ \emph{restricts to a cotorsion pair} in $\sE$
if the pair of classes $\sE\cap\sF$ and $\sE\cap\sC$ is a cotorsion pair
in~$\sE$.
 If $(\sF,\sC)$ is complete in $\sK$ and $(\sE\cap\sF\;\sE\cap\sC)$ is
complete in $\sE$, we say that $(\sF,\sC)$ \emph{restricts to
a complete cotorsion pair} in~$\sE$.

\begin{lem} \label{restricting-hereditary-cotorsion}
 Let\/ $\sK$ be an exact category, $\sE\sub\sK$ be a full subcategory
inheriting an exact category structure, and $(\sF,\sC)$ be
a hereditary cotorsion pair in\/~$\sK$.
 Assume that $(\sF,\sC)$ restricts to a cotorsion pair in\/~$\sE$,
the class\/ $\sE\cap\sF$ is generating in\/ $\sE$, and the class\/
$\sE\cap\sC$ is cogenerating in\/~$\sE$.
 Then the cotorsion pair $(\sE\cap\sF\;\sE\cap\sC)$ is hereditary
in\/~$\sE$.
\end{lem}

\begin{proof}
 This is essentially~\cite[Lemma~1.6]{Pal}.
\end{proof}

\begin{lem} \label{restricting-cotorsion-pairs-lemma}
 Let\/ $\sK$ be an exact category, $\sE\sub\sK$ be a full additive
subcategory, and $(\sF,\sC)$ be a complete cotorsion pair in\/~$\sK$.
 Assume that one of the following two conditions holds: \par
\textup{(a)} $\sF\sub\sE$, and\/ $\sE$ is closed under extensions and
kernels of admissible epimorphisms in\/~$\sK$; or \par
\textup{(b)} $\sC\sub\sE$, and\/ $\sE$ is closed under extensions and
cokernels of admissible monomorphisms in\/~$\sK$. \par
 Then $(\sF,\sC)$ restricts to a complete cotorsion pair in\/~$\sE$.
\end{lem}

\begin{proof}
 This is~\cite[Lemma~1.5]{Pal}.
\end{proof}

\begin{exs} \label{cotorsion-pairs-examples}
 (1)~For any associative ring $R$, the pair of classes of flat
$R$\+modules $\sF=R\modl_\fl$ and cotorsion $R$\+modules
$\sC=R\modl^\cot$ is a hereditary complete cotorsion pair in
the abelian category of $R$\+modules $\sE=R\modl$.
 This assertion summarizes the discussion in the beginning of
Section~\ref{cotorsion-modules}, including Theorem~\ref{flat-cover-thm}.
 It is helpful to keep in mind
Lemma~\ref{cotorsion-pair-direct-summands-lemma}.

 (2)~For any commutative ring $R$, the pair of classes of very flat
$R$\+modules $\sF=R\modl_\vfl$ and contraadjusted $R$\+modules
$\sC=R\modl^\cta$ is a hereditary complete cotorsion pair in
the abelian category of $R$\+modules $\sE=R\modl$.
 This assertion summarizes much of the discussion in
Section~\ref{very-eklof-trlifaj-subsect}.

 (3)~To give a more complicated example, one can observe that, for
any commutative ring $R$, the pair of classes of very flat $R$\+modules
$\sF=R\modl_\vfl$ and flat contraadjusted $R$\+modules
$\sC=R\modl_\fl^\cta$ is a hereditary complete cotorsion pair in
the exact category of flat $R$\+modules $\sE=R\modl_\fl$.
 This assertion can be obtained from~(2) using
Lemmas~\ref{restricting-hereditary-cotorsion}
and~\ref{restricting-cotorsion-pairs-lemma}(a).

 (4)~Similarly, the pair of classes of flat contraadjusted
$R$\+modules $\sF=R\modl_\fl^\cta$ and cotorsion $R$\+modules
$\sC=R\modl^\cot$ is a hereditary complete cotorsion pair in the exact
category of contraadjusted $R$\+modules $\sE=R\modl^\cta$.
 This follows from~(1) in view of
Lemmas~\ref{restricting-hereditary-cotorsion}
and~\ref{restricting-cotorsion-pairs-lemma}(b).
\end{exs}

 Given an additive category $\sA$, we denote by $\Com(\sA)$
the additive category of unbounded complexes in~$\sA$.
 The notation $\Hot(\sA)$ stands for the triangulated homotopy
category of complexes in~$\sA$.

 When $\sA$ is an exact category, we endow the additive category
$\Com(\sA)$ with the exact category structure in which the admissible
short exact sequences in $\Com(\sA)$ are the termwise admissible short
exact sequences of complexes over~$\sA$, i.~e., the short sequences of
complexes that are admissible exact in every degree.
 We will call this exact structure on $\Com(\sA)$ the \emph{termwise
exact structure}.
 If $\sA$ is an abelian category with the abelian exact structure,
then the category $\Com(\sA)$ is abelian, too, and the termwise exact
structure on $\Com(\sA)$ coincides with the abelian exact structure. 

\begin{lem} \label{G-plus-G-minus-complexes-Ext-adjunction}
 Let\/ $\sE$ be an exact category and $E\in\sE$ be an object.
 Denote by $D^{n,n+1}(E)$ the contractible two-term complex
$\dotsb\rarrow 0\rarrow E\overset=\rarrow E\rarrow0\rarrow\dotsb$
situated in the cohomological degrees $n$ and~$n+1$.
 Then, for any complex $C^\bu$ in\/ $\sE$ and any integer $i\ge0$,
the following natural isomorphisms of abelian groups hold:
\begin{align*}
 \Ext^i_{\Com(\sE)}(D^{n,n+1}(E),C^\bu)&\simeq\Ext^i_\sE(E,C^n); \\
 \Ext^i_{\Com(\sE)}(C^\bu,D^{n,n+1}(E))&\simeq\Ext^i_\sE(C^{n+1},E).
\end{align*}
\end{lem}

\begin{proof}
 This is~\cite[Lemma~3.1(5,6)]{Gil}, \cite[Lemma~4.15]{Sto},
or~\cite[Lemma~1.8]{Pal}.
 Let $\sE^\boZ$ denote the category of $\boZ$\+graded objects in $\sE$,
i.~e., the Cartesian product of $\boZ$ copies of $\sE$, endowed with
the natural (componentwise) exact category structure.
 The desired natural isomorphisms of the Ext groups are special cases
of the Ext\+adjunction lemma (\cite[Lemma~1.7(a)]{Pal}
or~\cite[Lemma~6.1]{PS6}) applied to the exact forgetful functor
$\Com(\sE)\rarrow\sE^\boZ$ and its exact adjoint functors on both sides.
\end{proof}

\begin{lem} \label{Ext-1-as-homotopy-Hom}
 Let\/ $\sE$ be an exact category, and let $F^\bu$ and $C^\bu$ be two
complexes in\/~$\sE$.
 Assume that\/ $\Ext^1_\sE(F^n,C^n)=0$ for all $n\in\boZ$.
 Then there is a natural isomorphism of abelian groups
$$
 \Ext^1_{\Com(\sE)}(F^\bu,C^\bu)\simeq\Hom_{\Hot(\sE)}(F^\bu,C^\bu[1]).
$$
\end{lem}

\begin{proof} \hbadness=1500
 More generally, for any two complexes $F^\bu$ and $C^\bu$ in $\sE$,
the group $\Hom_{\Hot(\sE)}(F^\bu,C^\bu[1])$ is naturally isomorphic
to the subgroup in $\Ext^1_{\Com(\sE)}(F^\bu,C^\bu)$ consisting of
the classes of degreewise split extensions $0\rarrow C^\bu\rarrow
E^\bu\rarrow F^\bu\rarrow0$.
 This well-known observation goes back, at least,
to~\cite[Lemma~2.1]{Gil}; see also~\cite[Section~1.3]{Bec},
\cite[proof of Lemma~5.1]{Sto}, \cite[Lemma~5.1]{PS4},
\cite[Lemma~1.6]{BHP}, or~\cite[Lemma~1.9]{Pal}.
 A further generalization can be found in~\cite[Lemma~9.41]{Pedg}.
\end{proof}

\begin{lem} \label{F-tilde-orthogonal-to-C-tilde}
 Let $(\sF,\sC)$ be a cotorsion pair in an exact category\/~$\sE$.
 Endow\/ $\sF$ and\/ $\sC$ with the exact category structures inherited
from\/~$\sE$.
 Let $F^\bu$ be an acyclic complex in\/ $\sF$ and $C^\bu$ be
an acyclic complex in\/~$\sC$.
 Then any morphism of complexes $F^\bu\rarrow C^\bu$ is homotopic to
zero.
 In other words, $\Ext^1_{\Com(\sE)}(F^\bu,C^\bu)=0$.
\end{lem}

\begin{proof}
 This is essentially~\cite[Lemma~3.9]{Gil}.
 The two assertions stated in the lemma are equivalent by
Lemma~\ref{Ext-1-as-homotopy-Hom}.
 To prove them, notice that any exact complex $F^\bu$ in $\sF$ is
the termwise stabilizing inductive limit of its finite subcomplexes of
canonical and silly truncation.
 In fact, this construction represents $F^\bu$ as a transfinitely
iterated extension, in the sense of the direct limit, of one-term
complexes corresponding to objects from~$\sF$ (in the sense of
the definition in Section~\ref{small-object-argument-subsect} below).
 It is also clear that, for any object $F\in\sF$ and any acyclic
complex $C^\bu$ in $\sC$, any morphism of complexes $F\rarrow C^\bu$
is homotopic to zero.
 Now one can refer to Lemma~\ref{telescope} for the conclusion that
the complex of abelian groups $\Hom_\sE(F^\bu,C^\bu)$ is acyclic, or
to Lemma~\ref{eklof-lemma-general} for the conclusion that
the Ext group $\Ext^1_{\Com(\sE)}(F^\bu,C^\bu)$ vanishes.
\end{proof}

 We refer to Section~\ref{finite-resolutions-subsect} for
the definition of the (co)resolution dimension with respect to
a (co)resolving subcategory in a weakly idempotent-complete
exact category.

\begin{lem} \label{finite-homol-dim-finite-coresol-dim}
 Let $(\sF,\sC)$ be a hereditary complete cotorsion pair in
a weakly idempotent-complete exact category\/~$\sE$ and $d\ge0$
be an integer.
 Then the following conditions are equivalent:
\begin{enumerate}
\item all objects of\/ $\sE$ have finite\/ $\sC$\+coresolution
dimensions not exceeding~$d$;
\item all objects of\/ $\sF$ have finite\/
$(\sF\cap\nobreak\sC)$\+coresolution dimensions not exceeding~$d$;
\item $\Ext^n_\sE(F,E)=0$ for all objects $F\in\sF$, $E\in\sE$,
and integers $n>d$ (in other words, all the objects of\/ $\sF$
have projective dimensions not exceeding~$d$ in\/~$\sE$);
\item the exact category\/ $\sF$ has finite homological dimension
not exceeding~$d$.
\end{enumerate}
\end{lem}

\begin{proof}
 Notice first of all that, for any complete cotorsion pair $(\sF,\sC)$
in an exact category $\sE$, the intersection $\sF\cap\sC$ is
the full subcategory of injective objects in $\sF$ and the full
subcategory of projective objects in $\sC$, and there are enough of both.

 (1)~$\Longrightarrow$~(2) The $\sC$\+coresolution dimensions of
the objects of $\sF$ in $\sE$ coincide with their
$(\sF\cap\nobreak\sC)$\+coresolution dimensions in $\sF$ by the dual
version of Corollary~\ref{fdim-subcategory-cor}.

 (2)~$\Longleftrightarrow$~(4) Follows from the first paragraph of
this proof.

 (1)~$\Longleftrightarrow$~(3) Obvious.
 
 (4)~$\Longrightarrow$~(3) Arguing as in~\cite[Lemma~8.3]{PS5}
and~\cite[Lemmas~3.4\+-3.5]{Pgen}, one can show that, for any
self-resolving full subcategory $\sF$ in an exact category $\sE$
(in the sense of Section~\ref{fully-faithful-subsect}) and any objects
$F\in\sF$ and $E\in\sE$, the groups $\Ext^*_\sE(F,E)$ can be computed
using resolutions of $F$ by objects from~$\sF$.
 So any given class $\xi\in\Ext^n_\sE(F,E)$ can be represented by
a cocycle $c\in\Hom_\sE(F_n,E)$ in the complex $\Hom_\sE(F_\bu,E)$
for some exact complex $F_\bu\rarrow F\rarrow0$ in~$\sF$.
 In other words, this means that there exist an object $G\in\sF$,
a class $\eta\in\Ext^n_\sF(F,G)$, and a morphism $g\:G\rarrow E$
such that $\xi=g\eta$.
 Now (4)~tells us that $\Ext^n_\sF(F,G)=0$, and it follows that
$\Ext^n_\sE(F,E)=0$.
\end{proof}

\subsection{Small object argument} \label{small-object-argument-subsect}
 The classical small object argument was invented by
Quillen~\cite[Lemma~II.3.3]{Quil} as a part of his model
category theory; see~\cite[Theorem~2.1.14]{Hov-book},
\cite[Proposition~1.3]{Bek}, or~\cite[Proposition~2.1]{SS}
for a general formulation.
 The interpretation of the Eklof--Trlifaj theorem~\cite{ET} as
an application of the small object argument was suggested by
Rosick\'y~\cite{Ros} and Hovey~\cite{Hov}.

 Two far-reaching generalizations of the Eklof--Trlifaj theorem,
in two different directions, can be found in the contemporary
literature: the version for efficient exact
categories~\cite[Theorem~2.13 and Corollary~2.15]{SS},
\cite[Theorem~5.16 and Corollary~5.17]{Sto-ICRA}
and the version for locally presentable abelian
categories~\cite[Corollary~3.6 and Theorem~4.8]{PR},
\cite[Theorems~3.3 and~3.4]{PS4}.
 The advantage of the former is that it is applicable to exact 
categories and not only to abelian ones; but we cannot use it for
our purposes, as the contraherent cosheaf categories are not efficient
(they do not even have infinite direct sums).
 The latter approach is more relevant to us in that it is applicable
to the abelian categories of contramodules alongside with the usual
modules, sheaves, and complexes.
 So we include the formulations of these results from~\cite{PR,PS4}
here, starting from some general definitions.

 Let $\sK$ be a category and $\alpha$~be an ordinal.
 An \emph{$\alpha$\+indexed smooth chain} (of objects and morphisms)
in $\sK$ is a direct system $(K_i\to K_j)_{0\le i<j<\alpha}$ in
$\sK$ indexed by the ordered set~$\alpha$ (i.~e., by the ordinals
$0\le i<\alpha$) such that, for every limit ordinal $j<\alpha$, one has
$K_j=\varinjlim_{i<j}K_i$ in~$\sK$.
 
 Let $\sE$ be an exact category and $F\in\sE$ be an object.
 An \emph{$\alpha$\+indexed filtration} on $F$ is defined as
an $(\alpha+\nobreak1)$\+indexed smooth chain
$(F_i\to F_j)_{0\le i<j\le\alpha}$ in $\sE$ such that $F_0=0$,
\,$F_\alpha=F$, and the morphism $F_i\rarrow F_{i+1}$ is
an admissible monomorphism in $\sE$ for all $0\le i<\alpha$.
 If this is the case, denote by $S_i\in\sE$ the cokernel of
the admissible monomorphism $F_i\rarrow F_{i+1}$.
 The object $F$ is said to be \emph{filtered by} the objects~$S_i$.
 In a different terminology, one says that the object $F\in\sE$ is
a \emph{transfinitely iterated extension} (\emph{in the sense of
the direct limit}) of the objects $S_i\in\sE$, \ $0\le i<\alpha$.

 One should keep in mind the possibilities of complicated behavior
of filtrations in exact or abelian categories where the direct limit
functors are not exact.
 See~\cite[Example~4.4]{PR} for an example of a nontrivial filtration
of a zero object.

 Given a class of objects $\sS\sub\sE$, we denote by $\Fil(\sS)\sub\sE$
the class of all objects of $\sE$ filtered by (objects isomorphic to) objects from~$\sS$.
 A class of objects $\sF\sub\sE$ is said to be \emph{deconstructible}
if there exists a \emph{set} of objects $\sS\sub\sE$ such that
$\sF=\Fil(\sS)$.

 The following result is known as the \emph{Eklof
lemma}~\cite[Lemma~1]{ET}, \cite[Lemma~6.2]{GT},
\cite[Proposition~2.12]{SS}, \cite[Proposition~5.7]{Sto-ICRA},
\cite[Lemma~4.5]{PR}, \cite[Proposition~1.3]{PS4},
\cite[Lemma~1.1]{BHP}, \cite[Lemma~7.5]{PS6}.

\begin{lem} \label{eklof-lemma-general}
 Let\/ $\sE$ be an exact category and\/ $\sC\sub\sE$ be a class of
objects.
 Then the class of objects\/ ${}^{\perp_1}\sC$ is closed under
transfinitely iterated extensions in\/~$\sE$.
 In other words, $\Fil({}^{\perp_1}\sC)={}^{\perp_1}\sC$.
\end{lem}

\begin{proof}
 The argument from~\cite[Lemma~4.5]{PR} is applicable.
\end{proof}

 For any class of objects $\sS$ in an exact category $\sE$, the pair
of classes $\sF={}^{\perp_1}(\sS^{\perp_1})$ and $\sC=\sS^{\perp_1}$
is a cotorsion pair in~$\sE$.
 The cotorsion pair $(\sF,\sC)$ is said to be \emph{generated by}
the class of objects~$\sS$.

 The notation $\sA^\oplus$ for a class of objects $\sA\subset\sE$ was
introduced in Section~\ref{cotorsion-prelim-subsect}.

 We refer to the book~\cite[Definition~1.17 and Theorem~1.20]{AR}
for the definition of a \emph{locally presentable category}.
 We are interested in locally presentable \emph{abelian} categories
(studied in the papers~\cite{PR,PS4} and the preprint~\cite{Pper}).
 In this context, it is helpful to keep in mind that, for abelian
categories, there is no difference between the concepts of a generator
and a strong generator (any generator is strong).
 All Grothendieck abelian categories are locally
presentable~\cite[Lemma~A.1]{Sto-hill}.
 The additive category of complexes $\Com(\sA)$ is locally presentable
for any locally presentable additive category~$\sA$
\,\cite[Lemma~6.3]{PS4}.

 The following theorem, going back to Eklof and
Trlifaj~\cite[Theorems~2 and~10]{ET}, is a far-reaching generalization
of the discussion in Section~\ref{very-eklof-trlifaj-subsect}.

\begin{thm} \label{eklof-trlifaj-general}
 Let\/ $\sK$ be a locally presentable abelian category, $\sS\sub\sK$ be
a \emph{set} of objects, and $(\sF,\sC)$ be the cotorsion pair generated
by\/ $\sS$ in\/~$\sK$.
 In this setting: \par
\textup{(a)} if the class\/ $\sF$ is generating and the class\/ $\sC$
is cogenerating in\/ $\sK$, then the cotorsion pair $(\sF,\sC)$
is complete; \par
\textup{(b)} if the class\/ $\Fil(\sS)$ is generating in\/ $\sK$,
then\/ $\sF=\Fil(\sS)^\oplus$.
\end{thm}

\begin{proof}
 Part~(a) is~\cite[Corollary~3.6]{PR} or~\cite[Theorem~3.3]{PS4}.
 Part~(b) is~\cite[Theorem~4.8(d)]{PR} or~\cite[Theorem~3.4]{PS4}.
\end{proof}

\begin{cor} \label{eklof-trlifaj-corollary}
 Let\/ $\sK$ be a locally presentable abelian category and\/
$\sF\sub\sK$ be a class of objects.
 Assume that the class\/ $\sF$ is generating, deconstructible, and
closed under direct summands in\/~$\sK$.
 Then the class\/ $\sF$ is a part of a cotorsion pair $(\sF,\sC)$
in\/~$\sK$.
 If the class\/ $\sC$ is cogenerating in\/ $\sK$ (e.~g., this always
holds for a Grothendieck abelian category~$\sK$), then the cotorsion
pair $(\sF,\sC)$ is complete.
\end{cor}

\begin{proof}
 Let\/ $\sS\sub\sK$ be a set of objects such that $\sF=\Fil(\sS)$.
 Then $\sS^{\perp_1}=\sF^{\perp_1}\sub\sK$ by
Lemma~\ref{eklof-lemma-general}.
 It remains to apply Theorem~\ref{eklof-trlifaj-general} to deduce
both assertions of the corollary.
\end{proof}

\begin{lem} \label{kaplansky-for-filtrations}
 Let\/ $\sA$ be a Grothendieck abelian category.
 Then, for any deconstructible class of objects\/ $\sG\sub\sA$,
the class of objects\/ $\sG^\oplus\sub\sA$ is deconstructible, too.
\end{lem}

\begin{proof}
 This is~\cite[Lemma~7.12]{GT} (for categories of modules over
associative rings) or~\cite[Proposition~2.9(1)]{Sto-hill}
(for Grothendieck categories).
\end{proof}

\begin{exs} \label{deconstructible-examples}
 For any associative ring $R$, the class of flat $R$\+modules
$R\modl_\fl$ is deconstructible in $R\modl$
\,\cite[Lemma~1 and Proposition~2]{BBE}.

 For any commutative ring $R$, the class of very flat $R$\+modules
$R\modl_\vfl$ is deconstructible in $R\modl$.
 This follows from Corollary~\ref{very-flat-transfinite} and
Lemma~\ref{kaplansky-for-filtrations} (see~\cite[Lemma~2.2]{ST}
for further details).
\end{exs}

 Let $\sF$ be an exact category with exact functors of infinite
direct sum.
 The full subcategory $\sD(\sF)^\lh$ in the exact category $\sD(\sF)$
was defined in Section~\ref{homotopy-adjusted}.
 Denote by $\Com(\sF)^\lh\sub\Com(\sF)$ the class of all complexes in
$\sE$ which, viewed as objects of the derived category $\sD(\sF)$,
belong to $\sD(\sF)^\lh$.
 Given an exact category $\sE$, we denote by $\Acycl(\sE)\sub\Com(\sE)$
the class of all acyclic complexes in~$\sE$.

 The following result goes back to~\cite[Theorem~4.12]{Gil1}
and~\cite[Lemma~4.9]{Gil3}.

\begin{prop} \label{dg-F-acycl-C-cotorsion-pair}
 Let\/ $\sA$ be a Grothendieck abelian category and $(\sF,\sC)$ be
a hereditary complete cotorsion pair generated by a set of objects
in\/~$\sA$.
 Endow the full subcategories\/ $\sF$ and\/ $\sC\sub\sA$ with the exact
category structures inherited from the abelian exact structure
of\/~$\sA$.
 Then the pair of classes\/ $\sF'=\Com(\sF)^\lh\sub\Com(\sF)\sub
\Com(\sA)$ and $\sC'=\Acycl(\sC)\sub\Com(\sC)\sub\Com(\sA)$ is
a hereditary complete cotorsion pair in the Grothendieck abelian
category\/ $\Com(\sA)$.
\end{prop}

\begin{proof}
 Let\/ $\sS_0\sub\sA$ be a set of objects generating the cotorsion pair
$(\sF,\sC)$ in\/~$\sA$.
 Since this cotorsion pair is complete by assumption, the class $\sF$
is generating in $\sA$, so we can choose a set of generators $\sS_1$
of the category $\sA$ such that $\sS_1\sub\sF$.
 Put $\sS_2=\sS_0\cup\sS_1$.
 Then the class $\Fil(\sS_2)$ is generating in $\sA$, so
Theorem~\ref{eklof-trlifaj-general}(b) tells us that
$\sF=\Fil(\sS_2)^\oplus$.
 By Lemma~\ref{kaplansky-for-filtrations}, there even exists a set of
objects $\sS_3\sub\sA$ such that $\sF=\Fil(\sS_3)$ (but we will not
need to use this fact).

 Denote by $\sS'\sub\Com(\sA)$ the set of all one-term complexes
$\dotsb\rarrow0\rarrow S\rarrow0\rarrow\dotsb$ with $S\in\sS_2$.
 Let $\sF'={}^{\perp_1}(\sS'{}^{\perp_1})$ and $\sC'=\sS'{}^{\perp_1}$
be the cotorsion pair in $\Com(\sA)$ generated by the set~$\sS'$.
 By Lemma~\ref{eklof-lemma-general}, we have $\Fil(\sS')\sub\sF'$.
 It is clear from the construction that the class $\Fil(\sS')$
is generating in $\Com(\sA)$, while the class $\sC'$ is cogenerating
since $\Com(\sA)$ is a Grothendieck abelian category (so it has enough
injectives).
 Thus both parts of Theorem~\ref{eklof-trlifaj-general} are applicable,
and we can conclude that $\sF'=\Fil(\sS')^\oplus$ and the cotorsion pair
$(\sF',\sC')$ is complete in $\Com(\sA)$.

 It is clear from Lemma~\ref{G-plus-G-minus-complexes-Ext-adjunction}
(for $i=1$) that $\sC'\sub\Com(\sC)$.
 In order to show that $\sC'=\Acycl(\sC)$, we use
Lemma~\ref{Ext-1-as-homotopy-Hom}.
 Indeed, let $C^\bu$ be a complex belonging to~$\sC'$.
 Denote by $Z^n$ the kernel of the differential
$C^n\rarrow C^{n+1}$ taken in the abelian category~$\sA$
(for some $n\in\boZ$).
 By assumption, there exists an object $F\in\sF$ together with
an epimorphism $F\rarrow Z^n$ in~$\sA$.
 This morphism corresponds to a degree~$n$ cocycle in the complex
$\Hom_\sA(F,C^\bu)$.
 If this cocycle is a coboundary, then the map $C^{n-1}\rarrow Z^n$
is an epimorphism in~$\sA$.
 So the complex $C^\bu$ is acyclic in~$\sA$.
 The condition of acyclicity of the complexes of abelian groups
$\Hom_\sA(F,C^\bu)$ for all objects $F\in\sF$ further implies that
the kernel $Z^{n-1}$ of the morphism $C^{n-1}\rarrow Z^n$
belongs to~$\sC$.
 Thus $C^\bu\in\Acycl(\sC)$.
 Conversely, it is obvious that the complex $\Hom_\sA(F,C^\bu)$ is
acyclic for all $F\in\sF$ and $C^\bu\in\Acycl(\sC)$.

 Clearly, $\sF'=\Fil(\sS')^\oplus\sub\Com(\sF)$ and
$\sF=\Fil(\sS_2)^\oplus\sub\sF'$.
 By Lemma~\ref{F-tilde-orthogonal-to-C-tilde}, all acyclic complexes
in the exact category $\sF$ belong to $\sF'$, that is
$\Acycl(\sF)\sub\sF'$.
 It is also clear that the class $\sF'$ is closed under shifts, cones,
and infinite direct sums in $\Com(\sF)$.
 By Lemma~\ref{Ext-1-as-homotopy-Hom}, the class $\sF'$ is closed
under homotopy equivalences of complexes in~$\sF$.
 Thus $\Com(\sF)^\lh\sub\sF'$.
 In order to prove the converse inclusion, let $G^\bu$ be a complex
from the class~$\sF'$.
 By Proposition~\ref{spal-for-exact}, there exists a complex
$F^\bu\in\Com(\sF)^\lh$ together with a quasi-isomorphism
$F^\bu\rarrow G^\bu$ of complexes in~$\sA$.
 The cone $H^\bu$ of this quasi-isomorphism is an acyclic complex
in $\sA$ belonging to~$\sF'$.

 Now, given an object $C\in\sC$, choose its injective coresolution
$J^\bu$ in $\sA$, and denote by $K^\bu$ the augmented coresolution
$C\rarrow J^\bu$.
 Then the complex of abelian groups $\Hom_\sA(H^\bu,J^\bu)$ is
acyclic, because the complex $H^\bu$ is acyclic in $\sA$, while
$J^\bu$ is a bounded below complex of injectives.
 On the other hand, the complex of abelian groups
$\Hom_\sA(H^\bu,K^\bu)$ is acyclic by
Lemma~\ref{Ext-1-as-homotopy-Hom}, because $H^\bu\in\sF'$ and
$K^\bu\in\Acycl(\sC)=\sC'$ (here we use the assumption that the class
$\sC$ is closed under cokernels of monomorphisms in~$\sA$).
 Therefore, the complex $\Hom_\sA(H^\bu,C)$ is acyclic, too.
 One can easily conclude that $\Ext_\sA^1(Z^n,C)=0$ for every
object of cocycles $Z^n$ of the complex $H^\bu$ in~$\sA$ (since
the terms of the complex $H^\bu$ belong to~$\sF$).
 Hence $Z^n\in\sF$.
 We have shown that the complex $H^\bu$ is acyclic in the exact
category $\sF$, so the morphism $F^\bu\rarrow G^\bu$ is
an isomorphism in $\sD(\sF)$.
 Thus $F^\bu\in\Com(\sF)^\lh$ implies $G^\bu\in\Com(\sF)^\lh$.

 Finally, the cotorsion pair $(\sF',\sC')$ is hereditary in
$\Com(\sA)$, because the class $\sC'=\Acycl(\sC)$ is closed under
cokernels of monomorphisms in $\Com(\sA$).
\end{proof}

\subsection{Local classes of rings, modules, and complexes}
\label{local-classes-subsect}
 Here and below, we will assume all our classes or full subcategories
of rings, modules, and complexes to be closed under isomorphisms.
 This section builds upon~\cite[Section~2]{Pal}.

 A class of commutative rings $\sR$ is said to be \emph{local} if it
satisfies the following two conditions:
\begin{description}
\item[Ascent for rings] For any ring $R\in\sR$ and any element $r\in R$,
the ring $R[r^{-1}]$ belongs to~$\sR$.
\item[Descent for rings] Let $R$ be a ring and $r_1$,~\dots, $r_d\in R$
be a finite set of elements generating the unit ideal in $R$, that is,
$(r_1,\dotsc,r_d)=R$.
 Assume that $R[r_j^{-1}]\in\sR$ for every $j=1$,~\dots,~$d$.
 Then $R\in\sR$.
\end{description}

\begin{lem} \label{rings-affine-communication}
 Let\/ $\sR$ be a local class of commutative rings.  Then \par
\textup{(a)} for any affine open subscheme $V$ in an affine scheme $U$,
one has\/ $\O(V)\in\sR$ whenever\/ $\O(U)\in\sR$; \par
\textup{(b)} for any affine open covering $U=\bigcup_\alpha V_\alpha$
of an affine scheme $U$, one has\/ $\O(U)\in\sR$ whenever\/
$\O(V_\alpha)\in\sR$ for every~$\alpha$.
\end{lem}

\begin{proof}
 This is (a particular case of what is) called \emph{Affine
Communication Lemma} in~\cite[Lemma~5.3.2]{Vak};
see also~\cite[Section Tag~01OO]{SP}.
 To prove part~(a), cover $V$ by affine open subschemes that are
principal affine in~$U$ (hence also in~$V$).
 For part~(b), cover each $V_\alpha$ by affine open subschemes
that are principal affine in~$U$.
\end{proof}

 Given a local class of commutative rings $\sR$, we will say that
a scheme $X$ is \emph{locally\+$\sR$} if one has $\O_X(U)\in\sR$
for every affine open subscheme $U\sub X$.
 It is clear from Lemma~\ref{rings-affine-communication} that it
suffices to check this condition for affine open subschemes $U_\alpha$
belonging to any given affine open covering $X=\bigcup_\alpha U_\alpha$
of the scheme~$X$.

 The theory of (co)locality and gluing below is developed in two
versions: either for modules/sheaves/cosheaves, or for complexes of
modules/sheaves/cosheaves.
 Hence the following notation with the two options built into it.

 Given a ring $R$, we denote by $\sK_R$ either the abelian category
of $R$\+modules $\sK_R=R\modl$, or the abelian category of complexes of
$R$\+modules $\sK_R=\Com(R\modl)$.
 This is our big ambient abelian category of modules/complexes.

 Let $R$ be a local class of rings.
 Suppose given a class of objects $\sL_R\sub\sK_R$ defined for every
ring $R\in\sR$.
 By abuse of terminology and notation, we will speak of ``a~class
of modules/complexes~$\sL$\,'' (over all rings $R\in\sR$).
 A class $\sL$ is called \emph{local} if it satisfies the following
two conditions:
\begin{description}
\item[Ascent] For any ring $R\in\sR$, any element $r\in R$, and any
module/complex $M\in\sL_R$, the module/complex $M[r^{-1}]=
R[r^{-1}]\ot_RM$ belongs to~$\sL_{R[r^{-1}]}$.
\item[Descent] Let $R\in\sR$ be a ring and $r_1$,~\dots, $r_d\in R$
be a finite set of elements generating the unit ideal in~$R$.
 Let $M\in\sK_R$ be a module/complex such that
$M[r_j^{-1}]\in\sL_{R[r_j^{-1}]}$ for every $j=1$,~\dots,~$d$.
 Then $M\in\sL_R$.
\end{description}

\begin{lem} \label{locality-affine-communication}
 Let\/ $\sR$ be a local class of commutative rings and $U$ be
a locally\+$\sR$ affine scheme.
 Let\/ $\sL$ be a local class of modules/complexes over rings $R\in\sR$,
and let\/ $\M$ be a quasi-coherent sheaf or complex of quasi-coherent
sheaves over~$U$.
 Then \par
\textup{(a)} for any affine open subscheme $V\sub U$ in the affine
scheme $U$, one has\/ $\M(V)\in\sL_{\O(V)}$ provided that\/
$\M(U)\in\sL_{\O(U)}$; \par
\textup{(b)} for any affine open covering $U=\bigcup_\alpha V_\alpha$
of the affine scheme $U$, one has\/ $\M(U)\in\sL_{\O(U)}$ provided
that\/ $\M(V_\alpha)\in\sL_{\O(V_\alpha)}$ for every~$\alpha$.
\end{lem}

\begin{proof}
 Similar to Lemma~\ref{rings-affine-communication}.
\end{proof}

 Let $\sR$ be a local class of commutative rings and $\sL$ be
a local class of modules or complexes over rings $R\in\sR$.
 Let $X$ be a locally\+$\sR$ scheme and $\M$ be a quasi-coherent
sheaf/complex of quasi-coherent sheaves over~$X$.
 We will say that the sheaf/complex $\M$ is \emph{locally\+$\sL$} if
one has $\M(U)\in\sL_{\O_X(U)}$ for every affine open subscheme
$U\subset X$.
 It is clear from Lemma~\ref{locality-affine-communication} that it
suffices to check this condition for affine open subschemes belonging
to any given affine open covering of the scheme~$X$.

 More generally, given a locally\+$\sR$ scheme $X$ and an open
subscheme $Y\subset X$, we will colloquially say that a sheaf/complex
$\M$ on $X$ is ``locally\+$\sL$ over~$Y$'' if the sheaf/complex
$\M|_Y$ on $Y$ is locally\+$\sL$.

 For any scheme $X$, put $\sK_X=X\qcoh$ if we are dealing with modules
(i.~e., if $\sK_R=R\modl$), or $\sK_X=\Com(X\qcoh)$ if we are dealing
with complexes of modules (i.~e., if $\sK_R=\Com(R\modl)$).
 Given a local class of commutative rings $\sR$, a local class of
modules/complexes $\sL$ over rings $R\in\sR$, and a locally\+$\sR$
scheme $X$, we will denote by $\sL_X\subset\sK_X$ the full subcategory
of locally\+$\sL$ sheaves or complexes.

 We will say that a class of modules/complexes $\sL$ over rings
$R\in\sR$ is \emph{very local} if $\sL$ is local and satisfies
the following
\begin{description}
\item[Direct image condition] Let $R\rarrow S$ be a homomorphism of
commutative rings such that the induced morphism of affine schemes
$\Spec S\rarrow\Spec R$ is an open embedding.
 Assume that $R\in\sR$ (then $S\in\sR$
by Lemma~\ref{rings-affine-communication}(a)); and let $N\in\sL_S$
be a module/complex.
 Then the object $N$, viewed as a module/complex over the ring $R$,
belongs to~$\sL_R$.
\end{description}

 A similar direct image condition was stated in~\cite[Section~2]{Pal}
for \emph{principal} affine open subschemes $\Spec S\sub\Spec R$,
that is, $S=R[r^{-1}]$ for some element $r\in R$.
 We do \emph{not} know whether the two versions of the direct image
condition are equivalent.
 Let us call the direct image condition restricted to principal
affine open subschemes $\Spec S=\Spec R[r^{-1}]\subset\Spec R$
the \emph{weak direct image condition}.
 The following lemma provides a partial comparison result.

\begin{lem} \label{weak-direct-image-for-local}
 Let\/ $\sR$ be a local class of commutative rings, $\sE$ be a very
local class of modules/complexes over rings $R\in\sR$, and\/
$\sL\sub\sE$ be a subclass such that the class\/ $\sL$ is local.
 Assume that, for every $R\in\sR$, the class of modules/complexes\/
$\sL_R$ is closed under finite direct sums and contains the kernel
of every epimorphism in\/ $\sK_R$ between objects from\/ $\sL_R$
such that the kernel belongs to\/~$\sE_R$.
 Assume further that the class\/ $\sL$ satisfies the weak direct image
condition; so, for any ring $R\in\sR$, any element $r\in R$, and any
module/complex $N\in\sL_{R[r^{-1}]}$, the underlying
$R$\+module/complex of $N$ belongs to~$\sL_R$.
 Then the class\/ $\sL$ satisfies the direct image condition; so,
it is very local.
\end{lem}

\begin{proof}
 Let $f\:R\rarrow S$ be a homomorphism of commutative rings such that
the induced morphism $\Spec S\rarrow\Spec R$ is an open embedding.
 Then the open subset $\Spec S$ can be covered by a finite collection
of principal affine open subschemes in $\Spec R$.
 So there exists a finite set of elements $r_1$,~\dots, $r_d\in R$
such that the localization homomorphism $R\rarrow R[r_j^{-1}]$
factorizes as $R\rarrow S\rarrow R[r_j^{-1}]$ for every
$j=1$,~\dots,~$d$, and the elements $f(r_1)$,~\dots, $f(r_d)$ generate
the unit ideal in~$S$.
 Obviously, the ring map $f\:R\rarrow S$ induces isomorphisms
$R[r_j^{-1}]\rarrow S[f(r_j)^{-1}]$ for all~$j$.

 Now assume that $R\in\sR$, and let $N\in\sL_S$ be a module or complex.
 Since the class $\sL$ is local by assumption, we have
$N[r_{i_1}^{-1},\dotsc,r_{i_k}^{-1}]\in\sL_{S[f(r_{i_1})^{-1},\dotsc,
f(r_{i_k})^{-1}]}$ for every subset of indices
$1\le i_1<\dotsb<i_k\le d$.
 When $k\ge1$, we have $S[f(r_{i_1})^{-1},\dotsc,f(r_{i_k})^{-1}]=
R[r_{i_1}^{-1},\dotsc,r_{i_k}^{-1}]$; so
$N[r_{i_1}^{-1},\dotsc,r_{i_k}^{-1}]\in\sL_{R[r_{i_1}^{-1},\dotsc,
r_{i_k}^{-1}]}$.
 By the weak direct image assumption of the lemma, it follows that
$N[r_{i_1}^{-1},\dotsc,r_{i_k}^{-1}]\in\sL_R$.

 It remains to observe that the \v Cech
coresolution~\eqref{cech-modules} of the $S$\+module/complex $N$,
\begin{multline} \label{cech-principal-over-S} \textstyle
 0\lrarrow N \lrarrow \bigoplus_j N[r_j^{-1}] \lrarrow
 \bigoplus_{i<j} N[r_i^{-1},r_j^{-1}] \\ \lrarrow
 \dotsb\lrarrow N[r_1^{-1},\dotsc,r_d^{-1}] \lrarrow 0,
\end{multline}
splits naturally as a complex (of complexes) of
$S[f(r_l)^{-1}]$\+modules after the functor $S[f(r_l)^{-1}]\ot_S{-}$
is applied to it for any index $l=1$,~\dots,~$d$.
 In fact, the localization splits into a direct sum of contractible
two-term complexes with the terms $N[r_{i_1}^{-1},\dotsc,r_{i_k}^{-1}]$,
where $l\in\{i_1,\dotsc,i_k\}$.
 So it follows from the finite direct sum, locality, and weak direct
image assumptions about the class $\sL$ that all (the terms and)
the modules/complexes of cocycles of~\eqref{cech-principal-over-S}
belong to~$\sL_S$.
 (Let us emphasize that $\Spec S$, but not $\Spec R$, is the union of
$\Spec R[r_j^{-1}]$ over $1\le j\le d$.)

 Hence the modules/complexes of cocycles
of the \v Cech coresolution~\eqref{cech-principal-over-S} also belong to
$\sE_S$, and therefore to~$\sE_R$ (as the class $\sE$ is assumed to be
very local).
 Now closedness of the class $\sL_R$ under the kernels of those
epimorphisms in $\sK_R$ whose kernels belong to $\sE_R$ allows to prove
by induction, moving from the rightmost to the leftmost end of
the acyclic complex~\eqref{cech-principal-over-S}, that $N\in\sL_R$.
\end{proof}

\begin{lem} \label{very-local-open-direct-image}
 Let\/ $\sR$ be a local class of commutative rings, $X$ be
a locally\+$\sR$ scheme, and $Y\subset X$ be an open subscheme
whose open embedding morphism $j\:Y\rarrow X$ is affine.
 Let\/ $\sL$ be a very local class of modules or complexes over
rings $R\in\sR$.
 Then the direct image functor $j_*\:Y\qcoh\rarrow X\qcoh$ takes\/
locally\+$\sL$ quasi-coherent sheaves or complexes to locally\+$\sL$
quasi-coherent sheaves/complexes.
\end{lem}

\begin{proof}
 Obvious from the definitions.
\end{proof}

 The following lemma is useful for application of set-theoretical
methods.

\begin{lem} \label{locally-L-sheaves-deconstructible}
 Let\/ $\sR$ be a local class of commutative rings and\/ $\sL$ be
a local class of modules or complexes over rings $R\in\sR$.
 Assume that, for every ring $R\in\sR$, the class of modules/complexes\/
$\sL_R$ is deconstructible in the abelian category\/ $\sK_R$
(in the sense of Section~\ref{small-object-argument-subsect}).
 Then, for any locally\+$\sR$ scheme $X$, the class of quasi-coherent
sheaves or complexes\/ $\sL_X$ is deconstructible in the abelian
category\/~$\sK_X$.
\end{lem}

\begin{proof}
 More specifically, let $\kappa$~be a regular cardinal.
 An $R$\+module is said to be \emph{$\kappa$\+presented} if it has
less than~$\kappa$ generators and relations (in some suitable
presentation).
 One can check that the class of all $\kappa$\+presented modules is
local over all commutative rings (cf.~\cite[Lemma~3.1]{PS6}).
 For any ring $R\in\sR$, denote by $\sL_R^{<\kappa}$ the class of
all $\kappa$\+presented modules from $\sL_R$ or complexes from
$\sL_R$ whose terms are $\kappa$\+presented modules.
 Then $\sL^{<\kappa}$ is a local class.

 Assume that $\kappa$~is uncountable, and that for every $R\in\sR$,
all modules/complexes from $\sL_R$ are filtered by modules/complexes
from~$\sL_R^{<\kappa}$.
 Assume that the scheme $X$ is a union of less than~$\kappa$ affine
open subschemes whose pairwise intersections are also unions of less
than~$\kappa$ open affines.
 Then all quasi-coherent sheaves or complexes from $\sL_X$ are
filtered by locally\+$\sL^{<\kappa}$ quasi-coherent sheaves or
complexes.
 This is a straightforward corollary of the \emph{Hill lemma} for
filtered $R$\+modules~\cite[Theorem~7.10]{GT};
see~\cite[Theorem~2.1]{Sto-hill} for a Grothendieck category version.
\end{proof}

\begin{exs} \label{local-classes-examples}
 Let $\sR$ be the local class of all commutative rings.
 Obviously, over all commutative rings $R$, the class of all modules
$\sL_R=R\modl$ is very local.
 Similarly, over all commutative rings $R$, the class of flat modules
$\sL_R=R\modl_\fl$ is well-known to be very local
(cf.~\cite[Example~2.3]{Pal}).

 Over all commutative rings $R$, the class of very flat modules
$\sL_R=R\modl_\vfl$ is local by Lemma~\ref{very-open-covering}(a),
and it satisfies the direct image condition by
Lemmas~\ref{very-scalars-veryflat-case}(b)
and~\ref{very-open-embedding}.
 So the class of very flat modules is very local
(cf.~\cite[Example~2.5]{Pal}).

 The class of projective modules is local over all commutative rings,
but \emph{not} very local~\cite[Section~II.3.1]{RG}, \cite{Pe}
(cf.~\cite[Example~2.6]{Pal}).
 The classes of contraadjusted and cotorsion modules are \emph{not}
local~\cite[Example~2.7]{Pal}.

 The class of injective modules is local over Noetherian rings
(i.~e., for the local class $\sR$ of all Noetherian commutative rings)
\cite[Lemma~II.7.16 and Theorem~II.7.18]{Har}, but \emph{not}
over commutative rings in general~\cite[Example~2.8]{Pal}.
 Injective modules also satisfy the direct image condition over all
commutative rings (by Lemma~\ref{cotors-coexten}(b)), so over
Noetherian commutative rings the class of injective modules is
very local.

 The class of homotopy injective complexes of injective modules is
\emph{not} local even over Noetherian commutative
rings~\cite[Example~6.5]{N-bb}, \cite{Bel}, \cite[Example~2.9]{Pal}.

 A variety of examples of very local classes are mentioned
throughout~\cite[Section~6]{Pal}.
\end{exs}

\subsection{Gluing cotorsion pairs in quasi-coherent sheaves}
\label{gluing-cotorsion-in-qcoh-subsect}
 The following construction, which we use throughout much of
Chapter~\ref{becker-on-qcomp-qsep-sect}, is a common generalization of
Section~\ref{quasi-compact-quasi-coherent} and~\cite[Section~4]{Pal}.
 The idea of plugging whole complexes into the construction of
Section~\ref{quasi-compact-quasi-coherent} comes
from~\cite[Lemma~6.3]{ES}.

 Let $\sR$ be a local class of commutative rings.
 Suppose given, for every ring $R\in\sR$, a full subcategory of modules
or complexes $\sE_R\subset\sK_R$ such that $\sE_R$ inherits an exact
category structure from the abelian exact structure of~$\sK_R$ (in
the sense of the definition in Section~\ref{fully-faithful-subsect}).
 Assume further that $\sE$ is a local class of modules/complexes over
rings $R\in\sR$ (in the sense of Section~\ref{local-classes-subsect}).

\begin{lem} \label{locally-E-qcoh-inherit-exact}
 For any locally\+$\sR$ scheme $X$, the full subcategory\/ $\sE_X
\sub\sK_X$ inherits an exact category structure from the abelian exact
structure of\/~$\sK_X$.
\end{lem}

\begin{proof}
 Follows from the characterization of full subcategories inheriting
exact category structures given in~\cite[Theorem~2.6]{DS}
or~\cite[Lemma~4.21]{Pedg}.
 It is helpful to keep in mind that the functors of sections over
affine open subschemes $U\sub X$ are exact on the abelian
category~$\sK_X$.
\end{proof}

 For the rest of this section we will assume that the class of
modules/complexes $\sE$ is very local (as defined in
Section~\ref{local-classes-subsect}).

\begin{thm} \label{quasi-coherent-gluing-theorem}
 Assume that we are given, for every ring $R\in\sR$, a hereditary
complete cotorsion pair $(\sF_R,\sC(R))$ in the exact
category\/~$\sE_R$.
 Assume further that the class of modules/complexes\/ $\sF_R$ is
very local for rings $R\in\sR$.
 For any locally\+$\sR$ affine scheme $U$, denote by\/ $\sC(U)\sub\sK_U$
the full subcategory corresponding to\/ $\sC(R)\sub\sK_R$ under
the natural equivalence\/ $\sK_U\simeq\sK_R$ for $R=\O(U)$.

 Let $X$ be a locally\+$\sR$ quasi-compact semi-separated scheme
and $X=\bigcup_\alpha U_\alpha$ be its finite affine open covering.
 As above, we denote by\/ $\sF_X\sub\sK_X$ the class of all
locally\+$\sF$ quasi-coherent sheaves or complexes on~$X$.
 Then there exists a hereditary complete cotorsion pair
$(\sF_X,\sC(X))$ in\/~$\sE_X$.
 Moreover, a quasi-coherent sheaf or complex\/ $\cP\in\sK_X$ belongs
to\/ $\sC(X)$ if and only if\/ $\cP$ is a direct summand of a finitely
iterated extension of direct images $j_\alpha{}_*\cP_\alpha$ of
quasi-coherent sheaves/complexes\/ $\cP_\alpha\in\sC(U_\alpha)$
in the exact category\/~$\sE_X$.
 Here $j_\alpha\:U_\alpha\rarrow X$ are the open embedding morphisms.

 A quasi-coherent sheaf or complex\/ $\G\in\sK_X$ belongs to\/
$\sF_X\cap\sC(X)$ if and only if\/ $\G$ is a direct summand of
a finite direct sum of direct images $j_\alpha{}_*\G_\alpha$ of
quasi-coherent sheaves/complexes\/ $\G_\alpha\in\sF_{U_\alpha}
\cap\sC(U_\alpha)$.
\end{thm}

 Let us introduce an auxiliary notation $\sC'_{\{U_\alpha\}}(X)$ for
the class of all finitely iterated extensions in $\sE_X$ of direct
images $j_\alpha{}_*\cP_\alpha$ of quasi-coherent sheaves/complexes\/
$\cP_\alpha\in\sC(U_\alpha)$, as in the theorem.
 The proof of the first two assertions of the theorem is based on
a lemma and two propositions.

\begin{lem} \label{quasi-coherent-gluing-orthogonality}
 Assume that we are given, for every ring $R\in\sR$, a pair of classes
of objects $(\sF_R,\sC(R))$ in the exact category\/ $\sE_R$ such that\/
$\Ext_{\sE_R}^1(F,C)=0$ for all $F\in\sF_R$ and $C\in\sC(R)$.
 Assume further that the class of modules/complexes\/ $\sF$ is local.
 Let $X$ be a locally\+$\sR$ quasi-compact semi-separated scheme
and $X=\bigcup_\alpha U_\alpha$ be its finite affine open covering.
 Then, for any\/ $\F\in\sF_X$ and\/ $\cP\in\sC'_{\{U_\alpha\}}(X)$,
one has\/ $\Ext_{\sE_X}^1(\F,\cP)=0$.
\end{lem}

\begin{proof}
 It suffices to show that $\Ext_{\sE_X}^1(\F,j_*\cQ)=0$ for any
affine open subscheme $U\sub X$ with the open embedding morphism
$j\:U\rarrow X$ and any quasi-coherent sheaves/complexes
$\F\in\sF_X$ and $\cQ\in\sC(U)$.
 Quite generally, for any locally\+$\sR$ scheme $X$ and any open
subscheme $Y\sub X$ whose open embedding morphism $j\:Y\rarrow X$
is affine, a natural isomorphism $\Ext_{\sE_X}^1(\F,j_*\G)\simeq
\Ext_{\sE_Y}^1(j^*\F,\G)$ holds for all $\F\in\sE_X$ and $\G\in\sE_Y$
by~\cite[Lemma~1.7(a)]{Pal} or~\cite[Lemma~6.1]{PS6}.
\end{proof}

\begin{prop} \label{quasi-coherent-gluing-precover}
 Assume that we are given, for every ring $R\in\sR$, a pair of classes
of objects $(\sF_R,\sC(R))$ in the exact category\/ $\sE_R$ such that
an (admissible) short exact sequence~\eqref{special-precover-sequence}
with $F\in\sF_R$ and $C'\in\sC(R)$ exists in\/ $\sE_R$ for every
object $E\in\sE_R$.
 Assume further that the class of modules/complexes\/ $\sF$ is very
local, and the class\/ $\sF_R$ is closed under extensions and kernels
of admissible epimorphisms in\/ $\sE_R$ for every $R\in\sR$.
 Let $X$ be a locally\+$\sR$ quasi-compact semi-separated scheme
and $X=\bigcup_\alpha U_\alpha$ be its finite affine open covering.
 Then any object\/ $\M\in\sE_X$ can be included in a short exact
sequence\/ $0\rarrow\cP\rarrow\F\rarrow\M\rarrow0$ in\/ $\sE_X$
with\/ $\F\in\sF_X$ and\/ $\cP\in\sC'_{\{U_\alpha\}}(X)$.
\end{prop}

\begin{proof}
 This is a straightforward generalization of
Lemmas~\ref{quasi-very-flat-cover} and~\ref{quasi-flat-cover},
based on the same construction from~\cite[Section~A.1]{EP}.
 The argument proceeds by induction on the number of affine open
subschemes covering~$X$.
 Assume that for some open subscheme $V\sub X$ with the open
embedding morphism $h\:V\rarrow X$ there is a short exact sequence
$0\rarrow\cQ\rarrow\cH\rarrow\M\rarrow0$ in $\sE_X$ such that
the quasi-coherent sheaf/complex $h^*\cH$ on $V$ belongs to $\sF_V$,
while $\cQ\in\sC'_{\{U_\alpha\}}(X)$.
 Pick an index~$\beta$ such that $U=U_\beta$ is not contained in~$V$;
we will construct a short exact sequence $0\rarrow\cP\rarrow\F\rarrow
\M\rarrow0$ in $\sE_X$ having the same properties with respect to
the open subscheme $U\cup V\sub X$.

 By assumption, there exists a short exact sequence $0\rarrow\cR\rarrow
\G\rarrow j^*\cH\rarrow0$ in $\sE_U$ with $\G\in\sF_U$ and $\cR\in
\sC(U)$, where $j=j_\beta\:U\rarrow X$ is the open embedding.
 The direct image $0\rarrow j_*\cR\rarrow j_*\G\rarrow j_*j^*\cH
\rarrow0$ of this short exact sequence under the affine morphism~$j$
is a short exact sequence in $\sK_X$ with the terms in $\sE_X$
(by Lemma~\ref{very-local-open-direct-image}), so it is admissible
exact in~$\sE_X$.
 Taking the pull-back of this short exact sequence in $\sE_X$ under
the adjunction morphism $\cH\rarrow j_*j^*\cH$, we obtain a short
exact sequence $0\rarrow j_*\cR\rarrow\F\rarrow\cH\rarrow0$ in~$\sE_X$.

 By the definition, in order to show that the restriction of $\F$ to
$U\cup V$ is locally\+$\sF$, it suffices to check that so are
the restrictions of $\F$ to $U$ and~$V$.
 We have $j^*\F\simeq\G$, which is locally\+$\sF$ by construction;
so $\F$ is locally\+$\sF$ over~$U$.
 On the other hand, the sheaf/complex $\cH$ is locally\+$\sF$ over $V$,
and consequently over $V\cap U$, by the induction assumption.
 Hence the sheaf/complex $\cR$ is also locally\+$\sF$ over $V\cap U$ as
the kernel of an admissible epimorphism $\G\rarrow j^*\cH$ in $\sE_U$.
 We are using the assumptions of the ascent condition for the classes
$\sE$ and $\sF$ here.
 The open embedding $U\cap V\rarrow V$ is an affine morphism, so
the sheaf/complex $j_*\cR$ is locally\+$\sF$ over~$V$
by Lemma~\ref{very-local-open-direct-image}.
 Now the admissible short exact sequence $0\rarrow j_*\cR\rarrow
\F\rarrow\cH\rarrow0$ implies that the sheaf/complex $\F$ is
locally\+$\sF$ over~$V$.

 Finally, the kernel $\cP$ of the composition of admissible epimorphisms
$\F\rarrow\cH\rarrow\M$ in $\sE_X$ is an extension of
the sheaves/complexes $\cQ$ and $j_*\cR$ in~$\sE_X$.
 It remains to recall that $U=U_\beta$, \,$j=j_\beta$, and
$\cR\in\sC(U)$; while $\cQ\in\sC'_{\{U_\alpha\}}(X)$.
 Hence $\cP\in\sC'_{\{U_\alpha\}}(X)$.
\end{proof}

\begin{prop} \label{quasi-coherent-gluing-preenvelope}
 Assume that we are given, for every ring $R\in\sR$, a pair of classes
of objects $(\sF_R,\sC(R))$ in the exact category\/ $\sE_R$ such that
an (admissible) short exact sequence~\eqref{special-precover-sequence}
with $F\in\sF_R$ and $C'\in\sC(R)$ exists in\/ $\sE_R$ for every
object $E\in\sE_R$.
 Assume further that the class of modules/complexes\/ $\sF$ is very
local, the class\/ $\sF_R$ is closed under extensions and kernels
of admissible epimorphisms in\/ $\sE_R$ for every $R\in\sR$, and
the class\/ $\sC(R)$ is cogenerating in\/ $\sE_R$ for every $R\in\sR$.
 Let $X$ be a locally\+$\sR$ quasi-compact semi-separated scheme
and $X=\bigcup_\alpha U_\alpha$ be its finite affine open covering.
 Then any object\/ $\M\in\sE_X$ can be included in a short exact
sequence\/ $0\rarrow\M\rarrow\cP\rarrow\F\rarrow0$ in\/ $\sE_X$
with\/ $\F\in\sF_X$ and\/ $\cP\in\sC'_{\{U_\alpha\}}(X)$.
\end{prop}

\begin{proof}
 Let us show that the class of all quasi-coherent sheaves/complexes
of the form $\bigoplus_\alpha j_\alpha{}_*\cP_\alpha$, where
$\cP_\alpha\in\sC(U_\alpha)$, is cogenerating in~$\sE_X$.
 In view of Proposition~\ref{quasi-coherent-gluing-precover}
and Lemma~\ref{salce-lemma}(b), this will imply the assertion of
the proposition.

 Indeed, given an object $\M\in\sE_X$, pick admissible monomorphisms
$j_\alpha^*\M\rarrow\cP_\alpha$ in $\sE_{U_\alpha}$ with
$\cP_\alpha\in\sC(U_\alpha)$ for all indices~$\alpha$.
 Let us check that the induced morphism $\M\rarrow\bigoplus_\alpha
j_\alpha{}_*\cP_\alpha$ is an admissible monomorphism in~$\sE_X$.

 Both the morphisms $\M\rarrow\bigoplus_\alpha j_\alpha{}_*j_\alpha^*\M$
and $\bigoplus_\alpha j_\alpha{}_*j_\alpha^*\M\rarrow
\bigoplus_\alpha j_\alpha{}_*\cP_\alpha$ are monomorphisms in~$\sK_X$;
hence so is their composition.
 Denote by $\N$ the cokernel of the monomorphism
$\M\rarrow\bigoplus_\alpha j_\alpha{}_*\cP_\alpha$ in~$\sK_X$.
 It remains to show that $\N\in\sE_X$.

 By the definition of $\sE_X$, we only need to check that
$j_\beta^*\N\in\sE_{U_\beta}$ for every index~$\beta$.
 Now we have a short exact sequence $0\rarrow j_\beta^*\M\rarrow
\bigoplus_\alpha j_\beta^*j_\alpha{}_*\cP_\alpha\rarrow j_\beta^*\N
\rarrow0$ in $\sK_{U_\beta}$ with the leftmost and middle terms
in~$\sE_{U_\beta}$.
 The direct summand projection map $\bigoplus_\alpha
j_\beta^*j_\alpha{}_*\cP_\alpha\rarrow j_\beta^*j_\beta{}_*\cP_\beta
=\cP_\beta$ is a (split) admissible epimorphism in $\sE_{U_\beta}$,
while the composition $j_\beta^*\M\rarrow\bigoplus_\alpha
j_\beta^*j_\alpha{}_*\cP_\alpha\rarrow\cP_\beta$ is an admissible
monomorphism in $\sE_{U_\beta}$ by construction.
 Applying K\"unzer's axiom~\cite[Exercise~3.11(i)]{Bueh}, we
conclude that $j_\beta^*\M\rarrow\bigoplus_\alpha
j_\beta^*j_\alpha{}_*\cP_\alpha$ is an admissible monomorphism
in~$\sE_{U_\beta}$.
 Hence the cokernel of the latter morphism in the category
$\sK_{U_\beta}$ belongs to~$\sE_{U_\beta}$.
\end{proof}

\begin{proof}[Proof of
Theorem~\ref{quasi-coherent-gluing-theorem}]
 In order to prove the first two assertions of the theorem, we apply
Lemma~\ref{cotorsion-pair-direct-summands-lemma} to the exact
category $\sE=\sE_X$ with two classes of objects $\sF=\sF_X$ and
$\sC=\sC'_{\{U_\alpha\}}(X)$.
 Lemma~\ref{quasi-coherent-gluing-orthogonality} establishes
the $\Ext^1$\+orthogonality assumption of
Lemma~\ref{cotorsion-pair-direct-summands-lemma} in this context,
while Propositions~\ref{quasi-coherent-gluing-precover}\+-%
\ref{quasi-coherent-gluing-preenvelope} provide the short exact
sequences~(\ref{special-precover-sequence}\+-%
\ref{special-preenvelope-sequence}).
 Lemma~\ref{cotorsion-pair-direct-summands-lemma} tells us that the pair
of classes $\sF_X^\oplus$ and $\sC(X)=(\sC'_{\{U_\alpha\}}(X))^\oplus$
is a complete cotorsion pair in~$\sE_X$; and it remains to observe
that the class $\sF_X$ is closed under direct summands in $\sE_X$,
because the class $\sF_R$ is closed under direct summands in $\sE_R$
for every ring $R\in\sR$.
 The cotorsion pair $(\sF_X,\sC(X))$ in $\sE_X$ is hereditary,
because the class $\sF_X$ is closed under kernels of admissible
epimorphisms in~$\sE_X$ (because the class $\sF_R$ is closed under
kernels of admissible epimorphisms in $\sE_R$, since the cotorsion
pair $(\sF_R,\sC(R))$ in $\sE_R$ is hereditary by assumption).

 The proof of the last assertion of the theorem is similar;
cf.\ Lemmas~\ref{very-flat-contraadjusted-quasi}
and~\ref{flat-cotorsion-quasi}.
 The ``if'' assertion is provided by
Lemma~\ref{very-local-open-direct-image} for the class $\sF$ and
by the description of the class $\sC(X)$ that we have obtained already.
 Let us prove the ``only if''.

 Let $\F\in\sF_X$ be a locally\+$\sF$ quasi-coherent sheaf/complex
on~$X$.
 For every index~$\alpha$, pick a special preenvelope
sequence~\eqref{special-preenvelope-sequence}
\,$0\rarrow j_\alpha^*\F\rarrow\G_\alpha\rarrow\cH_\alpha\rarrow0$
in the complete cotorsion pair $(\sF_{U_\alpha},\sC(U_\alpha))$
in the exact category~$\sE_{U_\alpha}$.
 So $\G_\alpha\in\sC(U_\alpha)$ and $\cH_\alpha\in\sF_{U_\alpha}$.
 Since $j_\alpha^*\F\in\sF_{U_\alpha}$, it follows that
$\G_\alpha\in\sF_{U_\alpha}\cap\sC(U_\alpha)$.

 Now consider the induced morphism $\F\rarrow\bigoplus_\alpha
j_\alpha{}_*\G_\alpha$ in $\sF_X\sub\sE_X$.
 Similarly to the proof of
Proposition~\ref{quasi-coherent-gluing-preenvelope}, one shows that
this morphism is an admissible monomorphism in $\sE_X$, and in fact,
by the same argument, even in (the inherited exact category
structure on)~$\sF_X$.
 We recall that $\sF$ is a local class, just like~$\sE$.

 Thus we have a short exact sequence $0\rarrow\F\rarrow\bigoplus_\alpha
j_\alpha{}_*\G_\alpha\rarrow\cH\rarrow0$ in~$\sF_X$.
 If $\F\in\sF_X\cap\sC(X)$, then $\Ext^1_{\sE_X}(\cH,\F)=0$; so
this short exact sequence splits and $\F$ is a direct summand of
$\bigoplus_\alpha j_\alpha{}_*\G_\alpha$.
\end{proof}

\begin{rem} \label{quasi-coherent-gluing-nonuniversal-remark}
 Theorem~\ref{quasi-coherent-gluing-theorem} is formulated above
for ``universal'' cotorsion pairs, defined uniformly for all
commutative rings $R$ from a certain local class~$\sR$.
 In practice, one often encounters cotorsion pairs in $R\modl$
depending on an additional datum for a particular ring~$R$, such as
the datum of a specific $R$\+module; or cotorsion pairs in $X\qcoh$
depending on an additional datum over the scheme~$X$.
 It is easy to extend Theorem~\ref{quasi-coherent-gluing-theorem}
to such contexts.
 One needs to fix a quasi-compact semi-separated scheme $X$, and
assume given an exact subcategory of modules or complexes
$\sE_U\sub\sK_{\O(U)}$ for every affine open subscheme $U\sub X$
(where $\sK_R=R\modl$ or $\Com(R\modl)$, while $\sE_U$ depends on
the embedding of $U$ into $X$ and \emph{not} only on the abstract
affine scheme~$U$) together with a hereditary complete cotorsion
pair $(\sF_U,\sC(U))$ in~$\sE_U$.
 The classes $\sE_U$ and $\sF_U$ must be local with respect to
affine open coverings of affine open subschemes $U\sub X$ and
preserved by direct images with respect to identity embeddings of
affine open subschemes $V\rarrow U$, where $V\sub U\sub X$.
 Then the same arguments and constructions as above produce
an exact subcategory $\sE_X\sub\sK_X$, where $\sK_X=X\qcoh$ or
$\Com(X\qcoh)$, and a hereditary complete cotorsion pair
$(\sF_X,\sC(X))$ in $\sE_X$, with the same description of
the classes $\sE_X$, \,$\sF_X$, and~$\sC(X)$.
\end{rem}

\subsection{Colocal classes of modules and complexes}
\label{colocal-classes-subsect}
 In this section we largely follow~\cite[Section~3]{Pal}.
 Colocality is dual to locality in that one uses the colocalization
functors $\Hom_R(R[r^{-1}],{-})$ instead of the localization functors
$R[r^{-1}]\ot_R{-}$.
 As illustrated by~\cite[Example~3.1]{Pal}, any meaningful discussion
of colocality has to presume the modules under consideration to be
at least contraadjusted.

 Given a ring $R$, we denote by $\sK^R=\sK_R^\cta\sub\sK_R$ either
the exact category of contraadjusted $R$\+modules $\sK^R=R\modl^\cta$
(with the exact structure inherited from the abelian exact
structure of $R\modl$, as usual), or the exact category of complexes of
contraadjusted $R$\+modules $\sK^R=\Com(R\modl^\cta)$ (with the exact
structure inherited from the abelian structure of $\Com(R\modl)$, or
which is the same, the termwise exact structure of the category of
complexes over $R\modl^\cta$).
 This is our big ambient exact category of modules/complexes.

 Let $R$ be a local class of rings.
 Suppose given a class of objects $\sL^R\sub\sK^R$ defined for every
ring $R\in\sR$.
 As above, we will abuse terminology and notation by speaking of
``a~class of modules/complexes~$\sL$\,'' (over all rings $R\in\sR$).
 A class $\sL$ is called \emph{colocal} if it satisfies the following
two conditions:
\begin{description}
\item[Coascent] For any ring $R\in\sR$, any element $r\in R$, and any
module/complex $M\in\sL^R$, the module/complex $\Hom_R(R[r^{-1}],M)$
belongs to~$\sL^{R[r^{-1}]}$.
\item[Codescent] Let $R\in\sR$ be a ring and $r_1$,~\dots, $r_d\in R$
be a finite set of elements generating the unit ideal in~$R$.
 Let $M\in\sK^R=\sK_R^\cta$ be a module/complex such that
$\Hom_R(R[r_j^{-1}],M)\in\sL^{R[r_j^{-1}]}$ for every $j=1$,~\dots,~$d$.
 Then $M\in\sL^R$.
\end{description}

\begin{lem} \label{colocality-affine-communication}
 Let\/ $\sR$ be a local class of commutative rings and $U$ be
a locally\+$\sR$ affine scheme.
 Let\/ $\sL$ be a colocal class of modules/complexes over rings
$R\in\sR$, and let\/ $\gM$ be a contraherent cosheaf or complex of
contraherent cosheaves over~$U$.
 Then \par
\textup{(a)} for any affine open subscheme $V\sub U$ in the affine
scheme $U$, one has\/ $\gM[V]\in\sL^{\O(V)}$ provided that\/
$\gM[U]\in\sL^{\O(U)}$; \par
\textup{(b)} for any affine open covering $U=\bigcup_\alpha V_\alpha$
of the affine scheme $U$, one has\/ $\gM[U]\in\sL^{\O(U)}$ provided
that\/ $\gM[V_\alpha]\in\sL^{\O(V_\alpha)}$ for every~$\alpha$.
\end{lem}

\begin{proof}
 Similar to Lemmas~\ref{rings-affine-communication}
and~\ref{locality-affine-communication}.
\end{proof}

 Let $\sR$ be a local class of commutative rings and $\sL$ be
a colocal class of modules or complexes over rings $R\in\sR$.
 Let $X$ be a locally\+$\sR$ scheme, $\bW$ be an open covering of $X$,
and $\gM$ be a $\bW$\+locally contraherent cosheaf/complex of
$\bW$\+locally contraherent cosheaves on~$X$.
 We will say that the cosheaf/complex $\gM$ is \emph{locally\+$\sL$}
if one has $\gM[U]\in\sL^{\O_X(U)}$ for every affine open subscheme
$U\sub X$ subordinate to~$\bW$.
 It is clear from Lemma~\ref{colocality-affine-communication} that it
suffices to check this condition for affine open subschemes belonging
to any given affine open covering of $X$ subordinate to~$\bW$.

 For any scheme $X$ with an open covering $\bW$, put
$\sK^X_\bW=X\lcth_\bW$ if we are dealing with modules (i.~e., if
$\sK^R=R\modl^\cta$), or $\sK^X_\bW=\Com(X\lcth_\bW)$ if we are dealing
with complexes of modules (i.~e., if $\sK^R=\Com(R\modl^\cta)$).
 So $\sK^X_\bW$ is an exact category with the exact structure defined
in Section~\ref{locally-contraherent} in the case of
$\sK^X_\bW=X\lcth_\bW$; or in the case of $\sK^X_\bW=\Com(X\lcth_\bW)$,
the termwise exact structure on complexes induced by the exact structure
on $X\lcth_\bW$ defined in Section~\ref{locally-contraherent}.

 Given a local class of commutative rings $\sR$, a colocal class of
modules/complexes $\sL$ over rings $R\in\sR$, and a locally\+$\sR$
scheme $X$ with an open covering $\bW$, we will denote by
$\sL^X_\bW\subset\sK^X_\bW$ the full subcategory of locally\+$\sL$
\,$\bW$\+locally contraherent cosheaves or complexes.
 For an affine scheme $U$, we put $\sK^U=\sK^U_{\{U\}}\simeq\sK^R$
and $\sL^U=\sL^U_{\{U\}}\simeq\sL^R$, where $R=\O(U)$.
 The open covering $\{U\}$ of $U$ consisting of the single open
subscheme $U\sub U$ is presumed here.

 We will say that a class of modules/complexes $\sL$ over rings
$R\in\sR$ is \emph{very colocal} if $\sL$ is colocal and satisfies
the following condition very similar to the condition from
Section~\ref{local-classes-subsect} under the same name.
\begin{description}
\item[Direct image condition] Let $R\rarrow S$ be a homomorphism of
commutative rings such that the induced morphism of affine schemes
$\Spec S\rarrow\Spec R$ is an open embedding.
 Assume that $R\in\sR$ (then $S\in\sR$
by Lemma~\ref{rings-affine-communication}(a)); and let $N\in\sL^S$
be a module/complex.
 Then the object $N$, viewed as a module/complex over the ring $R$,
belongs to~$\sL^R$.
\end{description}

 As mentioned in Section~\ref{local-classes-subsect}, we do \emph{not}
know whether it suffices to check the direct image condition for
principal affine open subschemes, i.~e., for rings $S=R[r^{-1}]$,
where $r\in R$.
 We call the direct image condition restricted to rings $S=R[r^{-1}]$
the \emph{weak direct image condition}.
 The following lemma is a dual/colocal version of
Lemma~\ref{weak-direct-image-for-local}.

\begin{lem} \label{weak-direct-image-for-colocal}
 Let\/ $\sR$ be a local class of commutative rings, $\sE$ be a very
colocal class of modules/complexes over rings $R\in\sR$, and\/
$\sL\sub\sE$ be a subclass such that the class\/ $\sL$ is colocal.
 Assume that, for every $R\in\sR$, the class of modules/complexes\/
$\sL^R$ is closed under finite direct sums and contains the cokernel
of every admissible monomorphism in\/ $\sK^R$ between objects from\/
$\sL^R$ whose cokernel belongs to\/~$\sE^R$.
 Assume further that the class\/ $\sL$ satisfies the weak direct image
condition; so, for any ring $R\in\sR$, any element $r\in R$, and any
module/complex $N\in\sL^{R[r^{-1}]}$, the underlying
$R$\+module/complex of $N$ belongs to~$\sL^R$.
 Then the class\/ $\sL$ satisfies the direct image condition; so,
it is very colocal.
\end{lem}

\begin{proof}
 The argument is dual-analogous to the proof of
Lemma~\ref{weak-direct-image-for-local}.
 In the notation of that proof, consider a module/complex $N\in\sL^S$.
 Then $\Hom_R(R[r_{i_1}^{-1},\dotsc,r_{i_k}^{-1}],N)\in\sL^R$ for
all $k\ge1$.
 The \v Cech resolution~\eqref{cech-contra} of
the $S$\+module/complex~$N$,
\begin{multline} \label{cech-contra-principal-over-S} \textstyle
 0\lrarrow\Hom_S(S[r_1^{-1},\dotsc,r_d^{-1}],N) \lrarrow \dotsb
 \\ \textstyle \lrarrow\bigoplus_{i<j}\Hom_S(S[r_i^{-1},r_j^{-1}],N)
 \\ \textstyle \lrarrow\bigoplus_j \Hom_S(S[r_j^{-1}],N)
 \lrarrow N\lrarrow0
\end{multline}
is an exact complex in the exact category $\sK^S$ (hence also
in~$\sK^R$), splitting naturally as a complex (of complexes) of
$S[f(r_l)^{-1}]$\+modules after the functor
$\Hom_S(S[f(r_l)^{-1}],{-})$ is applied to it, etc.
\end{proof}

\begin{lem} \label{very-colocal-open-direct-image}
 Let\/ $\sR$ be a local class of commutative rings, $X$ be
a locally\+$\sR$ scheme with an open covering\/ $\bW$, and
$Y\subset X$ be an open subscheme whose open embedding morphism
$j\:Y\rarrow X$ is affine, endowed with the induced open covering\/
$\bT=\bW|_Y=\{Y\cap W\mid W\in\bW\}$.
 Let\/ $\sL$ be a very colocal class of modules or complexes over
rings $R\in\sR$.
 Then the direct image functor $j_!\:Y\lcth_\bT\rarrow X\lcth_\bW$
takes\/ locally\+$\sL$ locally contraherent cosheaves or complexes to
locally\+$\sL$ locally contraherent cosheaves/complexes.
\end{lem}

\begin{proof}
 Obvious from the definitions.
\end{proof}

\begin{exs} \label{colocal-classes-examples}
 Let $\sR$ be the local class of all commutative rings.
 Over all commutative rings $R$, the class of contraadjusted modules
$\sL^R=R\modl^\cta$ satisfies coascent by
Lemma~\ref{very-scalars-veryflat-case}(a), and it satisfies codescent
\emph{by the definition}.
 The class of contraadjusted modules also satisfies the direct image
condition by Lemma~\ref{very-scalars-always}(a).
 So, the class of contraadjusted modules is very
colocal~\cite[Example~3.2]{Pal}.

 Over all commutative rings $R$, the class of cotorsion modules
$\sL=R\modl^\cot$ satisfies coascent by
Lemma~\ref{cotors-coexten}(a) and codescent by
Lemma~\ref{cotors-inj-covering}(a).
 The class of cotorsion modules also satisfies the direct image
condition by Lemma~\ref{cotors-restrict}(a).
 So, the class of cotorsion modules is very
colocal~\cite[Example~3.8]{Pal}.

 Over all commutative rings $R$, the class of injective modules
$\sL=R\modl^\inj$ satisfies coascent by
Lemma~\ref{cotors-restrict}(b) and codescent by
Lemma~\ref{cotors-inj-covering}(b).
 The class of injective modules also satisfies the direct image
condition by Lemma~\ref{cotors-coexten}(b).
 So, the class of injective modules is very
colocal~\cite[Example~3.7]{Pal}.

 The class of homotopy injective complexes of injective modules
is likewise very colocal over all commutative
rings~\cite[Example~3.9]{Pal}.
 A variety of other examples of very colocal classes are mentioned
throughout~\cite[Section~7]{Pal}.
\end{exs}

\subsection{Gluing cotorsion pairs in locally contraherent cosheaves}
\label{gluing-cotorsion-in-lcth-subsect}
 The following construction, which we use throughout much of
Chapter~\ref{becker-on-qcomp-qsep-sect}, is a common generalization of
Section~\ref{clp-subsection} and~\cite[Section~5]{Pal}.
 For examples of this kind of construction applied to certain
non-semi-separated schemes $X$, see
Corollary~\ref{coflasque-flat-loc-cotorsion-cotorsion-pair}
and Lemmas~\ref{coflasque-acyclic-preenvelope}\+-%
\ref{coflasque-homotopy-projective-precover}, and particularly
Section~\ref{acyclic-complexes-of-flat-cosheaves-subsect} (which
also provides a discussion in the proofs).

 Let $\sR$ be a local class of commutative rings.
 Suppose given, for every ring $R\in\sR$, a full subcategory of modules
or complexes $\sE^R\sub\sK^R=\sK_R^\cta$ such that $\sE^R$ inherits
an exact category structure from~$\sK^R$.
 Assume further that $\sE$ is a colocal class of modules/complexes
over rings $R\in\sR$.
 The notation in the following lemma, as well as the exact category
structure on\/ $\sK^X_\bW$, was introduced in
Section~\ref{colocal-classes-subsect}.

\begin{lem} \label{locally-E-lcth-inherit-exact}
 For any locally\+$\sR$ scheme $X$ with an open covering\/ $\bW$,
the full subcategory\/ $\sE^X_\bW\sub\sK^X_\bW$ inherits an exact
category structure from\/~$\sK^X_\bW$.
\end{lem}

\begin{proof}
 Similarly to the proof of Lemma~\ref{locally-E-qcoh-inherit-exact},
the argument is based on the abstract category-theoretic
characterization of full subcategories inheriting exact category
structures.
 One should keep in mind that the functors of cosections over affine
open subschemes $U\sub X$ subordinate to $\bW$ are exact on the exact
category~$\sK^X_\bW$.
\end{proof}

 For the rest of this section we will assume that the class of
modules/complexes $\sE$ is very colocal.

\begin{thm} \label{loc-contraherent-gluing-theorem}
 Assume that we are given, for every ring $R\in\sR$, a hereditary
complete cotorsion pair $(\sF(R),\sC^R)$ in the exact
category\/~$\sE^R$.
 Assume further that the class of modules/complexes\/ $\sC^R$ is
very colocal for rings $R\in\sR$.
 For any locally\+$\sR$ affine scheme $U$, denote by\/
$\sF(U)\sub\sK^U$ the full subcategory corresponding to\/
$\sF(R)\sub\sK^R$ under the natural equivalence\/ $\sK^U\simeq\sK^R$
for $R=\O(U)$.

 Let $X$ be a locally\+$\sR$ quasi-compact semi-separated scheme
with an open covering\/ $\bW$ and $X=\bigcup_\alpha U_\alpha$ be
a finite affine open covering of $X$ subordinate to\/~$\bW$.
 As above, we denote by\/ $\sC^X_\bW\sub\sK^X_\bW$ the class of all
locally\+$\sC$ \,$\bW$\+locally contraherent cosheaves or complexes
on~$X$.
 Then there exists a hereditary complete cotorsion pair
$(\sF(X),\sC^X_\bW)$ in\/~$\sE^X_\bW$.
 Moreover, a\/ $\bW$\+locally contraherent cosheaf or complex\/
$\gF\in\sK_X$ belongs to\/ $\sF(X)$ if and only if\/ $\gF$ is a direct
summand of a finitely iterated extension of direct images
$j_\alpha{}_!\gF_\alpha$ of contraherent cosheaves/complexes\/
$\gF_\alpha\in\sF(U_\alpha)$ in the exact category\/~$\sE^X_\bW$.
 Here $j_\alpha\:U_\alpha\rarrow X$ are the open embedding morphisms.

 A $\bW$\+locally contraherent cosheaf or complex\/ $\Q\in\sK^X_\bW$
belongs to\/ $\sF(X)\cap\sC^X_\bW$ if and only if\/ $\Q$ is a direct
summand of a finite direct sum of direct images $j_\alpha{}_*\Q_\alpha$
of contraherent cosheaves/complexes\/ $\Q_\alpha\in\sF(U_\alpha)
\cap\sC^{U_\alpha}$.

 Consequently, all the objects of\/ $\sF(X)$ are (globally) contraherent
cosheaves or complexes of contraherent cosheaves on~$X$.
 Viewed as a full subcategory in $X\lcth$ or\/ $\Com(X\lcth)$,
the subcategory $\sF(X)$ does not depend on the choice of an open
covering\/~$\bW$ of the scheme~$X$.
 Similarly, the full subcategory $\sF(X)\cap\sC^X_\bW$ in $X\lcth$ or\/
$\Com(X\lcth)$ does not depend on the choice of an open covering\/~$\bW$.
\end{thm}

 It is convenient to introduce an auxiliary notation
$\sF'_{\{U_\alpha\}}(X)$ for the class of finitely iterated extensions
in $\sE^X_\bW$ of direct images $j_\alpha{}_!\gF_\alpha$ of
contraherent cosheaves/complexes $\gF_\alpha\in\sF(U_\alpha)$, as
in the theorem.
 Similarly to Section~\ref{gluing-cotorsion-in-qcoh-subsect},
the proof of the first two assertions of the theorem is based on
a lemma and two propositions.

\begin{lem} \label{loc-contraherent-gluing-orthogonality}
 Assume that we are given, for every ring $R\in\sR$, a pair of classes
of objects $(\sF(R),\sC^R)$ in the exact category\/ $\sE^R$ such that\/
$\Ext_{\sE^R}^1(F,C)=0$ for all $F\in\sF(R)$ and $C\in\sC^R$.
 Assume further that the class of modules/complexes\/ $\sC$ is colocal.
 Let $X$ be a locally\+$\sR$ quasi-compact semi-separated scheme with
an open covering\/ $\bW$ and $X=\bigcup_\alpha U_\alpha$ be a finite
affine open covering of $X$ subordinate to\/~$\bW$.
 Then, for any\/ $\gF\in\sF'_{\{U_\alpha\}}(X)$ and\/ $\P\in\sC^X_\bW$,
one has\/ $\Ext_{\sE^X_\bW}^1(\gF,\P)=0$.
\end{lem}

\begin{proof}
 It suffices to show that $\Ext^1_{\sE^X_\bW}(j_!\gG,\P)=0$ for any
affine open subscheme $U\sub X$ subordinate to $\bW$, any
$\bW$\+locally contraherent cosheaf $\P\in\sC^X$, and any
contraherent cosheaf $\gG\in\sF(U)$ (were $j\:U\rarrow X$ denotes
the open embedding morphism).
 Quite generally, for any locally\+$\sR$ scheme $X$ with an open
covering $\bW$, any open subscheme $Y\sub X$ whose open embedding
morphism $j\:Y\rarrow X$ is affine, and the induced open covering
$\bT=\bW|_Y=\{Y\cap W\mid W\in\bW\}$ of $Y$, a natural isomorphism
$\Ext^1_{\sE^X_\bW}(j_!\Q,\P)\simeq\Ext^1_{\sE^Y_\bT}(\Q,j^!\P)$
holds for all $\P\in\sE^X_\bW$ and $\Q\in\sE^Y_\bT$
by~\cite[Lemma~1.7(a)]{Pal} or~\cite[Lemma~6.1]{PS6}.
\end{proof}

\begin{prop} \label{loc-contraherent-gluing-preenvelope}
 Assume that we are given, for every ring $R\in\sR$, a pair of classes
of objects $(\sF(R),\sC^R)$ in the exact category\/ $\sE^R$ such that
an (admissible) short exact
sequence~\eqref{special-preenvelope-sequence} with $C\in\sC^R$ and
$F'\in\sF(R)$ exists in\/ $\sE^R$ for every object $E\in\sE^R$.
 Assume further that the class of modules/complexes\/ $\sC$ is very
colocal, and the class\/ $\sC^R$ is closed under extensions and
cokernels of admissible monomorphisms in\/ $\sE^R$ for every $R\in\sR$.
 Let $X$ be a locally\+$\sR$ quasi-compact semi-separated scheme with
an open covering\/ $\bW$ and $X=\bigcup_{\alpha=1}^N U_\alpha$ be
a finite affine open covering of $X$ subordinate to\/~$\bW$.
 Then any object\/ $\gM\in\sE^X_\bW$ can be included in a short exact
sequence\/ $0\rarrow\gM\rarrow\P\rarrow\gF\rarrow0$ in\/ $\sE^X_\bW$
with\/ $\P\in\sC^X$ and\/ $\gF\in\sF'_{\{U_\alpha\}}(X)$.
\end{prop}

\begin{proof}
 This is a straightforward generalization of Lemmas~\ref{lin-envelope}
and~\ref{lct-envelope} (for other special cases, see
Lemmas~\ref{complexes-loc-cot-preenvelope}
and~\ref{complexes-loc-hot-inj-of-inj-preenvelope}).
 As in the proof of Lemma~\ref{lin-envelope}, we argue by induction
on $1\le\beta\le N$, considering the open subscheme
$V=\bigcup_{\alpha<\beta}U_\alpha$ with the induced covering
$\bT=\bW|_V$ and the identity embedding $h\:V\rarrow X$.
 Assume that we have constructed a short exact sequence $0\rarrow
\gM\rarrow\gK\rarrow\gG\rarrow0$ in $\sE^X_\bW$ such that
the $\bT$\+locally contraherent cosheaf/complex $h^!\gK$ on $V$ is
locally\+$\sC$, while the object $\gG$ is a finitely iterated
extension of the objects $j_\alpha{}_!\gF_\alpha$ for some
$\gF_\alpha\in\sF(U_\alpha)$ with $\alpha<\beta$.
 When $\beta=1$, it suffices to take $\gK=\gM$ and $\gG=0$ for
the induction base.
 Set $U=U_\beta$ and denote by $j\:U\rarrow X$ the open embedding
morphism.

 Let $0\rarrow j^!\gK\rarrow\R\rarrow\gH\rarrow0$ be a short exact
sequence~\eqref{special-preenvelope-sequence} in $\sE^U$ with
$\R\in\sC^U$ and $\gH\in\sF(U)$.
 Applying the direct image functor~$j_!$, which is exact as
a functor $U\ctrh\rarrow X\ctrh\sub X\lcth_\bW$, we obtain a short
exact sequence $0\rarrow j_!j^!\gK\rarrow j_!\R\rarrow j_!\gH\rarrow0$
in $\sE^X_\bW$ (by Lemma~\ref{very-colocal-open-direct-image}).
 Taking the push-out with respect to the adjunction morphism
$j_!j^!\gK\rarrow\gK$, we produce a short exact sequence
$0\rarrow\gK\rarrow\Q\rarrow j_!\gH\rarrow0$ in $\sE^X_\bW$.
 We need to show that $\Q$ is locally\+$\sC$ in restriction to
$U\cup V$, and for this purpose, it suffices to check that
the restrictions of $\Q$ to $U$ and $V$ are locally\+$\sC$.

 Indeed, in the restriction to $U$ we have $j^!j_!j^!\gK\simeq j^!\gK$,
hence $j^!\Q\simeq j^!j_!\R\simeq\R$ is locally\+$\sC$.
 On the other hand, denoting by $j'\:U\cap V\rarrow V$ and
$h'\:U\cap V\rarrow U$ the open embedding morphisms, we have
$h^!j_!\gH\simeq j'_!h'{}^!\gH$ (see the end of
Section~\ref{direct-inverse-loc-contra}).
 The contraherent cosheaf/complex $\gK$ is locally\+$\sC$ in
restriction to $V$, hence also in restriction to $U\cap V$;
and it follows that the contraherent cosheaf/complex $h'{}^!\gH$
is locally\+$\sC$ as the cokernel of the admissible monomorphism
$h'{}^!j^!\gK\rarrow h'{}^!\R$ between two locally\+$\sC$ objects
in $\sE^{U\cap V}_{\{U\cap V\}}$.
 We are using the assumptions of the coascent condition for
the classes $\sE$ and $\sC$ here.
 By Lemma~\ref{very-colocal-open-direct-image}, the contraherent
cosheaf/complex $j'_!h'{}^!\gH$ is locally\+$\sC$, too.
 Now in the short exact sequence $0\rarrow h^!\gK \rarrow h^!\Q
\rarrow h^!j_!\gH\rarrow0$ in $\sE^V_\bT$, the middle term is
locally\+$\sC$, because so are the two other terms.

 Finally, the composition of admissible monomorphisms $\gM\rarrow
\gK\rarrow\Q$ in $\sE^X_\bW$ is an admissible monomorphism whose
cokernel is an extension of the objects $j_!\gH$ and $\gG$, hence also
a finitely iterated extension of the objects $j_\alpha{}_!\gF_\alpha$
for some $\gF_\alpha\in\sF(U_\alpha)$ with $\alpha\le\beta$
(take $\gF_\beta=\gH$).
 This finishes the induction step and the proof of the proposition. 
\end{proof}

\begin{prop} \label{loc-contraherent-gluing-precover}
 Assume that we are given, for every ring $R\in\sR$, a pair of classes
of objects $(\sF(R),\sC^R)$ in the exact category\/ $\sE^R$ such that
an (admissible) short exact
sequence~\eqref{special-preenvelope-sequence} with $C\in\sC^R$ and
$F'\in\sF(R)$ exists in\/ $\sE^R$ for every object $E\in\sE^R$.
 Assume further that the class of modules/complexes\/ $\sC$ is very
colocal, the class\/ $\sC^R$ is closed under extensions and cokernels
of admissible monomorphisms in\/ $\sE^R$ for every $R\in\sR$, and
the class\/ $\sF(R)$ is generating in\/ $\sE^R$ for every $R\in\sR$.
 Let $X$ be a locally\+$\sR$ quasi-compact semi-separated scheme with
an open covering\/ $\bW$ and $X=\bigcup_\alpha U_\alpha$ be a finite
affine open covering of $X$ subordinate to\/~$\bW$.
 Then any object\/ $\gM\in\sE^X_\bW$ can be included in a short exact
sequence\/ $0\rarrow\P\rarrow\gF\rarrow\gM\rarrow0$ in\/ $\sE^X_\bW$
with\/ $\P\in\sC^X$ and\/ $\gF\in\sF'_{\{U_\alpha\}}(X)$.
\end{prop}

\begin{proof}
 Let us show that the class of all contraherent cosheaves/complexes of
the form $\bigoplus_\alpha j_\alpha{}_*\gF_\alpha$, where
$\gF_\alpha\in\sF(U_\alpha)$, is generating in~$\sE^X_\bW$.
 In view of Proposition~\ref{loc-contraherent-gluing-preenvelope}
and Lemma~\ref{salce-lemma}(a), this will imply the assertion of
the proposition.

 Indeed, given an object $\gM\in\sE^X_\bW$, pick admissible epimorphisms
$\gF_\alpha\rarrow j_\alpha^!\gM$ in $\sE^{U_\alpha}$ with
$\gF_\alpha\in\sF(U_\alpha)$ for all indices~$\alpha$.
 Let us check that the induced morphism $\bigoplus_\alpha j_\alpha{}_!
\gF_\alpha\rarrow\gM$ is an admissible epimorphism in~$\sE^X_\bW$.

 Both the morphisms $\bigoplus_\alpha j_\alpha{}_!\gF_\alpha\rarrow
\bigoplus_\alpha j_\alpha{}_!j_\alpha^!\gM$ and $\bigoplus_\alpha
j_\alpha{}_!j_\alpha^!\gM\rarrow\gM$ are admissible epimorphisms
in~$\sK^X_\bW$ (see~\eqref{contraherent-cech}); hence so is
their composition.
 Denote by $\gN$ the kernel of the admissible epimorphism
$\bigoplus_\alpha j_\alpha{}_!\gF_\alpha\rarrow\gM$ in~$\sK^X_\bW$.
 It remains to show that $\gN\in\sE^X_\bW$.

 By the definition of $\sE^X_\bW$, we only need to check that
$j_\beta^!\gN\in\sE^{U_\beta}$ for every index~$\beta$.
 Now we have a short exact sequence $0\rarrow j_\beta^!\gN\rarrow
\bigoplus_\alpha j_\beta^!j_\alpha{}_!\gF_\alpha\rarrow
j_\beta^!\gM\rarrow0$ in $\sK^{U_\beta}$ with the middle and rightmost
terms in~$\sE^{U_\beta}$.
 The direct summand inclusion map $\gF_\beta=j_\beta^!j_\beta{}_!
\gF_\beta\rarrow\bigoplus_\alpha j_\beta^!j_\alpha{}_!\gF_\alpha$ is
a (split) admissible monomorphism in $\sE^{U_\beta}$,
while the composition $\gF_\beta\rarrow\bigoplus_\alpha
j_\beta^!j_\alpha{}_!\gF_\alpha\rarrow j_\beta^!\gM$ is an admissible
epimorphism in $\sE^{U_\beta}$ by construction.
 Applying K\"unzer's axiom (the dual version
of~\cite[Exercise~3.11(i)]{Bueh}), we conclude that
$\bigoplus_\alpha j_\beta^!j_\alpha{}_!\gF_\alpha\rarrow j_\beta^!\gM$
is an admissible epimorphism in~$\sE^{U_\beta}$.
 Hence the kernel of the latter morphism in the category
$\sK^{U_\beta}$ belongs to~$\sE^{U_\beta}$.
\end{proof}

\begin{proof}[Proof of
Theorem~\ref{loc-contraherent-gluing-theorem}]
 In order to prove the first two assertions of the theorem, we apply
Lemma~\ref{cotorsion-pair-direct-summands-lemma} to the exact
category $\sE^X_\bW$ with two classes of objects
$\sF=\sF'_{\{U_\alpha\}}(X)$ and $\sC=\sC^X_\bW$.
 Lemma~\ref{loc-contraherent-gluing-orthogonality} establishes
the $\Ext^1$\+orthogonality assumption of
Lemma~\ref{cotorsion-pair-direct-summands-lemma} in the situation
at hand, while
Propositions~\ref{loc-contraherent-gluing-preenvelope}\+-%
\ref{loc-contraherent-gluing-precover} provide the short exact
sequences~\eqref{special-preenvelope-sequence}
and~\eqref{special-precover-sequence}.
 Lemma~\ref{cotorsion-pair-direct-summands-lemma} tells us that
the pair of classes $\sF(X)=\sF'_{\{U_\alpha\}}(X)^\oplus$
and $(\sC^X_\bW)^\oplus$ is a complete cotorsion pair in~$\sE^X_\bW$.
 It remains to observe that the class $\sC^X_\bW$ is closed under
direct summands in $\sE^X_\bW$, because the class $\sC^R$ is
closed under direct summands in $\sE^R$ for every $R\in\sR$.
 The cotorsion pair $(\sF(X),\sC^X_\bW)$ in $\sE^X_\bW$ is hereditary,
because the class $\sC^X_\bW$ is closed under cokernels of admissible
monomorpisms in~$\sE^X_\bW$ (because the class $\sC^R$ is closed
under cokernels of admissible monomorphisms in~$\sE^R$).

 The proof of the description of the class $\sF(X)\cap\sC^X_\bW$ stated
in the theorem is similar; cf.\ Lemma~\ref{clp-lin}.
 The ``if'' assertion is provided by
Lemma~\ref{very-colocal-open-direct-image} for the class $\sC$ and
by the description of the class $\sF(X)$ that we have obtained already.
 Let us prove the ``only if''.

 Let $\P\in\sC^X_\bW$ be a locally\+$\sC$ \,$\bW$\+locally contraherent
cosheaf/complex on~$X$.
 For every index~$\alpha$, pick a special precover
sequence~\eqref{special-precover-sequence} \,$0\rarrow\Q_\alpha
\rarrow\gG_\alpha\rarrow j_\alpha^!\P\rarrow0$ in the complete
cotorsion pair $(\sF(U_\alpha),\sC^{U_\alpha})$ in the exact
category~$\sE^{U_\alpha}$.
 So $\gG_\alpha\in\sF(U_\alpha)$ and $\Q_\alpha\in\sC^{U_\alpha}$.
 Since $j_\alpha^!\P\in\sC^{U_\alpha}$, it follows that
$\gG_\alpha\in\sF(U_\alpha)\cap\sC^{U_\alpha}$.

 Now consider the induced morphism $\bigoplus_\alpha j_\alpha{}_!
\gG_\alpha\rarrow\P$ in $\sC^X_\bW\subset\sE^X_\bW$.
 Similarly to the proof of
Proposition~\ref{loc-contraherent-gluing-precover}, one shows that
this morphism is an admissible epimorphism in $\sE^X_\bW$, and in fact,
by the same argument, even in (the inherited exact category structure
on)~$\sC^X_\bW$.
 We recall that $\sC$ is a colocal class, just like~$\sE$.

 Thus we have a short exact sequence $0\rarrow\gK\rarrow
\bigoplus_\alpha j_\alpha{}_!\gG_\alpha\rarrow\P\rarrow0$
in~$\sC^X_\bW$.
 If $\P\in\sF(X)\cap\sC^X_\bW$, then $\Ext^1_{\sE^X_\bW}(\P,\gK)=0$;
so this short exact sequence splits and $\P$ is a direct summand of
$\bigoplus_\alpha j_\alpha{}_!\gG_\alpha$.

 The assertions in the last paragraph of the theorem follow from
the descriptions of the two categories in question stated in
the preceding paragraphs (cf.\ Corollary~\ref{clp-independence}).
 The point is that the class of objects $\sF'_{\{U_\alpha\}}(X)^\oplus$
depends on the affine open covering $X=\bigcup_\alpha U_\alpha$, but
\emph{not} on the open covering~$\bW$ (just as the notation suggests).
 Notice that the full subcategory $X\lcth_\bW$ is closed under
extensions in $X\lcth$ by Corollary~\ref{nonaffine-lctrh-extensions},
and it follows that the full subcategory $\sE^X_\bW$ is closed under
extensions in $\sE^X_\bV$ for any open covering $\bV$ of $X$
subordinate to~$\bW$.
 Furthermore, the full subcategory $\sE^X_\bW$ is closed under direct
summands in $\sE^X_\bV$, since $X\lcth_\bW$ is closed under direct
summands in~$X\lcth$.
 Now, given two open coverings $\bW'$ and $\bW''$ of the scheme $X$,
it remains to choose an affine open covering $X=\bigcup_\alpha U_\alpha$
subordinate to both $\bW'$ and~$\bW''$.
\end{proof}

\begin{rem} \label{loc-contraherent-gluing-nonuniversal-remark}
 Theorem~\ref{loc-contraherent-gluing-theorem} is formulated above
for ``universal'' cotorsion pairs, defined uniformly for all commutative
rings $R\in\sR$.
 Similarly to the discussion in
Remark~\ref{quasi-coherent-gluing-nonuniversal-remark},
the assertion of Theorem~\ref{loc-contraherent-gluing-theorem} can be
extended to nonuniversal cotorsion pairs, depending on an additional
datum specific to a particular scheme~$X$.
 One needs to fix a quasi-compact semi-separated scheme $X$ with
an open covering $\bW$, and assume given an exact subcategory of
contraadjusted modules or complexes $\sE^U\sub\sK^{\O(U)}$ for every
affine open subscheme $U\sub X$ subordinate to $\bW$ (where
$\sK^R=R\modl^\cta$ or $\Com(R\modl^\cta)$, while $\sE^U$ depends
on the embedding of $U$ into~$X$) together with a hereditary complete
cotorsion pair $(\sF(U),\sC^U)$ in~$\sE^U$.
 The classes $\sE^U$ and $\sC^U$ must be colocal with respect to
affine open coverings of affine open subschemes $U\sub X$ subordinate
to $\bW$, and preserved by direct images with respect to identity
embeddings of affine open subschemes $V\rarrow U$, where $V\sub U
\sub X$.
 Then the same arguments and constructions as above produce an exact
subcategory $\sE^X_\bW\sub\sK^X_\bW$, where $\sK^X_\bW=X\lcth_\bW$
or $\Com(X\lcth_\bW)$, and a hereditary complete cotorsion pair
$(\sF(X),\sC^X_\bW)$ in $\sE^X_\bW$, with the same description of
the classes $\sE^X_\bW$, \,$\sC^X_\bW$, and~$\sF(X)$.
\end{rem}

\subsection{Becker's coderived and contraderived categories}
\label{becker-subsect}
 In this section we discuss an alternative approach to coderived
and contraderived categories, as compared to the definitions in
Section~\ref{derived-second-kind}.
 This approach is named after Becker~\cite[Remark~9.2]{PS4},
\cite[Section~7]{Pksurv}, \cite[Section~4]{Pctrl},
\cite[Examples~2.5(3) and~2.6(3)]{Pps}, \cite[Appendix~A]{Psemten},
with the reference to his important paper~\cite{Bec}, even though
it has its roots in the earlier work of J\o rgensen~\cite{Jorg},
Krause~\cite{K-st}, and Neeman~\cite{N-f}.

 This section is included in this Appendix~\ref{cotorsion-pairs-appx}, 
rather than in Appendix~\ref{exact-derived}, because the cotorsion
pairs are the main technique used for working with derived categories
of second kind in the sense of Becker.
 Indeed, the possibility to use powerful techniques of
set-theoretical homological algebra, expressed in the small object
argument (see Section~\ref{small-object-argument-subsect}), is the main
advantage of Becker's approach as compared to the definitions
in Section~\ref{derived-second-kind}.
 Among the disadvantages of the Becker coderived and contraderived
categories, there are the facts that it may be difficult to prove
that all Becker-coacyclic or Becker-contraacyclic complexes are acyclic
(if one's exact category is not abelian)~\cite[Lemma~A.2]{Psemten} or
to show that an exact, (co)product-preserving functor preserves
co/contraacyclicity in the sense of Becker (if the functor does not
have an adjoint on the suitable side)~\cite[Lemma~A.5]{Psemten}.
 See the discussion below.

 Let $\sE$ be an exact category.
 A complex $A^\bu$ in $\sE$ is said to be \emph{Becker-coacyclic} if
the complex $\Hom_\sE(A^\bu,J^\bu)$ is acyclic for any complex $J^\bu$
of injective objects in~$\sE$.
 The \emph{Becker coderived category} $\sD^\bco(\sE)$ is defined as
the triangulated Verdier quotient category
$$
 \sD^\bco(\sE)=\Hot(\sE)/\Acycl^\bco(\sE)
$$
of the unbounded homotopy category $\Hot(\sE)$ by the thick subcategory
$\Acycl^\bco(\sE)$ of Becker-coacyclic complexes.

 Dually, a complex $B^\bu$ in $\sE$ is said to be
\emph{Becker-contraacyclic} if the complex $\Hom_\sE(P^\bu,B^\bu)$ is
acyclic for any complex $P^\bu$ of projective objects in~$\sE$.
 The \emph{Becker contraderived category} $\sD^\bctr(\sE)$ is defined
as the triangulated Verdier quotient category
$$
 \sD^\bctr(\sE)=\Hot(\sE)/\Acycl^\bctr(\sE)
$$
of the homotopy category $\Hot(\sE)$ by the thick subcategory
$\Acycl^\bctr(\sE)$ of Becker-contraacyclic complexes.

 It is helpful to keep in mind that, in any exact category with
infinite direct sums and enough injective objects, the infinite direct
sum functors are exact.
 Dually, in any exact category with infinite products and enough
projective objects, the infinite product functors are
exact~\cite[Remark~5.2]{Pedg}.

 To avoid ambiguity, we will call the coacyclic and contraacyclic
complexes in the sense of Section~\ref{derived-second-kind}
\emph{Positselski-coacyclic} and \emph{Positselski-contraacyclic}.

\begin{lem} \label{Positselski-trivial-are-Becker-trivial}
\textup{(a)} In any exact category, all absolutely acyclic complexes
(as defined in Section~\ref{derived-second-kind}) are Becker-coacyclic
and Becker-contraacyclic. \par
\textup{(b)} In an exact category with exact functors of infinite direct
sum, all Positselski-coacyclic complexes are Becker-coacyclic. \par
\textup{(c)} In an exact category with exact functors of infinite
product, all Positselski-contraacyclic complexes are
Becker-contraacyclic.
\end{lem}

\begin{proof}
 This is essentially explained in the proof of
Lemma~\ref{homotopy-inj-proj-fully-faithful}.
\end{proof}

 The following two lemmas form a generalization
of~\cite[Lemma~A.2]{Psemten}.

\begin{lem} \label{bounded-Becker-trivial-are-acyclic}
\textup{(a)} In an exact category with enough injective objects,
a bounded below complex is Becker-coacyclic if and only if
it is acyclic. \par
\textup{(b)} In an exact category with enough projective objects,
a bounded above complex is Becker-contraacyclic if and only if
it is acyclic.
\end{lem}

\begin{proof}
 Part~(b): the ``only if'' implication is~\cite[Exercise~11.11]{Bueh}.
 The claim is that a bounded above complex $B^\bu$ in $\sB$ is acyclic
whenever, for every projective object $P\in\sB$, the complex
$\Hom_\sB(P,B^\bu)$ is acyclic.
 The ``if'' implication is well-known (and can be deduced from
Lemma~\ref{telescope} applied to the projective system of bicomplexes
of morphisms from silly truncations of a given complex of projective
objects $P^\bu$ to the bounded above acyclic complex~$B^\bu$).
\end{proof}

\begin{lem} \label{with-co-kernels-Becker-trivial-are-acyclic}
\textup{(a)} Let\/ $\sA$ be an exact category with enough injective
objects.
 Assume that all morphisms have cokernels in\/~$\sA$.
 Then any Becker-coacyclic complex in\/ $\sA$ is acyclic. \par
\textup{(b)} Let\/ $\sB$ be an exact category with enough projective
objects.
 Assume that all morphisms have kernels in\/~$\sB$.
 Then any Becker-contraacyclic complex in\/ $\sB$ is acyclic.
\end{lem}

\begin{proof}
 Part~(b): let $B^\bu$ be a complex in\/~$\sB$.
 Assume that the complex $\Hom_\sB(P,B^\bu)$ is acyclic for every
projective object $P\in\sB$.
 Then the subcomplexes of canonical truncation of the complex $B^\bu$
also have this property (notice that the functor $\Hom_\sB(P,{-})$
preserves kernels).
 It remains to apply (the proof of) the ``only if'' implication of
Lemma~\ref{bounded-Becker-trivial-are-acyclic}(b).
\end{proof}

\begin{rem} \label{Becker-co-contra-acyclic-are-acyclic-remark}
 We do \emph{not} know whether a Becker-coacyclic complex in an exact
category with enough injectives needs to be acyclic in general, or
dually, whether a Becker-contraacyclic complex in an exact category
with enough projectives needs to be acyclic.
 However, notice that the argument in the proof of
Lemma~\ref{with-co-kernels-Becker-trivial-are-acyclic} did not use
the full definition of co/contraacyclicity; rather, the claim in
part~(b) was that a complex $B^\bu$ is acyclic whenever the complex
of abelian groups $\Hom_\sB(P,B^\bu)$ is acyclic for every projective
object $P\in\sB$.
 So we have only used one-term complexes of projectives, rather than
arbitrary complexes of projectives, in this argument.

 This is certainly \emph{not} sufficient in general.
 To give a simple example, let $k$~be a field, $\Lambda=
k[\epsilon]/(\epsilon^2)$ be the $k$\+algebra of dual numbers,
and $\sE$ be the additive category of free (equivalently, projective
or injective) $\Lambda$\+modules.
 Endow $\sE$ with the trivial (split) exact category structure; so
all the objects of $\sE$ are both projective and injective.
 Let $\dotsb\rarrow\Lambda\overset\epsilon\rarrow\Lambda
\overset\epsilon\rarrow\Lambda\rarrow\dotsb$ be the unbounded,
noncontractible, acyclic complex of free $\Lambda$\+modules with
one generator (cf.~\cite[formula~(2) in Section~7.4]{Pksurv}).
 Denote this complex by~$C^\bu$.
 Viewed as a complex in $\sE$, the complex $C^\bu$ is \emph{not}
acyclic.
 It is not Becker-coacyclic or Becker-contraacyclic, either; indeed,
it is a noncontractible complex of projective-injective objects.
 Still, the complexes of abelian groups $\Hom_\sE(C^\bu,J)$
and $\Hom_\sE(P,C^\bu)$ are acyclic for all objects $J$, $P\in\sE$.

 To give another example, consider the exact category $\sB=R\modl_\fl$
of flat modules over an associative ring~$R$.
 Then it is indeed true that every Becker-contraacyclic complex in
$\sB$ is acyclic (and vice versa; the classes of Becker-contraacyclic
and acyclic complexes in $\sB$ coincide).
 But this is a nontrivial theorem of Neeman (see
Proposition~\ref{flat-projective-periodicity-complements}(a)) which
seems to be a special property of the class of flat modules.

 Any Becker-coacyclic complex in $R\modl_\fl$ is also acyclic (and
vice versa; the classes of Becker-coacyclic and acyclic complexes
in $R\modl_\fl$ also coincide).
 See~\cite[Theorems~7.14 and~7.18]{Pphil}.
 Generalizations to flat contraadjusted quasi-coherent sheaves and
flat/antilocally flat contraherent cosheaves can be found in
Corollaries~\ref{becker-co-contraderived-of-fl-cta},
\ref{acycl=bctracycl=bcoacycl-in-alf},
and~\ref{finite-krull-derived-equivalences}(a).
 For generalizations to flat comodules and contramodules,
see~\cite[Theorems~6.5 and~13.2]{Pflcc} and~\cite[Theorems~5.1 and~6.1,
and Corollary~9.1]{Pbc}.
\end{rem}

 The next lemma is a generalization of~\cite[Lemma~A.5]{Psemten}.

\begin{lem} \label{exact-with-adjoint-preservation-lemma}
\textup{(a)} Let\/ $\sA'$, $\sA''$ be exact categories and
$F\:\sA'\rarrow\sA''$ be an exact functor having a right adjoint.
 Then the functor $F$ takes Becker-coacyclic complexes to
Becker-coacyclic complexes. \par
\textup{(b)} Let\/ $\sB'$, $\sB''$ be exact categories and
$G\:\sB'\rarrow\sB''$ be an exact functor having a left adjoint.
 Then the functor $G$ takes Becker-contraacyclic complexes to
Becker-contraacyclic complexes.
\end{lem}

\begin{proof}
 Part~(b): let $L$ be the left adjoint functor to~$G$.
 Then $L$ takes projective objects in $\sB''$ to projective objects
in $\sB'$ (since $G$ is exact).
 Now if $B^\bu$ is a Becker-contraacyclic complex in $\sB'$ and
$P^\bu$ is a complex of projective objects in~$\sB''$, then
$\Hom_{\sB''}(P^\bu,G(B^\bu))\simeq\Hom_{\sB'}(L(P^\bu),B^\bu)$ is
an acyclic complex of abelian groups.
\end{proof}

 The following three theorems describe three contexts in which
the Becker coderived and Becker contraderived categories are
well-behaved.
 In the former two of the three cases, they agree with
the suitable Positselski derived categories of the second kind
(which are consequently well-behaved, too).

\begin{thm} \label{finite-homol-dim-becker-co-contra-derived}
\textup{(a)} Let\/ $\sA$ be an exact category of finite homological
dimension with enough injective objects.
 Then the classes of absolutely acyclic, Becker-coacyclic, and
acyclic complexes in\/ $\sA$ coincide.
 The composition of triangulated functors
$$
 \Hot(\sA^\inj)\lrarrow\Hot(\sA)\lrarrow\sD^\abs(\sA)=
 \sD^\bco(\sA)=\sD(A)
$$
is a triangulated equivalence\/
$\Hot(\sA^\inj)\simeq\sD^\abs(\sA)=\sD^\bco(\sA)=\sD(A)$. \par
\textup{(b)} Let\/ $\sB$ be an exact category of finite homological
dimension with enough projective objects.
 Then the classes of absolutely acyclic, Becker-contraacyclic, and
acyclic complexes in\/ $\sB$ coincide.
 The composition of triangulated functors
$$
 \Hot(\sB_\prj)\lrarrow\Hot(\sB)\lrarrow\sD^\abs(\sB)=
 \sD^\bctr(\sB)=\sD(\sB)
$$
is a triangulated equivalence\/
$\Hot(\sB_\prj)\simeq\sD^\abs(\sB)=\sD^\bctr(\sB)=\sD(\sB)$.
\end{thm}

\begin{proof}
 This is~\cite[Proposition~7.5]{Pphil}.
 The assertion that all acyclic complexes are absolutely acyclic in
an exact category of finite homological dimension is
Lemma~\ref{psemi-remark21}.
 The argument proving that the functors $\Hot(\sA^\inj)\rarrow
\sD^\abs(\sA)$ and $\Hot(\sB_\prj)\rarrow\sD^\abs(\sB)$ are triangulated
equivalences in the respective assumptions goes back to~\cite[proof of
Theorem~3.6]{Pkoszul}; for a far-reaching generalization,
see~\cite[proof of Theorem~5.6]{Pedg}.
 For another generalization of part~(b), see
Proposition~\ref{finite-resolutions} above (take $\sE=\sB$ and
$\sF=\sB_\prj$).
 In view of Lemma~\ref{Positselski-trivial-are-Becker-trivial}(a),
the assertion that the classes of absolutely acyclic and
Becker co/contraacyclic complexes coincide in the respective cases
follows.
 For a further discussion, see~\cite[Theorem~7.8]{Pksurv}.
\end{proof}

\begin{thm} \label{positselski-becker-co-contra-derived}
\textup{(a)} Let\/ $\sA$ be an exact category with exact functors of
infinite direct sum and enough injective objects.
 Assume that countable direct sums of injective objects have finite
injective dimensions in\/~$\sA$.
 Then the classes of Positselski-coacyclic and Becker-coacyclic
complexes in\/ $\sA$ coincide.
 The composition of triangulated functors
$$
 \Hot(\sA^\inj)\lrarrow\Hot(\sA)\lrarrow\sD^\co(\sA)=\sD^\bco(\sA)
$$
is a triangulated equivalence\/
$\Hot(\sA^\inj)\simeq\sD^\co(\sA)=\sD^\bco(\sA)$. \par
\textup{(b)} Let\/ $\sB$ be an exact category with exact functors of
infinite product and enough projective objects.
 Assume that countable products of projective objects have finite
projective dimensions in\/~$\sB$.
 Then the classes of Positselski-contraacyclic and Becker-contraacyclic
complexes in\/ $\sB$ coincide.
 The composition of triangulated functors
$$
 \Hot(\sB_\prj)\lrarrow\Hot(\sB)\lrarrow\sD^\ctr(\sB)=\sD^\bctr(\sB)
$$
is a triangulated equivalence\/
$\Hot(\sB_\prj)\simeq\sD^\ctr(\sB)=\sD^\bctr(\sB)$.
\end{thm}

\begin{proof}
 This is~\cite[Proposition~7.6]{Pphil}.
 The argument proving the assertion that the functors
$\Hot(\sA^\inj)\rarrow\sD^\co(\sA)$ and $\Hot(\sB_\prj)\rarrow
\sD^\ctr(\sB)$ are triangulated equivalences in the respective
assumptions goes back to~\cite[Sections~3.7\+-3.8]{Pkoszul};
for a far-reaching generalization, see~\cite[Theorem~5.10]{Pedg}.
 Another generalization of part~(b) is
Corollary~\ref{finite-homol-dim-equivalence-cor}
(take $\sE=\sB$ and $\sF=\sB_\prj$).
 In view of Lemma~\ref{Positselski-trivial-are-Becker-trivial}(b\+c),
the assertion that the classes of Positselski co/contraacyclic
complexes coincide with the respective Becker classes follows.
 For a further discussion, see~\cite[Theorem~7.9]{Pksurv}.
\end{proof}

 For a different kind of argument showing that, for certain special
classes of abelian categories, Becker's co/contraderived categories
coincide with Positselski's ones, see
Theorems~\ref{discrete-mod-coderived}(d)
and~\ref{contramod-contraderived}(d), and
Remark~\ref{no-needed-adjoints-remark} below.

 In the following theorem, unlike in all the previous results in
this section, the ``coderived'' and ``contraderived'' assertions
\emph{cannot} be obtained from one another by inverting the arrows.
 They are two different results.

\begin{thm} \label{coderived-of-grothendieck-contraderived-of-lpacepo}
\textup{(a)} Let\/ $\sA$ be a Grothendieck abelian category.
 Then the composition of triangulated functors
$$
 \Hot(\sA^\inj)\lrarrow\Hot(\sA)\lrarrow\sD^\bco(\sA)
$$
is a triangulated equivalence\/ $\Hot(\sA^\inj)\simeq\sD^\bco(\sA)$.
\par
\textup{(b)} Let\/ $\sB$ be a locally presentable abelian category
with enough projective objects.
 Then the composition of triangulated functors
$$
 \Hot(\sB_\prj)\lrarrow\Hot(\sB)\lrarrow\sD^\bctr(\sB)
$$
is a triangulated equivalence\/ $\Hot(\sB_\prj)\simeq\sD^\bctr(\sB)$.
\end{thm}

\begin{proof}
 Part~(a) is~\cite[Corollary~9.5]{PS4}.
 Part~(b) is~\cite[Corollary~7.4]{PS4}.
 The proofs are based on the technique of complete cotorsion pairs
and the small object argument (see
Section~\ref{small-object-argument-subsect}).
 Specifically, in part~(a) one proves that the pair of classes
 of coacyclic complexes and complexes of injective objects
($\Acycl^\bco(\sA)$, $\Com(\sA^\inj)$) is a hereditary complete
cotorsion pair in $\Com(\sA)$.
 In part~(b), one shows that the pair of classes of complexes of
projective objects and contraacyclic complexes
($\Com(\sB_\prj)$, $\Acycl^\bctr(\sB)$) is a hereditary complete
cotorsion pair in $\Com(\sB)$.

 In connection with part~(a), see also~\cite[Theorem~4.2]{Gil4}.
 For another approach to part~(a), see~\cite[Theorem~2.13]{N-i}
or~\cite[Corollary~5.13]{K-au}.
\end{proof}

 The following two results form a partial version of
Proposition~\ref{finite-resolutions} for Becker's contraderived
categories.

\begin{prop} \label{becker-contraderived-finite-resolutions}
 Let\/ $\sE$ be a weakly idempotent-complete exact category and\/
$\sF\sub\sE$ be a resolving full subcategory closed under direct
summands.
 Assume that the\/ $\sF$\+resolution dimensions of all objects of\/
$\sE$ do not exceed a fixed constant~$n$.
 Then the triangulated functor\/ $\sD^\bctr(\sF)\rarrow\sD^\bctr(\sE)$
induced by the embedding functor\/ $\sF\rarrow\sE$ is an equivalence
of triangulated categories.
\end{prop}

\begin{proof}
 Under the assumptions of the proposition, all the projective objects
of $\sE$ belong to~$\sF$.
 Furthermore, any object $P\in\sF$ that is projective in $\sF$ is also
projective in~$\sE$.
 Indeed, it suffices to show that any admissible epimorphism
$E\rarrow P$ in $\sE$ is split.
 Pick an admissible epimorphism $F\rarrow E$ in $\sE$ with $F\in\sF$;
then the composition $F\rarrow P$ is a split epimorphism, and
it follows that so is the morphism $E\rarrow P$.

 We have shown that a complex in $\sF$ is Becker-contraacyclic in $\sF$
if and only if it is Becker-contraacyclic in~$\sE$.
 In view of Lemma~\ref{pkoszul-lemma16}(a), it remains to check that
for every complex $E^\bu$ in $\sE$ there exists a complex $F^\bu$
in $\sF$ together with a morphism of complexes $F^\bu\rarrow E^\bu$
with Becker-contraacyclic cone.
 It was explained in the beginning of the proof of
Proposition~\ref{finite-resolutions} how to construct such
a complex $F^\bu$ and a morphism $F^\bu\rarrow E^\bu$ so that
the cone would be even absolutely acyclic.
 Any absolutely acyclic complex is Becker-contraacyclic by
Lemma~\ref{Positselski-trivial-are-Becker-trivial}(a).
\end{proof}

\begin{prop} \label{becker-contraderived-finite-coresolutions}
 Let\/ $\sE$ be a weakly idempotent-complete exact category with enough
projective objects and\/ $(\sF,\sC)$ be a hereditary complete cotorsion
pair in\/~$\sE$.
 Assume that the\/ $\sC$\+coresolution dimensions of all objects of\/
$\sE$ do not exceed a fixed constant~$n$.
 Then a complex in\/ $\sC$ is Becker-contraacyclic in\/ $\sC$ if and
only if it is Becker-contraacyclic as a complex in\/~$\sE$,
$$
 \Acycl^\bctr(\sC)=\Com(\sC)\cap\Acycl^\bctr(\sE).
$$
 The triangulated functor\/ $\sD^\bctr(\sC)\rarrow\sD^\bctr(\sE)$
induced by the embedding functor\/ $\sC\rarrow\sE$ is an equivalence
of triangulated categories.
\end{prop}

\begin{proof}
 By Lemma~\ref{finite-homol-dim-finite-coresol-dim}, the homological
dimension of the exact category $\sF$ is finite.
 Furthermore, there are enough projective and injective objects
in~$\sF$: in fact, one has $\sF_\prj=\sE_\prj$ and
$\sF^\inj=\sF\cap\sC=\sC_\prj$.
 So both parts of
Theorem~\ref{finite-homol-dim-becker-co-contra-derived} are applicable
to~$\sF$, telling us that for any complex $F^\bu$ in $\sF$ there exist
two complexes $P^\bu\in\Com(\sE_\prj)$ and $Q^\bu\in\Com(\sF\cap\nobreak
\sC)$ and two morphisms of complexes $P^\bu\rarrow F^\bu\rarrow Q^\bu$
such that the cones of both the morphisms are absolutely acyclic
in~$\sF$.

 In particular, for any complex $P^\bu$ in $\sE_\prj$ there exists
a complex $Q^\bu$ in $\sF\cap\sC$ together with a morphism of complexes
$P^\bu\rarrow Q^\bu$ with the cone absolutely acyclic in~$\sF$.
 Conversely, for any complex $Q^\bu$ in $\sF\cap\sC$ there exists
a complex $P^\bu$ in $\sE_\prj$ together with a morphism of complexes
$P^\bu\rarrow Q^\bu$ with the cone absolutely acyclic in~$\sF$.

 Furthermore, for any short exact sequence $0\rarrow K^\bu\rarrow
L^\bu\rarrow M^\bu\rarrow 0$ of complexes in $\sF$ and any complex
$B^\bu$ in $\sC$, the short sequence of complexes of abelian groups
$0\rarrow\Hom_\sE(M^\bu,B^\bu)\rarrow\Hom_\sE(L^\bu,B^\bu)\rarrow
\Hom_\sE(K^\bu,B^\bu)\rarrow0$ is exact.
 Consequently, the complex $\Hom_\sE(A^\bu,B^\bu)$ is acyclic for any
absolutely acyclic complex $A^\bu$ in~$\sF$.
 Thus, in the context of the previous paragraph, the complex
$\Hom_\sE(P^\bu,B^\bu)$ is acyclic if and only if the complex
$\Hom_\sE(Q^\bu,B^\bu)$ is acyclic.
 This proves the first assertion of the proposition.

 Finally, for any complex $E^\bu$ in $\sE$ there exists a complex
$C^\bu$ in $\sC$ together with a morphism of complexes $E^\bu\rarrow
C^\bu$ whose cone is absolutely acyclic in~$\sE$.
 A dual version of this construction was spelled out the beginning of
the proof of Proposition~\ref{finite-resolutions}.
 By Lemma~\ref{Positselski-trivial-are-Becker-trivial}(a), any
absolutely acyclic complex is Becker-contraacyclic in~$\sE$.
 By virtue of Lemma~\ref{pkoszul-lemma16}(b), the second assertion
of the proposition follows.
\end{proof}

 The next proposition is a version of
Proposition~\ref{infinite-resolutions}(b) for Becker's contraderived
categories.

\begin{prop} \label{becker-contraderived-infinite-resolutions}
 Let\/ $\sE$ be an exact category with countable products and\/
$\sF\sub\sE$ be a resolving full subcategory closed under countable
products.
 Then the triangulated functor\/ $\sD^\bctr(\sF)\rarrow\sD^\bctr(\sE)$
induced by the embedding functor\/ $\sF\rarrow\sE$ is an equivalence
of triangulated categories.
\end{prop}

\begin{proof}
 Let us first show that $\sF$ is closed under direct summands in~$\sE$.
 Indeed, if $E\oplus E'=F$ in $\sE$ and $F\in\sF$, then $E$ is
the kernel of a suitable split epimorphism $\prod_{n=1}^\infty F
\rarrow\prod_{n=1}^\infty F$ in~$\sE$ (``the cancellation trick'').
 As $\sF$ is closed under countable products and kernels of admissible
epimorphisms in $\sE$, it follows that $E\in\sF$.
 Now the first paragraph of the proof of
Proposition~\ref{becker-contraderived-finite-resolutions} is applicable,
and we can conclude that the classes of projective objects in $\sE$
and $\sF$ coincide.
 Thus a complex in $\sF$ is Becker-contraacyclic in $\sF$ if and only if
it is Becker-contraacyclic in~$\sE$.

 It remains to check that for every complex $E^\bu$ in $\sE$ there
there exists a complex $F^\bu$ in $\sF$ together with a morphism of
complexes $F^\bu\rarrow E^\bu$ with Becker-contraacyclic cone.
 It was explained in the proof of
Proposition~\ref{infinite-resolutions}(b) how to construct such
a complex and a morphism so that the cone would be
Positselski-contraacyclic.
 Any Positselski-contraacyclic complex is Becker-contraacyclic by
Lemma~\ref{Positselski-trivial-are-Becker-trivial}(c).
 This argument presumes existence and exactness of all products
in $\sE$ and preservation of $\sF\sub\sE$ by all products, but one can
see following the steps that countable products are sufficient and
their exactness is not used.
\end{proof}

 Our final proposition is important in the context of contraherent
cosheaves.

\begin{prop} \label{becker-contraacyclicity-for-cta-cot}
\textup{(a)} Let $R$ be a commutative ring and $B^\bu$ be a complex of
contraadjusted $R$\+modules.
 Then $B^\bu$ is Becker-contraacyclic as a complex in $R\modl^\cta$ if
and only if $B^\bu$ is Becker-contraacyclic as a complex in $R\modl$.
 Moreover, the embedding of exact/abelian categories
$R\modl^\cta\rarrow R\modl$ induces a triangulated equivalence\/
$\sD^\bctr(R\modl^\cta)\simeq\sD^\bctr(R\modl)$. \par
\textup{(b)} Let $R$ be an associative ring and $B^\bu$ be a complex of
cotorsion left $R$\+modules.
 Then $B^\bu$ is Becker-contraacyclic as a complex in $R\modl^\cot$ if
and only if $B^\bu$ is Becker-contraacyclic as a complex in $R\modl$.
 Moreover, the embedding of exact/abelian categories
$R\modl^\cot\rarrow R\modl$ induces a triangulated equivalence\/
$\sD^\bctr(R\modl^\cot)\simeq\sD^\bctr(R\modl)$.
\end{prop}

\begin{proof}
 Part~(a): the first assertion is~\cite[Theorem~7.28]{Pphil}, and
the second one is~\cite[Corollary~7.30]{Pphil}.
 The whole part~(a) is a particular case of
Proposition~\ref{becker-contraderived-finite-coresolutions}
(for $n=1$), and the proof of part~(b) is somewhat similar.

 Part~(b): the first assertion is~\cite[Theorem~7.19]{Pphil}, and
the second one is~\cite[Corollary~7.21]{Pphil}.
 We will only present a sketch of the argument, and mostly
for the first assertion.
 To begin with, recall that the projective objects of the exact category
$R\modl^\cot$ are the flat cotorsion $R$\+modules (see
Section~\ref{contraadjusted-exact-cat}).

 The point is that for any complex of flat cotorsion $R$\+modules
$Q^\bu$ there exists a complex of projective $R$\+modules $P^\bu$
together with a morphism of complexes $P^\bu\rarrow Q^\bu$ whose cone
is an acyclic complex of flat $R$\+modules with flat $R$\+modules
of cocycles.
 This is a particular case of
Proposition~\ref{flat-projective-periodicity-complements}(b).

 Conversely, for any complex of projective $R$\+modules $P^\bu$
there exists a complex of flat cotorsion $R$\+modules $Q^\bu$
together with a morphism of complexes $P^\bu\rarrow Q^\bu$ whose cone
is an acyclic complex of flat $R$\+modules with flat $R$\+modules
of cocycles.
 More generally, this assertion holds for any complex of flat
$R$\+modules $G^\bu$ in place of~$P^\bu$.
 Indeed, by~\cite[Definition~3.3 and Corollary~4.10]{Gil}
(see Lemma~\ref{acyclic-of-vfl-flat-arbitrary-of-cta-cot-pairs}
for a quasi-coherent sheaf version), for any complex of
$R$\+modules $G^\bu$ there exists a special preenvelope short exact
sequence of complexes of $R$\+modules $0\rarrow G^\bu\rarrow C^\bu
\rarrow F^\bu\rarrow0$, where $C^\bu$ is a complex of cotorsion
$R$\+modules and $F^\bu$ is a complex of flat $R$\+modules with
flat $R$\+modules of cocycles.
 The argument finishes as spelled out in~\cite[proof of
Theorem~7.18]{Pphil}.

 It remains to refer to Proposition~\ref{dw-cotorsion-are-dg-cotorsion}
in order to finish the proof of the first assertion of part~(b).
 The proof of the second assertion can be based on
Lemma~\ref{pkoszul-lemma16}(b) and the same construction of the special
preenvelope sequence from the previous paragraph (for an arbitrary
complex of $R$\+modules~$G^\bu$).
 Here one has to refer to Theorem~\ref{flat-projective-periodicity}(b)
to the effect that all acyclic complexes of flat $R$\+modules with
flat $R$\+modules of cocycles are Becker-contraacyclic in $R\modl$.
\end{proof}

\Section{Co-Contra Correspondence over a Flat Coring}
\label{over-flat-coring-appx}

 The aim of this appendix is to extend the assertion of
Theorem~\ref{co-contra-dualizing} from semi-separated Noetherian
schemes to semi-separated Noetherian stacks.
 These are the stacks that can be represented by groupoids with
affine schemes of vertices and arrows (see~\cite[Section~2.1]{AB}).
 We avoid explicit use of the stack language by working with what
Kontsevich and Rosenberg call ``finite covers''~\cite[Section~2]{KR}
(see also~\cite{KR2} and~\cite[Sections~2\+-6]{Pflcc}) instead.
 This naturally includes the noncommutative geometry situation.
 The corresponding algebraic language is that of flat corings over
noncommutative rings; thus this appendix provides a bridge
between the results of Chapter~\ref{dualizing-complex-sect} and
those of~\cite[Chapter~5]{Psemi}.

\subsection{Contramodules over a flat coring}  \label{contra-over-flat}
 Let $A$ be an associative ring with unit.
 We refer to the memoir and monographs~\cite{EM,BW,Psemi} for
the definitions of a (coassociative) \emph{coring}
(with counit) $\C$ over $A$, a \emph{left comodule} $\M$
over $\C$, and a \emph{right comodule} $\N$ over~$\C$.
 The definition of a \emph{left contramodule} $\P$ over $\C$ can be
found in~\cite[Section~III.5]{EM}
or~\cite[Section~0.2.4 or~3.1.1]{Psemi}
(see also~\cite[Section~2.5]{Prev}).
 We denote the abelian groups of morphisms in the additive category
of left $\C$\+comodules by $\Hom_\C(\L,\M)$ and the similar groups
related to the additive category of left $\C$\+contramodules by
$\Hom^\C(\P,\Q)$.

 The category $\C\comodl$ of left $\C$\+comodules is abelian \emph{and}
the forgetful functor $\C\comodl\rarrow A\modl$ is exact if and only
if $\C$ is a flat right $A$\+module~\cite[Sections~18.6 and~18.14]{BW}, 
\cite[Section~1.1.2]{Psemi}, \cite[Proposition~2.12(a)]{Prev},
\cite[Lemma~2.1]{Pflcc}.
 The category $\C\contra$ of left $\C$\+contramodules is abelian
\emph{and} the forgetful functor $\C\contra\rarrow A\modl$ is exact if
and only if $\C$ is a projective left
$A$\+module~\cite[Section~3.1.2]{Psemi},
\cite[Proposition~2.12(b)]{Prev}.
 The following counterexample shows that the category $\C\contra$ may be
not abelian even though $\C$ is a flat left and right $A$\+module.

\begin{ex}
 Let us consider corings $\C$ of the following form.
 The coring $\C$ decomposes into a direct sum of $A$\+$A$\+bimodules
$\C=\C_{11}\oplus\C_{12}\oplus\C_{22}$; the counit map $\C\rarrow A$
annihilates $\C_{12}$ and the comultiplication map takes $\C_{ik}$
into the direct sum $\bigoplus_j\C_{ij}\ot_A\C_{jk}$ for all
$i$, $j$, $k=1$,~$2$.
 Assume further that the restrictions of the counit map to
$\C_{11}$ and $\C_{22}$ are both isomorphisms $\C_{ii}\simeq A$.
 Notice that the datum of an $A$\+$A$\+bimodule $\C_{12}$ determines
the coring $\C$ in this case.
 A left $\C$\+contramodule $\P$ is the same thing as a pair of left
$A$\+modules $\P_1$ and $\P_2$ endowed with an $A$\+module morphism
$\Hom_A(\C_{12},\P_1)\rarrow\P_2$.

 The kernel of a morphism of left $\C$\+contramodules $(f_1,f_2)\:
(\P_1,\P_2)\rarrow(\Q_1,\Q_2)$ (taken in the additive category
$\C\contra$) is the $\C$\+contramodule $(\ker f_1\;\ker f_2)$.
 The cokernel of the morphism~$(f_1,f_2)$ can be computed as
the $\C$\+contramodule $(\coker f_1\;\gL)$, where $\gL$ is
the cokernel of the morphism from $\P_2$ to the fibered coproduct
$\Hom_A(\C_{12}\;\coker f_1)\sqcup_{\Hom_A(\C_{12}\;\Q_1)}\Q_2$
(where $\ker f$ and $\coker f$ denote the kernel and cokernel
of a morphism of $A$\+modules~$f$).

 Now setting $A$ to be the ring of integers $\boZ$ and $\C_{12}$ to be
the (bi)module of rational numbers $\boQ$ over~$\boZ$, one can check
that the category $\C\contra$ is not abelian.
 The discussion in~\cite[Example~B.1.1]{Pcosh} contained a mistake
at this point: the specific counterexample suggested there
does not work.
 Let us expand the discussion a bit in order to explain how to
correct the counterexample from~\cite{Pcosh} so as to make it valid.

 Firstly, consider two left $\C$\+contramodules $\P$ and $\Q$ with
$\P_1=\bigoplus_{n\in\mathbb N}\boZ$, \ $\P_2=0$, and $\Q_1=\Q_2=\boQ$,
the structure morphism $\Hom_\boZ(\boQ,\Q_1)\rarrow\Q_2$ being
the identity isomorphism.
 Choose the map $f_1\:\P_1\rarrow\Q_1$ to be surjective,
while of course the map $f_2\:\P_2\rarrow\Q_2$ is zero.
 Then the kernel of the $\C$\+contramodule morphism $f=(f_1,f_2)\:
\P\rarrow\Q$ is the $\C$\+contramodule $(\ker f_1\;0)$.
 Surprisingly (and counterintuitively), the cokernel of the morphism
of $\C$\+contramodules $(\ker f_1\;0)\rarrow(\P_1,0)=\P$ is
the $\C$\+contramodule $(\boQ,\boQ)$ isomorphic to~$\Q$.
 The cokernel of the morphism~$f$ is the contramodule $(0,0)$, and
the kernel of the cokernel of~$f$ is $(\boQ,\boQ)$.
 So the abelian category axiom holds for the morphism~$f$.

 In order to construct a morphism in $\C\contra$ for which the abelian
category axiom does \emph{not} hold, let us modify the morphism~$f$
as follows.
 Put $\P'=(\P_1,\boQ)=\P\oplus(0,\boQ)$.
 Let the map $f'_1\:\P'_1=\P_1\rarrow\Q_1$ be equal to~$f_1$,
while $f'_2\:\P'_2\rarrow\Q_2$ is the identity map.
 Then the kernel of $\C$\+contramodule morphism $f'\:\P'\rarrow\Q$
coincides with the kernel of the morphism~$f$, and the cokernel of~$f'$
coincides with the cokernel of~$f$.
 Thus the kernel of the cokernel of~$f'$ is $(\boQ,\boQ)=\Q$,
while the cokernel of the kernel of~$f'$ is $(\boQ,\boQ^2)=
\Q\oplus(0,\boQ)$.

 Similarly one can construct a nonflat coring $\C$ for which
the category $\C\comodl$ is not abelian.
 For a coring $\C=\C_{11}\oplus\C_{12}\oplus\C_{22}$ with
$\C_{ii}\simeq A$ as above, a left $\C$\+comodule is the same
thing as a pair of left $A$\+modules $\M_1$ and $\M_2$ endowed
with an $A$\+module morphism $\M_1\rarrow\C_{12}\ot_A\M_2$.
 The cokernel and the kernel of an arbitrary morphism of left
$\C$\+comodules can be computed in the way dual-analogous
to the above computation for $\C$\+contramodules.

 Specifically, the cokernel of a morphism of left $\C$\+comodules
$(g_1,g_2)\:(\L_1,\L_2)\rarrow(\M_1,\M_2)$ (taken in the additive
category $\C\comodl$) is the $\C$\+comodule $(\coker g_1\;\allowbreak
\coker g_2)$.
 The kernel of the morphism $(g_1,g_2)$ can be computed as
the $\C$\+comodule $(\K\;\ker g_2)$, where $\K$ is the kernel of
the morphism to $\M_1$ from the fibered product $(\C_{12}\ot_A\ker g_2)
\sqcap_{(\C_{12}\ot_A\L_2)}\L_1$.

 Setting $A=\boZ$ and $\C_{12}=\boZ/n$ with any $n\ge2$,
let us first consider the morphism of left $\C$\+comodules
$(g_1,g_2)\:(\L_1,\L_2)\rarrow(\M_1,\M_2)$ with $\L_1=\L_2=\boZ/n$,
the structure morphism $\L_1\rarrow\L_2/n\L_2$ being the identity map,
$\M_1=0$, $\M_2=\boZ/n^2$, and an injective map $g_2\:\L_2\rarrow\M_2$.
 Then the cokernel of~$g$ is the $\C$\+comodule $(0,\boZ/n)$ and,
surprizingly, the kernel of the cokernel of~$g$ is $(\boZ/n,\boZ/n)
\simeq\L$; while the kernel of~$g$ is $(0,0)$ and the cokernel of
the kernel of~$g$ is $(\boZ/n,\boZ/n)$.
 So the abelian category axiom holds for the morphism~$g$.

 In order to construct a morphism in $\C\comodl$ for which the abelian
category axiom does \emph{not} hold, we modify the morphism~$g$
as follows.
 Put $\M'=(\boZ/n,\M_2)=\M\oplus(\boZ/n,0)$; so the structure morphism
$\M'_1\rarrow\M'_2/n\M'_2$ is the zero map.
 Let the map $g'_2\:\L_2\rarrow\M_2=\M'_2$ be equal to~$g_2$, while
$g'_1\:\L_1\rarrow\M'_1$ is the identity map.
 Then the cokernel of the $\C$\+comodule morphism $g'\:\L\rarrow\M'$
coincides with the cokernel of the morphism~$g$, and the kernel of~$g'$
coincides with the kernel of~$g$.
 Thus the cokernel of the kernel of~$g'$ is $(\boZ/n,\boZ/n)=\L$,
while the kernel of the cokernel of~$g'$ is $(\boZ^2/n,\boZ/n)=
\L\oplus(\boZ/n,0)$.
\end{ex}

 It is clear from the example above that the category of all left
contramodules over a flat coring $\C$ does not have good homological
properties in general---at least, unless one restricts the class
of exact sequences under consideration to those whose exactness
is preserves by the functor $\Hom_A(\C,{-})$.
 Our preference is to restrict the class of contramodules instead
(or rather, at the same time).
 So we will be interested in the category $\C\contra^\Acot$ of
left $\C$\+contramodules whose underlying left $A$\+modules are
cotorsion modules (see Section~\ref{cotorsion-modules}).

 Assuming that $\C$ is a flat left $A$\+module, the category
$\C\contra^\Acot$ has a natural exact category structure where a short
sequence of contramodules is exact if and only if its underlying short
sequence of $A$\+modules is exact in the abelian category $A\modl$,
or equivalently, in the exact category $A\modl^\cot$.
 We denote the $\Ext$ groups computed in the exact category
$\C\contra^\Acot$ by $\Ext^{\C,*}(\P,\Q)$.
 Assuming that $\C$ is a flat right $A$\+module, so the category
$\C\comodl$ is abelian, we denote the $\Ext$ groups computed in
the abelian category $\C\comodl$ by $\Ext_\C^*(\L,\M)$.

 Given a left $A$\+module $U$, the left $\C$\+comodule $\C\ot_A U$
is said to be \emph{coinduced} from~$U$.
 For any left $\C$\+comodule $\L$, there is a natural isomorphism
$\Hom_\C(\L\;\C\ot_A\nobreak U)\simeq\Hom_A(\L,U)$.
 Given a left $A$\+module $V$, the left $\C$\+contramodule
$\Hom_A(\C,V)$ is said to be \emph{induced} from~$V$.
 For any left $\C$\+contramodule $\Q$, there is a natural
isomorphism $\Hom^\C(\Hom_A(\C,V),\Q)\simeq\Hom_A(V,\Q)$
\cite[Sections~1.1.2 and~3.1.2]{Psemi}.
 By Lemma~\ref{cotors-hom-nc}(a), the induced $\C$\+contramodule
$\Hom_A(\C,V)$ is $A$\+cotorsion whenever the coring $\C$ is
a flat left $A$\+module and the $A$\+module $V$ is cotorsion.

 Notice that, assuming $\C$ to be a flat right $A$\+module, the direct
summands of $\C$\+comodules coinduced from injective left $A$\+modules
are the injective objects of the abelian category $\C\comodl$, and
there are enough of them.
 We will denote the full subcategory of injective objects in
the abelian category $\C\comodl$ by $\C\comodl^\inj$.

 Similarly, assuming $\C$ to be a flat left $A$\+module, the direct
summands of $\C$\+con\-tra\-mod\-ules induced from flat cotorsion left
$A$\+modules are the projective objects of the exact category
$\C\contra^\Acot$, and there are enough of them
(as one can see using Theorem~\ref{flat-cover-thm}(b)).
 We will denote the full subcategory of projective objects in
the exact category $\C\contra^\Acot$ by $\C\contra^\Acot_\prj$.

 When $\C$ is a projective left $A$\+module, the direct summands of
$\C$\+con\-tra\-mod\-ules induced from projective left $A$\+modules
are the projective objects of the abelian category $\C\contra$, and
there are enough of them.
 We will denote the full subcategory of projective objects in
the abelian category $\C\contra$ by $\C\contra_\prj$.

 The following two theorems describe the Becker and Positselski
coderived and contraderived categories of comodules and contramodules.

\begin{thm} \label{comodules-contramodules-becker-co-contra-derived}
\textup{(a)} Let\/ $\C$ be a coring over a ring $A$ such that the right
$A$\+module\/ $\C$ is flat.
 Then the composition of triangulated functors\/
$$
 \Hot(\C\comodl^\inj)\lrarrow\Hot(\C\comodl)\lrarrow\sD^\bco(\C\comodl)
$$
is a triangulated equivalence\/ $\Hot(\C\comodl^\inj)\simeq
\sD^\bco(\C\comodl)$. \par
\textup{(b)} Let\/ $\C$ be a coring over a ring $A$ such that the left
$A$\+module\/ $\C$ is projective.
 Then the composition of triangulated functors\/
$$
 \Hot(\C\contra_\prj)\lrarrow\Hot(\C\contra)\lrarrow\sD^\bctr(\C\contra)
$$
is a triangulated equivalence\/ $\Hot(\C\contra_\prj)\simeq
\sD^\bctr(\C\contra)$.
\end{thm}

\begin{proof}
 Part~(a) is a special case of
Theorem~\ref{coderived-of-grothendieck-contraderived-of-lpacepo}(a).
 Part~(b) is a special case of
Theorem~\ref{coderived-of-grothendieck-contraderived-of-lpacepo}(b).
\end{proof}

\begin{qst}
 Can one prove a version of
Theorem~\ref{comodules-contramodules-becker-co-contra-derived}(b)
for the exact category $\C\contra^\Acot$, under the assumption that
$\C$ is a flat left $A$\+module?
 The set-theoretical argument proving
Theorem~\ref{coderived-of-grothendieck-contraderived-of-lpacepo}(b)
does not seem to be easily applicable in this case.
 In particular, the exact category $\C\contra^\Acot$, and even
the exact category $A\modl^\cot$, is usually \emph{not} ``efficient''
in the sense of~\cite{Sto-ICRA}, if only because it does not have
infinite direct sums.
 All we can prove in the direction of this question is the following
Theorem~\ref{comods-contramods-positselski-co-contra-derived}(c).
\end{qst}

\begin{thm} \label{comods-contramods-positselski-co-contra-derived}
\textup{(a)} Let\/ $\C$ be a coring over a ring $A$ such that the right
$A$\+module\/ $\C$ is flat.
 Assume additionally that countable direct sums of injective left
$A$\+modules have finite injective dimensions.
 Then the classes of Positselski-coacyclic and Becker-coacyclic
complexes in the abelian category\/ $\C\comodl$ coincide.
 The natural functor\/ $\Hot(\C\comodl^\inj)\rarrow \sD^{\co=\bco}
(\C\comodl)$ is an equivalence of triangulated categories. \par
\textup{(b)} Let\/ $\C$ be a coring over a ring $A$ such that the left
$A$\+module\/ $\C$ is projective.
 Assume additionally that countable products of projective left
$A$\+modules have finite projective dimensions (for example, this holds
if the ring $A$ is right coherent and flat left $A$\+modules have
finite projective dimensions).
 Then the classes of Positselski-contraacyclic and Becker-contraacyclic
complexes in the abelian category\/ $\C\contra$ coincide.
 The natural functor\/ $\Hot(\C\contra_\prj)\rarrow\sD^{\ctr=\bctr}
(\C\contra)$ is an equivalence of triangulated categories. \par
\textup{(c)} Let\/ $\C$ be a coring over a ring $A$ such that the left
$A$\+module\/ $\C$ is flat.
 Assume additionally that countable products of flat left $A$\+modules
have finite flat dimensions (for example, this holds if the ring $A$
is right coherent).
 Then the classes of Positselski-contraacyclic and Becker-contraacyclic
complexes in the exact category\/ $\C\contra^\Acot$ coincide.
 The natural functor\/ $\Hot(\C\contra^\Acot_\prj)\rarrow
\sD^{\ctr=\bctr}(\C\contra^\Acot)$ is an equivalence of triangulated
categories.
\end{thm}

\begin{proof}
 Part~(a) is a special case of
Theorem~\ref{positselski-becker-co-contra-derived}(a).
 Parts~(b\+c) are special cases of
Theorem~\ref{positselski-becker-co-contra-derived}(b).
 In particular, in the case of part~(c) it is helpful to observe that
any cotorsion module of finite flat dimension has a finite resolution
by flat cotorsion modules (in view of Theorem~\ref{flat-cover-thm}(b)).
\end{proof}

 Recall that the \emph{contratensor product} $\N\ocn_\C\P$ of a right
$\C$\+comodule $\N$ and a left $\C$\+contramodule $\P$ is an abelian
group constructed as the cokernel of the natural pair of maps
$\N\ot_A\Hom_A(\C,\P)\birarrow\N\ot_A\P$, one of which is induced
by the contraaction map $\Hom_A(\C,\P)\rarrow\P$, while the other
one is the composition of the maps induced by the coaction map
$\N\rarrow\N\ot_A\C$ and the evaluation map $\C\ot_A\Hom_A(\C,\P)
\rarrow\P$.
 For any right $\C$\+comodule $\N$ and any left $A$\+module $V$,
there is a natural isomorphism $\N\ocn_\C\Hom_A(\C,V)\simeq
\N\ot_A V$ \cite[Sections~0.2.6 and~5.1.1\+-2]{Psemi}.

 Given two corings $\C$ and $\E$ over associative rings $A$ and $B$,
and a $\C$\+$\E$\+bicomodule $\K$, the rules $\P\mpsto\K\ocn_\E\P$
and $\M\mpsto\Hom_\C(\K,\M)$ define a pair of adjoint functors
between the categories of left $\E$\+contramodules and left
$\C$\+comodules.
 In the particular case of $\C=\E=\K$, the corresponding functors
are denoted by $\Phi_\C\:\C\contra\rarrow\C\comodl$ and
$\Psi_\C\:\C\comodl\rarrow\C\contra$, so $\Phi_\C(\P)=\C\ocn_\C\P$
and $\Psi_\C(\M)=\Hom_\C(\C,\M)$.
 The functors $\Phi_\C$ and $\Psi_\C$ transform the induced left
$\C$\+contramodule $\Hom_A(\C,U)$ into the coinduced left
$\C$\+comodule $\C\ot_A U$ and back, inducing an equivalence
between the full subcategories of comodules and contramodules of
this form in $\C\comodl$ and $\C\contra$ \cite{BBW},
\cite[Section~5.1.3]{Psemi}, \cite[Section~3.4]{Prev}.

\subsection{Base rings of finite weak dimension}
 Let $\C$ be a coring over an associative ring~$A$.
 In this section we assume that $\C$ is a flat left and right
$A$\+module and, additionally, that $A$ is a ring of finite
weak dimension (i.~e., the functor $\Tor^A({-},{-})$ has finite
homological dimension).
 Equivalently, the latter condition means that all $A$\+modules have
finite flat dimensions.
 Notice that the injective dimension of any cotorsion $A$\+module
is also finite in this case.

 The content of this section is very close to that
of~\cite[Chapter~5]{Psemi}, and, partly, \cite[Section~9.1]{Psemi}.
 The main difference is that the left $A$\+module projectivity
assumption on the coring~$\C$ used in~\cite{Psemi} is weakened here
to the flatness assumption.
 That is why the duality-analogy between comodules and contramodules
is more obscure here than in \emph{loc.\ cit}.
 Also, the form of the presentation below may be more in line with
the main body of this book than with~\cite{Psemi}.

\begin{lem}  \label{finite-weak-dim-flat-inj-resolutions}
\textup{(a)} Any\/ $\C$\+comodule can be presented as the quotient
comodule of an $A$\+flat\/ $\C$\+comodule by a finitely iterated
extension of\/ $\C$\+comodules coinduced from cotorsion $A$\+modules.
\par
\textup{(b)} Any $A$\+cotorsion\/ $\C$\+contramodule has
an admissible monomorphism into an $A$\+injective\/ $\C$\+contramodule
such that the cokernel is a finitely iterated extension of\/
$\C$\+contramodules induced from cotorsion $A$\+modules.
\end{lem}

\begin{proof}
 Part~(a) is provable by the argument of~\cite[Lemma~1.1.3]{Psemi}
used together with the result of Theorem~\ref{flat-cover-thm}(b).
 Part~(b) is similar to~\cite[Lemma~3.1.3(b)]{Psemi}.
\end{proof}

 A left $\C$\+comodule $\M$ is said to be \emph{cotorsion} if
the functor $\Hom_\C({-},\M)$ takes short exact sequences of
$A$\+flat $\C$\+comodules to short exact sequences of abelian
groups.
 In particular, any $\C$\+comodule coinduced from a cotorsion
$A$\+module is cotorsion.

 An $A$\+cotorsion left $\C$\+contramodule $\P$ is said to be
\emph{projective relative to\/~$A$} (\emph{$\C/A$\+projective})
if the functor $\Hom^\C(\P,{-})$ takes short exact sequences
of $A$\+injective $\C$\+contramodules to short exact sequences
of abelian groups.
 In particular, any $\C$\+contramodule induced from a cotorsion
$A$\+module is $\C/A$\+projective.

\begin{cor}  \label{finite-weak-dim-cotors-rel-proj-characterizations}
\textup{(a)} A left\/ $\C$\+comodule\/ $\M$ is cotorsion if and only
if\/ $\Ext_\C^{>0}(\L,\M)\allowbreak=0$ for any $A$\+flat left\/
$\C$\+comodule\/~$\L$.
 In particular, the functor\/ $\Hom_\C(\L,{-})$ takes short exact
sequences of cotorsion left\/ $\C$\+comodules to short exact
sequences of abelian groups.
 The class of cotorsion\/ $\C$\+comodules is closed under extensions
and the passage to cokernels of injective morphisms. \par
\textup{(b)} An $A$\+cotorsion left\/ $\C$\+contramodule\/ $\P$ is\/
$\C/A$\+projective if and only if $\Ext^{\C,>0}(\P,\Q)=0$ for any
$A$\+injective left\/ $\C$\+contramodule\/~$\Q$.
 In particular, the functor\/ $\Hom^\C({-},\Q)$ takes short exact
sequences of\/ $\C/A$\+projective $A$\+cotorsion\/ $\C$\+contramodules
to short exact sequences of abelian groups.
 The class of\/ $\C/A$\+projec\-tive $A$\+cotorsion\/ $\C$\+contramodules
is closed under extensions and the passage to kernels of admissible
epimorphisms in\/ $\C\contra^\Acot$.
\emergencystretch=1.5em \hbadness=1300
\end{cor}

\begin{proof}
 Follows from there being enough $A$\+flat $\C$\+comodules in
$\C\comodl$ and $A$\+in\-jective $\C$\+contramodules in
$\C\contra^\Acot$, i.~e., weak forms of the assertions of
Lemma~\ref{finite-weak-dim-flat-inj-resolutions}
(cf.~\cite[Lemma~5.3.1]{Psemi}).
\end{proof}

 It follows, in particular, that the full subcategories of cotorsion
$\C$\+comodules and $\C/A$\+projective $A$\+cotorsion $\C$\+contramodules
in $\C\comodl$ and $\C\contra^\Acot$ can be endowed with the induced
exact category structures.
 We denote these exact categories by $\C\comodl^\cot$ and
$\C\contra^\Acot_\CApr$, respectively.

\begin{lem}  \label{finite-weak-dim-cotors-proj-resolutions}
\textup{(a)} Any\/ $\C$\+comodule admits an injective morphism into
a finitely iterated extension of\/ $\C$\+comodules coinduced from
cotorsion $A$\+modules such that the quotient\/ $\C$\+comodule
is $A$\+flat. \par
\textup{(b)} For any $A$\+cotorsion\/ $\C$\+contramodule there exists
an admissible epimorphism onto it from a finitely iterated extension
of\/ $\C$\+contramodules induced from cotorsion $A$\+modules such that
the kernel is an $A$\+injective\/ $\C$\+contramodule.
\end{lem}

\begin{proof}
 The assertions follow from
Lemma~\ref{finite-weak-dim-flat-inj-resolutions} together with
the existence of enough comodules coinduced from cotorsion modules
in $\C\comodl$ and enough contramodules induced from cotorsion
modules in $\C\contra^\Acot$ by virtue of the Salce lemma
(Lemma~\ref{salce-lemma}).
 An alternative argument (providing somewhat weaker assertions) can be
found in~\cite[Lemma~9.1.2]{Psemi}.
\end{proof}

\begin{cor}  \label{cotors-rel-proj-co-contra-cor}
\textup{(a)} A\/ $\C$\+comodule is cotorsion if and only if it is 
a direct summand of a finitely iterated extension of\/ $\C$\+comodules
coinduced from cotorsion $A$\+modules. \par
\textup{(b)} An $A$\+cotorsion\/ $\C$\+contramodule is\/
$\C/A$\+projective if and only if it is a direct summand of
a finitely iterated extension of\/ $\C$\+contramodules induced from
cotorsion $A$\+modules.
\end{cor}

\begin{proof}
 Follows from Lemma~\ref{finite-weak-dim-cotors-proj-resolutions} and
Corollary~\ref{finite-weak-dim-cotors-rel-proj-characterizations}.
\end{proof}

\begin{cor}
\textup{(a)} A\/ $\C$\+comodule is simultaneously cotorsion and
$A$\+flat if and only if it is a direct summand of a\/ $\C$\+comodule
coinduced from a flat cotorsion $A$\+module. \par
\textup{(b)} An $A$\+cotorsion\/ $\C$\+contramodule is simultaneously\/
$\C/A$\+projective and $A$\+in\-jective if and only if it is a direct
summand of a\/ $\C$\+contramodule induced from an injective $A$\+module.
\end{cor}

\begin{proof}
 Both the ``if'' assertions are obvious.
 To prove the ``only if'' in part~(a), consider an $A$\+flat cotorsion
$\C$\+comodule~$\M$.
 Using Theorem~\ref{flat-cover-thm}(a), pick an injective morphism
$\M\rarrow P$ from $\M$ into a cotorsion $A$\+module $P$ such that
the cokernel $P/\M$ is a flat $A$\+module.
 Clearly, $P$ is a flat cotorsion $A$\+module.
 The cokernel of the composition of $\C$\+comodule monomorphisms
$\M\rarrow\C\ot_A\M\rarrow\C\ot_A P$, being an extension of two
$A$\+flat $\C$\+comodules, is also $A$\+flat.
 (Notice that the natural injective $\C$\+comodule morphism $\M\rarrow
\C\ot_A\M$ is naturally split as an $A$\+module map, with the splitting
induced by the counit of~$\C$.)
 According to
Corollary~\ref{finite-weak-dim-cotors-rel-proj-characterizations}(a),
it follows that the $\C$\+comodule $\M$ is a direct summand of
$\C\ot_A P$.
 
 To prove the ``only if'' in part~(b), consider a $\C/A$\+projective
$A$\+injective $\C$\+contra\-module~$\P$.
 The natural morphism of $\C$\+contramodules $\Hom_A(\C,\P)\rarrow\P$
is surjective with an $A$\+injective kernel.
 (Notice that the natural surjective $\C$\+contramodule morphism
$\Hom_A(\C,\P)\rarrow\P$ is naturally split as an $A$\+module map,
with the splitting induced by the counit of~$\C$.)
 It remains to apply
Corollary~\ref{finite-weak-dim-cotors-rel-proj-characterizations}(b)
in order to conclude that the extension splits in $\C\contra^\Acot$.
\end{proof}

\begin{thm}  \label{finite-weak-dim-exact-co-contra}
\textup{(a)} For any cotorsion left\/ $\C$\+comodule\/ $\M$, the left\/
$\C$\+contramodule\/ $\Psi_\C(\M)$ is a\/ $\C/A$\+projective
$A$\+cotorsion\/ $\C$\+contramodule. \par
\textup{(b)} For any\/ $\C/A$\+projective $A$\+cotorsion left\/
$\C$\+contramodule\/ $\P$, the left\/ $\C$\+comodule\/ $\Phi_\C(\P)$
is cotorsion. \par
\textup{(c)} The functors\/ $\Psi_\C$ and\/ $\Phi_\C$ restrict to
mutually inverse equivalences between the exact subcategories\/
$\C\comodl^\cot\sub\C\comodl$ and\/ $\C\contra^\Acot_\CApr\sub
\C\contra^\Acot$.
\end{thm}

\begin{proof}
 The functor $\Psi_\C$ takes cotorsion left $\C$\+comodules $\M$ to
$A$\+cotorsion left $\C$\+con\-tramodules, since the functor
$\Hom_A(F,\Hom_\C(\C,\M))\simeq\Hom_\C(\C\ot_A F\;\M)$ is exact
on the exact category of flat left $A$\+modules~$F$.
 Furthermore, the functor $\Psi_\C\:\M\mpsto\Hom_\C(\C,\M)$ takes
short exact sequences of cotorsion $\C$\+comodules to short
exact sequences of $A$\+cotorsion $\C$\+contramodules, since
$\C$ is a flat left $A$\+module.

 Similarly, the functor $\Phi_\C$ takes short exact sequences of
$\C/A$\+projective $A$\+cotorsion $\C$\+contramodules to short
exact sequences of $\C$\+comodules, since $\Hom_\boZ(\C\ocn_\C\P\;
\allowbreak\boQ/\boZ)\simeq\Hom^\C(\P,\Hom_\boZ(\C,\boQ/\boZ))$
and the left $\C$\+contramodule $\Hom_\boZ(\C,\boQ/\boZ)$ is
$A$\+injective ($\C$ being a flat right $A$\+module).

 In view of these observations, all the assertions follow from
Corollary~\ref{cotors-rel-proj-co-contra-cor}.
 Alternatively, one could proceed along the lines of the proof
of~\cite[Theorem~5.3]{Psemi}, using the facts that any cotorsion
$\C$\+comodule has finite injective dimension in $\C\comodl$ and
any $\C/A$\+projective $A$\+cotorsion $\C$\+contramodule
has finite projective dimension in $\C\contra^\Acot$ (see the proof
of the next Theorem~\ref{finite-weak-dim-co-contra-derived}).
\end{proof}

\begin{thm} \label{finite-weak-dim-co-contra-derived}
\textup{(a)} The triangulated functor\/ $\sD^\abs(\C\comodl^\cot)
\rarrow\sD^\bco(\C\comodl)$ induced by the embedding of exact
categories\/ $\C\comodl^\cot\rarrow\C\comodl$ is an equivalence
of triangulated categories. \par
\textup{(b)} Assume additionally that countable direct sums of injective
left $A$\+modules have finite injective dimensions.
 Then the triangulated functor\/ $\sD^\abs(\C\comodl^\cot)
\rarrow\sD^{\co=\bco}(\C\comodl)$ induced by the embedding of exact
categories\/ $\C\comodl^\cot\rarrow\C\comodl$ is an equivalence
of triangulated categories. \par
\textup{(c)} The triangulated functor\/ $\sD^\abs(\C\contra^\Acot_\CApr)
\rarrow\sD^{\ctr=\bctr}(\C\contra^\Acot)$ induced by the embdding of
exact categories\/ $\C\contra^\Acot_\CApr\rarrow\C\contra^\Acot$ is
an equivalence of triangulated categories.
\end{thm}

\begin{proof}
 In parts~(a\+b), we notice that the injective dimension of any
cotorsion $\C$\+comodule $\M$ (as an object of $\C\comodl$) does not
exceed the weak homological dimension of the ring~$A$.
 Indeed, one can compute the functor $\Ext_\C^*({-},\M)$ using
$A$\+flat resolutions of the first argument (which exist
by Lemma~\ref{finite-weak-dim-flat-inj-resolutions}(a)).
 Now part~(a) can be obtained by comparing
Theorem~\ref{comodules-contramodules-becker-co-contra-derived}(a)
with Theorem~\ref{finite-homol-dim-becker-co-contra-derived}(a)
for the exact category $\sA=\C\comodl^\cot$.

 Part~(b) can be deduced by comparing part~(a) with
Theorem~\ref{comods-contramods-positselski-co-contra-derived}(a).
 Alternatively, given the description of injective $\C$\+comodules in
Section~\ref{contra-over-flat}, it follows from the assumptions
of part~(b) that countable direct sums of cotorsion $\C$\+comodules
have finite injective dimensions.
 It remains to apply the dual version of 
Corollary~\ref{finite-homol-dim-equivalence-cor} (or the remark
in the paragraph after it).

 Similarly, to prove part~(c) one first notices that the projective
dimension of any $\C/A$\+projective $A$\+cotorsion $\C$\+contramodule
$\P$ does not exceed the supremum of the injective dimensions of
cotorsion $A$\+modules, i.~e., the weak homological dimension of
the ring~$A$.
 Indeed, one can compute the functor $\Ext^{\C,*}(\P,{-})$ using
$A$\+injective coresolutions of the second argument (which exist
by Lemma~\ref{finite-weak-dim-flat-inj-resolutions}(b)).
 It remains to compare
Theorem~\ref{comods-contramods-positselski-co-contra-derived}(c)
with Theorem~\ref{finite-homol-dim-becker-co-contra-derived}(b)
for the exact category $\sB=\C\contra^\Acot_\CApr$.

 Alternatively, one can see from the description of projective
$A$\+cotorsion $\C$\+con\-tra\-mod\-ules
in Section~\ref{contra-over-flat} that infinite products of
$\C/A$\+projective $A$\+cotorsion $\C$\+contramodules have finite
projective dimensions in $\C\contra^\Acot$.
 Hence Corollary~\ref{finite-homol-dim-equivalence-cor} (or
the subsequent remark) applies.

 One could also argue in the way similar to the proof
of~\cite[Theorem~5.4]{Psemi} (cf.\ 
Theorem~\ref{gorenstein-base-co-contra-derived} below).
\end{proof}

\begin{cor}
\textup{(a)} The Becker coderived category of the abelian category
of left\/ $\C$\+comodules\/ $\sD^\bco(\C\comodl)$ and the contraderived
category of the exact category of $A$\+cotorsion left\/
$\C$\+contramodules\/ $\sD^{\ctr=\bctr}(\C\contra^\Acot)$
are naturally equivalent.
 The equivalence is provided by the derived functors of co-contra
correspondence\/ $\boR\Psi_\C$ and\/ $\boL\Phi_\C$ constructed in terms
of cotorsion coresolutions of comodules and relatively projective
resolutions of contramodules. \par
\textup{(b)}
 Assume additionally that countable direct sums of injective
left $A$\+modules have finite injective dimensions.
 Then the coderived category of left\/ $\C$\+comodules\/
$\sD^{\co=\bco}(\C\comodl)$ and the contraderived category of
$A$\+cotorsion left\/ $\C$\+contramodules\/
$\sD^{\ctr=\bctr}(\C\contra^\Acot)$ are naturally equivalent.
 The equivalence is provided by the derived functors of co-contra
correspondence\/ $\boR\Psi_\C$ and\/ $\boL\Phi_\C$.
\end{cor}

\begin{proof}
 Follows from Theorems~\ref{finite-weak-dim-exact-co-contra}
and~\ref{finite-weak-dim-co-contra-derived}.
\end{proof}

\subsection{Gorenstein base rings}
 Let $\C$ be a coring over an associative ring~$A$.
 We assume the ring $A$ to be left Gorenstein, in the sense that
the classes of left $A$\+modules of finite injective dimension and
of finite flat dimension coincide.
 In this case, it is clear that both kinds of dimensions are uniformly
bounded by a constant for those modules for which they are finite.

 More generally, it will be sufficient to require that the classes of
cotorsion left $A$\+modules of finite injective dimension and
of finite flat dimension coincide, countable direct sums of injective
left $A$\+modules have finite injective dimensions, and flat cotorsion
$A$\+modules have uniformly bounded injective dimensions.
 Then it follows that the flat and injective dimensions of
cotorsion $A$\+modules of finite flat/injective dimension
are also uniformly bounded.

\begin{thm}  \label{gorenstein-base-co-contra-derived}
\textup{(a)} Assume that the coring\/ $\C$ is a flat right $A$\+module.
 Then the coderived category\/ $\sD^\co(\C\comodl)$ of
the abelian category of left\/ $\C$\+comodules is equivalent to
the quotient category of the homotopy category of the additive
category of left\/ $\C$\+comodules coinduced from cotorsion
$A$\+modules of finite flat/injective dimension by its
minimal thick subcategory containing the total complexes of
short exact sequences of complexes of\/ $\C$\+comodules that
at every term of the complexes are short exact sequences of\/
$\C$\+comodules coinduced from short exact sequences of
cotorsion $A$\+modules of finite flat/injective dimension. \par
\textup{(b)} Assume that the coring\/ $\C$ is a flat left $A$\+module.
 Then the contraderived category\/ $\sD^\ctr(\C\contra^\Acot)$
of the exact category of $A$\+co\-torsion left\/ $\C$\+contramodules
is equivalent to the quotient category of the homotopy category of
the additive category of left\/ $\C$\+contramodules induced from
cotorsion $A$\+modules of finite flat/injective dimension by its
minimal thick subcategory cointaining the total complexes of
short exact sequences of complexes of\/ $\C$\+contramodules that
at every term of the complexes are short exact sequences of\/
$\C$\+contramodules induced from short exact sequences of
cotorsion $A$\+modules of finite flat/injective dimension.
\end{thm}

\begin{proof}
 The argument proceeds along the lines of the proof
of~\cite[Theorem~5.5]{Psemi} (see also~\cite[Question~5.4]{Psemi}
and~\cite[Sections~3.9\+-3.10]{Pkoszul}).

 Part~(a): the relative cobar-coresolution
$$
 0\lrarrow\M^\bu\lrarrow\C\ot_A\M^\bu\lrarrow\C\ot_A\C\ot_A\M^\bu
 \lrarrow\dotsb,
$$
totalized by taking infinite direct sums along the diagonals, provides
a closed morphism with a coacyclic cone from any complex of left
$\C$\+comodules $\M^\bu$ into a complex of coinduced $\C$\+comodules.
 The construction from the proof of~\cite[Theorem~5.5]{Psemi} provides
a closed morphism from any complex of coinduced left $\C$\+comodules
into a complex of $\C$\+comodules termwise coinduced from injective
$A$\+modules such that this morphism is coinduced from an injective
morphism of $A$\+modules at every term of the complexes.

 This allows to obtain a closed morphism with a coacyclic cone from
any complex of coinduced left $\C$\+comodules into a complex of
$\C$\+comodules termwise coinduced from injective $A$\+modules (using
the assumption that countable direct sums of injective $A$\+modules
have finite injective dimensions).
 The same construction from~\cite{Psemi} can be also used to obtain
a closed morphism from any complex of left $\C$\+comodules with
the terms coinduced from cotorsion $A$\+modules of finite injective
dimension into a complex of $\C$\+comodules termwise coinduced from
injective $A$\+modules such that the cone is homotopy equivalent to
a complex obtained from short exact sequences of complexes of
$\C$\+comodules termwise coinduced from short exact sequences of
cotorsion $A$\+modules of finite injective dimension using
the operation of cone repeatedly.
 
 The assertion of part~(a) follows from these observations by
a semi-orthogonal decomposition argument from the proofs
of~\cite[Theorems~5.4\+-5.5]{Psemi}.
 
 Part~(b): the contramodule relative bar-resolution
$$
 \dotsb\lrarrow\Hom_A(\C\ot_A\C\;\P^\bu)\lrarrow\Hom_A(\C\;\P^\bu)
 \lrarrow\P^\bu\lrarrow0,
$$
totalized by taking infinite products along the diagonals, provides
a closed morphism with a contraacyclic cone onto any complex of
$A$\+cotorsion left $\C$\+contramodules $\P^\bu$ from a complex of
$\C$\+contramodules termwise induced from cotorsion $A$\+modules.
 The construction dual to the one elaborated in the proof
of~\cite[Theorem~5.5]{Psemi} provides a closed morphism onto any complex
of $\C$\+contramodules termwise induced from cotorsion $A$\+modules from
a complex of $\C$\+contramodules termwise induced from flat cotorsion
$A$\+modules such that this morphism is induced from an admissible
epimorphism of cotorsion $A$\+modules at every term of the complexes.

 This allows to obtain a closed morphism with a contraacyclic cone
onto any complex of left $\C$\+contramodules termwise induced from
cotorsion $A$\+modules from a complex of $\C$\+contramodules
termwise induced from flat cotorsion $A$\+modules (since an infinite
product of flat cotorsion $A$\+modules, being a cotorsion $A$\+module
of finite injective dimension, has finite flat dimension).
 The same construction can be also used to obtain a closed morphism
onto any complex of left $\C$\+contramodules with the terms induced
from cotorsion $A$\+modules of finite flat dimension from a complex
of $\C$\+contramodules termwise induced from flat cotorsion $A$\+modules
such that the cone is homotopy equivalent to a complex obtained from
short exact sequences of complexes of $\C$\+contramodules termwise
induced from short exact sequences of cotorsion $A$\+modules of finite
flat dimension using the operation of cone repeatedly.

 Now the dual version of the same semi-orthogonal decomposition
argument from~\cite{Psemi} implies part~(b).
\end{proof}

\begin{cor}
 Assume that the coring\/ $\C$ is a flat left and right $A$\+module.
 Then the coderived category of left\/ $\C$\+comodules\/
$\sD^\co(\C\comodl)$ and the contraderived category of $A$\+cotorsion
left\/ $\C$\+contramodules\/ $\sD^\ctr(\C\contra^\Acot)$ are
naturally equivalent.
 The equivalence is provided by the derived functors of co-contra
correspondence\/ $\boR\Psi_\C$ and\/ $\boR\Phi_\C$ constructed in terms
of coresolutions by complexes of\/ $\C$\+comodules termwise
coinduced from cotorsion $A$\+modules of finite injective/flat
dimension and resolutions by complexes of\/ $\C$\+contramodules termwise
induced from cotorsion $A$\+modules of finite flat/injective dimension.
\end{cor}

\begin{proof}
 Follows from Theorem~\ref{gorenstein-base-co-contra-derived}
and the remarks at the end of Section~\ref{contra-over-flat}.
\end{proof}

\subsection{Corings with dualizing complexes}  \label{corings-dualizing}
 Let $A$ and $B$ be associative rings.
 We call a finite complex of $A$\+$B$\+bimodules $D^\bu$
a \emph{dualizing complex} for $A$ and~$B$ \,\cite{Yek,Miy,CFH} if
\begin{enumerate}
\renewcommand{\theenumi}{\roman{enumi}}
\item $D^\bu$ is simultaneously a complex of injective left
$A$\+modules and a complex of injective right $B$\+modules
(in the one-sided module structures obtained by forgetting
the other module structure on the other side);
\item as a complex of left $A$\+modules, $D^\bu$ is quasi-isomorphic
to a bounded above complex of finitely generated projective
$A$\+modules, and similarly, as a complex of right $B$\+modules,
$D^\bu$ is quasi-isomorphic to a bounded above complex of finitely
generated projective $B$\+modules;
\item the ``homothety'' maps $A\rarrow\Hom_{B^\rop}(D^\bu,D^\bu)$
and $B\rarrow\Hom_A(D^\bu,D^\bu)$ are quasi-isomorphisms of
complexes.
\end{enumerate}

 For similar definitions with less restrictive conditions,
see~\cite[Sections~3 and~4]{Pfp}, \cite[Sections~3 and~7]{Pps},
and Section~\ref{ind-affine-co-contra-subsect} below.

\begin{lem} \label{noncomm-ring-dualizing-lemma}
 Let $D^\bu$ be a dualizing complex for associative rings $A$ and~$B$.
 Then \par
\textup{(a)} if the ring $A$ is left Noetherian and $F$ is a flat
left $B$\+module, then the natural homomorphism of finite complexes
of left $B$\+modules $F\rarrow\Hom_A(D^\bu\;D^\bu\ot_B F)$ is
a quasi-isomorphism; \par
\textup{(b)} if the ring $B$ is right coherent and $J$ is an injective
left $A$\+module, then the natural homomorphism of finite complexes
of left $A$\+modules $D^\bu\ot_B\Hom_A(D^\bu,J)\rarrow J$ is
a quasi-isomorphism.
\end{lem}

\begin{proof}
 This is a particular case of~\cite[Lemma~4.2]{Pfp}
and~\cite[Lemma~3.2]{Pps}.

 Part~(a): first of all, $D^\bu\ot_B F$ is a complex of injective
left $A$\+modules by Lemma~\ref{coherent-tensor-hom-lemma}(a).
 Let ${}'\!\.D^\bu\rarrow D^\bu$ be a quasi-isomorphism of
complexes of left $A$\+modules between a bounded above complex
of finitely generated projective $A$\+modules ${}'\!\.D^\bu$ and
the complex~$D^\bu$.
 Then it suffices to show that the induced morphism of complexes
of abelian groups $F\rarrow\Hom_A({}'\!\.D^\bu\;D^\bu\ot_B F)$
is a quasi-isomorphism.
 Now the complex $\Hom_A({}'\!\.D^\bu\;D^\bu\ot_B F)$ is
isomorphic to $\Hom_A({}'\!\.D^\bu\;D^\bu)\ot_B F$, and
it remains to use the condition~(iii) for the morphism
$B\rarrow\Hom_A(D^\bu,D^\bu)$ together with the flatness
condition on the $B$\+module~$F$.

 Part~(b): the complex $\Hom_A(D^\bu,J)$ is a complex of flat
left $B$\+modules by Lemma~\ref{coherent-tensor-hom-lemma}(b).
 Let ${}''\!\.D^\bu\rarrow D^\bu$ be a quasi-isomorphism of
complexes of right $B$\+modules between a bounded above complex
of finitely generated projective $B$\+modules ${}''\!\.D^\bu$ and
the complex~$D^\bu$.
 It suffices to show that the induced morphism of complexes of
abelian groups ${}''\!\.D^\bu\ot_B\Hom_A(D^\bu,J)\rarrow J$
is a quasi-isomorphism.
 The complex ${}''\!\.D^\bu\ot_B\Hom_A(D^\bu,J)$ being isomorphic
to $\Hom_A(\Hom_{B^\rop}({}''\!\.D^\bu,D^\bu),J)$, it remains
to use the condition~(iii) for the morphism 
$A\rarrow\Hom_{B^\rop}(D^\bu,D^\bu)$ together with the injectivity
condition on the $A$\+module~$J$.
\end{proof}

 The following result is due to Christensen, Frankild, and
Holm~\cite[Proposition~1.5]{CFH}.
 For a generalization, see~\cite[Proposition~4.3]{Pps}.

\begin{cor}  \label{christensen-frankild-holm}
 Let $D^\bu$ be a dualizing complex for associative rings $A$ and~$B$.
 Assume that the ring $A$ is left Noetherian.
 Then the projective dimension of any flat left $B$\+module
does not exceed the length of~$D^\bu$.
\end{cor}

\begin{proof}
 Assume that the complex $D^\bu$ is concentrated in the cohomological
degrees from~$i$ to~$i+d$.
 It suffices to show that $\Ext^{d+1}_B(F,G)=0$ for any flat
left $B$\+modules $F$ and~$G$.
 Let $P_\bu$ be a projective left resolution of the $B$\+module~$F$.
 By Lemma~\ref{noncomm-ring-dualizing-lemma}(a), the natural
map of complexes of abelian groups $\Hom_B(P_\bu,G)\rarrow
\Hom_B(P_\bu\;\Hom_A(D^\bu\;D^\bu\ot_B G))$ is a quasi-isomorphism.
 The right-hand side is isomorphic to the complex
$\Hom_A(D^\bu\ot_B P_\bu\;D^\bu\ot_B G)$, which is quasi-isomorphic
to $\Hom_A(D^\bu\ot_B F\;D^\bu\ot_B G)$, since $D^\bu\ot_B G$ is
a finite complex of injective $A$\+modules.
 The corollary is proved.
 Notice that we have only used ``a~half of'' the conditions (i\+iii)
imposed on a dualizing complex~$D^\bu$.
\end{proof}

 Let $\C$ be a coring over an associative ring $A$ and $\E$ be a coring
over an associative ring~$B$.
 We assume $\C$ to be a flat right $A$\+module and $\E$ to be a flat
left $B$\+module.
 A \emph{dualizing complex} for $\C$ and $\E$ is defined as a triple
consisting of a finite complex of $\C$\+$\E$\+bicomodules $\D^\bu$, 
a finite complex of $A$\+$B$\+bimodules $D^\bu$, and a morphism of
complexes of $A$\+$B$\+bimodules $\D^\bu\rarrow D^\bu$ with
the following properties:
\begin{enumerate}
\renewcommand{\theenumi}{\roman{enumi}}
\setcounter{enumi}{3}
\item $D^\bu$ is a dualizing complex for the rings $A$ and~$B$;
\item $\D^\bu$ is a complex of injective left $\C$\+comodules
(forgetting the right $\E$\+comodule structure) and a complex
of injective right $\E$\+comodules (forgetting the left
$\C$\+comodule structure);
\item the morphism of complexes of left $\C$\+comodules
$\D^\bu\rarrow\C\ot_A D^\bu$ induced by the morphism of
complexes of left $A$\+modules $\D^\bu\rarrow D^\bu$ is
a quasi-isomorphism;
\item the morphism of complexes of right $\E$\+comodules
$\D^\bu\rarrow D^\bu\ot_B\E$ induced by the morphism of
complexes of right $B$\+modules $\D^\bu\rarrow D^\bu$ is
a quasi-isomorphism.
\end{enumerate}

\begin{lem}  \label{co-contra-inj-to-proj-lemma}
\textup{(a)} Assume that the ring $A$ is left Noetherian.
 Then for any\/ $\C$\+injective\/ $\C$\+$\E$\+bicomodule\/ $\K$
and any left\/ $\E$\+contramodule $\gF$ which is a projective object
of\/ $\E\contra^\Bcot$, the left\/ $\C$\+comodule\/ $\K\ocn_\E\gF$
is injective. \par
\textup{(b)} Assume that the ring $B$ is right coherent.
 Then for any\/ $\E$\+injective\/ $\C$\+$\E$\+bico\-module\/ $\K$
and any injective left\/ $\C$\+comodule\/ $\J$, the left\/
$\E$\+contramodule\/ $\Hom_\C(\K,\J)$ is a projective object of\/
$\E\contra^\Bcot$.
\end{lem}

\begin{proof}
 Part~(a): one can assume the left $\E$\+contramodule $\gF$ to be
induced from a flat (cotorsion) $B$\+module~$F$; then
$\K\ocn_\E\gF\simeq\K\ot_B F$.
 Hence it suffices to check, e.~g., that the class of injective left
$\C$\+comodules is preserved by filtered inductive limits.
 Let us show that $\C\comodl$ is a locally Noetherian Grothendieck
abelian category.

 The following is a standard argument (cf.~\cite[Section~18.16]{BW}).
 Let $\L$ be a left $\C$\+comodule.
 For any finitely generated $A$\+submodule $U\sub\L$, the full
preimage $\L_U$ of $\C\ot_A U\sub\C\ot_A\L$ with respect to
the $\C$\+coaction map $\L\rarrow\C\ot_A \L$ is
a $\C$\+subcomodule in $\L$ contained in~$U$.
 Since the left $\C$\+comodule $\C\ot_A\L$ is a filtered inductive
limit of its $\C$\+subcomodules $\C\ot_A U$, it follows that
the $\C$\+comodule $\L$ is a filtered inductive limit of
its $\C$\+subcomodules~$\L_U$.
 Given that $A$ is a left Noetherian ring, we can conclude that
any left $\C$\+comodule is the union of its $A$\+finitely
generated $\C$\+subcomodules.

 Part~(b): one can assume the left $\C$\+comodule $\J$ to be coinduced
from a left $A$\+module $J$, and the $A$\+module $J$ to have
the form $\Hom_\boZ(A,X)$ for a certain injective abelian group~$X$.
 Then we have $\Hom_\C(\K,\J)\simeq\Hom_A(\K,J)\simeq\Hom_\boZ(\K,X)$.
 The $\E$\+contramodule $\Hom_\boZ(\K,X)$ only depends on
the right $\E$\+comodule structure on $\K$, so one can assume $\K$
to be the right $\E$\+comodule coinduced from an injective 
right $B$\+module~$I$.
 Now $\Hom_\boZ(I\ot_B\E\;X)\simeq\Hom_B(\E,\Hom_\boZ(I,X))$ is
the left $\E$\+contramodule induced from the left $B$\+module
$\Hom_\boZ(I,X)$.
 The latter is a flat cotorsion $B$\+module by
Lemmas~\ref{cotors-hom-nc}(b) and~\ref{coherent-tensor-hom-lemma}(b).
\end{proof}

\begin{lem}  \label{coring-dualizing-lemma}
 Let\/ $\D^\bu\rarrow D^\bu$ be a dualizing complex for corings\/ $\C$
and\/~$\E$.
 Then \par
\textup{(a)} assuming that the ring $A$ is left Noetherian, for any
left\/ $\E$\+contramodule\/ $\gF$ that is a projective object of\/
$\E\contra^\Bcot$ the adjunction morphism\/ $\gF\rarrow
\Hom_\C(\D^\bu\;\D^\bu\ocn_\E\gF)$ is a quasi-isomorphism of
finite complexes over the exact category\/ $\E\contra^\Bcot$; \par
\textup{(b)} assuming that the ring $B$ is right coherent, for any
injective left\/ $\C$\+comodule\/ $\J$ the adjunction morphism\/
$\D^\bu\ocn_\E\Hom_\C(\D^\bu,\J)\rarrow\J$ is a quasi-isomorphism of
finite complexes over the abelian category\/ $\C\comodl$.
\end{lem}

\begin{proof}
 Part~(a): one can assume the $\E$\+contramodule $\gF$ to be induced
from a flat (cotorsion) left $B$\+module~$F$; then $\D^\bu\ocn_\E\gF
\simeq\D^\bu\ot_B F$.
 It follows from (the proof of)
Lemma~\ref{co-contra-inj-to-proj-lemma}(a) that
the morphism $\D^\bu\ot_B F\rarrow\C\ot_A D^\bu\ot_B F$ induced
by the quasi-isomorphism $\D^\bu\rarrow\C\ot_A D^\bu$ from~(vi) is
a quasi-isomorphism of complexes of injective left $\C$\+comodules.
 Any quasi-isomorphism of finite complexes of injective objects is
a homotopy equivalence.
 Hence we have a quasi-isomorphism $\Hom_\C(\D^\bu\;\D^\bu
\ot_B F)\rarrow\Hom_A(\D^\bu\;D^\bu\ot_B F)$.
 Furthermore, by~(vii) there is a quasi-isomorphism
$\Hom_A(D^\bu\ot_B\E\;D^\bu\ot_B F)\rarrow\Hom_A(\D^\bu\;D^\bu\ot_B F)$
and a natural isomorphism $\Hom_A(D^\bu\ot_B\E\;D^\bu\ot_B F)\simeq
\Hom_B(\E\;\Hom_A(D^\bu\;D^\bu\ot_B F))$.
 Finally, the natural morphism $\Hom_B(\E,F)\rarrow\Hom_B(\E\;
\Hom_A(D^\bu\;D^\bu\ot_B F))$ is a quasi-isomorphism by
Lemmas~\ref{cotors-hom-nc}(b)
and~\ref{noncomm-ring-dualizing-lemma}(a).
 Here the point is that the functor $\Hom_B(\E,{-})$ from a flat left
$B$\+module $\E$ preserves quasi-isomorphisms of (finite) complexes
of cotorsion left $B$\+modules.

 Part~(b): one can assume the $\C$\+comodule $\J$ to be induced from
an injective left $A$\+module~$J$; then $\Hom_\C(\D^\bu,\J)\simeq
\Hom_A(\D^\bu,J)$.
 It follows from (the proof of)
Lemma~\ref{co-contra-inj-to-proj-lemma}(b) that
the morphism $\Hom_A(D^\bu\ot_B\E\;J)\rarrow\Hom_A(\D^\bu\;J)$
induced by the quasi-isomorphism $\D^\bu\rarrow D^\bu\ot_B\E$ 
from~(vii) is a quasi-isomorphism of complexes of projective
objects in $\E\contra^\Bcot$.
 Any quasi-isomorphism of finite complexes of projective objects is
a homotopy equivalence.
 Hence we have a quasi-isomorphism $\D^\bu\ocn_\E\Hom_A(\D^\bu,J)
\larrow\D^\bu\ocn_\E\Hom_A(D^\bu\ot_B\E\;J)\simeq
\D^\bu\ocn_\E\Hom_B(\E,\Hom_A(D^\bu,J))\simeq
\D^\bu\ot_B\Hom_A(D^\bu,J)$.
 Furthermore, by~(vi) and Lemma~\ref{coherent-tensor-hom-lemma}(b)
there is a quasi-isomorphism $\D^\bu\ot_B\Hom_A(D^\bu,J)\rarrow
\C\ot_A D^\bu\ot_B\Hom_A(D^\bu,J)$.
 Now it remains to apply Lemma~\ref{noncomm-ring-dualizing-lemma}(b).
\end{proof}

 The following theorem does not depend on the existence of any
dualizing complexes.

\begin{thm}  \label{noetherian-coherent-co-contra-resolutions}
\textup{(a)} Let\/ $\C$ be a coring over an associative ring $A$;
asssume that\/ $\C$ is a flat right $A$\+module and\/ $A$ is
a left Noetherian ring.
 Then the classes of Becker-coacyclic and Positselski-coacyclic
objects in the abelian category\/ $\C\comodl$ coincide.
 The coderived category\/ $\sD^{\co=\bco}(\C\comodl)$ of the abelian
category of left\/ $\C$\+comodules is equivalent to the homotopy
category of complexes of injective left\/ $\C$\+comodules. \par
\textup{(b)} Let\/ $\E$ be a coring over an associative ring $B$;
assume that\/ $\E$ is a flat left $B$\+module and\/ $B$ is
a right coherent ring.
 Then the classes of Becker-contraacyclic and
Positselski-contraacyclic objects in the exact category\/
$\C\contra^\Bcot$ coincide.
 The contraderived category\/ $\sD^\ctr(\E\contra^\Bcot)$ of
the exact category of $B$\+cotorsion left\/ $\E$\+contramodules
is equivalent to the homotopy category of complexes of\/
projective objects in\/ $\E\contra^\Bcot$.
\end{thm}

\begin{proof}
 Part~(a) is a particular case of
Theorem~\ref{comods-contramods-positselski-co-contra-derived}(a).
 Essentially, it holds because there are enough injective objects in
$\C\comodl$ and the class of injective objects is closed under infinite
direct sums (the class of injective left $A$\+modules being closed
under infinite direct sums).
 Part~(b) is a particular case of
Theorem~\ref{comods-contramods-positselski-co-contra-derived}(c).
 Essentially, it is true because there are enough projective objects
in $\E\contra^\Bcot$ and the class of projective objects is preserved by
infinite products (the class of flat left $B$\+modules
being preserved by infinite products).
\end{proof}

 More generally, the assertion of part~(a) holds if $\C$ is a flat
right $A$\+module and countable direct sums of injective left
$A$\+modules have finite injective dimensions.
 Similarly, the assertion of~(b) is true if $\E$ is a flat
left $B$\+module and countable products of flat cotorsion left
$B$\+modules have finite flat dimensions
(cf.\ Corollary~\ref{finite-homol-dim-equivalence-cor}).

\begin{cor}  \label{noncomm-flat-dualizing-co-contra-cor}
 Let\/ $\C$ be a coring over an associative ring $A$ and\/
$\E$ be a coring over an associative ring~$B$.
 Assume that\/ $\C$ is a flat right $A$\+module, $\E$ is
a flat left $B$\+module, the ring $A$ is left Noetherian, and
the ring $B$ is right coherent.
 Then the datum of a dualizing complex\/ $\D^\bu\rarrow D^\bu$
for the corings\/ $\C$ and\/ $\E$ induces an equivalence of
triangulated categories\/ $\sD^\co(\C\comodl)\simeq
\sD^\ctr(\E\contra^\Bcot)$, which is provided by the derived functors\/
$\boR\Hom_\C(\D^\bu,{-})$ and\/ $\D^\bu\ocn_\E^\boL\nobreak{-}$.
\end{cor}

\begin{proof}
 Follows from Lemmas~\ref{co-contra-inj-to-proj-lemma}\+-%
\ref{coring-dualizing-lemma} and
Theorem~\ref{noetherian-coherent-co-contra-resolutions}.
\end{proof}

 Now let us assume that a coring $\E$ over an associative ring $B$
is a projective left $B$\+module.
 Then, as mentioned in Section~\ref{contra-over-flat}, the category
$\E\contra$ is abelian with enough projective objects; the latter are
the direct summands of the left $\E$\+contramodules induced from
projective left $B$\+modules.

\begin{lem}  \label{contra-cotorsion-contraflat}
 Assume that any left $B$\+module has finite cotorsion dimension
(or equivalently, any flat left $B$\+module has finite projective
dimension).  Then \par
\textup{(a)} any left\/ $\E$\+contramodule can be embedded into
a $B$\+cotorsion left\/ $\E$\+contramodule in such a way that
the cokernel is a finitely iterated extension of\/
$\E$\+contramodules induced from flat left $B$\+modules; \par
\textup{(b)} any left\/ $\E$\+contramodule can be presented as
the quotient contramodule of a finitely iterated extension of\/
$\E$\+contramodules induced from flat left $B$\+modules by
a $B$\+cotorsion left\/ $\E$\+contramodule.
\end{lem}

\begin{proof}
 The proof of part~(a) is similar to that
of~\cite[Lemma~3.1.3(b)]{Psemi} and uses
Theorem~\ref{flat-cover-thm}(a).
 The proof of part~(b) is similar to that of
Lemma~\ref{finite-weak-dim-cotors-proj-resolutions} and based on
part~(a) and the fact that there are enough contramodules induced
from flat (and even projective) $B$\+modules in $\E\contra$
(see Lemma~\ref{salce-lemma}(a), cf.\ Lemma~\ref{clf-cover}).
\end{proof}

 Let us call a left $\E$\+contramodule $\gF$ \emph{strongly contraflat}
if the functor $\Hom^\E(\gF,{-})$ takes short exact sequences of
$B$\+cotorsion $\E$\+contramodules to short exact sequences of
abelian groups (cf.\ Section~\ref{clf-subsection}
and~\cite[Section~5.1.6 and Question~5.3]{Psemi}).

 Assuming that any flat left $B$\+module has finite projective
dimension, it follows from Lemma~\ref{contra-cotorsion-contraflat}(a)
that the Ext groups computed in the exact category $\E\contra^\Bcot$
and in the abelian category $\E\contra$ agree.
 We denote these by $\Ext^{\E,*}({-},{-})$.

\begin{cor}  \label{strongly-contraflat-cor}
 Assume that any flat left $B$\+module has finite projective dimension.
 Then \par
\textup{(a)} A left\/ $\E$\+contramodule\/ $\gF$ is strongly contraflat
if and only if\/ $\Ext^{\E,>0}(\gF,\Q)=0$ for any $B$\+cotorsion left\/
$\E$\+contramodule\/~$\Q$.
 The class of strongly contraflat\/ $\E$\+contramodules is closed
under extensions and the passage to kernels of surjective morphisms
in\/ $\E\contra$. \par
\textup{(b)} A left\/ $\E$\+contramodule is strongly contraflat if and
only if it is a direct summand of a finitely iterated extension of\/
$\E$\+contramodules induced from flat $B$\+modules. \par
\end{cor}

\begin{proof}
 Part~(a) follows from (a weak version of)
Lemma~\ref{contra-cotorsion-contraflat}(a).
 Part~(b) follows from Lemma~\ref{contra-cotorsion-contraflat}(b)
and part~(a).
\end{proof}

\begin{thm}  \label{projective-coring-contra-resolutions}
\textup{(a)} Assuming that any flat left $B$\+module has finite
projective dimension, for any symbol\/ $\bst=\b$, $+$, $-$, $\empt$,
$\abs+$, $\abs-$, $\bctr$, $\ctr$, or~$\abs$ the triangulated functor\/
$\sD^\st(\E\contra^\Bcot)\rarrow\sD^\st(\E\contra)$ induced by
the embedding of exact categories $\E\contra^\Bcot\rarrow\E\contra$
is an equivalence of triangulated categories. \par
\textup{(b)} Assuming that countable products of projective left
$B$\+modules have finite projective dimensions (in particular, if
the ring $B$ is right coherent and flat left $B$\+modules have
finite projective dimensions), the contraderived category\/
$\sD^{\ctr=\bctr}(\E\contra)$ of the abelian category of left\/
$\E$\+contramodules is equivalent to the homotopy category of
complexes of projective left\/ $\E$\+contramodules. \par
\textup{(c)} Assuming that countable products of flat left
$B$\+modules have finite projective dimensions, the contraderived
category\/ $\sD^{\ctr=\bctr}(\E\contra)$ is equivalent to the absolute
derived category of the exact category of strongly contraflat left\/
$\E$\+contramodules.
\end{thm}

\begin{proof}
 Part~(a) for all the symbols~$\bst$ except $\bctr$
follows from Lemma~\ref{contra-cotorsion-contraflat}(a)
and the dual version of Proposition~\ref{finite-resolutions}.
 In the case $\bst=\bctr$, the argument is based on
Proposition~\ref{becker-contraderived-finite-coresolutions} applied
to the exact category $\sE=\E\contra$ with the hereditary complete
cotorsion pair $(\sF,\sC)$, where $\sF$ is the class of all strongly
contraflat left $\E$\+contramodules and $\sC=\E\contra^\Bcot$.
 By Lemma~\ref{contra-cotorsion-contraflat} and
Corollary~\ref{strongly-contraflat-cor} (see also
Lemma~\ref{cotorsion-pair-direct-summands-lemma}), \,$(\sF,\sC)$
is a hereditary complete cotorsion pair in~$\sE$.
 Since all left $B$\+modules have finite cotorsion dimensions, it
follows that all $\E$\+contramodules have finite $\sC$\+coresolution
dimensions; so
Proposition~\ref{becker-contraderived-finite-coresolutions}
is applicable.
 
 Part~(b) is a restatement of
Theorem~\ref{comods-contramods-positselski-co-contra-derived}(b);
see also Corollary~\ref{finite-homol-dim-equivalence-cor}.
 Part~(c) can be deduced from part~(b) together with the fact that
strongly contraflat left $\E$\+contramodules have finite projective
dimensions in $\E\contra$ and Proposition~\ref{finite-resolutions}.
 Alternatively, notice that countable products of strongly contraflat
$\E$\+contramodules have finite projective dimensions
by (the proof of) Corollary~\ref{strongly-contraflat-cor}(b),
so Corollary~\ref{finite-homol-dim-equivalence-cor} applies.
\end{proof}

\begin{cor}
 In the situation of
Corollary~\textup{\ref{noncomm-flat-dualizing-co-contra-cor}},
assume additionally that the coring\/ $\E$ is a projective
left $B$\+module.
 Then the datum of a dualizing complex\/ $\D^\bu\rarrow D^\bu$ for
the corings\/ $\C$ and\/ $\E$ induces an equivalence of
triangulated categories\/ $\sD^\co(\C\comodl)\simeq
\sD^\ctr(\E\contra)$, which is provided by the derived functors\/
$\boR\Hom_\C(\D^\bu,{-})$ and\/ $\D^\bu\ocn_\E^\boL\nobreak{-}$. 
\end{cor}

\begin{proof}
 In addition to what has been said in
Corollaries~\ref{christensen-frankild-holm},
\ref{noncomm-flat-dualizing-co-contra-cor} and
Theorem~\ref{projective-coring-contra-resolutions},
we point out that in our present assumptions the assertions of
Lemmas~\ref{co-contra-inj-to-proj-lemma}(a)
and~\ref{coring-dualizing-lemma}(a) apply to any strongly
contraflat left $\E$\+contramodule~$\gF$.
 Indeed, in view of Corollary~\ref{strongly-contraflat-cor}(b)
one only has to check that the functor $\N\ocn_\E{-}$ takes
short exact sequences of strongly contraflat left
$\E$\+contramodules to short exact sequences of abelian groups
for any right $\E$\+comodule~$\N$.
 This follows from part~(a) of the same Corollary, as
$\Hom_\boZ(\N\ocn_\E\gF\;\boQ/\boZ)\simeq
\Hom^\E(\gF,\Hom_\boZ(\N,\boQ/\boZ))$ and the left
$\E$\+contramodule $\Hom_\boZ(\N,\boQ/\boZ)$ is $B$\+cotorsion
by Lemma~\ref{cotors-hom-nc}(b).
 Therefore, one can construct the derived functor
$\D^\bu\ocn_\E^\boL{-}$ using strongly contraflat resolutions.
\end{proof}

\subsection{Base ring change}
 Let $\C$ be a coring over an associative ring~$A$.
 Given a right $\C$\+comodule $\N$ and a left $\C$\+comodule $\M$,
their \emph{cotensor product} $\N\oc_\C\M$ is an abelian group
constructed as the kernel of the natural pair of maps
$\N\ot_A\M\birarrow\N\ot_A\C\ot_A\M$.
 Given a left $\C$\+comodule $\M$ and a left $\C$\+contramodule $\P$,
their group of \emph{cohomomorphisms} $\cohom_\C(\M,\P)$ is
constructed as the cokernel of the natural pair of maps
$\Hom_A(\C\ot_A\M\;\P)\birarrow\Hom_A(\M,\P)$
\,\cite[Sections~1.2.1 and~3.2.1]{Psemi},
\cite[Sections~2.5\+-2.6]{Prev}.

 Let $\C$ be a coring over an associative ring $A$ and $\E$ be
a coring over an associative ring~$B$.
 A \emph{map of corings\/ $\C\rarrow\E$ compatible with
a ring homomorphism $A\rarrow B$} is an $A$\+$A$\+bimodule morphism
such that the maps $A\rarrow B$, \ $\C\rarrow\E$, and the induced
map $\C\ot_A\C\rarrow\E\ot_B\E$ form commutative diagrams with
the comultiplication and counit maps in $\C$ and~$\E$
\,\cite[Section~7.1.1]{Psemi}. 

 Let $\C\rarrow\E$ be a map of corings compatible with a ring
map $A\rarrow B$.
 Given a left $\C$\+comodule $\M$, one defines a left
$\E$\+comodule ${}_B\M$ by the rule ${}_B\M=B\ot_A\M$,
the coaction map being constructed as the composition
$B\ot_A\M\rarrow B\ot_A\C\ot_A\M\rarrow B\ot_A\E\ot_A\M
\rarrow\E\ot_A\M\simeq\E\ot_B(B\ot_A\M)$.
 Similarly, given a right $\C$\+comodule $\K$, there is 
a natural right $\E$\+comodule structure on the tensor product
$\K_B=\K\ot_A B$.
 Given a left $\C$\+contramodule $\P$, one defines a left
$\E$\+contramodule ${}^B\P$ by the rule ${}^B\P=\Hom_A(B,\P)$,
the contraaction map being constructed as the composition
$\Hom_B(\E,\Hom_A(B,\P))\simeq\Hom_A(\E,\P)\rarrow\Hom_A(\E\ot_A B\;\P)
\rarrow\Hom_A(\C\ot_A B\;\P)\simeq\Hom_A(B,\Hom_A(\C,\P))\rarrow
\Hom_A(B,\P)$.

 Assuming that $\C$ is a flat right $A$\+module, the functor
$\M\mpsto{}_B\M\:\C\comodl\rarrow\E\comodl$ has a right adjoint
functor, which is denoted by $\N\mpsto{}_\C\N\:\E\comodl\rarrow
\C\comodl$ and constructed by the rule ${}_\C\N=\C_B\oc_\E\M$.
 In particular, the functor $\N\mpsto{}_\C\N$ takes a coinduced
left $\E$\+comodule $\E\ot_B U$ into the coinduced left
$\C$\+comodule $\C\ot_A U$ (where $U$ is an arbitrary left
$B$\+module) \cite[Section~7.1.2]{Psemi}.

 Given a coring $\C$ over an associative ring $A$ and an associative
ring homomorphism $A\rarrow B$, one can define a coring ${}_B\C_B$
over the ring $B$ by the rule ${}_B\C_B=B\ot_A\C\ot_A B$.
 The counit in ${}_B\C_B$ is constructed as the composition
$B\ot_A \C\ot_A B\rarrow B\ot_A A\ot_A B\rarrow B$, and
the comultiplication is provided by the composition
$B\ot_A\C\ot_A B\rarrow B\ot_A\C\ot_A\C\ot_A B\simeq
B\ot_A\C\ot_A A\ot_A\C\ot_A B\rarrow B\ot_A\C\ot_A B\ot_A\C\ot_A B
\simeq(B\ot_A\C\ot_A B)\ot_B(B\ot_A\C\ot_A B)$.
 There is a natural map of corings $\C\rarrow {}_B\C_B$ compatible
with the ring map $A\rarrow B$.

 Assuming that $\C$ is a flat right $A$\+module and $B$ is
a faithfully flat right $A$\+module, the functors $\M\rarrow {}_B\M$
and $\N\mpsto{}_\C\N$ are mutually inverse equivalences between
the abelian categories of left $\C$\+comodules and left
${}_B\C_B$\+comodules.
 One proves this by applying the Barr--Beck monadicity theorem to
the conservative exact functor $\C\comodl\rarrow B\modl$ taking
a left $\C$\+comodule $\M$ to the left $B$\+module $B\ot_A\M$
(see~\cite[Section~2.1.3]{KR} or~\cite[Section~7.4.1]{Psemi}).

 Now let us assume that $\C$ is a flat left $A$\+module, $\E$ is
a flat left $B$\+module, and $B$ is a flat left $A$\+module.
 Then the functor $\P\mpsto {}^B\P$ takes $A$\+cotorsion left
$\C$\+contramodules to $B$\+cotorsion left $\E$\+contramodules
(by Lemma~\ref{cotors-hom-nc}(a)) and induces an exact functor
between the exact categories $\C\contra^\Acot\rarrow
\E\contra^\Bcot$.

 The left adjoint functor $\Q\mpsto {}^\C\Q$ to this exact functor
may not be everywhere defined, but one can easily see that it is
defined on the full subcategory of left $\E$\+contramodules
induced from cotorsion left $B$\+modules in the exact category
$\E\contra^\Bcot$, taking an induced contramodule $\Q=\Hom_B(\E,V)$
over $\E$ to the induced contramodule ${}^\C\Q=\Hom_A(\C,V)$
over~$\C$ (where $V$ is an arbitrary cotorsion left $B$\+module;
see Lemma~\ref{cotors-restrict}(a)).
 These words mean that the adjunction isomorphism $\Hom^\E(\Q\;{}^B\P)
\simeq\Hom^\C({}^\C\Q\;\P)$ holds for every $A$\+cotorsion left
$\C$\+contramodule $\P$ and every left $\E$\+contramodule $\Q$
induced from a cotorsion left $B$\+module.
 In such a context, we will say that the functor $\Q\mpsto {}^\C\Q$
(defined on a full subcategory of $\E\contra^\Bcot$) is a ``partial
left adjoint'' to the functor $\P\mpsto {}^B\P$.

 Finally, we assume that $\C$ is flat left $A$\+module and $B$ is
a faithfully flat left $A$\+module.
 Here a flat left $A$\+module $G$ is said to be \emph{faithfully flat}
if the abelian group $N\ot_AG$ is nonzero for any nonzero right
$A$\+module~$N$ (in other words, the exact functor ${-}\ot_AG\:\modr A
\rarrow\boZ\modl$ is faithful).

 It is clear from Example~\ref{loc-contraherent-counterexample}
(in light of the interpretation of quasi-coherent sheaves on schemes
as comodules~\cite[Section~2]{KR}, \cite[Example~2.5]{Pflcc} and
the dual-analogous point of view on contraherent cosheaves)
that the functor $\P\mpsto{}^B\P\:\C\contra^\Acot\allowbreak
\rarrow{}_B\C_B\contra^\Bcot$ is \emph{not} an equivalence of
exact categories in general (and not even in the case when $\C=A$).
 The aim of this section is to prove the following slightly weaker
result (cf.\ Corollaries~\ref{ctrh-lcth-cor}\+-\ref{lin-ctrh-lcth-cor}
and~\ref{finite-krull-ctrh-lcth-derived}).

\begin{thm}  \label{change-of-covering-contraderived-equivalence}
 Let\/ $\C$ be a coassociative coring over an associative ring $A$
and $A\rarrow B$ be an associative ring homomorphism.
 Assume that\/ $\C$ is a flat left $A$\+module and $B$ is a faithfully
flat left $A$\+module.
 Then the triangulated functor\/ $\sD^\ctr(\C\contra^\Acot)\allowbreak
\rarrow\sD^\ctr({}_B\C_B\contra^\Bcot)$ induced by the exact functor\/
$\P\mpsto{}^B\P$ is an equivalence of triangulated categories.
\end{thm}

\begin{cor}
 In the assumptions of the previous Theorem, the exact functor\/
$\P\mpsto{}^B\P\:\C\contra^\Acot\rarrow {}_B\C_B\contra^\Bcot$
is fully faithful and induces isomorphisms between the Ext
groups computed in the exact categories\/ $\C\contra^\Acot$
and\/ ${}_B\C_B\contra^\Bcot$.
 Any short sequence (or, more generally, bounded above complex)
in\/ $\C\contra^\Acot$ which this functor transforms into
an exact sequence in\/ ${}_B\C_B\contra^\Bcot$ is exact in\/
$\C\contra^\Acot$.
\end{cor}

\begin{proof}
 See~\cite[Section~4.1]{Psemi}
or Lemma~\ref{co-contra-bounded-fully-faithful}(b).
\end{proof}

 The proof of Theorem~\ref{change-of-covering-contraderived-equivalence}
is based on the following technical result, which provides adjusted
resolutions for the construction of the left derived functor of
the partial left adjoint functor $\Q\mpsto{}^\C\Q$ to the exact functor
$\P\mpsto{}^B\P$.
 For completeness, we formulate three versions of the (co)induced
resolution theorem; it is the third one that will be used
in the argument below.

\begin{thm}  \label{co-induced-resolutions}
\textup{(a)} Let\/ $\C$ be a coring over an associative ring $A$;
assume that\/ $\C$ is a flat right $A$\+module.
 Then the coderived category\/ $\sD^\co(\C\comodl)$ of the abelian
category of left\/ $\C$\+comodules is equivalent to the quotient
category of the homotopy category of complexes of coinduced\/
$\C$\+comodules by its minimal triangulated subcategory containing
the total complexes of short exact sequences of complexes of\/
$\C$\+comodules termwise coinduced from short exact sequences of
$A$\+modules and closed under infinite direct sums. \par
\textup{(b)} Let\/ $\C$ be a coring over an associative ring $A$;
assume that\/ $\C$ is a projective left $A$\+module.
 Then the contraderived category\/ $\sD^\ctr(\C\contra)$ of the abelian
category of left\/ $\C$\+contramodules is equivalent to the quotient
category of the homotopy category of complexes of induced\/
$\C$\+contramodules by its minimal triangulated subcategory containing
the total complexes of short exact sequences of complexes of\/
$\C$\+contramodules termwise induced from short exact sequences of
$A$\+modules and closed under infinite products. \par
\textup{(c)} Let\/ $\C$ be a coring over an associative ring $A$;
assume that\/ $\C$ is a flat left $A$\+module.
 Then the contraderived category\/ $\sD^\ctr(\C\contra^\Acot)$ of
the exact category of $A$\+cotorsion left\/ $\C$\+contramodules is
equivalent to the quotient category of the homotopy category of
complexes of\/ $\C$\+contramodules termwise induced from cotorsion
$A$\+modules by its minimal triangulated subcategory containing
the total complexes of short exact sequences of complexes of\/
$\C$\+contramodules termwise induced from short exact sequences of
cotorsion $A$\+modules and closed under infinite products. \par
\end{thm}

\begin{proof}
 This is yet another version of the results of~\cite[Proposition~1.5
and Remark~1.5]{EP}, \cite[Theorem~4.2.1]{Pweak}, \cite[Theorems~7.9
and~7.11]{Pedg}, the above Propositions~\ref{fully-faithful-prop}
and~\ref{infinite-resolutions}(b), etc.
 The difference is that the categories of coinduced comodules or
induced contramodules are not even exact.
 This does not change much, however.
 Let us sketch a proof of part~(c); the proofs of parts~(a\+b)
are similar (but simpler).
 
 The contramodule relative bar-resolution, totalized by taking
infinite products along the diagonals, provides a closed morphism
with a contraacyclic cone onto any complex of $A$\+cotorsion left
$\C$\+contramodules from a complex of $\C$\+contramodules termwise
induced from cotorsion $A$\+modules.
 This proves that the natural functor from the homotopy category
of complexes of left $\C$\+contramodules termwise induced from
cotorsion $A$\+modules to the contraderived category of
$A$\+cotorsion left $\C$\+contramodules is a Verdier localization
functor (see Lemma~\ref{pkoszul-lemma16}).

 In order to prove that the natural functor from the quotient
category of the homotopy category of termwise induced complexes
in the formulation of the theorem to the contraderived category
$\sD^\ctr(\C\contra^\Acot)$ is fully faithful, one shows
that any morphism in $\Hot(\C\contra^\Acot)$ from a complex of
contramodules termwise induced from cotorsion $A$\+modules to
a contraacyclic complex over $\C\contra^\Acot$ factorizes through
an object of the minimal triangulated subcategory containing
the total complexes of the short exact sequences of complexes of
$\C$\+contramodules termwise induced from short exact sequences of
cotorsion $A$\+modules and closed under infinite products.

 Let us only explain the key step, as the argument is mostly
similar to the ones spelled out in the references above.
 Let $\gM^\bu$ be the total complex of a short exact sequence
$\gU^\bu\rarrow\gV^\bu\rarrow\gW^\bu$ of complexes of $A$\+cotorsion
$\C$\+contramodules, and let $\gE^\bu$ be a complex of
$\C$\+contramodules termwise induced from cotorsion $A$\+modules.
 The construction dual to the one explained in~\cite[Theorem~5.5]{Psemi}
provides a closed morphism onto the complex $\gE^\bu$ from a complex of
$\C$\+contramodules $\gF^\bu$ termwise induced from flat cotorsion
$A$\+modules such that this morphism is induced from an admissible
epimorphism of cotorsion $A$\+modules at every term of the complexes
(cf.\ the proof of Theorem~\ref{gorenstein-base-co-contra-derived}(b)).

 Let $\R^\bu$ denote the kernel of the morphism of complexes
$\gF^\bu\rarrow\gE^\bu$.
 Then the cone $\P^\bu$ of the admissible monomorphism $\R^\bu
\rarrow\gF^\bu$ maps naturally onto $\gE^\bu$ with the cone
isomorphic to the total complex of a short exact sequence of
complexes of $\C$\+contramodules termwise induced from a short
exact sequence of cotorsion $A$\+modules.
 As a morphism of graded $\C$\+contramodules, the composition
$\P^\bu\rarrow\gE^\bu$ factorizes through the graded
$\C$\+contramodule $\gF^\bu$.
 The latter being a projective graded object of $\C\contra^\Acot$,
it follows that the composition $\P^\bu\rarrow\gE^\bu\rarrow\gM^\bu$
is homotopic to zero by Lemma~\ref{lifted-contractible}.
\end{proof}

\begin{proof}[Proof of
Theorem~\ref{change-of-covering-contraderived-equivalence}]
 First of all, we notice that the functor $\Q\mpsto {}^\C\Q$ acting
between the full subcategories of left ${}_B\C_B$\+contramodules
induced from cotorsion $B$\+modules and left $\C$\+contramodules
induced from cotorsion $A$\+modules in ${}_B\C_B\contra^\Bcot$
and $\C\contra^\Acot$ is fully faithful.
 Indeed, one has
\begin{gather*}
\Hom^{{}_B\C_B}(\Hom_B({}_B\C_B,U)\;\Hom_B({}_B\C_B,V))
\simeq\Hom_B(U,\Hom_B({}_B\C_B,V)) \\ \simeq
\Hom_B({}_B\C_B\ot_B U\;V)\simeq\Hom_A(\C\ot_A U\;V)\simeq
\Hom^\C(\Hom_A(\C,U),\Hom_A(\C,V)).
\end{gather*}
 Furthermore, the composition of functors $\Q\mpsto{}^\C\Q\mpsto
{}^B({}^\C\Q)$ is naturally isomorphic to the identity endofunctor
on the full subcategory of ${}_B\C_B$\+contramodules induced from
cotorsion $B$\+modules in ${}_B\C_B\contra^\Bcot$, as
$\Hom_A(B,\Hom_A(\C,V))\simeq\Hom_A(\C\ot_A B\;V)\simeq
\Hom_B(B\ot_A\C\ot_A B\;V)$.

 Clearly, the functor $\Q^\bu\mpsto{}^\C\Q^\bu$, acting between
the homotopy categories of complexes of ${}_B\C_B$\+contramodules
termwise induced from complexes of cotorsion $B$\+modules and
complexes of $\C$\+contramodules termwise induced from cotorsion
$A$\+modules, is ``partially left adjoint'' to the functor
$\P^\bu\mpsto{}^B\P^\bu\:\Hot(\C\contra^\Acot)\allowbreak\rarrow
\Hot({}_B\C_B\contra^\Bcot)$.
 By Theorem~\ref{co-induced-resolutions}(c), the former functor
can be used to construct a left derived functor
$\Q^\bu\mpsto\nobreak{}^\C_\boL\Q^\bu\:\sD^\ctr({}_B\C_B\contra^\Bcot)
\rarrow\sD^\ctr(\C\contra^\Acot)$.
 It is easy to see that the functor
$\Q^\bu\mpsto{}^\C_\boL\Q^\bu$ is left adjoint to the functor
$\P^\bu\mpsto\nobreak{}^B\P^\bu\:\sD^\ctr(\C\contra^\Acot)
\rarrow\sD^\ctr({}_B\C_B\contra^\Bcot)$.

 According to the above, the composition of the adjoint triangulated
functors $\Q^\bu\mpsto{}^B({}^\C_\boL\Q^\bu)\:\sD^\ctr
({}_B\C_B\contra^\Bcot)\rarrow\sD^\ctr({}_B\C_B\contra^\Bcot)$
is isomorphic to the identity functor.
 Hence the functor $\Q^\bu\mpsto{}^\C_\boL\Q^\bu\:
\sD^\ctr({}_B\C_B\contra^\Bcot)\rarrow\sD^\ctr(\C\contra^\Acot)$
is fully faithful, while the functor $\P^\bu\mpsto{}^B\P^\bu\:
\sD^\ctr(\C\contra^\Acot)\allowbreak\rarrow
\sD^\ctr({}_B\C_B\contra^\Bcot)$ is a Verdier localization functor.
{\hfuzz=1.2pt\par}

 It remains to show that the triangulated functor
$\Q^\bu\mpsto{}^\C_\boL\Q^\bu$ is essentially surjective.
 It is straightforward from the constructions that both functors
$\P^\bu\mpsto{}^B\P^\bu$ and $\Q^\bu\mpsto{}^\C_\boL\Q^\bu$
preserve infinite products.
 So it suffices to prove that the contraderived category
$\sD^\ctr(\C\contra^\Acot)$ coincides with its minimal triangulated
subcategory containing all complexes of $\C$\+contramodules
termwise induced from left $A$\+modules obtained by restriction of
scalars from cotorsion left $B$\+modules and closed under
infinite products.
 In view of the argument of
Lemmas~\ref{second-kind-complex-resolution}\+-\ref{telescope}
and the final paragraph of the proof of
Proposition~\ref{infinite-resolutions}(b),
we only need to check that there is an admissible epimorphism
onto any object of $\C\contra^\Acot$ from an induced
left $\C$\+contramodule of the desired type.

 The contraaction morphism $\Hom_A(\C,\P)\rarrow\P$ being an admissible
epimorphism in $\C\contra^\Acot$ for any $A$\+cotorsion left
$\C$\+contramodule $\P$, the problem reduces to constructing
an admissible epimorphism in $A\modl^\cot$ onto any cotorsion
left $A$\+module from an $A$\+module obtained by restriction of
scalars from a cotorsion left $B$\+module.
 Here we are finally using the assumption that $B$ is
a \emph{faithfully} flat left $A$\+module.
 According to~\cite[Section~7.4.1]{Psemi}, the latter is equivalent to
the ring homomorphism $A\rarrow B$ being injective and its cokernel
being a flat left $A$\+module.
 Assuming these conditions, the morphism $\Hom_A(B,P)\rarrow
\Hom_A(A,P)=P$ induced by the map $A\rarrow B$ is an admissible 
epimorphism in $A\modl^\cot$ for any cotorsion left $A$\+module $P$
by Lemma~\ref{cotors-hom-nc}(a).
 It remains to say that the left $B$\+module $\Hom_A(B,P)$ is
cotorsion by Lemma~\ref{cotors-coexten}(a). 
\end{proof}

\Section{Affine Noetherian Formal Schemes}

 This appendix is a continuation of~\cite[Appendix~B]{Pweak}.
 Its goal is to construct the derived co-contra correspondence
between quasi-coherent torsion sheaves and contraherent cosheaves
of contramodules on an affine Noetherian formal scheme with
a dualizing complex.
 The affineness condition allows to speak of (co)sheaves in terms of
(contra)modules over a ring, simplifying the task.
 For a general discussion of the philosophy of co-contra correspondence
in the context of a commutative ring with a finitely generated ideal,
see~\cite[Introduction and Remark~4.10]{Pmgm}.

\subsection{Torsion modules and contramodules}  \label{torsion-subsect}
 Let $R$ be a Noetherian commutative ring and $I\sub R$ be an ideal.
 Given an $R$\+module $E$, we denote by ${}_{(n)}E\sub E$ its
$R$\+submodule consisting of all the elements annihilated by~$I^n$.
 An $R$\+module $\M$ is said to be \emph{$I$\+torsion} if
$\M=\bigcup_{n\ge1}\,{}_{(n)}\M$, i.~e., for any element $m\in \M$
there exists an integer $n\ge1$ such that $I^nm=0$ in~$\M$.
 An $R$\+module $\P$ is called an \emph{$(R,I)$\+contramodule} if
$\Ext^*_R(R[s^{-1}],\P)=0$ for all $s\in I$, or equivalently,
the system of equations $q_n=p_n+sq_{n+1}$, \ $n\ge0$, is uniquely
solvable in $q_n\in \P$ for any fixed sequence $p_n\in \P$
\,\cite[Lemma~2.1(b)]{Pcta} (cf.\ the beginning of
Section~\ref{very-eklof-trlifaj-subsect}).

 Just as the category of $I$\+torsion $R$\+modules, the category
of $(R,I)$\+contramodules is abelian~\cite[Proposition~1.1]{GL},
\cite[Theorem~1.2(a)]{Pcta} and depends only on the $I$\+adic
completion of the ring~$R$ \,\cite[Theorem~B.1.1]{Pweak},
\cite[Proposition~1.5 and Corollary~3.7]{Pdc}
(cf.\ the ending part of Section~\ref{cotorsion-modules} above).
 Other relevant references include the papers and
preprint~\cite{Pmgm,Pcta,Pper}.
 We denote the category of $I$\+torsion $R$\+modules by
$(R,I)\tors$ and the category of $(R,I)$\+contramodules by
$(R,I)\contra$.
 Both categories are full subcategories in the abelian category of
$R$\+modules $R\modl$, closed under the kernels and cokernels;
the former subcategory is also closed under infinite direct sums,
while the latter one is closed under infinite products.

 Given a sequence of elements $r_n\in R$ converging to zero in
the $I$\+adic topology of $R$ and a sequence of elements $p_n\in \P$
of an $(R,I)$\+contramodule $\P$, the result of the ``infinite
summation operation'' $\sum_n r_n p_n$ is well-defined as
an element of~$\P$ \cite[Section~1.2]{Pweak}, \cite[Proposition~1.5
and Corollary~3.7]{Pdc}.
 One can define the \emph{contratensor product} $\P\ocn_{(R,I)} \M$ of
an $(R,I)$\+contramodule $\P$ and an $I$\+torsion $R$\+module $\M$
as the quotient $R$\+module of the tensor product $\P\ot_\boZ \M$ by
the relations $(\sum_n r_n p_n)\ot m=\sum_n p_n\ot r_n m$ for any
sequence of elements $r_n\in R$ converging to zero in $R$, any sequence
of elements $p_n\in \P$, and any element $m\in \M$.
 Here all but a finite number of summands in the right-hand side
vanish due to the conditions imposed on the sequence $r_n$
and the module~$\M$.

 In fact, however, for any $(R,I)$\+contramodule $\P$ and
$I$\+torsion $R$\+module $\M$ there is a natural isomorphism
$\P\ocn_{(R,I)} \M\simeq \P\ot_R \M$.
 This can be deduced from the fact that the functor
$(R,I)\contra\rarrow R\modl$ is fully faithful, or alternatively,
follows directly from the observation that, for any fixed integer
$n\ge0$, any sequence of elements of the ideal $I^n\sub R$ converging
to zero in the $I$\+adic topology of $R$ is a linear combination of
a finite number of sequences converging to zero in $R$ with
the coefficients belonging to~$I^n$
(cf.~\cite[proof of Proposition~B.9.1]{Pweak}).

 Clearly, the tensor product $E\ot_R \M$ of any $R$\+module $E$ and
any $I$\+torsion $R$\+module $\M$ is an $I$\+torsion $R$\+module.
 Similarly, the $R$\+module $\Hom_R(\M,E)$ of $R$\+linear maps
from an $I$\+torsion $R$\+module $\M$ into any $R$\+module $E$ is
an $(R,I)$\+contramodule.
 Finally, the $R$\+module $\Hom_R(E,\P)$ of homomorphisms from
any $R$\+module $E$ into an $(R,I)$\+contramodule $\P$ is
an $(R,I)$\+contramodule~\cite[Sections~1.5 and~B.2]{Pweak},
\cite[Lemma~6.1]{Pcta}.

 An $(R,I)$\+contramodule $\gF$ is said to be \emph{$(R,I)$\+contraflat}
if the functor $\M\mpsto \gF\ocn_{(R,I)}\M$ is exact on the abelian
category of $I$\+torsion $R$\+modules.
 According to the above, this is equivalent to the functor
$\M\mpsto \gF\ot_R \M$ being exact on the full abelian subcategory
$(R,I)\tors\sub R\modl$.

 By the Artin--Rees lemma, an object $\K\in (R,I)\tors$ is injective
if and only if it is injective in $R\modl$, and if and only if
its submodules ${}_{(n)}\K$ of elements annihilated by $I^n$ are
injective $R/I^n$\+modules for all $n\ge1$.
 By~\cite[Lemma~B.9.2]{Pweak} or~\cite[Corollary~10.3(a)]{Pcta},
an $(R,I)$\+contramodule $\gF$ is $(R,I)$\+contraflat if and only if
it is a flat $R$\+module, and if and only if its reductions $\gF/I^n\gF$
are flat $R/I^n$\+modules for all $n\ge1$.
 Furthermore, an object $\gF\in (R,I)\contra$ is projective
if and only if it is $(R,I)$\+contraflat and its reduction $\gF/I\gF$
is a projective $R/I$\+module, and if and only if all the reductions
$\gF/I^n\gF$ are projective
$R/I^n$\+modules~\cite[Corollary~B.8.2]{Pweak},
\cite[Theorem~10.5]{Pcta}.

\begin{thm}  \label{tors-contra-co-contra-derived}
\textup{(a)} The classes of Positselski-coacyclic and Becker-coacyclic
complexes in the abelian category $(R,I)\tors$ coincide.
 The coderived category\/ $\sD^{\co=\bco}((R,I)\tors)$ of the abelian
category of $I$\+torsion $R$\+modules is equivalent to the homotopy
category of complexes of injective $I$\+torsion $R$\+modules. \par
\textup{(b)} The contraderived category\/ $\sD^\ctr((R,I)\contra)$
of the abelian category of $(R,I)$\+contramodules is equivalent to
the contraderived category of the exact category of $R$\+flat
(i.~e., $(R,I)$\+contraflat) $(R,I)$\+contramodules. \par
\textup{(c)} The Becker contraderived category\/
$\sD^\bctr((R,I)\contra)$ of the abelian category of
$(R,I)$\+contramodules is equivalent to the Becker contraderived
category of the exact category of $R$\+flat $(R,I)$\+contramodules. \par
\textup{(d)} The Becker contraderived category\/
$\sD^\bctr((R,I)\contra)$ is equivalent to the homotopy category of
complexes of projective $(R,I)$\+contramodules. \par
\textup{(e)} Assume that the Noetherian ring $R/I$ has finite
Krull dimension.
 Then the classes of Positselski-contraacyclic and Becker-contraacyclic
complexes in the abelian category $(R,I)\contra$ coincide.
 The contraderived category\/ $\sD^{\ctr=\bctr}((R,I)\contra)$
is equivalent to the absolute derived category of the exact
category of $R$\+flat $(R,I)$\+contramodules. \hbadness=1075\par
\textup{(f)} Assume that the Noetherian ring $R/I$ has finite
Krull dimension.
 Then the contraderived category\/ $\sD^{\ctr=\bctr}((R,I)\contra)$
is equivalent to the homotopy category of complexes of projective
$(R,I)$\+contramodules.
\end{thm}

\begin{proof}
 Part~(a) holds by
Theorem~\ref{positselski-becker-co-contra-derived}(a), since there are
enough injectives in the abelian category $(R,I)\tors$ and the class
of injective objects is closed under infinite direct sums.
 Part~(b) is true, because there are enough $(R,I)$\+contraflat
(and even projective) objects in $(R,I)\contra$ and the class of
$(R,I)$\+contraflat $(R,I)$\+contramodules (or even flat
$R$\+modules) is closed under infinite products
(see Proposition~\ref{infinite-resolutions}(b)).
 Similarly, part~(c) holds by 
Proposition~\ref{becker-contraderived-infinite-resolutions}.
 Part~(d) is a particular case of
Theorem~\ref{coderived-of-grothendieck-contraderived-of-lpacepo}(b).

 Parts~(e\+f) hold, since in their assumptions any $(R,I)$\+contraflat
$(R,I)$\+contra\-module $\gF$ has finite projective dimension in
$(R,I)\contra$.
 Indeed, consider a projective resolution $\P_\bu$ of $\gF$ in
$(R,I)\contra$ and apply the functor $R/I\ot_R{-}$ to it.
 Being a bounded above exact complex of flat $R$\+modules,
the complex $\P_\bu\rarrow\gF$ will remain exact after taking
the tensor product.
 By Theorem~\ref{raynaud-gruson-flat-thm}, the $R/I$\+modules
of cycles in the complex $R/I\ot_R\P_\bu$ are eventually projective,
and it follows that the $(R,I)$\+contramodules of cycles in
$\P_\bu$ are eventually projective, too.
 Now it remains to apply
Corollary~\ref{finite-homol-dim-equivalence-cor},
Proposition~\ref{finite-resolutions}, and/or
Lemma~\ref{psemi-remark21} for Positselski's contraderived categories;
or Theorems~\ref{finite-homol-dim-becker-co-contra-derived}(b) and
\ref{positselski-becker-co-contra-derived}(b) for Becker's ones.
\end{proof}

 For noncommutative generalizations of the next two lemmas
(as well as the theorem following them),
see Sections~\ref{noncomm-noetherian-subsect}
and~\ref{ind-affine-co-contra-subsect}.

\begin{lem}  \label{hom-tensor-reduction}
\textup{(a)} For any $R$\+module $M$, injective $R$\+module $J$,
and\/ $n\ge1$, there is a natural isomorphism of $R/I^n$\+modules\/
$\Hom_R(M,J)/I^n\Hom_R(M,J)\simeq\Hom_{R/I^n}({}_{(n)}M\;{}_{(n)}J)$.
\par
\textup{(b)} For any $R$\+module $M$, flat $R$\+module $F$, and\/
$n\ge1$, there is a natural isomorphism of $R/I^n$\+modules\/
${}_{(n)}(M\ot_R F)\simeq\,{}_{(n)}M\ot_{R/I^n} F/I^n F$.
\end{lem}

\begin{proof}
 For any finitely generated $R$\+module $E$ there are natural
isomorphisms $\Hom_R(\Hom_R(E,M),J)\simeq E\ot_R\Hom_R(M,J)$
and $\Hom_R(E\;M\ot_R F)\simeq\Hom_R(E,M)\allowbreak\ot_R F$.
 It remains to take $E=R/I^n$.
\end{proof}

 A finite complex $\D^\bu$ of injective $I$\+torsion $R$\+modules is
said to be a \emph{dualizing complex} for the pair $(R,I)$ if for
every integer $n\ge1$ the complex of $R/I^n$\+modules ${}_{(n)}\D^\bu$
is a dualizing complex for the commutative Noetherian ring $R/I^n$.
 If a finite complex of injective $R$\+modules $D^\bu$ is
a dualizing complex for the ring $R$, then its subcomplex
$\D^\bu=\bigcup_{n\ge0}\,{}_{(n)}D^\bu$ consisting of all the elements
annihilated by some powers of $I\sub R$ is a dualizing complex
for $(R,I)$.

\begin{lem}  \label{dualizing-reduction}
 Let\/ $\D^\bu$ be a finite complex of injective $I$\+torsion
$R$\+modules such that the complex of $R/I$\+modules\/
${}_{(1)}\D^\bu$ is a dualizing complex for the ring~$R/I$.
 Then\/ $\D^\bu$ is a dualizing complex for the pair $(R,I)$.
\end{lem}

\begin{proof}
 Clearly, ${}_{(n)}\D^\bu$ is a finite complex of injective
$R/I^n$\+modules.
 Since any $R/I^n$\+module $K$ with finitely generated
$R/I$\+module ${}_{(1)}K$ is finitely generated itself,
and the functors $\Ext^i_{R/I^n}(R/I,{-})$ take finitely
generated $R/I^n$\+modules to finitely generated $R/I$\+modules,
one can check by increasing induction on the cohomological degree
that the $R/I^n$\+modules of cohomology of the complex ${}_{(n)}\D^\bu$
are finitely generated.

 It remains to show that the natural map $R/I^n\rarrow
\Hom_{R/I^n}({}_{(n)}\D^\bu\;{}_{(n)}\D^\bu)$ is a quasi-isomorphism.
 Both the left-hand and the right-hand sides are finite complexes
of flat $R/I^n$\+modules.
 By Lemma~\ref{hom-tensor-reduction}(a), the functor $P\mpsto P/IP$
transforms this morphism into the morphism $R/I\rarrow
\Hom_{R/I}({}_{(1)}\D^\bu\;{}_{(1)}\D^\bu)$, which is
a quasi-isomorphism by assumption.
 Hence the desired assertion follows by Nakayama's lemma.
\end{proof}

\begin{thm}  \label{tors-contra-correspondence}
 The datum of a dualizing complex\/ $\D^\bu$ for a pair $(R,I)$ with
a Noetherian commutative ring $R$ and an ideal $I\sub R$ induces
an equivalence of triangulated categories\/ $\sD^\co((R,I)\tors)
\simeq\sD^\ctr((R,I)\contra)$, which is provided by the derived
functors\/ $\boR\Hom_R(\D^\bu,{-})$ and\/ $\D^\bu\ot_R^\boL{-}$.
\end{thm}

\begin{proof}
 Notice first of all that $\sD^\co((R,I)\tors)=\sD^\bco((R,I)\tors)$
by Theorem~\ref{tors-contra-co-contra-derived}(a) and
$\sD^\ctr((R,I)\contra)=\sD^\bctr((R,I)\contra)$ by
Theorem~\ref{tors-contra-co-contra-derived}(e). {\hbadness=1450\par}

 The constructions of the derived functors are also based on
Theorem~\ref{tors-contra-co-contra-derived}(a,e).
 Applying the functor $\Hom_R(\D^\bu,{-})$ to a complex of
injective $I$\+torsion $R$\+modules produces a complex of
$R$\+flat $(R,I)$\+contramodules.
 Applying the functor $\D^\bu\ot_R{-}$ to a complex of $R$\+flat
$(R,I)$\+contramodules produces a complex of injective $I$\+torsion
$R$\+modules, and for any absolutely acyclic complex of flat
$R$\+modules (or an acyclic complex of flat $R$\+modules with
flat modules of cocycles) $F^\bu$, the complex of injective
($I$\+torsion) $R$\+modules $\D^\bu\ot_R F^\bu$ is contractible.

 It remains to show that the morphism of finite complexes
$\D^\bu\ot_R\Hom_R(\D^\bu,\J)\rarrow\J$ is a quasi-isomorphism
for any $I$\+torsion $R$\+module $\J$, and the morphism of
finite complexes $\gF\rarrow\Hom_R(\D^\bu\;\D^\bu\ot_R\gF)$
is a quasi-isomorphism for any $(R,I)$\+contramodule~$\gF$.
 Both assertions follow from Lemma~\ref{hom-tensor-reduction},
which implies that the morphisms become quasi-isomorphisms
after applying the functors $\K\mpsto{}_{(n)}\K$ and
$\P\mpsto\P/I^n\P$.
 To prove that the morphism $\gF\rarrow\Hom_R(\D^\bu\;\D^\bu\ot_R\gF)$
is a quasi-isomorphism, one also has to use part~(b) of the following
lemma.
\end{proof}

\begin{lem} \label{all-contramods-complete-flat-ones-separated}
\textup{(a)} For any $(R,I)$\+contramodule\/ $\P$, the natural
map \textup{(}$(R,I)$\+con\-tra\-mod\-ule morphism\textup{)}
$\P\rarrow\varprojlim_n\P/I^n\P$ is surjective. \par
\textup{(b)} For any $R$\+flat $(R,I)$\+contramodule\/ $\P$, the natural
map\/ $\P\rarrow\varprojlim_n\P/I^n\P$ is an isomorphism.
\end{lem}

\begin{proof}
 Part~(a) is~\cite[Lemma~A.2.3 and Remark~A.3]{Psemi},
\cite[Theorem~5.6]{Pcta}, or~\cite[Lemma~6.3(b)]{PR}.
 Part~(b) is~\cite[proof of Lemma~B.9.2]{Pweak},
\cite[Corollary~10.3(b)]{Pcta}, or~\cite[Corollary~6.15]{PR};
see also Section~\ref{flat-contramod-subsect}.
\end{proof}

\subsection{Contraadjusted and cotorsion contramodules}
\label{cta-cot-contramod-subsect}
 For another discussion related to the material of
Sections~\ref{cta-cot-contramod-subsect}\+-%
\ref{veryflat-contra-subsect} below, see~\cite[Section~5]{PSl1}.

 Let $R$ be a Noetherian commutative ring and $s\in R$ be an element.
 In the spirit of the definitions from
Sections~\ref{very-eklof-trlifaj-subsect} and~\ref{torsion-subsect},
let us say that an $R$\+module $P$ is \emph{$s$\+contraadjusted} if
$\Ext_R^1(R[s^{-1}],P)=0$, and that $P$ is an \emph{$s$\+contramodule}
if $\Ext_R^*(R[s^{-1}],P)=0$ \,\cite[Section~2]{Pcta}.
 The property of an $R$\+module $P$ to be $s$\+contraadjusted or
an $s$\+contramodule only depends on the abelian group $P$ with
the operator $s\:P\rarrow P$.

 An $R$\+module $P$ is an $s$\+contramodule if and only if it is
an $(R,(s))$\+contramodule, where $(s)\sub R$ denotes the principal
ideal generated by $s\in R$.
 More generally, given an ideal $I\sub R$, an $R$\+module $\P$ is
an $(R,I)$\+contramodule if and only if it is an $s$\+contramodule
for every $s\in I$, and it suffices to check this condition for
any given set of generators of the ideal~$I$ 
\cite[proof of Theorem~B.1.1]{Pweak}, \cite[Theorem~5.1]{Pcta}.

 Recall that any quotient $R$\+module of an $s$\+contraadjusted
$R$\+module is $s$\+contra\-adjusted.
 Given an $R$\+module $L$, let us denote by $L(s)\sub L$ the image
of the morphism $\Hom_R(R[s^{-1}],L)\rarrow L$ induced by
the localization map $R\rarrow R[s^{-1}]$.
 Equivalently, the $R$\+submodule $L(s)\sub L$ can be defined
as the maximal $s$\+divisible $R$\+submodule in $R$, i.~e.,
the sum of all $R$\+submodules (or even all $s$\+invariant
abelian subgroups) in $L$ in which the element $s$~acts surjectively.

 Therefore, if $P=L/L(s)$ denotes the corresponding quotient
$R$\+module, then one has $P(s)=0$.
 Notice also that one has $\Hom_R(R[s^{-1}],P)=0$ for any
$R$\+module $P$ for which $P(s)=0$.
 It follows that the $R$\+quotient module $L/L(s)$ is
an $s$\+contramodule whenever an $R$\+module $L$ is
$s$\+contraadjusted.

 Now let $s$, $t\in R$ be two elements; suppose that an $R$\+module 
$L$ is a $t$\+contramodule.
 Then the $R$\+module $\Hom_R(R[s^{-1}],L)$ is also
a $t$\+contramodule~\cite[Section~B.2]{Pweak},
\cite[Lemma~6.1(b)]{Pcta}, as is the image $L(s)$ of the morphism of
$t$\+contramodules $\Hom_R(R[s^{-1}],L)\rarrow L$.
 Hence the quotient $R$\+module $L/L(s)$ is also
a $t$\+con\-tra\-mod\-ule.
 Assuming additionally that the $R$\+module $L$ was
$s$\+contraadjusted, the quotient module $L/L(s)$ is both
a $t$\+ and an $s$\+contramodule (i.~e., it is
an $(R,(t,s))$\+contramodule).

 Recall that an $R$\+module $L$ is said to be \emph{contraadjusted}
if it is $s$\+contraadjusted for every element $s\in R$.
 Given a contraadjusted $R$\+module $L$ and an ideal $I\sub R$,
one can apply the above construction of a quotient $R$\+module
successively for all the generators of the ideal~$I$.
 Proceeding in this way, one in a finite number of steps (equal to
the number of generators of~$I$) obtains the unique maximal
quotient $R$\+module $\P$ of the $R$\+module $L$ such that $\P$
is an $s$\+contramodule for every $s\in I$, i.~e.,
$\P$ is an $(R,I)$\+contramodule.
 We denote this quotient module by $\P=L/L(I)$.

 More generally, all one needs in order to apply the above construction
of the maximal quotient $(R,I)$\+contramodule $\P=L/L(I)$ to
a given $R$\+module $L$ is that the ideal $I\sub R$ should have a set
of generators $s_1$,~\dots, $s_k$ such that $L$ is
an $s_j$\+contramodule for every $j=1$,~\dots,~$k$.

\begin{lem}  \label{no-divisible-torsion}
 Let $R$ be a Noetherian commutative ring, $s\in R$ be an element,
$F$ be a flat $R$\+module, and $M$ be a finitely generated $R$\+module.
 Then the natural $R$\+module homomorphism\/
$\Hom_R(R[s^{-1}]\; M\ot_R F)\rarrow M\ot_R F$ is injective.
 In other words, one has\/ $\Hom_R(R[s^{-1}]/R\;M\ot_R F)=0$.
\end{lem}

\begin{proof}
 The point is that the $s$\+torsion in $M\ot_R F$ is bounded, i.~e.,
there exists an integer $n\ge1$ such that every $s$\+torsion element
in $M\ot_R F$ is annihilated by~$s^n$.
 It follows immediately that any $R$\+linear map $R[s^{-1}]/R\rarrow
M\ot_R F$ vanishes.

 To show that such an integer~$n$ exists, consider the sequence of
annihilator submodules of the elements $s^m\in R$, \ $m\ge 1$, in
the $R$\+module~$M$.
 This is an increasing sequence of $R$\+submodules in~$M$.
 Let $N$ be the maximal submodule in this sequence and $n$~be
a positive integer such that $N$ is the annihilator of~$s^n$ in~$M$.
 Then any element of the $R$\+module $M\ot_R F$ annihilated by
$s^{n+1}$ is also annihilated by~$s^n$.
 
 Indeed, the element~$s$ acts by an injective endomorphism of 
an $R$\+module $M/N$; hence so does the element~$s^{n+1}$.
 Since $F$ is flat, it follows that $s^{n+1}$ must act injectively
in the tensor product $(M/N)\ot_R F\simeq(M\ot_R F)/(N\ot_R F)$.
 Hence any element of $M\ot_R F$ annihilated by $s^{n+1}$ belongs
to $N\ot_R F$ and is therefore annihilated by~$s^n$.
\end{proof}

\begin{lem}  \label{cta-cot-preserved}
 Let $R$ be a Noetherian commutative ring, $s\in R$ be an element,
and $I\sub R$ be an ideal.  Then \par
\textup{(a)} whenever an $R$\+module $L$ is $s$\+contraadjusted,
the $R$\+modules $L(s)$ and $L/L(s)$ are also $s$\+contraadjusted
(and $L/L(s)$ is even an $s$\+contramodule); \par
\textup{(b)} whenever an $R$\+module $L$ is contraadjusted,
the $R$\+modules $L(s)$ and $L/L(s)$ are also contraadjusted; \par
\textup{(c)} whenever an $R$\+module $L$ is contraadjusted,
the $R$\+modules $L(I)$ and $L/L(I)$ are also contraadjusted
(and $L/L(I)$ is in addition an $(R,I)$\+contramodule); \par
\textup{(d)} whenever an $R$\+module $L$ is cotorsion and\/
$\Hom_R(R[s^{-1}]/R\;L)=0$, the $R$\+modules $L(s)$ and $L/L(s)$ are
also cotorsion.
\end{lem}

\begin{proof}
 The parenthesized assertions in~(a) and~(c) have been explained above.
 Since the class of cotorsion $R$\+modules is closed under
cokernes of injective morphisms, while the classes of contraadjusted
and $s$\+contraadjusted $R$\+modules are even closed under quotients,
it suffices to check the assertions related to the $R$\+modules
$L(s)$ and~$L(I)$.
 Furthermore, there is a natural surjective morphism of
$R$\+modules $\Hom_R(R[s^{-1}],L)\rarrow L(s)$, which in
the assumptions of~(d) is an isomorphism.

 Now part (d)~is a particular case of Lemma~\ref{cotors-hom}(a),
part~(b) is provided by Lemma~\ref{very-tensor-hom}(b), and part~(a)
follows from the similar claim that the $R$\+module
$\Hom_R(R[s^{-1}],P)$ is $s$\+contraadjusted for any
$s$\+contraadjusted $R$\+module~$P$.
 Finally, (b) implies~(c) in view of the above recursive construction
of the $R$\+module $L(I)$ and the fact that the class of
constraadjusted $R$\+modules is closed under extensions.
\end{proof}

\begin{cor}  \label{flat-preserved}
 Let $R$ be a Noetherian commutative ring, $s\in R$ be an element,
and $I\sub R$ be an ideal.
 In this setting \par
\textup{(a)} if $F$ is an $s$\+contraadjusted flat $R$\+module,
then the $R$\+modules $F(s)$ and $F/F(s)$ are also flat and
$s$\+contraadjusted (and $F/F(s)$ is even an $s$\+contramodule); \par
\textup{(b)} if $F$ is a contraadjusted flat $R$\+module,
then $F(I)$ is a contraadjusted flat $R$\+module and $F/F(I)$
is an $R$\+flat $R$\+contraadjusted $(R,I)$\+contramodule; \par
\textup{(c)} if $F$ is a flat cotorsion $R$\+module, then $F(I)$ is
a flat cotorsion $R$\+module and $F/F(I)$ is an $R$\+flat
$R$\+cotorsion $(R,I)$\+contramodule.
\end{cor}

\begin{proof}
 Part~(a): in order to prove that the $R$\+module $F/F(s)$ is flat,
let us check that the map $M\ot_R F(s)\rarrow M\ot_R F$ induced by
the natural embedding $F(s)\rarrow F$ is injective for any
finitely generated $R$\+module~$M$.
 By Lemma~\ref{no-divisible-torsion}, one has
$F(s)\simeq\Hom_R(R[s^{-1}],F)\sub F$ and $(M\ot_R F)(s)\simeq
\Hom_R(R[s^{-1}]\;M\ot_R F)\sub M\ot_R F$.
 It remains to point out that the natural morphism
$M\ot_R\Hom_R(R[s^{-1}],F)\rarrow\Hom_R(R[s^{-1}]\;M\ot_R F)$
is an isomorphism by Lemma~\ref{coherent-lemma}.

 Part~(b): the above recursive construction of
the $(R,I)$\+contramodule $F/F(I)$ together with part~(a) imply
the assertion that $F/F(I)$ is a flat $R$\+module.
 Now it remains to use Lemma~\ref{cta-cot-preserved}(c).

 Part~(c) can be proved in the way similar to part~(b), using
the fact that the class of cotorsion $R$\+modules is closed
under extensions and Lemma~\ref{cta-cot-preserved}(d).
 Alternatively, (c)~can be deduced from 
Theorem~\ref{flat-cotorsion-classification}.
 Indeed, if $F\simeq\prod_\p \gF_\p$, where $\p\sub R$ are
prime ideals and $\gF_\p$ are $(R,\p)$\+contramodules, then
one easily checks that $F/F(I)$ is the product of $\gF_\p$ over
the prime ideals $\p$ containing~$I$, while $F(I)$ is
the product of $\gF_\q$ over all the other prime ideals~$\q$.
\end{proof}

\begin{lem}  \label{cta-cot-contra-envelope-lemma}
\textup{(a)} Let\/ $0\rarrow P\rarrow K\rarrow F\rarrow 0$ be
a short exact sequence of $R$\+modules, where $P$ is
an $s$\+contramodule, $K$ is a contraadjusted $R$\+module,
and $F$ is a flat $R$\+module.
 Then\/ $0\rarrow P\rarrow K/K(s)\rarrow F/F(s)\rarrow 0$ is
a short exact sequence of $R$\+modules in which $K/K(s)$
is a contraadjusted $R$\+module, $F/F(s)$ is a flat $R$\+module,
and all the three modules are $s$\+contramodules. \par
\textup{(b)} Let\/ $0\rarrow P\rarrow K\rarrow F\rarrow 0$ be
a short exact sequence of $R$\+modules, where $P$ is
an $s$\+contramodule, $K$ is a cotorsion $R$\+module,
and $F$ is a flat $R$\+module.
 Then\/ $0\rarrow P\rarrow K/K(s)\rarrow F/F(s)\rarrow 0$ is
a short exact sequence of $R$\+modules in which $K/K(s)$
is a cotorsion $R$\+module, $F/F(s)$ is a flat $R$\+module,
and all the three modules are $s$\+contramodules. 
\end{lem}

\begin{proof}
 Part~(a): first of all let us show that the morphism of
$R$\+modules $K\rarrow F$ restricts to an isomorphism of their
submodules $K(s)\rarrow F(s)$.
 Indeed, we have $\Hom_R(R[s^{-1}]/R\;F)=0$ by
Lemma~\ref{no-divisible-torsion} and $\Hom_R(R[s^{-1}]/R\;P)
\sub\Hom_R(R[s^{-1}],P)=0$ by the definition of
an $s$\+contramodule, hence $\Hom_R(R[s^{-1}]/R\;\allowbreak K)=0$.
 Therefore, there are isomorphisms $K(s)\simeq
\Hom_R(R[s^{-1}],K)$ and $F(s)\simeq\Hom_R(R[s^{-1}],F)$.
 Using the condition that $P$ is an $s$\+contramodule, that is
$\Ext^*_R(R[s^{-1}],P)=0$ again, we conclude that $K(s)\simeq F(s)$.

 It follows that the sequence $0\rarrow P\rarrow K/K(s)\rarrow
F/F(s)\rarrow0$ is exact.
 The $R$\+module $K$ being contraadjusted by assumption, its
quotient $R$\+module $F$ is contraadjusted, too.
 Now the quotient module $K/K(s)$ is contraadjusted by
Lemma~\ref{cta-cot-preserved}(b), the quotient module $F/F(s)$
is flat by Corollary~\ref{flat-preserved}(a), and both are
$s$\+contramodules by Lemma~\ref{cta-cot-preserved}(a).

 In part~(b), it only remains to prove that $K/K(s)$ is
a cotorsion $R$\+module.
 This follows from the argument above and
Lemma~\ref{cta-cot-preserved}(d).
\end{proof}

\begin{cor}  \label{cta-cot-contra-envelope-cor}
\textup{(a)} Let\/ $0\rarrow P\rarrow K\rarrow F\rarrow 0$ be
a short exact sequence of $R$\+modules, where $P$ is
an $(R,I)$\+contramodule, $K$ is a contraadjusted $R$\+module,
and $F$ is a flat $R$\+module.
 Then\/ $0\rarrow P\rarrow K/K(I)\rarrow F/F(I)\rarrow 0$ is
a short exact sequence of\/ $(R,I)$\+contramodules in which $K/K(I)$ is
a contraadjusted $R$\+module and $F/F(I)$ is a flat $R$\+module. \par
\textup{(b)} Let\/ $0\rarrow P\rarrow K\rarrow F\rarrow 0$ be
a short exact sequence of $R$\+modules, where $P$ is
an $(R,I)$\+contramodule, $K$ is a cotorsion $R$\+module,
and $F$ is a flat $R$\+module.
 Then\/ $0\rarrow P\rarrow K/K(I)\rarrow F/F(I)\rarrow 0$ is
a short exact sequence of\/ $(R,I)$\+contramodules in which $K/K(I)$
is a cotorsion $R$\+module and $F/F(I)$ is a flat $R$\+module.
\end{cor}

\begin{proof}
 Follows by recursion from Lemma~\ref{cta-cot-contra-envelope-lemma}.
\end{proof}

\begin{lem} \label{flat-contra-cover-lemma}
\textup{(a)} Let\/ $0\rarrow K\rarrow F\rarrow P\rarrow 0$ be
a short exact sequence of $R$\+modules, where $K$ is
a contraadjusted $R$\+module, $F$ is a flat $R$\+module,
and $P$ is an $s$\+contramodule.
 Then\/ $0\rarrow K/K(s)\rarrow F/F(s)\rarrow P\rarrow 0$ is
a short exact sequence of $R$\+modules in which $K/K(s)$
is a contraadjusted $R$\+module, $F/F(s)$ is a flat $R$\+module,
and all the three modules are $s$\+contramodules. \par
\textup{(b)} Let\/ $0\rarrow K\rarrow F\rarrow P\rarrow 0$ be
a short exact sequence of $R$\+modules, where $K$ is a cotorsion
$R$\+module, $F$ is a flat $R$\+module, and $P$ is
an $s$\+contramodule.
 Then\/ $0\rarrow K/K(s)\rarrow F/F(s)\rarrow P\rarrow 0$ is
a short exact sequence of $R$\+modules in which $K/K(s)$
is a cotorsion $R$\+module, $F/F(s)$ is a flat $R$\+module,
and all the three modules are $s$\+contramodules.
\end{lem}

\begin{proof}
 Part~(a): first we prove that the morphism of $R$\+modules
$K\rarrow F$ restricts to an isomorphism of their submodules
$K(s)\rarrow F(s)$.
 Indeed, $\Hom_R(R[s^{-1}]/R\;K)\sub \Hom_R(R[s^{-1}]/R\;F) = 0$
by Lemma~\ref{no-divisible-torsion}, hence $K(s)\simeq
\Hom_R(R[s^{-1}],K)$ and $F(s)\simeq\Hom_R(R[s^{-1}],F)$.
 Since $\Hom_R(R[s^{-1}],P)=0$ by assumption, we conclude
that $K(s)\simeq F(s)$, and it follows that the sequence
$0\rarrow K/K(s)\rarrow F/F(s)\rarrow P\rarrow 0$ is exact.

 The $R$\+module $F$, being an extension of $s$\+contraadjusted
$R$\+modules $K$ and $P$, is $s$\+contraadjusted, too.
 The rest of the argument coincides with the respective part
of the proof of Lemma~\ref{cta-cot-contra-envelope-lemma},
and so does the proof of part~(b).
\end{proof}

\begin{cor}  \label{flat-contra-cover-cor}
\textup{(a)} Let\/ $0\rarrow K\rarrow F\rarrow P\rarrow 0$ be
a short exact sequence of $R$\+modules, where $K$ is
a contraadjusted $R$\+module, $F$ is a flat $R$\+module,
and $P$ is an $(R,I)$\+contramodule.
 Then\/ $0\rarrow K/K(I)\rarrow F/F(I)\rarrow P\rarrow 0$ is
a short exact sequence of $(R,I)$\+contramodules in which $K/K(I)$ is
a contraadjusted $R$\+module and $F/F(I)$ is a flat $R$\+module. \par
\textup{(b)} Let\/ $0\rarrow K\rarrow F\rarrow P\rarrow 0$ be
a short exact sequence of $R$\+modules, where $K$ is
a cotorsion $R$\+module, $F$ is a flat $R$\+module,
and $P$ is an $(R,I)$\+contramodule.
 Then\/ $0\rarrow K/K(I)\rarrow F/F(I)\rarrow P\rarrow 0$ is
a short exact sequence of $(R,I)$\+contramodules in which $K/K(I)$
is a cotorsion $R$\+module and $F/F(I)$ is a flat $R$\+module.
\end{cor}

\begin{proof}
 Follows by recursion from Lemma~\ref{flat-contra-cover-lemma}.
\end{proof}

 Recall that, according to~\cite[Theorem~B.8.1]{Pweak}, the $\Ext$
groups/modules computed in the abelian categories $R\modl$ and
$(R,I)\contra$ agree (see also~\cite[Theorem~2.9]{Pmgm}).
 
\begin{cor}  \label{r-i-contra-cor}
 Let $R$ be a Noetherian commutative ring and $I\sub R$ be an ideal.
 Then \par
\textup{(a)} any $(R,I)$\+contramodule can be embedded into
an $R$\+cotorsion $(R,I)$\+contra\-module in such a way that
the quotient $(R,I)$\+contramodule is $R$\+flat; \par
\textup{(b)} any $(R,I)$\+contramodule admits a surjective morphism
onto it from an $R$\+flat $(R,I)$\+contramodule such that the kernel
is an $R$\+cotorsion $(R,I)$\+contramodule; \par
\textup{(c)} an $(R,I)$\+contramodule\/ $\Q$ is $R$\+cotorsion if
and only if one has\/ $\Ext_R^1(\gF,\Q)=0$ for any $R$\+flat
$(R,I)$\+contramodule\/~$\gF$; \par
\textup{(d)} an $(R,I)$\+contramodule\/ $\gF$ is $R$\+flat if and
only if one has\/ $\Ext_R^1(\gF,\Q)=0$ for any $R$\+cotorsion
$(R,I)$\+contramodule\/~$\Q$.
\end{cor}

\begin{proof}
 Parts~(a\+b) follow from Theorem~\ref{flat-cover-thm} together with
Corollaries~\ref{cta-cot-contra-envelope-cor}(b)
and~\ref{flat-contra-cover-cor}(b).
 Part~(c) is deduced from~(a) and part~(d) deduced from~(b) easily.
\end{proof}

 Let us denote by $(R,I)\contra^\cta$ the full exact subcategory of
$R$\+contraadjusted $(R,I)$\+contramodules and by $(R,I)\contra^\cot$
the full exact subcategory of $R$\+cotorsion $(R,I)$\+contramodules
in the abelian category $(R,I)\contra$.

 Notice that the full subcategory $(R,I)\contra^\cot\sub(R,I)\contra$
depends only on the $I$\+adic completion of the ring~$R$, as one
can see from Corollary~\ref{r-i-contra-cor}(c).
 The similar assertion for the full subcategory
$(R,I)\contra^\cta\sub(R,I)\contra$ will follow from the results
of Section~\ref{veryflat-contra-subsect} below.

\begin{thm}
\textup{(a)} Let $R$ be a Noetherian commutative ring and $I\sub R$
be an ideal.
 Then for any symbol\/ $\bst=\b$, $+$, $-$, $\abs+$, $\abs-$, $\bctr$,
$\ctr$, or~$\abs$, the triangulated functor\/
$\sD^\st((R,I)\contra^\cta)\rarrow\sD^\st((R,I)\contra)$ induced by
the embedding of exact categories $(R,I)\contra^\cta\rarrow(R,I)\contra$
is an equivalence of triangulated categories.  \par
\textup{(b)} Let $R$ be a Noetherian commutative ring and $I\sub R$
be an ideal such that the quotient ring $R/I$ has finite Krull
dimension.
 Then for any symbol\/ $\bst=\b$, $+$, $-$, $\abs+$, $\abs-$, $\bctr$,
$\ctr$, or~$\abs$, the triangulated functor\/
$\sD^\st((R,I)\contra^\cot)\rarrow\sD^\st((R,I)\contra)$ induced by
the embedding of exact categories $(R,I)\contra^\cot\rarrow(R,I)\contra$
is an equivalence of triangulated categories.
\end{thm}

\begin{proof}
 Part~(a) for all the symbols~$\bst$ except $\bctr$ follows
from Corollary~\ref{cta-cot-contra-envelope-cor}(a)
or~\ref{r-i-contra-cor}(a) together with the opposite version 
of Proposition~\ref{finite-resolutions}.
 In the case $\bst=\bctr$, the argument is based on
Proposition~\ref{becker-contraderived-finite-coresolutions} applied
to the abelian category $\sE=(R,I)\contra$ with the hereditary
complete cotorsion pair $(\sF,\sC)$, where $\sF$ is the class of
all very flat $(R,I)$\+contramodules (as defined in
Section~\ref{veryflat-contra-subsect} below) and
$\sC=(R,I)\contra^\cta$.
 By Corollary~\ref{contramod-veryflat-cor}, \,$(\sF,\sC)$ is a complete
cotorsion pair in~$\sE$.
 Since all $R$\+modules have contraadjusted dimensions at most~$1$,
it follows that all $(R,I)$\+contramodules have $\sC$\+coresolution
dimensions at most~$1$, too; so 
Proposition~\ref{becker-contraderived-finite-coresolutions}
is applicable.

 To prove part~(b) in the similar way, one needs to know that
in its assumptions any $(R,I)$\+contramodule has finite coresolution
dimension with respect to the coresolving subcategory
$(R,I)\contra^\cot \sub (R,I)\contra$.
 In view of Corollary~\ref{r-i-contra-cor}(c), this follows from
the fact that any $R$\+flat $(R,I)$\+contramodule has finite
projective dimension in $(R,I)\contra$, which was established in
the proof of Theorem~\ref{tors-contra-co-contra-derived}
(using~\cite[Corollary~B.8.2]{Pweak} or~\cite[Theorem~10.5]{Pcta}
together with Theorem~\ref{raynaud-gruson-flat-thm}).
 In the case $\bst=\bctr$, one can apply
Proposition~\ref{becker-contraderived-finite-coresolutions}
to the abelian category $\sE=(R,I)\contra$ with the hereditary
complete cotorsion pair $(\sF,\sC)$, where $\sF$ is the class of
all $R$\+flat $(R,I)$\+contramodules and $\sC=(R,I)\contra^\cot$.
 By Corollary~\ref{r-i-contra-cor}, \,$(\sF,\sC)$ is a complete
cotorsion pair in~$\sE$.
\end{proof}

\begin{rem}
 Let $X=\Spec R$ be the Noetherian affine scheme corresponding to
the ring $R$ and $U=X\setminus Z$ be the open complement to the closed
subscheme $Z=\Spec R/I\sub X$.
 Denote by $j\:U\rarrow X$ the natural open embedding morphism. 
 Then the exact category $(R,I)\contra^\cta$ is equivalent to
the full exact subcategory in the exact category $X\ctrh$ of
contraherent cosheaves on $X$ consisting of all the contraherent
cosheaves $\P\in X\ctrh$ with vanishing restrictions $j^!\P=\P|_U$
to the open subscheme~$U$.
 Similarly, the exact category $(R,I)\contra^\cot$ is equivalent
to the full exact subcategory in the exact category $X\ctrh^\lct$
consisting of all the (locally) cotorsion contraherent cosheaves
$\P$ on $X$ for which $j^!\P=0$.
 Indeed, a contraadjusted $R$\+module $P$ is an $s$\+contramodule
if and only if the corresponding contraherent cosheaf $\P=\widecheck P$
on $X$ vanishes in the restriction to $\Spec R[s^{-1}]\sub\Spec R$
(see Sections~\ref{contraherent-definition}\+-%
\ref{cohom-subsection} for the definitions and notation).
 In particular, when the ring $R/I$ is Artinian, the exact category
$\ker(j^!\:X\ctrh\to U\ctrh)=\ker(j^!\:X\ctrh^\lct\to U\ctrh^\lct)$
is abelian and equivalent to the abelian category
$(R,I)\contra=(R,I)\contra^\cta=(R,I)\contra^\cot$.
 The point here is that, over a complete Noetherian local ring, all
contramodules are cotorsion modules (see
Proposition~\ref{contramodules-cotorsion}(a)).
\end{rem}

\subsection{Very flat contramodules}  \label{veryflat-contra-subsect}
 Unlike the flatness, cotorsion, and contraadjustedness properties of
$(R,I)$\+contramodules, their very flatness property is \emph{not}
defined in terms of the similar property of $R$\+modules.
 Instead, it is described in terms of the reductions modulo $I^n$
and the very flatness properties of $R/I^n$\+modules.

 For a (straightforward) generalization of the next lemma,
see Lemma~\ref{mod-contramod-flat-ext-lemma}.

\begin{lem}  \label{vfl-cta-reduction-orthogonality}
\textup{(a)} Let $F$ be a flat $R$\+module and\/ $\Q$ be
an $(R,I)$\+contramodule such that the $R/I$\+module $F/IF$
is very flat, the natural $(R,I)$\+contramodule morphism\/
$\Q\rarrow\varprojlim_n\Q/I^n\Q$ is an isomorphism, and
the $R/I$\+module\/ $\Q/I\Q$ is contraadjusted.
 Then one has\/ $\Ext^{>0}_R(F,\Q)=0$. \par
\textup{(b)} Let $F$ be a flat $R$\+module and\/ $\Q$ be
an $R$\+flat $(R,I)$\+contramodule such that
the $R/I$\+module\/ $\Q/I\Q$ is cotorsion.
 Then one has\/ $\Ext^{>0}_R(F,\Q)=0$.
\end{lem}

\begin{proof}
 Part~(a): by (the proof of)
Lemma~\ref{fin-gen-nilpotent-quotient-scalars}(a),
the $R/I^n$\+modules $\Q/I^n\Q$ and the $R/I$\+modules
$I^{n-1}\Q/I^n\Q$ are contraadjusted, while by part~(b) of the same
lemma the $R/I^n$\+mod\-ules $F/I^nF$ are very flat.
 Now we follow the argument in the proof
of~\cite[Proposition~B.10.1]{Pweak}.
 The $R$\+module $F$ being flat by assumption, one has $\Ext^i_R(F\;
\Q/I^n\Q)\simeq\Ext^i_{R/I^n}(F/I^nF\;\Q/I^n\Q)=0$ for all $i>0$
and any~$n\ge1$.
 The natural map $\Hom_R(F\;\Q/I^n\Q)\rarrow\Hom_R(F\;\Q/I^{n-1}\Q)$
is surjective, since $\Ext^1_{R/I^n}(F/I^nF\;I^{n-1}\Q/I^n\Q)=0$.
 By~\cite[Lemma~B.10.3]{Pweak}, we have $\Ext_R^i(F\;\varprojlim_n
\Q/I^n\Q)\simeq\varprojlim_n^i\Hom_R(F\;\Q/I^n\Q)=0$ for all $i>0$.

 Part~(b): by (the proof of)
Lemma~\ref{fin-gen-nilpotent-quotient-scalars}(c),
the $R/I^n$\+modules $\Q/I^n\Q$ and the $R/I$\+modules $I^{n-1}\Q/I^n\Q$
are cotorsion, while by
Lemma~\ref{all-contramods-complete-flat-ones-separated}(b)
the natural map $\Q\rarrow\varprojlim_n\Q/I^n\Q$ is an isomorphism.
 The argument continues exactly the same as in part~(a).
\end{proof}

 For generalizations of respective parts of the following
corollary, see Lemma~\ref{contraadjusted-reduction-lemma} 
and Corollary~\ref{flat-cotors-mod-contramod-cor}.

\begin{cor}  \label{contraadjusted-reduction}
\textup{(a)} For any contraadjusted $R$\+module $Q$,
the $R/I$\+module $Q/IQ$ is contraadjusted.
 If the natural map $Q\rarrow\varprojlim_n Q/I^nQ$ is an isomorphism
for an $R$\+module $Q$ and the $R/I$\+module $Q/IQ$ is contraadjusted,
then the $R$\+module $Q$ is contraadjusted. \par
\textup{(b)} A flat $(R,I)$\+contramodule\/ $\Q$ is $R$\+cotorsion
if and only if the $R/I$\+module\/ $\Q/I\Q$ is cotorsion.
\end{cor}

\begin{proof}
 Part~(a): the first assertion is
Lemma~\ref{quotient-scalars-contraadjusted}(b).
 Since the $R/I$\+module $F/IF$ is very flat for any very flat
$R$\+module $F$ (see Lemma~\ref{very-scalars-always}(b)),
the second assertion follows from
Lemma~\ref{vfl-cta-reduction-orthogonality}(a).
 Part~(b) is similarly a consequence of
Lemma~\ref{quotient-scalars-cotorsion}(a) and
Lemma~\ref{vfl-cta-reduction-orthogonality}(b).
\end{proof}

 For a generalization of the next corollary, see
Corollary~\ref{cta-mod-contramod-cor}.

\begin{cor}  \label{veryflat-reduction}
 Let $F$ be a flat $R$\+module for which the $R/I$\+module
$F/IF$ is very flat.
 Then\/ $\Ext^{>0}_R(F,\P)=0$ for any $R$\+contraadjusted\/
$(R,I)$\+contramodule\/~$\P$.
\end{cor}

\begin{proof}
 By Corollary~\ref{flat-contra-cover-cor}(a) or~\ref{r-i-contra-cor}(b),
there exists a short exact sequence of $(R,I)$\+contramodules
$0\rarrow\Q\rarrow\gG\rarrow\P\rarrow0$, where the $R$\+module $\Q$
is contraadjusted and the $R$\+module $\gG$ is flat.
 Then the $R$\+module $\gG$ is also contraadjusted.
 By Corollary~\ref{contraadjusted-reduction}(a), it follows
that so are the $R/I$\+modules $\Q/I\Q$ and $\gG/I\gG$.
 Furthermore, according to
Lemma~\ref{all-contramods-complete-flat-ones-separated}(b)
one has $\gG\simeq\varprojlim_n\gG/I^n\gG$.
 The $(R,I)$\+contramodule $\Q$ being a subcontramodule of $\gG$,
it follows that $\bigcap_n I^n\Q=0$, hence also
$\Q\simeq\varprojlim_n\Q/I^n\Q$
by Lemma~\ref{all-contramods-complete-flat-ones-separated}(a).
 Now it remains to apply Lemma~\ref{vfl-cta-reduction-orthogonality}(a)
to (the $R$\+module $F$ and)
the $(R,I)$\+contramodules $\Q$ and~$\gG$.
\end{proof}

\begin{lem}  \label{vfl-cta-contra-envelope-cover-lemma}
\textup{(a)} Any $(R,I)$\+contramodule\/ $\gM$ can be included in
a short exact sequence of $(R,I)$\+contramodules\/ $0\rarrow\gM\rarrow
\P\rarrow\gF\rarrow0$, where the $R/I^n$\+modules\/ $\gF/I^n\gF$
are very flat, while the $R$\+module\/ $\P$ is contraadjusted. \par
\textup{(b)} Any $(R,I)$\+contramodule\/ $\gM$ can be included in
a short exact sequence of $(R,I)$\+contramodules\/ $0\rarrow \P\rarrow
\gF\rarrow\gM\rarrow0$, where the $R/I^n$\+modules\/ $\gF/I^n\gF$
are very flat, while the $R$\+module\/ $\P$ is contraadjusted.
\end{lem}

\begin{proof}
 First of all, let us show that any $R$\+flat $(R,I)$\+contramodule
$\gG$ can be included in a short exact sequence of
$(R,I)$\+contramodules $0\rarrow\gE\rarrow\gF\rarrow\gG\rarrow0$,
where the $R/I^n$\+modules $\gF/I^n\gF$ are very flat, while
the $R$\+module $\gE$ is (flat and) contraadjusted.
 By Theorem~\ref{eklof-trlifaj-very}(b), there exists
a short exact sequence of $R$\+modules $0\rarrow E\rarrow F\rarrow
\gG\rarrow 0$ such that the $R$\+module $F$ is very flat,
while the $R$\+module $E$ is contraadjusted.
 Then the $R$\+module $E$ is also flat.
 
 The short sequence of $R/I^n$\+modules $0\rarrow E/I^n E\rarrow
F/I^n F\rarrow\gG/I^n\gG\rarrow 0$ is exact for every $n\ge1$.
 Set $\gE=\varprojlim_n E/I^nE$ and $\gF=\varprojlim_n F/I^nF$;
recall that the natural map $\gG\rarrow\varprojlim_n\gG/I^n\gG$ is
an isomorphism by
Lemma~\ref{all-contramods-complete-flat-ones-separated}(b).
 Passing to the projective limit, we obtain a short exact sequence
of $(R,I)$\+contramodules $0\rarrow\gE\rarrow\gF\rarrow\gG\rarrow0$.
 The $R/I^n$\+modules $\gF/I^n\gF\simeq F/I^n F$ are very flat
by Lemma~\ref{very-scalars-always}(b), and the $R$\+module $\gE$
is contraadjusted by Corollary~\ref{contraadjusted-reduction}(a).

 Now it is easy to obtain part~(a) from 
Corollary~\ref{r-i-contra-cor}(a) and part~(b) from
Corollary~\ref{r-i-contra-cor}(b).
 Specifically, let $0\rarrow\gM\rarrow\Q\rarrow\G\rarrow0$ be
a short exact sequence of $(R,I)$\+contramodules with an
$R$\+cotorsion $(R,I)$\+contramodule $\Q$ and an $R$\+flat
$(R,I)$\+contramodule $\gG$, as in
Corollary~\ref{r-i-contra-cor}(a).
 Let $0\rarrow\gE\rarrow\gF\rarrow\gG\rarrow0$ be the short exact
sequence constructed above.
 Denote by $\P$ the pullback of the pair of surjective morphisms
$\Q\rarrow\gG$ and $\gF\rarrow\gG$.
 Then the $R$\+module $\P$ is contraadjusted as an extension of two
contraadjusted $R$\+modules $\Q$ and~$\gE$.

 Similarly, let $0\rarrow\Q\rarrow\gG\rarrow\gM\rarrow0$ be
a short exact sequence of $(R,I)$\+contramodules with an
$R$\+cotorsion $(R,I)$\+contramodule $\Q$ and an $R$\+flat
$(R,I)$\+contramodule $\gG$, as in
Corollary~\ref{r-i-contra-cor}(b).
 Let $0\rarrow\gE\rarrow\gF\rarrow\gG\rarrow0$ be the short exact
sequence constructed above.
 Then the kernel $\P$ of the composition of surjective morphisms
$\gF\rarrow\gG\rarrow\gM$ is contraadjusted as an extension of two
contraadjusted $R$\+modules $\Q$ and~$\gE$.

 Alternatively, in order to obtain the assertions of Lemma
one can simply apply the constructions of
Corollaries~\ref{cta-cot-contra-envelope-cor}(a)
and~\ref{flat-contra-cover-cor}(a), taking $F$ to be a very flat
$R$\+module and using Lemma~\ref{very-scalars-always}(b) together with
the observation that one has $(L/L(I))/I^n(L/L(I))\simeq L/I^nL$
for every $R$\+module $L$ for which the ideal $I\sub R$ admits a set of
generators~$s_j$ such that $L$ is $s_j$\+contraadjusted for every~$j$.
\end{proof}

 Let us call an $(R,I)$\+contramodule $\gF$ \emph{very flat} if 
the functor $\Hom_R(\gF,{-})$ takes short exact sequences of
$R$\+contraadjusted $(R,I)$\+contramodules to short exact sequences
of abelian groups.
 The next corollary says that this definition is equivalent
to the more familiar formulation is terms of the $\Ext^1$ vanishing
(cf.\ the definitions of very flat modules and contraadjusted
sheaves in Sections~\ref{very-eklof-trlifaj-subsect}
and~\ref{fHom-subsection}).

\begin{cor}  \label{contramod-veryflat-cor}
 Let $R$ be a Noetherian commutative ring and $I\sub R$ be an ideal.
 Then \par
\textup{(a)} any $(R,I)$\+contramodule can be embedded into
an $R$\+contraadjusted $(R,I)$\+con\-tramodule in such a way that
the quotient $(R,I)$\+contramodule is very flat; \par
\textup{(b)} any $(R,I)$\+contramodule admits a surjective morphism
onto it from a very flat $(R,I)$\+contramodule such that the kernel
is an $R$\+contraadjusted $(R,I)$\+contramodule; \par
\textup{(c)} an $(R,I)$\+contramodule\/ $\Q$ is $R$\+contraadjusted
if and only if\/ $\Ext_R^1(\gF,\Q)=0$ for any very flat
$(R,I)$\+contramodule\/~$\gF$; \par
\textup{(d)} an $(R,I)$\+contramodule\/ $\gF$ is very flat if and
only if\/ $\Ext_R^1(\gF,\Q)=0$ for any $R$\+contraadjusted
$(R,I)$\+contramodule\/~$\Q$; \par
\textup{(e)} an $(R,I)$\+contramodule\/ $\gF$ is very flat if and only
if the $R/I^n$\+module\/ $\gF/I^n\gF$ is very flat for every $n\ge1$.
\end{cor}

\begin{proof}
 The ``if'' assertion in part~(e) follows from
Corollary~\ref{veryflat-reduction} (together with the result
of~\cite[Lemma~B.9.2]{Pweak} or~\cite[Corollary~10.3(a)]{Pcta}
implying that the $R$\+module $\gF$ is flat).
 Parts~(a\+b) are provided by the respective parts of 
Lemma~\ref{vfl-cta-contra-envelope-cover-lemma} together with
the ``if'' assertion in~(e).
 The ``if'' assertions in parts~(c) and~(d) are deduced from
parts (a) and~(b), respectively; and to prove the ``only if'',
one only needs to know that any $(R,I)$\+contramodule can be
embedded into an $R$\+contraadjusted one.
 Finally, the ``only if'' assertion in part~(e) follows from
the construction in Lemma~\ref{vfl-cta-contra-envelope-cover-lemma}
on which the proof of part~(b) is based. \hbadness=1425
\end{proof}

\subsection{Affine geometry of $(R,I)$-contramodules}
 All rings in this section are presumed to be commutative and
Noetherian.
 Let $f\:R\rarrow S$ be a ring homomorphism, $I$ be an ideal in $R$,
and $J$ be an ideal in~$S$.
 We denote by $Sf(I)$ the extension of the ideal $I\sub R$ in
the ring~$S$.

\begin{lem}  \label{contramod-restrict}
\textup{(a)} If $f(I)\sub J$, then any $(S,J)$\+contramodule is
also an $(R,I)$\+con\-tramodule in the $R$\+module structure obtained
by restriction of scalars via~$f$. \par
\textup{(b)} If $J\sub Sf(I)$, then an $S$\+module is
an $(S,J)$\+contramodule whenever it is an $(R,I)$\+contramodule
in the $R$\+module structure obtained by restriction of
scalars via~$f$.
\end{lem}

\begin{proof}
 Part~(a) holds, since an $S$\+module $Q$ is an $(S,J)$\+contramodule
if and only if the system of equations $q_n=p_n+tq_{n+1}$, \ $n\ge0$,
is uniquely solvable in $q_n\in Q$ for any fixed sequence $p_n\in Q$
and any $t\in J$ \,\cite[Lemma~2.1(b)]{Pcta}.
 Part~(b) is true, because it suffices to check the previous
condition for the elements~$t$ belonging to any given set of
generators of the ideal $J\sub S$ only~\cite[Sections~B.1
and~B.7]{Pweak}, \cite[Theorem~5.1]{Pcta}.
\end{proof}

\begin{lem}  \label{contramod-exten}
 Assume that $J\sub Sf(I)$, and let\/ $\P$ be an $(R,I)$\+contramodule.
Then \par
\textup{(a)} the $S$\+module\/ $\Hom_R(S,\P)$ is
an $(S,J)$\+contramodule; \par
\textup{(b)} the $S$\+module $S\ot_R\P$ is an $(S,J)$\+contramodule
whenever $f$~is a finite morphism.
\end{lem}

\begin{proof} \hbadness=1300
 Both assertions follow from Lemma~\ref{contramod-restrict}(b).
 Indeed, the $R$\+module $\Hom_R(M,\P)$ is an $(R,I)$\+contramodule
for any $R$\+module $M$, because the contramodule infinite summation
operations can be defined on it (see the beginning of
Section~\ref{torsion-subsect}) or by~\cite[Lemma~6.1(b)]{Pcta}.
 The $R$\+module $M\ot_R\P$ is an $(R,I)$\+contramodule for any finitely
generated $R$\+module $M$, because $M\ot_R\P$ is the cokernel of
a morphism between two finite direct sums of copies of the
$(R,I)$\+contramodule~$\P$.
\end{proof}

 Denote by $\R$ the $I$\+adic completion of the ring $R$ and
by $\S$ the $J$\+adic completion of the ring $S$, both viewed as
topological rings.
 By~\cite[Theorem~B.1.1]{Pweak} or~\cite[Proposition~1.5 and
Corollary~3.7]{Pdc}, the full subcategories
$(R,I)\contra\sub R\modl$ and $(S,J)\contra\sub S\modl$ in
the categories of $R$- and $S$\+modules are equivalent
to the categories $\R\contra$ and $\S\contra$ of contramodules
over the topological rings $\R$ and~$\S$.

 Assume that $f(I)\sub J$; then the ring homomorphism $f\:R\rarrow S$
induces a continuous homomorphism of topological rings
$\phi\:\R\rarrow\S$.
 According to~\cite[Section~1.8]{Pweak} or~\cite[Section~2.9]{Pproperf}, 
there is a pair of adjoint functors of ``contrarestriction'' and
``contraextension'' of scalars $R^\phi\:\S\contra\rarrow\R\contra$ and
$E^\phi\:R\contra\rarrow\S\contra$.
 While the functor $R^\phi$ is easily identified with the functor
of restriction of scalars from Lemma~\ref{contramod-restrict},
the functor $E^\phi$, which is left adjoint to $R^\phi$, is defined
by the rules that $E^\phi$ is right exact and takes the free
$\R$\+contramodule $\R[[X]]=\varprojlim_n R/I^n[X]$ to the free
$\S$\+contramodule $\S[[X]]=\varprojlim_nS/J^n[X]$ for any set~$X$.
 Here $A[x]=A^{(X)}$ is a notation for the direct sum of $X$ copies
of an abelian group~$A$.

 When $f$~is a finite morphism and $J=Sf(I)$, the functor $E^\phi$
is simply the functor of extension of scalars from
Lemma~\ref{contramod-exten}(b).

\begin{lem}  \label{contramod-contraexten}
 For any $R$\+flat $(R,I)$\+contramodule\/ $\gF$ there is
a natural isomorphism of $(S,J)$\+contramodules $E^\phi(\gF)\simeq
\varprojlim_n(S\ot_R\gF)/\allowbreak J^n(S\ot_R\gF)$.
 In particular, the functor $E^\phi$ takes $R$\+flat
$(R,I)$\+contramodules to $S$\+flat $(S,J)$\+contramodules and
very flat $(R,I)$\+contramodules to very flat
$(S,J)$\+contramodules.
\end{lem}

\begin{proof}
 Since $F\mpsto(S\ot_RF)/J^n(S\ot_RF)=(S\ot_RF)/(J^n\ot_RF)=
S/J^n\ot_RF$ is an exact functor on the category of flat
$R$\+modules $F$, so is the functor
$F\mpsto\varprojlim_n(S\ot_RF)/J^n(S\ot_RF)$.
 Hence it suffices to compute this functor for free
$\R$\+con\-tramodules, and indeed we have $S/J^n\ot_R\R[[X]]
= S/J^n\ot_{R/I^n}(R/I^n\ot_R\R[[X]]) = S/J^n\ot_{R/I^n}R/I^n[X]
= S/J^n[X]$ and $\varprojlim_n S/J^n\ot_R\R[[X]] =
\varprojlim_n S/J^n[X] = \S[[X]]$, as desired.
 The remaining assertions follow by the way
of~\cite[Lemma~B.9.2]{Pweak} or~\cite[Corollary~10.3(a)]{Pcta}
and the above Corollary~\ref{contramod-veryflat-cor}(e)
(see also Lemma~\ref{projlim-reduction-lemma} below).
\end{proof}

\begin{lem}
\textup{(a)} If the map $\bar f\:R/I\rarrow S/J$ is surjective, then
the functor $E^\phi$ takes $R$\+contraadjusted $(R,I)$\+contramodules
to $S$\+contraadjusted $(S,J)$\+contra\-modules. \par
\textup{(b)} If the morphism $\bar f\:R/I\rarrow S/J$ is finite,
then the functor $E^\phi$ takes $R$\+flat $R$\+cotorsion
$(R,I)$\+contramodules to $S$\+flat $S$\+cotorsion
$(S,J)$\+contramodules.
\end{lem}

\begin{proof}
 Part~(a): by Corollary~\ref{r-i-contra-cor}(b)
or Corollary~\ref{contramod-veryflat-cor}(b), any
$R$\+contraadjusted $(R,I)$\+contramodule is a quotient
$(R,I)$\+contramodule of an $R$\+flat $R$\+contraadjusted
$(R,I)$\+contramodule.
 The functor $E^\phi$ being right exact and the class of contraadjusted
$S$\+modules being closed under quotents, it suffices to show that
the $S$\+module $E^\phi(\P)$ is contraadjusted for any $R$\+flat
$R$\+contraadjusted $(R,I)$\+contramodule~$\P$.
 Then, by Lemma~\ref{contramod-contraexten}, one has
$E^\phi(\P) = \varprojlim_n E^\phi(\P)/J^nE^\phi(\P)$ and
$E^\phi(\P)/JE^\phi(\P) \allowbreak= S/J\ot_R\P$, so it remains
to apply Corollary~\ref{contraadjusted-reduction}(a) together with
Lemma~\ref{quotient-scalars-contraadjusted}(b).
 Part~(b) follows directly from Lemma~\ref{contramod-contraexten},
Corollary~\ref{contraadjusted-reduction}(b), and 
Lemma~\ref{quotient-scalars-cotorsion}(a) in the similar way.
\end{proof}

 We keep assuming that $f(I)\sub J$.
 Let $g\:R\rarrow T$ be another ring homomorphism and $K=Tg(I)
\sub T$ be the extension of the ideal $I\sub R$ in the ring~$T$.
 Suppose that the commutative ring $H=S\ot_RT$ is Noetherian,
denote by $f'\:T\rarrow H$ and $g'\:S\rarrow H$ the related
ring homomophisms, and set $L=Hg'(J)\sub H$.

 Let $\gT$ and $\gH$ denote the adic completions of the rings
$T$ and $H$ with respect to the ideals $K$ and $L$, and let
$\psi\:\R\rarrow\gT$, \ $\phi'\:\gT\rarrow\gH$, and $\psi'\:
\S\rarrow\gH$ be the induced homomorphisms of topological rings.

\begin{lem}
 \textup{(a)} For any $(T,K)$\+contramodule\/ $\gN$ there is
a natural isomorphism of $(S,J)$\+contramodules
$E^\phi R^\psi(\gN)\simeq R^{\psi'}E^{\phi'}(\gN)$. \par
 \textup{(b)} Assume that the ring homomorphism $g\:R\rarrow T$
induces an open embedding of affine schemes\/ $\Spec T\rarrow
\Spec R$, while the morphism $\bar f\:R/I\rarrow S/J$ is finite.
 Then for any $R$\+flat $R$\+contraadjusted $(R,I)$\+contramodule\/
$\P$ there is a natural isomorphism of $(H,L)$\+contramodules
$E^{\phi'}(\Hom_R(T,\P))\simeq\Hom_S(H,E^\phi(\P))$.
\end{lem}

\begin{proof}
 Part~(a): the functor of contraextension of scalars $E^\phi\:
(R,I)\contra\rarrow(S,J)\contra$ is left adjoint to the functor of
(contra)restriction of scalars $R^\phi\:(S,J)\allowbreak\contra\rarrow
(R,I)\contra$, while the functor $R^\psi\:(T,K)\contra\rarrow
(R,I)\contra$ is left adjoint to the functor of coextension of scalars
taking an $(R,I)$\+contramodule $\P$ to the $(T,K)$\+contramodule
$\Hom_R(T,\P)$.
 To obtain the desired isomorphism of functors, one can start with
the obvious functorial isomorphism of $(T,K)$\+contramodules
$\Hom_R(T,R^\phi(\Q))\simeq R^{\phi'}\Hom_S(H,\Q)$ for any
$(S,J)$\+contramodule $\Q$, and then pass to the left adjoint functors.

 Part~(b): the $T$\+module $\Hom_R(T,\P)$ being flat by
Corollary~\ref{coherent-flat-local}(a), one has
$E^{\phi'}(\Hom_R(T,\P))\simeq\varprojlim_n\. H/L^n\ot_T\Hom_R(T,\P)
\simeq\varprojlim_n\. S/J^n\ot_R\Hom_R(T,\P)$ and
$\Hom_S(H,E^\phi(\P))\simeq\Hom_R(T,E^\phi(\P))\simeq
\Hom_R(T\;\varprojlim_n\. S/J^n\ot_R\P)\simeq
\varprojlim_n\Hom_R(T\;S/J^n\ot_R\P)$.
 Finally, one has $S/J^n\ot_R\Hom_R(T,\P)\simeq\Hom_R(T\;S/J^n
\allowbreak\ot_R\P)$ by Corollary~\ref{coherent-very}(c), since
the $R$\+modules $S/J^n$ are finitely generated.
\end{proof}

\begin{lem}
 Let $R\rarrow S_\alpha$ be a collection of morphisms of commutative
rings for which the corresponding collection of morphisms of affine
schemes\/ $\Spec S_\alpha\rarrow\Spec R$ is a finite open covering.
 Let $I$ be an ideal in $R$ and $J_\alpha$ be its extensions
in~$S_\alpha$.
 Then a contraadjusted $R$\+module $P$ is an $(R,I)$\+contramodule
if and only if the $S_\alpha$\+modules\/ $\Hom_R(S_\alpha,P)$ are
$(S_\alpha,J_\alpha)$\+contramodules for all~$\alpha$.
\end{lem}

\begin{proof}
 The ``only if'' is a particular case of
Lemma~\ref{contramod-exten}(a).
 To prove the ``if'', one can use the \v Cech
sequence~\eqref{cech-contra} together with
Lemma~\ref{contramod-restrict}(a) and the facts that the class of
$(R,I)$\+contramodules is preserved by the functors $\Hom_R$ from
any $R$\+module as well as the passages to the cokernels of $R$\+module
morphisms (and, actually, kernels and extensions, too).
\end{proof}

\subsection{Noncommutative Noetherian rings}
\label{noncomm-noetherian-subsect}
 The aim of this section is to generalize the main results
of~\cite[Appendix~B]{Pweak} and the above
Section~\ref{torsion-subsect} to noncommutative rings.
 Our exposition is somewhat sketchy with details of
the arguments omitted when they are essentially the same
as in the commutative case.
 For a further generalization, see
Sections~\ref{flat-contramod-subsect}\+-%
\ref{ind-affine-co-contra-subsect} below.

 Let $R$ be a right Noetherian associative ring, and let $m\sub R$
be an ideal generated by central elements in~$R$.
 Denote by $\R=\varprojlim_n R/m^n$ the $m$\+adic completion of
the ring $R$, viewed as a topological ring in the projective limit
($=$~$m$\+adic) topology.
 We refer to~\cite[Remark~A.3]{Psemi}, \cite[Section~1.2]{Pweak},
\cite[Sections~2.1\+-2.2]{Prev}, \cite[Sections~1.2 and~5]{PR},
\cite[Sections~2.5\+-2.7]{Pproperf}, \cite[Section~2.7]{Pcoun},
\cite[Section~6]{PS1} for the definitions of the abelian category
$\R\contra$ of \emph{left\/ $\R$\+contramodules} and the forgetful
functor $\R\contra\rarrow R\modl$ (see also the beginning of
Section~\ref{flat-contramod-subsect}).

\begin{thm}  \label{noncomm-contramod-fully-faithful}
 The forgetful functor\/ $\R\contra\rarrow R\modl$ is fully faithful.
 If $s_1$,~\dots, $s_k$ are some set of central generators of
the ideal $m\sub R$, then the image of the forgetful functor
consists precisely of those left $R$\+modules $P$ for which one has\/
$\Ext^*_R(R[s_j^{-1}],P)=0$ for all\/ $1\le j\le k$.
\end{thm}

 In other words, extending the terminology of
Section~\ref{cta-cot-contramod-subsect} to the noncommutative
situation, one can say that a left $R$\+module is an $\R$\+contramodule
if and only if it is an $s_j$\+contramodule for every~$j$.
 
\begin{proof}
 The argument follows the proof of~\cite[Theorem~B.1.1]{Pweak}.
 To show that $\Ext_R^*(R[s^{-1}],\P)=0$ for any left $\R$\+contramodule
$\P$, one can simply (contra)restrict the scalars to the subring
$K\sub R$ generated by $s$ over $\boZ$ in $R$, completed adically at
the ideal $m\cap K$ or $(s)\sub K$, reducing to the commutative case,
where the quoted theorem from~\cite[Appendix~B]{Pweak} is
directly applicable (cf.~\cite[Section~B.2]{Pweak}).

 To prove the fully faithfullness and the sufficiency of the condition
on an $R$\+module $P$, consider the ring of formal power series
$\gT=R[[t_1,\dotsc,t_k]]$ with central variables~$t_j$ and endow it with
the adic topology of the ideal generated by $t_1$~\dots,~$t_k$.
 There is a natural surjective open homomorphism of topological
rings $\gT\rarrow\R$ forming a commutative diagram with the ring
homomorphisms $R\rarrow\gT$ and $R\rarrow\R$ and taking~$t_j$ to~$s_j$.
 Consider also the polynomial ring $T=R[t_1,\dotsc,t_k]$; there are
natural ring homomorphisms $R\rarrow T\rarrow\gT$.
 The argument is based on two lemmas.

\begin{lem}
 The kernel\/ $\gJ$ of the ring homomorphism\/ $\gT\rarrow\R$ is
generated by the central elements $t_j-s_j$ as an ideal in an abstract,
nontopological ring\/~$\gT$.
 Moreover, any family of elements in\/ $\gJ$ converging to zero in
the topology of\/ $\gT$ can be presented as a linear combination of
$k$~families of elements in\/ $\gT$, converging to zero in the topology
of\/ $\gT$, with the coefficients $t_j-s_j$.
\end{lem}

\begin{proof}
 The proof is similar to that in~\cite[Sections~B.3\+-B.4]{Pweak};
the only difference is that one has to use the noncommutative versions
of Hilbert basis theorem~\cite[Theorem~1.9 and Exercise~1ZA(c)]{GW}
and Artin--Rees lemma~\cite[Theorem~13.3]{GW}.
\end{proof}

\begin{lem}
 The forgetful functor\/ $\gT\contra\rarrow T\modl$ identifies
the category of left contramodules over the topological ring\/ $\gT$
with the full subcategory in the category of left\/ $T$\+modules
consisting of all those modules which are $t_j$\+contramodules
for all~$j$.
\end{lem}

\begin{proof}
 This assertion is true for any associative ring~$R$;
the argument is the same as in~\cite[Sections~B.5\+-B.6]{Pweak}.
\end{proof}

 The proof of the theorem finishes similarly
to~\cite[Section~B.7]{Pweak}.
 The category $\R\contra$ is identified with the full subcategory
in $\gT\contra$ consisting of those left $\gT$\+contramodules $\P$ in
which the elements $t_j-s_j$ act by zero.
 The latter category coincides with the category of $T$\+modules
in which the variables~$t_j$ act the same as the elements~$s_j$ and
which are also $t_j$\+contramodules for all~$j$.
\end{proof}

\begin{prop}  \label{noncomm-r-flat-contramod}
 A left\/ $\R$\+contramodule\/ $\P$ is a flat $R$\+module if and only
if the $R/m^n$\+module\/ $\P/m^n\P$ is flat for every $n\ge1$.
 The natural map\/ $\P\rarrow\varprojlim_n\P/m^n\P$ is an isomorphism
in this case.
\end{prop}

\begin{proof}
 The same as in~\cite[Lemma~B.9.2]{Pweak}
or~\cite[Corollary~10.3(a)]{Pcta}.
 One computes the functor $M\mpsto M\ot_R\P$ on the category of finitely
generated right $R$\+modules $M$ and uses (the noncommutative version
of) the Artin--Rees lemma for such modules~$M$.
\end{proof}

\begin{prop}  \label{noncomm-mod-contramod-flat-ext}
 For any flat $R$\+module $F$ such that the $R/m$\+module $F/mF$ is
projective, and any\/ $\R$\+contramodule\/ $\Q$, one has\/
$\Ext_R^{>0}(F,\Q)=0$.
\end{prop}

\begin{proof}
 This assertion, provable in the same way
as~\cite[Proposition~B.10.1]{Pweak} (see also
Corollary~\ref{projective-contramod}(b) below),
does not depend on any Noetherianity assumptions.
\end{proof}

\begin{cor}  \label{noncomm-r-proj-contramod}
\textup{(a)} A left\/ $\R$\+contramodule\/ $\gF$ is projective if
and only if it is a flat $R$\+module and the $R/m$\+module $\gF/m\gF$
is projective.  \par
\textup{(b)} The forgetful functor\/ $\R\contra\rarrow R\modl$
induces isomorphisms of all the Ext groups.
\end{cor}

\begin{proof}
 To prove the ``only if'' assertion in part~(a), it suffices to
notice the isomorphism of $R/m^n$\+modules $\R[[X]]/m^n\R[[X]]\simeq
\R/m^n[X]$, which holds for any set $X$ and any $n\ge1$, and
use Proposition~\ref{noncomm-r-flat-contramod}.
 The ``if'' follows from the fully faithfulness assertion in
Theorem~\ref{noncomm-contramod-fully-faithful} and
Proposition~\ref{noncomm-mod-contramod-flat-ext}, and
part~(b) follows from the same together with part~(a).
 (Cf.~\cite[Section~B.8]{Pweak} and
Section~\ref{flat-contramod-subsect} below.)
\end{proof}

\begin{cor} \label{noncomm-contramod-contraderived}
\textup{(a)} The contraderived category\/ $\sD^\ctr(\R\contra)$ of
the abelian category of left\/ $\R$\+contramodules is equivalent to
the contraderived category of the exact category of $R$\+flat
left\/ $\R$\+contramodules. \par
\textup{(b)} The Becker contraderived category\/ $\sD^\bctr(\R\contra)$
the abelian category of left\/ $\R$\+contramodules is equivalent to
the Becker contraderived category of the exact category of $R$\+flat
left\/ $\R$\+contramodules. \par
\textup{(c)} The Becker contraderived category\/ $\sD^\bctr(\R\contra)$
is equivalent to the homotopy category of complexes of projective
left\/ $\R$\+contramodules. \par
\textup{(d)} Assume that all flat left $R/m$\+modules have finite
projective dimensions.
 Then the classes of Positselski-contraacyclic and Becker-contraacyclic
complexes in the abelian category $\R\contra$ coincide.
 The contraderived category\/ $\sD^{\ctr=\bctr}(\R\contra)$ is
equivalent to the absolute derived category of the exact
category of $R$\+flat left\/ $\R$\+contramodules. \par
\textup{(e)} Assume that all flat left $R/m$\+modules have finite
projective dimensions.
 Then the contraderived category\/ $\sD^{\ctr=\bctr}(\R\contra)$ is
equivalent to the homotopy category of complexes of projective
left\/ $\R$\+contramodules.
\end{cor}

\begin{proof}
 The argument is similar to the proof of
Theorem~\ref{tors-contra-co-contra-derived}(b\+f).
 In particular, Propositions~\ref{infinite-resolutions}(b) and
Proposition~\ref{becker-contraderived-infinite-resolutions}, and
Theorem~\ref{coderived-of-grothendieck-contraderived-of-lpacepo}(b)
are applicable.
 Clearly, the class of flat left $R$\+modules is closed under infinite
products in our assumptions, so it remains to notice that projective
left $\R$\+contramodules are $R$\+flat and, in the assumption of
parts~(d\+e), any $R$\+flat left $\R$\+contramodule has finite
projective dimension in $\R\contra$.
 The latter assertions follow from
Corollary~\ref{noncomm-r-proj-contramod}(a).
 An analogue of Theorem~\ref{tors-contra-co-contra-derived}(a) also
holds, but it requires Noetherianity on the other side; see the proof
of Theorem~\ref{noncomm-tors-contra-correspondence} below.
 (Cf.\ Theorems~\ref{discrete-mod-coderived}\+-%
\ref{contramod-contraderived}.)
\end{proof}

 Let $K$ be a commutative ring, $I\sub K$ be a finitely generated
ideal, and $A$ and $B$ be associative $K$\+algebras.
 For any $K$\+module $L$, we will denote by ${}_IL\sub L$
the submodule of elements annihilated by $I$ in~$L$.

\begin{lem}  \label{noncomm-hom-tensor-reduction}
\textup{(a)} For any left $A$\+module $M$ and injective left
$A$\+module $J$ there is a natural isomorphism of $K/I$\+modules
$\Hom_A(M,J)/I\Hom_A(M,J)\simeq\Hom_{A/IA}({}_IM,{}_IJ)$. \par
\textup{(b)} For any right $B$\+module $N$ and flat left
$B$\+module $F$ there is a natural isomorphism of $K/I$\+modules
${}_I(N\ot_BF)\simeq {}_IN\ot_{B/IB}F/IF$.
\end{lem}

\begin{proof}
 For any finitely presented $K$\+module $E$ there are natural
isomorphisms $\Hom_A(\Hom_K(E,M),J)\simeq E\ot_K\Hom_A(M,J)$ and
$\Hom_K(E\;N\ot_BF)\simeq\Hom_K(E,\allowbreak N)\ot_BF$.
 So it remains to take $E=K/I$.
 (Cf.\ Lemmas~\ref{hom-tensor-reduction}
and~\ref{ind-affine-hom-tensor-reduction}.)
\end{proof}

 Let $D^\bu$ be a finite complex of $A$\+injective and $B$\+injective
$A$\+$B$\+bimodules over~$K$ (i.~e., it is presumed that the left
action of $A$ and the right action of $B$ restrict to one and the same
action of $K$ in~$D^\bu$).
 We are using the definition of a dualizing complex for a pair of
noncommutative rings from Section~\ref{corings-dualizing}.

\begin{lem}  \label{noncomm-dualizing-reduction}
\textup{(a)} Assume that the ring $A$ is left coherent and
the ring $B$ is right coherent.
 Then ${}_ID^\bu$ is a dualizing complex for the rings $A/IA$ and $B/IB$
whenever $D^\bu$ is  dualizing complex the rings $A$ and~$B$. \par
\textup{(b)} Assume that the ideal $I$ is nilpotent, the ring $A$
is left Noetherian, and the ring $B$ is right Noetherian.
 Then $D^\bu$ is a dualizing complex for the rings $A$ and $B$
whenever ${}_ID^\bu$ is a dualizing complex for the rings $A/IA$
and $B/IB$.
\end{lem}

\begin{proof}
 Part~(a): clearly, ${}_ID^\bu$ is a complex of injective left
$A/IA$\+modules.
 The homothety map $B\rarrow\Hom_A(D^\bu,D^\bu)$ is a quasi-isomorphism
of finite complexes of flat left $B$\+modules (see
Lemma~\ref{coherent-tensor-hom-lemma}(b)) and therefore remains
a quasi-isomorphism after taking the tensor product with $B/IB$
over $B$ on the left, i.~e., reducing modulo~$I$.
 By Lemma~\ref{noncomm-hom-tensor-reduction}(a), it follows that
the map $B/IB\rarrow\Hom_{A/IA}({}_ID^\bu\;{}_ID^\bu)$ is
a quasi-isomorphism.

 A bounded above complex of left modules over a left coherent ring
is quasi-isomorphic to a bounded above complex of finitely generated
projective modules if and only if its cohomology modules are finitely
presented.
 To show that the $A/IA$\+mod\-ules of cohomology of the complex
${}_ID^\bu$ are finitely presented, we notice that the complex of
left $A/IA$\+modules ${}_ID^\bu$ is quasi-isomorphic to the complex
$\Hom_{B^\rop}(L_\bu,D^\bu)$ of right $B$\+module homomorphisms from
a resolution $L_\bu$ of the right $B$\+module $B/IB$ by finitely
generated free $B$\+modules into the complex~$D^\bu$.
 Locally in the cohomological grading, the complex
$\Hom_{B^\rop}(L_\bu,D^\bu)$ is a finitely iterated cone of morphisms
between shifts of copies of the complex $D^\bu$, so its $A$\+modules
of cohomology are finitely presented.
 It remains to point out that an $A/IA$\+module is finitely presented
if and only if it is finitely presented as an $A$\+module.

 Part~(b): whenever the ideal $I$ is nilpotent, a morphism of finite
complexes of flat $B$\+modules is a quasi-isomorphism if and only if
it becomes one after taking the tensor products with $B/IB$ over~$B$.
 This, together with the above argument, proves that the homothety
map $B\rarrow\Hom_A(D^\bu,D^\bu)$ is a quasi-isomorphism.
 To show that the $A$\+modules of cohomology of the complex $D^\bu$
are finitely generated, one can consider the spectral sequence
$E_2^{pq}=\Ext^p_{B^\rop}(B/IB,H^qD^\bu)\Longrightarrow
H^{p+q}\Hom_{B^\rop}(B/IB,D^\bu)$ converging to the cohomology
$A/IA$\+modules $H^n({}_ID^\bu)$. 
 Then one argues by increasing induction on~$q$, using the fact
that an $A$\+module $M$ is finitely generated provided that so
is the $A/IA$\+module ${}_IM$.
 (Cf.\ Lemmas~\ref{dualizing-reduction}
and~\ref{ind-affine-dualizing-reduction}.)
\end{proof}

 Let $K$ be a commutative Noetherian ring, $I\sub K$ be an ideal,
and $A$ and $B$ be associative $K$\+algebras such that the ring $A$
is left Noetherian and the ring $B$ is right Noetherian.
 A finite complex of $A$\+$B$\+bimodules $\D^\bu$ over $K$ is said to be
a \emph{dualizing complex for $A$ and $B$ over\/~$(K,I)$} if
\begin{enumerate}
\renewcommand{\theenumi}{\roman{enumi}}
\item $\D^\bu$ is a complex of injective left $A$\+modules
and a complex of injective right $B$\+modules;
\item any element in $\D^\bu$ is annihilated by some power of
the ideal~$I$;
\item for any (or, equivalently, some) $n\ge1$, the complex
${}_{(n)}\D^\bu={}_{I^n}\D^\bu$ is a dualizing complex for the rings
$A/I^nA$ and $B/I^nB$ (cf.\ Lemma~\ref{noncomm-dualizing-reduction}).
\end{enumerate}

 A $K$\+module (or $A$\+module) is said to be \emph{$I$\+torsion} if
every its element is annihilated by some power of the ideal $I\sub K$.
 A \emph{left\/ $(B,I)$\+contramodule} is a left contramodule over
the $I$\+adic completion of the ring $B$, or, equivalently, a left
$B$\+module that is a $(K,I)$\+contramodule in the $K$\+module
structure obtained by restriction of scalars
(see Theorem~\ref{noncomm-contramod-fully-faithful} and
Section~\ref{torsion-subsect}).

 We denote the abelian category of $I$\+torsion left $A$\+modules by
$(A,I)\tors$ and the abelian category of left $(B,I)$\+contramodules
by $(B,I)\contra$.

\begin{thm} \label{noncomm-tors-contra-correspondence}
 The choice of a dualizing complex\/ $\D^\bu$ for the rings $A$ and $B$
over $(K,I)$ induces an equivalence between the coderived category
of $I$\+torsion left $A$\+modules\/ $\sD^\co((A,I)\tors)$ and
the contraderived category of left $(B,I)$\+contra\-modules\/
$\sD^\ctr((B,I)\contra)$.
 The equivalence is provided by the derived functors\/
$\boR\Hom_A(\D^\bu,{-})$ and\/ $\D^\bu\ot_B^\boL{-}$.
\end{thm}

\begin{proof}
 Notice first of all that an object $\J\in (A,I)\tors$ is injective
if and only if the $A/I^nA$ modules ${}_{(n)}\J={}_{I^n}\J$ are
injective for all~$n$ (obviously), and if and only if it is
an injective $A$\+module (by the Artin--Rees lemma for centrally
generated ideals in noncommutative Noetherian
rings~\cite[Theorem~13.3]{GW}).
 The dual result for $B$\+flat $(B,I)$\+contramodules is provided by
Proposition~\ref{noncomm-r-flat-contramod}.

 Furthermore, we identify the coderived category $\sD^\co((A,I)\tors)$
with the homotopy category of (complexes of) injective $I$\+torsion
$A$\+modules (based on the facts that there are enough injectives and
the class of injectives is closed with respect to infinite direct sums
in $(A,I)\tors$) and the contraderived category $\sD^\ctr((B,I)\contra)$
with the absolute derived category of $B$\+flat $(B,I)$\+contramodules
(based on Corollary~\ref{noncomm-contramod-contraderived}(d)
and Lemma~\ref{christensen-frankild-holm}).
 The rest of the argument is no different from the proof of
Theorem~\ref{tors-contra-correspondence} (see also
Lemma~\ref{noncomm-ring-dualizing-lemma}) and based on
Lemma~\ref{noncomm-hom-tensor-reduction}.

 Similarly to Theorem~\ref{tors-contra-correspondence}, in the context
of the present theorem one also has $\sD^\co((A,I)\tors)=
\sD^\bco((A,I)\tors)$ (by the argument in the previous paragraph) and
$\sD^\ctr((B,I)\contra)=\sD^\bctr((B,I)\contra)$
(by Corollary~\ref{noncomm-contramod-contraderived}(d)).
\end{proof}

\Section{Ind-Affine Ind-Schemes}

 The aim of this appendix is to lay some bits of preparatory groundwork
for the definition of contraherent cosheaves of contramodules on
ind-schemes.
 We construct enough very flat and contraadjusted contramodules on
an ind-affine ind-scheme represented by a sequence of affine schemes
and their closed embeddings with finitely generated defining ideals,
and also enough cotorsion contramodules on an ind-Noetherian
ind-affine ind-scheme of totally finite Krull dimension.
 For a more abstract and general approach, see the paper~\cite{PR}.
 A version of co-contra correspondence for ind-coherent ind-affine
ind-schemes with dualizing complexes (and their noncommutative
generalizations) is also worked out.

\subsection{Flat and projective contramodules}
\label{flat-contramod-subsect}
 Let $R_0\larrow R_1\larrow R_2\larrow R_3\larrow\dotsb$ be a projective
system of associative rings, indexed by the ordered set of positive
integers, with surjective morphisms between them.
 Denote by $\R$ the projective limit $\varprojlim_n R_n$, viewed as
a topological ring in the topology of projective limit of discrete
rings~$R_n$.
 Clearly, the ring homomorphisms $\R\rarrow R_n$ are surjective; let
$\gI_n\sub\R$ denote their kernels.
 Then the open ideals $\gI_n$ form a base of neighborhoods of zero
in the topological ring~$\R$.

 We are interested in left contramodules over the topological ring~$\R$
(see~\cite[Remark~A.3]{Psemi}, \cite[Section~1.2]{Pweak},
\cite[Section~2.1]{Prev}, \cite[Sections~1.2 and~5]{PR},
\cite[Sections~2.5\+-2.7]{Pproperf}, \cite[Section~2.7]{Pcoun},
or~\cite[Section~6]{PS1} for the definition).
 They form an abelian category $\R\contra$ having enough projective
objects and endowed with an exact and faithful forgetful functor
$\R\contra\rarrow\R\modl$ preserving infinite products.
 Here $\R\modl$ denotes the abelian category of left modules over
the ring $\R$ (viewed as an abstract ring without any topology).

 The projective $\R$\+contramodules are precisely the direct summands
of the free $\R$\+contramodules $\R[[X]]$.
 Here $X$ is an arbitrary set of generators, and $\R[[X]]=\varprojlim_n
R_n[X]$ is the set of all maps $X\rarrow\R$ converging to zero in
the topology of~$\R$.

 For any left $\R$\+contramodule $\P$ and any closed two-sided ideal
$\gJ\sub\R$, we will denote by $\gJ\tim\P\sub\P$ the image of
the contraaction map $\gJ[[\P]]\rarrow\P$.
 As usually, for any left module $M$ over a ring $R$ and any ideal
$J\sub R$ the notation $JM$ (or, as it may be sometimes convenient,
$J\.{\cdot}M$) stands for the image of the action map $J\ot_RM\rarrow M$.
 For an $\R$\+contramodule $\P$, there is a (generally speaking, proper)
inclusion $\gJ\P=\gJ\.{\cdot}\.\P\sub\gJ\tim\P$.
 Of course, $\gJ\tim\P$ is an $\R$\+subcontramodule in $\P$, while
$\gJ\P$ is in general only an $\R$\+submodule.

\begin{lem}  \label{contramod-projlim-surjectivity}
 For any left\/ $\R$\+contramodule\/ $\P$, the natural map to
the projective limit\/ $\P\rarrow\varprojlim_n\P/(\gI_n\tim\P)$
is surjective.
\end{lem}

\begin{proof}
 This assertion is true for any topological ring $\R$ with a base of
topology consisting of open right ideals, and any decreasing sequence
of closed abelian subgroups $\R\supset\gI_0\supset\gI_1\supset
\gI_2\supset\dotsb$ converging to zero in the topology of~$\R$
(meaning that for any neighborhood of zero $\gU\sub\R$ there exists
$n\ge0$ such that $\gI_n\sub\gU$).
 Indeed, it suffices to show that for any sequence of elements
$p_n\in\gI_n\tim\P$, \ $n\ge0$, there exists an element $p\in\P$
such that $p-(p_0+\dotsb+p_n)\in\gI_{n+1}\tim\P$ for all~$n\ge0$.

 Now each element~$p_n$ can be obtained as an infinite linear
combination of elements of $\P$ with the family of coefficients
converging to zero in~$\gI_n$.
 The countable sum of such expressions over all $n\ge0$ is again
an infinite linear combination of elements of $\P$ with
the coefficients converging to zero in~$\R$.
 The value of the latter, understood in the sense of the contramodule
infinite summation operations in $\P$, provides the desired
element $p\in\P$.
 Another (and perhaps more illuminating) argument can be found
in~\cite[Lemma~A.2.3]{Psemi} (while counterexamples showing that
the map in question may not be injective are provided
in~\cite[Section~A.1]{Psemi} or~\cite[Section~1.5]{Prev}).
\end{proof}

 The following result is a version of Nakayama's lemma for
$\R$\+contramodules (see~\cite[Lemma~1.3.1]{Pweak} for a somewhat
more familiar formulation).
 A more general version can be found in~\cite[Lemma~6.14]{PR} and
another version in~\cite[Lemma~6.2]{Pproperf}; see
also~\cite[Lemmas~2.1 and~3.22]{Prev}.

\begin{lem} \label{generalized-contramodule-nakayama}
 For any left\/ $\R$\+contramodule\/ $\P$, the equations\/
$\gI_n\tim\P=\P$ for all\/ $n\ge0$ imply\/ $\P=0$.
\end{lem}

\begin{proof}
 This assertion holds for any (complete and separated) topological ring
$\R$ with a base of topology consisting of open right ideals, and any
sequence of closed abelian subgroups $\gI_0$, $\gI_1$,
$\gI_2$,~\dots~$\sub\R$ with the property that the sequence of
subgroups $\gI_0$, $\gI_0\gI_1$, $\gI_0\gI_1\gI_2$,~\dots\ converges to
zero in the topology of $\R$ (i.~e., for any neighborhood of zero
$\gU\sub\R$ there exists $n\ge0$ such that $\gI_0\dotsm\gI_n\sub\gU)$.
 Indeed, assume that the restrictions $\pi_n\:\gI_n[[\P]]\rarrow\P$
of the contraaction map $\pi\:\R[[\P]]\rarrow\P$ are surjective
for all $n\ge0$.
 Let $p\in\P$ be an element.

 Notice that for any surjective map of sets $f\:X\rarrow Y$ and any
$n\ge0$, the induced map $\gI_n[[f]]\:\gI_n[[X]]\rarrow\gI_n[[Y]]$ is
also surjective.
 Given a set $X$, define inductively the sets $\gI^{(-1)}[[X]]=X$
and $\gI^{(n)}[[X]]=\gI^{(n-1)}[[\gI_n[[X]]]]$ for all $n\ge0$.
 Let $p_0\in\gI_0[[\P]]$ be a preimage of $p\in\P$ under the map
$\pi_0\:\gI_0[[\P]]\rarrow\P$.
 Furthermore, let $p_n\in\gI^{(n)}[[\P]]$ be a preimage of
$p_{n-1}$ under the map $\gI^{(n-1)}[[\pi_n]]\:\gI^{(n)}[[\P]]
\rarrow\gI^{(n-1)}[[\P]]$.

 For any set $X$, the abelian group $\R[[X]]$ is complete in its
natural topology with the base of neighborhoods of zero formed by
the subgroups $\gU[[X]]$, where $\gU\sub\R$ are open right ideals.
 Besides, the ``opening of parentheses''/monad multiplication map
$\rho_X\:\R[[\R[[X]]]]\rarrow \R[[X]]$ is continuous, as is the map
$\R[[f]]\:\R[[X]]\rarrow\R[[Y]]$ induced by any map of sets
$f\:X\rarrow Y$.
 For every $n\ge0$, let $\rho^{(n)}_X\:\gI^{(n)}[[X]]\rarrow
\R[[X]]$ denote (the restriction of) the iterated monad
multiplication map.

 Set $q_n=\rho^{(n-1)}_{\gI_n[[\P]]}(p_n)\in\R[[\gI_n[[\P]]]]\subset
\R[[\R[[\P]]]]$ for all $n\ge1$.
 Due to our convergence condition on the products of the subgroups
$\gI_n\sub\R$, the sum $\sum_{n=1}^\infty q_n$ converges in the topology
of $\R[[\R[[\P]]]]$.
 Now we have $\R[[\pi_n]](q_n)=\rho_\P(q_{n-1})$ for all $n\ge2$ and
$\gI_0[[\pi_1]](q_1)=p_0$.
 Hence
$$\textstyle
 \R[[\pi]]\big(\sum_{n=1}^\infty q_n)-\rho_\P(\sum_{n=1}^\infty q_n)
 = p_0
$$
in $\R[[\P]]$ and $p=\pi(p_0)=0$ by the contraassociativity equation.
\end{proof}

\begin{lem}  \label{projlim-reduction-lemma}
 Let $P_0\larrow P_1\larrow P_2\larrow\dotsb$ be a projective system of
left $R_n$\+modules in which the morphism $P_{n+1}\rarrow P_n$
identifies $P_n$ with $R_n\ot_{R_{n+1}}P_{n+1}$ for every $n\ge0$.
 Let\/ $\P$ denote the\/ $\R$\+contramodule\/ $\varprojlim_n P_n$.
 Then the natural map\/ $\P\rarrow P_n$ identifies $P_n$ with\/
$\P/\gI_n\tim\P$.
 Conversely, for any left\/ $\R$\+contramodule $\P$ the projective
system $P_n=\P/\gI_n\tim\P$ satisfies the condition above.
\end{lem}

\begin{proof}
 Since $P_n$ is an $R_n$\+module and $\P\rarrow P_n$ is a morphism of
$\R$\+contramodules, the kernel of this morphism contains
$\gI_n\tim\P$.
 To prove the inverse inclusion, consider an element $p\in\P$
belonging to the kernel of the morphism $\P\rarrow P_n$.
 The image of the element~$p$ in $P_{n+1}$ belongs to
$(\gI_n/\gI_{n+1})P_{n+1}$ (where $\gI_n/\gI_{n+1}$ is an ideal
in~$R_{n+1}$).
 Let us decompose this image accordingy into a finite linear
combination of elements of $P_{n+1}$ with coefficients from
$\gI_n/\gI_{n+1}$, lift all the entering elements of $P_{n+1}$ to $\P$
and all the elements of $\gI_n/\gI_{n+1}$ to $\gI_n$, and subtract
from~$p$ the corresponding finite linear combination of elements
of $\P$ with coefficients in~$\gI_n$.

 The image of the resulting element $p'\in\P$ in $P_{n+1}$ vanishes,
so its image in $P_{n+2}$ belongs to $(\gI_{n+1}/\gI_{n+2})P_{n+2}$.
 Continuing this process indefinitely, we obtain an expression of
the original element~$p$ in the form of a countable linear
combination of elements from $\P$ with the coefficient sequence
converging to zero in~$\gI_n$.
 This proves the first assertion; the second one is straightforward.
\end{proof}

 A left $\R$\+contramodule $\gF$ is called \emph{flat} if
the map $\gF\rarrow\varprojlim_n\gF/\gI_n\tim\gF$ is an isomorphism
and the $R_n$\+modules $\gF/\gI_n\tim\gF$ are flat for all $n\ge0$.
 We will see below in Corollary~\ref{flat-reductions-enough} that
the former condition, in fact, follows from the latter one.

\begin{lem}  \label{contramod-kernel-flat}
 If\/ $\gG\rarrow\gF$ is a surjective morphism of flat\/
$\R$\+contramodules then its kernel\/ $\gH$ is also a flat\/
$\R$\+contramodule.
 Moreover, the sequences of $R_n$\+modules\/ $0\rarrow
\gH/\gI_n\tim\gH\rarrow\gG/\gI_n\tim\gG\rarrow
\gF/\gI_n\tim\gF\rarrow0$ are exact.
\end{lem}

\begin{proof}
 Clearly, for any short exact sequence of $\R$\+contramodules
$0\rarrow\gH\rarrow\gG\rarrow\gF\rarrow0$ there are short exact
sequences of $R_n$\+modules $0\rarrow\gH/(\gH\cap(\gI_n\tim\gG))
\rarrow\gG/\gI_n\tim\gG\rarrow\gF/\gI_n\tim\gF\rarrow0$ (because
the map $\gI_n\tim\gG\rarrow\gI_n\tim\gF$ induced by
a surjective map $\gG\rarrow\gF$ is surjective).
 If the $R_n$\+module $\gF/\gI_n\tim\gF$ is flat, then
the tensor product with $R_{n-1}$ over $R_n$ transforms this
sequence into the similar sequence corresponding to the ideal
$\gI_{n-1}\sub\R$.

 On the other hand, if the maps $\gF\rarrow\varprojlim_n\gF/
\gI_n\tim\gF$ and $\gG\rarrow\varprojlim_n\gG/\gI_n\tim\gG$ are
isomorphisms, then the passage to the projective limit of the short
exact sequences above allows to conclude that the map $\gH\rarrow
\varprojlim_n\gH/(\gH\cap(\gI_n\tim\gG))$ is an isomorphism.
 By Lemma~\ref{projlim-reduction-lemma}, it follows from these
observations that $\gH\cap(\gI_n\tim\gG)=\gI_n\tim\gH$,
so the sequences $0\rarrow\gH/\gI_n\tim\gH\rarrow
\gG/\gI_n\tim\gG\rarrow\gF/\gI_n\tim\gF\rarrow0$ are exact.
 Finally, now if the $R_n$\+modules $\gG/\gI_n\tim\gG$ are also flat,
then so are the $R_n$\+modules $\gH/\gI_n\tim\gH$.
\end{proof}

\begin{lem}  \label{contramod-extension-flat}
 If\/ $0\rarrow\gH\rarrow\gG\rarrow\gF\rarrow0$ is a short exact
sequence of\/ $\R$\+contramodules and the\/ $\R$\+contramodules\/ $\gH$
and\/ $\gF$ are flat, then so is the\/ $\R$\+contramodule\/~$\gG$.
\end{lem}

\begin{proof}
 In view of the proof of Lemma~\ref{contramod-kernel-flat}, we only
need to show that the map $\gG\rarrow\varprojlim_n\gG/\gI_n\tim\gG$
is an isomorphism.
 Choose a termwise surjective map onto the short exact sequence
$0\rarrow\gH\rarrow\gG\rarrow\gF\rarrow0$ from a short exact
sequence of free $\R$\+contramodules $0\rarrow\gW\rarrow\gV\rarrow
\gU\rarrow0$ (e.~g., $\gW=\R[[\gH]]$, \ $\gU=\R[[\gF]]$ or
$\R[[\gG]]$, and $\gV=\gW\oplus\gU$).
 Let $0\rarrow\gM\rarrow\gL\rarrow\gK\rarrow0$ be the corresponding
short exact sequence of kernels.

 Passing to the projective limit of short exact sequences $0\rarrow
\gL/\gL\cap(\gI_n\tim\gV)\rarrow\gV/\gI_n\tim\gV\rarrow\gG/\gI_n\tim\gG
\rarrow0$, we obtain a short exact sequence $0\rarrow\varprojlim_n
\gL/\gL\cap(\gI_n\tim\gV)\rarrow \varprojlim_n\gV/\gI_n\tim\gV\rarrow
\varprojlim\gG/\gI_n\tim\gG\rarrow0$.
 The $\R$\+contramodule $\gV$ being free, the intersection
$\bigcap_n\gI_n\tim\gV\sub\gV$ vanishes; so both the maps $\gV\rarrow
\varprojlim_n\gV/\gI_n\tim\gV$ and $\gL\rarrow\varprojlim_n\gL/
\gI_n\tim\gL$ are isomorphisms, while the map $\gL\rarrow
\varprojlim_n\gL/\gL\cap(\gI_n\tim\gV)$ is, at least, injective.
 It remains to show that the latter map is surjective; equivalently,
it means vanishing of the derived projective limit
$\varprojlim^1_n\gL\cap(\gI_n\tim\gV)$.

 Similarly, the assumptions that the maps $\gH\rarrow\varprojlim_n
\gH/\gI_n\tim\gH$ and $\gF\rarrow\varprojlim_n\gF/\gI_n\tim\gF$ are
isomorphisms are equivalently expressed as the vanishing of
derived projective limits $\varprojlim^1_n\gM\cap(\gI_n\tim\gW)$
and $\varprojlim^1_n\gK\cap(\gI_n\tim\gU)$.

 The short sequences $0\rarrow\gI_n\tim\gW\rarrow\gI_n\tim\gV\rarrow
\gI_n\tim\gU\rarrow0$ are exact, because the short sequence
$0\rarrow\gW\rarrow\gV\rarrow\gU\rarrow0$ splits.
 Passing to the fibered product of two short exact sequences
$0\rarrow\gI_n\tim\gW\rarrow\gI_n\tim\gV\rarrow\gI_n\tim\gU\rarrow0$
and $0\rarrow\gM\rarrow\gL\rarrow\gK\rarrow0$ over the short exact
sequence $0\rarrow\gW\rarrow\gV\rarrow\gU\rarrow0$ (into which both
of them are embedded), we obtain an exact sequence
$0\rarrow\gM\cap(\gI_n\tim\gW)\rarrow\gL\cap(\gI_n\tim\gV)\rarrow
\gK\cap(\gI_n\tim\gU)$.

 The $\R$\+contramodules $\gU$ and $\gF$ being flat, by
Lemma~\ref{contramod-kernel-flat} we have
$\gK\cap(\gI_n\tim\gU)=\gI_n\tim\gK$.
 It follows that the map $\gL\cap(\gI_n\tim\gV)\rarrow
\gK\cap(\gI_n\tim\gU)$ is surjective, so the whole sequence
$0\rarrow\gM\cap(\gI_n\tim\gW)\rarrow\gL\cap(\gI_n\tim\gV)\rarrow
\gK\cap(\gI_n\tim\gU)\rarrow0$ is exact.
 Applying the right exact functor $\varinjlim^1_n$ (which is, in
particular, exact in the middle), we obtain the desired vanishing.
\end{proof}

\begin{lem}  \label{contramod-proj-flat}
 Any projective\/ $\R$\+contramodule is flat.
\end{lem}

\begin{proof} \emergencystretch=1em
 It suffices to consider the case of a free left\/ $\R$\+contramodule
$\R[[X]]$.
 For any set $X$ and any closed ideal $\gJ\sub\R$ one has $\gJ
\tim(\R[[X]]) = \gJ[[X]]$ and $\R[[X]]/\gJ[[X]]\simeq(\R/\gJ)[X]$.
 So in particular $\R[[X]]/(\gI_n\tim\R[[X]])=R_n[X]$ is a flat
(and even free and projective) left $R_n$\+module and the natural
$\R$\+contramodule morphism $\R[[X]]\rarrow\varprojlim_n
\R[[X]]/(\gI_n\tim\R[[X]])$ is an isomorphism.
\end{proof}

\begin{cor}  \label{flat-reductions-enough}
 A left\/ $\R$\+contramodule\/ $\gF$ is flat if and only if
the $R_n$\+modules $\gF/\gI_n\tim\gF$ are flat for all $n\ge0$.
\end{cor}

\begin{proof}
 Set $\overline\P_n=\P/\gI_n\tim\P$ for any left $\R$\+contramodule
$\P$, and denote by $\P\longmapsto\boL_i\overline\P_n$ the left
derived functors of the right exact functor
$\P\longmapsto\overline\P_n$ constructed in terms of projective
resolutions of left $\R$\+contramodules~$\P$.
 It follows from Lemmas~\ref{contramod-kernel-flat}
and~\ref{contramod-proj-flat} that flat $\R$\+contramodules are
acyclic for the derived functors $\P\longmapsto\boL_*\overline\P_n$,
i.~e., one has $\boL_{>0}\overline\gF_n=0$ for all flat left
$\R$\+contramodules~$\gF$.

 Now let $\gG$ be a left $\R$\+contramodule such that all
the $R_n$\+modules $\overline\gG_n$ are flat.
 Set $\gH=\bigcap_n(\gI_n\tim\gG)$ and $\gF=\varprojlim_n
\overline\gG_n=\gG/\gH$;
according to Lemma~\ref{projlim-reduction-lemma}, one has
$\overline\gG_n\simeq\overline\gF_n$.
 Since the $\R$\+contramodule $\gF$ is flat, the long exact
sequence of derived functors $\P\longmapsto\boL_*\overline\P_n$
applied to the short exact sequence of left $\R$\+contramodules
$0\rarrow\gH\rarrow\gG\rarrow\gF\rarrow0$ proves that the short
sequence $0\rarrow\overline\gH_n\rarrow\overline\gG_n\rarrow
\overline\gF_n\rarrow0$ is exact.
 We have shown that $\overline\gH_n=0$ for all~$n\ge0$;
by Lemma~\ref{generalized-contramodule-nakayama}, it follows
that $\gH=0$.  (Cf.~\cite[Lemma~B.9.2]{Pweak}.)
\end{proof}

 Given two $\R$\+contramodules $\P$ and $\Q$, we will denote by
$\Hom^\R(\P,\Q)$ and $\Ext^{\R,*}(\P,\Q)$ the $\Hom$ and $\Ext$
groups in the abelian category $\R\contra$.

\begin{lem}  \label{contramod-contramod-flat-ext-lemma}
 Let\/ $\gF$ be a flat left\/ $\R$\+contramodule and\/ $\P$ be
a left\/ $\R$\+contra\-module for which the natural map\/
$\P\rarrow\varprojlim_n\P/\gI_n\tim\P$ is an isomorphism.
 Set $F_n=\gF/\gI_n\tim\gF$ and $P_n=\P/\gI_n\tim\P$, and assume
that one has\/ $\Ext^1_{R_0}(F_0,P_0)=0=\Ext^1_{R_{n+1}}(F_{n+1}\;
\allowbreak\ker(P_{n+1}\to P_n))$ for all\/ $n\ge0$.
 Then\/ $\Ext^{\R,1}(\gF,\P)=0$.
\end{lem}

\begin{proof}
 Let $\gG\rarrow\gF$ be a surjective morphism onto $\gF$ from
a projective $\R$\+contramod\-ule $\gG$ with the kernel~$\gH$.
 Set $G_n=\gG/\gI_n\tim\gG$ and similarly for~$H_n$.
 By Lemmas~\ref{contramod-proj-flat} and~\ref{contramod-kernel-flat},
the $\R$\+contramodules $\gG$ and $\gH$ are flat, and the short
sequences $0\rarrow H_n\rarrow G_n\rarrow F_n\rarrow0$ are exact.
 Let us show that any $\R$\+contramodule morphism $\gH\rarrow\P$
can be extended to an $\R$\+contramodule morphism $\gG\rarrow\P$.

 The datum of a morphism of $\R$\+contramodules $\gH\rarrow\P$ is
equivalent to that of a morphism of projective systems of
$R_n$\+modules $H_n\rarrow P_n$.
 Let us construct by induction an extension of this morphism to
a morphism of projective systems $G_n\rarrow P_n$ for $n\ge0$.
 Since $\Ext_{R_0}^1(F_0,P_0)=0$, the case of $n=0$ is clear.
 Assuming that the morphism $G_n\rarrow P_n$ has been obtained already,
we will proceed to construct a compatible morphism
$G_{n+1}\rarrow P_{n+1}$.

 Since $G_{n+1}$ is a projective $R_{n+1}$\+module, the composition
$G_{n+1}\rarrow G_n\rarrow P_n$ can be lifted to an $R_{n+1}$\+module
morphism $G_{n+1}\rarrow P_{n+1}$.
 (Notice that what is actually used here is the vanishing of
$\Ext^1_{R_{n+1}}(G_{n+1}\;\ker(P_{n+1}\to P_n))$.)
 The composition of a morphism so obtained with the embedding
$H_{n+1}\rarrow G_{n+1}$ differs from the given map $H_{n+1}\rarrow
P_{n+1}$ by an $R_{n+1}$\+module morphism $H_{n+1}\rarrow
\ker(P_{n+1}\to P_n)$.
 Given that $\Ext^1_{R_{n+1}}(F_{n+1}\;\ker(P_{n+1}\to P_n))=0$,
the latter map can be extended from $H_{n+1}$ to $G_{n+1}$ and added to
the previously constructed map $G_{n+1}\rarrow P_{n+1}$.
\end{proof}

 Suppose $R$ is an associative ring endowed with a descending sequence
of two-sided ideals $R\supset I_0\supset I_1\supset I_2\supset\dotsb$
such that the projective system of quotient rings $R/I_n$ is
isomorphic to our original projective system $R_0\larrow R_1\larrow R_2
\larrow\dotsb$.
 For example, one can always take $R=\R$ and $I_n=\gI_n$; sometimes
there may be other suitable choices of a ring $R$ with the ideals
$I_n$ as well.

 Then there is a natural ring homomorphism $R\rarrow\R$ inducing
the isomorphisms $R/I_n\simeq\R/\gI_n$.
 The restriction of scalars provides a forgetful functor
$\R\contra\rarrow R\modl$, which is exact and faithful, and
preserves infinite products.
  Given two $R$\+modules $L$ and $M$, we denote, as usually, by
$\Hom_R(L,M)$ and $\Ext_R^*(L,M)$ the Hom and Ext groups in
the abelian category of $R$\+modules.

\begin{lem}  \label{mod-contramod-flat-ext-lemma}
 Let $F$ be a flat left $R$\+module and\/ $\P$ be a left\/
$\R$\+contramodule for which the natural map\/
$\P\rarrow\varprojlim_n\P/\gI_n\tim\P$ is an isomorphism.
 Set $P_n=\P/\gI_n\tim\P$, and assume that that one has\/
$\Ext^1_{R_0}(F/I_0F\;P_0)=0=\Ext^1_{R_{n+1}}(F/I_{n+1}F\;\allowbreak
\ker(P_{n+1}\to\nobreak P_n))$ for all\/ $n\ge0$.
 Then\/ $\Ext^1_R(F,\P)=0$.

 If, moreover, $\Ext^i_{R_0}(F/I_0F\;P_0)=0=
\Ext^i_{R_{n+1}}(F/I_{n+1}F\;\allowbreak
\ker(P_{n+1}\to\nobreak P_n))$ for all\/ $n\ge0$ and $i>0$, then\/
$\Ext^i_R(F,\P)=0$ for all $i>0$.
\end{lem}

\begin{proof}
 The proof is similar to that of
Lemma~\ref{vfl-cta-reduction-orthogonality}.
 The $R$\+module $F$ being flat by assumption, one has
$\Ext^i_R(F,P_0)\simeq\Ext^i_{R_0}(F/I_0F\;P_0)$
and $\Ext^i_R(F\;\allowbreak\ker(P_{n+1}\to\nobreak P_n))\simeq
\Ext^1_{R_{n+1}}(F/I_{n+1}F\;\allowbreak\ker(P_{n+1}\to\nobreak P_n))$
for all $n\ge0$ and $i\ge0$.
 It remains to refer to the dual Eklof lemma~\cite[Proposition~18]{ET},
\cite[Lemma~6.37]{GT} in order to prove the first assertion.
 For the second assertion, one can refer to~\cite[Lemma~B.10.3]{Pweak}.
\end{proof}

\begin{cor}  \label{projective-contramod}
\textup{(a)} A left\/ $\R$\+contramodule\/ $\gF$ is projective
if and only if the $R_n$\+modules\/ $\gF/\gI_n\tim\gF$
are projective for all\/ $n\ge0$. \par
\textup{(b)} For any flat left $R$\+module $F$ such that
the $R_n$\+modules $F/I_nF$ are projective for all\/ $n\ge0$,
and any left\/ $\R$\+contramodule\/ $\Q$, one has\/
$\Hom^\R(\varprojlim_n F/I_nF\;\Q)\simeq\Hom_R(F,\Q)$ and\/
$\Ext^{\R,>0}(\varprojlim_n F/I_nF\;\Q)=\Ext_R^{>0}(F,\Q)=0$.
\end{cor}

\begin{proof}
 Part~(a): the ``only if'' assertion holds by (the proof of)
Lemma~\ref{contramod-proj-flat}.
 To prove the ``if'', notice first of all that the left
$\R$\+contramodule $\gF$ is flat by
Corollary~\ref{flat-reductions-enough}.
 Consider a short exact sequence of
$\R$\+contramodules $0\rarrow\gH\rarrow\gG\rarrow\gF\rarrow0$
with a projective $\R$\+contramodule~$\gG$.
 The natural map $\gH\rarrow\varprojlim_n\gH/\gI_n\tim\gH$ being
an isomorphism because the map $\gG\rarrow\varprojlim_n\gG/
\gI_n\tim\gG$ is (or by Lemma~\ref{contramod-kernel-flat}),
one has $\Ext^{\R,1}(\gF,\gH)=0$ by
Lemma~\ref{contramod-contramod-flat-ext-lemma}.
 Hence the short exact sequence splits and the $\R$\+contramodule
$\gF$ is a direct summand of the $\R$\+contramodule~$\gG$.

 Part~(b): a natural map $\Hom^\R(\varprojlim_nF/I_nF\;\Q)\rarrow
\Hom_R(F,\Q)$ for any $R$\+module $F$ and $\R$\+contramodule $\Q$
is induced by the $R$\+module morphism $F\rarrow\varprojlim_n F/I_nF$.
 For an $\R$\+contramodule $\P$ such that the map $\P\rarrow
\varprojlim_n\P/\gI_n\tim\P$ is an isomorphism, one has
$\Hom_R(F,\P)\simeq\varprojlim_n\Hom_{R_n}(F/I_nF\;\P/\gI_n\tim\P)$,
which is isomorphic to $\Hom^\R(\varprojlim_nF/I_nF\;\P)$
by Lemma~\ref{projlim-reduction-lemma}.

 When $F$ is also a flat $R$\+module with projective $R_n$\+modules
$F/I_nF$, one has $\Ext^{>0}_R(F,\P)=0$ by
Lemma~\ref{mod-contramod-flat-ext-lemma} and
$\Ext^{\R,>0}(\varprojlim_n F/I_nF\;\P)=0$ by part~(a).
 Now to prove the assertion of part~(b) in the general case, it
suffices to present an $\R$\+contramodule $\Q$ as the cokernel
of an injective morphism of $\R$\+contramodules $\gK\rarrow\P$
with $\P=\varprojlim_n\P/\gI_n\tim\P$ (and, consequently,
the same for~$\gK$).
\end{proof}

 When the ideals\/ $\gI_n/\gI_{n+1}\sub R_{n+1}$ are nilpotent
(i.~e., for every $n\ge0$ there exists $N_n\ge1$ such that
$(\gI_n/\gI_{n+1})^{N_n}=0$), it follows from
Corollary~\ref{projective-contramod}(a) and~\cite[Lemma~B.10.2]{Pweak}
that a left $\R$\+contramodule $\gF$ is projective if and only if
it is flat and the $R_0$\+module $\gF/\gI_0\tim\gF$ is projective.
 Similarly, it suffices to require that the $R_0$\+module $F/I_0F$
be projective in Corollary~\ref{projective-contramod}(b)
in this case.

\subsection{Co-contra correspondence}
\label{ind-affine-co-contra-subsect}
 In this section we consider a pair of projective systems $R_0\larrow
R_1\larrow R_2\larrow\dotsb$ and $S_0\larrow S_1\larrow S_2\larrow\dotsb$
of associative rings and surjective morphisms between them.
 Set $\R=\varprojlim_n R_n$ and $\S=\varprojlim_n S_n$, and denote by
$\gI_n\sub\R$ and $\gJ_n\sub\S$ the kernels of the natural surjective
ring homomorphisms $\R\rarrow R_n$ and $\S\rarrow S_n$.
 We assume that the rings $S_n$ are left coherent, the rings $R_n$
are right coherent, the kernels $\gJ_n/\gJ_{n+1}\sub S_{n+1}$ of
the ring homomorphisms $S_{n+1}\rarrow S_n$ are finitely generated as
left ideals, and the kernels $\gI_n/\gI_{n+1}\sub R_{n+1}$ of
the ring homomorphisms $R_{n+1}\rarrow R_n$ are finitely generated
as right ideals (cf.\ the discussion of reasonable ind-schemes and
reasonable topological rings in~\cite[Section~7.11.1]{BD2}
and~\cite[Sections~2.1 and~2.4(5)]{Psemten}).

 Given a left $\S$\+module $M$, we denote by ${}_{S_n}\!M\sub M$ its
submodule consisting of all the elements annihilated by~$\gJ_n$;
so ${}_{S_n}\!M$ is the maximal left $S_n$\+submodule in~$M$.
 A left $\S$\+module $\M$ is said to be \emph{discrete} if every its
element is annihilated by the ideal $\gJ_n$ for $n$~large enough.
 In other words, a left $\S$\+module $\M$ is discrete if its increasing
filtration ${}_{S_0}\M\sub{}_{S_1}\M\sub{}_{S_2}\M\sub\dotsb$ is
exhaustive, or equivalently, if the left action map $\S\times\M\rarrow
\M$ is continuous in (the projective limit topology of $\S$ and)
the discrete topology of~$\M$.
 We denote the full abelian subcategory of discrete left $\S$\+modules
by $\S\discr\sub\S\modl$.
 Clearly, a discrete left $\S$\+module $\J$ is an injective object
in $\S\discr$ if and only if all the left $S_n$\+modules
${}_{S_n}\J$ are injective.

 A discrete left $\S$\+module $\M$ is called \emph{finitely presented}
if there exists $n\ge0$ such that $\M={}_{S_n}\M$ and $\M$ is a finitely
presented left $S_n$\+module.
 It follows from the finite generatedness condition on the kernel of
the ring homomorphism $S_{n+1}\rarrow S_n$ that any finitely presented
left $S_n$\+module is at the same time a finitely presented left
$S_{n+1}$\+module.
 Furthermore, the functor $\M\mpsto{}_{S_n}\M\:\S\discr\rarrow S_n\modl$
commutes with filtered inductive limits.
 A discrete left $\S$\+module $\M$ is finitely presented if and only if
the functor of discrete $\S$\+module homomorphisms $\Hom_\S(\M,{-})$
from $\M$ preserves filtered inductive limits in the abelian
category of discrete left $\S$\+modules.
 The full subcategory of finitely presented discrete left $\S$\+modules
is closed under kernels, cokernels, and extensions in $\S\discr$; so it
is an abelian category.
 It follows that $\S\discr$ is a locally coherent Grothendieck category
in the sense of~\cite[Section~2]{Roos0}, \cite[Section~13]{PS3},
\cite[Section~9.5]{Pedg}, \cite[Section~8.2]{PS5}.
 So the topological ring $\S$ is topologically left coherent in
the sense of~\cite[Section~4]{Roos0}, \cite[Section~13]{PS3},
\cite[Section~7]{PS8}.

 We refer to~\cite{Sten0}, \cite[Section~1]{Pfp}, \cite[Section~6 and
Appendix~B]{Sto}, and~\cite[Section~2]{BHP} for discussions of
fp\+injective modules over coherent rings and fp\+injective objects
in locally coherent Grothendieck categories.
 A discrete left $\S$\+module $\J$ is called \emph{fp\+injective} if
the functor $\Hom_\S({-},\J)$ takes short exact sequences of finitely
presented discrete left $\S$\+modules to short exact sequences of
abelian groups.
 Given two discrete left $\S$\+modules $\L$ and $\M$, we denote by
$\Ext_\S^*(\L,\M)$ the Ext groups in the abelian category $\S\discr$.

\begin{lem} \label{fp-injective-discrete}
\textup{(a)} A discrete left\/ $\S$\+module\/ $\J$ is fp\+injective
if and only if the $S_n$\+module ${}_{S_n}\J$ is fp\+injective for
every $n\ge0$. \par
\textup{(b)} A discrete left\/ $\S$\+module\/ $\J$ is fp\+injective
if and only if\/ $\Ext_\S^1(\L,\J)=0$ for any finitely presented
discrete left\/ $\S$\+module\/ $\L$, and if and only if
$\Ext_\S^{>0}(\L,\J)=0$ for any such\/~$\L$. \par
\textup{(c)} The full subcategory of fp\+injective discrete left
$\S$\+modules is closed under extensions, the passages to the cokernels
of injective morphisms, and infinite direct sums (as well as filtered
inductive limits) in\/ $\S\discr$. \par
\textup{(d)} For any finitely presented discrete left\/ $\S$\+module\/
$\L$, the functor\/ $\Hom_\S(\L,{-})$ is exact on the exact category of
fp\+injective left\/ $\S$\+modules. \par
\textup{(e)} The functors\/ $\J\mpsto{}_{S_n}\J$ take short exact
sequences of fp\+injective discrete left\/ $\S$\+modules to short
exact sequences of fp\+injective left\/ $S_n$\+modules.
\end{lem}

\begin{proof}
 Part~(a) is straightforward.
 To prove part~(b), one notices that through any surjective morphism
$\E\rarrow\L$ from a discrete left $\S$\+module $\E$ onto a finitely
presented discrete left $\S$\+module $\L$ one can factorize
a surjective morphism $\K\rarrow\L$ from a finitely presented
discrete left $\S$\+module $\K$ onto~$\L$.
 The first two assertions of part~(c) follow from part~(b), and
the third one from part~(a).
 Part~(d) follows from part~(b).
 Part~(e) follows from part~(d) applied to the discrete left
$\S$\+modules $\L=S_n$.
\end{proof}

\begin{lem}  \label{finite-fp-injective-flat-projective-dimension}
\textup{(a)} The supremum of injective dimensions of fp\+injective
discrete left\/ $\S$\+modules (viewed as objects of the abelian
category\/ $\S\discr$) does not exceed the supremum of injective
dimensions of fp\+injective left $S_n$\+modules (viewed as objects
of the abelian category $S_n\modl$), taken over all\/ $n\ge0$. \par
\textup{(b)} Assume that, for a given $n\ge0$, there is an integer
$N_n\ge0$ such that every left ideal in the ring $S_n$ admits a set of
generators of the carginality not exceeding~$\aleph_{N_n}$.
 Then the injective dimension of any fp\+injective discrete left\/
$S_n$\+module is not greater than $N_n+1$.
 In particular, if all the left ideals in the rings $S_n$, for all
$n\ge0$, admit at most countable sets of generators, then the injective
dimension of any fp\+injective discrete left\/ $\S$\+module does not
exceed\/~$1$. \par
\textup{(c)} The supremum of projective dimensions of flat left\/
$\R$\+contramodules (viewed as objects of the abelian category\/
$\R\contra$) does not exceed the supremum of projective dimensions
of flat left $R_n$\+modules (viewed as objects of the abelian
category $R_n\modl$), taken over all\/ $n\ge0$.
\end{lem}

\begin{proof}
 Part~(a) follows from Lemma~\ref{fp-injective-discrete}(a,c,e) by
applying the functors $\J\mpsto{}_{S_n}\J$ to an injective coresolution
of a given fp\+injective discrete left $\S$\+module.
 Part~(b) is provided by~\cite[Proposition~2.3]{Pfp} together with
part~(a).
 Part~(c) holds in the less restrictive assumptions of
Section~\ref{flat-contramod-subsect}, and follows from
Lemmas~\ref{contramod-kernel-flat}, \ref{contramod-proj-flat},
and Corollary~\ref{projective-contramod}(a) by applying the functors
$\gF\mpsto{}\gF/\gI_n\tim\gF$ to a projective resolution of a given
flat left $\R$\+contramodule.
\end{proof}

\begin{lem}  \label{contramod-product-reduction}
 For any family of left\/ $\R$\+contramodules\/ $\P_\alpha$ and
any $n\ge0$, the two\/ $\R$\+subcontramodules\/
$\gI_n\tim\prod_\alpha\P_\alpha$ and\/ $\prod_\alpha\gI_n\tim\P_\alpha$
coincide in\/ $\prod_\alpha\P_\alpha$. 
\end{lem}

\begin{proof}
 The former subcontramodule is obviously contained in the latter one;
we have to prove the converse inclusion.
 For every $m\ge0$, pick a finite set of generators $\bar r_m^\gamma$
of the right ideal $\gI_m/\gI_{m+1}\sub \R/\gI_{m+1}$, and lift them
to elements $r_m^\gamma\in\R$.
 Then, for any left $\R$\+contramodule $\P$, any element of
the subcontramodule $\gI_n\tim\P\sub\P$ can be expressed in the form
$\sum_{m\ge n}^\gamma r_m^\gamma p_m^\gamma$ with some
$p_m^\gamma\in\P$.

 Indeed, the point is that, for any element $s\in\gI_l$, there exists
a family of elements $t_m^\gamma\in\R$ such that $s=\sum_{m\ge0}^\gamma
r_m^\gamma t_m^\gamma$ in $\R$ and $t_m^\gamma=0$ whenever $m<l$.
 Consequently, for any family of elements $s^\beta\in\gJ_n$ converging
to zero in the topology of $\R$ there exists a family of elements
$t_m^{\gamma,\beta}\in\R$ such that $s^\beta=\sum_{m\ge n}^\gamma
r_m^\gamma t_m^{\gamma,\beta}$ for every~$\beta$ and the family of
elements~$t_m^{\gamma,\beta}$ converges to zero in $\R$ as
$\beta$~varies for every fixed~$m$ and~$\gamma$.
 In fact, one can choose the elements~$t_m^{\gamma,\beta}$ so that
for every fixed~$m$ and~$\gamma$ one has $t_m^{\gamma,\beta}=0$
for all but a finite number of indices~$\beta$.
 Now, for every family of elements $p^\beta\in\P$ one has
$\sum_\beta s^\beta p^\beta=\sum_{m\ge n}^\gamma r_m^\gamma
(\sum_\beta t_m^{\gamma,\beta}p^\beta)$ in~$\P$.

 In particular, any element $p$ of the product $\prod_\alpha\gI_n\tim
\P_\alpha$ can be presented in the form $p=\big(\sum_{m\ge n}^\gamma
r_m^\gamma p_m^{\gamma,\alpha}\big)_\alpha$ with some elements
$p_m^{\gamma,\alpha}\in\P_\alpha$.
 Now one has $\big(\sum_{m\ge n}^\gamma
r_m^\gamma p_m^{\gamma,\alpha}\big)_\alpha =
\sum_{m\ge n}^\gamma r_m^\gamma((p_m^{\gamma,\alpha})_\alpha)$, which 
is an infinite sum of elements of $\prod_\alpha\P_\alpha$ with 
the coefficient family still converging to zero in $\gI_n\sub\R$,
proving that $p$~belongs to $\gI_n\tim\prod_\alpha\P_\alpha$.
 (Cf.~\cite[Lemmas~1.3.6\+-1.3.7]{Pweak}.)
\end{proof}


 Let us denote the exact subcategory of flat left $\R$\+contramodules by
$\R\contra_\fl\sub\R\contra$ and the exact subcategory of fp\+injective
discrete left $\S$\+modules by $\S\discr^\fpi\sub\S\discr$.
 Similarly, given a left coherent ring $S$, we let $S\modl^\fpi\sub
S\modl$ denote the exact subcategory of fp\+injective left $S$\+modules.

\begin{thm} \label{discrete-mod-coderived}
\textup{(a)} The coderived category\/ $\sD^\co(\S\discr)$ of
the abelian category of discrete left\/ $\S$\+modules is equivalent
to the coderived category of the exact category of fp\+injective
discrete left\/ $\S$\+modules. \par
\textup{(b)} The Becker coderived category\/ $\sD^\bco(\S\discr)$ of
the abelian category of discrete left\/ $\S$\+modules is equivalent
to the Becker coderived category of the exact category of fp\+injective
discrete left\/ $\S$\+modules. \par
\textup{(c)} The Becker coderived category\/ $\sD^\bco(\S\discr)$ is
equivalent to the homotopy category of complexes of injective
discrete left\/ $\S$\+modules. \par
\textup{(d)} Assume that, for every $n\ge0$, every fp\+injective
left $S_n$\+module has finite injective dimension in $S_n\modl$.
 Then the classes of Positselski-coacyclic and Becker-coacyclic
complexes in the abelian category\/ $\S\discr$ coincide.
 The coderived category\/ $\sD^{\co=\bco}(\S\discr)$ is equivalent
to the homotopy category of complexes of injective discrete left\/
$\S$\+modules. \par
\textup{(e)} Assume that every fp\+injective discrete left\/
$\S$\+module has finite injective dimension in\/ $\S\discr$.
 Then the classes of Positselski-coacyclic and Becker-coacyclic
complexes in the abelian category\/ $\S\discr$ coincide.
 The coderived category\/ $\sD^{\co=\bco}(\S\discr)$ is equivalent
to the absolute derived category of the exact category of fp\+injective
discrete left\/ $\S$\+modules. \par
\textup{(f)} Assume that every fp\+injective discrete left\/
$\S$\+module has finite injective dimension in\/ $\S\discr$.
 Then the coderived category\/ $\sD^{\co=\bco}(\S\discr)$ is equivalent
to the homotopy category of complexes of injective discrete left\/
$\S$\+modules.
\end{thm}

\begin{proof}
 Part~(a) holds, because the category $\S\discr$, being
a Grothendieck abelian category, has enough injective objects;
so it remains to use Lemma~\ref{fp-injective-discrete}(c)
together with the dual version of
Propositions~\ref{infinite-resolutions}(b).
 Part~(b) follows similarly from Lemma~\ref{fp-injective-discrete}(c)
and the dual version of
Proposition~\ref{becker-contraderived-infinite-resolutions}.
 Part~(c) is a special case of
Theorem~\ref{coderived-of-grothendieck-contraderived-of-lpacepo}(a).
 To deduce parts~(e\+f), one can apply (the dual versions of)
Corollary~\ref{finite-homol-dim-equivalence-cor},
Proposition~\ref{finite-resolutions}, and/or
Lemma~\ref{psemi-remark21} for Positselski's coderived categories;
or Theorems~\ref{finite-homol-dim-becker-co-contra-derived}(a) and
\ref{positselski-becker-co-contra-derived}(a) for Becker's ones.

 The most nontrivial assertion is~(d).
 Notice that, in the assumption of~(d), for every $n\ge0$ there exists
an integer $d_n\ge0$ such that the injective dimension of any
fp\+injective left $S_n$\+module does not exceed~$d_n$ (since the class
of all fp\+injective modules is closed under infinite direct sums and
products in $S_n\modl$); but the numbers~$d_n$ may grow to infinity
as $n\to\infty$.
 So the injective dimensions of fp\+injective discrete left
$\S$\+modules \emph{need not} be finite.

 Nevertheless, in view of parts~(a\+c) and
Lemma~\ref{Positselski-trivial-are-Becker-trivial}(b), we only need to
prove that every Becker-coacyclic complex $\cA^\bu$ in $\S\discr^\fpi$
is Positselski-coacyclic in $\S\discr^\fpi$ (or in $\S\discr$).
 Then the equivalences of categories $\sD^\co(\S\discr^\fpi)\rarrow
\sD^\co(\S\discr)$ and $\sD^\co(\S\discr^\fpi)\rarrow\sD^\bco
(\S\discr^\fpi)\rarrow\sD^\bco(\S\discr)$ would imply the desired
equivalence $\sD^\co(\S\discr)\rarrow\sD^\bco(\S\discr)$.

 Now we apply~\cite[Theorem~6.12]{Sto}, which tells us that $\cA^\bu$ is
an acyclic complex in the exact category $\S\discr^\fpi$.
 By Lemma~\ref{fp-injective-discrete}(a), it follows that
${}_{S_n}\cA^\bu$ is an acyclic complex in $S_n\modl^\fpi$ for every
$n\ge0$.
 By assumption and in view of Lemma~\ref{psemi-remark21}, we can
conclude that ${}_{S_n}\cA^\bu$ is an absolutely acyclic (hence
Positselski-coacyclic) complex in $S_n\modl^\fpi$.
 Finally, it follows that the complex $\cA^\bu=\varinjlim_{n\ge0}
{}_{S_n}\cA^\bu$ is Positselski-coacyclic in $\S\modl$ (or equivalently,
in $\S\modl^\fpi$), because the class of Positselski-coacyclic complexes
is preserved by the inductive limits of sequences indexed by
the natural numbers.
\end{proof}

\begin{thm} \label{contramod-contraderived}
\textup{(a)} The contraderived category\/ $\sD^\ctr(\R\contra)$ of
the abelian category of left\/ $\R$\+contramodules is equivalent
to the contraderived category of the exact category of flat
left\/ $\R$\+contramodules. \par
\textup{(b)} The Becker contraderived category\/ $\sD^\bctr(\R\contra)$
of the abelian category of left\/ $\R$\+contramodules is equivalent
to the Becker contraderived category of the exact category of flat
left\/ $\R$\+contramodules. \par
\textup{(c)} The Becker contraderived category\/ $\sD^\bctr(\R\contra)$
is equivalent to the homotopy category of complexes of projective
left\/ $\R$\+contramodules. \par
\textup{(d)} Assume that, for every $n\ge0$, every flat left
$R_n$\+module has finite projective dimension in $R_n\modl$.
 Then the classes of Positselski-contraacyclic and Becker-contraacyclic
complexes in the abelian category\/ $\R\contra$ coincide.
 The contraderived category\/ $\sD^{\ctr=\bctr}(\R\contra)$ is
equivalent to the homotopy category of complexes of projective left\/
$\R$\+contramodules. \par
\textup{(e)} Assume that every flat left\/ $\R$\+contramodule has finite
projective dimension in\/ $\R\contra$.
 Then the classes of Positselski-contraacyclic and Becker-contraacyclic
complexes in the abelian category\/ $\R\contra$ coincide.
 The contraderived category\/ $\sD^{\ctr=\bctr}(\R\contra)$ is
equivalent to the absolute derived category of the exact category of
flat left\/ $\R$\+contramodules. \par
\textup{(f)} Assume that every flat left\/ $\R$\+contramodule has finite
projective dimension in\/ $\R\contra$.
 Then the contraderived category\/ $\sD^{\ctr=\bctr}(\R\contra)$ is
equivalent to the homotopy category of complexes of projective left\/
$\R$\+contramodules.
\end{thm}

\begin{proof}
 Any left $\R$\+contramodule is a quotient contramodule of a flat (and
even projective) one, and the class of flat left $\R$\+contramodules
is closed under extensions and the passage to the kernels of
surjective morphisms
(see Lemmas~\ref{contramod-kernel-flat}\+-\ref{contramod-proj-flat}).
 Hence, in order to prove the assertions~(a\+b), it only remains to
show that the class of flat left $\R$\+contramodules is preserved by
infinite products in $\R\contra$
(see Propositions~\ref{infinite-resolutions}(b)
and~\ref{becker-contraderived-infinite-resolutions}).
 The latter follows from the definition of flatness for left
$\R$\+contramodules, the coherence condition on the rings $R_n$,
and Lemma~\ref{contramod-product-reduction}.
 Part~(c) is a special case of
Theorem~\ref{coderived-of-grothendieck-contraderived-of-lpacepo}(b).
 To deduce parts~(e\+f), one can apply
Corollary~\ref{finite-homol-dim-equivalence-cor},
Proposition~\ref{finite-resolutions}, and/or
Lemma~\ref{psemi-remark21} for Positselski's contraderived categories;
or Theorems~\ref{finite-homol-dim-becker-co-contra-derived}(b) and
\ref{positselski-becker-co-contra-derived}(b) for Becker's ones.

 The most nontrivial assertion is~(d).
 Once again, we notice that, in the assumption of~(d), for every
$n\ge0$ there exists an integer $d_n\ge0$ such that the projective
dimension of any flat left $R_n$\+module does not exceed~$d_n$ (since
the class of all flat modules is closed under infinite direct sums
and products in $R_n\modl$); but the numbers~$d_n$ may grow to infinity
as $n\to\infty$.
 So the projective dimensions of flat left $\R$\+contramodules
\emph{need not} be finite.

 Nevertheless, in view of parts~(a\+c) and
Lemma~\ref{Positselski-trivial-are-Becker-trivial}(c), we only need to
prove that every Becker-contraacyclic complex $\gB^\bu$ in
$\R\contra_\fl$ is Positselski-contraacyclic in $\R\contra_\fl$.
 Then the equivalences of categories $\sD^\ctr(\R\contra_\fl)\rarrow
\sD^\ctr(\R\contra)$ and $\sD^\ctr(\R\contra_\fl)\rarrow
\sD^\bctr(\R\contra_\fl)\rarrow\sD^\bctr(\R\contra)$ would imply
the desired equivalence $\sD^\ctr(\R\contra)\rarrow
\sD^\bctr(\R\contra)$.

 Now we apply a result of the preprint~\cite[Theorem~6.1]{Pbc}, which
tells us that $\gB^\bu$ is an acyclic complex in the exact category
$\R\contra_\fl$.
 By Lemma~\ref{contramod-kernel-flat}, it follows that
$\gB^\bu/\gI_n\tim\gB^\bu$ is an acyclic complex in $R_n\modl_\fl$
for every $n\ge0$.
 By assumption and in view of Lemma~\ref{psemi-remark21}, we can
conclude that $\gB^\bu/\gI_n\tim\gB^\bu$ is an absolutely acyclic
(hence Positselski-contraacyclic) complex in $R_n\modl_\fl$.
 Finally, it follows that the complex $\gB^\bu=\varprojlim_{n\ge0}
\gB^\bu/\gI_n\tim\gB^\bu$ is Positselski-contraacyclic in
$\R\contra$ (hence in $\R\contra_\fl$), because the class of
Positselski-contraacyclic complexes in $\R\contra$ is preserved by
the projective limits of sequences of surjective morphisms.
\end{proof}

\begin{rem} \label{no-needed-adjoints-remark}
 A more natural approach to proving
Theorems~\ref{discrete-mod-coderived}(d)
and~\ref{contramod-contraderived}(d) would be to show that
the functor $\M\mpsto{}_{S_n}\M$ takes Becker-coacyclic complexes
in $\S\discr^\fpi$ to Becker-coacyclic complexes in $S_n\modl^\fpi$,
while the functor $\gF\mpsto\gF/\gI_n\tim\gF$ takes Becker-contraacyclic
complexes in $\R\contra_\fl$ to Becker-contraacyclic complexes in
$R_n\modl_\fl$.
 Then it would remain to refer to
Theorem~\ref{finite-homol-dim-becker-co-contra-derived} to the effect
that all Becker-coacyclic complexes in $S_n\modl^\fpi$ are absolutely
acyclic and all Becker-contraacyclic copmlexes in $R_n\modl_\fl$ are
absolutely acyclic.
 Unfortunately, we are \emph{not} aware of any direct proof of
preservation of the Becker co/contraacyclicity by the respective
functors (while the preservation of the Positselski co/contraacyclicity
follows directly from the facts that the two functors are exact and
preserve infinite direct sums or products, respectively).
 The problem is that the two functors have no adjoints on the sides
where they are needed, so
Lemma~\ref{exact-with-adjoint-preservation-lemma} is not applicable.
\end{rem}

 A right $\R$\+module $\N$ is said to be discrete if every its element
is annihilated by the ideal $\gI_n=\ker(\R\to R_n)$ for $n\gg0$, i.~e.,
in the notation similar to the above, if $\N=\bigcup_n\N_{R_n}$.
 The abelian category of discrete right $\R$\+modules is denoted by
$\discrR\R$.
 For any associative ring $S$, any $\R$\+discrete $S$\+$\R$\+bimodule
$\K$, and any left $S$\+module $U$, the abelian group $\Hom_S(\N,U)$
is naturally endowed with a left $\R$\+contramodule structure as
the projective limit of the sequence of left $R_n$\+modules
$\Hom_S(\N,U)=\varprojlim_n\Hom_S(\N_{R_n},U)$.
 For any discrete right $\R$\+module $\N$ and any left
$\R$\+contramodule $\P$, their \emph{contratensor product} $\N\ocn_\R\P$
is defined as the inductive limit of the sequence of abelian groups
$\varinjlim_n\N_{R_n}\ot_{R_n}\P/(\gI_n\tim\P)$.

 To give a more fancy definition, the contratensor product $\N\ocn_\R\P$
is the cokernel (of the difference) of the natural pair of abelian group
homomorphisms $\N\ot_\boZ\R[[\P]]\birarrow\N\ot_\boZ\P$.
 Here one map is induced by the left contraaction map
$\R[[\P]]\rarrow\P$ and the other one is the composition
$\N\ot_\boZ\R[[\P]]\rarrow\N[\P]\rarrow\N\ot_\boZ\P$, where $\N[\P]$
is the group of all finite formal linear combinations of elements
of $\P$ with the coefficients in $\N$, the former map to be composed
is induced by the discrete right action map $\N\times\R\rarrow\N$,
and the latter map is just the obvious one.
 The contratensor product is a right exact functor of two
arguments $\ocn_\R\:\discrR\R\times\R\contra\rarrow\boZ\modl$.

 For any $\R$\+discrete $S$\+$\R$\+bimodule $\K$, any left
$\R$\+contramodule $\P$, and any left $S$\+module $U$, there are
natural isomorphisms of abelian groups
$\Hom_S(\K\ocn_\R\P\;U)\simeq\varprojlim_n\Hom_S(\K_{R_n}\ot_{R_n}
\P/\gI_n\tim\P\;U)\simeq\varprojlim_n\Hom_{R_n}(\P/\gI_n\tim\P\;
\Hom_S(\K_{R_n},U))\simeq\Hom^\R(\P,\Hom_S(\K,U))$
(cf.~\cite[Sections~3.1.2 and~5.1.1]{Psemi}, \cite[Section~5]{PR},
\cite[Section~2.8]{Pproperf}, \cite[Section~2.8]{Pcoun},
\cite[Section~7.2]{PS1}).

\medskip

 Given surjective ring homomorphisms $g\:S\rarrow {}'\!S$ and
$f\:R\rarrow {}'\!R$, a left $S$\+module $M$, and a right $R$\+module
$N$, we denote by ${}_{{}'\!S}M$ the submodule of all elements
annihilated by the left action of $\ker(g)$ in $M$ and by $N_{\.{}'\!R}$
the submodule of all elements annihilated by the right action of
$\ker(f)$ in~$N$.
 So ${}_{{}'\!S}M$ is the maximal left ${}'\!S$\+submodule in $M$
and $N_{\.{}'\!R}$ is the maximal right ${}'\!R$\+submodule in~$N$.
 Similarly, for any left $R$\+module $P$ we denote by ${}^{{}'\!R\!}P$
its maximal quotient left ${}'\!R$\+module $P/\ker(f)P$.

\begin{lem}  \label{ind-affine-hom-tensor-reduction}
 Assume that\/ $\ker(g)$ is finitely generated as a left ideal in $S$
and\/ $\ker(f)$ is finitely generated as a right ideal in~$R$.
 Let $K$ be an $S$\+$R$\+bimodule such that its submodules\/
${}_{{}'\!S}K$ and $K_{\.{}'\!R}$ coincide; denote this\/ 
${}'\!S$\+${}'\!R$\+subbimodule in $K$ by\/~${}'\!K$.
 Then \par
\textup{(a)} for any injective left $S$\+module $J$ there is 
a natural isomorphism of left\/ ${}'\!R$\+modules\/ ${}^{{}'\!R\!}
\Hom_S(K,J)\simeq\Hom_{\.{}'\!S}({}'\!K,\.{}_{{}'\!S}J)$; \par
\textup{(b)} for any flat left $R$\+module $F$ there is
a natural isomorphism of left\/ ${}'\!S$\+modules\/
${}_{{}'\!S}(K\ot_R F)\simeq {}'\!K\ot_{\,{}'\!R} {}^{{}'\!R\!}F$.
\end{lem}

\begin{proof}
 Part~(a): for any finitely presented right $R$\+module $E$,
there is a natural isomorphism
$\Hom_S(\Hom_{R^\rop}(E,K),J)\simeq E\ot_R\Hom_S(K,J)$.
 Taking $E={}'\!R$, we get $\Hom_{\.{}'\!S}({}_{{}'\!S}K,\.{}_{{}'\!S}J)
\simeq\Hom_S({}_{{}'\!S}K,J) = \Hom_S(K_{\.{}'\!R},\.J)\simeq
{}^{{}'\!R\!}\Hom_S(K,J)$.

 Part~(b): for any finitely presented left $S$\+module $E$
there is a natural isomorphism
$\Hom_S(E\;K\ot_RF)\simeq\Hom_S(E,K)\ot_RF$.
 Taking $E={}'\!S$, we get $K_{\.{}'\!R}\ot_{\,{}'\!R} {}^{{}'\!R\!}F
\simeq K_{\.{}'\!R}\ot_R F = {}_{{}'\!S}K \ot_R F \simeq
{}_{{}'\!S}(K \ot_R F)$.
\end{proof}

 The following definition of a dualizing complex for a pair of
coherent noncommutative rings is a slight generalization of
the one from~\cite[Section~4]{Pfp} (which is, in turn, a generalization
of the one from the above Section~\ref{corings-dualizing}).

 Let $S$ be a left coherent ring and $R$ be a right coherent ring.
 A complex of $S$\+$R$\+bimodules $D^\bu$ is called
a \emph{dualizing complex} for the rings $S$ and $R$ if it satisfies
the following conditions:
\begin{enumerate}
\renewcommand{\theenumi}{\roman{enumi}}
\item the terms of the complex $D^\bu$ are fp\+injective left
$S$\+modules and fp\+injective right $R$\+modules;
\item the complex $D^\bu$ is isomorphic, as an object of the coderived
category of the exact category of $S$\+fp\+injective and
$R$\+fp\+injective $S$\+$R$\+bimodules, to a finite complex of
$S$\+fp\+injective and $R$\+fp\+injective $S$\+$R$\+bimodules;
\item the $S$\+$R$\+bimodules of cohomology $H^*(D^\bu)$ of
the complex $D^\bu$ are finitely presented left $S$\+modules and
finitely presented right $R$\+modules; and
\item the homothety maps $S\rarrow\Hom_{\sD(\modr R)}(D^\bu,D^\bu)$ and
$R^\rop\rarrow\Hom_{\sD(S\modl)}\allowbreak(D^\bu,D^\bu[*])$ are
isomorphisms of graded rings.
\end{enumerate}

 A dualizing complex in the sense of this definition satisfies
the conditions of the definition of a dualizing complex in
the sense of~\cite[Section~4]{Pfp} if and only if it is a finite
complex of $S$\+$R$\+bimodules.
 Any dualizing complex in the sense of the above definition is
isomorphic, as an object of the coderived category of
$S$\+fp\+injective and $R$\+fp\+injective $S$\+$R$\+bimodules,
to a dualizing complex in the sense of~\cite{Pfp}.
 A dualizing complex in the sense of~\cite{Pfp} satisfies
the conditions of the definition of a dualizing complex of bimodules
in the sense of our Appendix~\ref{over-flat-coring-appx} if and only
if it is a complex of $S$\+injective and $R$\+injective (rather than
just fp\+injective) bimodules.
 This is what is called a ``strong dualizing complex''
in~\cite[Section~3]{Pfp}.

\begin{lem}  \label{ind-affine-dualizing-reduction}
 Let $D^\bu$ be a finite complex of $S$\+fp\+injective and
$R$\+fp\+injective $S$\+$R$\+bimodules.
 Suppose that the subcomplexes\/ ${}_{{}'\!S}D^\bu$ and $D^\bu_{\.{}'\!R}$
coincide in~$D^\bu$, and denote this complex of\/
${}'\!S$\+${}'\!R$\+bimodules by\/~${}'\!D^\bu$.  \par
\textup{(a)} Assume that the ring $S$ is left coherent, the ring $R$
is right coherent, all fp\+injective left $S$\+modules have finite
injective dimensions, all fp\+injective right $R$\+modules have finite
injective dimensions, $\ker(g)$ is finitely generated as a left ideal
in $S$, and\/ $\ker(f)$ is finitely generated as a right
ideal in~$R$.
 Then\/ ${}'\!D^\bu$ is a dualizing complex for the rings\/ ${}'\!S$
and\/ ${}'\!R$ whenever $D^\bu$ is dualizing complex for the rings $S$
and~$R$. \par
\textup{(b)} Assume that the ring $S$ is left Noetherian, the ring
$R$ is right Noetherian, and the ideals\/ $\ker(g)\sub S$ and\/
$\ker(f)\sub R$ are nilpotent.
 Then $D^\bu$ is a dualizing complex for the rings $S$ and $R$
whenever\/ ${}'\!D^\bu$ is a dualizing complex for the rings\/
${}'\!S$ and\/~${}'\!R$.
\end{lem}

\begin{proof}
 Part~(a): clearly, ${}'\!D^\bu={}_{{}'\!S}D^\bu$ is a finite complex of
fp\+injective left ${}'\!S$\+modules.
 Let $E^\bu$ be a finite complex of injective left $S$\+modules
endowed with a quasi-isomorphism of complexes of $S$\+modules
$D^\bu\rarrow E^\bu$.
 The natural map $R\rarrow\Hom_S(D^\bu,E^\bu)$ is a quasi-isomorphism
of finite complexes of flat left $R$\+modules
(see~\cite[Lemma~4.1(b)]{Pfp}) and therefore remains a quasi-isomorphism
after taking the tensor product with ${}'\!R$ over $R$ on the left,
i.~e., reducing modulo $\ker(f)$.
 By Lemma~\ref{ind-affine-hom-tensor-reduction}(a), it follows that
the map ${}'\!R\rarrow\Hom_{\.{}'\!S}({}'\!D^\bu,\.{}_{{}'\!S}E^\bu)$
is an quasi-isomorphism.
 It remains to notice that ${}_{{}'\!S}E^\bu$ is a complex of injective
left ${}'\!S$\+modules and the induced map ${}_{{}'\!S}D^\bu\rarrow
{}_{{}'\!S}E^\bu$ is a quasi-isomorphism of complexes of left
${}'\!S$\+modules, as the functor $J\mpsto {}_{{}'\!S}J^\bu$ is exact
on the exact category of fp\+injective left $S$\+modules.

 To show that the left ${}'\!S$\+modules of cohomology of the complex
${}'\!D^\bu=D^\bu_{\.{}'\!R}$ are finitely presented, we notice that it
is quasi-isomorphic, as a complex of left ${}'\!S$\+modules, to
the complex $\Hom_{R^\rop}(L_\bu,D^\bu)$ of right $R$\+module
homomorphisms from a resolution $L_\bu$ of the right $R$\+module
${}'\!R$ by finitely generated free $R$\+modules into
the complex~$D^\bu$.
 The argument finishes similarly to the proof of 
Lemma~\ref{noncomm-dualizing-reduction}(a).
 The proof of part~(b) is also similar to that of 
Lemma~\ref{noncomm-dualizing-reduction}(b).
\end{proof}

 A complex of $\S$\+$\R$\+bimodules $\D^\bu$ is called a \emph{dualizing
complex} for the projective systems of rings $(S_n)$ and~$(R_n)$ if
\begin{enumerate}
\renewcommand{\theenumi}{\roman{enumi}}
\item $\D^\bu$ is a complex of fp\+injective discrete left $\S$\+modules
and a complex of fp\+injective discrete right $\R$\+modules;
\item for every $n\ge0$, the two subcomplexes ${}_{S_n}\!\D^\bu$ and
$\D^\bu_{R_n}$ coincide in $\D^\bu$; and
\item the latter subcomplex, denoted by ${}_{S_n}\!\D^\bu = D^\bu_n
= \D^\bu_{R_n}\sub\D^\bu$, is a dualizing complex for the rings
$S_n$ and~$R_n$.
\end{enumerate}
 The complex $\D^\bu$ does not have to be bounded from either side,
and neither do its subcomplexes $D_n^\bu$; but the latter has to be
isomorphic, as an object of the coderived category of
$S_n$\+fp\+injective and $R_n$\+fp\+injective $S_n$\+$R_n$\+bimodules,
to a finite (dualizing) complex of $S_n$\+fp\+injective and
$R_n$\+fp\+injective $S_n$\+$R_n$\+bimodules.

\begin{thm}  \label{ind-affine-derived-co-contra}
 Assume that, for every $n\ge0$, all fp\+injective left $S_n$\+modules
have finite injective dimensions.
 Then the choice of a dualizing complex\/ $\D^\bu$ for the projective
systems of rings $(S_n)$ and $(R_n)$ induces an equivalence between
the coderived category of discrete left\/ $\S$\+modules\/
$\sD^{\co=\bco}(\S\discr)$ and the contraderived category of left\/
$\R$\+contramodules\/ $\sD^{\ctr=\bctr}(\R\contra)$.
 The equivalence is provided by the derived functors\/
$\boR\Hom_\S(\D^\bu,{-})$ and\/ $\D^\bu\ocn_\R^\boL{-}$ of the functor
of discrete\/ $\S$\+module homomorphisms from\/ $\D^\bu$ and
the functor of contratensor product with\/ $\D^\bu$ over\/~$\R$.
\end{thm}

\begin{proof}
 First of all, one has $\sD^\co(\S\discr)=\sD^\bco(\S\discr)$ by
Theorem~\ref{discrete-mod-coderived}(d).
 Furthermore, all flat left $R_n$\+modules have finite projective
dimensions by~\cite[Proposition~4.3]{Pfp}.
 Hence one has $\sD^\ctr(\R\contra)=\sD^\bctr(\R\contra)$ by
Theorem~\ref{contramod-contraderived}(d).

 The constructions of the derived functors are based on
Theorems~\ref{discrete-mod-coderived}(d)
and~\ref{contramod-contraderived}(a) or~(b).
 By Lemmas~\ref{ind-affine-hom-tensor-reduction}(a) and
\ref{projlim-reduction-lemma} together with~\cite[Lemma~4.1(b)]{Pfp},
applying the functor $\J^\bu\mpsto\Hom_\S(\D^\bu,\J^\bu)\simeq
\varprojlim_n\Hom_{S_n}({}_{S_n}\!\D^\bu,\.{}_{S_n}\J^\bu)$ to
a complex of injective discrete left $\S$\+modules $\J^\bu$
produces a complex of flat left $\R$\+contramodules. {\hfuzz=3pt\par}

 By Lemma~\ref{ind-affine-hom-tensor-reduction}(b)
and~\cite[Lemma~4.1(a)]{Pfp}, applying the functor
$\gF^\bu\mpsto\D^\bu\ocn_\R\gF^\bu\simeq
\varinjlim_n\D^\bu_{R_n}\ot_{R_n}\gF^\bu/(\gI_n\tim\gF^\bu)$
to a complex of flat left $\R$\+contramodules $\gF^\bu$
produces a complex of fp\+injective discrete left $\S$\+modules.
 Let us show that the complex $\D^\bu\ocn_\R\gF^\bu$ is coacyclic
whenever the complex $\gF^\bu$ is contraacyclic.

 Coacyclicity being preserved by inductive limits of sequences, it
suffices to show that the complexes of left $S_n$\+modules
$D^\bu_n\ot_{R_n}\gF^\bu/\gI_n\tim\gF^\bu$ are coacyclic.
 The functor $\gF\mpsto\gF/\gI_n\tim\gF$ preserves infinite products
of left $\R$\+contramodules by Lemma~\ref{contramod-product-reduction}
and short exact sequences of flat left $\R$\+contramodules by
Lemma~\ref{contramod-kernel-flat}.
 Hence it takes contraacyclic complexes of flat left
$\R$\+contramodules $\gF^\bu$ to contraacyclic complexes of
flat left $R_n$\+modules $\gF^\bu/\gI_n\tim\gF^\bu$.

 Now the exact category of flat left $R_n$\+modules has finite
homological dimension by~\cite[Proposition~4.3]{Pfp}, so any
contraacyclic complex of flat left $R_n$\+modules is absolutely acyclic
by Lemma~\ref{psemi-remark21} (see also
Section~\ref{finite-homol-dim-subsect}), and the functor
$D_n^\bu\ot_{R_n}\nobreak{-}$ clearly takes absolutely acyclic
complexes of flat left $R_n$\+modules to absolutely acyclic complexes
of left $S_n$\+modules.
 (Concerning the apparent ambiguity of our terminology, notice that
the classes of complexes of flat objects contraacyclic or absolutely
acyclic with respect to the abelian categories of arbitrary objects
coincide with those of flat objects contraacyclic or absolutely
acyclic with respect to the exact categories of flat objects by
Proposition~\ref{fully-faithful-prop}.)

 Furthermore, given a complex of injective discrete left $\S$\+modules
$\J^\bu$, in order to check that the adjunction morphism
$\D^\bu\ocn_\R\Hom_\S(\D^\bu,\J^\bu)\rarrow\J^\bu$ has a coacyclic cone,
it suffices to show that so does the morphism 
$D_n^\bu\ot_{R_n}\Hom_{S_n}(D_n^\bu,\.{}_{S_n}\J^\bu)
\rarrow{}_{S_n}\J^\bu$
(in view of Lemma~\ref{ind-affine-hom-tensor-reduction}, and
because the class of coacyclic complexes is closed with respect to
inductive limits of sequences).
 Given a complex of flat $\R$\+contramodules $\gF^\bu$, choose
a complex of injective discrete left $\S$\+modules $\J^\bu$ endowed
with a morphism of complexes of discrete left $\S$\+modules
$\D^\bu\ocn_\R\gF^\bu\rarrow\J^\bu$ with a coacyclic cone.
 This cone is a complex of fp\+injective discrete left $\S$\+modules
and, by Theorem~\ref{discrete-mod-coderived}(a), it is also coacyclic
with respect to the exact category of fp\+injective discrete left
$\S$\+modules.
 By Lemma~\ref{fp-injective-discrete}(e), it follows that the morphisms
of complexes of left $S_n$\+modules $D_n^\bu\ot_{R_n}\gF^\bu/
(\gI_n\tim\gF^\bu)\rarrow{}_{S_n}\J^\bu$ have coacyclic cones, too.
 In order to check that the morphism of complexes of left
$\R$\+contramodules $\gF^\bu\rarrow\Hom_\S(\D^\bu,\J^\bu)$
has a contraacyclic cone, it suffices to check that so does the morphism
$\gF^\bu/\gI_n\tim\gF^\bu\rarrow\Hom_{S_n}(D_n^\bu,\.{}_{S_n}\J^\bu)$
(because the class of contraacyclic complexes is closed with respect
to projective limits of sequences of surjective morphisms).

 Finally, the pair of adjoint functors $\boR\Hom_{S_n}(D_n^\bu,{-})$
and $D_n^\bu\ot_{R_n}^\boL{-}$ between the coderived category of
left $S_n$\+modules and the contraderived category of left
$R_n$\+modules is connected by a chain of isomorphisms of pairs
of adjoint functors with a similar pair of adjoint functors
corresponding to a finite dualizing complex of
$S_n$\+$R_n$\+bimodules.
 The latter pair of functors are mutually inverse equivalences
by~\cite[Theorem~4.5]{Pfp}.
\end{proof}

\subsection{Cotensor product of complexes of discrete modules}
 Now we specialize from the setting of the previous section to the case
when the two projective systems $R_0\larrow R_1\larrow R_2\larrow\dotsb$
and $S_0\larrow S_1\larrow S_2\larrow\dotsb$ coincide and the rings
$R_n=S_n$ are commutative.
 The morphisms $R_{n+1}\rarrow R_n$ are still assumed to be surjective,
their kernel ideals to be finitely generated, and the rings $R_n$
to be coherent.
 As above, we set $\R=\varprojlim_n R_n$ and denote by $\gI_n\sub\R$
the kernels of the surjective ring homomorphisms $\R\rarrow R_n$.

 A complex of $\R$\+modules $\D^\bu$ is called a \emph{dualizing
complex} for the topological ring $\R$ if it is a dualizing complex
for the two coinciding projective systems of rings $(R_n)=(S_n)$ in
the sense of the definition from the previous section, that is
\begin{enumerate}
\renewcommand{\theenumi}{\roman{enumi}}
\item $\D^\bu$ is a complex of fp\+injective discrete $\R$\+modules;
\item for every $n\ge0$, the maximal subcomplex of $R_n$\+modules
$D_n^\bu={}_{R_n}\!\D^\bu$ in $\D^\bu$ is isomorphic, as an object of
the coderived category of the exact category of fp\+injective
$R_n$\+modules, to a finite complex of fp\+injective $R_n$\+modules;
\item the cohomology modules $H^*(D_n^\bu)$ of the complex $D_n^\bu$
are finitely presented $R_n$\+modules; and
\item the homothety map $R_n\rarrow\Hom_{\sD(R_n\modl)}(D_n^\bu,D_n^\bu)$
is an isomorphism of graded rings.
\end{enumerate}
 One can conclude from Lemma~\ref{ind-affine-dualizing-reduction}(a)
that the property of a complex of $\R$\+modules to be a dualizing
complex does not in fact depend on the choice of the projective
system of discrete rings $(R_n)$ satisfying the above-listed
conditions and representing a given topological ring~$\R$.

 The construction of the functor of \emph{cotensor product} of
complexes of discrete $\R$\+modules below in this section is to be
compared to the constructions of the cotensor product of complexes
of quasi-coherent sheaves on a Noetherian scheme with a dualizing
complex in~\cite[Section~B.2.5]{EP} and of the semitensor product of
complexes of modules in a relative situation in~\cite[Section~6]{Pfp},
as well as to the construction of the cotensor product of complexes
of quasi-coherent torsion sheaves on an ind-Noetherian ind-scheme with
a dualizing complex in~\cite[Chapter~5]{Psemten}.

 The definition of the (\emph{contramodule}) \emph{tensor product}  
$\P\ot^\R\Q$ of two $\R$\+contramodules $\P$ and $\Q$ can be found
in~\cite[Section~1.6]{Pweak}.
 The contramodule tensor product $\ot^\R$ is a right exact functor
of two arguments defining an associative and commutative tensor
category structure on the category $\R\contra$ with the unit
object~$\R$.
 The functor $\ot^\R$ also commutes with infinite direct sums in
$\R\contra$.
 The tensor product of free $\R$\+contramodules is described by
the rule $\R[[X]]\ot^\R\R[[Y]]\simeq\R[[X\times Y]]$ for any
sets $X$ and~$Y$.

\begin{lem} \label{flat-contramodules-tensor-product}
 The tensor product of two flat\/ $\R$\+contramodules\/ $\gF$
and\/ $\gG$ is a flat\/ $\R$\+contramodule that can be constructed
explicitly as\/ $\gF\ot^\R\gG=\varprojlim_n (\gF/\gI_n\tim\gF)
\ot_{R_n}(\gG/\gI_n\tim\gG)$.
 The contramodule tensor product functor is exact on the exact
category of flat\/ $\R$\+contramodules.
\end{lem}

\begin{proof}
 The assertions of this Lemma hold true for any topological ring
$\R$ isomorphic to the topological projective limit of a sequence
of surjective morphisms of commutative rings.
 Indeed, the functor defined by the above explicit formula takes
flat $\R$\+contramodules to flat $\R$\+contramodules by
Lemma~\ref{projlim-reduction-lemma}, and is exact on the exact
category of flat $\R$\+contramodules by
Lemma~\ref{contramod-kernel-flat} and because the projective limit
of sequences of surjective maps is an exact functor.
 This functor is isomorphic to the functor $\ot^\R$ for flat
$\R$\+contramodules $\gF$ and $\gG$, since both functors are right
exact for flat $\R$\+contramodules and they are isomorphic
for the free ones.

 Moreover, there is a natural isomorphism of $R_n$\+modules
$(\P\ot^\R\Q)/(\gI_n\tim(\P\ot^\R\nobreak\Q))\allowbreak\simeq
(\P/\gI_n\tim\P)\ot_{R_n}(\Q/\gI_n\tim\Q)$ for any
$\R$\+contramodules $\P$ and $\Q$, because all the functors involved
are right exact and there is such a natural isomorphism for all
free or flat $\R$\+contramodules $\P$ and~$\Q$.
\end{proof}

 The tensor product of two complexes of $\R$\+contramodules
$\P^\bu\ot^\R\Q^\bu$ is constructed by totalizing the bicomplex
$\P^i\ot^\R\Q^j$ by taking infinite direct sums (in the category
$\R\contra$) along the diagonals.
 Since the reduction functor $\P\mpsto\gI_n\tim\P$, being left
adjoint to the embedding functor $R_n\modl\rarrow\R\contra$,
commutes with infinite direct sums, we know from
Corollary~\ref{flat-reductions-enough} that infinite direct sums
of flat $\R$\+contramodules are flat and from
Lemma~\ref{flat-contramodules-tensor-product} that the tensor product
of two complexes of flat $\R$\+contramodules $\gF^\bu$ and $\gG^\bu$
can be constructed as the projective limit of the complexes of
$R_n$\+modules $(\gF^\bu/\gI_n\tim\gF^\bu)\ot_{R_n}
(\gG^\bu/\gI_n\tim\gG^\bu)$.

 Recall that all flat $R_n$\+modules have finite projective dimensions
provided that the ring $R_n$ admits a dualizing complex and all
fp\+injective $R_n$\+modules have finite injective
dimensions~\cite[Proposition~4.3]{Pfp}, and that all fp\+injective
$R_n$\+modules have injective dimensions~$\le1$ provided that
all ideals in the rings $R_n$ are at most countably generated
(Lemma~\ref{finite-fp-injective-flat-projective-dimension}).

\begin{lem} \label{contramod-tensor-product-contraacyclic}
 Assume that, for every $n\ge0$, all flat $R_n$\+modules have finite
projective dimensions.
 Then the tensor product\/ $\gF^\bu\ot^\R\gG^\bu$ of two complexes of
flat $\R$\+contramodules is contraacyclic whenever one of
the complexes\/ $\gF^\bu$ and\/ $\gG^\bu$ is contraacyclic.
\end{lem}

\begin{proof}
 First of all we recall that a complex of flat $\R$\+contramodules is
contraacyclic as a complex in the abelian category of
$\R$\+contramodules if and only if it is contraacyclic as a complex
in the exact category of flat $\R$\+contramodules
(Theorem~\ref{contramod-contraderived}(a)).
 Furthermore, if the complex of flat $\R$\+contramodules $\gF^\bu$
is contraacyclic, then, by Lemmas~\ref{contramod-kernel-flat}
and~\ref{contramod-product-reduction}, so is the complex of flat
$R_n$\+modules $\gF^\bu/\gI_n\tim\gF^\bu$.
 In view of the assumption of the lemma, by Lemma~\ref{psemi-remark21}
it then follows that the latter complex of flat $R_n$\+modules is
absolutely acyclic.
  Hence the tensor product $(\gF^\bu/\gI_n\tim\gF^\bu)\ot_{R_n}
(\gG^\bu/\gI_n\tim\gG^\bu)$ is an absolutely acyclic complex of
flat $R_n$\+modules, too.
 Since the projective limit of sequences of surjective morphisms
preserves contraacyclicity, we can conclude that the complex
$\gF^\bu\ot^\R\gG^\bu$ is contraacyclic.
\end{proof}

 The result of Lemma~\ref{contramod-tensor-product-contraacyclic}
together with Theorem~\ref{contramod-contraderived}(a) allows to define
the left derived tensor product functor
$$
 \ot^{\R,\.\boL}\:\sD^\ctr(\R\contra)\times\sD^\ctr(\R\contra)\lrarrow
 \sD^\ctr(\R\contra)
$$
on the contraderived category of $\R$\+contramodules
$\sD^\ctr(\R\contra)$, providing it with the structure of
an (associative and commutative) tensor triangulated category.
 In order to compute the contraderived category object $\P^\bu
\ot^{\R,\.\boL}\Q^\bu$ for two complexes of $\R$\+contramodules $\P^\bu$
and $\Q^\bu$, one picks two complexes of flat $\R$\+contramodules
$\gF^\bu$ and $\gG^\bu$ endowed with morphisms of complexes of
contramodules with contraacyclic cones $\gF^\bu\rarrow\P^\bu$ and
$\gG^\bu\rarrow\Q^\bu$, and puts $\P^\bu\ot^{\R,\.\boL}\Q^\bu=
\gF^\bu\ot^\R\gG^\bu$.
 The one-term complex $\R$ is the unit object of this
tensor category structure.

 Furthermore, there is a natural associativity isomorphism of
discrete $\R$\+modules $\M\ocn_\R(\P\ot^\R\Q)\simeq
(\M\ocn_\R\P)\ocn_\R\Q$ for any $\R$\+contramodules $\P$ and $\Q$
and any discrete $\R$\+module~$\M$
(as one can see, e.~g., from the second paragraph of the proof
of Lemma~\ref{flat-contramodules-tensor-product} together with
the explicit construction of the contratensor product functor
$\ocn_\R$ in Section~\ref{ind-affine-co-contra-subsect}).
 The following lemma shows how to construct the left derived
contratensor product functor
$$
 \ocn_\R^\boL\:\sD^\co(\R\discr)\times\sD^\ctr(\R\contra)\lrarrow
 \sD^\co(\R\discr)
$$
providing the coderived category $\sD^\co(\R\discr)$ with
the structure of a triangulated module category over
the triangulated tensor category $\sD^\ctr(\R\contra)$.
 
\begin{lem}
 Assume that, for every $n\ge0$, all flat $R_n$\+modules have finite
projective dimensions.
 Then the contratensor product\/ $\M^\bu\ocn_\R\gF^\bu$ of a complex of
discrete\/ $\R$\+modules\/ $\M^\bu$ and a complex of flat\/
$\R$\+contramodules\/ $\gF^\bu$ is coacyclic provided that either
the complex\/ $\M^\bu$ is coacyclic, or the complex\/ $\gF^\bu$ is
contraacyclic.
\end{lem}

\begin{proof}
 If the complex of $\R$\+contramodules $\gF^\bu$ is contraacyclic then,
according to the proof
of Lemma~\ref{contramod-tensor-product-contraacyclic}, the complexes of
flat $R_n$\+modules $\gF^\bu/\gI_n\tim\gF^\bu$ are absolutely acyclic.
 It follows that the complexes of $R_n$\+modules
${}_{R_n}\M^\bu\ot_{R_n}(\gF^\bu/\gI_n\tim\gF^\bu)$ are absolutely
acyclic, too; hence their inductive limit $\M^\bu\ocn_\R\gF^\bu$
is a coacyclic complex of discrete $\R$\+modules.
 The contratensor product of a complex of flat $\R$\+contramodules and
a coacyclic complex of discrete $\R$\+modules is coacyclic since
the functor of contratensor product with a complex of flat
$\R$\+contramodules preserves infinite direct sums and exactness of
short sequences of complexes of discrete $\R$\+modules.
\end{proof}

 In order to compute the coderived category object
$\M^\bu\ocn_\R^\boL\P^\bu$ for a complex of discrete $\R$\+modules
$\M^\bu$ and a complex of $\R$\+contramodules $\P^\bu$, one picks
a complex of flat $\R$\+contramodules $\gF^\bu$ endowed with a morphism
of complexes of $\R$\+contramodules $\gF^\bu\rarrow\P^\bu$ with
a contraacyclic cone, and puts $\M^\bu\ocn_\R^\boL\P^\bu=
\M^\bu\ocn_\R\gF^\bu$.

 Now assume that we are given a dualizing complex $\D^\bu$ for
the topological ring $\R$, and that, for every $n\ge0$, all
fp\+injective $R_n$\+modules have finite injective dimensions.
 Then it is clear from the natural isomorphism $\D^\bu\ocn_\R
(\gF^\bu\ot^\R\gG^\bu)\simeq(\D^\bu\ocn_\R\gF^\bu)\ocn_\R\gG^\bu$
that the equivalence of triangulated categories
$\sD^\ctr(\R\contra)\simeq\sD^\co(\R\discr)$ from
Theorem~\ref{ind-affine-derived-co-contra} transforms the derived
contramodule tensor product functor $\ot^{\R,\.\boL}$ into the derived
contratensor product functor~$\ocn_\R^\boL$.

 The functor of \emph{cotensor product} of complexes of discrete
$\R$\+modules
$$
 \oc_{\D^\bu}\:\sD^\co(\R\discr)\times\sD^\co(\R\discr)\lrarrow
 \sD^\co(\R\discr)
$$
is obtained from the functor $\ocn_\R^\boL$ or $\ot^{\R,\.\boL}$ by
identifying one or both of the coderived categories $\sD^\co(\R\discr)$
in the arguments of the functors with the contraderived category
$\sD^\ctr(\R\contra)$ using the equivalence of 
Theorem~\ref{ind-affine-derived-co-contra}.
 Explicitly, this means that in order to compute the cotensor product
of two object of the coderived category of discrete $\R$\+modules,
one has to represent them by two complexes of injective discrete
$\R$\+modules $\I^\bu$ and $\J^\bu$ and apply one of the formulas
\begin{gather*}
 \I^\bu\oc_{\D^\bu}\J^\bu = \D^\bu\ocn_\R
 (\Hom_\R(\D^\bu,\I^\bu)\ot^\R\Hom_\R(\D^\bu,\J^\bu)) \\
 \simeq \I^\bu\ocn_\R\Hom_\R(\D^\bu,\J^\bu)\simeq
 \J^\bu\ocn_\R\Hom_\R(\D^\bu,\I^\bu).
\end{gather*}
 Using the latter two formulas, it suffices that only one of
the complexes of discrete $\R$\+modules $\I^\bu$ and $\J^\bu$ be
a complex of injective discrete $\R$\+modules.

 The functor of cotensor product of complexes of discrete
$\R$\+modules $\oc_{\D^\bu}$ defines the structure of an (associative
and commutative) tensor triangulated category on the coderived
category $\sD^\co(\R\discr)$.
 The dualizing complex $\D^\bu$ is the unit object of this tensor
category structure.

\subsection{Very flat and contraadjusted contramodules}
 For the rest of the appendix we stick to a single projective
system of associative rings and their surjective morphisms
$R_0\larrow R_1\larrow R_2\larrow\dotsb$ with the projective
limit $\R=\varprojlim_n R_n$.
 In this section we assume that the rings $R_n$ are commutative.

 Then an $\R$\+contramodule $\gF$ is said to be \emph{very flat} if
the $R_n$\+modules $\gF/\gI_n\tim\gF$ are very flat for all~$n$
(cf.\ Section~\ref{veryflat-contra-subsect}).
 It is clear from Corollary~\ref{flat-reductions-enough} that any
very flat $\R$\+contramodule is flat.
 In the case when the ideals $\gI_n/\gI_{n+1}\sub R_{n+1}$ are nilpotent
and finitely generated, it follows from
Lemma~\ref{fin-gen-nilpotent-quotient-scalars}(b)
that a flat $\R$\+contramodule $\gF$ is very flat if and only if
the $R_0$\+module $\gF/\gI_0\tim\gF$ is very flat.

\begin{cor}  \label{very-flat-contramod-cor}
 The class of very flat\/ $\R$\+contramodules contains the projective\/
$\R$\+contramodules and is closed under extensions and the passage
to the kernels of surjective morphisms in\/ $\R\contra$.
 The projective dimension of any very flat\/ $\R$\+contramodule
(as an object of\/ $\R\contra$) does not exceed\/~$1$.
\end{cor}

\begin{proof}
 Follows from
Lemmas~\ref{contramod-kernel-flat}\+-\ref{contramod-extension-flat}
and Corollary~\ref{projective-contramod}(a).
\end{proof}

 An $\R$\+contramodule $\Q$ is said to be \emph{contraadjusted}
if the functor $\Hom^\R({-},\Q)$ takes short exact sequences of
very flat $\R$\+contramodules to short exact sequences of abelian
groups (or, equivalently, of $\R$\+modules, or of $\R$\+contramodules).
 It is clear from the first assertion of
Corollary~\ref{very-flat-contramod-cor} that an $\R$\+contramodule
$\Q$ is contraadjusted if and only if $\Ext^{\R,1}(\gF,\Q)=0$ for any
very flat $\R$\+contramodule $\gF$, and if and only if
$\Ext^{\R,>0}(\gF,\Q)=0$ for any very flat~$\gF$.

 It follows that the class of contraadjusted $\R$\+contramodules is
closed under extensions and the passages to the cokernels of injective
morphisms.
 Moreover, the second assertion of
Corollary~\ref{very-flat-contramod-cor} implies that any quotient
$\R$\+contramodule of a contraadjusted $\R$\+contramodule is
contraadjusted.

 For the rest of this section we assume that (the rings $R_n$ are
commutative and) the ideals $\gI_n/\gI_{n+1}=\ker(R_{n+1}\to R_n)
\sub R_{n+1}$ are finitely generated.

\begin{lem}  \label{ideal-contraadjusted-lemma}
 Let $R$ be a commutative ring and $I\sub R$ be a finitely generated
ideal.
 Then the $R$\+module $IQ$ is contraadjusted for any contraadjusted
$R$\+module~$Q$.
\end{lem}

\begin{proof}
 The $R$\+module $IQ$ is a quotient module of a finite direct sum of
copies of the $R$\+module~$Q$, and the class of contraadjusted
$R$\+modules is preserved by the passages to finite direct sums and
quotients.
\end{proof}

\begin{lem}  \label{contraadjusted-reduction-lemma}
 Let\/ $\P$ be an\/ $\R$\+contramodule for which the natural map\/
$\P\rarrow\varprojlim_n\P/\gI_n\tim\P$ is an isomorphism.
 Then the\/ $\R$\+contramodule\/ $\P$ is contraadjusted whenever
the $R_n$\+modules\/ $\P/\gI_n\tim\P$ are contraadjusted
for all\/~$n$.
\end{lem}

\begin{proof}
 Applying Lemma~\ref{ideal-contraadjusted-lemma} to the ring
$R=R_{n+1}$ with the ideal $I=\gI_n/\gI_{n+1}$ and the $R$\+module
$Q=P_{n+1}=\P/\gI_{n+1}\tim\P$, we conclude that the $R_{n+1}$\+module
$\ker(P_{n+1}\to P_n) = (\gI_n/\gI_{n+1})P_{n+1}$
is contraadjusted.
 Then it remains to make use of
Lemma~\ref{contramod-contramod-flat-ext-lemma}.
\end{proof}

 When the ideals $\gI_n/\gI_{n+1}$ are nilpotent, in view of
Lemma~\ref{fin-gen-nilpotent-quotient-scalars}(a) it suffices to
require that the $R_0$\+module $\P/\gI_0\tim\P$ be contraadjusted
in Lemma~\ref{contraadjusted-reduction-lemma}.

 Let us assume that we are given the datum of a ring $R$ with ideals
$I_n\sub R$, as in Section~\ref{flat-contramod-subsect}.
 Assume additionally the ring $R$ to be also commutative.
 The following construction plays a key role in our approach.

\begin{lem}  \label{contramod-contraadjusted-main-lemma}
 Let\/ $\P$ be an\/ $\R$\+contramodule for which the natural map\/
$\P\rarrow\varprojlim_n\P/\gI_n\tim\P$ is an isomorphism.
 Then\/ the\/ $\R$\+contramodule\/ $\P$ can be embedded into
a contraadjusted\/ $\R$\+contramodule\/ $\Q$ in such a way that
the quotient\/ $\R$\+contramodule\/ $\Q/\P$ is very flat.
\end{lem}

\begin{proof}
 Let us consider $\P$ as an $R$\+module and embed it into
a contraadjusted $R$\+module $K$ in such a way that the quotient
$R$\+module $F=K/\P$ is very flat.
 Then, the $R$\+module $F$ being, in particular, flat, there are
exact sequences of $R_n$\+modules $0\rarrow \P/I_n\P\rarrow
K/I_nK\rarrow F/I_nF\rarrow0$.
 Furthermore, we have surjective morphisms of $R_n$\+modules
$\P/I_n\P\rarrow\P/\gI_n\tim\P$ and the induced short exact
sequences $0\rarrow\P/\gI_n\tim\P\rarrow K/(I_nK+\gI_n\tim\P)
\rarrow F/I_nF\rarrow0$.

 Passing to the projective limits of these systems of short exact
sequences, we obtain a natural morphism from the short exact
sequence $0\rarrow\varprojlim_n\P/I_n\P\rarrow\varprojlim_nK/I_nK
\rarrow\varprojlim_n F/I_nF\rarrow0$ to the short exact sequence
$0\rarrow\varprojlim_n\P/\gI_n\tim\P\allowbreak\rarrow\varprojlim_n
K/(I_nK+\gI_n\tim\P)\rarrow\varprojlim_n F/I_nF\rarrow0$.
 The map $\varprojlim_n\P/I_n\P\rarrow\varprojlim_n\P/\gI_n\tim\P$
is always surjective, since a surjective map $\P\rarrow
\varprojlim_n\P/\gI_n\tim\P$ (see
Lemma~\ref{contramod-projlim-surjectivity}) factorizes through it.
 Hence the map $\varprojlim_n K/I_nK\rarrow\varprojlim_n
K/(I_nK+\gI_n\tim\P)$ is also surjective.
 Furthermore, the projective system $K/(I_nK+\gI_n\tim\P)$ satisfies
the condition of Lemma~\ref{projlim-reduction-lemma}, because
$\gI_n\tim\P = I_n\P+\gI_{n+1}\tim\P$.

 Now put $\gL = \varprojlim_n K/I_nK$ and $\Q = \varprojlim_n
K/(I_nK+\gI_n\tim\P)$, and also $\gG=\varprojlim_n F/I_nF$.
 By Lemma~\ref{projlim-reduction-lemma}, we have
$\gL/\gI_n\tim\gL = K/I_nK$ and $\gG/\gI_n\tim\gG=F/I_nF$.
 Therefore, the $R$\+contramodule $\gG$ is very flat.
 Besides, $\Q/\gI_n\tim\Q = K/(I_nK+\gI_n\tim\P)$, and
the natural map
$\Q\rarrow\varprojlim_n\Q/\gI_n\tim\Q$ is an isomorphism.
 In addition, the $R_n$\+modules $K/I_nK$ are contraadjusted by
Lemma~\ref{quotient-scalars-contraadjusted}(b).
 Finally, the $R_n$\+module $\Q/\gI_n\tim\Q$ is a quotient module
of a contraadjusted $R_n$\+module $\gL/\gI_n\tim\gL=K/I_nK$,
and consequently is also contraadjusted.
 It remains to use Lemma~\ref{contraadjusted-reduction-lemma}.
\end{proof}

\begin{cor}  \label{contramod-vfl-covering}
 Any\/ $\R$\+contramodule can be presented as the quotient contramodule
of a very flat\/ $\R$\+contramodule by a contraadjusted\/
$\R$\+subcontramodule.
\end{cor}

\begin{proof}
 Present an arbitrary $\R$\+contramodule $\P$ as the quotient
contramodule of a free $\R$\+contramodule $\gH$; apply
Lemma~\ref{contramod-contraadjusted-main-lemma} to embed the kernel 
$\gK$ of the morphism $\gH\rarrow\P$ into a contraadjusted
$\R$\+contramodule $\Q$ so that the quotient contramodule $\gF =
\Q/\gK$ is very flat; consider the induced extensions
$0\rarrow\Q\rarrow\gG\rarrow\P\rarrow0$ and
$0\rarrow\gH\rarrow\gG\rarrow\gF\rarrow0$, and
use Lemma~\ref{contramod-extension-flat}.
 It is important here that the map $\gK\rarrow\varprojlim_n
\gK/\gI_n\tim\gK$ is an isomorphism, since so is the map
$\gH\rarrow\varprojlim_n\gH/\gI_n\tim\gH$.
 (This is the same construction as in Lemmas~\ref{salce-lemma}(a)
and~\ref{eklof-trlifaj-cta}.)
\end{proof}

\begin{cor}  \label{contramod-cta-embedding}
 Any\/ $\R$\+contramodule can be embedded into a contraadjusted\/
$\R$\+contramodule in such a way that the quotient contramodule is
very flat.
\end{cor}

\begin{proof}
 Use Corollary~\ref{contramod-vfl-covering} to present an arbitrary
$\R$\+contramodule $\P$ as the quotient $\R$\+contramodule of
a very flat $\R$\+contramodule $\gG$ by a contraadjusted
$\R$\+contra\-module~$\gK$; apply
Lemma~\ref{contramod-contraadjusted-main-lemma} to embed
the $\R$\+contramodule $\gG$ into a contraadjusted $\R$\+contramodule
$\gL$ so that the quotient contramodule $\gF = \gL/\gG$ is very flat;
and denote by $\gN$ the cokernel of the composition of injective
morphisms of $\R$\+contramodules $\gK\rarrow\gG\rarrow\gL$.
 Now the $\R$\+contramodule $\gN$ is contraadjusted as the quotient
contramodule of a contraadjusted $\R$\+contramodule (by
a contraadjusted $\R$\+subcontramodule), and the $\R$\+contramodule
$\P$ is embedded into $\gN$ with the cokernel~$\gF$.
 (This construction comes from~\cite[proof of Lemma~9.1.2(a)]{Psemi}.)
\end{proof}

\begin{lem}  \label{vfl-cta-reduction-nilp-ortho}
 Let $F$ be a flat $R$\+module and\/ $\Q$ be an\/ $\R$\+contramodule
such that the $R_n$\+modules $F/I_nF$ are very flat, the map\/
$\Q\rarrow\varprojlim_n\Q/\gI_n\tim\Q$ is an isomorphism, and
the $R_n$\+modules\/ $\Q/\gI_n\tim\Q$ are contraadjusted.
 Then one has\/ $\Ext^{>0}_R(F,\Q)=0$. \par
\end{lem}

\begin{proof}
 This is a particular case of Lemma~\ref{mod-contramod-flat-ext-lemma}.
 To show that its conditions are satisfied, one only has to recall
that the ideals $\gI_n/\gI_{n+1}\sub R_{n+1}$ are finitely generated by
our assumptions, and use Lemma~\ref{ideal-contraadjusted-lemma}.
\end{proof}

\begin{cor}  \label{cta-mod-contramod-cor}
 An\/ $\R$\+contramodule is contraadjusted if and only if it is
a contraadjusted\/ $R$\+module.
 In particular, the $R_n$\+modules\/ $\P/I_n\P$ and\/ $\P/\gI_n\tim\P$
are contraadjusted for any contraadjusted\/ $\R$\+contramodule\/~$\P$.

 Furthermore, one has\/
$\Hom^\R(\varprojlim_n F/I_nF\;\P)\simeq\Hom_R(F,\P)$ and\/
$\Ext^{\R,>0}(\varprojlim_n F/\allowbreak I_nF\;\P)=\Ext_R^{>0}(F,\P)=0$
for any contraadjusted\/ $\R$\+contramodule\/ $\P$ and any flat
$R$\+module $F$ for which the $R_n$\+modules $F/I_nF$ are very flat.
\end{cor}

\begin{proof}
 According to the constructions of
Lemma~\ref{contramod-contraadjusted-main-lemma}
and Corollaries~\ref{contramod-vfl-covering}\+-%
\ref{contramod-cta-embedding}, any contraadjusted $\R$\+contramodule
$\P$ can be obtained from the $\R$\+contramodules $\Q$ satisfying
the conditions of Lemma~\ref{vfl-cta-reduction-nilp-ortho} using
the operations of the passage to the cokernel of an injective
morphism and the passage to a direct summand.
 This proves the equation $\Ext^{>0}_R(F,\P)=0$, and consequently also
the ``only if'' part of the first assertion; then the second
assertion follows.

 Conversely, suppose that an $\R$\+contramodule $\P$ is a contraadjusted
$R$\+module.
 Let us present $\P$ as the quotient contramodule of a (very) flat
$\R$\+contramodule $\gF$ by a contraadjusted $\R$\+subcontramodule~$\Q$.
 As we have proved, $\Q$ is a contraadjusted $R$\+module, and
therefore so is the $R$\+module~$\gF$.
 Hence the $R_n$\+modules $\gF/\gI_n\tim\gF$ are contraadjusted, and
by Lemma~\ref{contraadjusted-reduction-lemma} it follows that
$\gF$ is a contraadjusted $\R$\+contramodule.
 Now the $\R$\+contramodule $\P$ is contraadjusted as the quotient
$\R$\+contramodule of a contraadjusted one.

 The equation $\Ext^{\R,>0}(\varprojlim_n F/I_nF\;\P)=0$ holds since
the $\R$\+contramodule in the first argument is very flat.
 To prove the isomorphism $\Hom^\R(\varprojlim_n F/I_nF\;\P)
\simeq\Hom_R(F,\P)$, one argues in the same way as in
Corollary~\ref{projective-contramod}(b), representing $\P$ as
a quotient contramodule of a very flat contramodule by a contraadjusted
subcontramodule.
\end{proof}

\subsection{Cotorsion contramodules}
 A left $\R$\+contramodule $\Q$ is said to be \emph{cotorsion} if
the functor $\Hom^\R({-},\Q)$ takes short exact sequences of flat
left $\R$\+contramodules to short exact sequences of abelian groups.
 It follows from Lemmas~\ref{contramod-kernel-flat}\+-%
\ref{contramod-proj-flat} that a left $\R$\+contramodule $\Q$ is
cotorsion if and only if $\Ext^{\R,1}(\gF,\Q)=0$ for any flat
left $\R$\+contramodule $\gF$, and if and only if
$\Ext^{\R,>0}(\gF,\Q)=0$ for any flat~$\gF$.
 Hence the class of cotorsion left $\R$\+contramodules is closed
under extensions and the passages to the cokernels of injective
morphisms and direct summands.

 The \emph{cotorsion dimension} of a left $\R$\+contramodule $\P$
is defined as the supremum of the set of all integers~$d$ for which
there exists a flat left $\R$\+contramodule $\gF$ such that
$\Ext^{\R,d}(\gF,\P)\ne0$.
 Clearly, a nonzero $\R$\+contramodule is cotorsion if and only if
its cotorsion dimension is equal to zero.
 Given a short exact sequence of $\R$\+contra\-modules $0\rarrow\gK
\rarrow\gL\rarrow\gM\rarrow0$ of the cotorsion dimensions $k$, $l$,
and~$m$, respectively, one has $l\le\max(k,m)$, \ $k\le\max(l,m+1)$,
and $m\le\max(l,k-1)$.

\begin{lem}  \label{cotorsion-reduction-lemma}
 Let\/ $\Q$ be a left\/ $\R$\+contramodule for which the natural
map\/ $\Q\rarrow\varprojlim_n\Q/\gI_n\tim\Q$ is an isomorphism.
 Set $Q_n=\Q/\gI_n\tim\Q$, and assume that the left $R_0$\+module $Q_0$
is cotorsion and the left $R_n$\+modules\/ $\ker(Q_{n+1}\to Q_n)$ are
cotorsion for all\/ $n\ge0$.
 Then\/ $\Q$ is a cotorsion\/ $\R$\+contramodule and a cotorsion
left $R$\+module.
\end{lem}

\begin{proof}
 This is essentially a special case of the dual Eklof
lemma~\cite[Proposition~18]{ET}, \cite[Lemma~6.37]{GT},
\cite[Lemma~4.5]{PR}.
 The equations $\Ext^{\R,1}(\gF,\Q)=0=\Ext^1_R(F,\Q)$ for any flat
left $\R$\+contramodule $\gF$ and any flat left $R$\+module $F$ follow
from Lemmas~\ref{contramod-contramod-flat-ext-lemma}
and~\ref{mod-contramod-flat-ext-lemma}, respectively.
\end{proof}

 For the rest of this section we assume that the rings $R_n$ are
commutative and Noetherian.
 The ring $R$ is also presumed to be commutative.
 Our approach is inspired by~\cite[Section~4.3]{Xu}.

\begin{cor}  \label{flat-cotorsion-reduction-cor}
 Let\/ $\gE$ be a flat\/ $\R$\+contramodule.
 Assume that the flat $R_n$\+modules\/ $E_n=\gE/\gI_n\tim\gE$ are
cotorsion for all\/ $n\ge0$.
 Then\/ $\gE$ is a flat cotorsion\/ $\R$\+contra\-module and
a cotorsion $R$\+module.
\end{cor}

\begin{proof}
 Since $\ker(E_{n+1}\to E_n)=(\gI_n/\gI_{n+1})E_{n+1}\simeq
\gI_n/\gI_{n+1}\ot_{R_{n+1}}E_{n+1}$, the assertions follow from
Lemma~\ref{cotorsion-reduction-lemma} in view of
Lemma~\ref{coherent-cotors}(a).
\end{proof}

\begin{prop}  \label{contramod-flat-cotorsion-embedding}
 Any flat\/ $\R$\+contramodule can be embedded into a flat
cotorsion\/ $\R$\+contramodule in such a way that the quotient\/
$\R$\+contramodule is flat.
\end{prop}

\begin{proof}
 In fact, we will even present a functorial construction of
such an embedding of flat $\R$\+contramodules.
 It is based on the functorial construction of an embedding
$G\rarrow\FC_S(G)$ of a flat module $G$ over a commutative
Noetherian ring $S$ into a flat cotorsion module $\FC_S(G)$
with a flat quotient $S$\+module $\FC_S(G)/G$, which was explained
in Lemma~\ref{noetherian-flat-cotorsion-embedding}.
 In addition, we will need the following lemma.

\begin{lem} \label{flat-cotorsion-embedding-quotient-scalars}
 Let $T$ be a commutative Noetherian ring, $S=T/I$ be its quotient
ring by an ideal $I\sub T$, and $H$ be a flat $T$\+module.
 Then there is a natural isomorphism of $S$\+modules\/ $\FC_T(H)/
I\FC_T(H)\simeq\FC_S(H/IH)$ compatible with the natural morphisms
$H\rarrow\FC_T(H)$ and $H/IH\rarrow\FC_S(H/IH)$.
\end{lem}

\begin{proof}
 Since the ideal $I\sub T$ is finitely generated, the reduction
functor $S\ot_T\nobreak{-}$ preserves infinite direct products.
 For any prime ideal $\q\sub T$ that does not contain $I$ and
any $T_\q$\+module $Q$ one has $Q/IQ=0$; in particular, this
applies to any $T_\q$\+contramodule.
 Now consider a prime ideal $\q\supset I$ in $T$, and let
$\p=\q/I$ be the related prime ideal in~$S$.
 It remains to construct a natural isomorphism of $S$\+modules
$\widehat H_\q/I\widehat H_q\simeq\widehat{(H/IH)}_\p$.
 Indeed, according to~\cite[proof of Lemma~B.9.2]{Pweak}
applied to the Noetherian ring $T$ with the ideal $\q\sub T$,
the complete ring $\gT=\widehat T_\q$, and the finitely generated
$T$\+module $S$, one has $S\ot_T\widehat H_\q \simeq
\varprojlim_n S\ot_T (H_\q/\q^nH_q) \simeq
\varprojlim_n (S\ot_TH)_\p/\p^n(S\ot_TH)_\p$, as desired.
\end{proof}

 So let $\gG$ be a flat $\R$\+contramodule; set
$G_n=\gG/\gI_n\tim\gG$.
 According to Lemma~\ref{flat-cotorsion-embedding-quotient-scalars},
the flat cotorsion $R_n$\+modules $\FC_{R_n}(G_n)$ naturally form
a projective system satisfying the condition of
Lemma~\ref{projlim-reduction-lemma}.
 Set $\gE=\varprojlim_n\FC_{R_n}(G_n)$; by
Corollary~\ref{flat-cotorsion-reduction-cor},
the flat $\R$\+contramodule $\gE$ is cotorsion.
 The projective system of $R_n$\+modules $\FC_{R_n}(G_n)/G_n$ also
satisfies the condition of Lemma~\ref{projlim-reduction-lemma},
being the cokernel of a morphism of two projective systems that do;
in view of the second assertion of
Lemma~\ref{noetherian-flat-cotorsion-embedding}, it follows that
the $\R$\+contramodule $\gF=\varprojlim_n\FC_{R_n}(G_n)/G_n$ is flat.
 Finally, the short sequence $0\rarrow\gG\rarrow\gE\rarrow\gF\rarrow0$
is exact, since $\varprojlim_n^1G_n=0$.
\end{proof}

\begin{cor}  \label{flat-cotors-mod-contramod-cor}  \hbadness=1100
 A flat\/ $\R$\+contramodule\/ $\gG$ is cotorsion if and only if
the $R_n$\+modules\/ $\gG/\gI_n\tim\gG$ are cotorsion for all\/ $n\ge0$.
 Any flat cotorsion\/ $\R$\+contramodule is a cotorsion $R$\+module.
\end{cor}

\begin{proof}
 In view of Corollary~\ref{flat-cotorsion-reduction-cor}, we only have
to show that the $R_n$\+modules $\gG/\gI_n\tim\gG$ are cotorsion for
any flat cotorsion $\R$\+contramodule~$\gG$.
 For this purpose, it suffices to apply the construction of
Proposition~\ref{contramod-flat-cotorsion-embedding} to obtain
a short exact sequence of $\R$\+contramodules $0\rarrow\gF\rarrow
\gE\rarrow\gG\rarrow0$, where the $\R$\+contramodule $\gF$ is flat
and the $\R$\+contramodule $\gE$, by construction, has the desired
property.
 Now $\gG$ is an $\R$\+contramodule direct summand of $\gE$, since
$\Ext^{\R,1}(\gG,\gF)=0$ by assumption.
\end{proof}
 
 From this point on and until the end of this section we assume that
the Krull dimensions of the Noetherian commutative rings $R_n$ are
uniformly bounded by a constant~$D$.
 Then it follows from the results of
Section~\ref{flat-contramod-subsect} and
Theorem~\ref{raynaud-gruson-flat-thm} that the projective dimension
of any flat $\R$\+contramodule does not exceed~$D$.
 Therefore, the cotorsion dimension of any $\R$\+contramodule also
cannot exceed~$D$.

\begin{lem}  \label{contramod-flat-cotorsion-theory}
\textup{(a)} Any\/ $\R$\+contramodule can be embedded into
a cotorsion\/ $\R$\+contra\-module in such a way that the quotient\/
$\R$\+contramodule is flat.  \par
\textup{(b} Any\/ $\R$\+contramodule admits a surjective map onto it
from a flat\/ $\R$\+contramodule with the kernel being a cotorsion\/
$\R$\+contramodule.
\end{lem}

\begin{proof}
 Part~(a): we will argue by decreasing induction in $0\le d\le D$,
showing that any $\R$\+contramodule $\P$ can be embedded into
an $\R$\+contramodule $\Q$ of cotorsion dimension~$\le d$ in such
a way that the quotient $\R$\+contramodule $\Q/\P$ is flat.

 Assume that we already know this for the cotorsion dimension~$d+1$.
 Let $\P$ be an $\R$\+contramodule; pick a surjective morphism
$\gL\rarrow\P$ onto $\P$ from a projective $\R$\+contramodule~$\gL$.
 Denote the kernel of the morphism $\gL\rarrow\P$ by $\gK$ and
use the induction assumption to embed it into an $\R$\+contramodule
$\gN$ of cotorsion dimension~$\le d+1$ so that the cokernel
$\gH=\gN/\gK$ is flat.

 Let $\gG$ be the fibered coproduct $\gN\sqcup_\gK\gL$; then
the $\R$\+contramodule $\gG$ is flat as an extension of two
flat $\R$\+contramodules $\gH$ and~$\gL$.
 Now we use Proposition~\ref{contramod-flat-cotorsion-embedding}
to embed $\gG$ into a flat cotorsion $\R$\+contramodule $\gE$ so that
the cokernel $\gF=\gE/\gG$ is a flat $\R$\+contramodule.
 Then the cokernel $\Q$ of the composition $\gN\rarrow\gG\rarrow\gE$
has cotorsion dimension~$\le d$.
 The $\R$\+contramodule $\P$ embeds naturally into $\Q$ with
the cokernel $\gF$, so the desired assertion is proved.

 We have also proved part~(b) along the way (it suffices to notice
that there is a natural surjective morphism $\gG\rarrow\P$ with
the kernel~$\gN$); see also Lemma~\ref{salce-lemma}(a) or the proof
of Lemma~\ref{eklof-trlifaj-cta}.
\end{proof}

\begin{cor}
 Any cotorsion\/ $\R$\+contramodule\/ $\P$ is also a cotorsion
$R$\+module.
 Besides, for any flat $R$\+module $F$ one has\/
$\Hom^\R(\varprojlim_n F/I_nF\;\P)\simeq\Hom_R(F,\P)$
and\/ $\Ext^{\R,>0}(\varprojlim_n F/I_nF\;\P)=\Ext_R^{>0}(F,\P)=0$.
\end{cor}

\begin{proof}
 Applying the construction of
Lemma~\ref{contramod-flat-cotorsion-theory}(a) to
the $\R$\+contramodule $\P$, we obtain a short exact sequence of
$\R$\+contramodules $0\rarrow\P\rarrow\Q\rarrow\gF\rarrow0$
with a flat $\R$\+contramodule $\gF$ and an $\R$\+contramodule
$\Q$ obtained from flat cotorsion $\R$\+contramodules using
the operation of passage to the cokernel of an injective morphism
at most~$D$ times.
 Then $\P$ is a direct summand of $\Q$, and since the class of
cotorsion $R$\+modules is preserved by the passages to the cokernels
of injective morphisms and direct summands, the first assertion
follows from Corollary~\ref{flat-cotors-mod-contramod-cor}.
 The Ext-vanishing assertions now hold by the definitions, and
the Hom isomorphism is obtained in the same way as in
Corollaries~\ref{projective-contramod}(b)
and~\ref{cta-mod-contramod-cor}.
\end{proof}

\bigskip
\addtocontents{toc}{\smallskip}

\end{document}